\documentclass[a4paper,10pt]{book}
\usepackage[utf8]{inputenc}
\title{La structure des courbes analytiques}
\DeclareUnicodeCharacter{00A0}{}
\author{\sc Antoine Ducros\\ \small Institut de mathématiques de Jussieu \\ \small Projet {\em Topologie et géométrie algébriques}\\ \small 4 place Jussieu 75005 Paris ~~FRANCE}
\date{}
\usepackage{mathrsfs}

\usepackage[francais]{babel}
\usepackage{amssymb}
\usepackage{hyperref} 
\newcommand{\conv}{\mathsf {Conv}\;}
\date{}
\input xy
\xyoption{all} 
\input{xypic} 
\makeatletter
\let\old@subsection\subsection
\newcommand{\subsection@star}[1]{\old@subsection*{#1}\addcontentsline{toc}{subsection}{\hspace{.2cm}#1}}
\renewcommand{\subsection}{\@ifstar{\subsection@star}{\old@section}}
\makeatother
\newcommand{\ti}{^\times}
\newcommand{\abs}[1]{|#1|}

\newcommand{\Aff}{{\Bbb A}}
\newcommand{\an}{^{\rm an}}
\newcommand{\arb}[2]{\mathsf{Arb}(#1,#2)}
\newcommand{\bbfer}[2]{{\mathbb B}(#1,#2)}
\newcommand{\bbouv}[2]{{\mathbb B}\zero (#1,#2)}
\newcommand{\bfer}[2]{B(#1,#2)}
\newcommand{\bnd}[1]{{\cal  #1}}
\newcommand{\bouv}[2]{B\zero (#1,#2)}
\newcommand{\br}[2]{\mathsf {br}(#1,#2)}
\newcommand{\brdan}{\partial^{\rm an}}
\newcommand{\CC}{{\Bbb C}}
\newcommand{\card}{\mathsf{card}}

\newcommand{\ctd}[1]{_{\subset #1}}
\newcommand{\DD}{{\Bbb D}}

\newcommand{\dtr}{_{[2,3]}}
\newcommand{\deux}[1]{\refstepcounter{subsection}\label{#1}\medskip\noindent {\bf (\thesubsection)}\hspace{.1cm}}

\newcommand{\et}{_{\rm \acute{e}t}}
\newcommand{\FF}{{\Bbb F}}
\newcommand{\form}[2]{\widehat{(#1,\got #2)}}
\newcommand{\fort}[2]{\widehat{(#1, #2)}}

\newcommand{\geom}{_{[0,2,3]}}
\newcommand{\gm}{{\Bbb G}_{m}}

\newcommand{\gmk}{{\Bbb G}_{m,k}}
\newcommand{\got}[1]{{\mathfrak #1}}
\newcommand{\gred}[1]{\red{#1}_{\rm gr}}
\newcommand{\grot}{_{\rm G}}

\newcommand{\hotimes}{\wid{\otimes}}
\newcommand{\hres}{{\sch H}}

\newcommand{\inter}[2]{\mathsf{Interv}(#1,#2)}
\newcommand{\intera}[2]{\mathsf{Interv}_{\rm a}(#1,#2)}
\newcommand{\interac}[2]{\mathsf{Interv}_{\rm a,c}(#1,#2)}
\newcommand{\interc}[2]{\mathsf{Interv}_{\rm c}(#1,#2)}
\newcommand{\inv}{^{-1}}
\newcommand{\kparf}{\widehat{k^{\rm parf}}}

\newcommand{\itb}{\medskip \item[$\bullet$]}
\newcommand{\JJ}{{\mathbb J}}

\newcommand{\kk}{\red k^a}
\newcommand{\KK}{\wid{k^a}}
\newcommand{\kum}{\mathsf {Kum}_\ell}

\newcommand{\NN}{{\Bbb N}}
\newcommand{\pk}{\PP^{1,{\rm an}}_k}
\newcommand{\pkk}{\PP^{1,{\rm an}}_{\KK}}
\newcommand{\PP}{{\Bbb P}}

\newcommand{\QQ}{{\Bbb Q}}

\newcommand{\red}{\widetilde}
\newcommand{\RR}{{\Bbb R}}

\newcommand{\rnorm}{\mathsf r_{\rm norm}}

\newcommand{\sbr}[1]{\mathsf {Sec}\;#1}
\newcommand{\sch}[1]{\mathscr #1}

\newcommand{\sgb}[2]{\mathsf{Secbr}(#1,#2)}
\newcommand{\shil}[2]{#2_{[#1]}}
\newcommand{\skel}[1]{\mathsf S(#1)}
\newcommand{\skelan}[1]{\mathsf S\an(#1)}

\newcommand{\spec}{\mathsf{Spec}\;}
\newcommand{\spf}{\mathsf{Spf}\;}
\newcommand{\stel}[2]{\mathsf{St}(#1,#2)}
\newcommand{\stela}[2]{\mathsf{St}_{\rm a}(#1,#2)}
\newcommand{\stelac}[2]{\mathsf{St}_{\rm a,c}(#1,#2)}
\newcommand{\stelc}[2]{\mathsf{St}_{\rm c}(#1,#2)}
\newcommand{\trois}[1]{\refstepcounter{subsubsection}\label{#1}\medskip\noindent {\bf
    (\thesubsubsection)}\hspace{.1cm}}

\newcommand{\wid}{\widehat}
\newcommand{\tp}{_{top}}
\newcommand{\typ}[1]{_{[#1]}}

\newcommand{\val}[1]{\langle{#1}\rangle}
\newcommand{\wbeth}[1]{\widehat {\beth(#1)}}
\newcommand{\widd}[2]{\wid #1^{#2}}

\newcommand{\zero}{^{\mbox{\tiny o}}}
\newcommand{\zeroo}{^{\mbox{\tiny oo}}}
\newcommand{\ZZ}{{\Bbb Z}}

\renewcommand{\top}{_{\rm top}}
\renewcommand{\leq}{\leqslant}
\renewcommand{\geq}{\geqslant}

\renewcommand{\phi}{\varphi}
\renewcommand{\dim}[1]{\mbox{\rm dim}_{#1}\;}
\renewcommand{\Bbb}{\mathbb}    
\renewcommand{\H}{\mbox{\rm H}}
\renewcommand{\epsilon}{\varepsilon}

\begin{document}
\maketitle
\pagestyle{myheadings}
\tableofcontents

\newpage
\chapter*{(Trop brève) introduction}\addcontentsline{toc}{chapter}{(Trop brève) introduction}
Le but de cette monographie, qui est très loin d'être sous sa forme définitive, est de procéder à une étude systématique 
des courbes analytiques \textit{au sens de Berkovich} sur un corps ultramétrique complet. Nous considèrerons comme connues
et utiliserons librement, 
sans références les définitions et propriétés générales des espaces de Berkovich. Nous renvoyons plus précisément aux textes fondateurs 
\cite{berkovich1990} et \cite{berkovich1993}, ainsi qu'aux travaux de l'auteur \cite{ducros2007} et \cite{ducros2009} sur la théorie de la dimension,
la notion d'espace quasi-lisse ou celle de composante irréductible. Nous ferons également un grand usage de la théorie de la réduction (graduée)
des germes d'espaces analytiques due à Temkin (\cite{temkin2004}). 

Au chapitre 4 de \cite{berkovich1990}, Berkovich donne une première description des courbes analytiques en se fondant sur le théorème de
réduction semi-stable. Dans le présent texte, nous changeons de point de vue : nous nous livrons à une étude directe locale et globale des courbes
analytiques, ce qui nous permet de montrer que toute courbe admet une \textit{triangulation} (\ref{deftri}, théorème \ref{theotri}), 
et nous en \textit{déduisons} le théorème de réduction semi-stable (\ref{semstabpot}).

Nous avions déjà
introduit et utilisé la notion de triangulation dans des travaux antérieurs sur la comparaison entre cohomologies étale et topologique
d'une courbe analytique (\cite{ducros2008}), mais nous avions alors adopté l'approche classique consistant à construire une triangulation
à partir d'un modèle semi-stable et non l'inverse comme nous le faisons ici. 

Insistons sur le fait que ce travail est encore en chantier, et que nous sommes conscients qu'il n'est pas publiable en
l'état (il pèche notamment par manque criant de références à la littérature existante). Mais comme il a d'ores et déjà été utilisé et cité
(\cite{colmez-d-n2022}, \cite{cubides-p2021}, \cite{gubler-j-r2021}), il nous a semblé préférable de le rendre accessible de manière pérenne
en le déposant sur \verb+arXiv+ .

\chapter{Topologie : la théorie des graphes réels}

\markboth{Graphes réels }{Graphes réels}

\section{Généralités}\label{RAP}

{\em Nous allons rappeler quelques faits élémentaires que nous utiliserons librement par la suite}.

\deux{boradadh} Si~$U$ est une partie d'un espace topologique~$X$, on notera~$\overline U^X$ (resp.~$\partial_XU$) l'adhérence (resp. le bord) de~$U$ dans~$X$ ; on écrira simplement~$\overline U$ et~$\partial U$ s'il n'y a pas d'ambiguïté sur l'espace ambiant. 

\deux{compprop} Nous dirons qu'une application continue~$\phi : Y\to X$ entre espaces topologiques est {\em compacte} si~$\phi\inv({\cal K})$ est compact pour tout compact~$\cal K$ de~$Y$ ; une application compacte dont le but est localement compact est fermée. 

\deux{actgconv} Si~$\mathsf G$ est un groupe topologique agissant sur un espace topologique~$X$ il sera toujours sous-entendu, sauf mention  expresse du contraire, que~$\mathsf G\times X\to X$ est continue.

\deux{appcont} Soit~$\phi : Y\to X$ une application continue et compacte entre espaces topologiques séparés et localement connexes. 

\trois{fprop} Soit~$x$ un point de~$X$ et soit~$y$ un antécédent de~$x$ isolé dans~$\phi\inv(x)$. Si~$U$ est un voisinage ouvert connexe de~$x$, on notera~$\phi\inv(U)_y$ la composante connexe de~$y$ dans~$\phi\inv(U)$ ; on dira que~$\phi\inv(U)$ {\em isole~$y$} si~$y$ est le seul antécédent de~$x$ sur~$\phi\inv(U)_y$ ; si~$\phi\inv(x)$ est discrète, on dira que~$\phi\inv(U)$ {\em sépare les antécédents de~$x$} s'il isole chacun d'eux. 

\medskip
Les ouvert de la forme~$\phi\inv(U)$, où~$U$ est un voisinage ouvert connexe de~$x$, forment une base de voisinages de la fibre~$\phi\inv(x)$ ; par conséquent, si~$U$ est un voisinage ouvert suffisamment petit de~$x$ et si~$\phi\inv(x)$ est finie alors~$\phi\inv(U)$ sépare les antécédents de~$x$. 

\trois{imouvferm} Supposons que~$\phi$ est de plus {\em ouverte}. Si~$U$ est un ouvert de~$X$ et si~$V$ est une composante connexe de~$\phi\inv(U)$ alors~$\phi(V)$ est une partie, connexe non vide, ouverte et fermée de~$U$ ; c'en est donc une composante connexe ; en particulier, si~$U$ est lui-même connexe on a~$\phi(V)=U$. 

\deux{compxmoinse} Soit~$X$ un espace topologique connexe et localement connexe, soit~$E$ un sous-ensemble fermé et non vide de~$X$, et soit~$U$ un ouvert contenu dans~$X\setminus E$.  Nous utiliserons systématiquement et de façon implicite l'équivalence des propriétés suivantes :

i)~$U$ est une composante connexe de~$X\setminus E$ ;

ii)~$U$ est connexe, non vide et~$\partial U\subset E$ ;

iii)~$U$ est connexe, non vide, et~$\partial U$ est non vide et contenu dans~$E$. 

\medskip
Justifions-la brièvement : si i) est vraie alors~$U$ est connexe, non vide, et fermée dans~$X\setminus E$, d'où ii) ; si ii) est vrai le bord de~$U$ ne peut être vide car sinon l'ouvert non vide~$U$ de~$X$ serait également fermé, contredisant la connexité de~$X$ ; et si iii) est vraie alors~$U$ est une partie ouverte, connexe, non vide et fermée de~$X\setminus E$, d'où i). 

\deux{compquoti} Soit~$\phi: Y\to X$ une surjection continue entre espaces topologiques et soit~$\sch R$ la relation d'équivalence sur~$Y$ dont les classes sont les fibres de~$\phi$. Si~$Y$ est quasi-compact et si~$X$ est séparé la flèche naturelle~$Y/\sch R\to X$ est un homéomorphisme. 

\deux{condflechecomp} Soit~$\phi: Y\to X$ une application continue entre espaces topologiques séparés. Les conditions suivantes sont équivalentes : 

\medskip
i)~$X$ est localement compact et~$\phi$ est compacte ; 

ii) tout point de~$X$ possède un voisinage compact dont l'image réciproque est compacte. 

\medskip
En effet, il est clair que i)$\Rightarrow$ ii). Réciproquement, supposons que ii) soit satisfaite ; cela entraîne immédiatement la locale compacité de~$X$. Soit~$\sch K$ un compact de~$X$ ; en vertu de l'hypothèse ii), o, peut le recouvrir par un nombre fini de compacts~$\sch K_i$ tels que~$\phi\inv(\sch K_i)$ soit compact pour tout~$i$. L'image réciproque de~$\sch K$ est alors fermée dans le compact~$\bigcup p\inv(\sch K_i)$, et est donc elle-même compacte. 

\medskip
Supposons que ces conditions sont satisfaites, et que~$\phi$ est de surcroît surjective. Soit~$X'$ une partie de~$X$ dont l'image réciproque~$Y'$ est fermée. Si~$V$ est un compact de~$X$ alors~$V\cap X'$ est égal à~$\phi(Y'\cap \phi\inv(V))$ ; comme~$Y'$ est fermé et~$\phi\inv(V)$ compact,~$V\cap X'$ est compact. Joint à la compacité locale de~$X$, ceci entraîne que~$X'$ est fermé ; par conséquent,~$X$ s'identifie au quotient de~$Y$ par la relation d'équivalence dont les classes sont les fibres de~$\phi$. 

\deux{stabgcomp} Soit~$X$ un espace topologique localement connexe et soit~$\mathsf G$ un groupe topologique agissant sur~$X$. Si~$U$ est une composante connexe de~$X$ le stabilisateur de~$U$ est un sous-groupe ouvert de~$\mathsf G$ : pour le voir, il suffit de choisir un point~$x$ sur~$U$ et de remarquer que le stabilisateur en question est l'image réciproque de~$U$ par l'application~$g\mapsto g.x$. 

\deux{quotsep} Soit~$X$ un espace topologique séparé et localement compact et soit~$\mathsf G$ un groupe topologique compact agissant sur~$X$. Nous allons montrer que~$X/\mathsf G$ est séparé et localement compact, et que la flèche~$X\to X/\mathsf G$ est ouverte et compacte. 

\trois{xsurgsep} {\em Le quotient~$X/\mathsf G$ est séparé.} Soient~$x$ et~$y$ deux points distincts de~$X/\mathsf G$. Ils correspondent à deux orbites distinctes~$\mathsf G.\xi$ et~$\mathsf  G.\eta$ de~$\mathsf G$. Par compacité de~$\mathsf G$ et par séparation et locale compacité de~$X$, il existe un voisinage ouvert~$V$ de~$\xi$ tel que~$\overline V$ soit compact et ne rencontre pas~$\mathsf G.\eta$. Le sous-ensemble~$\mathsf G.\overline V$ de~$X$ est l'image de la partie compacte~$\mathsf G\times \overline V$ de~$\mathsf G\times X$, et est donc compact ; son complémentaire~$W$ dans~$X$ est un ouvert~$\mathsf G$-invariant qui contient~$\mathsf G.\eta$ ; par ailleurs,~$\mathsf G.V$ est un ouvert~$\mathsf G$-invariant qui contient~$\mathsf G.x$ et ne rencontre pas~$W$. 

Si l'on note~$V'$ (resp.~$W'$) l'image de~$V$ (resp.~$W$) sur~$X/\mathsf G$ alors~$V'$ (resp.~$W'$) est un voisinage ouvert de~$x$ (resp.~$y$), et~$V'\cap W'=\emptyset$ ; par conséquent,~$X/\mathsf G$ est séparé.

\trois{xsurgloccomp} {\em Le quotient~$X/\mathsf G$ est localement compact et~$X\to X/\mathsf G$ est compacte.} Soit~$x\in X/\mathsf G$ ; il correspond à une orbite~$\mathsf G.\xi$ de~$\mathsf G$. Choisissons un voisinage compact~$\Omega$ de~$\xi$ ; le sous-ensemble~$\mathsf G.\Omega$ de~$X$ est compact (en tant qu'image de~$\mathsf G\times \Omega$), stable sous~$\mathsf G$, et est un voisinage de~$\mathsf G.\xi$ ; son image~$\Omega'$ sur~$X/\mathsf G$ est alors un voisinage de~$x$ qui est compact puisque~$\mathsf G.\Omega$ est compact et~$X/\mathsf G$ séparé ; de plus, l'image réciproque de~$\Omega'$ sur~$X$ est par construction égale à~$\mathsf G.\Omega$, et est en particulier compacte. On déduit alors de~\ref{condflechecomp} que~$X/\mathsf G$ est localement compact et que~$X\to X/\mathsf G$ est compacte. 

\trois{flechequotouv} {\em La flèche quotient~$X\to X/\mathsf G$ est ouverte.} Soit~$U$ un ouvert de~$X$ et soit~$V$ son image sur~$X/\mathsf G$. L'image réciproque de~$V$ sur~$X$ coïncidant avec l'ouvert~$\mathsf G.U$ de~$X$, le sous-ensemble~$V$ de~$X/\mathsf G$ en est une partie ouverte.

\trois{xsurhxsurg} {\em Remarque.} Soit~$\mathsf H$ un sous-groupe compact de~$\mathsf G$ ; de ce qui précède (appliqué aux flèches quotients~$X\to X/\mathsf G$ et~$X\to X/\mathsf H$) et de la surjectivité de~$X\to X/\mathsf H$, on déduit que~$X/\mathsf H\to X/\mathsf G$ est compacte et ouverte.

\deux{interconfin} Soit~$X$ un espace topologique localement connexe, soit~$V$ une partie fermée de~$X$ dont le bord est fini, et soit~$Y$ une partie connexe de~$X$ dont l'intersection avec~$V$ est localement connexe. Sous ces hypothèses,~$V\cap Y$ a un nombre fini de composantes connexes. 

En effet, si~$Y\subset V$ c'est évident ; supposons maintenant que~$Y$ n'est pas contenu dans~$V$, et soit~$Z$ une composante connexe de~$Y\cap V$ ; comme~$Z$ est fermée dans~$Y\cap V$, elle est fermée dans~$Y$. 

Supposons que~$Z\cap \partial V=\emptyset$ ; dans ce cas,~$Z$ serait un ouvert de~$Y\cap (V\setminus \partial  V)$, donc un ouvert de~$Y$ ; par connexité de~$Y$, on aurait~$Z=Y$ et donc~$Y\subset V$, ce qu'on a exclu. Il s'ensuit que~$Z$ rencontre~$\partial V$ ; comme celui-ci est fini,~$Y\cap V$ a un nombre fini de composantes connexes. 

\deux{connarcinj} Si~$X$ est un espace topologique séparé et connexe par arcs, il est connexe par arcs {\em injectifs}. 

\section{L'arbre des boules d'un espace ultramétrique~$E$}\label{EXWB}

\subsection*{Définition de l'espace~$\aleph(E)$ des boules de~$E$}

\deux{espmetrult} Soit~$(E,d)$ un espace métrique ; on suppose que la distance~$d$ est {\em ultramétrique}. Pour tout~$a\in E$ et tout~$r>0$ on note~$\bfer a r~$ (resp.~$\bouv a r$) la boule fermée (resp. ouverte) de~$E$ de centre~$a$ et de rayon~$r$. 

\trois{relboules} On se donne un ensemble de symboles~$\zeta_{a,r}$ indexé par~$E\times\RR_+$ ; on vérifie que la relation~$\sch R$ pour laquelle~$\zeta_{a,r}\sch R\zeta_{b,s}$ si et seulement si~$r=s$ et ~$d(a,b)\leq r$ est une relation d'équivalence, et l'on note~$\aleph(E)$ l'ensemble quotient. La relation «~$\zeta_{a,r}\leq \zeta_{b,s}$ si et seulement si~$r\leq s$ et~$d(a,b)\leq s$ » passe au quotient et induit sur~$\aleph (E)$ une relation d'ordre (partiel).

\medskip
Soit~$x\in \aleph(E)$ ; on peut l'écrire~$\zeta_{a,r}$ pour un certain~$a$ et un certain~$r$. Il résulte de la définition que le réel~$r$ ne dépend que de~$x$, et pas de l'écriture choisie. On l'appelle le {\em rayon} de~$x$,  et on le note~$\rho(x)$ ; l'application~$\rho$ est croissante. Si~$a\in E$ alors~$a\leq x$ si et seulement si~$x=\zeta_{a,\rho(x)}$ ; si~$a\leq x$ alors l'ensemble des éléments de~$E$ majorés par~$x$ est égal à~$\bfer a {\rho(x)}$.

\trois{einclaleph} L'application~$a\mapsto \zeta_{a,0}$ définit une bijection de~$E$ sur l'ensemble des points de rayon nul de~$\aleph(E)$. 

\trois{cofindesc} Soit~$a\in E$ et soit~$r\geq 0$. Si~$y$ est un point de~$\aleph (E)$ qui majore~$\zeta_{a,r}$ alors~$y=\zeta_{a,\rho(y)}$. Il s'ensuit que si 
$F$ est une {\em chaîne}, c'est-à-dire une partie totalement ordonnée de~$\aleph(E)$ alors~$\rho$ induit une injection croissante de~$F$ dans~$\RR_+$. Nous dirons qu'une chaîne~$F$ est {\em ouverte saturée} (resp. {\em fermée saturée}) si~$\rho(F)$ est de la forme~$]r;+\infty[$ (resp.~$[r;+\infty[$) pour un certain~$r\geq 0$.

\trois{lienvraiesboules}
{\em Lien avec les vraies boules de~$E$.} 
Si~$a\in E$ et si~$r\geq 0$, il résulte des définitions que la boule fermée~$B(a,r)$ ne dépend que de~$\zeta_{a,r}$. On vérifie immédiatement
que l'application~$\zeta_{a,r}\mapsto B(a,r)$ est croissante, mais elle n'est pas injective en général. En effet, 
on a
$B(a,r)=B(b,s)$ si et seulement si
$${\rm(i)}\;\;\;d(a,b)\leq \min(r,s)$$
et
$${\rm (ii)}\;\;\;\forall c\in E, \; (d(a,c)\leq r\iff d(a,c)\leq s).$$
 
La condition~(i)
équivaut à la conjonction des
deux égalités~$\zeta_{a,r}=\zeta_{b,r}$ et~$\zeta_{a,s}=\zeta_{b,s}$, et implique que~$\zeta_{a,r}$ et~$\zeta_{b,s}$ 
sont comparables ; ils sont
alors égaux si et seulement si~$r=s$. Mais la condition~(ii) signifie seulement que l'application~$d(a,.)$
ne prend aucune valeur dans~$\left]\min(r,s); \max(r,s)\right]$ ; cela peut se produire
même si~$r\neq s$, c'est-à-dire même si~$\min(r,s)<\max(r,s)$. Notons toutefois que si~$\max(r,s)$ appartient
à l'image de~$d(a,.)$ alors dire que~$d(a,.)$ ne rencontre pas $\left]\min(r,s); \max(r,s)\right]$ signifie 
que $\min(r,s)=\max(r,s)$, c'est-à-dire que~$r=s$. 

En conséquence, si~$B(a,r)=B(b,s)$
{\em et si~$\max(r,s)$ appartient
à l'image de~$d(a,.)$}
alors~$\zeta_{a,r}=\zeta_{b,s}$.

\subsection*{Une compactification partielle : l'espace~$\beth(E)$ des chaînes ouvertes saturées de~$\aleph(E)$}

\medskip
On note~$\beth(E)$ l'ensemble constitué des chaînes ouvertes saturées de~$\aleph(E)$, muni de la relation d'ordre {\em opposée} à celle de l'inclusion. Pour des raisons psychologiques, on préfèrera souvent penser à~$\beth(E)$ comme à un ensemble de points limite de~$\aleph(E)$ plutôt que comme à un sous-ensemble de~$\sch P(\aleph(E))$ ; il arrivera donc que l'on parle de la chaîne ouverte saturée~$F$ de~$\aleph(E)$ qui {\em correspond} à un point~$x$ de~$\beth(E)$. 

\medskip
Si~$x=\zeta_{a,r}$ est un élément de~$\aleph(E)$ l'ensemble des majorants stricts~$\aleph(E)^{>x}$ de~$x$ dans~$\aleph(E)$ est égal à~$\{\zeta_{a,s}\}_{s>r}$ ; c'est une chaîne ouverte saturée. 

\medskip
Soient~$x$ et~$y$ deux points de~$\aleph(E)$. Si~$x\leq y$ alors~$\aleph(E)^{>y}\subset \aleph(E)^{>x}$. Réciproquement, supposons que~$\aleph(E)^{>y}$ soit contenue dans~$\aleph(E)^{>x}$. Écrivons~$x=\zeta_{a,r}$ et~$y=\zeta_{b,s}$. L'inclusion~$\aleph(E)^{>y}\subset \aleph(E)^{>x}$ assure que pour tout~$t>s$, on a~$\zeta_{a,r}\leq \zeta_{b,t}$, c'est-à-dire~$r\leq t$ et~$d(a,b)\leq t$ ; il s'ensuit que~$r\leq s$ et~$d(a,b)\leq s$ ; par conséquent,~$x\leq y$.

\medskip
Si~$x\in \aleph(E)$, vérifions que le point~$x$ est la borne inférieure de~$\aleph(E)^{>x}$ (ce qui montrera en particulier que~$x$ est uniquement déterminé par~$\aleph(E)^{>x}$). C'en est clairement un minorant. Donnons-nous maintenant un minorant~$y$ de~$\aleph(E)^{>x}$ ; comme ce dernier est une chaîne ouverte saturée, il n'a pas de plus petit élément et ne contient donc pas~$y$ ; il s'ensuit que~$\aleph(E)^{>x}\subset \aleph(E)^{>y}$, et donc que~$y\leq x$ ({\em cf. supra})  ; on a donc bien ainsi~$x=\inf \aleph(E)^{>x}$. 

\medskip
Il résulte de ce qui précède que~$x\mapsto \aleph(E)^{>x}$ définit un isomorphisme d'ensembles ordonnés de~$\aleph(E)$ sur un sous-ensemble de~$\beth(E)$, sous-ensemble que l'on {\em identifiera} désormais par ce biais à~$\aleph(E)$ ; en vertu de~\ref{einclaleph}, on considérera également~$E$ lui-même comme un sous-ensemble de~$\beth(E)$. 

\medskip
Si~$F$ est une chaîne ouverte saturée de~$\aleph(E)$ et si~$\rho(F)=]r;+\infty[$ on définit le {\em rayon} de~$F$ comme étant égal à~$r$ ; si~$F$ est de la forme~$\aleph(E)^{>x}$ pour un certain~$x\in \aleph(E)$, le rayon de~$F$ est égal à celui de~$x$ : on a ainsi prolongé la fonction croissante~$\rho:  \aleph(E)\to \RR_+$ en une fonction strictement croissante~$\beth(E)\to \RR_+$ que l'on note encore~$\rho$.

\medskip
Il existe une notion de chaîne ouverte saturée (resp. fermée saturée) de~$\beth(E)$ : il suffit de décalquer {\em verbatim} les définitions données sur~$\aleph(E)$ ; toute chaîne ouverte saturée (resp. fermée saturée) de~$\aleph(E)$ est encore une chaîne ouverte saturée (resp. fermée saturée) de~$\beth(E)$. 

\trois{lemchainesat} Soit~$F$ une chaîne saturée de~$\aleph(E)$ et soit~$y\in F$. Les sous-ensembles~$F^{> y}$ et~$\aleph(E)^{>y}$ sont deux chaînes ouvertes saturées dont l'image par~$\rho$ est égale à~$]\rho(y);+\infty[$ ; comme~$F^{>y}\subset \aleph(E)^{>y}$, on a~$F^{>y}= \aleph(E)^{>y}$. 

\trois{chainesat} Soit~$x\in \beth(E)$ et soit~$F$ la chaîne ouverte saturée de~$\aleph(E)$ qui lui correspond. Nous allons montrer que~$\beth(E)^{>x}$ coïncide avec~$F$, et est en particulier une chaîne ouverte saturée de~$\beth(E)$ contenue dans~$\aleph(E)$. 

\medskip
{\em Montrons que~$F\subset \beth(E)^{>x}$}. Soit~$y\in F$ ; d'après le~\ref{lemchainesat}, l'ensemble~$F^{>y}$ est égal à~$\aleph(E)^{>y}$ ; on a donc~$\aleph(E)^{>y}\subset F$, et cette inclusion est stricte puisque~$y$ n'appartient pas à~$\aleph(E)^{>y}$. Comme les points de~$\beth(E)$ correspondant aux chaînes saturées~$\aleph(E)^{>y}$ et~$F$ sont respectivement~$y$ et~$x$, il vient~$y>x$ ; par conséquent,~$F\subset \beth(E)^{>x}$.

\medskip
{\em Montrons que~$\beth(E)^{>x}\subset F$}. Soit~$y\in \beth(E)^{>x}$ et soit~$G$ la chaîne ouverte saturée de~$\aleph(E)$ qui lui correspond. Comme~$y>x$, la chaîne~$G$ est contenue strictement dans~$F$, et possède de ce fait une borne inférieure~$z$ dans~$F$. La chaîne~$G$ est alors égale à~$F^{>z}$, et donc à~$\aleph(E)^{>z}$ en vertu de~\ref{lemchainesat}. Le point de~$\beth(E)$ qui correspond à~$G$ est par conséquent le point~$z$ de~$\aleph(E)$, ce qui signifie que~$y=z$, entraîne que~$y\in F$, et achève la démonstration. 

\trois{chainedesc} Soit~$F$ une chaîne ouverte saturée de~$\beth(E)$ et soit~$y\in F$. La chaîne~$F$ n'ayant pas de plus petit élément,  il existe~$z\in F$ tel que~$z<y$ ; par conséquent,~$y\in \beth(E)^{>z}$, lequel est contenu dans~$\aleph(E)$ d'après le~\ref{chainesat}. Ainsi,~$F\subset \aleph(E)$ ; il s'ensuit que~$F$ est une chaîne ouverte saturée de~$\aleph(E)$. Si~$x$ est le point qui lui correspond dans~$\beth(E)$, il résulte de ~\ref{chainesat} que~$F=\beth(E)^{>x}$. 

\medskip
{\em Le point~$x$ est égal à la borne inférieure de~$F=\beth(E)^{>x}$.} En effet,~$x$ est un minorant de~$\beth(E)^{>x}$. Soit~$y$ un minorant de~$\beth(E)^{>x}$ ; comme~$F$ n'a pas de plus petit élément, ~$y\notin F$ ; il s'ensuit que~$\beth(E)^{>x}\subset \beth(E)^{>y}$, ce qui signifie, compte-tenu de~\ref{chainesat}, que~$y\leq x$ et achève la preuve.

\deux{recapchaines} Notons quelques conséquences de ce qui précède. 

\trois{bijchain} L'application ~$x\mapsto \beth(E)^{>x}$ (resp.~$x\mapsto \beth(E)^{\geq x}$) établit un isomorphisme d'ensembles ordonnés entre~$\beth(E)$ et l'ensemble de ses chaînes saturées ouvertes (resp. fermées), dont la réciproque envoie une chaîne sur sa borne inférieure (resp. son plus petit élément) ; les chaînes saturées ouvertes de~$\beth(E)$ sont par ailleurs toutes contenues dans~$\aleph(E)$. 

\trois{satsat} Si~$F$ est une chaîne saturée de~$\beth(E)$ et si~$y\in F$ alors~$F^{>y}=\beth(E)^{>y}$ : on le voit en remarquant que ce sont deux chaînes saturées comparables pour l'inclusion et ayant même image par~$\rho$.  

\trois{minim} Si~$z\in \beth(E)-\aleph(E)$ ou si~$\rho(z)=0$, alors~$z$ est un élément minimal de~$\beth(E)$. En effet, soit~$t\leq z$. 

\medskip
La chaîne ouverte saturée~$\beth(E)^{>t}$ est contenu dans~$\aleph(E)$, et son image par~$\rho$ est égale à~$]\rho(t);+\infty[$ qui ne contient pas~$0$. Il s'ensuit que~$z$ ne peut appartenir à~$\beth(E)^{>t}$, et l'on a donc~$t=z$. 

\trois{remysupx} Soient~$x$ et~$y$ deux éléments de~$\beth(E)$ tels que~$y>x$. Le point~$y$ appartient alors à la chaîne ouverte saturée~$\beth(E)^{>x}$, laquelle est contenue dans~$\aleph(E)$ ; on a en particulier~$y\in \aleph(E)$. De l'inégalité~$\rho(y)>\rho(x)\geq 0$ il vient~$\rho(y)>0$ ; par conséquent,~$y$ est de la forme~$\zeta_{a,r}$ où~$a\in E$ et où~$r>0$.

\deux{notinterv} Il résulte de~\ref{bijchain} que si~$x$ et~$y$ sont deux éléments de~$\beth(E)$, il existe une chaîne saturée les contenant si et seulement si ils sont comparables. Supposons que ce soit le cas, et que l'on ait par exemple~$x\leq y$ ; et soit~$F$ une chaîne saturée les contenant. En vertu de~\ref{satsat}, l'ensemble des points~$z$ de~$F$ tels que~$x\leq z\leq y$ est égal à l'ensemble des points~$z$ de~$\beth(E)$ tels que~$x\leq z\leq y$ ; il ne dépend donc pas de~$F$, et nous le noterons~$[x;y]$. On utilisera de même les notations~$]x;y[, ]x;y], [x;y[$ dans un sens évident. 

\medskip
Si~$x\in \beth (E)$ on désignera par~$[x;\infty[$ (resp.~$]x;\infty[$) la chaîne saturée~$\beth(E)^{\geq x}$ (resp.~$\beth(E)^{>x}$). 

\medskip
Remarquons que si~$x$ et~$y$ sont deux points de~$\beth(E)$ et si~$]y;\infty[\subset [x;\infty[$ alors~$y\in [x;\infty[$ ; c'est en effet évident si~$y=x$ ; et sinon, on a~$]y;\infty[\subset ]x;\infty[$, ce qui équivaut d'après~\ref{bijchain} à l'inégalité~$x\leq y$. 

\deux{interchsat} On dira qu'une
chaîne~$F$ de~$\beth(E)$ est {\em convexe} si~$[x;y]\subset F$ pour tout couple~$(x,y)$ d'éléments de~$F$. Toute chaîne saturée est convexe, et l'intersection de deux chaînes convexes est une chaîne convexe. 

\deux{wedgearb} Soient~$x$ et~$y$ deux points de~$\beth(E)$.  Nous allons montrer que l'intersection des chaînes saturées fermées~$[x;\infty[$ et~$[y;\infty[$ est une chaîne saturée fermée. 

\medskip
{\em Cette intersection est non vide.} Choisissons en effet~$x'\in ]x;\infty[$ et~$y'\in ]y;\infty[$ ; les points~$x'$ et~$y'$ peuvent respectivement s'écrire~$\zeta_{a,r}$ et~$\zeta_{b,s}$ pour~$a,b,r$
 et~$s$ convenables (\ref{remysupx}). Si l'on pose~$R=\max(r,s,,d(a,b))$ alors~$\zeta_{a,R}$ majore~$x'$ et~$y'$, et {\em a fortiori}~$x$ et~$y$ ; par conséquent, ~$[x;\infty[\cap[y;\infty[\neq \emptyset$. 

\medskip
Si~$z$ est un point de~$\beth(E)$ majorant~$x$ et~$y$ alors~$[z;\infty[\subset [x;\infty[\cap [y;\infty[$, ce qui montre que~$[x;\infty[\cap [y;\infty[$ est une chaîne saturée. Enfin, si~$z$ est un point de~$\beth(E)$ tel que~$]z;\infty[\in [x;\infty[\cap [y;\infty[$ alors~$z\in [x;\infty[\cap [y;\infty[$ (fin du~\ref{notinterv}) ; il en découle que la chaîne saturée~$[x;\infty[\cap [y;\infty[$ est fermée, ce qui achève la preuve. 

\medskip
On notera~$x\wedge y$ le plus petit élément de~$[x;\infty[\cap [y;\infty[$ ; on peut également le définir comme la borne supérieure de~$\{x,y\}$.

\deux{introtopoarbree} Si~$a\in E$ et si~$r>0$, on note~$\bbfer a r$ l'ensemble~$\beth(E)^{\leq \zeta_{a,r}}$. Remarquons que~$\bbfer a r$ ne dépend que de~$\zeta_{a,r}$, et que~$\zeta_{a,r}$ se retrouve à partir de~$\bbfer a r$ : c'est son plus grand élément. Le réel~$r$ est donc égal au maximum de~$\rho$ sur~$\bbfer a r$. 

\trois{boulenaiv} Il découle immédiatement de la définition que si~$a\in E$ et si~$r>0$ alors~$\bbfer a r\cap E=\bfer a r$. 

\trois{incluboule} Soient~$a$ et~$b$ dans~$E$ et~$r$ et~$s$ deux réels strictement positifs. L'inclusion~$\bbfer b s \subset \bbfer a r$ équivaut à l'inégalité~$\zeta_{b,s}\leq \zeta_{a,r}$, qui est vérifiée si et seulement si~$s\leq r$ et~$d(a,b)\leq r$ ; on a donc~$\bbfer b s= \bbfer a r$ si et seulement si~$r=s$ et~$d(a,b)\leq r$, c'est-à-dire encore si et seulement si~$\bbfer b s \subset \bbfer a r$ et~$r=s$.

\trois{comparboules} Soient~$a$ et~$b$ deux éléments de~$E$ et soit~$r$ et~$s$ deux réels strictement positifs. {\em Supposons que~$\bbfer a r\cap \bbfer b s$ soit non vide}. Il existe alors~$x$ tel que~$x\leq \zeta_{a,r}$ et~$x\leq \zeta_{b,s}$. Les points~$\zeta_{a,r}$ et~$\zeta_{b,s}$ appartenant dès lors tous deux à la chaîne~$[x;\infty[$, ils sont comparables, d'où il découle que~$\bbfer a r\subset \bbfer b s~$ ou~$\bbfer b s\subset \bbfer a r$.

\deux{boulouv}  Si~$a\in E$ et si~$r>0$, on note~$\bbouv a r$ la réunion des~$\bbfer a t$ pour~$0<t<r$. Remarquons que~$r$ est égal à~$\sup\limits_{x\in \bbouv a r} \rho(x)$. 

\trois{boulenaivouv} Il découle immédiatement de la définition que si~$a\in E$ et si~$r>0$ alors~$\bbouv a r\cap E=\bouv a r$. 

\trois{incluboulouv} Soient~$a$ et~$b$ dans~$E$ et~$r$ et~$s$ deux réels strictement positifs. On déduit de~\ref{incluboule} que l'inclusion~$\bbouv b s \subset \bbouv a r$ est vérifiée si et seulement si~$s\leq r$ et~$d(a,b)<r$ ; on a donc~$\bbouv b s= \bbouv a r$ si et seulement si~$r=s$ et~$d(a,b)<r$, c'est-à-dire encore si et seulement si~$\bbouv b s \subset \bbouv a r$ et~$r=s$.  

\trois{comparboulouv} Soient~$a$ et~$b$ deux éléments de~$E$ et soit~$r$ et~$s$ deux réels strictement positifs. {\em Supposons que~$\bbouv a r\cap \bbouv b s$ soit non vide}. On déduit de~\ref{comparboules} que~$\bbouv  a r\subset \bbouv b s~$ ou~$\bbouv b s\subset \bbouv a r$.

\deux{wedgeboulouv} Soient~$x$ et~$y$ deux points de~$\beth(E)$ tels que~$(x\wedge y)>x$ et~$(x\wedge y)>y$. D'après~\ref{remysupx}, il existe~$a\in E$ et~$r>0$ tels que~$x\wedge y=\zeta_{a,r}$. Choisissons un point~$x'$ quelconque sur~$]x;x\wedge y[$. Comme~$x'>x$, il est égal à~$ \zeta_{b,s}$ pour un certain~$b\in E$ et un certain~$s>0$ (toujours en vertu de~\ref{remysupx}). Comme~$\zeta_{b,s}< \zeta_{a,r}$ on a~$b\in \bfer a r$,~$s<r$, et~$\zeta_{a,r}=\zeta_{b,r}$. Si~$z\in [x;x\wedge y[$ il existe un élément de~$[z;x\wedge y[$ qui appartient aussi à~$[x';x\wedge y[$ et est donc égal à~$\zeta_{b,t}$ pour un certain~$t\in [s;r[$ ; par conséquent,~$z\in \bbouv b r$. Ainsi,~$[x; x\wedge y[\subset \bbouv b r$. Et il existe de même~$c\in \bfer a r$ tel que~$[y;x\wedge y[\subset \bbouv c r$. 

\medskip
{\em L'intersection~$\bbouv b r\cap \bbouv c r$ est vide.} En effet, il suffit de vérifier par symétrie que l'on ne peut avoir~$\bbouv c r\subset \bbouv b r$. Supposons que ce soit le cas. Les points~$x$ et~$y$ appartenant tous deux à~$\bbouv  b r$, il existerait~$r'$ et~$r''$ dans~$]0;r[$ tels que~$x\leq \zeta_{b,r'}$ et~$y\leq \zeta_{b,r''}$ ; mais cela impliquerait que~$\zeta_{b,\max(r',r'')}$ majore~$x$ et~$y$, ce qui est absurde puisque~$\zeta_{b,\max(r',r'')}<\zeta_{b,r}=x\wedge y$. 

\deux{interconvbb} Soit~$F$ une chaîne convexe de~$\beth(E)$, soit~$a\in E$ et soit~$r>0$. 

\trois{complbf} Posons~$G=F\cap \bbfer a r$ et~$H=F\setminus G$ ; nous allons montrer que l'on est dans l'un (et un seul) des trois cas suivants : 

\medskip
$\bullet$~$G=F$ et~$H=\emptyset$ ; 

$\bullet$~$G=\emptyset$ et~$H=F$ ; 

$\bullet$~$G$ et~$H$ sont tous deux non vides,~$\zeta_{a,r}$ appartient à~$F$, et~$G$ (resp.~$H$) est l'ensemble des éléments~$z$ de~$F$ tels que~$z\leq \zeta_{a,r}$ (resp.~$z>\zeta_{a,r}$).

\medskip
Il suffit de vérifier que si~$G$ et~$H$ sont non vides alors~$\zeta_{a,r}\in F$. Supposons donc qu'il existe~$x\in G$ et~$y\in H$. Comme~$x\leq \zeta_{a,r}$ et comme~$y\notin G$, on ne peut avoir~$y\leq x$. Par conséquent~$y>x$ et~$[x;y]\subset [x;\infty[$. D'autre part,~$ \zeta_{a,r}$ appartient aussi à~$[x;\infty[$. Puisque~$y\notin G$, on ne peut avoir~$y\leq \zeta_{a,r}$ ; dès lors~$\zeta_{a,r}\leq y$ et~$\zeta_{a,r}\in [x;y]\subset F$, ce qu'on souhaitait établir.

\trois{complbo} Posons~$G'=F\cap \bbouv a r$ et~$H'=F\setminus G'$ ; nous allons montrer que l'on est dans l'un (et un seul) des trois cas suivants : 

\medskip
$\bullet$~$G'=F$ et~$H'=\emptyset$ ; 

$\bullet$~$G'=\emptyset$ et~$H'=F$ ; 

$\bullet$~$G'$ et~$H'$ sont tous deux non vides,~$\zeta_{a,r}$ appartient à~$F$, et~$G'$ (resp.~$H'$) est l'ensemble des éléments~$z$ de~$F$ tels que~$z< \zeta_{a,r}$ (resp.~$z\geq \zeta_{a,r}$). 

\medskip
En effet, supposons que~$G'$ et~$H'$ soient non vides. Vérifions tout d'abord que~$\zeta_{a,r}\in F$. Choisissons~$x\in G'$ et~$y\in H'$. Comme~$x\leq \zeta_{a,t}$ pour un certain~$t<r$ et comme~$y\notin G'$, on ne peut avoir~$y\leq x$. Par conséquent~$y>x$ et~$[x;y]\subset [x;\infty[$. D'autre part,~$ \zeta_{a,s}$ appartient aussi à~$[x;\infty[$ pour tout~$s\geq t$ ; en particulier,~$\zeta_{a,r}\in [x;\infty[$. Puisque~$y\notin G'$, on ne peut avoir~$y\leq \zeta_{a,s}$ si~$s<r$ ; dès lors~$\zeta_{a,s}\in [x;y[$ pour tout~$s\in [t;r[$, d'où il découle que~$\zeta_{a,r}\in [x;y]\subset F$, ce qu'on souhaitait établir. 

\medskip
Il reste à s'assurer que~$G'$ est exactement l'ensemble des éléments~$z$ de~$F$ tels que~$z< \zeta_{a,r}$. Il est clair que tout élément de~$G'$ est strictement majoré par~$\zeta_{a,r}$ ; réciproquement, si~$z$ appartient à~$F$ et est strictement inférieur à~$\zeta_{a,r}$ il existe~$s\in [t;r[$ tel que~$\zeta_{a,s}\in ]z;\zeta_{a,r}[$. On par conséquent,~$z\leq \zeta_{a,s}$ ; dès lors,~$z\in \bbouv a r$, et partant à~$G'$. 

\deux{defthetaar} Soit~$a\in E$ et soit~$r>0$. On notera~$\Theta(a,r)$ le quotient de~$\bfer a r$ obtenu en identifiant deux points dont la distance est strictement inférieure à~$r$. Si~$b$ et~$c$ sont deux points de~$\bfer a r$ alors~$\bbouv b r=\bbouv c r$ si et seulement si ils ont même image dans~$\Theta(a,r)$. On appellera {\em système complet de centres} de~$\bfer ar$ tout système de représentants de la relation d'équivalence définissant~$\Theta(a,r)$. 

\medskip
Il découle immédiatement des définitions que les assertions suivantes sont équivalentes :

\medskip
i)~$\bfer a r=\bouv a r$ ; 

ii)~$\Theta(a,r)$ est un singleton ; 

iii)~$\{a\}$ est un système complet de centres de~$\bfer a r$.

\deux{topoarbree} On munit~$\beth (E)$ de la topologie engendrée (par unions quelconques et intersections finies) par les parties~$U$ de l'une des formes suivantes : 

\medskip
$\bullet$ ~$U=\beth(E)-\bbfer a r$ avec~$a\in E$ et~$r>0$ ; 

$\bullet$~$U=\bbouv a r~$ avec~$a\in E$ et~$r>0$. 

\medskip
En vertu de~\ref{boulenaiv} et~\ref{boulenaivouv}, la topologie de~$E$ induite par celle de~$\beth(E)$ est sa topologie d'espace métrique.

\deux{basevoisal} Soit~$x\in E$. 

\trois{basevoisgen} Il découle de~\ref{comparboules} et~\ref{comparboulouv} que si~$x\in \beth(E)$ il possède une base de voisinages de la forme~$\bbouv a r- \coprod\limits_{i\in I} \bbfer {a_i}{r_i}$, où~$I$ est fini et où~$\bbfer {a_i}{r_i}$ est contenue dans~$\bbouv a r$ pour tout~$i$.

\trois{voiszetaar} Supposons que~$x=\eta_{b,s}$ pour un certain~$b\in E$ et un certain~$s>0$. Si~$a\in E$ et si~$r$ est un réel strictement positif tel que~$x\in \bbouv a r$ alors~$r>s$ et~$\bbouv a r=\bbouv b r$ ; par conséquent,~$x$ possède une base de voisinages de la forme~$\bbouv b r- \coprod\limits_{i\in I} \bbfer {a_i}{r_i}$, où~$r>s$, où~$I$ est fini et où~$\bbfer {a_i}{r_i}\subset \bbouv b r$ pour tout~$i$.

\medskip
Faisons maintenant deux remarques. 

\medskip
$\bullet$ Soit~$c\in \bbfer b s$, soit~$n\geq 1$ et soient~$c_1,\ldots, c_n\in \bbfer b s$ tels que~$d(c,c_i)<s$ pour tout~$i$. Si~$s_1,\ldots, s_n$ sont des réels appartenant
à~$]0;s[$ et si~$$R=\max (s_1,\ldots, s_n,d(c,c_1),\ldots, d(c,c_n))$$ alors~$R<s$,~$\bbfer c R\subset \bbfer b s$ et~$\bbfer c R$ contient chacune des~$\bbfer {c_i} {s_i}$. 

$\bullet$ Soit~$r>ss$, soit~$n\geq 1$, soient~$c_1,\ldots, c_n\in \bbouv b r-\bbfer b s~$ et soient~$r_1,\ldots, r_n$ des réels strictement
supérieurs à~$s$ tels que~$d(c_i,b)>r_i$ pour tout~$i$. Posons~$R=\min r_i$ ; on a alors~$r>R>s$, et~$\bbouv b R$ ne contient aucune des~$\bbfer {c_i} {r_i}$. 

\medskip
Fixons un système complet~$\sch S$ de centres de~$\bbfer b s$. On déduit des deux remarques ci-dessus que le point~$x=\eta_{b,s}$ possède une base de voisinages de la forme ~$\bbouv b r- \coprod\limits_{i\in I} \bbfer {a_i}{r_i}$ où~$r>s$, où~$I$ est fini, où~$r_i<s$ pour tout~$i$, et où les~$a_i$ sont des éléments deux à deux distincts de~$\sch S$. 

\medskip
Notons un cas particulier, celui où~$\bfer b s=\bouv b s$, c'est-à-dire où~$\{b\}$ est un système complet de centres de~$\bfer b s$ : le point~$x=\eta_{b,s}$ possède alors une base de voisinages de la forme ~$\bbouv b r-  \bbfer b{r'}$ où~$r>s>r'$.

\trois{voisefeuille} Supposons que~$\rho(x)=0$ ou que~$x\in \beth(E)-\aleph(E)$, et soient~$a\in E$ et~$r>0$ tels que~$x\notin \bbfer a r$. 

\medskip
Posons~$y=x\wedge \zeta_{a,r}$. Comme~$y\geq \zeta_{a,r}$ il est égal à~$\zeta_{a,R}$ pour un certain~$R\geq r$. Comme~$y\geq x$ et comme~$x\notin \bbfer a r$, on a nécessairement~$R>r$ ; et l'hypothèse faite sur~$x$ assure d'autre part que~$x\neq \zeta_{a,R}$ et donc que~$x<\zeta_{a,R}$. On déduit alors de~\ref{wedgeboulouv} qu'il existe~$b\in E$ tel que~$x\in \bbouv b R$ et tel que 
$\bbouv b R\cap \bbouv a R=\emptyset$ ; l'intersection de~$\bbouv b R$ et~$\bbfer a r$ est {\em a fortiori} vide. 

\medskip
Il découle de ce qui précède, de~\ref{comparboules} et de~\ref{comparboulouv} que~$x$ possède une base de voisinages de la forme~$\bbouv bR$. 

\deux{separarb} {\em L'espace topologique~$\beth(E)$ est séparé.} En effet, soient~$x$ et~$y$ deux points distincts de~$\beth(E)$ ; nous allons chercher à les séparer par deux ouverts disjoints. 

\medskip
{\em  Le cas où~$x$ et~$y$ sont 
comparables.} Supposons par exemple que~$x<y$ et choisissons~$z\in ]x;y[$. Comme~$z>x$ il est de la forme~$\zeta_{a,r}$ avec~$a\in E$ et~$r>0$ ; on a dès lors~$y=\zeta_{a,R}$ pour un certain~$R>r$. Soit~$s\in ]r;R[$. Les ouverts~$\bbouv a s$ et~$\beth(E)-\bbfer a s$ sont disjoints, le premier contient~$x$ et le second~$y$. 

\medskip
{\em Le cas où~$x$ et~$y$ ne sont pas comparables.} Posons~$z=x\wedge y$. On a alors~$z>x$ et~$z>y$.  On déduit de~\ref{wedgeboulouv} qu'il existe~$a$ et~$b$ dans~$E$ et~$r>0$ tels que~$x\in \bbouv a r$, tels que~$y\in \bbouv b r$, et tels que~$\bbouv a r\cap \bbouv b r=\emptyset$.

\deux{topchaines} Soit~$x\in \beth(E)$ ; l'application~$\rho$ induit une bijection~$\rho_x$ entre~$\beth(E)^{\geq x}$ et~$[\rho(x);+\infty[$ ; nous allons montrer qu'il s'agit d'un homéomorphisme.

\trois{rhoxcont} {\em L'application~$\rho_x$  est continue.} Soit~$r\in ]\rho(x);+\infty[$ ; il suffit de vérifier que~$\rho_x^{-1}(]r;+\infty[)$ et~$\rho_x^{-1}([\rho(x);r[)$ sont ouverts. Choisissons~$s\in ]\rho(x);r[$. On peut écrire~$s=\rho(y)$ pour un certain~$y\in ]x;\infty[$. Appartenant à~$]x;\infty[$, le point~$y$ est égal à~$\zeta_{a,s}$ pour un certain~$a\in E$. Il s'ensuit que~$\zeta_{a,r}$ appartient à~$[x;\infty[$, et c'est nécessairement l'unique antécédent de~$r$ par~$\rho_x$. La chaîne~$[x;\infty[$ rencontre~$\bbouv a r$ (elle contient~$\zeta_{a,s}$) et le complémentaire de~$\bbfer a r$ (elle contient~$\rho_x^{-1}(]r;+\infty[)$ qui est constitué de majorants stricts de~$\zeta_{a,r}$). Il découle alors de~\ref{complbf} et~\ref{complbo} que~$$\rho_x^{-1}(]r;+\infty[)= [x;\infty[\setminus \bbfer a r\;{\rm et}\; \rho_x^{-1}([\rho(x); r[)= [x;\infty[\cap \bbouv a r,$$ d'où l'assertion. 

\trois{rhoxouv}  {\em L'application~$\rho_x$  est ouverte.} Soient~$a\in E$ et~$r>0$. Il suffit de s'assurer que~$\rho_x([x;\infty[\cap \bbouv a r)$ et~$\rho_x([x;\infty[-\bbfer a r)$ sont ouverts dans~$[\rho(x); +\infty[$ ; mais c'est une conséquence immédiate de ~\ref{complbf} et~\ref{complbo}.

\deux{arbreecompact} Soit~$a\in E$ et soit~$r>0$ ; nous allons montrer que~$\bbfer a r$ est compact. Pour ce faire, on se donne un ensemble~$I$ muni d'un ultrafiltre~$\sch U$, et une suite~$(x_i)$ d'éléments de~$\bbfer a r$ indexée par~$I$ ; on va prouver qu'elle converge le long de~$\sch U$, ce qui permettra de conclure. Pour tout~$i \in I$, la chaîne fermée saturée~$[x_i;\infty[$ contient~$[\zeta_{a,r};\infty[$. Il s'ensuit que pour tout~$J$ appartenant à~$\sch U$, le sous-ensemble~$F_J:=\bigcap\limits_{i\in J}[x_i; \infty[$ de~$\beth(E)$ est une chaîne fermée saturée contenant~$[\zeta_{a,r}; \infty[$. Comme~$(F_J)_{J\in \sch U}$ est une famille filtrante de chaînes fermées saturées contenant~$[\zeta_{a,r}; \infty[$, la réunion~$F$ des~$F_J$ est une chaîne saturée (ouverte ou fermée) contenant~$\zeta_{a,r}$. Elle possède une borne inférieure (\ref{chainedesc}) ; appelons-la~$x$. Comme~$\zeta_{a,r}\in F$, on a~$x\in \bbfer a r$ ; nous allons maintenant vérifier que~$x$ est une valeur d'adhérence de~$(x_i)$. Soit~$b\in E$ et soit~$s>0$. 

\trois{xappbouv} {\em Supposons que~$x\in \bbouv b s$.} Il existe~$t<s$ tel que~$x\leq \zeta_{b,t}$. Quitte à rapprocher~$t$ de~$s$, on peut supposer que~$x\neq \zeta_{b,t}$. Comme~$F$ est égale ou bien à~$[x;\infty[$ ou bien à~$]x; \infty[$, on a~$\zeta_{b,t}\in F$. Il existe dès lors~$J\in \sch U$ tel que~$\zeta_{b,t}\in F_J$. Par définition de ce dernier, cela signifie que~$x_i\leq \zeta_{b,t}$ pour tout~$i\in J$ ; par conséquent,~$x_i\in \bbouv b s$ pour tout~$i\in J$. 

\trois{xappcompbf} {\em Supposons que~$x\in\beth(E)- \bbfer b s$.} Soit~$J$ l'ensemble des indices~$i\in I$ tels que~$x_i\leq \zeta_{b,s}$. Si~$J$ appartenait à~$\sch U$ la chaîne~$[\zeta_{b,s}; \infty[$ serait contenue dans~$F_J$, et {\em a fortiori} dans~$F$, ce qui entraînerait l'inégalité~$x\leq \zeta_{b,s}$, contredisant notre hypothèse sur~$x$. Par conséquent~$I\setminus J\in \sch U$ ; et par définition de~$J$, on a~$x_i\in 
\beth(E)- \bbfer b s$ pour tout~$i\in I\setminus J$, ce qui termine la preuve que~$x_i\to x$ le long de~$\sch U$. 

\deux{bethloccomp} Soit~$x\in \beth(E)$. Choisissons un point~$y\in ]x;\infty[$. Il est de la forme~$\zeta_{a,r}$ avec~$a\in E$ et~$r>0$ (\ref{remysupx}). Comme~$x\leq \zeta_{a,r}$, on a~$x\in \bbouv a{r+1}\subset \bbfer a {r+1}$ ; en vertu de~\ref{arbreecompact} {\em et sq.}, ~$\bbfer a {r+1}$ est compact, et~$\beth(E)$ est donc localement compact. 

\deux{compbferr} Soit~$a\in E$ et soit~$\sch K$ une partie compacte de~$\beth(E)$. On peut recouvrir~$\sch K$ par un nombre fini d'ouverts du type décrit au~\ref{basevoisal} ; il existe {\em a fortiori} une famille finie~$(a_i)$ d'éléments de~$E$ et, pour tout~$i$, un réel~$r_i>0$, tels que~$\sch K\subset \bigcup \bbouv {a_i} {r_i}$. Soit~$r$ un réel strictement supérieur pour tout~$i$ à~$r_i$ ainsi qu'à~$d(a,a_i)$ ; on a alors~$\sch K\subset \bbfer a r$.

\subsection*{Le compactifié d'Alexandrov~$\wbeth E$ de~$\beth(E)$}

\deux{remvide} Remarquons que si~$E$ est vide alors~$\beth(E)$ est vide aussi. {\em On suppose à partir de maintenant que~$ E$ est non vide.} L'espace~$\beth(E)$ n'est alors pas compact. En effet, choisissons~$a\in E$. Si~$\beth(E)$ était compact, il serait d'après~\ref{compbferr} contenu dans~$\bbfer a r$ pour un certain~$r>0$, ce qui est absurde puisque~$\zeta_{a,r+1}$ (par exemple) n'appartient pas à ~$\bbfer a r$. 

\deux{defwbeth} Désignons par~$\wbeth E~$ le compactifié d'Alexandrov de l'espace localement compact et non compact~$\beth(E)$. On peut écrire~$\wbeth E=\beth(E)\cup\{\infty\}$ ; une base de voisinages de~$\infty$ est constituée des parties de la forme~$\wbeth E\setminus \sch K$, où~$\sch K$ est un compact de~$\beth(E)$. 

\medskip
Il résulte de~\ref{compbferr} que si l'on se donne un élément~$a$ de~$E$ alors les parties de la forme~$\wbeth E-\bbfer a r$ où~$r>0$, forment une base de voisinages ouverts de~$\infty$. 

\deux{intervxinfini} Prolongeons la relation d'ordre sur~$\beth(E)$ à~$\wbeth E$ en décrétant que~$\infty$ est son plus grand élément. Si~$x\in \wbeth E$ l'ensemble~$[x;\infty]$ des éléments~$y$ de~$\wbeth E$ tels que~$x \leq y\leq \infty$ est alors égal au singleton~$\{\infty\}$ si~$x=\infty$, et à~$[x;\infty[\bigcup \{\infty\}$ sinon (où~$[x;\infty[$ a le sens donné au~\ref{notinterv}). 

Toute paire ~$\{x,y\}$ d'éléments de~$\wbeth E$ admet une borne supérieure, à savoir~$\infty$ si~$x$ ou~$y$ est égal à~$\infty$, et~$x\wedge y$ sinon ; on pourra donc noter cette borne~$x\wedge y$ dans tous les cas sans risque de confusion. 

\deux{rhohomeowh} Si l'on prolonge~$\rho$ à~$\wbeth E$ en posant~$\rho(\infty)=+\infty$ alors pour tout~$x\in \wbeth E$ l'application~$\rho$ induit une bijection~$[x;\infty]\simeq [\rho(x);+\infty]$. 

\medskip
{\em Cette bijection est un homéomorphisme.} Il suffit, par compacité  de~$[\rho(x);+\infty]$ et par séparation de~$\wbeth E$, de montrer que sa réciproque est continue ; le cas~$x=\infty$ est évident, ce qui permet de supposer que~$x\neq \infty$. On sait que la bijection naturelle~$[\rho(x);+\infty[\simeq [x;\infty[$ est continue ; il reste donc à s'assurer de la continuité en~$+\infty$. Soit~$r$ un réel strictement supérieur à~$\rho(x)$. Le point~$y$ de~$[x;\infty[$ qui lui correspond est strictement supérieur à~$x$, et est donc égal à~$\zeta_{a,r}$ pour un certain~$r>0$. L'ensemble des ouverts de~$\wbeth E$ de la forme~$\wbeth E-\bbfer a R$, pour~$R>r$, est une base de voisinages ouverts de~$\infty$. Or si~$R>0$ l'image réciproque de~$\wbeth E-\bbfer a R$ sur~$[\rho(x);+\infty]$ est précisément~$]R;\infty]$ qui est ouvert ; ceci achève de prouver notre assertion. 

\deux{definterxywb} La notation~$[x;y]$ a un sens dès que~$x$ et~$y$ sont deux points comparables de~$\wbeth E$. On l'étend à {\em tout} couple~$(x,y)$ de points de~$\wbeth E$ en posant~$[x;y]=[x;x\wedge y]\cup [y;x\wedge y]$ (notons que cette formule est trivialement vraie lorsque~$x$ et~$y$ sont comparables, puisque~$x\wedge y$ est alors égal à l'un des deux).  

\medskip
Pour tout couple~$(x,y)$ de points de~$\wbeth E$, le sous-ensemble~$[x;y]$ de~$\wbeth E$ est homéomorphe à un intervalle compact d'extrémités~$x$ et~$y$ (on l'a déjà vu si~$x$ et~$y$ sont comparables ; le cas général en résulte en vertu de~\ref{compquoti}). Nous étendons la définition de la convexité (qui ne s'appliquait jusqu'ici qu'aux chaînes) de façon naturelle : nous dirons qu'une partie~$X$ de~$\wbeth E$ est {\em convexe} si~$[x;y]\subset X$ dès que~$x$ et~$y$ appartiennent à~$X$. Toute partie convexe est connexe par arcs, ce qui s'applique notamment à~$\wbeth E$ lui-même ; l'intersection d'une famille de parties convexes est convexe. 

\medskip
Soit~$a\in E$ et soit~$r>0$. 

\trois{connectbfer} {\em Le fermé~$\bbfer a r$ est convexe.} Soient en effet~$x$ et~$y$ deux points de~$\bbfer a r$. Ils sont tous deux majorés par~$\zeta_{a,r}$ ; par conséquent,~$\zeta_{a,r}\geq x\wedge y$, et tout élément de~$[x;y]$ est donc majoré par~$\zeta_{a,r}$, ce qu'il fallait établir. 

\trois{connectbouv} {\em L'ouvert~$\bbouv a r$ est convexe.} En effet, si~$x$ et~$y$ appartiennent à~$\bbouv a r$, il existe~$s<r$ tel que tous deux soient situés sur~$\bbfer a s$. On a alors d'après le~\ref{connectbfer} ci-dessus~$[x;y]\subset \bbfer a s\subset \bbouv a r$. 

\trois{connectcompbfer} {\em L'ouvert~$\wbeth E-\bbfer ar$ est convexe.} En effet, soient~$x$ et~$y$ deux points de~$\wbeth E-\bbfer ar$, et soit~$z\in [x;y]$. Par définition de~$[x;y]$ le point~$z$ majore~$x$ ou~$y$, ce qui exclut que l'on ait~$z\leq \zeta_{a,r}$ ; ainsi,~$[x;y]\subset \wbeth E-\bbfer ar$. 

\deux{locconvexe} Il découle de~\ref{connectbouv} et~\ref{connectcompbfer} que tout point de~$\wbeth E$ possède une base de voisinages convexes ; en particulier,~$\wbeth E$ est localement connexe par arcs. 

\deux{compconncompzet} Soit~$a\in E$ et soit~$r>0$. Si~$b$ et~$c$ sont deux points de~$\bfer a r$ alors~$\bbouv b r=\bbouv c r$ si et seulement si~$d(b,c)<r$, c'est-à-dire si et seulement si ils ont même image dans l'ensemble~$\Theta(a,r)$ défini au~\ref{defthetaar} ; si~$\beta\in \Theta(a,r)$ on peut donc définir sans ambiguïté~$\bbouv \beta r$ comme étant égale à~$\bbouv b r$ pour {\em n'importe quel} antécédent~$b$ de~$\beta$ dans~$\bbfer a r$. 

\trois{descbbfermszeta} Soit~$x$ un élément de~$ \bbfer a r$ différent de~$\zeta_{a,r}$ ; choisissons~$y\in ]x;\zeta_{a,r}[$. Étant strictement supérieur à~$x$, le point~$y$ est égal à~$\zeta_{b,s}$ pour un certain~$b\in E$ et un certain~$s>0$. Comme~$y<\zeta_{a,r}$ on a~$b\in \bbfer a r$ et~$s<r$ ; par conséquent,~$x\in \bbfer b s\subset \bbouv b r$. 

\medskip
On peut ainsi écrire~$\bbfer a r\setminus \{\zeta_{a,r}\}=\coprod\limits_{\beta \in \Theta(a,r)} \bbouv \beta r$, et partant~$$\wbeth E\setminus\{\zeta_{a,r}\}=\left(\coprod\limits_{\beta \in \Theta(a,r)} \bbouv \beta r\right)\coprod \left(\wbeth E-\bbfer a r\right).$$

\trois{conclwbmsz} Les différents termes dont le membre de droite est la réunion disjointe sont des ouverts non vides et connexes de~$\wbeth E$ (\ref{connectbouv} et~\ref{connectcompbfer}); ce sont donc les composantes connexes de~$\wbeth E \setminus\{\zeta_{a,r}\}$. Comme~$\wbeth E$ est lui-même connexe, chacune de ces composantes a pour bord~$\{\zeta_{a,r}\}$. 

\deux{unisegm} Soient~$x$ et~$y$ deux points de~$\wbeth E$. Nous allons montrer que~$[x;y]$ est {\em l'unique} fermé de~$\wbeth E$ homéomorphe à un segment d'extrémités~$x$ et~$y$. Soit~$F$ un tel fermé ; il suffit de montrer que~$F$ {\em contient}~$[x;y]$. 

\trois{unisegmwedge} {\em Le fermé~$F$ contient~$x\wedge y$.} En effet, si~$x\wedge y=x$ ou~$x\wedge y=y$ c'est évident. Dans le cas contraire, ni~$x$ ni~$y$ ne sont égaux à~$\infty$, et il découle de~\ref{wedgeboulouv} qu'il existe~$b$ et~$c$ dans~$E$ et~$r>0$ tels que~$x\wedge y=\zeta_{b,r}=\zeta_{c,r}$, tels que~$x\in \bbouv b r~$ et~$c\in \bbouv c r$, et tels que~$\bbouv b r \cap \bbouv c r=\emptyset$ ; en vertu du~\ref{conclwbmsz} ci-dessus,~$\bbouv b r$ et~$\bbouv c r$ sont deux composantes connexes distinctes de~$\wbeth E\setminus \{x\wedge y\}$. Par conséquent,~$x\wedge y \in F$. 

\trois{unisegmconclu} Il suffit maintenant de vérifier que~$]x;x\wedge y[\subset F$ et~$]y;x\wedge y[\subset F$. Par symétrie, on peut se contenter de traiter le cas de~$]x;x\wedge y[$; soit donc~$z\in ]x;x\wedge y[$. Comme~$z>x$ il existe~$a\in E$ et~$r>0$ tel que~$z=\zeta_{a,r}$. Comme~$x<z$ il appartient à~$\bbfer a r\setminus\{\zeta_{a,r}\}$ ; comme~$x\wedge y>z$ il appartient à~$\wbeth E-\bbfer a r$. En en vertu du~\ref{conclwbmsz} ci-dessus,~$x$ et~$y$ sont situés sur deux composantes connexes différentes de~$\wbeth E\setminus\{\zeta_{a,r}\}$ ; par conséquent,~$z=\zeta_{a,r}$ appartient à~$F$, ce qui achève la démonstration. 

\deux{wbethearbr} On déduit du~\ref{unisegm} ci-dessus, de~\ref{connectbouv},~\ref{connectcompbfer}, et de~\ref{conclwbmsz}, que tout point de~$\wbeth E$ possède une base de voisinages ouverts~$U$ possédant les deux propriétés suivantes : 

\medskip
$\bullet$~$\partial U$ est fini ; 

$\bullet$ pour tout couple~$(x,y)$ d'éléments de~$U$, il existe un et un seul fermé de~$U$ homéomorphe à un segment d'extrémités~$x$ et~$y$. 

\deux{pizerowbe} Soit~$x$ un point de~$\wbeth E$ ; nous nous proposons d'étudier le cardinal de~$\pi_0(\wbeth E\setminus\{x\})$. 

\trois{pizerowbezar} Si~$x$ est de la forme~$\zeta_{a,r}$ avec~$r>0$, autrement dit si~$x\in \aleph(E)$ et si~$\rho(x)>0$ il résulte de~\ref{conclwbmsz} que~$$\pi_0(\wbeth E\setminus\{x\})=\sharp \Theta(a,r)+1.$$ Il s'ensuit que le cardinal de~$\pi_0(\wbeth E\setminus\{x\})$ vaut au moins deux, et qu'il est infini si et seulement si~$\Theta(a,r)$ est infini. 

Notons que~$\sharp\pi_0(\wbeth E\setminus\{x\})= 2$ si et seulement si~$\sharp \Theta(a,r)=1$, et que cela équivaut à l'égalité~$\bfer a r=\bouv a r$. 

\trois{pizerowbeextr} Supposons que~$x\in \beth(E)-\aleph(E)$ ou que~$x\in \aleph(E)$ et que~$\rho(x)=0$. et soient~$y$ et~$z$ deux éléments de~$\wbeth E\setminus\{x\}$. Si~$t\in [y;z]$ alors par définition de ce dernier on a~$t\geq y$ ou~$t\geq z$ ; si~$t=y$ ou si~$t=z$ il est différent de~$x$ par hypothèse, et si~$t>y$ ou~$t>z$ il résulte de la minimalité de~$x$ (\ref{minim}) que~$t\neq x$ ; par conséquent,~$[y;z]\subset \wbeth E\setminus\{x\}$ et ce dernier est connexe. 

\trois{pizerowbeinf} Supposons que~$x=\infty$. Si~$y$ et~$z$ sont deux points appartenant à~$\wbeth E\setminus\{x\}=\beth(E)$ il résulte de la définition de~$[y;z]$ que~$[y;z]\in \beth(E)$ ; par conséquent~$\wbeth E \setminus\{x\}$ est connexe.

\subsection*{Complément: boules emboîtées de~$E$ et espace~$\beth (E)$}

\deux{nastequiv} Soit~$\sch B$ une chaîne de boules fermées de~$E$. On dira que~$\sch B$ est  {\em évanescente} si l'intersection des boules appartenant à~$\sch B$ est vide ; comme~$E$ est non vide, cela force~$\sch B$ elle-même à être non vide ; cela la force également à ne contenir aucun singleton. 

\trois{defequiv} Si~$\sch B$ et~$\sch B'$ sont deux chaînes évanescentes de boules fermées de~$E$, on dira qu'elles sont {\em équivalentes} si~$B\cap B'\neq \emptyset$ pour tout~$(B,B')\in \sch B\times \sch B'$. Cette terminologie est raisonnable dans la mesure où nous définissons bien ainsi une relation d'équivalence. Pour le voir, il suffit d'en vérifier la transitivité, la réflexivité et la symétrie étant évidentes. 

\medskip
Soient donc~$\sch B, \sch B'$ et~$\sch B''$ trois chaînes évanescentes de boules fermées de~$E$ telles que ~$\sch B$ et~$\sch B'$ soient équivalentes, et telles que~$\sch B'$ et~$\sch B''$ soient équivalentes. Soit~$B$ une boule appartenant à~$\sch B$ et soit~$B''$ une boule appartenant à~$\sch B''$. Comme~$\sch B'$ est évanescente, il existe~$B'_0\in\sch B'$ ne contenant pas~$B$, et il existe~$B'_1\in \sch B'$ ne contenant pas~$B''$. Soit~$B'$ l'intersection de~$B'_0$ et~$B'_1$ ; c'est une boule appartenant à~$\sch B'$. Par construction,~$B'$ ne contient pas~$B$ ; comme~$\sch B$ et~$\sch B'$ sont équivalentes, on a alors nécessairement ~$B'\subset B$. Par le même raisonnement~$B'\subset B''$. Ainsi~$B\cap B''$ contient~$B'$ et est en particulier non vide, ce qu'il fallait établir. 

\trois{ultnast} Soit~$\sch B$ une chaîne évanescente de boules fermées de~$E$ et soit~$\overline {\sch B}$ l'ensemble des boules fermées de~$E$ qui contiennent une boule appartenant à~$\sch B$. Il est immédiat que~$ \overline {\sch B}$ est une chaîne évanescente de boules fermées qui est équivalente à~$\sch B$ et contient cette dernière. 

\medskip
{\em La chaîne~$\overline {\sch B}$ est la plus grande chaîne évanescente de boules fermées de~$E$ équivalente à~$\sch B$.} En effet, soit~$\sch B'$ une chaîne évanescente de boules fermées de~$E$ équivalent à~$\sch B$ et soit~$B'\in \sch B'$. Comme~$\sch B$ est évanescente, il existe~$B\in \sch B$ tel que~$B$ ne contienne pas~$B'$ ; comme les chaînes~$\sch B$ et~$\sch B'$ sont équivalentes, on a~$B\cap B'\neq \emptyset$, et partant~$B\subset B'$ ; ainsi,~$B'\in \overline{ \sch B}$, et l'on a bien~$\sch B'\subset \overline {\sch B}$, ce qu'il fallait démontrer.  

\medskip
On en déduit que l'ensemble des classes d'équivalence de chaînes évanescentes de boules fermées de~$E$ est en bijection naturelle avec l'ensemble des chaînes évanescentes maximales de boules fermées de~$E$.

\deux{chainebethchaine} Soit~$x\in \beth(E)$. On note~$E^{\leq x}$ l'ensemble des éléments de~$E$ majorés par~$x$. Si~$x\in \aleph(E)$ il est égal à~$\zeta_{a,r}$ pour un certain~$a\in E$ et un certain~$r\geq 0$, et~$E^{\leq x}$ est alors la boule fermée~$\bfer a r$ de~$E$ ; si~$x\notin \aleph(E)$ c'est un élément minimal de~$\beth(E)$ en vertu de~\ref{minim}, et~$E^{\leq x}$ est donc vide. Le point~$x$ étant par ailleurs la borne inférieure de~$]x;\infty[$, on a~$$E^{\leq x}=\bigcap\limits_{y\in ]x;\infty[} E^{\leq y}.$$

\deux{chainebethmax} Si~$F$ est une chaîne ouverte saturée de~$\aleph(E)$ dont la borne inférieure n'appartient pas à~$\aleph(E)$, il découle de ce qui précède que~$\got b(F):=\{ E^{\leq y}\}_{y\in F}$ est une chaîne évanescente de boules fermées de~$E$. Cette chaîne est maximale : en effet, soit~$B\in \got b(F)$ et soit~$B'$ une boule fermée de~$E$ contenant~$B$. La boule~$B$ est de la forme~$E^{\leq y}$ pour un certain~$y\in F$ ; écrivons~$y=\zeta_{a,r}$ avec~$a\in E$ et~$r\geq 0$ ; on a alors~$B=\bfer a r$. 

Comme~$B'\supset B$, on a~$a\in B'$ et~$B'$ est de ce fait égal à~$\bfer a s$ pour un certain~$s\geq 0$. Si~$s\leq r$ alors~$B'\subset B$, d'où l'égalité~$B'=B$ ; si~$s>r$ alors~$\zeta_{a,s}\in ]y;\infty[\subset F$ et~$B'=\bfer a s =E^{\leq \zeta_{a,s}}\in \got b(F)$. On a donc~$B'\in \got b (F)$ dans tous les cas ; par conséquent~$\got b(F)=\overline{\got b(F)}$, ce qui équivaut d'après~\ref{ultnast} à la maximalité de~$\got b(F)$.

\deux{chainemax} Soit~$\sch B$ une chaîne évanescente maximale de boules fermées de~$E$, et soit~$\got f(\sch B)$ l'ensemble des éléments~$y$ de~$\aleph(E)$ tels que~$E^{\leq y}\in \sch B$. 

\medskip
{\em L'ensemble~$\got f(\sch  B)$ est une chaîne ouverte saturée de~$\aleph(E)$ dont la borne inférieure n'appartient pas à~$\aleph(E)$, et~$\got b(\got f(\sch B))=\sch B$.} En effet, la chaîne~$\sch B$ est non vide et contient donc~$\bfer a r$ pour un certain~$a\in E$ et un certain~$r\geq 0$ ; l'élément~$\zeta_{a,r}$ de~$\aleph(E)$ appartient alors à~$\got f(\sch B)$, et il résulte immédiatement de la définition de~$\got f(\sch B)$ et de la maximalité de~$\sch B$ que~$\got f(\sch B)$ contient~$[\zeta_{a,r};\infty[$. 

\medskip
Soient~$y$ et~$z$ deux éléments de~$\got f(\sch B)$. Les boules~$E^{\leq y}$ et~$E^{\leq z}$ appartenant à~$\sch B$, leur intersection est non vide ; il en va {\em a fortiori} de même de~$\beth(E)^{\leq y}\cap \beth(E)^{\leq z}$, ce qui entraîne que ~$\beth(E)^{\leq y}$ et ~$\beth(E)^{\leq z}$ sont comparables pour l'inclusion. Par conséquent,~$y$ et~$z$ sont comparables, et~$\got f(\sch B)$ est une chaîne ; comme elle contient~$[\zeta_{a,r};\infty[$, elle est saturée. 

\medskip
Par définition, l'ensemble~$\{E^{\leq y}\}_{y\in \got f(\sch B)}$ est égal à la chaîne~$\sch B$ ; celle-ci étant évanescente,~$\bigcap\limits_{y\in\got f(\sch B)} E^{\leq y}=\emptyset$. Cela exclut que~$\got f(B)$ puisse avoir un plus petit élément~$x$ dans~$\aleph(E)$ (car sinon,~$E^{\leq x}$ serait contenue dans l'intersection précédente), et elle est en conséquence ouverte ; il découle dès lors de~\ref{chainebethchaine} que sa borne inférieure n'appartient pas à~$\aleph(E)$. 

\medskip
L'égalité~$\{E^{\leq y}\}_{y\in \got f(\sch B)}=\sch B$ se récrit sous la forme~$\got b(\got f(\sch B))=\sch B$, ce qui achève de prouver l'assertion requise. 

\deux{nastbij} Soit~$F$ une chaîne ouverte saturée de~$\aleph(E)$ dont la borne inférieure n'appartient pas à~$\aleph(E)$ ; nous allons montrer que~$\got f(\got b(F))=F$. Notons~$x$ la borne inférieure de~$F$, et~$x'$ celle de~$\got f(\got b(F))$, qui est une chaîne ouverte saturée d'après~\ref{chainemax}. 

Si~$y\in F$ alors~$E^{\leq y}\in \got b(F)$ ; il résulte donc des définitions que~$F\subset \got f(\got b(F))$, ce qui signifie que~$x'\leq x$. Comme~$x\notin \aleph(E)$, il est minimal en vertu de~\ref{minim}. Par conséquent,~$x'=x$, ce qui entraîne que~$F=\got f(\got b(F))$. 

\deux{recapbijnast} {\bf Conclusion.} Il découle de~\ref{chainebethmax},~\ref{chainemax} et~\ref{nastbij} que~$F\mapsto \got b(F)$ et~$\sch B\mapsto \got f(\sch B)$ établissent une bijection entre l'ensemble des chaînes ouvertes saturées de~$\aleph(E)$ dont la borne inférieure n'appartient pas à~$\aleph(E)$ et l'ensemble des chaînes évanescentes maximales de boules fermées de~$E$. 

\medskip
On peut reformuler cette assertion de la façon suivante : les applications~$x\mapsto \got b(]x;\infty[)$ et~$\sch B\mapsto \inf \got f(\sch B)$ établissent une bijection entre~$\beth(E)-\aleph(E)$  et l'ensemble des chaînes évanescentes maximales de boules fermées de~$E$ ; rappelons que ce dernier peut lui-même être canoniquement identifié à l'ensemble des classes d'équivalence de chaînes évanescentes de boules fermées de~$E$ (\ref{ultnast}).

\section{Arbres et graphes : propriétés de base}

\subsection*{Définitions et premiers exemples}

\deux{defarb} Nous dirons qu'un espace topologique~$X$ est un {\em graphe} s'il possède les propriétés suivantes :

\begin{itemize}
\medskip
\item[i)]~$X$ est séparé et localement compact ; 

\item[ii)]~$X$ possède une base d'ouverts~$U$ satisfaisant les conditions suivantes :

\medskip

\begin{itemize}
\item[$\alpha)$] pour tout couple~$(x,y)$ de points de~$U$ il existe un et un seul fermé~$[x;y]$ de~$U$ homéomorphe à un segment d'extrémités~$x$ et~$y$ ; 

\item[$\beta)$] le bord de~$U$ dans~$X$ est fini. 

\end{itemize}

\end{itemize}

\medskip
Si de plus~$X$ lui-même vérifie la propriété~$\alpha)$, nous dirons que c'est un {\em arbre.} Si~$X$ est un graphe, nous qualifierons de {\em sous-graphe} (resp. {\em sous-arbre}) de~$X$ tout sous-ensemble localement fermé de~$X$ qui est un graphe (resp. un arbre) pour la topologie induite. Il résulte des définitions qu'un graphe (resp. un arbre) est localement connexe par arcs (resp. connexe par arcs). Tout ouvert d'un graphe en est un sous-graphe ; un espace topologique est un graphe si et seulement si il est séparé et si chacun de ses points possède un voisinage ouvert qui est un arbre. 

\deux{premexarbr} {\bf Exemples.}

\trois{defgraphelocfin} Soit~$E$ un sous-ensemble fini de ~$\RR/2\pi \ZZ$ et soit~$\got G_E$ le sous-ensemble~$$\{0\}\cup\{z\in \CC^\times,|z|<1\;{\rm et}\;{\rm Arg}\;z\in E\}$$ de~$\CC$ ; l'espace topologique~$\got G_E$ est un arbre. Si~$X$ est un graphe et si~$x$
appartient à~$X$, on dira que~$X$ est {\em fini en~$x$} si~$x$ possède dans~$X$ un voisinage admettant,
pour un certain~$E$, 
un homéomorphisme sur~$\got G_E$
qui envoie~$x$ sur~$0$. Le cardinal de l'ensemble~$E$ en question est alors bien déterminé ; on l'appelle la {\em valence} de~$(X,x)$ ; on parlera parfois plus simplement de la valence de~$x$, s'il n'y a pas d'ambiguïté sur~$X$. 

\medskip
On dira qu'un graphe est {\em localement fini} s'il est fini en chacun de ses points. Si~$X$ est un graphe localement fini, on appellera {\em sommet topologique} de~$X$ (on omettra souvent l'épithète «topologique» s'il n'y a pas d'ambiguïté), tout point de~$X$ qui ne possède pas de voisinage homéomorphe à un intervalle ouvert, c'est-à-dire encore tout point~$x$ de~$X$ tel que la valence de~$(X,x)$ soit différente de~$2$ ; on dira d'un graphe qu'il est {\em fini} s'il est localement fini et n'a qu'un nombre fini de sommets et de composantes connexes ; tout graphe compact et localement fini est fini. 

\medskip
 Si~$X$ est un graphe fini et compact, tout sous-graphe fermé de~$X$ est fini et compact ; on en déduit que tout sous-graphe d'un graphe localement fini est localement fini. 

\trois{locfinfac} Soit~$X$ un graphe dont tout point~$x$ a un voisinage~$U$ qui est un arbre vérifiant l'égalité~$U=\bigcup [y_i;x]$,  où~$\{y_i\}$ est un ensemble {\em fini} de points de~$U$ et où~$[y_i;x]\cap[y_j;x]=\{x\}$ pour tout~$i\neq j$. Il résulte du~\ref{compquoti} que le graphe~$X$ est localement fini. 

\trois{wbetheexarbre} Soit~$E$ un espace métrique ultramétrique. L'espace topologique~$\wbeth E$ défini au~\ref{defwbeth} est un arbre compact en vertu de~\ref{wbethearbr}. 

\subsection*{Quelques propriétés locales des graphes}

\deux{voisptgraphecon} Soit~$X$ un graphe {\em connexe}. 

\trois{prsquette} Soit~$x\in X$ et soit~$U$ un voisinage ouvert de~$x$ dans~$X$ ; pour presque toute composante connexe~$V$ de~$X\setminus\{x\}$, le sous-ensemble~$V\cup\{x\}$ de~$X$ est compact,~$V$ est un arbre et~$V\subset U$ (ce qui entraîne que~$V$ est une composante connexe de~$U\setminus\{x\}$). 

\medskip
En effet, il existe un voisinage ouvert~$W$ de~$x$ dans~$X$ qui est contenu dans~$U$, relativement compact et satisfait~$\alpha)$ et~$\beta)$ ; comme~$\partial W$ est fini, il ne rencontre qu'un nombre fini de composantes connexes de~$X\setminus\{x\}$. Soit~$V$ une composante connexe de~$X\setminus\{x\}$ ne rencontrant pas~$\partial W$ ; on a~$\partial V=\{x\}$, ce qui entraîne que~$V$ rencontre~$W$ ; puisque~$V$ ne rencontre pas~$\partial W$, l'intersection~$V\cap W$, qui est ouverte et non vide, est fermée dans~$V$ ; par connexité de~$V$, elle est égale à~$V$, ce qui signifie que~$V$ est contenue dans~$W$, et {\em a fortiori} dans~$U$. Comme~$W$ est relativement compact,~$\overline V=V\cup\{x\}$ est compact ; comme~$W$ est un arbre,~$V$ est un arbre.

\trois{intercompvois} Soit~$x\in X$, soit~$V$ une composante connexe de~$X\setminus\{x\}$ et soit~$U$ un voisinage ouvert connexe de~$x$. Toute composante connexe de~$V\cap U$ est un ouvert connexe et non vide de~$U$ de bord contenu dans~$\{x\}$ ; c'est donc une composante connexe de~$U\setminus\{x\}$. 

\medskip
Le cardinal de~$\pi_0(U\cap V)$ est par ailleurs non nul et fini. En effet, il est non nul puisque~$x$ adhère à~$V$ ; et il est fini car sinon il existerait une infinité de composantes connexes de~$U\setminus\{x\}$ qui ne soient pas des composantes connexes de~$X\setminus\{x\}$, en contradiction avec~\ref{prsquette}. 

\medskip
Soit~$W$ une composante connexe de~$U\setminus\{x\}$ qui est relativement compacte dans~$U$ et est contenue dans~$V$. Comme~$W\cup\{x\}$ est compacte, c'est une partie fermée de~$V\cup \{x\}$ ; il s'ensuit que l'ouvert non vide~$W$ de~$V$ en est également un fermé. On a donc~$W=V$, et~$V$ est {\em a fortiori} relativement compacte dans~$X$. 

\subsection*{Parties convexes d'un arbre}

\deux{intervxycomp} Soit~$X$ un arbre. On dira qu'une partie~$E$ de~$X$ est {\em convexe} si~$[x;y]\subset E$ pour tout couple~$(x,y)$ de points de~$E$.

\trois{recolsegm} Si~$x,y$ et~$z$ sont trois points de~$X$ et si~$[x;z]\cap [z;y]=\{z\}$ alors~$[x;z]\cup[y;z]$ est un segment d'extrémités~$x$ et~$y$ (\ref{compquoti}), et est donc égal à~$[x;y]$. En général,~$[x;y]\cap [y;z]$ est convexe et est donc de la forme~$[z;t]$ pour un certain~$t\in X$ ; on a~$[x;t]\cap [t;y]=\{t\}$ et donc~$[x;y]=[x;t]\cup[t;y]\cup[y;z]$.

Il s'ensuit aisément (sans utiliser~\ref{connarcinj}) que si~$U$ est un ouvert de~$X$ et si~$x\in U$ alors l'ensemble des points~$y$ de~$U$ tels que~$[x;y]\subset U$ est ouvert et fermé dans~$U$ ; par conséquent,~$U$ est connexe si et seulement si il est convexe, c'est-à-dire si et seulement si c'est un sous-arbre de de~$X$.

\trois{interconn} L'intersection de deux parties convexes de~$X$ est convexe ; en particulier, l'intersection de deux sous-arbres ouverts de~$X$ est un sous-arbre ouvert de~$X$. La réunion de deux parties convexes de~$X$ qui se rencontrent est une partie convexe de~$X$ (on le déduit du~\ref{recolsegm} ci-dessus). 

\trois{arbreferme} Soit~$E$ une partie localement fermée et convexe de~$X$ ; c'est un espace séparé et localement compact, et pour tout couple~$(x,y)$ de points de~$E$ il existe un unique fermé de~$E$ homéomorphe à un segment d'extrémités~$x$ et~$y$, à savoir~$[x;y]$. Par ailleurs, si~$x\in E$ il possède une base de voisinages dans~$E$ satisfaisant les conditions~$\alpha)$ et~$\beta)$ du~\ref{defarb}. En effet, il possède une telle base dans~$X$, et si~$U$ est un voisinage de~$x$ dans~$X$ satisfaisant ~$\alpha)$ et~$\beta)$ alors~$\partial_E(U\cap E)\subset \partial_XU$, et est donc fini ; de plus,~$U\cap E$ est une partie convexe de l'arbre~$U$, et satisfait dès lors~$\beta)$ en vertu de ce qui précède. Par conséquent,~$E$ est un arbre. 

\medskip
Il s'ensuit que si~$F$ et~$G$ sont deux sous-arbres fermés de~$X$ d'intersection non vide alors alors~$F\cup G$ est un sous-arbre fermé de~$X$ ; si de plus 
$F$ et~$G$ sont localement finis
 alors~$F\cup G$ est localement fini : cela résulte de~\ref{recolsegm} et de~\ref{locfinfac}. 

\trois{arbrebijcom} Si~$x\in X$, si~$U$ est un voisinage ouvert et connexe de~$x$, et si~$V$ est une composante connexe de~$X\setminus\{x\}$ alors~$V\cap U$ est connexe, et est donc une composante connexe de~$U\setminus\{x\}$ (\ref{intercompvois}) ; il s'ensuit que l'application naturelle~$\pi_0(U\setminus\{x\})\to \pi_0(X\setminus\{x\})$ est bijective. 

\trois{descxycomp} Si~$x$ et~$y$ sont deux points distincts de~$X$ alors~$]x;y[$ est exactement l'ensemble des~$z\in X\setminus\{x,y\}$ tels que~$x$ et~$y$ soient situés sur deux composantes connexes distinctes de~$X\setminus\{x\}$. En effet, soit~$z\in X$. Si~$x$ et~$y$ sont situés sur la même composante connexe de~$X\setminus\{z\}$, celle-ci contient nécessairement~$[x;y]$, et par conséquent~$z$ n'appartient pas à~$]x;y[$. Supposons maintenant que~$x$ et~$y$ soient situés sur deux composantes connexes différentes, respectivement notées~$U$ et~$V$, de~$X\setminus\{z\}$. L'intervalle ouvert~$[x;z[$ (resp.~$[y;z[$) est alors contenu dans~$U$ (resp.~$V$) ; il s'ensuit que~$[x;z]\cap [z;y]=\{z\}$ ; par conséquent~$[x;z]\cup[z;y]=[x;y]$ et l'on a bien~$z\in ]x;y[$. 

\trois{intervcoinc} Soient~$x, y$ et~$z$ appartenant à~$X$. Si~$x$ et~$y$ sont situés sur deux composantes connexes différentes de de~$X\setminus\{z\}$ alors~$z\in ]y;x[$ d'après le~\ref{descxycomp} ci-dessus, et l'on a donc~$[y;z]\cap [x;z]=\{z\}$. Si~$x$ et~$y$ sont situés sur la même composante connexe~$U$ de~$X\setminus\{z\}$ alors~$[y;z]\cap [x;z]$ est de la forme~$[t;z]$ pour un certain~$t$ qui appartient à~$U$ (et est donc différent de~$z$). En effet~$[y;z]\cap [x;z]$ est de la forme~$[t;z]$ avec~$t\in [x;z]$ (\ref{recolsegm}) ; comme~$[x;z[\subset U$, il reste à s'assurer que~$t\neq z$. Si l'on avait~$t=z$, l'on aurait~$[y;z]\cap [x;z]=\{z\}$, et donc~$[y;z]\cup[x;z]=[y;x]$, d'où~$z\in ]y;x[$ ; mais cela contredirait l'hypothèse que~$y$ et~$x$ sont situés sur la même composante connexe de~$X\setminus\{z\}$.

\trois{adhcompinterv} Soit~$E$ une partie non vide et convexe de~$X$ ; un point~$x$ de~$X$ appartient à~$\overline E$ si et seulement~$[y;x[\subset E$ pour tout~$y\in E$. La condition est en effet clairement suffisante. Pour voir qu'elle est nécessaire, on suppose que~$x$ appartient à ~$\overline E$. Soit~$y\in E$. L'intersection de~$E$ et~$[y;x]$ est convexe et contient~$y$, c'est donc un intervalle non vide de~$[y;x]$ ; soit~$t$ sa borne supérieure ; nous allons montrer par l'absurde que~$t=x$. On suppose que ce n'est pas le cas ; la composante connexe de~$X\setminus\{t\}$ qui contient~$]t;x]$ est alors un voisinage ouvert de~$x$, qui par conséquent rencontre~$E$ ; soit~$z$ un point de~$E$ situé sur cette composante. L'intersection de~$[t;z]$ et~$[t;x]$ est de la forme~$[t;t']$ avec~$t'\neq t$ ; on a~$[t;t']\subset [t;x]$ et~$[y;t']\subset E$ puisque ce dernier est convexe; dès lors~$$[y;t']\subset [y;x]\cap E,$$ ce qui contredit la définition de~$t$. 

Fixons~$x$ dans~$E$. Si~$y$ et~$z$ sont deux points de~$\overline E$ alors~$$[y;z]\subset [x;y]\cup[x;z]\subset \overline E$$ en vertu de ce qui précède. Par conséquent,~$\overline E$ est convexe et est donc un sous-arbre fermé de~$X$ (\ref{arbreferme}).  Remarquons que cette dernière assertion s'étend trivialement au cas où~$E=\emptyset$. 

\trois{unecompdeuxpts} Soit~$S$ un sous-ensemble fermé de~$X$ et soit~$U$ une composante connexe de~$X\setminus S$ ; supposons que~$\partial U$ contienne deux éléments {\em distincts}~$x$ et~$y$ de~$S$. Choisissons~$z\in U$ ; en vertu du~\ref{adhcompinterv}, les intervalles~$[z;x[$ et~$[z;y[$ sont contenus dans~$U$. L'intersection~$[z;x]\cap [z;y]$ est convexe, compacte, et contient~$z$ ; elle est donc de la forme~$[z;t]$ pour un certain~$t\in U$. On a alors~$[x;t]\cap[t;y]=\{t\}$ et donc~$[x;y]=[x;t]\cup[t;y]$. Par conséquent,~$]x;y[\subset U$, et ne rencontre {\em a fortiori} pas~$S$ ; l'ouvert~$U$ est alors nécessairement {\em la} composante connexe de~$X\setminus S$ contenant~$]x;y[$. 

\trois{arbreestadm} Appliquons ce qui précède au cas où~$S$ est un sous-arbre fermé et non vide de~$X$. Soit~$U$ une composante connexe de~$X\setminus S$ ; l'ensemble~$\partial U$ est une partie non vide de~$S$, et nous allons montrer par l'absurde que c'est un singleton. Supposons que ce ne soit pas le cas. Il existerait alors deux points distincts~$x$ et~$y$ dans~$\partial_X U$. On aurait alors~$]x;y[\subset U\subset X\setminus S$ d'après le~\ref{unecompdeuxpts}, et~$]x;y[\subset S$ puisque~$S$ est un sous-arbre fermé ; on aboutit ainsi à une contradiction. 

\medskip
Soit~$x$ l'unique point de~$\partial U$, soit~$y\in U$ et soit~$z\in S$. L'intersection~$[y;z]\cap S$ est alors égale à~$[y;x]$ ; en effet, c'est une partie fermée et convexe de~$[y;z]$ qui contient~$z$ ; elle est donc de la forme~$[t;z]$ pour un certain~$t\in S$. L'intervalle ouvert~$[y;t[$ ne rencontrant pas~$S$, il est contenu dans une composante connexe de~$X\setminus S$ qui est nécessairement celle de~$y$, à savoir~$U$ ; dès lors~$t\in \partial U$ et~$t=x$.

\subsection*{Autour des sous-graphes d'un graphe donné}

\deux{cercle} Soit~$X$ un graphe connexe et soit~$x\in X$. On déduit de~\ref{descxycomp} et~\ref{intervcoinc} que l'ensemble des points pouvant être joints à~$x$ par un arc injectif est à la fois ouvert et fermé dans~$X$ ; on en déduit, sans avoir à utiliser~\ref{connarcinj}, que~$X$ est connexe par arcs injectifs.

\deux{sousgrapheinter} Soit~$X$ un graphe et soient~$E$ et~$F$ deux sous-graphes de~$X$. Soit~$V$ un sous-arbre ouvert de~$X$. Les intersections~$E\cap V$ et~$F\cap V$ sont des sous-graphes de~$V$ ; les composantes connexes de chacune d'elles sont donc des sous-arbres de~$V$, fermés (resp. ouverts, resp. localement finis) si~$E$ et~$F$ sont fermés (resp. ouverts resp. localement finis). On en déduit, à l'aide de~\ref{arbreferme} : que~$E\cap F\cap V$ est un sous-graphe de~$V$, localement fini si~$E$ et~$F$ sont localement finis ; et que~$(E\cup F)\cap V$ est un sous-graphe de~$V$ si~$E$ et~$F$ sont tous deux fermés ou bien tous deux ouverts, qui est là encore localement fini si~$E$ et~$F$ sont localement finis. 

\medskip
Par conséquent,~$E\cap F$ est un sous-graphe de~$X$, localement fini si~$E$ et~$F$ sont localement finis ; et si~$E$ et~$F$ sont tous deux fermés ou bien tous deux ouverts,~$E\cup F$ est un sous-graphe de~$X$, localement fini si~$E$ et~$F$ sont localement finis. 

\deux{bouclepasarbre} Soit~$X$ un graphe. Si~$x$ et~$y$ sont deux points de~$X$, et si~$I$ et~$J$ sont deux arcs injectifs joignant~$x$ à~$y$, la réunion~$I\cup J$ est un sous-graphe compact et fini de~$X$ (\ref{sousgrapheinter}). Il s'ensuit que s'il existe deux segments distincts~$I$ et~$J$ joignant~$y$ à~$x$, leur réunion contient une boucle ({\em i.e.} un compact homéomorphe à~$S^1$). 

\medskip
Par conséquent, un graphe connexe est un arbre si et seulement si il ne contient aucune boucle.

\deux{sousgraphead} Soit~$X$ un graphe et soit~$E$ un sous-graphe connexe de~$X$. Soit~$x\in \overline E$. Il existe un voisinage~$V$ de~$x$ dans~$X$ qui est un arbre compact dont le bord est fini. L'intersection~$E\cap V$ est, en vertu du~\ref{sousgrapheinter} ci-dessus, un sous-graphe de~$X$ ; elle est en particulier localement connexe. Il découle alors de~\ref{interconfin} que~$E\cap V$ a un nombre fini de composantes connexes ; appelons-les~$E_1,\ldots,E_r$. Si~$V'$ désigne l'intérieur de~$V$ dans~$X$, chacune des intersections~$E_i\cap V'$ est un sous-arbre de~$V'$. L'intersection~$\overline E^X\cap V'$ est égale à l'adhérence dans~$V'$ de~$\bigcup E_i\cap V'$, c'est-à-dire à la réunion des~$\overline{E_i\cap V'}^{V'}$ ; par conséquent,~$\overline E^X\cap V'$ est réunion disjointe de sous-arbres fermés de~$V'$ ; en particulier,~$x$ possède un voisinage dans~$\overline E^X$ qui est un arbre, et~$\overline E^X$ est donc un sous-graphe fermé de~$X$. 

\medskip
De plus, le point~$x$ appartient à~$\overline{E_i\cap V'}^{V'}$ pour au moins un~$i$ ; cela force~$E_i\cap V'$ à être non vide. Choisissons~$y\in E_i\cap V'$ ; il découle du~\ref{adhcompinterv} que~$[y;x[\subset E_i\cap V'\subset E$. 

On a ainsi démontré l'existence, pour tout~$x\in \overline E^X$, d'un point~$y$ de~$E$ et d'un segment~$I$ d'extrémités~$x$ et~$y$ tel que~$I\setminus\{x\}\subset E$. Il est clair, réciproquement, que si~$x$ est un point de~$X$ pour lequel existent un tel~$y$ et un tel~$I$, alors~$x\in \overline E^X$. 

\deux{intervalsousgr} Soit~$X$ un graphe, soit~$E$ un sous-graphe connexe de~$X$ et soient~$x$ et~$y$ deux points distincts de~$\overline E\setminus E$. En vertu du~\ref{sousgraphead} ci-dessus, il existe un segment~$I_x$ aboutissant à~$x$, non réduit à~$\{x\}$ et tel que~$I_x\setminus\{x\}\subset E$ ; on note~$x'$ sa deuxième extrémité, et l'on définit de même~$I_y$ et~$y'$. Par connexité de~$E$, il existe un segment~$J$ d'extrémités~$x'$ et~$y'$ tracé sur~$E$  ; la réunion~$\Gamma$ de~$I_x, I_y$ et~$\Gamma$ est un sous-arbre compact et connexe de~$X$ tel que~$\Gamma\setminus\{x,y\}\subset E$. Comme~$\Gamma$ est connexe, il existe un segment~$I$ tracé sur~$\Gamma$ et reliant~$x$ à~$y$ ; on a alors~$I\setminus\{x,y\}\subset E$. 

\subsection*{Composantes connexes de certains ouverts d'un graphe}

\deux{compxmscomp} Soit~$X$ un graphe connexe et soit~$\sch K$ une partie compacte de~$X$. Soit~$\Pi$ un sous-ensemble de~$\pi_0(X\setminus \sch K)$ tel que toutes les composantes connexes appartenant à~$\Pi$ soient relativement compactes. La réunion~$Y:=\sch K\cup \bigcup\limits_{V\in \Pi}V$ est alors une partie compacte de~$X$. 

En effet, donnons-nous une famille ~$(U_i)_{i\in I}$ d'ouverts de~$X$ qui recouvrent~$Y$. Nous allons prouver que l'on peut en extraire un sous-recouvrement fini ; on se ramène, quitte à raffiner~$(U_i)$, au cas où les~$U_i$ sont à bord fini. Comme~$\sch K$ est compact, il existe un sous-ensemble fini~$I'$ de~$I$ tel que~$\sch K$ soit contenu dans~$U:=\bigcup\limits_{i\in I'}U_i$ ; notons que~$\partial U$ est fini. 

Soit~$V\in \Pi$. Comme~$\partial V$ est un sous-ensemble non vide de~$\sch K$, l'intersection~$U\cap V$ est un ouvert non vide de~$V$ ; si de plus~$V$ ne rencontre pas~$\partial U$, cet ouvert est également fermé dans~$V$, et donc égal à~$V$, ce qui signifie que~$V\subset U$. 

Soient~$V_1,\ldots,V_r$ les composantes connexes appartenant à~$\Pi$ qui rencontrent~$\partial U$. Comme chacune d'elles est relativement compacte, la réunion~$\bigcup \overline {V_j}$ est une partie compacte de~$X$ contenue dans~$Y$. Il existe donc un sous-ensemble fini~$I''$ de~$I$ tel que~$\bigcup \overline {V_j}$ soit contenue dans~$\bigcup\limits_{i\in I''}U_i$.

Par construction, ~$Y\subset \bigcup\limits_{i\in I'\cup I''}U_i$, ce qui achève la démonstration.

\deux{graphemoinss} Soit~$X$ un graphe connexe et soit~$S$ un sous-ensemble fini de~$X$.

\trois{compgraphemoinss} Presque toute composante connexe~$W$ de~$X\setminus S$ est un arbre relativement compact et à bord singleton. 

Pour le voir, on procède par récurrence sur le cardinal de~$S$ ; si~$S=\emptyset$ il n'y a rien à démontrer. On suppose~$S\neq \emptyset$ et le résultat établi en cardinal strictement inférieur à celui de~$S$ ; on choisit un point~$s$ dans l'ensemble non vide~$S$, et l'on note~$U$ l'ouvert de~$X$ complémentaire de~$S\setminus\{s\}$ ; c'est un voisinage ouvert de~$s$. Soit~$\Pi$ le sous-ensemble de~$\pi_0(X\setminus\{s\})$ constitué des composantes~$W$ telles que~$W\subset U$, telles que~$W$ soit un arbre, et telles que~$W\cup\{s\}$ soit compact ; en vertu du~\ref{prsquette}, le complémentaire de~$\Pi$ dans~$\Pi_0(X\setminus\{s\})$ est fini ; soient~$W_1,\ldots,W_r$ ses éléments. 

Fixons~$i$. L'ouvert~$W_i$ de~$X$ est un graphe connexe, et le cardinal de~$W_i\cap S$ est strictement inférieur à celui de~$S$ (puisque~$s\notin W_i$). On peut donc appliquer l'hypothèse de récurrence : presque toutes les composantes connexes de~$W_i\setminus S$ sont des arbres relativement compacts et à bord singleton dans~$W_i$, et sont {\em a fortiori} relativement compacts et à bord singleton dans~$X$ ; soient~$W_{i,1},\ldots,W_{i,r_i}$ les composantes connexes de~$W_i\setminus S$ qui ne sont pas des arbres relativement compacts et à bord singleton dans~$W_i$. 

Par construction, si~$W$ est ou bien une composante connexe de~$X\setminus\{s\}$ qui appartient à~$\Pi$ ou bien une composante connexe de~$W_i\setminus S$ pour l'un des~$W_i$, alors~$W$ est une partie ouverte, connexe et non vide de~$X\setminus S$ telle que~$\partial_X W\subset S$ ; par conséquent,~$W$ est une composante connexe de~$X\setminus S$. 

Réciproquement, si~$W$ est une composante connexe de~$X\setminus S$ et si elle n'appartient pas à~$\Pi$, elle est nécessairement contenue dans~$W_i$ pour un certain~$i$ ; c'est une partie connexe, ouverte et non vide de~$W_i\setminus S$ telle que~$\partial _{W_i}W\subset \partial_X W\subset S$, et partant une composante connexe de~$W_i\setminus S$. 

Les composantes connexes de~$X\setminus S$ sont donc exactement les éléments de~$\Pi$ et les composantes connexes des~$W_i\setminus S$ pour~$i\in \{1,\ldots,r\}$ ; par construction, ce sont toutes des arbres relativement compacts et à bord singleton dans~$X$, à l'exception des~$W_{i,j}$ qui sont en nombre fini, ce qui achève la démonstration. 

\trois{boncompac} Soit~$\Pi$ un sous-ensemble de~$\pi_0(X\setminus S)$. Le  sous-ensemble~$$Y:=\bigcup\limits_{W\in \Pi}\overline W$$ de~$X$ en est un sous-graphe fermé, qui est compact si chacune des composantes appartenant à~$\Pi$ est relativement compacte. 

\medskip
En effet, soit~$T$ la réunion des~$\partial W$ pour~$W$ parcourant~$\Pi$. Toute composante connexe appartenant à~$\Pi$ est alors une composante connexe de~$X\setminus T$ ; comme~$Y$ est réunion de~$T$ et des~$W$ pour~$W\in \Pi$, son complémentaire~$X-Y$ est réunion de composantes connexes de~$X\setminus T$ et est donc ouvert ; ainsi,~$Y$ est un sous-ensemble fermé de~$X$. Il hérite à ce titre du caractère séparé et localement compact de ce dernier.

Soit~$y\in Y$ ; si~$y\notin T$ il est situé sur une certaine composante connexe~$V$ de~$X\setminus T$ qui est contenue dans~$Y$, et il possède alors un voisinage dans~$V$, et {\em a fortiori} dans~$Y$, qui est un arbre. 

Supposons maintenant que~$y\in T$ ; choisissons un voisinage~$U$ de~$y$ dans~$X$ qui est un arbre et qui ne rencontre aucun autre élément de~$T$. Si~$V$ est une composante connexe de~$X\setminus T=(X-(T\setminus\{y\}))\setminus\{y\}$ alors~$V\cap U$ est réunion finie de composantes connexes de~$U\setminus\{y\}$ (\ref{intercompvois}). Par conséquent,~$U\cap Y$ est réunion de~$\{y\}$ et de composantes connexes de~$U\setminus\{y\}$ ; si~$W$ est l'une d'elles, alors pour tout~$x\in W$ l'intervalle~$[x;y]$ est contenu dans~$W\cup\{y\}$ (\ref{adhcompinterv}) et donc dans~$U\cap Y$ ; il en découle que le fermé~$U\cap Y$ de~$U$ (qui est un voisinage ouvert de~$y$ dans~$Y$) est convexe ; c'est par conséquent un arbre (\ref{arbreferme}). 

On déduit de ce qui précède que~$Y$ est un graphe.

\medskip
Supposons maintenant que chacune des composantes appartenant à~$\Pi$ est relativement compacte ; il résulte alors de~\ref{compxmscomp}, appliqué avec~$\sch K=T$, que~$Y$ est compact. 

\subsection*{Définition générale de la valence ; points isolés, unibranches et pluribranches}
\deux{defunipluribr} Si~$X$ est un graphe, si~$x\in X$ et si~$U$ est un voisinage ouvert de~$x$ dans~$X$ qui est un arbre, le cardinal de~$\pi_0(U\setminus\{x\})$ ne dépend que de~$x$, et pas de~$U$ (\ref{arbrebijcom}) ; nous l'appellerons la {\em valence} de~$(X,x)$ ; cette définition est compatible avec la précédente (donnée au~\ref{defgraphelocfin} ) lorsque~$X$ est fini en~$x$. On dira que~$x$ est un point {\em unibranche} (resp. {\em pluribranche}) de~$X$ si la valence de~$(X,x)$ est égale à~$1$ (resp. strictement supérieure à~$1$). Remarquons que la valence de~$(X,x)$ peut être nulle ; cela se produit si et seulement si~$x$ est un point isolé de~$X$. 

\medskip
Notons que si~$x$ appartient à un intervalle ouvert~$I$ tracé sur~$X$, il est nécessairement pluribranche : en effet, supposons que ce soit le cas fixons un voisinage ouvert~$U$ de~$x$ qui est un arbre. Soit~$J$ la composante connexe de~$x$ dans~$U\cap I$ ; c'est un intervalle ouvert. Choisissons deux points~$y$ et~$z$ sur~$J\setminus\{x\}$ situés de part et d'autre de~$x$ ; on a alors~$x\in [y;z]$, et les composantes connexes de~$U\setminus\{x\}$ contenant~$y$ et~$z$ sont donc distinctes, d'où notre assertion. 

\medskip
Donnons un corollaire de cette remarque : si~$X$ est un graphe connexe non réduit à un point, les points pluribranches sont denses dans~$X$ : cela provient du fait que si~$U$ est un ouvert connexe et non vide de~$X$ alors~$U$ n'est lui-même pas réduit à un point et contient donc un intervalle ouvert non vide. 

\deux{bordpluribr}  Soit~$X$ un graphe et soit~$x\in X$. Nous allons montrer qu'il possède un voisinage ouvert~$V$ dans~$X$ possédant les propriétés suivantes : 

\medskip
i)~$V$ est un arbre ; 

ii)~$\overline V$ est un arbre compact ; 

iii)~$\partial V$ est fini, constitué de points pluribranches, et est un singleton si~$x$ est unibranche.

\medskip
Si~$x$ est isolé, l'on peut prendre~$V=\{x\}$ ; on suppose maintenant que~$x$ n'est pas isolé. On se ramène immédiatement, quitte à restreindre~$X$, au cas où ce dernier est un arbre, nécessairement non réduit à~$\{x\}$ ; le point~$x$ possède un voisinage ouvert~$U$ à bord fini et relativement compact dans~$X$ qui n'est pas~$X$ tout entier ; le bord de~$U$ est alors non vide. Le compact~$\Gamma:=\bigcup\limits_{y\in \partial U} [x;y]$ est un sous-arbre fini de~$X$ (\ref{arbreferme}) qui contient~$x$. Soit~$S$ l'ensemble des sommets de~$\Gamma$ autres que~$x$. Pour tout ~$y\in \partial U$ on a~$[x;y[ \subset U(\ref{adhcompinterv})$ ; par conséquent,~$\Gamma\setminus \partial  U$ est égal à~$\Gamma\cap U$ qui est convexe ; on en déduit que~$\partial U\subset S$. La composante connexe~$\Delta$ de~$x$ dans~$\Gamma\setminus S$ est de la forme~$\bigcup\limits_{1\leq i\leq r}[x;y_i[$ où~$r$ est un entier non nul, où~$y_i\in \overline U$ pour tout~$i$, et où~$[x;y_i[\cap [x;y_j[=\{x\}$ dès que~$i\neq j$. Si de plus~$x$ est unibranche l'ouvert~$X\setminus\{x\}$ est connexe, et il résulte alors de~\ref{intervcoinc} que~$r=1$. 

\medskip
Choisissons pour tout~$i$ un un point~$z_i$ sur~$]x;y_i[$ ; chacun des~$z_i$ appartient à un intervalle ouvert tracé sur~$X$ est est de ce fait pluribranche (\ref{defunipluribr}). Soit~$V$ la composante connexe de~$x$ dans~$X\setminus\{z_1,\ldots,z_r\}$. Comme~$[x;z_i[\subset V$ pour tout~$i$, le bord de~$V$ est exactement~$\{z_i\}_{1\leq i\leq r}$ ; il est donc fini, constitué de points pluribranches, et est un singleton dès que~$r=1$, ce qui est notamment le cas dès que~$x$ est unibranche. 

Il suffit pour conclure de prouver que~$V$ est relativement compact ; pour ce faire, on va démontrer qu'il est contenu dans l'ouvert relativement compact~$U$. On raisonne par l'absurde en supposant qu'il existe~$v\in V-U$. Le segment~$[v;x]$ est inclus dans~$V$ par convexité de celui-ci et contient nécessairement un point~$y$ de~$\partial U$ ; l'ensemble~$S\cap [y;x]$ est fini et non vide (car~$y\in S$), et son élément le plus proche de~$x$ est égal à~$y_i$ pour un certain~$i$. On a alors~$[y_i;x]\subset [v;x]\subset V$ et donc~$z_i\in V$, ce qui est contradictoire. 

\subsection*{Existence de points unibranches dans le cas compact}

\deux{infexist} Soit~$X$ un arbre et soit~$x\in X$. Pour tout couple~$(y,z)$ d'éléments de~$X$, la relation~$\leq$ définie par la condition~$z\leq y$ si et seulement si~$y\in [z;x]$ est une relation d'ordre partiel pour laquelle~$x$ est le plus grand élément de~$X$. Si~$y\in X$ l'ensemble~$X_{\leq y}$ des éléments~$z$ tels que~$z\leq y$ est la réunion de~$\{y\}$ et des composantes connexes de~$X\setminus\{y\}$ qui ne contiennent pas~$x$ ; c'est donc un sous-arbre fermé de~$X$ (\ref{boncompac}). 

\medskip
Supposons de plus que~$X$ soit {\em compact}. Pour tout~$y\in X$, l'ensemble~$X_{\leq y}$ est alors compact ; ceci entraîne immédiatement (compte-tenu du fait que~$X$ est non vide) que toute partie totalement ordonnée de~$X$ admet un minorant ; autrement dit,~$X$ est inductif et possède, par le lemme de Zorn, un élément minimal~$x_0$.

Comme~$x$ est le plus grand élément de~$X$, on a~$x_0\neq x$ dès que~$X$ n'est pas réduit à~$\{x\}$. Plaçons-nous sous cette dernière hypothèse. Le point ~$x_0$ est alors unibranche : en effet dans le cas contraire~$X\setminus\{x_0\}$ aurait au moins deux composantes connexes, et donc au moins une composante connexe~$U$ ne contenant pas~$x$ ; mais on aurait alors~$y<x_0$ pour tout~$y\in U$, ce qui contredirait la minimalité de~$x_0$. 

\subsection*{La notion de partie convexe d'un graphe} 

\deux{soussat} Soit~$X$ un graphe ; nous dirons qu'un sous-ensemble~$Y$ de~$X$ est {\em convexe} si pour tout couple~$(x,y)$ de points de~$Y$, tout segment d'extrémités~$x$ et~$y$ tracé sur~$X$ est contenu dans~$Y$ ; dans le cas où~$X$ est un arbre, on retrouve la notion usuelle de convexité. Il découle de la définition que l'intersection d'une famille de parties convexes de~$X$ est convexe. 

\trois{satimplgraphe} Soit~$Y$ une partie localement fermée et convexe de~$X$. Si~$V$ est un sous-arbre de~$X$ alors~$Y\cap V$ est une partie localement fermée et convexe de~$V$ ; c'est donc un sous-arbre de~$X$. 

\medskip
Soit~$y\in Y$. Il existe un voisinage ouvert~$V$ de~$y$ dans~$X$ qui est un arbre ; l'intersection~$V\cap Y$ est un voisinage de~$y$ dans~$Y$, et est un arbre d'après ce qui précède ; en conséquence,~$Y$ est un sous-graphe de~$X$.

\trois{remskgrad} Soit~$Y$ une partie de~$X$ telle que~$X-Y$ soit réunion d'ouverts de~$X$ dont le bord a au plus un élément ; le sous-ensemble~$Y$ de~$X$ est alors un sous-graphe fermé et convexe de~$X$. 

En effet,~$Y$ est fermé en vertu de notre hypothèse ; pour voir que c'est un sous-graphe convexe, il suffit en vertu de~\ref{satimplgraphe} d'en vérifier la convexité.

On raisonne par l'absurde. On suppose donc qu'il existe deux points~$x$ et~$y$ de~$Y$, un segment~$I$ les joignant, et un point~$z$ sur~$J:=I\setminus\{x,y\}$ qui n'appartient pas~$Y$. Par hypothèse, il existe un voisinage ouvert~$U$ de~$z$ dans~$X$ qui ne rencontre pas~$Y$ et tel que~$\partial U$ ait au plus un élément. Soit~$J_0$ la composante connexe de~$z$ dans~$U\cap J$ ; le bord de~$U$ ayant au plus un élément, il existe une extrémité~$t$ de~$J_0$ qui est située sur~$U$. Comme~$U$ ne rencontre pas~$Y$, le point~$t$ n'est égal ni à~$x$, ni à~$y$ ; il appartient donc à~$U\cap J$, ce qui contredit le fait que le composante~$J_0$ est fermée dans~$U\cap J$.

\deux{bordcomplconvgr} Soit ~$X$ un graphe, soit~$\Gamma$ un sous-graphe convexe et fermé de~$X$, et soit~$U$ une composante connexe de~$X-\Gamma$. Le bord de~$U$ comprend alors au plus un élément. En effet, s'il existait deux points distincts~$x$ et~$y$ sur~$\partial U$, l'on pourrait en vertu de~\ref{intervalsousgr} tracer un segment~$I$ reliant~$x$ à~$y$ et tel que~$I\setminus\{x,y\}\subset U$, contredisant ainsi la convexité de~$\Gamma$.

\section{La compactification arboricole} 

\subsection*{L'espace des bouts d'un graphe connexe}

\deux{limindks} Soit~$X$ un graphe connexe. Pour tout sous-ensemble fini~$S$ de~$X$, on note~$\sch K_X(S)$ la réunion de~$S$ et des composantes connexes de~$X\setminus S$ qui sont relativement compactes. C'est un sous-graphe compact de~$X$, et~$X\setminus \sch K_X(S)$ est réunion disjointe d'un nombre fini de composantes connexes de~$X\setminus S$ (\ref{boncompac} et~\ref{compgraphemoinss}). 

\trois{kxscof} La famille des~$\sch K_X(S)$, où~$S$ parcourt l'ensemble des parties finies de~$X$ constituée de points pluribranches, est cofinale dans celle de toutes les parties compactes de~$X$. En effet, soit~$Y$ un compact de~$X$ ; on peut le recouvrir par un nombre fini d'ouverts relativement compacts et à bord fini constitué de points pluribranches (\ref{bordpluribr}). La réunion~$U$ des ouverts en question est-elle même relativement compacte, et son bord~$S$ est fini et constitué de points pluribranches. Toute composante connexe de~$U$ est une composante connexe relativement compacte de~$X\setminus S$ ; par conséquent,~$\sch K(S)\supset U\supset Y$.

\medskip
\trois{kxsunion} Remarquons que si~$S$ et~$S'$ sont deux sous-ensembles finis de~$X$ avec~$S\subset S'$ et si~$U$ est une composante connexe de~$X\setminus S$ alors~$U\setminus S'$ est réunion de composantes connexes de~$X\setminus S'$, qui sont relativement compactes si~$U$ est relativement compacte ; il s'ensuit que~$\sch K_X(S)\subset \sch K_X(S')$. 

\trois{defboutsx} On désignera par~$\got d X$ l'ensemble~$\lim\limits_{\leftarrow}\pi_0(X\setminus\sch K_X(S))$, où~$S$ parcourt l'ensemble des parties finies de~$X$ ; on appellera {\em bout} de~$X$ tout élément de~$\got d X$, et~$\wid X$ la réunion disjointe de~$X$ et~$\got d X$. Remarquons qu'en vertu de~\ref{kxscof},
on pourrait tout aussi bien définir~$\got d X~$ comme étant égal à~$\lim\limits_{\leftarrow}\pi_0(X\setminus \sch K)$ où~$\sch K$ parcourt l'ensemble des parties compactes de~$X$ 

\medskip
Si~$U$ un ouvert de~$X$, il résulte de~\ref{boncompac} que les propositions suivantes sont équivalentes : 

\medskip
i) il existe un ensemble fini~$S\subset X$ tel que~$U$ soit une réunion de composantes connexes (ou encore un ouvert fermé) de~$X\setminus S$ ;

ii)~$\partial  U$ est fini. 

\medskip
Si ces conditions sont vérifiées, le sous-ensemble~$S$ de i) peut toujours être pris égal à~$\partial U$.

\medskip
Soit~$U$ un ouvert de~$X$ à bord fini, et soit~$S$ comme dans i) ; remarquons que l'ouvert~$U$ contient une composante connexe de~$X\setminus \sch K_X(S)$ si et seulement si il n'est pas relativement compact ; on peut voir~$\pi_0(U)$ comme un sous-ensemble de~$\pi_0(X\setminus \partial U)$. L'image réciproque de~$\pi_0(U)$ par la flèche~$$\pi_0(X\setminus S)\to \pi_0(X\setminus \partial  U)$$ est égale à~$\pi_0(U)$ ; par conséquent, l'image réciproque de~$\pi_0(U)$ par la flèche canonique~$\got d X \to \pi_0(X\setminus \sch K_X(S))\subset \pi_0(X\setminus S)$ est égale à  l'image réciproque de~$\pi_0(U)$ par la flèche canonique~$\got d X \to \pi_0(X\setminus \sch K_X(\partial _XU))\subset \pi_0(X\setminus \partial  _XU)$ ; elle ne dépend donc que de~$U$, et pas de~$S$. La réunion de~$U$ et de cette image réciproque est un sous-ensemble de~$\wid X$ que l'on notera~$\widd U X$ ; si~$U$ et~$V$ sont deux ouverts de~$X$ à bord fini et si~$V\subset U$ (resp.~$V\cap U=\emptyset$) alors~$\widd V X\subset \widd U X$ (resp.~$\widd U X \cap \widd V X=\emptyset$) : on prouve la première affirmation en considérant~$V$ comme un ouvert fermé de~$X\setminus (\partial V\cup \partial  U)$, et la seconde en considérant~$V$ et~$U$ comme deux ouverts fermés disjoints de~$X\setminus (\partial  V\cup \partial U)$.

\trois{varpiunonvide} Si~$U$ est un ouvert à bord fini de~$X$ les assertions suivantes sont équivalentes : 

a)~$\widd U  X=U$ ;

b)~$\widd U X \cap \got d X=\emptyset$ ;

c)~$U$ est relativement compact. 

\medskip
En effet, il résulte des définitions que a)$\iff$ b) et que c)~$\Rightarrow$ b). Supposons maintenant que l'ouvert~$U$ n'est pas relativement compact, et montrons que~$\widd U X \cap \got  d X \neq \emptyset$. Comme~$U$ n'est pas relativement compact, il existe une composante connexe de~$U$ qui n'est pas contenue dans~$\sch K_X(\partial U)$, c'est-à-dire qui n'est elle-même pas relativement compacte ; quitte à remplacer~$U$ par cette composante, on peut supposer que~$U$ est connexe et non vide. Soit~$S$ un sous-ensemble fini de~$X$ contenant~$\partial U$. Comme~$U$ n'est pas relativement compact, il n'est pas contenu dans~$\sch K_X(S)$ et appartient donc à l'image de l'application~$\pi_0(X\setminus \sch K_X(S))\to \pi_0(X\setminus \sch K_X(\partial  U))$ ; ceci valant quelle que soit la partie finie~$S$ de~$X$ contenant~$\partial_X U$, il existe un élément de~$\got d X$ dont l'image dans~$\pi_0(X\setminus \sch K_X(\partial U))$ est égale à~$U$, c'est-à-dire un élément de~$\got d X \cap\widd U X$. 

\medskip
On en déduit que~$X$ est compact si et seulement si~$\got d X=\emptyset$ : il suffit d'appliquer ce qui précède avec~$U=X$.

\trois{defxchapeau} On munit~$\wid X$ de la topologie engendrée par les parties de la forme~$\widd U X~$, où~$U$ est ouvert à bord fini de~$X$ ; on définit la même topologie en se restreignant aux ouverts~$U$ à bord fini qui sont connexes et non vides. Par construction,~$X$ est un ouvert de~$\wid X$, et tout point de~$\got d X$ adhère à~$X$ ; par conséquent,~$X$ est dense dans~$\wid X$ ; son bord dans~$\wid X$ est précisément~$\got d X$. La restriction de la topologie de~$\wid X$ à~$\got d X$ coïncide avec la topologie de limite projective sur ce dernier ; dès lors, le sous-ensemble~$\got d X$ de~$\wid X$ est profini. 

\trois{varpiuubarre} Soit~$U$ un ouvert à bord fini de~$X$. L'intersection~$\widd U X\cap \got d X~$ est précisément l'ensemble des points de~$\got d X$ qui adhèrent à~$U$, et également l'ensemble des points de~$\got d X$ qui adhèrent à~$\widd U X$. En effet, soit~$x$ un point de~$\got d X$ ; il existe une unique composante connexe~$V$ de~$X\setminus \partial  U$ telle que~$x\in \widd V X$.

Si~$V\not\subset U$ alors~$x\notin \widd U X$, et~$\widd V X$ est un voisinage ouvert de~$x$ qui ne rencontre pas~$\widd U X$ ; par conséquent,~$x$ n'adhère pas à~$\widd U X$, et {\em a fortiori} pas à~$U$. 

Si~$V\subset U$ alors~$x\in \widd U X$ ; si~$W$ est un ouvert à bord fini de~$X$ tel que~$x\in \widd W X$ alors~$\widd W X \cap \widd U X\neq \emptyset$ ; par conséquent,~$W\cap U\neq \emptyset$, et~$\widd W X$ rencontre ~$U$ ; le point~$x$ adhère donc à~$U$, et {\em a fortiori} à~$\widd U X$, ce qui achève de prouver les équivalences requises. 

\medskip
On déduit de ce qui précède que~$U$ est dense dans~$\widd U X$, ce qui entraîne que~$\widd U X$ est connexe dès que~$U$ est connexe ; on en déduit aussi que~$\partial_{\wid X}\widd U X=\partial_X U$ ; il s'ensuit que~$\partial_{\wid X}\widd U X$ est fini. 

\deux{xchapxomp} {\bf Proposition.} {\em Soit~$X$ un graphe connexe. L'espace topologique~$\wid X$ est compact.} 

\medskip
{\em Démonstration.} Nous procédons en deux temps. 

\trois{xchapsep} {\em L'espace~$\wid X$ est séparé.}  Soient~$x$ et~$y$ deux points distincts  de~$\wid X$ ; nous allons exhiber deux ouverts disjoints de~$\wid X$ contenant respectivement~$x$ et~$y$. 

Si~$x$ et~$y$ appartiennent à~$X$, on le fait en utilisant la séparation de~$X$. 

Supposons que ~$x\in \got d X$ et~$y\in X$ ; soit~$U$ un voisinage ouvert de~$y$ dans~$X$, relativement compact et dont le bord~$S$ est fini, et soit~$V$ la composante connexe de~$X\setminus \sch K_X(S)$ égale à l'image de~$x$. Les ouverts~$\widd UX=U$ et~$\widd V X$ de~$\wid X$ sont disjoints, le premier contient~$x$ et le second~$y$. 

Supposons maintenant que~$x$ et~$y$ appartiennent tous deux à~$\got d X$. Comme ils sont distincts, il existe une partie finie~$S$ de~$X$ telles que les images respectives~$U$ et~$V$ de~$x$ et~$y$ dans~$\pi_0(X\setminus \sch K_X(S))$ soient distinctes ; les ouverts~$\widd U X$ et~$\widd V X$ de~$\wid X$  sont alors disjoints ; le premier contient~$x$ et le second~$y$. 

\trois{xchapqc} {\em L'espace~$\wid X$ est quasi-compact.} Soit~$\sch U$ un recouvrement ouvert de~$\wid X$ ; nous allons en extraire un sous-recouvrement fini. On peut supposer que les ouverts qui constituent~$\sch U$ appartiennent sont tous de la forme~$\widd U X$, où~$U$ est un ouvert à bord fini de~$X$. Comme~$\got d X$ est compact, il existe un ensemble fini~$\{U_1,\ldots,U_n\}$ de tels ouverts de~$X$ tels que les~$\widd {U_i}X$ appartiennent à~$\sch U$ et recouvrent~$\got d X$.

Soit~$S$ la réunion des~$\partial_X U_i$ ; si~$V$ est une composante connexe de~$X\setminus \sch K_X(S)$ alors~$\widd V X$ rencontre~$\got d X$ (\ref{varpiunonvide}) ; par conséquent,~$\widd V X$ rencontre~$\widd {U_i}X$ pour un certain~$i$ ; dès lors,~$V$ rencontre~$U_i$ (\ref{defboutsx}) ; comme~$U_i$ est un ouvert fermé de~$X\setminus \partial _X U_i$, on a~$V\subset U_i$. 

On peut de ce fait écrire~$X=\sch K_X(S)\cup \bigcup U_i$ et donc~$$\wid X=X\cup \got d X =\sch K_X(S)\cup \bigcup \widd {U_i}X\;;$$  comme~$\sch K_X(S)$ est compact, il existe une famille finie~$(V_j)$ d'ouverts de~$\sch U$ qui le recouvrent ; les~$\widd{U_i}X$ et les~$V_j$ constituent alors un sous-recouvrement fini du recouvrement~$\sch U~$ de~$\wid X$.~$\Box$

\deux{lemmeboutsfin} {\bf Lemme.} {\em Soit~$X$ un graphe connexe tel que~$\got d X$ soit fini. Il existe un sous-ensemble fini~$S_0$ de~$X$ tel que pour tout ensemble fini~$S$ de~$X$ vérifiant la condition~$\sch K_X(S)\supset \sch K_X(S_0)$, la flèche naturelle de~$\got d X$ vers~$\pi_0(X\setminus \sch K_X(S))$ soit bijective.}

\medskip
{\em Démonstration.} Comme~$\got d X$ est fini, il existe un sous-ensemble fini~$S_0$ de~$X$ tel que l'application naturelle~$\got d X \to \pi_0(X\setminus \sch K_X(S_0))$ soit injective. Soit~$S$ une partie finie de~$X$ tel que~$\sch K_X(S)\supset \sch K_X(S_0)$. Comme  l'application~$\got d X \to \pi_0(X\setminus \sch K_X(S_0))$ se factorise par~$\pi_0(X\setminus \sch K_X(S))$, la flèche de~$\got d X$ vers ~$\pi_0(X\setminus \sch K_X(S))$ est injective aussi. Elle est par ailleurs surjective en vertu du~\ref{varpiunonvide}, et donc bijective.~$\Box$ 

\subsection*{Propriété universelle de la compactification d'un graphe} 

\deux{propunivwidx} {\bf Théorème.} {\em Soit~$X$ un espace topologique compact et soit~$Y$ un sous-ensemble localement fermé de~$X$ qui est un graphe connexe. Supposons que tout point~$x\in \overline Y-Y$ possède dans~$X$ une base de voisinages ouverts à bord fini, et dont l'intersection avec~$Y$ est connexe. Il existe alors un homéomorphisme canonique~$\overline Y \simeq \wid Y$.}

\medskip
{\em Démonstration.} Pour tout~$x\in \overline Y-Y$ on notera~$\sch V_x$ l'ensemble des voisinages ouverts et à bord fini de~$x$ dont l'intersection avec~$Y$ est connexe ; c'est par hypothèse une base de voisinages de~$x$. 

\medskip
Soit~$x\in \overline Y-Y$ et soit~$S$ un sous-ensemble fini de~$Y$. Il existe~$V\in \sch V_x$ tel que~$V\cap \sch K_Y(S)=\emptyset$. L'intersection~$V\cap Y$ étant connexe et non vide, elle est contenue dans une unique composante connexe~$\varpi(V,S)$ de~$X\setminus \sch K_X(S)$. Celle-ci ne dépend en réalité pas de~$V$ : si~$V'$ est un autre élément de~$\sch V_x$ évitant~$\sch K_Y(S)$ il existe~$V''$ appartenant à~$\sch V_x$ et contenu dans~$V\cap V'$, et il est clair que~$\varpi(V,S)=\varpi(V'',S)=\varpi(V',S)$ ; on peut donc écrire~$\varpi(S)$ au lieu de~$\varpi(V,S)$. L'inclusion~$(V\cap Y)\subset \varpi(S)$ pour tout~$V\in \sch V_x$ permet de caractériser~$\varpi(S)$ comme la seule composante connexe de~$Y\setminus \sch K_Y(S)$ à laquelle~$x$ adhère. 

\medskip
La donnée pour tout~$S$ de la composante~$\varpi(S)$ définit un point de~$\got d Y$ que l'on notera~$\phi(x)$. En posant~$\phi(x)=x$ pour tout~$x\in Y$, on étend~$\phi$ en une application de~$\overline Y$ vers~$\wid Y$, dont nous allons montrer que c'est un homéomorphisme. Il suffit, par compacité, de s'assurer qu'elle est continue et bijective. 

\trois{phiestcont} {\em L'application~$\phi$ est continue.} Comme~$Y$ est localement fermé,~$Y$ est ouvert dans~$\overline Y$ ; par conséquent,~$\phi$ est continue en tout point de~$Y$. Il reste à montrer que~$\phi$ est continue en tout point de~$\overline Y -Y$ ; soit donc~$x$ un tel point. 

\medskip
Soit~$S$ un sous-ensemble fini de~$Y$ et soit~$U$ la composante connexe de~$Y\setminus S$ telle que~$\phi(x)\in \widd U Y$ ; il suffit de prouver que~$\phi\inv(\widd U Y)$ est un voisinage de~$x$ dans~$\overline Y$.

\medskip
Soit~$V\in \sch V_x$ tel que~$V\cap \sch K_Y(S)=\emptyset$ ;  par définition de~$\phi$, on a~$V\cap Y \subset U$. 

Soit~$y\in V\cap \overline Y$. Si~$y\in Y$ alors~$\phi(y)=y\in V\subset U\subset \widd U Y$. Si~$y\in \overline Y-Y$ alors l'inclusion~$(V\cap Y)\subset U$ implique que la seule composante de~$Y\setminus \sch K_Y(S)$ à laquelle~$y$ adhère est nécessairement~$U$ ; il s'ensuit, là encore par définition de~$\phi$, que~$\phi(y)\in \widd U Y$. 

\medskip
Par conséquent, le voisinage~$V\cap \overline Y$ de~$x$ dans~$\overline Y$ est contenu dans~$\widd U Y$, ce qui achève d'établir la continuité de~$\phi$. 

\trois{phiestinj} {\em L'application~$\phi$ est injective.} Par la définition même de~$\phi$, la seule chose qu'il y a à vérifier est l'injectivité de~$\phi_{|\overline Y-Y}$. 

\medskip
Donnons-nous donc deux points distincts~$x$ et~$y$ sur~$\overline Y-Y$. En vertu du caractère séparé de~$Y$ il existe~$V\in \sch V_x$ et~$W\in \sch V_y$ tels que~$V\cap W=\emptyset$. Soit~$S$ l'ensemble fini~$Y\cap (\partial V \cup \partial W)$. L'intersection~$V\cap Y$ est connexe, non vide, ouverte et fermée dans~$Y\setminus S$ ; c'est donc une composante connexe de~$Y\setminus S$ ; il en va de même de~$W\cap Y$. Comme~$V\cap W=\emptyset$, les deux composantes connexes~$Y\cap V$ et~$Y\cap W$ de~$Y\setminus S$ sont distinctes ; par définition de~$\phi$, il vient~$\phi(x)\neq \phi(y)$, ce qui achève de prouver l'injectivité de~$\phi_{|\overline Y-Y}$, et partant celle de~$\phi$. 

\trois{phiestsurj} {\em L'application~$\phi$ est surjective.} Soit~$x\in \wid Y$ ; si~$x\in Y$ alors~$\phi(x)=x$ ; supposons maintenant que~$x\in \got d Y$. Pour toute famille finie~$S$ de points de~$Y$, notons~$\varpi(S)$ la composante connexe de~$Y\setminus \sch K_Y(S)$ égale à l'image de~$x$ ; l'intersection des~$\overline {\varpi(S)}$ est non vide par compacité de~$X$ ; elle est contenue dans~$\overline Y$.  

Soit~$y\in Y$ ; il existe un voisinage ouvert~$U$ de~$y$ dans~$Y$ qui est relativement compact et à bord fini ; comme~$\varpi(\partial_Y U)$ est une composante connexe de~$Y\setminus \partial _Y U$ qui n'est pas relativement compacte,~$y\notin \overline {\varpi(\partial_Y U)}$ ; il s'ensuit que l'ensemble non vide~$\bigcap\limits_S \overline {\varpi(S)}$ est contenu dans~$\overline Y-Y$ ; si~$z$ désigne un point de cet ensemble on a par construction~$\phi(z)=x$ ; par conséquent,~$\phi$ est surjective.~$\Box$ 

\deux{ensboutsfini} {\bf Corollaire.} {\em Soit~$X$ un graphe et soit~$Y$ un sous-graphe connexe de~$X$. Supposons que~$\got d Y$ est fini et que~$\overline Y-Y$ a un cardinal supérieur ou égal à celui de~$\got d Y$. Il existe alors un homéomorphisme canonique~$\overline Y\simeq \wid Y$. }

\medskip
{\em Démonstration.} Quitte à remplacer~$X$ par l'adhérence de~$Y$, on peut supposer que~$Y$ est un ouvert dense de~$X$. Soit~$N$ le cardinal de~$\got d Y$. Soit~$r\geq N$ et soit~$(x_1,\ldots,x_r)$ une famille de points deux à deux distincts de~$\partial Y$. Il existe une famille~$(V_1,\ldots,V_r)$ d'ouverts de~$X$ qui satisfait les conditions suivantes : 

\medskip
$\alpha)$ pour tout~$i$, l'ouvert~$V_i$ est un voisinage connexe et relativement compact de~$x_i$ dont le bord est fini ; 

$\beta)$ les~$V_i$ sont deux à deux disjoints. 

\medskip
Soit~$S$ l'ensemble fini~$Y\cap \left(\bigcup  \partial  V_i\right)$. On a~$Y=\sch K_Y(S)\coprod U_1\coprod\ldots\coprod U_s$, où les~$U_j$ sont les composantes connexes non relativement compactes (dans~$Y$) de~$Y\setminus S$. 

\medskip
Fixons~$i\in\{1,\ldots,r\}$. Le point~$x_i$ adhère à~$Y$ ; comme~$\sch K_Y(S)$ est compact,~$x_i$ n'adhère pas à ~$\sch K_Y(S)$ ; il existe donc un indice~$j(i)$ tel que~$x_i$ adhère à~$U_{j(i)}$ ; la composante~$U_{j(i)}$ rencontre alors~$V_i$. Comme~$Y\cap \partial V_i\subset S$, l'intersection~$U_{j(i)}\cap V_i$ est fermée dans~$U_{j(i)}$ ; étant par ailleurs ouverte et non vide, elle coïncide avec~$U_{j(i)}$ ; autrement dit,~$U_{j(i)}\subset V_i$. 

Ainsi, chacun des ouverts~$V_i$ contient une composante~$U_{j(i)}$. Comme les~$V_i$ sont deux à deux disjoints, l'application~$i\mapsto j(i)$ est injective, d'où il découle que~$s\geq r\geq N$. D'autre part, l'application naturelle~$\got d Y \to \pi_0(Y\setminus \sch K_Y(S))$ est surjective (\ref{varpiunonvide}), d'où la majoration~$s\leq N$, puis les égalités~$s=r=N$. Cela implique que le cardinal de~$\partial Y$ est égal à~$N$, et que~$\{x_1,\ldots,x_r\}=\partial Y$ ; on renumérote les~$U_j$ de sorte que~$U_i$ soit contenu dans~$V_i$ pour tout~$i$ . 

\medskip
Fixons~$i$. Comme~$V_i\cap V_j=\emptyset$ pour tout~$j\neq i$, l'ouvert~$V_i$ ne rencontre pas~$S$. Si~$U$ est une composante connexe de~$Y\setminus S$ qui diffère de~$U_i$ l'intersection~$U\cap V_i$ est un ouvert fermé strict de~$V_i$, et est donc vide par connexité de ce dernier ; par conséquent,~$V_i\cap X$ coïncide avec~$U_i$, et est en particulier connexe. 

\medskip
Si~$W$ est un voisinage ouvert connexe à bord fini de~$x_i$ dans~$V$, la famille~$(V_1,\ldots,V_{i-1},W,V_{i+1},\ldots,V_N)$ satisfait encore les conditions~$\alpha)$ et~$\beta)$ ci-dessus ; par conséquent,~$W\cap Y$ est connexe. Ainsi,~$x_i$ possède une base de voisinages à bord fini dont l'intersection avec~$Y$ est connexe. 

\medskip
Le graphe~$X$ est la réunion de~$\sch K_Y(S)$, des~$U_i$ et de~$\partial Y=\{x_1,\ldots,x_N\}$ ; puisque~$U_i\subset V_i\subset \overline{V_i}$ pour tout~$i$, il vient~$X=\sch K_Y(S)\cup \coprod \overline{V_i}$ ; comme~$\sch K_Y(S)$ et les~$\overline {V_i}$ sont compacts,~$X$ est compact. Par ce qui précède, chacun des points de~$\partial Y$ a une base de voisinages à bord fini dont l'intersection avec~$Y$ est connexe ; le théorème~\ref{propunivwidx} assure alors qu'il existe un homéomorphisme canonique entre~$\overline Y$ (qui ici n'est autre que~$X$) et~$\wid Y$.~$\Box$

\deux{adhouvdensecpct} {\bf Corollaire.} {\em Soit~$X$ un arbre et soit~$Y$ un sous-arbre de~$X$. La compactification~$\overline Y^{\wid X}$ de~$Y$ s'identifie à~$\wid Y$.} 

\medskip
{\em Démonstration.} Soit~$x\in \partial_{\wid X}Y$. Soit~$S$ un sous-ensemble fini de~$X$ et soit~$U$ la composante connexe de~$X\setminus S$ telle que~$x\in \widd U X$. On a les égalités~$$\widd U X \cap Y=(\widd U X\cap X)\cap Y=U\cap Y.$$ Par convexité,~$U\cap Y$ est connexe ; par conséquent,~$x$ possède dans~$\wid X$ une base de voisinages ouverts et à bord fini dont l'intersection avec~$Y$ est connexe, et le corollaire découle alors aussitôt du théorème~\ref{propunivwidx} . ~$\Box$

\deux{corollbordprof} {\bf Corollaire.} {\em Soit~$X$ un graphe et soit~$U$ un sous-graphe ouvert et connexe de~$X$. Le bord de~$U$ dans~$X$ possède une base d'ouverts compacts ; il est en particulier totalement discontinu, et profini s'il est compact.}

\medskip
{\em Démonstration.} Soit~$x\in \partial U$ ; il s'agit de montrer qu'il possède une base de voisinages ouverts compacts dans~$\partial U$. Soit~$V$ un voisinage de~$x$ dans~$X$ qui est un arbre compact et à bord fini. En vertu de~\ref{interconfin}, l'ensemble~$\pi_0(U\cap V)$ est fini ; soient~$U_1,\ldots,U_r$ les composantes connexes de~$V\cap U$ dont l'adhérence contient~$x$. La question étudiée étant locale, il suffit de démontrer que~$x$ possède une base de voisinages ouverts compacts dans~$\bigcup \partial_V U_i$. 

\medskip
Soit~$W$ une composante connexe de~$V\setminus\{x\}$. Soient~$i$ et~$j$ deux indices tels que~$U_i$ et~$U_j$ soient incluses dans~$W$. Soit~$y\in U_i$ et soit~$z\in U_j$ ; comme~$x$ adhère à~$U_i$ et~$U_j$, les deux intervalles~$[y;x[$ et~$[z;x[$ sont respectivement inclus dans~$U_i$ et~$U_j$ ; leur intersection étant non vide (elle est de la forme~$[t;x[$ avec~$t\in W$), on a nécessairement~$i=j$. 

\medskip
Ainsi chaque composante connexe de~$V\setminus\{x\}$ contient {\em au plus} un ouvert~$U_i$ ; on en déduit que~$\partial_V U_i\cap \partial_V U_j=\{x\}$ pour tout couple~$(i,j)$ avec~$i\neq j$. 

\medskip
Fixons~$i$. Il résulte du corollaire~\ref{adhouvdensecpct} que~$\overline {U_i}^V$ s'identifie à~$\wid {U_i}$ ; par conséquent,~$\partial_V U_i$ est homéomorphe à~$\got d U_i$, et est en particulier profini ; il s'ensuit que~$x$ possède une base de voisinages ouverts compacts dans~$\partial_V U_i$. 

\medskip
Soit~$\Omega$ un voisinage de~$x$ dans~$\bigcup \partial_V U_i$. Par ce qui précède, il existe pour tout~$i$ un sous-ensemble compact~$\sch K_i$ de~$\partial_V U_i$ dont le complémentaire~$\sch K'_i$ dans~$\partial_V U_i$ est compact, contient~$x$ et est contenu dans~$\Omega$ ; si~$i$ et~$j$ sont deux indices distincts alors~$\sch K_i \cap \partial_V U_j=\emptyset$ puisque~$\partial_V U_i\cap \partial_V U_j=\{x\}$. 

 La réunion des~$\sch K_i$ est un compact de~$\bigcup \partial_V U_i$ dont le complémentaire dans~$\bigcup \partial_V U_i$ est le compact~$\bigcup \sch K'_i$, qui contient~$x$ et est contenu dans~$\Omega$ ; ceci achève la démonstration.~$\Box$

\deux{defabout} Si~$X$ est un arbre, si~$Y$ est un sous-arbre de~$X$, et si~$\omega\in \got d X$, nous dirons que~$Y$ {\em aboutit} à~$\omega$ si ce dernier appartient à~$\got d Y$ modulo l'identification entre~$\wid Y~$ et~$\widd Y {\wid X}$ fournie par le corollaire~\ref{adhouvdensecpct}.  

\subsection*{L'enveloppe convexe d'une partie compacte d'un arbre}

\deux{segmlong} On dira qu'un espace topologique est une {\em droite éventuellement longue} s'il est connexe, non vide, non compact, et localement homéomorphe à~$\RR$. Une droite éventuellement longue est une «vraie» droite ({\em i.e.}, est homéomorphe à~$\RR$) si et seulement si elle est paracompacte, ou encore si et seulement si elle est métrisable ; toute droite éventuellement longue est un arbre à deux bouts. 

\medskip
Si~$x$ et~$y$ sont deux points d'un espace topologique~$X$, on dira que~$X$ est un {\em segment éventuellement long d'extrémités~$x$ et~$y$} si~$X=\{x\}=\{y\}$ ou s'il existe une droite éventuellement longue~$I$ et un homéomorphisme~$X\simeq \wid I$ envoyant~$\{x,y\}$ sur~$\got d I$. 

\medskip
Si~$x$ est un point d'un espace topologique~$X$, on dira que~$X$ est une {\em demi-droite éventuellement longue issue de~$x$} s'il existe une droite éventuellement longue~$I$ et un plongement topologique de~$X$ dans~$\wid I$ envoyant~$x$ sur l'un des deux points de~$\got d I$ et identifiant~$X\setminus\{x\}$ à~$I$. 

\deux{propuniondrl} {\bf Proposition.} {\em  Soit~$X$ un arbre et soit~$S$ un sous-ensemble fermé et discret de~$X$. Soit~$(I_j)_{j\in J}$ une famille de sous-arbres fermés de~$X$ tels que chacun des~$I_j$ soit ou bien une droite éventuellement longue, ou bien une demi-droite éventuellement longue issue d'un point de~$S$. La réunion~$\Gamma$ des~$I_i$ est alors un sous-graphe fermé localement fini de~$X$.}

\medskip
{\em Démonstration.} Soit~$V$ un sous-arbre ouvert de~$X$ relativement compact, à bord fini, et rencontrant au plus un point de~$S$ ; nous allons montrer que~$\Gamma\cap V$ est un sous-graphe fini et fermé de~$V$, ce qui suffira à conclure puisque~$X$ admet un recouvrement par de tels~$V$. 

\medskip
Soit~$\sch I$ l'ensemble des sous-arbres~$I$ de~$V$ tels que l'une des deux propriétés suivantes soit satisfaite : 

\medskip
$\bullet$ il existe deux point distincts~$x$ et~$y$ de~$\partial V$ tels que I=$]x;y[$ ; 

$\bullet$~$V\cap S$ est un singleton~$\{z\}$ et il existe~$t\in \partial V$ tel que~$I=[z;t[$. 

\medskip
L'ensemble~$\sch I$ est fini, et chacun de ses éléments est un sous-arbre fermé et fini de~$V$. 

\medskip
Soit~$j\in J$. L'intersection~$I_j\cap V$ est un ouvert convexe et relativement compact de~$I_j$, ce qui implique que~$\overline{V\cap I_j}$ est ou bien vide, ou bien un segment tracé sur le fermé~$I_j$. Supposons que~$V\cap I_j$ est non vide, et distinguons deux cas. 

{\em Le cas où~$I_j$ est une droite éventuellement longue.} L'intersection~$I_j\cap V$ est alors de  la forme~$]x;y[$ où~$x$ et~$y$ sont deux points distincts de~$I_j$ situés sur~$\partial V$, et l'on a donc~$(I_j\cap V)\in \sch I$. 

{\em Le cas où~$I_j$ est une demi-droite éventuellement longue issue d'un point~$z$ de~$S$.} L'intersection~$I_j\cap V$ peut sous cette hypothèse : ou bien être de la forme~$]x;y[$ où~$x$ et~$y$ sont deux points de~$I_j$ situés  sur~$\partial V$, et l'on a alors~$(I_j\cap V)\in \sch I$ ; ou bien être de la forme~$[z;t[$ où~$t$ est un point de~$I_j$ situé sur~$\partial V$, et l'on a encore~$(V\cap I_j)\in \sch I$. 

\medskip
Ainsi~$V\cap I_j$ appartient-il à~$\sch I$ dès qu'il est non vide. Par conséquent, ~$V\cap \Gamma$ est réunion d'éléments de~$\sch I$. Comme~$\sch I$ est un ensemble fini de sous-arbres fermés et finis de~$V$, l'intersection~$V\cap \Gamma$ est un sous-graphe fermé et fini de~$V$, ce qui achève la démonstration.~$\Box$  

\deux{defenvconvcomp} Soit~$X$ un arbre et soit~$E$ une partie de~$X$. La réunion~$\conv E$ des segments~$[x;y]$, où~$(x,y)$ parcourt~$E^2$, est visiblement {\em l'enveloppe convexe} de~$E$, c'est-à-dire la plus petite partie convexe de~$X$ contenant~$E$. 

\deux{theoenvconv} {\bf Théorème.} {\em Soit~$X$ un arbre et soit~$E$ une partie compacte de~$X$. L'ensemble~$\conv E$ est un sous-arbre compact de~$X$. Si~$x$ est un point de~$\conv E$ tel que~$x\notin E$ ou tel que~$x$ soit un point isolé de~$E$, l'arbre~$\conv E$ est fini en~$x$.}

\medskip
{\em Démonstration.} On procède en plusieurs étapes. 

\trois{envconvredcomp} {\em Réduction au cas où~$X$ est compact.} Tout point de~$X$ ayant un voisinage qui est un arbre compact, on peut recouvrir le compact~$E$ par une famille finie~$V_1,\ldots,V_n$ d'arbres compacts non vides. Pour tout~$i$, on choisit un point~$x_i$ sur~$V_i$. La réunion des~$V_i$ et des segments~$[x_i;x_j]$ où~$(i,j)$ parcourt~$\{1,\ldots,n\}$  est un compact convexe, et donc un sous-arbre compact de~$X$ ; par construction, il contient~$E$, et~$\conv E$ ; en le substituant à~$X$, on se ramène au cas où ce dernier est compact.

\trois{debutenvconv} Posons~$U=X\setminus E$ ; c'est une réunion disjointe d'arbres. Si~$V$ désigne une composante connexe de~$U$ alors~$V\cap \conv E$ est une réunion d'intervalles ouverts de la forme~$]y;z[$ avec~$y\in E$ et~$z\in E$ ; ils sont tous fermés dans~$V$. Il résulte de la proposition~\ref{propuniondrl} que~$V\cap \conv E$ est un sous-graphe fermé et localement fini de~$V$. Par conséquent,~$U\cap \conv E$ est un sous-graphe fermé et localement fini de~$U$ ; cela entraîne que~$E\cup(U\cap \conv E)=\conv E$ est une partie fermée, et partant compacte, de~$X$. Comme elle est convexe par construction, c'en est un sous-arbre compact, dont on vient de voir qu'il est fini en tout point de~$U$.

\trois{envconvpointiso} Soit~$x\in E$ tel que~$E\setminus\{x\}$ soit compact ; posons~$E'=E\setminus\{x\}$, et soit~$U'$ la composante connexe de~$x-$ dans~$X\setminus E'$. L'intersection~$U'\cap \conv E$ est réunion d'intervalles ouverts de la forme~$]y;z[$ avec~$y\in E'$ et~$z\in E'$ et d'intervalles semi-ouverts de la forme~$[x;t[$ avec~$t\in E'$ ; chacun d'eux est fermé dans~$U'$. En vertu de la propostion ~\ref{propuniondrl},~$U'\cap \conv E$ est un graphe localement fini ; en particulier,~$\conv E$ est fini en~$x$.~$\Box$ 

\deux{compdeninf} {\bf Proposition.} {\em Soit~$X$ un arbre compact et soit~$U$ un sous-arbre ouvert de~$X$. L'espace topologique~$U$ est dénombrable à l'infini.}

\medskip
{\em Démonstration.} Si~$U$ est vide il n'y a rien à démontrer. On suppose~$U$ non vide, on choisit~$x\in U$, et l'on munit~$X$ de la relation d'ordre partiel définie par~$x$ (\ref{infexist}). Le fermé~$\partial U$ est une partie compacte et totalement discontinue de~$X$ (cor.~\ref{corollbordprof}). L'enveloppe convexe~$\conv (\partial U \cup\{x\})$ coïncide avec la réunion des~$[y;x]$ où~$y$ parcourt~$\partial U$. Comme~$x$ est un point isolé de~$\partial U\cup\{x\}$, l'intersection~$$\Gamma:=\conv {(\partial U\cup\{x\})}\cap U$$ est un arbre localement fini en vertu du théorème~\ref{theoenvconv} ci-dessus ; il contient~$]y;x]$ pour tout~$y\in \partial U$ (\ref{adhcompinterv}). On note~$S$ l'ensemble des sommets de~$\Gamma$ ; c'en est une partie fermée et discrète.

\trois{finsomgamau} Soit~$y\in \partial  U$. Nous allons définir une suite~$(\sigma_i(y))_i$ de points de~$]y;x]$ par le procédé récursif suivant. 

On pose~$\sigma_0(y)=x$. Si~$i$ est un entier tel que~$\sigma_i(y)$ soit défini, on distingue deux cas : ou bien~$S\cap ]y;\sigma_i(y)[$ est non vide, et l'on définit alors~$\sigma_{i+1}(y)$ comme son plus grand élément (qui existe en raison du caractère fermé et discret de~$S$) ; ou bien il est vide, et l'on interrompt alors la construction de la suite, qui ne sera donc définie que jusqu'au rang~$i$. 

\medskip
Par construction,~$(\sigma_i(y))_i$ est une suite strictement décroissante de points de~$]y;x]$, et~$\{\sigma_i(y)\}_{i>0}$ est égal à~$S\cap ]y;x[$. Il résulte alors du caractère fermé et discret de~$S$ que si la suite~$(\sigma_i(y))_i$ est infinie,~$\sigma_i(y)$ tend vers~$y$ quand~$i$ tend vers l'infini. 

\medskip
Supposons que~$(\sigma_i(y))_i$ soit finie, et soit~$j$ le plus grand indice de son domaine de définition. L'intervalle~$[\sigma_j(y);y[$ ne comporte alors aucun sommet de~$\Gamma$ ; on choisit une suite arbitraire et strictement décroissante de points de~$]\sigma_j(y);y[$ qui tend vers~$y$, et que l'on numérote~$\sigma_{j+1}(y),\sigma_{j+2}(y)$, etc. On prolonge ainsi la suite finie~$(\sigma_i(y))_i$ en une suite définie sur~$\NN$. 

\trois{indepmarque} Supposons qu'un point~$z$ de~$\Gamma$ soit égal à~$\sigma_j(y)$ pour un certain entier~$j$ et un certain~$y\in \partial U$, et qu'il existe un point~$y'$ de~$\partial  U$ distinct de~$y$ et tel que~$z\in [y';x]$ ; nous allons montrer que~$z=\sigma_j(y')$. 

\medskip
Si~$z=x$ on a~$j=0$ et~$z=\sigma_0(y')$ ; si~$z\neq x$ alors~$j\neq 0$. Plaçons-nous sous cette dernière hypothèse ; l'intersection de~$[y;x]$ et~$[y';x]$ contient~$z$, elle est donc de la forme~$[t;x]$ avec~$t\neq x$ ; comme~$y\neq y'$, comme~$]y;x]\subset \Gamma$ et comme~$]y';x]\subset \Gamma$, on a~$t\in \Gamma$.

Par construction,~$[y;t]\cap[y';t]=[y;t]\cap [t;x]=[y';t]\cap [t;x]=\{t\}$ ; par conséquent, la	valence de~$(\Gamma,t)$ est au moins 3, et~$t$ est un sommet de~$\Gamma$. L'intersection~$S\cap [t;x[$ est finie et contient~$t$ ; elle est donc de la forme~$\{t_1,\ldots t_r\}$ où~$(t_i)$ est une suite strictement décroissante telle que~$t_r=t$. On a par construction~$\sigma_i(y)=t_i$ et~$\sigma_i(y')=t_i$ pour tout~$i\in\{1,\ldots,r\},$ et~$\sigma_i(y)<t$ et~$\sigma_i(y')<t$ pour tout~$i>r$ ; comme~$z\in [t;x[$ il est égal à~$t_s$ pour un certain~$s$ compris entre~$1$ et~$r$ ; on a par conséquent,~$j=s$ et~$z=\sigma_j(y')$, comme annoncé. 

\trois{defptmarque} Qualifions de {\em point marqué} de~$\Gamma$ tout point de la forme~$\sigma_i(y)$ pour un certain~$i$ et un certain~$y\in \partial U$ ; en vertu de~\ref{indepmarque}, l'entier~$i$ ne dépend pas du choix de~$y$ ; on l'appellera {\em l'indice} du point marqué en question. Si~$y\in \partial U$, il découle de~\ref{indepmarque} que l'ensemble des points marqués de~$\Gamma$ appartenant à~$[y;x]$ est exactement l'ensemble~$\{\sigma_i(y)\}_i$. Par construction, tout sommet de~$\Gamma$ est marqué. 

\medskip
Nous allons montrer par récurrence sur~$i$ que l'ensemble des points marqués de~$\Gamma$ d'indice~$i$ est fini. C'est évident pour~$i=0$, puisque~$x$ est le seul point marqué d'indice~$0$. Supposons que ce soit vrai à un certain rang~$i>0$, et montrons-le au rang~$i+1$. Si~$z$ est un point marqué d'indice~$i+1$, l'ensemble des points marqués de l'intervalle~$[z;x]$ est une suite strictement décroissante~$x=z_0,z_1,\ldots, z_{i+1}=z$, où chacun des~$z_j$ est d'indice~$j$. L'intervalle~$]z;z_i[$ ne contient aucun point marqué de~$\Gamma$. 

\medskip
Soit maintenant~$t$ un point marqué d'indice~$i$ de~$\Gamma$ et soient~$w$ et~$w'$ deux points marqués distincts d'indice~$i+1$ tels que~$]w;t[$ et~$]w';t[$ ne contiennent aucun point marqué de~$\Gamma$. L'intersection~$[w;t]\cap [w';t]$ est de la forme~$[w'';t]$. Comme~$w$ et~$w'$ sont tous deux marqués et distincts, on ne peut avoir~$w\in [w';t[$ ou~$w'\in [w;t[$ ; par conséquent,~$w''\neq w$ et~$w''\neq w'$. Si~$w''$ était distinct de~$t$, on aurait alors~$[w;w'']\cap [w';w'']=[w;w'']\cap [w'';t]=[w';w'']\cap [w'';t]=\{w''\}$, et la valence de~$(\Gamma,w'')$ serait au moins égale à~$3$, ce qui est absurde puisque~$]w;t[$ ne contient aucun point marqué de~$\Gamma$, et {\em a fortiori} aucun sommet. On a dès lors~$w''=t$ ; par conséquent,~$[w;t]\cap [w';t]=\{t\}$. 

Il s'ensuit que l'ensemble des points marqués~$w$ d'indice~$i+1$ tels que~$]w;t[$ ne contienne aucun point marqué de~$\Gamma$ est fini, et plus précisément de cardinal majoré par la valence de~$(\Gamma,t)$. 

\medskip
On déduit de ce qui précède et de la finitude du nombre de points marqués d'indice~$i$ de~$\Gamma$, qui constitue notre hypothèse de récurrence, que le nombre de points marqués d'indice~$i+1$ de~$\Gamma$ est fini, ce qu'on souhaitait établir. 

\trois{famdencomp} Pour tout~$i>0$, notons~$E_i$ l'ensemble fini des points marqués d'indice~$i$ de~$\Gamma$, et désignons par~$U_i$ la composante connexe de~$U\setminus E_i$ qui contient~$x$. 

\medskip
{\em L'adhérence~$\overline {U_i}$ est contenue dans~$U$.} Pour le voir, on raisonne par l'absurde : si ce n'était pas le cas, il existerait~$y\in \partial U$ qui soit adhérent à~$U_i$. L'intervalle~$]y;x]$ serait alors contenu dans l'ouvert~$U_i$ (\ref{adhcompinterv}) ; un tel intervalle rencontrant~$E_i$, on aboutirait à une contradiction. 

\medskip
{\em Tout point de~$U$ appartient à~$U_i$ pour un certain~$i$.} En effet, soit~$y\in U$. Nous allons montrer que~$[y;x]$ ne contient qu'un nombre fini de points marqués de~$\Gamma$ ; cela entraînera l'existence d'un entier~$i$ tel que~$[y;x]\cap E_i=\emptyset$, c'est-à-dire tel que~$y\in U_i$. On distingue deux cas. 

\medskip
Si~$y\in \Gamma$ il est situé sur un intervalle~$[z;x]$ pour un certain~$x\in \partial U$, et l'intervalle~$[y;x]\subset [z;x]$ ne contient alors qu'un nombre fini de points marqués de~$\Gamma$. 

Supposons maintenant que~$y\notin \Gamma$, et soit~$V$ la composante connexe de~$U-\Gamma$ contenant~$y$. Le bord de~$V$ dans~$U$ est un singleton~$\{t\}$ pour un certain~$t\in \Gamma$, et~$[y;x]\cap \Gamma=[t;x]$ (\ref{arbreestadm}) ; par le raisonnement suivi ci-dessus,~$[t;x]$ ne contient qu'un nombre fini de points marqués de~$\Gamma$, et l'intervalle~$[y;x]$ ne contient donc lui-même qu'un nombre fini de points marqués de~$\Gamma$. 

\medskip
{\em Conclusion.} Pour tout entier~$i$, le compact~$\overline {U_i}$ de~$X$ est contenu dans~$U$ ; comme~$U$ est égal à la réunion des~$U_i$, il est {\em a fortiori} égal à la réunion des~$\overline {U_i}$, et est de ce fait dénombrable à l'infini.~$\Box$

\deux{paracompgen} {\bf Théorème.} {\em Soit~$X$ un graphe paracompact et soit~$U$ un ouvert de~$X$ ; l'espace~$U$ est paracompact.} 

\medskip
{\em Démonstration.} On se ramène immédiatement au cas où les graphes~$U$ et~$X$ sont connexes ; l'espace~$X$ est alors dénombrable à l'infini, et il s'agit de montrer qu'il en va de même de~$U$. 

\medskip
Tout point de~$X$ possède un voisinage qui est un arbre compact à bord fini (\ref{bordpluribr}). Par conséquent,~$X$ est réunion dénombrable de sous-arbres compacts à bord fini. Soit~$Y$ un tel sous-arbre. L'intersection~$Y\cap U$ est un sous-graphe de~$X$ (\ref{sousgrapheinter}) et est donc localement connexe. D'après le~\ref{interconfin}, le sous-graphe~$Y\cap U$ a un nombre fini de composantes connexes. En vertu de la proposition~\ref{compdeninf}, chacune d'elle est dénombrable à l'infini ; il s'ensuit que~$Y\cap U$ est dénombrable à l'infini. Ceci valant quel que soit~$Y$, l'ouvert~$U$ est dénombrable à l'infini.~$\Box$

\deux{xchapeaupresquarb} {\bf Théorème.} {\em Soit~$X$ un arbre et soient~$x$ et~$y$ deux points de~$\wid X$. Il existe un unique fermé~$[x;y]$ de~$\wid X$ qui soit un segment éventuellement long extrémités~$x$ et~$y$, et~$]x;y[$ est le sous-ensemble de~$X$ formé des points~$z$ tels que~$x$ et~$y$ appartiennent à deux composantes connexes distinctes de~$\wid X\setminus\{z\}$. Si~$U$ est un ouvert connexe et à bord fini de~$X$ tel que~$x$ et~$y$ appartiennent à~$\widd U X$ alors~$[x;y]\subset \widd U X$.}

\medskip
{\em Démonstration.} Si~$x=y$ le théorème est évident ; on suppose maintenant que~$x\neq y$.

\trois{constxy} Soit~$I$ l'ensemble des points~$z$ de~$X\setminus\{x,y\}$ tels que~$x$ et~$y$ soient situés sur deux composantes connexes différentes de~$\wid X\setminus\{z\}$. Rappelons que si~$U$ est un ouvert à bord fini de~$X$, on a~$\partial_{\wid X} \widd U X=\partial_X U$ (\ref{varpiuubarre}). 

\medskip
Soient~$U$ et~$V$ deux ouverts connexes, non vides et à bord fini de~$X$, tels que~$x\in \widd U X~$ et~$y\in \widd V X~$, et tels que~$(\widd U X \cup \partial_X U)\cap ( \widd V X\cup \partial_X V)=\emptyset$. Choisissons~$z\in U$ et~$t\in V$. Comme~$t\notin U$, l'intersection de~$[z;t]$ et de~$U$ est de la forme~$[z,u[$ avec~$u\in \partial_X U$ ; de même,~$[z;t]\cap V$ est de la forme~$]v;t]$ avec~$v\in \partial_X V$. Comme~$U\cup \partial_X U$ et~$V\cup \partial_X V$ sont disjoints,~$u<v$ si l'on oriente~$[z;t]$ de~$z$ vers~$t$ ; nous allons montrer que l'intersection de~$I$ avec~$X-(U\cup V)$ est égale à~$[u;v]$ (cela entraînera que~$u$ et~$v$ ne dépendent pas de~$z$ et~$t$, ce que l'on pourrait vérifier directement). 

\medskip
On procède par double inclusion. Soit~$w$ un élément de~$X-(U\cup V)$ qui n'appartient pas à~$[u;v]$. Comme~$u\in \partial_X U$ et comme~$v\in \partial_X V$, la réunion~$$\widd U X \cup[u;v]\cup \widd V X$$ est une partie connexe de~$\wid X\setminus\{w\}$, qui contient~$x$ et~$y$ ; par conséquent,~$w\notin I$ ; il s'ensuit que~$I\cap (X-(U\cup V))\subset [u;v]$. 

\medskip
Réciproquement, soit~$w\in [u;v]$ et soit~$U'$ (resp.~$V'$) la composante connexe de~$X\setminus\{w\}$ qui contient~$U$ (resp.~$V$). Comme~$w\in [z;t]$ ces deux composantes sont distinctes. Il s'ensuit que~$\widd {U'}X$ et~$\widd {V'} X$ sont disjoints. Or~$\wid X\setminus\{w\}$ est réunion disjointe des~$\widd W X$, où~$W$ parcourt~$\pi_0(X\setminus\{w\})$ ; par conséquent,~$\widd {U'} X$ et~$\widd {V'}X$ sont deux composantes connexes disjointes de~$\wid X\setminus\{w\}$ ; la première contient~$x$ et la seconde~$y$, ce qui montre que~$w\in I$, et partant que~$[u;v]\subset I\cap (X-(U\cup V))$. 

\medskip
Comme~$\wid X$ est compact, les ouverts de la forme~$\widd U X \coprod \widd V X$, où~$U$ et~$V$ sont comme ci-dessus, forment une base de voisinages de~$\{x,y\}$ ; compte-tenu de ce qui précède, on en déduit les faits suivants : 

\medskip

$\bullet$~$I$ est non vide ; 

$\bullet$ si~$w$ et~$w'$ sont deux points de~$I$ il existe un ouvert de~$I$ contenant~$w$ et~$w'$ et homéomorphe à~$\RR$ ;

$\bullet$ le compact~$\overline I ^{\wid X}$ contient~$x$ et~$y$. 

\medskip
Il en découle que~$I$ est une variété topologique connexe et non vide de dimension~$1$ ; comme son bord dans~$\wid X$ est non vide (il contient~$x$ et~$y$), l'espace topologique~$I$ n'est pas compact ; c'est donc une droite éventuellement longue. 

\medskip
L'ensemble~$\got d I$ comprend deux éléments, et~$\partial_{\wid X}I$ comprend au moins deux éléments, à savoir~$x$ et~$y$ ; en vertu de la proposition~\ref{ensboutsfini},~$\overline I ^{\wid X}$ s'identifie alors à~$\wid I$, et est donc un segment éventuellement long d'extrémités~$x$ et~$y$. 

\trois{uniqui} Montrons maintenant que~$\overline I ^{\wid X}$ est le seul fermé de~$\wid X$ qui soit un segment éventuellement long d'extrémités~$x$ et~$y$. Soit~$J$ un tel fermé et soit~$w$ un point de~$X\setminus\{x,y\}$ qui n'appartient pas à~$J$. Comme~$J$ est connexe et contient~$\{x,y\}$, les points~$x$ et~$y$ sont situés sur la même composante connexe de~$\wid X\setminus\{w\}$ ; par conséquent,~$w\notin I$. Il s'ensuit que~$I\subset J$ ; par compacité de~$J$, on a~$\overline I ^{\wid X}\subset J$. Si~$z$ est un point de~$J$ différent de~$x$ et~$y$ alors~$J\setminus\{z\}$ a deux composantes connexes ; il en découle que toute partie connexe de~$J$ qui contient~$x$ et~$y$ coïncide avec~$J$ ; en particulier,~$\overline I ^{\wid X}=J$.

\trois{xytracesuru} Soit~$U$ un ouvert connexe et à bord fini de~$X$ tel que~$x$ et~$y$ appartiennent à~$\widd U X$. Si~$w$ est un point de~$X-U$ et si~$U'$ désigne la composante connexe de~$X\setminus\{w\}$ contenant~$U$ (celui-ci est non vide puisque~$x$ et~$y$ appartiennent à~$\widd U X$) alors~$x$ et~$y$ appartiennent à~$\widd {U'}X$ ; il s'ensuit que ~$w\notin]x;y[$ ; comme on a par ailleurs~$]x;y[\subset X$, cela signifie que~$]x;y[\subset U$ ; par conséquent,~$[x;y]\subset \widd U X$, ce qui achève la démonstration.~$\Box$ 

\deux{remordrelong} {\em Remarque.} Soient~$X,x$ et~$y$ comme dans le théorème ci-dessus, et soit~$z$ un point de~$\got d X\setminus\{x,y\}$. Comme~$]x;y[\subset X$, on a~$[x;y]\subset \wid X\setminus\{z\}$ ; le segment éventuellement long~$[x;y]$ étant connexe,~$x$ et~$ y$ sont situés sur la même composante connexe de~$\wid X\setminus\{z\}$ ; ainsi,~$]x;y[$ peut être caractérisé comme l'ensemble des points~$z$ de~$\wid X$ (et pas seulement de~$X$) tels que~$x$ et~$y$ soient situés sur deux composantes connexes distinctes de~$\wid X\setminus\{z\}$ : les résultats du~\ref{descxycomp} s'étendent donc {\em mutatis mutandis} à~$\wid X$.

\medskip
Il en va de même, de façon évidente, de ceux énoncés au début du~\ref{infexist} : si~$x\in \wid X$ la relation~$\leq$  définie par la condition~$z\leq y$ si et seulement si~$y\in [z;x]$ est une relation d'ordre partiel pour laquelle~$x$ est le plus grand élément de~$\wid X$ ; si~$y\in \wid X$ l'ensemble~$\wid X_{\leq y}$ des éléments~$z$ tels que~$z\leq y$ est la réunion de~$\{y\}$ et des composantes connexes de~$\wid X\setminus\{y\}$ qui ne contiennent pas~$x$ ; par conséquent,~$\wid X-\wid X_{\leq y}$ est réunion de composantes connexes de~$\wid X\setminus\{y\}$, et~$\wid X_{\leq y}$ est ainsi une partie fermée de~$\wid X$.  

\deux{lemextremchap} {\bf Lemme.} {\em Soit~$X$ un arbre, soient~$x$ et~$y$ deux points de~$X$ et soit~$z\in \got d X$ ; l'intersection~$[x;z]\cap [y;z]$ est de la forme~$[t;z]$ pour un certain~$t\in X$.}

\medskip
{\em Démonstration.} Comme~$[x;z]$ est compact,~$[x;z[$ est un fermé de~$X$ ; son intersection avec~$[y;x]$ est donc de la forme~$[y;t]$ pour un certain~$t$ ; la réunion de~$[y;t]$ et de~$[t;z]$ est un segment éventuellement long joignant~$y$ à~$z$, c'est donc~$[y;z]$ et l'on a alors~$[y;z]\cap [x;z]=[t;z]$.~$\Box$ 

\deux{lemmabout} {\bf Lemme.} {\em Soit~$X$ un arbre, soit~$Y$ un sous-arbre de~$X$ et soit~$\omega \in \got d X$. Les assertions suivantes sont équivalentes : 

\medskip
i) il existe~$x\in X$ tel que~$[x;\omega[\subset Y$ ; 

ii)~$Y$ aboutit à~$\omega$ ; 

iii)~$Y$ est non vide et~$[x;\omega[\subset Y$ pour tout~$x\in Y$.}

\medskip
{\em Démonstration}. Il est clair que i)$\Rightarrow$ii) et que iii)$\Rightarrow$i). Il reste à s'assurer que ii)$\Rightarrow$i). Supposons donc que~$Y$ aboutit à~$\omega$. Dans ce cas,~$Y\neq \emptyset$ puisque~$\got d Y$ contient~$\omega$. Soit~$x\in Y$ ; comme~$\got d Y$ contient~$\omega$, il existe un segment éventuellement long joignant~$x$ à~$\omega$ sur~$\wid Y$ ; ce dernier coïncide nécessairement avec l'unique segment éventuellement long joignant~$x$ à~$\omega$ sur~$X$, c'est-à-dire avec~$[x;\omega]$. Par conséquent,~$[x;\omega]\subset \wid Y$ ; il s'ensuit que~$[x;\omega[\subset Y$.~$\Box$

\deux{comparbor} {\bf Proposition.} {\em Soit~$X$ un arbre. Les quatre assertions suivantes sont équivalentes : 

\medskip

i) il existe un plus petit arbre compact contenant~$X$ comme sous-arbre ; 

ii) il existe un arbre compact contenant~$X$ comme sous-arbre ; 

iii)~$X$ est paracompact ; 

iv)~$\wid X$ est un arbre. 

\medskip
Si elles sont satisfaites,~$\wid X$ est le plus petit arbre compact contenant~$X$ comme sous-arbre.}

\medskip
{\em Démonstration.} L'implication i)~$\Rightarrow$ ii) est évidente. Supposons que ii) soit vraie, et soit~$Y$ un arbre compact dont~$X$ est un sous-arbre. Étant localement fermé dans~$Y$, l'arbre~$X$ est ouvert dans~$\overline X^Y$, qui est un arbre compact (\ref{adhcompinterv}) ; la proposition~\ref{compdeninf} assure alors que~$X$ est paracompact. 

Supposons maintenant que iii) soit vraie, et prouvons iv). Nous allons tout d'abord vérifier que les propriétés~$\alpha)$ et~$\beta)$ de~\ref{defarb} sont satisfaites par tout ouvert de la forme~$\widd U X$, où~$U$ est un ouvert connexe et à bord fini de~$X$. 

\medskip
Soit donc~$U$ un tel ouvert ; comme~$\partial_{\wid X}\widd U X$ est égal à~$\partial_X U$ (\ref{varpiuubarre}), il est fini et~$\widd U X$ satisfait~$\beta)$. 

Montrons maintenant qu'il satisfait~$\alpha)$. Soient~$x$ et~$y$ dans~$\widd U X$ ; le théorème~\ref{xchapeaupresquarb} assure l'existence d'un unique segment éventuellement long joignant~$x$ à~$y$ et tracé sur~$\widd U X$. 

Il suffit maintenant de s'assurer que~$[x;y]$ est un «vrai» segment, autrement dit que~$]x;y[$ est une «vraie» droite. Mais l'on a~$]x;y[=[x;y]\cap (X\setminus\{x,y\})$. Comme~$X$ est paracompact,~$X\setminus\{x,y\}$ est paracompact (\ref{paracompgen}) ; son fermé~$]x;y[$ est donc paracompact, et est de ce fait une vraie droite, ce qui achève de prouver que~$\widd U X$ satisfait~$\alpha)$. 

\medskip
L'espace~$\wid X$ est compact, et chacun des ses points a une base de voisinages satisfaisant~$\alpha)$ et~$\beta)$ ; de plus,~$\wid X$ lui-même satisfait~$\alpha)$ : il suffit d'appliquer ce qui précède avec~$U=X$. Par conséquent,~$\wid X$ est un arbre. 

\medskip
Enfin, supposons que iv) soit vraie, et soit~$Y$ un arbre compact contenant~$X$ comme sous-arbre. D'après le théorème~\ref{adhouvdensecpct},~$\overline X^Y$ s'identifie à~$\wid X$, qui est un arbre par hypothèse ; c'est donc le plus petit arbre compact contenant~$X$ comme sous-arbre.~$\Box$ 

\deux{corolladhsgr} {\bf Corollaire.} {\em Soit~$X$ un graphe et soit~$U$ un sous-arbre ouvert de~$X$ ; supposons que~$\got d U$ est fini. Les assertions suivantes sont équivalentes : 

\medskip
i)~$\overline U^X$ s'identifie à~$\wid U$ ; 

ii)~$\overline U^X$ est un arbre compact ; 

iii) le cardinal de~$\partial_X U$ est égal à celui de~$\got d U$ ; 

iv)  le cardinal de~$\partial_X U$ est supérieur ou égal à celui de~$\got d U$.}

\medskip
{\em Démonstration.} Supposons que i) soit vraie. On sait que~$\overline U^X$ est un sous-graphe compact de~$X$ (\ref{sousgraphead}), qui est connexe puisque~$U$ est connexe. Si~$x$ et~$y$ sont deux points de~$\overline U^X$ il existe donc un segment d'extrémités~$x$ et~$y$ tracé sur~$\overline U^X$ ; on déduit alors du théorème~\ref{xchapeaupresquarb} ci-dessus qu'un tel segment est unique ; par conséquent,~$\overline U^X$ est un arbre compact. 

Supposons que ii) soit vraie. Dans ce cas,~$U$ est un ouvert dense de l'arbre compact~$\overline U^X$ ; par conséquent,~$\overline U^X$ s'identifie à~$\wid U$ en vertu du théorème~\ref{adhouvdensecpct} ; en particulier,~$\partial_X U$ s'identifie à~$\got d U$, d'où iii).

Il est clair que iii)$\Rightarrow~$ iv). Si iv) est vraie, on déduit i) de la proposition~\ref{ensboutsfini}.~$\Box$

\deux{rembout} Soit~$X$ un graphe et soit~$U$ un sous-arbre ouvert de~$X$. Dans la suite, nous utiliserons abondamment, sans la rappeler explicitement, l'équivalence des propositions ci-dessous, qui résulte du corollaire~\ref{corolladhsgr} : 

a)~$U$ a exactement un bout et est relativement compact ;

b)~$U$ a exactement un bout et~$\partial_X U$ est non vide ; 

c)~$\partial_X U$ est un singleton et~$\overline U^X$ est un arbre compact. 

\medskip
Ainsi, on peut par exemple reformuler le~\ref{compgraphemoinss} comme suit : si~$X$ est un graphe connexe et si~$S$ est un sous-ensemble fini de~$X$, presque toute composante connexe de~$X\setminus S$ est un arbre à un bout relativement compact ; ou déduire de~\ref{bordpluribr} que si~$X$ est un graphe et~$x$ un point unibranche de~$X$, alors~$x$ possède un voisinage dans~$X$ qui est un arbre à un bout relativement compact. 

\medskip
\deux{lemdmbouts} {\bf Corollaire.} {\em Soit~$X$ un arbre et soit~$Y$ un sous-arbre fermé et non vide de~$X$. Les assertions suivantes sont équivalentes : 

\medskip
i)~$\got d Y =\got d X$ ; 

ii) toute composante connexe de~$X-Y$ est relativement compacte ; 

iii) toute composante connexe de~$X-Y$ est un arbre à un bout relativement compact.}

\medskip
{\em Démonstration.} Comme~$Y$ est un sous-arbre fermé et non vide de l'arbre~$X$, le bord de toute composante connexe de~$X-Y$ est un singleton ; par conséquent, ii)$\iff$iii).

\medskip
{\em Supposons que i) soit vraie}. Soit~$U$ une composante connexe de~$X-Y$ ; nous allons montrer que~$\overline U^X$ est compacte, ce qui revient à démontrer qu'elle est fermée dans~$\wid X$, c'est-à-dire encore que~$\partial_{\wid X}U$ ne rencontre pas~$\got d X$. On raisonne par l'absurde, en supposant que~$\partial_{\wid X}U$ contient un bout~$\omega$ de~$X$.

Choisissons~$x\in U$. Comme~$\omega\in  \partial_{\wid X}U$, l'on a~$[x;\omega[\subset U$. Par ailleurs, on déduit du lemme~\ref{lemmabout} et de i) qu'il existe~$y\in Y$ tel que~$[y;\omega[\subset Y$. L'intersection~$[x;\omega[\cap [y;\omega[$ est de la forme~$[t;\omega[$ pour un certain~$t\in X$ et est en particulier non vide ; mais elle est d'autre part contenue dans~$U\cap Y$ qui est vide, ce qui est contradictoire.

\medskip
{\em Supposons que ii) soit vraie.} Soit~$\omega\in \got d X$ ; nous allons montrer que~$\omega\in \got d Y$, en raisonnant par l'absurde ; on suppose donc que~$\omega\notin \got d Y$. Choisissons un point~$y$ sur l'arbre non vide~$Y$.  Comme~$Y$ n'aboutit pas à~$\omega$, il existe~$z\in [y;\omega[$ tel que~$z\notin Y$ ; comme~$y\in Y\cap [y;\omega[$ et comme~$Y\cap [y;\omega[$ est convexe,~$[z;\omega[$ ne rencontre pas~$Y$. Si~$U$ désigne la composante connexe de~$X-Y$ contenant~$[z;\omega[$ alors~$\omega\in \overline U^{\wid X}$ ; ceci entraîne que 
$\overline U^X$ n'est pas fermée dans~$\wid X$, et partant pas compacte, ce qui contredit ii). ~$\Box$

\section{Sous-graphes admissibles et squelettes}

\subsection*{Sous-graphes admissibles}

\deux{defarbadm} Soit~$X$ un graphe. On dira qu'un sous-graphe~$\Gamma$ de~$X$ est {\em admissible} s'il est fermé et si toute composante connexe de~$X-\Gamma$ est un arbre à un bout relativement compact. 

\trois{admrencontrecomp} Si~$\Gamma$ est un sous-arbre admissible de~$X$ alors toute composante connexe de~$X-\Gamma$ a un bord non vide ; il s'ensuit que~$\Gamma$ rencontre toutes les composantes connexes de~$X$ ; il est notamment non vide dès que~$X$ est non vide. 

\trois{admimplcon} Si~$\Gamma$ est un sous-arbre admissible de~$X$, il est convexe en vertu de~\ref{remskgrad}. 

\trois{compladmconv} Si~$\Gamma$ est un sous-arbre admissible de~$X$ et si~$V$ est une composante connexe de~$X-\Gamma$ alors~$\overline V$ est un sous-graphe convexe de~$X$. En effet, soient~$x$ et~$y$ deux points de~$V$, et soit~$I$ un segment les joignant.  Pour montrer que~$I\subset \overline V$, on raisonne par l'absurde ; on suppose donc que~$I$ rencontre~$X-\overline V$, et l'on choisit une composante connexe~$J$ de~$(X-\overline V)\cap I$. Comme~$J$ ne contient ni~$x$ ni~$y$, son bord dans~$I$ compte deux éléments ; mais ce bord est contenu dans~$\partial V$ qui est un singleton, d'où une contradiction. 

\deux{arbradmarbr} Soit~$X$ un arbre non vide. Il résulte du lemme~\ref{lemdmbouts}, de~\ref{admrencontrecomp} et de~\ref{admimplcon} qu'un sous-graphe fermé~$Y$ de~$X$ est admissible si et seulement si~$Y$ est un arbre non vide tel que~$\got d Y=\got d X$.

\trois{admzerobout} Supposons que l'arbre~$X$ soit compact, autrement dit que~$\got d X=\emptyset$. On déduit du~\ref{arbradmarbr} ci-dessus que si~$x\in X$ alors~$\{x\}$ est un sous-graphe admissible de~$X$. Comme tout sous-graphe admissible de~$X$ est non vide, il s'ensuit que les sous-graphes admissibles minimaux de~$X$ sont exactement les singletons ; si~$X$ lui-même n'est pas un singleton, il ne possède donc pas de plus petit sous-graphe admissible. 

\trois{admunbout} Supposons que l'arbre~$X$ ait exactement un bout, que l'on note~$\omega$.  On déduit du~\ref{arbradmarbr} ci-dessus que  pour tout~$x\in X$, le sous-arbre~$[x;\omega[$ de~$X$ est admissible. 

\medskip
Par ailleurs, si~$\Gamma$ un sous-graphe admissible de~$X$ contenant~$x$, il aboutit à~$\omega$, là encore en vertu de~\ref{arbradmarbr}, et contient dès lors~$[x;\omega[$ ; pour tout~$t$ différent de~$x$ et situé sur~$[x;y[$, l'intervalle~$[t;y[$ est encore un sous-graphe admissible de~$X$, qui est contenu {\em strictement} dans~$\Gamma$. 

\medskip
Il découle de ce qui précède, et du fait que tout sous-graphe admissible de~$X$ est non vide, que~$X$ ne possède pas de sous-graphe admissible minimal. 

\subsection*{Le squelette d'un graphe}

\deux{defexsquel} Soit~$X$ un graphe. Le sous-ensemble de~$X$ formé des points admettant un voisinage qui est un arbre ayant au plus un bout est un ouvert de~$X$ ; son fermé complémentaire est noté~$\skel X$ et est appelé le {\em squelette} de~$X$. Il résulte immédiatement de la définition : que si~$U$ est un ouvert (resp. une composante connexe) de~$X$ alors~$\skel X\cap U\subset \skel U$ (resp.~$\skel X\cap U=\skel U$) ; et que~$\mathsf S(X)$ est contenu dans tout sous-graphe admissible de~$X$. 

\deux{cerclesquel} Soit~$X$ un graphe. Si~$C$ est une boucle de~$X$, elle est contenue dans~$\skel X$. Pour le voir, on raisonne par l'absurde : on suppose qu'il existe~$x\in C$ et un voisinage ouvert~$U$ de~$x$ qui soit un arbre ayant au plus un bout. Comme~$U$ est un arbre,~$C$ n'est pas contenu dans~$U$. La composante connexe~$I$ de~$x$ dans~$C\cap U$ est alors un intervalle ouvert, dont chacune des deux extrémités appartient à~$\partial U$. Ce dernier ayant au plus un élément, les deux extrémités en question coïncident, ce qui signifie que~$\overline I=C$, et implique que~$C\subset \overline U$ ; mais~$\overline U$ est un arbre (compact), et l'on aboutit ainsi à une contradiction.

\deux{squelinterarb} Soit~$X$ un graphe ; il découle de~\ref{remskgrad} que~$\skel X$ est un sous-graphe fermé et convexe de~$X$. 

\deux{lemdrinclskel} {\bf Lemme.}{ \em Soit~$X$ un graphe, soit~$U$ un ouvert de~$X$ tel que~$\partial U\subset \skel X$, et soit~$I$ un fermé de~$U$ qui est une droite éventuellement longue. Le sous-graphe~$I$ de~$X$ est contenu dans~$\skel X$.}

\medskip
{\em Démonstration.} Soit~$x\in I$. Pour montrer que~$x\in \skel X$, on raisonne par l'absurde en supposant que ce n'est pas le cas. Il existe alors un voisinage ouvert~$V$ de~$x$ dans~$X$ qui est un arbre ayant au plus un bout. Soit~$J$ la composante connexe de~$x$ dans~$I\cap V$ ; c'est une droite éventuellement longue. L'adhérence de~$J$ dans~$\wid V$ s'identifie à~$\wid J$ ; par conséquent,~$\partial_{\wid V}J$ compte exactement deux éléments. Comme~$\wid V-V$ contient par hypothèse au plus un élément,~$\partial_V J$ est non vide ; soit~$y\in \partial_V J$. La composante~$J$ étant fermée dans~$I\cap V$, le point~$y$ ne peut appartenir à~$I$ ; il est donc situé sur~$\partial I\subset \partial U$, et par conséquent sur~$\skel X$ en vertu de notre hypothèse. Mais cela contredit le fait que~$y$ appartient à l'ouvert~$V$, qui est un arbre ayant au plus un bout.~$\Box$ 

\deux{corolldrinclskel} {\bf Corollaire} { \em Soit~$X$ un graphe et soit~$I$ un fermé de~$X$ qui est une droite éventuellement longue. Le sous-graphe~$I$ de~$X$ est contenu dans~$\skel X$.}

\medskip
{\em Démonstration.} C'est simplement l'énoncé du lemme~\ref{lemdrinclskel} ci-dessus, dans le cas particulier où~$U=X$.~$\Box$ 

\deux{squelarb} {\bf Proposition.} {\em Soit~$X$ un graphe et soit~$U$ un sous-arbre ouvert de~$X$ tel que~$\partial U\subset \skel X$. L'intersection~$\skel X\cap U$ est alors égale à la réunion des~$]x;y[$, où~$(x,y)$ parcourt~$(\got d U)^2$.}

\medskip
{\em Démonstration.} On procède par double inclusion. Posons~$\Gamma=\bigcup\limits_{(x,y)\in (\got d U)^2} ]x;y[$.

\medskip
{\em Montrons que~$\skel X\cap U\subset \Gamma$.} Nous allons en fait établir que~$\skel U \subset \Gamma$ ce qui suffira à conclure puisque~$\skel X\cap U \subset \skel U$.

Il résulte de~\ref{propuniondrl} que~$\Gamma$ est un sous-graphe fermé (et localement fini) de~$U$ ; comme~$\Gamma$ est convexe par contruction, c'est même un sous-arbre fermé de~$U$ dont l'ensemble des bouts, là encore par construction, coïncide avec~$\got d U$. On déduit alors de~\ref{arbradmarbr} que si~$\got d U$ est non vide,~$\Gamma$ est un sous-arbre admissible de~$U$, contenant de ce fait~$\mathsf S(U)$ ; et si~$\got d U=\emptyset$ alors~$U$
 est un arbre compact et l'on a~$\mathsf S(U)=\emptyset=\Gamma$.

\medskip
{\em Montrons que~$\Gamma\subset \skel X\cap U$.} Cela découle directement du lemme~\ref{lemdrinclskel} ci-dessus.~$\Box$ 

\deux{corollsquelarb} {\bf Corollaire.} {\em Soit~$X$ un arbre. Le squelette~$\skel X$ coïncide avec la réunion des~$]x;y[$, où~$(x,y)$ parcourt~$(\got d X)^2$.}

\medskip
{\em Démonstration.} C'est simplement l'énoncé de la proposition~\ref{squelarb} ci-dessus, dans le cas particulier où~$U=X$.~$\Box$

\deux{squelvide} {\bf Proposition.} {\em Soit~$X$ un graphe connexe. Son squelette est vide si et seulement si~$X$ est un arbre ayant au plus un bout.}

\medskip
{\em Démonstration.} Si~$X$ est un arbre ayant au plus un bout, son squelette est vide par définition. Réciproquement, supposons que le squelette de~$X$ soit vide. Comme toute boucle de~$X$ est contenue dans son squelette,~$X$ n'a pas de boucle et est un arbre ; on déduit alors du corollaire~\ref{corollsquelarb} que~$\got d X$ a au plus un élément.~$\Box$ 

\deux{compxmskl}{\bf Proposition.} {\em Soit~$X$ un graphe connexe dont le squelette est non vide. Toute composante connexe de~$X-\skel X$ est un arbre à un bout relativement compact. }

\medskip
{\em Démonstration.} Soit~$U$ une composante connexe de~$X-\skel X$ ; comme~$\skel X$ est non vide,~$\partial U$ est un sous-ensemble non vide de~$\skel X$.

Toute boucle de~$X$ étant contenue dans~$\skel X$, l'ouvert~$U$ ne contient aucune boucle et est donc un arbre. L'inclusion~$\partial U\subset \skel X$ implique, d'après la proposition~\ref{squelarb}, que~$\skel X\cap U$ est égale à la réunion des~$]x;y[$ pour~$(x,y)$ parcourant~$(\got d U)^2$. Mais~$S(X)\cap U=\emptyset$ par définition de~$U$ ; il s'ensuit que~$U$ a au plus un bout. Si~$\got d U$ était vide,~$U$ serait compact et son bord dans~$X$ serait vide, ce qui est absurde ; dès lors~$U$ est un arbre à un bout, relativement compact puisque~$\partial U$ est non vide. ~$\Box$ 

\deux{theosquel} {\bf Théorème.} {\em Soit~$X$ un graphe. Le graphe~$\skel X$ est localement fini, et la valence de~$(\skel X,x)$ est au moins égale à~$2$ pour tout~$x\in \skel X$.}

\medskip
{\em Démonstration.} Soit~$U$ un ouvert de~$X$ qui est un arbre à bord fini dont l'adhérence est un arbre compact. Le squelette~$\skel U$ s'identifie, d'après le corollaire~\ref{corollsquelarb}, à~$\bigcup ]x;y[$ pour~$(x,y)$ parcourant~$(\got d U)^2$. Comme~$\overline U$ est un arbre compact,~$\got d U\simeq \partial U$ et est en particulier fini ; il s'ensuit que~$\skel U$ est un arbre fini. L'arbre~$\skel X\cap U$ est contenu dans~$\skel U$, et est par conséquent un arbre fini ; ceci valant quel que soit~$U$, le squelette~$\skel X$ est un graphe localement fini. 

\medskip
Soit maintenant~$x\in \skel X$, et soit~$\Gamma$ un voisinage ouvert de~$x$ dans~$\skel X$ qui est un arbre à bord fini ; l'arbre~$\Gamma$ est une composante connexe de~$\skel X\setminus \partial  \Gamma$.

Soit~$U$ la composante connexe de~$X\setminus \partial  \Gamma$ contenant~$\Gamma$. On a~$U\cap \skel X=\Gamma$. En effet, supposons qu'il existe un point~$y$ de~$U\cap \skel X$ n'étant pas situé sur~$\Gamma$ ; par connexité de~$U$, il existerait un segment~$I$ d'extrémités~$x$ et~$y$ tracé sur~$U$, et qui ne rencontrerait dès lors pas~$\partial  \Gamma$. En vertu de la convexité de~$\skel X$, le segment~$I$ serait contenu dans~$\skel X$ ; mais il n'existe pas de segment tracé sur~$\mathsf S(X)$ et reliant~$x$ à un point de~$\skel X-\Gamma$ sans passer par~$\partial \Gamma$, d'où une contradiction.  

\medskip
Comme~$U\cap \skel X=\Gamma$, toute boucle de~$U$ est contenue dans~$\Gamma$ qui est un arbre ; il s'ensuit que~$U$ n'a pas de boucle et est donc un arbre. Son bord est par construction contenu dans~$\skel X$. La proposition~\ref{squelarb} assure alors que~$\skel X\cap U$ est égal à la réunion des~$]x;y[$ où~$(x,y)$ parcourt~$(\got d U)^2$ ; on en déduit que la valence de~$(\skel X\cap U,y)$ est au moins égale à 2 pour tout~$y\in \skel X\cap U$ ; en particulier, la valence de~$(\skel X,x)$ est au moins égale à~$2$, ce qui achève la démonstration.~$\Box$

\deux{propexistadm} {\bf Proposition.} {\em Soit~$X$ un graphe. Il existe un sous-graphe admissible et localement fini de~$X$.}

\medskip
{\em Démonstration.} En raisonnant composante connexe par composante connexe, on se ramène au cas où~$X$ est connexe et non vide. Si~$\skel X$ est non vide, il répond à la question ; sinon~$X$ est un arbre ayant au plus un bout, et l'on conclut à l'aide de~\ref{admzerobout} et~\ref{admunbout}.~$\Box$ 

\subsection*{La rétraction canonique sur un sous-graphe admissible}

\deux{defretcan} Soit~$X$ un graphe et soit~$\Gamma$ un sous-graphe admissible de~$X$. Soit~$r$ l'application de~$X$ vers~$\Gamma$ que l'on définit comme suit : 

\medskip
$\bullet$ si~$x\in \Gamma$ alors~$r(x)=x$ ; 

$\bullet$ si~$x\notin \Gamma$ alors~$r(x)$ est l'unique point du bord de la composante connexe de~$x$ dans~$X-\Gamma$ . 

\medskip
Par construction, l'application~$r$ est une rétraction de l'inclusion~$\Gamma\hookrightarrow X$ ; on l'appelle la {\em rétraction canonique} de~$X$ sur~$\Gamma$. 

\deux{introxrx} Soit~$X$ un graphe, soit~$\Gamma$ un sous-graphe admissible de~$X$ et soit~$r$ la rétraction canonique de~$X$ sur~$\Gamma$. 

\trois{xrx} Soit~$x\in X$. Il existe un et un seul segment~$I$ tracé sur~$X$ et joignant~$x$ à~$r(x)$. Si~$x\notin \Gamma$ alors~$I\setminus\{r(x)\}$ est contenu dans la composante connexe de~$x$ dans~$X-\Gamma$. 

\medskip
En effet, c'est évident si~$x=r(x)$. Sinon, soit~$U$ la composante connexe de~$X-\Gamma$ contenant~$x$.  Son adhérence~$\overline U$ est un arbre, et~$\partial U=\{r(x)\}$ ; par conséquent, il existe un unique segment~$I$ tracé sur~$\overline U$ et joignant~$x$ à~$r(x)$, et l'on a~$I\setminus\{r(x)\}\subset U$. D'autre part,~$\overline U$ est une partie convexe de~$X$ (\ref{compladmconv}) ; il s'ensuit que~$I$ est l'unique segment tracé {\em sur~$X$} et joignant~$x$ à~$r(x)$. 

\medskip
Nous nous permettrons de désigner cet unique segment joignant~$x$ à~$r(x)$ par~$[x;r(x)]$ même lorsque~$X$ n'est pas un arbre ; notons que~$[x;r(x)]\cap \Gamma=\{r(x)\}$. 

\trois{concatena} Soient~$x$ et~$y$ deux points de~$X$, et soit~$I$ un segment tracé sur~$X$ et joignant~$x$ à~$y$. Le graphe~$\Gamma$ étant convexe,~$J:=I\cap \Gamma$ est un intervalle. 

\medskip
S'il existe une composante connexe~$U$ de~$X-\Gamma$ qui contient~$x$ et~$y$ alors~$I$ est contenu dans~$U$, puisque ce dernier est un arbre ; par conséquent,~$J=\emptyset$. 

\medskip
Si ce n'est pas le cas alors~$J$ est non vide. Soit~$z$ (resp.~$t$) la borne de~$J$ située du côté de~$x$ (resp. de~$y$) ;  on a ou bien~$z=x$, ou bien~$[x;z[\subset X-\Gamma$, et de même avec~$y$ et~$t$. Il s'ensuit que~$z=r(x)$ et~$t=r(y)$ ; autrement dit,~$I=[x;r(x)]\cup J\cup[r(y);y]$, et~$J$ est {\em un} segment tracé sur~$\Gamma$ et joignant~$r(x)$ à~$r(y)$ ; on a~$J\cap [x;r(x)]=\{x\}$ et~$J\cap [r(y);y]=\{y\}$. 

\trois{convtestcomp} Si~$\Delta$ est une partie de~$X$ contenant~$\Gamma$ telle que~$\Delta\cap \overline V$ soit un sous-arbre compact de~$\overline V$ pour toute composante connexe~$V$ de~$X-\Gamma$ alors~$\Delta$ est un sous-graphe fermé convexe de~$X$. 

\medskip
Pour le voir, commençons par remarquer que les hypothèses faites sur~$\Delta$ entraînent que~$[x;r(x)]\subset \Delta$ pour tout~$x\in \Delta$. 

\medskip
{\em Le sous-ensemble~$\Delta$ de~$X$ est fermé.} Comme~$\Gamma\subset \Delta$, on peut écrire~$X-\Delta$ comme la réunion des~$\overline V-(\Delta\cap \overline V)$ pour~$V$ parcourant~$\pi_0(X-\Gamma)$. Si~$V\in \pi_0(X-\Gamma)$ le singleton~$\partial V$ est contenu dans~$\Gamma$, d'où l'égalité~$\overline V-(\Delta\cap \overline V)=V-(\Delta\cap \overline V)$ ; et ce dernier est, par compacité de~$\Delta\cap \overline V$, un ouvert de~$X$. Il s'ensuit que~$\Delta$ est un fermé de~$X$. 

\medskip
{\em Le sous-ensemble~$\Delta$ de~$X$ est convexe.} Soient~$x$ et~$y$ deux points de~$\Delta$ et soit~$I$ un segment les joignant. En vertu du~\ref{concatena} ci-dessus deux cas peuvent se présenter : 

\medskip
$\bullet$ ou bien~$x$ et~$y$ sont situés sur une même composante connexe~$U$ de~$X-\Gamma$, et~$I\subset U$ ; comme~$\Delta\cap \overline U$ est un sous-arbre de~$\overline U$, on a alors~$I\subset \Delta$ ; 

$\bullet$ ou bien~$I$ est de la forme~$[x;r(x)]\cup J\cup [y;r(y)]$ où~$J$ est un segment tracé sur~$\Gamma$ et joignant~$r(x)$ à~$r(y)$ ; étant alors réunion de trois intervalles contenus dans~$\Delta$, le segment~$I$ est lui-même contenu dans~$\Delta$. 

\medskip
Ainsi,~$\Delta$ est un sous-ensemble à la fois fermé et convexe de~$X$ ; c'en est donc un sous-graphe fermé et convexe.

\deux{thretcan} {\bf Théorème.} {\em Soit~$X$ un graphe, soit~$\Gamma$ un sous-graphe admissible de~$X$, soit~$U$ un ouvert connexe, non vide et à bord fini de~$\Gamma$, et soit~$\Delta$ un sous-graphe de~$\Gamma$. Soit~$r$ la rétraction canonique de~$X$ vers~$\Gamma$. 

\medskip
i ) L'image réciproque~$r\inv(U)$ est la composante connexe de~$X\setminus \partial _\Gamma U$ contenant~$U$. 

ii) L'application~$r:X\to \Gamma$ est continue et compacte. 

iii) L'image réciproque~$r\inv(\Delta)$ est un sous-graphe de~$X$.

iv) Le sous-graphe~$\Delta$ de~$r\inv(\Delta)$ est admissible, et la rétraction canonique de~$r\inv(\Delta)$ sur~$\Delta$ est la restriction de~$r$. 

v) L'on a~$\partial_X(r\inv(\Delta))=\partial_\Gamma\Delta$, et~$r\inv(\Delta)$ est relativement compact si et seulement si~$\Delta$ est relativement compact.

vi) Le sous-graphe~$\overline \Delta$ de~$\overline{r\inv(\Delta)}$ est admissible, et la rétraction canonique de~$\overline{r\inv(\Delta)}$ sur~$\overline \Delta$ est la restriction de~$r$.}

\medskip
{\em Démonstration.} On prouve chacune des assertions séparément.

\trois{rmnsunu} {\em Preuve de i).} On appelle~$V$ la composante connexe de~$X\setminus \partial _\Gamma U$ qui contient~$U$ et l'on procède par double inclusion. 

\medskip
Soit~$x\in r\inv(U)$. Si~$x\in \Gamma$ alors~$r(x)=x$ et l'on a donc~$x\in U\subset V$. Sinon, soit~$W$ la composante connexe de~$X-\Gamma$ contenant~$x$ ; c'est un arbre ouvert relativement compact de~$X$ de bord~$\{r(x)\}$. Comme~$r(x)$ adhère à~$W$ et comme~$V$ est un voisinage de~$r(x)$,  l'ouvert~$V\cap W$ de~$W$ est non vide ; le bord de~$V$ est contenu dans~$\partial_\Gamma U$, donc dans~$\Gamma$, et ne rencontre de ce fait pas~$W$ ; par conséquent,~$V\cap W$ est également fermé dans~$W$, et partant égal à~$W$ tout entier. Autrement dit ,~$W\subset V$ ; en particulier,~$x\in V$. 

\medskip
Soit~$x\in V$. Si~$x\in U$ alors~$x\in r\inv(U)$. Sinon, soit~$W$ la composante connexe de~$V-U$ contenant~$x$ ; son bord dans~$V$ est une partie non vide de~$U$, au sein de laquelle on choisit un point~$y$. 

{\em La composante~$W$ ne rencontre pas~$\Gamma$}. En effet s'il existait un point~$z$ de~$\Gamma$ situé sur~$W$, l'on pourrait, en vertu de la connexité de~$V$, tracer sur~$V$ un segment d'extrémités~$y$ et~$z$ ; ce segment serait contenu dans~$\Gamma$ par convexité de ce dernier (\ref{remskgrad}), et relierait le point~$y$ de~$U$ au point~$z$ de~$\Gamma-U$ sans passer par~$\partial_\Gamma U$ (puisque~$V\cap \partial_\Gamma U=\emptyset$), ce qui est absurde. 

La composante connexe~$W$ de~$V-U$ ne rencontre pas~$\Gamma$, et son bord est contenu dans~$(U\cup\partial V)\subset\Gamma$ ; par conséquent,~$W$ est une composante connexe de~$X-\Gamma$ ; on a alors nécessairement~$y=r(x)$, et donc~$x\in r\inv(U)$. 

\trois{rcompact} {\em Preuve de ii.} La continuité découle de i), et du fait que les ouverts 
connexes, non vides et à bord fini forment une base de la topologie de~$\Gamma$. 

Soit~$\sch K$ une partie compacte de~$\Gamma$. Toute composante connexe de~$X-\Gamma$ dont l'unique point du bord appartient à~$\sch K$ est une composante connexe relativement compacte de~$X\setminus \sch K$ ; on déduit alors de la définition de~$r$ et de~\ref{compxmscomp} que~$r\inv(\sch K)$ est compact.

\trois{rmoinsdelta}{\em Preuve de iii).} Comme~$r$ est continue,~$r\inv(\Delta)$ est une partie localement fermée de~$X$. Nous allons montrer que c'est un graphe. 

\medskip
{\em Étape technique intermédiaire.} On suppose provisoirement que~$\Delta$ est un arbre, et l'on se donne deux points ~$x$ et~$y$ sur~$r\inv(\Delta)$  ; nous allons vérifier qu'il existe un unique segment d'extrémités~$x$ et~$y$ tracé sur~$r\inv(\Delta)$. 

\medskip
$\bullet$ {\em Existence}. Si~$x$ et~$y$ sont situés sur une même composante connexe~$V$ de~$X-\Gamma$, on peut prendre le segment joignant~$x$ à~$y$ sur l'arbre~$V$ ; sinon, on peut prendre la réunion~$[x;r(x)]\cup I\cup [r(y):y]$, où~$I$ est le segment joignant~$r(x)$ à~$r(y)$ sur l'arbre~$\Delta$.

$\bullet$ {\em Unicité.} Si~$x$ et~$y$ sont situés sur une même composante connexe~$V$ de~$X-\Gamma$, la convexité de~$\overline V$ (\ref{compladmconv}) assure que le seul segment joignant~$x$ à~$y$ sur~$X$ est celui qui les relie sur~$\overline V$ (et qui est en fait contenu dans~$V$) ; sinon, tout segment joignant~$x$ à~$y$ sur~$r\inv(\Delta)$ est, en vertu de~\ref{concatena}, de la forme~$[x;r(x)]\cup J\cup [r(y); y]$, où~$J$ est un segment joignant~$r(x)$ à~$r(y)$ tracé sur~$\Gamma\cap r\inv(\Delta)=\Delta$ ; l'unicité souhaitée découle alors du fait qu'il y a, sur l'arbre~$\Delta$, un unique segment joignant~$r(x)$ à~$r(y)$.

\medskip
{\em Fin de la preuve de iii).} On ne suppose plus que~$\Delta$ est un arbre. Il suffit maintenant de montrer que tout point de~$r\inv(\Delta)$ a un voisinage ouvert dans~$r\inv(\Delta)$ qui est un arbre. Soit donc~$x\in r\inv(\Delta)$, et soit~$\Gamma_0$ un voisinage ouvert de~$r(x)$ dans~$\Gamma$ qui est un arbre tel que~$\Gamma_0\cap \Delta$ soit un sous-graphe fermé de~$\Gamma_0$. Soit~$\Delta_0$ la composante connexe de~$r(x)$ dans~$\Delta\cap \Gamma_0$ ; c'est un ouvert de~$\Delta$, et un sous-arbre fermé de~$\Gamma_0$. D'après l'étape technique ci-dessus, l'ouvert~$r\inv(\Gamma_0)$ de~$X$ est un arbre, et son fermé~$r\inv(\Delta_0)$ en est une partie convexe, et est donc elle-même un arbre. Or~$r\inv(\Delta_0)$ est un ouvert de~$r\inv(\Delta)$ qui contient~$x$, ce qui achève de prouver iii).

\trois{presqueconclurmoinsdelta} {\em Preuve de iv).} Par la définition même de~$r$, toute composante connexe de~$r\inv(\Delta)-\Delta$ est une composante connexe de~$X-\Gamma$, et est en particulier un arbre à un bout relativement compact ; comme on a par ailleurs, là encore par définition de~$r$, l'égalité~$\Delta=\Gamma\cap r\inv(\Delta)$, le graphe~$\Delta$ est fermé dans~$r\inv(\Delta)$ et est donc un sous-graphe admissible de ce dernier. Quant à la rétraction canonique de~$r\inv(\Delta)$ sur~$\Delta$, il résulte immédiatement de sa définition qu'elle est égale à~$r_{|r\inv(\Delta)}$.

\trois{conclurmoinsdelta} {\em Preuve de v).} Montrons tout d'abord l'égalité~$\partial_X r\inv(\Delta)=\partial_\Gamma \Delta$ ; on procède par double inclusion. 

\medskip
Soit ~$x\in \partial_\Gamma\Delta$. Si~$U$ est un voisinage ouvert de~$x$ dans~$X$ alors~$U$ rencontre~$\Delta\subset r\inv(\Delta)$, et~$U$ rencontre aussi~$\Gamma-\Delta\subset r\inv(\Gamma-\Delta)=X-r\inv(\Delta)$. Par conséquent,~$x\in \partial_X r\inv(\Delta)$. 

\medskip
Réciproquement, soit~$x\in \partial_X r\inv(\Delta)$ ; on a alors~$x\in r\inv(\overline \Delta^X)$. Montrons tout d'abord que~$x\in \overline \Delta^X$. Si ce n'était pas le cas, la composante connexe de~$x$ dans~$X-\Gamma$ serait un ouvert de~$X$ contenu dans~$r\inv(r(x))$ et serait donc ou bien contenue dans~$r\inv(\Delta)$, ou bien contenue dans son complémentaire, contredisant ainsi l'hypothèse que~$x\in \partial_X r\inv(\Delta)$. 

Montrons maintenant que~$x\in \partial_\Gamma \Delta$. Si ce n'était pas le cas, il existerait un ouvert~$\Gamma_0$ de~$\Gamma$ contenu dans~$\Delta$ et contenant~$x$, et~$r\inv(\Gamma_0)$ serait alors un voisinage ouvert de~$x$ inclus dans~$r\inv(\Delta)$, contredisant là encore l'appartenance de~$x$ à~$\partial_X r\inv(\Delta)$.

\medskip
On donc bien démontré que~$\partial_X r\inv(\Delta)=\partial_\Gamma \Delta$. Il reste à s'assurer que~$r\inv(\Delta)$ est relativement compact si et seulement si il en va de même pour~$\Delta$ ; mais c'est une conséquence formelle de la compacité et de la surjectivité de~$r$.

\trois{adhretradm} {\em Preuve de vi.} En vertu de v), on a~$\overline{r\inv(\Delta)}=r\inv(\Delta)\cup\partial \Delta$ ; par conséquent, toute composante connexe de~$\overline{r\inv(\Delta)}-\overline \Delta$ est une composante connexe de~$r\inv(\Delta)-\Delta$, et vi) est alors une conséquence immédiate de iv).~$\Box$

\deux{propretcan} Soit~$X$ un graphe, soit~$\Gamma$ un sous-graphe admissible de~$X$ et soit~$r$ la rétraction canonique de~$X$ sur~$\Gamma$. 

\trois{pizeroiso} La continuité de~$r$ établie par le théorème~\ref{thretcan} ci-dessus assure que l'application naturelle~$\pi_0(\Gamma)\to \pi_0(X)$ est injective. Mais~$\Gamma$ rencontre par ailleurs toutes les composantes connexes de~$X$ ; il en découle que~$\pi_0(\Gamma)\simeq \pi_0(X)$.

\trois{gammxarbeq} Si~$X$ est connexe (cela équivaut par ce qui précède à la connexité de~$\Gamma$) alors~$X$ est un arbre si et seulement si~$\Gamma$ est un arbre : cela résulte du fait que~$\Gamma$ contient~$\skel X$, qui lui-même contient toutes les boucles de~$X$.

\trois{boutsxboutsgamma} Si~$\sch K$ est une partie compacte de~$X$ alors~$r(\sch K)$ est une partie compacte de~$\Gamma$, et~$r\inv(r(\sch K))$ est, par compacité de~$r$, un compact de~$X$ contenant~$\sch K$. L'ensemble des parties de~$X$ de la forme~$r\inv(\Delta)$, où~$\Delta$ est un compact de~$\Gamma$, est donc cofinal dans l'ensemble des compacts de~$X$. On en déduit, en se fondant sur~\ref{pizeroiso}, que les espaces topologiques~$\got b(\Gamma)$ et~$\got b(X)$ sont naturellement homéomorphes. Dans le cas où~$X$ et~$\Gamma$ sont des arbres, cela signifie plus précisément que l'inclusion naturelle~$\wid \Gamma\subset \wid X$ identifie~$\got d\Gamma$ et~$\got d X$, ce que l'on savait déjà (\ref{arbradmarbr}). 

\deux{lemmcommret} {\bf Lemme.} {\em Soit~$f:Y\to X$ une application continue entre graphes, et soit~$\Gamma$ un sous-graphe admissible de~$X$ tel que~$f\inv(\Gamma)$ soit un sous-graphe admissible de~$Y$ ; soient~$r$ et~$\rho$ les rétractions canoniques respectives de~$X$ sur~$\Gamma$ et de~$Y$ sur~$f\inv(\Gamma)$. On a l'égalité~$f\circ \rho=r\circ f$.}

\medskip
{\em Démonstration.} Soit ~$y\in Y$. Si~$y\in f\inv(\Gamma)$ alors~$f(y)\in \Gamma$ et l'on a donc~$f(\rho(y))=f(y)=r(f(y))$. Supposons maintenant que~$y\notin f\inv(\Gamma)$, soit~$V$ la composante connexe de~$y$ dans~$Y-f\inv(\Gamma)$ et soit~$U$ la composante connexe de~$X-\Gamma$ contenant~$f(V)$. Le point~$\rho(y)$ est par définition l'unique point de~$\partial V$ ; son image~$f(\rho(y))$ est donc un point de~$\Gamma$ adhérant à~$U$ ; c'est de ce fait l'unique point de~$\partial U$, et il vient ~$r(f(y))=f(\rho(y))$.~$\Box$ 

\deux{extendadm} {\bf Proposition.} {\em Soit~$X$ un graphe, soit~$\Gamma$ un sous-graphe admissible de~$X$, et soit~$\Delta$ un sous-graphe fermé de~$X$ contenant~$\Gamma$. Les propositions suivantes sont équivalentes : 

\medskip
i)~$\Delta$ est admissible ; 

ii)~$\pi_0(\Gamma)\to \pi_0(\Delta)$ est une bijection ; 

iii) pour toute composante connexe~$V$ de~$X-\Gamma$, l'intersection~$\Delta\cap \overline V$ est un arbre ; 

iv)~$\Delta$ est convexe.}

\medskip
{\em Démonstration.} Si i) est vraie alors ii) est vraie en vertu de~\ref{pizeroiso}. 

\medskip
Supposons que ii) est vraie, et soit~$V$ une composante connexe de~$X-\Gamma$ ; soit~$x$ l'unique point de~$\partial V$. L'intersection~$\Delta\cap \overline V$ est un sous-graphe compact de~$\overline V$ contenant~$x$. Si ce graphe n'était pas un arbre, il possèderait une composante connexe~$\Delta'$ ne contenant pas~$x$. On aurait donc~$\Delta'\subset V$, et~$\Delta'$ serait par conséquent ouverte dans~$\Delta\cap V$, et partant dans~$\Delta$ ; étant par ailleurs compacte, connexe et non vide,~$\Delta'$ serait une composante connexe de~$\Delta$ ne rencontrant pas~$\Gamma$, ce qui contredirait ii) ; par conséquent, iii) est vraie. 

\medskip
Si iii) est  vraie, la propriété iv) découle de~\ref{convtestcomp}. 

\medskip
Supposons que iv) soit vraie et soit~$W$ une composante connexe de~$X-\Delta$. Nous allons montrer que~$W$ est un arbre à un bout relativement compact, ce qui achèvera la démonstration. La composante connexe~$W'$ de~$X-\Gamma$ qui contient~$W$ est un arbre à un bout relativement compact ; il s'ensuit que~$W$ est un arbre relativement compact et à bord non vide ; la convexité de~$\Delta$ garantit alors, en vertu de~\ref{bordcomplconvgr}, que~$\partial W$ est un singleton.~$\Box$ 

\deux{pluspetitadm} {\bf Proposition.} {\em Soit~$X$ un graphe, soit~$\Gamma$ un sous-graphe admissible de~$X$ et soit~$\Delta$ un sous-graphe fermé de~$X$. Pour toute composante connexe~$V$ de~$X-\Gamma$, désignons par~$T(V)$ l'enveloppe convexe de~$(\Delta \cap\overline V)\cup \partial V$ dans l'arbre~$\overline V$ ; la réunion~$\Theta:=\Gamma\cup \bigcup\limits_{V\in \pi_0(X-\Gamma)}T(V)$ est alors l'enveloppe convexe de~$\Gamma\cup \Delta$ ; c'est aussi le plus petit sous-graphe admissible de~$X$ contenant~$\Gamma\cup \Delta$, et~$\Theta$ est localement fini dès que~$\Gamma$ et~$\Delta$ le sont.}

\medskip
{\em Démonstration.} On a par construction~$\Theta\supset \Delta\cup \Gamma$. Il découle de~\ref{convtestcomp} que~$\Theta$ est un sous-graphe fermé et convexe de~$X$ ; comme il est par ailleurs clair que toute partie convexe de~$X$ contenant~$\Gamma$ et~$\Delta$ contient~$\Theta$, celui-ci est bien l'enveloppe convexe de~$\Gamma \cup \Delta$.  On déduit alors de la proposition~\ref{extendadm} que~$\Theta$ est le plus petit sous-graphe admissible de~$X$ contenant~$\Gamma\cup \Delta$.

\medskip
Supposons maintenant que~$\Gamma$ et~$\Delta$ sont localement finis. Dans ce cas, si~$V\in \pi_0(X-\Gamma)$ alors~$\Delta\cap \overline  V$ est un graphe compact et fini, d'où il découle que~$T(V)$ est lui-même compact et fini : en effet, si l'on note~$y$ l'unique point de~$\partial V$ et si l'on choisit un sous-ensemble fini~$S$ de~$\Delta\cap \overline V$ rencontrant chacune de ses composantes connexes, l'on a~$T(V)=\bigcup\limits_{z\in S}[z;y]$. 

\medskip
Soit maintenant~$x\in \Theta$ et soit~$\Gamma_0$ un voisinage de~$r(x)$ dans~$\Gamma$ qui est un sous-graphe compact et fini de ce dernier. L'intersection de~$\Delta$ et du compact~$ r\inv(\Gamma_0)$ est un graphe compact et fini ; par conséquent,~$(\Delta\cap r\inv(\Gamma_0))-\Gamma_0$ a un nombre fini de composantes connexes ; il s'ensuit que l'ensemble~$\Pi$ des composantes connexes de~$r\inv(\Gamma_0)-\Gamma_0$ qui rencontrent~$\Delta$ est fini.

On a par construction~$\Theta\cap r\inv(\Gamma_0)=\Gamma_0\cup \bigcup\limits_{V\in \Pi}T(V)$. Comme~$\Pi$ est fini, et comme~$T(V)$ est pour tout~$V\in \pi_0(X-\Gamma)$ un arbre compact et fini ({\em cf. supra}), l'intersection de~$\Theta$ et de~$ r\inv(\Gamma_0)$ est un graphe fini. Ainsi,~$x$ possède un voisinage dans~$\Theta$ qui est un graphe fini ; par conséquent,~$\Theta$ est localement fini.~$\Box$ 

\deux{lemmadmsquel} {\bf Lemme.} {\em Soit~$X$ un graphe et soit~$\Gamma$ un sous-arbre admissible de~$X$. Les propositions suivantes sont équivalentes : 

\medskip

i)~$\Gamma=\skel X$ ; 

ii)~$\Gamma$ ne possède aucun point isolé ou unibranche.}

\medskip
{\em Démonstration.} L'implication i)$\Rightarrow$ ii) est une partie du théorème~\ref{theosquel}. Prouvons maintenant que ii)$\Rightarrow i)$ ; nous allons en réalité établir l'implication contraposée ; on suppose donc que~$\Gamma$ contient strictement~$\skel X$ (ce qui entraîne que~$X\neq \emptyset$) et l'on va établir qu'il possède au moins un point isolé ou unibranche. 

\medskip
En raisonnant composante par composante, on se ramène au cas où~$X$ est connexe. 

Si~$\skel X$ est non vide c'est un sous-arbre admissible de~$X$ strictement contenu dans~$\Gamma$. 

Si~$\skel X$ est vide alors~$X$ est ou bien un arbre à un bout, auquel cas~$\Gamma$ contient strictement un sous-arbre admissible de~$X$ (\ref{admunbout}), ou bien un arbre compact et non vide ; dans ce dernier cas,~$\Gamma$ peut ou bien être un singleton, auquel cas il consiste en un point isolé et la démonstration est terminée, ou bien compter au moins deux points, et le choix d'un point quelconque de~$\Gamma$ définit alors un sous-arbre admissible de~$X$ (\ref{admzerobout}), qui est strictement contenu dans~$\Gamma$. 

\medskip
On peut donc faire l'hypothèse que~$\Gamma$ contient strictement un sous-arbre admissible~$\Gamma'$ de~$X$. Il existe dès lors une composante connexe~$V$ de~$X-\Gamma'$ et un point~$x$ sur~$V\cap \Gamma$ ; le bord de~$V$ est un singleton~$\{y\}$ pour un certain~$y\in \Gamma'$. 

En vertu de la proposition~\ref{extendadm}, l'intersection~$\Gamma\cap \overline V$ est un arbre compact, qui n'est pas réduit à~$\{y\}$ puisqu'il contient~$x$. D'après~\ref{infexist}, l'arbre~$\Gamma\cap \overline V$ possède un point~$z$ distinct de~$y$ et unibranche. Comme~$z$ appartient à l'ouvert~$V\cap\Gamma$ de~$\Gamma$, c'est un point unibranche de~$\Gamma$.~$\Box$

\deux{propwadm} {\bf Théorème.} {\em Soit~$X$ un graphe et soit~$\Delta$ un sous-graphe de~$X$. Les assertions suivantes sont équivalentes :

\medskip
i) il existe un sous-graphe admissible~$\Gamma$ de~$X$ dont~$\Delta$ est un ouvert ; 

ii) il existe un ouvert~$U$ de~$X$ dont~$\Delta$ est un sous-graphe admissible. 

\medskip
De plus, si elles sont satisfaites : 

\medskip

- le graphe~$\Gamma$ de i) peut être choisi localement fini si~$\overline \Delta$ est localement fini ; 

- l'ouvert~$U$ de ii) est unique ; 

- si~$\Gamma$ est comme dans i) et si~$r$ désigne la rétraction canonique de~$X$ sur~$\Gamma$ alors l'unique ouvert~$U$ de ii) est égal à~$r\inv(\Delta)$ (qui ne dépend donc pas du choix de~$\Gamma$).}

\medskip
{\em Démonstration.} On procède par double implication, en établissant en cours de preuve les assertions supplémentaires énoncées en fin de théorème.

\trois{rempreal}{\em Remarques préliminaires.} Soit~$U$ un ouvert dont~$\Delta$ est un sous-graphe admissible et soit~$V$ une composante connexe de~$U-\Delta$ ; c'est un arbre à un bout relativement compact dans~$U$, dont l'unique point~$x$ du bord est situé sur~$\Delta$ ; c'est {\em a fortiori} un arbre à un bout relativement compact dans~$X$ de bord~$\{x\}$. 

Soit~$\Gamma$ un sous-graphe admissible de~$X$ dont~$\Delta$ est un ouvert. On a alors~$V\cap \Gamma=\emptyset$ : en effet, s'il existait~$y\in \Gamma\cap V$ le segment~$[y;x]$ de l'arbre~$\overline V$ serait contenu dans~$\Gamma$ (\ref{admimplcon}) ; mais alors tout voisinage de~$x$ dans~$\Gamma$ rencontrerait~$[y;x[$ et donc~$\Gamma-\Delta$, contredisant le fait que~$\Delta$ est ouvert dans~$\Gamma$. 

Par conséquent,~$V$ est une composante connexe de~$X-\Gamma$.

\trois{ouvgamimplu} {\em On suppose que i) est vraie.} Soit~$r$ la rétraction canonique de~$X$ sur~$\Gamma$. On sait que~$r\inv(\Delta)$ est un ouvert de~$X$ dont~$\Delta$ est un sous-graphe admissible. Soit maintenant~$U$ un ouvert de~$X$ dont~$\Delta$ est un sous-graphe admissible ; nous allons montrer que~$U=r\inv(\Delta)$.

\medskip
{\em Montrons que~$U\subset r\inv(\Delta)$.} Soit~$V$ une composante connexe de~$U-\Delta$. C'est en vertu de~\ref{rempreal} une composante connexe de~$X-\Gamma$ dont l'unique point~$x$ du bord est situé sur~$\Delta$ ; par définition de~$r$, on a~$r(y)=x\in \Delta$ pour tout~$y\in V$. Il s'ensuit que~$r(U-\Delta)\subset \Delta$ ; comme par ailleurs~$r(\Delta)\subset \Delta$, on a~$r(U)\subset \Delta$, ce que nous voulions prouver.

\medskip
{\em Montrons que~$r\inv(\Delta)\subset U$.} Soit~$x\in r\inv(\Delta)$. Si~$x\in \Gamma$, on a~$r(x)=x$ et~$x$ appartient donc à~$\Delta$, et en particulier à~$U$. Sinon, soit~$V$ la composante connexe de~$X-\Gamma$ contenant~$x$ ; son bord est~$\{r(x)\}$. Comme~$U$ est un voisinage de~$r(x)$, il rencontre~$V$. Soit~$W$ une composante connexe de~$U-\Delta$ telle que~$W\cap V\neq \emptyset$. D'après~\ref{rempreal},~$W$ est une composante connexe de~$X-\Gamma$ ; l'on a donc~$V=W$ et~$V\subset U$ ; en particulier,~$x\in U$.

\trois{uimlouvgam} {\em On suppose que ii) est vraie.} Choisissons un sous-graphe admissible localement fini~$\Theta$ de~$X$,  et soit~$\Xi$ le plus petit sous-graphe admissible de~$X$ contenant~$\overline \Delta$ et~$\Theta$ (prop.~\ref{pluspetitadm}) ; si~$\overline \Delta$ est localement fini,~$\Xi$ est localement fini. 

\medskip
Soit~$\Gamma$ le fermé~$\Xi-(U-\Delta)$ de~$\Xi$ ; on a~$\Gamma\cap U=\Delta$, et~$\Delta$ est donc un ouvert de~$\Gamma$. Nous allons montrer que~$\Gamma$ est un sous-graphe admissible de~$X$ ; en tant que sous-graphe fermé de~$\Xi$, il sera alors automatiquement localement fini si~$\Xi$ est localement fini, et donc si~$\overline \Delta$ est localement fini. 

\medskip
Soit~$V$ une composante connexe de~$X-\Gamma$ ; le but de ce qui suit est de prouver que~$V$ est un arbre à un bout relativement compact, ce qui permettra de conclure. On distingue deux cas. 

\medskip
{\em Supposons que~$V\cap U=\emptyset$.} Dans ce cas~$V$ ne rencontre pas~$\Xi$ (puisque~$\Xi\cap (X-U)=\Gamma\cap (X-U)$), et son bord est contenu dans~$\Gamma\subset \Xi$ ; par conséquent,~$V$ est une composante connexe de~$X-\Xi$, et est donc un arbre à un bout relativement compact. 

\medskip
{\em Supposons que~$V$ rencontre~$U$.} Soit~$x\in V\cap U$. Comme~$\Gamma\cap U=\Delta$, le point~$x$ n'est pas situé sur~$\Delta$, et la composante connexe~$U'$ de~$x$ dans~$U-\Delta$ ne rencontre pas~$\Gamma$ ; elle est donc contenue dans~$V$. L'ouvert~$U'$ de~$X$ est un arbre à un bout relativement compact dont l'unique point du bord est situé sur~$\Delta\subset \Gamma$ (\ref{rempreal}).Par conséquent,~$U'$ est fermé dans~$V$ ; il s'ensuit que~$V$ est égal à~$U'$, et est donc bien un arbre à un bout relativement compact.~$\Box$ 

\deux{deltabem} Soit~$X$ un graphe. On dira d'un sous-graphe~$\Delta$ de~$X$ qui satisfait les conditions équivalentes du théorème~\ref{propwadm} ci-dessus qu'il est {\em faiblement admissible} ; si c'est le cas, l'unique ouvert de~$X$ dont~$\Delta$ est un sous-graphe admissible sera noté~$\Delta^\flat$. 

\trois{wadmadm} Si~$\Delta$ est un sous-graphe admissible de~$X$ alors~$\Delta$ est faiblement admissible et~$\Delta^\flat=X$. 

\trois{wadmouv} Si~$U$ est un ouvert de~$X$ et si~$\Delta$ est un sous-graphe faiblement admissible de~$U$, il existe un ouvert~$V$ de~$U$ dont~$\Delta$ est un sous-graphe admissible ; par conséquent,~$\Delta$ est un sous-graphe faiblement admissible de~$X$, et l'on~$\Delta^\flat=V$ indépendamment du fait que l'on voie~$\Delta$ comme un sous-graphe faiblement admissible de~$U$ ou de~$X$ ; nous utiliserons ce fait implicitement à plusieurs reprises. 

\trois{wadmbord} Soit~$\Delta$ un sous-graphe faiblement admissible de~$X$. On déduit alors des assertions v) et vi) du théorème~\ref{thretcan} les faits suivants :~$\partial_X\Delta^\flat=\overline \Delta-\Delta$, et~$\Delta^\flat$ est relativement compact si et seulement~$\overline \Delta$ est compact ; de plus,~$\overline \Delta$ est un sous-graphe admissible de~$\overline{\Delta^\flat}$.

\subsection*{Dimension topologique d'un graphe paracompact}

\deux{propdimtop} {\bf Proposition.} {\em Soit~$X$ un graphe paracompact non vide. La dimension topologique de~$X$ vaut~$0$ si~$X$ est discret, et 1 sinon.}

\medskip
{\em Démonstration.} Si le graphe~$X$ est discret, sa dimension topologique est nulle. Réciproquement, si sa dimension topologique est nulle, tout recouvrement ouvert de~$X$ peut être raffiné en un recouvrement constitué d'ouverts deux à deux disjoints ; cela interdit à deux points de~$X$ de se trouver sur la même composante connexe de~$X$, et ce dernier est dès lors discret. 

\medskip
Il suffit donc maintenant de démontrer (sans plus s'inquiéter de savoir si~$X$ est discret ou non) que la dimension topologique de~$X$ est majorée par 1. Soit~$(U_i)$ un recouvrement ouvert de~$X$ ; nous allons montrer qu'il peut être raffiné en un recouvrement localement fini par des ouverts dont les intersections {\em trois à trois} sont vides, ce qui permettra de conclure. 

\medskip
Comme~$X$ est paracompact, et comme il possède une base d'ouverts connexes à bord fini, on peut raffiner~$(U_i)$ en un recouvrement localement fini~$(V_j)$ où les~$V_j$ sont des ouverts connexes, non vides et à bord fini. Pour tout~$j$, choisissons un point~$x_j$ sur~$V_j$. Comme~$(V_j)$ est localement fini, la réunion des~$\partial V_j$ et des~$\{x_j\}$ est une partie fermée et discrète de~$X$ ; celle-ci est contenue dans un sous-graphe admissible et localement fini~$\Gamma$ de~$X$, et l'on note~$r$ la rétraction canonique de~$X$ sur~$\Gamma$. 

\medskip
Fixons~$j$. L'ouvert~$V_j$ est une composante connexe de~$X\setminus \partial  V_j$. Si l'on désigne par~$\Delta_j$ la composante connexe de~$\Gamma\setminus \partial  V_j$ qui contient~$x_j$ alors~$r\inv(\Delta_j)$ est une composante connexe de~$X\setminus \partial  V_j$ contenant~$x_j$ ; par conséquent,~$r\inv(\Delta_j)=V_j$. Remarquons qu'on a alors~$r(V_j)=\Delta_j$, ce qui entraîne que les~$\Delta_j$ recouvrent~$\Gamma$. 

\medskip
En tant que fermé de~$X$, le graphe~$\Gamma$ est paracompact. Comme il est localement fini, il est localement de dimension topologique majorée par 1 ; il s'ensuit que la dimension topologique de~$\Gamma$ est majorée par~$1$. 

\medskip
On peut de ce fait raffiner~$(\Delta_j)$ en un recouvrement localement fini~$(\Theta_\ell)$ de~$\Gamma$ constitué d'ouverts dont les intersections trois à trois sont vides. La famille~$(r\inv(\Theta_\ell))$ constitue alors un recouvrement localement fini de~$X$ par des ouverts dont les intersections trois à trois sont vides ; par construction, ce recouvrement raffine~$(V_j)$, et {\em a fortiori}~$(U_i)$.~$\Box$

\subsection*{Sous-graphes admissibles et revêtements topologiques}

\deux{covtop} Si~$X$ est un graphe, on notera~$\mathsf{Revtop}\;X$ la catégorie des revêtements topologiques de~$X$.  

\deux{remrecouvgr} Si~$X$ est un graphe (resp. un graphe localement fini) et si~$Y\to X$ est un revêtement topologique alors~$Y$ est séparé, et tout point de~$Y$ a un voisinage ouvert qui est un arbre (resp. un arbre localement fini) ; par conséquent,~$Y$ est un graphe (resp. un graphe localement fini). 

\deux{arbresimplco} {\bf Proposition.} {\em Soit~$X$ un arbre. Tout revêtement topologique de~$X$ est trivial.}

\medskip
{\em Démonstration.} Soit~$Y\to X$ un revêtement topologique de~$X$ ; nous allons montrer qu'il est trivial. On peut raisonner composante par composante, et partant supposer que~$Y$ est connexe et non vide ; il s'agit alors de vérifier que les fibres de~$Y\to X$ sont des singletons. 

\medskip
On raisonne par l'absurde. On suppose donc qu'il existe un point~$x$ de~$X$ et deux antécédents distincts~$y_1$ et~$y_2$ de~$x$ sur~$Y$. Soit~$I$ un segment tracé sur~$Y$ et d'extrémités~$y_1$ et~$y_2$, et soit~$\Gamma$ l'image de~$I$ sur~$X$. 

Comme~$Y\to X$ est un revêtement,~$I\to X$ est localement injective sur sa source. Joint à la compacité de~$I$, cela entraîne que~$\Gamma$ est une réunion finie de segments ; étant par ailleurs connexe,~$\Gamma$ est un sous-arbre compact et fini de~$X$. 

\medskip
L'arbre compact et fini~$\Gamma$ contient~$x$, sans être réduit à~$\{x\}$ (la flèche~$I\to X$ est localement injective, et en particulier non constante) ; par conséquent,~$\Gamma$ possède au moins un point unibranche~$\xi$ qui est différent de~$x$ ; choisissons un antécédent~$\eta$ de~$\xi$ sur~$I\setminus\{y_1,y_2\}$. Comme~$I\to \Gamma$ est injective au voisinage de~$\eta$, il existe un intervalle ouvert tracé sur~$\Gamma$ et contenant~$\xi$, ce qui contredit le caractère unibranche de ce dernier et achève la démonstration.~$\Box$ 

\deux{equivrevsquel} {\bf Théorème.} {\em Soit~$X$ un graphe et soit~$\Gamma$ un sous-graphe admissible de~$X$. Soit~$j$ l'inclusion de~$\Gamma$ dans~$X$ et soit~$r$ la rétraction canonique de~$X$ sur~$\Gamma$. Les foncteurs~$j^*$ et~$r^*$ établissent une équivalence entre~$\mathsf{Revtop}\; X$ et~$\mathsf{Revtop}\;\Gamma$.}

\medskip
{\em Démonstration.} Nous procédons en deux temps. 

\trois{descrevgraphe} {\em Généralités sur les revêtements d'un graphe.} Soit~$Y$ un graphe, et soit~$\sch F$ un ensemble de sous-arbres ouverts {\em non vides} de~$Y$ qui recouvrent~$Y$ ; notons~$\sch N(\sch F)$ l'ensemble des triplets~$(Z,Z',U)$ où~$Z$ et~$Z'$ appartiennent à~$\sch F$ et où~$U$ est une composante connexe de~$Z\cap Z'$. Il résulte de la proposition~\ref{arbresimplco} que~$\mathsf{Revtop}\;Y$ est équivalente de manière naturelle à la catégorie~$\mathsf C_{\sch F}$ dont la définition suit. 

\medskip
\begin{itemize}

\item[$\bullet$] Un objet de~$\mathsf C_{\sch F}$ consiste en : 

\medskip

\begin{itemize}

\item[$\diamond$] pour tout~$Z\in \sch F$, un ensemble~$\mathsf E_Z$ ; 

\item[$\diamond$] pour tout~$(Z,Z',U)\in \sch N(\sch F)$, un isomorphisme~$\iota_{Z,Z',U} :\mathsf E_Z\simeq \mathsf E_{Z'}$ ; 

\medskip
\end{itemize}

ces données étant sujettes aux conditions de cocycle usuelles. 
\medskip

\item[$\bullet$] Une flèche entre deux objets~$\left((\mathsf E_Z),( \iota_{Z,Z',U})\right)$ et~$\left((\mathsf F_Z), (\upsilon_{Z,Z',U})\right)$ consiste en une famille d'applications~$(\mathsf E_Z\to \mathsf F_Z)$ compatible aux isomorphismes~$\iota_{Z,Z',U}$ et~$\upsilon_{Z,Z',U}$. 

\end{itemize}

\medskip
\trois{applicdeescrev} {\em Retour à la preuve proprement dite}. Soit~$\sch F$ un ensemble de sous-arbres ouverts non vides de~$\Gamma$ qui recouvrent~$\Gamma$. Notons~$\sch G$ l'ensemble des~$r\inv(\Delta)$ pour~$\Delta$ parcourant~$\sch F$ ; c'est un ensemble de sous-arbres ouverts non vides de~$X$ qui recouvrent~$X$. 

\medskip
Pour tout ouvert~$\Omega$ de~$\Gamma$, on a~$\Omega=\Gamma\cap r\inv(\Omega)$, et~$\pi_0(\Omega)\simeq \pi_0(r\inv(\Omega))$. Il s'ensuit : 

\medskip
$\bullet$ que~$\Delta\mapsto r\inv(\Delta)$ et~$Z\mapsto \Gamma\cap Z$ mettent~$\sch F$ et~$\sch G$ en bijection ; 

$\bullet$ que~$(\Delta,\Delta',\Omega)\mapsto (r\inv(\Delta),r\inv(\Delta'),r\inv(\Omega))$~$${\rm et}\;(Z,Z',U)\mapsto (\Gamma\cap Z,\Gamma\cap Z',\Gamma \cap U)$$ mettent~$\sch N(\sch F)$ et~$\sch N(\sch G)$ en bijection. 

\medskip
Ces bijection permettent de définir deux foncteurs~$ \mathsf C_{\sch F}\to \mathsf C_{\sch G}$ et~$\mathsf C_{\sch G}\to \mathsf C_{\sch F}$ quasi-inverses (et même {\em strictement} inverses) l'un de l'autre, et partant deux foncteurs~$\mathsf{Revtop}\;\Gamma\to \mathsf{Revtop}\;X$ et~$\mathsf{Revtop}\;X\to \mathsf{Revtop}\;\Gamma$ quasi-inverses l'un de l'autre. Par construction, ces deux derniers sont respectivement isomorphes à~$r^*$ et~$j^*$, ce qui achève la démonstration.~$\Box$ 

\deux{corollh1} {\bf Corollaire.} {\em Soit~$X$ un graphe et soit~$\ell$ un entier au moins égal à~$2$. Les assertions suivantes sont équivalentes : 

\medskip
i) tout revêtement topologique de~$X$ est trivial ; 

ii)~$\H^1(X,\ZZ/\ell \ZZ)=0$ ; 

iii) les composantes connexes de~$X$ sont des arbres.

}

\medskip
{\em Démonstration.} Il est clair que i)$\Rightarrow$ii). Montrons que ii)$\Rightarrow$iii) ; nous allons plus précisément établir la contraposée. On fait donc l'hypothèse que l'une des composantes connexes de~$X$ n'est pas un arbre ; il existe alors une boucle~$C$ sur~$X$. Choisissons un sous-graphe admissible et localement fini~$\Gamma$ de~$X$ ; il contient~$\mathsf S(X)$, et partant~$C$. Nous allons montrer que~$\H^1(\Gamma,\ZZ/\ell)\neq 0$, ce qui entraînera, en vertu du théorème~\ref{equivrevsquel}, que~$\H^1(X,\ZZ/\ell)\neq 0$ et achèvera ainsi la preuve. 

\medskip
Comme~$C$ est un cercle, il existe un~$\ZZ/\ell\ZZ$-torseur {\em connexe}~$C'\to C$. Le graphe~$\Gamma$ étant localement fini, il existe un intervalle ouvert~$I$ non vide tracé sur~$C$ ne rencontrant aucun sommet de~$\Gamma$ ; l'intervalle~$I$ est alors un ouvert {\em de~$\Gamma$} ; fixons ~$x\in I$.

\medskip
Le torseur~$C'$ est déployé au-dessus de~$C\setminus\{x\}$ ; donnons-nous un isomorphisme~$\iota : C'\times_C(C\setminus\{x\})\simeq (C\setminus\{x\})\times \ZZ/\ell \ZZ$. Soit~$\Gamma'$ le~$\ZZ/\ell\ZZ$-torseur sur~$\Gamma$ obtenu en recollant~$(\Gamma\setminus\{x\})\times \ZZ/\ell \ZZ$ et~$C'\times_CI$ {\em via} la restriction de~$\iota$ à~$I\setminus\{x\}$. Par construction,~$\Gamma'\times_\Gamma C\simeq C'$ ; en conséquence~$\Gamma'\to \Gamma$ n'est pas trivial et~$\H^1(\Gamma,\ZZ/\ell \ZZ)\neq 0$, ce qu'on souhaitait prouver. 

\medskip
Enfin l'implication iii)$\Rightarrow$i) est une conséquence directe du lemme~\ref{arbresimplco}.~$\Box$ 

%
%
%
%
%

\section{Toises et homotopies}

\subsection*{Toises sur un graphe}

\deux{deflong} Soit~$X$ un graphe et soit~$Y$ une partie de~$X$ ; soit~$\sch I$ l'ensemble des segments tracés sur~$Y$. On appellera {\em toise} sur~$Y$ toute application~$l$ de~$\sch I$ dans~$\RR_+$ possédant les propriétés suivantes : 

\medskip
$\bullet$ si~$x\in Y$ et si~$I$ est un segment issu de~$x$ et tracé sur~$Y$ alors~$t\mapsto l([x;t])$ induit un homéomorphisme~$I\simeq [0:l(I)]$ ;

$\bullet$ si~$I\in \sch I$ et si~$I'$ et~$I''$ sont deux segments tracés sur~$I$ tels que~$I=I'\cup I''$ alors~$l(I)=l(I')+l(I'')-l(I\cap I'')$. 

\medskip
On parlera du réel~$l(I)$ comme de la {\em longueur} de~$I$ (relative à la toise~$l$). On dira qu'une toise~$l$ sur~$Y$ est bornée par un réel~$M$ si~$l(I)\leq M$ pour tout~$I\in \sch I$. 

\deux{remlocfinmetr} {\em Remarque.} Soit~$X$ un graphe muni d'une toise~$l$. Pour tout~$x\in X$ et tout~$r>0$ notons~$\mathsf B(x,r)$ l'ensemble des points~$y$ de~$X$ tels qu'il existe un segment tracé sur~$X$ joignant~$y$ à~$x$ et de longueur strictement inférieure à~$r$. Les~$\mathsf B(x,r)$ constituent, pour~$x$ parcourant~$X$ et~$r$ parcourant~$\RR\ti_+$, une base d'une topologie sur~$X$. On prendra garde que celle-ci est en général {\em strictement plus fine} que la topologie de~$X$ (elle permet de sectionner simultanément toutes les branches issues d'un point donné) ; c'est pour cette raison que nous avons choisi d'utiliser le terme «toise» plutôt que «métrique», qui nous a semblé être implicitement supposé définir la topologie de l'espace considéré. 

\medskip
Notons toutefois un cas dans lequel la topologie définie par~$l$ coïncide avec la topologie de~$X$ : c'est celui où~$X$ est localement fini. Dans cette situation, une toise n'est donc rien d'autre qu'une métrique au sens classique, sans ambiguïté ni risque de confusion ; un graphe localement fini muni d'une toise (resp. admettant une toise) est simplement un graphe localement fini métrisé (resp. métrisable). 
 
\deux{longcaslocfin} {\bf Proposition.} {\em Soit~$X$ un graphe localement fini et soit~$S$ un sous-ensemble fermé et discret de~$X$. Les assertions suivantes sont équivalentes :

\medskip
i)~$X$ admet une toise ; 

ii)~$X\setminus S$ admet une toise ; 

iii)~$X$ est paracompact.}

\medskip
{\em Démonstration.} Il est clair que i)$\Rightarrow$ ii). 

\medskip
\trois{longxs}{\em Montrons que ii) entraîne i) et iii).} On fait l'hypothèse qu'il existe une toise~$l_0$ sur~$X\setminus S$. Quitte à agrandir~$S$, on peut supposer que celui-ci contient tous les sommets de~$X$, ainsi qu'au moins un point sur chaque composante connexe de~$X$ qui est une boucle. Soit~$I$ une composante connexe de~$X\setminus S$ ; c'est une droite éventuellement longue. Si~$x\in I$, la toise~$l_0$ permet d'immerger chacune des deux demi-droites de~$I$ issues de~$x$ dans~$\RR_+$ ; par conséquent,~$I$ est paracompact et est donc une vraie droite. Il existe dès lors un homéomorphisme~$\iota_I : I\simeq ]0;1[$. 

\medskip
Les homéomorphismes~$\iota_I$ permettent de définir la longueur~$l_1(J)$ de tout intervalle (compact ou non) tracé sur~$X\setminus S$ ; c'est un réel appartenant à~$[0;1]$. 

\medskip
Soit ~$I$ un segment tracé sur~$X$ ; l'ensemble~$I\cap S$ est fini, et~$I\setminus S$ a donc un nombre fini de composantes connexes. On note alors~$l(I)$ la somme des des~$l_1(J)$, où~$J$ parcourt~$\pi_0(I\setminus S)$. Il est immédiat que l'application~$l$ ainsi construite est une toise sur~$X$. 

\medskip
Nous allons montrer que~$X$ est paracompact ; soit~$X'$ une composante connexe de~$X$, et soit ~$x\in X'$. Pour tout~$N\geq 1$, notons~$X'_N$ la réunion des segments~$I$ tracés sur~$X'$ et issus de~$x$ possédant les propriétés suivantes : 

\medskip
i)~$I\cap (S\cup\{x\})$ possède au plus~$N$ éléments ; 

ii) si~$y$ désigne l'autre extrémité de~$I$ et si~$z$ désigne le point de~$I\cap (S\cup\{x\})$ le plus proche de~$y$ alors~$l([y;z])\leq 1-1/N$. 

\medskip
Il résulte de la construction de~$l$ que chacun des~$X'_N$ est compact, et que les~$X'_N$ recouvrent~$X'$. Par conséquent,~$X'$ est dénombrable à l'infini et~$X$ est paracompact. 

\trois{parlong}{\em Montrons que iii) entraîne i).} Supposons~$X$ paracompact, et soit~$\Sigma$ un sous-ensemble fermé et discret de~$X$ contenant ses sommets ainsi qu'au moins un point par composante connexe de~$X$ qui est une boucle. Chaque composante connexe de~$X-\Sigma$ est un ouvert de~$X$ et est donc paracompacte ; comme une telle composante est par ailleurs une droite éventuellement longue, les composantes connexes de~$X-\Sigma$ sont finalement de vraies droites. Il existe dès lors une toise sur~$X-\Sigma$ ; en vertu de l'implication ii)$\Rightarrow$i) déjà établie, cela entraîne l'existence d'une toise sur~$X$.~$\Box$ 

\deux{equivparac} {\bf Corollaire.} {\em Soit~$X$ un graphe. Les trois propositions suivantes sont équivalentes : 

\medskip
i) tout sous-graphe localement fini de~$X$ admet une toise ; 

ii) il existe un sous-graphe localement fini et admissible de~$X$ admettant une toise ; 

iii)~$X$ est paracompact.}

\medskip
{\em Démonstration.} Il est clair que i)$\Rightarrow$ii). 

\medskip
Supposons que ii) soit vraie et soit~$\Gamma$ un sous-graphe admissible et localement fini de~$X$ admettant une toise. En vertu de la proposition~\ref{longcaslocfin} ci-dessus,~$\Gamma$ est paracompact. La rétraction canonique de~$X$ sur~$\Gamma$ étant compacte, iii) est vérifiée. 

\medskip
Supposons que iii) soit vraie et soit~$\Delta$ un sous-graphe localement fini de~$X$. Comme~$X$ est paracompact,~$\Delta$ est paracompact ; il résulte alors de la proposition~\ref{longcaslocfin} ci-dessus que~$\Delta$ admet une toise ; par conséquent, i) est vérifiée, ce qui achève la démonstration.~$\Box$

\deux{paramfini} Soit~$X$ un graphe et soit~$\sch E$ un ensemble fini de sous-graphes fermés de~$X$ recouvrant~$X$ ; supposons que chacun des graphes appartenant à~$\sch E$ admette une toise ; il existe alors une toise sur~$X$ lui-même. Pour le voir, on choisit pour tout~$Y\in \sch E$ une toise~$l_Y$ sur ~$Y$ ; on étend par additivité la définition de~$l_Y(I)$ à tout compact~$I$ de~$Y$ qui est une réunion finie disjointe de segments. 

L'application~$I\mapsto \sum\limits_{Y\in \sch E} l_Y(I\cap Y)$ est alors une toise sur~$X$. 

\deux{paramtoise} Soit~$X$ un arbre, soit~$Y$ une partie convexe de~$X$ et soit~$x\in Y$. Une {\em toise partielle sur~$Y$ basée en~$x$} est une application~$\delta : Y\to \RR_+$ telle que~$\delta$ induise pour tout~$y\in Y$ un homéomorphisme~$[x;y]\simeq [0;\delta(y)]$. 

\medskip
Soit~$\delta$ une toise partielle sur~$Y$ basée en~$x$. Soit~$I$ un segment tracé sur~$Y$, soient~$y$ et~$z$ ses deux extrémités et soit~$t$ le point tel que~$[y;x]\cap [z;x]=[t;x]$. Posons~$l(I)=\delta(y)+\delta(z)-\delta(t)$ ; on vérifie aussitôt : que~$l$ est une toise sur~$Y$, que l'on dira {\em induite par~$\delta$} ; que~$l([y;x])=\delta(y)$ pour tout~$y\in Y$ ; et que si~$\delta$ est bornée par un réel~$M$ alors~$l$ est bornée par 2M. 

\deux{yconvmborn} Soit~$X$ un arbre et soit~$Y$ un sous-ensemble convexe de~$X$. Supposons que~$Y$ admet une toise~$l$ ; il admet alors une toise~$l'$ {\em bornée par~$1$}. 

\medskip
En effet, si~$Y$ est vide c'est évident ; supposons maintenant~$Y\neq \emptyset$ et choisissons~$x$ sur~$Y$. L'application~$\delta$ de~$Y$ dans~$\RR_+$ qui envoie~$y$ sur~$\displaystyle{\frac{l([y;x])}{2(1+l([y;x]))}}$ est une toise partielle basée en~$x$, qui est par construction bornée par~$1/2$ ; on peut alors prendre pour~$l'$ la toise sur~$Y$ induite par~$\delta$. 

\deux{toiseconcat} Soit~$X$ un graphe et soit~$\Gamma$ un sous-graphe admissible de~$X$. Supposons données : 

\medskip
$\bullet$ une toise~$l_\Gamma$ sur~$\Gamma$ ; 

$\bullet$ pour toute composante connexe~$U$ de~$X-\Gamma$, une toise~$l_U$ sur~$\overline U$. 

\medskip
Si~$x\in X$ on pose~$\delta(x)=0$ si~$x\in \Gamma$, et~$\delta (x)=l_U([x;r(x)])$ sinon, où~$U$ est la composante connexe de~$X-\Gamma$ contenant~$x$. Soit~$I$ un segment tracé sur~$X$ et soient~$x$ et~$y$ ses extrémités. On définit le réel~$l(I)$ comme suit : 

\medskip
- si~$I\cap \Gamma=\emptyset$, on pose~$l(I)=l_U(I)$ où~$U$ est la composante connexe de~$X-\Gamma$ contenant~$I$ ; 

- sinon, on sait ({\em cf.}~\ref{concatena}) que~$I$ est de la forme~$[x;r(x)]\cup J\cup[r(y);y]$, où~$J$ est un segment tracé sur~$\Gamma$ joignant~$r(x)$ à~$r(y)$, et l'on pose alors~$$l(I)=\delta(x)+\delta(y)+l_\Gamma(J).$$

\medskip
On vérifie aussitôt que~$l$ est une toise sur~$X$, dont on dira qu'elle est obtenue par {\em concaténation} de~$l_\Gamma$ et des~$l_U$. 

\deux{xadmtoise} {\bf Lemme.} {\em Soit~$X$ un graphe. Les assertions suivantes sont équivalentes :

\medskip
i)~$X$ admet une toise ; 

ii) tout sous-graphe compact de~$X$ admet une toise et tout sous-graphe localement fini de~$X$ admet une toise ;

iii)~$X$ est paracompact et tout sous-graphe compact de~$X$ admet une toise.}

\medskip
{\em Démonstration.} Il est clair que i)$\Rightarrow$ii), et ii)$\iff$iii) se déduit du corollaire~\ref{equivparac}.

\medskip
Supposons que ii) soit vérifiée. Choisissons un sous-graphe localement fini et admissible~$\Gamma$ de~$X$. Par hypothèse,~$\Gamma$ admet une toise~$l_\Gamma$. Si~$U$ est une composante connexe de~$X-\Gamma$, son adhérence~$\overline U$ est compacte et admet donc une toise~$l_U$. La concaténation des~$l_U$ et de~$l_\Gamma$ fournit alors une toise sur~$X$ ; ainsi, ii)$\Rightarrow$i), ce qui achève la démonstration.~$\Box$

\subsection*{Homotopies}

\deux{paramhom} Soit~$X$ un graphe, soit~$\Gamma$ un sous-graphe admissible de~$X$ et soit~$r$ la rétraction canonique de~$X$ sur~$\Gamma$. Si~$U$ est une composante connexe de~$X-\Gamma$, l'unique point de son bord sera noté~$\gamma_U$. 

\medskip
Supposons donnée pour toute composante connexe~$U$ de~$X-\Gamma$ une toise~$l_U$ {\em bornée par 1} sur l'arbre~$\overline U$ ; nous allons construire explicitement, en fonctions des~$l_U$, une homotopie reliant~$r$ à l'identité. 

Si~$x\in X$ on posera~$\delta(x)=0$ si~$x\in \Gamma$ et~$\delta(x)=l_U([x;\gamma_U])$ si~$x\notin \Gamma$ et si~$U$ est la composante connexe de~$X-\Gamma$ contenant~$x$ ; par construction,~$\delta$ est majorée par~$1$ et induit pour tout~$x$ un homéomorphisme entre~$[r(x);x]$ et~$[0;\delta(x)]$.

\deux{theoparhom}{\bf Proposition.} {\em On reprend les notations du~\ref{paramhom} ci-dessus. L'application~$h$ de~$[0;1]\times X$ vers~$X$
 qui envoie un couple~$(x,t)$ sur~$\delta_{|[r(x);x]}\inv(t)$ si~$t\leq \delta(x)$ et sur~$x$ sinon est continue ; on a~$h(0,.)=r$ et~$h(1,.)=\mathsf{Id}_X$.}
 
 \medskip
 {\em Démonstration.} Seule la continuité de~$h$ nécessite une preuve ; celle-ci requiert un certain nombre de préliminaires. 
 
 \trois{imrecxleqy} {\em Images réciproques par~$h$ de sous-ensembles particuliers.} Fixons une composante connexe~$U$ de~$X-\Gamma$ et soit~$x$ un point de~$U$. 
 
\medskip
Soit~$V$ une composante connexe de~$\overline U\setminus\{x\}$ ne contenant pas~$\gamma_U$, c'est-à-dire encore une composante connexe relativement compacte de~$U\setminus\{x\}$ ; il résulte de la définition de~$h$ que~$h\inv(V)\subset [0;1]\times U$. Le complémentaire de~$V$ dans~$\overline U$ est un compact qui est réunion de~$\{x\}$ et de  composantes connexes de~$\overline U\setminus\{x\}$ ; il est donc convexe. Il s'ensuit que si~$z\in \overline U-V$ alors~$[z;\gamma_U]\subset \overline U-V$, et l'on a donc~$h(t,z)\in \overline U-V$ pour tout~$t$ ; par conséquent,~$h\inv(V)\subset [0;1]\times V$. 

Si~$z\in V$ alors~$[z;\gamma_U]$ contient~$x$, et~$[z;\gamma_U]\cap V=[z;x[$ ; dès lors~$h(t,z)\in V$ si et seulement si~$t>\delta(x)$. 

En conclusion, on a démontré que~$h\inv(V)=]\delta(x); 1]\times V$ ; c'est un ouvert de~$[0;1]\times X$.

\medskip
Soit~$W$ la composante connexe de~$\gamma_U$ dans~$\overline U\setminus\{x\}$ ; son fermé complémentaire dans~$\overline U$ est le compact~$\sch K_U(x)$. Il résulte de la définition de~$h$ que~$$h\inv(\sch K_U(x))\subset [0;1]\times U.$$ Si~$z\in W$ alors~$[z;\gamma_v]\subset W$, et l'on a donc~$h(t,z)\in W$ pour tout~$t$ ; par conséquent,~$h\inv(\sch K_U(x))\subset [0;1]\times \sch K_U(x).$

Si~$z\in \sch K_U(x)$ alors~$[z;\gamma_U]$ contient~$x$, et~$[z;\gamma_U]\cap \sch K_U(x)=[z;x]$ ; dès lors~$h(t,z)\in \sch K_U(x)$ si et seulement si~$t\geq\delta(x)$. 

En conclusion, on a démontré que~$h\inv(\sch K_U(x))=[\delta(x); 1]\times \sch K_U(x)$ ; c'est un compact de~$[0;1]\times X$.

\trois{bordfiniubarre} On désigne toujours par~$U$ une composante connexe de~$X-\Gamma$. Soit~$V$ un ouvert connexe, non vide et à bord fini de~$\overline U$ pose~$S=\partial_{\overline U}V$. Si~$W$ est une composante connexe de~$\overline U\setminus S$ autre que~$V$, c'est une composante connexe de~$\overline U-\overline V$, et comme~$\overline V$ est un arbre, le bord de~$W$ est nécessairement un singleton, dont l'unique élément appartient à~$S$. 

\medskip
Soit~$x\in S$. On peut écrire~$\overline U$ comme la réunion disjointe de~$\{x\}$, des composantes connexes de~$\overline U-\overline V$ de bord~$\{x\}$, de~$V$, de~$S\setminus\{x\}$, et des composantes connexes de~$\overline U-\overline V$ dont l'unique point du bord diffère de~$x$. La réunion de~$V$, de~$S\setminus\{x\}$, et des composantes connexes de~$\overline U-\overline V$ dont l'unique point du bord diffère de~$x$ étant convexe, c'est une composante connexe de~$\overline U\setminus\{x\}$.

\trois{declinschk} Les faits énoncés ci-dessus se déclinent comme suit dans deux cas particuliers. 

\medskip
{\em Supposons que~$V$ contient~$\gamma_U$.}  On a alors~$$\overline U-V=\coprod_{x\in S}\sch K_U(x).$$

\medskip
{\em Supposons que~$V\subset U$ et que~$\gamma_U\notin S$ (autrement dit,~$V$ est relativement compact dans~$U$)}. Soit~$x$ l'unique point du bord de la composante connexe de~$\overline U-\overline V$ contenant~$\gamma_U$, et soit~$W$ la composante connexe de~$\overline U\setminus\{x\}$ contenant~$V$. On a alors~$$V=W-\coprod_{y\in S\setminus\{x\}} \sch K_U(y).$$

\trois{conth}{\em Continuité de~$h$.} Il s'agit de montrer que si~$V$ est un ouvert de~$X$ alors~$h\inv(V)$ est ouvert ; il suffit bien entendu se limiter au cas où~$V$ appartient à une base donnée de la topologie de~$X$. Il est donc loisible de supposer que~$V$ est connexe, non vide, et à bord fini ; on peut de surcroît faire l'hypothèse que si~$V$ ne rencontre pas~$\Gamma$, il est relativement compact dans~$X-\Gamma$. 

\medskip
{\em Le cas où~$V$ est un ouvert relativement compact de~$X-\Gamma$.} Soit~$U$ la composante connexe de~$X-\Gamma$ contenant~$V$ ; l'ouvert~$V$ est connexe, non vide, relativement compact et à bord fini dans~$U$. D'après~\ref{declinschk} il existe~$x\in \partial V$ et une composante connexe relativement compacte~$W$ de~$V\setminus\{x\}$ tels que~$V=W-\coprod\limits_{y\in \partial V \setminus\{x\}} \sch K_U(y).$ Il découle alors de~\ref{imrecxleqy} que~$h\inv(V)$ est un ouvert de~$[0;1]\times X$. 

\medskip
{\em Le cas où~$V$ rencontre~$\Gamma$.} Soit~$V_0$ l'intersection de~$V$ et de~$\Gamma$. C'est un ouvert non vide de~$\Gamma$ dont le bord est l'ensemble fini~$\partial V\cap \Gamma$ ; de plus,~$V_0$ est connexe en vertu de~\ref{admimplcon}. Le théorème~\ref{thretcan} assure alors que~$r\inv(V_0)$ est la composante connexe de~$X-(\partial V\cap \Gamma)$ contenant~$V_0$ ; il s'ensuit que~$V\subset r\inv(V_0)$. 

\medskip
Soit~$U$ une composante connexe de~$r\inv(V_0)-V_0$ ;  en d'autres termes,~$U$ est une composante connexe de~$X-\Gamma$ telle que~$\gamma_U\in V_0$. Comme~$\gamma_U\in V$, l'intersection~$V\cap U$ est un ouvert non vide de~$U$. 

Si~$U$ ne rencontre pas~$\partial V$, cet ouvert est fermé dans~$U$ et coïncide donc avec~$U$, ce qui signifie que~$U\subset V$.

Supposons que~$U$ rencontre~$\partial V$ (comme ce dernier est fini, cela ne se produit que pour un ensemble fini de composantes connexes de ~$r\inv(V_0)-V_0$),  et soit~$S$ l'ensemble fini~$U\cap \partial V$. Si~$x\in U\cap V$ la connexité de~$V$ entraîne l'existence d'un segment~$I$ tracé sur~$V$ et joignant~$x$ à~$\gamma_U$. Soit~$J$ la composante connexe de~$x$ dans~$I\cap U$ ; sa borne supérieure appartient à~$\partial U$ et est donc égale à~$\gamma_U$ ; par conséquent,~$I$ est le segment~$[x;\gamma_U]$ de l'arbre~$\overline U$, et~$[x;\gamma_U]\subset V$. Il s'ensuit que l'ouvert~$V\cap \overline U$ de~$\overline U$ est connexe ; son bord dans~$\overline U$ est par construction égal à~$S$. 

Il découle alors de~\ref{declinschk} que~$\overline U-V= \coprod\limits_{x\in S}\sch K_U(x).$ 

\medskip
Pour tout~$x\in \partial  V- \Gamma$, notons~$U_x$ la composante connexe de~$X-\Gamma$ qui contient~$x$. Il résulte de ce qui précède que~$$V=r\inv(V_0)-\coprod_{x\in \partial V-\Gamma} \sch K_{U_x}(x).$$ On voit immédiatement à l'aide de la définition de~$h$ que~$h\inv(r\inv(V_0))$ est égal à~$[0;1]\times r\inv(V_0)$ ; et d'après~\ref{imrecxleqy},~$h\inv( \sch K_{U_x}(x))$ est compact pour tout~$x\in \partial V s-\Gamma$ ; par conséquent,~$h\inv(V)$ est un ouvert de~$[0;1]\times X$.~$\Box$

\deux{fortcont} Rappelons qu'un espace topologique~$X$ est dit {\em fortement contractile} s'il est non vide et si~$X$ se rétracte par déformation sur {\em chacun} de ses points.

\deux{theoadmtoise} {\bf Théorème.} {\em Soit~$X$ un graphe admettant une toise. 

\medskip
\begin{itemize}

\item[1)] Soit~$\Gamma$ un sous-graphe admissible de~$X$ et soit~$r$ la rétraction canonique de~$X$ sur~$\Gamma$. L'inclusion de~$\Gamma$ dans~$X$ est une équivalence homotopique ; plus précisément il existe une application continue~$h: [0;1]\times X\to X$ telle que~$h(0,.)=r, h(1,.)=\mathsf{Id}_X$ et~$h(t,x)=x$ pour tout~$t\in [0;1]$ et tout~$x\in \Gamma$.

\item[2)] Le graphe~$X$ est paracompact. 

\item[3)] Le graphe~$X$ est localement fortement contractile, et a le type d'homotopie d'un graphe localement fini métrisable. 

 \item[4)] Les assertions suivantes sont  équivalentes  :
 
 \medskip
 \begin{itemize}
 \item[$i)$]~$X$ est fortement contractile ; 
 
 \item[$ii)$]~$X$ est contractile ; 
 
 \item[$iii)$]~$X$ est un arbre non vide.
\end{itemize}
\end{itemize}}

\medskip
{\em Démonstration}. Choisissons un sous-graphe localement fini et admissible~$\Delta$ de~$X$.

\trois{rephrasehomo} {\em Preuve de 1}. Soit~$U$ une composante connexe de~$X-\Gamma$ ; par restriction, l'arbre compact~$\overline U$ admet par hypothèse une toise~$l_U$, et il résulte de~\ref{yconvmborn} qu'on peut choisir~$l_U$ bornée par 1 ; la proposition~\ref{theoparhom} fournit alors l'homotopie requise. 

\trois{xtoiseparac} {\em Preuve de 2) et d'une partie de 3)}. Par restriction, le graphe localement fini~$\Delta$ admet une toise (ce qui revient à dire qu'il est métrisable), et le corollaire~\ref{equivparac} assure alors que~$X$ est paracompact. D'après l'assertion 1) déjà prouvée,~$X$ est homotopiquement équivalent à~$\Delta$, et a donc bien le type d'homotopie d'un graphe localement fini métrisable. 

\trois{toiseglobparacomp} {\em Preuve de 4) et conclusion.}.

\medskip
{\em Supposons que~$X$ soit fortement contractile} ; il est alors clairement contractile.

{\em Supposons que~$X$ soit contractile.} Ceci entraîne que~$X$ est connexe, non vide, et que chacun de ses revêtements topologiques est trivial ; il découle dès lors du corollaire~\ref{corollh1} que~$X$ est un arbre non vide.

{\em Supposons que~$X$ soit un arbre non vide.} Soit~$x\in X$ ; il existe un sous-graphe localement fini et admissible~$\Delta'$ de~$X$ contenant~$x$ ; étant admissible dans~$X$, le graphe~$\Delta'$ est lui-même un arbre non vide. Pour tout~$y\in \Delta'$ et tout~$t\in [0;1]$ notons~$h(t,y)$ l'unique point~$z$ de~$[x;y]$ tel que~$l([x;z])=tl([x;y])$ ; l'application~$h$ définit alors une rétraction par déformation de~$\Delta'$ sur~$\{x\}$ ; comme~$X$ admet une rétraction par déformation sur~$\Delta'$ (assertion 1) ), il admet par concaténation une rétraction par déformation sur~$\{x\}$ ; en conséquence,~$X$ est fortement contractile. 

\medskip
Enfin tout point de~$X$ a un voisinage ouvert qui est un arbre, admettant par restriction une toise, et est donc par ce qui précède fortement contractile ; ainsi,~$X$ est localement fortement contractile, ce qui achève la preuve de 3) et la démonstration.~$\Box$

\section{Branches d'un graphe}

\deux{branchegr} Soit~$X$ un graphe et soit~$x\in X$ ; on notera~$\arb X x$ l'ensemble des sous-arbres ouverts de~$X$ contenant~$x$. Une {\em branche issue de~$x$} sera un élément de~$\lim\limits_{\stackrel \leftarrow V}\pi_0(V\setminus \{x\})$ pour~$V$ parcourant~$\arb X x$ ; on notera~$\br X x$ l'ensemble des branches de~$X$ issues de~$x$ ; la flèche naturelle~$\br X x \to \pi_0(V\setminus\{x\})$ est une bijection pour tout~$V\in \arb X x$. 

\trois{brgraphefin} Le cardinal de~$\br X x$ est égal à la valence de~$(X,x)$. 

\trois{bdev} Si~$b$ est une branche issue de~$x$ et si~$V\in \arb X x$, on notera~$b(V)$ l'image de~$b$ dans~$\pi_0(V\setminus\{x\})$ ; si~$V\in \arb X x$ et si~$V'\in \arb V x$ alors~$b(V')=b(V)\cap V'$.  

\trois{defsec} Si~$b$ est une branche issue de~$x$, une {\em section} de~$b$ sera un ouvert de~$X$ de la forme~$b(V)$ pour un certain~$V$ ; on désignera par~$\sbr b$ l'ensemble des sections de~$b$ ; remarquons que deux sections de~$b$ se rencontrent toujours.

On notera~$\sgb X x$  l'ensemble des ouverts de~$X$ qui appartiennent à~$\sbr b$ pour au moins une branche~$b$ issue de~$x$. 

\trois{finibranchecontr} Si~$\{b_1,\ldots,b_r\}$ est un ensemble fini de branches de~$X$ issues de~$x$, si~$Z_i$ est pour tout~$i$ une section de~$b_i$, et si~$U$ est un voisinage ouvert de~$x$, il existe~$V\in \arb X x$ que~$b_i(V)\subset Z_i$ pour tout~$i$ : il suffit de se donner pour tout~$i$ un élément~$W_i$ de~$\arb X x$ tel que~$b_i(W_i)=Z_i$, et de prendre pour~$V$ n'importe quel élément de~$\arb {U\cap \bigcap W_i} x$.

\trois{etresecb}{\bf Lemme.} {\em Soit~$U$ un ouvert connexe et non vide de~$X$. 

\medskip
1) Si~$U\in \sgb X x$ il existe une {\em unique} branche~$b\in \br X x$ telle que~$U\in \sbr b$.

2) Si~$U\in \sgb X x$ alors~$x$ est un point isolé de~$\partial U$. 

3) Si~$x$ est un point isolé de~$\partial U$ et si~$\overline U$ est un arbre alors~$U\in \sgb X x$.}

\medskip
{\em Démonstration.} Supposons que~$U\in \sbr b$ pour une certaine~$b\in \br X x$, et soit~$V\in \arb X x$ tel que~$b(V)=U$. Si~$\beta$ est une branche issue de~$x$ et différente de~$b$, la composante~$\beta(V)$ de~$V\setminus\{x\}$ est disjointe de~$b(V)$ ; c'est donc une section de~$\beta$ qui ne rencontre pas~$U$, ce qui exclut que~$U$ lui-même soit une section de~$\beta$, d'où 1). 

\medskip
Supposons que~$U\in \sgb X x$ ; il existe alors~$V\in \arb X x$ tel que~$U$ soit une composante connexe de~$V\setminus\{x\}$. On a alors~$V\cap \partial U=\{x\}$, et~$x$ est bien un point isolé de~$\partial U$, d'où 2). 

\medskip
Supposons maintenant que~$x$ soit un point isolé de~$\partial U$ et que~$\overline U$ soit un arbre. Comme~$x$ est isolé dans~$\partial U$, il existe un voisinage compact~$W$ de~$x$ dans~$X$ qui est un arbre et ne rencontre pas~$\partial U\setminus\{x\}$. Nous allons tout d'abord montrer que le voisinage~$W\cup \overline U$ de~$x$ est un arbre. 

\medskip
Par construction,~$W\cup \overline U$ est un sous-graphe fermé et connexe de~$X$ ; pour s'assurer que c'est un arbre, on peut se contenter de vérifier qu'il ne contient aucune boucle. On raisonne par l'absurde : supposons qu'il existe une boucle~$C$ tracée sur~$W\cup \overline U$. Comme~$W$ est un arbre, il existe un point~$y$ sur~$C-W$ ; comme~$\overline U$ est un arbre, il existe un point~$z$ sur~$C-\overline U$. Comme~$\partial_{\overline U\cup W}\overline U=\{x\}$, chacun des deux intervalles ouverts tracés sur~$C$ et joignant~$y$ à~$z$ contient~$x$, ce qui est contradictoire. 

\medskip
Soit~$V$ la composante connexe de~$x$ dans l'intérieur de~$W\cup\overline U$ ; c'est un sous-arbre ouvert de~$X$ qui contient~$U\cup\{x\}$. Comme~$W$ est par ailleurs un compact ne rencontrant pas~$\partial U\setminus\{x\}$, aucun point de~$\partial U\setminus\{x\}$ n'appartient à l'intérieur de~$W\cup\overline U$ ; par conséquent~$\partial U\cap V=\{x\}$. Il s'ensuit que~$U$ est une partie ouverte, fermée, connexe et non vide de~$V\setminus\{x\}$, c'est-à-dire une composante connexe de~$V\setminus\{x\}$ ; ainsi,~$U\in \sgb X x$.~$\Box$ 

\trois{defbrdefisec} Si~$U\in \sgb  X x$, l'unique branche issue de~$x$ dont~$U$ soit une section ({\em cf.} lemme~\ref{etresecb} ci-dessus) sera appelée la branche {\em définie par~$U$.}  

\trois{ouvdefsec} Soit~$b$ une branche de~$X$ issue de~$x$. Soit~$V\in \arb X x$ et soit~$U$ un élément de~$\sgb X x$ contenu dans~$b(V)$. Le lemme~\ref{etresecb} assure que~$x$ est un point isolé de~$\partial U$ ; c'est {\em a fortiori} un point isolé de~$\partial_V U$. Comme~$\overline U^V$ est nécessairement un arbre puisque~$V$ en est un, on déduit de l'assertion 3) du lemme~\ref{etresecb} l'existence de~$V_0\in \arb V x$ tel que~$U$ soit une composante connexe de~$V_0\setminus\{x\}$ ; il découle alors de l'inclusion~$U\subset b(V)$ que~$U=b(V_0)$. 

\medskip
Notons qu'en vertu de ce qui précède, deux éléments de~$\sgb X x$ qui sont comparables pour l'inclusion définissent la même branche. 

\trois{ouvcontient} Si~$b$ est une branche de~$X$ issue de~$x$, nous dirons par abus qu'un ouvert~$U$ de~$X$ {\em contient}~$b$ s'il en contient une section. Si~$U$ est un ouvert de~$X$ contenant~$b$ alors~$x\in \overline U$ : cela résulte du fait que~$x$ adhère à toute section de~$b$.

\trois{compcontient} Si~$V$ est un voisinage ouvert connexe de~$x$ dans~$X$ et si~$W$ est une composante connexe de~$X\setminus\{x\}$ alors~$W$ contient au moins une branche issue de~$x$ : cela résulte immédiatement du fait que si~$U\in \arb V x$ alors~$U\cap W$ est une réunion de composantes connexes de~$U\setminus\{x\}$ qui est non vide puisque~$x$ adhère à~$W$. 

\trois{defbasesecb} Si~$b$ est une branche de~$X$ issue de~$x$, on appellera {\em base} de sections de~$b$ tout sous-ensemble~$\sch B$ de~$\sbr b$ tel que pour tout~$U\in \sbr b$, il existe~$V$ appartenant à~$\sch B$ et contenue dans~$U$.

\trois{branchegen} Si~$b$ est une branche de~$X$ issue de~$x$, on peut retrouver~$x$ et~$b$ elle-même à partir de~$\sbr b$ : le point~$x$ est l'unique point de~$X$ adhérent à tout ouvert élément de~$\sbr b$ ; et l'assertion 1) du lemme~\ref{etresecb} assure qu'une branche issue de~$x$ est uniquement déterminée par {\em n'importe laquelle} de ses sections. 

\medskip
Cette remarque permet de définir une branche de~$X$ comme un sous-ensemble de~${\sch P}(X)$ qui est l'ensemble des sections d'une branche~$\beta$ issue d'un point~$\xi$ de~$X$ ; dans ce cas~$\beta$ et~$\xi$ sont uniquement déterminés, et l'on dira parfois que~$\xi$ est {\em l'origine} de~$\beta$. 

\medskip
Si~$E$ est un sous-ensemble de~$X$, et si~$U$ est un ouvert de~$X$, on notera~$\br X E$ (resp.~$\br X E\ctd U$ l'ensemble des branches de~$X$ dont l'origine appartient à~$E$ (resp. dont l'origine appartient à~$E$ et qui sont contenues dans~$U$). 

\subsection*{Images directe et réciproque d'une branche}
 
\deux{imdir} Soit~$\phi: Y\to X$ une application continue entre graphes, soit~$x$ un point de~$X$ dont la fibre est finie, et soit~$y$ un antécédent de~$x$.

Soit~$U$ un élément de~$\arb X x$ et soit~$V$ un élément de~$\arb {\phi\inv(U)}y$ qui ne rencontre pas~$\phi\inv(x)\setminus\{y\}$. 
 Soit~$W$ une composante connexe de~$V\setminus\{y\}$ ; son image est contenue dans une composante connexe de~$U\setminus\{x\}$. Modulo les bijections~$$\br X x\simeq \pi_0 (U\setminus\{x\})\;{\rm et}\;\br Y y \simeq \pi_0(V\setminus\{y\}),$$ on définit ainsi une application~$\br Y y\to \br X x$ qui ne dépend pas du choix de~$U$ et 
$V$ et sera notée~$b\mapsto \phi(b)$. Si~$a\in \br X x$ on désignera par~$\phi\inv(a)$ l'ensemble des branches~$b$ de~$Y$ telles que~$\phi(b)=a$. 

\deux{imrec} Soit~$\phi: Y\to X$ une application continue {\em compacte et ouverte} entre graphes, soit~$x$ un point de~$X$ dont la fibre est finie et soit~$a\in \br X x$. Soit~$U_0$ un élément de~$\arb X x$ tel que~$\phi\inv(U_0)$ soit une réunion disjointe d'arbres séparant les antécédents de~$x$. 

\trois{descimrec} L'ouvert~$\phi\inv(U_0)$ est la réunion disjointe des~$\phi\inv(U_0)_y$ où~$y$ parcourt la fibre~$\phi\inv(x)$. Soit~$W$ une composante connexe de~$\phi\inv(a(U_0))$ ; elle est contenue dans~$\phi\inv(U_0)_{y}$ pour un certain~$y$ uniquement déterminé, et est plus précisément une composante connexe de~$\phi\inv(U_0)_y\setminus\{y\}$. Elle définit donc une branche~$b_W$ de~$Y$ issue de~$y$ ; il résulte de la définition de l'image d'une branche que~$W\mapsto b_W$ établit une bijection entre~$\pi_0(\phi\inv(a(U_0)))$ et~$\phi\inv(a)$ ; nous noterons~$b\mapsto b(U_0)$ la bijection réciproque. Pour tout~$y\in \phi\inv(x)$ la composante~$\phi\inv(U_0)_y$ se surjecte sur~$U_0$, ce qui entraîne l'existence d'une branche~$b$ appartenant à~$\phi\inv(a)$ telle que~$b(U_0)\subset \phi\inv(U_0)_y$, c'est-à-dire encore telle que l'origine de~$b$ soit égale à~$y$. 

\medskip
Soit~$U\in \arb {U_0} x$. On a~$\phi\inv(U)_y=\phi\inv(U_0)_y\cap \phi\inv(U)$ pour tout élément~$y$ de~$\phi\inv(x)$ ; il s'ensuit que si~$b\in \phi\inv(a)$ alors~$$b(U)=b(U_0)\times_{U_0} U=b(U_0)\times_{a(U_0)}a(U).$$

\trois{basesecimrec} Soit~$b\in \phi\inv(a)$ et soit~$y$ l'origine de~$b$. Soit~$\Omega\in \arb {\phi\inv(U_0)_y}y$. Il existe~$U\in \arb {U_0} x$ tel que~$\phi\inv(U)_y\subset \Omega$ ; on a alors~$b(U)\subset \Omega$ ; la composante connexe de~$\Omega\setminus\{y\}$ qui contient~$b(U)$ coïncide nécessairement avec~$b(\Omega)$ (\ref{ouvdefsec}) ; ainsi,~$b(U)\subset b(\Omega)$. Par conséquent,~$\{b(U)\}_{U\in \arb {U_0}x}$ est une base de sections de~$b$.  

\trois{imrecboncomp} Soit~$V$ une section de~$a$ telle que~$\phi\inv(V)$ soit de la forme~$\coprod W_i$ où~$W_i$ est pour tout~$i$ un ouvert de~$Y$ possédant les propriétés suivantes : 

\medskip
$\bullet$~$\overline {W_i}$ contient un unique antécédent~$y_i$ de~$x$ ; 

$\bullet$~$W_i\in \sgb Y {y_i}.$ 

\medskip
Si~$b_i$ désigne pour tout~$i$ la branche issue de~$y_i$ et définie par~$W_i$, alors~$\{b_i\}_i$ coïncide avec~$\phi\inv(a)$. En effet, soit~$U\in \arb {U_0} x$ tel que~$a(U)\subset V$ ; on a alors~$\phi\inv(a(U))=\coprod\limits_{b\in \phi\inv(a)} b(U)$. Pour tout~$i$, le produit fibré~$W_i\times_V a(U)$ est un ouvert fermé de~$\phi\inv(a(U))$, qui est non vide puisque~$W$ se surjecte sur~$V$ ; il s'écrit donc ~$\coprod\limits_{b\in \sch B_i}b(U)$, où~$\sch B_i$ est un sous-ensemble non vide de~$\phi\inv(a)$ ; on a~$\phi\inv(a)=\coprod \sch B_i$. 

\medskip
Fixons~$i$. Si~$b\in \sch B_i$ l'origine de la branche~$b$ adhère à~$b(U)$, et {\em a fortiori} à~$W_i$ ; elle est donc égale à~$y_i$. Les ouverts~$b(U)$ et~$W_i$ appartiennent alors tous deux à~$ \sgb Y {y_i}$, et~$b(U)\subset W_i$. Il s'ensuit qu'ils définissent la même branche issue de~$y_i$ (\ref{ouvdefsec}) , et donc que~$b=b_i$. Ainsi,~$\sch B_i=\{b_i\}$ pour tout~$i$, et l'on a finalement bien~$\phi\inv(a)= \{b_i\}_i$. 

\trois{imdirboncomp} Soit~$y$ un antécédent de~$x$, et soit~$b$ une branche issue de~$y$ appartenant à~$\phi\inv(a)$. Soit~$W$ une section de~$b$ telle que~$\phi(W)\in \sgb X x$ ; la branche définie par~$\phi(W)$ coïncide alors avec~$a$. En effet, on déduit de~\ref{basesecimrec} que~$W$ contient une section~$W'$ de~$b$ telle que~$\phi(W')$ soit une section de~$a$ ; comme~$\phi(W')\subset \phi(W)$, il résulte de~\ref{ouvdefsec} que~$\phi(W)$ est une section de~$a$.

\deux{imdombr} Soit~$X$ un graphe, soit~$Y$ un sous-graphe de~$X$ et soit~$x\in Y$ ; nous allons montrer que~$\br Y x\to \br X x$ est {\em injective}. Soit 
$X'\in \arb X x$ et soit~$Y'\in \arb {Y\cap X'} x$ ; soit~$U$ une composante connexe de~$X'\setminus\{x\}$ , et soient~$V$ et~$W$ deux composantes connexes de~$Y'\setminus\{x\}$ contenues dans~$U$. Choisissons~$v$ dans~$V$ et choisissons~$w$ dans~$W$. Les intervalles~$[v;x[$ et~$[w;x[$ sont respectivement tracés sur~$V$ et~$W$ ; mais ils sont également tous deux tracés sur~$U$, puisque~$v$ et~$w$ appartiennent à~$U$ ; leur intersection est donc non vide (elle est de la forme~$[t;x[$ avec~$t\in U$), ce qui entraîne que~$V=W$ ; ainsi,~$\br Y x$ s'injecte dans~$\br X x$, comme annoncé. 

Pour cette raison, nous identifierons le plus souvent implicitement~$\br Y x$ à un sous-ensemble de~$\br X x$. Modulo cette convention, on peut par exemple écrire que si~$X$ et~$Y$ sont des arbres et si~$b\in \br Y x$ alors~$b(Y)=b(X)\cap Y$. 

\medskip
Notons par ailleurs une conséquence triviale de l'injectivité de la flèche~$\br Y x\to\br X x$ : si~$Y$ est localement fini, le cardinal de~$\br X x$ est minoré par la valence de~$Y$ en~$x$.

\subsection*{Branches et étoiles}

\deux{aboutpropet} Soit~$X$ un graphe et soit~$x\in X$. 

\trois{aboutprop} Nous dirons qu'un intervalle ouvert tracé sur~$X$ aboutit  {\em proprement} à~$x$ si~$x\in \partial I$ et si~$\overline I$ est un arbre ; cela revient à demander que~$\overline I$ soit ou bien un segment dont~$x$ est l'une des extrémités, ou bien une demi-droite dont~$x$ est l'extrémité. Si~$I$ aboutit proprement à~$x$ alors~$I\cup\{x\}$ est un graphe localement fini dont~$x$ est un point unibranche ; lorsqu'on verra l'unique branche de~$I\cup\{x\}$ comme un élément de~$\br X x$,  on dira que c'est la branche {\em définie par~$I$}. 

\trois{notint} L'ensemble des intervalles ouverts tracés sur~$X$ et aboutissant proprement à~$X$ sera noté~$\inter X x$ ; le sous-ensemble de~$\inter X x$ formé des intervalles faiblement admissibles (resp. relativement compacts, resp. faiblement admissibles et relativement compacts) dans~$X$ sera noté~$\intera X x$ (resp.~$\interc X x$, resp.~$\interac X x$). Si~$b\in \br X x$, le sous-ensemble de~$\inter X x$ formé des intervalles définissant une branche~$b$ donnée sera désigné par~$\inter X b$ ; si~$I\in \inter X x$, le sous-ensemble de~$\inter X x$ formé des intervalles contenus dans~$I$ sera noté~$\inter X x \ctd I$.

\medskip
On emploiera également, avec un sens évident, des notations telles~$\intera X b$,~$\intera X x \ctd I$,~$ \inter X b \ctd I$,~$ \intera X b\ctd I$, etc.

\trois{iflatidefb} Soit~$b\in \br X x$ et soit~$I\in \intera X b$. Comme~$\partial I^\flat=\partial I$, lequel contient au plus deux éléments dont~$x$, le point~$x$ est isolé dans~$\partial I^\flat$ ; comme~$\overline I$ est un arbre,~$\overline {I^\flat}$ est un arbre. On déduit alors du lemme~\ref{etresecb} que~$I^\flat\in \sgb X x$ ; comme~$I\subset I^\flat$, la branche définie par~$I^\flat$ est égale à~$b$. 

\deux{defetoile} Une {\em étoile} est la donnée d'un arbre fini~$\Gamma$ et d'un point~$x$ distingué sur~$\Gamma$, que l'on appelle son {\em sommet}, et qui est tel que~$\Gamma\setminus\{x\}$ soit réunion finie disjointe d'intervalles ouverts, que l'on appelle les {\em arêtes} de~$\Gamma$ ; le nombre de ces intervalles est égal à la valence de~$(\Gamma,x)$ ; on l'appellera simplement la valence de~$\Gamma$. 

\medskip
Si~$X$ est un graphe, si~$x\in X$ et si~$\Gamma$ est un sous-arbre de~$X$ contenant~$x$, on dira par abus que~$\Gamma$ {\em est} une étoile de sommet~$x$ si le {\em couple}~$(\Gamma,x)$ définit une étoile de sommet~$x$, c'est-à-dire si et seulement si~$\Gamma\setminus\{x\}$ est réunion finie disjointe d'intervalles ouverts ; on dira aussi d'un tel~$\Gamma$ qu'il est une étoile de sommet~$x$ tracée sur~$X$. 

\trois{exempet} Si~$\Gamma$ est une étoile et si~$x$ est son sommet, la valence de~$\Gamma$ est nulle (resp. égale à~$1$, resp. égale à~$2$) si et seulement si~$\Gamma=\{x\}$ (resp.~$\Gamma$ est un intervalle ouvert, resp.~$\Gamma$ est une demi-droite d'origine~$x$). 

\trois{sometintrins} Si~$\Gamma$ est une étoile et si~$x$ est son sommet, la valence de~$(\Gamma,y)$ est égale à~$2$ pour tout~$y\in \Gamma\setminus\{x\}$ ; si la valence de~$\Gamma$ est différente de~$2$, le point~$x$ ne dépend donc que du graphe localement fini~$\Gamma$ et peut être caractérisé comme son unique sommet {\em topologique}. 

Notons par contre que si~$I$ est un intervalle ouvert et si~$y$ est {\em n'importe quel} point de~$I$ alors~$I$ peut être considéré comme une étoile de sommet~$y$. 

\trois{souset} Si~$\Gamma$ est une étoile et si~$x$ est son sommet, on appellera {\em sous-étoile} de~$\Gamma$ tout voisinage ouvert convexe de~$x$ dans~$\Gamma$ ; les sous-étoiles de~$\Gamma$ sont exactement les étoiles de sommet~$x$ tracées sur~$\Gamma$. 
 
\deux {etbasex} Soit~$X$ un graphe et soit~$x\in X$. On notera~$\stel X x$ l'ensemble des étoiles de sommet~$x$ tracées sur~$X$ dont l'adhérence dans~$X$ est un arbre. Un sous-graphe~$\Gamma$ de~$X$ appartient à~$\stel X x$ si et seulement si il est de la forme~$\{x\}\cup \bigcup\limits_{i\in E}^N I_i$ où~$E$ est un ensemble fini, où chacun des~$I_i$ est un intervalle ouvert aboutissant proprement à~$x$ et où~$\overline {I_i}\cap \overline {I_j}=\{x\}$ dès que~$i\neq j$. Si c'est le cas, les~$I_i$ sont alors les arêtes de~$\Gamma$.

\medskip
Le sous-ensemble de~$\stel X x$ formé des étoiles dont toutes les arêtes sont faiblement admissibles dans~$X$ sera noté~$\stela X x$ ; celui formé des étoiles relativement compactes dans~$X$ sera noté~$\stelc X x$ ; on désignera par~$\stelac X x$ l'intersection de~$\stela X x$ et~$\stelc X x$. 

\trois{netdefbr} Pour tout~$\Gamma\in \stel X x$, le sous-ensemble ~$\br \Gamma x$ de~$\br X x$ est constitué des branches définies par les arêtes de~$\Gamma$. 

\trois{locfinibaseet} Si~$\Delta$ est un sous-graphe de~$X$, on notera~$\stel X x \ctd \Delta$ le sous-ensemble de~$\stel X x$ formé des étoiles contenues dans~$\Delta$ ; on utilisera à l'occasion les notations~$\stela X x \ctd \Delta$, etc., dans un sens évident. 

\medskip
Si~$\Delta$ contient~$x$ et est fini en~$x$, et si~$N$ désigne la valence de~$(\Delta,x)$ alors les étoiles de valence~$N$ appartenant à~$\stel X x \ctd \Delta$ constituent une base de voisinages ouverts de~$x$ dans~$\Delta$.

\deux{etoilebranche} Nous utiliserons librement dans la suite le fait élémentaire suivant : soit~$X$ un graphe, soit~$x\in X$ et soit~$\sch B$ un sous-ensemble fini de~$\br X x$ ; il existe alors~$\Gamma\in \stel X x$ telle que~$\br \Gamma x=\sch B$. En effet, choisissons~$V\in \arb X x$ ; donnons-nous pour chaque~$b\in \sch B$ un point~$x_b$ de~$b(V)$. Pour tout~$b\in \sch B$,  l'intervalle~$[x_b;x[$ est tracé sur~$b(V)$ ; il s'ensuit que~$\Gamma:=\{x\}\cup\bigcup\limits_{b\in \sch B} ]x_b;x[$ est une étoile de sommet~$x$ tracée sur~$X$ et que~$\br \Gamma x=\sch B$.

\deux{unionetbasex} Soit~$X$ un graphe, soit~$x\in X$ et soient~$\Gamma$ et~$\Gamma'$ deux étoiles tracées sur~$X$ de sommet~$x$. 

\trois{existdeltacommun} Il existe~$\Delta\in \stel X x$ possédant les propriétés suivantes : 

\medskip
$\bullet$~$\Delta$ est un ouvert de~$\Gamma\cup \Gamma'$ ; 

$\bullet$~$\br \Delta x=\br \Gamma x\cup \br {\Gamma'}x$. 

$\bullet$ la réunion de~$\{x\}$ et des~$b(\Delta)$ où~$b$ parcourt~$\br \Gamma x$ (resp.~$\br {\Gamma'}x$) est égale à l'ouvert~$\Delta \cap \Gamma$ (resp.~$\Delta\cap \Gamma'$) de~$\Gamma$ (resp.~$\Gamma'$). 

\medskip
Pour le voir, on choisit~$V\in \arb X x$ ; la composante connexe~$\Gamma_0$ (resp.~$\Gamma'_0$) de ~$x$ dans~$V\cap \Gamma$ (resp.~$V\cap \Gamma'$) appartient à~$\stel V x$ ; il existe un voisinage ouvert convexe~$W$ de~$x$ dans~$V$ tel que~$W\cap \Gamma\subset \Gamma_0$ et~$W\cap \Gamma'\subset \Gamma'_0$ ; par convexité,~$W\cap \Gamma$ (resp.~$W\cap \Gamma'$) appartient à~$\stel W x$ ; quitte à remplacer~$X$ par~$W$ et~$\Gamma$ et~$\Gamma'$ par leurs  intersections avec~$W$, il est loisible de supposer que~$X$ est un arbre. 

\medskip
Il est alors clair que l'on peut prendre pour~$\Delta$ la réunion de~$\{x\}$ et de~$$\bigcup_{b\in \br \Gamma x -\br {\Gamma'}x} b(\Gamma)\;\cup\; \bigcup _{b\in\br{\Gamma'}x-\br \Gamma x} b(\Gamma')\;\cup\; \bigcup_{b\in \br \Gamma x\cap \br {\Gamma '} x} b(\Gamma)\cap b(\Gamma')\;.$$ 

\trois{propdeltacommun} Soit~$\Delta\in \stel X x$ satisfaisant les conditions énoncées au~\ref{existdeltacommun} ci-dessus ; toute sous-étoile de~$\Delta$ les satisfait encore. L'ouvert~$\Delta\cap \Gamma\cap \Gamma'$ de~$\Gamma\cap \Gamma'$ est égal à~$\{x\}\bigcup\limits_{b\in \br \Gamma x\cap \br {\Gamma'}x} b(\Delta)$ ; c'est une étoile tracée sur~$X$ dont l'ensemble des branches issues de~$x$ est égal à~$\br \Gamma x \cap \br{\Gamma'}x$. 

\trois{consdeltacommun} Comme~$\Gamma$ et~$\Gamma'$ sont deux sous-graphes localement finis, l'intersection~$\Gamma\cap \Gamma'$ est un graphe localement fini. Comme~$\overline \Gamma$ et~$\overline {\Gamma'}$ sont deux sous-graphes localement finis et {\em fermés} de~$X$, leur réunion est un sous-graphe fermé et localement fini de~$X$. Le complémentaire de~$\Gamma\cup \Gamma'$ dans celui-ci est contenu dans~$\partial \Gamma\cup \partial \Gamma'$, et partant fini ; il est en particulier fermé. Par conséquent,~$\Gamma\cup \Gamma'$ est un ouvert de~$\overline \Gamma\cup \overline {\Gamma'}$ est de ce fait un sous-graphe localement fini de~$X$.   

\medskip
Il découle de~\ref{existdeltacommun} et~\ref{propdeltacommun} que ~$\br {\Gamma\cap \Gamma'}x=\br \Gamma x \cap \br {\Gamma'} x$, et que~$\br {\Gamma \cup\Gamma'}x=\br \Gamma x\cup \br {\Gamma' }x~$; 

\deux{exaradm} Soit~$X$ un graphe, soit~$x\in X$ et soit~$\Gamma$ une étoile de somme~$x$ tracée sur~$X$. Pour tout~$b\in \br \Gamma x$, donnons-nous une section~$U_b$ de~$b$. Nous allons montrer qu'il existe~$\Delta\in \stelac X x\ctd \Gamma$ telle que pour tout~$b\in \br \Gamma x$ l'arête~$b(\Delta)$ de~$\Delta$ soit contenue dans~$U_b$.  

\medskip
Il existe~$V\in \arb X x$ tel que~$b(V)\subset U_b$ pour tout~$b\in \br \Gamma x$ (\ref{finibranchecontr}) ; la composante connexe~$\Gamma_0$ de~$x$ dans~$\Gamma\cap V$ est une sous-étoile de~$\Gamma$ ; il existe un sous-graphe admissible et localement fini~$\Gamma'$ de~$V$ qui contient~$\Gamma_0$. 

Soit~$b\in \br \Gamma x$. Comme~$\Gamma'$ est localement fini, il existe un point~$x_b$ sur~$b(\Gamma_0)$ tel que~$]x_b;x[$ ne contienne aucun sommet de~$\Gamma'$ ; l'intervalle~$]x_b;x[$ est  dès lors ouvert dans~$\Gamma'$, et partant faiblement admissible dans~$V$. 

Comme~$\partial \left(]x_b;x[^\flat\right)=\partial\left(]x_b;x[\right)=\{x_b,x\}$, l'ouvert~$]x_b;x[^\flat$ est inclus dans~$V\setminus\{x\}$ ; comme il est connexe et rencontre~$b(V)$ (puisqu'il contient~$]x_b;x[$), il est contenu dans~$b(V)$, et {\em a fortiori} dans~$U_b$. 

La réunion~$\Delta:=\{x\}\cup\bigcup]x_b;x[$ satisfait alors les conditions requises.

\medskip
{\em Remarque.} Si~$\Delta$ est une sous-étoile de~$\Gamma$ possédant les propriétés énoncées ci-dessus, il en va de même de toute sous-étoile de~$\Delta$.

\deux{etadmcomp} Soit~$X$ un graphe, soit~$x\in X$ et soit~$\Gamma\in \stela X x$. Soit~$V$ la composante connexe de~$x$ dans~$X\setminus \partial  \Gamma$. 

\trois{iflatar} Soit~$I$ une arête de~$\Gamma$ et soit~$b$ la branche qu'elle définit. L'ouvert~$I^\flat$ appartient à~$\sgb X x$ et la branche qu'il définit est précisément~$b$ (\ref{iflatidefb}) ; il est connexe, non vide, et son bord est contenu dans~$\{x\}\cup \partial \Gamma$ ; il est donc fermé dans~$V\setminus\{x\}$ et en est de ce fait une composante connexe. Soit~$J$ une arête de~$\Gamma$ distincte de~$I$ et soit~$\beta$ la branche définie par~$J$ ; comme l'ouvert~$J^\flat$ appartient à~$\sgb X x$ et définit~$\beta$, il diffère de~$I^\flat$ ; les composantes connexes~$I^\flat$ et~$J^\flat$ de~$V\setminus\{x\}$ sont dès lors disjointes. 

\trois{etnarnbr} Soit~$N$ la valence de~$\Gamma$ et soient~$I_1,\ldots,I_N$ ses arêtes ; les~$I_i^\flat$ sont en vertu de ce qui précède des composantes connexes deux à deux distinctes de~$V\setminus\{x\}$. Comme toute composante connexe de~$V\setminus\{x\}$ contient une  branche issue de~$x$ (\ref{compcontient}), les assertions suivantes sont équivalentes : 

\medskip
i) le cardinal de~$\br X x$ est égal à~$N$ ; 

ii)~$V\setminus\{x\}=\coprod I_i^\flat$.

\medskip
Supposons qu'elles soient satisfaites et soit~$W$ une composante connexe de~$V-\Gamma$. Si~$\partial W=\{x\}$ alors~$W$ est une composante connexe de~$V\setminus\{x\}$ qui évite~$\Gamma$ et ne peut donc être l'un des~$I_i^\flat$, ce qui contredit ii) ; par conséquent, il existe~$i$ et un point~$y\in I_i$ tel que~$y\in \partial W$. L'intersection~$W\cap I_i^\flat$ est alors non vide, ouverte et fermée dans~$W$ (puisque~$\partial I_i^\flat\subset \{x\}\cup \partial \Gamma$) ; elle coïncide donc avec~$W$, ce qui signifie que~$W\subset I_i^\flat$. Dès lors,~$W$ est une composante connexe de~$I_i^\flat-I_i$ ; comme~$I_i$ est un sous-graphe admissible de~$I_i^\flat$, la composante~$W$ est un arbre à un bout relativement compact. 

\medskip
Il s'ensuit que le sous-graphe fermé~$\Gamma$ de~$V$ est admissible ; par conséquent,~$\Gamma$ est faiblement admissible et~$V=\Gamma^\flat$. 

\deux{brxcardn} Soit~$X$ un graphe et soit~$x\in X$. Supposons que~$\br X x$ soit fini, et notons~$N$ son cardinal. 

\medskip
Soit~$\Gamma$ un élément de~$\stel X x$ de valence~$N$ et soient~$I_1,\ldots,I_n$ ses arêtes.

\medskip
Les assertions suivantes sont alors équivalentes : 

\medskip
i)~$\Gamma\in \stela X x$, {\em i.e.} les~$I_i$ sont faiblement admissibles ; 

ii)~$\Gamma$ est faiblement admissible. 

\medskip
En effet, i)$\Rightarrow$ ii) est une conséquence du~\ref{etnarnbr} ci-dessus, et ii)$\Rightarrow$i) découle du fait que chacun des~$I_i$ est un ouvert de~$\Gamma$. Si ces conditions sont satisfaites, il découle de~\ref{etnarnbr} que les~$I_i^\flat$ sont exactement les composantes connexes de~$\Gamma^\flat\setminus\{x\}$.

\trois{baseetadm} Soit~$\Gamma$ un élément de~$\stel X x$ de valence~$N$ (ce qui revient à demander que~$\br \Gamma x=\br X x$). Nous allons montrer que~$x$ possède une base de voisinages ouverts de la forme~$\Delta^\flat$, où~$\Delta\in \stelac X x \ctd \Gamma$. 

\medskip
Soit~$V$ un voisinage ouvert de~$x$ ; il contient une sous-étoile~$\Gamma_0$ de~$\Gamma$, qui est une~$x$-étoile de~$V$. On déduit de~\ref{exaradm} qu'il existe une sous-étoile~$\Delta$ de~$\Gamma_0$ dont les arêtes sont faiblement admissibles dans~$V$ ; il découle alors de~\ref{etnarnbr} que~$\Delta$ est faiblement admissible dans~$V$ ; cela équivaut à dire que~$\Delta$ est faiblement admissible dans~$X$ et que~$\Delta^\flat\subset V$, ce qui achève notre preuve. 

\trois{basenbouts} Si~$\Gamma$ est un élément de~$\stela X x$ de valence~$N$, l'ouvert~$\Gamma^\flat$ est un arbre à~$N$ bouts relativement compact dont~$\Gamma$ est un sous-arbre admissible, et est même le squelette si~$N\geq 2$, puisque la valence de~$(\Gamma,y)$ est alors au moins égale à~$2$ pour tout point~$y$ de~$\Gamma$. Ainsi, le point~$x$ possède en particulier une base de voisinages qui sont des arbres à~$N$ bouts. 

\deux{caractfinibrx} Soit~$X$ un graphe et soit~$x\in X$. Nous allons montrer que~$\br X x$ est fini si et seulement si il existe un entier~$n$ tel que~$x$ possède une base de voisinages formée d'arbres à~$n$ bouts, et que si c'est le cas, le cardinal de~$\br X x$ est le plus petit entier~$n$ à satisfaire cette propriété. Compte-tenu du~\ref{basenbouts} ci-dessus, il suffit de vérifier que si~$x$ possède une base de voisinages formée d'arbres à~$n$ bouts, alors~$\br X x$ a pour cardinal au plus~$n$.

On raisonne par contraposition. On suppose donc qu'il existe~$n+1$ éléments distincts~$b_1,\ldots,b_{n+1}$ dans~$\br X x$. En vertu de~\ref{exaradm}, il existe~$\Gamma\in \stela X x$ telle que~$\br \Gamma x=\{b_1,\ldots,b_{n+1}\}$. Si~$V$ désigne la composante connexe de~$x$ dans~$X\setminus \partial  \Gamma$, 
il découle de~\ref{iflatar} que les ouverts~$b_1(\Gamma)^\flat, \ldots, b_{n+1}(\Gamma)^\flat$ sont des composantes connexes deux à deux distinctes de~$V\setminus\{x\}$. 

\medskip
Soit~$W$ un sous-arbre ouvert de~$V$ contenant~$x$ ; nous allons monter que~$W$ a au moins~$n+1$ bouts, ce qui permettra de conclure. Il suffit de s'assurer que~$W\setminus\{x\}$ a au moins~$n+1$ composantes connexes non compactes. Fixons~$i$, et soit~$W_i$ une composante connexe de~$W\setminus\{x\}$ contenue dans~$b_i(\Gamma)^\flat$ (il en existe au moins une car~$W$ rencontre~$b_i(\Gamma)^\flat$, à laquelle~$x$ adhère). Si~$W_i$ était compacte, ce serait un ouvert fermé non vide de~$I_i^\flat$, et donc~$I_i^\flat$ elle-même, qui n'est pas compacte (c'est un arbre à deux bouts) ; par conséquent~$W_i$ n'est pas compacte, ce qui achève la démonstration. 

\deux{declun} Soit~$X$ un graphe et soit~$x\in X$. Les résultats énoncés ci-dessus s'appliquent notamment au cas des étoiles de sommet~$x$ dont la valence~$1$, c'est-à-dire au cas des demi-droites issues de~$x$ ; ils se déclinent alors comme suit. 

\medskip
Si~$b\in \br X x$, il existe~$I\in \inter X b$ (\ref{etoilebranche} ; nous utiliserons librement ce fait dans la suite) ; si~$I$ et~$I'$ sont deux intervalles appartenant à~$\inter X b$, il existe un intervalle~$I''$ appartenant à~$\inter X b$ qui est contenu dans~$I\cap I'$ (\ref{existdeltacommun}) ; et si~$I\in \inter X b$ les ouverts~$J^\flat$, où~$J\in \interac X x\ctd J$, forment une base de sections de~$b$ (\ref{exaradm}).

\medskip
Le point~$x$ est unibranche si et seulement si il est non isolé et possède une base de voisinages ouverts constituée d'arbres à un bout (\ref{caractfinibrx}) ; si c'est le cas et si~$I\in \inter X x$ alors : 

\medskip
$\bullet$ pour tout~$J\in \intera X x \ctd I$ la demi-droite~$J\cup\{x\}$ est faiblement admissible dans~$X$ (\ref{etnarnbr}) ; 

$\bullet$ les arbres à un bout~$(J\cup \{x\})^\flat$ constituent, lorsque l'intervalle~$J$ parcourt ~$\intera X x \ctd I$, une base de voisinages ouverts de~$x$ (\ref{baseetadm}). 

\deux{orbranchegr} Soit~$X$ un graphe, soit~$x\in X$ et soit~$b\in \br X x$. Pour tout intervalle non vide~$I$, 
on note~${\rm Or}(I)$ l'ensemble
des orientations de~$I$. Soient~$I$ et~$I'$ deux intervalles appartenant à~$\inter X b$. D'après
le~\ref{declun} ci-dessus, il existe un intervalle~$I''$ appartenant à~$\inter X b$ qui est contenu dans~$I\cap I'$. 
La famille~$\{{\rm Or}(I)\}_{I\in \inter X b}$ apparaît ainsi comme un système filtrant dont toutes
 les flèches de transition sont des isomorphismes. On note~${\rm Or}(b)$ sa limite (à la fois inductive et projective,
 selon le sens dans lequel on considère les isomorphismes) ; c'est un ensemble à deux éléments appelées {\em orientations}
 de la branche~$b$. Pour tout~$I\in \inter X b$, on dispose d'un isomorphisme naturel~${\rm Or}(I)\simeq {\rm Or}(b)$. 
 
 On dit que~$b$ est orientée si on a choisi une orientation de~$b$ ; on dit que~$b$ est orientée {\em vers~$x$} si on 
 a choisi l'orientation induisant sur tout~$I\in \inter X b$ l'orientation vers~$x$, et
 {\em depuis~$x$}
 sinon.

\subsection*{Branches et bouts}

\deux{xmoinssbr}{\bf Proposition.} {\em Soit~$X$ un graphe et soit~$S$ un sous-ensemble fini de~$X$. Soit~$U$ une composante connexe de~$X\setminus S$, soit~$\sch F$ l'ensemble des fermés de~$X$ contenus dans~$U$ et soit~$\Pi$ l'ensemble~$\lim\limits_{\stackrel {\leftarrow}{Y\in \sch F}} \pi_0(U-Y)$. 

\medskip
i) Pour tout~$x\in S$ l'ensemble~$\br X  x \ctd U$ est fini, et non vide si et seulement si~$x\in \partial U$. 

ii) Si~$Z_b$ est pour tout~$b\in \br X S \ctd U$ une section de~$b$ contenue dans~$U$ alors~$U-\bigcup Z_b$ est un fermé de~$X$.

iii) Pour tout~$Y\in \sch F$, il existe~$Y' \in \sch F$ tel que~$Y\subset Y'$ et tel que~$U-Y'$ s'écrive~$$\coprod_{b\in \br X S \ctd U} Z_b,$$ où~$Z_b$ est pour tout~$b$ une section de~$b$ ; on dispose donc d'une bijection naturelle~$\br X S \ctd U \simeq \Pi$. 

iv) La flèche naturelle~$\Pi\to \got d U$ est injective, et bijective si et seulement si~$U$ est relativement compacte. On dispose donc en vertu de iii) d'une injection~$\br X S\ctd U\hookrightarrow \got d U$ qui est bijective si et seulement si~$U$ est relativement compacte. }

\medskip
{\em Démonstration.} On procède en plusieurs temps. 

\trois{segalborddeu}{\em Réduction au cas où~$S=\partial U$.} Le bord de~$U$ est contenu dans~$S$, et~$U$ est une composante connexe de~$X\setminus \partial  U$. Par ailleurs, si~$x$ est un élément de~$S$ tel qu'il existe une branche de~$X$ issue de~$x$ et contenue dans~$U$ alors~$x\in \partial U$ ; on peut donc, pour établir la proposition, remplacer~$S$ par~$\partial U$, c'est-à-dire supposer que~$S=\partial U$.

\trois{interuvs} {\em Remarques générales et preuve de i)}. Soit~$x$ un point de~$S$ et soit~$V\in \arb X x$ tel que~$V\cap S=\{x\}$. Soit~$X'$ la composante connexe de~$x$ dans~$X-(S\setminus\{x\})$ ; on a~$V\subset X'$ ; comme~$x\in \partial U$, l'ouvert~$U$ est contenu dans~$X'$ et~$U$ est ainsi une composante connexe de~$X'\setminus\{x\}$. Par conséquent~$V\cap U$ est une réunion finie de composantes connexes de~$V\setminus\{x\}$ (\ref{intercompvois}), qui est non vide puisque~$x\in \partial U$ ; les composantes en question sont nécessairement celles qui sont de la forme~$b(V)$ où~$b\in \br X x \ctd U$. On en déduit notamment que~$\br X x \ctd U$ est fini et non vide, d'où l'assertion i).

\trois{fermesux} Soit~$V$ un voisinage ouvert de~$S$. Comme~$S=\partial U$, le fermé~$U-V$ de~$U$ est également fermé dans~$X$. 

\trois{bilanvoiss} L'ensemble fini~$S$ possède une base de voisinages de la forme~$\coprod V_x$, où~$V_x$ est pour tout~$x\in S$ un élément relativement compact et à bord fini  de~$\arb X x$, et où les~$V_x$ sont deux à deux disjoints. Si~$V=\coprod V_x$ est un tel voisinage et si l'on pose~$T=U\cap \partial V$ les faits suivants sont vérifiés : 

\medskip
$\bullet$~$V\cap U=\coprod\limits_{x\in S,b\in \br X x \ctd U}b(V_x)$ (cela résulte de~\ref{interuvs}) ; 

$\bullet$~$U-V$ est fermé dans~$X$ (\ref{fermesux}) ; 

$\bullet$ chacune des~$b(V_x)$ est une composante connexe de~$U\setminus T$, qui n'est pas relativement compacte dans~$U$ puisque~$x\in \partial b(V_x)$ ; 

$\bullet$ le compact~$\sch K_U(T)$ est contenu dans~$U-V$ (c'est une conséquence directe de l'assertion précédente).

\trois{cupzbferme} {\em Preuve de ii}. Choisissons~$V=\coprod V_x$ comme au~\ref{bilanvoiss} ci-dessus, en exigeant de surcroît que~$b(V_x)\subset Z_b$ pour tout~$x$ et toute~$b\in \br X x \ctd U$ (\ref{finibranchecontr}). Le fermé ~$U-\bigcup Z_b$ de~$U$ est alors contenu dans~$U-V$, lequel est fermé dans~$X$ ; par conséquent,~$U-\bigcup Z_b$ est fermé dans~$X$.

\trois{brxsufini} {\em Preuve de iii).} Soit~$Y\in \sch F$. Choisissons~$V=\coprod V_x$ comme au~\ref{bilanvoiss} ci-dessus, en exigeant de surcroît que chacun des~$V_x$ évite~$Y$. Le sous-ensemble~$U-V$ de~$U$ est alors un fermé de~$X$ qui possède la propriété requise. 

\trois{gotdbrxs}{\em Assertion iv) : preuve de l'injectivité.} Soient~$\varpi$ et~$\varpi'$ deux éléments de ~$\Pi$ dont les images dans~$\got d U$ coïncident ; pour tout~$Y\in \sch F$, on notera~$\varpi(Y)$ et~$\varpi'(Y)$ les images respectives de~$\varpi$ et~$\varpi'$ dans~$\pi_0(U-Y)$. 

Soit~$Y\in \sch F$. Choisissons~$V=\coprod V_x$ comme au~\ref{bilanvoiss} ci-dessus, en exigeant de surcroît que chacun des~$V_x$ évite~$Y$ ; posons~$T=U\cap \partial V$ et~$Z=U-V$. Le sous-ensemble~$Z$ de~$U$ est fermé dans~$X$ et contient~$Y$ ; on a~$\partial_UZ=\partial_XZ=T$ et~$\sch K_U(T)\subset Z$. 

\medskip
Les ouverts~$\varpi(Z)$ et~$\varpi'(Z)$ étant deux composantes connexes de~$U-Z$, ce sont deux composantes connexes de~$U\setminus T$ qui ne sont pas contenues dans~$\sch K_U(T)$ ; autrement dit, ce sont deux composantes connexes de~$U\setminus \sch K_U(T)$. Comme les images de~$\varpi$ et~$\varpi'$ dans~$\got d U$ coïncident par hypothèse, on a l'égalité~$\varpi(Z)=\varpi'(Z)$, d'où~$\varpi(Y)=\varpi'(Y)$. Ceci valant pour tout~$Y\in \sch F$, l'application~$\Pi\to \got d U$  est bien injective. 

\trois{gotdbrxrel} {\em Fin de la preuve de iv)}. Si~$U$ est relativement compacte alors~$\sch F$ est l'ensemble des compacts de~$U$, et~$\Pi\to \got d U$ est  clairement bijective. 

\medskip
Supposons maintenant que~$U$ ne soit pas relativement compacte. Choisissons~$V=\coprod V_x$ comme au~\ref{bilanvoiss} ci-dessus ; posons~$Z=U-V$ et~$T=U\cap \partial V$. On a~$\sch K_U(T)\subset Z$. Comme~$U$ n'est pas relativement compact, il ne peut être contenu dans~$\sch K_U(T)\cup \overline V$ ; il existe donc une composante connexe~$W$ de~$U\setminus T$ qui n'est pas relativement compacte dans~$U$ et n'est pas l'une des~$b(V_x)$, ce qui veut dire que~$W\subset Z$. Comme~$W$ n'est pas relativement compacte, il existe~$y\in \wid U$ dont l'image dans~$\pi_0(U\setminus T)$ est égale à~$W$. Si~$y$ était l'image d'un certain élément~$\varpi$ de~$\Pi$, la composante~$\varpi(Z)$ de~$U-Z$ devrait être contenue dans~$W$, ce qui est absurde puisque~$W\subset Z$.~$\Box$ 

\deux{boutconv} On conserve les notations de la proposition~\ref{xmoinssbr} ci-dessus. Soit~$\omega\in \got d U$. Si~$\omega$ est l'image d'une branche~$b\in \br X S \subset U$ et si~$x$ désigne l'origine de~$b$ on dira alors que~$\omega$ {\em converge} vers~$x$.  

\deux{consboutbr} Illustrons la proposition~\ref{xmoinssbr} dans deux cas très simples, qui seront constamment utilisés par la suite.

\trois{contboutbr1} Soit~$X$ un graphe et soit~$U$ un ouvert de~$X$ qui est un graphe connexe relativement compact ayant exactement un bout (le cas le plus utilisé sera celui où~$U$ est un arbre à un bout relativement compact). Comme~$U$ a exactement un bout, il n'est pas compact ; comme il est relativement compact,~$\partial U$ est non vide. La proposition~\ref{ensboutsfini} assure alors que~$\partial U$ est un singleton ; si~$x$ désigne son unique élément, on déduit de la proposition~\ref{xmoinssbr} ci-dessus que~$\br X x\ctd U$ est un singleton. L'unique bout de~$U$ converge vers~$x$. 

\trois{constboutbr2} Soit~$X$ un graphe et soit~$U$ un ouvert de~$X$ qui est un arbre à deux bouts. Son squelette~$\skel U$ est alors un sous-graphe faiblement admissible de~$X$, son adhérence~$\overline{\skel U}$ est égale~$ \skel U\cup \partial U$ et est un sous- graphe admissible de~$\overline U$. Trois cas peuvent se présenter. 

\medskip
\begin{itemize}

\item[$\bullet$] {\em L'ensemble~$\br X {\partial U}\ctd U$ est vide.} Cela se produit si et seulement si~$\partial U$ est vide, et aucun des deux bouts de~$U$ ne converge alors. 

\item[$\bullet$] {\em L'ensemble~$\br X {\partial U}\ctd U$ est un singleton}. Dans ce cas ~$U$ n'est pas relativement compact, et~$\partial U$ est nécessairement lui-même un singleton~$\{x\}$. Comme~$\skel U\cup\{x\}$ est un sous-graphe admissible de~$\overline U$, il est connexe et non compact ; par conséquent, c'est une demi-droite éventuellement longue issue de~$x$, et~$\overline U$ est un arbre à un bout. L'ensemble~$\br X x\ctd U$ est un singleton. L'un des bouts de~$U$ converge vers~$x$, l'autre ne converge pas. 

\medskip
\item[$\bullet$] {\em L'ensemble~$\br X  {\partial U}\ctd U$ a deux éléments.} Dans ce cas~$\overline U$ est compact, et son sous-graphe admissible~$\skel U\cup \partial U$ est compact et connexe.

\medskip
\begin{itemize}
\item[$\diamond$] Si~$\partial U$ a deux éléments~$x$ et~$y$ alors~$\overline U$ s'identifie à la compactification arboricole de~$U$, et~$\skel U\cup\{x,y\}$ est un segment. Chacun des ensembles~$\br X x\ctd U$ et~$\br X y \ctd U$ est un singleton. L'un des bouts de~$U$ converge vers~$x$, l'autre vers~$y$. 

\item[$\diamond$] Si~$\partial U$ est un singleton~$\{x\}$ alors le graphe compact~$\skel U\cup \{x\}$ est un cercle, et~$\overline U$ n'est donc pas un arbre. Comme~$\skel U\cup\{x\}$ est de valence 2 en chacun de ses points, c'est le squelette de~$\overline U$. L'ensemble~$\br X x\ctd U$ a deux éléments, et chacun des deux bouts de~$U$ converge vers ~$x$.  

\end{itemize}

\end{itemize}

\section{Cochaînes harmoniques sur un graphe 
localement
fini}

\subsection*{Définition et premières propriétés}

\deux{comp-conn-orient}
Soit~$X$ une variété topologique séparée purement de dimension~$1$, connexe et non vide. C'est
un cercle ou une droite éventuellement longue, et elle possède donc exactement deux orientations.

Si~$U$ est un ouvert connexe et non vide de~$X$ et si~$o$ est une orientation 
sur~$U$ (resp.~$X$) on se permettra de noter encore~$o$ l'orientation induite
sur~$X$ (resp.~$U$), si cela ne prête pas à confusion.

\deux{pre-def-harm}
Soit~$\Gamma$ un graphe localement fini et soit~$S$ un sous-ensemble
fermé et discret de~$\Gamma$, contenant tous ses sommets. L'ouvert~$\Gamma \setminus S$ 
de~$\Gamma$ est alors une variété topologique purement de dimension~$1$. 
On note~$\mathsf{Or}(\Gamma,S)$ l'ensemble des couples~$(I,o)$ où~$I$ est un ouvert connexe 
non vide de~$\Gamma\setminus S$ et où~$o$ est une orientation sur celui-ci ; on note~$\mathsf{Or}^{\max}(\Gamma,S)$
le sous-ensemble de~~$\mathsf{Or}(\Gamma,S)$
formé des couples~$(I,o)$ tels que~$I$ soit une composante connexe de~$\Gamma\setminus S$.

\deux{def-coch-a}
Soit~$A$ un groupe abélien. Une
{\em $A$-cochaîne}
sur~$(\Gamma,S)$ est 
une application~$\phi : \mathsf{Or}(\Gamma, S)\to A$ satisfaisant les deux axiomes suivants : 

\medskip
$\bullet$ pour tout~$(I,o)\in  \mathsf{Or}(\Gamma, S)$ eut tout ouvert connexe~$J$
de~$\Gamma\setminus S$ contenant~$I$, on a~$\phi(I,o)=\phi(J,o)$ ; 

$\bullet$ pour tout ouvert connexe non vide~$I$ de~$\Gamma\setminus S$ dont on note~$o$
et~$o'$ les deux orientations, on a ~$\phi(I,o)+\phi(I,o')=0$. 

\trois{coch-harm-gp} L'ensemble~$\mathsf{Coch}(\Gamma,S,A)$
des~$A$-cochaînes sur~$(\Gamma,S)$ constituent de façon naturelle
un groupe
abélien, et un~$A$-module si~$A$ est le (groupe abélien sous-jacent à)
un anneau commutatif unitaire. 

\trois{coch-harm-compconn}
On peut également définir une~$A$-cochaîne sur~$(\Gamma,S)$
comme une application~$\phi$
de~$\mathsf{Or}^{\rm max}(\Gamma,S)$
vers~$A$ telle que~$\phi(I,o)+\phi(I,o')=0$ pour toute composante connexe~$I$
de~$\Gamma \setminus S$ dont on note~$o$ et~$o'$ les deux orientations. En effet, une telle~$\phi$ étant
donnée, on vérifie aussiôt qu'elle admet un unique prolongement~$\psi$ à~$\mathsf{Or}(\Gamma,S)$ satisfaisant les axiomes requis ;
si~$(I,o)\in \mathsf{Or}(\Gamma, S)$ et si~$J$ désigne la composante connexe de~$\Gamma\setminus S$
contenant~$I$, alors~$\psi (I,o)=\phi(J,o)$.

\trois{def-somm-en-x}
Soit~$x\in S$ et soit~$\phi$ une~$A$-cochaîne sur~$(\Gamma,S)$. 
Comme~$\Gamma$ est localement fini et comme~$S$ est discret, il existe un voisinage ouvert~$\Delta$
de~$x$ dans~$\Gamma$ qui est une étoile de sommet~$x$
et est tel que~$\Delta\cap S=\{x\}$. 
Soient~$I_1,\ldots, I_n$ les composantes connexes 
de~$\Delta\setminus\{x\}$; pour tout~$s$, on note~$o_i$
l'orientation de~$I_i$ dans la direction de~$x$. 
La somme~$\sum_j \phi(I_i,o_i)$ ne dépend pas de~$\Delta$, et est notée~$\partial_x(\phi)$.

\trois{def-coch-harm}
La
cochaîne~$\phi$ est
dite
{\em harmonique}
si~$\partial_x(\phi)=0$ pour tout
point~$x$
de~$S$. 
L'ensemble des~$A$-cochaînes harmoniques
sur~$(\Gamma,S)$ est un sous-groupe 
de~$\mathsf{Coch}(\Gamma,S,A)$, et un sous-module
de ce dernier lorsque~$A$ est un anneau. 

\deux{coch-harm-restr}
Soit~$T$ un sous-ensemble fermé et discret de~$\Gamma$
contenant~$S$. 
L'ensemble $\mathsf{Or}(\Gamma,T)$ est un sous-ensemble
de~$\mathsf{Or}(\Gamma,S)$, et la restriction induit un morphisme
$\mathsf{Coch}(\Gamma,S,A)\to \mathsf{Coch}(\Gamma,T,A)$.

\trois{inj-harm-restr}
{\em Ce morphisme est injectif.}
En effet, soit~$\phi \in
\mathsf{Coch}(\Gamma,S,A)$ une cochaîne
induisant la cochaîne nulle sur~$(\Gamma,T)$, et soit~$(I,o)\in \mathsf{Or}(\Gamma,S)$. 
Comme~$T$ est discret, il existe un intervalle ouvert non vide~$J$ 
contenu dans~$I$ et ne rencontrant pas~$T$ ; on a alors
$\phi(I,o)=\phi(J,o)=0$, et~$\phi$ est bien nulle. 

\trois{partial-x-indeps} Soit~$\phi \in \mathsf{Coch}(\Gamma,S,A)$ et soit~$\psi$
la cochaîne induite sur~$(\Gamma,T)$. Soit~$x\in S$. 
Soit~$\phi$ une~$A$-cochaîne sur~$(\Gamma,S)$ et soit~$\psi$
la~$A$-cochaîne induite sur~$(\Gamma, T)$. Soit~$x\in S$. Choisissons un voisinage
ouvert~$\Delta$ de~$x$ dans~$\Gamma$ qui est une étoile telle que~$\Delta\cap T=\{x\}$. On a
{\em a  fortiori}
$\Delta\cap S=\{x\}$, et il résulte alors immédiatement des
définitions que~$\partial_x(\psi)=\partial_x(\phi)$. 

\deux{coch-harm-indep}
Nous allons montrer que
$\mathsf{Coch}(\Gamma,S,A)\hookrightarrow \mathsf{Coch}(\Gamma,T,A)$
induit un {\em isomorphisme} entre les groupes de
$A$-cochaînes harmoniques 
sur~$(\Gamma,S)$ et~$(\Gamma, T)$. 

\trois{iso-coch-1}
Soit~$\phi$ une~$A$-cochaîne harmonique sur~$(\Gamma,S)$,
et soit~$\psi$ son image dans~$\mathsf{Hom}(\mathsf{Or}(\Gamma,T),A)$. 
Nous allons vérifier que~$\psi$ est une cochaîne harmonique sur~$(\Gamma,T)$. Soit~$x\in T$. 

\medskip
{\em Supposons que~$x\in S$.} On a d'après~\ref{partial-x-indeps}
l'égalité~$\partial_x(\psi)=\partial_x(\phi)=0$. 

\medskip
{\em Supposons que~$x\notin S$.}
Le graphe~$\Gamma$
est de valence~$2$ en~$x$ ; il existe donc un voisinage ouvert~$I$
de~$x$ dans~$\Gamma$ qui est un intervalle 
dont l'intersection avec~$T$ est égale à~$\{x\}$. 
Soient~$I'$ et~$I''$
les deux composantes connexes de~$I\setminus \{x\}$ ; on note $o'$ (resp.~$o''$)
l'orientation de~$I'$ (resp.~$I''$)
dans la direction de~$x$. 
On a alors 

$$\partial_x(\psi)=\psi(I',o')+\psi(I'',o'')=\phi(I',o')+\phi(I'',o'')=\phi(I,o)+\phi(I,o'')=0$$
puisque les deux orientations~$o'$
et~$o''$ de~$I$ sont opposées. Ainsi, $\psi$ est harmonique. 

\trois{iso-coch-2}
Soit~$\psi$ une~$A$-cochaîne harmonique sur~$(\Gamma, T)$, soit~$(I,o)\in \mathsf{Or}(\Gamma,S)$, et soit~$J$ un
ouvert connexe et non vide
de~$I$ ne rencontrant pas~$T$. Nous allons montrer que~$\psi(J,o)$ ne dépend
que de~$(I,o)$, et pas du choix de~$J$.

\medskip
Pour ce faire, on choisit un second ouvert connexe et non vide~$J'$ 
de~$I$ ne rencontrant pas~$T$. On se donne un point~$a$ sur~$J$, un point~$b$
sur~$J'$, et un segment joignant~$a$ à~$b$ que  l'on se permet de noter~$[a;b]$ (même si~$I$ est un cercle). 
L'intervalle~$J$ contient un intervalle
ouvert non vide inclus dans~$]a;b[$, et de même pour~$J'$ ; on est donc ramené au cas où~$J$
et~$J'$ sont contenus dans~$]a;b[$ puis, quitte à les agrandir, à celui où chacun d'eux est une composante connexe
de~$]a;b[\setminus T$. 

\medskip
L'intersection~$]a;b[\cap T$ est un ensemble fini~$x_1<x_2<\ldots x_n$. Il suffit dès lors
de démontrer que pour tout~$i$ compris entre~$1$ et~$n$, 
on a
$$\psi(]x_{i-1};x_i[,o)=\psi(]x_i;x_{i+1}[,o),$$ avec les conventions~$x_0=a$ et~$x_n=b$. 
Mais c'est une conséquence immédiate du fait que~$\delta_{x_i}(\psi)=0$. 

\medskip
Il est en conséquence
licite de poser~$\phi(I,o)=\psi(J,o)$ où~$J$
est {\em n'importe quel}
ouvert connexe et non vide de~$I$ ne rencontrant
par~$T$. Par construction, $\phi$ est une~$A$-cochaîne sur~$(\Gamma,S)$, qui
induit la~$A$-cochaîne~$\psi$ sur~$(\Gamma,T)$. Comme~$\psi$ est harmonique, 
on déduit de~\ref{partial-x-indeps}
que~$\phi$ est harmonique, ce qui achève la démonstration. 

\deux{def-harm}
Le groupe des~$A$-cochaînes harmoniques sur~$(\Gamma,S)$ est donc canoniquement
indépendant de~$S$ ; on le note~$\mathsf{Harm}(\Gamma,A)$, et on dira simplement de ses éléments
que ce sont des $A$-cochaînes harmoniques
{\em sur~$\Gamma$.} 

\subsection*{Interprétation
cohomologique}
\deux{lien-coh-compact}
Supposons
que~$A$ est un anneau
commutatif unitaire, et soit~$\lambda$ une
application~$A$-linéaire de~$\H^1_{\rm c}(\Gamma, A)$ vers~$A$.

Soit~$S_0$
l'ensemble des sommets de~$\Gamma$, et soit~$I$ 
un ouvert connexe et non vide
de~$\Gamma\setminus _0$. 
Soit~$u : I\hookrightarrow \Gamma$
l'immersion ouverte, 
et soit~$o$ une orientation sur~$I$. Le choix de~$o$ fournit une
{\em classe fondamentale}~$h$
appartenant à~$\H^1_{\rm c}(I,A)$ ; posons~$\phi(I,o)=\lambda(u_*h)$. 
L'application~$\phi$ est alors une~$A$-cochaîne harmonique sur~$(\Gamma,S_0)$,
et l'on définit ainsi un morphisme de~$A$-modules
$$ \mathsf{Hom}_A(\H^1_{\rm c}(\Gamma,A), A)\to \mathsf{Harm}(\Gamma,A).$$

\deux{iso-coho-harm}
On munit~$A$ de la topologie discrète, et on fait de~$\mathsf{Hom}_A(\H^1_{\rm c}(\Gamma,A), A)$
et~$\mathsf{Harm}(\Gamma,A)$
des~$A$-modules topologiques comme suit.

\medskip
$\bullet$ On considère~$\mathsf{Hom}_A(\H^1_{\rm c}(\Gamma,A), A)$ 
comme la limite projective des~$A$-modules discrets~$\mathsf{Hom}_A(M,A)$, où~$M$ parcourt l'ensemble
des sous-modules de type fini de~$\H^1_{\rm c}(\Gamma,A)$. 

\medskip
$\bullet$ On munit~$\mathsf{Harm}(\Gamma,A)$
de la topologie induite par la topologie
produit de~$A^{\mathsf{Or}^{\max}(\Gamma,S_0)}$. 

\medskip
La flèche~$\mathsf{Hom}_A(\H^1_{\rm c}(\Gamma,A), A)\to \mathsf{Harm}(\Gamma,A)$
construite au~\ref{lien-coh-compact}
ci-dessus est alors un {\em isomorphisme
de~$A$-modules topologiques.}

\section{Quotients d'un graphe}

Dans tout ce paragraphe, on désigne par~$\mathsf G$ un groupe {\em profini}.

\deux{defaret} Si~$\Gamma$ est un arbre localement fini et si~$S$ est l'ensemble de ses sommets, on appellera {\em arête} de~$\Gamma$ toute composante connexe de~$\Gamma\setminus S$. Toute arête de~$\Gamma$ est une droite éventuellement longue.

Si~$\phi$ est une application dont la source~$\Gamma$ est un graphe localement fini, on dira que~$\phi$ est {\em injective par morceaux} s'il existe un sous-ensemble fermé et discret~$S$ de~$\Gamma$ tel que la restriction de~$\phi$ à toute composante connexe de~$\Gamma\setminus S$ soit injective.

\subsection*{Action de~$\mathsf G$ sur un graphe : premières propriétés}

\deux{remgenactggraph} Commençons par quelques remarques générales. Soient~$X$ et~$Y$ deux graphes et soit~$\phi$ un homéomorphisme~$Y\simeq X$. Soit~$\Delta$ un sous-graphe de~$Y$ et soit~$y\in Y$ ; les propriétés qui suivent sont conséquences immédiates des définitions. 

\medskip
$\bullet$~$\phi(\skel Y) =\skel X$ ;

$\bullet$~$\phi(\Delta)$ est admissible si et seulement si~$\Delta$ est admissible et si c'est le cas,~$\phi$ commute aux rétractions canoniques respectives de~$Y$ sur~$\Delta$ et~$X$ sur~$\phi(\Delta)$ ; 

$\bullet$~$\phi(\Delta)$ est faiblement admissible si et seulement si~$\Delta$ est faiblement admissible et si c'est le cas,~$\phi(\Delta^\flat)=\phi(\Delta)^\flat$.

\medskip
De plus, si~$y\in Y$ et si~$x=\phi(y)$ alors~$\phi$ induit une bijection~$\br Y y \simeq \br X x$ ; cela résulte de la définition même des branches, et peut également être vu également comme un cas particulier des constructions de~\ref{imdir} {\em et sq.}. Si~$y$ est situé sur~$\Delta$, on vérifie aussitôt que les bijections~$\br Y  y\simeq \br X x$ et~$\br \Delta y \simeq \br {\phi(\Delta)} x$ sont compatibles aux inclusions~$\br \Delta y\subset \br Y y~$ et~$\br{ \phi(\Delta) }x\subset \br X x$.

\medskip
Nous aurons surtout l'occasion d'appliquer les remarques qui précèdent en travaillant avec un graphe~$X$ muni d'une action de~$\mathsf G$, sur lequel on considèrera des automorphismes de la forme~$x\mapsto g(x)$ avec~$g\in \mathsf G$ ; elles ont notamment, dans ce contexte, les conséquences suivantes : 

\medskip
$\bullet$ si~$\Delta$ est un sous-graphe faiblement admissible de~$X$ stable sous~$\mathsf G$ alors~$\Delta^\flat$ est stable sous~$\mathsf G$ ; 

$\bullet$ si~$x\in X$, le stabilisateur~$\mathsf H$ de~$x$ agit sur~$\br X x$.

\deux{prolactprof} {\bf Proposition.} {\em Soit~$X$ un graphe connexe sur lequel~$\mathsf G$ agit. L'action de~$\mathsf G$ sur~$X$ s'étend d'une unique manière en une action de~$\mathsf G$ sur~$\wid X$.}

\medskip
{\em Démonstration.} Il résulte de la définition de~$\wid X$ et de la densité de~$X$ dans~$\wid X$ que tout automorphisme de~$X$ s'étend de manière unique en un automorphisme de~$\wid X$. Par conséquent, l'action de~$\mathsf G$ sur~$X$ s'étend d'une unique manière en une opération {\em ensembliste} de~$\mathsf G$ sur~$\wid X$ qui est continue par rapport à la seconde variable ; il reste à s'assurer que~$\mathsf G\times \wid X\to \wid X$ est continue. 

\medskip
Soit~$S$ un sous-ensemble fini de~$X$ et soit~$U$ une composante connexe de~$X\setminus S$. Soit~$(g,x)$ un élément de~$\mathsf G\times \wid X$ tel que~$g(x)\in \widd U X$ ; nous allons établir l'existence d'un voisinage de~$(g,x)$ dans~$\mathsf G\times \wid X$ dont l'image est contenue dans~$\widd U X$, ce qui permettra de conclure. 

\trois{casxinx} Si~$x\in X$ alors~$g(x)\in \widd U X\cap X=U$, et il résulte de la continuité de~$G\times X\to X$ qu'il existe un voisinage ouvert de~$(g,x)$ dans l'ouvert~$\mathsf G\times X$ de~$\mathsf G\times \wid X$ dont l'image est contenue dans~$U$, et {\em a fortiori} dans~$\widd U X$.

\trois{casxnotinx} Supposons maintenant que~$x\in \got d X$ ; choisissons un voisinage compact~$\mathsf G_0$ de~$g$ dans~$\mathsf G$. Le sous-ensemble~$\mathsf G_0\inv.S$ de~$X$ en est une partie compacte ; il existe donc un sous-ensemble fini~$T$ de~$X$ tel que~$\mathsf G_0\inv.S\subset \sch K_X(T)$. Si~$h\in \mathsf G_0$ on a~$S\subset h(\sch K_X(T))=\sch K_X(h(T))$.

\medskip
Comme~$x\in \got d X$ il appartient à~$\widd V X$ pour une certaine composante connexe~$V$ de~$X\setminus \sch K_X(T)$. Si~$h$ appartient à~$\mathsf G_0$ alors~$h(V)$ est une composante connexe de~$X\setminus \sch K_X(h(T))$, et est donc contenue dans une composante connexe de~$X\setminus S$. 

En particulier,~$g(V)$ est contenu dans une composante connexe~$W$ de~$X\setminus S$ ; il s'ensuit que~$g(\widd V X)\subset \widd W X$ ; comme~$g(x)\in \widd U X$, on a~$W=U$ et donc~$g(V)\subset U$.

\medskip
Choisissons~$y\in V$. Par continuité de~$\mathsf G\times X\to X$, il existe un voisinage ouvert~$\mathsf H$ de~$g$ dans~$\mathsf G_0$ tel que~$\mathsf H.y\subset U$. Si~$h\in \mathsf H$ alors~$h(V)$ rencontre~$U$ par construction ; par conséquent, la composante connexe de~$X\setminus \sch K_X(S)$ contenant~$h(V)$ n'est autre que~$U$ ; il s'ensuit que~$h(\widd V X)\subset \widd U X$. 

\medskip
Ainsi, l'image réciproque de~$\widd U X$ par la flèche~$\mathsf G\times \wid X\to \wid X$ contient le voisinage~$\mathsf H\times \widd V X$ de~$(g,x)$, ce qui achève la démonstration.~$\Box$ 

\deux{lemmeprof} {\bf Lemme.} {\em Soit~$X$ un arbre sur lequel~$\mathsf G$ agit. Soit~$g\in \mathsf G$ et soit~$x$ un point de~$\wid X$ tel que~$g(x)=x$ ; munissons~$\wid X$ de l'ordre partiel défini par~$x$ (rem.~\ref{remordrelong}). Si~$y\in \wid X$ est tel que~$y$ et~$g(y)$ soient comparables alors~$g$ fixe~$[y;x]$ point par point.}

\medskip
{\em Démonstration.} Comme~$g(x)=x$ l'automorphisme~$g$ est croissant. 

\trois{gyegaly} {\em Montrons tout d'abord que~$g(y)=y$.} Quitte à remplacer~$g$ par~$g\inv$ et~$y$ par~$g(y)$, on peut supposer que~$y\leq g(y)$. Par une récurrence immédiate~$g^n(y)\geq y$ pour tout~$n$. Comme~$\mathsf G$ est profini,~$g^{n ! }$ tend vers l'élément neutre~$e$ de~$\mathsf G$ quand~$n$ tend vers l'infini ; par conséquent,~$g^{n!-1}(y)\to g\inv(y)$ quand~$n$ tend vers l'infini.

Par ce qui précède,~$g^{n!-1}(y)\geq y$ pour tout~$n\geq 2$ ; on en déduit par passage à la limite (licite du fait que~$\wid X _{\leq z}~$ est fermé pour tout~$z\in \wid X$, {\em cf.} rem.~\ref{remordrelong}) que~$g\inv(y)\geq y$ ; en appliquant~$g$, il vient~$y\geq g(y)$ ; comme par ailleurs~$y\leq g(y)$, on a bien~$g(y)=y$.

\trois{gzegalz} {\em Conclusion}. Comme~$g(y)=y$ on a~$g([y;x])=[y;x]$. Soit~$z\in [y;x]$. Par ce qui précède,~$g(z)$ appartient à~$[y;x]$ et est de ce fait comparable à~$z$ ; on déduit alors du~\ref{gyegaly}, appliqué à~$z$ au lieu de~$y$, que~$g(z)=z$.~$\Box$

\deux{corolcroisdl} {\bf Corollaire.} {\em Soit~$I$  une droite éventuellement longue sur laquelle~$\mathsf G$ agit. Si~$g$ est un élément de~$\mathsf G$ qui ne permute pas les orientations de~$I$ alors~$g$ agit trivialement sur~$I$.}

\medskip
{\em Démonstration.} Soient~$x$ et~$y$ les deux bouts de~$I$. Dire que~$g$ n'échange pas les orientations de~$I$ signifie que~$g(x)=x$  (ou encore que~$g(y)=y$). Il n'y a plus alors qu'à appliquer le lemme précédent en remarquant que l'ordre induit par~$x$ sur~$\wid I$ est total.~$\Box$

\deux{corollprof} {\bf Corollaire.}  {\em Soit~$I$ une droite éventuellement longue. Si~$\mathsf G$ agit sur~$I$ de façon non triviale, il existe~$x\in I$ et une  symétrie~$\sigma$ de centre~$x$ telle que l'action de~$\mathsf G$ se factorise par une surjection continue~$\mathsf G\to \{\mathsf{Id},\sigma\}$.}

\medskip
{\em Démonstration.} Donnons-nous une action non triviale de~$\mathsf G$ sur~$I$ ; par hypothèse, il existe~$g\in \mathsf G$ et~$y\in I$ tels que~$g(y)\neq y$. L'action de~$g^2$ respecte nécessairement les orientations de~$I$ (c'est le cas pour le carré de n'importe quel homéomorphisme de~$I$) ; il s'ensuit, en vertu du corollaire~\ref{corolcroisdl} ci-dessus, que l'action de~$g^2$ est triviale ; on a en particulier~$g^2(y)=y$, et~$g$ stabilise donc le segment~$[y;g(y)]$. Il induit de ce fait un automorphisme de~$[y;g(y)]$ qui en permute les extrémités, et possède donc un point fixe~$x$ sur~$]y;g(y)[$. Si~$I_1$ (resp.~$I_2$) désigne la composante connexes de~$I\setminus\{x\}$ contenant~$y$ (resp.~$g(y)$) alors~$g$ induit un homéomorphisme~$\tau : I_1\simeq I_2$, dont la réciproque est également induite par~$g$, l'action de~$g^2$ étant triviale. L'application~$\sigma$ de~$I$ dans~$I$ qui fixe~$x$ et envoie~$z$ sur~$\tau(z)$ (resp.~$\tau\inv(z)$) si~$z\in I_1$ (resp.~$I_2$) est une symétrie de centre~$x$, et~$g(z)=\sigma(z)$ pour tout~$z\in I$. 

\medskip
Si~$h$ est maintenant un élément quelconque de~$\mathsf G$ alors ou bien~$h$ agit en préservant les orientations de~$I$, c'est-à-dire trivialement en vertu du lemme~\ref{lemmeprof} ci-dessus, ou bien il les inverse ; mais dans ce dernier cas~$hg$ les préserve et agit  trivialement, ce qui signifie que~$h(z)=g\inv(z)=\sigma(z)$ pour tout~$z\in I$.~$\Box$

\deux{preservor} Si~$\Gamma$ est un arbre localement fini sur lequel~$\mathsf G$ agit, le groupe~$\mathsf G$ agit sur l'ensemble des arêtes de~$\Gamma$. On dira que~$\mathsf G$ {\em préserve les orientations des arêtes de~$\Gamma$} si pour toute arête~$I$ de~$\Gamma$ et tout élément~$g$ du stabilisateur de~$I$, l'action de~$g$ sur~$I$ préserve les orientations de~$I$ -- ce qui revient à demander qu'elle soit triviale (cor.~\ref{corolcroisdl}). 

\deux{quotarbrefin} {\bf Lemme.} {\em Soit~$\Gamma$ un arbre compact, fini et non vide sur lequel~$\mathsf G$ agit. Il existe un point de~$\Gamma$ qui est fixe sous~$\mathsf G$ ; si~$\mathsf G$ agit sans permuter les orientations des arêtes de~$\Gamma$ ce point peut être choisi parmi les sommets de~$\Gamma$.}

\medskip
{\em Démonstration.} On procède en deux temps.

\trois{ptfixegen} {\em Montrons que~$\mathsf G$ possède un point fixe sur~$\Gamma$.} On procède par récurrence sur le nombre~$N$ de sommets de~$\Gamma$. 

\medskip
Si ~$N=1$ alors~$\Gamma$ est réduit à un point, et l'assertion est trivialement vraie. 

Si~$N=2$ alors~$\Gamma$ est un segment (non trivial), et l'assertion se déduit alors immédiatement du corollaire~\ref{corollprof}.

On suppose maintenant que~$N\geq 3$ et que l'assertion est vraie en rang inférieur ou égal à~$N-1$. Soit~$S$ l'ensemble des sommets unibranches de~$\Gamma$. Si~$x\in S$ il existe une unique arête de~$\Gamma$ issue de~$x$ ; l'autre extrémité de cette arête sera notée~$\rho(x)$. C'est un sommet de~$\Gamma$ ; comme~$N\geq 3$ le graphe~$\Gamma$ n'est pas un segment, et la valence de~$\rho(x)$ est donc au moins 3 ; en particulier,~$\rho(x)\notin S$. 

\medskip
Soit~$\Gamma'$ le sous-ensemble de~$\Gamma$ égal à~$\Gamma-\bigcup\limits_{x\in S}[x;\rho(x)[$ ; c'est un sous-arbre compact et non vide de~$\Gamma$ qui est stable par~$\mathsf G$ et possède strictement moins de~$N$ sommets. En vertu de l'hypothèse de récurrence, il existe un point de~$\Gamma'$ qui est fixe sous~$\mathsf G$, ce qui achève de prouver notre assertion.

\trois{ptfixecroiss} {\em Montrons que si~$\mathsf G$ ne permute pas les orientations des arêtes de~$\Gamma$ il existe un sommet de~$\Gamma$ qui est fixe sous~$\mathsf G$.} D'après le~\ref{ptfixegen} ci-dessus,~$\mathsf G$ possède un point fixe~$x$ sur~$\Gamma$. Si~$x$ est un sommet, la démonstration est terminée. Sinon, l'arête~$I$ de~$\Gamma$ qui contient~$x$ est stable sous~$\mathsf G$ ; comme~$\mathsf G$ ne permute pas les orientations des arêtes de~$\Gamma$, l'action de~$\mathsf G$ sur~$I$ est triviale et si~$y$ désigne l'un des deux points de~$\partial I$ alors~$y$ est un sommet de~$\Gamma$ qui est fixe sous~$\mathsf G$.~$\Box$ 

\deux{orbfinlocfin} {\bf Lemme.} {\em Soit~$\Gamma$ un graphe localement fini sur lequel~$\mathsf G$ agit. Les orbites de~$\mathsf G$ sur~$\Gamma$ sont finies.}

\medskip
{\em Démonstration.} Soit~$S$ l'ensemble des sommets de~$\Gamma$. 

\trois{orbfindl} Si~$\Gamma$ est une droite éventuellement longue, la finitude des orbites de~$\mathsf G$ est une conséquence immédiate du corollaire~\ref{corollprof}.

\trois{orbfinxs} Si~$x\in S$, son orbite est un compact contenu dans l'ensemble des sommets de~$\Gamma$, lequel est discret ; par conséquent,~$\mathsf G.x$ est fini.

\trois{orbfinpass} Il reste à s'assurer que les orbites de~$\mathsf G$ sur~$\Gamma\setminus S$ sont finies. Une telle orbite est un compact contenu dans~$\coprod\limits_{I\in \pi_0(\Gamma\setminus S)} I$ et ne rencontre donc qu'un nombre fini de composantes connexes.  Il suffit dès lors de démontrer que si~$I$ est une composante connexe de~$\Gamma\setminus S$ et si~$\mathsf H$ est son stabilisateur, les orbites de~$\mathsf H$ sur~$I$ sont  finies ; autrement dit on s'est ramené, quitte à remplacer~$\mathsf G$ par~$\mathsf H$ et~$\Gamma$ par~$I$, au cas où~$\Gamma$ est connexe, non vide, et sans sommets. On distingue deux cas. 

\medskip
{\em Premier cas :~$\Gamma$ n'est pas compact.} C'est alors une droite éventuellement longue, et la conclusion découle de~\ref{orbfindl}. 

\medskip
{\em Second cas :~$\Gamma$ est compact.} C'est alors un cercle. Soit~$x\in \Gamma$. L'orbite~$\mathsf G.x$ étant totalement discontinue, elle n'est pas égale à~$\Gamma$ tout entier ; choisissons un point~$y$ sur~$\Gamma$ qui ne rencontre pas~$\mathsf G.x$ ; le point~$x$ n'appartient alors pas à~$\mathsf G.y$. L'orbite de~$x$ est un compact contenu dans~$\coprod \limits_{I\in \pi_0(\Gamma-\mathsf G.y)} I$ ; elle ne rencontre donc qu'un nombre fini de composantes connexes de~$\Gamma-\mathsf G.y$ ; chacune de ces composantes est un ouvert connexe et non vide de~$\Gamma\setminus\{y\}$, et est de ce fait un intervalle ouvert. La finitude de~$\mathsf G.x$ découle alors de~\ref{orbfindl}.~$\Box$

\subsection*{Quotient d'un arbre sous l'action de~$\mathsf G$}

\deux{propprel} {\bf Proposition.} {\em Soit~$X$ un arbre sur lequel~$\mathsf G$ agit et soient~$x$ et~$y$ deux points de~$X$ ; on désigne par~$p$ la flèche quotient~$X\to X/\mathsf G$. 

\medskip
i) L'image~$p([x;y])$ est un arbre compact et fini ayant au plus trois arêtes, et~$p_{|[x;y]}$ est injective par morceaux. 

ii) Il existe un unique fermé de~$X/\mathsf G$ homéomorphe à un segment d'extrémités~$p(x)$ et~$p(y)$.}

\medskip
{\em Démonstration.} Notons tout d'abord les faits suivants, que nous utiliserons implicitement tout au long de la preuve : l'espace quotient~$X/\mathsf G$ est séparé et l'application~$p$ est compacte (\ref{quotsep} {\em et sq.}). 

\medskip
Si~$x=y$ la proposition est évidente ; on suppose à partir de maintenant que~$x\neq y$. Soit~$Z$ l'enveloppe convexe du compact~$\mathsf G.x\cup \mathsf G.y$ de~$X$ ; c'est un sous-arbre compact de~$X$ ; l'ouvert~$Z':=Z-\left(\mathsf G.x\cup \mathsf G.y\right)$ de~$Z$ est un graphe localement fini.

\medskip
Comme~$x\neq y$ le segment~$[x;y]$ n'est pas réduit à un singleton ; il ne peut donc être entièrement contenu dans le compact totalement discontinu~$\mathsf G.x\cup \mathsf G.y$, et il s'ensuit que~$Z'$ est non vide. Choisissons un point~$z$ sur~$Z'$. Son orbite est finie (lemme~\ref{orbfinlocfin}) ; la réunion des segments~$[w;w']$, où~$(w,w')$ parcourt~$(\mathsf G.z)^2$, est un sous-arbre non vide, compact et fini de~$Z$ qui est stable sous l'action de~$\mathsf G$. Il possède donc en vertu du lemme~\ref{quotarbrefin} un point~$t$ invariant sous~$\mathsf G$. 

\trois{ytyprime} Nous allons montrer que~$p_{|[y;t]\cup[x;t]}$ est injective par morceaux et que~$p([y;t]\cup[x;t])$ est un arbre compact et fini ayant au plus trois arêtes ; cela entraînera immédiatement  l'assertion i). 

\medskip
{\em Les applications~$p_{|[x;t]}$ et~$p_{|[y;t]}$ sont injectives.} En effet il suffit pour le voir, les points~$x$ et~$y$ et jouant des rôles analogues, de montrer que~$p_{|[x;t]}$ est injective. Soient~$x_1$ et~$x_2$ deux points de~$[x;t]$ ayant même image par~$p$. Cela signifie qu'il existe~$g\in \mathsf G$ tel que~$g(x_1)=x_2$. Étant tous deux situés sur~$[x;t]$, les points~$x_1$ et~$x_2$ sont comparables pour l'ordre défini par~$t$, lequel est fixe sous~$\mathsf G$ ; le lemme~\ref{lemmeprof} assure alors que~$x_1=g(x_1)=x_2$, et~$p_{|[x;t]}$ est bien injective. 

\medskip
{\em L'image~$p([x;t]\cup[y;t])$ est un arbre compact fini ayant au plus trois arêtes}. Par ce qui précède,~$p([x;t])$ (resp.~$p([y;t])$ est un segment d'extrémités~$p(x)$ et~$p(t)$ (resp.~$p(y)$ et~$p(t)$) ; posons~$I=p([x;t])\cap p([y;t])$. Nous allons démontrer que~$I$ est un segment dont~$p(t)$ est l'une des extrémités, ce qui entraînera aussitôt que~$p([x;t]\cup[y;t])=p([x;t])\cup p([y;t])$ est un arbre fini et compact ayant au plus trois arêtes.

\medskip
On sait déjà que~$I$ est compact ; pour établir qu'il s'agit d'un segment dont~$p(t)$ est l'une des extrémités, il suffit dès lors de vérifier que si~$\omega \in I$ alors le segment~$[\omega;p(t)]$ de~$p([x;t])$ est contenu dans~$I$. 

Soit donc~$\omega \in I$ ; il existe un (unique) point~$\xi\in [x;t]$ tel que~$\omega=p(\xi)$ et un (unique) point~$\eta\in [y;t]$ tel que~$\omega=p(\eta)$. Les points~$\xi$ et~$\eta$ ayant même image par~$p$, il existe~$g\in \mathsf G$ tel que~$\xi=g(\eta)$. Le point~$t$ étant invariant sous~$\mathsf G$ on a alors~$g([\eta;t])=[\xi;t]$ ; par conséquent,~$p([\xi;t])$ est contenu dans~$I$. Or par injectivité de~$p_{|[x;t]}$ l'image~$p([\xi;t])$ est exactement le segment~$[\omega;p(t)]$ de~$p([x;t])$, ce qui achève la preuve de l'assertion requise. 

\trois{uniquesegmquot} Nous allons maintenant montrer ii). Commençons par une remarque : le point~$t$ étant fixe sous~$\mathsf G$, le compact~$Z$ peut s'écrire comme la réunion des~$[g(x);t]\cup [g(y);t]$ où~$g$ parcourt~$\mathsf G$ ; par conséquent,~$p(Z)$ est égal à~$p([x;t]\cup[y;t])$, et est en particulier un arbre. Pour établir ii), il suffit dès lors de vérifier que tout segment d'extrémités~$p(x)$ et~$p(y)$ tracé sur~$X/\mathsf G$ est  contenu dans~$p(Z)$. 

\medskip
Soit~$U$ une composante connexe de~$X/\mathsf G-p(Z)$ et soit~$V$ une composante connexe de~$p\inv(U)$ ; l'ouvert~$V$ de~$X$ est une composante connexe de~$X-Z$ ; comme~$Z$ est un arbre,~$\partial V$ est un singleton~$\{v\}$ pour un certain~$v\in Z$. La projection~$p$ étant fermée,~$p(\overline V)=p(V)\cup\{p(v)\}=U\cup \{p(v)\}$ est un fermé de~$X/\mathsf G$ ; par conséquent,~$\partial U$ est le singleton~$\{p(v)\}$. 

Il s'ensuit, en vertu de~\ref{remskgrad}, que tout segment tracé sur~$X/\mathsf G$ et joignant deux points de~$p(Z)$ est contenu dans~$p(Z)$, ce qui achève la démonstration.~$\Box$ 

\deux{theoquot} {\bf Théorème.} {\em Soit~$X$ un graphe sur lequel~$\mathsf G$ agit et soit~$p:X\to X/\mathsf G$ la flèche quotient.

\medskip
i) L'espace quotient~$X/\mathsf G$ est un graphe ; si~$X$ est un arbre,~$X/\mathsf G$ est un arbre.

ii) Tout point de~$X/\mathsf G$ possède une base de voisinages~$V$ satisfaisant la condition suivante :~$V$ est un arbre compact et~$p\inv(V)$ est une réunion disjointe finie d'arbres compacts.

iii) Si~$\Gamma$ est un sous-graphe fermé de~$X$ alors~$p(\Gamma)$ est un sous-graphe fermé de~$X/\mathsf G$ ; si de plus ~$\Gamma$ est localement fini alors~$p(\Gamma)$ est localement fini, et~$p_{|\Gamma}$ est injective par morceaux.

}

\medskip
{\em Démonstration.}  Remarquons pour commencer que l'espace quotient~$X/\mathsf G$ est séparé et localement compact,  et que l'application~$p$ est compacte et ouverte  (\ref{quotsep} {\em et sq.}). 

\trois{quotbien} {\em Preuve de i) et ii).} Soit~$x\in X/\mathsf G$, soit~$\xi$ un antécédent de~$x$ sur~$X$, et soit~$U$ un voisinage ouvert de~$x$ dans~$X/\mathsf G$. Il existe un voisinage ouvert connexe~$\Omega$ de~$\xi$ dans~$p\inv(U)$ tel que~$\overline \Omega$ soit compact et contenu dans~$p\inv(U)$. Comme~$\mathsf G.\xi$ est compact et rencontre toutes les composantes connexes de~$\mathsf G.\Omega$ (l'ouvert~$g(\Omega)$ étant connexe pour tout~$g\in \mathsf G$), le nombre de composantes connexes de~$\mathsf G.\Omega$ est fini. 

Le sous-ensemble~$\mathsf G.\overline \Omega$ de~$X$ est compact (c'est l'image de~$\mathsf G\times \overline \Omega$), contient~$\mathsf G.\Omega$ et est contenu dans~$\overline{\mathsf G.\Omega}$ ; en conséquence,~$\overline{\mathsf G.\Omega}$ coïncide avec~$\mathsf G.\overline \Omega$ et est en particulier compact ; par ailleurs, il découle de la finitude de~$\pi_0(\mathsf G.\Omega)$ que~$ \overline{\mathsf G.\Omega}=\mathsf G.\overline \Omega$ est un sous-graphe de~$X$. 

 Le graphe compact~$\mathsf G.\overline \Omega$ comprend un nombre fini de boucles, puisqu'elles sont toutes situées sur son squelette, lequel est un arbre compact et fini. Soit~$(\xi_1,\ldots,\xi_r)$ une famille finie de points de~$\mathsf G.\Omega$ telle que chaque boucle de~$\mathsf G.\overline \Omega$ contienne l'un des~$\xi_i$ et telle que~$\mathsf G.\xi$ ne rencontre pas~$\{\xi_1,\ldots,\xi_r\}$ (c'est possible car le compact totalement discontinu~$\mathsf G.\xi$ ne contient aucune boucle). Il existe un voisinage ouvert connexe~$\Omega_0$ de~$\xi$ dans~$\Omega$ qui possède les propriétés  suivantes : 
 
 \medskip
~$\bullet$ le compact~$\overline \Omega_0$ est contenu dans~$\Omega-\mathsf G.\{\xi_1,\ldots,\xi_r\}$ ; 
 
~$\bullet$~$\partial \Omega_0$ est fini.	
 
\medskip
En appliquant à~$\Omega_0$ le raisonnement tenu ci-dessus à propos de~$\Omega$, on voit que le sous-graphe ouvert~$\mathsf G.\Omega_0$ de~$X$ a un nombre fini de composantes connexes, que son adhérence est compacte et qu'elle coïncide avec~$\mathsf G.\overline {\Omega_0}$. Par construction, cette adhérence ne contient aucun des~$\xi_i$, et donc aucune des boucles tracées sur~$\mathsf G.\overline \Omega$. Ainsi,~$\mathsf G.\overline {\Omega_0}$ apparaît comme un graphe compact sans boucle, c'est-à-dire comme une réunion disjointe finie d'arbres compacts. 

On note~$\Omega_0'$ (resp.~$\Omega''_0$) la composante connexe de~$\xi$ dans~$\mathsf G.\Omega_0$ (resp.~$\mathsf G.\overline{\Omega_0}$), et~$\mathsf G'$ (resp.~$\mathsf G''$) le stabilisateur de~$\Omega_0'$ (resp.~$\Omega_0''$). Par construction,~$\Omega'_0$ (resp.~$\Omega_0''$) est un sous-arbre ouvert (resp. un sous-arbre compact) de~$X$.

\medskip
Posons~$U_0=p(\Omega_0)$. Comme~$p$ est ouverte,~$U_0$ est un voisinage ouvert de~$x$ dans~$U$. Comme~$p$ est fermée,~$p(\overline \Omega_0)$ est fermé dans~$X$, ce qui implique que~$\partial U_0\subset p(\partial \Omega_0)$ ; comme~$\partial \Omega_0$ est fini,~$\partial U_0$ est fini.

On a~$U_0=p(\Omega_0)=p(\mathsf G.\Omega_0)=p(\Omega'_0)$. L'ouvert~$U_0$ s'identifie au quotient~$\Omega'_0/\mathsf G'$ ; la proposition~\ref{propprel} assure alors que pour tout couple~$(z,t)$ de points de~$U_0$, il existe un unique segment d'extrémités~$z$ et~$t$ tracé sur~$U_0$ ; ceci achève de prouver que~$X/\mathsf G$ est un graphe. 

\medskip
Cela étant établi, on déduit immédiatement de la proposition ~\ref{propprel} que~$X/\mathsf G$ est un arbre si~$X$ est un arbre ; ainsi, i) est établie. 

\medskip
Posons~$V=p(\overline \Omega_0)=p(\mathsf G.\overline {\Omega_0})=p(\Omega''_0)$ ; c'est un voisinage compact de~$x$ contenu dans~$U$. Par compacité, ~$V$ s'identifie au quotient~$\Omega''_0/\mathsf G''$ ; comme~$\Omega''_0$ est un arbre,~$V$ est un arbre d'après l'assertion i) déjà établie. Par construction,~$p\inv(V)=\mathsf G.\overline {\Omega_0}$ est une réunion finie disjointe d'arbres compacts, d'où ii). 

\trois{imsousgra}{\em Preuve de iii)}. Soit~$\Gamma$ un sous-graphe fermé de~$X$. Les assertions à établir pouvant se vérifier localement sur~$X/\mathsf G$ l'on peut, grâce à l'assertion ii) prouvée ci-dessus, supposer que~$X$ est une réunion finie disjointe d'arbres compacts, et que~$X/\mathsf G$ est un arbre. 

\medskip
Une réunion finie de sous-graphes fermés (resp. de sous-graphes fermés et localement finis) de~$X/\mathsf G$ est encore un sous-graphe fermé (resp. un sous-graphe fermé et localement fini) de~$X/\mathsf G$ ; cela permet de se ramener, en raisonnant composante par composante, au cas où le graphe compact~$\Gamma$ est connexe et non vide puis, quitte à remplacer~$X$ par la composante connexe de~$X$ contenant~$\Gamma$, et~$\mathsf G$ par le stabilisateur correspondant, à celui où~$X$ est lui-même également connexe. Autrement dit, il suffit de traiter le cas où~$X$ est un arbre compact et où~$\Gamma$ est un sous-arbre fermé de~$X$. 

\medskip
L'image de~$\Gamma$ sur~$X/\mathsf G$ est une partie compacte de~$X/\mathsf G$. Soient~$x$ et~$y$ deux éléments de~$p(\Gamma)$, et soit~$\xi$ (resp.~$\eta$) un antécédent de~$x$ (resp.~$y$) sur~$\Gamma$. L'assertion i) de la proposition~\ref{propprel} assure que~$p([\xi;\eta])$ est un arbre ; il contient donc~$[x;y]$ ; par conséquent, le compact~$p(\Gamma)$ de~$X/\mathsf G$ est convexe, et en est dès lors un sous-arbre fermé. 

\medskip
Si le sous-arbre compact~$\Gamma$ de~$X$ est localement fini, il est fini. Il découle alors de l'assertion i) de la proposition ~\ref{propprel} : que  le sous-arbre~$p(\Gamma)$ de~$X/\mathsf G$ est réunion finie de sous-arbres compacts et finis, ce qui entraîne qu'il est lui-même compact et fini ; et que~$p_{|\Gamma}$ est injective par morceaux.~$\Box$ 

\deux{quoetflat} {\bf Lemme.} {\em Soit~$X$ un graphe sur lequel~$\mathsf G$ agit, soit~$p:X\to X/\mathsf G$ la flèche quotient et soit~$\Gamma$ un sous-graphe admissible de~$X$ stable sous~$\mathsf G$ ; l'image~$p(\Gamma)$ est alors un sous-graphe admissible de~$X/\mathsf G$.}

\medskip
{\em Démonstration.} Comme~$\Gamma$ est un sous-graphe fermé de~$X$, son image~$p(\Gamma)$ est un sous-graphe fermé de~$X/\mathsf G$ (th.~\ref{theoquot},  iii) ).

\trois{pgammafad} Soit~$U$ une composante connexe de~$X/\mathsf G-p(\Gamma)$ et soit~$V$ une composante connexe de~$p\inv(U)$. L'ouvert~$V$ est une composante connexe de~$X-\Gamma$ et est donc un arbre à un bout relativement compact dans~$X$ ; on note~$\mathsf H$ le stabilisateur de~$V$ dans~$\pi_0(p\inv(U))$, et~$x$ l'unique point de~$\partial V$ ; le point~$x$ appartient à~$\Gamma$ et est fixe sous~$\mathsf H$, et son image sur~$X$ n'appartient pas à~$U$. 

\medskip
On a par compacité de~$p$ l'égalité~$\overline U=p(\overline V)=U\cup\{p(x)\}$. Comme~$V\simeq U/\mathsf H$, comme~$p(x)\notin U$ et comme~$x$ est fixe sous~$\mathsf H$, la flèche~$\overline V\to p(\overline V)$ se factorise par une bijection continue de~$(\overline V/\mathsf H)$ sur~$p(\overline V)$ ; le graphe~$\overline V$ étant compact, il vient~$p(\overline V)\simeq \overline V/\mathsf H$ ; en vertu du théorème~\ref{theoquot},  le compact~$\overline V/\mathsf H$ est un arbre.

\medskip
Ainsi,~$\overline U=p(\overline V)\simeq \overline V/\mathsf H$ est un arbre compact et~$\overline U-U$ est le singleton~$\{p(x)\}$ ; l'ouvert~$U$ est de ce fait un arbre à un bout relativement compact dans~$X/\mathsf G$. Par conséquent, le sous-graphe fermé~$p(\Gamma)$ de~$X/\mathsf G$ est admissible.~$\Box$

\deux{gprofpluribr} {\bf Lemme.} {\em Soit~$X$ un graphe sur lequel~$\mathsf G$ agit, soit~$x\in X$ et soit~$\mathsf H$ le stabilisateur de~$x$.  

\medskip
A) Si~$x$ est pluribranche,~$\mathsf H$ est un sous-groupe ouvert de~$\mathsf G$. 

\medskip
B) Soit~$\Gamma\in \stel X x$. Il existe~$\Delta\in \stelac X x$, stable sous~$\mathsf H$, et possédant les propriétés suivantes : 

\medskip
i) ~$\br \Delta x\supset \br \Gamma x$ ; 

ii) la réunion~$\{x\}\cup\bigcup\limits_{b\in \br \Gamma x}b(\Delta)$ est une sous-étoile de~$\Gamma$.

\medskip
C) Soit~$\Delta$ comme en B) et soit~$b\in \br \Delta x$. Le stabilisateur de~$b$ dans~$\mathsf H$ coïncide avec le stabilisateur de la composante connexe~$b(\Delta)$ de~$\Delta\setminus\{x\}$ et fixe~$b(\Delta)$ point par point.}

\medskip
{\em Démonstration.} D'après l'assertion ii) du théorème~\ref{theoquot}, il existe un voisinage compact~$V$ de~$x$ tel que~$\mathsf G.V$ soit une réunion disjointe finie d'arbres compacts. Quitte à remplacer~$x$ par sa composante connexe dans l'intérieur de~$\mathsf G. V$, et~$\mathsf G$ par le stabilisateur de la composante en question (qui est un ouvert de~$\mathsf G$), on peut supposer que~$X$ est un arbre. 

\trois{stabpluriouv} {\em Preuve de A)}. Supposons que~$x$ est pluribranche ; il existe alors deux points~$y$ et~$z$ de~$X$ tels que~$x\in ]y;z[$. L'orbite~$\mathsf G.x$ étant un compact totalement discontinu, elle ne contient ni~$]y;x[$ ni~$]z;x[$ ; par conséquent on peut supposer, quitte à rapprocher~$y$ et~$z$ de~$x$, que ni~$y$ ni~$z$ n'appartiennent à~$\mathsf G. x$. Soit~$E$ l'enveloppe convexe du compact~$\mathsf G. y\cup \mathsf G.z$ ; c'est un sous-arbre compact de~$X$ qui contient~$x$. Ni~$y$ ni~$z$ n'appartenant à~$\mathsf G.x$, le point~$x$ n'est pas situé sur ~$\mathsf G.y\cup \mathsf G.z$ ; or le graphe~$E-\left(\mathsf G.y\cup \mathsf G.z\right)$ est localement fini, et stable sous~$\mathsf G$. Il découle alors du lemme~\ref{orbfinlocfin} que~$\mathsf G.x$ est finie, ce qui achève de prouver A). 

\trois{etstabh} {\em Preuve de B)}. En vertu de~\ref{exaradm},~$\stelac X x \ctd \Gamma$ est non vide ; soit~$\Gamma_0$ un élément de cet ensemble. Choisissons pour tout~$b\in \br \Gamma x$ un point~$x_b$ sur l'arête~$b(\Gamma_0)$ ; étant situé sur un intervalle ouvert, le point~$x_b$ est pluribranche ; en vertu de l'assertion A) déjà prouvée, l'orbite~$\mathsf H.x_b$ est finie. 

\medskip
Il s'ensuit que le compact~$\Gamma_1:=\mathsf G.\left(\bigcup [x_b;x]\right)=\bigcup\limits_{g\in \mathsf H}[g(x_b);x]$ est un arbre compact et fini, évidemment invariant par~$\mathsf H$. Soit~$S$ l'ensemble de ses sommets et soit~$\Delta$ la composante connexe de~$x$ dans~$\Gamma_1\setminus S$ ; c'est étoile de sommet~$x$ tracée sur~$X$, relativement compacte et stable sous~$\mathsf H$ ; nous allons montrer que~$\Delta$ satisfait i) et ii), puis que ses arêtes sont faiblement admissibles. 

\medskip
{\em La condition i).} Elle résulte du fait que~$\Delta$ contient pour tout~$b\in \br \Gamma x$ l'intervalle ouvert~$]x_b;x[$ qui aboutit proprement à~$x$ et définit~$b$.

\medskip
{\em La condition iii).} Soit~$b\in \br \Gamma x$. Comme~$]x_b;x[$ et~$b(\Delta)$ sont toutes deux tracées sur~$b(X)$, leur intersection est un intervalle ouvert~$I$ tracé sur~$b(\Delta)$ et aboutissant à~$x$. Nous allons montrer par l'absurde que~$I=b(\Delta)$ ; il en découlera que~$b(\Delta)\subset ]x_b;x[\subset b(\Gamma)$, ce qui permettra de conclure. 

Supposons donc que~$I\neq b(\Delta)$ ; l'intervalle~$I$ aboutit dans ce cas à un point~$y$ de~$b(\Delta)$. Mais il y a alors au moins trois branches de~$\Gamma_1$ issues de~$y$, à savoir celles respectivement définies par les intervalles~$]y;x[$,~$]x_b;y[$ et~$b(\Delta)-[y;x[$ ; dès lors~$y$ est un sommet de~$\Gamma_1$, ce qui contredit son appartenance à~$b(\Delta)$.

\medskip
{\em Les arêtes de~$\Delta$ sont faiblement admissibles.} Soit~$I$ une arête de~$\Delta$ ; elle aboutit à un sommet~$y$ de~$\Gamma_1$ qui diffère de~$x$ ; ce sommet est un point de~$[g(x_b);x[$ pour un certain~$g\in \mathsf H$ et un certain~$b\in \br \Gamma x$ ; par conséquent,~$I$ est un ouvert de~$]g(x_b);x[=g(]x_b;x[)$. Or~$]x_b;x[$ est lui-même un ouvert de~$b(\Gamma_0)$, qui est faiblement admissible par choix de~$\Gamma_0$ ; par conséquent,~$I$ est faiblement admissible.

\medskip
{\em Preuve de C)}. Soit~$g\in \mathsf H$. Comme~$b$ est la branche issue de~$x$ définie par~$b(\Delta)$, il résulte des définitions que~$g$ envoie la branche~$b$ sur la branche issue de~$x$ définie par l'arête~$g(b(\Delta))$ de~$\Delta$ ; par conséquent, le stabilisateur de~$b$ dans~$\mathsf H$ coïncide avec le stabilisateur de la composante connexe~$b(\Delta)$ de~$\Delta\setminus\{x\}$.

Si~$h$ appartient au stabilisateur en question et si~$y\in b(\Delta)$ alors~$h(y)$ appartient aussi à~$b(\Delta)$, et est donc comparable à~$y$ pour l'ordre induit par~$x$ sur l'arbre~$X$ ; compte-tenu du fait que~$h(x)=x$, cela entraîne en vertu du lemme~\ref{lemmeprof} que~$h(y)=y$, ce qui achève la démonstration.~$\Box$ 

\deux{xunixpluri} Soit~$X$ un graphe sur lequel~$\mathsf G$ agit et soit~$x\in X$ ; soit~$p$ la flèche quotient~$X\to  X/\mathsf G$ , et soit~$\xi$ l'image de~$x$ sur~$X/\mathsf G$. 

\trois{pluribiendef} Supposons que~$x$ est pluribranche.  D'après l'assertion A) du lemme~\ref{gprofpluribr} ci-dessus, la fibre de~$p$ en~$\xi$ est finie ; comme~$p$ est par ailleurs compacte et ouverte, les résultats généraux de~\ref{imdir},~\ref{imrec} et~\ref{descimrec}-\ref{imdirboncomp} s'appliquent ici ; on dispose notamment d'une application naturelle~$\br X x \to \br {X/\mathsf G}\xi$. 

\trois{uniuni} Supposons que~$x$ est unibranche ; nous allons montrer que~$\xi$ l'est aussi. L'assertion ii) du théorème~\ref{theoquot} fournit un voisinage ouvert~$V$ de~$\xi$ dans~$X/\mathsf G$ qui est un arbre, et qui est tel que~$p\inv(V)$ soit une réunion finie d'arbres. Quitte à remplacer~$X$ par la composante connexe de~$x$ dans~$p\inv(V)$, et~$\mathsf G$ par le stabilisateur de la composante en question, on peut supposer que~$X$ et~$X/\mathsf G$ sont des arbres. 

\medskip
Comme~$x$ est unibranche il n'est pas isolé et~$X$ n'est pas réduit à~$\{x\}$ ; par conséquent, le compact totalement discontinu~$\mathsf G.x$ n'est pas égal à~$X$ tout entier, et~$X-\mathsf G.x$ est un ouvert non vide de~$X$. Il est convexe : en effet, si~$y$ et~$z$ sont deux points de~$X-\mathsf G.x$, alors tout point de~$]y;z[$ est pluribranche (puisque situé sur un intervalle ouvert) et ne peut donc appartenir à~$\mathsf G.x$. 

\medskip
Il s'ensuit que l'image de~$X-\mathsf G.x$ sur~$X/\mathsf G$, qui n'est autre que~$X/\mathsf G\setminus\{\xi\}$, est connexe et non vide. Par conséquent,~$\xi$ est un point unibranche de l'arbre~$X/\mathsf G$. 

\trois{isolisol} Supposons que~$x$ est isolé. Dans ce cas c'est une composante connexe de~$X$ ; son stabilisateur est alors nécessairement ouvert, et~$X$ est donc un ensemble fini discret. Il en va dès lors de même de~$X/\mathsf G$, et~$\xi$ est en particulier isolé.

\subsection*{Compactification et quotient}

\deux{xchapeausurg}{\bf Proposition.} {\em Soit~$X$ un graphe connexe sur lequel~$\mathsf G$ agit. Il existe un homéomorphisme naturel~$\wid X/\mathsf G\simeq\wid{ \left(X/\mathsf G\right)}$ qui identifie~$(\got d X)/\mathsf G$ et~$\got d(X/\mathsf G)$.}

\medskip
{\em Démonstration.} On désigne par~$p$ la projection~$X\to X/\mathsf G$. Soit~$E$ l'ensemble des points pluribranches de~$X$ ; c'est une partie stable sous l'action de~$\mathsf G$ ; si~$x\in E$ l'orbite~$\mathsf G.x$ est finie par la proposition~\ref{gprofpluribr}. Soit~$\sch S$ l'ensemble des parties finies de~$E$ stable sous~$\mathsf G$. 

\medskip
L'ensemble des compacts de la forme~$\sch K_X(S)$, où~$S$ est une partie finie de~$E$, est cofinal dans l'ensemble des compacts de~$X$ (\ref{kxscof}). Si~$S$ est une partie finie de~$E$, l'ensemble~$\mathsf G.S$ appartient à~$\sch S$ ;  comme~$S\subset \mathsf G.S$, le compact~$\sch K_X(\mathsf G.S)$ contient~$\sch K_X(S)$ ; par conséquent, l'ensemble des compacts de la forme~$\sch K_X(S)$, où~$S\in \sch S$, est cofinal dans l'ensemble des compacts de~$X$. 

\medskip
Soit~$S\in \mathsf S$ ; le compact~$\sch K_X(S)$ est invariant sous~$\mathsf G$ ; son image~$p(\sch K_X(S))$ est une partie compacte de~$X/\mathsf G$. L'invariance de~$S$ sous~$\mathsf G$ implique que~$$p(\sch K_X(S))-p(S)=p(\sch K_X\setminus S)\;;$$ de ce fait,~$\sch K_X(S)-p(S)$ est un ouvert de~$X/\mathsf G$ ; comme c'est aussi un fermé de~$X/\mathsf G-p(S)$, c'est une réunion de composantes connexes de~$X/\mathsf G-p(S)$, qui sont toutes relativement compactes. 

Soit~$U$ une composante connexe de~$X/\mathsf G-p(\sch K_X(S))$ ; son image réciproque est une réunion disjointe et non vide de composantes connexes de~$X\setminus \sch K_X(S)$ ; si~$U$ était relativement compacte dans~$X/\mathsf G$, chacune de ces composantes le serait dans~$X$ par compacité de~$p$, contredisant la définition de~$\sch K_X(S)$ ; par conséquent,~$U$ n'est pas relativement compacte. 

Il résulte de tout ceci que~$p(\sch K_X(S))$ est égal à~$\sch K_{X/\mathsf G}(p(S))$. 

\medskip
Par ailleurs, si~$\sch K$ est un compact de~$X/\mathsf G$, son image réciproque~$p\inv(\sch K)$ est compacte, et donc contenue dans~$\sch K_X(S)$ pour un certain~$S\in \sch S$ ; par conséquent,~$\sch K\subset \sch K_{X/\mathsf G}(p(S))$ ; ainsi, l'ensemble des compacts de~$X/\mathsf G$ de la forme~$\sch K_{X/\mathsf G}(p(S))$ est cofinal dans l'ensemble des compacts de~$X/\mathsf G$.

\medskip
On dispose donc de deux identifications canoniques~$$\got dX\simeq \lim_{\stackrel{S\in \sch S}\leftarrow} \pi_0\left(X\setminus \sch K_X(S)\right)\;{\rm et}\;\got d (X/\mathsf G)\simeq \lim_{\stackrel{S\in \sch S}\leftarrow} \pi_0\left(X/\mathsf G\setminus \sch K_X(p(S))\;\right).$$

\medskip
Pour tout~$S$ appartenant à~$\sch S$, le groupe~$\mathsf G$ agit sur~$\pi_0(X\setminus \sch K_X(S))$ et le quotient ~$\pi_0(X\setminus \sch K_X(S))/\mathsf G$ s'identifie {\em via}~$p$ à~$\pi_0(X/\mathsf G\setminus \sch K_X(p(S))\;).$ On en déduit, compte-tenu de ce qui précède, l'existence d'une application naturelle~$(\got d X)/\mathsf G\to \got d(X/\mathsf G)$. 

\medskip
{\em Cette application est injective.} Soient~$x$ et~$y$ deux points de~$\got d X$ ayant même image dans~$\got d(X/\mathsf G)$. Pour tout~$S\in \sch S$, notons~$x_S$ et~$y_S$ les images de~$x$ et~$y$ dans~$\pi_0(X\setminus \sch K_X(S))$, et~$\mathsf G_S$ l'ensemble des~$g\in \mathsf G$ qui envoient~$x_S$ sur~$y_S$ ; c'est une partie compacte (et ouverte) de~$\mathsf G$. 

Si~$S\in \sch S$ il résulte de notre hypothèse que~$x_S$ et~$y_S$ ont même image dans~$\pi_0(X/\mathsf G\setminus \sch K_X(p(S))\;)= \pi_0(X\setminus \sch K_X(S))/\mathsf G$ ; par conséquent,~$\mathsf G_S$ est non vide. Par compacité, l'intersection des~$\mathsf G_S$ pour~$S$ parcourant~$\sch S$ est non vide ; si~$g$ est  un élément quelconque de celle-ci alors~$y=g(x)$ ; il s'ensuit que~$y$ et~$x$ ont même image dans~$(\got d X)/\mathsf G$ et l'application étudiée est bien injective.

\medskip
{\em Cette application est surjective.} Soit~$x$ un point de~$\got d(X/\mathsf G)$. Pour tout~$S\in \sch S$, soit~$x_S$ l'image de~$x$ dans~$\pi_0(X/\mathsf G\setminus \sch K_X(p(S))\;)$ et soit~$\xi_S$ l'image réciproque de~$x_S$ dans~$\pi_0(X\setminus \sch K_X(S))$. La famille des~$(\xi_S)$ est un système projectif filtrant d'ensembles finis et non vide, sa limite projective est donc non vide et si~$\xi$ désigne un point de cette dernière, on a par construction~$p(\xi)=x$ ; l'application étudiée est donc bien surjective.

\medskip
Ainsi, ~$(\got d X)/\mathsf G$ s'identifie à~$\got d(X/\mathsf G)$ ; on définit par ce biais une bijection {\em ensembliste}~$\wid X/\mathsf G\to \wid {X/\mathsf G}$ ; il reste à s'assurer que c'est un homéomorphisme. Ses source et but étant compacts, on peut se contenter d'en vérifier la continuité, c'est-à-dire de prouver que la flèche composée~$\wid X\to \wid X/\mathsf G \to \wid {X/\mathsf G}$ est continue.

\medskip
Il s'agit de démontrer que si~$U$ est un ouvert de~$\wid {X/\mathsf G}$, son image réciproque sur~$X$ est ouverte ; il suffit de le faire lorsque~$U$ parcourt une base donnée d'ouverts de~$\wid {X/\mathsf G}$. Par cofinalité de l'ensemble des~$\sch K_{X/\mathsf G}(p(S))$, où~$S$ parcourt~$\sch S$, dans celui de tous les compacts de~$\mathsf G$, on se ramène ainsi à traiter les deux cas suivants : 

\medskip
$\bullet$ le cas où~$U$ est un ouvert de~$X/\mathsf G$, qui est trivial ; 

$\bullet$ le cas où~$U$ est de la forme~$\widd V {X/\mathsf G}$ où~$V$ est une composante connexe de~$ X/\mathsf G\setminus \sch K_X(p(S))$ pour un certain~$S\in \sch S$. L'ouvert~$p\inv(V)$ de~$X$ est réunion finie de composantes connexes de~$X\setminus \sch K_X(S)$, et il résulte de la construction de notre bijection que l'image réciproque de~$\widd V {X/\mathsf G}$ est égale à~$\widd V X$, qui est un ouvert de~$\wid X$ ; ceci achève la démonstration.~$\Box$ 

\deux{exarbreunbout} {\bf Exemple  : l'action de~$\mathsf G$ sur un arbre à un bout.} Soit~$X$ un arbre  un bout sur lequel~$\mathsf G$ agit ; on note~$\omega$ l'unique élément de~$\got d X$, qui est fixe sous~$\mathsf G$. Il résulte du théorème~\ref{theoquot} et de la proposition~\ref{xchapeausurg} que~$X/\mathsf G$ est un arbre à un bout.

\trois{unboutglobinv} Soit~$y$ un point pluribranche de~$X$ ; son orbite sous~$\mathsf G$ est finie d'après la proposition~\ref{gprofpluribr}). Posons~$I=\bigcap\limits_{z\in \mathsf G.y}[y;\omega]$. Il découle de la finitude de~$\mathsf G.y$ et le lemme~\ref{lemextremchap} que~$I$ est de la forme~$[x;\omega]$ pour un certain~$x\in X$. 

\medskip
Le point~$\omega$ est invariant sous~$\mathsf G$ ; si~$t\in [x;\omega]$ et~$g\in \mathsf G$ on a~$g(t)\in [x;\omega]$, et~$g(t)$ est donc comparable à~$t$ pour l'ordre défini par~$\omega$. En vertu du lemme~\ref{lemmeprof}, le groupe~$\mathsf G$ fixe~$[x;\omega]$ point par point. 

\trois{unboutquotinj} Si~$\xi\in X$ l'application~$[\xi;\omega]\to \wid{X/\mathsf G}$ est injective. En effet, soient~$\eta$ et~$\zeta$ deux points de~$[\xi;\omega]$ ayant même image sur~$\wid{X/\mathsf G}\simeq \wid X/\mathsf G$ ; il existe alors~$g\in \mathsf G$ tel que~$g(\eta)=\zeta$ ; comme ils sont situés sur~$[\xi;\omega]$, les points~$\eta$ et~$\zeta$ sont comparables pour l'ordre induit par~$\omega$. Ce dernier étant fixe par~$g$, on conclut à l'aide du lemme~\ref{lemmeprof} que~$\eta=g(\eta)=\zeta$, d'où notre assertion. 

\medskip
Soit~$x$ l'image de~$\xi$ sur~$\wid {X/\mathsf G}$, et soit~$w$ celle de~$\omega$. On en déduit de ce qui précède, en vertu de la compacité de~$[\xi;\omega]$, que~$\wid X\to \wid{X/\mathsf G}$ induit un homéomorphisme entre~$[\xi;\omega]$ et son image ; cette dernière est donc  un segment éventuellement long d'extrémités~$x$ et~$w$ ; elle coïncide dès lors avec~$[x;w]$. 

\medskip
Si~$\xi$ est pluribranche il existe une composante connexe de~$X\setminus\{x\}$ qui n'est pas celle contenant~$]x;\omega[$ ; par conséquent, il existe~$x'\in X$ tel que~$x\in ]x';\omega[$. En vertu de ce qui précède, l'image de~$]x';\omega[$ sur~$X/\mathsf G$ est un intervalle ouvert, qui contient évidemment~$\xi$. On en déduit que~$\xi$ est pluribranche. 

\trois{anteccroissants} Soit~$x$ un point de~$X/\mathsf G$ ayant un nombre fini~$N$ d'antécédents sur~$X$. Pour tout~$t\in [x;w[$, le point~$t$ a au plus~$N$ antécédents sur~$X$. En effet, soient~$\xi_1,\ldots,\xi_N$ les antécédents de~$x$ sur~$X$. Pour tout~$g\in \mathsf G$ et tout indice~$i$, l'image de~$[x_i;\omega[$ par~$g$ est égale à~$[g(\xi_i);\\omega[$ et donc à~$[\xi_j;\omega[$ pour un certain~$j$. Il s'ensuit que la réunion des~$[\xi_i;\omega[$ est stable sous~$\mathsf G$. Comme son image est égale à~$[x;w[$ en vertu du~\ref{unboutquotinj} ci-dessus,~$\bigcup [\xi_i;\omega[$ est égale à l'image réciproque de~$[x;w[$. Toujours d'après~\ref{unboutquotinj}, la flèche~$[\xi_i;\omega[\to [x;w[$ est un homéomorphisme pour tout~$i$. Par conséquent, tout point de~$[x;w[$ a au plus~$N$ antécédents sur~$X$, comme annoncé. 

\deux{exarbdeuxbouts} {\bf Exemple : l'action de~$\mathsf G$ sur un arbre à deux bouts.} Soit~$X$ un arbre à deux bouts sur lequel~$\mathsf G$ agit. Le squelette de~$X$ est la droite éventuellement longue~$]\xi;\eta[$ où~$\xi$ et~$\eta$ sont les deux bouts de~$X$. 

\trois{gpermutedeux} {\em Le cas où~$\mathsf G$ échange~$\xi$ et~$\eta$.} Il découle alors de la proposition~\ref{xchapeausurg} ci-dessus que~$X/\mathsf G$ est un arbre à un bout ; si~$\omega$ désigne l'unique élément de~$\got d (X/\mathsf G)$, le corollaire~\ref{corollprof} entraîne que l'image de~$\skel X$ sur~$X/\mathsf G$ est de la forme~$[t;\omega[$ pour un certain~$t\in X/\mathsf G$, égal à l'image de l'unique point fixe de~$\mathsf G$ sur~$\skel X$. La demi-droite~$[t;\omega[$ est un sous-arbre admissible de~$\mathsf X/\mathsf G$ (\ref{admunbout} ou lemme~\ref{quoetflat}) mais n'en est pas le squelette -- ce dernier est vide. 

\trois{gnepermutepas} {\em Le cas où~$\mathsf G$ n'échange pas~$\xi$ et~$\eta$.} Il découle alors de la proposition~\ref{xchapeausurg} ci-dessus que~$X/\mathsf G$ est un arbre à deux bouts, qui correspondent aux images de~$\xi$ et~$\eta$ sur~$\wid {X/\mathsf G}$, que l'on notera respectivement~$x$ et~$y$. En vertu du corollaire~\ref{corolcroisdl}, l'action de~$\mathsf G$ sur~$[\xi;\eta]$ est triviale ; l'image de~$\skel X$ sur~$X/\mathsf G$ s'identifie donc à~$ ]x;y[$, c'est-à-dire à~$\mathsf S(X/\mathsf G)$. 

Par conséquent, l'image réciproque de~$\mathsf S(X/\mathsf G)$ sur~$X$ est égale à~$\skel X$, et~$\skel X\to \mathsf S(X/\mathsf G)$ est un homéomorphisme. Le lemme~\ref{lemmcommret} asssure que~$X\to X/\mathsf G$ commute aux rétractions canoniques de ses source et but sur leurs squelettes respectifs. 

\medskip

Si~$\mathsf H$ est un sous-groupe fermé de~$\mathsf G$, on déduit de ce qui précède, appliqué d'une part directement et d'autre part à la flèche quotient~$X\to \mathsf X/\mathsf H$, que~$X/\mathsf H$ est un arbre à deux bouts, que l'image réciproque de~$\mathsf S(X/\mathsf G)$ sur~$X/\mathsf H$ est égale à~$ \mathsf S(X/\mathsf H)$, que~$\mathsf S(X/\mathsf H)\to \mathsf S(X/\mathsf G)$ est un homéomorphisme et que~$X/\mathsf H\to X/\mathsf G$ commute aux rétractions canoniques de ses source et but sur leurs squelettes respectifs. 

\subsection*{Étude locale autour d'un point fixe : sections et voisinages stables}

\deux{corollstabbr} {\bf Proposition.} {\em Soit~$X$ un graphe sur lequel~$\mathsf G$ agit, soit~$p: X\to X/\mathsf G$ la flèche quotient et soit~$x$ un point de~$X$ invariant par~$\mathsf G$ ; soit~$b$ une branche issue de~$x$.

\medskip
\begin{itemize}
\item[1)] Le stabilisateur~$\mathsf H$ de~$b$ dans~$\mathsf G$ est un sous-groupe ouvert de~$\mathsf G$. 

\item[2)]  Soit~$\sch I$ le sous-ensemble de~$\interac X b~$ constitué des intervalles~$I$ possédant les propriétés suivantes : 

\medskip
\begin{itemize}

\item[$\alpha)$]~$I$ est fixé point par point par~$\mathsf H$ ; 

\item[$\beta)$] si~$g\in \mathsf G-\mathsf H$ l'on a~$g(I^\flat)\cap I^\flat=\emptyset$ ;

\item[$\gamma)$] l'application~$p$ induit un homéomorphisme~$\overline I\simeq p(\overline I)$. 

\end{itemize}

\medskip
Les ouverts de~$X$ de la forme~$I^\flat$ avec~$I\in \sch I$ forment une base de sections de~$b$. 

\item[3)] Si~$I\in \sch I$ l'intervalle ouvert~$p(I)$ aboutit proprement à~$p(x)$, est faiblement admissible dans~$X/\mathsf G$ et l'on a~$p(I)^\flat=p(I^\flat)$ ; ce dernier s'identifie à~$I^\flat/\mathsf H$ et est une section de~$p(b)$. 

\item[4)]  Le groupe~$\mathsf G$ agit transitivement sur~$p\inv(p(b))$.

\end{itemize} }

\medskip
{\em Démonstration.} Choisissons une demi-droite~$\Gamma$ issue de~$x$ telle que~$\Gamma\setminus\{x\}$ définisse~$b$.  

\trois{stabropen} {\em Preuve de 1)}. On peut voir~$\Gamma$ comme une étoile tracée sur~$X$, de sommet~$x$ et de valence 1, et l'on a~$\br \Gamma x=\{b\}$. L'assertion B) du lemme~\ref{gprofpluribr} affirme alors l'existence d'une étoile~$\Delta\in \stela X x$ qui est table sous~$\mathsf G$ et satisfait les conditions  i) et ii) de son énoncé. La branche~$b$ appartient à~$\br \Delta x$ (condition i) ), et~$\mathsf H$ s'identifie au stabilisateur de la composante connexe~$b(\Delta)$ de~$\Delta\setminus\{x\}$ ({\em loc. cit.}, assertion C)) ; c'est donc automatiquement un sous-groupe ouvert de~$\mathsf G$. 

\trois{secinvar} {\em Preuve de 2)}. Les ouverts de la forme~$I^\flat$, où~$I\in \intera X b$, forment une base de sections de~$b$ (\ref{declun}). Par ailleurs, si~$I\in \sch I$, il en va de même de tout intervalle ouvert de~$I$ aboutissant à~$x$ ; il suffit donc pour établir 2) d'exhiber {\em un} intervalle~$I$ appartenant à~$\sch I$. 

\medskip
L'intervalle~$b(\Delta)$ aboutit à~$x$ et est relativement compact dans~$X$ ; il est faiblement admissible. Le stabilisateur de~$b(\Delta)$ dans~$\Delta\setminus\{x\}$ est exactement~$\mathsf H$, et ce dernier fixe~$b(\Delta)$ point par point (lemme~\ref{gprofpluribr}, assertion C). Si~$g\in \mathsf G-\mathsf H$ alors~$g$ ne stabilise pas l'arête~$b(\Delta)$ de~$\Delta$, ce qui signifie que l'arête~$g(b(\Delta))$ est distincte de~$\b(\Delta)$ ; ceci entraîne que les ouverts~$b(\Delta)^\flat$ et~$g(b(\Delta))^\flat$ sont disjoints (\ref{iflatar}). 

\medskip
La restriction de~$p$ à~$\overline {b(\Delta)}$ étant injective par morceaux en vertu de l'assertion iii) du théorème~\ref{theoquot}, il existe~$I\in \inter X x \ctd{b(\Delta)}$ tel que~$p$ induise un homéomorphisme~$\overline I\to p(\overline I)$ ; par construction,~$I\in \sch I$, ce qui achève de prouver 3).

\trois{imegalsec} {\em Preuve de 3)}. Soit~$I\in \sch I$ ; il découle directement de~$\gamma)$ que~$p(I)$ est un intervalle ouvert aboutissant proprement à~$x$. En vertu de~$\beta)$ 'image~$p(I^\flat)$ s'identifie à~$I^\flat/\mathsf H$. Comme~$I^\flat$ est un arbre à deux bouts, et comme~$\mathsf H$ agit trivialement sur~$I$ en vertu de~$\alpha)$, l'ouvert~$p(I^\flat)\simeq I^\flat/\mathsf H$ est un arbre à deux bouts de squelette~$p(I)$ (\ref{gnepermutepas} ). En particulier,~$p(I)$ est un sous-arbre admissible de~$p(I^\flat)$ ; il s'ensuit que~$p(I)$ est un intervalle faiblement admissible de~$X/\mathsf G$ et que~$p(I)^\flat=p(I^\flat)$. 

\medskip
L'intervalle~$I$ aboutit proprement à~$x$. Par conséquent, l'ouvert~$p(I)^\flat$ appartient à~$\sgb {X/\mathsf G}{p(x)}$ ; comme~$I^\flat$ est une section de~$b$, la branche définie par~$p(I)^\flat=p(I^\flat)$ est nécessairement~$p(b)$ (\ref{imdirboncomp}). 

\trois{gagtransbr} {\em Preuve de 4)}. Choisissons~$I\in \sch I$. Pour tout~$g\in \mathsf G$ l'ouvert~$g(I^\flat)$ est une section de~$g(b)$ ; il résulte par ailleurs de~$\alpha)$ et~$\beta)$ que~$p\inv(p(I^\flat))$ est la réunion disjointe des~$g(I^\flat)$ où~$g$ parcourt un système de représentants de~$\mathsf G/\mathsf H$ ; on peut donc écrire comme une réunion disjointe~$\coprod\limits_{c \in \mathsf G.b} U_c$ où~$U_c$ est pour tout~$c$ une section de~$c$. 

D'après 3), l'ouvert~$p(I^\flat)$ est une section de~$p(b)$. Le point~$x$ étant le seul antécédent de~$p(x)$ sur~$X$ (il est fixe sous~$\mathsf G$) il découle de ce qui précède et de~\ref{imrecboncomp} que~$p\inv(p(b))=\mathsf G.b$, ce qui prouve 4) et achève la démonstration.~$\Box$

\deux{corollimreci} {\bf Corollaire.} {\em Soit~$X$ un graphe sur lequel~$\mathsf G$ agit, soit~$p:X\to X/\mathsf G$ la flèche quotient et soit~$\xi$ un point de~$X/\mathsf G$ tel que~$p\inv(\xi)$ soit finie. Soit~$a$ appartenant à~$\br {X/\mathsf G}\xi$ ; il existe~$J_0\in  \intera {X/\mathsf G} a$ tel que pour tout intervalle ouvert~$J$ de~$J_0$ aboutissant à~$\xi$ les propriétés suivantes soient vérifiées : 

\medskip
$a)$~$p\inv(J)\to J$ est un revêtement topologique (évidemment trivial) dont les feuillets sont en nombre finis et sont tous des intervalles ouverts admissibles de~$X$. 

$b)$~$p\inv(J^\flat)=\coprod\limits_{I\in \pi_0(p\inv(I))} I^\flat$.}

\medskip
{\em Démonstration.} Il existe un voisinage ouvert~$U$ de~$\xi$ tel que~$p\inv(U)$ sépare les antécédents de~$\xi$. Pour tout~$x\in p\inv(\xi)$, le stabilisateur de la composante connexe~$p\inv(U)_x$ de~$p\inv(U)$ coïncide avec le stabilisateur~$\mathsf G_x$ de~$x$, et~$U$ s'identifie à~$p\inv (U)_x/\mathsf G_x$ ; par conséquent on peut, en raisonnant composante connexe par composante connexe, se ramener au cas où~$p\inv(\xi)$ est de la forme~$\{x\}$ pour un certain~$x$ fixe sous~$\mathsf G$. Choisissons une branche~$b$ de~$X$ issue de~$x$ et située au-dessus de~$a$. 

Soit~$\sch I$ comme dans le corollaire~\ref{corollstabbr} ci-dessus et soit~$I_0\in \sch I$ (son existence étant assurée par l'assertion 2) de {\em loc. cit.}) ; posons~$J_0=p(I_0)$. Par définition de~$\sch I$ et en vertu de l'assertion 3) de {\em loc. cit}, l'intervalle~$J_0$ satisfait les conditions souhaitées.~$\Box$

\deux{corollimrecih} {\bf Corollaire.} {\em Soit~$X$ un graphe sur lequel~$\mathsf G$ agit et soit~$\mathsf H$ un sous-groupe fermé de~$\mathsf G$. Soit~$\Gamma$ un sous-graphe fermé
 de~$X/\mathsf G$ et soit~$\Gamma'$ son image réciproque sur~$X/\mathsf H$ ; on suppose que les fibres sur~$X$ des points de~$\Gamma$ sont toutes finies.
 
 \medskip
 1) Le sous-ensemble fermé~$\Gamma'$ de~$X/\mathsf H$ en est un sous-graphe.
 
 2) Si de plus~$\Gamma$ est localement fini alors~$\Gamma'$ est localement fini, et~$\Gamma'\to \Gamma$ est injective par morceaux. }

\medskip
{\em Démonstration.} Soit~$\Gamma''$ l'image réciproque de~$\Gamma$ sur~$X/\mathsf G$. On procède en deux temps, en commençant par établir 2). 

\trois{imrecihloc} {\em Preuve de 2)}. On suppose que~$\Gamma$ est localement fini. En vertu du corollaire~\ref{corollimrecih} ci-dessus (et sans avoir recours à la propriété b) de son énoncé), le fermé~$\Gamma''$ de~$X$ en est un sous-graphe localement fini, et~$\Gamma''\to \Gamma$ est injective par morceaux. Comme~$\Gamma'$ s'identifie à~$\Gamma''/\mathsf H$, on déduit de l'énoncé iii) du théorème~\ref{theoquot} que~$\Gamma'$ est un sous-graphe fermé localement fini de~$X/\mathsf H$, et que~$\Gamma''\to \Gamma'$ est injective par morceaux. L'assertion requise s'ensuit aussitôt.

\trois{imrecihgen} {\em Preuve de 1)}. On ne fait plus d'hypothèse sur~$\Gamma$. Comme~$\Gamma'\simeq \Gamma''/\mathsf H$ il suffit, grâce à l'énoncé i) du théorème~\ref{theoquot}, de démontrer que~$\Gamma''$ est un sous-graphe de~$X$. La propriété requise étant locale, l'énoncé ii) de {\em loc. cit.} permet de supposer que~$X/\mathsf G$ est un arbre et que~$X$ est une réunion finie disjointe d'arbres ; en raisonnant composante par composante sur~$X$ (et en remplaçant à chaque fois~$\mathsf G$ par le stabilisateur de la composante considérée), on se ramène ensuite au cas où~$X$ et~$X/\mathsf G$ sont des arbres ; on peut enfin, en raisonnant encore composante par composante, mais cette fois-ci sur~$\Gamma$, supposer que celui-ci est un sous-arbre de~$X/\mathsf G$. Soit~$x\in \Gamma''$ et soit~$\xi$ son image sur~$X/\mathsf G$.

\medskip
Soit~$\Delta$ un sous-arbre fermé et localement fini de~$\Gamma$ qui contient~$\xi$. L'image réciproque de~$\Delta$ sur~$X$ est, en vertu de l'assertion 2) déjà établie, un sous-graphe fermé et localement fini de~$X$ ; notons~$\sch C_x(\Delta)$ la composante connexe de~$x$ dans ce dernier ; c'est un sous-arbre fermé de~$X$. Soit~$\Theta(x)$ la réunion des~$\sch C_x(\Delta)$, où~$\Delta$ parcourt l'ensemble des sous-arbres fermés et localement finis de~$\Gamma$ contenant~$\xi$. Par sa construction même,~$\Theta(x)$ est une partie convexe de~$X$ qui est contenue dans~$\Gamma''$. 

Nous allons montrer que~$\Theta(x)$ est un fermé de~$\Gamma''$. Pour ce faire, considérons un point~$x'$ de~$\Gamma$ qui n'appartient pas à~$\Theta(x)$. {\em Le segment~$[x';x]$ n'est alors pas contenu dans~$\Gamma''$.} En effet, s'il l'était, son image sur~$X/\mathsf G$ serait, par l'énoncé iii) de {\em loc. cit.}, un sous-arbre compact et fini~$\Delta$ de~$\Gamma$, contenant~$\xi$ ; mais on aurait alors~$[x';x]\subset \sch C_x(\Delta)\subset \Theta(x)$, ce qui est absurde. 

\medskip
Il existe donc un point~$z\in ]x';x[$ qui n'appartient pas à~$\Gamma''$. Soit~$U$ la composante connexe de~$x'$ dans~$X\setminus\{z\}$ et soit~$t\in U\cap \Gamma''$. Le segment~$[t;x]$ contient~$z$, et n'est donc pas inclus dans~$\Gamma''$ ; par convexité de~$\Theta(x)$, il s'ensuit que~$t\notin \Theta(x)$. 

Ainsi,~$U\cap \Gamma''$ est un voisinage de~$x'$ dans~$\Gamma''$ qui ne rencontre pas~$\Theta(x)$ ; ce dernier est donc bien fermé dans~$\Gamma''$, et partant dans~$X$ ; étant par ailleurs convexe,~$\Theta(x)$ est un sous-arbre fermé de~$X$. 

\medskip
Soit~$y$ un antécédent de~$\xi$, et supposons que~$\Theta(x)\cap \Theta(y)\neq \emptyset$. Si~$t$ désigne un point de cette intersection alors il existe deux sous-arbres fermés et localement finis~$\Delta$ et~$\Delta'$ de~$\Gamma$ qui contiennent~$x$ et qui sont tels que~$t\in \sch C_x(\Delta)\cap \sch C_y (\Delta')$ ; il s'ensuit que~$t\in \sch C_x(\Delta\cup\Delta')$ et que~$t\in \sch C_y(\Delta\cup\Delta')$ ; par conséquent,~$\sch C_x(\Delta\cup\Delta')$ et~$\sch C_y(\Delta\cup\Delta')$, coïncident, ce qui entraîne que~$\Theta(x)=\Theta(y)$. Ainsi {\em l'ensemble} des fermés~$\Theta(y)$, où~$y$ parcourt la fibre de~$\xi$, est-il un ensemble de parties fermées deux à deux disjointes de~$\Gamma''$. 

Par ailleurs, ces fermés recouvrent~$\Gamma''$. En effet, soit~$x'\in \Gamma''$ et soit~$\xi'$ son image sur~$\Gamma$. La composante connexe de~$x'$ dans l'image réciproque de~$[\xi';\xi]$ sur~$X$ se surjecte sur~$[\xi';\xi]$ et contient donc un antécédent~$y$ de~$\xi$ ; on a alors~$x'\in \sch C_y([\xi';\xi])\subset \Theta(y)$. 

Il s'ensuit que chacun des~$\Theta(y)$ est également un ouvert de~$\Gamma''$ ; c'est en particulier le cas de~$\Theta(x)$, qui est par ailleurs un arbre. Ainsi,~$x$ possède un voisinage ouvert dans~$\Gamma''$ qui est un arbre ; ceci valant pour tout~$x\in \Gamma''$, le fermé~$\Gamma''$ de~$X$ en est bien un sous-arbre.~$\Box$

\deux{corollstabet} {\bf Proposition.} {\em Soit~$X$ un graphe sur lequel~$\mathsf G$ agit et soit~$x$ un point de~$X$ fixe sous~$\mathsf G$. Supposons que~$\br X x$ est fini, et soit~$\Gamma$ une étoile de sommet~$x$ tracée sur~$X$ dont la valence est égale au cardinal de~$\br X x$. Les ouverts~$\mathsf G$-invariants de la forme~$\Delta^\flat$, où~$\Delta$ est une étoile appartenant à~$\stelac X x \ctd \Gamma$ et stable sous~$\mathsf G$, forment une base de voisinages de~$x$.}

\medskip
{\em Démonstration.} Soit~$V$ un voisinage ouvert de~$X$ ; il s'agit de montrer qu'il contient un voisinage de~$x$ de la forme~$\Delta^\flat$, où~$\Delta$ appartient à~$\stelac X x \ctd \Gamma$ et est stable sous~$\mathsf G$.

L'ouvert~$V$ contient une sous-étoile~$\Gamma_0$ de~$\Gamma$ ; et il existe une sous-étoile~$\Gamma_1$ de~$\Gamma_0$ qui appartient à~$\stela V x$ (\ref{exaradm}) ; on peut donc, quitte à remplacer~$\Gamma$ par~$\Gamma_1$, supposer que~$\Gamma$ appartient à~$\stela V x$. L'assertion B) du lemme~\ref{gprofpluribr} affirme l'existence d'une étoile~$\Delta$ appartenant à~$\stelac X x$, stable sous~$\mathsf G$ et satisfaisant les conditions i) et ii) de son énoncé. 

Comme~$\br \Delta x\supset \br \Gamma x$ (condition i) ) et comme~$\br X x$ a pour cardinal~$N$, on a~$\br \Delta x=\br \Gamma x$ ; il découle alors de la condition ii) que~$\Delta$ est une sous-étoile de~$\Gamma$. Comme les arêtes de~$\Gamma$ sont faiblement admissibles dans~$V$ et comme leur nombre est égal au cardinal de~$\br V x$, l'étoile~$\Gamma$ est elle-même faiblement admissible dans~$V$ (\ref{etnarnbr}), et il en va de même de son ouvert~$\Delta$ ; cela entraîne que~$\Delta^\flat\subset V$ et achève la démonstration.~$\Box$

\subsection*{Quotient d'un graphe par une relation d'équivalence d'un type particulier}

\deux{quotrelfin} {\bf Proposition.} {\em Soit~$X$ un graphe, soit~$S$ un sous-ensemble fermé et discret de~$X$ et soit~$\sch R$ une relation d'équivalence sur~$X$ dont toute classe non singleton est une partie finie de~$S$. L'espace quotient~$X/\sch R$ est alors un graphe, la flèche quotient~$p:X\to X/\sch R$ est compacte, et si~$X$ est localement fini alors~$X/\sch R$ est localement fini. De plus si~$\Gamma$ est un sous-graphe de~$X/\sch R$ alors~$p\inv(\Gamma)$ est un sous-graphe de~$X$.}

\medskip
{\em Démonstration.} On procède en plusieurs étapes. 

\trois{xsurrsep} {\em L'espace~$X/\sch R$ est séparé.} Soient~$x$ et~$y$ deux points distincts de~$X/\sch R$ ; les fibres~$p\inv(x)$ et~$p\inv(y)$ sont deux ensembles finis disjoints. Le sous-ensemble~$S$ étant fermé et discret, il existe deux ouverts~$U$ et~$V$ de~$X$ possédant les propriétés suivantes : 

\medskip
$\bullet$~$U\cap V=\emptyset$ ; 

$\bullet$~$U$ contient~$p\inv(x)$ et~$V$ contient~$p\inv(y)$ ; 

$\bullet$~$(U\cap S)\subset p\inv(x)$ et~$(V\cap S)\subset p\inv(y)$. 

\medskip
Par construction,~$U$ et~$V$ sont saturés pour~$\sch R$, et~$p(U)\cap p(V)=\emptyset$ ; dès lors~$p(U)$ et~$p(V)$ sont deux ouverts disjoints de~$X/\sch R$ dont le premier contient~$x$ et le second~$y$ ; ainsi,~$X/\sch R$ est séparé.

\trois{xsurrloccomp} {\em L'espace~$X/\sch R$ est localement compact et~$p$ est compacte.} Nous allons établir que tout point de~$\sch X/\sch R$ possède un voisinage compact dont l'image réciproque est compacte, ce qui permettra de conclure. Soit~$x\in X/\sch R$. La fibre~$p\inv(\xi)$ étant finie et~$S$ étant discret, il existe un voisinage ouvert~$U$ de~$p\inv(x)$ qui est relativement compact et tel que~$\overline U\cap S\subset p\inv(x)$. L'ouvert~$U$ est saturé pour~$\sch R$, et~$p(U)$ est dès lors un voisinage ouvert de~$x$. Il est contenu dans~$p(\overline U)$ qui est compact ; remarquons que~$\overline U$ est lui aussi saturé, ce qui entraîne que~$p\inv(p(\overline U))$ coïncide avec le compact~$\overline U$. 

\trois{xsurrlocarb} {\em L'espace~$X/\sch R$ est un arbre.} Soit~$x\in X/\sch R$ et soit~$U$ un voisinage ouvert de~$X$ dans~$X/\sch R$. Soient~$\xi_1,\ldots,\xi_r$ les antécédents de~$x$ sur~$\xi$. Il existe~$r$ sous-arbres ouverts~$V_1,\ldots,V_r$ de~$X$ possédant les propriétés suivantes :

\medskip
i) pour tout~$i$ l'ouvert~$V_i$ est un voisinage de~$\xi_i$ dont le bord est fini, et qui est réunion finie d'intervalles ouverts si~$X$ est localement fini ; 

ii) les~$V_i$ sont deux à deux disjoints ;

iii)~$S\cap (\bigcup V_i)\subset \{\xi_1,\ldots,\xi_r\}$. 

\medskip
Posons~$V=p(\bigcup V_i)$ ; par construction,~$\bigcup V_i$ est un ouvert de~$p\inv(U)$ saturé pour~$\sch R$ ; il s'ensuit que~$V$ est un voisinage ouvert de~$x$ dans~$U$, dont nous allons montrer qu'il satisfait les conditions~$\alpha)$ et~$\beta)$ du~\ref{defarb}. 

Comme~$p$ est fermée, la finitude de~$\partial V_i$ pour tout~$i$ entraîne celle de~$\partial V$. 

\medskip
Soient maintenant~$z$ et~$t$ deux points de~$V$ ; soit~$\zeta$ (resp.~$\tau$) un antécédent de~$z$ (resp.~$t$) et soit~$i_0$ (resp.~$i_1$)  l'entier tel que~$\zeta\in V_{i_0}$ (resp.~$\tau \in V_{i_1}$). Nous allons établir l'existence et l'unicité d'un segment d'extrémités~$z$ et~$t$ tracé sur~$V$. 

\medskip
{\em Existence}. Nous allons construire un segment~$I$ joignant~$z$ à~$t$. On distingue deux cas. Supposons d'abord que~$i_0=i_1$. La restriction de~$p$ à~$V_{I_0}$ est injective ; par conséquent, l'image~$I$ par~$p$ du segment~$[\zeta ; \tau]$ de l'arbre~$V_{i_0}$ est alors un segment d'extrémités~$z$ et~$t$.

Supposons maintenant que~$i_0\neq i_1$. La restriction de~$p$ à~$V_{i_0}$ (resp.~$V_{i_1}$) étant injective, l'image par~$p$ du segment~$[\zeta;\xi_{i_0}]$ de l'arbre~$V_{i_0}$ (resp. du segment~$[\tau;\xi_{i_1}]$ de l'arbre~$V_{i_1}$) est un segment~$I_0$ (resp.~$I_1$) d'extrémités~$z$ et~$ x$ (resp.~$t$ et~$x$) ; de plus, on a par construction~$I_0\cap I_1=\{x\}$ ; par conséquent, la réunion~$I:=I_0\cup I_1$ est un segment d'extrémités~$z$ et~$t$. 

\medskip
{\em Unicité.} Nous allons montrer que~$I$ est le seul segment tracé sur~$V$ et qui joigne~$z$ à~$t$. Il suffit, en vertu de~\ref{remskgrad}, de prouver que le bord dans~$V$ de toute composante connexe de~$V-I$ est un singleton. Remarquons pour commencer que si~$\Omega$ est un ouvert saturé de~$X$ alors~$(X/\sch R)-p(\Omega)=p(X-\Omega)$ ; l'application~$p$ étant fermée,~$p(X-\Omega)$ est fermé et~$p(\Omega)$ est en conséquence ouvert. 

\medskip
On définit le sous-ensemble~$\Pi$ de~$\sch P(\bigcup V_i)$ de la façon suivante, en distinguant trois cas. 

\medskip
$\bullet$  {\bf Le cas où~$i_0=i_1$ et où~$\xi_{i_0}\notin [\zeta;\tau]$.} On appelle~$W$ la composante connexe de~$V_{i_0}-[\zeta;\tau]$ contenant~$\xi_{i_0}$, et~$W'$ la réunion~$W\cup \bigcup\limits_{i\neq i_0} V_i$ ; on pose alors~$\Pi=\pi_0(V_{i_0}-[\zeta;\tau])\setminus\{W\}\cup\{W'\}$. 

$\bullet$  {\bf Le cas où~$i_0=i_1$ et où~$\xi_{i_0}\in [\zeta;\tau]$.}  On pose alors~$$\Pi=\pi_0(V_{i_0}-[\zeta;\tau])\cup \bigcup\limits_{i\neq i_0} \pi_0(V_i\setminus\{\xi_i\}).$$

$\bullet$ {\bf Le cas où~$i_0\neq i_1$.} On pose alors~$$\Pi=\pi_0(V_{i_0}-[\zeta;\xi_{i_0}])\cup \pi_0(V_{i_1}-[\xi_{i_1};\tau])\cup \bigcup\limits_{i\notin\{i_0,i_1\}} \pi_0(V_i\setminus\{\xi_i\}).$$

\medskip
L'ouvert~$V-I$ est par construction la réunion disjointes des~$p(\Omega)$, où~$\Omega$ parcourt~$\Pi$ ; ce sont tous des sous-ensembles connexes et non vides de~$V$ (en ce qui concerne~$p(W')$ dans le premier cas, cela résulte du fait qu'il est réunion de parties connexes contenant toutes~$x$). Soit~$\Omega\in \Pi$ ; c'est par sa définition même un ouvert saturé de~$V$ ; les~$p(\Omega)$ pour~$\Omega$ parcourant~$\Pi$ sont dès lors des ouverts de~$V$, et ce sont ainsi nécessairement les composantes connexes de~$V-I$. 

\medskip
On déduit du fait qu'un sous-arbre fermé d'un arbre est admissible que~$\partial_{\bigcup V_i}\Omega$ est pour tout~$\Omega\in \Pi$ un singleton ; l'application~$\bigcup V_i\to V$ induite par~$p$ est encore compacte, et partant fermée ; il s'ensuit que~$\partial_V p(\Omega)$ est un singleton pour tout~$\Omega\in \Pi$, ce qui achève la preuve de l'unicité ; ainsi,~$X/\sch R$ est un arbre.

\trois{xsurrlocfin} Si~$X$ est localement fini, il existe dans chacun des~$V_i$ ci-dessus un voisinage compact~$W_i$ de~$\xi_i$ qui est une réunion finie de segments ; soit~$W'_i$ l'intérieur de~$W_i$ dans~$X$. La flèche~$p_{|W_i}$ étant injective pour tout~$i$, l'image par~$p$ de~$\bigcup W_i$ est une réunion finie de segments ; l'ouvert~$\bigcup W'_i$ étant saturé,~$p(\bigcup W'_i)$ est un ouvert. Par conséquent,~$p(\bigcup W_i)$ est un voisinage compact de~$x$ dans le graphe~$X/\sch R$ ; ce voisinage étant réunion finie de segments,~$X/\sch R$ est localement fini.

\trois{imrecgrxsurr} Soit~$\Gamma$ un sous-graphe de~$X/\sch R$. Son image réciproque~$p\inv(\Gamma)$ est localement fermée ; pour montrer que c'est un sous-graphe de~$X$, il suffit de s'assurer que chacun de ses points possède un voisinage qui est un arbre. 

\medskip
Soit~$\xi \in p\inv(\Gamma)$ ; posons~$x=p(\xi)$, et notons~$\xi=\xi_1,\ldots,\xi_n$ les antécédents de~$\xi$ sur~$X$. Il existe une famille finie~$(V_1,\ldots,V_n)$ de sous-arbres ouverts deux à deux disjoints de~$X$ tels que~$V_i\cap S\subset \{\xi\}$ pour tout~$i$. La réunion~$V=\bigcup V_i$ est un ouvert saturé de~$X$ ; par conséquent,~$U:=p(V)$ est un voisinage ouvert de~$x$ dans~$X/\sch R$. Soit~$U'$ un voisinage compact de~$x$ dans~$U$ qui est un arbre et qui est tel que~$U'\cap \Gamma$ soit un sous-arbre compact de~$U'$ ; posons~$V'_i=p\inv(U)\cap V_i$ pour tout~$i$. Comme~$p$ est compacte, chacun des~$V'_i$ est un voisinage compact de~$\xi_i$ dans~$V_i$. Comme~$\Gamma\cap U'$ est compact,~$p\inv(\Gamma)\cap V'_i$ est compact pour tout~$i$. 

\medskip
Soient~$\eta$ et~$\tau$ deux points de~$p\inv(\Gamma)\cap V'_1$ et soient~$y$ et~$t$ leurs images respectives sur~$X$. Par définition des~$V_i$, la restriction de~$p$ à~$V_1$ est injective ; elle induit donc un homéomorphisme entre le segment~$[\eta;\tau]~$ de l'arbre~$V_1$ et son image sur~$U'$, qui coïncide dès lors avec~$[z;t]$. Puisque~$\Gamma\cap U'$ est un sous-arbre de~$U'$, le segment~$[z;t]$ est contenu dans~$\Gamma\cap U'$ ; par conséquent,~$[\eta;\tau]\subset p\inv(\Gamma)\cap V'_1$. Ainsi, le voisinage~$p\inv(\Gamma)\cap V'_1$ de~$\xi=\xi_1$ dans~$p\inv(\Gamma)$ apparaît-il comme une partie convexe et compacte de l'arbre~$V_1$ ; c'est donc un arbre, ce qui achève la démonstration.~$\Box$ 

\subsection*{Quotients et toises}

\deux{transfparfacile} {\bf Lemme.} {\em Soit~$\phi : Y\to X$ une application continue et compacte entre graphes, telle que la restriction de~$\phi$ à tout sous-graphe localement fini de~$Y$ soit injective par morceaux. Supposons que tout sous-graphe compact de~$X$ admette une toise ; il en va alors de même de tout sous-graphe compact de~$Y$.}

\medskip
{\em Démonstration.} Soit~$V$ un sous-graphe compact de~$Y$. Son image étant compacte, elle est contenue dans un sous-graphe compact~$U$ de~$X$. Par hypothèse, il existe une toise~$\lambda$ sur~$U$.

\medskip
Soit~$I$ un segment tracé sur~$V$ et soient~$x$ et~$y$ ses deux extrémités. Orientons~$I$ de~$x$ vers~$y$. La restriction de~$\phi$ à~$I$ étant injective par morceaux, il existe une suite finie~$x_0=x<x_1<\ldots<x_n=y$ de points de~$I$ tels que~$\phi$ soit injective sur~$[x_i;x_{i+1}]$ pour tout~$i\in \{0,\ldots,n-1\}$ ; si~$0\leq i\leq n-1$ alors~$\phi([x_i;x_{i+1}])$ est un segment tracé sur~$U$. 

\medskip
Le réel~$\sum\limits_i  \lambda(\phi([x_i;x_{i+1}]))$ ne dépend pas du choix de la subdivision~$(x_i)$ : il suffit en effet de vérifier qu'il est invariant par raffinement de cette dernière, ce qui est immédiat ; il est donc licite de le noter~$l(I)$. 

\medskip
Par sa construction même,~$l$ est une toise sur~$V$.~$\Box$

\deux{transflong}{\bf Proposition.} {\em Soit~$X$ un graphe sur lequel~$\mathsf G$ agit.  Les assertions suivantes sont équivalentes : 

\medskip
i) tout sous-graphe compact de~$X$ admet une toise ; 

ii) tout sous-graphe compact de~$X/\mathsf G$ admet une toise. }

\medskip
{\em Démonstration.} On note~$p$ la flèche quotient~$X\to X/\mathsf G$. Comme la restriction de~$p$ à tout sous-graphe localement fini de~$X$ est, en vertu du théorème~\ref{theoquot}, injective par morceaux, l'implication ii)$\Rightarrow$i) est une conséquence directe du lemme~\ref{transfparfacile}. Il reste à montrer que i)$\Rightarrow$ ii) ; on suppose donc que tout sous-graphe compact de~$X$ admet une toise. Soit~$U$ un sous-graphe compact de~$X$. 

\trois{longimpllongquot} {\em Modification du problème à traiter.} Comme~$U$ est compact, il existe une famille finie~$(X'_i,X_i)_i$ où~$X_i$ est pour tout~$i$ un sous-arbre ouvert de~$X/\mathsf G$, où~$X'_i$ est pour tout~$i$ un sous-arbre compact de~$X'_i$, et où les~$X'_i$ recouvrent~$U$. Il suffit de montrer que chacun des~$X_i$ admet une toise : en effet si c'est le cas alors par restriction chacun des~$X'_i\cap U$ admettra une toise, et il n'y aura plus qu'à appliquer~\ref{paramfini} au recouvrement fini de~$U$ par les~$U\cap X'_i$. 

\medskip
L'application~$p$ étant compacte,~$p\inv(X_i)$ est pour tout~$i$ un sous-graphe ouvert relativement compact de~$X$ ; il admet donc une toise d'après l'hypothèse faite sur~$X$.

\medskip
On s'est ainsi ramené à établir l'assertion suivante : {\em soit~$X$ un graphe sur lequel~$\mathsf G$ agit ; supposons que~$X$ admet une toise et que~$X/\mathsf G$ est un arbre ; le quotient~$X/\mathsf G$ admet alors une toise.} 

\medskip
C'est elle que nous entendons désormais démontrer ; on fixe une toise~$l$ sur~$X$ ; on note~$X_{\rm pluri}$ (resp.~$(X/\mathsf G)_{\rm pluri}$) l'ensemble des points pluribranches de~$X$ (resp. de~$X/\mathsf G$).

\trois{toisepluribr} {\em Construction d'une toise~$l'$ sur~$(X/\mathsf G)_{\rm pluri}$.} Soit~$I$ un segment tracé sur ~$(X/\mathsf G)_{\rm pluri}$ et soient~$x$ et~$y$ les extrémités de~$I$. En vertu de~\ref{uniuni} et~\ref{isolisol},~$p\inv(I)\subset X_{\rm pluri}$. On déduit dès lors de la compacité de~$I$ et du corollaire~\ref{corollimreci} l'existence d'une suite finie~$x=x_0<x_1<\ldots<x_r=y$ d'éléments de~$I$ (orienté de~$x$ vers~$y$) telle que pour tout~$i$ compris entre~$0$ et~$r-1$, l compact~$p\inv([x_i;x_{i+1}])$ s'écrive comme une réunion finie~$\bigcup\limits _ j I_{i,j}$, où~$I_{i,j}$ est pour tout~$(i,j)$ un segment s'envoyant homéomorphiquement sur~$[x_i;x_{i+1}]$, et où les intérieurs des~$I_{i,j}$ sont deux à deux disjoints. 

\medskip
Pour tout~$i\in \{0,\ldots,r-1\}$ notons~$l_i$ la moyenne des~$l(Y_{i,j})$ pour~$j$ variable. La somme des~$l_i$ ne dépend que de~$I$, et pas de la subdivision~$(x_i)$ : il suffit en effet de s'assurer qu'elle est insensible à un raffinement de cette dernière, ce qui est immédiat ; il est dès lors licite de poser~$l'(I)=\sum l_i$. Par sa construction même,~$l'$ définit une toise sur~$(X/\mathsf G)_{\rm pluri}$. 

\trois{conclutoise} {\em Conclusion.} Comme tout point de~$X/\mathsf G$ situé sur un intervalle ouvert appartient à~$(X/\mathsf G)_{\rm pluri}$, ce dernier ensemble est une partie convexe de l'arbre~$X/\mathsf G$. On déduit alors du~\ref{toisepluribr} ci-dessus et de~\ref{yconvmborn} qu'il existe une toise~$l''$ sur~$(X/\mathsf G)_{\rm pluri}$ qui est bornée par~$1$. 

\medskip
Soit~$I$ un segment tracé sur~$X/\mathsf G$. Notons~$\lambda(I)$ la borne supérieure des~$l''(J)$, où~$J$ parcourt l'ensemble des segments contenus dans l'intérieur de~$I$ (qui est constitué de points pluribranches) ; l'application~$I\mapsto \lambda(I)$ est alors une toise sur~$X/\mathsf G$.~$\Box$ 

\deux{corolltoisexsurg} {\bf Corollaire.} {\em Soit~$X$ un graphe sur lequel~$\mathsf G$ agit. Le quotient~$X/\mathsf G$ admet une toise si et seulement si~$X$ admet une toise}.

\medskip
{\em Démonstration.} Comme~$X\to X/\mathsf G$ est compacte,~$X$ est paracompact si et seulement si~$X/\mathsf G$ est paracompact. Le corollaire est dès lors une conséquence immédiate de la proposition~\ref{transflong} ci-dessus et du lemme~\ref{xadmtoise}.~$\Box$

\deux{transfr}{\bf Proposition.} {\em Soit~$X$ un graphe, soit~$S$ un sous-ensemble fermé et discret de~$X$ et soit~$\sch R$ une relation d'équivalence sur~$X$ dont toute classe non singleton est une partie finie de~$S$. Les assertions suivantes sont équivalentes : 

\medskip
i) tout sous-graphe compact de~$X$ admet une toise ; 
ii) tout sous-graphe compact du graphe~$X/\sch R$ ({\em cf.} prop.\ref{quotrelfin}) admet une toise. }

\medskip
{\em Démonstration.} Soit~$p$ la projection~$X\to X/\sch R$. Il résulte de la forme même de~$\sch R$ que la restriction de~$p$ à tout sous-graphe localement fini de~$X$ est injective par morceaux. Par conséquent, ii)$\Rightarrow$i) résulte du lemme~\ref{transfparfacile}. Supposons maintenant que i) est vraie, et montrons ii). 

\medskip
Soit~$U$ un sous-graphe compact de~$X$. Par compacité de~$U$ et de~$p$, il existe une famille finie~$(X_i)$ de sous-arbres compacts de~$X/\sch R$ possédant les propriétés suivantes : 

\medskip
$\bullet$ les~$X_i$ recouvrent~$U$ ; 

$\bullet$ pour tout~$i$, l'image réciproque de~$X_i$ s'écrit comme une union finie disjointe~$\coprod X'_{ij}$, où~$X'_{ij}$ est pour tout~$j$ un compact comprenant au plus un élément de~$S$. 

\medskip
Fixons~$i$ et~$j$. Comme~$p\inv(X_i)$ est un sous-graphe compact de~$X$ (prop.~\ref{quotrelfin}), son ouvert fermé~$X'_{ij}$ est un sous-graphe compact de~$X$. Comme~$X'_{ij}$ contient au plus un élément de~$S$, l'application~$p$ induit un homéomorphisme de~$X'_{ij}$ sur son image~$X_{ij}$, qui est donc un sous-graphe compact de~$X/\sch R$. 

\medskip
En vertu de notre hypothèse, le sous-graphe compact~$X'_{ij}$ de~$X$ admet une toise ; par conséquent,~$X_{ij}$ qui lui est homéomorphe en admet une aussi. Le sous-graphe~$U$ de~$X$ est recouvert par ses sous-graphes fermés~$U\cap X_{ij}$, qui admettent chacun une toise par restriction. On déduit alors de~\ref{paramfini} que~$U$ admet une toise.~$\Box$ 

\deux{corolltoisexsurr} {\bf Corollaire.} {\em Soient~$X$ et~$\sch R$ comme dans la proposition~\ref{transfr} ci-dessus. Le quotient~$X/\sch R$ admet une toise si et seulement si~$X$ admet une toise}.

\medskip
{\em Démonstration.} Comme~$X\to X/\sch R$ est compacte d'après la proposition~\ref{quotrelfin}, le graphe~$X$ est paracompact si et seulement si~$X/\sch R$ est paracompact. Le corollaire est dès lors une conséquence immédiate de la proposition~\ref{transfr} ci-dessus et du lemme~\ref{xadmtoise}.~$\Box$

\section{Espaces dont un quotient d'un certain type est un graphe}

\deux{theofondquot} {\bf Théorème.} {\em Soit~$\mathsf G$ un groupe {\em fini} et soit~$Y$ un espace topologique séparé, localement connexe par arcs et localement compact sur lequel~$\mathsf G$ agit. Supposons que les deux conditions suivantes soient satisfaites : 

\medskip
i) il existe un entier~$\ell\geq 2$ tel que~$\H^1(Y,\ZZ/\ell\ZZ)$ soit fini ;  

ii)~$Y/\mathsf G$ est un graphe. 

\medskip
L'espace~$Y$ est alors lui-même un graphe.}

\medskip
{\em Démonstration.} On peut supposer, en raisonnant composante par composante, que~$Y$ est connexe. Désignons par~$X$ l'espace quotient~$Y/\mathsf G$  (qui est alors lui aussi connexe), par~$\phi$ la flèche quotient~$Y\to X$ (qui est ouverte, propre et à fibres finies) et par~$d$ le cardinal maximal d'une famille libre du~$\ZZ/\ell\ZZ$-module~$\H^1(Y,\ZZ/\ell\ZZ)$ (qui est fini d'après i) ).

\trois{boucle} Soit~$C$ un cercle tracé sur~$Y$ et soit~$J$ un ouvert de~$C$ possédant la propriété suivante :

$(*)$ {\em il existe un intervalle ouvert~$I$ non vide tracé sur~$X$ tel que~$J$ soit un ouvert fermé de~$\phi\inv(I)$ s'envoyant homéomorphiquement sur~$I$.}

Nous allons démontrer qu'il existe alors un~$\ZZ/\ell\ZZ$-revêtement topologique~$Z$ de~$Y$ déployé sur~$Y-J$ et tel que la classe de~$Z\times_Y C$ dans~$\H^1(C,\ZZ/\ell \ZZ)$ engendre ce dernier. 

\medskip
{\em Simplification du contexte}. Pour cela, choisissons un point~$y$ sur~$J$ ; son image~$x$ sur~$I$ est un point pluribranche de~$X$. Il existe un voisinage ouvert~$U$ de~$x$ dans~$X$ qui possède les propriétés suivantes : 

\medskip
$\bullet$~$U$ est un arbre et~$U\cap I$ est un intervalle ouvert ; 

$\bullet$ si~$V$ désigne la composante connexe de~$y$ dans~$\phi\inv(U)$, alors~$V\cap C\subset J$ et~$V\cap \phi\inv(I)\subset J$. 

\medskip
Choisissons un intervalle ouvert~$J'$ de~$J$ contenant~$y$ et inclus dans~$V$, et notons~$I'$ l'image de~$J'$ sur~$I$. Il existe un intervalle ouvert~$I''$ de~$I'$ aboutissant à~$x$ qui est faiblement admissible dans~$X$ et tel que l'arbre à deux bouts~$\Omega:=(I'')^\flat$ soit contenu dans~$U$ (\ref{declun}). Notons que~$\Omega\cap I=\Omega\cap U\cap I$ est par convexité un intervalle ouvert, dont~$I''$ est une partie fermée puisque~$\partial I''$ ne rencontre pas~$\Omega$ ; par conséquent,~$\Omega\cap I=I''$. 

 Soit~$J''$ l'image réciproque de~$I''$ sur~$J$ et soit~$W$ la composante connexe de~$\phi\inv(\Omega)$ qui contient~$J''$. Comme~$J''\subset J'\subset V$, l'on a~$W\subset V$,  et donc~$W\cap C\subset J$ et~$W\cap \phi\inv(I)\subset J$ ; comme~$\Omega\cap I=I''$, il s'ensuit que~$W\cap C=J''$ et~$W\cap \phi\inv(I'')=J''$.

\medskip
Quitte à remplacer~$I$ par~$I''$ et~$J$ par~$J''$, l'on peut donc supposer {\em qu'il existe un ouvert~$\Omega$ de~$X$ qui est un arbre à deux bouts  dont~$I$ est le squelette et une composante connexe~$W$ de~$\phi\inv(\Omega)$ telle que~$W\cap C=J$ et~$W\cap \phi\inv(I)=J$}. On fixe deux intervalles ouverts~$I_0$ et~$I_1$ de~$I$, disjoints et aboutissant chacun à l'une de ses extrémités ; l'on note~$J_0$ et~$J_1$ les images réciproques respectives de~$I_0$ et~$I_1$ sur~$J$, et~$\Omega_0$ et~$\Omega_1$ les images réciproques respectives de~$I_0$ et~$I_1$ par la rétraction canonique~$\Omega\to I$ ; enfin, l'on désigne par~$W_0$ et~$W_1$ les images réciproques respectives de~$\Omega_0$ et~$\Omega_1$ sur~$W$. 

\medskip
{\em Construction d'un~$\ZZ/\ell\ZZ$-torseur~$Z\to Y$.} Soit~$\sch K$ l'image réciproque sur~$W$ du compact~$\Omega-(\Omega_0\cup \Omega_1)$ ; c'est un compact de~$Y$. L'on construit alors un~$\ZZ/\ell\ZZ$-revêtement topologique~$Z$ de~$Y$ en recollant~$(Y-{\sch K})\times \ZZ/\ell\ZZ$ et~$W\times \ZZ/\ell\ZZ$ comme suit : au-dessus de~$W_0$, l'on identifie pour tout~$i$ le feuillet~$W_0\times \{i\}$ du premier revêtement au feuillet~$W_0\times \{i\}$ du second ; au-dessus de~$W_1$, l'on identifie pour tout~$i$ le feuillet~$W_1\times\{i\}$ du premier revêtement au feuillet~$W_1\times\{i+1\}$ du second.

\medskip
{\em Étude du torseur~$Z\times_yC$}. Il découle de notre construction, du fait que~$W\cap C$ est égal à~$J$, et de la description explicite des revêtements du cercle que la classe de~$Z\times_YC$ dans~$\H^1(C,\ZZ/\ell\ZZ)$ engendre ce dernier. 

\medskip
{\em Trivialisation de~$Z$ au-dessus de~$Y-J$}. Soit~$W'$ une composante connexe de~$W-J$. Comme~$J$ est l'image réciproque de~$I$ sur~$W$, l'ouvert~$W'$ est une composante connexe de~$\phi\inv(\Omega')$ pour une certaine composante connexe~$\Omega'$ de~$\Omega-I$. L'espace~$\Omega$ étant un arbre de squelette~$I$, le bord de~$\Omega'$ dans~$\Omega$ consiste en un et un seul point de~$I$ ; la flèche~$J\to I$ étant un homéomorphisme, le bord de~$W'$ dans~$W$ consiste en un et un seul point de~$J$. Soit~$i\in\{0,1\}$. Puisque~$\Omega'$ est contenue dans~$\Omega_i$ si et seulement si l'unique point de son bord appartient à~$I_i$, la composante~$W'$ est contenue dans~$W_i$ si et seulement si l'unique point de son  bord appartient à~$J_i$.

\medskip
Soit~$\sigma$ l'application de~$Y-J$ dans~$Z$ définie comme suit : si~$y$ appartient à~$Y-J\setminus \sch K$ l'on pose~$\sigma(y)=(y,0)\in (Y-{\sch K})\times \ZZ/\ell \ZZ\subset Z$ ; si~$y\in W-J$ alors~$y$ appartient à une composante connexe~$W'$ de~$W-J$ ; si l'unique point du bord de~$W'$ dans~$W$ appartient (resp. n'appartient pas) à~$J_0$, l'on définit~$\sigma(y)$ comme l'élément~$(y,0)$ (resp.~$(y,1)$) de~$W\times \ZZ/\ell\ZZ\subset Z$.  Par construction,~$\sigma$ est une section continue de~$Z\times_Y(Y-J)\to Y-J$ ; dès lors,~$Z$ se trivialise au-dessus de~$Y-J$. 

\trois{remapresboucle} {\em Remarque à propos de la propriété (*) ci-dessus.} Soit~$C$ un cercle tracé sur~$Y$ ; supposons qu'il existe un graphe compact et fini~$\Gamma$ sur~$X$ tel que~$\phi\inv(\Gamma)$ soit un graphe compact et fini contenant~$C$ ; tout ouvert non vide de~$C$ contient alors un ouvert~$J$ satisfaisant~$(*)$ : cela résulte du fait que la flèche~$\phi\inv(\Gamma)\to \Gamma\simeq \phi\inv(\Gamma)/\mathsf G$ est alors injective par morceaux d'après le théorème~\ref{theoquot}. 

\trois{imrecarb} Soit~$\Gamma$ un graphe compact et fini tracé sur~$X$, soit~$\Gamma_0$ le complémentaire d'un nombre fini de points dans~$\Gamma$ et soit~$N$ un entier non nul. Supposons que les fibres de~$\phi$ en les points de~$\Gamma_0$ sont toutes de cardinal~$N$ ; nous allons montrer que~$\phi\inv(\Gamma_0)\to \Gamma_0$ est un revêtement à~$N$ feuillets, puis que~$\phi\inv(\Gamma)$ est un graphe fini. On peut supposer pour ce faire que~$\Gamma$ est sans point isolé, auquel cas il coïncide avec l'adhérence de~$\Gamma_0$. 

\medskip
Soit~$x\in \Gamma_0$ et soient~$x_1,\ldots,x_N$ ses antécédents. Soient~$V_1,\ldots,V_N$ des voisinages ouverts respectifs de~$x_1,\ldots,x_N$ dans~$\phi\inv(\Gamma)$, que l'on choisit deux à deux disjoints, et soit~$W$ un voisinage {\em compact} de~$x$ dans~$\Gamma_0$ qui est contenu dans~$\bigcap \phi(V_i)$. Désignons pour tout~$i$ par~$W_i$ l'image réciproque de~$W$ dans~$V_i$ ; c'est un fermé de~$\phi\inv(W)$ et donc un compact ; par choix de~$W$, les flèches~$W_i\to W$ sont toutes surjectives ; comme tout point de~$\Gamma_0$ a exactement~$N$ antécédents, elles sont toutes bijectives ; par compacité de leurs sources et buts, ce sont des homéomorphismes, et~$\phi\inv(\Gamma_0)\to \Gamma_0$ est bien un revêtement à~$N$ feuillets.

\medskip
Pour montrer que~$\phi\inv(\Gamma)$ est un arbre on peut, quitte à amputer~$\Gamma_0$ d'un ensemble fini de points, supposer que l'adhérence dans~$\Gamma$ de toute composante connexe de~$\Gamma_0$ est un arbre. Soit~$\Delta$ une composante connexe de~$\phi\inv(\Gamma_0)$ ; l'application~$\Delta\to \phi(\Delta)$ est un homéomorphisme, et~$\phi(\Delta)$ est une composante connexe de~$\Gamma_0$ ; dès lors,~$\phi(\Delta)$ et~$\overline{\phi(\Delta)}$ sont deux arbres. Comme~$\phi(\overline \Delta)$ est un fermé de~$X$ contenu dans le compact~$\overline{\phi(\Delta)}$, c'est exactement~$\overline{\phi(\Delta)}$. Soit~$a\in \overline{\phi(\Delta)}-\phi(\Delta)$. Il existe un voisinage ouvert~$U$ connexe de~$a$ dans~$X$ tel que les composantes connexes de~$\phi\inv(U)$ séparent les antécédents de~$a$ ; comme~$\phi(\Delta)$ est un arbre,~$\overline{\phi(\Delta)}\setminus\{a\}$ est connexe, et~$\overline{\phi(\Delta)}$ est donc de valence~$1$ en~$a$. Il est par conséquent loisible de supposer, quitte à restreindre~$U$, que~$U\cap \overline{\phi(\Delta)}$ est de la forme~$[a;x[$ pour un certain~$x$ tel que~$]a;x[\subset\phi(\Delta)$. L'intersection de~$\Delta$ et~$\phi\inv(U)$ est dès lors égale à~$\phi_{|\Delta}\inv(]a;x[)$, qui est homéomorphe à~$]a;x[$ et est de ce fait contenu dans une composante connexe de~$\phi\inv(U)$ ; par conséquent,~$\overline\Delta$ contient au plus un antécédent de~$a$. Ceci étant vrai quel que soit le point~$a$ de~$\overline{\phi(\Delta)}-\phi(\Delta)$, la surjection~$\overline \Delta \to \overline{\phi(\Delta)}=\phi(\overline \Delta)$ est une bijection continue entre compacts, et donc un homéomorphisme ; on en déduit que~$\overline \Delta$ est un arbre compact ; ceci vaut pour tout~$\Delta\in \pi_0(\phi\inv(\Gamma_0))$.

\medskip
Comme~$\phi\inv(\Gamma)\to \Gamma$ est ouverte et comme~$\Gamma_0$ est dense dans~$\Gamma$, l'ouvert~$\phi^{1}(\Gamma_0)$ de~$\phi\inv(\Gamma)$ en est une partie dense. Autrement dit, l'on a~$$\phi\inv(\Gamma)=\overline{\phi\inv(\Gamma_0)}=\bigcup_{\Delta \in\pi_0(\phi\inv(\Gamma_0))}  \overline \Delta.$$ Étant compact,~$\phi\inv(\Gamma)$ est donc homéomorphe au quotient de~$\coprod\limits_{\Delta \in\pi_0(\phi\inv(\Gamma_0))}  \overline \Delta$ par une relation d'équivalence convenable, qu'engendre un ensemble fini de couples de points ; il s'ensuit que~$\phi\inv(\Gamma)$ est un graphe fini (prop.~\ref{quotrelfin}).

\trois{imrecinterv} Supposons donnés un entier~$N$ et, pour tout~$i$ compris entre~$1$ et~$N$, un intervalle~$[a_i; b_i]$ tracé sur~$X$ (avec~$a_i\neq b_i$) ; nous faisons les hypothèses suivantes : 

\medskip
$\bullet$ les intervalles~$[a_i;b_i]$ sont deux à deux disjoints ; 

$\bullet$ pour tout~$i$, les fibres de~$\phi$ en les points de~$]a_i;b_i[$ ont toutes même cardinal, et ce dernier majore strictement les cardinaux de~$\phi\inv(a_i)$ et~$\phi\inv(b_i)$.

\medskip
Nous allons montrer que sous ces conditions~$\phi\inv([a_i;b_i])$ est pour tout~$i$ un graphe fini ayant au moins une boucle, puis que~$N$ est nécessairement inférieur ou égal à~$d$. 
 
\medskip
Fixons~$i$. Il résulte du~\ref{imrecarb} ci-dessus que ~$\phi\inv(]a_i;b_i[)\to ]a_i;b_i[$ est un revêtement, et que~$\phi\inv([a_i;b_i])$ est un graphe compact, que l'on note~$\Gamma$. 

\medskip
Nous allons montrer que~$\Gamma$ possède au moins une boucle. Le cardinal de~$\phi\inv(a_i)$ étant par hypothèse strictement inférieur au nombre de feuillets du revêtement~$\phi\inv(]a_i;b_i[)\to ]a_i;b_i[$, il existe un antécédent de~$a_i$ auquel aboutissent au moins deux branches de~$\Gamma$. Par homogénéité (l'on a~$[a_i;b_i]\simeq \Gamma/\mathsf G)$, c'est    le cas de {\em tous} les antécédents de~$a_i$ ; il en va de même de ceux de~$b_i$. Soit~$\Gamma_0$ une composante connexe de~$\Gamma$ et soit~$n$ le nombre de feuillets de~$\phi\inv(]a_i;b_i[)$ qu'elle contient. À chaque point de~$\Gamma_0$ situé au-dessus de~$a_i$ aboutissent au moins deux de ces feuillets. Il y a donc au plus~$n/2$ tels points ; de même, il y a au plus~$n/2$ antécédents de~$b_i$ sur~$\Gamma_0$. Il en résulte que la caractéristique d'Euler de~$\Gamma_0$ est majorée par~$-n+2n/2=0$ ; par conséquent,~$h^1(\Gamma_0)$ est non nul, et~$\Gamma$ possède ainsi au moins une boucle. 

\medskip
Pour tout~$i$, choisissons une boucle~$C_i$ sur le graphe~$\phi\inv([a_i;b_i])$. On déduit de la remarque~\ref{remapresboucle} et de la construction faite au~\ref{boucle} qu'il existe pour chaque~$i$ un~$\ZZ/\ell\ZZ$-revêtement topologique~$Z_i$ de~$Y$ trivial en dehors de~$C_i$ tel que la classe de~$Z_i\times_YC_i$ engendre~$\H^1(C_i,\ZZ/\ell \ZZ)$, lequel est isomorphe à ~$\ZZ/\ell \ZZ$ ; il en découle aisément que les classes des~$Z_i$ dans~$\H^1(Y,\ZZ/\ell\ZZ)$ forment une famille libre, et donc que~$N\leq d$.

\trois{imrecgrfini} Soit~$\Gamma$ un graphe fini et compact tracé sur~$X$ ; nous allons établir que~$\phi\inv(\Gamma)$ est un graphe fini et que~$h^1(\phi\inv(\Gamma))\leq d$.

\medskip
Soit~$N$ le cardinal maximal d'une fibre de~$\phi$ au-dessus d'un point de~$\Gamma$ ; montrons par récurrence sur~$N$ que~$\phi\inv(\Gamma)$ est un graphe fini. 

Si~$N=1$ alors~$\phi\inv(\Gamma)\to \Gamma$ est une bijection continue entre compacts, et donc un homéomorphisme ; par conséquent,~$\phi\inv(\Gamma)$ est un graphe fini. 

Supposons que~$N>1$ et que la propriété est vraie pour tous les entiers compris entre~$1$ et~$N-1$. Soit~$U$ l'ensemble des points de~$\Gamma$ dont les fibres ont exactement~$N$ points ; il résulte du caractère ouvert et propre de~$\phi$ que~$U$ est un ouvert. Sa trace sur chaque arrête ouverte de~$\Gamma$ est une réunion (au plus) dénombrable d'intervalles ouverts ; en vertu de~\ref{imrecinterv}, cette réunion est en réalité finie. Par conséquent :

\medskip
- le complémentaire de~$U$ dans~$\Gamma$ est un graphe fini ;

- l'adhérence de~$U$ dans~$\Gamma$ est un graphe fini~$\Delta$, et~$\Delta-U$ est fini.

\medskip
Il découle de la définition de~$U$ et de~\ref{imrecarb} que~$\phi\inv(\Delta)$ est un graphe fini. Par ailleurs les fibres de~$\phi$ en les points de~$\Gamma-U$ sont toutes de cardinal strictement inférieur à~$N$ ; par l'hypothèse de récurrence,~$\phi\inv(\Gamma-U)$ est un graphe fini. Comme~$\phi\inv(\Gamma)$ est compact, il s'identifie au quotient de~$\phi\inv(\Gamma-U)\coprod \phi\inv(\Delta)$ par une relation d'équivalence convenable, qu'engendre un ensemble fini de couples de points ; en conséquence,~$\phi\inv(\Gamma)$ est un graphe fini (prop.~\ref{quotrelfin}). 

\medskip
Posons~$r=h^1(\phi\inv(\Gamma))$. Il existe~$r$ boucles ~$C_1,\ldots,C_r$ sur~$\phi\inv(\Gamma)$ et, pour tout~$i\in \{1,\ldots r\}$. On déduit de la remarque~\ref{remapresboucle} et du~\ref{boucle} qu'il existe pour tout~$i$ un intervalle ouvert~$I_i$ non vide tracé sur~$C-\bigcup\limits_{j\neq i}C_j$ et un~$\ZZ/\ell\ZZ$-revêtement~$Z_i$ de~$Y$ qui est trivial au-dessus de~$Y-I_i$, et tel que~$Z_i\times _Y C_i$ engendre~$\H^1(C_i,\ZZ/\ell\ZZ)$. Il s'ensuit que les classes des~$Z_i$ dans~$\H^1(Y,\ZZ/\ell \ZZ)$ forment une famille libre ; en conséquence,~$r\leq d$. 

\trois{graphstab} Soient~$\Gamma$ et~$\Gamma'$ deux graphes finis compacts et connexe tracés sur~$X$ tels que~$\Gamma$ soit contenu dans~$\Gamma'$. L'on a~$d\geq h^1(\phi\inv(\Gamma)')\geq h^1(\phi\inv(\Gamma))$ ; comme toute composante connexe de~$\phi\inv(\Gamma')$ se surjecte sur~$\Gamma'$, l'application naturelle de~$\pi_0(\phi\inv(\Gamma))$ vers~$\pi_0(\phi\inv(\Gamma'))$ est surjective ; notons que sa source et sont buts sont finis et de cardinal majoré par celui de~$\mathsf G$. 

\medskip
Il en résulte qu'il existe un graphe fini compact et connexe~$\Gamma_0$ tracé sur~$X$ tel que pour tout couple~$(\Gamma,\Gamma')$ de graphes finis, compacts et connexes tracés sur~$X$, contenant~$\Gamma_0$ et tels que~$\Gamma\subset\Gamma'$, les flèches~$$\pi_0(\phi\inv(\Gamma))\to \pi_0(\phi\inv(\Gamma)')\;{\rm et}\;\H_1(\phi\inv(\Gamma'),\ZZ)\to \H_1(\phi\inv(\Gamma),\ZZ)$$ soient bijectives ; cela signifie que~$\phi\inv(\Gamma)$ est un sous-graphe admissible de~$\phi\inv(\Gamma')$.

\medskip
Notons~$\got G$ l'ensemble des graphes finis compacts et connexes tracés sur~$X$ et contenant~$\Gamma_0$ ; il est naturellement ordonné et filtrant ; l'on a~$X=\bigcup\limits_{\Gamma \in \got G}\Gamma$. 

\trois{conngraphe} Soient~$y$ et~$y'$ deux points distincts d'un ouvert connexe~$\Omega$ de~$Y$. Nous allons démontrer qu'il existe~$\Gamma\in \got G$ tel que~$y$ et~$y'$ soient situés sur la même composante connexe de~$\phi\inv(\Gamma)\cap \Omega$. 

\medskip
Comme~$Y$ est localement connexe par arcs, il existe un chemin continu~$[0;1]\to \Omega$ joignant~$y$ à~$y'$ ; comme~$Y$ est séparé, l'on peut supposer que ce chemin est {\em injectif}, autrement dit que c'est un homéomorphisme sur son image ; soit~$I$ cette dernière. Nous allons montrer que pour tout point~$t$ de~$I$ il existe~$\Gamma_t\in \got G$ et un voisinage ouvert~$I_t$ de~$t$ dans~$I$ tel que~$I_t\subset \phi\inv(\Gamma_t)$. La compacité de~$I$ permettra alors de le recouvrir par un nombre fini de~$I_t$, et si~$\Gamma$ désigne un élément de~$\got G$ contenant les~$\Gamma_t$ correspondant, l'on aura~$I\subset \phi\inv(\Gamma)\cap \Omega$ et la conclusion requise. 

\medskip
Soit~$t\in I$ et soit~$x$ son image sur~$X$. Il existe un voisinage ouvert~$U$ de~$x$ qui est un arbre, et qui est tel que la composante connexe~$V$ de~$t$ dans~$\phi\inv(U)$ ne contienne aucun autre antécédent de~$x$. Fixons une branche de~$I$ issue de~$t$ (il y en a au plus deux) et soit~$\tau$ un point de cette branche, tel que~$[t;\tau]\subset V$ ; posons~$\xi=\phi(\tau)$. Comme~$\phi([t;\tau])$ est un sous-ensemble connexe de l'arbre~$U$, il contient l'intervalle~$[x;\xi]$. La composante connexe~$\Delta$ de~$t$ dans le graphe~$\phi\inv([x;\xi])$ est contenue dans~$V$, et coïncide avec~$\phi_{|V}\inv([x;\xi])$. Soit~$U'$ une composante connexe de~$U-[x;\xi]$ aboutissant à un point~$\zeta$ de~$]x;\xi[$ et soit~$J$ une composante connexe de~$\phi_{|]t;\tau[}\inv(U')$ ; c'est un intervalle ouvert tracé sur~$]t;\tau[$ ; notons~$\alpha$ et~$\beta$ ses deux bornes ; l'on a~$\phi(\alpha)=\phi(\beta)=\zeta$. La concaténation de~$J$ et de l'un des intervalles reliant~$\alpha$ à~$\beta$ sur~$\Delta$ est un cercle~$C$. L'on déduit du~\ref{boucle} qu'il existe un~$\ZZ/\ell\ZZ$-torseur~$Z\to Y$ qui est trivial sur~$Y-J$ et tel que la classe de~$Z\times_Y C$ engendre~$\H^1(C,\ZZ/\ell \ZZ)$. 

\medskip
Soit~$r\in\NN$. Supposons qu'il existe pour tout~$i\in \{1,\ldots, r\}$ une composante connexe~$U_i$ de~$U-[x;\xi]$ aboutissant à un point de~$]x;\xi[$ et une composante connexe~$J_i$ de~$\phi_{|]t;\tau[}\inv(U_i)$, les~$J_i$ étant en outre supposées deux à deux disjointes. Le procédé ci-dessus fournit pour tout~$i$ un~$\ZZ/\ell\ZZ$-torseur~$Z_i\to Y$ trivial au-dessus de~$Y-J_i$ et un cercle~$C_i$ ne rencontrant aucun des~$J_j$ pour~$j\neq i$ et tel que~$Z_i\times_ YC_i$ engendre~$\H^1(C_i,\ZZ/\ell \ZZ)$. Les classes des~$Z_i$ dans~$H^1(Y,\ZZ/\ell\ZZ)$ sont alors linéairement indépendantes ; par conséquent,~$r\leq d$. 

\medskip
Ceci implique l'existence de~$\xi_0\in ]x;\xi[$ tel que~$\phi(]t;\tau[)$ ne rencontre aucune composante connexe de~$U-[x;\xi]$ aboutissant à un point de~$]x;\xi_0[$. Soit~$\tau_0$ un antécédent de~$\xi_0$ sur~$]t;\tau[$. Comme~$x$ et~$\xi_0$ adhèrent tous deux à~$\phi(]t;\tau_0[)$, la composante connexe~$U_0$ de~$U\setminus\{x;\xi_0\}$ qui contient~$\phi(]t;\tau_0[)$ est celle qui contient~$]x;\xi_0[$ ; par ce qui précède,~$\phi(]t;\tau_0[)$ ne rencontre aucune composante connexe de~$U_0-]x;\xi_0[$ ; dès lors,~$\phi(]t;\tau_0[)\subset ]x;\xi_0[$ et l'on a donc~$\phi([t;\tau_0[)\subset [x;\xi_0[$. 

\medskip
En raisonnant de même le cas échéant sur la seconde branche de~$I$ issue de~$t$, on obtient l'existence d'un voisinage ouvert~$I_t$ de~$t$ dans~$I$ et d'un graphe fini compact et connexe~$\Gamma_t$ de~$X$, que l'on peut toujours supposer appartenir à~$\got G$ quitte à l'agrandir, tel que~$I_t\subset \phi\inv(\Gamma_t)$ ; c'est ce que l'on souhaitait démontrer. 

\trois{rgamma} Fixons~$\Gamma\in \got G$. Soit~$y\in Y$. Si~$\Gamma'$ est un graphe appartenant à~$\got G$ et contenant~$\Gamma$ tel que~$y\in \phi\inv(\Gamma')$, alors~$\phi\inv(\Gamma)$ est un sous-graphe admissible de~$\phi\inv(\Gamma')$ ; l'image de~$y$ par la rétraction canonique de~$\phi\inv(\Gamma')$ vers~$\phi\inv(\Gamma)$ ne dépend pas du choix de~$\Gamma'$ (il suffit de vérifier qu'elle ne change pas lorsqu'on agrandit~$\Gamma'$, ce qui est immédiat) ; il est donc licite de la noter~$r_\Gamma(y)$. Notons que si~$y\in  \phi\inv(\Gamma)$ alors~$r_\Gamma(y)=y$. 

\medskip
Supposons que~$y\notin \phi\inv(\Gamma)$ et soit~$U$ sa composante connexe dans~$Y-\phi\inv(\Gamma)$. Montrons que~$\partial U=\{r_\Gamma(y)\}$. Par définition de~$r_\Gamma(y)$, il existe un intervalle fermé~$I$ tracé sur~$Y$ d'extrémités~$y$ et~$r_\Gamma(y)$, et tel que~$I\cap \phi\inv(\Gamma)=\{r_\Gamma(y)\}$. Cette dernière condition garantit que~$I\setminus\{r_\Gamma(y)\}$ est tracé sur~$U$, et donc que~$r_\Gamma(y)$ est adhérent à~$U$. 

\medskip
Supposons qu'il existe~$z\in \phi\inv(\Gamma)$ qui soit adhérent à~$U$ et soit distinct de~$r_\Gamma(y)$. Il existe un voisinage ouvert et connexe~$V$ de~$z$ qui ne contient pas~$r_\Gamma(y)$. Puisque~$z$ adhère à~$U$, son voisinage~$V$ rencontre~$U$, et~$V\cup U$ est par conséquent connexe. Il existe donc, en vertu du~\ref{conngraphe}, un graphe~$\Gamma'$ appartenant à~$\got G$, que l'on peut toujours supposer contenir~$\Gamma$, et tel que~$y$ et~$z$ soient situés sur la même composante connexe de~$\phi\inv(\Gamma')\cap (U\cup V)$. Ceci entraîne l'existence d'un intervalle tracé sur~$\phi\inv(\Gamma')$ qui relie~$y$ à~$z$ sans passer par~$r_\Gamma(y)$, ce qui contredit le fait que ce dernier est l'image de~$y$ par la rétraction canonique~$\phi\inv(\Gamma')\to \phi\inv(\Gamma)$ ; par conséquent,~$\partial U=\{r_\Gamma(y)\}$ ; remarquons que ceci entraîne que~$r_\Gamma(z)=r_\Gamma(y)$ pour tout~$z\in U$. 

\medskip
\trois{rgamcont} Nous nous proposons maintenant d'établir la continuité de l'application~$r_\Gamma$. Soit~$\Delta$ un ouvert de~$\Gamma$. Pour montrer que~$r_\Gamma\inv(\Delta)$ est ouvert, on peut supposer que~$\Delta$ est connexe et non vide. Soit~$\Theta$ le compact~$\Gamma-\Delta$ et soit~$U$ la composante connexe de~$Y-\Theta$ qui contient~$\Delta$ ; nous allons démontrer que~$r_\Gamma\inv(\Delta)$ est égal à~$U$. 

\medskip
Soit~$y\in r_\Gamma\inv(\Delta)$. Si~$y\in \phi\inv(\Gamma)$ alors~$y\in \Delta\subset U$ ; supposons que~$y\notin  \phi\inv(\Gamma)$ et soit~$V$ la composante connexe de~$y$ dans~$Y-\phi\inv(\Gamma)$. La composante connexe~$W$ de~$y$ dans~$Y-\Theta$ est fermée dans ce dernier, et contient~$V$. Comme~$\partial V$ est égal à~$\{r_\Gamma(y)\}$ et comme~$r_\Gamma(y)\notin \Theta$, on a~$r_\Gamma(y)\in W$. Par conséquent,~$W$ rencontre~$\Delta$ et coïncide de ce fait avec~$U$ ; ainsi,~$y\in U$. 

\medskip
Soit~$y\in U$. Si~$y\in \phi\inv(\Gamma)$ alors~$y\in U\cap  \phi\inv(\Gamma)=\Delta$, et l'on a donc~$r_\Gamma(y)=y\in \Delta$. Supposons que~$y\notin \phi\inv(\Gamma)$, et soit~$V$ sa composante connexe dans~$Y-\phi\inv(\Gamma)$. Le sous-ensemble~$V$ de~$U$ en est une partie ouverte, connexe et non vide ; comme elle diffère de~$U$ (elle ne rencontre pas~$\phi\inv(\Gamma)$), elle n'est pas fermée dans~$U$. L'unique point du bord de~$V$, qui n'est autre que~$r_\Gamma(y)$, appartient donc à~$U$ ; comme~$U\cap \phi\inv(\Gamma)=\Delta$, on a~$r_\Gamma(y)\in \Delta$ ; ainsi,~$y$ appartient à~$r_\Gamma\inv(\Delta)$, ce qu'il fallait démontrer. 

\trois{rgamarb} Soit~$\Delta$ un sous-arbre ouvert de~$\Gamma$, et soit~$U$ l'ouvert~$r_\Gamma\inv(\Delta)$.  Soient~$x$ et~$y$ deux points de~$U$ ; nous allons démontrer qu'il existe un et un seul fermé de~$U$ homéomorphe à un intervalle fermé d'extrémités~$x$ et~$y$. Tout tel fermé est, en vertu du~\ref{conngraphe} contenu dans~$\phi\inv(\Gamma')\cap U$ pour un certain graphe~$\Gamma'$ appartenant à~$\got G$, que l'on peut toujours supposer contenir~$\Gamma$. Il suffit donc de démontrer que pour tout~$\Gamma'$ appartenant à~$\got G$, contenant~$\Gamma, \phi(x)$ et~$\phi(y)$, il existe un et un seul fermé de~$\phi\inv(\Gamma')\cap U$  homéomorphe à un intervalle fermé d'extrémités~$x$ et~$y$. Mais si l'on se donne~$\Gamma'$ satisfaisant ces conditions, ~$\phi\inv(\Gamma')\cap U$  est l'image réciproque de~$\Delta$ par la rétraction canonique de~$\phi\inv(\Gamma')$ sur~$\phi\inv(\Gamma)$, et est donc un arbre, d'où l'assertion. 

\medskip
Montrons maintenant que~$\partial U$ est fini. Comme~$r_\Gamma$ est continue,~$r_\Gamma\inv(\overline \Delta)$ est un fermé de~$Y$, qui contient donc le bord de~$U$. Soit~$x$ un élément de~$r_\Gamma\inv(\overline \Delta)-U$. Supposons que~$x$ n'appartienne pas à~$\phi\inv(\Gamma)$ et soit~$V$ la composante connexe de~$x$ dans~$Y-\phi\inv(\Gamma)$. On a pour tout~$z\in V$ l'égalité~$r_\Gamma(z)=r_\Gamma(x)$ ; or~$r_\Gamma(x)\notin \Delta$ (puisque~$x\notin U$) ; par conséquent,~$V$ ne rencontre pas~$U$ et~$x$ n'appartient pas au bord de~$U$. Il s'ensuit que~$\partial U\subset \phi\inv(\Gamma)$. Compte-tenu du fait que~$\partial U\subset r_\Gamma\inv(\overline \Delta)$, il vient~$\partial U\subset \partial \Delta$, et~$\partial U$ est dès lors fini (l'inclusion réciproque étant triviale, on a en fait~$\partial U=\partial \Delta$). 

\trois{conclygraphe}{\em Conclusion :~$Y$ est un graphe.} L'espace~$Y$ est séparé et localement compact. Soit~$y\in Y$ et soit~$V$ un voisinage ouvert de~$y$ dans~$Y$ ; on va montrer que~$V$ contient un voisinage ouvert de~$y$ satisfaisant les conditions i) et ii) du~\ref{defarb}, ce qui prouvera que~$Y$ est un graphe. Quitte à restreindre~$V$, on peut le supposer relativement compact ; son bord est alors une partie compacte de~$Y$, qui ne rencontre de ce fait qu'un nombre fini de composantes connexes de~$Y\setminus\{y\}$. Choisissons~$\Gamma\in \got G$ tel que~$\phi\inv(\Gamma)$ contienne~$y$ et rencontre toutes les composantes connexe de~$Y\setminus\{y\}$ qui intersectent~$\partial V$, c'est-à-dire qui ne sont pas contenues dans~$V$. Soit~$\Delta$ un arbre ouvert de~$\phi\inv(\Gamma)$ contenant~$y$ et tel que~$\overline \Delta\subset V$. Le bord de~$V$ étant compact, il ne rencontre qu'un nombre fini de composantes connexes de~$Y-\overline \Delta$ ; par conséquent,~$V$ contient toutes les autres. En particulier,~$V$ contient presque toutes les composantes connexes de~$Y-\phi\inv(\Gamma)$ dont l'unique point du bord appartient à~$\Delta$ ; notons par ailleurs qu'une composante connexe de~$Y-\phi\inv(\Gamma)$ dont l'unique point du bord est~$y$ est une composante connexe de~$Y\setminus\{y\}$ qui ne rencontre pas~$\phi\inv(\Gamma)$, et est donc contenue dans~$V$. 

\medskip
Il découle de ce qui précède que l'on peut supposer, quitte à restreindre~$\Delta$, que~$V$ contient {\em toutes} les composantes connexes de~$Y-\phi\inv(\Gamma)$ dont l'unique point du bord appartient à~$\Delta$ ; comme~$\Delta\subset V$, on a~$r_\Gamma\inv(\Delta)\subset V$ ; et d'après~\ref{rgamarb} ci-dessus,~$r_\Gamma\inv(\Delta)$ satisfait les conditions i) et ii) du ~\ref{defarb}, ce qui achève la démonstration.~$\Box$

\deux{contrexboucles} {\em Remarque.} Dans l'énoncé du théorème~\ref{theofondquot} ci-dessus, l'hypothèse que~$\H^1(Y,\ZZ/\ell\ZZ)$ est fini pour au moins un entier~$\ell$ supérieur ou égal à~$2$ est indispensable. 

Pour le voir, désignons par~$Y$ le quotient de~$[0;1]\times\{0,1\}$ par la relation d'équivalence obtenue en identifiant~$(1/n,0)$ à~$(1/n,1)$ pour tout~$n>0$ ainsi que~$(0,0)$ à~$(0,1)$. On vérifie aisément que l'involution de ~$[0;1]\times\{0,1\}$ induite par la transposition de~$\{0,1\}$ passe au quotient et définit ainsi une action de~$\ZZ/2\ZZ$ sur~$Y$, qui est telle que~$Y/(\ZZ/2\ZZ)$ soit homéomorphe à~$[0;1]$. 

Et pourtant l'espace compact~$Y$ n'est pas un graphe, puisqu'il contient une infinité de boucles qui s'accumulent au voisinage de~$(0,0)=(0,1)$ ; mais la présence de cette infinité de boucles implique justement (par exemple {\em via} le~\ref{boucle}) que~$\H^1(Y,\ZZ/\ell \ZZ)$ est infini pour tout~$\ell\geq 2$. 

\chapter{Algèbre commutative}\markboth{Algèbre commutative}{Algèbre commutative}

\section{Morphismes étales et anneaux locaux
henséliens}

\subsection*{Schémas étales sur une base normale} 

\deux{remexttorseur} Soit~$\sch X$ un schéma intègre,
soit~$\sch Y$
 un~$\sch X$-schéma et soit~$\sch U$ un ouvert non
 vide de~$\sch X$. 
 
 \trois{locconnexeppf} Supposons $\sch Y$ plat et localement
 de présentation finie sur~$\sch X$. Il est alors localement connexe. 
 En effet, la propriété étant locale sur~$\sch Y$, on peut le supposer
 de présentation finie. La flèche~$\sch Y\to \sch X$ étant ouverte, tout point
 de~$\sch Y$ est adhérent à la fibre générique de~$\sch Y\to \sch X$. 
 Comme celle-ci a un nombre fini de composantes
 irréductibles, le schéma~$\sch Y$ est réunion
 {\em finie} de fermés irréductibles, d'où l'assertion. 
  
\trois{etalenormal} On suppose à partir 
de maintenant que~$\sch X$
est normal et que 
le morphisme~$\sch Y\to \sch X$ est étale. Le schéma~$\sch Y$ est alors
normal,
et ses composantes connexes
sont dès lors irréductibles. En conséquence, si~$\sch Z$ est une
composante connexe de~$\sch Y$, son ouvert~$\sch Z\times_{\sch X}\sch U$,
qui est non vide puisque~$\sch Y\to \sch X$ est ouvert,
est irréductible. 
Il en résulte que
l'application~$\pi_0(\sch Y\times_{\sch X}\sch U)\to \pi_0(\sch Y)$ 
est bijective. 

\trois{finietalenormal} On fait l'hypothèse supplémentaire
que le morphisme étale~$\sch Y\to \sch X$ est fini. L'ensemble 
des sections de~$\sch Y$
au-dessus de~$\sch X$ (resp. $\sch U$) est alors en bijection avec l'ensemble
des composantes connexes de~$\sch Y$ (resp.~$\sch Y\times_{\sch X}\sch U$)
dont le degré sur~$\sch X$ (resp.~$\sch U$) est égal à~$1$. ll s'ensuit,
compte-tenu de~\ref{etalenormal}, que la restriction
induit une bijection~$\sch Y(\sch X)\simeq \sch Y(\sch U)$. 

\medskip
\trois{finietalehombij} On suppose
toujours que~$\sch Y\to \sch X$ est fini étale ; pour
tout~$\sch X$-schéma~$\sch Z$, on note~$\sch Z_{\sch U}$ le
produit fibré~$\sch Z\times_{\sch X}\sch U$. 
Soit~$\sch Z$ un~$\sch X$-schéma étale. Il est normal,
et en 
appliquant~\ref{finietalenormal} à chacune des composantes
connexes de~$\sch Z$, on voit que
la flèche~$\sch Y(\sch Z)\to \sch Y (\sch Z_{\sch U})$ est bijective. 

\medskip
Autrement dit,
$$\mathsf{Hom}_{\sch X} (\sch Z,\sch Y)\to \mathsf{Hom}_{\sch U}(\sch Z_{\sch U}, \sch Y_{\sch U})$$ est bijective.

\trois{applitorseurs} Soit~$\sch G$ un schéma en groupes fini étale sur~$\sch X$. Supposons donnée
une action de~$\sch G_{\sch U}$ sur~$\sch Y_{\sch U}$. Il
résulte de~\ref{finietalehombij} 
qu'elle s'étend d'une unique manière en une action de~$\sch G$ sur~$\sch Y$. 

Supposons de plus que~$\sch Y_{\sch U}$ soit 
un~$\sch G_{\sch U}$-torseur. La flèche induite
$$\sch G_{\sch U}\times_{\sch U}\sch Y_{\sch U}\to \sch Y_{\sch U}\times_{\sch U}\sch Y_{\sch U}$$
est alors un isomorphisme. En utilisant à nouveau~\ref{finietalehombij},
on en déduit que~$\sch G\times_{\sch X}\sch Y\to \sch Y\times_{\sch X}\sch Y$
est un isomorphisme ; en conséquence,~$\sch Y$ est un~$\sch G$-torseur. 

\medskip
Donnons-nous maintenant un (autre)~$\sch G$-torseur~$\sch T$ sur~$\sch X$,
et supposons qu'il existe un isomorphisme
de~$\sch G_{\sch U}$-torseurs~$\sch T_{\sch U}\simeq \sch Y_{\sch U}$.
 En utilisant une fois encore~\ref{finietalehombij}, on voit 
que cet isomorphisme se prolonge (d'une unique manière) en un
isomorphisme de~$\sch G$-torseurs~$\sch T\simeq \sch Y$. 

\subsection*{Généralités sur les~$\mu_\ell$-torseurs}

\deux{introhilb90} Soit~$\ell$ un entier. Soit $(\mathsf C,\sch O)$ un site annelé, 
et soit~$\mu_\ell$
le faisceau des racines~$\ell$-ièmes de l'unité sur~$\mathsf C$. On suppose que
sur le site annelé~$\mathsf C$, 
l'endomorphisme~$z\mapsto z^\ell$
du faisceau~${\mathbb G}_m$ est surjectif ; 
le cas à avoir en tête est celui où~$\mathsf C$ est le (petit)
site étale d'un schéma ou d'un espace
analytique sur lequel~$\ell$ est inversible. On dispose alors sur le site~$\mathsf C$
d'une
suite exacte dite {\em de Kummer}
 $$\diagram 1\rto &\mu_\ell \rto &\gm \rrto^{z\mapsto z^\ell}&&\gm \rto&1\enddiagram,$$ 
 qui induit pour tout objet~$X$ de~$\mathsf C$
 une injection~$\sch O(X)\ti/(\sch O(X)\ti)^\ell
 \hookrightarrow \H^1(X,\mu_\ell)$. Si~$f\in \sch O(X)\ti$, et s'il n'y a pas d'ambiguïté
 sur~$X$ ni sur~$\ell$, on notera~$(f)$ l'image de~$f$ dans~$\H^1(X,\mu_\ell)$.
 La classe de cohomologie~$(f)$ correspond
 à une classe d'isomorphie de~$\mu_\ell$-torseurs sur~$X$,
 que l'on peut décrire explicitement : c'est la classe du~$\mu_\ell$-torseur~$U\mapsto
 \{z\in \sch O(U)\ti, z^\ell=f_{|U}\}$. 
 
 \trois{kummscheme} Soit~$\sch X$ un schéma sur lequel~$\ell$ est inversible et soit~$f$ une
 fonction inversible sur~$\sch X$. Tout faisceau étale localement constant sur~$\sch X \et$ est représentable, 
 et c'est en particulier le cas de tout~$\mu_\ell$-torseur. En conséquence, 
$\H^1(\sch X\et,\mu_\ell)$ classifie les~$\mu_\ell$-torseurs étales {\em schématiques}
 sur~$X$ à isomorphisme près ;  la classe d'isomorphie correspondant
 par ce biais à~$(f)$ est celle de~$\spec \sch O_{\sch X}[T]/(T^n-f)$. 
 
 \medskip
 Bien qu'il ne soit défini qu'à isomorphisme près,
 nous nous permettrons souvent de parler {\em du} torseur,
 et parfois simplement {\em du} revêtement étale (si l'on oublie
 l'action de~$\mu_\ell$) défini par~$(f)$ ; dans le cas d'une base affine, 
 nous évoquerons tout aussi bien l'algèbre finie étale définie par~$(f)$. 
 
  \trois{hilb90corps} Soit~$A$ un anneau local
  dans lequel~$\ell$ est inversible. Le
 théorème 90 de Hilbert assure que~$\H^1(A,\gm)=0$ (nous
 commettons l'abus usuel consistant à écrire~$\H^\bullet(A,.)$ au lieu
 de~$\H^\bullet((\spec A)\et,.)$). Il s'ensuit que
 la flèche~$f\mapsto (f)$ induit un 
 isomorphisme~$A\ti/(A\ti)^\ell\simeq \H^1(A,\mu_\ell)$. 
 
 \trois{muellegalzsurell} Soit~$K$ un corps dans lequel~$\ell$ est 
 inversible ; on suppose 
 que~$T^\ell-1$ est scindé dans~$K$.  
 
 \medskip
 Soit~$y\in K\ti$, soit~$d$ l'ordre
 de~$(y)$ et soit~$z$ une racine~$\ell$-ième de~$y$
 dans une clôture algébrique~$K^a$ de~$K$. Le 
 polynôme~$T^n-y$ s'écrit~$\prod\limits_{\omega \in \mu_\ell(K)} T-\omega z$.  
 Soit~$P$ un diviseur unitaire irréductible de~$T^n-\ell$ dans~$K[T]$, soit~$r$ 
 son degré, et soit~$s$ le PGCD de~$n$ et~$r$. Le polynôme~$P$ 
 s'écrit~$\prod \limits_{\omega \in \Omega} T-\omega z$ pour un certain
 sous-ensemble~$\Omega$ de~$\mu_\ell(K)$ de cardinal~$r$. Son terme
 constant est égal à~$z^r\prod \limits_{\omega\in \Omega} \omega$ ; comme il
 appartient à~$K$, on a~$z^r\in k\ti$ ; puisque~$z^\ell\in k\ti$, on a~$z^s\in k\ti$. Par
 conséquent, $y^s=(z^s)^\ell\in (k\ti)^\ell$, et~$d$ divise donc~$s$, et {\em a fortiori}~$r$. 
 
 Soit~$u\in k\ti$ tel que~$y^d=u^\ell$ ; il existe alors~$\omega \in \mu_\ell(K)$ tel
 que~$y=(\omega u)^{\ell/d}$ (en effet~$\mu_d(K)=(\mu_\ell(K))^{\ell/d}$ puisque~$T^\ell-1$
 est scindé dans~$K$). Posons~$ v=\omega u$. On a
 les égalités 
 $$T^\ell-y=T^\ell-v^{\ell/d}=(T^d)^{\ell/d}-v^{\ell/d}=\prod_{\zeta\in \mu_{\ell/d}(K)}T^d-\zeta v.$$ 
 Chacun des facteurs de la décomposition ci-dessus 
 appartient à~$K[T]$ et est de degré~$d$ ; par ce qui précède, ces facteurs
 sont nécessairement irréductibles. Il y en a exactement~$\ell/d$ ; en particulier, $T^\ell-y$
 est irréductible si et seulement si$~d=\ell$.

 \medskip
 En termes géométriques, le~$\mu_\ell$-torseur défini par~$(y)$ a~$\ell/d$ composantes connexes ;
 il est connexe 
 si et seulement si~$d=\ell$. 
 
 \subsection*{La notion de résidu} 
 
 \deux{cohcourbeproj} Soit~$F$ un corps, soit~$\sch C$
une~$F$-courbe projective, irréductible et normale, et soit~$\ell$ un entier inversible
dans~$F$. 

\trois{lieuetaletorseur} Soit~$h\in \H^1(F(\sch C),\mu_\ell)$ ; cette classe
définit une~$F(\sch C)$-algèbre étale~$A$. On note~$\sch D$ la~$F$-courbe
projective et normale, finie sur~$\sch C$, dont la fibre générique s'identifie
à~$\spec A$. Soit~$\sch C'$ l'ouvert dense de~$\sch C$ au-dessus
duquel~$\sch D$ est étale. En vertu de~\ref{applitorseurs}, 
la structure
de~$\mu_\ell$-torseur
sur~$\spec A\to \spec F(\sch C)$
définie
par la classe~$h$
se prolonge d'une unique
manière à~$\sch D\times_{\sch C}\sch C'\to \sch C'$. 

\medskip
Par ailleurs, soit~$\sch V$ un ouvert non vide de~$\sch C$ au-dessus
duquel le morphisme~$\spec A\to \spec F(\sch C)$ admet
un prolongement fini étale~$\sch W\to \sch V$. Comme~$\sch W$ est normal,
il s'identifie à~$\sch D\times_{\sch C}\sch V$ ; en conséquence,~$\sch V\subset \sch C'$. 

\medskip
Ainsi, un point~$P$
de~$\sch C$ appartient à 
$\sch C'$ si et seulement si 
le~$\mu_\ell$-torseur~$\spec A\to \spec F(\sch C)$
s'étend en un~$\mu_\ell$-torseur au-dessus de~$\spec \sch O_{\sch C,\sch P}$.

\trois{defresiduh} Fixons 
une fonction~$f\in F(\sch C)\ti$ telle
que~$h=(f)$. Soit~$P$ un point fermé de~$\sch C$. La
classe modulo~$\ell$ de la valuation~$P$-adique de~$f$ ne dépend 
que de~$h$ ; on la note~$\delta_P(h)$. On dit que~$\delta_P(h)$
est le {\em résidu de~$h$ en~$P$}. 

\medskip
Il résulte de~\ref{lieuetaletorseur} et~\ref{hilb90corps}
que le point~$P$ appartient à~$\sch C'$ si et seulement
si~$f$ appartient à~$\sch O_{\sch C,P}\ti$
modulo~$(F(\sch C)\ti)^\ell$, soit encore
si et seulement si~$\delta_P(h)=0$.

\subsection*{La suite exacte des résidus} 

{\em À partir de maintenant, on suppose~$F$ {\em algébriquement clos}}.

\deux{suiteexacteresidus} Soit~$\sch J$ 
la jacobienne de la courbe~$\sch C$,
et soit~$_\ell \sch J$ le sous-groupe
de~$\ell$-torsion de~$\sch J$. Nous allons démontrer qu'il existe un
morphisme
naturel~$\iota : \;\;_\ell\sch J(F) \to \H^1(F(\sch C),\mu_\ell)$ tel que
la suite

$$\diagram 0\rto & _\ell \sch J(F)\rrto^\iota &&\H^1(F(\sch C),\mu_\ell) \rto^{\prod \partial_P}
&\bigoplus\limits_{P\in \sch C(F)}\ZZ/\ell \ZZ\rto^{\;\;\;\;\;\;\sum}&\ZZ\ell \ZZ\rto&0\enddiagram$$
soit exacte. 

\medskip
\begin{itemize}

\item[$\bullet$] Soit~$x\in \;_\ell \sch J(F)$ et soit~$D$ un diviseur de degré~$0$ 
sur~$\sch C$ dont la classe dans le groupe de Picard est égale à~$x$. Comme
le point~$x$ est de~$\ell$-torsion, le diviseur~$\ell D$ est égale à~${\rm div}(f)$
pour une certaine fonction~$f\in F(\sch C)\ti$, bien déterminée à un élément
de~$F\ti$ près ; comme~$F$ est algébriquement clos, la classe~$(f)$
est bien définie. {\em Elle ne dépend pas du
choix de~$D$.} En effet, soit~$D'$ un (autre) représentant de~$x$ ;
choisissons~$g\in F(\sch C)\ti$ telle que~${\rm div}(g)=\ell D'$. On a
alors~${\rm div}(f/g)=\ell(D-D')$ ; puisque~$D$ et~$D'$ représentent~$x$, 
il existe~$h\in F(\sch C)\ti$ telle que~$D-D'={\rm div}(h)$. En conséquence,
${\rm div}(f/g)={\rm div}(h^\ell)$, ce qui signifie qu'il existe~$\lambda\in F\ti$
tel que~$f=\lambda h^\ell g$ ; comme~$F$ est algébriquement clos, 
il vient~$(f)=(g)$, comme annoncé.

On définit ainsi un morphisme~$\iota : \;_\ell\sch J(F)\to \H^1(F(\sch C),\mu_\ell)$. Montrons qu'il est
injectif. Soit~$x$ tel que~$\iota(x)=0$,
et soient~$D$ et~$f$ comme
ci-dessus. Puisque~$\iota(x)=0$, il existe~$g\in F(\sch C)\ti$ telle
que~$f=g^\ell$ ; comme le
diviseur de~$f$ est égal à~$\ell D$, on a~${\rm div}(g)=D$ et~$x=0$. 

\medskip
\item[$\bullet$] La surjectivité de~$\Sigma$ est évidente. 

\medskip
\item[$\bullet$]  La composée~$\Sigma\circ \prod \partial_P$ est nulle
parce que le diviseur d'une fonction rationnelle inversible
a un degré nul.

\medskip
Soit~$(P_i)$ une famille finie de points fermés
deux à deux distincts de~$\sch C$ et soit~$(e_i)$ une famille d'entiers
de somme nulle modulo~$\ell$. Soit~$D$
le diviseur~$\sum e_i P_i$. Choisissons un point
fermé~$Q$ de~$\sch C$, et soit~$N$ tel
que~${\rm deg}\;D=\sum e_i=N\ell$. Le diviseur~$D-N\ell Q$ est de degré~$0$, 
et définit donc un point~$x\in \sch J(F)$. 

Par divisibilité
du groupe~$\sch J(F)$, il existe~$y\in \sch J(F)$ tel que~$\ell y=x$. En d'autres termes, il existe un 
diviseur~$E$ de degré~$0$ sur~$\sch C$ et une fonction rationnelle
inversible~$f$ sur~$\sch C$
tels
que~$${\rm div}(f)=D-N\ell Q-\ell E=D-\ell(NQ+E).$$
Par construction, 
l'image de~$(f)$ dans~$\bigoplus\limits_{P\in \sch C(F)}\ZZ/\ell \ZZ$ est égale à~$\sum \overline{e_i} P_i$. Ainsi, 
la suite étudiée est exacte en~$\bigoplus\limits_{P\in \sch C(F)}\ZZ/\ell \ZZ$.

\item[$\bullet$] Il reste à montrer l'exactitude de la suite en~$\H^1(F(\sch C),\mu_\ell)$. 
Si~$x\in\; _\ell\sch J(F)$, son image~$h$ dans~$\H^1(F(\sch C),\mu_\ell)$ est par construction
la classe d'une fonction dont le diviseur est multiple de~$\ell$ ; il s'ensuit que~$\partial_P(h)=0$ pour
tout~$P\in \sch C(F)$. 

Réciproquement, soit~$f\in F(\sch C)\ti$ telle que~$(f)$ ait
un résidu nul en tout point
~$P$ de~$\sch C(F)$. 
Cette dernière condition signifie que~${\rm div}(f)=\ell D$
pour un certain diviseur~$D$, dont le degré
est nécessairement nul ; 
si l'on appelle~$x$ la classe de~$D$ dans~$\sch J(F)$, l'égalité~${\rm div}(f)=\ell D$ entraîne
que~$x\in \;_\ell \sch J(F)$ ; on
a alors~$\iota(x)=(f)$
par définition de~$\iota$,
ce qui achève la
démonstration.
\end{itemize} 

\deux{commenth1u} Soit~$\sch U$ un ouvert de~$\sch C$. On
peut reformuler les résultats de~\ref{cohcourbeproj} {\em et sq.} 
comme suit : la restriction au point générique induit un isomorphisme
entre~$\H^1(\sch U\et,\mu_\ell)$ et le sous-groupe
de~$\H^1(F(\sch C),\mu_\ell)$ formé des
classes de résidu nul en tout point fermé de~$\sch U$. 

\trois{comptorscourbes} {\em Remarque.} 
Soit~$h\in \H^1(\sch U\et,\mu_\ell)$, soit~$d$ son
ordre et soit~$\sch V\to \sch U$ le~$\mu_\ell$-torseur
qu'elle définit. Comme~$F$ est algébriquement clos,
le polynôme~$T^\ell-1$ est scindé dans~$F$, et {\em a fortiori}
dans~$F(\sch C)$. 

D'après~\ref{muellegalzsurell}, 
le schéma~$\sch V\times_{\sch C}\spec F(\sch C)$ a~$\ell/d$
composantes connexes. Comme~$\sch U$ est normal,~$\sch V$
est normal
et ses composantes connexes sont dès lors 
irréductibles. En 
conséquence,~$\sch V$ a lui aussi~$\ell/d$ composantes
connexes.

\trois{h1genreg} En vertu 
de~\ref{commenth1u},
le groupe~$\H^1(\sch C\et,\mu_\ell)$
s'identifie au
sous-groupe de~$\H^1(F(\sch C),\mu_\ell)$ formé des
classes de résidu nul en tout point fermé de~$\sch C$ ; 
on déduit 
alors de la suite exacte des résidus~(\ref{suiteexacteresidus})
que~$\H^1(\sch C\et,\mu_\ell)$ est canoniquement
isomorphe à~$_\ell\sch J(F)$. 

\medskip
Soit~$g$ le genre de~$\sch C$. La variété abélienne~$\sch J$
est de dimension~$g$, et~$_\ell\sch J(F)$ 
est donc isomorphe
(non canoniquement si~$g>0$ et~$\ell>1$) 
à~$(\ZZ/\ell \ZZ)^{2g}$. Ainsi~$\H^1(\sch C,\mu_\ell)\simeq (\ZZ/\ell\ZZ)^{2g}$.

\deux{remmuellezsurell} {\em Remarque.} Comme~$F$ est 
algébriquement clos, les faisceaux étales $\mu_\ell$ 
et~$\ZZ/\ell\ZZ$ sont isomorphes
(non canoniquement si~$\ell>2$)
sur~$F$. 

Une fois fixé un tel isomorphisme, on
peut récrire~\ref{suiteexacteresidus}, \ref{commenth1u}, \ref{comptorscourbes}
et~\ref{h1genreg} 
en remplaçant partout~$\H^1(.,\mu_\ell)$ par~$\H^1(.,\ZZ\ell\ZZ)$. 

\subsection*{Anneaux locaux henséliens} 

\deux{prehensel} Soit~$A$ un anneau local d'idéal maximal~$\got m$,
et soit~$k$ son corps résiduel.  

\trois{etesssurj} Soit~$B_0$ une~$k$-algèbre finie étale. 
Il existe alors une~$A$-algèbre finie étale~$B$ telle que~$B_0\simeq B/\got m$ ; 
en effet, pour le voir, on se ramène immédiatement au cas où~$B_0$ est une extension
finie séparable de~$k$. Le théorème de l'élément primitif assure l'existence
d'un polynôme unitaire, irréductible et séparable~$P_0\in k[T]$ tel
que~$B\simeq k[T]/(P_0)$. Si~$P$ désigne un relevé unitaire
de~$P_0$ dans~$A[T]$ on peut alors prendre~$B=A[T]/P$. 

\trois{relefactun} Soient~$Q$ un polynôme
unitaire de~$k[T]$ et soit~$P$
un polynôme unitaire de~$A[T]$ relevant~$Q$. 

\medskip
Supposons donnés
deux polynômes unitaires~$Q_1$ et~$Q_2$
de~$k[T]$ tels que~$Q_1$
et~$Q_2$ soient premiers entre et tels que~$Q=Q_1Q_2$. 
Il existe alors {\em au plus}
un couple~$(P_1,P_2)$ de polynômes unitaires de~$A[T]$ relevant
respectivement~$Q_1$ et~$Q_2$ et tels que~$P=P_1P_2$. 

En effet, supposons qu'un tel couple~$(P_1,P_2)$
existe. L'anneau
quotient~$A[T]/(P_1, P_2)$ est un~$A$-module fini, 
dont la
réduction modulo~$\got m$ s'identifie 
à~$k[T]/(Q_1,Q_2)$ qui est nul par hypothèse ; 
par le lemme
de Nakayama, il vient~$A[T]/(P_1, P_2)=0$. Ainsi,~$(P_1)+(P_2)=A[T]$, et
le lemme chinois fournit alors un isomorphisme~$i: A[T]/P\simeq A[T]/P_1\times A[T]/P_2$,
qui relève l'isomorphisme~$j:k[T]/Q\simeq k[T]/Q_1\times k[T]/Q_2$. 
Soit~$e$ l'idempotent de~$A[T]/P$ correspondant à
~$(1,0)$ {\em via}~$i$. 
Il relève l'idempotent~$f$ de~$k[T]/Q$ qui correspond
à~$(1,0)$ {\em via}~$j$. 
Or~$f$ admet {\em au plus} un relevé idempotent
dans~$A[T]/P$. Pour le voir, considérons un élément~$x$
de~$\got mA[T]/P$ tel que~$e':=e+x$ soit idempotent. On 
a alors~$$x^3=(e-e')^3 =e-3ee'+3ee'-e'=e-e'=x,$$ et 
donc~$x(1-x^2)=0$. Mais~$x$ appartient à~$\got  mA[T]/P$,
qui est le radical de Jacobson de~$A[T]/P$ ; en 
conséquence,~$1-x^2$ est inversible et~$x=0$, 
d'où l'assertion requise. Il s'ensuit que~$e$ ne dépend pas
du couple~$(P_1,P_2)$. 

Or~$P_1$ peut se caractériser comme l'annulateur de~$e$, et~$P_2$ comme
celui de~$1-e$. En conséquence, le couple~$(P_1,P_2)$ est
uniquement déterminé, comme annoncé. 

\medskip
Si~$\alpha\in k$ alors~$\alpha$ est racine simple
de~$Q$ si et seulement si~$Q$ s'écrit~$(X-\alpha)R$
avec~$R$ premier à~$X-\alpha$ ; et si~$a\in A$, on
a~$P(a)=0$ si et seulement si~$X-a$ divise~$P$. 

On déduit alors de ce qui précède qu'une racine simple
de~$Q$ admet~{\em au plus} un relèvement en une racine
de~$P$. 

\deux{rapphensel} Soit~$A$ un anneau local
 d'idéal maximal~$\got m$
 et de corps résiduel~$k$.  On note~$a\mapsto \overline a$
 la réduction modulo~$\got m$. 
 Les propriétés suivantes sont équivalentes. 
 
 \medskip
 i) Toute~$B$-algèbre finie est un produit fini d'anneaux locaux. 
 
 ii) Pour tout polynôme unitaire~$P\in A[T]$ et toute 
 factorisation~$\overline P=Q_1Q_2$ où~$Q_1$
 et~$Q_2$ sont deux polynômes
 unitaires et
 premiers entre eux de~$k[T]$, il existe un (unique)
 couple~$(P_1,P_2)$ de polynômes
 unitaires de~$A[T]$ tels que~$P=P_1P_2$ et~$\overline {P_i}=Q_i$ 
 pour tout~$i$. 
 
 iii) Pour tout polynôme unitaire~$P\in A[T]$ et toute 
racine simple~$\alpha$
de~$\overline P$ dans~$k$, il existe
un (unique) élément~$a\in A$ tel que~$\overline a=\alpha$
et~$P(a)=0$. 

\medskip
Lorsqu'elles sont satisfaites, on dit que~$A$ est~{\em hensélien.} 

\deux{complimplhens} Soit~$A$
un anneau local d'idéal maximal~$\got m$. 
Si~$A$ est complet pour la topologie~$\got m$-adique, il est hensélien
(on prouve~iii) par la méthode de Newton, c'est le «lemme de Hensel»
originel).

\deux{proplochens} Soit~$A$ un anneau local hensélien. 

\trois{stablefinihens} Toute~$A$-algèbre finie locale est hensélienne, comme on 
le voit à l'aide de la propriété~i). 

\trois{eqetalehens1} Soit~$B$ une~$A$-algèbre finie étale. C'est un produit
de~$A$-algèbres finies étales locales ; il résulte aisément de la forme
explicite de ces dernières 
et de~iii) que~${\rm Hom}_A(B,A)\to {\rm Hom}_k(B\otimes_Ak,k)$ est bijective.

\trois{eqetalehens2}
Soit~$C$ une~$A$-algèbre finie  étale. Elle s'écrit comme un produit~$\prod C_i$
de~$A$-algèbres finies, étales, locales et henséliennes ; pour tout~$i$, 
l'idéal~$\got mC_i$ de~$C_i$ est son idéal 
maximal, puisque~$C_i/\got mC_i$ est une~$k$-algèbre
locale finie étale, et partant un corps.

En appliquant~\ref{eqetalehens1} à la~$C_i$-algèbre
finie étale~$B\otimes_AC_i$ pour tout~$i$, 
on voit 
que~${\rm Hom}_A(B,C)\to {\rm Hom}_k(B\otimes_Ak,C\otimes_Ak)$ est bijective.

\medskip
Le foncteur~$B\mapsto B\otimes_Ak$, de la catégorie des~$A$-algèbres finies
étales vers celle des~$k$-algèbres finies étales, est ainsi pleinement fidèle. 
En vertu de~\ref{etesssurj}, c'est une équivalence de catégories. 

\subsection*{Hensélisé d'un anneau local} 

\deux{henselise} Soit~$A$ un anneau local, soit~$\got m$
son idéal maximal et soit~$k$ son corps résiduel. On note~$\sch X$ et~$\bf x$ 
les spectres respectifs de~$A$ et~$k$, et~$j$ l'immersion
fermée~${\bf x}\hookrightarrow \sch X$. On considère~$\sch O_{\sch X}$ comme
un faisceau sur le {\em petit} site étale~$\sch X\et$ ; on
pose~$A^h=\H^0({\bf x}\et, j\inv \sch O_{\sch X})$. 

\trois{prophenlocet} Soit~$\mathsf F$ l'ensemble des (classes
d'isomorphie de) couples~$(\sch Y,\pi)$ où~$\sch Y$ est
un~$\sch X$-schéma étale et de présentation finie,
et où~$\pi : {\bf x}\to \sch Y$ est un morphisme
tel que~$j$ soit la flèche
composée~$$\diagram {\bf x}\rto^\pi&\sch Y\rto &\sch X\enddiagram.$$

L'anneau~$A^h$ est alors la limite
inductive des anneaux locaux~$\sch O_{\sch Y, \pi({\bf x})}$,
pour~$(\sch Y,\pi)$ parcourant~$\mathsf F$; c'est donc une~$A$-algèbre
locale. 

Remarquons que les morphismes de transition du
système inductif considéré sont des morphismes locaux
et plats entre anneaux locaux ; ils sont dès lors
fidèlement plats, et en particulier injectifs. 

Soit~$(\sch Y,\pi)\in \mathsf F$. Le schéma~$\sch Y\times_{\sch X}{\bf x}$
est un~$k$-schéma étale de type fini, c'est donc le spectre
d'un produit fini d'extension finies
séparables de~$k$. En particulier,~$\sch O_{\sch Y,\pi({\bf x})}/\got m \sch O_{\sch Y,\pi({\bf x})}$
est une extension finie séparable de~$k$. Cette extension 
se plonge par construction dans~$\kappa({\bf x})=k$ ; 
il vient~$\sch O_{\sch Y,\pi({\bf x})}/\got m 
\sch O_{\sch Y,\pi({\bf x})}=k$. 
Il s'ensuit que l'idéal maximal de~$A^h$ est~$\got mA^h$
et que son corps résiduel est~$k$.

\medskip
{\em L'anneau~$A^h$ est hensélien.} En effet, soit~$P$ un polynôme
unitaire de~$A^h[T]$, et soit~$\alpha$ une racine simple
de l'image~$\overline P$ de~$P$ dans~$k[T]$. Le polynôme~$P$ 
provient d'un polynôme (que nous notons encore~$P$) à coefficient
dans~$\sch O_{\sch Y,\pi({\bf x})}$ pour un certain~$(\sch Y,\sch \pi)\in \mathsf F$ ; 
quitte à restreindre~$\sch Y$, on peut supposer qu'il est affine
et que~$P\in \sch O_{\sch Y}(\sch Y)[T]$. Le schéma~$\sch Z:=\spec \sch O_{\sch Y}(\sch Y)[T]$ 
est fini et plat sur~$\sch Y$ ; la racine simple~$\alpha$ de~$\overline P$
définit un point~${\bf z}$ de~$\sch Z$ situé au-dessus de~$\pi({\bf x})$, et de corps
résiduel~$k$. Le schéma~$\sch Z$ est étale en~$\bf z$
au-dessus de~$\sch Y$, et {\em a fortiori}
au-dessus de~$\sch X$ ; comme~$\bf z$ a pour corps résiduel~$k$,
l'immersion~$\pi$ se factorise par une flèche~${\bf x}\to \sch Z$ d'image~${\bf z}$.

Il existe donc un morphisme naturel de~$\sch O_{\sch Y,{\bf y}}$-algèbres 
de~$\sch O_{\sch Z,{\bf z}}$ vers~$A^h$. Par construction,
ce morphisme envoie~$\overline T$ sur une racine de~$P$ 
qui relève~$\alpha$, ce qui permet de conclure. 

\medskip
On dit que~$A^h$ est le~{\em hensélisé} de~$A$.

\trois{idealnalpha} Nous allons maintenant donner
une description un peu plus explicite de~$A^h$. Soit~$\mathsf E$ l'ensemble des
couples~$(P,\alpha)$ où~$P$ est un polynôme unitaire
à coefficients dans~$A[T]$ et~$\alpha$ une racine simple
de~$\overline P$ dans~$k$.

\medskip
Soit~$(P,\alpha)\in \mathsf E$. Soit~$\sch Y$ 
l'ouvert de~$\spec A[T]/P$ formé des points en lesquels
celui-ci est étale sur~$\sch X$. L'évaluation en~$\alpha$ définit
un morphisme
surjectif~$A[T]/P\to k$, dont
on notera~$\got n(\alpha)$ le noyau ; c'est un idéal maximal
de~$A[T]/P$. Le point correspondant
de~$\spec A[T]/P$ est situé sur~$\sch Y$, et son corps résiduel est~$k$ ; 
c'est donc l'image d'un~$\sch X$-morphisme~$\pi : {\bf x}\to \sch Y$. 
Le couple~$(\sch Y,\pi)$ appartient à~$\mathsf F$, et la~$A$-algèbre locale~$(A[T]/P)_{\got  n(\alpha)}$
est égale à~$\sch O_{\sch Y,\pi({\bf x}})$. 

Réciproquement, il résulte de la description locale des morphismes étales que
pour tout~$(\sch Y,\pi)\in \mathsf F$, il existe~$(P,\alpha)\in \mathsf E$
tel que~$\sch O_{\sch Y,\pi({\bf x})}\simeq (A[T]/P)_{\got n(\alpha)}$.
On en déduit aussitôt les assertions suivantes. 

\medskip
$\bullet$ La famille des anneaux~$(A[T]/P)_{\got n(\alpha)}$
constitue, lorsque~$(P,\alpha)$ parcourt~$\mathsf E$, un système inductif
filtrant de~$A$-algèbres locales, de corps résiduel~$k$ et d'idéal maximal
engendré par~$\got m$ ; les morphismes de transition 
de ce système sont fidèlement 
plats, et en particuliers injectifs. 

$\bullet$ L'anneau~$A^h$ est la limite inductive du système précédent. 

\trois{propunivhensa} Soit~$B$ un anneau local hensélien
et soit~$A\to B$ un morphisme
local. Soit~$(P,\alpha)\in \mathsf E$ ; comme~$B$ est hensélien,
le polynôme~$P$ possède dans~$B$ une unique racine 
relevant~$\alpha$. On en déduit immédiatement, à l'aide
de la description de~$A^h$ donnée au~\ref{idealnalpha} ci-dessus,
qu'il existe un unique morphisme 
local de~$A$-algèbres de~$A^h$ dans~$B$. 

\deux{henselnormal} Supposons que~$A$ est normal, et 
soit~$(P,\alpha)\in \mathsf E$. En tant que localisé
d'une~$A$-algèbre étale, l'anneau local
~$(A[T]/P)_{\got n(\alpha)}$
est normal, et en particulier intègre. 
La fibre générique de~$\spec (A[T]/P)_{\got n(\alpha)}\to \spec A$
est intègre et de dimension nulle ; elle est donc réduite au point
générique de~$\spec (A[T]/P)_{\got n(\alpha)}$. 

Soit~$B$ une~$A$-algèbre locale hensélienne et
intègre telle
que~$A\to B$ soit injective. Le noyau de~$(A[T]/P)_{\got n(\alpha)}\to B$
est alors un idéal premier de~~$(A[T]/P)_{\got n(\alpha)}$ situé
au-dessus de l'idéal~$(0)$ de~$A$. Par ce qui précède, ce noyau est nécessairement
trivial ; autrement
dit, la flèche canonique~$(A[T]/P)_{\got n(\alpha)}\to B$
est injective. 

\medskip
Ceci valant pour
tout~$(P,\alpha)\in \mathsf E$, l'anneau~$A^h$ 
est normal, et s'injecte dans~$B$.

\section{Algèbre commutative graduée}

Dans son article publié dans Israel Journal of Maths., Temkin 
a développé une théorie de la {\em réduction des germes 
d'espaces analytiques}, 
extrêmement
efficace 
pour l'étude locale des espaces de Berkovich, et dont nous 
ferons un usage abondant dans la suite. Elle repose 
sur le formalisme de l'{\em algèbre commutative graduée}, 
et sur la {\em réduction graduée} des algèbres normées. 

\medskip
Par ailleurs, l'algèbre commutative graduée 
et la réduction graduée se sont révélées 
particulièrement adaptées à l'étude 
des extensions de corps valués. D'une part, 
elles permettent souvent une présentation unifiée 
de raisonnements qui jusqu'alors nécessitaient 
de considérer séparément le corps résiduel et le groupe des valeurs ; 
d'autre part, elles fournissent une description particulièrement 
agréable des extensions modérément ramifiées, lesquelles 
jouent un rôle absolument crucial dans notre mise au jour 
de la structure locale des courbes analytiques 
 : en quelque sorte, ces extensions sont à la réduction graduée 
 ce que les extensions non ramifiées 
 sont à la réduction classique.
 
 \medskip
En conséquence, nous allons, dans cette section, présenter 
succinctement le formalisme
de l'algèbre commutative graduée, et les avatars 
gradués d'un certain nombre de théorèmes 
d'algèbre commutative usuelle, 
la plupart du temps sans démonstrations -- celles-ci étant en général 
des retranscriptions~{\em mutatis mutandis} des preuves classiques. 
Puis nous referons en détail, en utilisant 
systématiquement la réduction graduée, toute la 
théorie des extensions algébriques de corps valués : 
groupes de décomposition, d'inertie et de 
ramification, hensélisation, extensions non ramifiées et modérément ramifiées. 

\medskip
Mentionnons
que nous avons choisi 
une approche un peu inhabituelle de l'algèbre
graduée, qui
diffère notamment de celle de Temkin. Celui-ci
considère un anneau gradué comme un
anneau~$A$ muni d'une décomposition en somme directe
~$A=\bigoplus A_r$, où~$r$ parcourt un groupe abélien. 
Mais comme en pratique il n'arrive {\em jamais},
pour le type de problème que nous aurons à considérer ici,
que
l'on ait à additionner deux éléments (non nuls)
appartenant
à des sommandes distincts, il nous a semblé plus simple
de remplacer les sommes directes par des {\em réunions disjointes} :
cela ne change sur le fond
strictement rien
aux énoncés ni à leurs preuves, et permet
d'éviter 
les répétitions fastidieuses de l'adjectif «homogène», indispensables
lorsqu'on s'en tient à la définition traditionnelle. 

De ce fait, les objets qui nous tiendront lieu
d'anneaux gradués ne sont pas des anneaux 
(l'addition n'y est pas partout définie). Nous avons donc
modifié la terminologie.

\subsection*{Généralités}

On fixe un groupe abélien divisible~$D$, noté multiplicativement
(son élément neutre est en particulier noté~$1$) ; en pratique,~$D$ sera souvent égal à~$\RR\ti_+$.  

\deux{rappgrad} Un~{\em $D$-annéloïde} (sous-entendu, commutatif et unitaire)
est un ensemble~$A$ muni d'une décomposition~$A=\coprod \limits_{r\in D}A_r$,
que l'on appelle {\em graduation}, et des
données supplémentaires suivantes : 

\medskip
$\bullet$ pour tout~$r\in D$, une loi interne~$+$ sur~$A_r$ faisant de ce
dernier un groupe abélien, d'élément neutre noté~$0_r$ ;  

$\bullet$ une loi
interne commutative et associative~$\times$
sur~$A$
induisant pour tout~$(r,s)\in D^2$
une application bilinéaire~$A_r\times A_s\to A_{rs}$ ; 

$\bullet$ un élément~$1\in A_1$ tel que~$1\cdot a=a$ pour tout~$a\in A$.

\deux{diversdgrad} Soit~$A$ un~$D$-annéloïde. 

\trois{deganngrad} Si~$a\in A$, on 
appelle {\em degré}
de~$a$
l'unique élément~$r$ de~$D$ tel
que~$a$ appartienne à~$A_r$.

\trois{convzerooide} {\em Une convention.} Il arrivera
fréquemment, lorsque~ la valeur de~$r$ est
définie sans ambiguïté par le contexte, que
l'on utilise le symbole~$0$ au lieu de~$0_r$. Par exemple, 
on se permettra le plus souvent d'écrire~$a=0$
au lieu de~$a=0_{\deg a}$. 

\medskip
On  notera~$A^{\neq 0}$ l'ensemble des éléments non nuls
de~$A$, c'est-à-dire la réunion disjointe des~$A_r\setminus\{0_r\}$. 
Le sous-ensemble~$\deg (A^{\neq 0})$ de~$D$ est non vide si et seulement
si~$A$ est non nul, c'est-à-dire si et seulement si~$A\neq \coprod \{0\}_r$ ; cela
revient à demander que~$1\neq 0$. 

\trois{def-ati-anneloide}
Un élément~$a$ de~$A$ est dit {\em inversible}
s'il possède un inverse pour la multiplication ; cet inverse
est alors nécessairement unique, et sera en général noté~$a^{-1}$. On désignera par~$A\ti $ l'ensemble des éléments
inversibles de~$A$ ; il est stable sous
la loi~$\times$ qui en fait un groupe abélien. 

\trois{moranngrad} Si~$B$ est un~$D$-annéloïde, 
un morphisme de~$D$-annéloïdes 
de~$A$ vers~$B$
est une application~$f$ de~$A$ vers~$B$ 
qui préserve le degré, qui est telle que l'application induite~$A_r\to B_r$
soit un morphisme de groupes pour tout~$r$, qui
commute au produit et qui envoie~$1$
sur~$1$. On note~$D$-$\mathsf{Annel}$ la catégorie
des~$D$-annéloïdes. 

\trois{changegroupe} {\em Changement de groupes.} 
Soit~$D'$ un groupe divisible contenant~$D$. 
Si~$B$ est un~$D$-annéloïde, on note~$B\odot _DD'$ 
le~$D'$-annéloïde~$B\coprod (\coprod\limits_{r\in D'\setminus D} \{0_r\})$ ; 
si~$C$ est un~$D'$-annéloïde, on note~$C_D$
 le~$D$-annéloïde~$\coprod\limits_{r\in D}C_r$.
 Le 
 couple~$$(B\mapsto B\odot_D D', C\mapsto C_D)$$ 
 est un couple de foncteurs adjoints. Il induit une
 équivalence entre $D$-$\mathsf{Annel}$
 et la sous-catégorie pleine de
 ~$D'$-$\mathsf{Annel}$ constituée des objets~$C$ tels
 que~$\deg(C^{\neq 0})$
soit contenu dans~$D$.

\trois{commentextgroup} En pratique, les propriétés que nous allons considérer dans ce texte
seront invariantes par cette équivalence de catégories ; il pourra donc arriver
que l'on identifie subrepticement un~$D'$-annéloïde dont toutes 
les composantes de degré appartenant à~$D'\setminus D$ sont nulles
au~$D$-annéloïde correspondant.

%
 
 \trois{liensalgclass} Si~$A$ est un~$D$-annéloïde, 
 le sommande~$A_1$ est un anneau ordinaire. On définit par ce biais
 une équivalence entre la catégorie des anneaux
 et celle des~$\{1\}$-annéloïdes. La théorie des annéloïdes
 contient donc l'algèbre commutative ordinaire
 comme cas particulier. 
 
\medskip
{\em À partir de maintenant, nous dirons «annéloïde» au lieu de «~$D$-annéloïde». }

\deux{idgrad} Soit~$A$ un annéloïde. 
Une~{\em $A$-algèbre} est
un annéloïde~$B$
muni d'un morphisme~$A\to B$. Un {\em idéal} de~$A$ est une partie~$I$ de~$A$
stable par multiplication
par tous les éléments de~$A$, et telle que~$I\cap A_r$ soit pour tout~$r$ un sous-groupe de~$A_r$. 

On définit de façon évidente le quotient de~$A$ par un idéal.  

\deux{locoide} Soit~$A$ un annéloïde. Si~$S$ est une partie 
multiplicative de~$A$ ({\em i.e.}
$S$ contient~$1$ et est stable par multiplication),
la catégorie des~$A$-algèbres
dans lesquelles les éléments de~$S$ sont inversibles admet
un objet initial, noté~$S\inv A$. Pour tout~$r\in D$, le
sommande~$(S\inv A)_r$ est constitué de fractions~$a/s$
avec~$a\in A, s\in S$, et~$\deg a/\deg s=r$. Si~$a/s$ et~$b/t$ sont deux telles
fractions alors~$(a/s)=(b/t)$ si et seulement si il existe~$\sigma \in S$
vérifiant l'égalité~$\sigma(at-bs)=0$.
Si~$S=\{f^n\}_n$
 pour un certain~$f\in A$
on écrira~$A_f$ au lieu de~$S^{-1}A$.

\deux{noyimoid} Soit~$f: A\to B$ un morphisme d'annéloïdes. L'image
de~$f$ est un sous-annéloïde de~$B$, en un sens évident. Son
{\em noyau}~$I$ est l'ensemble
des éléments~$a$ de~$A$ tels que~$f(a)=0$. C'est
un idéal de~$A$, et~$f$ induit
un isomorphisme
entre~$A/I$ et~$f(A)$.  

\deux{annintgrad} Un annéloïde~$A$ est dit {\em intègre} s'il est 
non nul et si~$ab=0$ si et seulement
si~$a=0$ ou~$b=0$. Si~$A$ est un
annéloïde intègre,~$\deg(A^{\neq 0})$ est un sous-monoïde unitaire
de~$D$.

\deux{corspgrad} Un {\em corpoïde} est un annéloïde non nul dans lequel tout élément 
non nul est inversible. 
Si~$K$ est un corpoïde,~$\deg(K\ti)=\deg(K^{\neq 0})$ est un
sous-groupe de~$D$. 

Un corpoïde est intègre ; si~$A$ est un annéloïde
intègre et si~$S$ désigne l'ensemble
des éléments non nuls de~$A$, 
le localisé~$K:=S\inv A$ est un corpoïde, 
appelé {\em corpoïde des fractions} de~$A$, 
dont les éléments sont de la forme~$a/b$, où~$a$ et~$b$ sont deux éléments 
de~$A$, et où~$b$ est non nul. 
Tout plongement de~$A$ dans un corpoïde
se factorise d'une unique manière par~$K$. 

\medskip
Tout morphisme d'un corpoïde dans un annéloïde non nul
est injectif. 

\deux{idmaxoide} Soit~$A$ un annéloïde. Un idéal~$I$ de~$A$ est 
dit {\em premier} (resp. {\em maximal}) si~$A/I$ est intègre (resp. est 
un corpoïde). Si~$I$ est un idéal strict de~$A$, il est contenu dans un
idéal maximal ; en particulier, tout annéloïde non nul
possède un idéal maximal (appliquer ce qui précède à
l'idéal nul~$\coprod \{0_r\}$). 

\deux{anneloideloc} Un annéloïde est dit {\em local}
s'il possède un et un seul
idéal maximal. Soit~$A$ un annéloïde et soit~$\got p$
un idéal premier de~$A$ ; posons~$S=A\setminus \got p$. 
Le localisé~$S\inv A$ est un annéloïde local, qui sera le plus
souvent noté~$A_{\got p}$. 

\deux{spectroide}
{\bf Spectre d'un annéloïde}.
Soit~$A$ un annéloïde. On note~$\spec A$ l'ensemble
des idéaux premiers de~$A$, muni de la topologie dite de Zariski, dont les fermés
sont exactement les~$V(I)=:\{\got p,\got p\subset I\}$ où~$I$ est un idéal de~$A$. On vérifie aussitôt 
qu'une base d'ouverts
de~$\spec A$ est constituée des parties de la forme~$D(f):=\{\got p, f\notin \got p\}$ où~$f\in A$. 
On a~$\spec A=\varnothing$ si et seulement si~$A=\{0\}$. 

\trois{notation-fx-spectroide}
Comme en théorie de schémas classiques, on préfère penser
à~$\spec A$ comme à un espace topologique dont chaque point
{\em correspond}
à un idéal premier de~$A$. Soit~$x\in \spec A$ et soit~$\got p$ l'idéal premier
correspondant. On note~$\kappa(x)$ le corpoïde
${\rm Frac}\;A_{\got p}$, qui est appelé le {\em corpoïde résiduel}
de~$x$ ; et~$f\mapsto f(x)$ le morphisme
naturel de~$A$ dans~$\kappa(x)$. 

On peut dès lors écrire 
$V(I)=\{x\in \spec A, f(x)=0\;\forall\;x\in I\}$ pour tout idéal~$I$ de~$A$,
et~$D(f)=\{x\in \spec A, f(x)\neq 0\}$ pour toute~$f\in A$. 

\trois{fonct-specoide}
Si~$B$ est un annéloïde et si~$\phi \colon A\to B$ est un morphisme, la 
formule~$\got q\mapsto \phi^{-1}(\got q)$ définit une application continue de~$\spec B$ vers~$\spec A$ ; si~$y\in \spec B$
et si~$x$ désigne son image sur~$\spec A$, on dispose d'un plongement naturel~$\kappa(x)\hookrightarrow \kappa(y)$.

\medskip
Pour tout idéal~$I$ de~$A$ la flèche naturelle~$\spec A/I\to \spec A$ identifie~$\spec A/I$ au fermé~$V(I)$ de~$\spec A$,
en préservant les corpoïdes résiduels : 
pour tout~$f\in A$ la flèche naturelle~$\spec A_f\to \spec A$ identifie~$\spec A_f$ à l'ouvert~$D(f)$ de~$\spec A$, 
en préservant les corpoïdes résiduels.

\trois{cond-nullite-specoide}
Soit~$f\in A$. 

\medskip
Dire que~$f(x)\neq 0$ pour tout~$x\in A$ signifie que~$f$ n'appartient à aucun
idéal premier de~$A$ ; cela revient à demander qu'elle n'appartienne à aucun idéal strict de~$A$, 
c'est-à-dire encore que~$(f)=A$. Autrement dit, $D(f)=\spec A$ si et seulement si~$f\in A\ti$. 

\medskip
Dire que~$f(x)=0$ pour tout~$x\in \spec A$, c'est-à-dire que~$f$ appartient à tous les idéaux
premiers de~$A$, revient à demander que~$\spec A_f\simeq D(f)=\varnothing$ et donc que~$A_f=\{0\}$ ; mais la
description explicite de~$A_f$ au moyen de fractions assure que c'est le cas si et seulement si~$f$ est nilpotent. 

\trois{spec-annel-0}
Si~$A$ est intègre l'idéal nul de~$A$ est premier, et le point
correspondant de~$\spec A$ est dense ; en particulier, $\spec A$ est un espace topologique
irréductible.

\deux{polgrad} Soit~$A$ un annéloïde
et soit~${\bf r}=(r_i)$ une famille
d'éléments de~$D$. 
Posons~${\bf T}=(T_i)$. 
On note~$A[{\bf r}\inv{\bf T}]$ l'annéloïde
défini comme suit : 
pour tout~$s$ appartenant à~$D$, le sommande~$A[{\bf r}\inv{\bf T}]_s$ 
est constitué des sommes finies (formelles)~$\sum a_I{\bf T}^I$ 
avec~$a_I\in A_{s{\bf r}^{-I}}$ pour tout~$I$ ; les
opérations sont définies de façon évidentes. 
Lorsqu'on voudra évoquer le degré usuel d'un élément de~$A[{\bf r}\inv{\bf T}]$, 
on parlera de son degré {\em monomial}. 

\trois{univ-algpol-oide}
Soit~$(b_i)$ une famille d'éléments 
d'une~$A$-algèbre~$B$, chaque~$b_i$ étant de degré~$r_i$. Il existe un 
unique morphisme de~$A$-algèbres
de~$A[{\bf r}\inv {\bf T}]$ vers~$B$ qui
envoie~$T_i$ sur~$b_i$ pour tout~$i$. Ce morphisme
 est noté~$P\mapsto P(b_i)_i$, et est évidemment
donné
par la formule que cette notation suggère. Son image 
est la plus petite sous-$A$-algèbre de~$B$ contenant les~$b_i$ et sera notée~$A[b_i]_i$ ; 
on dit que c'est la sous-algèbre de~$B$ {\em engendrée par les~$b_i$}. 

\trois{oide-typefini}
On dit qu'une~$A$-algèbre~$B$ est {\em de type fini}
s'il existe une famille finie~$(b_i)$ d'éléments 
de~$B$ telle que~$B=A[b_i]_i$. 

\deux{ext-corpoides}
Soit~$K$ un corpoïde. Une
{\em extension de~$K$} est une~$K$-algèbre
qui est un corpoïde, c'est-à-dire encore un corpoïde~$L$
muni d'un plongement~$K\hookrightarrow L$.

\medskip
Si~$L$ est une extension de~$K$ et si~$(x_i)_i$ est une famille
d'éléments de~$L$, la plus-petite sous-extension de~$L$
contenant les~$x_i$ s'identifie au corps des
fractions de la~$K$-algèbre~$K[x_i]_i$, et sera noté~$K(x_i)_i$. 

\deux{krcorpsgrad} Soit~$K$ un corpoïde
et soit~$\bf r$ un polyrayon.
La~$K$-algèbre~$K[{\bf r}\inv{\bf T}]$ est intègre, et son corpoïde
des fractions sera noté~$K({\bf r}\inv{\bf T})$. On dira que~$\bf r$ est
 {\em~$K$}-libre s'il constitue une famille libre de~$\QQ\otimes_{\ZZ}(D/\deg(K\ti))$. 
Si~$\bf r$ est~$K$-libre 
le corpoïde~$K({\bf r}\inv{\bf T})$ est égal à~$K[{\bf r}\inv{\bf T}, {\bf r}{\bf T}\inv]$,
c'est-à-dire au quotient de~$K[{\bf r}\inv{\bf T}, {\bf r}{\bf S}]$
par l'idéal~$({\bf T}\cdot {\bf S}-{\bf 1})$. 

\deux{principalgrad} Soit~$K$ un corpoïde et soit~$r\in D$. 
L'annéloïde intègre~$K[T/r]$
est euclidien : si~$P\in K[T/r]$ et si~$S$ est un élément
non nul de~$K[T/r]$, il existe un unique
couple~$(Q,S)$ d'éléments de~$K[T/r]$ tel
que~$P=SQ+R$ et tel que le degré monomial de~$R$ soit 
strictement inférieur à celui de~$P$. On en déduit que
l'annéloïde intègre~$K[T/r]$ est principal,
c'est-à-dire
que ses idéaux sont monogènes. À ce titre, il satisfait la déclinaison dans notre
cadre
des propriétés arithmétiques usuelles des anneaux principaux :
existence d'un PGCD, lemme de Gauß, existence et unicité de
la décomposition d'un élément non nul en produit d'irréductibles
(qui s'étend au corpoïde~$K(T/r)$, à condition d'accepter les exposants négatifs).

\medskip
Si~$P\in K[T/r]$, une {\em racine}
de~$P$ dans~$K$ est un élément~$x$
de~$K_r$ tel que~$P(x)=0$. On vérifie à l'aide
de la division euclidienne qu'un élément~$x$ 
de~$K_r$ est une racine
de~$P$ si et seulement si~$T-x$ divise~$P$ ; on en déduit
que si~$P$ est non nul et de degré monomial~$n$, il a
au plus~$n$ racines dans~$K$, comptées avec multiplicités
-- la multiplicité d'une racine~$x$ de~$P$ est le plus
grand entier~$e$ tel que~$(T-x)^e$ divise~$P$. 

\deux{corpoide-rationnel-tf}
Soit~$K$ un corpoïde. 

\trois{tf-polyrayon-1}
Soit~$r\in D$. Le corpoïde~$K(T/r)$ est de type
fini comme~$K$-algèbre si et seulement si~$r$ est~$K$-libre. 

\medskip
En effet, si~$r$ est~$K$-libre
le corpoïde~$K(T/r)$ est égal à~$K[T/r, rT^{-1}]$
et il est donc engendré par~$T$ et~$T^{-1}$ comme~$K$-algèbre.

\medskip
Supposons maintenant que~$r$ est d'ordre fini~$n$ dans~$D/\deg(K\ti)$, et soit~$\lambda\in K\ti$
tel que~$r^n=\deg \lambda$. 
Nous allons montrer par l'absurde
que~$K(T/r)$ n'est pas de type fini. On suppose donc qu'il est de type fini. Il est
alors engendré par un nombre fini de fractions, et il existe donc une famille finie
~$(P_1,\ldots,P_n)$ d'éléments non nuls de~$K[T/r]$  tels que tout élément de~$K[T/r]$
puisse s'écrire comme une fraction ayant un dénominateur de la forme~$\prod P_i^{n_i}$. 
Il s'ensuite, en vertu de la théorie de la décomposition des éléments de~$K(T/r)\ti$ en produit d'éléments irréductibles 
(avec exposants éventuellement négatifs), que~$K[T/r]$ n'a qu'un nombre fini d'éléments irréductibles. 

\medskip
Soit~$\sch P$ l'ensemble des polynômes irréductibles unitaires appartenant
à~$K_1[\tau]$. On vérifie aussitôt que
$$\{P(T^n/\lambda)\}_{P\in \sch P}$$ est une ensemble d'éléments irréductibles 
deux à deux non associés de~$K[T/r]$ ; il est infini puisque~$\sch P$ l'est, 
et l'on aboutit ainsi à une contradiction. 

\trois{tf-polyrayon-qque}
Soit maintenant~${\bf r}=(r_1,\ldots, r_n)$ une famille finie d'éléments de~$D$. Nous allons
montrer que le corpoïde~$K({\bf T}/{\bf r})$ est de type fini comme~$K$-algèbre
si et seulement si~$\bf r$ est~$K$-libre. 

\medskip
Si~$\bf r$ est~$K$-libre, la~$K$-algèbre $K({\bf T}/{\bf r})$
est de type fini : on le voit en combinant~\ref{tf-polyrayon-1}
et une récurrence triviale sur~$n$. 

\medskip
Supposons que~$\bf r$ ne soit pas~$K$-libre. Quitte à renuméroter les~$r_i$, on peut 
supposer que~$r_n$ n'est pas~$K(T_1/r_1,\ldots, T_{n-1}/r_{n-1})$-libre ; il résulte alors 
de~\ref{tf-polyrayon-1}
que le corpoïde~$K({\bf T}/{\bf r})$ n'est pas de type fini comme~$K(T_1/r_1,\ldots, T_{n-1}/r_{n-1})$-algèbre ; 
il n'est {\em a fortiori}
pas de type fini comme~$K$-algèbre.

\deux{espvectgrad} {\em Modules et espaces vectoriels.} Soit~$A$ un annéloïde.  
On appelle {\em~$A$-module} un ensemble~$M$
muni d'une décomposition $M=\coprod \limits_{r\in D}  M_r$, d'une structure de groupe
abélien sur chacun des~$M_r$, et d'une famille d'applications
bilinéaires~$A_r\times M_s\to M_{rs}$ 
telles que~$a(bm)=(ab)m$ et~$1m=m$ 
pour tout~$(a,b,m)$. 
On dispose d'une notion évidente d'application linéaire graduée
de degré donné~$r$ entre deux~$A$-modules
(une telle application
multiplie le degré par~$r$), de sous-module,
de somme directe de modules,
de produit tensoriel de modules, etc. 

\medskip
Si $r\in D$ on note~$A(r)$ le~$A$-module
déduit de~$A$ par décalage de la graduation
de sorte que~$A(r)_s= A_{rs}$ pour tout
$s$. Si $M$ est un~$A$-module
et si $(m_i)$
est une famille d'éléments
de~$M$ de degrés respectifs~$r_i$, 
il existe une unique application linéaire
de degré~$1$ de
$\bigoplus  A(r_i^{-1})$ dans~$M$ qui envoie~$1_i$ sur~$m_i$
pour tout~$i$ ; on dit que la famille~$(m_i)$ est libre (resp. génératrice, resp. une base)
si cette application est injective (resp. surjective, resp. bijective).
Un~$A$-module est dit de type fini s'il admet une famille génératrice finie. 

\medskip
Si~$D'$ est un groupe abélien divisible contenant~$D$ et
si~$M$ est un~$A$-module, 
le~$A\odot_D D'$-module~$M\odot_D D'$ est défini
de façon évidente. 

\medskip
Si~$K$ est un corpoïde,
on parlera de~$K$-espace vectoriel
plutôt que de~$K$-module gradué. Si~$M$ est un~$K$-espace vectoriel, 
il possède une base, et toutes ses bases ont même cardinal ; ce dernier est appelé la {\em dimension} de~$M$. 

\deux{diffgrad} Soit~$ A$ un annéloïde, soit~$B$ une~$A$-algèbre et soit~$M$
un $B$-module. Une {\em $ A$-dérivation} de $B$ dans $M$ 
est une application $A$-linéaire  ${\rm d} :  B\to  M$ de degré $1$ telle que ${\rm d}a=0$ pour tout $a\in  A$
et $${\rm d}(bb')=b\;{\rm d}b'+b'{\rm d}b$$ pour tout couple $(b,b')$ d'éléments de $B$. 
La catégorie des $B$-modules munis d'une $A$-dérivation de source $B$
admet un objet initial que l'on note $\Omega_{B/ A}$ ; c'est le {\em module des différentielles de $B$ sur $ A$}. 

\deux{nakagrad} 
On dispose d'un avatar du lemme de Nakayama : si~$A$ est un annéloïde
local d'idéal maximal~$\got m$ et si le~$A$-module~$\got m$ est de type fini,
une famille finie $(m_i)$ d'éléments de $\got m$ engendre~$\got m$ si et seulement si 
les classes des $m_i$ engendrent le $k$-espace 
vectoriel $\got m/\got m^2$, où $ k= A/\got m$. Le cardinal minimal d'une famille
génératrice de~$\got m$ est donc égal à la dimension de $\got m/\got m^2$. 

\deux{clotintegr} Soit $A$ un annéloïde et soit $B$ une $A$-algèbre. 
Soit~$b\in B$ et soit~$r$ son degré. La plus petite sous-$A$-algèbre
de~$B$ contenant~$b$ est l'ensemble~$A[b]$ des éléments
de la forme~$P(b)$ où~$P\in A[T/r]$. 
On dit que $b$ est {\em entier sur~$A$} s'il existe un entier~$n$
et un élément $P\in  A[T/r]_{r^n}$ unitaire et de degré monomial égal à~$n$
 tel que~$P(b)=0$ ; si~$b$ n'est pas nilpotent,
 un tel~$P$ n'est pas un monôme, ce qui force~$r$ à 
 appartenir à~$\deg(A^{\neq 0})^{\QQ}$ (écrire l'égalité
 des degrés de deux monômes distincts de~$P$). Si~$b$
 est entier,
 la~$A$-algèbre~$A[b]$ est {\em finie},
 c'est-à-dire finie comme~$A$-module. 
 Réciproquement, on voit aisément,
 {\em via} par exemple grâce à un 
 avatar gradué du théorème de Cayley-Hamilton,
que si~$b$ vit dans
dans une sous-$A$-algèbre de~$B$ qui est
finie, $b$ est entier sur~$A$. 
Lorsque~$A$ est un corpoïde, on emploiera
l'adjectif «algébrique» de préférence
à «entier». 
L'ensemble des éléments de~$B$ entiers sur~$A$ 
est une sous-$A$-algèbre de~$B$
que l'on appelle la {\em fermeture intégrale} de~$A$ 
dans~$B$. 

\deux{going-upoide}
Soit~$f: A\to B$ un morphisme
d'annéloïdes. 

\trois{going-upoinde-facile}
Comme dans le cas classique, on démontre que si~$f$ est injective, si~$A$ et~$B$ sont intègres et si la~$A$-algèbre$B$
est entière  ({\em i.e.} constituée d'éléments entiers), alors~$A$
est un corpoïde si et seulement si~$B$ est un corpoïde. 

\trois{going-upoide-general}
Par des méthodes standard, 
on déduit de~\ref{going-upoinde-facile}
le lemme de {\em going-up}, dont voici l'énoncé.
Supposons que la~$A$-algèbre~$B$ est entière. Soit $\got p$ un idéal premier de $A$, et soit~$\sch P$
l'ensemble des idéaux premiers~$\got q$ de l'annéloide~$B$ tels que~$f^{-1}(\got q)=\got p$.

1) Si~$\got q\in \sch P$ alors~$\got q$ est un idéal maximal de~$B$
si et seulement si $\got p$ est un idéal maximal de~$A$. 

2) Les éléments de~$\sch P$ sont deux à deux non comparables pour l'inclusion. 

3) Si $\got p'$ est un idéal premier de~$A$ contenant $\got p$ et si~$\got q\in \sch P$,
il existe un idéal premier
$\got q'$ de $ B$ contenant $\got q$ tel que $f\inv(\got q')=\got p'$. 

\subsection*{Extensions de corpoïdes, degré de transcendance et théorie de Galois}

\deux{algindep-grad} Soit~$K$
un corpoïde et soit~$A$ une~$K$-algèbre. 
Soit~$(x_i)$ une famille d'éléments de~$A$. Pour tout~$i$,
on note~$r_i$ le degré de~$x_i$. 
On dit que les~$x_i$ sont {\em algébriquement indépendants sur~$K$} 
si~$P(x_i)\neq 0$ pour tout élément non nul~$P$
de~$K[{\bf T}/{\bf r}]$ ; notons que la famille
vide est algébriquement indépendante sur~$K$ si et seulement si~$A$ est non nulle. 
Si~$x$ est un élément de~$F$ on dira qu'il est 
{\em transcendant } sur~$K$ si la famille singleton~$\{x\}$ est algébriquement 
indépendante ; l'élément~$x$ est transcendant
si et seulement si il n'est pas algébrique. 
Il est algébrique si et seulement  si~$K[x]$ est de dimension
finie sur~$K$. 

\deux{theotransalg} Soit~$K\hookrightarrow F$ une extension de corpoïdes. 

\trois{algoiderel} Un élément~$x$ de~$F$ est algébrique sur~$K$
si et seulement si~$K[x]$ est un corpoïde. L'ensemble
des éléments de~$F$ algébriques sur~$K$ est un sous-corpoïde
de~$F$. On dit que~$F$ est une extension
{\em algébrique} de~$K$ si tous les éléments de~$F$ sont
algébriques sur~$K$ ; notons
que si c'est le cas,~$\deg(F\ti)/\deg(K\ti)$ est de torsion.

\medskip
Si~$F$ est une extension algébrique de~$K$, on 
notera~$\mathsf {Gal}(F/K)$
le groupe des~$K$-automorphismes de~$F$.

\trois{pominoide} Soit~$x$ un élément de~$F$ algébrique
sur~$K$, et soit~$r$ son degré. On appelle~{\em polynôme minimal}
de~$x$ sur~$K$ l'unique générateur unitaire de
l'annulateur de~$x$ dans~$k[T/r]$. Il est irréductible ; son degré monomial
est simplement appelé le
{\em degré de $x$ sur $K$}, et il coïncide avec la dimension
du~$K$-espace vectoriel~$K[x]$. 

\trois{degtroide} Une {\em base de transcendance}
de~$F$ sur~$K$ est une famille maximale d'éléments de~$F$-algébriquement indépendants sur~$K$. 

\medskip
Si~$(x_i)_{i\in I}$ est une famille d'éléments de~$F$ telle que~$F$ soit algébrique sur~$K(x_i)_i$ et si~$J$ est un sous-ensemble
de~$I$ tel que les~$x_i$ pour~$i$ parcourant~$J$ soient algébriquement indépendants sur~$K$, il existe par le lemme de Zorn
un sous-ensemble~$J$' de~$I$ contenant~$J$, tel que les~$x_i$ pour~$i$ parcourant~$J'$ soient algébriquement indépendants sur~$K$,
et qui est maximal pour cette propriété. On vérifie que~$(x_i)_{i\in J'} $ est une base de transcendance de~$L$ sur~$K$. 

Il existe donc des bases de transcendance de~$L$ sur~$K$ : il suffit d'appliquer ce qui précède à n'importe quelle
famille~$(x_i)_{i\in I}$ énumérant tous les éléments
de~$F$ et au sous-ensemble~$J=\varnothing$ de~$I$. 

\medskip
On démontre que toutes les bases de transcendance de~$F$ sur~$K$ ont  le
même cardinal, appelé le {\em degré de transcendance de~$F$ sur~$K$}.  

\deux{corpoide-alg-clos}
Soit~$K$ un corpoïde. Les assertions suivantes sont équivalentes. 

\medskip
$\bullet$ Le corpoïde~$K$ n'admet aucune extension algébrique stricte. 

$\bullet$ Pour tout~$r\in D$, tout polynôme non nul
appartenant à~$K[T/r]$ est scindé dans~$K$. 

$\bullet$ Pour tout~$r\in D$,  
tout polynôme appartenant à~$K[T/r]$ et de degré monomial strictement positif
a une racine dans~$K$. 

\medskip
Lorsqu'elles sont satisfaites, on dit que~$K$ est {\em algébriquement clos}. 

\deux{nullstoide}
{\bf Le {\em Nullstellensatz}
gradué.}
Soit~$K$ un corpoïde et soit~$F$ une extension de~$K$. Les assertions suivantes
sont équivalentes : 

\medskip
i) il existe une famille finie et~$K$-libre~$\bf r$ d'éléments de~$D$ telle que~$F$ 
soit isomorphe (comme extension de~$K$) à une extension finie de~$K({\bf r}^{-1}{\bf T})$ ; 

ii) la~$K$-algèbre~$F$ est de type fini. 

\medskip
En effet, si~i) est vraie alors comme~$K({\bf r}^{-1}{\bf T})$ est une~$K$-algèbre de type fini
(\ref{tf-polyrayon-qque}), la~$K$-algèbre~$F$ est de type fini et~ii) est vraie.

\medskip
Supposons~ii) vraie et choisissons une famille génératrice finie~$(t_1,\ldots, t_n)$ de la~$K$-algèbre~$F$. 
Quitte à renuméroter les~$t_i$, on peut supposer qu'il existe~$m$ tel que les~$t_i$ pour~$i\leq m$ soient 
algébriquement indépendants sur~$K$, et tel que pour tout~$j\leq m$, l'élément~$t_j$ soit algébrique
sur~$K(t_i)_{1\leq i\leq m}$. Pour tout~$i\leq m$, posons~$r_i=\deg t_i$.
Le corpoïde~$F$ est alors fini sur~$K(t_i)_{i\leq m}\simeq K({\bf r}^{-1}{\bf T})$ ; il suffit pour conclure de démontrer que~$\bf r$ est~$K$-libre,
c'est-à-dire encore en vertu de~\ref{tf-polyrayon-qque}
que la~$K$-algèbre~$K(t_i)_{i\leq m}\simeq K({\bf r}^{-1}{\bf T})$
est de type fini. 

\medskip
Pour tout~$j>m$, l'élément~$t_j$ de~$F$ est algébrique sur~$K(t_i)_{i\leq m}$, et est donc entier sur un localisé
$(K[t_i]_{i\leq m})_{f_j}$ de~$K[t_i]_{i\leq m}$ pour un certain~$f_j$ non nul de~$K[t_i]_{i\leq m}$. 
Si l'on note~$f$ le produit des~$f_j$
chacun des~$t_j$ est entier sur~$(K[t_i]_{i\leq m})_f$,
et~$F$ est donc entière sur~$(K[t_i]_{i\leq m})_f$. Comme~$F$ est un corpoïde, $(K[t_i]_{i\leq m})_f$
est un corpoïde (\ref{going-upoinde-facile}). Puisque
$$K[t_i]_{i\leq m}\subset (K[t_i]_{i\leq m})_f\subset K(t_i)_{i\leq m},$$
il vient~$$ (K[t_i]_{i\leq m})_f=K(t_i)_{i\leq m}$$ et ce dernier est donc
une~$K$-algèbre de type fini, ce qui achève la  démonstration. 

\deux{geom-nullst}
Soit~$A$ une~$K$-algèbre de type fini. 
Notons~$(\spec A)_\diamond $ le sous-ensemble de~$\spec A$ constitué
des points~$x$ tels que~$\kappa(x)$ soit une extension finie de
$K({\bf r}^{-1}{\bf t})$ pour une certaine famille finie~$K$-libre $\bf r$
d'éléments de~$D$. Soit~$E$ un sous-ensemble de~$\spec A$ qui est une combinaison
booléenne finie d'ouverts et fermés de~$\spec A$ ; il est alors non vide si et seulement si
son intersection avec~$(\spec A)_\diamond$ est non vide.

En effet, supposons que~$E$ soit non vide ; nous allons montrer que~$E$
rencontre~$(\spec A)_\diamond$. L'ensemble~$E$ est réunion finie de parties de la forme
$U\cap F$, où~$U$ et~$F$ sont respectivement ouverts et fermés dans~$\spec A$ ; l'une de ces 
parties est donc non vide, ce qui permet de 
supposer que~$E$ lui-même est de la forme~$U\cap F$
avec~$U$ et~$F$ comme ci-dessus. 

On peut écrire~$F=V(I)$ pour un certain idéal~$I$ de~$A$, et~$F$ s'identifie dès lors au spectre
de la~$K$-algèbre de type fini~$A/I$, avec préservation des corpoïdes résiduels. En remplaçant~$A$
par~$A/I$, on se ramène au cas où~$F=\spec A$ et où~$E$ est ouvert. Il est dès lors
réunion de
parties de la forme~$D(f)$  avec~$f\in A$ ; l'une de ces parties est donc non vide, 
ce qui permet de supposer que~$E$ lui-même
est 
de la forme~$D(f)$ avec~$f\in A$. Il s'identifie dès lors au spectre
de la~$K$-algèbre de type fini~$A_f$, avec préservation des corpoïdes résiduels. On peut
maintenant, quitte à remplacer~$A$ par~$A_f$, se ramener au cas où~$E=\spec A$. 

Puisque~$E$ est non vide, $A$ est non nul et possède donc un idéal maximal. Si~$x$
désigne le point correspondant de~$\spec A$ alors
le corpoide $\kappa(x)$ est un quotient de~$A$, et donc une~$K$-algèbre
de type fini. Le {\em Nullstellensatz} établi au~\ref{corpoide-alg-clos}
ci-dessus assure alors que~$x\in (\spec A)_\diamond$.

\deux{defdeltad} On fixe un sous-groupe divisible~$\Delta$
de~$D$. Dans la pratique, ce qui suit sera le plus souvent appliqué
lorsque~$\Delta=\{1\}$. 

\deux{extfiniresoide} Soit~$K\hookrightarrow L$ une extension 
de corpoïdes. On note~$E_\Delta(L/K)$ le groupe~$\deg(L\ti)/\deg (K\ti)\cdot (\Delta\cap \deg(L\ti))$. 
Remarquons que~$E_{\{1\}}(L/K)$ est simplement
le groupe~$\deg (L\ti)/\deg (K\ti)$. 

\trois{interpdimoide} Soit~$\bf a$
un sous-ensemble de~$L\ti$ tel que~$\deg({\bf a})$ constitue
un système de représentants de la flèche quotient~$(\deg L\ti) \to E_\Delta(L/K)$,
et soit~$\bf b$ une base de~$L_\Delta$ sur~$K_\Delta$. Nous allons montrer que~${\bf a}\cdot{\bf b}$ 
est une base de~$L$ sur~$K$, ce qui impliquera 
que
$$[L:K]=\card E_\Delta(L/K)\cdot [L_ \Delta : K_\Delta], $$ et en particulier que
$$[L:K]=(\deg (L\ti):\deg (K\ti))\cdot[L_1:K_1].$$ 

\medskip
{\em La famille~${\bf a}\cdot {\bf b}$ est libre.}
Soit~$r\in D$, 
et soient~$(\lambda_{a,b})_{a\in {\bf a},b\in {\bf b}}$ une
famille
d'éléments de~$K$ presque tous nuls tels que~$\deg(\lambda_{a,b})=r(\deg a)\inv(\deg b)\inv$
pour tout~$(a,b)$. 
Supposons que~$\sum \lambda_{a,b}ab=0$. Comme les degrés des éléments
de~$\bf a$ sont deux à deux distincts modulo~$\deg(K\ti)\cdot (\Delta\cap \deg(L\ti))$, 
on déduit de la famille
d'égalités~$\deg(\lambda_{a,b})=r(\deg a)\inv(\deg b)\inv $ qu'il existe~$a_0\in \bf a$
tel que~$\lambda_{a,b}$
soit nul dès que~$a\neq a_0$. Pour tout~$b$,
posons~$\mu_b=\lambda_{a_0,b}a_0$. C'est un élément de~$K$
de degré~$(\deg b)\inv$, et l'on a~$\sum \mu_b b=0$ ; comme~$\bf b$
est libre sur le corps~$K_D$, il vient~$\mu_b=0$ pour tout~$b$,
et donc~$\lambda_{a_0,b}=0$ pour tout~$b$. Finalement,
on a bien~$\lambda_{a,b}=0$ pour tout~$(a,b)$, ce qu'il fallait démontrer.

\medskip
{\em La famille~${\bf a}\cdot {\bf b}$ est génératrice.} 
Soit~$x\in L\ti$. 
Il existe~$a\in {\bf a}$ tel que~$\deg x$ soit égal à~$\deg a$ modulo~$\deg(K\ti)\cdot (\Delta\cap \deg(L\ti))$. 
On peut en
conséquence écrire~$x=\lambda ay$ avec~$y\in L_\Delta$ et~$\lambda \in K\ti$. Comme~$\bf b$ est une base
de~$L_\Delta$ sur~$K_\Delta$, on peut écrire~$y=\sum \mu_bb$ où les~$\mu_b$
appartiennent à~$K_\Delta$. Il vient~$x=\sum_b \lambda \mu_b ab$, ce qui termine la
preuve. 

\trois{casdegldegk} Soit~$\cal S$ un système de représentants
de~$D/\Delta$. Supposons que~$\deg (L\ti)=\deg(K\ti)\cdot (\Delta\cap \deg(L\ti))$ et soit~$F$
une sous-extension de~$L$. On a alors~$\deg (F\ti)=\deg (K\ti)\cdot (\Delta\cap \deg(F\ti))$.
Il s'ensuit, en vertu
de~\ref{interpdimoide} (ou par un calcul explicite direct)
que~$F=\coprod_{r\in \cal S} F_\Delta\cdot K_r.$ Les flèches
$$F\mapsto F_\Delta,\Lambda \mapsto \coprod_{r\in S}\Lambda\cdot K_r$$ 
établissent ainsi une bijection entre l'ensemble des sous-extensions de~$L$
et celui des sous-extensions de~$L_\Delta$.

\trois{interpdegtroide} On ne 
suppose plus que~$\deg(L\ti)=\deg(K\ti)\cdot (\Delta\cap \deg(L\ti))$. 
Soit~${\bf u}$ une base de transcendance de~$L_\Delta$ sur~$K_\Delta$,
et soit~${\bf v}$ une partie de~$L\ti$ tels que~$\deg({\bf v})$ soit une base de ~$\QQ\otimes_{\ZZ}E_\Delta(L/K)$.
Nous allons montrer que~${\bf u}\cup{\bf v}$ est
une base de transcendance de~$L$ sur~$K$, ce qui impliquera que
$$\deg {\rm tr.}(L/K)=\deg {\rm tr.} (L_\Delta/K_\Delta) +\dim {\QQ} \QQ\otimes_{\ZZ}E_\Delta(L/K),$$
et en particulier que
$$\deg {\rm tr.}(L/K)=\deg {\rm tr.} (L_1/K_1) +\dim {\QQ} \QQ\otimes_{\ZZ}(\deg(L\ti)/\deg(K\ti)).$$

\medskip
{\em La famille~${\bf u}\cup {\bf v}$ est algébriquement indépendante.} 
Soit~$P\in K[|{\bf u}|\inv{\bf T}, |{\bf v}|\inv{\bf S}]$ tel
que~$P({\bf u},{\bf v})=0$. Nous allons montrer que~$P$ est nul. 

Les monômes qui constituent~$P$ ont tous même degré~$r$ ; 
comme~$\deg({\bf v})$ est une~$\QQ$-base de~$\deg(L\ti)$
 modulo~$\deg(K\ti)\cdot (\Delta\cap \deg(L\ti))$, 
il existe un monôme unitaire~$Q\in K[\deg({\bf v})\inv{\bf S}]$  et
un élément~$R$ de~$K[|{\bf u}|\inv{\bf T}]$ 
tels que~$P=QR$. Comme~$Q$ est un monôme unitaire, 
on a~$Q({\bf v})\neq 0$ ; par conséquent,~$R({\bf u})=0$. 
On peut écrire~$R=\alpha R^\sharp$, où~$\alpha$ est un élément 
non nul de~$K$ et où~$R^\sharp$ est
somme de monômes dont tous les coefficients appartiennent à~$K_\Delta$. 
On a~$R^\sharp({\bf u})=0$ ; comme~$\bf u$ est une base de transcendance de~$L_\Delta$
sur~$K_\Delta$, il vient~$R^\sharp=0$ et finalement~$P=0$, ce qu'il fallait démontrer. 

\medskip
{\em Le corpoïde~$L$ est algébrique sur~$K({\bf u}\cup{\bf v})$.} Soit~$\lambda\in L\ti$. 
Comme~$\deg({\bf v})$ est une base de~$\QQ\otimes_{\ZZ}E_\Delta(L/K)$,
il existe un entier~$n>0$, un multi-indice~$J$ à coordonnées dans~$\ZZ$, 
un élément~$\alpha\in K\ti$ et un élément~$\beta\in L_\Delta\ti$ 
tels que
l'on ait~$\lambda^n=\alpha\beta\cdot {\bf v}^J$. 
Comme~$\beta \in L_\Delta$, 
il est algébrique sur~$K_\Delta({\bf u})$ ; par conséquent,~$\lambda$ 
est algébrique sur~$K({\bf u}\cup {\bf v})$, ce qui achève la démonstration

\trois{genoide} Soit~$\bf w$ une partie de~$L_\Delta$ l'engendrant
(comme corppoïde) sur~$K_\Delta$, et soit~$\bf g$ une partie de~$L\ti$ telle
que~$\deg({\bf g})$ engendre~$E_\Delta(L/K)$. Le
sous-ensemble~${\bf w}\cup{\bf g}$ de~$L$ engendre alors ce dernier sur~$K$.

\medskip
En effet, soit~$\lambda\in L\ti$. Comme~$\deg({\bf g})$
engendre~$E_\Delta(L/K)$, 
il existe un multi-indice~$J$ à coordonnées dans~$\ZZ$, 
un élément~$\alpha$ de~$K\ti$ et un élément~$\beta$
de~$L_\Delta\ti$. 
tels que~$\lambda=\alpha\beta\cdot {\bf g}^J$. 
Comme~$\beta \in L_\Delta$, 
il appartient à~$K_\Delta({\bf w})$ ; par conséquent,~$\lambda$ 
appartient~$K({\bf w}\cup {\bf g})$, ce qui achève la démonstration. 

\medskip
Notons le cas particulier où~$\Delta=\{1\}$ : si~$\bf w$
est une 
partie de~$L_1$ l'engendrant
(comme corps) sur~$K_1$, et si~$\bf g$ une partie de~$L\ti$ telle
que~$\deg({\bf g})$ engendre~$\deg (L\ti)$ 
modulo~$\deg (K\ti)$, alors~${\bf w}\cup{\bf g}$
engendre~$L$ sur~$K$. 

\trois{equivseparoide} Supposons que~$L$ soit une extension
algébrique de~$K$, et soit~$p$ l'exposant 
caractéristique de~$K$. On a alors équivalence entre les assertions suivantes :

\medskip
i) $L$ est séparable sur~$K$ ; 

ii) $L_\Delta$ est séparable sur~$K_\Delta$, et~$E_\Delta(L/K)$ est sans~$p$-torsion.

\medskip
En effet, supposons que~i) soit vraie ; 
le polynôme minimal sur $K$ de tout élément 
de~$L$ est alors séparable. 
C'est en particulier vrai pour les éléments dont
le degré appartient à~$\Delta$, 
et~$L_\Delta$ est dès lors séparable
sur~$K_\Delta$. Par ailleurs, soit~$x\in L^{*}$ tel
que
$\deg(x^p)$ appartienne à $~\deg(K\ti)\cdot (\Delta\cap \deg (L\ti))$. 
Il existe alors~$\ell \in L_\Delta$ et~$\omega\in K$ tel que~$x^{p}=\ell\omega$.

On en déduit que~$x$ est purement inséparable 
sur le sous-corpoïde 
$L_\Delta\cdot K$ de~$L$. Le corpoïde~$L$
est par hypothèse séparable sur~$K$, 
il l'est {\em a fortiori} sur~$L_\Delta\cdot K$. En conséquence,~$x\in L_\Delta\cdot K$ ; ceci implique que ~
le degré de~$x$
appartient à~$\deg(K\ti)\cdot (\Delta\cap \deg (L\ti))$, et achève de prouver~ii).  

\medskip
Supposons réciproquement 
que~ii) soit vraie. Le
sous-corpoïde ~$L_\Delta\cdot K$ de~$L$
est alors séparable sur~$K$. 
Soit~$x\in L\ti$. Il existe un entier~$m$
premier 
à~$p$ 
tel que~$\delta(x)^{m}\in \deg(K\ti)\cdot (\Delta\cap \deg(L\ti))$.
Il existe
donc ~$\ell \in L_\Delta$ et $\omega\in K$ tel que~$x^{m}=\ell\omega$ ; en conséquence,~$x$ est séparable 
sur~$L_\Delta\cdot K$. 
Le corpoïde~$L$ est donc séparable sur~$L_\delta\cdot K$, lequel est lui-même séparable 
sur~$K$ ; dès lors,~$L$ est séparable sur~$K$, ce qui termine la preuve de~i). 

\deux{alg-clos-oide}
Soit~$K$ un corpoïde. Les assertions suivantes sont équivalentes : 

\medskip
i) le corps~$K_1$ est algébriquement clos et~$\deg(K\ti)$ est divisible :

ii) le corpoïde~$K$ est algébriquement clos. 

\medskip
En effet, supposons~i) vraie, et soit~$L$ une extension finie de~$K$. 
D'après~\ref{interpdimoide}, le corps~$L_1$ est une extension finie de~$K_1$, 
et est donc égal à~$K_1$ puisque celui-ci est algébriquement clos ; et le groupe~$\deg (L\ti)/\deg(K\ti)$ est fini,
et donc trivial puisque $\deg~K\ti$ est divisible. En utilisant à nouveau {\em loc. cit.}, il vient~$[L:K]=1$ ; ainsi, $L=K$,
d'où~ii). 

\medskip
Supposons que~ii) soit vraie. Soit~$P\in K_1[T]$ un polynôme unitaire non constant. On peut le voir comme
un élément de~$K[T]$ ; l'hypothèse ii) assure qu'il a une racine dans~$K$, laquelle appartient
à~$K_1$ par définition (le degré de~$T$ a été choisi égal à~$1$) ; en conséquence, $K_1$ est algébriquement 
clos. 

Soit~$x\in K\ti$, soit~$r$ son degré et soit~$n>0$. Comme~$K$ est algébriquement clos,
le polynôme~$T^n-x\in K[T/r^{1/n}]$ a une racine~$y$ dans~$K$. Par construction, $y\in K\ti$ 
et~$\deg y=\deg x/n$. En conséquence, $\deg(K\ti)$ est divisible.

\deux{etarrel} Soit~$K$
un corpoïde, soit~$(T_i)$ une famille
d'indéterminées et soit~$(r_i)$ une famille
d'éléments de~$D$.

\trois{etarelalgclos} {\em Le corpoïde~$K$ est algébriquement clos dans~$K({\bf r}\inv{\bf T})$.}
  Par un argument fondé sur le lemme de Zorn
  (ou une récurrence transfinie), il suffit 
  de traiter
  le cas où
 ~$(T_i)$ est une famille
 singleton~$\{T\}$, et~$(r_i)$ une
 famille singleton~$\{r\}$. Soit~$f\in K(r\inv T)$ 
 un élément algébrique sur~$K$ ; il s'agit 
  de montrer que~$f\in K$. Si~$f=0$ c'est évident,
  on peut donc supposer~$f$ non nul ; soit $s$ son degré. 
  Il existe
  un polynôme unitaire~$P\in K[X/s]$
  tel que~$P(f)=0$ ; écrivons~$P=\sum a_iT^i$. Comme~$f\neq 0$
  et comme~$P$ est irréductible, $a_0\neq 0$. 
  On note~$\delta$ le degré 
  {\em monomial} de~$P$ en la variable~$T$ ; il est nécessairement au moins égal à~$1$. 
 
 \medskip
 L'élément~$f$ est de la forme~$Q_0/Q_1$ 
 où les~$Q_i$ 
 appartiennent à~$K[T/s]$ et où~$Q_1\neq 0$ ; 
 on peut supposer que les~$Q_i$ sont premiers entre eux. 
 On a~$P(f)=0$ ; 
 il vient~$$Q_0^n+a_{n-1}Q_0^{n-1}Q_1+\ldots +a_0Q_1^n=0.$$ 
 
On en déduit que~$Q_0$ divise~$a_0Q_1^n$, et donc que~$Q_0$ 
 divise~$a_0$ 
 par le lemme de Gauß ; autrement dit,~$Q_0$ est constant. 
 De même~$Q_1$ divise~$Q_0^n$, donc divise~$1$
 et est lui aussi constant, ce qui prouve que~$f\in  K$. 
 
 \trois{abhyankaroide}
 Supposons que la famille~$(T_i)$ est finie. Le
 groupe
 quotient~$\deg(K({\bf r}\inv T)\ti)/\deg(K\ti)$ est alors de type fini,
 et l'extension de corps~$K_\Delta\hookrightarrow K({\bf r}\inv T)_\Delta$
 est de type fini, et est transcendante pure. 
 
 \medskip
 Pour le voir, on se ramène, en raisonnant par récurrence, au cas
 où ~$(T_i)$ est une famille
 singleton~$\{T\}$, et~$(r_i)$ une
 famille singleton~$\{r\}$. Il découle immédiatement
 de la définition de~$K(r\inv T)$ que~$\deg(K(r\inv T)\ti)$
 est égal à~$\deg(K\ti)\cdot r^\ZZ$, d'où la première assertion. 
 
 \medskip
 En ce qui concerne la seconde assertion, on distingue deux cas. 
 
 \medskip
 {\em Le cas où~$r$
n'appartient pas au groupe divisible~$\deg(K\ti)^\QQ\cdot \Delta$.} Le
degré d'un monôme de~$K[T/r]$ de la forme~$aT^i$ avec~$a$ et~$i$ non
nuls ne peut alors appartenir à~$\Delta$ ; il s'ensuit que~$K(T/r)_\Delta =K_\Delta$.

\medskip
{\em Le cas où~$r\in \deg(K\ti)^\QQ\cdot \Delta$}. Soit~$m$ l'ordre de~$r$ 
modulo~$\deg(K\ti)\cdot \Delta$, et 
soit~$\lambda$ un élément de~$K$
de degré~$r^m$ modulo~$\Delta$. 
Posons~$S=T^m/\lambda$ ; c'est un élément
de~$K(T/r)_\Delta$. Soit~$f$ un élément
non nul de~$K(T/r)_\Delta$. Il s'écrit~$P/Q$,
où~$P$ et~$Q$ sont deux éléments de~$K[T/r]$, avec~$Q\neq 0$. 
Soit~$p$ le degré de~$P$ et
soit~$s$ celui de~$Q$. Écrivons~$p=p'r^i$
et~$s=s'r^j$, où~$p'$ et~$s'$ appartiennent à~$\deg(K\ti)$ et~$i$
et~$j$ à~$\ZZ$. Soient~$a$ et~$b$
deux éléments
de~$K\ti$ de degrés respectifs~$p'$ et~$s'$. On peut alors écrire~$P=aT^iP_0$
et~$Q=bT^jQ_0$ avec~$P_0$ et~$Q_0$ de degré~$1$. Ainsi,~$f=(a/b)T^{i-j}P_0/Q_0$. 

Le degré de~$f$ est égal à~$r^{i-j}\deg(a/b)$,
et il appartient à~$\Delta$ par
définition de~$f$. Par ailleurs, 
$m$ est l'ordre de~$r$ modulo~$\deg(K\ti)\cdot \Delta$. Il s'ensuit : 

\medskip
$\bullet$ que $P_0$ s'écrit~$\sum a_i T^{mi}$, avec~$a_i$ de degré~$r^{-im}$ pour tout~$i$,
ou 
encore~$\sum \alpha_i (T^m/\lambda)^i$ avec~$\alpha_i$ de degré~$1$ pour tout~$i$ ; 

$\bullet$ de même, que $Q_0$ s'écrit~$\sum \beta_i (T^m/\lambda)^i$ avec~$\beta_i$ de degré~$1$ pour tout~$i$ ; 

$\bullet$ que~$i-j$ s'écrit~$m\ell$ pour un certain~$\ell$, et que~$(a/b)T^{i-j}=\lambda^\ell(a/b)(T^m/\lambda)^\ell$.

\medskip
Comme le degré de~$f$ 
et celui de~$T^m/\lambda$ appartiennent à~$\Delta$, il en va de même du degré de~$\lambda^{\ell}(a/b)$ ; en 
conséquence,~$f\in K_\Delta(T^m/\lambda)$ et l'on a finalement~$K(T/r)_\Delta=K_\Delta(T^m/\lambda)$.

\trois{interpfinioide} Soit~$L$ une extension de type fini de~$K$, et soit~$(t_1,\ldots,t_n)$
une base de transcendance de~$L$ sur~$K$. Posons~$F=K(t_1,\ldots, t_n)$ ; le
corpoïde~$L$ est alors fini sur~$F$. 

En vertu du~\ref{abhyankaroide} ci-dessus,~$F_\Delta$ est de type fini sur~$K_\Delta$, et~$\deg(F\ti)$ est
de type fini modulo~$\deg(K\ti)$. Il résulte par ailleurs de~\ref{interpdimoide} que~$L_\Delta$ est fini
sur~$F_\Delta$, et que~$\deg(L\ti)$ est fini modulo~$\deg(F\ti)$. Par conséquent,~$L_\Delta$
est de type fini sur~$K_\Delta$, et~$\deg(L\ti)/\deg(K\ti)$ est de type fini. 

\medskip
Réciproquement, supposons que~$L_\Delta$ est de type fini sur~$K_\Delta$, et
que~$\deg(L\ti)/\deg(K\ti)$ est de type fini, ou même seulement que~$E_\Delta(L/K)$
est de type fini. Il résulte alors de~\ref{genoide}
que~$L$ est de type fini sur~$K$. 

\trois{alg-clos-ext-kr}
On suppose maintenant que~$K$ est algébriquement clos, c'est-à-dire que~$K_1$ est algébriquement
clos et que~$\deg (K\ti)$ est divisible (\ref{alg-clos-oide}), et l'on 
désigne toujours par~$L$ est une extension de type fini 
de~$K$. Les assertions suivantes sont 
alors équivalentes : 

\medskip
i) il existe une famille finie et~$k$-libre~$\bf r$ d'éléments de~$D$
telle que~$L$ soit une extension finie de~$K({\bf r}^{-1}{\bf T})$ ; 

ii) $L_1=K_1$ ; 

iii) $\dim{\QQ}\QQ\otimes_{\ZZ}\deg (L\it)/\deg(K\ti)={\rm deg. tr}(L/k)\;;$

iv)  il existe une famille finie et~$k$-libre~$\bf s$ d'éléments de~$D$
telle que~$L$ soit isomorphe à~$K({\bf s}^{-1}{\bf T})$. 

\medskip
Il en effet clair que~i)$\Rightarrow$ii), et~ii)$\iff$iii) en vertu 
de~\ref{interpdegtroide} et du fait que~$K_1$ est algébriquement clos. 

\medskip
Supposons maintenant que les propriétés équivalentes~iii) et~ii) soient satisfaites.  
Le quotient~$\deg (L\ti)/\deg (K\ti)$ est un 
 groupe abélien de type fini d'après~\ref {interpfinioide}, 
 et est {\em sans torsion}
 puisque~$\deg(K\ti)$ est divisible. Il est donc libre de rang fini. Choisissons
 un~$n$-uplet $(s_1,\ldots, s_n)$ de~$\deg (L\ti)^n$ qui soit une base de~$\deg(L\ti)$ modulo
 $\deg (K\ti)$ ; pour tout~$i$, soit~$t_i\in L\ti$ tel que~$\deg t_i=s_i$. Par construction, $\bf s$ est~$K$-libre, 
et le sous-corps~$F:=K(t_i)$ de~$L$ s'identifie à~$K({\bf s}^{-1}{\bf T})$. 

Il vient~$F_1=K_1=L_1$, et~$\deg (F\ti)=\deg(K\ti)\cdot s_1^{\ZZ}\cdot\ldots \cdot s_n^{\ZZ}=\deg(L\ti)$. On déduit
alors de~\ref{interpdimoide}
que~$[L:F]=1$, et donc que
$$L=F\simeq K({\bf s}^{-1}{\bf T}),$$
d'où~iv). L'implication~iv)$\Rightarrow$i) est triviale. 

\deux{geom-integre-oide}
Soit~$K$ un corpoïde algébriquement clos et soient~$A$ et~$B$ deux~$K$-algèbres
réduites (resp. intègres). Nous allons montrer que~$A\otimes_KB$ est réduit (resp. intègre). 

\trois{geom-integre-caspart-oide}
{\em Preuve dans un cas particulier}.
On suppose que~$B$ est de la forme~$K({\bf r}^{-1}{\bf T})$ pour une certaine famille
finie et~$K$-libre ${\bf r}=(r_1,\ldots, r_n)$ d'éléments de~$D$. Le produit
tensoriel~$(A\otimes_K B)_r$ s'identifie pour tout~$r\in D$ à l'ensemble
des polynômes de la forme~$\sum_{I\in \ZZ^n} a_I {\bf T}^I$ tels que~$\deg(a_I){\bf r}^I=r$ pour tout~$I$. 
Un calcul immédiat montre alors que~$A\otimes_K B$ est intègre (resp. réduit). 

\trois{geom-integre-casgen-oide}
{\em Preuve dans le cas général}. Un raisonnement standard
fondé sur le passage à la limite inductive permet de
se ramener au cas où les~$K$-algèbres~$A$ et~$B$ sont de type fini. 

\medskip
Rappelons que~$(\spec B)_\diamond$ désigne l'ensemble des
points~$x\in \spec B$ tels que~$\kappa(x)$ soit une extension finie de
$K({\bf r}^{-1}{\bf T})$ pour une certaine famille
finie et~$K$-libre $\bf r$ d'éléments de~$D$ ; mais 
comme~$K$ est algébriquement clos, il revient au même,
en vertu de~\ref{alg-clos-ext-kr},
de demander que~$\kappa(x)$ soit isomorphe
à~$K({\bf r}^{-1}{\bf T})$ pour une certaine famille
finie et~$K$-libre $\bf r$ d'éléments de~$D$. 
Fixons une base~$(a_i)$ de~$A$
comme~$K$-espace vectoriel. 

\medskip
{\em Supposons~$A$ et~$B$ réduits}. Soit~$f$ un élément
nilpotent de~$A\otimes_K B$ ; écrivons~$f=\sum a_i \otimes b_i$
avec les~$b_i\in B$. Pour tout~$x\in (\spec B)_\diamond$, l'élément
$\sum a_i \otimes b_i(x)$ de~$A\otimes_K \kappa(x)$ est nilpotent ; 
comme~$A\otimes_K \kappa(x)$ est réduit d'après le cas particulier
traité au~\ref{geom-integre-caspart-oide}, 
on a~$\sum a_i \otimes b_i(x)=0$, et donc~$b_i(x)=0$ pour tout~$i$
puisque les~$a_i\otimes 1$ forment une base du~$\kappa(x)$-espace
vectoriel~$A\otimes_K \kappa(x)$. 

On a ainsi montré que~$b_i(x)=0$ pour tout~$i$ et tout~$x\in (\spec B)_\diamond$. Il
s'ensuit d'après~\ref{geom-nullst}
que~$V(b_i)=\spec B$ pour tout~$i$, c'est-à-dire que chacun des~$b_i$ est nilpotent. 
Comme~$B$ est réduit les~$b_i$ sont tous nuls, 
$f=0$ et~$A\otimes_K B$ est réduit. 

\medskip
{\em Supposons~$A$ et~$B$ intègres}. Ils sont alors
en particulier non nuls, et~$A\otimes_K B$ est donc non nul. 
Soient maintenant~$f$ et~$g$ deux éléments de~$A\otimes_K B$
tels que~$fg=0$. Écrivons~$f=\sum a_i\otimes b_i$ et~$g=\sum a_i \otimes c_i$
où les~$b_i $ les~$c_i$ appartiennent à~$B$. Posons~$F=V(b_i)_i\subset \spec B$, 
et~$G=V(c_i)_i$. 

Soit~$x\in (\spec B)_\diamond$. On a
$(\sum a_i \otimes b_i(x))(\sum a_i \otimes c_i(x))=0$ dans~$A\otimes_K \kappa(x)$ ; 
comme~$A\otimes_K \kappa(x)$ est intègre d'après le cas particulier
traité au~\ref{geom-integre-caspart-oide}, 
on a~$\sum a_i \otimes b_i(x)=0$ ou~$\sum a_i \otimes c_i(x)=0$ ; 
puisque les~$a_i\otimes 1$ forment une base du~$\kappa(x)$-espace
vectoriel~$A\otimes_K \kappa(x)$, il s'ensuit que l'une au moins 
des deux familles~$(b_i(x))_i$ et~$(c_i(x))_i$ est nulle. Autrement dit, 
$x\in F\cup G$. 

Il en résulte que
le fermé~$F\cup G$ de~$\spec B$ contient~$(\spec B)_\diamond$ ; il 
est de ce fait égal à~$\spec B$
en vertu de~\ref{geom-nullst}. Puisque~$B$ est intègre, $\spec B$ est irréductible, et
il est donc égal à~$F$ ou~$G$. S'il est égal à~$F$ on a~$V(b_i)=\spec B$ pour tout~$i$, 
et chacun des~$b_i$ est dès lors nilpotent et partant nul (puisque~$B$ est intègre) ; il vient~$f=0$. 
De même, on a~$g=0$ si~$\spec B=G$, ce qui achève la démonstration. 

\deux{ruptdecomp} Soit $K$ un corpoïde, soit~$r\in D$
et soit~$P$ un élément irréductible de $K[T/r]$. 
Un {\em corpoïde de rupture} de~$P$ est une extension~$F$ de~$K$ 
({\em}i.e. une~$K$-algèbre qui est un corpoïde) engendrée
 par  une racine de~ $P$. 
Si~$P$ est un élément non nul de~$K[T/r]$, 
un {\em corpoïde de décomposition} de~$P$ est une 
extension~de~$K$ dans laquelle~$P$ est scindé, 
et qui est engendrée par les racines de~$P$. 

\medskip
Si~$P$ est un élément irréductible de $K[T/r]$, il admet un corpoïde
de rupture,
à savoir~$K[T/r]/(P)$ ; 
on en déduit que si~$P$ est un élément non nul de~$K[T/r]$, 
il admet un corps de décomposition. 
Combiné avec un raisonnement fondé sur le lemme de Zorn, 
ceci entraîne l'existence d'une {\em clôture algébrique} de $K$, 
c'est-à-dire d'une extension 
algébrique~$L$ de $K$ telle que pour 
tout~$r\in D$, tout élément non nul~$P$ de $L[T/r]$ soit scindé dans~$L$. 

\deux{dedekgrad} Soit~$K$ un corpoïde,
soit~$A$ une $K$-algèbre et soit~$L$ une extension
de $K$. 
L'ensemble~$E$
des applications $K$-linéaires de $A$ dans 
$L$ est un~$K$-espace vectoriel, gradué par le degré
des applications. 
Si $\phi_1,\ldots, \phi_n$ sont des morphismes
deux à deux distincts de $A$ dans $L$, 
ils forment une famille libre dans~$E$ : c'est l'avatar
gradué du lemme d'indépendance
des caractères. 
On en déduit, en appliquant ce lemme sur le corpoïde~$L$
à la famille des application $\psi_i :  L\otimes_K A\to  L$
induites par les $\phi_i$, que
si~$A$ est de dimension finie $d$ sur $ k$, 
il existe au plus $d$ morphismes de~$A$ dans~$L$.

\trois{borneplong} En particulier, si~$F$ est une extension finie de~$K$, 
l'ensemble des~$K$-plongements de~$F$ dans~$L$ est fini, 
de cardinal majoré par~$[F: K]$. 

\trois{bornecorpsdcec} Supposons de plus que~$F$ soit un corpoïde de décomposition 
d'un élément $P$ non nul de $K[T/r]$ pour un certain $r$. 
Dans ce cas, si $P$ est scindé (resp. scindé à racines simples) 
dans~$L$ il y a au moins un (resp. exactement~$[F: K]$) plongement(s) de~$F$ dans~$L$.
Cela peut se démontrer par récurrence sur~$[F:K]$, 
le point clef étant que si $Q$ est un élément irréductible de $K[T/r]$
alors $T\mapsto \phi(\overline T)$ établit une bijection 
entre $\mathsf{Hom}_K(K[T/\gamma]/Q,L)$ et l'ensemble des racines de $Q$ dans $L$. 

\trois{automorphgrad} Il s'ensuit que deux corpoïdes 
de décomposition d'un même élément
non nul~$P$ de $ K[T/r]$ sont isomorphes ; on en déduit, à l'aide du lemme de Zorn, 
que deux clôtures algébriques de~$K$ sont isomorphes. 

\deux{algetale} {\em Algèbres étales}. Soit~$K$ un corpoïde et soit $\overline K$ 
une clôture algébrique de~$K$. 

\trois{polsep} Si $r\in D$, un élément non nul $P$ de~$\in  K[T/r]$ est dit {\em séparable} 
s'il est premier à son polynôme dérivé $P'$ 
(qui est de degré~$s/r$ si~$P$ est de degré $s$) ; 
il revient au même de demander 
que les racines de~$P$ dans $\overline K$ 
soient toutes simples, ou 
encore que $\overline K[T/r]/P$ soit réduit. 
Si $P$ est irréductible, 
il est séparable si et seulement si $P'\neq 0$ ; 
en général, $P$ est séparable si et seulement si
il n'a que des facteurs irréductibles simples
et eux-mêmes séparable. 

Soit $L$ une extension algébrique de
$K$ et 
soit~$x\in L$. 
On dit que~$x$ est {\em séparable} sur~$K$ 
si son polynôme minimal est séparable. 

\medskip
On dit que~$L$ est séparable sur~$K$ si
tous ses éléments sont séparables sur~$K$. 

\trois{intralget} Si $ A$ est une $K$-algèbre finie,
 elle est isomorphe à un produit fini de $K$-algèbres locales artiniennes ({\em i.e.} 
 dont l'idéal maximal 
 est constitué d'éléments nilpotents). 
 Si~$A$ est locale, si~$\got m$ désigne son idéal
 maximal et si $ A/\got m= K$ alors $\Omega_{A/ K}$ 
 s'identifie au dual de~$\got m/\got m^2$. 

\trois{algetkbarre} On en déduit à l'aide du lemme de Nakayama gradué 
que pour une $\overline K$-algèbre finie locale (resp. finie quelconque)~$A$ 
les assertions suivantes sont équivalentes : 

1) $A$ est réduite ; 

2) $ A$ est égale à $\overline K$ (resp. est le produit d'un nombre fini de copies de $\overline K$) ; 

3) $\Omega_{ A/\overline K}=0$.

\trois{algetgen} Si~$A$ est une $K$-algèbre finie, 
on dit que~$A$ est {\em étale} si~$\Omega_{ A/ K}=0$. 

Par ce qui précède, il revient au même de demander 
que~$\overline K\otimes_K A$ 
soit réduite, ou encore que  $\overline K \otimes_K A$ 
soit un produit fini de copies de $\overline K$. 

\medskip
Indiquons quelques propriétés qui résultent immédiatement de la définition. 

Si~$A$ est une $K$-algèbre étale, elle est réduite 
(puisqu'elle se plonge dans $\overline K\otimes_K A$) et 
est donc un produit de corps. 

Si~$A$ et~$B$ sont deux~$K$-algèbres 
finies,~$A
\times  B$ est étale si et seulement si~$A$ et~$B$ sont étales.

Si $ A$ est une~$K$-algèbre étale et si~$B$ est une sous-algèbre 
de~$A$ alors~$B$ est étale : cela résulte du fait que~$\overline K\otimes_K A$ 
est alors une sous-algèbre de~$\overline K\otimes_K B$, 
et est en particulier réduite.

\trois{algetsep} Il résulte de ce qui précède qu'une~$K$-algèbre finie est étale 
si et seulement si elle s'écrit comme un produit d'extensions finies de~$K$ 
qui sont étales. Pour comprendre ce que sont les~$K$algèbres étales, 
il suffit donc de comprendre à quelle condition une extension 
finie~$L$ de~$K$ est étale. Soit donc~ $L$ une extension
finie de~$K$. 

\medskip
Supposons que~$L$ soit étale, soit~$x$ un élément de~$L$
dont on note~$r$ le degré et soit~$P$ son polynôme minimal. 
En tant que sous-algèbre de~$L$, la~$K$-algèbre graduée~$K[x]\simeq  K[T/r]/P$ est étale, 
c'est-à-dire qu'elle reste réduite après extension des scalaires à $\overline K$  ; 
par conséquent,~$P$ est séparable. On en conclut que~$L$ est séparable. 

\medskip
Réciproquement, soit~$L$ une extension finie de~$K$ engendrée comme~$K$-algèbre 
par une famille finie $(x_1,\ldots, x_n)$ d'éléments  séparables sur~$K$. 
Pour tout~$i$, notons~$P_i$ le polynôme minimal 
de~$x_i$. Comme~$x_i$ 
est séparable,~$P'_i$ est premier à~$P_i$ et~$P'_i(x_i)$
 est donc non nul ; l'égalité~$P_i(x_i)=0$ entraîne que~$P'_i(x_i){\rm d}x_i=0$ 
 et donc que~${\rm d}x_i=0$ ; 
 ceci valant pour tout~$i$, on a~$\Omega_{ L/K}=0$, 
 et
$L$ est étale. 

\medskip
Ainsi, une~$K$-algèbre finie est étale 
si et seulement si elle s'écrit comme un produit fini d'extensions 
séparables de~$K$. 
Et pour qu'une extension finie~$L$ 
de~$K$ soit séparable, 
il suffit qu'elle soit engendrée 
par une famille finie d'éléments séparables.

\trois{fermesp} Énonçons maintenant quelques conséquences de ce qui précède. 

\medskip
Si~$L$ est une extension finie séparable 
de~$K$ et si~$F$ est une extension finie séparable 
de~$L$ alors~$F$ est une extension finie séparable de~$K$.

Si~$L$ est une extension quelconque de $K$, 
l'ensemble des éléments de~$L$ séparables sur~$K$
est une extension algébrique de~$K$ appelée
{\em fermeture séparable} de~$K$
dans~$L$ ; 
celle-ci ne possède aucune extension séparable stricte dans~$L$.

Si~$L$ est une extension algébrique
de~$K$, la fermeture séparable de~$K$
dans~$L$
est réduite à~$K$
si et seulement si~$L$
est une extension {\em radicielle} de~$K$, 
c'est-à-dire si et seulement si le polynôme minimal 
de tout élément de~$L$ est de la forme $T^{p^n}-a$ 
où $a\in K$ et où~$p$ est l'exposant caractéristique de~$K$.

\medskip
La fermeture séparable $K^s$ de~$K$ dans $\overline K$ est une extension séparable de~$K$
n'admettant aucune extension finie séparable stricte. 
Une telle extension est appelée une {\em clôture séparable de~$K$} ; 
deux clôtures séparables de~$K$ sont isomorphes. 

\deux{autfixe} Soit~$L$ une extension de~$K$ et soit~$\mathsf G$ un groupe fini d'automorphismes de~$L$ ; 
supposons que~$L^{\mathsf G}=K$. 
Soit~$x\in L$ et soit~$r$ son degré. Soit~$P$ l'élément~$\prod\limits_{y\in \mathsf G.x} (T-y)$
de~$L[T/r]$ ; il est de degré~$r^{\card(\mathsf G.x})$. 
Le polynôme~$P$ est invariant sous~$\mathsf G$
et appartient donc à~$K[T/r]$ ; il est scindé à racines simples dans~$L$
et annule~$x$ ; par conséquent,~$x$ est séparable sur~$K$, 
et~$\mathsf G.x$ engendre une sous-extension de~$L$ finie et séparable sur~$K$
 qui est stable sous $\mathsf G$. 

Ainsi,~$L$ est réunion de ses sous-extensions finies et séparables sur~$K$ stables sous~$\mathsf G$. 
Soit~$F$ l'une d'elle et soit~$n$ sa dimension sur~$K$. Le groupe~$\mathsf G$ agit 
sur $\overline K \otimes_K F\simeq \overline K^n$ ; 
il découle de l'égalité~$L^{\mathsf G}=K$ que l'ensemble des éléments 
de~$\overline K \otimes_K F\simeq \overline K^n$ 
invariants sous~$\mathsf G$ est la diagonale $\overline K\cdot(1,\ldots, 1)$. 

Par ailleurs, le groupe des automorphismes de la $\overline K$-algèbre~$\overline K^n$
s'identifie à $\got S_n$, agissant par permutation des facteurs : cela résulte (par exemple) 
du fait que les éléments de la forme $(0,\ldots,0,1,0,\ldots, 0)$ de $\overline K^n$ sont ses idempotents 
non nuls minimaux, et sont donc permutés par tout automorphisme. 
Comme la sous-algèbre de~$\overline K^n$ formée des éléments invariants sous~$\mathsf G$ 
est égale à la diagonale,~$\mathsf G$ agit transitivement sur~$\{1,\ldots, n\}$,
ce qui oblige~$n$ à être inférieur ou égal au cardinal de~$\mathsf G$ ; ainsi, $[F:K]\leq \card(\mathsf G)$. 

\medskip
Par conséquent,~$L$ est une extension finie séparable de $K$,
et l'on a la majoration $[L:K]\leq \card(\mathsf G)$. 
Comme par ailleurs~$\mathsf G$ se  plonge dans le groupe des~$K$-automorphismes de~$L$, 
lequel a un cardinal borné par~$[L:K]$, on en déduit que $[L:K]=\card(\mathsf G)$ 
et que~$\mathsf G$ est égal au groupe des~$K$-automorphismes de $L$. 

\deux{thogalgrad} Les ingrédients sont désormais réunis pour disposer d'une théorie de Galois graduée,
qui se déduit formellement de tout ce qui précède.
Nous allons en donner une formulation classique, et une plus catégorique, fidèle au point de vue de Grothendieck. 

\trois{theogalgradclass} {\em La formulation classique pour les extensions finies.} Soit~$K$ un corpoïde
et soit~$L$ une extension de~$K$. Les assertions suivantes sont équivalentes : 

\medskip
i) il existe un groupe fini~$\mathsf G$ d'automorphismes de~$L$ tel que~$K=L^{\mathsf G}$ ; 

ii) l'extension graduée~$L/K$ est finie, séparable, et normale 
(ce qui signifie que si~$x$ est un élément de~$L$, son polynôme minimal~$P$ est scindé dans~$L$) ; 

iii) il existe $r\in D$ et un élément non nul et séparable $P$ de~$K[T/r]$ 
tel que~$L$ soit un corpoïde de décomposition de~$P$ sur~$K$.

\medskip
De plus si i) est vraie alors $\mathsf G=\mathsf {Gal}(L/K)$ et~$[L:K]=\card( \mathsf G)$, et si ii) ou iii) est vraie
alors $K=L^{\mathsf{Gal}(L/K)}$. 

\medskip
Lorsqu'une extension~$L$ de~$K$ satisfait ces conditions équivalentes, on dit qu'elle est {\em finie galoisienne}. 
Si~$L$ est une extension finie galoisienne de~$K$, on dispose sur~$L$ 
de la {\em correspondance de Galois} : les 
applications~$\mathsf H\mapsto L^{\mathsf H}$ et $F\mapsto \mathsf{Gal}(L/F)$ 
établissent une bijection décroissante entre l'ensemble des sous-groupes de~$\mathsf{Gal}(L/K)$
et celui des sous-extensions de~$L/K$. 

Si~$\mathsf H$ est un  sous-groupe de~$\mathsf{Gal}(L/k)$ 
et si~$F$ est la sous-extension qui lui correspond,~$F$ 
est une extension galoisienne de~$K$ si et seulement si~$\mathsf H$ est un sous-groupe distingué
de~$\mathsf{Gal}(L/K)$, ce qui revient à demander que~$F$ soit stable 
sous $\mathsf{Gal}(L/K)$ ; si c'est le cas, la restriction des automorphismes
induit un isomorphisme~$\mathsf{Gal}(L/K)/\mathsf H\simeq \mathsf{Gal}(F/K)$. 

\trois{galinfoide} {\em La formulation classique, 
pour les extensions quelconques}. Soit~$L$ une extension
algébrique de~$K$ et soit~$\mathsf G$ le groupe~$\mathsf {Gal}(L/K)$. Soit~$x\in L$. 
Son orbite sous~$\mathsf G$ est constituée d'éléments annulés par le polynôme
minimal de~$x$, et elle est donc finie ; il s'ensuit que~$L$ est réunion
de ses sous-extensions finies stables sous~$\mathsf G$,
puis que~$\mathsf G$ est isomorphe
à la limite projective des~$\mathsf{Gal}(F/K)$ où~$F$ parcourt l'ensemble
des sous-extensions finies de~$K$. C'est donc un groupe profini. 

\medskip
On dit que~$L$ est galoisienne si elle est séparable et normale. On 
suppose à partir de maintenant 
que c'est le cas. Soient~$x_1,\ldots,x_r$ des éléments de~$L$, soit~$\sch P$
{\em l'ensemble} (et non la famille)
des polynômes minimaux des~$x_i$ et soit~$P$ le produit
des éléments de~$\sch P$. Il résulte de notre hypothèse sur~$L$ que~$P$
est scindé à racines simples dans~$L$ ; la sous-extension de~$L$
engendrée par les racines de~$P$ est par construction finie galoisienne, et
stable sous l'action de~$\mathsf G$.

\medskip
Les faits suivants se démontrent par passage à la limite à 
partir du cas des extensions finies galoisiennes. 
Les 
applications~$\mathsf H\mapsto L^{\mathsf H}$ et $F\mapsto \mathsf{Gal}(L/F)$ 
établissent une bijection décroissante entre l'ensemble des sous-groupes 
{\em fermés} de~$\mathsf{Gal}(L/K)$
et celui des sous-extensions de~$L/K$. 

Si~$\mathsf H$ est un  sous-groupe 
fermé de~$\mathsf{Gal}(L/k)$ 
et si~$F$ est la sous-extension qui lui correspond,~$F$ 
est une extension galoisienne de~$K$ si et seulement si~$\mathsf H$ est un sous-groupe distingué
de~$\mathsf{Gal}(L/K)$, ce qui revient à demander que~$F$ soit stable 
sous $\mathsf{Gal}(L/K)$ ; si c'est le cas, la restriction des automorphismes
induit un isomorphisme~$\mathsf{Gal}(L/K)/\mathsf H\simeq \mathsf{Gal}(F/K)$.

\trois{theogalcat} {\em La formulation catégorique de
la théorie de Galois.} Soit~$\mathsf G$ le
groupe~$\mathsf {Gal}(K^s/K)$ ; notons qu'il coïncide
avec~$\mathsf{Gal}\; (\overline K/K)$
(cela provient du fait qu'une extension radicielle n'a pas d'automorphismes non triviaux). 

\medskip
Le foncteur~$A\mapsto \mathsf{Hom}_K(A,K^s)$ établit une équivalence entre la catégorie 
des $K$-algèbres étales et celle des~$\mathsf G$-ensembles discrets finis (un~$\mathsf G$-ensemble
 discret est un ensemble muni d'une action 
 continue de~$\mathsf G$ pour laquelle les stabilisateurs sont ouverts).

\deux{remelprem}  {\em Remarque.} Ainsi, la théorie de Galois classique survit 
à peu près mot pour mot dans le contexte gradué.
Mentionnons toutefois un théorème qui n'est plus vrai : celui de l'élément primitif  ; 
nous allons en effet construire un contre-exemple à celui-ci. 

\medskip
Soit $F$ un corps, et soit~$K$ le corpoïde
tel que~$K_1=F$ et~$K_r=\{0_r\}$ pour tout~$r\neq 1$.
Soit~$(r,s)$ une famille libre du $\QQ$-espace vectoriel $D$,
et soit~$L$ le corpoïde~$K[T/r,r T^{-1}, S/s, s S^{-1}]$.
Soit~$n$ un entier strictement supérieur à $1$ et inversible dans~$F$,
et soit $L'$ la~$L$-algèbre~$L[U/r^{1/n}, V/s^{1/n}]/(U^n-T, V^n-S)$.

On vérifie aisément 
que~$L'$ est une extension graduée finie et séparable de degré~$n^2$
 de~$L$. Soit $x\in L'$  ; 
 il est de la forme~$aU^mV^{m'}$ où $a\in L$ ; par conséquent,~$x^n\in L$
 et~$x$ ne peut engendrer~$L'$ qui est de degré~$n^2$ sur~$L$.

 \deux{descgaloideclassetgp} {\bf Lemme.} {\em Soit~$K$ un corpoïde
 et soit~$L$ une extension normale de~$K$. 
 
 \medskip
 i) Le corps~$L_1$
 est alors une extension normale de~$K_1$,
 qui est galoisienne si~$L$ est galoisienne.
 
 ii) La restriction~$\mathsf{Gal}(L/K)\to \mathsf{Gal}(L_1/K_1)$ est surjective. 
 
 iii) Si~$g\in \mathsf{Gal}(L/K)$ agit trivialement sur~$L_1$ alors 
 pour tout~$x\in L\ti$, l'élément~$g(x)/x$ de~$L_1\ti$
 ne dépend que de~$\deg x$, et l'on construit par ce biais un 
 isomorphisme
 $$\mathsf{Ker}(\mathsf{Gal}(L/K)\to \mathsf{Gal}(L_1/K_1))\simeq \mathsf{Hom}(\deg(L\ti)/\deg(K\ti),L_1\ti).$$}
 
 \medskip
 {\em Démonstration.} Supposons~$L$ est normale (resp. normale
et séparable sur~$K$), et soit~$x\in L_1$. Le polynôme minimal
de~$x$ appartient
par définition à~$K_1[T]$. En vertu de notre hypothèse sur~$L$,
il est scindé (resp. scindé à racines simples) dans~$L$. Comme
ses racines
appartiennent nécessairement à~$L_1$, il est scindé 
(resp. scindé à racines simples) dans~$L_1$, et~$L_1$
est une extension normale (resp. normale et séparable) de~$K_1$,
ce qui achève de prouver~i).

\medskip
Le sous-corpoïde~$L_1K$ de~$L$ est égal
à~$\coprod_{r\in D} L_1K_r$. Il en résulte que tout~$K$-automorphisme
de~$L$ stabilise~$L_1K$, et que tout~$K_1$-automorphisme
de~$L_1$ s'étend d'une unique manière en un~$K$-automorphisme de~$L_1K$. 

Autrement dit,~$L_1K$ est une sous-extension galoisienne de~$L$, 
et~$\mathsf{Gal}(L_1K/K)$ s'identifie {\em via} la restriction
à~$\mathsf{Gal}(L_1/K_1)$. 
L'assertion~ii)
découle dès lors de la surjectivité de~$\mathsf{Gal}(L/K)\to \mathsf{Gal}(L_1K/K)$.

Nous allons maintenant   prouver~iii). 
Pour tout~$d\in \deg(L\ti)/\deg(K\ti)$, posons~$$\Lambda_d=\coprod_{r\in \deg(L\ti), r=d\;{\rm mod} \;\deg(K\ti)}L_r.$$

\medskip
Pour tout~$d\in \deg(L\ti)/\deg(K\ti)$, le groupe~$(L_1K)\ti$ agit
sur~$\Lambda_d$, et l'ensemble des éléments non nuls de~$\Lambda_d$ est un
torseur sous~$(L_1K)\ti$. Remarquons par ailleurs
que~$L^\times \subset \coprod \Lambda_d$. 

\medskip
Soit $\psi$ un morphisme de groupes de~$\deg(L\ti)/\deg(K\ti)$ 
dans~$L_1\ti$. Il existe un unique automorphisme~$g_\psi$ de~$L$
qui agit pour tout~$d$ sur~$\Lambda_d$ par multiplication par~$\psi(d)$ ; 
l'action de~$g_\psi$ sur~$\Lambda_1$ étant triviale, l'automorphisme~$g_\psi$
appartient à~$\mathsf {Gal}(L/L_1K)$.

La flèche~$\psi \mapsto g_\psi$ est visiblement
un morphisme de groupes injectif
de~ $\mathsf{Hom}(\deg(L\ti)/\deg(K\ti),L_1\ti)$
vers~$\mathsf{Gal}(L/L_1K)$. 
Il suffit pour conclure de montrer que ce
morphisme est surjectif. 

Soit $\phi\in \mathsf {Gal}(L/L_1K)$ et soit $d\in \deg(L\ti)/\deg(K\ti)$. 
Si $x$ et $y$ sont deux éléments non nuls de $\Lambda_d$, leur quotient appartient à $L_1K$ et l'on a
donc~$\phi(x)/\phi(y)=x/y$ ; ceci implique l'existence d'un élément~$\psi(d)$ de $(L_1K)\ti$
tel que $\phi(x)=\psi(d)x$ pour tout $x\in \Lambda_d$ ; 
cette égalité appliquée à n'importe quel élément non nul~$x$
de~$\Lambda_d$ force le degré de $\psi(d)$
à être égal à $1$, et l'on a donc~$\psi(d)\in L_1\ti$. 
Comme~$\phi$ commute à la multiplication, on a
nécessairement~$\psi\in\mathsf{Hom}(\deg(L\ti)/\deg(K\ti),L_1\ti)$ ; 
par construction,~$\phi=g_\psi$.~$\Box$ 

\deux{remlunstarmu} On conserve les hypothèses et notations
du lemme~\ref{descgaloideclassetgp}
ci-dessus. 

\trois{actionout} Le noyau de la
surjection~$\mathsf{Gal}(L/K)\to \mathsf{Gal}(L_1/K_1)$
étant abélien, l'action extérieure de~$\mathsf{Gal}(L_1/K_1)$
sur celui-ci est une vraie action ; il résulte alors de sa construction
que l'isomorphisme de~iii) est~$\mathsf{Gal}(L_1/K_1)$-équivariant. 

\trois{lunstarmu} Comme $L$
est une extension algébrique de $K$, le
groupe $\deg(L^\times)/\deg(K^\times)$ est de torsion. 
Il s'ensuit que
$$\mathsf{Hom}(\deg(L^\times)/\deg(K^\times),L_1^\times)
=\mathsf{Hom}(\deg(L^\times)/\deg(K^\times),\mu(L_1)),$$
où~$\mu$ désigne le foncteur des racines de l'unité. 

\trois{assezracine} Soit~$p$ l'exposant
caractéristique de~$K$. Désignons 
par~$E$ le groupe quotient~$\deg(L^\times)/\deg(K^\times)$. 
Pour tout
entier~$n>0$, on note~$n_{p'}$ le plus
grand diviseur premier à~$p$ de~$n$ ; si~$A$
est un groupe abélien de torsion, on note~$A_{p'}$
sa partie de torsion première à~$p$. 

\medskip
Supposons que~$L$ est une extension finie de~$K$. 
On pose~$c=\card (\mathsf{Gal}(L/K))$
et~$c_1=\card (\mathsf{Gal}(L_1/K_1))$, et l'on 
note~$d$
le cardinal du
groupe~$\mathsf{Hom}
(E,\mu(L_1))$. Le groupe
multiplicatif de~$L_1$ étant sans~$p$-torsion,
l'entier~$d$ est premier à~$p$,
et~$\mathsf{Hom}(E,\mu(L_1))=
\mathsf{Hom}(E_{p'},\mu(L_1))$. 
Comme
l'extension~$L$ est une
extension normale
de~$K$, 
on a
~$[L:K]_{p'}=c_{p'}$ ; et puisque~$L_1$ 
est une extension normale de~$K_1$
en vertu du lemme~\ref{descgaloideclassetgp}, 
il vient aussi~$c_{1,p'}=[L_1:K_1]_{p'}$. On déduit
par ailleurs 
du même lemme et de~\ref{lunstarmu}
que~$c=c_1d$, puis que~$c_{p'}=c_{1,p'}d$. 
On a également~$[L:K]=[L_1:K_1]\cdot \card E$ (\ref{interpdimoide}). 
En conséquence,
$$[L:K]_{p'}=[L_1:K_1]_{p'} \cdot\card E_{p'}=c_{1,p'} \cdot \card E_{p'}.$$
Comme on a d'autre part~$[L:K]_{p'}=c_{p'}=c_{1,p'}d$,
il vient~$d=\card E_{p'}$.

Autrement dit,
on obtient {\em tous}
les caractères de~$E_{p'}$ 
par plongement dans~$L_1^\times$. 
Cela signifie que si~$a$ désigne l'exposant
de~$E_{p'}$,
le polynôme~$T^a-1$ est scindé
dans~$L_1$. Tout plongement
de~$\mu(L_1)$ dans~$\QQ/\ZZ$ induit un
isomorphisme entre~$\mathsf{Hom}(E,\mu(L_1))=
\mathsf{Hom}(E_{p'},\mu(L_1))$
et le dual~$E_{p'}^\vee$ de~$E_{p'}$
(qui est lui-même isomorphe, non canoniquement
sauf cas triviaux,  à~$E_{p'}$).

\trois{assezracinepasfini} On ne suppose plus que~$L$ est finie sur~$K$. 
Par passage à la limite, on déduit
 les faits suivant du~\ref{assezracine} ci-dessus : 

\medskip
$\bullet$ si~$a$ divise l'exposant généralisé de~$E_{p'}$,
le polynôme~$T^a-1$ est scindé dans~$L_1$ ; 

$\bullet$ tout plongement de~$\mu(L_1)$ dans~$\QQ/\ZZ$
induit un isomorphisme entre~$\mathsf{Hom}(E,\mu(L_1))$
et le dual de Pontrjagin~$E_{p'}^\vee$ de~$E_{p'}$. 

\subsection*{Valuations}

\deux{dominvaloide} Soit~$K$ un 
corpoïde, et soient~$A$ et $B$ deux
sous-annéloïdes locaux de~$K$. On dit que~$B$ 
{\em domine}~$A$ si~$A\subset B$ et si l'idéal maximal
de~$A$ est contenu dans celui de~$B$. 

\medskip
Soit~$A$ un sous-annéloïde~$K$. Nous allons
montrer qu'il est maximal
pour la relation de domination si et seulement si il possède la propriété
suivante : {\em pour tout élément non nul~$x$ de~$K$, 
on a $x\in A$ ou~$x\inv \in A$}
(notons
qu'elle entraîne que~$K$ est le corpoïde des fractions
de~$A$). 

\trois{dommaxfacile} Supposons que pour tout élément non nul~$x$ de~$K$, 
on a $x\in A$ ou~$x\inv \in A$, et soit~$B$ un sous-annéloïde local
de~$K$ dominant~$A$. Soit~$x\in B$. Supposons
que~$x\notin A$. On a alors par hypothèse~$x\inv\in A$. 

Comme~$x\notin A$, l'élément~$x\inv$ de~$A$
est non inversible. Il appartient
donc à l'idéal maximal de~$A$, et partant à
celui de~$B$ ; mais cela contredit l'appartenance
de~$x=(x\inv)\inv$ à~$B$.
Ainsi,~$x\in A$. Il s'ensuit que~$B\subset A$, et donc que~$B=A$. En
conséquence,~$A$ est maximal
pour la relation de domination. 

 \trois{dommaxdifficile} Supposons que~$A$ est maximal pour la relation
 de domination. Soit~$x\in K$, soit~$r$ son degré,
 et soit~$\got m$ 
 l'idéal maximal de~$A$. Nous allons montrer que~$x\in A$, 
 ou que~$x\neq 0$ et~$x\inv \in A$. 
 
 \medskip
 {\em Supposons que~$A[x]/\got m A[x]$ est non nul}. Il possède
 alors un idéal maximal, dont l'image réciproque dans~$A[x]$
 est un idéal premier~$\got n$ de ce dernier, contenant~$\got m$. 
 Le localisé~$A[x]_{\got n}$ est un sous-annéloïde local
 de~$K$ dominant~$A$ ; par maximalité, il est égal à~$A$ et l'on 
 a donc~$x\in A$. 
 
 \medskip
 {\em Supposons que~$A[x]/\got m A[x]=0$}. On a
 en particulier~$1\in \got mA[x]$ ,
 ce qui veut dire que l'on peut
 écrire~$1=a_0+\sum\limits_{i\geq 1} a_i x^i$ où
 les~$a_i$ appartiennent à~$\got m$, et où
 chaque~$a_i$ est de degré~$r^{-i}$. 
 Comme~$a_0\in \got m$, l'élément~$1-a_0$
 de~$A$ est inversible ; en faisant passer~$a_0$
 à gauche puis en multipliant par~$(1-a_0)\inv$
 les deux membres de l'égalité,
 on se ramène
 au cas où~$a_0=0$. 
 
 Comme~$1=\sum \limits_{i\geq 1}a_i x^i$ on a~$x\neq 0$. Si~$N$
 est le plus grand exposant de~$x$ apparaissant dans le membre de droite, 
 on voit en divisant l'égalité par~$x^N$ que~$1/x$ est entier
sur~$A$. Le {\em going-up}
 assure alors que~$A[1/x]$ possède un idéal maximal~$\got n$ dont l'intersection
 avec~$A$ est égale à~$\got m$. Le localisé~$A[1/x]_{\got n}$ est un sous-annéloïde local
 de~$K$ dominant~$A$, c'est donc~$A$ lui-même et il s'ensuit que~$(1/x)\in A$.

\deux{valgrad} Soit~$K$ un corpoïde.  Soit~$G$
un groupe abélien totalement ordonné, noté multiplicativement. On désigne
par~$G_0$ le monoïde ordonné obtenu en adjoignant formellement à~$G$
un plus petit élément absorbant~$0$.
Une {\em valuation sur~$K$ à valeurs dans~$G_0$}
est une application~$\abs.$ de~$K$ vers~$G_0$,
où~$0$ est un 
\medskip

$\bullet$~$\abs 1=1$,~$\abs {0_r}=0$ pour tout~$r$,
 et~$\abs {ab}=\abs a \cdot \abs b$ pour tout~$(a,b)$  (ce qui implique que~$\abs a \neq 0$ dès que~$a\neq 0$) ; 
 
~$\bullet$  pour tout~$r\in D$ et tout couple~$(a,b)$ d'éléments de~$K_r$ on a l'inégalité
ultramétrique~$\abs {a+b}\leq \max(\abs a, \abs b)$. 
 
 \medskip
Deux valuations~$\abs .: K\to G_0$ et~$\abs .' :K\to G'_0$ 
sont dites
{\em équivalentes} s'il existe une valuation
~$\abs .'': K\to G''_0$ et deux morphismes strictement
croissants~$G''\hookrightarrow G$ et~$G''\hookrightarrow G'$
 tels que le diagramme
 $$\diagram &&G_0\\K\urrto^{\abs.}\drrto_{\abs .'}\rrto^{\abs.''}&&G''_0\uto \dto&\\
 &&G'_0\enddiagram$$ commute.

\trois{defannvaloide} Soit~$\abs .~$ une valuation sur~$K$. 
L'ensemble des éléments ~$x$ de~$K$ tels que~$\abs x\leq 1$ est un
sous-annéloïde~$\sch O_{\abs .}$ de~$K$, qui est appelé {\em l'annéloïde de~$\abs.$.}
L'annéloïde~$\sch O_{\abs .}$ est un annéloïde local ; son idéal maximal
est l'ensemble des éléments~$x$ de~$K$ tels que~$\abs x <1$. 

Il résulte immédiatement
de la définition d'une valuation que
si~$x$ est un élément non nul de~$K$ 
alors~$x\in \sch O_{\abs .}$ ou~$x\inv\in \sch O_{\abs .}$. 

Réciproquement, soit~$A$ un sous-annéloïde
local de~$K$ tel que l'on
ait~$x\in A$ ou~$x\inv \in A$ pour
tout~$x\in K\ti$. Le groupe quotient~$G:=K\ti/A\ti$
possède alors un ordre dont les éléments~$\leq 1$
sont exactement
les classes d'éléments non nuls de~$A$, et la flèche~$K\to G_0$ 
obtenue en prolongeant la flèche quotient
$K\ti\to G$
est une valuation d'annéloïde~$A$. 

\medskip
Il est immédiat que ces constructions établissent
une bijection entre l'ensemble des classes d'équivalence
de valuations sur~$K$ et celui des sous-annéloïdes locaux~$A$
de~$K$ tels que~$x\in A$ ou~$x\inv \in A$ pour
tout~$x\in k\ti$, c'est-à-dire encore celui des sous-annéloïdes
locaux de~$K$  maximaux pour la relation de domination
({\em cf.}~\ref{dominvaloide} {\em et sq.}).

\trois{defvaltrivoide} La valuation {\em triviale} sur~$K$ est celle dont l'annéloïde
est~$K$ tout entier ; elle envoie tout élément non nul sur~$1$.

\trois{extvaloide} Soit~$K\hookrightarrow F$ une extension de corpoïdes
et soit~$\abs .$ une valuation sur~$K$.  Il résulte
du lemme de Zorn que tout sous-annéloïde local
de~$F$ est dominé par un sous-annéloïde local
maximal pour la relation de domination, donc par l'annéloïde 
d'une valuation. C'est en particulier le cas
de l'annéloïde de~$\abs .$ ; en conséquence, 
la valuation~$\abs.$ admet un prolongement
à~$F$.

\trois{exgauss} Soit~$K$ un corpoïde valué et soit~$G$ un groupe ordonné
contenant~$\abs {K\ti}$. Soit~${\bf d}=(d_i)_{i\in I}$ une famille d'éléments
de~$D$, et soit~${\bf r}=(r_i)_{i\in I}$ une famille d'éléments de~$G$ (indexée
par le même ensemble que~$\bf d$). La formule

$$\sum a_I {\bf T}^I\mapsto \max |a_I|\cdot {\bf r}^I$$ définit une valuation 
(dite {\em de Gauß}) sur~$K({\bf d}\inv{\bf T})$ ; on la notera~$\eta_{K,{\bf d},{\bf r}}$, 
ou parfois simplement~$\eta_{K{,\bf r}}$ ou même~$\eta_{\bf r}$ si le contexte est suffisamment
clair.

\trois{extvaloidealg} Soit~$K\hookrightarrow F$ une extension
de corpoïdes
et soit~$\abs .$ une valuation sur~$F$.
Soit~$t$ un élément 
de~$F\ti$ algébrique sur~$K$, et 
soit~$P=\sum a_iT^i$ son polynôme
minimal. L'égalité~$\sum a_it^i=0$ entraîne l'existence
de deux entiers~$i$ et~$j$ distincts tels que~$\abs {a_i t^i}=\abs{a_jt^j}$ ; il s'ensuit qu'il 
existe un entier~$n>0$, majoré par le degré de~$t$ sur~$K$, tel
que~$\abs t^n\in \abs {K\ti}$. Il s'ensuit aussi que si la restriction
de~$\abs .$ à~$K$ est triviale, alors
~$\abs .$ est elle-même triviale. 

\trois{extvaloidettrans} Soit~$K\hookrightarrow F$ une extension
de corpoïdes et soit~$\abs .$ une valuation sur~$K$.
Soit~$G$ un groupe ordonné contenant~$\abs{K\ti}$, soit~$(t_i)$ une 
famille d'éléments de~$F$ algébriquement indépendants sur~$K$
et soit~$(r_i)$ une famille d'éléments de~$G$. Il existe alors une valuation
sur~$F$ à valeurs dans~$G_0^\QQ$, prolongeant~$\abs .$ et envoyant~$t_i$
sur~$r_i$ pour tout~$i$. 

\medskip
En effet, quitte à agrandir les familles~$(t_i)$ et~$(r_i)$, on peut supposer
que~$(t_i)$ est une famille algébriquement indépendante maximale
de~$F$ sur~$K$, et donc que~$F$ est algébrique sur son
sous-corpoïde~$K(t_i)_i$. On choisit 
alors un prolongement quelconque à~$F$ de la
valuation de Gauß~$\eta_{\bf r}$ sur~$K(t_i)_{i\in I}$ ; il résulte
du~\ref{extvaloidealg} ci-dessus que ledit prolongement
est à valeurs dans~$G_0^\QQ$. 

\medskip
En particulier,~$\abs .$ admet un prolongement à~$F$ à valeurs dans~$\abs K^\QQ$ : il suffit
d'appliquer ce qui précède en prenant pour~$(t_i)$ et~$(r_i)$ les familles {\em vides}.

\deux{zr} Soit~$F$ un corpoïde. 
On note~$\PP_F$ l'espace de Riemann-Zariski de~$F$, 
c'est-à-dire l'ensemble des classes d'équivalence de 
valuations de~$K$. Pour toute partie~$X$ de~$\PP_F$ et 
tout sous-ensemble~$A$ de~$F$, on note~$X\{A\}$ (resp.~$X\{\{A\}\}$) 
le sous-ensemble de~$X$ formé des valuations~$\abs .$ dont 
l'annéloïde (resp. l'idéal maximal) contient~$A$ ou, si l'on 
préfère, telles que~$\abs a\leq 1$ (resp.~$\abs a <1$) pour tout~$a\in E$.
On munit~$\PP_K$ de la topologie engendrée 
par les parties de la forme~$\PP_K\{A\}$ où~$A$ est {\em fini}.

Soient~$A$ et~$B$ deux sous-ensembles 
de~$F$ ; 
posons~$X=\PP_F\{A\}\{\{B\}\}$, et munissons-le de la topologie induite. 
Il est quasi-compact : cela résulte du fait qu'il est fermé dans~$\PP_F\{A\}$, et que ce dernier 
est quasi-compact. On qualifiera d'{\em affine} tout ouvert de~$X$ de la forme~$X\{C\}$ où~$C$ 
est un sous-ensemble fini de~$K$ ; remarquons que~$X$ est lui-même affine (prendre~$C=\emptyset$).

\medskip
Soit~$K$ un sous-corpoïde
de~$F$ et soit~$\abs .$ une valuation  de~$K$. Soit~$A$ (resp.~$\got m$) 
 l'annéloïde de~$\abs .$ (resp. l'idéal maximal de~$\abs .$). 
 Le sous-ensemble~$\PP_F\{A\}\{\{\got m\}\}$ de~$\PP_K$ est alors l'ensemble des valuations qui 
 prolongent~$\abs .$. On le notera le plus souvent~$\PP_{F/(K,\abs .)}$,  ou simplement~$\PP_{F/K}$
 s'il n'y a pas d'ambiguïté sur la valuation de~$K$ ; lorsque celle-ci n'est pas précisée,~$K$ 
 sera implicitement considéré comme muni de la valuation triviale.

\deux{caspartvalclass} {\bf Valuations classiques.} Soit~$K$ un corps. On peut le voir
comme un~$\{1\}$-corpoïde, et dès lors définir une valuation sur~$K$
au sens de la section précédente. Il n'y a pas de conflit de terminologie : une telle
valuation est précisément une valuation (de Krull) usuelle
sur~$K$. Les espaces~$\PP_K$ et~$\PP_{K/(K_0,\abs .)}$ (où~$K_0$ est
un sous-corps de~$K$ muni d'une valuation de Krull~$\abs .$) 
sont eux-mêmes les espaces de Riemann-Zariski usuels ; si~$K$ est une
 extension de type fini de~$K_0$, l'espace topologique~$\PP_{K/(K_0,\abs .)}$ s'identifie 
 à la limite projective des fibres spéciales de tous les modèles intègres, 
 propres et plats du corps~$K$ sur~$\spec \sch O_{\abs .}$. 

\subsection*{Réduction à la Temkin} 

On désigne toujours par~$D$ un groupe abélien divisible ; on fixe par ailleurs un groupe abélien divisible {\em ordonné}~$G$. Dans ce qui suit, et sauf mention expresse du contraire, 
les annéloïdes seront des~$D$-annéloïdes, et les valuations des valuations à valeurs dans~$G_0$.

\deux{ineg} Soit~$A$ un
annéloïde. Une {\em semi-norme $||.||$} sur~$A$ est 
une application
de~$A$ dans~$G_0$ qui satisfait les axiomes suivants : 

$\bullet$ $||0||=0$ et~$||ab||\leq ||a||\cdot ||b||$ pour tout~$(a,b)\in A^2$ ; 

$\bullet$ $||a+b||\leq \max(||a||,||b||)$ pour tout~$r\in D$ et tout~$(a,b)\in A_r^2$.

\trois{exseminorme} {\em Exemples.} Si~$K$ est 
un corpoïde, toute valuation sur~$K$ est 
une semi-norme, et même une norme ({\em i.e} son noyau est trivial). 
Plus généralement,
donnons-nous une famille finie~$(K_i)$
de corpoïdes, chacun des~$K_i$ étant muni d'une valuation. 
L'application~$(x_1,\ldots, x_n)\mapsto \max |x_i|$ est alors 
une semi-norme sur~$\prod K_i$. 

\deux{defredgradoide} Soit~$A$ un annéloïde muni d'une semi-norme~$||.||$. 
Pour tout $(r,g)$
appartenant à~$D\times G$, on note $A_{r,\leq g}$ (resp. $A_{r,<g}$) le sous-groupe 
additif de $A_r$ formé des éléments $a$ tels que $||a||\leq g$ (resp. $||a||<g$). 
Nous désignerons par $\red A$ la réduction de~$A$
au sens de Temkin, c'est-à-dire
le~$(D\times G$)-annéloïde~$\coprod\limits_{(r,g)\in D\times G}A_{r,\leq g}/A_{r,<g}$. 
Si $a\in A$ et si $g$ est un élément de~$G$ tel que $||a||\leq g$, on notera $\red{a}_{g}$ l'image de $a$
dans $\red{A}_{\deg a, g}$ ; si $g=||a||$, on écrira simplement $\red{a}$. 


\medskip
La formation de~$\red A$ est fonctorielle en~$A$
pour les
morphismes contractants d'annéloïdes semi-normés ; notons
qu'un morphisme isométrique d'annéloïdes semi-normés induit
une injection entre les annéloïdes résiduels correspondants. 

\deux{cascorspval} Soit~$K$ un corpoïde muni d'une valuation~$\abs .$, soit~$\sch O_{\abs .}$ 
l'annéloïde correspondant, et soit~$\got m_{\abs .}$ son idéal maximal. 

\trois{corpsresid} Le~$(D\times G)$-annéloïde~$\red K$ est alors un corpoïde, appelé {\em corpoïde résiduel} de $K$ ; le corpoïde quotient~$\sch O_{\abs .}/\got m_{\abs .}$ s'identifie à~$\red K_{D\times \{1\}}$. 
Le groupe~$\deg(\red K\ti)$ est égal à~$(\deg,\abs .)(K\ti)$.

\trois{baseresoide} Soit~$A$ une~$K$-algèbre munie d'une semi-norme~$||.||$ compatible
avec la valuation de~$K$. Soit $(a_{i})$ une famille d'éléments 
de $A$ telle que ~$||a_i||\neq 0$ pour tout~$i$. Les propositions suivantes sont équivalentes :

\medskip
i) on a $||\sum \lambda_{i}a_{i}||=\max |\lambda_{i}|\cdot ||a_{i}||$ pour toute famille 
$(\lambda_{i})$ d'éléments de $K$
presque tous nuls et tels que les~$\lambda_i a_i$ aient tous même degré (dans la suite, nous dirons simplement qu'une telle famille~$(\lambda_i)$ est {\em adaptée} à~$(a_i)$). 

ii) $(\red{a_{i}})$ est une famille libre du $\red K$-espace vectoriel $\red A$. 

\medskip
Pour le voir, posons~$g_i=||a_i||$ pour tout~$i$. La propriété~i) est {\em fausse} si et seulement si il
existe~$h\in G$, et une famille~$(\lambda_i)$ de scalaires adaptée à~$(a_i)$
telle
que~$|\lambda_i|\leq hg_i^{-1}$ pour tout~$i$, 
telle que~$|\lambda_i|=hg_i^{-1}$ pour au moins un~$i$, et telle que~$||\sum \lambda_i a_i||<h$. Mais
cela revient à demander qu'il existe
une famille~$(\mu_i)$ d'éléments de~$\red K$ 
adaptée à~$(\red{a_i})$ telle que l'un au moins des~$\mu_i$ soit  non nul, et 
telle que l'élément~$\sum \mu_i \red{a_i}$ de~$\red A$ soit nul (prendre~$\mu_i=\red {\lambda_i}_{hg_i^{-1}}$). 
Autrement dit, cela revient à nier que~$(\red{a_i})$ soit libre sur~$\red K$. 
Ainsi, i) est fausse si et seulement si~ii) est fausse, ce qu'il fallait démontrer.

\trois{dimresoide} Soit~$(x_i)$ une base de~$\red A$ sur~$\red K$. 
Chaque~$x_i$
est alors non nul, donc de la forme~$\red{a_i}$ pour un certain~$a_i$ de norme
non nulle dans~$A$. Comme~$(a_i)$ satisfait ii) par construction,
elle satisfait i) et est dès lors libre.
Il s'ensuit que~$\dim {\red K}\red A\leq \dim KA$. 

\trois{algindepetgauss} Soit~$(t_i)$ une famille d'éléments de~$A$. Pour tout~$i$, on pose~$r_i=||t_i||$
et l'on suppose que~$r_i\neq 0$. Les assertions suivantes sont alors équivalentes :

\medskip
a) on a~$\left\|\sum \lambda_I {\bf t}^I\right\|=\max |\lambda_I|\cdot {\bf r}^I$ pour toute famille~$(\lambda_I)$ adaptée
à~$({\bf t})^I$ ; 

b) la famille~$(\red{\bf t})$ est algébriquement indépendante sur~$\red K$.

\medskip
En effet, supposons que~a) soit vérifiée. On en déduit pour commencer que~$||{\bf t}^I||={\bf r}^I$
pour tout~$I$, et donc que~$||\sum \lambda_I {\bf t}^I||=\max |\lambda_I|\cdot ||{\bf t}^I||$ pour toute
famille~$(\lambda_I)$ adaptée
à~${\bf t}^I$. Il découle alors de~\ref{baseresoide} que~$(\red {\bf t}^I)$ est une famille libre sur~$\red K$ ; autrement dit,~b) est vérifiée.

\medskip
Réciproquement, supposons que~b) soit vérifiée. On en déduit pour commencer que~$\red {\bf t}^I$ est non nul
pour tout~$I$, c'est-à-dire que~$||{\bf t}^I||={\bf r}^I$ pour tout~$I$. 

Comme~b) est vérifiée,~$(\red {\bf t}^I)$ est libre sur~$\red K$; on déduit alors de~\ref{baseresoide} et de ce qui précède
que~$||\sum \lambda_I {\bf t}^I||=\max |\lambda_I|\cdot ||{\bf t}^I||=\max |\lambda_I|\cdot {\bf r}^I$ pour toute famille~$(\lambda_I)$ adaptée
à~$({\bf t}^I)$, d'où~a).

\trois{algindepgaussuite} Remarquons que la propriété~a) peut se reformuler en demandant que les~$t_i$ soient 
algébriquement indépendants sur~$K$ et que~$||.||_{|K[{\bf t}]}$ soit égale à (la restriction de) la valuation
de Gauß~$\eta_{K, \deg({\bf t}), {\bf r}}$. 

Supposons que~a), et partant~b), soient satisfaites. Soit~$f=\sum \lambda_I {\bf t}^I$ un élément
non nul de~$K[{\bf t}]$. Soit~$\sch J$ l'ensemble des indices~$I$ tels que~$||f||=|\lambda_I|{\bf r}^I$. 
Posons~$g=\sum_{I\in \sch J} \lambda_I {\bf t}^I$. On a alors en vertu de~a) 
l'inégalité~$||f-g||<||f||$ ; 
en conséquence,~$\red f=\red g=\sum_{I\in \sch J}
 \red {\lambda_I}\red{{\bf t}^I}$. Ainsi,~$\red{K[{\bf t}]}$ est 
égal à~$\red K[\red {\bf t}]$. Comme~$\red {\bf t}$ est une famille algébriquement indépendante
sur~$\red K$, la~$\red K$-algèbre~$\red K[\red {\bf t}]$ s'identifie à
l'algèbre de polynômes~$\red K[(\deg ({\bf t}), {\bf r})\inv {\bf T}]$. 

\trois{exprecouronne} Supposons maintenant donnés deux éléments~$r$ et~$s$
 de~$G$ avec~$r\leq s$, et un élément~$t$ inversible
de~$A$ tel que~$\|\sum a_i t ^i||=\max \{|a_i|r^i, |a_i|s^i\}_i$ pour toute famille~$(a_i)_{i\in \ZZ}$
adaptée à~$(t^i)$. On a alors~$||t|^i ||=s^i$ et~$||t^{-i}||=r^{-i}$ pour
tout~$i\geq 0$,  et~$\|\sum a_i t ^i||=\max |a_i|\cdot ||t^i||$
pour toute famille~$(a_i)$ adaptée à~$(t^i)$.

Posons~$\tau=\red t $ et~$\sigma=\red {t\inv}$ ; on a
par ce qui précède~$\red {t ^i}=\tau ^i$ et~$\red {t^{-i}}=\sigma^i$ 
pour tout~$i\geq 0$. 

\medskip
On déduit alors de~\ref{baseresoide} que la famille~$\bf b$
formée des~$\sigma^i$ pour~$i<0$,
de~$\{1\}$ 
et des~$\tau^i$ pour~$i> 0$ est libre sur~$\red K$. De plus, soit~$f=\sum a_i t^i$ un élément
non nul de~$K[t,t\inv]$, et soit~$\sch J$ l'ensemble des indices~$i$ tels que~$||f||=|	a_i|\cdot ||t^i||$. 
Posons~$g=\sum_{i\in \sch J} a_i t^i$. On a alors en vertu de~a) 
l'inégalité~$||f-g||<||f||$ ; 
en conséquence,
$$\red f=\red g=\sum_{i \in \sch J, i<0}
 \red {a_i}\sigma^{-i}+\sum_{i\in \sch J, i\geq 0} \red{a_i} \tau^i.$$ Ainsi,~$\bf b$
 est une base de~$\red{K[t,t\inv]}$ sur~$\red K$. 

\medskip
Le produit~$\tau \sigma$ est égal par définition à~$\red{tt\inv}_{1, s/r}=\red 1_{1, s/
r}$ ; 
il est donc égal à~$1$ si~$r=s$, et à~$0_{1, s/r}$ sinon. Si l'on pose~$g=(s\inv ,\deg t\inv)$ 
et~$h=(r,\deg t)$ on dispose donc d'un morphisme naturel de~$\red K[g\inv S, h\inv T]/(TS-\epsilon)$ 
vers~$\red{K[t,t\inv]}$ envoyant~$T$ sur~$\tau$ et~$S$ sur~$\sigma$, où~$\epsilon$
est égal à~$1$ si~$r=s$ et à~$0_{1, r/s}$ sinon. Comme~$\bf b$ est une base de~$\red{K[t,t\inv]}$, 
ce morphisme est un isomorphisme.   

\deux{descextgrad} Soit~$K\hookrightarrow L$ une extension de corpoïdes valués.
Lorsqu'on
parlera de {\em sous-extension} de~$\red L$, c'est toujours à sa structure
naturelle d'extension de~$\red K$ que l'expression fera référence.

\trois{algalgtilde} Si~$L$ est une extension finie
de~$K$ alors~$\red L$ 
est une extension finie de~$K$ en vertu de~\ref{dimresoide}. Il s'ensuit, 
par un argument de passage à la limite inductive, que
si~$L$ est algébrique sur~$K$, alors~$\red L$ est algébrique sur~$\red K$.

\trois{descresgrad} Posons~$\delta =(\deg,\abs.) : L\ti \to D\times G$. Le groupe quotient

$$\deg (\red L\ti)/\left[(\deg (\red K\ti)\cdot ((D\times\{1\})\cap \deg(\red L\ti))\right] $$
est égal à
$$\delta(L\ti)/\left[ \delta(K\ti)\cdot 
\left((D\times\{1\})\cap\delta(L\ti)\right)\right].$$

Le sous-groupe~$ \delta(K\ti)\cdot 
\left((D\times\{1\})\cap\delta(L\ti)\right)$
de~$D\times G$ est l'ensemble
des éléments de la forme~$(\deg z, \abs z)$ avec~$z\in L\ti$ et tels
que~$\abs z\in \abs{K\ti}$ ; en conséquence, le groupe quotient considéré
ci-dessus s'identifie à~$|L\ti|/|K\ti|$. 

\trois{interpretef} Soit~$\bf a$
un sous-ensemble de~$L\ti$ tel que~$|{\bf a}|$ constitue un système de représentants de~$|L\ti|/|k\ti|$
et soit~$\bf b$ une base de~$\red L_{D\times \{1\}}$ sur~$\red K_{D\times \{1\}}$. 
Il résulte
de~\ref{interpdimoide} et de~\ref{descresgrad}
que~$\red{\bf a}\cdot \{\bf b\}$ est une base de~$\red L$ sur~$\red K$ ; on a
donc~$[\red L : \red K]=[\red L_{D\times \{1\}}:\red K_{D\times \{1\}}]\cdot (|L\ti|:|k\ti|)$. 

\medskip
Compte-tenu de~\ref{dimresoide}, il vient~
$$[L:K]\geq [\red L:\red K]=[\red L_{D\times \{1\}}:\red K_{D\times \{1\}}]\cdot (|L\ti|:|k\ti|).$$ 

\trois{interpdegtrres} Soient~$d$ le degré de transcendance de~$\red L_{D\times\{1\}}$ sur~$\red K_{D\times\{1\}}$, et ~$\rho$
la dimension du~$\QQ$-espace vectoriel~$\QQ\otimes_{\ZZ}(|L\ti|/|k\ti|)$.
Soit~${\bf u}$ une base de transcendance de~$\red L_{D\times\{1\}}$ sur~$\red K_{D\times\{1\}}$,
et soit~${\bf v}$ une partie de~$L\ti$ tels que~$|{\bf v}|$ soit une base de ~$\QQ\otimes_{\ZZ}(|L\ti|/|k\ti|)$.
Il résulte de~\ref{interpdegtroide} et de~\ref{descresgrad}
que~${\bf u}\cup\red{\bf v}$ est
une base de transcendance de~$\red L$ sur~$\red K$. En conséquence, le degré de transcendance de~$\red L$ sur~$\red K$ est égal à~$d+\rho$.

\trois{transeteta} Soit~$(t_i)$ une famille d'éléments de~$L\ti$. 
Il résulte de~\ref{algindepetgauss}
que les assertions suivantes sont équivalentes : 

\medskip
a) les~$t_i$ sont algébriquement indépendants sur~$K$, et
~$\abs ._{|K({\bf t})}=\eta_{K,\abs {\bf t}}$ ; 

b) les éléments~$\red{t_i}$ de~$\red L$ sont algébriquement indépendants sur~$\red K$.

\medskip
En conséquence, le degré de transcendance
de~$\red L$ sur~$\red K$ est majoré par celui
de~$L$ sur~$K$.

 \trois{remresidgauss} {\em On suppose
 de maintenant que les conditions équivalentes~a) et~b) ci-dessus 
 sont satisfaites}. 
 En vertu de~\ref{algindepetgauss}, l'annéloïde~$\red{K[{\bf t}]}$ s'identifie
 à~$\red K[\red{\bf t}]$ ; en conséquence, le corpoïde résiduel~$\red{K({\bf t})}$ 
 est égal au sous-corpoïde de~$\red L$ engendré par~$\red K$ et~$\red {\bf t}$. 
 Notons que ledit sous-corpoïde
 est par ailleurs
 naturellement isomorphe à~$\red K((\deg({\bf t}, |{\bf t}|)\inv{\bf T})$ puisque~b) est vraie par hypothèse.

 \deux{abhyankarfini} Soit~$K\hookrightarrow L$ une extension de corpoïdes valués,
 de type fini comme extension
 de corpoïdes. On a vu au~\ref{transeteta} ci-dessus que le degré de transcendance
 de~$\red L$ sur~$\red K$ est majoré par celui de~$L$ sur~$K$. Supposons
 que l'on ait égalité, et soit~$(t_1,\ldots, t_n)$ une base
 de transcendance de~$\red L$ sur~$\red K$. Chaque~$t_i$
 est alors non nul, et est donc de la forme~$\red{a_i}$ pour un
 certain~$a_i\in k\ti$. En vertu de~\ref{transeteta}, les~$a_i$ sont algébriquement
 indépendants sur~$K$ ; d'après notre hypothèse sur les degrés de transcendance, 
 ils forment une base de transcendance de~$L$ sur~$K$. Le corps~$L$ est donc
 fini sur~$F:=K(a_1,\ldots, a_n)$. 
 
 \medskip
Il résulte de~\ref{remresidgauss} que~$\red F=\red K(t_1,\ldots, t_n)$, 
et de~\ref{dimresoide} que~$\red L$ est fini sur~$\red F$. En conséquence, 
le corpoïde~$\red L$ est de type fini sur~$\red K$. En d'autres
termes,~$\red L_{D\times\{1\}}$ est de type fini sur~$\red K_{D\times\{1\}}$, et~$|L\ti|/|k\ti|$ est 
de type fini (\ref{interpfinioide}). 

\deux{uniqueprolgauss} Soit~$K\hookrightarrow L$ une extension 
{\em aglébrique} de corpoïdes
valués. Donnons-nous une
famille~${\bf d}=(d_i)_{i\in I}$
d'éléments de~$D$ et une famille${\bf r}=(r_i)_{i\in I}$ d'éléments de~$G$,
les deux familles étant indexées par le même ensemble~$I$. 
La
valuation~$\eta_{L,{\bf d}, {\bf r}}$
est alors l'unique prolongement à~$L({\bf d}\inv{\bf T})$
de la valuation
$\eta_{K,{\bf d},{\bf r}}$ du corps~$K({\bf d}\inv {\bf T})$. 

\medskip
En effet, fixons un prolongement~$\abs .'$ de~$\eta_{K,{\bf d},{\bf r}}$
à~$L({\bf d}\inv {\bf T})$. La famille~$\red{\bf T}$
est algébriquement indépendante sur~$\red K$ d'après~\ref{transeteta}, et
comme~$L$ est algébrique sur~$K$, le corpoïde~$\red L$ est
algébrique sur~$\red  K$ (\ref{algalgtilde}) ; il s'ensuit 
que~$\red{{\bf T}}$
est algébriquement indépendante sur~$\red L$. 
Comme on a par ailleurs~$\abs {\bf T}'={\bf r}$ il découle
de~\ref{transeteta}
que~$\abs .'=\eta_{L,{\bf d},{\bf r}}$. 

\subsection*{Un théorème de Chevalley pour les espaces de Zariski-Riemann gradués}

\deux{lemrac} {\bf Lemme.} {\em Soit~$K$ un corpoïde, soit~$r\in D$ et soit~$P=\sum a_i T^i$ un élément de~$K[T/r]$ ; soit~$\abs .$ une valuation de~$K$ et soit~$L$ une extension
de~$K$ dans laquelle~$P$ est scindé. Les assertions suivantes sont équivalentes : 

\medskip
i) il existe
une racine~$\alpha$ de~$P$
telle que~$\abs  \alpha\leq 1$ ; 

ii) il existe~$i\geq 1$ tel que~$\abs  {a_i}\geq \abs  {a_0}$.}

\medskip
{\em Démonstration.} Supposons que i) soit vraie ; on a ~$a_0=-\sum\limits_{i\geq 1}a_i\alpha^i$, d'où l'inégalité~$$\abs  {a_0}\leq \max_{i\geq 1} \abs  {a_i}\abs  {\alpha^i}\leq \max_{i\geq 1} \abs  {a_i}$$ (puisque~$\abs  \alpha \leq 1$) et ii) est établie. 

Si ii) est vraie, écrivons~$P=a_n \prod (T-\alpha_i)$ avec~$\alpha_i\in L^r$ pour tout~$i$. Si chacun des~$\abs  {\alpha_i}$ était strictement supérieur à~$1$, les relations coefficients-racines entraîneraient immédiatement que~$\abs  {a_0}>\abs  {a_i}$ pour tout~$i>1$, ce qui contredirait l'hypothèse ii) ; par conséquent, i) est vraie.~$\Box$ 

\deux{imvalqc} {\bf Proposition.} {\em Soit~$F$ un corpoïde, et soit~$K\hookrightarrow L$ une extension finie de corpoïdes au-dessus de~$F$. Soit~$\mathsf U$ un ouvert quasi-compact de~$\PP_{L/F}$ ; l'image de~$\mathsf U$ sur~$\PP_{K/F}$ est alors un ouvert quasi-compact de ce dernier.}

\medskip
{\em Démonstration.} On peut supposer que~$\mathsf U$ est de la forme~$\PP_{L/F}\{f_1,\ldots,f_n\}$ où chacun des~$f_i$ est un élément de~$L$ dont on note~$r_i$ le degré. 

\trois{casunseulfi} {\em Le cas où~$n=1$} ; on pose alors~$f=f_1$ et~$r=r_1$. Soit~$P\in K[T/r]$ le polynôme minimal de~$f$ sur~$K$ ; écrivons~$P=\sum a_i T^i~$ ; notons que comme~$f\neq 0$ on a~$a_0\neq 0$ ; soit~$J$ l'ensemble des indices~$i$ tels que~$a_i\neq 0$. Soit~$\abs  .$ une valuation de~$K$ et soit~${\mathbb K}$ une clôture algébrique de~$K$ ; choisissons un prolongement~$\abs  .'~$ de~$\abs  .$ à~${\mathbb K}$. Dire que la valuation~$\abs  .$ appartient à l'image de~$\PP_{L/F}$ signifie qu'elle admet une extension à~$L$ dont l'anneau contient~$f$. Comme toute extension de~$\abs  .$ à~$L$ se déduit de~$\abs  .'$ {\em via} un~$K$-plongement de~$L$ dans~${\mathbb K}$, la valuation~$\abs  .$ appartient à l'image de~$\PP_{L/F}$  si et seulement si~$P$ possède dans~$\mathbb K$ une racine~$\alpha$ telle que~$\abs  \alpha '\leq 1$ ; en vertu du lemme~\ref{lemrac} ci-dessus, cela revient à demander qu'il existe~$i$ avec~$\abs  {a_i}\geq \abs  {a_0}$. Le scalaire~$a_0$ étant non nul, l'image de~$\mathsf U$ est finalement égale à la réunion des~$\PP_{K/F}\{a_0/a_i\}$ pour~$i$ parcourant~$J$ ; c'est donc bien un ouvert quasi-compact de~$\PP_{K/F}$, ce qui achève la preuve dans le cas où~$n=1$. 

\trois{casplusfi} {\em Le cas général}. Si~$n=0$ alors~$\PP_{L/F}\{f_1,\ldots,f_n\}=\PP_{L/F}$, et son image est égale à~$\PP_{K/F}$ tout entier, ce qui termine la démonstration.

\medskip
Supposons
maintenant que
$n>0$.


%


\medskip
Soit~$f$ l'élément~$\sum f_iT_i$ de~$L(r_1\inv T_1,\ldots,r_n\inv T_n)$ ; remarquons que~$f$ est de degré~$1$. Si~$\abs  .'$ est une valuation sur~$L$,
il est immédiat que~$\abs .'\in \PP_{L/F}\{f_1,\ldots,f_n\}$ si et seulement si~$\eta_{(L,\abs .'),{\bf r},{\bf 1}}\in \PP_{L({\bf r}\inv{\bf T})/F}\{f\}$. 

\medskip
Par ailleurs, il découle de~\ref{uniqueprolgauss} que si~$\abs  .\in \PP_{K/F}$ et si~$\abs .'$ est un prolongement
de~$\abs .$ à~$L$,
alors~$\eta_{(L,\abs .'), {\bf r}, {\bf 1}}$
est l'unique prolongement de~$\eta_{(K,\abs .), {\bf r},{\bf 1}}$
à~$L({\bf r}\inv {\bf T})$. 

\medskip
Il s'ensuit que
la valuation~$\abs .$ appartient à l'image de l'ouvert quasi-compact~$\PP_{L/F}\{f_1,\ldots,f_n\}$ si et seulement si~$\eta_{(K,\abs .),{\bf r}, {\bf 1}}$
 appartient à l'image~$\mathsf V$ de~$ \PP_{L ({\bf r}\inv {\bf T})/F}\{f\}$ sur~$\PP_{K ({\bf r}\inv{\bf T})/F}$.

Or en vertu du cas~$n=1$ traité au~\ref{casunseulfi} ci-dessus,~$\mathsf V$ est un ouvert quasi-compact de ~$\PP_{K({\bf r}\inv{\bf T})/F}$ ; il en résulte, en vertu de la
description explicite d'une valuation de Gauß, 
que~l'image récirproque de~$\mathsf V$ par l'application
$$\abs . \mapsto \eta_{(K,\abs .),{\bf r},{\bf 1}}$$ est un ouvert quasi-compact de~$\PP_{K/F}$, ce qui achève la démonstration.~$\Box$ 

\subsection*{Exemples d'homéomorphismes entre espaces de Zariski-Riemann gradués} 

\deux{torsionhomeozr} {\bf Lemme.} {\em Soit~$F$ un corpoïde, 
soit~$K$ une
extension de~$F$
et soit~$L$ une extension de~$K$. Supposons que~$L\ti /K\ti$ est 
de torsion. L'application continue~$\PP_{L/F}\to \PP_{K/F}$ est alors 
un homéomorphisme qui préserve les ouverts affines.} 

\medskip
{\em Démonstration.}  Soit~$\abs .$ une valuation sur~$L$ 
et soit~$f\in L$. Par hypothèse, il existe un entier~$n$ tel
que~$f^n\in K$. Comme~$f^n$ appartient à l'annéloïde
de~$\abs. $ si et seulement si c'est le cas de~$f$, on en déduit : 

$\bullet$ que~$\abs .$ est entièrement déterminée par sa restriction
à~$K$, et donc que~$\PP_{L/F}\to \PP_{K/F}$ est injective ; 

$\bullet$ que si~$(f_1,\ldots, f_r)$ sont des éléments de~$L$, 
et si les~$n_i$ sont des entiers tels que~$f_i^{n_i}\in K$
pour tout~$i$ alors l'image de~$\PP_{L/F}\{f_1,\ldots, f_r\}$
sur~$\PP_{K/F}$ est égale à $\PP_{L/F}\{f_1^{n_1},\ldots, f_r^{n_r}\}$ ; 
la bijection continue $\PP_{L/F}\to \PP_{K/F}$ est donc bien 
un homéomorphisme qui préserve les ouverts affines.~$\Box$ 

\deux{bijvalexem} {\em Exemples.} Mentionnons deux cas particuliers dans
lesquels ce lemme s'applique. 

\trois{bijvalinsep} Si~$L$ est une extension algébrique purement inséparable
de~$K$ alors~$L\ti/K\ti$ est de torsion ($p$-primaire). 

\trois{bijvalgroupetors} S'il exixte un sous-groupe~$\Delta$
de~$D$ tel que~$D/\Delta$ soit de torsion et tel que~$K=L_{\Delta}$
alors~$L\ti/K\ti$ est de torsion. 

\deux{memegroupebijval} {\bf Lemme.} {\em Soit~$K$ un corpoïde
et soit~$L$ une extension de~$K$. 
Si~$\deg (L\ti)\deg (K\ti)$ est
de torsion alors~$\PP_{L/K}\to \PP_{L_1/K_1}$ est un homéomorphisme
préservant les ouverts affines. } 

\medskip
{\em Démonstration.} Posons~$\Delta=\deg (K\ti)$. Le lemme~\ref{torsionhomeozr}
assure que
l'application continue~$\PP_{L/K}\to \PP_{L_\Delta/K}$ est un homéomorphisme préservant
les ouverts affines. Comme~$(L_\Delta)_1=L_1$ on peut remplacer~$L$
par~$L_\Delta$ et ainsi supposer que~$\deg(L\ti)=\deg (K\ti)$. 

\medskip
Soit~$\abs .\in \PP_{L_1/K_1}$ et soit~$f\in L$. Comme~$\deg (L\ti)=\deg (K\ti)$
il existe~$g\in L_1$ et~$h\in K\ti $ tels que~$f=gh$. La trivialité de~$\abs ._{|K_1}$ 
assure que~$\abs g$ ne dépend que de~$f$, et pas du choix de~$g$ et~$h$ ; 
il est donc licite de le noter~$\abs f'$. On vérifie immédiatement que~$\abs .'$
est une valuation,
dont la restriction à~$K$ est triviale et donc la restriction à~$L_1$ coïncide avec~$\abs . $. 
De plus, on voit aussitôt à l'aide de la formule
qui définit~$\abs .'$ 
que celle-ci est l'unique prolongement
de~$\abs .$ à~$L$ dont la restriction à~$K$ soit triviale. 
Ainsi, l'application continue~$\PP_{L/K}\to \PP_{L_1/K_1}$ est bijective. 

\medskip
Soient~$f_1,\ldots, f_n$ des éléments de~$L$.
Comme~$\deg (L\ti)=\deg (K\ti)$, on peut écrire chacune des~$f_i$ s
comme un produit~$g_ih_i$ où~$g_i\in L_1$ et~$h_i\in K\ti $ ; on a
alors~$\PP_{L/K}\{f_1,\ldots, f_n\}=\PP_{L/K}\{g_1,\ldots, g_r\}$ ; ainsi, 
l'image de~$\PP_{L/K}\{f_1,\ldots, f_n\}$
sur~$\PP_{L_1/K_1}$ est 
l'ouvert affine
$\PP_{L/K}\{g_1,\ldots, g_r\}$, ce qui termine
la démonstration.~$\Box$

\section{Les corps henséliens et leurs extensions modérément ramifiées}

On désigne par~$D$ un groupe abélien divisible. Dans cette section, 
les valuations
seront à valeurs dans~$D$, et les annéloïdes seront des~$D$-annéloïdes.

\subsection*{Rappels sur les prolongements
d'une valuation à une extension algébrique} 

\deux{convultramkzero} {\bf Quelques conventions.} Sauf  mention expresse du contraire,
la valuation d'un corps valué~$K$ sera notée~$\abs .$, son anneau sera noté~$K\zero$
et l'idéal maximal de celui-ci sera noté~$K\zeroo$. 

La notation~$\red K$ désignera le corpoïde résiduel de~$K$ ; le
corps résiduel traditionnel sera vu comme un sommande de~$\red K$,
et noté en conséquence~$\red K_1$. 

\deux{bourbakithval} Soit~$K$ un corps valué, et soit~$L$ une extension algébrique de~$K$. 
Soit~$B$ la fermeture intégrale
de~$K\zero$ dans~$L$.  On démontre les faits suivants. 

\medskip
i) L'application~$\got m\mapsto B_{\got m}$ établit
une bijection entre l'ensemble des idéaux maximaux de~$B$
et l'ensemble des anneaux de valuation de~$L$ qui dominent~$K\zero$ ; 
ces ensembles sont finis dès que~$L$ est finie sur~$K$. 

ii) Si~$L$ est une extension galoisienne de~$K$,
le groupe~${\rm Gal}(L/K)$ agit transitivement sur l'ensemble
des idéaux maximaux de~$B$, et donc d'après~i) sur l'ensemble
des prolongements de~$|.|$ à~$L$. 

\deux{approxfaible} Soient~$|.|_1,\ldots, |.|_r$
des prolongements deux à deux distincts 
de~$|.|$ à~$L$, et soient~$\got m_1,\ldots, \got m_r$
 les idéaux maximaux de~$B$ correspondants. 

Soit~$I$ un sous-ensemble
de~$\{1,\ldots, r\}$. Il existe alors~$x\in B$ tel
que~$|x-1|_i<1$ pour tout~$i\in I$ 
et~$|x|_j<1$ pour tout~$j\notin I$. 

En effet, soit~$j\notin I$. Par le lemme
d'évitement des idéaux premiers,~$\got m_j$
n'est pas contenu dans~$\bigcup\limits_{i\in I}  \got m_i$ ;
choisissons~$y_j\in \got m_j$ tel que~$y_j$
n'appartienne à aucun des~$\got m_i$. 

Posons~$y=\prod\limits_{j\notin I} y_j$. 
On a alors~$y\in \bigcap\limits_{j\notin I} \got m_j$, 
et~$y\notin\bigcup\limits_{i\in I}  \got m_i$. 

Comme~$\got m_i+\got m_j=B$ si~$i\neq j$,
le lemme chinois assure que
$$B/\prod_{i\in I}\got m_i\simeq \prod_{i\in I} B/\got m_i.$$

Par construction, l'image de~$y$ dans le terme de droite est inversible ; 
en conséquence, 
il existe~$z\in B$ tel que~$zy-1\in \prod\limits_{i\in I}\got m_i$. 

L'élément~$x:=yz$ de~$B$ satisfait alors aux conditions requises.                          

\subsection*{La notion de corps valué hensélien}

\deux{premimcorpshens} Soit~$K$ un corps valué. Soit~$P=T^{n}+\sum_{i\leq n-1} a_{i}T^{i}$
un polynôme unitaire à coefficients dans $K$,
et soit $L$ une extension valuée de~$K$ dans laquelle~$P$ est scindé. 
Posons $\rho(P)=\sup_{0\leq i\leq n-1} |a_{i}|^{1/(n-i)}$. Si $x$ 
est un élément de~$L$ tel que $|x|>\rho(P)$, alors $|x|^{n}
>|a_{i}|\cdot|x|^{i}$ pour tout 
entier~$i\leq n-1$ ; en 
conséquence,~$x$ ne peut être racine de $P$.
Par ailleurs, si~$\rho(P)$ majorait strictement en valeur absolue
toutes les racines de~$P$ dans~$L$, 
l'expression des~$a_{i}$ en fonction de ces dernières 
fournirait l'inégalité absurde~ $\sup_{0\leq i\leq n-1} |a_{i}|^{1/(n-i)}<\rho(P)$ ; 
il en découle que~$\rho(P)$
coïncide avec la borne supérieure 
(et même le maximum dès que~$n>0$) 
des valeurs absolues des racines de $P$ dans $L$. 
Il s'ensuit que si~$Q$ divise~$P$, alors~$\rho(Q)\leq \rho(P)$.

\trois{notationptilde} Pour tout $r\in D$ majorant $\rho(P)$, on notera $\red{P}_{r}$ le polynôme $X^{n}+\sum \widetilde{a_{i,r^{n-i}}}T^{i}$ ; 
c'est un élément de $\red K[T/r]$ de degré $r^{n}$. Il est égal à~$T^n$ si et seulement si~$r>\rho(P)$.
Si~$P$ s'écrit $RS$ avec~$R$ et~$S$ unitaires, on a~$\red{P}_{r}=\red{R}_{r}\red{S}_{r}$ ; 
en particulier, si l'on écrit $P=\prod (T-x_{i})$ dans $L$, alors $\red{P}_{r}=\prod (T-\widetilde{x_{i,r}})$.
Notons une conséquence simple de ce fait : si~$\red P_r$ est séparable,~$P$ est séparable. 

\medskip
Si $\rho(P)>0$ (autrement dit si $P$ n'est pas une puissance de $X$), on écrira simplement $\red{P}$ au lieu de $\red{P}_{\rho(P)}$.  

\medskip
Si $r\in D$ et si $R$ est un élément unitaire de $\red K[T/r]$,
on appellera {\em relèvement admissible} de $R$ tout polynôme unitaire $\bnd R$ appartenant à $K[T]$, 
tel que $\rho({\bnd R})\leq r$ et tel que $\red{\bnd R}_{r}=R$ ; l'existence d'un tel relèvement est immédiate.

\trois{uniciterelevement} Soit~$r\in D$ et soit~$P$ un polynôme unitaire de~$K[T]]$
tel que~$\rho(P)\leq r$. Supposons
donnée une factorisation~$\red P_r=QR$ dans~$\red K[T/r]$, où~$Q$ et~$R$ sont unitaires et premiers entre eux,
Il y a alors {\em au plus} un
relèvement admissible~$(\bnd Q,\bnd R)$ de~$(Q,R)$ tel que~$P=\bnd Q\bnd R$. 

En effet, supposons donné un tel relèvement, choisissons un
corps de décomposition~$L$ de~$P$, et un prolongement de~$|.|$
à~$L$. Écrivons~$P=\prod (T-\lambda_i)$ avec~$\lambda_i\in L$
pour tout~$i$. 

Pour tout indice~$i$, la réduction~$\red{\lambda_i}_r$ est une
racine de~$\red P_r$. C'est donc ou bien une racine
de~$Q$, ou bien une racine de~$R$, les deux situations
étant exclusives l'une de l'autre puisque~$Q$ et~$R$ sont premiers
entre eux. Si l'on note~$I$ (resp.~$J$) l'ensemble des indices~$i$
tels que~$\red{\lambda_i}_r$ est une
racine de~$Q$ (resp. de~$R$), on voit que
l'on a nécessairement 
$$\bnd Q=\prod_{i\in I}(X-\lambda_i)\;{\rm et}\;\bnd R=\prod_{i\in J}(X-\lambda_i),$$
d'où notre assertion. 

Notons un cas particulier : si~$\alpha$ est une racine simple de~$\red P_r$, il existe
au plus un élément~$a$ de~$K$ tel que~$|a|\leq r, \red a_r=\alpha$ et~$P(a)=0$. 

\deux{lemalgclossepclosoide}
{\bf Lemme.}
{\em Soit~$K$ un corps valué. Si~$K$ est 
algébriquement clos (resp. séparablement clos), il en va de même de~$\red K$.}

\medskip
{\em Démonstration.} Soit~$P$ un polynôme unitaire 
appartenant à~$K[T/r]$ pour un certain~$r\in D$, et
soit~$\cal P$ un relèvement admissible de~$P$. On a
vu au~\ref{notationptilde}
ci-dessus que si~$P$ 
est séparable, il en va de même de~$\cal P$ ; et que si~$\cal P$ est scindé,
il en va de même de~$P$. La conclusion du lemme s'ensuit aussitôt.~$\Box$

\deux{prehensgrad} {\bf Lemme.} {\em Soit~$K$ un corps valué,
soit~$P$ un polynôme unitaire irréductible de~$K[T]$, et soit~$r$
un élément de~$D$ majorant~$\rho(P)$. Si il existe deux polynômes
unitaires~$R$ et~$S$ non constants et premiers entre eux dans~$\red K[T/r]$ 
tels que~$\red P=RS$ alors~$r=\rho(P)$ et il existe au moins
deux prolongements distincts de~$|.|$ à~$L$.}

\medskip
{\em Démonstration.} Comme $R$ et $S$ sont premiers
entre eux,~$\red{P}_{r}$ n'est pas
une puissance de~$T$
et l'on a donc $\rho(P)=r$. 
Fixons un prolongement de~$|.|$ à~$L$,
encore noté~$|.|$. 
Écrivons $P=\prod (T-x_{i})$ dans $L$. 
Puisque $\red P_r=RS$, 
il existe deux indices $i$ et $j$ distincts tels 
que~$\widetilde{x_{i,r}}$~(resp. $\widetilde{x_{j,r}}$)
soit une racine de $R$ (resp. de $S$) ; comme~$R$ et $S$~sont 
premiers entre eux, $\widetilde{x_{j,r}}$ n'annule pas $R$. 

\medskip
Soit $\bnd R$ un relevé admissible de $R$. Si $m$ désigne le degré monomial
de~$R$ 
on a alors~$|\bnd R(x_{i})|<r^{m}$ puisque $R(\widetilde{x_{i,r}})=0$, 
et~$|\bnd R(x_j)|=r^m$ puisque~$R(\red{x_{j,r})}\neq 0$. 

\medskip
Comme~$P$ est irréductible, 
il existe un~$k$-automorphisme~$g$ de~$L$
envoyant~$x_{i}$ sur~$x_{j}$ ; on a 
alors~$g(\bnd R(x_{i}))=\bnd R(x_j)$. Par
ce qui précède,
$$|g(\bnd R(x_{i}))|=r^m>|\bnd R(x_i)|.$$
Ainsi,~$|.|\circ g$ est un prolongement de~$|.|$
à~$L$ qui est distinct de~$|.|$ .~$\Box$

\deux{equivhens} {\bf Proposition.} {\em Soit~$K$ un corps valué. 
Les assertions suivantes sont équivalentes.

\medskip
i) L'anneau local~$K\zero$ est hensélien. 

ii) Pour toute extension algébrique~$L$ de~$K$, la valuation~$|.|$
admet un unique prolongement à~$L$. 

iii) Pour tout~$r\in D$, pour tout polynôme unitaire~$P\in K[T]$ tel
que~$\rho(P)\leq r$,  pour toute factorisation~$\red P_r=QR$
où~$Q$ et~$R$ sont des polynômes unitaires et premiers entre eux de~$K[T/r]$,
 il existe un (unique) relèvement
admissible~$(\bnd Q,\bnd R)$ de~$(Q,R)$ tel 
que~$P=\bnd Q\bnd R$. 

iv)  Pour tout~$r\in D$, pour
tout polynôme unitaire~$P\in K[T]$
et pour
toute racine
simple~$\alpha$ de~$\red P_r$
dans~$\red K$, il existe un (unique)
élément~$a$ de~$K$ tel que
l'on ait~$|a|\leq r$, $\red a_r=\alpha$
et~$P(a)=0$. }

\medskip
{\em Démonstration.} On procède par implications successives. 

\trois{equivhensi} Supposons~i) vraie, et soit~$L$ une extension algébrique
de~$K$. Soit~$B$ la fermeture intégrale de~$A$
dans~$L$, et soit~$\sch B$ 
l'ensemble des sous-algèbres de type fini de~$B$. 
Soit~$C\in \sch B$. L'algèbre~$C$ est 
un produit fini d'anneaux locaux puisque~$A$ est hensélien ; 
étant d'autre part intègre, elle est locale. 

Par ailleurs, si~$C$ et~$C'$ sont deux éléments
de~$\sch B$ avec~$C\subset C'$ alors~$C'$ est entière sur~$C$,
et~$C\hookrightarrow C'$ est dès lors un morphisme local. On en déduit
que~$B$ est un anneau local. En vertu de~\ref{bourbakithval}, cela équivaut
à dire que~$|.|$ admet un unique prolongement à~$L$, d'où~ii). 

\medskip
\trois{equivhensii} Supposons~ii) vraie, et montrons~iii). 
On raisonne par récurrence sur~$\deg P$. 
S'il est nul le résultat est trivial ;
on suppose donc qu'il est strictement
positif, et que le résultat a été 
prouvé en degré~$<\deg P$. 

\medskip
Si~$R$ ou~$S$ est constant, le résultat est  évident. On peut
donc les supposer non constants. Comme on 
fait l'hypothèse que~ii) est vraie, il résulte
du lemme~\ref{prehensgrad} que~$P$ {\em n'est pas}
irréductible. \'Ecrivons~$P=P_1P_2$, où~$P_1$
et~$P_2$ sont unitaires de degré strictement
inférieur à celui de~$P$. 

\medskip
Pour tout $i$ notons $Q_i$ (resp. $R_i$) le PGCD unitaire de~$Q$ (resp.~$R$)
et~$\red P_{i,r}$. Comme~$\red P_r=QR$
et comme~$Q$ et~$R$
sont premiers entre eux, 
il vient~$\red P_{1,r}=Q_1R_1$ et~$\red P_{2,r}=Q_2R_2$. Les 
polynômes~$Q_1$ et~$R_1$ étant premiers entre eux, 
l'hypothèse de récurrence assure l'existence
d'un relèvement admissible~$(\bnd Q_1,\bnd R_1)$ de~$(Q_1,R_1)$
tel que~$P_1=\bnd Q_1\bnd R_1$. 
De même, il existe un
relèvement admissible~$(\bnd Q_2,\bnd R_2)$ de~$(Q_2,R_2)$
tel que~$P_2=\bnd Q_1\bnd R_2$. 
Le couple~$(\bnd Q_1\bnd Q_2,\bnd R_1\bnd R_2)$ est alors un relèvement
admissible de~$(Q,R)$, 
et l'on
a~$P=P_1P_2=\bnd Q_1\bnd Q_2\bnd R_1\bnd R_2$, d'où~iii). 

\trois{equivhensiii} Si~iii) est vraie, alors~iv) est vraie : ce n'est qu'une reformulation
de~iii) dans le cas particulier où l'un des deux polynômes~$Q$ ou~$R$ est de degré
monomial égal à~$1$. 

\trois{equivhensiv} Supposons que~iv) est vraie. Dans le cas particulier
où~$r=1$, cette propriété signifie exactement que toute racine simple
modulo~$K\zeroo$ d'un polynôme unitaire~$P\in K[T]$ se relève en une racine
de~$P$. En conséquence,~$K\zero$ est hensélien.~$\Box$ 

\deux{defhensval} On dit qu'un corps valué~$K$ est {\em hensélien}
s'il satisfait les quatre propriétés équivalentes ci-dessus. 

\trois{remhentestsep} Si~$L\hookrightarrow L'$ est une
extension algébrique radicelle, alors pour 
tout~$x\in L'$ il existe un entier~$N>0$ tel que~$x^N\in L$. Il s'ensuit que
toute valuation de~$L$ admet
un unique prolongement à~$L'$ ; pour s'assurer qu'un corps~$K$ est 
hensélien, on peut donc se contenter de vérifier
propriété~ii) ci-dessus pour les extensions algébriques~{\em séparables}~$L$
de~$K$. 

\trois{exthenshens} Il résulte de la caractérisation
d'un corps hensélien par la propriété~ii) de la proposition~\ref{equivhens}
ci-dessus que  toute extension algébrique d'un corps
hensélien est encore hensélienne.

\subsection*{Groupe de décomposition, hensélisé d'un corps valué}

\deux{defgroupedec} Soit~$K\hookrightarrow L$ une extension
galoisienne de corps valués. 
On appelle {\em groupe de décomposition} de~$L$ sur~$K$
le sous-groupe de~$\mathsf{Gal}(L/K)$ formé
des automorphismes~$g$ tels que~$\abs . \circ g=\abs .$, et on
le note~$\mathsf D(L/K)$.

\deux{groupedecompcomp} Soit~$F$ une sous-extension galoisienne de~$L$. La restriction 
induit un morphisme~$\mathsf D(L/K)\to \mathsf D(F/K)$. 

\medskip
{\em Ce morphisme est surjectif.} En effet, soit~$g\in \mathsf D(F/K)$. Prolongeons
l'automorphisme~$g$ 
en un automorphisme de~$L$, noté encore~$g$. 
Les valuations~$\abs .$ et~$\abs . \circ g$ de~$L$ 
induisent toutes deux la valuation~$\abs .$ sur~$F$,
puisque~$g|_F\in \mathsf D(F/K)$. En conséquence, il existe
un élément~$h\in \mathsf{Gal}\;(L/F)$ tel que~$\abs .\circ g=\abs .\circ h$. 

Il vient~$\abs .\circ gh\inv=\abs .$ . Ainsi,~$gh\inv \in \mathsf D(L/K)$. 
Comme~$h\in \mathsf {Gal}\;(L/K)$, on a~$gh\inv|_F=g|_F$, 
ce qui achève la démonstration. 

\deux{henselimmediat} {\bf Proposition.} {\em Soit~$K\hookrightarrow L$ une
extension galoisienne
de corps valués, et soit~$F$ le corps~$L^{\mathsf D(L/K)}$. On a~$\red F=\red K$
(ou, ce qui revient au même,~$\red F_1=\red K_1$ et~$|F\ti|=|k\ti|$).}

\medskip
{\em Démonstration.} En vertu de~\ref{groupedecompcomp}, il suffit de montrer
l'assertion pour les sous-extensions finies galoisiennes de~$L$, ce qui autorise à
supposer que~$L$ est finie sur~$K$. Soient~$\abs .=\abs ._0,\abs ._1,\ldots, \abs ._r$ les
prolongements à~$L$
de la valuation~$\abs .$ du corps~$K$. 

\medskip
Comme deux prolongements
de la valuation~$|.|$ de~$F$ au
corps~$L$ sont conjugués sous 
l'action de~$\mathsf{Gal}(L/F)=\mathsf D(L/K)$, la
restriction de~$\abs._i$ à~$F$ est différente de~$\abs .$ dès
que~$i\neq 0$. Il résulte
alors de~\ref{approxfaible},
qu'il existe~$z\in F\ti$ tel que~$|z-1|<1$ et~$|z|_j<1$ pour tout~$j>0$. 

\medskip
Soit~$x\in F\ti$. Nous allons  montrer que~$\red x\in \red K$. Pour tout
entier~$n$, on a~$\red {z^n}=1$ 
et partant~$\red{z^nx}=\red x$. Par ailleurs,
comme~$|z|_j\neq 1$ pour tout~$j>0$, il existe
un entier~$n$ tel que~$|z|_j^n\cdot |x|_j$ soit
différent de~$|z^nx|=|x|$ pour tout~$j>0$. On peut
donc, quitte à remplacer~$x$ par~$z^nx$, supposer
que~$|x|_j\neq |x|$ dès que~$j>0$.

\medskip
Soit~$\mathsf E$ l'ensemble des conjugués de~$x$. 
Soit~$y\in\mathsf E$. 
Il est de la forme~$g(x)$ 
avec~$g\in\mathsf{Gal}(L/K)$. Si~$y\neq x$
alors~$g\notin\mathsf D(L/K)$,
 ce qui veut dire que~$\abs .\circ g=|.|_j$ pour
un certain~$j>0$ ; par conséquent,~$|y|=|x|_j\neq |x|$. 
On note~$\mathsf F$ le sous-ensemble de~$\mathsf E$ formé des
conjugués~$y$ de~$x$ tels que~$|y|>|x|$, et~$P$
le polynôme minimal de~$x$. On écrit~$P=T^m+a_1T^{m-1}+\ldots+a_m$,
et l'on pose~$a_0=1$. Soit~$N$ le cardinal de~$\mathsf F$. 
Pour toute partie~$\mathsf P$
de~$\mathsf E$ de cardinal~$N+1$
et différente de~$\mathsf F\cup\{x\}$, on a 
$$\left| \prod_{y\in \mathsf P}y\right|<\left| \prod_{y\in \mathsf F\cup\{x\}}y\right|.$$
Comme~$$a_{N+1}=(-1)^{N+1}\sum_{\mathsf P\subset \mathsf E, \;{\rm card}\; \mathsf P=N+1}
\;\;\;
\prod_{y\in \mathsf P}y,$$
il vient~$\red{a_{N+1}}= (-1)^{N+1}\prod\limits_{y\in \mathsf F\cup\{x\}}\red y$. 

Pour toute partie~$\mathsf P$
de~$\mathsf E$ de cardinal~$N$ 
et différente de~$\mathsf F$, on a 
$$\left| \prod_{y\in \mathsf P}y\right|<\left| \prod_{y\in \mathsf F}y\right|.$$
Comme~$$a_{N}=(-1)^N\sum_{\mathsf P\subset \mathsf E, \;{\rm card}\; \mathsf P=N}
\;\;\;\;\prod_{y\in \mathsf P}y,$$
il vient~$\red{a_N}= (-1)^N\prod\limits_{y\in \mathsf F}\red y$. 

\medskip
Il s'ensuit que $\red {a_N}\neq 0$,
et que~$\red x=-\red{a_{N+1}}/\red{a_N}$. 
En conséquence,~$\red x\in \red K$.~$\Box$ 

\deux{deschensvalth} {\bf Proposition.} {\em Soit~$K$ un corps valué, et soit~$K^s$ 
une clôture séparable de~$K$. On fixe un prolongement de~$\abs .$ à~$K^s$,
que l'on note encore~$\abs .$ . Soit~$F$ le sous-corps~$(K^s)^{\mathsf D(K^s/K)}$. 

\medskip
i) Le corps valué~$F$ (muni de la restriction de~$\abs .$) 
est hensélien, et l'anneau~$F\zero$ s'identifie au hensélisé de~$K\zero$.

ii) Pour toute extension valuée hensélienne~$L$ de~$K$, il existe un unique~$K$-plongement isométrique
de~$K$ dans~$L$.  }

\medskip
{\em Démonstration.} Comme le groupe de Galois de~$K^s$ sur~$F$
est le stabilisateur de~$\abs .$, la valuation~$\abs .$ de~$F$ admet une unique
extension à~$F^s$, et donc à toute extension séparable de~$F$. En 
conséquence,~$F$ est
hensélien (\ref{remhentestsep}). L'anneau local~$F\zero$
est donc lui-même hensélien. 
Il en résulte l'existence
d'un~$(K\zero)$-morphisme canonique de~$(K\zero)^h$ dans~$F\zero$.

L'anneau~$K\zero$ est par ailleurs normal, $F\zero$ est intègre, 
et~$K\zero$ s'injecte dans~$F\zero$. On
déduit alors de~\ref{henselnormal}
que~$(K\zero)^h$ est normal et s'injecte dans~$F\zero$. 

\medskip
Soit~$E$ le corps des fractions de~$(K\zero)^h\subset F$. Comme~$(K\zero)^h$ est 
normal, il contient la fermeture intégrale~$B$ de~$K\zero$ dans~$E$. D'autre part,
l'anneau local~$(K\zero)^h$ domine~$K\zero$ ; en conséquence, il domine~$B_{\got m}$
pour un certain idéal maximal~$\got m$ de~$B$. L'anneau~$B_\got m$ étant l'anneau d'une valuation
qui prolonge~$\abs .$, il coïncide par maximalité avec~$(K\zero)^h$. 

\medskip
L'inclusion~$(K\zero)^h\hookrightarrow F\zero$ est locale. L'idéal maximal
de~$(K\zero)^h$ est dès lors égal à~$E\zeroo$. 
Comme un anneau de valuation 
d'un corps est connu dès qu'on connaît
son idéal maximal, il vient~$(K\zero)^h=E\zero$. 

\medskip
L'anneau~$(K\zero)^h$ est hensélien ; le sous-corps~$E$ de ~$K^s$
est donc hensélien. En conséquence, la valuation~$\abs .$ de~$E$ admet un unique
prolongement à~$K^s$. Cela signifie que~$\mathsf{Gal}(K^s/E)\subset \mathsf D(K^s/K)$. 
Autrement dit,~$F\subset E$ ; compte-tenu du fait que~$E\subset F$, il
vient~$E=F$ et~$F\zero =(K\zero)^h$, ce qui achève de prouver~i). 

\medskip
Donnons-nous maintenant une extension valuée hensélienne~$L$ de~$K$. 
Il existe alors un unique~$K\zero$-morphisme 
local de~$F\zero=(K\zero)^h$ dans~$L\zero$, qui est injectif 
d'après la remarque~\ref{henselnormal}. 

\medskip
Or la donnée d'un~$K$-plongement isométrique de~$F$ 
dans~$L$ équivaut à celle d'un~$K\zero$-plongement local
de~$F\zero$ dans~$L\zero$ ; compte-tenu de ce qui précède,
on en déduit~ii).~$\Box$

\deux{convcorphens} Si~$K$ est un corps valué, 
son hensélisé~$K^h$ est canonique en vertu de l'assertion~ii)
de la proposition~\ref{deschensvalth}. Cette même assertion,
ou la définition directe de~$K^h$, assurent que~$K^h=K$ 
dès que~$K$ est hensélien. 

\medskip
Si~$L$ est une extension algébrique de~$K^h$, on la
considèrera comme implicitement munie de l'unique prolongement
de la valuation de~$K^h$ ; ce prolongement fait de~$L$ un corps valué
hensélien. 

\deux{lemextvalhens} {\bf Lemme.} {\em Soit~$K$ un corps valué et
soit~$L$ une extension algébrique~$K$. 

\medskip
i) La~$K^h$-algèbre~$K^h\otimes_{K}L$ est réduite ; c'est en particulier
un produit
fini d'extensions de~$K^h$ si~$L$ est finie sur~$K$. 

ii) Soit~$(K^h\otimes_K L\to L_i)_i$ la famille
des quotients de~$K^h\otimes_K L$ par ses idéaux maximaux. 
Pour tout~$i$, désignons par~$\abs ._i$ la restriction à~$L$ 
de la valuation de~$L_i$. Les~$\abs ._i$ sont alors
deux à deux distinctes, et sont exactement les prolongements
de~$\abs .$ à~$L$. De plus,~$L_i$ s'identifie
pour tout~$i$ à~$(L,\abs ._i)^h$.

iii) Supposons que~$L$ soit galoisienne. Pour tout~$i$,
l'injection naturelle~$\mathsf{Gal}(L_i/K^h)\hookrightarrow \mathsf{Gal}\;(L/K)$
a pour image le groupe de décomposition de la valuation~$\abs ._i$.} 

\medskip
{\em Démonstration.} Comme~$K^h$ est réunion d'extensions finies
séparables de~$K$,
la~$K^ h$-algèbre~$K^h\otimes_{K}L$
est réduite, d'où i).  

\medskip
Fixons~$i$. Le corps valué~$L_i$ est
hensélien ; en conséquence~$(L,\abs ._i)^h$
se plonge canoniquement de manière
isométrique dans~$L_i$ ; par ailleurs,~$K^h$
se plonge également canoniquement de manière isométrique
dans~$(L,\abs ._i)^h$. Le sous-corps~$(L,\abs ._i)^h$ de~$L_i$
contenant~$L$ et~$K^h$, il coïncide avec~$L_i$ (puisque
ce dernier est un quotient
de~$K^h\otimes_KL$) ; plus précisément,
ce sont les~{\em $K^h\otimes_KL$-algèbres}
$L_i$ et~$(L,\abs ._i)^h$, 
que l'on identifie par ce biais. 
Comme la~$K^h\otimes_KL$-algèbre~$(L,\abs ._i)^h$
ne dépend que de~$\abs ._i$, on
peut retrouver le facteur~$L_i$ à partir de~$\abs ._i$. De ce fait,
les~$\abs ._i$ sont deux à deux distinctes. 

Enfin, soit~$\abs .'$ une valuation de~$L$ prolongeant~$\abs .$ .
Comme~$K^h$ se plonge canoniquement de manière isométrique dans~$(L,\abs .')^h$,
il existe un morphisme~$K^h\otimes_KL\to (L,\abs .')^h$, qui se factorise nécessairement
par l'un des~$L_i$. On a alors par construction~$\abs .'=\abs ._i$, 
ce qui achève de prouver~ii).~$\Box$

\medskip
Supposons maintenant que~$L$ est galoisienne,
et fixons~$i$. Comme~$K^h$
est hensélien, tout~$K^h$-automorphisme de~$L_i$ en préserve
la valuation, et l'image de~$\mathsf{Gal}(L_i/K^h)$ est donc bien contenue
dans~$\mathsf D((L,\abs ._i)/K)$.

Pour établir l'inclusion réciproque, fixons un~$K$-plongement
de~$L$ dans
une clôture séparable~$K^s$ de~$K$, et munissons celle-ci
d'un prolongement arbitraire de~$\abs ._i$, encore noté~$\abs ._i$. 
On identifie~$K^h$ (resp.~$L_i$) 
au sous-corps
de~$K^s$ formé des éléments invariants
sous~$\mathsf D((K^s,\abs ._i)/K)$ (resp. $\mathsf D((K^s,\abs ._i)/L)$).  

Soit~$g$
un élément de~$\mathsf D((L,\abs ._i)/K)$. D'après~\ref{groupedecompcomp},
on peut relever~$g$ en un élément~$\gamma$
de~$\mathsf D((K^s,\abs ._i)/K)$. On a alors~$(\gamma_{|L_i})_{|L}=g$, 
et la démonstration est terminée.~$\Box$ 

\deux{ramenerhensel} Soit~$K\hookrightarrow L$ une extension
galoisienne de corps valués et soit~$g$ un
élément de~$\mathsf D(L/K)$. 
Comme~$g$ préserve~$\abs .$, il induit un~$\red K$-automorphisme
du corpoïde résiduel~$\red L$,
et le morphisme~$\mathsf D(L/K)\to \mathsf{Gal}(\red L/\red K)$
ainsi induit est continu (pour que l'image d'un
élément~$g$ de~$\mathsf D(L/K)$ fixe un ensemble fini~$E$
d'éléments de~$\red L\ti$, il suffit qu'il fixe un ensemble fini~$F\subset L\ti$
tel que~$\red F=E$). 

\medskip
En vertu du lemme~\ref{lemextvalhens} ci-dessus, 
l'extension~$K^h\hookrightarrow L^h$ 
est galoisienne, et son groupe
de Galois s'identifie naturellement à~$\mathsf D(L/K)$. 
De plus, la proposition~\ref{henselimmediat}
assure que~$\red{K^h}=\red K$ et~$\red {L^h}=\red L$. 

En conséquence, pour étudier le groupe~$\mathsf D$ et son
action sur~$\red K$, on peut
remplacer~$K$ par~$K^h$
et~$L$ par~$L^h$, supposer que~$K$ (et~$L$) sont
henséliens, et que~$\mathsf D=\mathsf{Gal}(L/K)$. 

\deux{defdefstable} Soit~$K$ un corps valué et soit~$L$
une extension finie de~$K$. Soient~$\abs ._1,\ldots, \abs ._n$ 
les prolongements de~$\abs .$ à~$L$. Pour tout~$i$, 
notons~$\red L_i$ le corpoïde résiduel de~$(L,\abs ._i)$. 

Le lemme~\ref{lemextvalhens} ci-dessus
assure que~$K^h\otimes_K L\simeq \prod (L,\abs ._i)^h$. 
On a~$\red{K^h}=\red K$ et~$\red{(L,\abs ._i)^h}=\red L_i$ ; 
en vertu de\ref{dimresoide}, 
il vient~$[\red L_i:\red K]\leq [(L,\abs ._i)^h:K^h]$. 

On en déduit que
$$\sum [\red L_i:\red K]\leq [L:K].$$ 

\trois{defstable} On dit que l'extension~$L$ de~$K$ est {\em sans défaut}
si~$\sum [\red L_i:\red K]=[L:K]$. Le corps
valué~$K$ est dit {\em stable}
si toute extension finie de~$K$ est sans défaut. 

\trois{trivstable} Si la valuation de~$K$ est triviale, le corps~$K$ est
stable. 

\trois{avdstable} Si la valuation de~$K$ est «discrète»,
c'est-à-dire si~$|K\ti|$ est libre de rang~$1$, le corps~$K$
est stable si et seulement si l'anneau de valuation
discrète~$K\zero$ est 
excellent. C'est automatiquement le cas dans chacune des situations
suivantes : 

\medskip
\begin{itemize}

\item[$\bullet$] $K$ est de caractéristique nulle ; 

\item[$\bullet$] $K$ est complet ; 

\item[$\bullet$] $K$ est le corps des fonctions d'une courbe
projective normale et irréductible~$\sch C$
sur un corps~$k$ et~$K\zero$ est l'anneau local en un 
point fermé de~$\sch C$. 

\end{itemize}

\trois{anticipe0} Nous verrons un peu plus bas que si~$\red K$
est de caractéristique nulle alors~$K$ est stable.

\deux{henselsansdef} Soit~$K$ un corps valué
hensélien. 
Une extension finie~$L$ de~$K$ est sans défaut si~$[\red L:\red K]=[L:K]$. 

Si~$F$ est une extension finie de~$K$ et si~$L$ est une extension finie de~$F$,
alors~$L$ est sans défaut sur~$K$ si et seulement si~$F$ est sans défaut sur~$K$
et~$L$ est sans défaut sur~$F$ : c'est une conséquence triviale de la transitivité des degrés et
des inégalités
$$[\red L:\red K]\leq [L:K], \;\;[\red L:\red F]\leq [L:F]\;\;{\rm et}\;[\red F:\red K]\leq [L:K].$$

\subsection*{Théorie de Galois d'un corps valué hensélien}

\deux{galsurj} {\bf Proposition.} {\em Soit~$K$ un corps
valué hensélien et soit~$L$ une extension galoisienne de~$K$. 
L'extension~$\red L$ de~$\red K$ est alors normale, 
et
la flèche~$\mathsf {Gal}(L/K)\to \mathsf{Gal}(\red L/\red K)$
est surjective.} 

\medskip
{\em Démonstration.} Soit $\xi$ un élément non nul de $\red{L}$, soit~$r$
son degré 
et soit~$x\in L$ tel que~$\red x=\xi$ ; soit $P$ le polynôme minimal de $x$ sur $L$. 
Écrivons $P=\prod (T-x_{i})$, où $x_{1},\ldots,x_{d}$ sont les éléments de l'orbite de $x$ sous $\mathsf {Gal}(L/K)$. 
Les~$x_i$ sont tous de valeur absolue égale 
à~$r$. 
Le polynôme $\red P_r$ est un élément 
unitaire de $\red K[T/r]$
qui annule $\xi$ et est scindé dans $\red L$ (il est égal à~$\prod (T-\red {x_i})$ ) ; 
le polynôme minimal $R$ de $\xi$ est donc scindé dans $\red L$,
ce qui montre que $\red L$ est une extension normale de $\red K$. 

\medskip
Soit $\eta$ une racine de $R$ distincte de $\xi$ ; elle coïncide avec $\red{x_{i}}$ pour un certain 
entier~$i$.
Il existe $g\in \mathsf {Gal}(L/k)$ tel que $g(x)=x_{i}$, et l'image de $g$ dans $\mathsf{Gal} (\red{L}/\red K)$ envoie 
par construction~$\xi$ sur $\eta$. 
On en déduit qu'un élément $\red L$ est invariant sous l'image de $\mathsf {Gal}(L/k)$ si et seulement s'il
coïncide avec tous ses conjugués, donc si et seulement s'il est invariant sous
$\mathsf{Gal} (\red{L}/\red K)$ ; par la correspondance de Galois, 
l'image de $\mathsf {Gal}(L/K)$ (qui est un sous-groupe
compact de~$\mathsf{Gal}(\red L/\red K)$)
est égale à $\mathsf{Gal} (\red{L}/\red K)$ tout entier.~$\Box$ 

\deux{descextresoide} On conserve les hypothèses et notations
de la proposition~\ref{galsurj} ci-dessus. 

\trois{galredlredk} Il découle 
de~\ref{descgaloideclassetgp} {\em et sq.} 
que~$\red L_1$ est une extension normale
de~$\red K_1$, 
que la restriction~$\mathsf{Gal}(\red L/\red K)\to \mathsf{Gal}(\red L_1/\red K_1)$
est surjective, et que la formule
$$g\mapsto (|x|\mapsto \red{g(x)/x})$$
définit sans ambiguïté un isomorphisme
$$\mathsf{Ker}(\mathsf{Gal}(\red L/\red K)\to \mathsf{Gal}(\red L_1/\red K_1))
\simeq \mathsf{Hom}(|L^\times|/|K^\times|,\mu(\red L_1)).$$

Comme~$\mathsf{Ker}(\mathsf{Gal}(\red L/\red K)\to \mathsf{Gal}(\red L_1/\red K_1))$
est abélien, il hérite d'une structure naturelle de~$\mathsf{Gal}(\red L_1/\red K_1)$module,
et l'isomorphisme ci-dessus est alors~$\mathsf {Gal}(\red L_1/\red F_1)$ équivariant. 

\trois{inertieetram} On note~$\mathsf W(L/K)$ 
le noyau de~$\mathsf{Gal}(L/K)\to \mathsf{Gal}(\red L/\red K)$, 
et~$\mathsf I(L/K)$ celui de la surjection 
composée~$\mathsf{Gal}(L/K)\to \mathsf{Gal}(\red L/\red K)
\to  \mathsf{Gal}(\red L_1/\red K_1)$. On dit que~$\mathsf I(L/K)$ (resp.~$\mathsf W(L/K)$)
est le~{\em groupe d'inertie} (resp. le~{\em groupe de ramification}) de~$L$
sur~$K$. 

\medskip
Soit~$g\in \mathsf {Gal}(L/K)$. Par définition,~$g$ appartient à~$\mathsf W(L/K)$ si et seulement
si~$|g(x)-x|<|x|$ pour tout~$x\in L\ti$ ; et il appartient à~$\mathsf I(L/K)$
si l'inégalité précédente est vraie pour tout~$x$ {\em de valeur absolue égale à~$1$}. 

\trois{quotinertram} Par définition,~$\mathsf I(L/K)$ 
est un sous-groupe distingué de~$\mathsf {Gal}(L/K)$,
et~$\mathsf {Gal}(L/K)/\mathsf I(L/K)$ s'identifie
à~$\mathsf{Gal}(\red L_1/\red K_1)$. 

\medskip
De même, $\mathsf W(L/K)$ 
est un sous-groupe distingué de~$\mathsf {Gal}(L/K)$,
et~$\mathsf {Gal}(L/K)/\mathsf W(L/K)$ s'identifie
à~$\mathsf{Gal}(\red L/\red K)$.

\medskip
Quant au quotient~$\mathsf I(L/K)/\mathsf W(L/K)$, 
il s'identifie d'après~\ref{galredlredk}
à~$\mathsf{Hom}(|L\ti|/|k\ti|, \mu(\red L_1))$, par la
formule~$g\mapsto (|x|\mapsto \red{g(x)/x})$. 

Comme $\mathsf I(L/K)/\mathsf W(L/K)$ est abélien, 
il hérite d'une vraie action du
groupe~$\mathsf{Gal}(\red L_1/\red K_1)$,
pour laquelle l'identification ci-dessus est équivariante.

\trois{assezracinetilde} Soit~$p$ l'exposant caractéristique de~$\red K$. Rappelons
que si~$A$ est un groupe abélien de torsion,
on désigne par~$A_{p'}$ sons sous-groupe
de torsion première à~$p$, et par~$A^\vee$ son dual de Pontrjagin. 

\medskip
D'après~\ref{assezracinepasfini}, le polynôme~$T^a-1$
est scindé dans~$\red L_1$ 
pour tout
entier~$a$ divisant l'exposant généralisé de~$(|L^\times|/|K^\times|)_{p'}$,
et tout plongement de~$\mu(\red L_1)$ dans~$\QQ/\ZZ$ induit un isomorphisme
$$\mathsf{Hom}(|L\ti|/|k\ti|, \mu(\red L_1))\simeq (|L\ti|/|K\ti|)_{p'}^\vee.$$

\deux{ietwcompres} {\bf Lemme.} {\em Soit~$K$ un corps valué hensélien, soit~$L$
une extension galoisienne de~$K$, et soit~$F$ une sous-extension galoisienne de~$L$. 
Les restrictions~$\mathsf W(L/K)\to \mathsf W(F/K)$ et~$\mathsf I(L/K)\to \mathsf I(F/K)$
sont surjectives.} 

\medskip
{\em Démonstration.}
Soit~$g\in \mathsf {Gal}(F/K)$ agissant
trivialement sur~$\red F$ (resp.~$\red F_1$). Relevons~$g$ 
en un élément~$h$ de~$\mathsf{Gal}(L/K)$. 
La restriction induisant une surjection 
de~$\mathsf {Gal}(L/F)$ vers le groupe de
Galois de~$\red L$ sur~$\red F$ (resp.
de~$\red L_1$ sur~$\red F_1$), il existe~$h'\in \mathsf{Gal}(L/F)$
tel que~$hh'$ agisse trivialement sur~$\red F$ (resp.~$\red F_1$). 

\medskip
Par construction,~$hh'$ est un relevé de~$g$ à~$\mathsf W(L/K)$
(resp.~$\mathsf I(L/K)$).~$\Box$

\deux{rampsylow} {\bf Lemme.} {\em Soit~$K$ un corps valué hensélien, soit~$L$
une extension galoisienne de~$K$, et soit~$p$ 
l'exposant caractéristique de~$\red K$. 
Le groupe~$\mathsf W(L/K)$ est trivial si~$p=1$,
et est l'unique pro-$p$-sous-groupe de Sylow de~$\mathsf I(L/K)$
sinon.}

\medskip
{\em Démonstration.}
Comme
le groupe quotient $\mathsf I(L/K)/\mathsf W(L/K)$
est isomorphe à~$(|L\ti|/|K\ti|)_{p'}^\vee$, il est d'ordre
généralisé premier à~$p$. Puisque~$\mathsf W(L/K)$ est distingué
dans~$\mathsf I(L/K)$, il suffit de démontrer que c'est un
pro-$p$-groupe si~$p\neq 1$, et qu'il est trivial sinon. 

\medskip
Le lemme~\ref{ietwcompres} ci-dessus permet de supposer que~$L$ est
une extension finie de~$K$ (le cas général s'en déduira par passage
la limite). Il suffit de
démontrer que tout élément d'ordre
premier à~$p$ de~$\mathsf W(L/K)$ est trivial. 
Soit $g$ un élément de $\mathsf W(L/K)$ dont l'ordre $q$ est premier à~$p$
et soit~$E$ le sous-corps de~$L$ formé des éléments invariants sous~$g$.
Soit $x$ appartenant à~$L$ ; l'élément~$y=x-\mathsf{Tr}_{L/E}(x)/q$
de~$L$ est de trace nulle sur~$E$. Autrement dit, on a~$\sum g^{i}(y)=0$, 
soit encore~$qy-\sum (g^{i}(y)-y)=0$. Comme~$q$ est premier à~$p$, on a $|qy|=|y|$ ;
et  puisque~$g$
 appartient à~$\mathsf W(L/k)$, on 
 a~$|g^{i}(y)-y|<|y|$ pour tout~$i$ si $y\neq 0$ ;
 par conséquent $y=0$ et $x\in E$. Ainsi, $E=L$, et  $g=\mathsf{Id}$.~$\Box$ 
 
 \deux{defsetd} Soit~$K$ un corps valué hensélien
 et soit~$L$ une extension algébrique de~$K$. 
 On note~$\sch S_f(L/K)$ l'ensemble des sous-extensions
 finies et séparables de~$\red L$, 
 et~$\sch D_f(L/K)$ l'ensemble des sous-extensions finies
 et sans défaut de~$L$ dont le corpoïde
 résiduel est séparable sur~$\red K$. 
 
 \medskip
 Si~$E$ est une sous-extension
 finie de~$L$ 
 et si~$F$ est une sous-extension
 de~$E$
 alors~$E\in \sch D_f(L/K)$ 
 si et seulement si~$F\in \sch D_f(L/K)$ et~$E\in \sch D_f(L/F)$
 (\ref{henselsansdef}).

\deux{uniquextsep} {\bf Lemme.} {\em La flèche~$E\mapsto \red E$
établit une bijection 
d'ensembles ordonnées
entre~$\sch D_f(L/K)$ et~$\sch S_f(L/K)$.} 

\medskip
 {\em Démonstration.} La flèche~$E\mapsto \red E$
 est clairement croissante. 
 Montrons maintenant
 qu'elle induit une bijection de~$\sch D_f(L/K)$
 vers~$\sch S_f(L/K)$. 
 Soit~$\Lambda$
 appartenant à
~$ \sch S_f(L/K)$ ; nous allons
 prouver qu'elle possède un unique antécédent
 dans~$\sch D_f(L/K)$. On pose~$d=[\Lambda:\red K]$, et
 l'on raisonne par récurrence sur~$d$. Si~$d=1$
 on a~$\Lambda=\red K$ et l'assertion
 à prouver est immédiate. 
Supposons~$d$
strictement supérieur à~$1$
et la propriété vraie pour les entiers~$<d$. 
 
 \medskip
 Comme~$d>1$ il existe~$x\in \Lambda\setminus \red K$. 
 Soit~$P$ son polynôme minimal, soit~$\delta$
 le degré de~$P$ (il est~$>1$
 par choix de~$x$)
 et soit~$\cal P$
 un relèvement admissible de~$P$ dans~$K[T]$. 
 Comme~$L$ est hensélien et
 comme~$P$ est séparable, 
 il existe une unique racine~$\xi$
 de~$\cal P$ dans~$L$ ; posons~$F=K[\xi]$. 
 
 Le polynôme~$\mathcal P$ étant de degré~$\delta$,
 on a~$[F:K]\leq \delta$. Par ailleurs,
 $$[F:K]\geq [\red F:\red K]\geq [\red K[x]:\red K]=\delta.$$ 
 Par conséquent, $[F:K]=\delta$ et~$\red F=\red K[x]$ ; ainsi,~$F\in \sch D_f(L/K)$. 
 
 \medskip
 Si~$E$ est une sous-extension 
 de~$L$ telle que~$x\in \red E$ le caractère
 hensélien de~$E$ assure qu'il
 existe une unique racine de~$\cal P$ dans~$E$
 relevant~$x$, qui coïncide nécessairement avec~$\xi$ ; 
 ainsi,~$F\subset E$. Puisque~$F\in \sch D_f(L/K)$, le
 corps~$E$ appartient à~$\sch D_f(L/K)$ si et seulement si il
 appartient à~$\sch D_f(L/F)$. 
 
 \medskip
 Comme~$\delta>1$, on a~$[\Lambda:\red K[x]]<d$. L'hypothèse
 de récurrence assure alors qu'il existe une unique~$E\in \sch D_f(L/F)$
 telle que~$\red E=\Lambda$ ; par ce qui 
 précède,~$E$ est l'unique extension appartenant
 à~$\sch D_f(L/K)$ telle que~$\red E=\Lambda$. Ainsi,~$\Lambda$
 a bien un unique antécédent dans~$\sch D_f(L/K)$, comme annoncé. 
 
 \medskip
 Soit maintenant~$\Lambda$ et~$\Lambda'$ deux extensions appartenant
 à~$\sch S_f(L/K)$ telle que~$\Lambda\subset \Lambda'$, et soient~$E$
 et~$E'$ leurs antécédents respectifs dans~$\sch D_f(L/K)$. Comme~$\Lambda\in \sch S_f(E'/K)$,
 elle possède un unique antécédent dans~$\sch D_f(E'/K)$ ; 
 par unicité, celui-ci coïncide 
 nécessairement avec~$E$ et l'on a en particulier~$E\subset E'$, 
 ce qui achève la démonstration.~$\Box$
 
 \deux{dssysfiltr1} Comme~$\sch S_f(L/K)$ est un ensemble ordonné filtrant, il
 en va de même en vertu du lemme~\ref{uniquextsep}
 ci-dessus de l'ensemble ordonné~$\sch D_f(L/K)$. 
 
Soit~$(F_i)$ une famille finie de sous-extensions de~$L$
appartenant à~$\sch D_f(L/K)$. Par ce qui précède, il existe
une extension~$F\in \sch D_f(L/K)$ contenant les~$F_i$ ; il
s'ensuit 
que la composée des~$F_i$ appartient à~$\sch D(L/K)$. 

 \deux{dssysfiltr2} Soit~$E$ une sous-extension de~$L$. 
 Il résulte de~\ref{dssysfiltr1}
 que les assertions suivantes sont équivalentes :

 \medskip
 i) toute sous-extension finie de~$E$ appartient à~$\sch D_f(L/K)$ ; 
 
 ii) $E$ est réunion d'extensions finies appartenant
 à~$\sch D_f(L/K)$. 
 
 \medskip
 On note~$\sch D(L/K)$ l'ensemble des sous-extensions
 de~$L$ satisfaisant ces conditions. Si~$E\in \sch D(L/K)$, il
 en va de même de toute sous-extension de~$E$. 
 
 \medskip
Soit~$\sch S(L/K)$ l'ensemble de toutes les
 sous-extensions
 séparables 
 de~$\red L$. Soit~$\Lambda\in \sch S(L/K)$. 
 Le
 sous-ensemble
 de~$\sch D_f(L/K)$ constitué des
 corps~$E$ 
 tels que~$\red E\subset \Lambda$ est filtrant ; la
 réunion
 des extensions appartenant audit sous-ensemble
 est un élément~$\Lambda^\diamond$
 de~$\sch D(L/K)$.  

\deux{uniquextsepgen} {\bf Lemme.} {\em Les flèches~$E\mapsto \red E$
et~$\Lambda\mapsto \Lambda^\diamond$ établissent une bijection
entre~$\sch D(L/K)$ et~$\sch S(L/K)$.}

\medskip
{\em Démonstration.} Si~$\Lambda\in \sch S(L/K)$,
il résulte immédiatement de la définition
que~$\red{ \Lambda^\diamond}=\Lambda$.

\medskip
Soit~$E\in \sch D(L/K)$. Si~$F$
est une sous-extension finie de~$E$ alors
~$F\in \sch D_f(L/K)$ et~$\red F\subset \red E$ ; 
en conséquence~$E\subset \red E^\diamond$. 

Inversement, 
soit~$F\in \sch D_f(L/K)$ 
telle que~$\red F\subset \red E$. Choisissons une partie génératrice
finie~$S$ de~$\red F$ sur~$\red K$. Pour tout~$\xi\in S$, 
il existe une sous-extension~$E_\xi$ de~$E$, appartenant
à~$\sch D(L/K)_f$ et telle que~$\xi \in \red {E_\xi}$. 
La composée des~$E_\xi$ est contenue dans~$E$,
et appartient à~$\sch D_f(L/K)$ ; par construction,
son corpoïde résiduel contient~$\red F$. On déduit alors
du lemme~\ref{uniquextsep} qu'elle contient~$F$, d'où
l'inclusion~$F\subset E$. Il vient~$\red E^\diamond \subset E$.~$\Box$

 \deux{plusgrandslk} L'ensemble~$\sch S(L/K)$
 admet un plus grand élément :
 la fermeture séparable~$\red L_{\rm sep}$
 de~$\red K$ dans~$\red L$ ; par
 conséquent,~$\sch D(L/K)$ admet un plus grand élément~$E$,
 à savoir~$\red L_{\rm sep}^\diamond$.

\deux{desmodram} {\bf Théorème.} {\em Supposons
que~$L$ est galoisienne sur~$K$. 
On a l'égalité~$L^\mathsf W(L/K)=\red L^\diamond_{\rm sep}$.}

\medskip
{\em Démonstration.}
Le lemme~\ref{ietwcompres} permet
de se ramener au cas où~$L$ est finie sur~$K$. 

Par définition 
de~$\mathsf W(L/K)$, le groupe~$\mathsf {Gal}(F/K)$
s'identifie naturellement à~$\mathsf{Gal}(\red L/\red K)$ ; cela signifie
que les surjections
$$\mathsf{Gal}(F/K)\to \mathsf{Gal}(\red F/\red K)\;\;
{\rm et}\; \;\mathsf{Gal}(\red L/\red K)\to \mathsf{Gal}(\red F/\red K)$$
sont bijectives. La bijectivité de la seconde
implique que~$\red L$ est purement inséparable sur~$\red F$. 

\medskip
On a par ailleurs
$$\card (\mathsf{Gal}(F/K))=[F:K]\geq [\red F:\red K]\geq \card(\mathsf{Gal}(\red F/\red K))
= \card(\mathsf{Gal}(F/K)).$$

On en déduit que~$[\red F:\red K]=[F:K]$ et que~$\red F$ est galoisienne 
sur~$\red K$. Il s'ensuit que~$F\in \sch D(L/K)$, 
et que~$\red F=\red L_{\rm sep}$. 
En conséquence,~$F=\red L_{\rm sep}^\diamond$.~$\Box$ 

\subsection*{Extensions modérément ramifiées}

\deux{defmodram} Soit~$K$ un corps
valué hensélien, soit~$p$ son
exposant caractéristique résiduel, 
et soit~$F$
une extension algébrique de~$K$. 
On dit que~$F$ est
{\em modérément ramifiée}
si~$F\in \sch D(F/K)$, 
c'est-à-dire si~$\red F$ 
est séparable sur~$\red K$
et si~$F=\red F^\diamond$ ; 
il revient
au même de demander que~$\red F$ soit
séparable sur~$\red K$ et 
que toute sous-extension finie de~$F$ soit
sans défaut. Notons qu'en vertu de~\ref{equivseparoide}, 
la séparabilité de~$\red F$
sur~$\red K$
peut elle-même se reformuler comme suit, dans un 
langage plus classique :~$\red F_1$ est séparable sur~$\red K_1$,
et~$F\ti|/|K\ti|$ est 
sans~$p$-torsion.
Si~$F$ est modérément ramifiée, il en va de même de chacune
de ses sous-extensions. 

\trois{modramfinitriv} Si~$F$ est finie sur~$K$, elle est modérément ramifiée
si et seulement si
elle est sans défaut et si~$\red F$ séparable sur~$\red K$. 

\trois{modrambij} Supposons que~$F$ soit modérément ramifiée. 
Toute sous-extension de~$F$ appartient alors
à~$\sch D(F/K)$, et toute sous-extension
de~$\red F$ appartient à~$\sch S(F/K)$.
On déduit dès lors du lemme~\ref{uniquextsepgen}
que l'application~$E\mapsto \red E$ établit une bijection
entre l'ensemble des sous-extensions de~$F$ et celui des
sous-extensions de~$\red F$.

\trois{modramwfixe}
Soit~$L$ une extension algébrique
de~$F$ qui est galoisienne sur~$K$. On déduit du
théorème~\ref{ietwcompres}
que~$F$ est modérément ramifiée
si et seulement si~$F\subset L^{\mathsf W(L/K)}$. 

\trois{modramgal}  Supposons que~$F$ soit galoisienne. 
D'après le~\ref{modramwfixe} ci-dessus, elle est modérément
ramifiée si et seulement si~$\mathsf W(F/K)$ est trivial. Plaçons-nous
sous cette hypothèse. 
On dispose alors par définition du groupe~$\mathsf W(F/K)$ 
d'un isomorphisme~$\mathsf{Gal}(F/K)\simeq \mathsf {Gal}(\red F/\red K)$. 

Si~$g\in \mathsf {Gal}(F/K)$ et si~$E$ est une sous-extension
de~$F$, on a~$\red{g(E)}=\red g(\red E)$, où~$\red g$
désigne l'image de~$g$ dans~$\mathsf {Gal}(\red F/\red K)$. Il s'ensuit
que la bijection~$E\mapsto \red E$ de~\ref{modrambij}
coïncide avec
celle déduite de 
l'isomorphisme~$\mathsf{Gal}(F/K)\simeq \mathsf {Gal}(\red F/\red K)$
et de la correspondance de Galois. 

\medskip
Le groupe~$\mathsf W(F/K)$
étant trivial, le groupe~$\mathsf I(F/K)$
est naturellement isomorphe
à~$\mathsf {Hom}(|F\ti|/|K\ti|),\mu(\red F_1))$  
(\ref{quotinertram}).
L'indice~$(|E\ti|:|K\ti|)$ est de
surcroît
premier
à~$p$ pour toute sous-extension finie~$E$
de~$F$ ({\em cf.} \ref{modramfinitriv}). En
conséquence,
il résulte de~\ref{assezracinetilde}
que
tout plongement de~$\mu(\red F_1)$ dans~$\QQ/\ZZ$
induit un isomorphisme~$\mathsf I(F/K)\simeq (|F\ti|/|K\ti|)^\vee$.

\deux{algetalemodram} Soit~$A$ une~$K$-algèbre étale ;
écrivons~$A=\prod K_i$, où~$(K_i)$ est une famille finie d'extensions
séparables de~$K$. On munit l'algèbre~$A$ de la
semi-norme~$(x_i)\mapsto \sup |x_i|$ ; notons
que~$\red A\simeq \prod \red K_i$. On dit que~$A$ est modérément
ramifiée si c'est le cas de chacune des~$K_i$. 

\medskip
Fixons une clôture séparable~$K^s$ de~$K$. La catégorie
des~$K$-algèbres finies étales s'identifie à celle
des~$\mathsf{Gal}(K^s/K)$-ensembles finis et 
discrets. Si~$A$ est une~$K$-algèbre étale,
on déduit de~\ref{modramwfixe}
que~$A$ est modérément ramifiée si et
seulement si~$\mathsf W(K^s/K)$ agit trivialement
sur le~$\mathsf{Gal}(K^s/K)$-ensemble correspondant ; 
on dispose donc d'une équivalence 
entre la catégorie des~$K$-algèbres finies étales modérément
ramifiées et celle des~$\mathsf {Gal}(\red{K^s}/\red K)$-ensembles
finis discrets. Comme~$\red K^s$ est séparablement clos en vertu
du lemme~\ref{lemalgclossepclosoide}, cette dernière catégorie est elle-même
équivalente à celle des~$\red K$-algèbres finies étales. 

\medskip
On a ainsi construit une équivalence 
entre la catégorie des~$K$-algèbres étales modérément ramifiées et celle
des~$\red K$-algèbres étales. On déduit du~\ref{modramgal}
que cette équivalence est donnée par la flèche~$A\mapsto \red A$.

\deux{prempmodram} {\bf Lemme}. {\em Soit~$K$ un corps valué hensélien et soit~$E$
une extension finie de~$K$ de degré premier à l'exposant caractéristique~$p$
de~$\red K$. L'extension~$E$ est alors modérément ramifiée.} 

\medskip
{\em Démonstration.} Soit~$L$ une extension 
finie de~$E$ galoisienne sur~$K$. Si~$p=1$ le groupe~$\mathsf W(L/K)$
est trivial (lemme~\ref{rampsylow}), et~$E$ est modérément ramifiée. 

Supposons~$p>1$,
et soit~$\mathsf H$ un~$p$-sous-groupe de Sylow de~$\mathsf{Gal}(L/E)$. 
Comme~$[E:K]$ est premier à~$p$, le groupe~$\mathsf H$ 
est aussi un~$p$-sous-groupe de Sylow de~$\mathsf{Gal}(L/K)$. 
Puisque le sous-groupe distingué~$\mathsf W(L/K)$ 
de~$\mathsf {Gal}(L/K)$ est un~$p$-groupe
d'après le lemme~\ref{rampsylow}, il est contenu dans~$\mathsf H$,
et~{\em a fortiori} dans~$\mathsf{Gal}(L/E)$, ce
qui termine la preuve.~$\Box$ 

\subsection*{Extensions non ramifiées} 

\deux{defnonram} Soit~$K$ un corps
valué hensélien, soit~$p$ son
exposant caractéristique résiduel, 
et soit~$F$
une extension algébrique de~$K$. 
On dit que~$F$ est
{\em non ramifiée}
si~$F$ est modérément
ramifiée et si~$|F\ti|=|K\ti|$. 
Si~$F$ est non
ramifiée, il en va de même de chacune
de ses sous-extensions. 

\trois{nonramfinitriv} Si~$F$ est finie sur~$K$, elle est non
ramifiée
si et seulement si~$\red F_1$ est une extension
séparable sur~$\red K_1$, de degré~$[F:K]$ : cela provient
de l'égalité
$$[\red F:\red K]=(|F\ti|:|K\ti|)\cdot[\red F_1:\red K_1]$$
et de~\ref{equivseparoide}.  

\trois{nonrambij} Supposons que~$F$ soit
non ramifiée. 
On déduit de~\ref{modrambij}
et de~\ref{casdegldegk} que~$E\mapsto \red E_1$ établit une bijection
entre l'ensemble des sous-extensions de~$F$ et celui des
sous-extensions de~$\red F_1$.

\trois{nonramwfixe}
Soit~$L$ une extension algébrique
de~$F$ qui est galoisienne sur~$K$. On déduit du
théorème~\ref{ietwcompres}
que~$F$ est modérément ramifiée
si et seulement si~$F\subset L^{\mathsf I(L/K)}$. 

\trois{nonramgal}  Supposons que~$F$ soit galoisienne. 
D'après le~\ref{modramwfixe} ci-dessus, elle est modérément
ramifiée si et seulement si~$\mathsf I(F/K)$ est trivial. Plaçons-nous
sous cette hypothèse. 
On dispose alors par définition du groupe~$\mathsf I(F/K)$ 
d'un isomorphisme~$\mathsf{Gal}(F/K)\simeq \mathsf {Gal}(\red F_1/\red K_1)$. 

Si~$g\in \mathsf {Gal}(F/K)$ et si~$E$ est une sous-extension
de~$F$, on a~$\red{g(E)_1}=\red g(\red E_1)$, où~$\red g$
désigne l'image de~$g$ dans~$\mathsf {Gal}(\red F_1/\red K_1)$. Il s'ensuit
que la bijection~$E\mapsto \red E_1$ de~\ref{nonrambij}
coïncide avec
celle déduite de 
l'isomorphisme~$\mathsf{Gal}(F/K)\simeq \mathsf {Gal}(\red F_1/\red K_1)$
et de la correspondance de Galois. 

\deux{algetalenonram} Soit~$A$ une~$K$-algèbre étale ;
écrivons~$A=\prod K_i$, où~$(K_i)$ est une famille finie d'extensions
séparables de~$K$. On munit l'algèbre~$A$ de la
semi-norme~$(x_i)\mapsto \sup |x_i|$ ; notons
que~$\red A_1\simeq \prod \red {K_{i,1}}$. On dit que~$A$ est non
ramifiée si c'est le cas de chacune des~$K_i$. 

\medskip
Fixons une clôture séparable~$K^s$ de~$K$. La catégorie
des~$K$-algèbres finies étales s'identifie à celle
des~$\mathsf{Gal}(K^s/K)$-ensembles finis et 
discrets. Si~$A$ est une~$K$-algèbre étale,
on déduit de~\ref{nonramwfixe}
que~$A$ est modérément ramifiée si et
seulement si~$\mathsf I(K^s/K)$ agit trivialement
sur le~$\mathsf{Gal}(K^s/K)$-ensemble correspondant ; 
on dispose donc d'une équivalence 
entre la catégorie des~$K$-algèbres finies étales modérément
ramifiées et celle des~$\mathsf {Gal}(\red{K^s_1}/\red K_1)$-ensembles
finis discrets. Comme~$\red K^s_1$ est séparablement clos en vertu
du lemme~\ref{lemalgclossepclosoide}, cette dernière catégorie est elle-même
équivalente à celle des~$\red K_1$-algèbres finies étales. 

\medskip
On a ainsi construit une équivalence 
entre la catégorie des~$K$-algèbres étales non
ramifiées et celle
des~$\red K_1$-algèbres étales. On déduit du~\ref{modramgal}
que cette équivalence est donnée par la flèche~$A\mapsto \red A_1$.

\deux{nonramalgcoval} {\bf Proposition.} {\em Soit~$A$ une~$K$-algèbre finie étale ;
notons~$A\zero$ la~$K\zero$-algèbre~$\{a\in A, ||a||\leq 1\}$. Les assertions suivantes
sont équivalentes : 

\medskip
i) $A$ est non ramifiée ; 

ii) $A\zero$ est finie étale sur~$K\zero$ ; 

iii) il existe une~$K\zero$-algèbre finie étale~$B$ telle
que~$A\simeq B\otimes_{K\zero} K$. 

\medskip
Si de plus si elles sont satisfaites alors la~$K\zero$-algèbre~$B$ de~iii)
est nécessairement égale à~$A\zero$, et~$\red A_1=A\zero\otimes_{K\zero} \red K_1$.}

\medskip
{\em Démonstration.} On peut raisonner facteur par facteur, et partant 
supposer que~$A$ est un corps. 

\medskip
{\em Supposons que~i) est vraie.} Choisissons un
élément~$\xi$ de~$\red A_1$
qui l'engendre sur~$\red K_1$, et soit~$P$ un relevé unitaire
dans~$K\zero[T]$ du polynôme minimal 
de~$\xi$ sur~$\red K_1$. Le corps~$A$ étant hensélien, le polynôme~$P$
possède une racine~$x$ dans~$A\zero$ située au-dessus de~$\xi$. 
Comme~$\red {K[x]}_1$ contient~$\xi$, son degré sur~$\red K_1$ est
au moins égal à~$[A:K]$, et il est par ailleurs majoré par~$[K[x]:K]$ ; 
il vient~$A=K[x]$. En particulier,~$x$ est de degré~$[A:K]=\deg P$
sur~$K$, et~$P$ est donc son polynôme minimal. Puisqu'il est unitaire,~$K\zero[x]
\simeq K\zero[T]/P$ ; le polynôme~$P$ se réduisant
par construction en un polynôme séparable de~$\red K_1[T]$, 
la~$K\zero$-algèbre~$K\zero[x]$ est finie étale. Elle est dès lors normale,
et coïncide de ce fait avec la fermeture intégrale de~$K\zero$ dans~$A$, qui n'est
autre que~$A\zero$, d'où~ii). 

\medskip
L'implication~ii)$\Rightarrow$~iii) est évidente. 

\medskip
{\em Supposons que~iii)
est vraie}.  Comme~$B$ est plate, elle est
sans~$K\zero$-torsion et se plonge donc dans~$A$ ; 
c'est une sous-$K\zero$-algèbre de~$A$ qui est entière sur~$K\zero$,
et normale car étale sur~$K\zero$. Elle
coïncide de ce fait avec la fermeture intégrale de~$K\zero$ dans~$A$, qui n'est
autre que~$A\zero$. 

\medskip
L'anneau~$K\zero$ étant hensélien, la~$K\zero$-algèbre finie et intègre~$A\zero$
est locale. La~$\red K_1$-algèbre~$A\zero\otimes_{K\zero}\red K_1$ étant à la fois
locale et étale, c'est un corps ; cela signifie que~$K\zeroo A\zero$ est l'idéal 
maximal de~$A\zero$, et~$A\zero\otimes_{K\zero}\red K_1$ s'identifie
en conséquence
au corps résiduel~$\red A_1$. Celui-ci est donc séparable et
de degré~$[A:K]$ sur~$\red K_1$, d'où~i).~$\Box$

\subsection*{Les~$\mu_\ell$-torseurs sur un corps valué hensélien}

\deux{muellehens} Soit~$K$ un corps valué hensélien, et soit~$\ell$ un
entier premier à l'exposant caractéristique de~$\red K$.

\trois{hilb90hensel} Le théorème 90 de Hilbert fournit trois
isomorphismes naturels
$$\H^1(K,\mu_\ell)\simeq K\ti/(K\ti)^\ell,\;\;\H^1(K\zero,\mu_\ell)\simeq K^{\circ, \times}/(K^{\circ, \times})^\ell$$
$${\rm et}
\;\H^1(\red K_1,\mu_\ell)\simeq \red K_1\ti/( \red K_1\ti)^\ell.$$

Il en résulte immédiatement, compte-tenu
du caractère hensélien de~$K\zero$, 
que la flèche naturelle~$\H^1(K\zero,\mu_\ell)\to \H^1(\red K_1,\mu_\ell)$
est un isomorphisme
(ce que l'on pouvait également déduire de~\ref{eqetalehens2}).

\trois{absksanstors} Comme~$|K\ti|$ est sans torsion, la suite exacte

$$1\to K^{\circ,  \times} \to K\ti \to |K\ti|\to 1$$ 
induit une suite exacte
$$1\to K^{\circ, \times}/(K^{\circ, \times})^\ell\to K\ti/(K\ti)^\ell\to |K\ti|/|K\ti|^\ell \to 1$$
que l'on peut récrire
$$\diagram 0\rto& \H^1(K\zero,\mu_\ell)\rto\dto_{\simeq}& \H^1(K,\mu_\ell)\rto^\partial &|K\ti|/|K\ti|^\ell\rto &1
\\&H^1(\red K_1,\mu_\ell)&&&\enddiagram.$$

La flèche~$\partial : \H^1(K,\mu_\ell)\to |K\ti|/|K\ti|^\ell$ est appelée {\em résidu}. Si~$h\in \H^1(K,\mu_\ell)$
alors~$\partial (h)=0$ si et seulement si~$h$ provient de~$\H^1(K\zero,\mu_\ell)$ ; cela revient à demander, 
en vertu de~\ref{applitorseurs} et de la proposition~\ref{algetalenonram}, que la~$K$-algèbre
étale définie par~$h$ soit non ramifiée. 

\trois{ktilde1sepclos} Supposons que~$\red K_1\ti$ soit~$\ell$-divisible, ce qui est notamment le
cas lorsque~$\red K_1$ est séparablement clos. On a alors~$\H^1(\red K_1,\mu_\ell)=0$ ; 
par ce qui précède, le résidu~$\partial$ induit un 
isomorphisme~$\H^1(K,\mu_\ell)\simeq |K\ti|/|K\ti|^\ell$. 

\trois{abskldiv} Supposons que~$|K\ti|$ soit~$\ell$-divisible. 
On a alors par ce qui précède un diagramme d'isomorphismes  
$$\diagram \H^1(\red K_1,\mu_\ell)&\H^1(K\zero,\mu_\ell)\lto_\simeq\rto^\simeq &\H^1(K,\mu_\ell)\enddiagram.$$

\subsection*{Extensions modérément ramifiées : le cas d'un corps résiduel séparablement clos}

\deux{khenssepclos} Soit~$K$ un corps valué hensélien dont le corps résiduel~$\red K_1$ est 
séparablement clos, soit~$K^s$ une clôture séparable de~$K$, et soit~$L$
le sous-corps~$(K^s)^{\mathsf W(K^s/K)}$
de~$K^s$. On note~$p$ l'exposant caractéristique de~$\red K$. 

\medskip
Comme~$\red K_1$ est séparablement clos, le groupe~$\mathsf I(K^s/K)$ est égal
à~$\mathsf{Gal}(K^s/K)$ tout entier, et~$\mu (\red {K^s_1})=\mu(\red K_1)$. Il s'ensuit
que~$\abs x \mapsto \red {g(x)/x}$ induit un isomorphisme
entre~$\mathsf {Gal}(L/K)$ et~$\mathsf{Hom}(|L\ti|/|K\ti|, \mu(\red K_1))$. 

\medskip
Soit~$F$ une sous-extension de~$L$. Le groupe de Galois de~$L$
sur~$K$ étant abélien, l'extension~$F$ est galoisienne. Elle est
également modérément ramifiée ; en conséquence,~$\abs x \mapsto
 \red {g(x)/x}$ induit un isomorphisme
entre~$\mathsf {Gal}(F/K)$ et~$\mathsf{Hom}(|F\ti|/|K\ti|, \mu(\red K_1))$. 

\medskip
La surjection~$\mathsf {Gal}(L/K)\to \mathsf{Gal}(F/K)$ correspond,
modulo les isomorphismes évoqués ci-dessus, à la flèche de restriction
de~$\mathsf{Hom}(|L\ti|/|K\ti|, \mu(\red K_1))$
vers~$\mathsf{Hom}(|F\ti|/|K\ti|, \mu(\red K_1))$. Le groupe~$\mathsf{Gal}(L/F)$
est donc égal, modulo l'identification entre~$\mathsf {Gal}(L/K)$
et~$\mathsf{Hom}(|L\ti|/|K\ti|, \mu(\red K_1))$, à l'orthogonal de~$|F\ti|/|K\ti|$. 

\deux{pontrjabijplong} Comme~$L$ est modérément 
ramifiée sur~$K$, le groupe~$|L\ti|/|K\ti|$ est sans~$p$-torsion. 
En combinant 
correspondance de Galois et dualité de Pontrjagin (comme~$\red F_1$ est 
séparablement clos,~$\mu(\red F_1)$ est isomorphe à la partie 
de torsion première à~$p$ de~$\QQ/\ZZ$, et 
l'on n'a pas besoin ici d'invoquer~\ref{assezracinepasfini}),
 on voit que pour tout sous-groupe~$\Delta$
de~$|L\ti|/|K\ti|$, il existe une et une seule sous-extension~$F$ de~$L$
telle que~$|L\ti|/|K\ti|=\Delta$.

\deux{pontrjabij} On en déduit que l'application~$F\mapsto |F\ti|/|K\ti|$
induit une bijection entre l'ensemble des classes d'isomorphie
d'extensions modérément
ramifiées de~$K$ et celui des sous-groupes de~$|L\ti|/|K\ti|$ 
(on considère~$|F\ti|$ et~$|L\ti|$ comme canoniquement
plongés dans~$|K\ti|^\QQ$). 

\deux{descexpl} On peut donner une description précise
du sous-groupe~$|L\ti|/|K\ti|$ de~$|K\ti|^\QQ/|K\ti|$ : c'est précisément
sa partie de torsion première à~$p$. En effet, on sait déjà que~$|L\ti|/|K\ti|$
est de torsion première à~$p$. Il suffit donc de montrer que si~$\Delta$ est un 
sous-groupe fini de~$|K\ti|^\QQ/|K\ti|$ qui est de torsion première à~$p$, 
il existe une extension modérément ramifiée~$F$ de~$K$ telle
que~$|F\ti|/|K\ti|=\Delta$. C'est ce que nous allons faire, en 
décrivant {\em explicitement}
une telle~$F$. 

\medskip
Choisissons une famille~$(d_1,\ldots, d_n)$ d'éléments
de~$|K\ti|^\QQ\setminus |K\ti|$
tels que~$\Delta$
soit somme directe des~$\overline {d_i}^\ZZ$. 
Pour tout~$i$, on note~$e_i$ l'ordre de~$\overline {d_i}$, et l'on choisit un
élément~$x_i$ de~$K$ tel que~$d_i^{e_i}=\abs {x_i}$. 

\medskip
Soit~$A$ la~$K$-algèbre finie~$K[T_1,\ldots, T_n]/(T_i^{e_i}-x_i)$. Sa dimension
est égale au produit des~$e_i$, c'est-à-dire au cardinal de~$\Delta$. Fixons un quotient~$F$
de~$A$ par un idéal maximal. Par construction,~$x_i$ a pour tout~$i$ une racine~$e_i$-ième
dans~$F$ ; il s'ensuit que
$$|F\ti|/|K\ti|\supset \langle\overline{ \abs {x_1}}^{1/e_i},\ldots,\overline{\abs{x_n}}^{1/e_n}\rangle=\langle \overline {d_1},\ldots, \overline {d_n}\rangle=\Delta.$$

On a ainsi~$$\card\; \Delta =[A:K]\geq [F:K]\geq (|F\ti|:|K\ti|)\geq \card\; \Delta.$$ Par
conséquent,~$F=A, [F:K]=\card \;\Delta$ et~$|F\ti|/|K\ti|=\;\Delta$. Notons que
l'égalité~$[F:K]=(|F\ti|:|K\ti|)$
implique que~$\red F_1=\red K_1$ (et que~$F$ est sans défaut).

\medskip
Il suffit pour conclure de montrer
que~$F$ est modérément ramifiée. 
Or comme le groupe~$\Delta$ est de cardinal premier à~$p$,
le quotient~$|F\ti|/|K\ti|$
est de cardinal premier à~$p$ ; on a vu par ailleurs
que~$\red F_1=\red K_1$, et que~$F$ est sans défaut. Elle est donc
modérément ramifiée.

\chapter{Géométrie analytique}\markboth{Géométrie analytique}{Géométrie analytique}

\section{Espaces de Berkovich : rappels et notations} 

\subsection*{Quelques conventions et notations}
On fixe {\em pour toute la suite du texte} un corps ultramétrique complet~$k$ (sa valeur absolue peut être triviale). On se donne une clôture séparable~$k^{\rm s}$ de~$k$, et l'on note~$\KK$ le complété de la clôture algébrique~$k^a$ de~$k^s$ ; on note par ailleurs~$\kparf$ le complété de la clôture parfaite~$k^{\rm parf}$ de~$k$. Si~$|k\ti|\neq\{1\}$ alors~$k^s$ est dense dans~$\KK$. Le corps résiduel de~$\KK$ est une clôture algébrique~$\kk_1$ de~$\red k_1$ ; on note par ailleurs~$p$ l'exposant caractéristique de~$\red k_1$. On pose~$\mathsf G={\mathsf{Gal}}\;(k^s/k)$. On appellera extension {\em presque algébrique} de~$k$ toute extension complète de~$k$ admettant un~$k$-plongement isométrique dans~$\KK$. 

\medskip
On appellera {\em polyrayon} une famille finie de réels strictement positifs ; si cette famille est libre lorsqu'on la voit dans~$\QQ\otimes_{\ZZ}(\RR\ti_+/|k\ti|)$, on dira que le polyrayon en question est {\em~$k$-libre.} Si~${\bf r}=(r_1,\ldots,r_n)$ est un polyrayon~$k$-libre, la~$k$-algèbre affinoïde~$k\{T_1/r_1,\ldots,T_n/r_n, r_1T_1\inv, \ldots, r_nT_n\inv\}$ est un corps que l'on notera~$k_{\bf r}$.

\deux{introb} Nous travaillerons avec la notion d'espace analytique {\em au sens de Berkovich}, en ayant soin de n'utiliser que les résultats de la théorie qui sont prouvés {\em sans recourir à l'existence de modèles formels ou algébriques particuliers}. Un espace de Berkovich est localement compact et localement connexe par arcs. Si~$\bnd X$ est une~$k$-variété algébrique, c'est-à-dire un~$k$-schéma séparé de type fini, l'on notera~$\bnd X\an$ son analytification. Si~$a\in k$ et si~$r\in \RR_+$, l'on notera~$\eta_{a,r}$ le point de~$\Aff^{1,\rm an}_k$ défini par la semi-norme~$\sum a_i(T-a)^i\mapsto \max |a_i|r^i$ ; lorsque~$a=0$, l'on écrira~$\eta_r$. Si l'on travaille sur une extension complète~$L$ de~$k$ on utilisera les notations~$\eta_{a,r,L}$ et~$\eta_{r,L}$.

\trois{changeb} Si~$X$ (resp.~$\sch A$) est un espace~$k$-analytique (resp. une~$k$-algèbre de Banach) et si~$L$ est une extension complète de~$k$, l'on notera~$X_L$ (resp.~$\sch A_L$) l'espace~$L$-analytique~$X\hotimes_k L$ (resp. la~$L$-algèbre de Banach~$\sch A\hotimes_kL$) ; si~$\bf r$ est un polyrayon~$k$-libre, on écrira souvent~$X_{\bf r}$ ou~$\sch A_{\bf r}$ au lieu de~$X_{k_{\bf r}}$ ou~$\sch A_{k_{\bf r}}$.

\trois{kaquot} Soit~$X$ un espace~$k$-analytique. Si~$L$ une extension complète de~$k$ l'on dispose d'une application continue naturelle~$X_L\to X$ qui est compacte et surjective. La flèche~$X_{\KK}\to X$ identifie topologiquement~$X$ au quotient~$X_{\KK}/\mathsf G$ ; il s'ensuit que~$X_L\to X$ est ouverte pour toute extension presque algébrique~$L$ de~$k$ (rem.~\ref{xsurhxsurg}). 

Si~$\bf r$ est un polyrayon~$k$-libre, la fibre de~$X_{\bf r}\to X$ en un point~$x$ de~$X$ est un espace~$\hres(x)$-affinoïde dont le bord de Shilov est un singleton~$\{\sigma(x)\}$ ; l'application~$\sigma$ est une section continue de~$X_{\bf r}\to X$, appelée sa {\em section de Shilov.} Si~$X$ est connexe alors~$X_{\bf r}$ l'est aussi. 

\trois{morphconst} On dira qu'un morphisme~$\phi : Y\to X$ entre espaces~$k$-analytiques est {\em constant} si~$Y$ est vide ou si~$\phi(Y)$ est un point rigide de~$X$ ; si~$Y$ est connexe et si~$\phi$ est G-localement constant alors~$\phi$ est constant. On dira qu'une fonction analytique sur~$Y$ est constante si le morphisme induit~$Y\to \Aff^{1,{\rm an}}_k$ est constant ; si~$Y$ est connexe (resp. connexe et quasi-compact), il revient au même de demander que l'image de~$f$ dans~$\sch O_{Y_{\rm red, G }}(Y_{\rm red})$ (resp. que~$f$) soit entière sur~$k$. 

\trois{sx} Soit~$X$ un espace~$k$-analytique. On désigne par~$\got s(X)$ l'algèbre des fonctions analytiques sur~$X$ qui annulent G-localement un polynôme {\em séparable} non nul à coefficients dans~$k$ ; les fonctions appartenant à~$\got s(X)$ sont G-localement constantes. 

Supposons de plus~$X$ connexe et non vide ; l'algèbre~$\got s(X)$ est alors une extension finie séparable de~$k$, qui est contenue dans l'algèbre des fonctions constantes sur~$X$ et coïncide avec celle-ci si~$X$ est géométriquement réduit ; l'espace ~$X_{\KK}$ est connexe si et seulement si~$\got s(X)=k$, et l'on dit dans ce cas que~$X$ est {\em géométriquement} connexe.

\trois{corpsres} Si~$X$ est un espace analytique et si~$x\in X$, on notera~$\hres(x)$ le corps résiduel complété de~$x$ ; si de plus~$X$ est bon, on notera~$\kappa(x)$ le corps résiduel de l'anneau local~$\sch O_{X,x}$. On notera~$\got s(x)$ la clôture séparable de~$k$ dans~$\hres(x)$. Pour tout domaine analytique~$V$ de~$X$ contenant~$x$, l'évaluation en~$x$
induit un plongement $\got s(V)\hookrightarrow \hres(x)$. Si~$[\got s(x)=k]$ est finie, le caractère hensélien de l'anneau local~$\sch O_{X,x}$ assure qu'il existe un voisinage analytique~$U$ de~$x$
dans~$X$ tel que~$\got s(U)\hookrightarrow \got s(x)$ soit un isomorphisme ; il est immédiat que~$\got s(V)\hookrightarrow \got s(x)$ est encore un isomorphisme
pour tout domaine analytique~$V$ de~$U$ contenant~$x$. 

\deux{dimespb} On dispose d'une bonne théorie de la dimension en théorie de Berkovich.
Si~$X$ est un espace~$k$-analytique, on a
$$\dim k X=\sup_{x\in X}\;{\rm deg. tr.} \; \red{\hres(x)}/\red k.$$

\deux{normalis} Nous utiliserons à plusieurs reprises implicitement le fait suivant : si~$X$ est un espace analytique et si~$\phi : X'\to X$ est son normalisé, alors~$X$ s'identifie topologiquement au quotient de~$X'$ par la relation dont les classes sont les fibres de~$\phi$ ; cela résulte directement de~\ref{condflechecomp}. 

\subsection*{Réduction à la Temkin : le cas affinoïde} 

\deux{temkglob} Soit~$\sch A$ une algèbre~$k$-affinoïde et soit~$X$ l'espace~$\sch M(\sch A)$. On désigne par~$\red{\sch A}$ l'annéloïde résiduel associé à~$\sch A$ {\em munie de sa semi-norme spectrale} (\ref{defredgradoide}) ; c'est une~$\red k$-algèbre de type fini, dont on note~$\red X$ le spectre, c'est-à-dire l'ensemble des idéaux premiers ;  on dispose d'une application de réduction~$\rho: X\to \red X$ qui est anticontinue. Si~$\sch A$ est strictement~$k$-affinoïde alors~$\red X$ s'identifie à la réduction usuelle de~$X$, à savoir~$\spec \red {\sch A_1}$. 

\trois{genshil} Si~$\xi$ est un point générique de~$\red X$ ({\em i.e.} le point générique d'une de ses composantes
irréductibles) alors~$\rho\inv(\xi)$ est un singleton ; le bord de Shilov de~$X$ est exactement~$\{\rho(\xi)\}$, où~$\xi$ parcourt l'ensemble des points génériques de~$\red X$. 

\trois{tubeadh} Si~$\xi$ et~$\eta$ sont deux points de~$\red X$ et si~$y\in \rho\inv(\eta)$ alors~$y\in \overline{\rho\inv(\xi)}$ si et seulement si~$\xi\in\overline {\{\eta\}}$. 

En effet, supposons que~$y\in  \overline{\rho\inv(\xi)}$. Soit~$\sch U$ l'ouvert de Zariski de~$\red X$ complémentaire du fermé~$\overline{\{\eta\}}$. Son image réciproque~$\rho\inv(\sch U)$ est un fermé de~$X$ qui ne contient pas~$y$ ; dès lors, il ne contient pas~$\rho\inv(\xi)$. Par conséquent,~$\xi\notin \sch U$ et l'on a donc~$x\in \overline{\{\eta\}}$. 

\medskip
Réciproquement, supposons que~$\xi\in \overline{\{\eta\}}$ et soit~$U$ un voisinage affinoïde de~$y$ dans~$X$ ; soit~$\sch B$ l'algèbre des fonctions analytiques sur~$U$. On désigne par~$\eta'$ l'image de~$y$ sur~$\red U$ et par~$\got p$ (resp.~$\got p'$) l'idéal premier de~$\red{ \sch A}$ (resp.~$\red{ \sch B}$) qui correspond à~$\eta$ (resp.~$\eta'$). Comme~$U\hookrightarrow X$ est par construction intérieur en~$y$, la~$\red{\sch A}/\got p$-algèbre~$\red{\sch B}/\got p'$ est finie ; en vertu du {\em going-up} gradué,~$\overline{\{\eta'\}}\to \overline{\{\eta\}}$ est surjective ; en particulier,~$\xi$ appartient à l'image de~$\overline{\{\eta'\}}\to \overline{\{\eta\}}$, et donc à celle de~$\red U\to \red X$ ; mais cela signifie qu'il existe~$x\in U$ dont l'image sur~$\red X$ est égale à~$\xi$ ; autrement dit,~$U$ rencontre~$\rho\inv(\xi)$. 

\trois{tubecon} Si~$\xi$ est  un point de~$\red X$ dont le corpoïde résiduel est fini sur~$\red k$, l'ouvert~$\rho\inv(\xi)$ est connexe.

\subsection*{Réduction à la Temkin : le cas local} 

\deux{defgerm} Nous allons commencer par rappeler la définition de la catégorie des germes (ponctuels) d'espaces~$k$-analytiques dont nous nous servirons abondammment. 

\trois{defisopointe} Un morphisme~$f: (Y,y)\to (X,x)$ d'espaces~$k$-analytiques pointés est appelé un {\em isomorphisme local} s'il induit un isomorphisme d'un voisinage ouvert de~$y$ dans~$Y$ sur un voisinage ouvert de~$x$ dans~$X$.

\trois{germlocalise} On définit la catégorie des {\em germes d'espaces~$k$-analytiques} 
comme la localisée de la catégorie des espaces analytiques pointés par la famille des isomorphismes locaux. On vérifie alors qu'un morphisme~$(Y,y)\to (X,x)$ entre espaces~$k$-analytiques pointés induit un isomorphisme de germes si et {\em seulement si} c'est un isomorphisme local. 

\deux{compldefgerm} Soit~$\phi: (Y,y)\to (X,x)$ un morphisme de germes d'espaces~$k$-analytiques. 
On appellera {\em représentant} de~$\phi$ la donnée d'un morphisme~$\psi : (Z,z)\to (T,t)$ d'espaces analytiques pointés et d'un diagramme commutatif de morphisme de germes~$$\diagram (Y,y)\dto_\phi\rto^{\sim} & (Z,z)\dto^\psi \\(X,x)\rto ^{\sim}& (T,t)\enddiagram.$$ Tout morphisme de germes admet un représentant. 

\medskip
On dira qu'un morphisme de germes est fini (resp. fini et plat, resp. à fibre finie, resp. sans bord) s'il admet un représentant fini (resp. fini et plat, resp. dont la fibre distinguée est finie, resp. sans bord). Un morphisme de germes~$(Y,y)\to (X,x)$ à fibre finie admet toujours un représentant~$(Z,z)\to (T,t)$ qui est compact (il suffit de prendre~$Z$ et~$T$ compacts) et dont la fibre en~$t$ est le singleton~$\{z\}$. 

\deux{bongerme} Nous dirons qu'un germe~$(X,x)$ d'espace~$k$-analytique est bon (resp. séparé, resp. sans bord) s'il existe un espace~$k$-analytique pointé~$(Y,y)$, où~$Y$ est bon (resp. séparé, resp. sans bord) et un isomorphisme de~$k$-germes~$(Y,y)\simeq (X,x)$. 

\deux{domgerme} Un {\em domaine analytique} d'un germe d'espace analytique~$(X,x)$ est une classe d'isomorphie de~$(X,x)$-germes de la forme~$(V,x)$ où~$V$ est un domaine analytique de~$X$ contenant~$x$.

\deux{temloc} Soit~$(X,x)$ un germe d'espace~$k$-analytique. Temkin a défini la {\em réduction graduée}~$\red{(X,x)}$ de~$(X,x)$. C'est un espace topologique quasi-compact non vide muni d'une application continue vers l'espace de Zariski-Riemann~$\PP_{\red{\hres(x)}/\red k}$ qui est un homéomorphisme local ; le diagramme~$\red{(X,x)}\to \PP_{\red{\hres(x)}/\red k}$ est fonctoriel en le germe~$(X,x)$. 

\trois{rappbijdom} La flèche$(V,x)\mapsto \red{(V,x)}$ établit une bijection entre l'ensemble des domaines analytiques de~$(X,x)$ et celui des ouverts quasi-compacts et non vides de~$\red{(X,x)}$. 

\trois{critredtemk} Le germe~$(X,x)$ est bon (resp. séparé, resp. sans bord) si et seulement si~$\red{(X,x)}$ est un ouvert affine de~$\PP_{\red{\hres(x)}/\red k}$ (resp. est un ouvert de~$\PP_{\red{\hres(x)}/\red k}$, resp. est égal à $\PP_{\red{\hres(x)}/\red k}$). 

\trois{constredtemk} Donnons une description de~$\red{(X,x)}$ dans le cas où~$(X,x)$ est bon. On peut alors supposer que~$X$ est affinoïde ; soit~$\sch A$ l'algèbre correspondante. On dispose d'un morphisme naturel d'annéloïdes~$f : \red {\sch A}\to \red{\hres(x)}$, induit par l'évaluation en~$x$. 
On a alors
$$\red{(X,x)}=\PP_{\red{\hres(x)}/\red k}\{f(\red {\sch A})\},$$ qui est bien un ouvert affine
de~$\PP_{\red{\hres(x)}/\red k}$ puisque~$\red {\sch A}$ est de type fini sur~$\red k$. 

\trois{fonct} On ne suppose plus que~$(X,x)$ est bon. Soient~$(f_1,\ldots,f_r$ des fonctions 
analytiques inversibles
au voisinage de~$x$ et soit~$(V,x)$ le domaine analytique de~$(X,x)$ défini par les inégalités~$|f_i|\leq |f_i(x)|$
pour tout~$i$. On a alors
$$\red {(V,x)}=\red{(X,x)}\times_{\PP_{\red{\hres(x)}/\red k}}
\PP_{\red{\hres(x)}/\red k}\{\red{f_1(x)},\ldots, \red{f_r(x)}\}.$$ 

\deux{morsansbord} La caractéristion des germes sans bord en termes de leur réduction de Temkin donnée au~\ref{critredtemk} 
ci-dessus admet la généralisation suivante au cas relatif : si~$(Y,y)\to (X,x)$ est un morphisme de germes
d'espaces~$k$-analytiques, il est sans bord si et seulement si l'application continue naturelle
$$\red{(Y,y)}\to \PP_{\red{\hres(y)}/\red k}\times_{ \PP_{\red{\hres(y)}/\red k}}\red{(X,x)}$$ est bijective 
(c'est alors un homéomorphisme).

\subsection*{Espaces quasi-lisses}

\deux{defql} Soit~$X$ un espace~$k$-analytique et soit~$x\in X$. On a~$$\dim {\hres(x)} \Omega^1_X\otimes\hres(x)\geq \dim x X\; ;$$on dit que~$X$ est {\em quasi-lisse} en~$x$ si l'on a égalité. Le lieu de quasi-lissité de~$X$ en est un ouvert de Zariski, sur lequel~$\Omega^1$ est localement libre ; il coïncide avec le lieu de régularité géométrique de~$X$. On dit que~$X$ est {\em quasi-lisse} s'il est quasi-lisse en chacun de ses points. 

\trois{qlform} Si~$V$ est un domaine analytique de~$X$ contenant~$x$ alors~$$\dim {\hres(x)} \Omega^1_X\otimes\hres(x)=\dim {\hres(x)} \Omega^1_X\otimes\hres(x)$$ et~$\dim x X=\dim x V$ ; par conséquent,~$V$ est quasi-lisse en~$x$ si et seulement si~$X$ est quasi-lisse en~$x$. 

\trois{qlsansbd} Si~$X$ est lisse en~$x$, il est quasi-lisse en~$x$. 

\trois{domlql} Il résulte de~\ref{qlform} et~\ref{qlsansbd} que tout domaine analytique d'un espace analytique lisse est quasi-lisse ; nous allons établir, dans le cas d'un bon espace, une réciproque locale de cette assertion. Soit donc~$V$ un bon espace~$k$-analytique et soit~$x\in V$. Supposons que~$V$ est  quasi-lisse en~$x$. Nous allons démontrer qu'il existe un voisinage affinoïde de~$x$ dans~$V$ qui s'identifie à un domaine affinoïde d'un espace~$k$-analytique lisse. 

\medskip
Posons~$d=\dim x V$ ; comme~$V$ est bon et comme~$\Omega^1_V$ est localement libre de rang~$d$ au voisinage de~$x$, il existe un voisinage affinoïde quasi-lisse~$W$ de~$x$ dans~$X$~$d$ fonctions analytiques~$f_1,\ldots,f_d$ sur~$V'$ telles que~$\Omega^1_W$ soit libre de base~${\rm d}f_1,\ldots,{\rm d}f_d$ ; soit~$\phi : W\to \Aff^{d,{\rm an}}_k$ le morphisme induit par les~$f_i$. Par choix des~$f_i$, le faisceaux~$\Omega_{W/\Aff^{d,{\rm an}}_k}$ est nul ; par conséquent,~$\phi$ est purement de dimension relative nulle ; en particulier~$\dim x \phi=0$.

\medskip
Il est dès lors loisible de restreindre~$W$ de sorte qu'il soit connexe et que le morphisme~$\phi$ admette une factorisation~$W\to W'\to T\to \Aff^{d,{\rm an}}_k$ où :

\medskip
$\bullet$~$T\to \Aff^{d,{\rm an}}_k$ est étale (cela entraîne que~$T$ est lisse et purement de dimension~$d$) ; 

$\bullet$~$W'$ est un domaine affinoïde connexe de~$T$ ; 

$\bullet$~$W\to W'$ est fini. 

\medskip
Soient~$\sch A$ et~$\sch B$ les algèbres~$k$-affinoïdes respectivement associées à~$W'$ et~$W$ ; les espaces~$W$ et~$W'$ sont connexes, non vides et quasi-lisses (le premier par hypothèse, le second en tant que domine affinoïde de l'espace lisse~$T$) ; ils sont donc tous deux réguliers, et dès lors irréductibles ; par construction, on a~$\dim k W=\dim k W'=d$ ; il s'ensuit que~$W\to W'$ est surjectif.  Le morphisme~$\spec \sch B\to \spec \sch A$ est un morphisme fini et dominant entre deux schémas intègres et réguliers ; il est donc plat. 

\medskip
On a~$\Omega_{W/\Aff^{d,{\rm an}}_k}=0$ ; faisceau~$\Omega^1_{W/W''}$ est {\em a fortiori} nul ; par conséquent, le morphisme fini et plat ~$\spec \sch B\to \spec \sch A$ est non ramifié, et partant étale. 

\medskip
Soit~$t$ l'image de~$x$ sur~$T$. Les catégories des germes finis étales sur~$(T,t)$ et~$(W',t)$ étant équivalentes, le bon germe~$(W,x)$ s'identifie à un germe de domaine affinoïde fermé d'un germe~$(X,x)$ qui est fini et étale sur le germe lisse~$(T,t)$, et est donc lui-même lisse ; ceci achève la démonstration. 

\trois{rremlql} Il découle de~\ref{qlsansbd} et~\ref{domlql} que si~$X$ est un espace~$k$-analytique et si~$x$ appartient à l'intérieur analytique de~$X$ alors~$X$ est quasi-lisse en~$x$ si et seulement si il est lisse en~$x$.

\section{Quelques compléments
sur les espaces analytiques
généraux}
\subsection*{Morphismes finis et plats}

\deux{remtopsep} Soit~$L$ une extension complète de~$k$, et soit~$X$ un espace~$k$-analytique tel que~$X_L$ soit topologiquement séparé ; l'espace~$X$ est alors topologiquement séparé. En effet, soient~$V$ et~$W$ deux domaines affinoïdes de~$X$ ; comme~$X_L$ est topologiquement séparé, l'intersection~$V_L\cap W_L$ est compacte ; son image sur~$X$ est dès lors quasi-compacte. Mais  l'image en question n'est autre que~$V\cap W$, qui est contenu dans l'espace topologique séparé~$V$ ; il s'ensuit que~$V\cap W$ est compact, ce qui achève ma démonstration.

\deux{finpropfin} Soit~$\phi : Y\to X$ un morphisme entre espaces~$k$-analytiques. Les propositions suivantes sont équivalentes : 

\medskip
i)~$\phi$ est fini ; 

ii)~$\phi$ est propre et à fibres (ensemblistement) finies. 

\medskip
En effet; i)$\Rightarrow ii)$ est clair. Supposons que ii) soit vraie, soit~$x\in X$ et soient~$y_1,\ldots,y_r$ les antécédents de~$x$ sur~$Y$. Pour tout~$i$, le point~$y_i$ est isolé dans~$\phi\inv(x)$ et appartient à l'intérieur analytique de~$\phi$ ; par conséquent,~$\phi$ est fini en~$y_i$, et induit donc un morphisme fini d'un voisinage ouvert~$V_i$ de~$y_i$ vers un voisinage ouvert~$U_i$ de~$x$ ; on peut, quitte à les restreindre, supposer les~$V_i$ deux à deux disjoints. Par compacité de~$\phi$, il existe un voisinage ouvert~$U$ de~$x$ contenu dans~$\bigcap U_i$ et tel que~$\phi\inv(U)\subset \coprod V_i$. On peut dès lors écrire~$$\phi\inv(U)=\coprod \phi_{|V_i}\inv(U)\;$$  par conséquent,~$\phi\inv(U)\to U$ est fini ; la finitude étant une propriété locale sur le but,~$\phi$ est fini. 

\deux{comirrdoman} Soit~$X$ un espace~$k$-analytique, soit~$V$ un domaine analytique fermé de~$X$ et soit~$Y$ une composante irréductible de~$V$. Si~$\partial \an Y=0$ alors~$Y$ est une composante irréductible de~$X$. En effet, dans ce cas,~$Y\hookrightarrow X$ est compacte, sans bord et à fibres finies, donc propre. Il s'ensuit que son image, qui n'est autre que~$Y$, est un fermé de Zariski de~$X$. Mais comme~$Y$ est une composante irréductible de~$V$, son adhérence pour la topologie de~$X$ est une composante irréductible de~$X$ ; or d'après ce qui précède, cette adhérence n'est autre que~$Y$ elle-même, d'où notre assertion. 

\deux{descfinplat} Soit~$\phi: Y\to X$ un morphisme entre espaces~$k$-analytiques. Supposons qu'il existe une extension complète~$L$ de~$k$ telle que~$\phi_L : Y_L\to X_L$ soit fini (resp. fini et plat). Sous ces hypothèses nous allons prouver que~$\phi$ est fini (resp. fini et plat). 

\trois{descfini} {\em Montrons que le morphisme~$\phi$ est fini.} On peut supposer, en raisonnant G-localement sur le but, que~$X$ est compact. L'espace~$Y_L$ est alors compact par compacité de~$\phi_L$ ; on en déduit, à l'aide de~\ref{remtopsep}, que~$Y$ est compact. La flèche~$\phi$ est donc compacte, et ses fibres sont finies puisque c'est le cas de celles de~$X_L$. 

Par ailleurs, le fait que~$\phi_L$ soit sans bord entraîne, en vertu d'un théorème de descente dû à Conrad et Temkin, que~$\phi$ est sans bord ; par conséquent,~$\phi$ est fini. 

\trois{descplat} {\em Le cas où~$\phi_L$ est plat.} Soit~$V$ un domaine affinoïde de~$X$ et soit~$W$ son image réciproque sur~$Y$ ; c'est un domaine affinoïde de~$Y$. On note~$\sch A$ et~$\sch B$ les algèbres respectivement associées à~$V$ et~$W$ ; l'algèbre~$\sch B$ est finie sur~$\sch A$. Par hypothèse,~$\sch B_L$ est finie et plate sur~$\sch A_L$ ; combiné à la platitude de~$\sch A_L$ sur~$\sch A$ et à la fidèle platitude de~$\sch B_L$ sur~$\sch B$, cela entraîne la platitude de~$\sch B$ sur~$\sch A$, qui est ce qu'on souhaitait établir. 

\trois{remdescfiniplat} {\em Remarques.} En pratique, si l'on sait {\em a priori} que~$\phi$ est sans bord ({\em e.g} si~$Y$ est sans bord), on n'a pas besoin d'utiliser le résultat de descente de Conrad et Temkin ; et si l'on sait {\em a priori} que l'espace~$Y$ est topologiquement séparé, il est inutile de se référer à~\ref{remtopsep}. 

Il peut également arriver que la platitude de~$\phi$ soit automatique, une fois connue sa finitude ; c'est notamment lorsque~$Y$ est une courbe normale et~$X$ une courbe réduite ({\em cf.}~\ref{nonconstfiniplat} {\em infra}.)   

\deux{introfl} Soit~$\phi:Y\to X$ un morphisme fini et {\em plat} entre espaces~$k$-analytiques.

\trois{ouvpropfin} L'application continue sous-jacente à~$\phi$, dont on a rappelé qu'elle est compacte et à fibres finies, est également ouverte, et partant justiciable des résultats généraux rappelés au début de l'article (\ref{appcont} et~\ref{imouvferm}-~\ref{fprop}).

\trois{degbut} Si~$U$ est un domaine affinoïde de~$X$ contenant~$x$, son image réciproque~$\phi\inv(U)$ est un domaine affinoïde de~$Y$, dont l'algèbre des fonctions est un~${\sch O}_U(U)$-module fini et localement libre ; si~$U$ est connexe, le rang de ce module est bien défini, et ne dépend pas du choix de~$U$ ; on l'appelle le {\em degré de~$\phi$ au-dessus de~$x$} et on le note~$\deg_x \phi$ ; lorsque~$X$ est bon,~$Y$ l'est aussi et~$\deg_x \phi$ est simplement le rang du~${\sch O}_{X,x}$-module libre~$\prod\limits_{y\in \phi\inv(x)} {\sch O}_{Y,y}$. 

\trois{degsource} Soit~$y\in \phi\inv(x)$ et soit~$U$ un domaine analytique de~$X$ contenant~$x$ et tel que~$\phi\inv(U)$ isole~$y$ ; le degré de~$\phi\inv(U)_y\to U$ au-dessus de~$x$ ne dépend que de~$x$, et pas de~$U$ ; on l'appelle le degré de~$\phi$ {\em en~$y$} et on le note~$\deg ^y \phi$ ; dans le cas où~$X$ et~$Y$ sont bons~$\deg^y \phi$ est simplement le rang du~${\sch O}_{X,x}$-module libre~${\sch O}_{Y,y}$.

\deux{divdeg} Si~$\phi: Y\to X$ est un morphisme fini et plat entre espaces analytiques~$x\mapsto \deg _x \phi$ est une fonction localement constante sur~$X$ ; sa valeur est donc bien définie dès que~$X$ est connexe et non vide, et est alors simplement appelée le degré de~$\phi$. Soit~$x\in X$ et soit~$U$ un voisinage analytique connexe de~$x$ tel que~$\phi\inv(U)$ sépare les antécédents de~$x$ ; l'égalité~$$\deg \; (\phi\inv(U)\to U)=\sum_{y\in \phi\inv(x)}\deg \;(\phi\inv(U)_y\to U)$$ peut se récrire~$$\deg_x \phi=\sum_{y\in \phi\inv(x)}\deg ^y \phi.$$ Il en résulte notamment que le cardinal de~$\phi\inv(x)$ est inférieur ou égal au degré de~$\phi$ au-dessus de~$x$. Si~$X$ est connexe et non vide, le cardinal de toute fibre de~$\phi$ est donc majoré par~$\deg \;\phi$. 

\deux{deggerm} Soit~$\phi : (Y,y)\to (X,x)$ un morphisme fini et plat entre {\em germes} d'espaces~$k$-analytiques. Il possède un représentant~$\psi : (Z,z)\to (T,t)$ qui est fini et plat ; le degré de~$\psi$ en~$z$ ne dépend alors que de~$\phi$, et pas du choix de~$\psi$ ; on l'appelle le degré de~$\phi$ et on le note~$\deg\;\phi$. 

\subsection*{Points d'Abhyankar}

\deux{bord-shilov-gen}
Soit~$\sch A$ une~$k$-algèbre de Banach. 
On suppose que la semi-norme spectrale de~$\sch A$
est multiplicative (c'est par exemple le cas si la norme 
de~$\sch A$ est multiplicative, auquel cas elle coïncide
avec sa semi-norme spectrale) ; par définition, cela revient
à demander qu'il existe~$x\in \sch M(\sch A)$ (nécessairement
unqiue) tel que~$|f(x)|=\sup_{y\in \sch M(A)}|f(y)|$
pour toute~$f\in \sch A$, 
c'est-à-dire encore tel que~$|f(y)|\leq |f(x)|$
pour toute~$f\in \sch M(A)$. Nous allons
montrer que
sous cette hypothèse, le singleton~$\{x\}$ est 
le bord de Shilov de~$\sch M(\sch A)$. On sait déjà
par hypothèse
que toute fonction~$f\in \sch A$ atteint son 
maximum en~$x$.

\medskip
Soit maintenant~$\sch K$ un compact
de~$\sch M(\sch A)$ ne contenant pas~$x$. Pour tout
point~$y$
de~$\sch K$, 
il existe~$f_y$
appartenant à~$A$ telle que~$|f_y(y)|\neq |f_y(x)|$, et donc 
telle
que $|f_y(y)|<|f_y(x)|$ (ceci implique que~$|f_y(x)|\neq 0$). 
L'inégalité
stricte~$|f_y|<|f_y(x)|$ reste vraie
sur un voisinage ouvert~$V_y$ de~$y$, et l'inégalité large correspondante est vraie sur tout~$\sch K$.

Par compacité, il existe~$y_1,\ldots, y_m\in \sch K$ tels que les~$V_{y_i}$ recouvrent~$\sch K$. 
Par construction, le maximum~$|\prod f_{y_i}(x)|$ 
de~$|\prod f_{y_i}|$ n'est pas atteint sur~$\sch K$. 
Il s'ensuit que~$x$
est bien le bord de Shilov de~$\sch M(\sch A)$, comme
annoncé. 

\deux{def-univmult}
Soit~$\sch A$
une~$k$-algèbre de Banach. On dit que la norme de~$\sch A$ est
{\em universellement multiplicative}
si la norme tensorielle de~$L\hotimes_k \sch A$ est multiplicative
pour tout extension complète~$L$ de~$k$. Supposons que ce soit le cas,
et soit~$L$ une extension complète de~$k$. 
Comme sa norme  est multiplicative, l'algèbre
$L\hotimes_k \sch A$ est intègre, et il en va
de même de~$L\otimes_k\sch A$, qui s'injecte dans~$L\hotimes_k\sch A$. 
On en déduit que~$\sch A$
est intègre et que~$k$ est algébriquement clos dans~${\rm Frac}\;\sch A$. 

\deux{affinoide-univmult}
Soit~$\sch A$ une algèbre~$k$-affinoïde dont la norme est universellement multiplicative, 
et soit~$x$ le point correspondant de~$\sch M(\sch A)$. La valeur absolue de la~$k$-algèbre
de Banach~$\hres(x)$ est alors universellement multiplicative, ce qui implique
notamment que~$k$ est algébriquement clos dans~$\hres(x)$ ({\em cf}.~\ref{def-univmult}
{\em supra}). 

\medskip
Soit~$L$ une extension complète de~$k$. La norme tensorielle
de~$L\hotimes_k \sch A$
étant multiplicative, elle définit un point~$y$ sur~$\sch M(L\hotimes_k \sch A)$. 
La flèche
canonique~$\sch M(L\hotimes_k \sch A)\to \sch M(\sch A)$ est surjective ; en conséquence,
$\sch A\to L\hotimes_k \sch A$ préserve les (semi)-normes spectrales, ce qui signifie que~$y$ est situé
sur la fibre~$\sch M(L\hotimes_k \hres(x))$ de~$\sch M(L\hotimes_k \sch A)$ en~$x$. 

\medskip
Soient~$a$ et~$b$ deux éléments de~$\sch A$ avec~$b(x)\neq 0$ (ce qui revient à demander
que~$b\neq 0$), et soit~$c\in L\hotimes_k \sch A$. 
Soit~$z\in \sch M(L\hotimes_k \hres(x))$. On a
$$|(a(z)b(z)^{-1}c(z)|=|a(x)b(x)^{-1}|\cdot  |c(z)|\leq |a(x)b(x^{-1})|\cdot |c(y)|=|a(y)b(y)^{-1}c(y)|.$$
Par densité, il vient~$|f(z)|\leq |f(y)|$ pour toute~$f\in L\hotimes_k \hres(x)$. Ainsi, 
$y$ est le point de~$\sch M(L\hotimes_k \hres(x))$ 
correspondant à la (semi)-norme spectrale de~$L\hotimes_k \sch A$,
et~$\{y\}$ est donc le bord de
Shilov de~$\sch M(L\hotimes_k \hres(x))$ d'après~\ref{bord-shilov-gen}. 

\deux{etar-univ-mult}
Soit~${\bf r}=(r_1,\ldots, r_n)$ un polyrayon. Pour toute extension complète~$L$
de~$k$, 
le produit tensoriel complété~$L\hotimes_k k\{{\bf r}^{-1}{\bf T}\}$
s'identifie à~$L\{{\bf r}^{-1}{\bf T}\}$, dont la norme est multiplicative et correspond
au point~$\eta_{{\bf r}, L}$ de~$\sch M(L\{{\bf r}^{-1}{\bf T}\})$. 

\medskip
Il résulte alors
de~\ref{affinoide-univmult}
que la valeur
absolue de~$\hres(\eta_{\bf r})$ est universellement multiplicative, que~$k$
est algébriquement clos dans~$\hres(\eta_{\bf r})$, et que~$\{\eta_{{\bf r},L}\}$ est pour
toute extension complète~$L$ de~$k$ le bord de Shilov
de~$L\hotimes_k \hres(\eta_{\bf r})$.

\deux{definition-pointabhyankar}
Soit~$X$ un espace~$k$-analytique
et soit~$x$
un point de~$X$. 

\trois{vraiment-def-abhyankar}
Soit~$Y\hookrightarrow X$
une immersion dont l'image contient~$x$. La dimension
de~$Y$ est alors minorée par~$d_k(x)$. On dit que le point~$x$ 
est {\em d'Abhyankar}
si l'on peut trouver une telle immersion~$Y\hookrightarrow X$
avec~$\dim {} Y=d_k(x)$ (on dit alors souvent que~$d_k(x)$ est 
le {\em rang}
de~$x$) ; si c'est le cas, on peut toujours supposer que~$Y$
s'identifie à un sous-espace analytique fermé d'un domaine affinoïde~$X'$ de~$X$. 

\medskip
Notons qu'un point d'Abhyankar de rang $0$ est simplement un point rigide.

\trois{definition-abhyankar-presentation}
On suppose à partir de maintenant que~$x$ est 
un point d'Abhyankar de~$X$ dont on note~$n$
le rang. Une
{\em présentation d'Abhyankar}
de~$x$ consiste en les données suivantes : 

\medskip
$\bullet$ un sous-espace analytique fermé~$Y$ d'un domaine affinoïde~$X'$ de~$X$ tel que~$x\in Y$
et tel que~$\dim {} Y=n$ ; 

$\bullet$ un morphisme~$\phi \colon Y\to \Aff^{d_k(x), {\rm an}}_k$ donné par
une famille finie~$(f_1,\ldots, f_n)$ de fonctions inversibles en~$x$ et telles que
les~$\red{f_i(x)}$ forment une base de transcendance de~$\red{\hres(x)}$ sur~$k$ ; 
il revient au même de demander que~$\phi(x)$ soit de la forme~$\eta_{\bf r}$ pour un certain
polyrayon~$\bf r$. 

\medskip
Il existe toujours une présentation d'Abhyankar de~$x$. 

\trois{conclusion-presentation-abhyankar}
Soit~$(Y,\phi)$ une présentation d'Abhyankar de~$x$, et soit~$\bf r$
le polyrayon tel que~$\phi(x)=\eta_{\bf r}$.

\medskip
Si~$t\in \phi\inv(\eta_{\bf r})$ on a
$n\geq d_k(t)=d_{\hres(\eta_{\bf r})}(t)+n$, et donc~$d_{\hres(\eta_{\bf r})}(t)=0$ ; ainsi, $\phi^{-1}(\eta_{\bf r})$ 
est de dimension nulle. Ceci entraîne que le corps~$\hres(x)$ est une extension finie
de~$\hres(\eta_{\bf r})$. Il en résulte les faits suivants. 

\medskip
$\bullet$ Le corpoïde~$\red{\hres(x)}$ est fini sur~$\red{\hres(\eta_{\bf r})}$ d'après~\ref{dimresoide} ; comme~$\red{\hres(\eta_{\bf r})}$
est transcendant pur sur~$\red k$ (là encore, en vertu de~\ref{remresidgauss}), le
corpoïde~$\red{\hres(x)}$ est de type fini sur~$\red k$. 

Cela signifie (\ref{interpfinioide}) que~$\red{\hres(x)_1}$ est de type fini sur~$\red k_1$, et que
le groupe~$|\hres(x)\ti|/|k\ti|$ 
est de type fini. 

$\bullet$ La clôture algébrique
de~$k$ dans~$\hres(x)$ est finie sur~$k$ en vertu de~\ref{etar-univ-mult} ; le corps~$\got s(x)$
est {\em a fortiori}
une extension finie de~$k$.  

\trois{suite-antecedent}
Soit~$L$ une extension complète de~$k$ et soit~$\pi : Y_L\to Y$
l'application canonique. Comme
la fibre~$\phi^{-1}(\eta_{\bf r})$  est de dimension nulle, 
l'ensemble compact~$\pi^{-1}(x)\cap \phi^{-1}(\eta_{{\bf r},L})$ est discret, et partant fini, 
et chacun de ses éléments~$y$ vérifie l'égalité
$d_L(y)=d_L(\eta_{{\bf r}, L})=n=\dim {} Y_L$ ; en conséquence, tous les points
de~$\pi^{-1}(x)\cap \phi^{-1}(\eta_{{\bf r},L})$ sont des points d'Abhyankar de~$X_L$ 
de même rang que~$x$.

\deux{prop-shilov-antec}
{\bf Proposition.}
{\em L'ensemble fini~$E:=\pi^{-1}(x)\cap \phi^{-1}(\eta_{{\bf r},L})$ est le bord de Shilov
de~$\pi^{-1}(x)=\sch M(\hres(x)\hotimes_k L)$ ; en particulier, il est non vide
et ne dépend pas
du choix de la présentation~$(Y,\phi)$
de~$x$.}

\medskip
{\em Démonstration.} 
Posons
$$A=\hres(\eta_{\bf  r})\hotimes_k L\;\;{\rm et}\;\;
B=\hres (x)\hotimes_k L\simeq A\otimes_{\hres(\eta_{\bf r})}\hres (x).$$
Par construction, $B$ est libre de rang fini sur~$A$. On note~$p$ la flèche naturelle
de~$\mathscr M(B)$ vers~$\mathscr M(A)$ ; on a~$E=p^{-1}(\eta_{{\bf r},L})$. 
Rappelons par ailleurs que le bord de Shilov de~$\sch M(\sch A)$ est le singleton
$\{\eta_{{\bf r}, L}\}$ (\ref{etar-univ-mult}). 

\trois{antec-shilov-inclu}
Soit~$b\in B$ et soit~$\chi=T^N+\sum_{i<N} a_iT^i\in A[T]$ son polynôme caractéristique. 
Soit~$y\in \mathscr M(A)$. Le maximum de~$|g|$ sur~$p^{-1}(y)$
est égal à la norme spectrale de l'image de~$b$ dans~$B\otimes_A \hres (y)$, c'est-à-dire
à~$\max |a_i(y)|^{1/{N-i}}$. Puisque~$\eta_{{\bf r}, L}$ est l'unique point du bord de Shilov de~$\mathscr M(A)$,
on a
$$\max |a_i(y)|^{1/{N-i}}\leq \max |a_i(\eta_{{\bf r}, L})|^{1/{N-i}}.$$ Par conséquent, le maximum de~$|b|$ 
sur~$p^{-1}(y)$ est majoré par son maximum sur~$E=p^{-1}(\eta_{{\bf r},L})$. Ainsi,~$|b|$ atteint son maximum
sur~$E$. 

\trois{antec-shilov-conclu}
Soit maintenant~$\sch K$ un compact de~$\mathscr M(B)$ ne contenant pas~$E$ ; nous allons construire une fonction~$b\in B$ 
dont le maximum n'est pas atteint sur~$\sch K$. 

\medskip
Soit~$e$
un point
de~$E$ n'appartenant
pas à~$\sch K$ ; il existe
une
fonction~$g$
appartenant à~$B\otimes_A \hres (\eta_{{\bf r},L})$ qui
vaut~$1$ sur~$e$ et~$0$ en tous les autres points de~$E$. Par densité
du corps des fractions de~$A$ dans~$\hres(\eta_{{\bf r}, L})$,
on peut trouver deux éléments
$a$ et~$\alpha$ de~$A$ avec~$\alpha\neq 0$ et un élément~$\beta$ de~$b$ tels que
la fonction~$|\beta a/\alpha|$ vaille~$1$ sur~$e$ et soit majorée strictement 
par~$1/2$ en tout point
de~$E\setminus \{e\}$.
Si l'on pose~$R=|(\alpha/a)(\eta_{{\bf r},L})|$,
la fonction~$|\beta|$ vaut~$R$ en~$e$ et est majorée
strictement
par~$R/2$
sur~$E\setminus \{e\}$. 
Posons~$$
U=\{z\in \mathscr M(B), |\beta(z)|<R/2\}\;\;{\rm et}\;\;V=\{z\in \mathscr M(B)\setminus \sch K, |\beta(z)|>3R/4\}.$$
Par construction, $U\cup V$ contient~$E$ ; par
compacité, il existe un voisinage ouvert~$W$ de~$\eta_{{\bf r}, L}$ 
dans~$\mathscr M(A)$ tel que~$p^{-1}(W)\subset U\cup V$ ; désignons
par~$\sch H$
le compact~$\mathscr M(A)\setminus W$. 

Comme~$\eta_{{\bf r}, L}$ 
est l'unique
point du bord
de Shilov de~$\mathscr M(A)$, il existe une fonction~$\lambda \in A$ telle que le maximum
$|\lambda(\eta_{{\bf r}, L})|$
de~$|\lambda|$
ne soit pas atteint
sur~$\sch H$. Il existe alors~$m>0$
tel que
$$\left(\sup_{z\in \sch H} |\lambda(z)|\right)^m \cdot \sup_{z\in p^{-1}(\sch H)} |\beta(z)|<|\lambda(\eta_{{\bf r}, L})|^m\cdot R.$$

Puisque~$|\beta(e)|=R$, le maximum de~$|\lambda^m\beta|$
sur~$\mathscr M(B)$ est au moins égal
à~$|\lambda(\eta_{{\bf r}, L})|^m\cdot R$ ; par ce qui précède, il n'est pas atteint sur~$p^{-1}(\sch H)$. Il n'est par ailleurs pas atteint sur~$U$ non plus : 
on a en effet pour tout~$z\in U$ les inégalités~$|\lambda(z)|\leq |\lambda(\eta_{{\bf r}, L})|$ et~$|\beta(z)|<R/2$. 

On a par construction~$\sch K\cap V=\emptyset$
et~$\mathscr M(B)\setminus p^{-1}(\sch H)=p^{-1}(W)\subset U\cup V$ ;
il vient~$\sch K\subset p^{-1}(\sch H)\cup U$. On en déduit que
le maximum de~$|\lambda^m\beta|$ n'est pas atteint sur~$\sch K$, ce qui achève la démonstration.~$\Box$

\deux{ex-ant-can}
Soit~$X$ un espace~$k$-analytique
et soit~$x$ un point d'Abhyankar de~$X$. 
Soit~$L$ une extension complète de~$k$.
Il résulte de~\ref{suite-antecedent}
et de la
proposition~\ref{prop-shilov-antec}
que~$\sch M(\hres(x)\hotimes_k L)\subset X_L$
possède un bord de Shilov fini dont tout point~$y$
est d'Abhyankar de même rang que~$x$. Ce bord de Shilov
sera noté~$\shil L x$. 

\deux{exemples-shil-lx}
{\bf Exemples.}

\trois{shil-l-etar}
Soit~${\bf r}=(r_1,\ldots, r_n)$ un polyrayon. On a alors~$\eta_{{\bf r},[L]}=\{\eta_{{\bf r}, L}\}$.
On peut voir ce fait comme une conséquence
triviale de la proposition~\ref{prop-shilov-antec}, appliquée en prenant
la présentation~$(\Aff^{n,\rm an}_k, \mathsf{Id})$. Mais on en a donné une preuve directe
au préalable en~\ref{suite-antecedent}, et on a d'ailleurs utilisé plusieurs fois ce résultat au cours
de la preuve de la proposition~\ref{prop-shilov-antec}.

\trois{shil-l-trans}
Soit~$X$ un espace~$k$-analytique, soit~$x$ un point
d'Abhyankar de~$X$, soit~$L$ une extension complète de~$k$ et soit~$L'$ une
extension complète de~$L$. On a alors
$$\shil {L'} x=\bigcup_{y\in \shil Lx}\shil {L'} y.$$
En effet, posons~$n=d_k(x)$ et choisissons
une présentation d'Abhyankar~$(Y,\phi)$ de~$x$. On désigne par~$\bf r$
le poyrayon tel que~$\phi(x)=\eta_{\bf r}$ et par
$$\pi \colon Y_L\to X, \pi'\colon Y_{L'}\to X\;\;{\rm et}\;\varpi \colon Y_{L'}\to L$$
les flèches canoniques. On a 

$$\shil {L'} x=(\pi')^{-1}(x)\cap \phi^{-1}(\eta_{{\bf r}, L'})\;\;\;\text{(prop.~\ref{prop-shilov-antec})}$$
$$=\bigcup_{y\in \pi^{-1}(x)} \varpi^{-1}(y)\cap \phi^{-1}(\eta_{{\bf r}, L'})$$
$$=\bigcup_{y\in \pi^{-1}(x)\cap \phi^{-1}(\eta_{\bf r}, L)} \varpi^{-1}(y)\cap \phi^{-1}(\eta_{{\bf r}, L'})\;\;\;\text{(si~}\phi(z)=\eta_{{\bf r}, L'}\;\text{alors}\;\phi
(\varpi(z))=\eta_{{\bf r}, L}\text{)}$$
$$=\bigcup_{y\in \shil L x}\shil {L'} y\;\;\text{(prop.~\ref{prop-shilov-antec})}.$$

\trois{shil-l-presquealg}
Soit~$X$ un espace~$k$-analytique de dimension finie,
soit~$x$
un point d'Abhyankar
de~$x$, 
et soit~$L$ une extension presque algébrique de~$k$ ; soit~$\pi \colon X_L\to X$
la flèche canonique. 
Le sous-ensemble~$\shil L x$ de~$X_L$ est alors égal à
la fibre~$\pi^{-1}(x)$ toute entière. 

\medskip
En effet, en considérant un~$k$-plongement isométrique
$L\hookrightarrow \KK$ et en appliquant~\ref{shil-l-trans}
on se ramène au cas où~$L=\KK$. Il
résulte immédiatement de la définition que le sous-ensemble
$\shil {\KK} x$ de~$\pi^{-1}(x)$ est invariant sous l'action de Galois. 
Comme~$\shil {\KK} x$  est non vide et
comme l'action de Galois
est transitive sur~$\pi^{-1}(x)$, 
il vient~$\shil {\KK} x=\pi^{-1}(x)$. 

\trois{conseq-shil-presqualg}
Conservons les hypothèses et notations de~\ref{shil-l-presquealg}. 
On y a établi l'égalité~$\shil L x=\pi^{-1}(x)$ ; elle entraîne en particulier
la finitude de~$\pi^{-1}(x)$, que l'on peut établir directement facilement. 

\medskip
En effet, on se ramène immédiatement au cas où~$L=\KK$ ; 
notons~$n$ le rang de~$x$, et donnons-nous une présentation
d'Abhyankar~$(Y,\phi)$ de~$x$ ; soit~$\bf r$ le polyrayon
tel que~$\phi(x)=\eta_{\bf r}$. 
Tout élément de~$\pi^{-1}(x)$
s'envoie {\em via}
$\phi$ au-dessus d'un antécédent de~$\eta_{\bf r}$ sur~$\Aff^{n, \rm an}_{\KK}$ ; mais 
la description explicite de~$\eta_{{\bf r}, \KK}$ assure que ce dernier est invariant
sous l'action de Galois, ce qui entraîne qu'il est le seul antécédent de~$\eta_{\bf r}$. 
Par conséquent,
la fibre~$\pi^{-1}(x)$ est contenue
dans~$\phi^{-1}(\eta_{{\bf r}, \KK})$, qui est compact
et de dimension nulle, et donc fini. 

\trois{shil-l-x-connexe}
Soit~$X$ un espace~$k$-analytique, soit~$x$ un point d'Abhyankar
de~$X$ et soit~$L$ une extension complète de~$K$. Soient~$L_1,\ldots, L_r$
les extensions composées de~$\got s(X)$ et~$L$ au-dessus de~$k$. On
a
$$X_L=\coprod_i X\times_{\got s(x)}L_i.$$ Lorsqu'on voit~$X$ comme un espace~$\got s(x)$-analytique,
le point~$x$ reste d'Abhyankar, ce qui permet de définir un sous-ensemble~$\shil {L_i}x$ de~$X\times_{\got s(x)}L_i$. 
On a alors~$\shil L x=\coprod \shil {L_i}x$. 

\medskip
En effet, on peut écrire
$$\hres(x)\hotimes_k L=\hres(x)\hotimes_{\got s(x)}(\got s(x)\otimes_k L)=\prod \hres(x)\hotimes_{\got s(x)}L_i.$$
Il est alors
immédiat que le bord de Shilov
de
$$\sch  M(\hres(x)\hotimes_k L)=\sch M(\prod \hres(x)\hotimes_{\got s(x)}L_i)=\coprod \sch M( \hres(x)\hotimes_{\got s(x)}L_i)$$
est la réunion des bords de Shilov des~$ \hres(x)\hotimes_{\got s(x)}L_i$, ce qui est précisément l'égalité annoncée.

%

\subsection*{Descente des immersions compactes}

\deux{defimmerintro} Soit~$\phi:Y\to X$ un morphisme entre espaces~$k$-analytiques. 

\trois{defimmer} On dira que~$\phi$ est une {\em immersion compacte} si pour tout domaine analytique compact~$U$ de~$X$ il existe un domaine analytique compact~$V$ de~$U$ tel que~$\phi\inv(U)\to U$ se factorise par une immersion fermée~$\phi\inv(U)\to V$ ; si c'est le cas,~$\phi$ est topologiquement une injection compacte. 

\trois{eximmer} S'il existe un domaine analytique fermé~$X'$ de~$X$ tel que~$Y\to X$ se factorise par une immersion fermée~$Y\to X'$ alors~$Y\to X$ est une immersion compacte. La réciproque est vraie si~$X$ est séparé et paracompact : en effet, plaçons-nous sous cette hypothèse et supposons que~$Y\to X$ soit une immersion compacte ; choisissons un G-recouvrement localement fini~$(U_i)$ de~$X$ par des domaines analytiques compacts ; pour tout~$i$, désignons par~$V_i$ un domaine analytique compact de~$U_i$ tel que~$\phi\inv(U_i)\to U_i$ se factorise par une immersion fermée~$\phi\inv(U_i)\to V_i$ ; la réunion~$V$ des~$V_i$ est alors un domaine analytique fermé de~$X$, et~$Y\to X$ se factorise par~$V$ ; étant G-localement sur son but une immersion fermée,~$Y\to V$ est une immersion fermée. 

\trois{immergloc} Supposons qu'il existe un G-recouvrement~$(U_i)$ de~$X$ par des domaines analytiques compacts et, pour tout indice~$i$, un domaine analytique compact~$V_i$ de~$U_i$ tel que~$\phi\inv(U_i)\to U_i$ se factorise par une immersion fermée~$\phi\inv(U_i)\to V_i$ ; dans ce cas~$\phi$ est une immersion compacte. 

En effet, soit~$U$ un domaine analytique compact de~$X$. Il existe un recouvrement affinoïde {\em fini}~$(W_j)$ de~$U$ raffinant le G-recouvrement~$(U\cap U_i)$ de~$U$. Fixons~$j$, et soit~$i_j$ un indice tel que~$W_j\subset U_{i_j}$. La flèche~$\phi\inv(W_j)\to W_j$ se factorise alors par une immersion fermée~$\phi\inv(W_j)\to W_j\cap V_{i_j}$. Appelons~$V$ la réunion des~$W_j\cap V_{i_j}$ ; c'est un domaine analytique compact de~$U$ ; la flèche~$\phi\inv(U)\to U$ se factorise par un morphisme~$\phi\inv(U)\to V$, qui est G-localement sur son but une immersion fermée ; par conséquent,~$\phi\inv(U)\to V$ est une immersion fermée, et~$\phi$ est bien une immersion compacte. 

\deux{descenteimmer} {\bf Proposition.} {\em Soit~$\phi: Y\to X$ un morphisme entre espaces~$k$-analytiques et soit~$L$ une extension complète de~$k$. Les assertions suivantes sont équivalentes :

\medskip
1)~$\phi$ est une immersion compacte ; 

2)~$\phi_L$ est une immersion compacte ;

3)~$\phi_L$ est un monomorphisme compact d'espaces~$L$-analytiques ;

4)~$\phi$ est un monomorphisme compact d'espaces~$k$-analytiques.}

\medskip
{\em Démonstration.} L'implication~$1)\Rightarrow 2)$ découle du~\ref{immergloc} ci-dessus ; l'implication~$2)\Rightarrow 3)$ est triviale. 

\medskip
Montrons que~$3)\Rightarrow 4)$ ; on suppose donc que~$\phi_L$ est un monomorphisme compact d'espaces~$L$-analytiques. Si~$T$ est une partie compacte de~$X$ et si~$T_L$ désigne son image réciproque sur~$X_L$ alors~$\phi_L\inv(T_L)\to \phi\inv(T)$ est surjective ; par compacité de~$\phi$ et de~$X_L\to X$, le sous-ensemble~$\phi_L\inv(Z_L)$ de~$Y_L$ est compact ; par conséquent,~$\phi\inv(Z)$ est compacte et~$\phi$ est elle-même compacte.

Soit~$Z$ un espace~$k$-analytique ; montrons que~${\rm Hom}\;(Z,Y)\to {\rm Hom}\;(Z,X)$ est injective, ce qui achèvera d'établir 4). Soient~$\psi$ et~$\psi'$ deux morphismes de~$Z$ vers~$Y$ dont les composées avec~$Y\to X$ coïncident. Comme~$\phi_L$ est un monomorphisme, on a~$\psi_L=\psi'_L$ ; ceci implique immédiatement, compte-tenu de la surjectivité de~$Z_L\to Z$, que les applications {\em ensemblistes} sous-jacentes à~$\psi$ et~$\psi'$ sont identiques. Soit~$V$ un domaine affinoïde de~$Z$ et soit~$U$ un domaine affinoïde de~$Y$ tel que~$\psi(V)=\psi'(V)$ soit contenu dans~$U$ ; soient~$\sch A$ et~$\sch B$ les algèbres affinoïdes correspondant respectivement à~$U$ et~$V$. Comme~$\psi_L=\psi'_L$, les deux morphismes~${\sch A}_L\to {\sch B}_L$ respectivement induits par~$\psi_L$ et~$\psi'_L$ sont égaux ; la flèche~${\sch B}\to {\sch B}_L$ étant injective, les deux morphismes~${\sch A}\to {\sch B}$ respectivement induits par~$\psi$ et~$\psi'$ sont égaux ; il s'ensuit que~$\psi=\psi'$, ce qu'on voulait démontrer. 

\medskip
Montrons enfin que 4)$\Rightarrow 1)$ ; on suppose donc maintenant que~$\phi$ est un monomorphisme compact d'espaces~$k$-analytiques ; nous allons montrer que c'est une immersion compacte. Par compacité de~$\phi$, l'espace~$\hres(x)$-analytique~$Y_x$ est compact. La flèche~$\phi$ étant un monomorphisme,~$Y\to Y\times_XY$ est un isomorphisme ; il en va alors de même de~$Y_x\to Y_x\times_{\hres(x)}Y_x$ , ce qui force l'espace compact~$Y_x$ à être de dimension nulle ; c'est en conséquence un espace~$\hres(x)$-affinoïde de dimension nulle, c'est-à-dire le spectre analytique d'une~$\hres(x)$-algèbre de Banach finie~$\sch A$. Celle-ci est, par ce qui précède, isomorphe à~${\sch A}\otimes_{\hres(x)}\sch A$ ; par un argument de dimension, on a~${\sch A}=\{0\}$ ou~${\sch A}=\hres(x)$. Autrement dit, le cardinal d'une fibre de~$\phi$ est au plus égal à~$1$ ; cela signifie que~$Y\to X$ est ensemblistement injective. 

Soit~$U$ un domaine affinoïde de~$X$. Comme~$\phi$ est compacte,~$\phi\inv(U)$ est un domaine analytique compact de~$Y$. Soit~$(V_i)$ un recouvrement affinoïde fini de~$\phi\inv(U)$. Fixons~$i$ ; les flèches~$V\to U$ et~$V_i\to V$ sont des monomorphismes ; par conséquent,~$V_i\to U$ est un monomorphisme. En vertu d'un théorème de Temkin et du~\ref{immergloc},~$V_i\to U$ est une immersion compacte ; il existe donc un domaine analytique compact~$U_i$ de~$U$ tel que~$V_i\to U$ se factorise par une immersion fermée~$V_i\to U_i$. Comme~$V_i\to U_i$ est sans bord,~$V_i\hookrightarrow \phi\inv(U_i)$ est sans bord ; par conséquent,~$V_i$ est une réunion de composantes connexes de~$\phi\inv(U_i)$. Comme~$\phi$ est injective,~$\phi(\phi\inv(U_i)-V_i)$ est un compact de~$U_i$ qui est disjoint de~$\phi(V_i)$ ; par conséquent l'on peut, quitte à remplacer~$U_i$ par un voisinage analytique compact convenable de~$\phi(V_i)$ dans celui-ci, supposer que~$V_i=\phi\inv(U_i)$. Posons~$V=\bigcup U_i$ ; c'est un domaine analytique compact de~$U$ par lequel~$\phi\inv(U)\to U$ se factorise ; il découle de notre construction que~$\phi\inv(U)\to V$ est G-localement sur son but une immersion fermée ; c'est donc une immersion fermée, ce qui achève la démonstration.~$\Box$

\deux{desciso} {\bf Corollaire.} {\em Soit~$\phi: Y\to X$ un morphisme entre espaces~$k$-analytiques et soit~$L$ une extension complète de~$k$. 

\medskip
i)~$\phi$ identifie~$Y$ à un domaine analytique fermé de~$X$ si et seulement si~$\phi_L$ identifie~$Y_L$ à un domaine analytique fermé de~$X_L$ ;

ii)~$\phi$ est un isomorphisme si et seulement si~$\phi_L$ est un isomorphisme.

}

\medskip
{\em Démonstration.} Comme~$Y\to X$ est surjective si et seulement si~$Y_L\to X_L$ est surjective, ii) sera une conséquence immédiate de i). Il suffit dès lors de montrer i) ; l'implication directe est évidente, et l'on peut donc se contenter de montrer l'implication réciproque. 

\medskip
On suppose que~$\phi_L$ identifie~$Y_L$ à un domaine analytique fermé de~$X_L$. Le morphisme~$\phi_L$ est en particulier une immersion compacte ; par la proposition~\ref{descenteimmer} ci-dessus,~$\phi$ est une immersion compacte. Soit~$U$ un domaine analytique compact de~$X$ ; il existe un domaine analytique compact~$V$ de~$U$ tel que~$\phi\inv(U)\to U$ se factorise par une immersion fermée~$\phi\inv(U)\to V$. Soit~$W$ un domaine affinoïde de~$V$ ; la flèche~$\phi\inv(W)\to W$ est une immersion fermée ; appelons~$\sch A$ et~$\sch B$ les algèbres affinoïdes respectivement associées à~$W$ et~$\phi\inv(W)$. La flèche~$\phi\inv(W_L)\to W_L$ identifie à la fois~$\phi\inv(W_L)$ à un fermé de Zariski et à un domaine analytique compact de~$W_L$ ; par conséquent, elle identifie~$\phi\inv(W_L)$ à une réunion de composantes connexes de~$W_L$, et ~$\spec {\sch B}_L\to \spec {\sch A}_L$ est de ce fait une immersion ouverte ; comme~$\phi\inv(W)\to W$ est une immersion fermée,~$\sch B$ est une~$\sch A$-algèbre de Banach finie, et~${\sch B}_L$ s'identifie donc à~${\sch B}\otimes_{\sch A}{\sch A}_L$. Par descente fidèlement plate schématique, l'immersion fermée~$\spec {\sch B}\to \spec \sch A$ est une immersion ouverte,~$\phi\inv(W)\to W$ identifie de ce fait~$\phi\inv(W)$ à une réunion de composantes connexes de~$W$, et en particulier à un domaine analytique compact de~$W$. Il s'ensuit que~$\phi\inv(U)\to V$ identifie~$\phi\inv(U)$ à un domaine analytique compact de~$V$, et donc de~$U$ ; la propriété que l'on souhaite établir étant G-locale sur le but, ceci achève la démonstration.~~$\Box$ 

\subsection*{Descente des domaines analytiques}

\deux{imetsansr} {\bf Proposition.} {\em Soit~$Y\to X$ un morphisme fini \'etale d'espaces~$k$-analytiques, et soit~$V$ un domaine analytique fermé de~$Y$. Son image sur~$X$ est un domaine analytique fermé de~$X$.}

\medskip
{\em Démonstration.} La question étant G-locale sur~$X$, on peut supposer~$X$ affinoïde et connexe. Dans ce cas,~$Y$ est affinoïde. En raisonnant composante connexe par composante connexe sur~$Y$, on se ramène au cas où~$Y$ est de surcroît connexe ; enfin, quitte à remplacer~$Y$ par un revêtement fini étale~$Y'\to Y$ convenable, et~$V$ par son image réciproque sur~$Y'$, on peut supposer que~$Y\to X$ est un revêtement galoisien ; le domaine analytique~$\bigcup\limits _{g\in \mathsf{Gal}\;(Y/X)} g(V)$ a même image que~$V$ sur~$X$, ce qui permet de le substituer à~$V$ et donc de supposer que~$V$ est stable sous ~$\mathsf{Gal}\;(Y/X)$. 

\medskip
Soit~$y\in V$ et soit~$x$ son image sur~$X$. L'extension graduée~$\red{\hres(y)}/\red{\hres(x)}$ est alors galoisienne, et l'ouvert quasi-compact~$\red{(V,y)}$ de~$\PP_{\red{\hres(y)}/\red{\hres(x)}}$ est stable sous l'action de Galois. 

La proposition~\ref{imvalqc}
assure que l'image de~$\red{(V,y)}$ sur~$\red{(X,x)}$ est un ouvert quasi-compact~$\mathsf U$ de ce dernier. Il correspond à un domaine analytique fermé~$(U,x)$ de~$(X,x)$. L'image réciproque de~$(U,x)$ sur~$(Y,x)$ est le domaine analytique fermé de~$(Y,y)$ correspondant à l'image réciproque de~$\mathsf U$ sur~$\red{(Y,y)}$ ; comme ~$\red{(V,y)}$ est stable sous Galois, cette image réciproque est égale à~$\red{(V,y)}$ -- et coïncide avec l'image réciproque de ~$\mathsf U$ sur~$\PP_{\red{\hres(y)}/\red{\hres(x)}}$ tout entier. Par conséquent,~$(V,y)\to (X,x)$ se factorise par~$(U,x)$, et~$(V,y)\to (U,x)$ est sans bord. Étant par ailleurs quasi-étale (puisque~$Y\to X$ est étale), il est étale. 

Il existe dès lors un voisinage analytique compact et connexe~$U'$ de~$x$ dans~$U$ et un voisinage analytique compact~$V'$ de~$y$ dans~$V$ tel que~$Y\to X$ induise un morphisme fini et étale de~$V'$ vers~$U'$ ; l'image de~$V'$ sur~$X$ est alors égale à~$U'$. 

\medskip
Ainsi, tout point de~$V$ a un voisinage dans~$V$ dont l'image sur~$X$ est un domaine analytique compact de~$X$. Par compacité de~$V$, l'image de~$V$ sur~$X$ est un domaine analytique compact de~$X$.~$\Box$ 

\deux{descgaldoman} {\bf Théorème.} {\em Soit~$X$ un espace~$k$-analytique et soit~$V$ un sous-ensemble de~$X$. Les assertions suivantes sont équivalentes : 

\medskip
i)~$V$ est domaine analytique de~$X$ ; 

ii) l'image réciproque de~$V$ sur~$X_{\KK}$ est un domaine analytique de~$X_{\KK}$.}

\medskip
{\em Démonstration.} L'implication ii)$\Rightarrow$i) est claire ; on suppose maintenant que ii) est vraie, et l'on va prouver i) ; pour toute extension complète~$L$ de~$k$, on notera~$V_L$ l'image réciproque de~$V$ sur~$X_{\KK}$. 

L'assertion à prouver étant locale pour la G-topologie de~$X$, on peut supposer que~$X$ est affinoïde et~$V$ compact. En vertu du théorème de Gerritzen-Grauert,~$V_{\KK}$ est alors une réunion finie de domaines rationnels de~$X_{\KK}$. 

\trois{domratdesfin}{\em Tout domaine rationnel de~$X_{\KK}$ est défini sur une extension finie séparable de~$k$ contenue dans~$\KK$.} En effet, soit~$W$ un tel domaine ; il est défini par une conjonction finie d'inégalités~$\{|f_i|\leq r_i |g|\}_i$ où les~$r_i$ sont des réels strictement positifs et où~$g$ et les~$f_i$ sont des fonctions analytiques sans zéro commun sur~$X_{\KK}$. Il s'ensuit immédiatement que~$|g|$ est minorée par un réel strictement positif sur~$W$, et donc que l'on peut rajouter une condition du type~$1\leq r |g|$ à la liste de celles qui décrivent~$W$. 

Si l'on se donne des fonctions analytiques~$f'_1,\ldots,f'_n, g'$ sur~$X_{\KK}$ telles que les réels~$||f'_1-f_1||, \ldots, ||f'_n-f_n||, ||g'-g||$ soient suffisamment petits, le domaine rationnel de~$X_{\KK}$ défini par les inégalités~$$\{|f'_i|\leq r_i |g'|\}_i\;{\rm et}\; 1\leq r |g'|$$ coïncide avec~$W$. Il existe donc une extension finie~$F$ de~$k$ contenue dans~$\KK$ et des fonctions analytiques~$f'_1,\ldots,f'_n, g'$ sur~$X_F$ telles que le domaine rationnel~$W$ puisse être défini par les inégalités~$$\{|f'_i|\leq r_i |g'|\}_i\;{\rm et}\; 1\leq r |g'|.$$ Si l'on désigne par~$F_{\rm sep}$ la fermeture séparable de~$k$ dans~$F$ il existe un entier non nul~$n$ tel que les~$(f'_i)^{p^n}$ et~$(g')^{p^n}$ proviennent de l'anneau des fonctions analytiques sur~$X_{F_{\rm sep}}$ ; comme~$W$ peut être tout aussi bien défini par les inégalités~$$\{|(f'_i)^{p^n}|\leq r_i^{p^n} |(g')^{p^n}|\}_i\;{\rm et}\; 1\leq r^{p^n} |(g')^{p^n}|,$$ il provient d'un domaine rationnel de~$X_{F_{\rm sep}}$.

\trois{suitedescalgdom} Il découle de ce qui précède qu'il existe une extension finie séparable~$L$ de~$k$ contenue dans~$\KK$ telle que~$V_{\KK}$ provienne d'un domaine analytique de~$X_L$, qui est nécessairement égal à l'image de~$V_{\KK}$ sur~$X_{\KK}$, c'est-à-dire à~$V_L$. Comme~$V$ est lui-même égal à l'image de~$V_L$ sur~$X$, la proposition~\ref{imetsansr} garantit alors que~$V$ est un domaine analytique compact de~$X$.~$\Box$

\subsection*{Descente du caractère affinoïde}  

\deux{remadmkl} Si~$L$ est une extension complète de~$k$ alors~$k\hookrightarrow L$ est une injection admissible ; il s'ensuit que pour toute~$k$-algèbre de Banach~$\sch A$ la flèche~${\sch A}\hookrightarrow {\sch A}_L$ est une injection admissible. 

\deux{rembanach} Soit~$X$ un espace~$k$-analytique {\em quasi-compact} ; choisissons un recouvrement affinoïde fini~$(X_i)$ de~$X$. Pour tout~$i$, l'algèbre des fonctions analytique sur~$X_i$ est~$k$-affinoïde ; elle est de ce fait munie d'une classe d'équivalence de normes (bien définie à équivalence près) ; choisissons une norme~$||.||_i$ dans cette classe. L'application~$f\mapsto \sup\limits_i||f_{|X_i}||_i$ est une norme sur l'algèbre~$\sch A$ des fonctions analytique sur~$X$, qui fait de celle-ci une~$k$-algèbre de Banach. À équivalence près,~$||.||$ ne dépend pas du choix des~$||.||_i$.

Elle ne dépend pas non plus du recouvrement~$(X_i)$. Pour le voir, donnons-nous un second recouvrement affinoïde fini ~$(X'_j)$ de~$X$, choisissons pour tout~$j$ une norme~$||.||'_j$ sur l'algèbre des fonctions analytiques sur~$X'_j$ et notons~$||.||'$ la norme sur~$\sch A$ déduite des~$||.||'_j$. 

Fixons~$i$ ; il existe un recouvrement affinoïde fini~$(X''_\ell)$ de~$X_i$ tel que pour tout~$\ell$ le domaine affinoïde~$X''_\ell$ soit contenu dans~$X'_{j_\ell}$ pour un certain indice~$j_\ell$ ; par ailleurs, on choisit pour tout~$\ell$ une norme~$||.||''_\ell$ sur l'algèbre des fonctions analytiques sur~$X''_\ell$. 

Pour tout~$\ell$, il existe un réel positif~$C_\ell$ tel que~$||g_{|X''_\ell}||''_\ell\leq C_\ell ||g||'_{j_\ell}$ pour toute fonction analytique~$g$ sur~$X'_{j_\ell}$.  Le théorème d'acyclicité de Tate assure d'autre part l'existence d'un réel positif~$D_i$ tel que~$||g||_i\leq D_i\sup_\ell ||g_{|X''_\ell}||''_\ell$  pour toute fonction analytique~$g$ sur~$X_i$. On en déduit que l'on a pour toute fonction~$f$ appartenant à~$\sch A$ les inégalités~$$||f_{|X_i}||\leq D_i\sup_\ell (C_\ell ||f_{|X'_{j_\ell}}||'_{j_\ell})\leq D_i(\sup_ \ell C_\ell )||f||'.$$ Il existe donc pour tout ~$i$ une constante positive~$\Delta_i$ telle que~$||f_{|X_i}||\leq \Delta_i||f||'$ pour toute~$f\in \sch A$ ; par conséquent,~$||f||\leq (\sup_i \Delta_i) ||f||'$ pour  toute~$f\in \sch A$ ; il s'ensuit, par symétrie de leurs rôles, que ~$||.||$ et~$||.||'$ sont équivalentes, comme annoncé.

\deux{banachimpl} Dans ce qui suit, on considèrera systématiquement l'algèbre des fonctions analytiques d'un espace~$k$-analytique quasi-compact~$X$ comme munie de la (classe d'équivalence de la) norme définie par un recouvrement affinoïde fini quelconque de~$X$.

\deux{prepxladm} Soit~$X$ un espace~$k$-analytique quasi-compact. 

\trois{alxladm} Soit~$L$ une extension complète de~$k$. Si~$\sch B$ désigne l'algèbre de Banach des fonctions analytique sur~$X_L$, il existe une flèche naturelle~${\sch A}_L\hookrightarrow \sch B$ ; on vérifie immédiatement, en le testant sur un recouvrement affinoïde (et à l'aide de~\ref{remadmkl}) que cette flèche est injective et admissible. 

\trois{alxliso} Supposons de plus que~$L$ possède une base topologique sur~$k$, c'est-à-dire qu'il existe une famille~$(e_i)$ d'éléments de~$L$ tel que~$L$ s'identifie, comme~$k$-espace de Banach, à la somme directe complétée~$\wid {\bigoplus k.e_i}$. L'injection admissible~${\sch A}_L\hookrightarrow \sch B$ est alors un isomorphisme. Pour le voir, donnons-nous une fonction~$f$ appartenant à~$\sch B$. Pour tout domaine affinoïde~$V$ de~$X$, la restriction de~$f$ à~$V_L$ possède une {\em unique} écriture de la forme~$\sum a_i e_i$ , où~$(a_i)$ est une famille de fonctions analytiques sur~$V$ satisfaisant les conditions de convergences usuelles. L'unicité évoquée assure le recollement coefficient par coefficient de cette décomposition ; il existe donc une famille~$(a_i)$ de fonctions analytiques sur~$X$ telles que pour tout domaine affinoïde~$V$ de~$X$ la famille~$(a_{i|V})$ satisfasse les conditions de convergence habituelles et vérifie l'égalité~$f_{|V}=\sum a_{i|V}e_i$. La considération d'un recouvrement affinoïde fini de~$X$ assure que la famille~$(a_i)$ satisfait elle-même les conditions de convergence habituelles, et l'on a alors (par vérification G-locale) l'égalité~$f=\sum a_ie_i$, d'où notre assertion.

\deux{affkr}{\bf Proposition.} {\em Soit~$X$ un espace~$k$-analytique et soit~$\bf r$ un polyrayon~$k$-libre. L'espace~$X_{\bf r}:=X_{k_{\bf r}}$ est~$k_{\bf r }$-affinoïde si et seulement si~$X$ est~$k$-affinoïde}. 

\medskip
{\em Démonstration.}  Seule l'implication réciproque requiert une preuve ; supposons donc que~$X_{\bf r}$ est~$k_{\bf r}$-affinoïde ; il est alors compact, ce qui entraîne la quasi-compacité de~$X$. Soit~$\sch A$ (resp.~$\sch B$) l'algèbre de Banach des fonctions analytiques~$X$ (resp.~$X_{\bf r}$). En vertu du~\ref{alxliso}, l'algèbre~$\sch B$ s'identifie à~${\sch A}_{\bf r}$ ; comme~$\sch B$ est~$k_{\bf r}$-affinoïde,~$\sch A$ est~$k$-affinoïde. 

On dispose alors d'un morphisme naturel~$\phi : X\to {\sch M}({\sch A})$ ; le morphisme~$\phi_{\bf r}$ est la composée de~$X_{\bf r}\to \sch M(\sch B)$ et~$\sch M(\sch B)\to \sch M(\sch A_{\bf r})$ ; or la première de ces flèches est un isomorphisme car~$X_{\bf r}$ est par hypothèse~$k_{\bf r}$-affinoïde, et on a vu que la seconde en est un aussi ; par conséquent,~$\phi_{\bf r}$ est un isomorphisme ; le corollaire~\ref{desciso} assure alors que~$\phi$ est un isomorphisme, et partant que~$X$ est~$k$-affinoïde.~$\Box$ 

\deux{affrad}{\bf Théorème.} {\em On fait l'hypothèse que~$k$ est de caractéristique non nulle (celle-ci coïncide alors avec l'exposant caractéristique~$p$ de~$\red k$). Soit~$X$ un espace~$k$-analytique. Supposons qu'il existe un entier~$n$ tel que~$X_{k^{1/p^n}}$ soit~$k^{1/p^n}$-affinoïde ; l'espace~$X$ est alors~$k$-affinoïde.}

\medskip
{\em Démonstration.} On se ramène par une récurrence triviale sur~$n$ au cas où~$n=1$ ; par ailleurs, la proposition~\ref{affkr} ci-dessus permet, quitte à étendre les scalaires à~$k_{\bf r}$ pour un certain polyrayon~$k$-libre~$\bf r$ bien choisi, de supposer que~$|k\ti|\neq 1$ et que~$X_{k^{1/p}}$ est strictement~$k^{1/p}$-affinoïde. Notons~$\sch A$ (resp.~$\sch B$) l'algèbre des fonctions analytiques sur~$X$ (resp.~$X_{k^{1/p}}$). Comme~$X_{k^{1/p}}$ est compact, il découle de~\ref{remtopsep} que~$X$ est compact. La définition de~$X$ ne fait dès lors intervenir qu'une quantité (au plus) dénombrable de paramètres ; il en va de même de l'isomorphisme entre~$X_{k^1/p}$ et le spectre (analytique) d'une algèbre strictement~$k^{1/p}$-affinoïde ; par conséquent, on peut supposer que~$k$ est topologiquement de type dénombrable sur son sous-corps complet premier (qui est~$\FF_p$ muni de la valuation triviale). 

Grâce à cette hypothèse, le corps~$k^{1/p}$ possède une base topologique sur~$k$ (que l'on peut choisir comme étant une~$p$-base topologique). Il résulte alors du~\ref{alxliso} que l'injection admissible~${\sch A}_{k^{1/p}}\to \sch B$ est un isomorphisme. 

\medskip
L'algèbre~$\sch B$ est strictement ~$k^{1/p}$-affinoïde ; on choisit une surjection admissible~$k^{1/p}\{T_1,\ldots,T_n\}\to \sch B$. Pour tout~$i$, on désigne par~$f_i$ l'image de~$T_i$ dans~$\sch B$. Si~$(g_1,\dots,g_n)$ est un~$n$-uplet d'éléments de~$\sch B$ tel que~$||g_i-f_i||$ soit suffisamment petit pour tout~$i$, le morphisme de ~$k^{1/p}\{T_1,\ldots,T_n\}$ dans~$\sch B$ qui envoie chacun des~$T_i$ sur~$g_i$ est encore une surjection admissible ; on peut donc, quitte à perturber un peu les~$f_i$, supposer qu'il existe une extension finie~$F$ de~$k$ contenue dans~$k^{1/p}$ telle que~$f_i$ appartienne pour tout~$i$ à~${\sch A}_F$, vue comme une sous-algèbre de~$\sch B$ {\em via} l'isomorphisme~${\sch A}_{k^{1/p}}\simeq \sch B$. 

\medskip
Le morphisme~$F\{T_1,\ldots,T_n\}\to {\sch A}_F$ qui envoie chaque~$T_i$ sur~$f_i$ devient par construction une surjection admissible après extension des scalaires de~$F$ à~$k^{1/p}$. Comme~$k$ est topologiquement de type dénombrable sur son sous-corps premier,~$k^{1/p}$ possède une base topologique sur~$F$ (que l'on peut choisir comme étant une~$p$-base topologique) ; il s'ensuit que~$F\{T_1,\ldots,T_n\}\to {\sch A}_F$ est lui-même une surjection admissible, et donc que~${\sch A}_F$ est strictement~$F$-affinoïde. 

\medskip
L'algèbre~${\sch A}_F$ étant strictement~$F$-affinoïde, elle est en particulier strictement~$k$-affinoïde. L'injection~${\sch A}\hookrightarrow {\sch A}_F$ est finie ; de plus, si~$a$ est un élément de~${\sch A}_F$ de norme spectrale majorée par~$1$, alors~$a^p$ est un élément de~${\sch A}$ dont la norme spectrale dans~$\sch A$ coïncide avec la norme spectrale dans~${\sch A}_F$ (\ref{remadmkl}) et est donc majorée par~$1$ ; il s'ensuit que~$\red {\sch A} \to \red{\sch A_F}$ est entier, puis que~$\sch A$ est elle-même strictement~$k$-affinoïde.

\medskip
On dispose alors d'un morphisme naturel~$\phi : X\to {\sch M}({\sch A})$ ; le morphisme~$\phi_{k^{1/p}}$ est la composée de~$X_{k^{1/p}}\to \sch M(\sch B)$ et~$\sch M(\sch B)\to \sch M(\sch A_{k^{1/p}})$ ; or la première de ces flèches est un isomorphisme car~$X_{k^{1/p}}$ est par hypothèse~$k^{1/p}$-affinoïde, et on a vu que la seconde en est un aussi ; par conséquent,~$\phi_{k^{1/p}}$ est un isomorphisme ; le corollaire~\ref{desciso} assure alors que~$\phi$ est un isomorphisme, et partant que~$X$ est~$k$-affinoïde.~$\Box$ 

\deux{affred} {\bf Théorème}. {\em Soit~$X$ un espace~$k$-analytique et soit~$\sch I$ un faisceau cohérent d'idéaux sur~$X$ dont  les éléments sont localement nilpotents ; soit~$Y$ le sous-espace analytique fermé de~$X$ défini par~$\sch I$. L'espace~$X$ est~$k$-affinoïde si et seulement si~$Y$ est~$k$-affinoïde.}

 \medskip
 {\em Démonstration.} La proposition~\ref{affkr} ci-dessus permet, quitte à étendre les scalaires à~$k_{\bf r}$ pour n'importe quel~$\bf r$ non vide, de supposer que~$|k\ti|\neq 1$. L'implication directe est évidente ; nous allons établir la réciproque, et nous supposons donc que~$Y$ est~$k$-affinoïde. Comme~$Y$ et~$X$ ont le même espace topologique sous-jacent,~$X$ est compact ; il existe donc un entier~$n\geq 1$ tel que~$\sch I^n=0$ ; on se ramène par une récurrence immédiate sur~$n$ au cas où~$n=2$. On appelle~$\sch A$ (resp.~$\sch B$) l'algèbre de Banach des fonctions analytiques sur~$X$ (resp.~$Y$).
 
 \trois{fascoh} {\em La cohomologie de tout faisceau cohérent sur~$X\grot$ est triviale en degrés strictement positifs.} En effet, si~$\sch F$ est un tel faisceau, il s'insère dans une suite exacte~$0\to \sch I.\sch F\to  \sch F\to \sch F/(\sch I.\sch F)\to 0$ ; les deux termes extrêmes de cette suite sont annulés par~$\sch I$, et peuvent donc être vus comme des faisceaux cohérents sur~$Y$ ; comme~$Y$ est affinoïde, leur cohomologie est triviale en degrés strictement positifs, et il en va donc de même de celle de~$\sch F$. 
 
 \trois{secglobai} {\em Si~$I$ désigne l'espace des sections globales de~$\sch I$ alors~$I$ est un idéal de type fini de~$\sch A$, et la flèche naturelle~$\sch A/I\to \sch B$ est un isomorphisme.} En effet, l'on dispose d'une suite exacte~$0\to \sch I\to \sch O_{X\grot}\to \sch O_{X\grot}/\sch I\to 0$ ; celle-ci induit, en vertu du~\ref{fascoh}, une suite exacte~$$0\to \H^0(X,\sch I)\to \H^0(X,\sch O_{X\grot})\to \H^0(X,\sch O_{X\grot}/\sch I)\to 0,~$$ que l'on peut récrire~$0\to I\to \sch A\to \sch B\to 0$, d'où notre seconde assertion.
 
 En ce qui concerne la première, nous allons en fait montrer la n\oe thérianité de~$\sch A$. Soit~$(f_m)$ une suite d'éléments de~$\sch A$ ; désignons pour tout~$m$ par~$J_m$ l'idéal de~$\sch A$ engendré par~$f_1,\ldots,f_m$ ; nous allons montrer que la suite~$(J_m)$ est stationnaire, ce qui suffira à conclure. Pour tout~$m$, notons~$\sch J_m$ l'image (faisceautique) du morphisme~${\sch O}_{X\grot}^m\to \sch O_{X\grot}$ défini par les~$f_i$ pour~$i$ variant de~$1$ à~$m$. Si~$V$ est un domaine affinoïde de~$X$, la nœthérianité de son algèbre de fonctions et la théorie des faisceaux cohérents sur~$V$ entraînent la stationnarité de la suite~$(\sch J_{m|V})$ ; comme~$X$ est compact, il peut être recouvert par un nombre fini de domaines affinoïdes ; la suite~$(\sch J_m)$ est donc stationnaire. Or il résulte de~\ref{fascoh} que l'on a pour tout~$m$ l'égalité~$J_m=\H^0(X,\sch J_m)$ ; par conséquent,~$(J_m)$ est stationnaire, ce qu'on souhaitait établir. 
 
 \trois{concluaffred} {\em Conclusion.} Comme~$|k\ti|\neq \{1\}$ la surjection~$\sch A\to \sch A/I\simeq \sch B$ est admissible et induit donc,~$I$ étant nilpotent, un homéomorphisme entre les {\em espaces topologiques}~$\sch M(\sch A)$ et~$\sch M(\sch B)$ ; elle préserve de ce fait les semi-normes spectrales. Choisissons une surjection admissible~$k\{T_1/r_1,\ldots,T_\ell/r_\ell\}\to \sch B$ ; pour tout~$i$, donnons-nous un antécédent~$f_i$ dans~$\sch A$ de l'image de~$T_i$ dans~$\sch B$ ; le rayon spectral de l'image en question est majoré par~$r_i$, et celui de~$f_i$ l'est donc aussi ; par conséquent,~$k\{T_1/r_1,\ldots,T_\ell/r_\ell\}\to \sch B$ se relève en un morphisme~$\lambda : k\{T_1/r_1,\ldots,T_\ell/r_\ell\}\to \sch A$. Soit~$(i_1,\ldots,i_m)$ une famille génératrice de~$I$ et soit~$a\in \sch A$. Comme~$k\{T_1/r_1,\ldots,T_\ell/r_\ell\}\to \sch B$ est surjective, il existe~$i\in I$ et~$a'\in k\{T_1/r_1,\ldots,T_\ell/r_\ell\}$ tels que~$a=\lambda(a')+i$ ; on peut écrire~$i=a_1i_1+\ldots+a_mi_m$, où les~$a_j$ sont éléments de~$\sch A$. Par le même raisonnement que précédemment, il existe pour tout~$j$ un élément~$a'_j$ de ~$k\{T_1/r_1,\ldots,T_\ell/r_\ell\}$ tel que~$a=\lambda(a'_j)$ modulo~$I$. Compte-tenu du fait que~$I$ est de carré nul, on a finalement~$a=\lambda(a')+\sum \lambda(a'_j)i_j$. 
 
 \medskip
 Chacun des~$i_j$ est nilpotent ; l'application~$\lambda$ se prolonge donc en un morphisme~$$k\{T_1/r_1,\ldots,T_\ell/r_\ell,S_1,\ldots,S_m\}\to \sch A$$ qui envoie~$S_j$ sur~$i_j$ pour tout~$j$ ; en vertu de ce qui précède, ce morphisme est une surjection, nécessairement admissible puisque~$|k\ti|\neq\{1\}$ ; la~$k$-algèbre de Banach~$\sch A$ est ainsi~$k$-affinoïde. 
 
 \medskip
 On dispose d'un morphisme naturel d'espaces~$k$-analytiques de~$X$ vers~$\sch M(\sch A)$, qui s'insère dans un diagramme commutatif~$$\diagram Y\rto \dto&\sch M(\sch B)\dto\\X\rto&\sch M(\sch A)\enddiagram$$ dans lequel les flèches verticales sont des immersions fermées données par un faisceau d'idéaux nilpotents, et dont la flèche horizontale supérieure est un isomorphisme puisque~$Y$ est affinoïde. Le morphisme composé de~$Y$ vers~$\sch M(\sch A)$ est sans bord ; en raison de la surjectivité de~$Y\to X$, le morphisme~$X\to \sch M(\sch A)$ est lui-même sans bord ; c'est par ailleurs un homéomorphisme puisque c'est le cas des trois autres flèches du diagramme. Par conséquent,~$X\to \sch M(\sch A)$ est fini (\ref{finpropfin}), ce qui implique que~$X$ est~$k$-affinoïde (et comme~$\sch A=\H^0(X,\sch O_{X\grot})$ la flèche~$X\to \sch M(\sch A)$ est  alors un isomorphisme).~$\Box$

\section{Les courbes analytiques : premières propriétés}

\subsection*{Classification des points} 

\deux{courban} Une {\em courbe} algébrique (resp.~$k$-analytique) sur~$k$ est
un~$k$-schéma séparé et localement de type fini (resp. un espace~$k$-analytique séparé) 
{\em purement}
de dimension 1. 

\deux{classepts} Soit~$x$ un point d'une courbe~$k$-analytique. Comme~$X$ 
est de dimension~$1$, le degré de transcendance
du corpoïde~$\red{\hres(x)}$ sur~$\red k$ est inférieur ou égal
à~$1$. De plus, s'il est égal à~$1$ alors~$\red{\hres(x)}$ est de type fini
sur~$\red k$, ce qui revient à dire que~$\red{\hres(x)}_1$ est de type fini 
sur~$\red k_1$, et que~$|\hres(x)\ti|/|k\ti|$ est de type fini
(\ref{interpfinioide}). 

\medskip
On se trouve donc dans l'un des quatre cas suivants, exclusifs l'un de l'autre. 

\medskip
$\bullet$ Le corps~$\hres(x)$ est une extension presque algébrique de~$k$ ; on dit alors que~$x$ est {\em de type 1} ; tout point rigide est de type 1, et la réciproque est vraie si~$k$ est algébriquement clos ou trivialement valué. 

$\bullet$ Le corps~$\red{\hres(x)}_1$ est de type fini et de degré de transcendance~$1$ sur~$\red k_1$, et~$|\hres(x)\ti|/|k\ti|$ est fini ; on dit alors que~$x$ est {\em de type 2}. 

$\bullet$ Le corps~$\red{\hres(x)}_1$ est une extension finie de~$\red k_1$, et ~$|\hres(x)\ti|/|k\ti|$ est de type fini et de rang rationnel égal à~$1$ ; on dit alors que~$x$ est {\em de type 3} ;

$\bullet$ Le corps~$\hres(x)$ n'est pas une extension presque algébrique de~$k$, mais admet un plongement isométrique dans une extension immédiate de~$\KK$ ; on dit alors que~$x$ est {\em de type 4}.

\medskip
Notons que~$x$ est de type~$1$ ou~$4$ (resp.~$2$ ou~$3$) si et seulement si le
corpoïde~$\red{\hres(x)}$ est algébrique (resp. de type fini et de
degré de transcendance~$1$) sur~$\red k$. 

\deux{extype} {\em Exemple.} Soit~$a\in  k\ti$ et soit~$r\in \RR\ti_+$.
Il résulte
de~\ref{remresidgauss}
que~$\red{T-a}$ est
transcendant sur~$\red k$ et
que~$\red{\hres(\eta_{a,r})}=\red k(\red{T-a})$ ; 
en conséquence, $\eta_{a,r}$ est de type 2 ou 3. 

\medskip
On a par ailleurs par définition
de~$\eta_{a,r}$ 
l'égalité~$|\hres(\eta_{a,r})\ti|=|k\ti|\cdot r^{\ZZ}$ ; le quotient
$|\hres(\eta_{a,r})\ti|/|k\ti|$ est 
donc engendré par~$r$. 

\trois{extype-2}
{\em Supposons que~$r\in |k\ti|^{\QQ}$.} Soit~$m$ l'ordre de~$r$
modulo~$|k\ti|$. 
Le groupe~$|\hres(\eta_{a,r})\ti|/|k\ti|$ est cyclique d'ordre~$m$, et~$\eta_{a,r}$
est donc de type 2. En vertu de~\ref{abhyankaroide},
$\tau:=\red{\frac{(T-a)^m}\lambda}$ est transcendant sur~$\red k_1$, et
~$\red{\hres(\eta_{a,r})}_1=\red k_1(\tau)$. 

\trois{extype-3}
{\em Supposons que~$r\notin |k\ti|^{\QQ}$.} 
Le groupe~$|\hres(\eta_{a,r})\ti|/|k\ti|$ est 
alors libre de rang~$1$, et~$\eta_{a,r}$
est donc de type 3. En vertu de~\ref{abhyankaroide},
~$\red{\hres(\eta_{a,r})}_1=\red k_1$. 

\deux{nottypepts} Si~$X$ est une courbe~$k$-analytique et si~$Y$ est un sous-ensemble de~$X$ on notera~$Y\typ 0~$ l'ensemble des points rigides de~$Y$ et~$Y\typ 1$ (resp.~$Y\typ 2$, resp. ....) l'ensemble de ses points de type~$1$ (resp.~$2$, resp. ...). On désignera par~$Y\dtr$ la réunion de~$Y\typ 2$ et~$Y\typ 3$, par~$Y\geom$ la réunion de~$Y\typ 0$,~$Y\typ 2$, et~$Y\typ 3$, etc.

\deux{corollgener} Soit~$X$ une courbe~$k$-analytique et soit~$x\in X\dtr$ ; nous allons dire quelques mots de~$\red{\hres(x)}$ et~$\PP_{\red{\hres(x)}/\red k}$. 

\trois{groupeval23gen} Commençons par une remarque générale. Soit~$\val .$
une valuation non triviale appartenant à~$\PP_{\red{\hres(x)}/\red k}$
et soit~$\kappa$ son corpoïde résiduel. Le 
groupe~$\val {\red{\hres(x)}\ti}$ est alors de rang rationnel au moins~$1$.
Comme~$\red{\hres(x)}$ est de type fini
et de degré de transcendance~$1$ sur~$\red k$, on déduit 
de~\ref{interpdegtrres} 
et~\ref{abhyankarfini}  
les faits suivants : 

\medskip
$\bullet$ le groupe~$\val {\hres(x)\ti}$ est libre de rang~$1$ ; 

$\bullet$ le corpoïde~$\kappa_{\RR\ti_+\times\{1\}}$ est une extension finie de~$\red k$.

\trois{valxideux}{\em Supposons que~$x\in X\typ 2$.} Le corps résiduel classique~$\red {\hres(x)}_1$ est alors de la forme~$\red k_1({\sch C})$, où~$\sch C$ est une~$\red k_1$-courbe projective, normale et intègre (mais pas nécessairement géométriquement intègre), qui est
uniquement déterminée et que l'on appellera {\em la courbe résiduelle} en~$x$ ; le groupe~$|\hres(x)\ti|/|k\ti|$ est fini.

\medskip
Le lemme~\ref{memegroupebijval} 
assure 
que~$\PP_{\red{\hres(x)}/\red k}\to \PP_{\red{\hres(x)}_1/\red k_1}$ est un homéomorphisme (préservant les ouverts affines). 
 L'espace topologique~$\PP_{\red{\hres(x)_1}/\red k_1}=\PP_{\red k_1 (\sch C)/\red k_1}$
est lui-même homéomorphe à l'espace topologique sous-jacent au schéma~$\sch C$ ; au point générique correspond la valuation triviale, et à chaque point fermé une valuation discrète.

%

\trois{valxitrois}{\em Supposons que~$x\in X\typ 3$}. Le groupe~$|\hres(x)\ti|/|k\ti|$ est alors de rang rationnel égal à~$1$, et~$\red {\hres(x)}_1$ est une extension finie de~$\red k$. Soit~$r$ un élément de~$|\hres(x)\ti|$ qui n'appartient pas à~$\sqrt{ |k\ti|}$, et soit~$\tau\in \red {\hres(x)}\ti$ un élément
de degré~$r$. Désignons
par~$F$
le sous-corpoïde de~$\red{\hres(x)}$ engendré par~$\red k$
et~$\red{\hres(x)}_1$ ; son groupe des degrés est égal à~$|k\ti|$, et l'annéloïde~$F[\tau, \tau^{-1}]$ est dès lors un corpoïde. 

Par construction,~$F$ est fini sur~$\red k$ ; il vient~$\PP_{\red{\hres(x)}/\red k}=\PP_{\red {\hres(x)}/F}$. Par ailleurs, le groupe~$\deg(\red{\hres(x)}\ti)/\deg(F\ti)$ est de 
torsion par choix de~$\tau$ ; on déduit alors du lemme~\ref{torsionhomeozr}
que
$$\PP_{\red{\hres(x)}/\red k}=\PP_{\red{\hres(x)}/F}\simeq \PP_{F[\tau,\tau\inv]/F}.$$

\medskip
L'annéloïde d'une valuation
sur~$F[\tau,\tau^{-1}]$ contient
nécessairement~$\tau$ ou~$\tau^{-1}$. Il s'ensuit aussitôt 
que~$\PP_{F[\tau,\tau^{-1}]}$ compte trois éléments : la valuation triviale, celle dont l'anneau est~$F[\tau]$, et celle dont l'anneau est~$F[\tau^{-1}]$. Les deux dernières sont des points fermés de~$\PP_{F[\tau,\tau^{-1}]}$, et la première en est un point ouvert et
dense.

\subsection*{Bonté des courbes analytiques} 

\deux{lemomega} {\bf Lemme.} {\em Soit~$X$ une courbe~$k$-analytique et soit~$x\in X\dtr$ ; soit~$\langle.\rangle$ une~$\red k$-valuation non triviale de~$\red {\hres(x)}$. Il existe un élément~$\omega$ de~$\red{\hres(x)}$ tel que~$\langle \omega \rangle>1$ et tel que~$\langle.\rangle$ soit la seule~$\red k$-valuation de~$\red{\hres(x)}$ à posséder cette propriété. }

\medskip
{\em Démonstration.} On distingue deux cas. 

\medskip
{\em Supposons que~$x\in X\typ 2$}.  On reprend les notations du~\ref{valxideux} ; la valuation ~$\langle .\rangle$ est induite par une valuation discrète de~$\red k_1({\sch C})$, elle-même définie par un point fermé~$\sch P$ de~$\sch C$. Le théorème de Riemann-Roch assure que~$N\sch P$ est engendré par ses sections globales pour~$N$ assez grand, ce qui fournit une fonction rationnelle sur~$\sch C$ ayant un et un seul pôle situé en~$\sch P$ ; on peut prendre pour~$\omega$ une telle fonction. 

\medskip
{\em Supposons que~$x\in X\typ 3$}.  On reprend les notations du~\ref{valxitrois} ; la valuation~$\langle .\rangle$ est ou bien la seule~$\red k$-valuation de~$\red{\hres(x)}$ dont l'anneau ne contient pas~$\tau$, ou bien la seule dont l'anneau ne contient pas~$\tau\inv$ ; on peut donc prendre~$\omega=\tau$ dans le premier cas, et~$\omega=\tau\inv$ dans le second. ~$\Box$ 

\deux{courbesbonnes} {\bf Proposition.} {\em Toute courbe~$k$-analytique est un bon espace.}

\medskip
{\em Démonstration.} Soit~$X$ une courbe~$k$-analytique et soit~$x$ un point de~$X$. Si~$x\in X\typ {1,4}$ il appartient à l'intérieur de~$X$, et~$(X,x)$ est donc bon. Supposons que~$x\in X\dtr$. Comme~$X$ est par définition séparée,~$\red{(X,x)}$ est un ouvert quasi-compact et non vide de~$\PP_{\red{\hres(x)}/\red k}$ ; il résulte de la forme explicite de ce dernier (\ref{valxideux} et~\ref{valxitrois}) que~$\red{(X,x)}$ est de la forme~$\PP_{\red{\hres(x)}/\red k}\setminus\{\langle.\rangle_1,\ldots,\langle.\rangle_r\}$ où les~$\langle.\rangle_i$ sont des valuations non triviales, ou encore fermées en tant que points de~$\PP_{\red{\hres(x)}/\red k}$. Pour chaque~$i$, il existe en vertu du lemme~\ref{lemomega} ci-dessus un élément~$\omega_i$ de~$\red{\hres(x)}$ tel que~$\langle\omega_i\rangle_i>1$ et telle que~$\langle.\rangle_i$ soit la seule~$\red k$-valuation de~$\red{\hres(x)}$ à posséder cette propriété. L'ouvert~$\red{(X,x)}$ peut dès lors se décrire comme étant égal à~$\PP_{\red{\hres(x)}/\red k}\{\omega_1,\ldots,\omega_r\}$ ; il est donc affine, ce qui équivaut à la bonté de~$(X,x)$.~$\Box$

\section {Étude détaillée de la droite projective}

\deux{ordrea1} On munit~$\Aff^{1,\rm an}_k$ de l'ordre partiel pour lequel~$x\leq y$ si~$|f(x)|\leq |f(y)|$ pour tout~$f\in k[T]$. Toute chaîne non vide~$(x_i)_i$ de points de~$\Aff^{1,\rm an}_k$ admet une borne inférieure, à savoir le point défini par la semi-norme~$f\mapsto \inf |f(x_i)|$. 

\deux{etaardar} Soit~$a\in k$ et soit~$r\geq 0$. Si~$x\in \Aff^{1,\rm an}_k$ les assertions suivantes sont équivalentes : 

\medskip
i)~$x\leq \eta_{a,r}$ ; 

ii)~$|(T-a)(x)|\leq r$. 

\medskip
En effet, i)$\Rightarrow$ ii) est clair. Supposons que ii) soit vraie, et soit~$f\in k[T]$ ; écrivons~$f=\sum a_i(T-a)^i$. On a alors~$$|f(x)|\leq \max |a_i|\cdot |(T-a)(x)|^i\leq \max |a_i|r^i=|f(\eta_{a,r})|,$$ d'où i).

\deux{compareta} Soient~$a$ et~$b$ deux éléments de~$k$ et~$r$ et~$s$ deux réels positifs. Les assertions suivantes sont équivalentes : 

\medskip
i)~$\eta_{b,s}\leq \eta_{a,r}$ ; 

ii)~$\max(|b-a|,s)\leq r$. 

\medskip
En effet, en vertu de l'égalité~$T-a=T-b+b-a$ on a~$$|(T-a)(\eta_{b,s})|=\max(|b-a|, s),$$ et l'équivalence requise provient dès lors du~\ref{etaardar} ci-dessus, que l'on applique avec~$x=\eta_{b,s}$. 

\medskip
Remarquons une conséquence de ce qui précède :~$\eta_{a,r}=\eta_{b,s}$ si et seulement~$r=s$ et~$|a-b|\leq r$. 

\deux{dar} Soit~$a\in k$ et soit~$r> 0$. On notera~$\DD(a,r)$ (resp.~$\DD\zero (a,r)$) le disque fermé (resp. ouvert) de Berkovich de centre~$a$ et de rayon~$r$,  c'est-à-dire le sous-ensemble de~$\Aff^{1,{\rm an}}_k$ défini par l'inégalité~$|T-a|\leq r$ (resp.~$|T-a|<r$). Si~$L$ est une extension complète de~$k$, l'ensemble des~$L$-points de~$\DD(a,r)$ (resp.~$\DD\zero(a,r)$) est en bijection avec le disque ouvert (resp. fermé) de centre~$a$ et de rayon~$r$ de l'espace métrique~$L$. 

\medskip
Le disque fermé~$\DD(a,r)$ est un domaine affinoïde de~$\Aff^{1,{\rm an}}_k$, d'algèbre associée~$k\{(T-a)/r\}$. En vertu de~\ref{etaardar}, il peut être défini comme l'ensemble des~$x$ majorés par~$\eta_{a,r}$. La norme de~$k\{(T-a)/r\}$ est induite par~$\eta_{a,r}$, et le bord de Shilov de~$\DD(a,r)$ est donc~$\{\eta_{a,r}\}$.  

\medskip
Le disque ouvert~$\DD(a,r)$ est un ouvert de~$\Aff^{1,{\rm an}}$, égal à la réunion des~$\DD(a,s)$ pour~$s\in ]0;r[$. 

\deux{includiscberk} Soient~$a$ et~$b$ dans~$k$ et soient~$r$ et~$s$ deux réels strictement positifs. 

\trois{incluclass} Comme~$\DD(a,r)$ (resp.~$\DD(b,s)$) est par ce qui précède l'ensemble des points de~$\Aff^{1,\rm an}$
majorés par~$\eta_{a,r}$ (resp.~$\eta_{b,s}$), on a ~$\DD(b,s)\subset  \DD(a,r)$ si et seulement si~$\eta_{b,s}\leq \eta_{a,r}$, c'est-à-dire encore, 
d'après~\ref{compareta}, si et seulement si~$\max(|b-a|,s)\leq r$.

\trois{interinclub} Si~$\DD(b,s)\cap  \DD(a,r)\neq \emptyset~$ alors~$\DD(b,s)$ et~$\DD(a,r)$ sont comparables pour l'inclusion. En effet, d'après notre hypothèse il existe une extension complète~$L$ de~$k$ telle que~$\DD(b,s)(L)\cap  \DD(a,r)(L)\neq \emptyset$. Choisissons~$\lambda$ dans cette intersection ; on a~$|\lambda-a|\leq r$ et~$|\lambda-b|\leq s$, d'où~$|b-a|\leq \max(r,s)$ ; on conclut à l'aide de~\ref{incluclass}. 

\trois{interdiscouvber} En écrivant un disque ouvert comme réunion croissante de disques fermés, on déduit de ce qui précède les faits suivants : 

\medskip
i) On a ~$\DD\zero(b,s)\subset  \DD\zero(a,r)$ si et seulement si~$|a-b|<r$ et~$s\leq r$. 

\medskip
ii) Si~$\DD\zero(b,s)\cap  \DD\zero(a,r)\neq \emptyset~$ alors~$\DD\zero(b,s)$ et~$\DD\zero(a,r)$ sont comparables pour l'inclusion. 

\deux{disquespascomp} Soient~$a$ et~$b$ deux éléments de~$k$ et soient~$r$ et~$s$ deux réels strictement positifs, et supposons que les boules~$\DD(a,r)$ et~$\DD(b,s)$ soient disjointes. Soit~$x\in \DD(a,r)$ et soit~$y\in \DD(b,s)$. Comme~$\DD(a,r)$ est l'ensemble des éléments majorés par~$\eta_{a,r}$, tout élément majoré par ~$x$ appartient enore à~$\DD(a,r)$ ; par conséquent,~$y$ n'est pas majoré par~$x$. Par symétrie,~$x$ n'est pas majoré par~$y$ et~$x$ et~$y$ sont ainsi incomparables. 

\medskip
L'assertion analogue pour les boules ouvertes s'en déduit immédiatement. 

\deux{invnormeconst} Soit~$f\in k[T]$, soit~$a\in k$ et soit~$r>0$. Si~$f$ ne s'annule pas sur~$\DD(a,r)$, elle est inversible sur ce dernier ; comme le bord de Shilov de~$\DD(a,r)$ est un singleton, cela entraîne que~$|f|$ est constante sur~$\DD(a,r)$. 

\medskip
Ce résultat s'étend immédiatement à~$\DD\zero(a,r)$, que l'on écrit comme réunion croissante de disques fermés. 

\deux{voisinfxi} Soit~$F$ une chaîne non vide de points de~$\Aff^{1,\rm an}_k$ ; on suppose que pour tout~$y\in F$ il existe~$a\in k$ et~$r>0$ tels que~$y=\eta_{a,r}$. Soir~$f$
un élément de~$k[T]$ ; posons~$x=\inf F$ et~$\rho=|f(x)|=\inf\limits_{y\in F} |f(y)|$. Soient~$\rho_-$ et~$\rho_+$ deux réels tels que~$\rho_-<\rho<\rho_+$.

\trois{voisinfxboulouv} Il existe~$y\in F$ tel que~$|f(y)|<\rho_+$. Soient~$a$ et~$r$ tels que le point~$y$ de~$F$ soit égal à~$\eta_{a,r}$ ; comme~$\eta_{a,r}=\inf\limits_{s>r}\eta_{a,s}$, il existe~$R>r$ tel que~$|f(\eta_{a,R})|<\rho_+$ ; comme on a l'encadrement~$x\leq y<\eta_{a,R}$, le point~$x$ appartient à~$\DD\zero(a,R)$. L'égalité~$|f(\eta_{a,R})|<\rho_+$ implique que~$|f|<\rho_+$ sur~$\DD(a,R)$.

\trois{voisinfxboulefer} Soit~$\rho'\in ]\rho_-;\rho[$. Nous allons maintenant montrer l'existence d'une famille finie~$((a_i,r_i))_i$ d'éléments de~$k\times \RR\ti_+$ possédant les propriétés suivantes : 

\medskip
$i)$ les boules fermées~$\DD(a_i,r_i)$ sont deux à deux disjointes et contenues dans~$\DD\zero(a,R)$ ; 

$ii)$ soit~$G$ une chaîne non vide de~$\Aff^{1,{\rm an}}_k$ constituée de points de la forme~$\eta_{b,s}$ avec~$b\in k$ et~$s\geq 0$, et soit~$\xi$ sa borne inférieure ; si~$\xi\in \DD\zero(a,R)-\coprod \DD(a_i,r_i)$ alors~$|f(\xi)|\geq \rho'$.

\medskip
L'assertion est évidente si~$f=0$ ; on peut donc supposer~$f\neq 0$. Commençons par remarquer qu'il suffit d'établir l'assertion ii) lorsque~$\xi$ est lui-même de la forme~$\eta_{b,s}$. Supposons en effet qu'elle ait été prouvée dans ce cas particulier, et soit~$G$ comme dans ii). Comme le point~$\xi$ appartient à~$\DD\zero(a,R)$, il existe~$t<R$ 
tel que~$|(T-a)(\xi)|<t$ ; par conséquent, il existe~$y\in G$ tel que~$|(T-a)(y)|<t$. Si~$z\in G$ alors~$\xi\leq z$, ce qui empêche~$z$ d'être majoré par l'un des~$\eta_{a_i,r_i}$. Ainsi, tout élément~$z$ de~$G$
majoré par~$y$ appartient à~$\DD\zero(a,t)-\coprod \DD(a_i,r_i)$, 
et {\em a fortiori}
à $\DD\zero(a,R)-\coprod \DD(a_i,r_i)$. Ceci entraîne, en vertu de notre hypothèse, que~$|f(z)|\geq \rho'$. Il s'ensuit que~$|f(\xi)|\geq \rho'$. 

\medskip
Il suffit donc de montrer l'existence de la famille~$((a_i,r_i))_i$ de sorte que ii) soit satisfaite lorsque~$\xi$ est lui-même de la forme~$\eta_{b,s}$. 

\medskip
On procède comme suit. Supposons tout d'abord que pour tout~$b\in k$ tel que~$|b-a|<R$, l'on ait~$|f(b)|\geq \rho'$. Comme~$s\mapsto |f(\eta_{b,s})|$ est croissante pour tout~$b$, on a~$|f(\eta_{b,s})|\geq \rho'$ pour tout~$s\geq 0$ et tout~$b\in \DD\zero(a,R)(k)$, et l'on peut prendre pour~$((a_i,r_i))_i$ la famille {\em vide}.  

\medskip
Supposons maintenant qu'il existe~$a_1\in k$ tel que~$|a_1-a|<R$ et tel que~$|f(a_1)|<\rho'$. L'application~$s\mapsto |f(\eta_{a_1,s})|$ étant croissante,  
et l'on a~$$|f(\eta_{a_1,R})|=|f(\eta_{a,R})|\geq \rho=|f(x)|.$$ Il existe donc~$r_1\in [0;R[$ tel que~$|f(\eta_{a_1,r_1})|=\rho'$ et qui est maximal pour cette propriété. 

On doit à nouveau distinguer deux cas. Le premier est celui dans lequel~$|f(b)|\geq \rho'$ pour tout~$b\in k$ tel que~$|b-a|<R$ et~$|b-a_1|>r_1$. Dans ce cas, soit~$(b,s)$ tel que~$\eta_{b,s}\in \DD\zero(a,R)-\DD(a_1,r_1)$. Si~$|b-a_1|\leq r_1$ alors~$s>r_1$, et l'on a donc~$|f(\eta_{b,s})|>\rho'$ par choix de~$r_1$ ; et si~$|b-a_1|>r_1$ on a alors ~$|f(\eta_{b,s})|\geq |f(\eta_{b,0})|=|f(b)|\geq \rho'$, et la famille à un élément~$(a_1,r_1)$ convient. 

Le second est celui où il existe~$a_2$ tel que~$|a_2-a_1|>r_1$,~$|a_2-a|<R$, et~$|f(a_2)|<\rho'$. Dans ce cas, il existe par le même raisonnement
 que ci-dessus un réel~$r_2\in ]0;R[$ tel que~$|f(\eta_{a_2,r_2})|=\rho'$ et qui est maximal pour cette propriété. 
 
 {\em Les boules~$\DD(a_1,r_1)$ et~$\DD(a_2,r_2)$ sont disjointes.} En effet, supposons le contraire. 
Comme~$|a_2-a_1|>r_1$, la boule~$\DD(a_2,r_2)$ n'est pas contenue dans~$\DD(a_1,r_1)$ : elle contient donc strictement~$\DD(a_1,r_1)$, ce qui implique que~$r_2>r_1$ et que~$\eta_{a_2,r_2}=\eta_{a_1,r_2}$ ; par conséquent,~$|f(\eta_{a_1,r_2})|=\rho'$, contredisant la définition de~$r_1$. 

\medskip
En poursuivant ce procédé, on construit de façon récursive une suite~$((a_i,r_i))_i$ d'éléments de~$k\times \RR\ti_+$  (finie ou non, {\em cf.~$\gamma)$} ) telle que : 

\medskip
$\alpha)$  les~$\DD(a_i,r_i)$ sont des boules fermées deux à deux disjointes contenues dans~$\DD(a,R)$ ; 

$\beta)$ pour tout~$i$, on a~$|f(a_i)|<\rho'$ et~$|f(\eta_{a_i,r_i})|=\rho'$ ;

$\gamma)$ la suite s'arrête au rang~$n$ si et seulement si~$|f(\eta_{b,s})|\geq a'$ pour tout~$(b,s)$ tel que~$\eta_{b,s}\in \DD\zero(a,R)-\coprod\limits_{i=1}^n\DD(a_i,r_i).$

\medskip
La condition~$\beta)$ entraîne que~$|f|$ n'est constante sur aucune des~$\DD(a_i,r_i)$, et partant que~$|f|$ n'est inversible sur aucune des~$\DD(a_i,r_i)$ (\ref{invnormeconst}). Chacune des boules~$\DD(a_i,r_i)$ contient donc au moins un point rigide en lequel~$f$ s'annule. Comme~$f\neq 0$, l'ensemble des points rigides en lesquels elle s'annule est fini. Par conséquent, la suite s'arrête à un certain rang~$n$, et la condition~$\gamma)$ garantit que la famille~$((a_1,r_1),\ldots, (a_n,r_n))$ satisfait les conditions requises.

\subsection*{L'arbre compact~$\wbeth k$ se plonge dans~$\pk$}

\deux{defplongwb} On déduit de~\ref{compareta} : que la formule~$\zeta_{a,r}\mapsto \eta_{a,r}$ définit sans ambiguïté une application injective~$\got t$ de~$\aleph(k)$ dans~$\Aff^{1,\rm an}_k$ (pour la définition de~$\aleph(k)$, {\em cf.}~\ref{relboules}) ;
et que~$\got t(x)\leq \got t(y)$ si et seulement si~$x\leq y$. Il s'ensuit que si~$F$ est une chaîne de~$\aleph(k)$ alors~$\got t (F)$ est une chaîne de~$\Aff^{1,\rm an}_k$ ; si de plus~$F$ est non vide,~$\inf\limits_{y\in F} \got t (y)$ est donc bien défini (\ref{ordrea1}). 

\medskip
Si~$a\in k$ et~$r\in \RR_+$ il est immédiat que~$\eta_{a,r}=\inf\limits_{R>r} \eta_{a,R}$. Autrement dit,~$\got t(\zeta_{a,r})=\inf\limits_{y\in ]\zeta_{a,r};\infty[}\got t(y)$. 

\medskip
On prolonge~$\got t$ à~$\beth(k)$ en posant~$\got t(x)=\inf\limits_{y\in ]x;\infty[}\got t (y)$ pour tout~$x$ appartenant à~$ \beth(E)-\aleph(E)$. Il résulte de la remarque précédente que l'égalité~$$\got t(x)=\inf\limits_{y\in ]x;\infty[}\got t(y)$$ vaut en réalité pour tout~$x\in \beth(k)$ ; ceci implique,~$\got t_{|\aleph(k)}$ étant croissante, que~$\got t$ est croissante. 

\deux{injgott} Nous allons montrer que l'application~$\got t : \beth(k)\to \pk$ est injective, et plus précisément qu'elle établit un {\em isomorphisme d'ensembles ordonnés} de~$\beth(k)$ sur son image. Sachant que~$\got t$ est croissante,  il reste à vérifier que si~$x$ et~$y$ sont deux points de~$\beth (k)$ tels que~$\got t(x)\leq \got t(y)$ alors~$x\leq y$. 

\medskip
On procède par l'absurde. Supposons tout d'abord que~$x>y$, et choisissons~$z\in ]y;x[$. Comme~$z$ et~$x$ sont strictement supérieurs à~$y$
ils appartiennent tous deux à~$\aleph(k)$ ; par conséquent,~$\got t(z)<\got t(x)$ ; comme~$\got t$ est croisssante,~$\got t(y)\leq \got t(z)$, d'où l'inégalité~$\got t(y)<\got t(x)$, contradictoire avec nos hypothèses. 

\medskip
Supposons maintenant que~$x$ et~$y$ soient non comparables ; cela signifie que~$x\wedge y>x$ et~$x\wedge y>y$. Il résulte alors de~\ref{wedgeboulouv} qu'il existe deux éléments~$a$ et~$b$ de~$k$ et deux réels strictement positifs~$r$ et~$s$ tels que~$x\in \bbouv a {r'}$, tels que~$y\in \bbouv b{ s'}$, et tels que~$\bbouv a {r'}\cap \bbouv b {s'}=\emptyset$. Comme~$\got t$ est croissante, les inégalités~$x\leq \zeta_{a,r'}$ et~$y\leq \zeta_{b,s'}$ impliquent que~$\got t(x)\leq \eta_{a,r'}$ et~$\got t(y)\leq \eta_{b,s'}$ ; autrement dit,~$\got t(x)\in \DD(a,r')$ et~$\got t(y)\in \DD(b,s')$. 

Puisque~$\bbfer a{r'}\cap \bbfer b{s'}=\emptyset$, on a~$|a-b|\geq \max(r',s')$ ; il s'ensuit que ~$\DD(a,r')\cap \DD(b,s')$ est également vide, et partant que~$\got t(x)$ et~$\got t(y)$ ne sont pas comparables (\ref{disquespascomp}). On aboutit ainsi à une contradiction. 

\deux{gottinvboul} On déduit de~\ref{injgott} que l'on a pour tout~$a\in k$ et tout~$r>0$ les égalités~$$\got t^{-1}(\DD(a,r))=\bbfer a r\;{\rm et}\; \got t^{-1}(\DD\zero(a,r))=\bbouv a r.$$

\deux{contgott} On étend~$\got t$ à~$\wbeth k$ en posant~$\got t (\infty)=\infty$ ; l'application~$\got t$ ainsi prolongée reste injective. 

\medskip
{\em L'application~$\got t$ est continue.} Sa continuité en tout point de~$\wbeth k$ provient de~\ref{voisinfxi} {\em et sq.} ainsi que du~\ref{gottinvboul} ci-dessus. Il reste à établir sa continuité en~$\infty$. Le point~$\infty$ de~$\pk$ admet une base de voisinages constitués des ouverts de la forme~$\pk-\DD(0,r)$ avec~$r>0$ ; or si~$r>0$ il résulte de~\ref{gottinvboul} que~$\got t^{-1}(\pk -\DD(0,r))=\wbeth k-\bbfer a r$, qui est un ouvert de~$\wbeth k$ ; ceci achève la démonstration. 

\medskip
Ainsi,~$\got t~$ apparaît comme une injection continue d'un espace topologique compact dans un espace topologique séparé. Par conséquent, elle établit un homéomorphisme entre l'arbre compact~$\wbeth k$ et son image~$\got t(\wbeth k)$, que l'on notera~$\got P(k)$. 
La fonction rayon~$\rho : \wbeth k \to [0;+\infty]$ induit {\em via	}
l'homéomorphisme~$\wbeth k\simeq
\got P(k)$
une fonction de~$\got P(k)$ dans~$[0;+\infty]$ que nous noterons encore~$\rho$.

\deux{remplongek} {\em Remarque.} L'espace métrique~$k$ se plonge 
topologiquement dans~$\wbeth k$ {\em via}
la flèche~$a\mapsto \zeta_{a,0}$. Composée avec~$\got t$, 
cette flèche induit un plongement topologique~$k\hookrightarrow \pk$, donné 
par la formule~$a\mapsto \eta_{a,0}$ ; il coïncide ainsi avec le plongement naturel~$k\hookrightarrow \PP^1(k)\hookrightarrow \PP^{1,\rm an}_k$.

\deux{propgottwb} {\bf Propriétés de l'arbre~$\got P(k)$ : le cas général}. 

\trois{gottwwbdesc} Il résulte de la définition de~$\got t$ que~$\got P(k)$ est la réunion de~$\PP^1(k)$, de l'ensemble des points de la forme~$\eta_{a,r}$ avec~$a\in k$ et~$r>0$, et enfin de l'ensemble des points de la forme~$\got t(x)$ avec~$x\in \beth(k)-\aleph(k)$. Ce dernier ensemble s'identifie, en vertu de~\ref{recapbijnast}, à celui des classes d'équivalence de chaînes évanescentes de boules fermées de~$k$. 

\trois{branchegottwb} Soit~$x\in \wbeth k$ et soit~$\xi$ le point~$\got t(x)$ de~$\pk$.  

\medskip
Supposons que~$x$ appartient à~$k\cup\{\infty\}$, c'est-à-dire que~$\xi\in \PP^1(k)$, ou que~$x\in \beth (k) -\aleph(k)$ ; on déduit alors de~\ref{pizerowbeextr} et~\ref{pizerowbeinf} que~$\got P(k)$ est unibranche en~$\xi$. 

\medskip
Supposons maintenant que~$x$ soit de la forme~$\zeta_{a,r}$ avec~$r>0$ ; on a alors~$\xi=\eta_{a,r}$. Distinguons maintenant deux cas : 

\medskip
$\bullet$ le cas où~$r\notin |k\ti|$ ; on a alors~$\bfer a r=\bouv a r$, et l'on déduit de~\ref{pizerowbezar} que~$\got P(k)$ est de valence 2 en~$\xi$ ; 

$\bullet$ le cas où~$r\in |k\ti|$ ; l'ensemble~$\Theta(a,r)$ (défini au~\ref{defthetaar}) est en bijection canonique avec~$\red k_r$, et 
donc en bijection (non canonique en général) 
avec~$\red k_1$. Ceci entraîne, en vertu de ~\ref{pizerowbezar}, que la valence de 
$\got P(k)$ en~$\xi$ est égale au cardinal de~$\PP^1(\red k_1)$. 

\deux{valdiscrwb} {\bf Un cas particulier important : celui d'un corps local.} Supposons que~$k$ est local, c'est-à-dire que~$|k\ti|$ est libre de rang~$1$, et que~$\red k$ est fini ; soit~$q$ le cardinal de ce dernier. 

\trois{valdiscrwb123} Comme~$|k\ti|$ est libre de rang~$1$, {\em  il n'existe pas de chaîne évanescente de boules fermées de~$k$} ; par conséquent,~$\got P(k)$ est la réunion de~$\PP^1(k)$ et de l'ensemble des points de la forme~$\eta_{a,r}$ avec~$a\in k$ et~$r>0$ (\ref{gottwwbdesc}). 

\trois{wbethmoinsp} Tout point de~$\PP^1(k)$ est un point unibranche de  ~$\got P(k)$ (~\ref{branchegottwb}). Comme par ailleurs~$\PP^1(k)$ est compact (puisque~$k$ est local), il s'ensuit que~$\got P(k)-\PP^1(k)$ est un sous-arbre ouvert de~$\got P(k)$. 

\medskip
Soit~$\xi\in \got P(k)-\PP^1(k)$ ; il s'écrit~$\eta_{a,r}$ pour un certain~$a\in k$ et un certain~$r>0$. D'après~\ref{branchegottwb}, la valence de~$(\got P(k),\xi)$ vaut~$2$ si~$r\notin |k\ti|$, et~$q+1$ sinon ; compte-tenu du fait que~$|k\ti|$ est libre de rang 1, ceci entraîne que~$\got P(k)-\PP^1(k)$ est un arbre localement fini, dont les sommets sont exactement les points de la forme~$\eta_{a,r}$ avec~$a\in k$ et~$r\in |k\ti|$, et sont tous de valence~$q+1$. 

\medskip
\trois{rembt} {\bf Remarque.} On déduit de sa construction que~$\got P(k)-\PP^1(k)$ s'identifie à {\em l'arbre de Bruhat-Tits} de~$\mathsf{GL}_2(k)$.

\subsection*{La courbe~$\pkk$ s'identifie à l'arbre~$\got P(\KK)$}

\deux{wbethgotpkk} Nous allons montrer que le compact~$\got P(\KK)$ de~$\pkk$ est égal à~$\pkk$ tout entier. Pour cela, donnons-nous un point~$x$ de~$\pkk$ ; 
nous allons vérifier qu'il appartient à l'image de~$\got t$. C'est évident si~$x=\infty$ ; on suppose
maintenant que~$x\in \Aff^{1,\rm an}_{\KK}$. Soit~$F$ l'ensemble des
majorants de~$x$ dans~$\Aff^{1,\rm an}_{\KK}$
qui sont de la forme~$\eta_{a,r}$, avec~$a\in \KK$ et~$r\geq 0$. 

\medskip
Soient~$y$ et~$z$ deux éléments de~$F$ ; écrivons~$y=\eta_{a,r}$ et~$z=\eta_{b,s}$. Par définition de~$F$, le point~$x$ appartient à l'intersection des boules~$\DD(a,r)$ et~$\DD(b,s)$, laquelle est en conséquence non vide ; il s'ensuit que ces deux boules sont comparables pour l'inclusion, ce qui implique que~$y$ et~$z$ sont comparables. L'ensemble~$F$ est donc une chaîne, qui contient par ailleurs~$\eta_r$ pour tout~$r>|T(x)|$. 

\medskip
Il s'ensuit que~$F=\got t(G)$ pour une certaine chaîne saturée~$G$ de~$\aleph(k)$. Soit~$x_0$ la borne inférieure de~$G$ ; nous allons montrer que~$x=\got t(x_0)$, ce qui permettra de conclure. Par définition de~$\got t$, on a~$\got t(x_0)=\inf \limits_{t\in G}\got t(t)=\inf\limits_{y\in F} y$ ; il s'agit donc de s'assurer que~$x=\inf\limits_{y\in F} y$.

\medskip
Soit~$a\in k$ ; posons~$r=|T(x)-a|$. On a~$x\in \DD(a,r)$ ; autrement dit,~$\eta_{a,r}\in F$. Comme~$|(T-a)(\eta_{a,r})|=r$, on a~$\inf\limits_{y\in F} |(T-a)(y)|\leq r=|(T-a)(x)|$. La définition de~$F$ assurant par ailleurs que~$|(T-a)(x)|\leq \inf\limits_{y\in F}|(T-a)(y)|$, il vient~$$|(T-a)(x)|=\inf_{y\in F} |(T-a)(y)|.$$ Ceci vaut pour tout~$a\in \KK$. Le corps~$\KK$ étant algébriquement clos, une semi-norme multiplicative sur~$\KK[T]$ est entièrement déterminée par ses valeurs sur les polynômes de la forme~$T-a$ avec~$a\in \KK$ ; il s'ensuit que~$x=\inf\limits_{y\in F} y$, ce qui achève la démonstration. 

\deux{pointspkkintro} Ainsi,~$\pkk$ est un arbre compact. Les faits suivants résultent de~\ref{propgottwb}. 

\trois{ptspkk} L'arbre~$\pkk$ est la réunion de~$\PP^1(\KK)$, de l'ensemble des points de la forme~$\eta_{a,r}$ avec~$a\in \KK$ et~$r>0$, et de~$\got t(\beth (\KK)-\aleph(\KK))$, lequel s'identifie à l'ensemble des classes d'équivalence de chaînes évanescentes de boules fermées de~$\KK$. Cet ensemble est vide si et seulement si~$\KK$ est sphériquement complet, ou encore si et seulement si il n'admet pas d'extension immédiate stricte. 

\trois{ptsextr} Si~$x\in \PP^1(\KK)$ ou si~$x$ appartient à~$\got t(\beth (\KK)-\aleph(\KK))$ alors~$x$ est un point unibranche de~$\pkk$. 

\trois{ptsetar} Si~$x$ est de la forme~$\eta_{a,r}$ avec~$r>0$ alors la valence de~$\pkk$ en~$x$ vaut~$2$ si~$r\notin |(\KK)\ti|$, et est sinon infinie et plus précisément égale au cardinal de~$\PP^1(\kk_1)$. 

\deux{typepointspkk} Soit~$x$ un point de~$\pkk$. 

\trois{pkkt1} Supposons que~$x\in \PP^1(\KK)$ ; le point~$x$ est alors de type 1. 

\trois{pkkt23} Supposons que~$x=\eta_{a,r}$ avec~$a\in \KK$ et~$r\in |(\KK)\ti|$. Le point~$x$ est alors de type 2 si~$r\in |\KK\ti|$, et de type~$3$ sinon 
(\ref{extype}
{\em et sq.}). 

\trois{pkkt4} Supposons que~$x\in \got t(\beth(\KK)-\aleph(\KK))$. Nous allons montrer qu'il est de type 4, en excluant les trois autres cas. Le point~$x$ n'est pas rigide puisqu'il n'appartient pas à~$\PP^1(\KK)$. 
Il ne peut pas non plus être de type 2 ou 3. 
Pour le voir, on raisonne par l'absurde en supposant donc
que c'est le cas. Le degré de transcendance
de~$\red {\hres(x)}$ sur~$({\kk})$ 
est alors égal à~$1$. Comme~$\KK$ est algébriquement clos,
$\KK(T)\ti$ est engendré par les monômes~$(T-a)$ où~$a$ parcourt~$\KK$, et il existe
dès lors un élément~$a$ dans~$\KK$ tel que~$\red{(T-a)(x)}$ soit transcendant sur~${\kk}$, ce qui signifie 
que~$x$ 
est égal à~$\eta_{a,|(T-a)(x)|}$
(\ref{transeteta}) et débouche ainsi sur une contradiction. 

\deux{basevoispkk} Nous allons maintenant donner une description des bases de voisinages des points de~$\pkk$, selon leur type ; elle se fonde sur les résultats établis aux~\ref{basevoisal} {\em et sq.} Soit donc~$x$ un point de~$\pkk$. 

\trois{voispkk14} Supposons que~$x$ est de type 1 ou 4 ; il possède alors une base de voisinages de la forme~$\DD\zero(a,r)$ avec~$a\in \KK$ et~$r>0$ (\ref{voisefeuille}). Notons une conséquence de ce fait : l'intersection~$\bigcap\limits_{(a,r), \;
x\in \DD\zero(a,r)} \DD\zero(a,r)=\{x\}$ (c'était {\em a priori} évident pour les points de type 1, mais pas pour ceux de type 4). 

\trois{voispkkt2} Supposons que~$x$ est de type 2, c'est-à-dire de la forme~$\eta_{a,r}$ avec~$a\in \KK$ et~$r\in |(\KK)\ti|$. Il possède alors (\ref{voiszetaar}) une base de voisinages de la forme~$\DD\zero(a,R)-\coprod\limits_{i\in I} \DD(a_i,r_i)$ où~$R>r$, où~$I$ est fini, où les~$r_i$ sont strictement inférieurs à~$r$, et où les~$a_i$ appartiennent à~$\bfer a r$ et sont tels que~$|a_i-a_j|=r$ dès que~$i\neq j$. 

\trois{voispkkt3} Supposons que~$x$ est de type 3, c'est-à-dire de la forme~$\eta_{a,r}$ avec~$a\in \KK$ et~$r\notin |(\KK)\ti|$. Il possède alors (\ref{voiszetaar}) une base de voisinages de la forme~$\DD\zero(a,R)-\DD(a,s)$ où~$R>r$ et où~$s<r$. 

\subsection*{L'arbre~$\pk$ et certains de ses sous-arbres finis} 

\deux{exretracan}
La courbe~$\PP^{1,\rm an}_{\KK}$ est un arbre compact d'après~\ref{wbethgotpkk} ; comme~$\pk$ est homéomorphe à~$\pkk/\mathsf G$,
c'est également un arbre compact, dont l'arbre~$\got P(k)$ est un sous-arbre compact.

\deux{chemina} Soient~$a$ et~$b$ deux éléments de~$k$, que l'on voit comme appartenant à~$\pk$ (ils s'identifient respectivement à~$\eta_{a,0}$ et~$\eta_{b,0}$). Il résulte 
de~\ref{defplongwb} {\em et sq.} 
que l'intervalle~$[a;b]$ est égal à~$$\{\eta_{a,r}\}_{0\leq r\leq |b-a||}\cup\{\eta_{b,r}\}_{0\leq r\leq |b-a|}, $$ et que l'intervalle~$[a;\infty[$ est égal à~$\{\eta_{a,r}\}_{0\leq r\leq +\infty}$, 
avec la convention~$\eta_{a,+\infty}=\infty$.  En particulier,~$[0;\infty]=\{\eta_r\}_{0\leq r\leq +\infty}$. 
 
\deux{retrainfini}
En tant que sous-arbre compact et non vide de~$\pk$, l'intervalle~$[0;\infty]$ en est un sous-arbre admissible.

\trois{trayonetaR} Comme~$|T(\eta_r)|=r$
pour tout~$r$, la restriction de~$|T|$ à~$[0;\infty]$ coïncide avec (la
restriction de) la fonction rayon~$\rho$ ; c'est en particulier une fonction 
strictement croissante. 

\trois{tconstanteetaR}
Nous allons démontrer
que~$|T|$ est localement constante en dehors de~$[0;\infty]$ et que
la rétraction canonique
$\pi: \pk\to [0;\infty]$ est égale à~$x\mapsto \eta_{|T(x)|}$. 

\medskip
Soit~$x\in \pk$. Si~$x=0$ ou~$x=\infty$ alors~$\rho(x)=x=\eta_{|T(x)|}$ ; supposons maintenant que~$x\notin\{0,\infty\}$. Soit~$L$ une
extension complète de~$k$
telle que~$x$ possède un antécédent~$L$-rationnel sur~$\PP^{1,\rm an}_L$, correspondant à un élément~$\lambda$ de~$L\ti$. 
Posons~$r=|T(x)|=|\lambda|$. L'intervalle~$I:=[\lambda; \eta_{\lambda, r,L}]$ est égal à~$\{\eta_{\lambda, t, L}\}_{0\leq t\leq r}$, et~$|T|$ est égale à~$r$ identiquement sur~$I$. 
L'image~$J$ de~$I$ sur~$\pk$ est une partie connexe contenant~$x$ et~$\eta_r$, et sur laquelle~$|T|=r$ identiquement. Il s'ensuit que~$J\cap[0;\infty]=\{\eta_r\}$, ce qui entraîne que~$\pi(x)=\eta_r$. 

Supposons que~$x$ n'appartienne pas à~$[0,\infty]$ et soit~$U$ sa composante connexe dans~$\pk\setminus [0;\infty]$. On a alors~$\pi(y)=\pi(x)$ pour tout~$y\in U$ ; compte-tenu de ce qui précède, 
il vient~$|T(y)|=|T(x)|$. Par conséquent, $|T|$ est localement constante en dehors de~$[0;\infty]$. 

\trois{tmoinsacommet} Soit~$a\in k$. Il résulte immédiatement
des définitions que l'automorphisme de translation par~$a$ de~$\pk$ 
induit un homéomorphisme de~$\got P(k)$ sur lui-même
qui préserve~$\rho$. On en déduit, au vu de ce qui précède,
les faits suivants : 

\medskip
$\bullet$ la restriction de~$|T-a|$ à~$[a;\infty]$ coïncide avec~$\rho$, et est
en particulier strictement croissante ; 

$\bullet$ la fonction~$|T-a|$ est localement constante en dehors de~$[a;\infty]$ et la rétraction
canonique de~$\pk$ sur~$[a;\infty]$ est égale à~$x\mapsto \eta_{a,|(T-a)(x)|}$.

\deux{intropolyvarp1} Soit~$f\in k(T)$ une fraction rationnelle non nulle et {\em scindée}, c'est-à-dire
quotient de deux polynômes scindés premiers entre eux. On note~$E$ l'ensemble des zéros et pôles de~$f$ ; 
pour tout~$a\in E$, on note~$m_a$ la multiplicité correspondante, qui appartient à~$\ZZ\setminus \{0\}$.

\trois{envconvinfini} Soit~$\Gamma$ la réunion des~$[a;\infty]$ pour~$a\in E$. Comme
la fonction~$f$ s'écrit~$\lambda \prod\limits _{a\in E}(T-a)^{m_a}$
pour un certain~$\lambda \in k^*$, on déduit du~\ref{tmoinsacommet}
ci-dessus que~$|f|$ est localement constante en dehors de~$\Gamma$. 

\trois{descfrho} Soit~$\Gamma'$ l'ensemble des points de valence 2 de l'arbre 
fini~$\Gamma$, et soit~$I$ 
un intervalle ouvert non vide de~$\Gamma'$. 
Par construction de~$\Gamma$, l'intervalle~$I$ 
est de la forme~$]\eta_{a,r} ; \eta_{a,R}[$ où~$a\in E$ et où~$0\leq r<R\leq +\infty$.  
De plus, si~$b\in E$ on est nécessairement dans l'un des deux cas suivants : 

\medskip
$\bullet$ $[b;\infty] \supset ]\eta_{a,r};\eta_{a,R}[$; c'est le cas si et seulement si~$b\leq \eta_{a,r}$, 
soit encore si et seulement si $|(T-a)(b)|\leq r$, c'est-à-dire
enfin si et seulement si~$b$
appartient
à la boule fermée~$B(a,r)$ du corps~$k$. 

$\bullet$ $[b;\infty] \cap  ]\eta_{a,r};\eta_{a,R}[=\varnothing$; c'est le cas si et seulement si~$b$
n'est majoré par aucun~$\eta_{a,s}$  pour~$s<R$, soit encore si et seulement si~$|(T-a)(b)|\geq R$, 
c'est-à-dire enfin 
si et seulement si~$b$ n'appartient pas à la boule ouverte~$B\zero(a,R)$ du corps~$k$. 

\medskip
Il s'ensuit, en vertu de~\ref{tmoinsacommet}, que la restriction de~$|f|$
à~$]\eta_{a,r}\;\eta_{a,R}[$ coïncide avec~$\rho^{\sum\limits_{b\in E\cap B(a,r)} m_b}$.

\trois{deflambdaa} Pour tout~$a\in E$, on note~$P_a$ l'ensemble des réels~$r\geq 0$ tels
que~$\sum\limits_{b\in E\cap B(a,r)}m_b=0$. On note~$E'$ l'ensemble des~$a\in E$ tels que~$P_a$ soit non vide, 
et~$E''$ l'ensemble~$E\setminus E'$. Si~$a\in E'$ on désigne par~$\lambda(a)$ la borne inférieure de~$P_a$ ; si~$a\in E''$
on pose~$\lambda(a)=+\infty$. 

\medskip
Soit~$a\in E'$, et soit~$r\in P_a$ ; notons que
comme~$m_a$ est non nulle, la boule~$B(a,r)$
contient au moins un élément de~$E$ distinct de~$a$. 
Si l'on pose~$r_0=\max\limits_{b\in B(a,r)\cap E, b \neq a} |a-b|$
alors~$\sum_{b\in B(a,r_0)\cap E}m_b=0$. Il s'ensuit 
que $\lambda(a)$ est de la forme~$|a-b|$ 
pour un certain~$b\in E\setminus \{a\}$ et qu'il appartient à~$P_a$.

\trois{predescgraphvar} Convenons de dire
qu'une boule fermée~$C$ de~$k$ est
{\em adaptée à~$f$}
si elle est de la forme~$B(a,\lambda(a))$ pour un certain~$a\in E'$. Si c'est le cas, 
on dira qu'un tel~$a$ est un {\em centre admissible}
de~$C$. Notons que comme~$\lambda(a)\in |k^*|$
il coïncide avec le diamètre de~$C$
(et ne dépend en particulier pas du centre admissible~$a$ considéré).

\trois{descgraphvar}
Soit~$C$ une boule fermé de~$k$ adaptée à~$f$ et soit~$r$ son diamètre. 
On note~$\Gamma_C$ la réunion des~$[a; \eta_{a,r}]$ où~$a$ parcourt l'ensemble
des centres admissibles de~$C$.
Par construction,~$\Gamma_C$ est un sous-arbre fermé et non vide
de~$\Gamma$, contenu dans~$\DD(a,r)$ pour tout centre admissible~$a$ de~$C$. 

\medskip
Soit~$D$ une boule fermée de~$k$ adaptée à~$f$ et distincte de~$C$. Nous allons
montrer que~$\Gamma_C\cap \Gamma_D =\varnothing$. Soit~$s$ le
diamètre
de~$D$, soit~$c$ un centre admissible de~$C$ et soit~$d$ un centre admissible de~$D$. 

\medskip
{\em Supposons que~$|c-d|>\max(r,s)$.}
On a alors~$\DD(c,r)\cap \DD(d,s)=\varnothing$, 
et {\em a fortiori}, $\Gamma_C\cap \Gamma_D=\varnothing$. 

\medskip
{\em Supposons que~$|c-d|\leq \max(r,s)$.}
On a alors une relation d'inclusion non triviale
entre~$C$ et~$D$, disons par exemple~$C\subset D$. Comme~$C$ et~$D$ sont
par hypothèse distinctes, 
l'inclusion est stricte, ce qui entraîne que~$r<s$. Supposons que~$\Gamma_C\cap \Gamma_D$
soit non vide. 
La réunion de ces deux sous-arbres
de~$\pk$ 
 est alors un arbre, et on a en particulier~$[c;\eta_{c,s}]\subset \Gamma_C\cup \Gamma_D$.
Comme~$\Gamma_C\subset \DD(c,r)$, on a~$]\eta_{c,r};\eta_{c,s}]\subset \Gamma_D$. 

Comme~$\Gamma_D$ est la réunion finie des~$[a;\eta_{a,s}]$ où~$a$
parcourt l'ensemble des centres admissibles de~$D$, il existe un tel~$a$ tel
que~$]\eta_{c,r};\eta_{c,s}]\subset [a;\eta_{a,s}]$ ; par compacité,
ce dernier intervalle contient aussi~$\eta_{c,r}$. Ceci entraîne que~$\eta_{c,r}=\eta_{a,r}$, ou encore que~$|c-a|\leq r$  ; puisque~$c$ 
est un centre admissible de~$C$, on a
$$0=\sum_{\alpha \in B(c,r)\cap E}m_\alpha=\sum_{\alpha \in B(a,r)\cap E}m_\alpha,$$ ce qui contredit le fait que~$a$ est un centre admissible
de la boule $D$, dont le diamètre~$s=\lambda(a)$ est~$>r$. 

\trois{compinfinigraphvar}
On note~$\Gamma_\infty$ la réunion des~$[a;\infty]=[a,\lambda(a)]$ pour~$a$
parcourant~$E''$ ; 
c'est un sous-arbre connexe de~$\pk$ (éventuellement vide). 

Soit~$C$ une 
boule fermée de~$k$ adaptée à~$f$ et soit~$r$ son
diamètre. 
L'intersection~$\Gamma_\infty\cap \Gamma_C$ est vide. En effet, supposons qu'elle ne le soit pas. Il existe
alors un
centre admissible~$c$
de~$C$ et un élément~$a$ de~$E''$ tel que~$[c;\eta_{c,r}]$ rencontre~$[a;\infty]$. Il existe donc~$s\leq r$ tel que~$\eta_{a,s}=\eta_{c,s}$ ; 
on a {\em a fortiori}
$\eta_{a,r}=\eta_{c,r}$ ; en conséquence,~$B(a,r)=B(c,r)$. Mais~$c$ est un centre admissible de~$C$, d'où
les égalités~$0=\sum_{\alpha \in D(c,r)\cap E}m_\alpha=\sum_{\alpha \in D(a,r)\cap E}m_\alpha,$ ce qui contredit
l'appartenance de~$a$ à~$E''$. 

\trois{gammafenfin} Notons~$\Gamma_f$ la réunion des~$[a; \eta_{a,\lambda(a)}]$
où~$a$ parcourt~$E$. C'est un sous-graphe compact de~$\Gamma$. 
Par ce qui précède, on peut l'écrire comme la réunion disjointe
de~$\Gamma_\infty$ et des~$\Gamma_C$ où~$C$ parcourt l'ensemble des boules fermées
de~$k$ adaptées à~$f$ ; ses composantes connexes
sont exactement les~$\Gamma_C$ ainsi
que~$\Gamma_\infty$ si celui-ci est non vide ; notons que par construction, aucune d'elles 
n'est réduite à un point.

\medskip
Soit~$I$ un intervalle ouvert non vide contenu
dans~$\Gamma'$. 
D'après~\ref{descfrho}, $I$ est de la forme~$]\eta_{a,r};\eta_{a,R}[$ avec~$a\in E$,
et la restriction de~$|f|$ à~$I$ est égale à~$\rho^{\sum\limits_{\alpha\in E\cap B(a,r)} m_\alpha}$. 
Nous allons montrer que les assertions suivantes
sont équivalentes : 

\medskip
1) la restriction de~$|f|$ à~$I$ est constante ; 

2) $\sum\limits_{\alpha \in E\cap B(a,r)} m_\alpha=0$ ; 

3) $I\cap \Gamma_f=\emptyset$ ; 

4) $I\not\subset \Gamma_f$.

\medskip
Il est clair que~1)$\iff$2). Supposons maintenant que~2) est vraie et soit~$b$
appartenant à~$E$. Supposons que~$[b;\infty]$ rencontre~$I$. Dans ce
cas, il le contient (car~$I\subset \Gamma'$) ; par compacité, il contient~$\eta_{a,r}$. Il vient
alors~$\eta_{a,r}=\eta_{b,r}$, puis~$B(a,r)=B(b,r)$.
En conséquence,
$\sum\limits_{\alpha \in E\cap B(b,r)}m_\alpha=0$, et donc~$\lambda(b)\leq r$.
Il s'ensuit que~$[b;\eta_{b,\lambda(b)}]$ 
ne rencontre pas~$I=]\eta_{b,r}; \eta_{b,R}[$.
Ceci valant pour tout~$b\in E$, on a bien~$I\cap \Gamma_f=\emptyset$, et~3) est vraie. 

\medskip
L'implication~3)$\Rightarrow$4) est évidente. Supposons enfin que~4) soit vraie et soit~$b$
appartenant à~$E\cap B(a,r)$. Comme~$I\not\subset \Gamma_f$, 
l'intervalle~$[b;\eta_{b,\lambda(b)}]$ ne contient pas~$I=]\eta_{b,r};\eta_{b,R}[$. Il s'ensuit trivialement que~$\lambda(b)<R$. Mais on a même~$\lambda(b)\leq r$.
En effet,
supposons~$\lambda(b)=s$ avec~$r<s<R$. Cela signifie que~$\sum\limits_{\alpha \in E\cap B(b,s)}m_\alpha=0$, 
et que~$\sum\limits_{\alpha \in E\cap B\zero(b,s)}m_\alpha\neq 0$. Il existe
donc un élément~$\alpha$
de~$E$ tel que~$|b-\alpha|=s$. L'arbre~$\Gamma$
contient~$]\eta_{b;r}\;\eta_{b,R}[$ et~$[\alpha ; \eta_{\alpha,s}]=[\alpha ; \eta_{b, s}]$, dont l'intersection est~$\{\eta_{b,s}\}$.
Le point~$\eta_{b,s}$ est 
donc de valence au moins 3 sur~$\Gamma$, contredisant le choix de~$I$. On a bien finalement~$\lambda (b)\leq r$. 

\medskip
Ceci entraîne que~$b\in E'$, et que~$B(b,\lambda(b))\subset B(a,r)$. Ainsi, $E\cap B(a,r)$ est
contenu dans une réunion 
(nécessairement disjointe) de boules adaptées à~$f$ et elles-mêmes
incluse dans~$B(a,r)$. Il s'ensuit que~$\sum\limits_{\alpha\in E\cap B(a,r)}m_\alpha=0$, 
ce qui montre~2) et achève la preuve de l'équivalence souhaitée. 

\trois{gammafgraphevar}
Soit~$U$ l'ouvert de~$\pk\setminus E$
formé des points au voisinage desquels~$|f|$ est constante. Nous allons montrer que~$U$
est le complémentaire de~$\Gamma_f$, par double inclusion. 

\medskip
{\em L'ouvert~$U$ ne rencontre pas~$\Gamma_f$.}
Soit~$x\in \Gamma_f$ et soit~$V$ un voisinage ouvert de~$x$. Comme~$x$
n'est pas un point isolé de~$\Gamma_f$, l'ouvert~$V$ contient un intervalle ouvert~$I$
aboutissant à~$x$. Quitte à restreindre~$I$, on peut supposer qu'il
est contenu dans~$\Gamma'$. 
Puisque~$I\subset \Gamma_f$, il résulte de~\ref{gammafenfin}
que la restriction de~$|f|$ à~$I$ est de la forme~$\rho^m$ 
pour un certain~$m$
{\em non nul.}
En conséquence,~$|f|_{|V}$ n'est pas constante, et~$x\notin U$. 

\medskip
{\em La fonction~$|f|$ est localement constante en dehors de~$\Gamma_f$.}
Soit~$x\in \pk\setminus \Gamma_f$. Si~$x\notin \Gamma$ alors~$|f|$ est
constante au voisinage de~$x$ (\ref{envconvinfini}). 

Supposons que~$x\in \Gamma$. Par compacité de~$\Gamma_f$, il existe un voisinage
ouvert~$\Delta$
de~$x$ dans~$\Gamma$ qui ne rencontre pas~$\Gamma_f$. Quitte à le restreindre, 
on peut supposer que~$\Delta\setminus \{x\}$ est une réunion finie d'intervalles ouverts contenus
dans~$\Gamma'$ et aboutissant à~$x$. Il résulte
de~\ref{gammafenfin}
que~$|f|$ est constante sur chacun de ces intervalles ouverts ; par continuité, elle est constante
de valeur~$|f(x)|$ sur~$\Delta$. 

Soit~$r$ la rétraction canonique de~$\pk$ sur~$\Gamma$. L'ouvert~$r^{-1}(\Delta)$ est un voisinage
de~$x$ dans~$\pk\setminus E$. La fonction~$|f|$ est
comme on l'a vu plus haut localement constante sur~$r^{-1}(\Delta)\setminus \Delta$. Par ailleurs, soit~$y\in
r^{-1}(\Delta)\setminus \Delta$ et soit~$V$ la composante connexe de~$y$ dans  $y\in
r^{-1}(\Delta)\setminus \Delta$. Par définition, $r(y)$ est le seul point de~$\partial V$. La valeur constante de~$|f|$ sur~$V$ est donc égale à~$|f(r(y))|$,
c'est-à-dire à~$|f(x)|$ puisque~$r(y)\in \Delta$. Ainsi, $|f|$ est constante de valeur~$|f(x)|$ sur~$r^{-1}(\Delta)$. 

\deux{gammafcasgen}
On désigne maintenant par~$f$ 
une fonction non nulle quelconque de~$k(T)$. 

\trois{gammafkbarre}
La fonction~$f$ est scindée dans~$\KK(T)$. D'après~\ref{intropolyvarp1}
{\em et sq.}, 
l'ouvert~$U$
de~$\pkk$ formé des points au voisinage desquels~$f$ est définie, inversible et de norme 
constante est le complémentaire d'un sous-arbre fini, non vide et non singleton~$\Gamma_{f, \KK}$
de~$\pkk$ ; si~$I$ désigne un intervalle ouvert
non vide de~$\Gamma_{f, \KK}$ constitué de points de valence~$2$
de ce dernier, la restriction de~$|f|$ à~$I$ est de la forme~$\rho^N$ pour un certain
entier relatif non nul~$N$ ; en particulier, cette restriction est strictement monotone. 
Il résulte de sa définition que~$U$ est stable sous~$\mathsf G$, et
il en va évidemment
de même de~$\Gamma_{f,\KK}$.

\trois{gammafkbarresurg}
La fonction~$f$ est définie et inversible sur~$U/\mathsf G$, et~$|f|$
y est localement constante car la
projection~$\pkk\to \pk\simeq \pkk/\mathsf G$ est ouverte. 
Par ailleurs, le quotient $\Gamma_{f,\KK}/\mathsf G$ est 
un arbre, et~$\Gamma_{f,\KK}\to \Gamma_{f,\KK}/\mathsf G$
est injective par morceaux 
(th.~\ref{theoquot}). 
On en déduit que~$\Gamma_{f,\KK}/\mathsf G$
n'est pas réduit à un singleton, et qu'il existe un sous-ensemble fini~$S$
de~$\Gamma_{f, \KK}/\mathsf G$ tel que pour tout intervalle ouvert non vide~$I$ contenu
dans~$(\Gamma_{f, \KK}/\mathsf G)\setminus S$, la restriction de~$|f|$ à~$I$ soit strictement monotone. 

\medskip
En conséquence, $U/\mathsf G$ est exactement l'ensemble des points en lesquels~$f$ est définie et 
au voisinage desquels~$|f|$ est constante. Si~$f$ est scindée, $\Gamma_{f,\KK}/\mathsf G$ coïncide
donc avec l'arbre~$\Gamma_f$ défini au~\ref{gammafenfin} ; on peut par conséquent
sans risque d'ambiguïté poser
en général~$\Gamma_f=\Gamma_{f,\KK}/\mathsf G$. Nous dirons que~$\Gamma_f$ est {\em l'arbre de variation}
de~$f$. 

\trois{exemvarypol} {\em Exemple : le cas d'un polynôme.}
Supposons que~$f\in k[T]$, et soit~$E$ l'ensemble de ses racines dans~$\KK$. Pour tout~$a\in E$ et tout~$r>0$, la boule
fermée de~$\KK$ de centre~$a$ et de rayon~$r$ ne contient aucun pôle de~$f$. Il s'ensuit que~$\lambda(a)=+\infty$,
que l'arbre~$\Gamma_{f,\KK}$ 
est la réunion des~$[a;\infty]$ pour~$a\in E$, et que~$|f|$ est strictement croissante sur~$[a;\infty]$ pour
tout~$a\in E$. 

\medskip
Soit~$a\in E$ et soit~$\alpha$ son image sur~$\pk$ ; c'est le point rigide défini par l'annulation
du polynôme minimal de~$a$ sur~$k$. Il résulte de~\ref{unboutquotinj} (appliqué en prenant pour~$X$
l'arbre à un bout~$\Aff^{1,\rm an}_{\KK}$, avec~$\omega=\infty$) que la flèche quotient
$\pkk\to \pk$ induit un homéomorphisme entre~$[a;\infty]$ et~$[\alpha; \infty]$. 

\medskip
Soit~$F$ l'ensemble des points rigides de~$\pk$ en lesquels~$f$ s'annule. Il résulte de ce qui précède que~$\Gamma_f$ est
égal à~$\bigcup_{\alpha \in F} [\alpha; \infty]$,  et que~$|f|$ est strictement croissante sur~$[\alpha;\infty]$ pour
tout~$\alpha\in E$.

\section{Toute courbe analytique est un graphe}

\subsection*{Courbes analytiques : allure globale}

\deux{courbegraphe} {\bf Théorème.} {\em Toute courbe analytique est un graphe.}

\medskip
{\em Démonstration.} Soit~$X$ une courbe~$k$-analytique. Comme~$X$ est topologiquement le quotient de~$X_{\KK}$ par l'action de Galois il suffit, pour montrer que~$X$ est un graphe, de s'assurer que~$X_{\KK}$ en est un ; on peut par conséquent supposer~$k$ algébriquement clos.  On sait que~$X$ est localement compact et localement connexe par arcs ; il suffit donc de vérifier {\em localement} que~$X$ est un graphe.

\medskip
Comme~$X$ est bon, on se ramène au cas où~$X$ est affinoïde. Il suffit, en vertu de la proposition~\ref{quotrelfin}, de démontrer que le normalisé de~$X$ est un graphe ; autrement dit,
on peut supposer que~$X$ est normal. 

\medskip
Comme~$k$ est algébriquement clos,~$X$ est quasi-lisse. Chacun de ses points possède donc un voisinage isomorphe à un domaine affinoïde de l'analytification d'une~$k$-courbe projective, irréductible et lisse ; on peut dès lors faire l'hypothèse que~$X$ est un domaine affinoïde de~${\sch X}\an$, pour une certaine ~$k$-courbe algébrique, projective, irréductible et lisse~$\sch X$. Le complémentaire de~$X$ dans~${\sch X}\an$ est une réunion de composantes connexes de~${\sch X}\an\setminus \partial  X$ ; comme~$\partial X$ est fini (c'est le bord de Shilov de~$X$), il suffit de démontrer que~${\sch X}\an$ est un graphe. 

\medskip
Il existe un morphisme fini et plat de~${\sch X}$ vers~$\PP^1_k$ ; soit~$\sch Y$ le normalisé de~${\sch X}$ dans une clôture normale de l'extension~$k(T)\subset \kappa({\sch X})$ induite par ce morphisme. Si~$\mathsf G$ désigne le groupe de Galois de cette clôture normale, et si~$\mathsf H$ est le sous-groupe de~$\mathsf G$ correspondant à~$\kappa({\sch X})$, on dispose d'homéomorphismes~$$\pk\simeq {\sch Y}\an/\mathsf G\;{\rm et}\;{\sch X}\an\simeq {\sch Y}\an/\mathsf H.$$

\medskip
Nous avons vu au~\ref{wbethgotpkk} ci-dessus que l'espace topologique~$\pk$ est un arbre. Soit~$\ell$ un entier strictement supérieur à~$1$ et inversible dans~$k$. Tout revêtement topologique de~${\sch Y}\an$ pouvant être vu comme un revêtement étale de cet espace, il résulte de GAGA et de la finitude du groupe~$\H^1({\sch Y}_{\rm \acute{e}t},\ZZ/\ell\ZZ)$ que~$\H^1({\sch Y}\an_{\rm top},\ZZ/\ell\ZZ)$ est fini. Le théorème~\ref{theofondquot} assure alors que~${\sch Y}\an$ est un graphe ; il s'ensuit, en vertu du théorème~\ref{theoquot}, que~${\sch X}\an$ est un graphe.~$\Box$

\deux{pointrigfinbr} { \bf Proposition.} {\em Soit~$X$ une courbe~$k$-analytique et soit~$x\in X\typ 0$. L'ensemble~$\br X x$ est fini, et son cardinal est égal au nombre d'antécédents de~$x$ sur le normalisé de~$X$.}

\medskip
{\em Démonstration.} Soit~$X'$ la normalisée de~$X$, soit~$\phi : X'\to X$ la flèche canonique et soient~$x_1,\ldots,x_r$ les antécédents de~$x$ sur~$X'$. Il existe un voisinage affinoïde~$V$ de~$x$ dans~$X$ possédant les propriétés suivantes : 

\medskip

$\bullet$~$V$ est un arbre et l'espace~$V_{\rm red}\setminus\{x\}$ est normal ; 

$\bullet$~$\phi\inv(V)$ est de la forme~$\coprod \limits_{1\leq i\leq r} V_i$ où~$V_i$ est pour tout~$i$ un voisinage affinoïde de~$x_i$ qui est un arbre.

\medskip
Il résulte des hypothèses faites sur~$V$ que~$\coprod \left(V_i\setminus\{x_i\}\right)\to V\setminus\{x\}$ est bijective ; étant par ailleurs finie et de ce fait compacte, cette flèche est un homéomorphisme. Il s'ensuit que~$\pi_0(V\setminus\{x\})$ est en bijection naturelle avec~$\coprod_i \pi_0(V_i\setminus\{x_i\})$ ; par conséquent, l'application canonique ~$$ \coprod_i \br {X'}{x_i}\to \br X x$$ est bijective.

\medskip
Cela permet de se ramener au cas où la courbe~$X$ est normale ; il s'agit alors de prouver que~$\br X x$ est un singleton. Soit~$F$ le complété de la clôture parfaite de~$k$ et soit~$Z$ la normalisée de~$X_F$. La flèche~$X_F\to X$ est un homéomorphisme. Soit~$V$ un ouvert connexe non vide de~$X_F$ ; il est de la forme~$U_F$, où~$U$ est un ouvert connexe et non vide de~$X$. Comme~$X$ est normale,~$U$ est irréductible ; il s'ensuit que~$V=U_F$ est irréductible ; par conséquent, l'image réciproque~$W$ de~$V$ sur~$Z$, qui s'identifie à la normalisée de~$V$, est connexe. Ceci valant pour tout ouvert connexe et non vide~$V$ de~$X_F$, la flèche~$Z\to X$ est un homéomorphisme ; on peut donc, pour établir notre assertion, remplacer~$X$ par~$Z$, c'est-à-dire finalement supposer~$X$ normale et~$k$ parfait ; la courbe~$X$ est alors quasi-lisse. Le corps~$\hres(x)=\kappa(x)$ est une extension finie séparable de~$k$, qui se plonge dans l'anneau local hensélien~$\sch O_{X,x}$  ; on peut donc, quitte à restreindre~$X$ et à remplacer~$k$ par~$\hres(x)$, supposer que~$x$ est un~$k$-point ; étant lisse, il possède un voisinage qui est un disque, et est de ce fait unibranche.~$\Box$

\subsection*{Morphismes finis et plats entre courbes analytiques : le degré d'une branche sur son image}

\deux{bracourbe} Soit~$X$ une courbe~$k$-analytique ; c'est un graphe (th.~\ref{courbegraphe}) ce qui autorise à parler de ses branches (\ref{branchegr}). Si~$Y\to X$ est un morphisme compact entre courbes~$k$-analytiques ({\em e.g.}~$\phi$ est fini, ou~$\phi$ identifie~$Y$ à un domaine analytique fermé de~$X$), si~$y\in Y$, si~$x$ désigne son image sur~$X$ et si~$\phi\inv(x)$ est fini, on dispose d'une application naturelle~$\br Y y \to \br X x$ (\ref{imdir}). On vérifie immédiatement que la flèche~$\br Y y \to \br X x$ reste bien définie lorsque~$(Y,y)\to (X,x)$ est un morphisme de {\em germes} de courbes~$k$-analytiques à fibre finie. 

\deux{degbr} Soit~$\phi : Y\to X$ un morphisme fini et plat entre {\em courbes}~$k$-analytiques ;  soit~$x\in X$ ; soit~$a$ une branche de~$X$ issue de~$x$. 

\trois{phimoinsfin} Si~$d$ désigne le degré de~$\phi$ au-dessus de~$x$ et si~$U$ est un voisinage ouvert de~$x$ dans~$X$ qui est un arbre, le cardinal de toute fibre de~$\phi$ en un point de~$U$ est majoré par~$d$ ; il s'ensuit que pour tout ouvert connexe~$V$ de~$U$ le cardinal de~$\pi_0(\phi\inv(V))$ est majoré par~$d$ ; par conséquent,~$\phi\inv(a)$ est fini de cardinal au plus égal à~$d$.

\trois{defdegbr} Soit~$b\in \phi\inv(a)$ ; soit~$U$ un voisinage ouvert de~$x$ dans~$X$ qui est un arbre et qui est tel que~$\phi\inv(U)$ soit une réunion disjointe d'arbres et sépare les antécédents de~$x$. Soit~$W$ la composante connexe de~$\phi\inv(a(U))$ qui correspond à~$b$. Le morphisme~$V\to a(U)$ est fini et plat, {\em et son degré ne dépend pas du choix de~$U$} : on le déduit du fait que si~$U'$ est un voisinage ouvert connexe de~$x$ dans~$U$, alors la composante connexe de~$\phi\inv(a(U'))$ qui correspond à~$b$ n'est autre que~$W\times_{a(U)}a(U')$. Ce degré est appelé {\em degré de~$b$ au-dessus de~$a$} et est noté~$\deg \;(b\to a)$.

\trois{egaldegbr} Soit~$y$ un antécédent de~$x$. On a~$$\deg^y\phi=\deg \;(\phi\inv(U)_y\to U)=\deg\;(\phi\inv(U)_y\times_Ua(U)\to a(U))$$~$$=\sum_{W\in \pi_0(f\inv(a(U))), W\subset \phi\inv(U)_y}\deg \;(W\to a(U))$$~$$=\sum_{b\in \br Y y \cap \phi\inv(a)}\deg \;(b\to a).$$

\subsection*{L'exemple des fonctions non constantes} 

 \deux{fincourbes} Soit~$\phi : Y\to X$ un morphisme entre espaces~$k$-analytiques ; supposons que~$Y$ est de dimension finie~$d$. Pour tout~$x\in X$, la dimension  de~$\phi\inv(x)$ est au plus~$d$ (on peut le vérifier après extension des scalaires, et donc supposer~$x$ rigide, auquel cas~$\phi\inv(x)$ est un fermé de Zariski de~$Y$, d'où l'assertion) ; la dimension de~$\phi$ en tout point de~$Y$ est donc majorée par~$d$, et le lieu~$Z$ des points où elle vaut~$d$ est par conséquent un fermé de Zariski de~$Y$. {\em Il est purement de dimension~$d$.} En effet, il s'agit de vérifier que si~$z\in Z$ alors~$\dim z Z=d$ ; on peut s'en assurer après extension des scalaires, et donc supposer~$z$ rigide ; dans ce cas, la fibre~$\phi\inv(\phi(z))$ est un fermé de Zariski de~$Y$, qui possède par définition de~$Z$  une composante irréductible~$T$ de dimension~$d$ passant par~$z$ ; mais  l'on a alors clairement~$T\subset Z$ et donc~$ \dim z Z=d$, d'où notre assertion. 

\deux{lemmefincourbes} {\bf Lemme.} {\em Soit~$d$ un entier et soit~$\phi : Y\to X$ un morphisme entre espaces~$k$-analytiques ; supposons~$Y$ purement de dimension~$d$. Soit~$\{Y_i\}_{i\in I}$ l'ensemble des composantes irréductibles de~$Y$, et soit~$J$ le sous-ensemble de~$I$ formé des indices~$i$ tels que~$\phi_{|Y_i}$ soit constante. Le lieu~$Z$ des points de~$Y$ en lesquels~$\phi$ est de dimension~$d$ est alors égal à~$\bigcup_{i\in J}Y_i$.}

\medskip
{\em Démonstration.} Soit~$i\in J$. Comme~$\phi_{|Y_i}$ est constante,~$\phi(Y_i)$ est égal à~$\{x\}$ pour un certain point rigide~$x$ de~$X$. L'inclusion~$Y_i\subset \phi\inv(x)$ implique que la dimension de~$\phi$ en tout point de~$Y_i$ est minorée par~$d$, et donc égale à~$d$ puisque~$\dim {}Y=d$ ; par conséquent,~$\bigcup_{i\in J}Y_i\subset Z$. 

\medskip
Montrons maintenant l'inclusion réciproque. Comme~$Z$ est purement de dimension~$d$ (\ref{fincourbes}), il est réunion de composantes irréductibles de~$Y$ ; il suffit dès lors de s'assurer que si~$i\in I$ est tel que~$Y_i\subset Z$ alors~$i\in J$.

\medskip
Soit~$i$ tel que~$Y_i\subset Z$. Choisissons un polyrayon~$k$-libre~$\bf r$ tel que~$Y_{i, \bf r}-\coprod\limits_{j\neq i} Y_{j,\bf r}$ ait un point~$k_{\bf r}$-rigide~$y$ ; l'image de~$y$ sur~$X_{\bf r}$ est un point~$k_{\bf r}$-rigide~$x$ ; comme~$\phi_{\bf r}: Y_{\bf r}\to X_{\bf r}$ est de dimension~$d$ en tout point de~$Y_{i, \bf r}$, la fibre~$\phi_{\bf r}\inv(x)$ est un fermé de Zariski de~$Y_{\bf r}$ qui contient~$y$ et est de dimension~$d$ en celui-ci ; puisque~$Y_{i,\bf r}$ est la seule composante irréductible de~$Y_{\bf r}$ qui contienne~$y$, on a~$Y_{i,\bf r}\subset \phi_{\bf r}\inv(x)$ ; autrement dit, ~$\phi_{\bf r}(Y_{i,\bf r})=\{x\}$.

Si~$t$ désigne l'image de~$x$ sur~$X$ alors~$\phi(Y_i)=\{t\}$ ; nous allons montrer que~$t$ est rigide, ce qui prouvera que~$\phi_{|Y_i}$ est constante et achèvera la démonstration. Si~$z$ est un antécédent de~$t$ sur~$X_{\bf r}$, la fibre de~$Y_{i,\bf r}$ en~$z$ est non vide ; par conséquent,~$z=x$. Il s'ensuit que~$\sch M(\hres(t)\hotimes_k k_{\bf r})$ est réduit à un point ; par conséquent,~$\bf r$ est~$\hres(t)$-libre et le corps~$\hres(x)$ s'identifie alors à~$\hres(t)_{\bf r}$. Par hypothèse,~$\hres(t)_{\bf r}=\hres(x)$ est une~$k_{\bf r}$-algèbre de Banach finie ; il en découle que~$\hres(t)$ est une~$k$-algèbre de Banach finie et le point~$t$ est bien rigide.~$\Box$

\deux{corollimpasrig} {\bf Corollaire.}
{\em Soit~$\phi : Y\to X$ un morphisme entre espaces~$k$-analytiques ; supposons~$Y$ purement de dimension~$d$. Si~$x$ est un point non rigide de~$X$, la dimension de~$\phi\inv(x)$ est strictement inférieure à~$d$.}

\medskip
{\em Démonstration.} Soit~$y$ un point de~$Y$ en lequel la dimension de~$\phi$ est~$d$. En vertu du lemme~\ref{lemmefincourbes} ci-dessus,~$y$ est situé sur une composante irréductible~$Y'$ de~$Y$ telle que~$\phi_{|Y'}$ soit constante ; par conséquent,~$\phi(y)$ est rigide, et diffère donc de~$x$.~$\Box$ 

\deux{consaltfinc} Si~$\phi : Y\to X$ est un morphisme entre espaces~$ k$-analytiques et si~$Y$ est une courbe, le lemme~\ref{lemmefincourbes} ci-dessus implique que les assertions suivantes sont équivalentes : 

\medskip
a)~$\phi$ est de dimension nulle en tout point de~$Y$ ;

b) la restriction de~$\phi$ à chacune des composantes irréductibles de~$Y$ est non constante. 

\medskip
Remarquons que si~$\phi$ satisfait ces deux propositions alors~$\phi$ est fini en tout point de~$Y$ en lequel il est intérieur, et en particulier en tout point de l'intérieur analytique de~$Y$. 

\deux{exfoncnoncons} Soit~$Y$ une courbe~$k$-affinoïde et soit~$\sch A$ l'algèbre associée ; nous allons montrer qu'il existe une fonction analytique~$f: Y\to \Aff^{1,{\rm an}}_k$ qui satisfait les conditions équivalentes a) et b) du~\ref{consaltfinc} ci-dessus. 

Il s'agit d'exhiber une fonction~$f\in \sch A$ dont la restriction à chacune des composantes irréductibles de~$Y$ soit non constante. En raisonnant composante connexe par composante connexe, on se ramène au cas où~$Y$ est connexe et non vide, et l'on distingue deux cas. 

\medskip
 {\em Le cas où~$Y$ est irréductible.} Comme~$\dim k Y=1$, il existe une fonction~$f\in \sch A$ qui n'est pas entière sur~$k$, et répond au problème posé. 
 
 {\em Le cas où l'ensemble~$\{Y_i\}$ des composantes irréductibles de~$Y$ comprend au moins deux éléments.}  On choisit alors pour tout~$i$ une fonction~$f_i$ dans~$\sch A$ dont la restriction à~$Y_i$ est génériquement inversible, et dont le lieu des zéros contient~$\bigcup\limits_{j\neq i}Y_j$ ; on pose~$f=\sum f_i$ ; la définition de~$f$ assure que pour tout indice~$i$ la restriction de~$f$ à~$Y_i$ est génériquement inversible, et qu'elle s'annule au moins en un point (ceci parce que~$Y$ est connexe et possède au moins deux composantes irréductibles, ce qui garantit que~$Y_i$ rencontre~$Y_j$ pour au moins un indice~$j\neq i$) ; par conséquent,~$f$ répond là encore au problème posé.

\deux{nonconstfiniplat} Soit~$\phi : Y\to X$ un morphisme entre courbes~$k$-analytiques qui satisfait les conditions équivalentes a) et b) du~\ref{consaltfinc} et soit~$y$ un point de ~$Y$ en lequel~$\phi$ est intérieur. Le morphisme de germes~$(Y,y)\to (X,x)$ est alors fini. Si de plus~$X$ est normale en~$x$ et~$Y$ est réduite en~$y$, alors~$(Y,y)\to (X,x)$ est {\em plat}. Pour le voir on peut, quitte à restreindre~$X$ et~$Y$ (opération qui préserve visiblement la propriété a) ), supposer que~$X$ et~$Y$ sont toutes deux affinoïdes, que~$Y$ est réduite, que~$X$ est normale et irréductible et que~$Y\to X$ est fini. Si~$Z$ est une composante irréductible de~$Y$ son image sur~$X$ est un fermé de Zariski, nécessairement égal à~$X$ pour des raisons de dimension. Soit~$\sch A$ (resp.~$\sch B$) l'algèbre des fonctions analytiques sur~$X$ (resp.~$Y$). 

\medskip
L'anneau~$\sch A$ est intègre, noethérien, normal et de dimension de Krull au  plus égale à~$1$.  Le schéma~$\sch B$ est réduit, et chaque composante irréductible de~$\spec \sch B$ se surjecte sur~$\spec \sch A$. Il s'ensuit que si~$f$ est un élément non nul de~$\sch A$ alors~$f$ n'est pas un diviseur de zéro dans~$\sch B$ ; par conséquent,~$\sch B$ est~$\sch A$-plat, ce qui prouve notre assertion. 

\deux{nonconstimnonrig} Soit~$\phi : Y\to X$ un morphisme entre courbes~$k$-analytiques, soit~$y\in Y$ et soit~$x$ son image sur~$X$ ; on suppose que~$x$ n'est pas rigide. La dimension de~$\phi$ en~$y$ est alors nulle d'après le corollaire~\ref{corollimpasrig}. On suppose de plus que~$\phi$ est intérieur en~$y$, et que~$X$ est réduite en~$x$. 

\trois{platcarcorps} Le morphisme fini de germes~$(Y,y)\to (X,x)$ est alors plat : en effet, comme~$X$ est réduite et~$x$ non rigide, l'anneau local~$\sch O_{X,x}$ est un corps, et~$\sch O_{Y,y}$ est dès lors un~$\sch O_{X,x}$-module plat. 

\trois{bondeviss} Supposons de plus que~$Y$ est réduite en~$y$ ; comme~$y$ n'est pas rigide, l'anneau local~$\sch O_{Y,y}$ est un corps ; par conséquent,~$Y$ est normale en~$y$, et~$X$ est normale en~$x$ pour le même motif. Il existe donc un voisinage affinoïde connexe et normal~$W$ de~$y$ dans~$Y$, et un voisinage affinoïde connexe et normal~$V$ de~$x$ dans~$X$, tels que~$\phi$ induise un morphisme fini et plat~$W\to V$ ; soient~$\sch A$ et~$\sch B$ les algèbres de fonctions respectives de~$V$ et~$W$ ; elles sont intègres,~$\spec \sch B\to \spec \sch A$ est fini et dominant, et~$y$ et~$x$ sont situés au-dessus des points génériques respectifs de~$\spec \sch B$ et~$\spec \sch A$. 

\medskip
Le dévissage de l'extension~$\kappa(\spec \sch A)\hookrightarrow \kappa(\spec \sch B)$ en ses parties séparables et radicielles induit un dévissage analogue de~$\spec \sch B\to \spec \sch A$ au-dessus d'un ouvert non vide~$\sch U$ de~$\spec \sch A$ ; soit~$V'$ un voisinage affinoïde connexe de ~$x$ dans~$\sch U\an$ et soit~$W'$ la composante connexe de~$y$ dans~$W\times_VV'$. Par construction,~$W'\to V'$ admet une factorisation~$W'\to Z\to V'$ où~$W'\to Z$ est fini, plat et radiciel, où~$Z$ est connexe et où~$Z\to V'$ est fini étale. 

\deux{pasplatmaisouv} Si~$\phi: Y \to X$ est un morphisme fini entre courbes~$k$-analytiques et si~$X$ est normale, alors~$\phi$ est ouvert. En effet, comme la question est purement topologique, on peut remplacer~$Y$ par~$Y_{\rm red}$ et donc supposer que~$Y$ est réduite ; mais en vertu du~\ref{nonconstfiniplat} ci-dessus,~$Y\to X$ est alors plat, et partant ouvert.

\section{Disques et couronnes déployés et virtuels} 

\subsection*{Disques et couronnes plongés}

\deux{defdi} Si~$I$ est un intervalle de~$\RR_+$ rencontrant~$\RR\ti_+$, nous noterons~$\DD_I$ le domaine analytique de~$\Aff^{1,{\rm an}}_k$ défini par la condition~$|T|\in I$. 
Il résulte de~\ref{retrainfini} 
que l'application~$r\mapsto \eta_r$ induit un homéomorphisme entre~$I$ et un intervalle admissible~$\eta_I$ de~$\DD_I$, et 
que la rétraction canonique de~$\DD_I$ vers~$\eta_I$ est donnée par la formule~$x\mapsto \eta_{|T(x)|}$ ; en particulier, $\DD_I$ est connexe et non vide. 

Si~$I$ est ouvert dans~$\RR$ alors~$\eta_I=\mathsf S(\DD_I)$, et le domaine~$\DD_I$ est affinoïde si et seulement si~$I$ est compact. Comme~$\Aff^{1,{\rm an}}_k$ est lisse,~$\DD_I$ est quasi-lisse. Si~$r\in \RR\ti_+$ nous écrirons~$\DD_r$ au lieu de~$\DD_{\{r\}}$. Si~$L$ est une extension complète de~$k$, nous utiliserons les notations~$\DD_{I,L}$,~$\DD_{r,L}$ et~$\eta_{I,L}$ dans un sens évident. Pour toute extension complète~$L$ de~$k$, l'espace~$L$-analytique~$\DD_{I,L}$ est connexe ; par conséquent,~$\DD_I$ est géométriquement connexe.

\deux{etainvgal} Soit~$I$ un intervalle de~$\RR_+$ rencontrant~$\RR\ti_+$ et soit~$L$ une extension complète de~$k$. Il résulte des définitions que la flèche~$\DD_{I,L}\to \DD_I$ induit un homéomorphisme~$\eta_{I,L}\simeq \eta_I$, et que~$\DD_{I,L}\to \DD_I$ commute aux rétractions canoniques de~$\DD_{I,L}$ sur~$\eta_{I,L}$ et de~$\DD_I$ sur~$\eta_I$. 

\medskip
Par ailleurs, le point~$\eta_{r,\KK}$ de~$\Aff^{1,{\rm an}}_{\KK}$ est par sa définition même invariant sous l'action de Galois ; il s'ensuit que si~$L$ est une extension presque algébrique de~$k$, l'image réciproque de~$\eta_I$ sur~$\DD_{I,L}$ est exactement~$\eta_{I,L}$. 

\deux{annfoncdi} Si~$I$ est un intervalle de~$\RR_+$ rencontrant~$\RR\ti_+$, l'anneau~$k_I$ des fonctions analytiques sur~$\DD_I$ s'identifie à la~$k$-algèbre des séries~$\sum\limits _{i\in \ZZ}a_iT^i~$ telles que~$|a_i|r^i$ tende vers zéro lorsque~$|i|$ tend vers l'infini pour tout~$r\in I\setminus\{0\}$ et telles que~$a_i=0$ pour tout~$i<0$ si~$0\in I$. 

\deux{borsshiv} Soit~$r\in \RR\ti_+$ et soit~$s>r$. Nous allons décrire les bords de Shilov des espaces~$k$-affnoïdes~$k_{[0;r]}$,~$k_r$ et~$k_{[r;s]}$ ainsi que les annéloïdes résiduels de leurs algèbres de fonctions.

\trois{borsshivdisc} Soit~$f=\sum a_iT^i$ une fonction analytique sur~$\DD_{[0;r]}$ ; pour tout~$x$ appartenant $\DD_{[0;r]}$, on a~$$|f(x)|\leq \max |a_i|\cdot r^i=|f(\eta_r)|\;;$$ 
par conséquent, le bord de Shilov de~$\DD_{[0;r]}$ est le singleton~$\{\eta_r\}$. 

\medskip
La (semi-) norme spectrale de
l'algèbre affinoïde~$k_{[0;r]}$
est donc égale à~$\sum a_i T^i\mapsto \max |a_i|r^i$.
En vertu
de~\ref{algindepgaussuite} et
par densité
de~$k[T]$ dans~$k_{[0;r]}$,
il s'ensuit que~$\tau \mapsto \red T$ induit 
un isomorphisme~$\red k[r\inv \tau]\simeq \red{k_{[0;r]}}$. 

\trois{borsshivcercl} Soit~$f=\sum a_iT^i$ une fonction analytique sur~$\DD_r$ ; pour tout~$x$ appartenant $\DD_r$, on a~$$|f(x)|\leq \max |a_i|\cdot r^i=|f(\eta_r)|\;;$$ 
par conséquent, le bord de Shilov de~$\DD_r$ est le singleton~$\{\eta_r\}$. 

\medskip
La (semi-) norme spectrale de~l'algèbre
affinoïde~$k_r$
est donc égale à~$\sum a_i T^i\mapsto \max |a_i|r^i$.
En vertu
de~\ref{exprecouronne} et
par densité
de~$k[T,T\inv ]$ dans~$k_r$, il s'ensuit
que~$\tau \mapsto \red T$ induit 
un isomorphisme~$\red k[r\inv \tau , r\tau \inv]\simeq \red{k_r}$.

\trois{borsshivcour} Soit~$f=\sum a_iT^i$ une fonction analytique sur~$\DD_{[r;s]}$ ; pour tout~$x$ appartenant $\DD_{[r;s]}$, on a~$$|f(x)|\leq \max \{|a_i|\cdot r^i,  |a_i|\cdot s^i\}_i=\max (|f(\eta_r)|, |f(\eta_s)|)\;;$$ 
par conséquent, le bord de Shilov de~$\DD_{[r;s]}$ est contenu dans~$\{\eta_r,\eta_s\}$. Comme~$|T(\eta_r)|<|T(\eta_s)|$ et~$|T\inv(\eta_r)|>|T\inv(\eta_s)|$, ce bord de Shilov
est en fait {\em exactement égal}
à la paire~$\{\eta_r,\eta_s\}$.

\medskip
La (semi-) norme spectrale de
l'algèbre affinoïde~$k_{[0;r]}$
est donc égale à
$\sum a_i T^i\mapsto \max \{|a_i|r^i,|a_i|s^i\}_i$.
En vertu
de~\ref{exprecouronne} et
et
par densité
de~$k[T,T\inv ]$ dans~$k_{[r;s]}$, il s'ensuit
que~$\tau\mapsto \red T, \sigma \mapsto \red{T\inv}$ induit 
un isomorphisme~$\red k[r\inv \tau , s\sigma \inv]/\sigma \tau\simeq \red{k_{[r;s]}}$.

\deux{foncdr} Soit~$r\in \RR\ti_+$ et soit~$f$ une fonction analytique non nulle sur~$\DD_r$ ; écrivons~$f=\sum a_iT^i$, et notons~$\sch E$ l'ensemble fini et non vide des indices~$j\in \ZZ$ tels que~$|a_j|r^j=\max\limits_i|a_i|r^i$. 

\trois{jsinglefinv} Supposons que l'ensemble~$\sch E$ soit un singleton~$\{j\}$ ; on a alors~$a_j\neq 0$ et~$f=a_jT^j (1+u)$, où~$u$ appartient à l'anneau ~$\sch O_{\DD_r}(\DD_r)\zeroo$ formé des fonctions analytiques sur~$\DD_r$ dont la norme est strictement majorée par~$1$ en tout point de~$\DD_r$ ; par conséquent,~$f$ est inversible. 

\trois{finvsingle} Réciproquement, supposons que~$f$ est inversible, et montrons que~$\sch E$ est un singleton. Le bord de Shilov de~$\DD_r$ étant un singleton, la valeur absolue de la fonction inversible~$f$ est constante sur~$\DD_r$. Il s'ensuit que~$\red f$ ne s'annule pas sur~$\red{\DD_r}$ ; par conséquent,~$\red f$ est un élément inversible de la~$\red k$-algèbre
~$\red {k_r}$. Or on a vu au~\ref{borsshivcour} 
que~$\red T$ est transcendant sur~$\red k$ et que~$\red {k_r}=\red k [\red T, \red T\inv]$. On en déduit
que~$\red f$ est de
la forme~$\red a \red T^j$ pour un certain~$a\in k\ti$ et un certain~$j\in \ZZ$ ; mais la définition de la flèche~$g\mapsto \red g$ assure alors que~$\sch E=\{j\}$. 

\trois{factf} En général, le lieu des zéros de~$f$ sur~$\DD_r$ est un ensemble fini de points rigides. On peut donc écrire~$f=gh$, où~$g$ est une fonction inversible et où~$h$ est un polynôme unitaire en~$T$ dont toutes les racines dans~$\KK$ sont de valeur absolue égale à~$r$.

\medskip
{\em Si~$d$ est le degré de~$h$ alors~$d=\max \sch E-\min \sch E$}. Pour le voir, on peut étendre les scalaires et donc supposer~$h$ de la forme~$\prod (T-\alpha_\ell)$ où les~$\alpha_\ell$ sont des éléments de~$k$ de valeur absolue égale à~$r$. Écrivons~$g=\sum b_iT^i$ ; comme~$g$ est inversible, il résulte de~\ref{finvsingle} qu'il existe un unique indice~$j$ tel que~$|b_j|r^j=\max |b_i|r^i$ ; un calcul immédiat montre alors que~$\sch E$ contient~$j$ et~$j+d$, et est contenu dans~$\{j,j+1,\ldots,j+d\}$, d'où notre assertion. Remarquons que celle-ci entraîne que~$k_r/f$ est une~$k$-algèbre finie de dimension~$\max \sch E-\min \sch E$. 

\deux {invgen} Soit~$I$ un intervalle de~$\RR_+$ rencontrant~$\RR\ti_+$, et soit~$f=\sum a_iT^i$ un élément de~$k_I$. 

\trois{casipaszero} Si~$0\notin I$, il résulte de~\ref{jsinglefinv} et de~\ref{finvsingle} que~$f$ est inversible si et seulement si pour tout~$r\in I$ il existe un unique entier~$j(r)$ tel que~$$|a_{j(r)}|r^{j(r)}=\max |a_i|r^i.$$ Si c'est le cas, la fonction~$r\mapsto j(r)$ est localement constante, et partant constante sur~$I$. Il s'ensuit que~$f$ est inversible si et seulement si il existe un entier~$j$ (nécessairement unique) tel que~$|a_j|r^j>\max |a_i|r^i$ pour tout~$i\neq j$ et {\em pour tout~$r\in I$.}

\trois{casizero} Si~$0\in I$, il résulte du~\ref{casipaszero} que~$f$ est inversible si et seulement si les deux conditions suivantes sont satisfaites : 

\medskip
i)~$a_0=f(0)\neq 0$;

ii) il existe un entier~$j$ (nécessairement unique) tel que~$|a_j|r^j> |a_i|r^i$ pour tout~$r\in I\setminus\{0\}$ et pour tout~$i\neq j$. 

\medskip
L'unique entier~$j$ de la condition ii) ne peut, compte-tenu de i), être égal qu'à~$0$. Par conséquent,~$f$ est inversible si et seulement si~$a_0\neq 0$ et~$|a_i|r^i<|a_0|$ pour tout~$r\in I$ et tout~$i>0$ ; si c'est le cas, la  norme de~$f$ est identiquement égale à~$|a_0|$ sur~$\DD_I$. 

\deux{courdiscnoniso} Soit~$I$ un intervalle de~$\RR_+$ rencontrant~$\RR\ti_+$. Il découle du~\ref{casizero} ci-dessus que si~$0\in I$ et si~$f$ est une fonction inversible sur~$\DD_I$, il existe un scalaire~$a\neq 0$ tel que~$|f-a|<|a|$ en tout point de~$\DD_I$ ; cette propriété est prise en défaut si~$0\notin I$, la fonction~$T$ constituant un contre-exemple. Il s'ensuit que si~$I$ et~$J$ sont deux intervalles de~$\RR_+$ rencontrant~$\RR\ti_+$ et si~$\DD_I\simeq \DD_J$ alors~$I$ et~$J$ ou bien contiennent tous deux~$0$ ou bien sont tous deux contenus dans~$\RR\ti_+$. 

\deux{finplatvali} Soit~$I$ un intervalle de~$\RR_+$ rencontrant~$\RR\ti_+$, et soit~$f=\sum a_iT^i$ un élément de~$k_I$. On note~$\phi$
le morphisme~$\DD_I\to \Aff^{1,\rm an}_k$ induit par~$f$. On suppose qu'il existe un entier~$j\in \ZZ$ (nécessairement unique) tel que~$|a_j|r^j>|a_i|r^i$ pour tout~$r\in I\setminus\{0\}$ et tout~$i\neq j$ ; cela équivaut à dire, en vertu de~\ref{casipaszero}, que~$f$ ne s'annule pas sauf peut-être à l'origine si~$0\in I$ ; nous dirons que~$a_jT^j$ est le {\em monôme strictement dominant} de~$f$.  On note~$\Lambda$ l'application de~$\RR_+$ dans lui-même qui envoie~$r$ sur~$|a_j|r^j$. 

Si~$0\in I$, notre hypothèse force~$j$ à être égal à la valuation~$T$-adique de la série~$f$ de~$k_I\subset k[[T]]$. On peut donc écrire dans tous les cas (que~$0$ appartienne ou non à~$I$) la fonction~$f$ comme un produit~$a_jT^j(1+u)$ où~$u$ est une fonction analytique sur~$\DD_I$ telle que~$|u(x)|<1$ pour tout~$x\in \DD_I$. 

\medskip
Notons deux conséquences de cette égalité : si~$y\in \DD_I$ et 
si~$x=\phi(y)$ alors~$|T(x)|=|a_j|\cdot |T(y)|^j=\Lambda(|T(y)|)$, et l'élément~$\red{f(y)}$ de~$\red{\hres(y)}$ est égal à~$\red{a_j}(\red{T(y)})^j$.

\medskip
{\em On suppose de plus que~$j\neq 0$.} L'application continue~$\Lambda$ induit alors un homéomorphisme de~$I$ sur un intervalle~$J$ de~$\RR_+$, et l'on peut voir~$\phi$, ce que nous faisons à partir de maintenant, comme un morphisme de~$\DD_I$ vers~$\DD_J$. Nous allons montrer : que~$\phi$ est fini et plat de degré~$|j|$ ; que~$\phi\inv(\eta_J)=\eta_I$ ; que~$\phi(\eta_r)=\eta_{\Lambda(r)}$ pour tout~$r\in I$ ; et que~$\hres(\eta_r)/\hres(\phi(\eta_r))$ est sans défaut pour tout~$r\in I$. 

\trois{phicompact} {\em Le morphisme~$\phi$ est compact.} Cela résulte du fait que si~$J'$ est un intervalle compact contenu dans~$J$ alors~$\phi\inv(\DD_{J'})=\DD_{\Lambda\inv(J')}$. 

\trois{phisansbord} {\em Le morphisme~$\phi$ est sans bord.} Soit~$y\in \DD_I$ ; posons~$r=|T(y)|$ et~$x=\phi(y)$. Si~$r=0$ ou si~$r$ n'est pas une extrémité de~$I$, alors~$y$ appartient à l'intérieur analytique de~$\DD_I$ et~$\phi$ est par conséquent intérieur en~$y$.

Supposons maintenant que~$r$ soit non nul et soit une extrémité de~$I$. Quitte à composer à la source par l'isomorphisme~$\DD_{I\inv}\simeq \DD_I$ donné par~$T\mapsto T\inv$, on peut supposer que~$r$ est le plus grand élément de~$I$ ; quitte à composer au but par l'isomorphisme~$\DD_J\simeq \DD_{J\inv}$ donné par~$T\mapsto T\inv$, on peut supposer que~$j>0$.

Comme~$(\DD_I,y)$ est le 
domaine analytique du germe sans bord~$(\Aff^{1,{\rm an}}_k,y)$ défini par l'inégalité~$|T|\leq r=|T(y)|$, la réduction~$\red{(\DD_I,y)}$ est l'ouvert~$\PP_{\red{\hres(y})/\red k}\{\red {T(y)}\}$. Comme~$\Lambda(r)$ est le plus grand élément de~$J$ (puisqu'on a supposé~$j>0$), on voit de même que la réduction~$\red{(\DD_J,x)}$ est l'ouvert~$\PP_{\red{\hres(x)}/\red k}\{\red {T(x)}\}$. Par la forme même de~$f$, le plongement de~$\red {\hres(x)}$ dans~$\red{\hres(y)}$ induit par~$\phi$ envoie~$\red {T(x)}$ sur~$\red{a_j}(\red{T(y)})^j$. L'entier~$j$ étant strictement positif, l'image réciproque de~$\red{(\DD_J,x)}=\PP_{\red{\hres(x)}/\red k}\{\red {T(x)}\}$ sur~$\PP_{\red{\hres(y)}/\red k}$ est égale à~$\PP_{\red{\hres(y)}/\red k}\{\red {T(y)}\}$, c'est-à-dire à~$\red{(\DD_I,y)}$ ; il s'ensuit que~$\phi$ est intérieur en~$y$. 

\trois{preconclfin} {\em Description des fibres de~$\phi$.} Soit~$x\in \DD_I$ ; nous allons montrer que~$\phi\inv(x)$ est de la forme~$\sch M(\sch A)$ où~$\sch A$ est une~$k$-algèbre de Banach finie de rang~$|j|$. 

\medskip
Prouvons-le tout d'abord lorsque~$x$ est~$k$-rationnel. Il est alors donné par une équation de la forme~$T=\alpha$, où~$\alpha$ est un élément de~$k$ tel que~$|\alpha|\in J$. Posons~$r=\Lambda\inv(\alpha)$. La fibre~$\phi\inv(x)$ est le lieu des zéros de~$f-\alpha$. 

 Si~$\alpha=0$ on a~$0\in J$ et donc~$0\in I$ ; l'entier~$j$ est alors strictement positif et la fonction~$f$ est, comme on l'a vu, égale à~$T^jg$ où~$g$ est une fonction inversible sur~$\DD_J$. Il s'ensuit que la fibre~$\phi\inv(0)$ est, en tant qu'ensemble, égale à~$\{0\}$ ; en tant qu'espace analytique, elle s'identifie par conséquent à~$\sch M(k_{[0;s]}/T^jg)$ pour n'importe quel réel~$s$ de~$I\setminus\{0\}$, et donc à~$\sch M(k[T]/T^j)$, d'où notre assertion. 
 
 Si~$\alpha\neq 0$ la fibre~$\phi\inv(x)$ est contenue dans~$\DD_r$ et s'identifie dès lors à~$\sch M(k_r/(f-\alpha))$. On peut écrire~$f-\alpha=\sum b_iT^i$, où~$b_i=a_i$ si~$i\neq 0$ et~$b_0=a_0-\alpha$ ; on a~$|b_j|r^j=|b_0|=|\alpha|$ et~$|b_i|r^i<|\alpha|$ pour tout~$i\notin\{0,j\}$. On déduit alors de~\ref{factf} que~$k_r/(f-\alpha)$ est une~$k$-algèbre de Banach finie de rang~$|j|$, ce qui prouve là encore notre assertion. 

\medskip
Supposons maintenant~$x$ quelconque. Il existe une extension complète~$L$ de~$k$ telle que~$x$ ait un antécédent~$L$-rationnel~$x'$ sur~$\DD_{J,L}$. Il résulte de ce qui précède que~$\phi_L\inv(x')$ est compacte et de dimension nulle ; la fibre~$\phi\inv(x)$ hérite de ces propriétés, ce qui revient à dire qu'elle est de la forme~$\sch M(\sch A)$ où~$\sch A$ est une~$k$-algèbre de Banach finie. On a dès lors~$\phi_L\inv(x')\simeq \sch M(\sch A\otimes_kL)$ ; comme le rang du~$L$-espace vectoriel~$\sch A\otimes_kL$ est égal à~$|j|$ en vertu du cas particulier déjà traité, le rang du~$k$-espace vectoriel~$\sch A$ est égal à~$|j|$, ce qui achève la démonstration. 

\trois{conclfini} {\em Le morphisme~$\phi$ est fini et plat de degré~$|j|$.} La flèche ~$\phi$ est compacte et sans bord, c'est-à-dire propre ; étant de surcroît à fibres finies, elle est finie (\ref{finpropfin}) ; comme~$\DD_I$ est réduit et comme~$\DD_J$ est normal,~$\phi$ est plate (\ref{nonconstfiniplat}) ; on déduit alors de l'étude des fibres menée ci-dessus que~$\deg \phi=|j|$. 

\trois{imsquel}{\em Image réciproque de~$\eta_J$ et étude de~$\phi(\eta_r)$}. Si~$0\in J$ on a vu ci-dessus que la fibre~$\phi\inv(\eta_0)=\phi\inv(0)$ est réduite à~$\{0\}=\{\eta_{0}\}$ ;  il découle des égalités~$\hres(\eta_0)=\hres(\phi(\eta_0))=k$ que l'extension~$\hres(\eta_0)/\hres(\phi(\eta_0))$ est triviale, et en particulier sans défaut. 

\medskip
Soit maintenant~$r$ un élément de~$I\setminus\{0\}$. Posons~$y=\eta_r$ et~$x=\phi(y)$ ; on a vu en début de preuve que~$|T(x)|=\Lambda(|T(y)|)=\Lambda(r)$. 

L'élément~$\red{T(y)}$ de~$\red{\hres(y)}$
est transcendant sur~$\red k$ ; par conséquent,~$\red{a_j}(\red{T(y)})^j$ est transcendant sur~$\red k$ (l'entier~$j$ est non nul), ce qui implique que~$x=\eta_{\Lambda(r)}$. 

\medskip
Le plongement~$\red{\hres(x)}\hookrightarrow \red{\hres(y)}$ induit par~$\phi$ envoie~$\red{T(x)}$ 
sur~$\red{a_j}(\red{T(y)})^j$. L'élément~$\tau^j-\red{T(x)}/\red{a_j}$
de~$\red {\hres(x)} [r\inv \tau]$
est irréductible, 
par exemple par l'avatar gradué du critère d'Eisenstein qui vaut pour l'annéloïde principal~$\red k[\red{T(x)}]$. 
Le degré de~$\red{\hres(y)}$ sur~$\red{\hres(x)}$ vaut donc~$|j|$ ; par conséquent, le degré de~$\hres(y)$ sur~$\hres(x)$ est au moins égal à~$|j|$ ; il est par ailleurs au plus égal à~$|j|$ puisque~$\phi$ est de degré~$j$. 

On a donc~$[\hres(y):\hres(x)]=|j|=\deg\phi$, ce qui implique que~$y$ est l'unique antécédent de~$x$ ; puisque~$[\red{\hres(y)}:\red{\hres(x)}]=|j|$, l'extension~$\hres(y)/\hres(x)$ est sans défaut. 

\medskip
On a bien établi que ~$\phi\inv(\eta_J)=\eta_I$, que~$\phi(\eta_r)=\eta_{\Lambda(r)}$ pour tout~$r\in I$, et que~$\hres(\eta_r)/\hres(\phi(\eta_r))$ est sans défaut pour tout~$r\in I$. 

\deux{remdegvalj} Soit~$I$ un intervalle de~$\RR_+$ rencontrant~$\RR\ti_+$ et soit~$f=\sum a_iT^i$ une fonction appartenant à~$k_I$ ; supposons que le monôme~$a_0$ de~$f$ soit strictement dominant. La fonction~$f$ peut alors s'écrire~$a_0(1+u)$, où~$u\in \sch O_{\DD_I}(\DD_I)\zeroo$. Si~$\phi$ désigne le morphisme de~$\DD_I$ vers~$\Aff^{1,{\rm an}}_k$ induit par~$f$, il résulte de l'écriture de~$f$ que~$|f|=|a_0|$ identiquement sur~$\DD_I$, 
et que~$\phi(\DD_I)$ est contenu dans le disque de centre~$a_0$ et de rayon~$|a_0|$ ; par conséquent,~$\phi(\DD_I)$ ne rencontre pas~$\eta_{[0;+\infty[}$. 

\deux{defcoorddi} Soit~$I$ un intervalle de~$\RR_+$ rencontrant~$\RR\ti_+$ ; on dira qu'une fonction analytique~$f\in \sch O_{\DD_I}(\DD_I)$ est une {\em fonction coordonnée} sur~$\DD_I$ si et seulement si elle induit un isomorphisme~$\DD_I\simeq \DD_J$ pour un certain intervalle~$J$ de~$\RR_+$ rencontrant~$\RR\ti_+$. 

\deux{isomcour} Soit~$I$ un intervalle non vide de~$\RR\ti_+$ et soit~$J$ un intervalle de~$\RR_+$ rencontrant~$\RR\ti_+$ ; soit~$f=\sum a_iT^i$ un élément de~$k_I$ et soit~$\phi: \DD_I\to \Aff^{1,{\rm an}}_k$ le morphisme induit par~$f$. 

\trois{fimplisomcour} Supposons qu'il existe~$j\in\{-1,1\}$ tel que le monôme~$a_jT^j$ de~$f$ soit strictement dominant, et que~$J=|a_j|.I^j$ (ce qui implique que~$0\notin J$) ; dans ce cas, il résulte de~\ref{finplatvali} que~$\phi$ induit un morphisme fini et plat de degré~$1$, c'est-à-dire un isomorphisme, entre~$\DD_I$ et~$\DD_J$. 

\trois{isomcourimplf} Réciproquement, supposons que~$\phi$ induise un isomorphisme entre~$\DD_I$ et~$\DD_J$. Dans ce cas~$0$ n'appartient pas à~$J$ (\ref{courdiscnoniso}), et~$f$ est donc inversible. Par conséquent, elle possède un monôme strictement dominant~$a_jT^j$. L'entier~$j$ ne peut être nul, car sinon~$\phi(\DD_I)$ ne rencontrerait pas~$\eta_J$ (\ref{remdegvalj}), contredisant par là la surjectivité de~$\phi$ ; par conséquent,~$\phi$ induit un morphisme fini et plat de degré~$|j|$ de~$\DD_I$ sur~$\DD_{|a_j|.I^j}$ (\ref{finplatvali}) ; il s'ensuit que~$|j|=1$, que~$J=|a_j|.I^j$, et que~$\phi$ induit l'homéomorphisme~$\eta_r\mapsto \eta_{|a|.^j}$ de~$\eta_I$ vers~$\eta_J$. 

\medskip
\trois{conclucoordcour} Ainsi,~$f$ est une fonction coordonnée sur~$\DD_I$ si et seulement~$a_jT^j$ est strictement dominant pour un certain~$j\in \{1,-1\}$. Si c'est le cas,~$f$ induit un isomorphisme~$\DD_I\simeq \DD_{|a_j|I^j}$ ; si de plus~$I$ est infini, et si l'on oriente~$\eta_I=\{\eta_r\}_{r\in i}$ dans le sens où~$r$ croît, alors~$j=1$ (resp.~$j=-1$) si et seulement si~$|f|_{|\eta_I}$ est croissante (resp. décroissante). 

\deux{isomdisc} Soit~$I$ un intervalle de~$\RR_+$ contenant~$0$ et rencontrant~$\RR\ti_+$, et soit~$J$ un intervalle de~$\RR_+$ rencontrant~$\RR\ti_+$. 

\trois{fimplisomdisc} Soit~$f=\sum a_iT^i$ un élément de~$k_I$. Supposons que le monôme~$a_1T$ de~$f$ soit strictement dominant et que~$J=|a_1|.I$ (ce qui implique que~$0\in J$) ; soit~$a$ un élément de~$k$ dont la valeur absolue appartient à~$J$. Il résulte de ~\ref{finplatvali} que~$f$ induit un morphisme fini et plat de degré~$1$, c'est-à-dire un isomorphisme, entre~$\DD_I$ et~$\DD_J$. La translation par~$a$ définissant un automorphisme de~$\DD_J$, la fonction~$f+a$ induit un isomorphisme entre~$\DD_I$ et~$\DD_J$. 

\trois{isomdiscimplf} Réciproquement, soit~$g\in k_I$ induisant un isomorphisme~$\psi$ entre~$\DD_I$ et~$\DD_J$ ; l'existence même de~$g$ force~$J$ à contenir~$0$ (\ref{courdiscnoniso}). L'image de l'origine par~$\phi$ est un point~$k$-rationnel de~$\DD_J$, c'est-à-dire un élément~$a\in k$ tel que~$|a|\in J$. La translation par~$(-a)$ définissant un automorphisme de~$\DD_J$, sa composée avec~$\psi$ est un isomorphisme~$\phi: \DD_I\to \DD_J$ qui est induit par~$f:=g-a$ et qui envoie l'origine sur l'origine ; par conséquent,~$f(0)=0$, et l'injectivité (ensembliste) de~$\phi$ implique que~$f$ est inversible sur~$\DD_I\setminus\{0\}$. On en déduit que~$f$ possède un monôme strictement dominant~$a_jT^j$ avec~$j>0$ ; en vertu de~\ref{finplatvali}, la fonction~$f$ induit un morphisme fini et plat de degré~$j$ de~$\DD_I$ sur~$\DD_{|a_j|I^j}$ ; comme~$\phi$ est un isomorphisme,~$j=1$ et~$J=|a_1|.I$ .

\trois{conclucoorddisc} Ainsi, une fonction~$g$ sur~$\DD_I$ est une fonction coordonnée si et seulement si elle est de la forme~$a+f$ où~$f=\sum a_iT^i$ est une fonction dont le monôme~$a_1T$ est strictement dominant et où~$a$ est un scalaire tel que~$|a|\in |a_1|I$, et si c'est le cas~$g$ induit un isomorphisme~$\DD_I\simeq \DD_{|a_1|I}$. 

\deux{remoy} Soit~$I$ un intervalle non vide de~$\RR\ti_+$. On a vu au~\ref{casipaszero} que si~$f$ est une fonction analytique sur~$\DD_I$, elle est inversible si et seulement si elle possède un monôme strictement dominant~$a_jT^j$ (\ref{finplatvali}). Si c'est le cas, elle s'écrit~$a_jT^j(1+u)$, où~$u\in \sch O_{\DD_I}(\DD_I)\zeroo$, et la fonction continue
$|f|$ sur~$\DD_I$ est égale à~$|a_j|.|T|^j$.

\trois{unpluszeroo} Posons~$\sch Z_I=\sch O_{\DD_I}(\DD_I)\ti/k\ti.(1+\sch O _{\DD_I}(\DD_I)\zeroo)$. On déduit de~\ref{remoy} que l'application qui envoie une fonction inversible sur le degré de son monôme quasi-dominant induit un isomorphisme de groupes~$\sch Z_I\simeq  \ZZ$ ; en vertu de~\ref{conclucoordcour}, une fonction analytique~$f$ sur~$\DD_I$ en est une fonction coordonnée si et seulement si elle est inversible et s'envoie sur un générateur de~$\sch Z_I$. 

\medskip
Soit~$J$ un intervalle non vide de~$I$. Il résulte des définitions que le diagramme~$$\diagram \sch Z_I\rto^\sim\dto &\ZZ\dto  \\ \sch Z_J\rto^\sim&\ZZ \enddiagram$$ est commutatif ; en particulier,~$\sch Z_I\to \sch Z_J$ est un isomorphisme. 

\medskip
Notons deux conséquences de ce qui précède : 

\medskip
$\bullet$ la restriction induit également un isomorphisme~$$\sch O_{\DD_I}(\DD_I)\ti/(1+\sch O _{\DD_I}(\DD_I)\zeroo)\simeq \sch O_{\DD_J}(\DD_J)\ti/(1+\sch O _{\DD_J}(\DD_J)\zeroo)\;;$$  

$\bullet$ si~$f\in \sch O_{\DD_I}(\DD_I)\ti$ alors~$f$ est une fonction coordonnée sur~$\DD_I$ si et seulement si~$f_{|\DD_J}$ est une fonction coordonnée sur~$\DD_J$. 

\trois{gautpask} Soit~$\phi$ un automorphisme de~$\DD_I$ {\em au-dessus d'un automorphisme~$\mathsf g$ de~$k$ qui n'est pas nécessairement trivial}. L'automorphisme~$\phi$ agit sur le groupe~$\sch Z_I$ par multiplication par un certain~$\epsilon\in\{-1,1\}$ ; on dira que~$\phi$ est {\em direct} (resp. indirect) si~$\epsilon=1$ (resp.~$-1$). 

La fonction~$\phi^*T$ est de la forme~$aT^\epsilon(1+u)$ avec~$a\in k\ti$ et~$u\in \sch O_{\DD_I}(\DD_I)\zeroo$. L'automorphisme~$\phi$ se décompose sous la forme~$\diagram \DD_I\rto ^\psi&\DD_I\rto ^\pi&\DD_I\enddiagram$, où~$\pi$ est le morphisme de changement de base associé à l'automorphisme~$\mathsf g$ de~$ k$ et où~$\psi$ est un~$k$-automorphisme, induit par la fonction coordonnée~$aT^\epsilon(1+u)$. 

La flèche~$\pi$ fixe~$\eta_I$ point par point (\ref{etainvgal}), et~$\psi$ induit un automorphisme de~$\eta_I$ qui est nécessairement de la forme~$\eta_r\mapsto \eta_{|a|.r^\epsilon}$ (\ref{isomcour} {\em et sq.}). Par conséquent,~$\phi(\eta_I)=\eta_I$ et l'homéomorphisme de~$\eta_I$ induit par~$\phi$ est égal à~$\eta_r\mapsto \eta_{|a|.r^\epsilon}$. Il s'ensuit : 

\medskip
$\bullet$ que si~$I$ est infini alors~$\phi$ est direct si et seulement si il préserve les deux orientations de~$\eta_I$ (cette condition est automatiquement satisfaite si~$I$ est semi-ouvert) ; 

$\bullet$ que si~$\phi$ est direct alors~$|a|.I=I$, ce qui entraîne que~$|a|=1$ (et donc que~$\phi$ fixe~$\eta_I$ point par point) sauf éventuellement si~$I=\RR\ti_+$ ; 

$\bullet$ que si~$\phi$ est direct et si~$\phi$ appartient à un groupe compact d'automorphismes agissant sur~$\DD_I$, alors~$\phi$ fixe~$\eta_I$ point par point (et ce, même si~$I$ est égal à~$\RR\ti_+$) : pour le voir, on remarque que sous notre hypothèse, l'ensemble~$\{\phi^n(\eta_r)\}_n=\{\eta_{|a|^nr}\}_n$ est contenu dans un compact de~$\DD_I$, ce qui signifie que~$\{|a|^nr\}_r$ est contenu dans un compact de~$\RR\ti_+$, et donc que~$|a|=1$. 

\trois{casdescnorm} {\em Remarque.} Si~$I$ est ou bien infini, ou bien de la forme~$\{r\}$ avec~$r\in \RR\ti_+-\sqrt{|k\ti|}$, l'application qui envoie une fonction inversible~$f$ sur~$|f|$ induit, en vertu de l'expression de~$|f|$ donnée au~\ref{remoy}, un isomorphisme~$$\sch O_{\DD_I}(\DD_I)\ti/(1+\sch O _{\DD_I}(\DD_I)\zeroo)\simeq |\sch O_{\DD_I}(\DD_I)\ti|.$$

On peut dès lors, au~\ref{unpluszeroo}, remplacer~$\sch O_{\DD_I}(\DD_I)\ti/(1+\sch O _{\DD_I}(\DD_I)\zeroo)$ par~$|\sch O_{\DD_I}(\DD_I)\ti|$ et~$\sch Z_I$ par ~$|\sch O_{\DD_I}(\DD_I)\ti|/|k\ti|$, et de même concernant~$J$ s'il est lui aussi de la forme requise.

\deux{extremcour} Soit ~$X$ une courbe~$k$-analytique et soit~$U$ un ouvert de~$X$ qui est isomorphe à~$\DD_I$ pour un certain intervalle ouvert non vide~$I$ de~$\RR\ti_+$. Soit~$x$ un point de~$\partial U$ soit~$\omega$ un bout de~$U$ convergeant vers~$x$ (il peut arriver que les deux bouts de~$U$ convergent vers~$x$, et~$\skel {\overline U}$ est alors un cercle, {\em cf.}~\ref{constboutbr2} ). Soit~$r$ la borne de~$I$ correspondant à~$\omega$. Nous allons établir les faits suivants. 

\medskip
a) Le point~$x$ est~$k$-rationnel si~$r\in \{0,+\infty\}$, de type 2 si~$r \in \sqrt {|k\ti|}$, et de type 3 si~$r\in \RR\ti_+-\sqrt{|k\ti|}$.

\medskip
b) Si~$x$ est de type 2 ou 3 alors~$|\hres(x)\ti|\subset |k\ti|.r^\ZZ$, et~$X$ est quasi-lisse en~$x$. 

\medskip
Pour ce faire, on peut toujours supposer (quitte à composer par~$T\mapsto T\inv$) que~$r$ est la borne supérieure de~$I$ ; on fixe un isomorphisme entre~$U$ et~$\DD_{]r_0;r[}$ pour un certain~$r_0\in [0;+\infty[$.

\trois{premremtypx} {\em Une première remarque.} Soit~$f$ une fonction définie sur un voisinage~$V$ de~$x$ et inversible sur~$V\setminus\{x\}$. Il existe~$s\in ]r_0;r[$ tel que l'intersection de~$V$ avec~$U=\DD_{]r_0;r[}$ contienne~$\DD_{]s;r[}$. La restriction de~$f$ à~$\DD_{]s;r[}$ est inversible ; elle s'écrit donc, en vertu de~\ref{casipaszero}, ~$a_iT^i +\sum\limits_{j\in \ZZ\setminus\{i\}}a_jT^j$ pour un certain~$i\in \ZZ$, avec~$a_i\neq 0$ et~$|a_j|t^j<|a_i|t^i$ pour tout~$j\neq i$ et tout~$t\in ]s;r[$. On a alors~$|f|=|a_i|.|T|^i$ sur~$\DD_{]s;r[}$ ; par passage à la limite,~$|f(x)|=|a_i|.r^i$. Nous dirons que~$a_iT^i$ est {\em le monôme strictement dominant de~$f$ le long de~$\omega$.}

\trois{ypaspast4} {\em Si~$x$ est de type 1 ou 4, il est~$k$-rationnel.} Comme~$x$ est de type 1 ou 4, toute fonction inversible au voisinage de~$x$ est de norme constante au voisinage de~$x$, et a donc un monôme strictement dominant le long de~$\omega$ qui est de degré nul. 

Soit~$\lambda\in \kappa(x)\ti$ et soit~$f$ une fonction définie et inversible au voisinage de~$x$ telle que~$f(x)=\lambda$ ; par ce qui précède, son monôme strictement dominant le long de~$\omega$ est de la forme~$a$ avec~$a\in k$. 

Si~$f-a$ était inversible en~$x$, son monôme strictement dominant le long de~$\omega$ serait par construction non nul, ce qui est contradictoire ; par conséquent,~$\lambda=f(x)=a$ et~$\kappa(x)=k$.

\trois{yrigomegainf} {\em Le cas où~$x$ est~$k$-rationnel, de type 2 ou de type 3}. Nous allons tout d'abord nous donner une fonction~$f$ définie sur un voisinage ouvert~$V$ de~$x$ et inversible sur~$V\setminus\{x\}$, par un procédé dépendant du type de~$x$. 

\medskip
{\em  Le cas où~$x$ est~$k$-rationnel.} Choisissons un voisinage affinoïde~$W$ de~$x$ dans~$X$, soit~$\sch A$ l'algèbre des fonctions analytiques sur~$W$, et soit~$\bf x$ le point fermé de~$\spec \sch A$ correspondant à~$x$. Il existe~$f\in \sch A$ tel que~$\bf x$ soit un point isolé du lieu des zéros (ensembliste) de~$f$ sur~$\spec \sch A$ ; par conséquent,~$f(x)=0$ et il existe un voisinage ouvert~$V$ de~$x$ dans~$X$ contenu dans~$W$ et tel que~$f$ soit inversible sur~$V\setminus\{x\}$. 

{\em Le cas où~$x$ est de type 2}. On choisit~$f$ de sorte que~$f(x)\notin k$. 

{\em Le cas où~$x$ est de type 3}. On choisit~$f$ de sorte que~$|f(x)|\in \RR\ti_+-\sqrt{|k\ti|}$. 

\medskip
On note~$aT^i$ le monôme strictement dominant de~$f$ le long de~$\omega$. 

\medskip
{\em Supposons que~$x$ est~$k$-rationnel.} On a alors~$0=|f(x)|=|a|.r^i$, ce qui entraîne que~$r=+\infty$ (et que~$i<0$). 

\medskip
{\em Supposons que~$x$ est de type 2.} Si~$i=0$ alors comme~$f(x)\notin k$, la fonction~$f-a$ est inversible au voisinage de~$x$, et si~$bT^j$ désigne son monôme strictement dominant on a nécessairement~$j\neq 0$ ; en remplaçant~$f$ par~$f-a$,  on se ramène au cas où~$i\neq 0$. 

On a~$|f(x)|=|a|.r^i$. Comme~$x$ est de type 2, on a~$|f(x)|\in \sqrt{|k\ti|}$ ; il s'ensuit,~$i$ étant non nul, que~$r\in \sqrt{|k\ti|}$.

\medskip
{\em Supposons que~$x$ est de type 3}. On a alors d'une part~$|f(x)|\in \RR\ti_+-\sqrt{|k\ti|}$, et d'autre part~$|f(x)|=|a|.r^i$. Par conséquent,~$r\in \RR\ti_+-\sqrt{|k\ti|}$ (et~$i\neq 0$). 

\trois{borddiscql} Supposons maintenant que~$x\in X\dtr$. Si~$f$ est une fonction inversible au voisinage de~$x$, il résulte de~\ref{premremtypx} que~$|f(x)|=|a|r^i$ pour un certain~$a\in k\ti$ et un certain~$i\in \ZZ$. Par conséquent,~$|\hres(x)\ti|\subset |k\ti|.r^\ZZ$. 

\medskip
Soit~$Y$ le lieu des points en lequel~$X$ n'est pas quasi-lisse : c'est un fermé de Zariski de~$X$. Le point~$x$ étant de type 2 ou 3, son adhérence de Zariski~$X'$ dans~$X$ n'est autre que sa composante irréductible, et c'est un voisinage de~$x$ dans~$X$. Comme~$U\simeq \DD_{]r_0;r[}$, l'ouvert~$U$ est lisse, et~$Y\cap U=\emptyset$ ; par conséquent,~$Y$ ne contient aucun voisinage de~$x$ ; en particulier,~$Y$ ne contient pas~$X'$, et ne contient dès lors pas~$x$, ce qu'on souhaitait établir. 

\subsection*{Disques et couronnes : le point de vue intrinsèque}

\deux{rappbemol} Rappelons que si~$X$ est un graphe et~$\Delta$ un sous-graphe faiblement admissible de~$X$, on désigne par~$\Delta^\flat$ l'unique ouvert de~$X$ dont~$\Delta$ soit un sous-graphe admissible. La plupart du temps, nous appliquerons cette définition lorsque~$X$ sera un arbre à un bout~$\omega$ et~$\Delta$ de la forme~$]x;\omega[$ avec~$x\in X$ (qu'un tel sous-arbre soit faiblement admissible découle de~\ref{admunbout}).

\deux{defdc} Un {\em~$k$-disque} est un espace~$k$-analytique {\em isomorphe} à un ouvert~$\DD_{[0;r[}$ de~$\Aff^{1,{\rm an}}_k$ pour un certain~$r>0$ ; une fonction induisant un tel isomorphisme sera qualifiée de {\em fonction coordonnée} sur le disque en question ; cette définition est, en vertu de~\ref{conclucoorddisc}, compatible avec celle donnée au~\ref{defcoorddi}. 

\medskip
Un $k$-disque est non vide, lisse, et géométriquement connexe. 

\deux{compp1msd} {\em Exemple.} Soit~$x$ un point de type 2 ou 3 de~$\pk$, et soit~$Z$ une composante connexe de~$\pk\setminus\{x\}$ possédant un~$k$-point~$z$. Supposons qu'il existe une composante connexe~$Z'$ de~$\pk\setminus\{z\}$ distincte de~$Z$ et possédant un~$k$-point~$z'$ ; dans ce cas,~$Z$ est un disque (et il en va de même de~$Z'$ par symétrie). 

\medskip
Pour le voir, on se ramène par une homographie au cas où~$z=0$ et où~$z'=\infty$. Étant situé sur le segment qui joint~$0$ à~$\infty$, le point~$x$ est égal à~$\eta_r$ pour un certain~$r>0$, et~$Z$ est alors nécessairement le disque ouvert de centre~$0$ et de rayon~$r$.

\deux{introdisc} Soit~$X$ un~$k$-disque ; choisissons une fonction coordonnée~$f$ sur~$X$ ; elle induit un isomorphisme~$X\simeq \DD_{[0;r[}$ de~$\Aff^{1,{\rm an}}_k$ avec~$0<r$.

\trois{topdisc} L'espace topologique~$X$ est un arbre paracompact à un bout ; son squelette est par conséquent vide.

\trois{radisque} Il résulte de~\ref{conclucoorddisc} : qu'une fonction~$g$ sur~$X$ en est une fonction coordonnée si et seulement si elle s'écrit~$a+\sum_{i>0} a_if^i$, où~$\sum a_iT^i$ est une série appartenant à~$k_{[0;r[}$ dont le monôme~$a_1T$ est strictement dominant, et où~$a$ est un scalaire tel que~$|a|<|a_1|.r$ ; et que si c'est le cas,~$g$ induit un isomorphisme~$X\simeq \DD_{[0;|a_1|.r[}$. 

\medskip
Le réel~$r$ de la définition~\ref{defdc} est donc bien déterminé modulo~$|k\ti|$. On pourra ainsi parler sans ambiguïté du {\em rayon modulo~$|k\ti|$} d'un~$k$-disque, que l'on peut définir comme la borne supérieure, en norme, de n'importe quelle de ses fonctions coordonnées ; il est clair que deux~$k$-disques sont isomorphes si et seulement si ils ont même rayon modulo~$|k\ti|$.

\deux{compatdisc} Un ouvert de~$\Aff^{1,{\rm an}}_k$ est un~$k$-disque si et seulement si il est défini par une inégalité de la forme~$|T-\alpha|<R$, avec~$\alpha\in k$ et~$R>0$ : cela vient du fait que tout ouvert {\em strict} de~$\Aff^{1,{\rm an}}_k$ ayant exactement un bout et possédant un~$k$-point est de cette forme, et que~$\Aff^{1,{\rm an}}_k$ lui-même n'est pas un~$k$-disque en raison du théorème de Liouville.  

\deux{defdcour} Si~$I$ est un intervalle non vide de~$\RR\ti_+$, une {\em~$k$-couronne de type~$I$} est un espace~$k$-analytique {\em isomorphe} à~$\DD_I$ ; une fonction induisant un tel isomorphisme sera qualifiée de {\em fonction coordonnée} sur la couronne en question ; cette définition est, en vertu de~\ref{conclucoordcour}, compatible avec celle donnée au~\ref{defcoorddi}. 

\medskip
Une $k$-couronne est non vide, quasi-lisse et géométriquement connexe. 

\trois{dependtype} Si~$X$ est une couronne de type~$I$, alors~$X_L$ est une~$L$-couronne de type~$I$ pour tout extension complète~$L$ de~$k$. Il résulte de~\ref{isomcour} {\em et sq.} que l'intervalle~$I$ est uniquement déterminé {\em à une transformation de~$\RR\ti_+$ près de la forme~$r\mapsto \lambda r^j$ avec~$\lambda \in |k\ti|$ et~$j\in\{-1,1\}$} ; lorsqu'on parlera {\em du} type de~$X$, il s'agira donc d'une classe d'intervalles.

\trois{nottypes} On appellera~$k$-couronne de type ~$]*,*[$ (resp.~$]0,*[$, resp.~$[*,+\infty[$, etc.) toute~$k$-couronne de type~$I$ pour un certain intervalle~$I$ de la forme~$]r;R[$  avec~$0<r<R$ (resp.~$]0;r[$ avec~$r>0$, resp.~$[r;+\infty[$ avec~$0<r$, etc.) ; notons que les notions de~$k$-couronne de type~$]0;*[$ et~$]*,+\infty[$ coïncident, et qu'une~$k$-couronne de type~$]0;+\infty[$ est un espace~$k$-analytique isomorphe à~$\Aff^{1,\rm an}_k$. On dira qu'une~$k$-couronne est {\em ouverte} si son type est ouvert. 

\deux{compp1msc} {\em Exemple.} Soient~$x$ et~$x'$ deux points distincts de~$\pk$, et soit~$U$ la composante connexe de~$\pk\setminus\{x,x'\}$ contenant~$]x;x'[$ (que l'on peut également décrire comme l'ouvert~$]x;x'[^\flat$). 

\trois{xxprimet23} Supposons que~$x$ et~$x'$ sont de type 2 ou 3 et qu'il existe deux composantes connexes~$Z$ et~$Z'$ de~$\pk\setminus\{x,x'\}$ possédant les propriétés suivantes : 

\medskip
$\bullet$~$Z(k)\neq\emptyset$ et~$Z'(k)\neq \emptyset$ ; 

$\bullet$~$\partial Z=\{x\}$ et~$\partial Z=\{x'\}$. 

\medskip
Sous ces hypothèses,~$U$ est une~$k$-couronne de type~$]*,*[$. Pour le voir, on choisit~$z\in Z(k)$ et~$z'\in Z'(k)$, puis l'on se ramène par une homographie au cas où~$z=0$ et où~$z'=\infty$. Comme~$x$ et~$x'$ sont tous deux situés sur le segment qui joint~$0$ à~$\infty$, il existe deux réels strictement positifs distincts~$r$ et~$r'$ tels que~$x=\eta_r$ et~$x'=\eta_{r'}$ ; notre assertion s'ensuit aussitôt. 

\trois{xkrat} Supposons que~$x$ est~$k$-rationnel, que~$x'$ est de type 2 ou 3 et qu'il existe~$Z'$ comme ci-dessus ; alors~$U$ est une~$k$-couronne de type~$]0;*[$ : pour le voir, on choisit~$z'\in Z'(k)$, puis l'on se ramène par une homographie au cas où~$x=0$ et où~$z'=\infty$. Comme~$x'$ est situé sur le segment qui joint~$0$ à~$\infty$, il existe ~$r'>0$ tel que~$x=\eta_{r'}$ ; notre assertion s'ensuit aussitôt.

\trois{deuxkrat} Supposons que~$x$ et~$x'$ soient~$k$-rationnels ; alors~$U$ est une~$k$-couronne de type~$]0;+\infty[$ : pour le voir, on se ramène par une homographie au cas où~$x=0$ et où~$x'=\infty$, et l'assertion est alors triviale.

\deux{oxmodulo} Soit~$X$ une~$k$-couronne ; le groupe~$$\sch Z(X):=\sch O_X(X)\ti/k\ti(1+\sch O_X(X)\zeroo)$$ est libre de rang~$1$, et une fonction~$f$ sur~$X$ est une fonction coordonnée si et seulement si elle est inversible et s'envoie sur un générateur de~$\sch Z(X)$ (cela résulte de~\ref{unpluszeroo}). Si le type de~$X$ est infini ou bien de la forme~$\{r\}$ avec~$r\notin\sqrt{|k\ti|}$ alors le groupe ~$\sch Z(X)$ s'identifie à~$|\sch O_X(X)\ti|/|k\ti|$ (rem.~\ref{casdescnorm}).

\deux{introcour} Soit~$X$ une~$k$-couronne ; choisissons une fonction coordonnée~$f$ sur~$X$ ; elle induit un isomorphisme~$X\simeq \DD_I$ pour un certain intervalle~$I$ non vide de~$\RR\ti_+$. 

\trois{topcour} En tant qu'espace topologique,~$X$ est un arbre paracompact ; l'isomorphisme~$X\simeq \DD_I$ identifie~$\eta_I$ à un sous-arbre admissible de~$X$, qui ne dépend pas du choix de~$f$ (\ref{isomcour} {\em et sq.}) ; on l'appelle le {\em squelette analytique} de~$X$ et on le note~$\skelan X$  ; si~$I$ est ouvert alors~$X$ est un arbre à deux bouts et~$\skelan X=\skel X$ ; si~$I$ n'est pas ouvert alors~$X$ est un arbre ayant au plus un bout et~$\skel X=\emptyset$. Si~$L$ est une extension complète de~$k$ alors~$X_L\to X$ induit un homéomorphisme~$\skelan {X_L}\to \skelan X$ (\ref{etainvgal}) .

 \trois{modcour} On déduit du~\ref{conclucoordcour}, qu'une fonction~$g$ sur~$X$ est une fonction coordonnée si et seulement si~$g$ s'écrit~$\sum a_if^i$, où~$\sum a_iT^i$ est une série de~$k_I$ dont~$a_jT^j$ est un monôme strictement dominant pour un certain~$j\in \{-1,1\}$ ; si c'est le cas alors~$g$ induit un isomorphisme~$X\simeq \DD_{|a_j|I^j}$ (cela découle encore de~\ref{conclucoordcour}), et si de plus~$I$ est infini alors le sens de variation de~$|g|$ sur~$\skelan X$ est égal (resp. opposé) à celui de~$|f|$ si et seulement si~$j=1$ (resp.~$j=(-1)$) ; notons que le sens de variation en question ne dépend
 que de la classe de~$|g|$ modulo $|k\ti|$, c'est-à-dire encore de l'image de~$g$ dans~$\sch Z(X)$.  
 
Ainsi, lorsque~$I$ est infini, choisir un générateur de~$\sch Z(X)$ revient à choisir une orientation sur~$\skelan X$ (celle pour laquelle la fonction coordonnée correspondante est croissante en norme). Pour cette raison, nous nous permettrons, même dans le cas où~$I$ est un singleton, d'appeler~{\em orientation de $X$} un générateur de~$\sch Z(X)$.

\medskip
Si~$g$ est une fonction coordonnée sur~$X$, le quotient des bornes supérieure et inférieure de~$|g|$ sur~$X$ est par ce qui précède un élément de~$]0;+\infty[$ égal au quotient des bornes supérieure et inférieure de~$I$ ; il ne dépend donc pas du choix de la fonction coordonnée~$g$, et est appelé le {\em module} de~$X$ et est noté~$\mathsf{Mod}(X)$ ; il est invariant par extension des scalaires, et est égal à~$1$ si et seulement si~$I$ est un singleton (ce qui revient à demander que~$\skelan X$ soit un singleton). 

\medskip
Si~$J$ est un intervalle non vide de~$\skelan X$, on notera~$J^\sharp$ l'image réciproque de~$J$ par la rétraction canonique de~$X$ sur~$\skelan X$ ; si~$J$ est ouvert dans~$\skelan X$ alors~$J^\sharp=J^\flat$  (notons que~$J^\flat$ n'est tout simplement pas défini si~$J$ n'est pas ouvert dans~$\skelan X$). Le sous-ensemble~$J^\sharp$ de~$X$ en est un domaine analytique, qui est une couronne de type~$|g(J)|$ pour n'importe quelle fonction coordonnée~$g$ sur~$X$ ; son squelette analytique coïncide avec~$J$. On qualifiera de {\em sous-couronne} de~$X$ tout domaine analytique de~$X$ de cette forme. 

\medskip
Si~$J$ est un intervalle non vide de~$\skelan X$ alors ~$\mathsf{Mod}(J^\sharp)=\mathsf{Mod}({\overline J}^\sharp)$. Et si~$J_1$ est un intervalle de~$\skelan X$ rencontrant~$J$ alors~$$\mathsf{Mod}((J\cup J_1)^\sharp)=\mathsf{Mod}(J^\sharp).\mathsf{Mod}(J_1^\sharp).(\mathsf{Mod}(J\cap J_1)^\sharp)\inv.$$

\trois{rextmodk} Supposons que~$\skelan X$ n'est pas un singleton, et soit~$\omega$ l'une des deux extrémités du segment~$\wid{\skelan X}$ de la compactification~$\wid X$ ; notons que si~$\skelan X$ est un intervalle ouvert alors~$\omega$ est simplement l'un des deux bouts de~$X$. Choisissons une fonction coordonnée sur~$X$ qui soit croissante en norme sur~$\skelan X$ lorsqu'on oriente celui-ci vers~$\omega$. Sa limite en~$\omega$ est, en vertu de la description des fonctions coordonnées donnée au~\ref{modcour}, un élément de~$]0;+\infty]$ qui est bien déterminé modulo~$|k\ti|$ et que l'on appellera {\em le rayon extérieur de~$X$ en~$\omega$ modulo~$|k\ti|$.} Si~$J$ est un intervalle non vide de~$\skelan X$ aboutissant à~$\omega$, le rayon extérieur de~$J^\sharp$ modulo~$|k\ti|$ en~$\omega$ est égal à celui de~$X$. 

Si~$L$ est une extension complète de~$k$ et si l'on appelle~$\omega_L$ le point de~$\wid{\skelan {X_L}}\simeq \wid {\skelan X}$ qui correspond à~$\omega$, il résulte des définitions que le rayon extérieur de~$X_L$ en~$\omega_L$ est égal à celui de~$X$ en~$\omega$ modulo~$|L\ti|$.

\trois{boutinf} Si~$X$ est de type~$]0,*[$ (ou~$]*,+\infty[$, ce qui revient au même), il y a exactement un bout de~$X$ en lequel son rayon extérieur modulo~$|k\ti|$ est infini ; on l'appellera le {\em bout infini} de~$X$. Si~$L$ est une extension complète de~$k$, il résulte de l'invariance du rayon extérieur par extension des scalaires que le bout infini de~$X_L$ est celui qui correspond au bout infini de~$X$ {\em via} l'homéomorphisme~$\skel {X_L}\simeq \skel X$. 

\trois{souscourtyp} Si~$Y$ est une sous-couronne de~$X$, la restriction des fonctions induit des isomorphismes $$|\sch O_X\ti|\simeq |\sch O_Y\ti|\;{\rm et}\;\sch Z(X)\simeq \sch Z(Y),$$ et une fonction analytique inversible~$g$ sur~$X$ en est une fonction coordonnée si et seulement si sa restriction à~$Y$ est une fonction coordonnée de~$Y$ (\ref{unpluszeroo} et rem.~\ref{casdescnorm}). 

\trois{zxextscal} Si~$L$ est une extension complète de~$k$, toute fonction coordonnée sur~$X$ est encore une fonction coordonnée sur~$X_L$ ; par conséquent, la flèche naturelle~$\sch Z(X)\to \sch Z(X_L)$ est un isomorphisme. 

\trois{dircourgen} Si~$\phi$ est un automorphisme de~$X$ {\em  au-dessus d'un automorphisme de~$k$ qui n'est pas nécessairement trivial}, on dira que~$\phi$ est {\em direct} (resp. {\em indirect}) s'il préserve (resp. permute) les deux orientations de~$X$. Si~$\skelan X$ est infini,~$\phi$ est direct si et seulement si il préserve les orientations de~$\skelan X$, ce qui est automatique si celui-ci est semi-ouvert (\ref{gautpask}).

\deux{morcourintersk} Soient~$X$ et~$Y$ deux~$k$-couronnes et soit ~$\phi: Y\to X$ un morphisme tel que~$\phi(Y)$ rencontre~$\skelan X$. Choisissons une fonction coordonnée~$g$ sur~$Y$ et une fonction coordonnée~$f$ sur~$X$. La fonction~$g$ identifie naturellement~$\skelan Y$ à~$\eta_{|g|(\skelan Y)}$. 

\medskip
L'image réciproque de~$f$ sur~$Y$ est une fonction inversible, elle est donc de la forme~$\sum a_i g^i$ où~$\sum a_iT^i$ est une série admettant un monôme strictement dominant~$a_dT^d$. Comme~$\phi(Y)\cap \skelan X\neq \emptyset$, il découle de~\ref{remdegvalj} que~$d\neq 0$. On déduit alors de~\ref{finplatvali} {\em et sq.} les faits suivants : 

\medskip
$\bullet$ il existe un intervalle non vide~$I$ de~$\skelan X$ tel que~$\phi$ induise un morphisme fini et plat de degré~$|d|$ de~$Y$ sur~$I^\sharp$ ; 

$\bullet$ la fonction~$f$ identifie~$I$ à~$\eta_{|a_d|.|g|(\skelan Y)^d}$ ;

$\bullet$ l'image réciproque de~$I$ sur~$Y$ est égale à~$\skelan Y$, et~$\skelan Y \to I$ est l'homéomorphisme qui correspond à~$\eta_r\mapsto \eta_{|a_d|.r^d}$ {\em via} les identifications évoquées ci-dessus ; 

$\bullet$ pour tout~$y\in I$ d'image~$x$ sur~$X$ l'extension~$\hres(y)/\hres(x)$ est sans défaut et de degré~$|d|$. 

\trois{effmodule} Il s'ensuit, en vertu du lemme~\ref{lemmcommret}, que~$Y\to I^\sharp$ commute aux rétractions canoniques de ses source et but sur leurs squelettes analytiques respectifs~$\skelan Y$ et~$I$. 

Soit~$J$ un intervalle non vide tracé sur~$\skelan Y$. Il découle de ce qui précède que~$J^\sharp=\phi^{-1}(\phi(J)^\sharp)$, et ~$J^\sharp \to \phi(J)^\sharp$ est par conséquent fini et plat de degré~$|d|$ ; le type de~$J^\sharp$ est égal à~$|g|(J)$, et celui de~$\phi(J)$ à~$|f|(\phi(J))=|a_d|.|g|(J)^d$, et donc à~$|g|(J)^d$ ({\em cf.}~\ref{dependtype}).  

\trois{platinddd} En vertu de~\ref{remoy}, l'image de~$\sch Z(X)\simeq \sch Z(I^\sharp)$ dans~$\sch Z(Y)$ est un sous-groupe d'indice~$|d|$ de ce dernier. 

\trois{caspartfinplat} Remarquons que si~$\phi$ est lui-même fini et plat,~$\phi(Y)$ est égal à~$X$ et rencontre donc~$\skelan X$ : ce qui précède est alors valable avec~$I=\skelan X$. 

\deux{raddisccour} Soit~$X$ un~$k$-disque et soit~$Y$ un domaine analytique de~$X$ qui est une~$k$-couronne aboutissant à l'unique bout~$\omega$ de~$X$. Choisissons une fonction coordonnée~$f$ sur~$X$ ; elle induit un isomorphisme~$X\simeq \DD_{[0;R[}$ pour un certain~$R>0$ ; soit~$I$ l'intervalle ouvert de~$X$ correspondant à~$\eta_{[0;R[}$ {\em via} cette identification. Puisque~$\skelan Y$ aboutit à~$\omega$, il existe un sous-intervalle ouvert~$J$ de~$I$ aboutissant à~$\omega$ et contenu dans~$\skelan Y$ ; l'ouvert~$I^\flat$ de~$X$ est une~$k$-couronne aboutissant à~$\omega$, et~$f_{|I^\flat}$ est une fonction coordonnée sur~$I^\flat$, qui est croissante en norme sur~$I$. On en déduit que le rayon extérieur modulo~$|k\ti|$ de~$Y$ en~$\omega$ (qui coïncide avec celui de~$I^\flat$, {\em cf. supra}) est égal à~$R$, c'est à dire au rayon de~$X$ modulo~$|k\ti|$. 

\deux{raynorm} Soit~$X$ un~$\KK$-disque, soit~$x\in X\dtr$ et soit~$\omega$ l'unique bout de~$X$. Choisissons une fonction coordonnée~$f$ sur~$X$ telle que la composante connexe de~$X\setminus\{x\}$ contenant l'unique élément de~$f\inv(0)$ soit relativement compacte (il revient au même de demander qu'elle ne contienne pas~$]x;\omega[$). Le morphisme~$f$ identifie~$X$ à l'ouvert de~$\Aff^{1,{\rm an}}_k$ défini par l'inégalité~$|T|<R$ pour un certain réel~$R>0$ ; par construction, il envoie~$x$ sur un point de l'intervalle joignant l'origine à~$\eta_{R,\KK}$ ; par conséquent,~$f(x)=\eta_{r,\KK}$ pour un certain~$r$ vérifiant~$0<r<R$. 

\medskip
{\em Le rapport~$r/R$ ne dépend que de~$x$, et pas du choix de~$f$.} En effet,~$f$ identifie~$]x;\omega[^\flat$ à l'ouvert de~$\Aff^{1,{\rm an}}_r$ défini par les inégalités~$r<|T|<R$ ; par conséquent,~$]x;\omega[^\flat$ est une couronne et~$r/R$ est égal à {\em l'inverse} de son module. 

\medskip
Nous dirons que~$r/R$ est le {\em rayon normalisé} de~$x$, et nous le noterons~$\rnorm (x)$. La définition de~$\rnorm$ s'étend comme suit à~$X$ tout entier : si~$y\in X\typ {1,4}$ alors~$]y;\omega[\subset X\dtr$ ; l'application~$\rnorm$ est donc une bien définie sur~$]y;x[$ ; elle est par sa construction même continue, strictement décroissante lorsqu'on se dirige vers~$y$, et positive ; elle admet donc une limite en~$y$, que l'on appelle encore le rayon normalisé de~$y$ et que l'on note~$\rnorm(y)$ ; pour tout~$y\in X$, on vérifie immédiatement que~$\rnorm(y)\in [0;1[$ et que~$\rnorm (y)=0$ si et seulement si~$y$ est de type 1, c'est-à-dire~$\KK$-rationnel. 

\medskip
Soit~$I$ un intervalle tracé sur~$X$ et aboutissant à~$\omega$ ; par définition de~$\rnorm$, la restriction de celle-ci à~$I$ est continue, strictement monotone, et tend vers~$1$ à l'approche de~$\omega$. 

\subsection*{Revêtements de Kummer des couronnes} 

\deux{notquot} Si $X$ est un espace~$k$-analytique, et si~$P$ est un polynôme non nul à coefficients dans~$\sch O_X$, 
nous nous permettrons de noter~$\sch M(\sch O_X[T]/P)$ le~$X$-espace 
analytique obtenu par recollement des espaces~$\sch M(\sch O_X(V)[T]/P)$, où~$V$ parcourt l'ensemble des domaines
affinoïdes de~$X$.  

\deux{rappdefkummer} Soit~$\ell$ un entier inversible {\em dans~$k$} et soit~$X$ un espace~$k$-analytique. 
On dispose
sur le site~$X\et$ de la {\em suite exacte
de Kummer} $$\diagram 1\rto& \mu_\ell\rto& \gm\rto^{z\mapsto z^\ell}& \gm\rto &1\enddiagram$$ 
qui induit une injection~$f\mapsto (f)$ de~$\sch O_X(X)\ti/(\sch O_X(X)\ti)^\ell\hookrightarrow \H^1(X\et,\mu_\ell)$, dont
l'image sera notée~$\kum(X)$. 

Tout faisceau étale localement constant sur~$X \et$ est représentable, 
 et c'est en particulier le cas de tout~$\mu_\ell$-torseur. En conséquence, 
$\H^1(X\et,\mu_\ell)$ classifie les~$\mu_\ell$-torseurs étales {\em analytiques}
 sur~$X$ à isomorphisme près ;  la classe d'isomorphie correspondant
 par ce biais à~$(f)$ est celle de~$\sch M(\sch O_X[T]/(T^\ell-f))$ ; nous dirons qu'un tel~$\mu_\ell$-torseur est {\em de Kummer.} 

Si~$Y$ est un espace~$k$-analytique et si~$\phi : Y\to X$ est un morphisme, la
flèche naturelle~$\H^1(X\et,\mu_\ell)\to \H^1(Y\et,\mu_\ell)$ envoie~$\kum(X)$ dans~$\kum(Y)$.

\deux{exemplekumm} {\em Exemple.} Soit~$I$ un intervalle non vide de~$\RR\ti_+$. 
La fonction~$T^\ell$ induit un revêtement fini et plat~$\DD_I\to \DD_{I^\ell}$ ; on vérifie immédiatement
que ce revêtement identifie~$\DD_I$ au torseur de Kummer~$\sch M(\sch O_{\DD_{I^\ell}}[T]/(T^\ell-S))$, où
nous notons~$S$ (et non~$T$) la coordonnée standard de $\DD_{I^\ell}$, pour éviter toute confusion.

\deux{kummercouronne} Soit~$X$ 
une~$k$-couronne. Comme~$X$ est géométriquement
connexe, l'application composée $$k\ti/(k\ti)^\ell\to \sch O_X(X)\ti/(\sch O_X(X)\ti)^\ell\simeq \kum (X)$$
est injective ; on note~$\kum'(X)$ son conoyau. {\em À partir de maintenant, on suppose que $\ell$ est premier à~$p$}. 

\trois{kummcourpremp} Le groupe~$1+\sch O_X(X)\zeroo$ est
alors~$\ell$-divisible. Il en résulte que~$\kum'(X)$ est canoniquement isomorphe
à~$\sch Z(X)/\ell \sch Z(X)$. Il est en conséquence isomorphe
à~$\ZZ/\ell\ZZ$, mais non canoniquement dès que $\ell \geq 3$ : il y a alors {\em deux} générateurs
distingués de $\kum'(X)$, opposés l'un de l'autre ; chacun correspond à une orientation de $X$. 

Si~$k\ti$ est~$\ell$-divisible, ce qui est par exemple le cas
si~$k$ est algébriquement clos, alors~$\kum(X)=\kum'(X)$ ; et ce qu'on vient
de voir concernant~$\kum'(X)$ vaut dès lors pour~$\kum(X)$. 

\trois{kummcourfonct} Soit~$Y$ une sous-couronne de~$X$. La flèche~$\kum(X)\to \kum(Y)$ induit
un morphisme $\kum'(X)\to \kum'(Y)$ qui, en vertu de ce qui précède et
de~\ref{souscourtyp}, est un isomorphisme ; il en résulte
que~$\kum(X)\simeq \kum(Y)$. 

Soit~$L$ une extension complète de~$k$. La flèche ~$\kum(X)\to \kum(X_L)$
induit un morphisme~$\kum'(X)\to \kum'(X_L)$, qui, en vertu de ce qui précède et de~\ref{zxextscal},
est un isomorphisme. 

\trois{kummpoint} Soit~$x$ un point de type 3 de $\skelan X$. Le singleton~$\{x\}$
est alors une sous-couronne de~$X$, et l'on dispose dès lors d'isomorphismes naturels
$$\kum(X)\simeq \kum(\{x\})\simeq \H^1(\hres(x),\mu_\ell),$$ le premier d'après
le~\ref{kummcourpremp} ci-dessus, et le second d'après le théorème 90 de Hilbert appliqué
au corps~$\hres(x)$. 

\deux{revkummcor} {\em Nature coronaire des revêtements de Kummer de~$X$.} On fixe une fonction coordonnée~$f$
sur~$X$. 

\trois{kummcorbase} L'espace~$k$-analytique~$\sch M(\sch O_X[T]/(T^\ell-f))$ est une~$k$-couronne dont~$T$ est une fonction
coordonnée : cela découle immédiatement
de~\ref{exemplekumm}. 

\trois{kummcordivi} Soit~$n$ un entier, et soit~$d$ le PGCD de~$n$ et~$\ell$ ; écrivons~$\ell=rd$ et~$n=sd$. On suppose
que~$T^d-1$ est scindé dans~$k$, et l'on note~$\mu_1,\ldots, \mu_d$ ses racines. Posons~$Y=\sch M(\sch O_X[T]/(T^\ell-f^n))$. 
On a 
$$T^\ell-f^n=T^{rd}-f^{sd}=\prod_i (T^r-\mu_i f^s)$$ (en tant que sections de~$\sch O_X[T]$). 
Si~$i$ et~$j$ sont deux entiers distincts compris entre~$1$ et~$d$,
la différence~$(T^r-\mu_i f^s)-(T^r-\mu_jf^s)=(\mu_i-\mu_j)f$ est une section inversible 
du faisceau d'algèbres~$\sch O_X[T]$, d'où par le lemme chinois
un isomorphisme naturel de~$\sch O_X$-algèbres
$$\sch O_X[T]/(T^\ell-f^n)\simeq \prod_i \sch O_X[T]/(T^r-\mu_if^s),$$ puis une décomposition 
$$Y\simeq \coprod_i \sch M(\sch O_X[T]/(T^r-\mu_if^s)).$$ 

Choisissons  deux entiers relatifs~$u$ et~$v$ tels que $ur+vs=1$ et fixons~$i$. On vérifie immédiatement que 

$$T\mapsto \mu_i^v\tau^s\;{\rm et}\;\tau \mapsto f^vT^u$$ définissent deux isomorphismes réciproques
l'un de l'autre entre les $\sch O_X$-algèbres~$\sch O_X[T]/(T^r-\mu_if^s)$ et $\sch O_X[\tau]/(\tau^r-\mu_i^uf).$ 

Dès lors,  $$\sch M(\sch O_X[T]/(T^r-\mu_if^s))\simeq \sch M(\sch O_X[\tau]/(\tau^r-\mu_i^uf)).$$ En vertu de~\ref{kummcorbase}, ce dernier est une~$k$-couronne
dont~$\tau$ est une fonction coordonnée.

\medskip
Il résulte de ce qui précède que toute composante connexe de~$Y$ est une~$k$-couronne dont $f^vT^u$ est une fonction coordonnée. 

\trois{kummcorpgen} Soit maintenant~$Y$ un~$\mu_\ell$-torseur de Kummer quelconque sur~$X$. D'après~\ref{kummcourpremp}, il existe~$n\in \NN$ et
~$a\in k\ti$ 
tels que~$Y\simeq \sch M(\sch O_X[T]/(T^\ell-af^n))$ ; fixons deux entiers~$u$ et~$v$ tels que~$u\ell+vn={\rm PGCD}(\ell,n).$ 

Soit~$Z$ une composante
connexe de~$Y$, et soit~$Z'$ une composante connexe de~$Z_{\KK}$. Soit~$\alpha$ une racine~$n$-ième
de~$a$ dans $\KK$. On a~$$Y_{\KK}\simeq  \sch M(\sch O_{X_{\KK}}[T]/(T^\ell-af^n))\simeq  \sch M(\sch O_{X_{\KK}}[T]/(T^\ell-(\alpha f)^n)).$$

Comme~$\alpha f$ est une fonction coordonnée de la~$\KK$-couronne~$X_{\KK}$, il résulte de~\ref{kummcordivi} que~$Z'$ est une~$\KK$-couronne
dont~$(\alpha f)^vT^u$ est une fonction coordonnée ; c'est évidemment aussi le cas de~$f^vT^u$. 

Soit~$\phi$ le morphisme d'espaces~$\got s(Z)$-analytiques de~$Z$ vers~$\Aff^{1,\rm an}_{\got s(Z)}$ induit par~$f^vT^u$. La composante~$Z'$ s'identifie à~$Z\hotimes_{\got s(Z)}\KK$
relativement à un certain plongement~$\got s(Z)\hookrightarrow \KK$. Le morphisme~$\phi\hotimes_{\got s(Z)}\KK$ n'est alors autre que le morphisme de~$Z'$ vers~$\Aff^{1,\rm an}_{\KK}$
fourni par la fonction par~$f^vT^u$. Par conséquent,~$\phi\hotimes_{\got s(Z)}\KK$ induit un isomorphisme~$Z'\simeq \DD_{I,\KK}$ pour un certain intervalle non vide~$I$ de~$\RR\ti_+$. 
Il s'ensuit que~$\phi(Z)\subset \DD_{I,\got s(Z)}$ ; par descente, le morphisme induit~$Z\to \DD_{I,\got s(Z)}$ est un isomorphisme (cor.~\ref{desciso}). 

\medskip
Ainsi,~$Z$ est une~$\got s(Z)$-couronne dont~$f^vT^u$ est une fonction coordonnée.

\subsection*{Disques et couronnes virtuels}

\deux{pseudisc} Soit~$X$ une courbe~$k$-analytique connexe et non vide. Nous la qualifierons de {\em disque virtuel sur~$k$} si~$X_{\KK}$ est un~$\KK$-disque, et si c'est le cas nous dirons qu'une extension complète~$L$ de~$k$ {\em déploie}~$X$ si~$X_L$ est un~$L$-disque. Si~$|k\ti|\neq\{1\}$ (resp. si~$|k\ti|=\{1\}$) nous qualifierons de {\em gentiment} virtuel tout disque virtuel sur~$k$ qui est déployé par une extension finie séparable (resp. finie) de~$k$. 

Nous dirons simplement que~$X$ est un disque virtuel (resp. gentiment virtuel), {\em sans mention du corps de base}, si c'est un disque virtuel {\em sur le corps~$\got s(X)$}.

\trois{psdlisse} Tout disque virtuel sur~$k$ est lisse, non vide et géométriquement connexe.

\trois{psdcbase} Si~$X$ est un disque virtuel sur~$k$ alors~$X_L$ est un disque virtuel sur~$L$ pour toute extension complète~$L$ de~$k$. 

\trois{constn} Si~$X$ est un disque virtuel sur~$k$, toute fonction inversible sur ~$X$
est de norme constante : pour le voir, on se ramène immédiatement au cas déployé par extension des scalaires, et l'assertion requise résulte alors de~\ref{casizero}. 

\deux{exempsd} Donnons maintenant quelques exemples de disques virtuels. 

\trois{discsv} Tout~$k$-disque est un disque virtuel sur~$k$. 

\trois{discmodram} Si~$X$ est un disque virtuel qui est déployé par une extension finie, séparable et modérément ramifiée de~$k$, alors~$X$ est un~$k$-disque. 

\trois{discnontriv}{\em Un exemple de disque virtuel non trivial.} Soit~$X$ l'ouvert de~$\Aff^{1,{\rm an}}_{\QQ_2}$ défini par l'inégalité~$|T^2-2|<|2|$. Si~$K$ est une extension complète de~$\QQ_2$ dans laquelle~$2$ est un carré alors~$X_K$ est un disque de~$\Aff^{1,{\rm an}}_K$, qui peut être décrit par l'inégalité~$|T-\sqrt 2|<|\sqrt 2|$. Il en résulte que~$X$ est un disque virtuel sur~$\QQ_2$ ;  mais~$X$ n'est pas un~$\QQ_2$-disque, puisqu'on vérifie immédiatement que~$X(\QQ_2)=\emptyset$. 

\trois{compshildisc} Soit~$a\in k$, soit~$r$ un réel strictement positif et soit~$Y$ une composante connexe de~$\pk\setminus\{\eta_{a,r}\}$. Comme ~$\eta_{a,r}$ est pluribranche, il existe une composante connexe~$Z$ de~$\pk\setminus\{\eta_{a,r}\}$ qui est distincte de~$Y$. Soit~$L$ une extension finie de~$k$ telle que~$Y(L)$ et~$Z(L)$ soient non vides (son existence est assurée par densité de~$k^a$ dans~$\KK$ si~$|k\ti|\neq\{1\}$, et par la description explicite de~$\pk$ sinon). Le point~$\eta_{a,r}$ ne possédant qu'un antécédent sur~$\PP^{1,{\rm an}}_L$, à savoir~$\eta_{a,r,L}$, l'espace~$Y_L$ (resp.~$Z_L$) s'écrit comme une réunion disjointe finie~$\coprod Y_i$ (resp.~$\coprod Z_j$) de composantes connexes de~$\PP^{1,{\rm an}}_L\setminus\{\eta_{a,r,L}\}$. Puisque~$Y$ et~$Z$ possèdent chacun un~$L$-point, il existe deux indices~$i_0$ et~$j_0$, et deux~$L$-points~$y$ et~$z$ de~$\PP^{1,{\rm an}}_L$ respectivement situés sur~$Y_{i_0}$ et~$Z_{j_0}$. Cela entraîne que~$Y_{i_0}$ est un~$L$-disque (\ref{compp1msd}). Il en résulte que~$Y$ est un disque virtuel, déployé par une extension finie de~$\got s(Y)$. 

\medskip
\deux{psdtop} Soit~$X$ un disque virtuel sur~$k$. Il s'identifie topologiquement au quotient du disque~$X_{\KK}$ par~$\mathsf G$ ; en conséquence,~$X$ est un arbre paracompact à un bout. Notons~$\omega$ (resp.~$\omega_{\KK}$) l'unique bout de~$X$ (resp.~$X_{\KK}$). 

\trois{intervinj} Si~$y\in X_{\KK}$ et si~$ x$ désigne son image sur~$X$, il découle de~\ref{unboutquotinj} que~$X_{\KK}\to X$ induit un homéomorphisme~$[y;\omega_{\KK}[\to [x;\omega[$. 

\trois{raynormarithm} Soit~$U$ un ouvert de~$X_{\KK}$ qui est une couronne et soit~$g\in \mathsf G$. L'ouvert~$g(U)$ s'identifiant à~$U\hotimes_g\KK$, c'est une couronne de même module que~$U$. Il s'ensuit que la fonction~$\rnorm$ sur~$X_{\KK}$ est invariante sous l'action de~$\mathsf G$ ; elle définit donc une application, que l'on note encore~$\rnorm$ et appelle encore le rayon normalisé, de~$X$ vers~$[0;1[$. En vertu du~\ref{intervinj} ci-dessus, si~$I$ est un intervalle tracé sur~$X$ et aboutissant à~$\omega$ alors la restriction de~$\rnorm$ à~$I$ est continue, strictement monotone, et tend vers~$1$ à l'approche de~$\omega$ ; si~$x\in X$, le rayon normalisé de~$x$ est nul si et seulement si~$x$ est de type 1. 

\deux{defpsc} Soit~$X$ une courbe~$k$-analytique connexe et non vide, et soit~$I$ un intervalle non vide de~$\RR\ti_+$. Nous qualifierons~$X$ de {\em couronne virtuelle sur~$k$} de type~$I$ si~$X_{\KK}$ est une~$\KK$-couronne de type~$I$ sur laquelle~$\mathsf G$ agit par automorphismes {\em directs} (\ref{dircourgen}) ; on dira qu'une extension complète~$L$ de~$k$ {\em déploie}~$X$ si~$X_L$ est une~$L$-couronne. Si~$|k\ti|\neq\{1\}$ (resp. si~$|k\ti|=\{1\}$) nous qualifierons de {\em gentiment} virtuelle toute couronne virtuelle sur~$k$ qui est déployée par une extension finie séparable (resp. finie) de~$k$.

\medskip
Nous dirons simplement que~$X$ est une couronne virtuelle (resp. gentiment virtuelle) {\em sans mention du corps de base} si c'est une couronne virtuelle {\em sur son corps des constantes~$\got s(X)$} ; nous appellerons couronne virtuelle {\em ouverte} toute couronne virtuelle de type ouvert. 

\trois{remtypcourvirt} Si~$X$ est une couronne virtuelle de type~$I$, alors~$I$ est bien déterminé à une transformation de la forme~$r\mapsto \lambda r^d$ près avec~$\lambda\in \sqrt{|k\ti|}$ et~$d\in \{-1,1\}$ ; lorsqu'on parlera {\em du} type de~$X$, il s'agira donc d'une classe d'intervalles. 

\trois{nottypes-virt} On appellera~$k$-couronne virtuelle de type ~$]*,*[$ (resp.~$]0,*[$, resp.~$[*,+\infty[$, etc.) toute~$k$-couronne virtuelle de type~$I$ pour un certain intervalle~$I$ de la forme~$]r;R[$  avec~$0<r<R$ (resp.~$]0;r[$ avec~$r>0$, resp.~$[r;+\infty[$ avec~$0<r$, etc.) ; notons que les notions de~$k$-couronne virtuelle de type~$]0;*[$ et~$]*,+\infty[$ coïncident. 

On dira qu'une couronne virtuelle est {\em ouverte} si son type est ouvert. 

\trois{psclisse} Si~$X$ est une couronne virtuelle sur~$k$ alors~$X$ est non vide, quasi-lisse et géométriquement connexe ; si~$X$ est ouverte alors~$X$ est sans bord. 

\trois{psc2} Soit~$X$ une courbe~$k$-analytique connexe et non vide, soit~$F$ le complété d'une extension galoisienne de~$k$ dont on note~$\mathsf H$ le groupe de Galois, et soit~$\mathsf G'$ le sous-groupe de~$\mathsf G$ formé des éléments qui stabilisent~$F$ point par point. Supposons que~$X_F$ soit une~$F$-couronne  (auquel cas~$X_{\KK}$ est {\em a fortiori} une~$\KK$-couronne) ; le groupe~$\mathsf H$ stabilise alors~$\skelan {X_F}$ (\ref{gautpask}). Comme (l'image réciproque de ) toute fonction fonction coordonnée de~$X_F$ est encore une fonction coordonnée de~$X_{\KK}$, le groupe~$\mathsf G'$ agit sur~$X_{\KK}$ par automorphismes directs, et~$\mathsf G$ agit sur~$X_{\KK}$ par automorphismes directs si et seulement si~$\mathsf H$ agit sur~$X_F$ par automorphismes directs ; autrement dit,~$X$ est une~$k$-couronne virtuelle si et seulement si~$\mathsf H$ agit sur~$X_F$ par automorphismes directs. 

\medskip
Si~$\skelan {X_F}$ n'est pas un singleton cette dernière condition revient à demander, d'après~\ref{dircourgen}, que~$\mathsf H$ agisse par homéomorphismes croissants sur~$\skelan {X_F}$ et même, en vertu de la compacité de~$\mathsf H$, qu'il fixe~$\skelan {X_F}$ point par point. 

Si~$\skelan {X_F}$ est ouvert (auquel cas il coïncide avec~$\skel{X_F}$),  cette dernière condition est satisfaite si et seulement si~$X$ est un arbre à deux bouts (\ref{exarbdeuxbouts} {\em et sq.}). 

\medskip
Remarquons que si~$\skelan {X_F}$ est un singleton, il est évidemment fixé point par point par~$\mathsf H$, que l'action soit ou non par automorphismes directs. 

\trois{psdsquel} Soit~$X$ une couronne virtuelle sur~$k$ et soit~$F$ une extension presque algébrique de~$k$. Identifions~$F$ à un sous-corps complet de~$\KK$, contenant~$k$, et notons~$\mathsf G'$ le sous-groupe de~$\mathsf G$  formé des éléments qui stabilisent~$F$ point par point. Les espaces~$X$ et~$X_F$ s'identifient respectivement à~$X_{\KK}/\mathsf G$ et~$X_{\KK}/\mathsf G'$. Comme~$X$ est une couronne virtuelle,~$\mathsf G$ agit par automorphismes directs sur~$X_{\KK}$ ; c'est {\em a fortiori} le cas de~$\mathsf G'$, et~$X_F$ est donc une~$F$-couronne virtuelle. 

\medskip
D'après le ~\ref{psc2}, le groupe~$\mathsf G$ agit trivialement sur~$\mathsf S(X_{\KK})$~\ref{psc2}), et son image sur~$X$ sera appelée le {\em squelette analytique} de~$X$ et sera notée~$\skelan X$ ; cette terminologie est compatible avec la précédente dans le cas où~$X$ est déjà déployée sur~$k$ (\ref{etainvgal}). Notons que~$\skelan {X_{\KK}}$ est l'image réciproque de~$\skelan X$ sur~$X_{\KK}$, et que~$\skelan{X_{\KK}}\to \skelan X$ est un homéomorphisme. 

\medskip
Comme~$\skelan X_{\KK}$ est un sous-graphe admissible de~$X_{\KK}$ stable sous~$\mathsf G$ et comme~$\skelan X$ est son image sur~$X$, le lemme~\ref{quoetflat} assure que~$\skelan X$ est un sous-graphe admissible de~$X$. 

\medskip
On déduit des résultats qui précèdent, appliqués d'une part directement, et d'autre part à la courbe~$F$-couronne virtuelle~$X_F$, que l'image réciproque de~$\skelan X$ sur~$X_F$ est égale à~$\skelan {X_F}$ et que~$\skelan {X_F}\to \skelan X$ est un homéomorphisme. Le lemme~\ref{lemmcommret} garantit alors que~$X_F\to X$ commute aux rétractions canoniques de~$X_F$ et~$X$ sur leurs squelettes analytiques respectifs.

\medskip
Si~$x\in \skelan X$, il a d'après ce qui précède un et un seul antécédent sur~$X_{\KK}$ ; le corps~$k$ est donc séparablement clos dans~$\kappa(x)$ et~$\hres(x)$.

\trois{psextqlc} Soit~$X$ une couronne virtuelle sur~$k$ et soit~$L$ une extension complète quelconque de~$k$. Soit~$F$ une extension complète composée de~$L$ et~$\KK$ ; la fermeture algébrique de~$L$ dans~$F$ est galoisienne sur~$L$, et son groupe de Galois~$\mathsf H$ agit sur~$F$. Comme~$X_{\KK}$ est une~$\KK$-couronne,~$X_F$ est une ~$F$-couronne, et~$X_F\to X_{\KK}$ induit un homéomorphisme entre~$\skelan {X_F}$ et ~$\skelan {X_{\KK}}$ (\ref{defdi}). Le groupe~$\mathsf H$ s'identifie naturellement à un sous-groupe de~$\mathsf G$ ; ce dernier agit par automorphismes directs sur~$X_{\KK}$ ; comme (l'image réciproque de) toute fonction fonction coordonnée de~$X_{\KK}$ est encore une fonction coordonnée de~$X_F$, le groupe~$\mathsf H$ agit par automorphismes directs sur~$X_F$ ; par conséquent,~$X_L$ est une~$L$-couronne virtuelle (\ref{psc2}). Chacune des applications~$X_F\to X_L, X_F\to X_{\KK}$ et~$X_{\KK}\to X$ induit un homéomorphisme entre le squelette analytique de sa source et celui de son but, et commute aux rétractions canoniques sur ces derniers (\ref{defdi} et~\ref{psdsquel}) ; il en résulte que~$X_L\to X$ induit un homéomorphisme entre ~$\skelan {X_L}$ et~$\skelan X$ et commute aux rétractions canoniques sur ces derniers. 

\trois{boutinfvirt} Soit~$X$ une~$k$-couronne virtuelle {\em ouverte}. En vertu de~\ref{exarbdeuxbouts} {\em et sq.}, il découle de l'égalité~$\skelan {X_{\KK}}=\skel {X_{\KK}}$ que ~$\skelan X=\skel X$. 

Supposons de plus que~$X$ est de type~$]0,*[$.  Si~$F$ est une extension complète de~$k$ déployant~$X$, le bout de~$X$ qui correspond au bout infini de~$X_F$ {\em via} l'homéomorphisme~$\skel {X_F}\simeq \skel X$ ne dépend pas de~$F$ : cela résulte de l'invariance du bout infini d'une couronne de type~$]0,*[$ par extension des scalaires. On l'appelle le {\em bout infini de~$X$.} Si~$L$ est une extension complète quelconque de~$k$, le bout de~$X$ qui correspond au bout infini de~$X_L$ {\em via} l'homéomorphisme~$\skel {X_L}\simeq \skel X$ n'est autre que le bout infini de~$X$. 

\trois{compcourvirt} Soit~$X$ une couronne virtuelle sur~$k$ et soit~$x\in \skelan X$. Soit~$y$ l'unique antécédent de~$x$ sur~$X_{\KK}$ ; on déduit de la description explicite de~$\skelan {X_{\KK}}$ (\ref{topcour}) que toutes les composantes connexes de~$X_{\KK}\setminus\{y\}$ sont des~$\KK$-disques, à l'exception de celles qui sont de la forme~$I^\flat$ pour~$I\in \pi_0(\skelan {X_{\KK}}\setminus\{y\}$ (il y en a au plus 2), qui sont des~$\KK$-couronnes ; on en déduit aussi que toutes les composantes connexes de~$X_{\KK}-\skelan{X_{\KK}}$ sont des disques. 

Il en résulte que toutes les composantes connexes de~$X\setminus\{x\}$ sont des disques virtuels, à l'exception de celles qui sont de la forme~$I^\flat$ pour~$I\in \pi_0(\skelan {X_{\KK}}\setminus\{y\}$ (il y en a au plus 2), qui sont des~$k$-couronnes virtuelles de~$X$ ; et que toutes les composantes connexes de~$X-\skelan X$ sont des disques virtuels. 

\medskip
Si~$\rho$ désigne la rétraction canonique de~$X$ sur~$\skelan X$, il résulte de ce qui précède et de~\ref{constn} que pour toute fonction inversible~$f$ sur~$X$ et tout~$x\in X$ on a~$|f(x)|=|f(\rho(x))|$. 

\trois{souspseudom} Soit~$X$ une couronne virtuelle sur~$k$ et soit~$I$ un intervalle non vide de~$\skelan X$ ; soit~$I_{\KK}$ l'image réciproque de~$I$ sur~$X_{\KK}$. L'image réciproque de~$I^\sharp$ sur~$X_{\KK}$ est la sous-couronne~$(I_{\KK})^\sharp$ de~$X_{\KK}$ ; on en déduit que~$I^\sharp$ est un domaine analytique de~$X$ (prop.~\ref{descgaldoman}), puis que c'est une~$k$-couronne virtuelle de squelette analytique égal à~$I$.

\trois{souspseu} Si~$X$ est une couronne virtuelle sur~$k$, on dira qu'un domaine analytique~$Y$ de~$X$ en est une {\em sous-couronne virtuelle} s'il est de la forme~$I^\sharp$ pour un certain un intervalle non vide~$I$ de~$\skelan X$. Dans ce cas,~$Y$ est une couronne virtuelle sur~$k$ de squelette analytique~$I$. Plus précisément,~$Y_L$ est une sous-couronne de~$X_L$ pour toute extension complète~$L$ de~$k$ déployant~$X$ : on le déduit de la commutation de~$X_L\to X$ aux rétractions canoniques sur les squelettes analytiques de ses source et but (\ref{psextqlc}).

\trois{modvir} Soit~$X$ une couronne virtuelle sur~$k$ et soit~$F$ une extension complète de~$k$ déployant~$X$ ; le module de la couronne~$X_F$ est visiblement indépendant du choix de~$F$ ; on l'appellera le module de~$X$.

\deux{discvirtcontcour} Soit~$X$ un disque virtuel sur~$k$, soit~$\omega$ son unique bout et soit~$\omega_{\KK}$ l'unique bout de~$X_{\KK}$. 

\trois{existifixpp} En vertu de~\ref{unboutglobinv}, il existe un intervalle ouvert tracé sur~$X_{\KK}$, aboutissant à~$\omega_{\KK}$, et invariant point par point sous l'action de~$\mathsf G$ ; son image sur~$X$ est un intervalle ouvert~$I$ qui aboutit à~$\omega$.

\trois{uniqueantec} Si~$x\in I$, il a par construction un unique antécédent~$x_{\KK}$ sur~$X_{\KK}$, ce qui signifie que~$\got s(X)=k$ ; on en déduit que tout ouvert de~$X$ rencontrant~$I$ est géométriquement connexe.

\trois{discvcourv} Soit~$x\in I$ et soit~$x_{\KK}$ son unique antécédent sur~$X_{\KK}$. Soit~$U$ (resp.~$V$) la composante connexe de~$X\setminus\{x\}$ contenant~$]x;\omega[$ (resp. une composante connexe relativement compacte de~$X\setminus\{x\}$).  Comme~$U$ est géométriquement connexe,~$U_{\KK}$ est une composante connexe de~$X_{\KK}\setminus\{x_{\KK}\}$ ; quant à~$V_{\KK}$, c'est une union finie et non vide de composantes connexes relativement compactes de~$X_{\KK}\setminus\{x_{\KK}\}$. 

\medskip
Identifions~$X_{\KK}$ à~$\DD_{[0;R[}$ pour un certain~$R>0$ {\em via} une fonction coordonnée. Le point~$x$ étant situé sur~$I$, il est pluribranche ; par conséquent,~$x_{\KK}$ est de type 2 ou 3, et donc de la forme~$\eta_{a,r}$ pour un certain élément~$a$ de~$\KK$ (de valeur absolue strictement inférieure à~$R$) et un certain~$r\in ]0;R[$. Il en découle que~$U_{\KK}$ est une couronne de type~$]*,*[$, et que toute composante connexe de~$V_{\KK}$ est un disque. 

\medskip
Il s'ensuit que~$x$ est de type 2 ou 3, que~$U$ est une couronne virtuelle de type~$]*,*[$ sur~$k$, et que~$V$ est un disque virtuel ; si~$V$ est de plus géométriquement connexe (c'est par exemple le cas si~$V$ contient~$I-[x;\omega[$) alors~$V$ est un disque virtuel sur~$k$.  

\trois{relcompdiscvirt} Si~$y\in X$ il existe un point~$x$ situé sur~$I\cap ]y;\omega[$ ; en vertu de ce qui précède, la composante connexe de~$y$ dans~$X\setminus\{x\}$ est un disque virtuel sur~$k$ qui est relativement compact dans~$X$.

\deux{bordcourt23} Soit~$X$ une courbe~$k$-analytique, soit~$U$ un ouvert de~$X$ qui est une couronne virtuelle ouverte, et soit~$x\in \partial U$. Soit~$\omega$ un bout de~$U$ convergeant vers~$x$. Choisissons une sous-couronne virtuelle ouverte stricte~$Z$ de~$U$ qui aboutit à~$\omega$ ; le bord de~$Z$ est de la forme~$\{x,y\}$ avec~$y\in U$. 

\medskip
Soit~$U'$ une composante connexe de~$U_{\KK}$. C'est un arbre à deux bouts ; l'ensemble~$\got d U'$ s'identifie naturellement à~$\got d U$. Soit~$\omega'$ le bout de~$U'$ correspondant à~$\omega$ ; c'est l'unique bout de~$U'$ auquel l'image réciproque~$Z'$ de~$Z$ sur~$U'$ aboutit. Comme~$Z'$ est un arbre à deux bouts, son bord dans~$X_{\KK}$ compte au plus deux points, et il se surjecte par ailleurs sur~$\partial Z$. Par conséquent,~$\partial Z'$ possède exactement un point~$x'$ au-dessus de~$x$, et~$\omega'$ converge vers~$x'$. 

Le point~$g(x')$ appartient au bord de~$g(U')$ pour tout~$g\in \mathsf G$, et~$\pi_0(U_{\KK})$ est fini ; l'orbite~$\mathsf G.x'$ est de ce fait finie. 

On fixe un isomorphisme~$U'\simeq \DD_{I,\KK}$ pour un certain intervalle ouvert~$I$ de~$\RR\ti_+$ (la couronne virtuelle~$U$ est donc de type~$I$), et l'on appelle~$r$ la borne de~$I$ qui correspond à~$\omega'$.  

\medskip
{\em Supposons que~$r\notin\{0,\infty\}$.} En vertu de~\ref{extremcour} {\em et sq.} le point~$x'$ est situé sur le lieu quasi-lisse de~$X_{\KK}$ et est de type 2 (resp. 3) si et seulement si~$r$ appartient (resp. n'appartient pas) à~$|(\KK)\ti|$. Par conséquent,~$x$ est situé sur le lieu quasi-lisse de~$X$ et est de type 2 (resp. 3) si et seulement si~$r$ appartient (resp. n'appartient pas) à~$\sqrt{|k\ti|}$. 

\medskip
{\em Supposons que~$r\in \{0,\infty\}$.} Il découle alors de {\em loc. cit.} que~$x'$ est un~$\KK$-point ; compte-tenu de la finitude de~$\mathsf G.x'$, ceci entraîne que~$x$ est un point {\em rigide} de~$X$. 

\deux{commadhcour} {\em Remarque.} Si les deux bouts de~$U$ convergent vers~$x$, il résulte de ce qui précède appliqué à chacun des deux bouts de~$U$ que les deux bornes de~$I$ appartiennent toutes deux à~$\sqrt{|k\ti|}$ si~$x$ est de type 2, appartiennent toutes deux à~$\RR\ti_+-\sqrt{|k\ti|}$ s'il est de type 3, et appartiennent toutes deux à~$\{0,\infty\}$ s'il est rigide. Par conséquent si~$U$ est de type~$]0;*[$, ou de type~$]r;R[$ avec~$r\in \sqrt{|k\ti|}$ et~$R\in  \RR\ti_+-\sqrt{|k\ti|}$, il existe un seul bout de~$U$ convergeant vers~$x$.  

\medskip
On a ainsi prouvé l'assertion suivante : si~$X$ est une courbe~$k$-analytique et si~$U$ est un ouvert de~$X$ qui est une couronne virtuelle de type~$]0;*[$ ou de type ~$]r;R[$ avec~$r\in \sqrt{|k\ti|}$ et~$R\in  \RR\ti_+-\sqrt{|k\ti|}$, alors~$\overline U$ est un arbre.

\deux{borddisct23} Soit~$X$ une courbe~$k$-analytique, et soit~$U$ un ouvert relativement compact de~$X$ qui est un disque virtuel. Comme~$U$ est un arbre à un bout,~$\partial U$ est un singleton~$\{x\}$.  Il existe en vertu de~\ref{discvirtcontcour} {\em et sq.} une couronne virtuelle~$V$ de type~$]*,*[$ contenue dans~$U$ et dont l'adhérence dans~$X$ contient~$x$. Il découle alors de~\ref{bordcourt23} que~$x$ est de type 2 ou 3 et que~$X$ est quasi-lisse en~$x$. 

\deux{morcourvirt} {\bf Proposition.} {\em Soient~$Y$ et~$X$ deux~$k$-courbes analytiques ; on suppose que chacune d'elles est une couronne virtuelle. Soit~$\phi : Y\to X$ un morphisme tel que~$\phi(Y)$ rencontre~$\skelan X$. 

\medskip
i) Il existe un intervalle non vide~$I$ de~$\skelan X$ tel que~$\phi$ induise un morphisme fini et plat~$Y\to I^\sharp$. 

ii) On a~$\phi\inv(I)=\skelan Y$, et~$\skelan Y\to I$ est un homéomorphisme.

iii) La flèche~$Y\to I^\sharp$ commute aux rétractions canoniques de ses source et but sur leurs squelettes analytiques respectifs~$\skelan Y$ et~$I$. 

iv) Si~$J$ est un intervalle non vide de~$\skelan Y$ alors~$J^\sharp=\phi\inv(\phi(J)^\sharp)$ ; si~$J_0$ désigne le type de~$J^\sharp$ alors le type de~$\phi(J)^\sharp$ est égal à~$J_0^{\frac {\deg\;(Y\to I^\sharp)}{[\got s(Y):\got s(X)]}}$, et l'on a en particulier~$$\mathsf{Mod}(\phi(J)^\sharp)=\mathsf{Mod}(J^\sharp)^{\frac {\deg\;(Y\to I^\sharp)}{[\got s(Y):\got s(X)]}}.$$}

\medskip
{\em Démonstration.} On va tout d'abord se ramener au cas où~$X$ et~$Y$ sont géométriquement connexes. 

\trois{xsxsy} Le morphisme~$\phi$ admet une factorisation canonique~$$Y\to X\otimes_{\got s(X)}\got s(Y)\to X.$$ La flèche~$\psi : X\otimes_{\got s(X)}\got s(Y)\to X$ est finie et plate de degré~$[\got s(Y):\got s(X)]$, le fermé~$\psi\inv(\skelan X)$ de~$X\otimes_{\got s(X)}\got s(Y)$ est égal à~$\skelan {X\otimes_{\got s(X)}\got s(Y)}$, et~$$\skelan {X\otimes_{\got s(X)}\got s(Y)}\to \skelan X$$ est un homéomorphisme ; le morphisme~$\psi$ commute aux rétractions canoniques de ses source et but sur leurs squelettes analytiques respectifs, et si~$J$ est un intervalle non vide tracé sur~$\skelan {X\otimes_{\got s(X)}\got s(Y)}$ alors les couronnes virtuelles~$J^\sharp$ et~$\psi(J)^\sharp$ ont même type.

On peut donc remplacer~$X$ par ~$X\otimes_{\got s(X)}\got s(Y)$, puis~$\phi$ par~$Y\to X\otimes_{\got s(X)}\got s(Y)$ et enfin~$k$ par~$\got s(Y)$, et par là se ramener au cas où~$Y$ et~$X$ sont géométriquement connexes.

\trois{deuxcourgeoconn}{\em Preuve lorsque~$Y$ et~$X$ sont géométriquement connexes}. Chacune des~$\KK$-courbes~$Y_{\KK}$ et~$X_{\KK}$ est alors une couronne ; si~$I$ est un intervalle tracé sur~$\skelan X$ (resp.~$\skelan Y$), on notera~$I_{\KK}$ son image réciproque sur~$\skelan {Y_{\KK}}$ (resp. ~$\skelan {X_{\KK}}$) ; on écrira simplement~$I_{\KK}^\sharp$ au lieu de~$(I^\sharp)_{\KK}\simeq (I_{\KK})^\sharp$. Par hypothèse,~$\phi(Y)$ rencontre~$\skelan X$ ; en conséquence,~$\phi_{\KK}(Y_{\KK})$ rencontre~$\skelan {X_{\KK}}$. Il découle alors de~\ref{morcourintersk} : qu'il existe un intervalle non vide~$I$ tracé sur~$\skelan X$ tel que~$\phi_{\KK}$ induise un morphisme fini et plat~$Y_{\KK}\to I_{\KK}^\sharp$ ; que~$\skelan Y_{\KK}$ est égal à l'image réciproque de~$I_{\KK}$ sur~$Y_{\KK}$ ; et que l'application continue~$\skelan Y_{\KK}\to I_{\KK}$ est un homéomorphisme. 

 Dès lors,~$\phi$ induit un morphisme de~$Y$ vers le domaine analytique~$I^\sharp$ de~$X$ ; et comme~$Y_{\KK}\to I_{\KK}^\sharp$ est finie et plate, la flèche~$Y\to I^\sharp$ est finie et plate (\ref{descfinplat} {\em et sq.} -- les remarques du~\ref{remdescfiniplat} s'appliquent ici), de degré égal à celui de~$Y_{\KK}\to I_{\KK}^\sharp$, d'où i). On déduit également de ce qui précède que~$\skelan Y$ est égal à l'image réciproque de~$I$ sur~$Y$, et que~$\skelan Y\to I$ est un homéomorphisme, d'où ii).  Il résulte alors aussitôt du lemme~\ref{lemmcommret} que~$\phi$ commute aux rétractions canoniques de ses source et but sur leurs squelettes analytiques respectifs, d'où iii). 
 
\medskip
Soit~$J$ un intervalle non vide de~$\skelan Y$ ; l'égalité~$J^\sharp=\phi\inv(\phi(J)^\sharp)$ est une conséquence immédiate de iii). Le type~$J_0$ de~$J^\sharp$ est égal à celui de~$J_{\KK}^\sharp$, et le type de~$\phi(J)^\sharp$ est égal à celui de~$\phi(J)_{\KK}^\sharp=\phi(J_{\KK})^\sharp$ ; mais ce dernier type est, en vertu de~\ref{effmodule}, égal à~$J_0^{\deg\;(Y_{\KK}\to I_{\KK}^\sharp)}$, soit encore à~$J_0^{\deg \;(Y\to I^\sharp)}.$~$\Box$ 

\deux{remmorcourfiplat} {\em Remarque.} Soit~$\phi :Y\to X$ un morphisme entre couronnes virtuelles. Si~$\phi$ est fini et plat alors~$\phi(Y)=X$, et la proposition ci-dessus s'applique donc, avec~$I=\skelan X$ tout entier. 

\deux{torskummpseudocour} Soit~$X$ une~$k$-couronne virtuelle et soit~$\ell$ un entier premier à~$p$. Soit~$Y$ un~$\mu_\ell$-torseur de Kummer sur~$X$, et soit~$Z$ une composante connexe de~$Y$. C'est alors une couronne virtuelle, et gentiment virtuelle si~$X$ est gentiment virtuelle. En effet, soit~$L$ une extension complète de~$k$ déployant~$X$. Si~$Z'$ est une composante connexe de~$Z_L$, alors~$Z'$ est une composante connexe de~$Y_L$, et est dès lors une~$\got s(Z')$-couronne en vertu de~\ref{kummcorpgen}, d'où notre assertion. 

En particulier, le morphisme~$Z\to X$ est  justiciable de la remarque~\ref{remmorcourfiplat} ci-dessus. Il s'ensuit notamment, compte-tenu de l'assertion ii) de la proposition~\ref{morcourvirt}, que tout point de~$\skelan X$ a un et un seul antécédent sur~$Z$. 
 
\subsection*{Existence de toises sur les sous-graphes compacts d'une courbe analytique}

\deux{toisedv} {\em La toise canonique sur la compactification arboricole d'un disque virtuel.} Soit~$X$ un disque virtuel sur~$k$, soit~$\wid X$ sa compactification arboricole (sans autre structure que celle d'un espace topologique) et soit~$\omega$ l'unique bout de~$X$. Soit~$\delta$ l'application de~$\widehat X$ dans~$\RR_+$ qui envoie~$\omega$ sur~$0$ et~$x$ sur~$1-\rnorm (x)$ si~$x\in X$ ; c'est une toise partielle sur~$\wid X$, basée en~$\omega$ et bornée par~$1$. La toise sur~$\wid X$ induite par~$\delta$ est invariante par tout automorphisme de~$X$ ; elle sera appelée la {\em toise canonique} sur~$\wid X$ ; elle est bornée par~$2$.

\medskip
Si~$Y$ est une courbe~$k$-analytique et si~$U$ est un ouvert de~$Y$ qui est un disque virtuel relativement compact,~$\overline U$ s'identifie naturellement à~$\wid U$ ; par conséquent, on pourra parler de la toise canonique sur~$\overline U$. 

\deux{toisep1} {\em La toise standard sur la droite projective.} Si~$U$ est une composante connexe de~$\pk\setminus\{\eta_1\}$, c'est un disque virtuel (\ref{compshildisc}) ; on notera~$l_U$ la toise canonique sur~$\overline U=U\cup\{\eta_1\}$. 

\medskip
Le sous-arbre~$\{\eta_1\}$ de~$\pk$ est admissible. La toise définie par concaténation des~$l_U$ et de la toise triviale sur~$\{\eta_1\}$ sera appelée la {\em toise standard} de~$\pk$ ; elle est bornée par 2. 

\deux{theotoise}{\bf Théorème.} {\em Soit~$X$ une courbe~$k$-analytique et soit~$Y$ un sous-graphe compact de~$X$. Le graphe~$Y$ admet une toise.}

\medskip
{\em Démonstration.} On commence par plusieurs réductions. 

\trois{xegalycomp}{\em On peut supposer que~$X$ est affinoïde et que~$Y=X$.} En effet, il existe par compacité de~$Y$ une famille finie~$(X_i)$ de domaines affinoïdes de~$X$ qui recouvrent~$Y$. Il suffit de montrer que~$\bigcup X_i$ admet une toise ; on peut pour ce faire, en vertu de~\ref{paramfini}, se contenter de montrer que chacun des~$X_i$ admet une toise ; on s'est donc bien ramené à la situation annoncée. 

\trois{suppxnorm} {\em On peut supposer que~$X$ est un domaine affinoïde de l'analytifiée  d'une~$k$-courbe algébrique projective, irréductible et lisse~$\sch X$.} En effet, le problème étant purement topologique, on peut étendre les scalaires au complété de la clôture parfaite de~$k$, et donc supposer~$k$ parfait. Grâce à la proposition~\ref{transfr} on peut remplacer~$X$ par sa normalisation et donc le supposer quasi-lisse. Il admet alors un recouvrement fini~$(X_i)$ où chaque~$X_i$ est un domaine affinoïde de~$X$ s'immergeant dans l'analytification d'une~$k$-courbe algébrique irréductible et lisse, que l'on peut compactifier en une courbe lisse puisque~$k$ est parfait. D'après le~\ref{paramfini}, il suffit de s'assurer que chacun des~$X_i$ admet une toise ; on s'est donc bien, là encore, ramené à la situation annoncée.

\trois{toisecourbes} Comme~$\sch X$ est lisse, il existe un morphisme fini et génériquement étale~$\sch X\to \PP^1_k$ ; soit~$\sch Y$ la normalisée de~$\sch X$ dans une clôture galoisienne de l'extension~$\kappa(\sch X)/k(T)$, soit~$\mathsf G$ le groupe de Galois de~$\kappa(Y)$ sur~$k(T)$, et soit~$\mathsf H$ le sous-groupe de~$\mathsf G$ correspondant à~$\kappa(\sch X)$. On dispose d'identifications topologiques naturelles~$$\sch X\an \simeq \sch Y\an/\mathsf H\;{\rm et}\;\pk\simeq \sch Y\an/\mathsf G.$$ 

\medskip
On a vu au~\ref{toisep1} que~$\pk$ admet une toise. On déduit alors du théorème~\ref{transflong} que~$\sch Y\an$ admet une toise, puis que~$\sch X\an$ en admet une. Par restriction,~$X$ admet une toise, ce qui achève la démonstration.~$\Box$ 

\section{Algébrisation} 

\subsection*{Algébrisation des courbes analytiques propres}

\deux{pascercle} {\bf Proposition.} {\em Soit~$X$ une courbe~$k$-analytique propre. L'ensemble des points rigides de~$X$ en est une partie dense pour la topologie de Zariski}.

\medskip
{\em Démonstration.} Si~$|k\ti|\neq\{1\}$, il résulte du fait que~$X$ est sans bord que l'ensemble de ses points rigides en est une partie dense pour la topologie usuelle, et {\em a fortiori} pour la topologie de Zariski. On peut donc supposer que~$|k\ti|=\{1\}$. 

\medskip
Soit~$\sch U$ un ouvert de Zariski non vide de~$X$ ; nous allons montrer qu'il contient un point rigide.  Comme~$\KK$ est une clôture algébrique de~$k$, il suffit de montrer que~${\sch U}_{\KK}$ contient un point rigide, ce qui permet de se ramener au cas où~$k$ est algébriquement clos. On peut, quitte à remplacer~$X$ par sa normalisée, supposer que~$X$ est lisse, puis qu'elle est irréductible.  Dans ce cas~$X\setminus \sch U$ est un ensemble fini de~$k$-points ; il suffit donc, pour conclure, de démontrer que~$X$ contient une infinité de~$k$-points. 

\trois{aumoinst21} {\em L'ensemble~$X\typ 2$ est non vide}. Supposons que ce ne soit pas le cas ; on a alors~$X=X\typ {1,3}$.  Si~$x\in X\typ 3$, le groupe~$|\hres(x)\ti|$ est libre de rang~$1$, et possède donc un unique générateur strictement plus petit que~$1$ ; notons-le~$r(x)$. Le point~$x$ possède un voisinage ouvert dans~$X$ qui est isomorphe à un ouvert de l'analytifiée~$\sch Y\an$ d'une~$k$-courbe projective,  lisse et irréductible~$\sch Y$ ; soit~$y$ le point de~$\sch Y\an$ correspondant à~$x$ {\em via} l'isomorphisme évoqué. 

Pour tout point fermé~$P$ de~$\sch Y$ et tout~$r \in ]0;1[$ notons~$\xi_{P,r}$ le point de~$\sch Y\an$ correspondant à la valeur absolue~$r^{v_P}$ de~$\kappa(\sch Y)$, où~$v_P$ est la valuation discrète associée à~$P$ ; notons que~$|\hres(\xi_{P,r})\ti|=r^{\ZZ}$. Comme~$y$ est de type 3 il est de la forme~$\xi_{P,r}$ avec~$r$ et~$P$ comme ci-dessus ; l'ensemble~$\{\xi_{P,s}\}_{0<s<1}$ est alors un voisinage ouvert de~$y$ dans~$\sch Y\an$. 

Il s'ensuit que~$x$ possède un voisinage dans~$X$ qui est contenu dans~$X\typ 3$, homéomorphe à un intervalle ouvert, et sur lequel l'application~$r$ est continue et strictement monotone. 

\medskip
Si~$x$ est un~$k$-point de~$X$, il possède un voisinage ouvert~$V$ qui est un disque ; quitte à le restreindre, on peut supposer celui-ci isomorphe à l'ouvert de~$\Aff^{1,{\rm an}}_k$ défini par l'inégalité~$|T|<s$ pour un certain~$s\in ]0;1[$ ; cet ouvert peut simplement se décrire comme la réunion de l'origine et de~$\{\eta_t\}_{0<t<s}$ ; il s'ensuit que tout point de~$V\setminus\{x\}$ est de type 3, et que~$r(z)\to 0$ quand le point~$z$ de~$V\setminus\{x\}$ tend vers~$x$. 

\medskip
On déduit de ce qui précède que~$r$ se prolonge en une application continue de~$X$ dans~$\RR_+$ telle que~$r(x)=0$ pour tout~$k$-point~$x$ de~$X$. Par compacité,~$r$ atteint son maximum en un point~$x$ de~$X$. Comme~$X$ n'est pas constitué uniquement de points rigides (c'est un espace de dimension~$1
$) et comme~$r$ ne s'annule par construction qu'en les points rigides de~$X$, le point~$x$ est de type 3. Mais on a vu plus haut qu'il possède alors un voisinage homéomorphe à un intervalle ouvert sur lequel~$r$ est strictement monotone, ce qui contredit le fait qu'elle atteint son maximum en~$x$. 

\trois{concluinfrig} {\em Conclusion}. On peut donc choisir un point~$x$ de type 2 sur~$X$. Par lissité, il existe un voisinage ouvert~$V$ de~$x$ dans~$X$ et un morphisme fini, étale et surjectif de~$\phi$ de~$V$ vers un ouvert~$U$ de~$\Aff^{1,{\rm an}}_k$. Comme~$\phi(x)$ est de type 2, c'est nécessairement le point~$\eta_1$ de~$\Aff^{1,{\rm an}}_k$ ; en tant que voisinage de~$\eta_1$, l'ouvert~$U$ contient une infinité de~$k$-points ; comme~$\phi$ est fini et surjectif,~$V$ contient une infinité de~$k$-points, ce qui achève la démonstration.~$\Box$

\deux{alganpropo} {\bf Théorème.} {\em Soit~$X$ une~$k$-courbe propre. Elle est isomorphe à l'analytification d'une~$k$-courbe projective.}

\medskip
{\em Démonstration.} Commençons par rappeler que comme~$X$ est propre, les espaces vectoriels de cohomologie d'un faisceau cohérent sur~$X$ sont tous de dimension finie ; ce fait sera implicitement utilisé dans la suite. 

\medskip
L'ensemble des points de~$X$ dont l'anneau local est de Cohen-Macaulay est un ouvert de Zariski de~$X$ contenant tous les points non rigides de~$X$ (l'anneau local en un tel point est artinien), et rencontrant par conséquent chaque composante irréductible de~$X$. Cela entraîne, en vertu de la proposition~\ref{pascercle} ci-dessus, qu'il existe un ensemble fini~$D$ de points rigides de~$X$, rencontrant toutes les composantes irréductibles de~$X$ et constitué de points dont l'anneau local est de Cohen-Macaulay. 

\medskip
Chaque point  de~$D$ peut alors être localement défini, {\em ensemblistement}, comme le lieu d'annulation d'une fonction~$f$ qui n'est pas diviseur de zéro ; il existe donc un diviseur de Cartier effectif~$\mathsf D$ de~$X$ dont  le support est égal à~$D$ ; soit~$\sch L$ le fibré en droites correspondant à~$\mathsf D$. 

\medskip
\trois{lsecglob}{\em Pour~$n$ assez grand le fibré~${\sch L}^{\otimes n}$ est engendré par ses sections globales.} Comme~${\sch L}_{|X-\mathsf D}$ est trivial, il suffit de vérifier cette assertion en chacun des points de~$\mathsf D$. Pour tout entier~$n$, la tensorisation par~${\sch L}^{\otimes n}$ de la suite exacte~$$0\to {\sch O}_X(-\mathsf D)\to {\sch O}_X\to {\sch O}_{\mathsf D}\to 0$$ fournit une suite exacte~$$0\to {\sch L}^{\otimes n-1}\to {\sch L}^{\otimes n}\to {\sch L}^{\otimes n}\otimes {\sch O}_{\mathsf D}\to 0,$$ d'où une suite exacte de cohomologie~$$\H^0(X,{\sch L}^{\otimes n})\to  \H^0(X,{\sch L}^{\otimes n}\otimes{\sch O}_{\mathsf D})\to \H^1(X,{\sch L}^{\otimes n-1})\to \H^1(X,{\sch L}^{\otimes n})\to 0,$$ le faisceau~${\sch L}^{\otimes n}\otimes {\sch O}_{\mathsf D}$ étant à support dans~$\mathsf D$ et donc à cohomologie triviale en rang~$>0$. 

\medskip
Par ce qui précède, la suite~$(\dim k \H^1(X,{\sch L}^{\otimes n}))_n$ est décroissante, et donc stationnaire à partir d'un certain rang. Par conséquent,~$$\H^1(X,{\sch L}^{\otimes n-1})\to \H^1(X,{\sch L}^{\otimes n})$$ est bijective pour~$n$ assez grand ; il s'ensuit que~$\H^0(X,{\sch L}^{\otimes n})\to  \H^0(X,{\sch L}^{\otimes n}\otimes{\sch O}_{\mathsf D})$ est surjective pour~$n$ assez grand, d'où notre assertion.

\trois{fabricproj} Quitte à multiplier~$\mathsf D$ (et~$\sch L$) par un entier convenable, on peut donc supposer que~$\sch L$ est engendré par ses sections globales. Soit~$\mathsf A$ la sous-$k$-algèbre graduée de~$\bigoplus\limits_n\H^0(X,{\sch L}^{\otimes n})$ engendrée par~$\H^0(X,{\sch L})$ ; c'est une~$k$-algèbre graduée de type fini ; comme~$\sch L$ est engendré par ses sections globales, on dispose d'un morphisme naturel d'espaces annelés de~$X$ vers~${\sch X}:={\rm Proj}\;\mathsf A$ ; il résulte de la propriété universelle de l'analytifiée que ce morphisme se factorise par une flèche~$\phi : X\to {\sch X}\an$. 

\medskip
Comme~$\mathsf D$ est effectif, la fonction~$1$ définit une section globale de~$\sch L$ qui s'annule exactement en les points de~$D$, lequel rencontre chacune des composantes irréductibles de~$X$ ; soit~$a$ l'élément de degré~$1$ de~$\mathsf A$ qui correspond à cette section. 

\medskip
Soit~$Y$ une composante irréductible de~$X$ et soit~$y$ l'un des points de~$D$ situé sur~$Y$. L'image~$x$ de~$y$ sur~${\sch X}\an$ est située sur le lieu des zéros~$Z$ de~$a$ (vue comme section d'un fibré ample sur ~${\sch X}\an$) ; si~$y'$ est un point de~$Y$ n'appartenant pas à~$D$ alors~$\phi(y)\notin Z$ ; en particulier,~$\phi(y)\neq x$ ; il s'ensuit que~$Y$ n'est pas contenue dans~$\phi\inv(x)$. 

\medskip
Si~$t$ est un point rigide de~$X$, la fibre~$\phi\inv(\phi(t))$ ne contient par ce qui précède aucune composante irréductible de~$X$ ; par conséquent,~$\phi\inv(\phi(t))$ est de dimension nulle. L'ensemble des points de~$X$ en lesquels la dimension relative de~$X$ est nulle en est un ouvert de Zariski ; comme il contient tous les points rigides de~$X$, c'est~$X$ elle-même en vertu de la proposition~\ref{pascercle}. 

\medskip
Comme~$X$ est propre,~$\phi$ est compacte ; comme~$\phi$ est de dimension relative~$0$ en tout point de~$X$, et est par ailleurs sans bord puisque~$X$ est sans bord,~$\phi$ est fini en tout point de~$X$ ; étant propre et à fibres finies,~$\phi$ est fini. Il découle alors du caractère projectif de~${\sch X}$ et de GAGA que~$X$ est algébrisable.~$\Box$

\deux{remnorm} {\em Remarque.} Si~$X$ est normale, la preuve du théorème précédent peut être simplifiée : on choisit un point rigide~$P$ sur~$X$ (grâce à la proposition~\ref{pascercle}) ; c'est un diviseur de Cartier, et le même argument que ci-dessus montre que~${\sch O}(nP)$ est engendré par ses sections globales pour~$n$ assez grand.

En particulier, il existe~$n$ tel que~${\sch O}(nP)$ admette une section globale ne s'annulant pas en~$P$ ; celle-ci peut s'interpréter comme une fonction méromorphe non constante, et fournit ainsi un morphisme~$X\to \pk$ ; on procède alors encore une fois comme ci-dessus : on montre que ce morphisme est fini, et on conclut par GAGA. 

\deux{remxalgcan} {\bf Remarque.} En vertu de GAGA l'analytification des variétés propres est un foncteur pleinement fidèle. On peut donc reformuler le théorème~\ref{alganpropo} en disant que~$\sch X\mapsto \sch X\an$ établit une équivalence entre la catégorie des~$k$-courbes algébriques propres et celle des courbes~$k$-analytiques propres. 

\subsection*{Algébrisation des courbes formelles propres}

\deux{defcourbesform} On appellera {\em~$k\zero$-courbe formelle} (resp. {\em~$k\zero$-courbe algébrique}) tout~$k\zero$-schéma formel (resp.~$k\zero$-schéma) plat, de type fini et purement de dimension (resp. de dimension relative) 1. 

\deux{algebrformcurv} {\bf Théorème.} {\em Soit~$\got X$ une~$k\zero$-courbe formelle propre. 

\medskip
i) Il existe une~$k\zero$-courbe algébrique projective~$\bnd X$ telle que~$\got X$ s'identifie à la complétion de~$\bnd  X$ le long de sa fibre spéciale. 

ii) Supposons que~$|k\ti|$ est libre de rang 1, que~$\got X_{\red k}:=\got X\otimes_{k\zero}\red k$ et~$\got X_\eta$ sont réduits, que les composantes irréductibles de~$
\got X_{\red k}$ sont génériquement lisses, et qu'il existe un sous-corps dense~$F$ de~$k$ et une~$F$-courbe algébrique~$\sch Y$ telle que~$$\got X_\eta\simeq( {\sch Y}\times_Fk)\an\;;$$ il existe alors une~$F\zero$ -courbe projective~$\bnd Y$ telle que~$\bnd Y_\eta\simeq \sch Y$ et tel que~$\got X$ s'identifie à la complétion de~$\bnd Y\otimes_{F\zero}k\zero$ le long de sa fibre spéciale.}

\medskip
{\em Démonstration.} Le théorème est évident si~$|k\ti|=\{1\}$ ; on suppose à partir de maintenant que ce n'est pas le cas, et l'on fixe un élément~$a$ de~$k$ tel que~$0<|a|<1$. 

\trois{algebgen} {\em Preuve de i).} Choisissons un fibré ample~$\sch L$ sur la~$\red k$-courbe projective~$\got X_{\red k}$. L'obstruction à relever~$\sch L$ en un fibré en droites sur~$\got X$ vit dans~$\H^2(\got X_{\red k},k\zeroo{\sch O}_{\got X})$ ; l'espace topologique~$\got X_{\red k}$ étant de dimension cohomologique~$1$, cette obstruction s'annule et~$\sch L$ se relève donc en un fibré en droites~$\got L$ sur~$\got X$. 

\medskip
Soit~$n$ un entier. Pour tout idéal de définition~$\bnd I$ de~$k\zero$, le~$k\zero/\bnd I$-module~$\H^0(\got X/\bnd I, \got L^{\otimes n}/\bnd I)$ est plat et de présentation finie, et par conséquent libre puisque~$k\zero/\bnd I$ est local ; comme le~${\sch O}_{\got X/\bnd I }$-module~$\got L^{\otimes n}/\bnd I$ est~$k\zero/\bnd I$-plat, la formation de ses sections globales commute aux changements de base. 

On déduit de ce qui précède que~$$\H^0(\got X/\bnd J, \got L^{\otimes n}/\bnd J)\simeq \H^0(\got X/\bnd I,\got L^{\otimes n}/\bnd I)/\bnd J$$ pour tout couple~$(\bnd I, \bnd J)$ d'idéaux de définition de~$k\zero$ tel que~$\bnd I\subset \bnd J$ ; il s'ensuit que~$\H^0(\got X,\got L^{\otimes n})$ est un~$k\zero$-module libre de rang fini. 

\medskip
Comme~$\sch L$ est ample, la~$\red k$-algèbre~$\bigoplus\limits_n \H^0(\got X_{\red k},\sch L^{\otimes n})$ est engendrée par ses~$N$ premiers sommandes pour un certain entier~$N$, ce qui revient à dire que pour tout entier~$n$, l'application naturelle~$$\bigoplus_{n_1,\ldots,n_r\leq N,\;\sum n_i=n}  \left(\bigotimes_i \H^0(\got X_{\red k},\sch  L^{\otimes n_i})\right)\to  \H^0(\got X_{\red k},\sch  L^{\otimes n})$$ est surjective. En vertu du lemme de Nakayama, il s'ensuit que l'application~$$\bigoplus_{n_1,\ldots,n_r\leq N,\;\sum n_i=n}  \left(\bigotimes_i \H^0(\got X,\got L^{\otimes n_i})\right)\to  \H^0(\got X,\got L^{\otimes n})$$ est surjective pour tout~$n$, ce qui siginifie que la~$k\zero$-algèbre~$\bigoplus \H^0(\got X,\got L^{\otimes n})$ est de type fini ; comme elle est somme directe de modules libres, elle est sans torsion et donc plate ; par conséquent, elle est de présentation finie. 

Posons~$\bnd X={\rm Proj}\; \bigoplus\limits_n \H^0(\got X,\got L^{\otimes n})$. Comme~$\sch L$ est engendré par ses sections globales, il existe un morphisme naturel de~$\got X$ vers le complété~$\wid {\bnd X}$ de~$\sch X$ le long de sa fibre spéciale. 

Pour tout idéal de définition~$\sch I$ de de~$k\zero$, le fibré en droites~$\got L/\sch I$ sur le~$k\zero/\sch I$-schéma~$\got X/\sch I$ est relativement ample ; par conséquent, ~$\got X/{\sch I}\to \wid {\sch X}/{\sch I}$ est un isomorphisme ; ceci valant pour tout~$\sch I$, le morphisme~$\got X\to \wid {\bnd X}$ est un isomorphisme, ce qui achève de prouver i). 

\trois{algebvaldisc} {\em Preuve de ii).} Plaçons-nous maintenant sous les hypothèses de ii). Choisissons sur chacune des composantes irréductibles de~$\got X_{\red k}$ un point fermé en lequel~$\got X$ est lisse (c'est possible par hypothèse) ; notons~$(P_i)_i$ la collection finie de points fermés ainsi obtenus. 

Pour tout~$i$, l'image réciproque~$U_i$ de~$P_i~$ sur~$\got X_\eta$ est un ouvert non vide. Il existe un morphisme fini et plat de~$\bnd Y$ vers~$\PP^1_F$ ; la flèche~$$\got X_\eta\simeq(\bnd Y_ k)\an \to \pk$$ est ouverte ; il s'ensuit, compte-tenu de la densité de~$k^s$ dans~$\KK$, de la densité de~$F$ dans~$k$, et du lemme de Krasner, que l'image de~$U_i$ sur~$\pk$ contient un point rigide provenant d'un point fermé de~$\PP^1_F$ ; ce dernier possède par construction un antécédent~$Q_i$ sur~$\bnd Y$ dont l'image réciproque~$Q'_i$ sur~$\sch Y_ k$ est, modulo l'identification entre ~$(\sch  Y_k)\an$ et~$\got X_\eta$, contenue dans~$U_i$. Comme~$P_i$ est un point lisse de~$\got X_s$, il est contenu dans le lieu lisse du schéma formel~$\got X$ ; par conséquent, le fermé~$Q'_i$ est contenu dans le lieu lisse de~$\got X_\eta\simeq( \sch  Y_k)\an$, et~$Q_i$ est dès lors contenu dans le lieu lisse de~$\sch Y$. 

\medskip
Soit~$\sch L$ le fibré en droites sur~$\sch Y$ correspondant au diviseur de Cartier~$\sum Q_i$ ; l'image réciproque~${\sch L}_k$ de~$\sch L$ sur~$\sch Y_ k$ est le fibré en droites associé au diviseur de Cartier~$\sum Q'_i$. Si~$P$ est un point quelconque du support de~$\sum Q'_i$, il induit, {\em via} l'identification~$(\sch Y_k)\an\simeq \got X_\eta$, un morphisme fini~${\rm Spf}\;\hres(P)\zero\to \got X$, dont l'image est un fermé de Zariski de~$\got X$ de support ~$\{P_i\}$ pour un certain~$i$ ; c'est donc un fermé de codimension 1 du lieu lisse de~$\got X$, et {\em a fortiori} de son lieu régulier ; c'est de ce fait un diviseur de Cartier sur~$\got X$. Il s'ensuit que~${\sch L}_k\an$, vu comme fibré en droites sur~$\got X_\eta$, s'identifie à la fibre générique d'un fibré en droites~$\got L$ sur~$\got X$, qui est tel que~$\got L_{\red k}$ soit associé à un diviseur de Cartier de la forme~$\sum n_iP_i$ où les~$n_i$ sont des entiers strictements positifs ; en particulier,~$\got L_{\red k}$ est ample ; il s'ensuit que~$\got L_\eta$ et~$\sch L$ sont amples. 

\medskip
Choisissons un recouvrement formel affine~$(\got  U_j)$ de~$\got X$ qui trivialise~$\got L$. Pour tout~$j$, l'espace analytique~$(\got U_j)_\eta$ est un domaine affinoïde de~$\got X_\eta$. Comme~$\got X_{\red k}$ est réduite, la norme spectrale sur~$\H^0((\got U_j)_\eta, {\sch O}_{(\got U_j)_\eta})$ prend ses valeurs dans~$|k|$, et~$$\H^0(\got U_j, {\sch O}_{\got U_j})=\H^0\left((\got U_j)_\eta, {\sch O}_{(\got U_j)_\eta}\right)\zero.$$

Fixons un entier~$n$. Le fibré~$\got L^{\otimes n}$ est trivialisé par le recouvrement~$(\got U_j)$. On définit comme suit une métrique sur~$\H^0(\got X_\eta,\got L_\eta ^{\otimes n})$ : une section globale~$f$ de~$\got L_\eta^{\otimes n}$ étant donnée, on choisit pour tout~$j$ une trivialisation de~$\got L^{\otimes n}$ sur~$\got U_j$ ; elle permet de voir la restriction~$f_{|(\got U_j)_\eta}$ comme  une fonction analytique sur ~$(\got U_j)_\eta$, dont la norme ne dépend pas de la trivialisation en question, et on définit~$||f||$ comme le maximum des normes ainsi obtenues lorsque~$j$ varie.  En vertu de ce qui précède,~$||f||\in k$ et~$f$ s'étend en une section globale de~$\got L^{\otimes n}$ sur~$\got X$ si et seulement si~$||f||\leq 1$.

\medskip
Le~$k$-espace vectoriel~$\H^0(\got X_\eta,\got L_\eta^{\otimes n})$ étant muni d'une norme qui prend ses valeurs dans~$|k|$, il possède une base orthonormale pour cette dernière. On peut par ailleurs écrire, en vertu des théorèmes GAGA,~$$\H^0(\got X_\eta,\got L_\eta^{\otimes n})\simeq\H^0(\sch Y_k\an,(\sch L_k\an)^{\otimes n})\simeq \H^0(\sch Y_k,\sch L_k^{\otimes n})\simeq\H^0(\sch Y,\sch L^{\otimes n})\otimes_Fk,$$ et partant identifier~$\H^0(\got X_\eta,\got L_\eta^{\otimes n})$ à~$\H^0(\sch Y,\sch L^{\otimes n})\otimes_Fk$. 

Par densité de~$F$ dans~$k$, il existe une base orthonormée de~$\H^0(\got X_\eta,\got L_\eta^{\otimes n})$ constituée d'éléments appartenant à~$\H^0(\sch Y,\sch L^{\otimes n})$. 

Soit~$\mathsf A_n$ le~$F\zero$-module~$\H^0(\sch Y,\sch L^{\otimes n})\zero$ ; on déduit de ce qui précède que ~$$\H^0(\got X,\got L^{\otimes n})=\H^0(\got X_\eta,\got L_\eta^{\otimes n})\zero\simeq \mathsf A_n\otimes_{F\zero}k\zero$$ et que~$\mathsf A_n\otimes_{F\zero}F\simeq \H^0(\sch Y,\sch L^{\otimes n})$. 

Par conséquent, si l'on appelle~$\mathsf B$ la~$F\zero$-algèbre graduée~$\bigoplus \mathsf A_n$, alors~$$\mathsf B\otimes_{F\zero}k\zero\simeq \bigoplus \limits_n \H^0(\got X,\got L^{\otimes n}).$$ Le fibré~$\got L_s$ étant ample, on a vu lors de la preuve de l'assertion i) au~\ref{algebgen} que la~$k\zero$-algèbre graduée~$\bigoplus \limits_n \H^0(\got X,\got L^{\otimes n})$ est de présentation finie, et que~$\got X$ s'identifie au complété formel le long de sa fibre spécial du~$k\zero$-schéma projectif associé. Par conséquent,~$\mathsf B$ est de présentation finie, et si l'on pose~$\bnd  Y={\rm Proj}\;\mathsf B$ alors~$\got X$ s'identifie au complété formel de~$\bnd Y\otimes_{F\zero}k\zero$ le long de sa fibre spéciale.

\medskip
On a par ailleurs~$\mathsf B\otimes_{F\zero} F\simeq \bigoplus\limits_n \H^0(\sch Y,\sch L^{\otimes n})$. On en déduit, le fibré~$\sch L$ étant ample, que~$\bnd Y_\eta\simeq \sch Y$, ce qui achève la démonstration.~$\Box$

\deux{remgotxalg} {\bf Remarque.} En vertu de GAGA, la restriction de~$\sch X\mapsto \widehat {\sch X}$ à la catégorie des~$k\zero$-schémas propres est pleinement fidèle. On peut donc reformuler l'assertion i) du théorème~\ref{algebrformcurv} en disant que~$\sch X\mapsto \widehat {\sch X}$ induit une équivalence entre la catégorie des~$k\zero$-courbes algébriques propres et celle des~$k\zero$-courbes formelles propres.

\chapter{Étude locale des courbes analytiques}
 
\markboth{Étude locale}{Étude locale}
\section{Chirurgie sur les courbes analytiques}\label{CHI}

\subsection*{Prolongement d'un revêtement étale}

\deux{applikrasner} {\bf Lemme.} {\em Soit~$x$ un point de~$\pk$ et soit~$U$ une composante connexe de~$\pk\setminus\{x\}$. Soit~$y\in U$, soit~$V\to ]y;x[^\flat$ un revêtement fini étale et soit~$t\in ]y;x[$. Il existe un espace~$k$-analytique~$W$ normal et lisse en dehors d'un ensemble fini de points rigides et un morphisme fini et plat~$W\to U$ tel que~$$W\times_U]t;x[^\flat\simeq_U V\times_U]t;x[^\flat.$$}

\medskip
{\em Démonstration.} Choisissons~$\tau\in ]y;t[$. Appartenant à l'intérieur d'un intervalle ouvert tracé sur~$\pk$, le point~$\tau$ est de type 2 ou 3 et est de ce fait situé au-dessus du point générique du schéma~$\PP^1_k$ ; le corps des fonctions~$F$ de ce dernier s'identifie donc naturellement à un sous-corps dense de~$\hres(\tau)$. La fibre~$V_\tau$ est isomorphe au spectre analytique d'une~$\hres(\tau)$-algèbre finie étale, laquelle provient par le lemme de Krasner d'une~$F$-algèbre finie étale. Il existe par conséquent un morphisme~${\sch Y}\to \PP^1_k$ qui est fini, plat et génériquement étale sur son but, et un isomorphisme~$V_\tau\simeq {\sch Y}\an_\tau$ ; en vertu des propriétés satisfaites par~${\sch Y}\to \PP^1_k$, la~$k$-courbe~$\sch Y$ est lisse en dehors d'un ensemble fini de points fermés ; par ailleurs l'on peut, quitte à remplacer~$\sch Y$ par sa normalisée, la supposer normale. 

\medskip
Le point~$\tau$ étant situé au-dessus du point générique de~$\pk$, le morphisme~${\sch Y}\an\to \pk$ est étale au-dessus d'un voisinage de~$\tau$ ; par conséquent, l'isomorphisme~$V_\tau\simeq {\sch Y}\an_\tau$ s'étend en un~$\Omega$ isomorphisme~$\iota$ entre~$V\times_U\Omega$ et~${\sch Y}\an\times_{\pk}\Omega$ pour un voisinage ouvert convenable~$\Omega$ de~$\tau$ dans~$U$, voisinage qui contient un ouvert de la forme~$]\tau';\tau[^\flat$ pour un certain~$\tau'\in ]y;\tau[$ ; par abus, l'on notera encore~$\iota$ l'isomorphisme déduit de~$\iota$ par le changement de base~$]\tau';\tau[^\flat\hookrightarrow \Omega$. 

\medskip
Soit~$Z$ la composante connexe de~$\pk\setminus\{\tau\}$ qui contient~$y$, et soit~$W$ l'espace~$k$-analytique obtenu par recollement de~${\sch Y}\an\times_{\pk}Z$ et~$V\times_U]\tau'; x[^\flat$ le long de l'isomorphisme~$\iota$ entre leurs ouverts respectifs~${\sch Y}\an\times_{\pk} ]\tau';\tau[^\flat$ et~$V\times_U]\tau';\tau[^\flat$. Puisque~$\iota$ commute aux projections sur~$]\tau';\tau[^\flat$, les flèches~$${\sch Y}\an\times_{\pk}Z\to Z\;{\rm et}\;V\times_U]\tau'; x[^\flat \to ]\tau'; x[^\flat~$$ induisent un morphisme~$ \psi: W\to U$. 

\medskip
En vertu de l'égalité~$Z\cap ]\tau'; x[^\flat=]\tau';\tau[^\flat$, le morphisme~$\psi$ est fini et plat localement sur son but, et partant globalement ; l'on déduit par ailleurs de l'inclusion~$]t;x[^\flat\subset ]\tau;x[^\flat$ que~$\psi\inv(]t;x[^\flat)$ est~$U$-isomorphe à~$V\times_U ]t;x[^\flat$. 

\medskip
La~$k$-courbe~${\sch Y}$ est normale et lisse en dehors d'un ensemble fini de points fermés ; par conséquent, l'ouvert~${\sch Y}\an\times_{\pk}Z$ de~${\sch Y}\an$ est normal, et lisse en dehors d'un ensemble fini de points rigides ; par ailleurs,~$V$ est muni d'un morphisme étale vers l'espace lisse~$U$, et est de ce fait lisse, et {\em a fortiori} normal. Comme~$W$ est construit par recollement de~${\sch Y}\an\times_{\pk}Z$ et~$V\times_U]\tau'; x[^\flat$, il est normal, et lisse en dehors d'un ensemble fini de points rigides.~$\Box$

\subsection*{La chirurgie proprement dite}

\deux{theochir}  {\bf Théorème.} {\em Soit~$X$ une courbe~$k$-analytique sans bord et génériquement lisse, soit~$S$ un sous-ensemble fini de points lisses de~$X$, et soit~$\Omega$ une réunion finie de composantes connexes de~$X\setminus S$ telle que~$X-\Omega$ soit compact ; on note~${\sch E}$ l'ensemble des composantes connexes de~$X\setminus S$ non incluses dans~$\Omega$, et~$\sch B$ l'ensemble fini~$\br X S\ctd \Omega$. Il existe : 

\medskip
\begin{itemize}

\item [-] un fermé~$T$ de~$X$ contenu dans~$\Omega$ ; 

\item[-]  une~$k$-courbe algébrique projective~$\sch X$ dont le lieu lisse est dense, et qui est normale si~$X$ est normale ;

\item[-] une famille~$(\Omega_b)_{b\in {\sch B}}$ d'ouverts connexes, non vides, normaux et deux à deux disjoints de~${\sch X}\an$ ;

\item[-] un isomorphisme entre~$X\setminus T$ et un ouvert de~${\sch X}\an$ modulo lequel : 

 \begin{itemize}

\medskip
\item[1)]~${\sch X}\an\setminus S=\coprod\limits_{W\in \sch E}W\coprod \left(\coprod\limits_b\Omega_b\right)$ ;

\medskip
\item[2)]~$\sch B$ s'identifie à~$\br {{\sch X}\an} S \ctd {\coprod\limits_b\Omega_b}$, et pour tout~$b\in \sch B$, la branche correspondante de~${\sch X}\an$ appartient à~$\br{ {\sch X}\an} S \ctd {\Omega_b}$ et en est le seul élément. 

\end{itemize}
\end{itemize}
}

\medskip
{\em Démonstration}. Pour tout~$b\in \sch B$ l'on note~$x_b$ l'origine de~$b$. L'on choisit un élément~$(Z_b)\in \prod\limits_{b\in \sch B}\sbr b$ tel que les sections~$Z_b$ soient deux à deux disjointes ; le sous-ensemble~$T:=\Omega-\coprod Z_b$ de~$\Omega$ est un fermé du graphe~$X$. 

\trois{constcompk} {\em La construction des~$\Omega_b$ et d'un espace analytique~$Y$ qui se révélera {\em a posteriori} être égal à~${\sch X}\an$}.

\medskip
{\em Fixons~$b\in \sch B$}. Comme~$x_b$ est un point lisse de~$X$, il existe un revêtement fini étale~$\phi$ de source un voisinage ouvert de~$x_b$ dans~$X$ et de but un ouvert de~$\pk$ ; on peut supposer que la source de~$\phi$ est un arbre et que~$x_b$ est le seul antécédent de~$\phi(x_b)$. 

\medskip
Soit~$U$ la composante connexe de~$\pk\setminus\{\phi(x_b)\}$ qui correspond à la branche~$\phi(b)$. On peut supposer, quitte à restreindre~$Z_b$, qu'il existe un élément~$y\in U$ tel que~$]y;\phi(x_b)[^\flat$ soit contenu dans le but de~$\phi$ et tel que~$Z_b$ soit une composante connexe de~$\phi\inv(]y;\phi(x_b)[^\flat)$ ; en vertu du lemme~\ref{applikrasner}, on peut de surcroît faire l'hypothèse qu'il existe un espace~$k$-analytique normal~$\Omega_b$, lisse en dehors d'un nombre fini de points rigides, muni d'un morphisme fini et plat~$\psi$ vers~$U$ et d'un~$U$-isomorphisme~$\iota_b : \psi\inv(]y;\phi(x_b)[^\flat)\simeq Z_b$. 

\medskip
Les ouverts de~$X$ de la forme~$Z_b\times_{]y;\phi(x_b)[^\flat}]t;\phi(x_b)[^\flat$, où~$t$ parcourt~$]y;\phi(x_b)[$, forment une base de sections de~$b$ ; chacun d'eux s'identifie, {\em via} l'isomorphisme~$\iota_b$, à l'ouvert~$\psi\inv(]t;\phi(x_b)[^\flat)$ de~$\Omega_b$. 

\medskip
Les parties de~$U$ de la forme~$U-]t;\phi(x_b)[^\flat$, où~$t$ parcourt~$]y;\phi(x_b)[$, forment un système cofinal de compacts de~$U$ ; il résulte alors de la finitude de~$\psi$ que les parties de~$\Omega_b$ de la forme~$\Omega_b-\psi\inv(]t;\phi(x_b)[^\flat)$, où~$t$ parcourt~$]y;\phi(x_b)[$, forment un système cofinal de compacts de~$\Omega_b$ ; cela peut se reformuler comme suit : {\em modulo l'isomorphisme~$\iota_b$, les parties de la forme~$\Omega_b-Z$, où~$Z$ est une section de~$b$ contenue dans~$Z_b$, forment un système cofinal de compacts de~$\Omega_b$.} 

\medskip
Soit~$Y$ l'espace~$k$-analytique obtenu en recollant~$X\setminus T$ et~$\coprod \Omega_b$ le long de~$\coprod Z_b$, que l'on peut voir {\em via}~$\coprod \iota_b$ comme un ouvert de chacun des deux espaces en jeu ; en vertu de ce qui précède, les propriétés suivantes sont satisfaites (on se permettra de voir à l'occasion certains sous-ensembles de~$X$ ne rencontrant pas~$T$ comme des sous-ensembles de~$Y$, sans mention explicite de cette identification, en espérant que le contexte évite toute ambiguïté) : 

\medskip
-~$Y$ est topologiquement séparé ;

- pour toute branche~$b$ appartenant à~$\sch B$, l'adhérence de~$\Omega_b$ dans~$Y$ est égale à~$\Omega_b\cup \{x_b\}$ et est compacte ; 

- pour toute branche~$b$ appartenant à~$\sch B$, l'on a~$\br Y S \ctd {\Omega_b}=\{b\}$. 

\medskip
Par construction, l'on peut écrire~$$Y=\left(\coprod \Omega_b\right)\coprod S \coprod \left(\coprod_{W\in \sch E} W\right)$$~$$=\left(\bigcup \Omega_b\cup\{x_b\}\right)\cup (X-\Omega).$$ Chacun de ces deux derniers termes est compact ; l'espace topologiquement séparé~$Y$ est donc compact. 

\medskip
{\em L'espace~$Y$ est lisse sur~$k$ en dehors d'un ensemble fini de points rigides, et normal si~$X$ est normale}. La courbe~$X$ est génériquement lisse ; par conséquent,~$X\setminus T$ est génériquement lisse, et normal si~$X$ est normale. On a vu plus haut que chacun des~$\Omega_b$ est normal et lisse en dehors d'un ensemble fini de points rigides ; comme~$Y$ est construit par recollement de~$X-Z$ avec~$\coprod \Omega_b$, et comme~$Y$ est compacte, l'assertion requise s'ensuit aussitôt. 

\medskip
\trois{yegalyan} {\em Algébrisation et conclusion.} L'ouvert~$X\setminus T$ de~$X$ est sans bord, puisque~$X$ est sans bord ; chacun des~$\Omega_b$ est fini sur un ouvert de~$\pk$ et est en particulier sans bord. Il en résulte que la courbe analytique~$Y$ est sans bord ; on a vu plus haut qu'elle est compacte. Elle est dès lors propre ; il s'ensuit qu'elle est isomorphe à l'analytification d'une~$k$-courbe algébrique projective~${\sch X}$. Comme~$Y$ est lisse en dehors d'un ensemble fini de points rigides, et normal si~$X$ est normale, le lieu lisse de la courbe~${\sch X}$ est dense, et celle-ci est normale si~$X$ est normale.~$\Box$

\section{Branches et valuations}

\deux{imtoutoucomp} {\bf Lemme.} {\em Soit~$\sch X$ une courbe algébrique projective intègre munie d'un morphisme fini et plat sur~$\PP^1_k$ dont on note~$\phi$ l'analytifié. Soit~$x \in {\sch X}\an$, soit~$V$ une composante connexe de~${\sch X}\an\setminus\{x\}$, et soit~$U$ l'ouvert~$\phi(V)$ de~$\pk$. 

\medskip
i) Si~$V$ ne rencontre pas~$\phi\inv(\phi((x))$, l'ouvert~$U$ est une composante connexe de~$\pk\setminus\{\phi(x)\}$, et~$V$ est une composante connexe de~${\sch X}\an-\phi\inv(\phi(x))$ ; c'est {\em a fortiori} une composante connexe de~$\phi\inv(U)$.

ii) Si~$V$ contient au moins un antécédent de~$\phi(x)$ alors~$U=\pk$.}

\medskip
{\em Démonstration.} Soit~$E$ l'ensemble fini~$V\cap \phi\inv(\phi(x))$. Si~$E$ est vide alors~$V\subset {\sch X}\an-\phi\inv(\phi(x))$ et en est  donc une composante connexe ; dès lors son image par la flèche finie et plate~$\phi$ est une composante connexe de~$\pk\setminus\{\phi(x)\}$, d'où i). 

\medskip
Montrons maintenant ii). Supposons que~$E$ est non vide. Si~$W$ est une composante connexe de~$V\setminus E$, c'est une composante connexe de~${\sch X}\an-\phi\inv(\phi(x))$, et~$\phi(W)$ est de ce fait une composante connexe de~$\pk\setminus\{\phi(x)\}$. Par conséquent,~$\Omega:=\phi(V\setminus E)$ est une réunion de composantes connexes de~$\pk\setminus\{\phi(x)\}$. L'ensemble~$E$ étant non vide,~$U=\Omega\cup \{x\}$. 

\medskip
L'application~$\phi$ est ouverte ; il s'ensuit que~$U$ est un voisinage de~$x$. Or le seul voisinage de~$x$ dans~$\pk$ qui soit réunion de~$\{x\}$ et de composantes connexes de~$\pk\setminus\{\phi(x)\}$ est~$\pk$ lui-même ; par conséquent,~$U=\pk$.~$\Box$

\subsection*{L'anneau associé à une branche}

\deux{anneauob} Soit~$X$ une courbe~$k$-analytique, soit~$x \in X$ et soit~$b$ une branche de~$X$ issue de~$x$. On notera~${\sch O}_X(b)$ la limite inductive des~${\sch O}_X(V)$ où~$V$ parcourt l'ensemble des sections de~$b$ ; il existe un morphisme naturel de~${\sch O}_{X,x}$ dans~${\sch O}_X(b)$.

\medskip
Nous dirons qu'une fonction~$f$ appartenant à~${\sch O}_X(b)$ est {\em modérée} si  pour tout~$\lambda\in \RR_+$ l'une des trois assertions suivantes est satisfaite : 

\medskip
A) il existe une section~$Z$ de~$b$ sur laquelle~$f$ est définie et telle que~$|f(z)|=\lambda$ pour tout~$z\in Z$ ;

B) il existe une section~$Z$ de~$b$ sur laquelle~$f$ est définie et telle que~$|f(z)|<\lambda$ pour tout~$z\in Z$ ;

C) il existe une section~$Z$ de~$b$ sur laquelle~$f$ est définie et telle que~$|f(z)|>\lambda$ pour tout~$z\in Z$.

\medskip
Notons que comme deux sections de~$b$ se rencontrent toujours, ces assertions sont exclusives l'une de l'autre. Si~$f$ satisfait A) (resp. B), resp. C)) nous dirons que~$|f|=\lambda$ (resp.~$|f|<\lambda$, resp.~$|f|>\lambda$) {\em le long de~$b$.}

\deux{obgerm} {\bf Lemme.} {\em Soit~$\phi:(Y,y)\to (X,x)$ un morphisme fini et plat entre germes ponctuels de courbes~$k$-analytiques. Soit~$b\in \br Y y$ et soit~$a$ son image dans~$\br X x$;

\medskip
i) Il existe un morphisme naturel de ~${\sch O}_{X,x}$-algèbres~$\phi^*:{\sch O}_X(a)\to {\sch O}_Y(b)$.

\medskip
ii) Si~$b\in \br Y y$ et si~$\phi\in {\sch O}_X(a)$ alors~$f$ est modérée si et seulement si~$\phi^*f$ est modérée, et dans ce cas le comportement de~$|f|$ le long de~$a$ est le même que celui de~$|\phi^*f|$ le long de~$b$. 

\medskip
iii) Le~${\sch O}_X(a)$-module~${\sch O}_Y(b)$ est de présentation finie et localement libre de rang~$\deg (b\to a)$ ; 

\medskip
iv) Si~$\phi$ est galoisien de groupe~$\mathsf G$ alors~$\mathsf G$ agit naturellement sur~$\br Y y$ ; si~$\mathsf H$ désigne le stabilisateur de~$b$ alors~$\mathsf H$ agit naturellement sur~${\sch O}_Y(b)$ et~${\sch O}_Y(b)^{\mathsf H}$ s'identifie à~${\sch O}_X(a)$. }

\medskip
{\em Démonstration}. On se ramène immédiatement au cas où~$\phi$ est induit par un « vrai » morphisme~$(Y,y)\to (X,x)$ entre courbes~$k$-analytiques pointées, noté encore~$\phi$ ; on peut supposer que~$\phi$ est fini et plat, que~$y$ est le seul antécédent de~$x$ sur~$Y$, que~$Y$ et~$X$ sont des arbres et, pour la preuve de iv), que~$\phi$ est galoisien de groupe~$\mathsf G$. Pour tout~$\beta\in \phi\inv(a)$ et toute section~$Z$ de~$a$, on note~$Z_\beta$ la composante connexe de~$\phi\inv(Z)$ qui correspond à~$\beta$. 

\medskip
Pour toute~$Z\in \sbr a$ il existe une flèche canonique~${\sch O}_X(Z)\to {\sch O}_Y(Z_b)$ ; par passage à la limite on obtient une flèche~$\lim\limits_{\stackrel \to Z} {\sch O}_X(Z)\to \lim\limits_{\stackrel \to S}{\sch O}_Y(Z_b)$ ; mais le terme de gauche s'identifie à~${\sch O}_X(a)$ et celui de droite à~${\sch O}_Y(b)$, d'où i). 

\medskip
Si~$Z\in \sbr a$, si~$\lambda\in \RR_+$ et si~$f\in {\sch O}_X(Z)$ il découle de la surjectivité de~$Z_b\to Z$ que~$|f|$ est égale  (resp. strictement inférieure , resp. strictement supérieure ) à~$\lambda$ en tout point de~$Z$ si et seulement si~$|\phi^*f|$ est égale à (resp. strictement inférieure à, resp. strictement supérieure à)~$\lambda$ en tout point de~$Z_b$, d'où ii).

\medskip
Prouvons maintenant iii) ; posons~$n=\deg \;\phi$. Il existe un voisinage ouvert~$U$ de~$x$ dans~$X$ tel que le~${\sch O}_X$-module~$\phi_*{\sch O}_Y$ soit libre de rang~$n$ au-dessus de~$U$ ; soit~$(f_1,\ldots,f_n)$ une base du~${\sch O}_X(U)$-module~${\sch O}_Y(\phi\inv(U))$. Pour toute section~$Z$ de~$a$ qui est contenue dans~$U$, le~${\sch O}_{X}(Z)$-module~$\prod\limits_\beta{\sch O}_{Y}(Z_\beta)={\sch O}_Y(\phi\inv(Z))$ est libre de base~$(f_i)$ ; par passage à la limite, le~${\sch O}_X(a)$-module~$\prod\limits_\beta {\sch O}_Y(\beta)$ est libre de base~$(f_i)$. 

\medskip
Fixons une section~$Z$ de~$a$ qui est contenue dans~$U$. Par ce qui précède, la flèche naturelle~$$\left(\prod  {\sch O}_Y(Z_b)\right)\otimes_{{\sch O}_X(Z)}{\sch O}_X(a)\to  \prod {\sch O}_Y(\beta)$$ est un isomorphisme ; par conséquent, ~$ {\sch O}_Y(Z_\beta)\otimes_{{\sch O}_X(Z)}{\sch O}_X(a)\to {\sch O}_Y(\beta)$ est un isomorphisme pour tout~$\beta$. Il suffit donc de montrer que le~${\sch O}_X(Z)$-module~${\sch O}_Y(Z_b)$ est projectif, de présentation finie et de rang égal à~$\deg \;(b\to a)$. Qu'il soit projectif et de présentation finie résulte du fait qu'il est facteur direct du~${\sch O}_X(Z)$-module libre de type fini~${\sch O}_Y(\phi\inv(Z))$ ; il reste à étudier son rang.

\medskip
Soit~$V$ un domaine affinoïde {\em non vide} de~$Z$. Le~${\sch O}_X(V)$-module ~${\sch O}_Y(\phi\inv(V))$ est libre de base~$(f_i)$. Par conséquent, la flèche naturelle~$${\sch O}_Y(\phi\inv(Z))\otimes_{{\sch O}_X(Z)}{\sch O}_X(V)\to {\sch O}_Y(\phi\inv(V))$$ est un isomorphisme. Comme sa source et son but sont respectivement égaux à~$\left(\prod   {\sch O}_Y(Z_\beta)\right)\otimes_{{\sch O}_X(Z)}{\sch O}_X(V)$ et~$\prod {\sch O}_Y(Z_\beta\cap \phi\inv(V))$, le morphisme~$${\sch O}_Y(Z_\beta)\otimes_{{\sch O}_X(Z)}{\sch O}_X(V)\to {\sch O}_Y(Z_\beta\cap \phi\inv(V))$$ est un isomorphisme pour tout~$\beta$. Le~${\sch O}(V)$-module~${\sch O}_Y(Z_b)\otimes_{{\sch O}_X(Z)}{\sch O}_X(V)$ est dès lors isomorphe à~${\sch O}_Y(Z_b\cap \phi\inv(V))$ qui est lui-même localement libre de rang~$\deg\;(b\to a)$ d'après la définition du degré d'une branche au-dessus d'une autre. 

\medskip
Comme~$Z$ est connexe,~${\sch O}_X(Z)$ n'a pas d'idempotents non triviaux. Le rang du~${\sch O}_X(Z)$-module projectif et de présentation finie~${\sch O}_Y(Z_b)$, vu comme fonction sur~$\spec {\sch O}_X(Z)$, est donc constant. Comme ledit module devient, après tensorisation par l'anneau non nul~${\sch O}(V)$, localement libre de rang~$\deg \;(b\to a)$, son rang est nécessairement~$\deg \;(b\to a)$, et iii) est démontré. 

\medskip
Prouvons iv). Que~$\mathsf G$ agisse sur~$\br Y y$ résulte de~\ref{bracourbe} ; que le stabilisateur~$\mathsf H$ de~$b$ agisse sur~${\sch O}_Y(b)$ est un cas particulier de l'assertion i) ci-dessus, déjà établie. Soit~$Z\in \sbr a$ et soit~$g\in \mathsf G$ ; on a pour toute branche~$\beta$ de~$\phi\inv(a)$ l'égalité~$g(Z_\beta)=Z_{g(\beta)}$ ; on peut en particulier caractériser~$\mathsf H$ comme le sous-groupe de~$\mathsf G$ qui laisse invariante la composante~$Z_b$. Par descente galoisienne l'on a~${\sch O}_X(Z)={\sch O}_Y(Z_b)^{\mathsf H}$ ; un passage à la limite sur~$Z$ fournit alors l'identification requise entre~${\sch O}_X(a)$ et~${\sch O}_X(b)^{\mathsf H}$.~$\Box$

\deux{normsurb} {\bf Proposition.} {\em Soit~$X$ une courbe~$k$-analytique, soit~$x$ un point de~$X$ et soit~$Y$ un domaine analytique fermé de~$X$ contenant~$x$. Soit~$f$ un élément de~${\sch O}_{X,x}$. 

\medskip
i) Pour toute branche~$b$ de~$X$ issue de~$x$ l'image de~$f$ dans~${\sch O}_X(b)$ est modérée ({\em cf.}~\ref{anneauob}).

\medskip
ii) Pour presque toute branche~$b$ de~$X$ issue de~$x$, on a~$|f|=|f(x)|$ le long de~$b$. 

iii) Il existe un sous-ensemble fini~$\sch B$ de~$\br X x$ et un voisinage ouvert~$V_0$ de~$x$ dans~$X$ qui est un arbre, et qui est tel que pour tout voisinage ouvert connexe~$V$ de~$x$ contenu dans~$V_0$ l'on ait~$V\cap Y=V-\coprod\limits_{b\in \sch B} b(V)$ ; l'image de~$\br Y x$ dans~$\br X x~$ par l'injection naturelle est égale à~$\br X x\setminus \sch B$ ; le domaine analytique~$Y$ est un voisinage de~$x$ dans~$X$ si et seulement si~${\sch B}=\emptyset$. }

\medskip
{\em Démonstration.} On procède en plusieurs étapes. 

\trois{pririg} {\em Preuve de i) et ii) lorsque~$x\in X\typ 0$.} Soit~$\lambda\in \RR_+$. Si~$|f(x)|>\lambda$ alors~$|f|>\lambda$ au voisinage de~$x$.
Si~$|f(x)|\leq \lambda$ et si~$\lambda >0$ alors le lieu de validité~$V$ de l'inégalité~$|f|\leq \lambda$ est un domaine analytique de~$X$, et comme~$x$ est un point rigide de~$V$ il appartient à l'intérieur topologique de~$V$ dans~$X$ ; par conséquent,~$|f|\leq \lambda$ au voisinage de~$x$. 

Si~$|f(x)|=\lambda=0$, choisissons un voisinage affinoïde~$X_0$ de~$x$ dans~$X$. Soit~$b$ une branche de~$X$ issue de~$x$ ;  elle possède une section~$Z$ contenue dans une composante irréductible~$X'_0$ de~$X'$. Si~$f_{|X'_0}$ est nilpotente alors~$|f|$  est identiquement nulle sur~$Z$ ; sinon l'ensemble des zéros de~$f$ sur~$X'_0$ est fini, et l'on peut donc restreindre~$Z$ de sorte que~$|f|$ ne s'y annule pas. 

\medskip
L'assertion i) est une conséquence immédiate des faits qui précèdent, et ii) découle de la finitude de~$\br X x~$ (rem.~\ref{pointrigfinbr}). 

\medskip
\trois{prisansbd} {\em Preuve de i) et ii) lorsque~$x\in X-X\typ 0 -\brdan X$.} La fonction~$f$ induit un morphisme de germes~$\phi$ de~$(X,x)$ vers~$(\Aff^{1,\rm an}_k,z)$ pour un certain~$z\in \Aff^{1,\rm an}_k$. Si~$z$ est rigide alors~$f$ est constante, et~$|f|$ est donc constante au voisinage de~$x$, d'où i) et ii). Sinon,~$\phi$ est fini et plat. Si~$\alpha$ appartient à~$\br {\Aff^{1,\rm an}}z$ alors~$\phi\inv(\alpha)$ est fini ; ceci permet en vertu de l'assertion ii) du lemme~\ref{obgerm}, de se ramener au cas où~$X=\Aff^{1,\rm an}_k$ et où~$f=T$. Si~$x$ n'appartient pas à~$]0;\infty[$ alors~$|T|$ est constante au voisinage de~$x$, d'où i) et ii). Si~$x\in ]0;\infty[$, posons~$r=|T(x)|$. Si~$\lambda>r$ (resp.~$\lambda <r$) alors~$|T|<\lambda$ (resp.~$|T|>\lambda$) au voisinage de~$x$. 

\medskip
Soit~$V^-$ (resp.~$V^+$) l'ouvert de~$\Aff^{1,\rm an}_k$ défini par la condition~$|T|<r$ (resp.~$|T|>r$). Les ouverts~$V^-$ et~$V^+$ de~$\Aff^{1,\rm an}_k$ sont deux composantes connexes de~$\Aff^{1,\rm an}_k\setminus\{x\}$ ; et si~$W$ est une composante connexe de~$\Aff^{1,\rm an}_k\setminus\{x\}$ distincte de~$V^-$ et~$V^+$ alors~$|T|$ est identiquement égale à~$r$ sur~$W$. Compte-tenu du fait que~$\br {\Aff^{1,\rm an}_k}x\to \pi_0(\Aff^{1,\rm an}_k\setminus\{x\})$ est bijective, i) et ii) s'ensuivent immédiatement.

\trois{coroldomferm} {\em Preuve de iii), ainsi que de i) et ii) dans le cas général.} Supposons tout d'abord que~$x\in X\typ{1,4}$ ; dans ce cas, i) et ii) sont satisfaites en vertu de~\ref{pririg} et~\ref{prisansbd} ; et iii) est vérifiée puisque~$Y$ est alors un voisinage de~$x$ dans~$X$. 

Supposons maintenant que~$x\in X\typ{2,3}$. Quitte à étendre les scalaires au complété de la clôture parfaite de~$k$ puis à réduire~$X$ (ce qui ne modifie pas les espaces topologiques en jeu), on peut supposer que~$X$ est quasi-lisse en~$x$ ; en remplaçant~$X$ par un voisinage ouvert convenable de~$x$, on se ramène au cas où~$X$ est lui-même un domaine analytique fermé d'une courbe lisse~$X'$.

\medskip
{\em Preuve de iii)}. Les deux dernières assertions de iii) sont des conséquences évidentes de la première. Pour établir celle-ci, il suffit de démontrer l'existence d'un voisinage~$V_0$ de~$x$ qui est un arbre et qui est tel que~$V_0\cap Y=V_0-\coprod\limits_{b\in \sch B} b(V_0)$ ; le fait que~$b(V_0)\cap V=b(V)$ pour tout voisinage ouvert connexe~$V$ de~$x$ dans~$V_0$ conduira alors à la conclusion souhaitée. 

\medskip
Il existe un voisinage ouvert~$W$ de~$x$ dans~$X'$, une famille finie~$(f_1,\ldots,f_n)$ de fonctions analytiques sur~$W$ et un entier~$m\geq n$ vérifiant la propriété suivante : le domaine analytique~$X\cap W$ (resp.~$Y\cap W$) de~$W$ est  égal à l'ensemble des points~$w\in W$ tels que~$|f_i(w)|\leq |f_i(x)|$ pour tout~$i\leq m$ (resp. pour tout~$i\leq n$) ; on peut choisir~$W$ suffisamment petit pour que ce soit un arbre, et pour que les bords topologiques de~$X\cap W$ et~$Y\cap W$ dans~$W$ soient tous deux contenus dans~$\{x\}$ ; cette dernière condition assure que toute composante connexe de~$W\setminus\{x\}$ qui rencontre~$X$ (resp.~$Y$) est contenue dans~$X$ (resp.~$Y$). 

Notons qu'en vertu du cas intérieur déjà traité (\ref{prisansbd}) chaque fonction~$f_i$ est modérée sur toute branche de~$X'$ issue de~$x$. 

\medskip
Soit~${\sch B}_0$ (resp.~${\sch B}_1$) le sous-ensemble de~$\br W x$ formé des branches~$b$ pour lesquelles il existe un entier~$j$ appartenant à~$\{1,\ldots, n\}$ (resp. à~$\{1,\ldots,m\}$) tel que~$|f_j|>|f_j(x)|$ le long de~$b$. Les ensembles~${\sch B}_0$ et~${\sch B}_1$ sont finis d'après le cas intérieur déjà traité (\ref{prisansbd}) ; on peut donc, quitte à raboter certaines des composantes connexes de~$W\setminus\{x\}$, supposer que pour toute~$b\in {\sch B}_0$ (resp.~${\sch B}_1$) il existe~$j$ appartenant à~$\{1,\ldots, n\}$ (resp. à~$\{1,\ldots,m\}$) tel que~$|f_j(z)|>|f_j(x)|$ pour tout~$z\in b(W)$. Soit~$b\in \br W x$. 

\medskip
Si~$b\notin {\sch B}_0$ alors~$|f_i|=|f_i(x)|$ le long de~$b$ pour tout~$i$ appartenant à~$\{1,\ldots, n\}$ ;  il existe donc une section de~$b$ contenue dans~$Y$ ; par conséquent,~$b(W)$ rencontre~$Y$ ce qui, en vertu de l'hypothèse faite sur~$W$, assure que~$b(W)\subset Y$. 

Si~$b\in {\sch B}_0$ il existe~$j\in  \{1,\ldots, n\}$ tel que~$|f_j(z)|>|f_j(x)|$ pour tout~$z\in b(W)$, ce qui assure que~$b(W)\cap Y=\emptyset$.

En utilisant {\em mutatis mutandis} les mêmes arguments, on voit que si~$b\notin {\sch B}_1$ alors~$b(W)\subset X$, et que si~$b\in {\sch B}_1$ alors~$b(W)\cap X=\emptyset$. 

\medskip
Il s'ensuit que~$W\cap X=W-\coprod\limits_{b\in {\sch B}_1} b(W)$, et que~$W\cap Y=W-\coprod\limits_{b\in {\sch B}_0} b(W)$. Si l'on pose~$V_0=W\cap X$ alors~$V_0$ est un voisinage ouvert de~$x$ dans~$X$  ; l'ensemble~$\br X x$ s'identifie à~$\br W x -{\sch B}_1$ ; et modulo cette identification, on a~$$V_0\cap Y=V_0-\coprod_{b\in {\sch B}_0-{\sch B}_1}b(V_0),$$ d'où iii). 

{\em Preuve de i) et ii).} Elles découlent immédiatement de l'assertion iii) que l'on vient d'établir, appliquée au domaine analytique fermé~$X'$ de~$X$, et des assertions i) et ii) dans le cas intérieur, prouvées au~\ref{prisansbd}.~$\Box$ 
 
\deux{lisseval}{\bf Théorème.} {\em Soit~$X$ une courbe~$k$-analytique, soit~$x\in X-X\typ 0 -\brdan X$,  et soit~$b\in \br X x$. Supposons qu'il existe une branche~$\beta$ de~$X$ issue de~$x$ et différente de~$b$. Le point~$x$ appartient alors à~$X\dtr$, et il existe une fonction~$f$ inversible dans~${\sch O}_{X,x}$ telle que~$|f|>|f(x)|$ le long de~$b$ et telle que~$|f|\leq f(x)$ le long de chacune des autres branches de~$X$ issues de~$x$.}

\medskip
{\em Démonstration.} On procède en plusieurs temps.

\trois{valbrlisse} {\em Réduction au cas où~$X$ est lisse et connexe.} Choisissons un voisinage affinoïde~$V$ de~$x$ dans~$X$ et une fonction analytique sur~$V$ dont la restriction à chacune des composantes irréductibles de~$V$ est non constante. Cette fonction définit un morphisme~$V\to \Aff^{1,\rm an}_k$ qui est de dimension relative nulle, et intérieur en~$x$ ; soit~$\xi$ l'image de~$x$ sur~$\Aff^{1,\rm an}_k$ ; comme~$x$ est non rigide,~$\xi$ est non rigide. 

La droite affine étant réduite, il existe en vertu de~\ref{bondeviss} un voisinage affinoïde~$V_0$ de~$x$ dans~$V_{\rm red}$ et un voisinage affinoïde~$U$ de~$\xi$ dans~$\Aff^{1,\rm an}_k$ tel que~$V\to \Aff^{1,\rm an}_k$ induise un morphisme fini et plat~$V_0\to U$ admettant une factorisation~$V_0\to V'_0\to U$ où~$V'_0\to U$ est fini étale, et où~$V_0\to V'_0$ est fini, radiciel et plat. On peut successivement remplacer :

\medskip
-~$X$ par le voisinage~$V$ de~$x$ ; 

-~$V$ par~$V_{\rm red}$ (le problème étudié est insensible au quotient par un idéal nilpotent) ; 

-~$V_{\rm red}$ par le voisinage~$V_0$ de~$x$ ; 

-~$V_0$ par~$V'_0$ et~$x$ par son image sur~$V'_0$ (en raison du fait que~$V_0\to V'_0$ induit un homéomorphisme entre les espaces topologiques sous-jacents).

\medskip
On se ramène ainsi au cas où~$X$ est lisse en~$x$, puis, quitte à restreindre encore~$X$, au cas où elle est lisse et connexe.

\medskip
\trois{valbrgen}{\em Preuve dans le cas où~$X$ est lisse et connexe.} Soit~$\Omega$ la réunion de~$b(X)$ et des composantes connexes de~$X\setminus\{x\}$ qui ne sont pas relativement compactes ; le théorème~\ref{theochir}, appliqué à~$X$, à son ouvert~$\Omega$, et à l'ensemble~$S=\{x\}$, permet de remplacer~$X$ par l'analytification~${\sch X}\an$ d'une~$k$-courbe projective normale, connexe et génériquement lisse~${\sch X}$ telle que si~$V$ désigne la composante connexe de~${\sch X}\an\setminus\{x\}$ contenant~$b$, alors~$\br {{\sch X}\an} x\ctd V=\{b\}$.

\medskip
Par hypothèse, il existe une branche~$\beta$ de~${\sch X}\an$ issue de~$x$ et distincte de~$b$ ; soit~$W$ la composante connexe de~${\sch X}\an\setminus\{x\}$ contenant~$\beta$ ; on a~$W\neq V$. Il existe un point rigide~$P$ sur~$V$ (c'est dû au {\em Nullstellensatz} si~$|k\ti|\neq \{1\}$, et à la description directe de~${\sch X}\an$ sinon). Par le théorème de Riemann-Roch le fibré en droites~${\sch O}(NP)$ sur la courbe~${\sch X}$ est engendré par ses sections globales pour~$N$ assez grand ; il existe donc une fonction méromorphe sur~$\sch X$ dont~$P$ est le seul pôle ; soit~$\psi: {\sch X}\an\to \pk$ le morphisme fini et plat qu'elle induit. Par construction,~$\infty\notin \psi(W)$ ; l'on déduit alors du lemme~\ref{imtoutoucomp} que~$W$ est une composante connexe de~${\sch X}\an-\psi\inv(\psi(x))$, et que~$\psi(W)$ une composante connexe de~$\pk\setminus\{\psi(x)\}$. L'infini n'appartient pas à~$\psi(W)$, et ne peut être égal à~$\psi(x)$ qui n'est pas rigide ;  l'ouvert~$\pk\setminus\{\psi(x)\}$ de~$\pk$ n'est donc pas connexe ; il en découle que~$\psi(x)$ est de type 2 ou 3, et il en va de même de~$x$. 

\medskip
Soit~$U$ la composante connexe de~$\pk\setminus\{\psi(x)\}$ qui contient~$\infty$. Comme ce dernier n'a qu'un antécédent sur~${\sch X}\an$, à savoir~$P$, l'ouvert~$\psi\inv(U)$ de~${\sch X}\an$ est connexe et contenu dans~$V$. Soit~$a$ l'unique élément de~$\br \pk {\psi(x)}\ctd U$. Toute branche de~${\sch X}\an$ située au-dessus de~$a$ est contenue dans~$\psi\inv(U)$, et {\em a fortiori} dans~$V$. On sait par ailleurs que l'ensemble~$\psi\inv(a)\cap \br {{\sch X}\an} x$ est non vide, et il est par ce qui précède contenu dans~$\br {{\sch X}\an}x\ctd V$ qui est le singleton~$\{b\}$ ; il en résulte que~$b$ est la seule branche issue de~$x$ et située au-dessus de~$a$. Ceci permet, en vertu de l'assertion ii) du lemme~\ref{obgerm}, de supposer que~${\sch X}=\PP^1_k$, et que~$V$ est la composante connexe de~$\pk\setminus\{x\}$ qui contient~$\infty$. 

\medskip
Choisissons un point rigide~$z$ sur~$\pk-V$ (il en existe car~$x$ est de type 2 ou 3) et soit~$\mathsf m$ le polynôme minimal de~$T(z)$ sur~$k$ ; posons~$f=\mathsf m(T)$. 
Il résulte de~\ref{exemvarypol}
que l'application~$|f|$ est strictement croissante sur l'intervalle~$[z; \infty]$ qui contient~$x$, et est localement constante sur~$\pk-[z; \infty]$ ; par conséquent,~$f$ satisfait la propriété requise.~$\Box$

\deux{brancheval} Soit~$X$ une courbe~$k$-analytique, soit~$x\in X$ et soit~$b\in \br X x$. Comme~$X$ est un bon espace,~$\kappa(x)$ est dense dans~$\hres(x)$ et l'on a donc~$\red{\kappa(x)}=\red{\hres(x)}$. 

Soit~$\lambda\in \RR\ti_+$. Si~$f$ est un élément de~${\sch O}_{X,x}$ tel que~$|f|\leq \lambda$ le long de~$b$ alors~$|f(x)|\leq \lambda$, puisque~$x$ adhère à toute section de~$b$. Cela a donc un sens de considérer l'élément~$\red{ f(x)}_\lambda$ de~$\red {\kappa(x)}_\lambda$ ; si~$g$ est un élément de~${\sch O}_{X,x}$ tel que~$ |g(x)|\leq \lambda$ et~$\red {g(x)}_\lambda=\red {f(x)}_\lambda$ alors~$|g(x)-f(x)|<\lambda$ et l'on a donc~$|g|\leq \lambda$ le long de~$b$.

\medskip
Le sous-ensemble de~${\sch O}_{X,x}$ formé des fonctions~$f$ telles que~$|f|\leq \lambda$ le long de~$b$ est stable par somme ; si~$f$ est un élément de~${\sch O}_{X,x}$ tel que~$|f(x)|=\lambda$ et~$|f|>\lambda$ le long de~$b$ alors~$|1/f(x)|=1/\lambda$ et~$|1/f|<1/\lambda$ le long de~$b$ ; et si~$\alpha$ est un élément de~$k$ tel que~$|\alpha|=\lambda$ alors~$|\alpha|=\lambda$ le long de~$b$. 

On déduit de ces faits que~$$\coprod_{\mu \in \RR\ti_+}\{\red {f(x)}_{\mu}\}_{f\in {\sch O}_{X,x}\;{\rm et}\;|f|\leq \mu\;{\rm le}\; {\rm long}\; {\rm de}\;b}$$ est l'annéloïde
d'une valuation~$\langle .\rangle_b$ de~$\red {\kappa(x)}=\red{\hres(x)}$ qui est triviale sur~$\red k$. Si~$f$ est une fonction définie au voisinage de~$x$ telle que~$|f(x)|=\lambda$ alors~$\langle \red{ f(x)}_\lambda\rangle_b\leq 1$ (resp.~$<1$, resp.~$>1$) si et seulement si~$|f|\leq \lambda$ (resp.~$<\lambda$, resp.~$>\lambda$) le long de~$b$. La restriction de~$\langle .\rangle_b$ à~$\red {\kappa(x)}_1=\red{\hres(x)}_1$ est une valuation classique, triviale sur~$\red k_1$ ; lorsqu'on la compose avec~$|.|$ on obtient une valuation~$|.|_b$ sur~$\hres(x)$, qui par construction raffine~$|.|$ et coïncide avec celle-ci sur~$k$ ; si~$f\in {\sch O}_{X,x}$ alors~$|f(x)|_b\leq 1$ (resp.~$<1$, resp.~$>1$) si et seulement si~$|f|\leq 1$ (resp.~$<1$, resp.~$>1$) le long de~$b$ ; le corps résiduel de~$|.|_b$ est celui de la restriction de~$\langle.\rangle_b$ à~$\red{\hres(x)}_1$.

\deux{kappahresmeme} Il découle des égalités ~$|\kappa(x)\ti|=|\hres(x)\ti|$ et~$\red{\kappa(x)}_1=\red{\hres(x)}_1$
et de la construction même de~$|.|_b$ que~$(\hres(x),|.|_b)$ a même groupe des valeurs et même corps résiduel que~$(\kappa(x),|.|_{b|\kappa(x)})$. 

\deux{remvalbt3} {\em Remarque.} Si~$\red {\hres(x)}_1$ est algébrique sur~$k$, autrement si~$x\in X\typ{1,3,4}$,  la valeur absolue de~$\hres(x)$ n'admet pas de raffinement strict coïncidant avec elle sur~$k$, et l'on a donc alors~$|.|_b=|.|$~.

\deux{remtransnat} {\em Remarque.} L'application~$b\mapsto \langle.\rangle_b$ (resp.~$b\mapsto |.|_b$) définit par sa construction même une {\em transformation naturelle} entre les foncteurs~$$(X,x)\mapsto \br X x\;{\rm et}\;(X,x)\mapsto \PP_{\red{\hres(x)}/\red k}\; {\rm (resp.}\;{\rm et}\;(X,x)\mapsto \PP_{\hres(x)}\;\;)$$ ayant pour source la catégorie dont les objets sont les germes de courbes~$k$-analytiques en un point non rigide, et les flèches les morphismes à fibre finie.

\deux{bijvalgrad} {\bf Théorème.} {\em Soit~$X$ une courbe~$k$-analytique et soit~$x\in X\typ{2,3}$ ; notons~$\langle.\rangle_0$ la valuation
triviale sur~$\red {\hres(x)}$. L'application~$b\mapsto \langle.\rangle_b$ établit une injection~$\br X x\hookrightarrow \PP_{\red{\hres(x)}/\red k}\setminus\{\langle.\rangle_0\}$ dont l'image~$\sch U$ est  de complémentaire fini ; cette injection est bijective si et seulement si~$x$ appartient à l'intérieur analytique de~$X$ ; l'ouvert quasi-compact~$\sch U\cup \{\langle.\rangle_0\}$ de~$\PP_{\red {\hres(x)}/\red k}$ est égal à~$\red {(X,x)}$.}
 
\medskip
{\em Démonstration.} On procède en plusieurs étapes.

\trois{surjbrval} {\em Preuve dans le cas particulier où~$x\notin \brdan X$.} On a alors l'égalité~$\red{(X,x)}=\PP_{\red{\hres(x)}/\red k}$ ; il suffit dès lors de démontrer que~$\br X x \to \PP_{\red{\hres(x)}/\red k}$ est injective, d'image égale à~$\PP_{\red{\hres(x)}/\red k}\setminus\{\langle.\rangle_0\}$. 

\medskip
{\em Montrons que toute~$\red k$-valuation non triviale de~$\red {\hres(x)}$ est de la forme~$\langle.\rangle_b$ pour une certaine~$b\in \br X x$.} Soit~$\langle.\rangle$ une~$\red k$-valuation non triviale de~$\red {\hres(x)}$. Le lemme~\ref{lemomega} fournit un élément~$\omega$ de~$\red{\hres(x)}$ tel que~$\langle \omega \rangle>1$ et tel que~$\langle.\rangle$ soit la {\em seule}~$\red k$-valuation de~$\red{\hres(x)}$ à posséder cette propriété ; remarquons que~$\omega$ est obligatoirement transcendant sur~$\red k$.

\medskip
Soit~$\lambda$ le degré de~$\omega$. Choisissons~$f\in {\sch O}_{X,x}\ti$ telle que~$\red{f(x)}=\omega$. Comme~$\red {f(x)}$ est transcendant et de degré~$\lambda$ sur~$\red k$ et comme~$X$ est sans bord, la fonction~$f$ induit un morphisme fini et plat~$\phi$ du germe~$(X,x)$ vers le germe~$(\PP^{1,\rm an}_k,\eta_\lambda)$. Soit~$a$ l'unique branche de~$\pk$ issue de~$\eta_\lambda~$ et contenue dans la composante connexe de~$\infty$ dans~$\pk\setminus\{\eta_\lambda\}$ ; il existe une branche~$b$ de~$X$ issue de~$x$ et située au-dessus de~$a$. Par construction,~$|T|>\lambda$ le long de~$a$ ; en conséquence,~$|f|>\lambda$ le long de~$b$; il s'ensuit que~$\langle \omega\rangle_b =\langle \red {f(x)}\rangle_b>1$, et par définition de~$\omega$ cette inégalité force~$\langle.\rangle_b$ à être égale à~$\langle.\rangle$.

\medskip
Il suffit pour conclure de démontrer que si~$b\in \br X x$ alors~$\langle.\rangle_b$ est non triviale, et différente de~$\langle.\rangle_\beta$ pour toute branche~$\beta$ de~$X$ issue de~$x$ et différente de~$b$. Soit donc~$b\in \br X x$. Il découle de~\ref{valxideux} et~\ref{valxitrois} qu'il existe {\em au moins deux}~$\red k$-valuations non triviales sur~$\red {\hres(x)}$ ; il y a par conséquent, en vertu du~\ref{surjbrval}, au moins deux branches de~$X$ issues de~$x$. Nous sommes dans les conditions d'application du théorème~\ref{lisseval} ; on en déduit l'existence d'une fonction inversible~$f$ dans~${\sch O}_{X,x}$ telle que~$|f|>|f(x)|$ le long de~$b$, et telle que~$|f|\leq |f(x)|$ le long de~$\beta$ pour toute branche~$\beta$ de~$\br X x$ différente de~$b$. Si l'on pose~$\lambda=|f(x)|$, l'on a donc~$\langle \red{ f(x)}_\lambda\rangle_b>1$, ce qui montre que~$\langle.\rangle_b$ est non triviale, et~$\langle\red {f(x)}_\lambda \rangle_\beta\leq 1$ pour toute~$\beta\in \br X x\setminus\{b\}$, ce qui montre que~$\langle.\rangle_\beta\neq\langle.\rangle_b$ quelle que soit~$\beta\in \br X x\setminus\{b\}$ et achève la démonstration dans le cas où~$x\notin \brdan X$. 

\trois{passclotpar} {\em Le cas général : une première réduction}. Soit~$F$ le complété de la clôture parfaite de~$k$. L'application~$X_F\to X$ est un homéomorphisme ; si~$x_F$ désigne l'unique antécédent de~$x$ sur~$X_F$ alors~$\br {X_F}{x_F}$ s'identifie à~$\br X x$ ; le corpoïde~$\red {\hres(x_F)}$ (resp.~$\red F$) est une extension radicielle de~$\red {\hres (x)}$ (resp.~$\red k$), ce qui entraîne que ~$\PP_{\red{\hres(x_F)}/\red F}\to  \PP_{\red{\hres(x)}/\red k}$ est bijective ; ce fait, combiné au critère de Temkin assure par ailleurs que~$x_F\in \brdan X_F$ si et seulement si~$x\in \brdan X$. 

\medskip
Ces remarques
(couplées à la compatibilité des différentes constructions en jeu à l'extension des scalaires de~$k$ à~$F$)
autorisent à supposer, quitte à étendre les scalaires à~$F$, que~$k$ est parfait. Les questions considérées sont insensibles aux phénomènes de nilpotence, ce qui permet de faire l'hypothèse que~$X$ est réduite, et partant quasi-lisse en~$x$. Il est alors loisible de la restreindre de sorte qu'elle s'identifie à un domaine analytique fermé d'une courbe lisse, et en particulier sans bord, que l'on notera~$X'$.

\trois{surjbrvalgen} {\em Preuve du cas général.} L'assertion iii) de la proposition~\ref{normsurb} assure que~$\br X x\to \br {X'}x$ est une injection dont l'image est de complémentaire fini ; par conséquent,~$\br X x\to \PP_{\red{\hres(x)}/\red k}$ induit une bijection entre~$\br X x$ et un sous-ensemble~$\sch U$ de ~$\PP_{\red{\hres(x)}/\red k}\setminus\{\langle.\rangle_0\}$ qui est de complémentaire fini. Soient~$\langle.\rangle_1,\ldots,\langle.\rangle_n$ les valuations appartenant à~$\PP_{\red{\hres(x)}/\red k}\setminus\{\langle.\rangle_0\}\setminus \sch U$ ; pour tout~$i$, le lemme~\ref{lemomega} assure l'existence d'un élément homogène~$\omega_i$ de~$\red{\hres(x)}$ tel que~$\langle\omega_i\rangle_i>1$ et tel que~$\langle.\rangle_i$ soit la seule~$\red k$-valuation de~$\red{\hres(x)}$ à posséder cette propriété ; pour tout~$i$, on choisit~$f_i$ dans~$\sch O_{X',x}\ti$ telle que~$\red {f_i(x)}=\omega_i$. 

\medskip
Par choix des~$\omega_i$, l'ouvert quasi-compact~$\PP_{\red{\hres(x)}/\red k}\{\omega_1,\ldots,\omega_n\}$ est égal à~$\sch U\cup \{\langle.\rangle_0\}$ ; par la théorie de Temkin, il s'identifie à~$\red{(Y,x)}$, où~$Y$ est le domaine analytique fermé de~$X'$ défini par les inégalités~$|f_i|\leq \|f_i(x)|$ . On déduit de cette description que le sous-ensemble~$\br Y x$ de~$\br {X'}x$ est obtenu par l'ablation pour tout~$i$ de l'unique branche issue de~$x$ le long de laquelle~$|f_i|>|f_i(x)|$, qui correspond à la valuation~$\langle.\rangle_i$ ; autrement dit,~$\br Y x=\br X x$ ; en vertu de l'assertion  iii) de la proposition~\ref{normsurb}, les germes~$(Y,x)$ et~$(X,x)$ coïncident ; par conséquent,~$\sch U\cup \{\langle.\rangle_0\}=\red{(X,x)}$, ce qu'il fallait démontrer.~$\Box$ 

\deux{apllbijbfval} Regardons maintenant plus précisément la façon dont le théorème ci-dessus se décline selon le type de~$x$.

\trois{casbij2} {\em Supposons que~$x\in X\typ 2$}. En vertu de~\ref{valxideux}, le théorème~\ref{bijvalgrad} se traduit alors comme suit : si~$\sch C$ désigne la courbe résiduelle en~$x$, l'application qui envoie~$b$ sur la restriction de~$\langle.\rangle_b$ à~$\red k_1({\sch C})$ induit une injection, dont l'image est de complémentaire fini, de~$\br X x$ dans l'ensemble des~$\red k_1$-valuations non triviales de~$\red k_1 ({\sch C})$ ; ce dernier s'identifiant à l'ensemble des points fermés de~$\sch C$, on a finalement construit une bijection entre~$\br X x$ et l'ensemble des points fermés d'un ouvert de Zariski non vide de~$\sch C$, qui est égal à~$\sch C$ tout entier si et seulement si~$x$ est un point intérieur de~$X$. 

\medskip
Soit~$b\in \br X x$ et soit~$\sch P$ le point fermé de~$\sch C$ correspondant. La valuation~$|.|_b$ de~$\hres(x)$ est  composée de~$|.|$ et de la valuation discrète associée à~$\sch P$ ; par conséquent, si~$\lambda \in \hres(x)$ est  tel  que~$|\lambda|=1$ et que~$\red \lambda$ ait un zéro d'ordre~$1$ en~$\sch P$, alors~$|\hres(x)\ti|_b=|\hres(x)\ti|\oplus|\lambda|_b^\ZZ$, et~$|\lambda|_b$ est infiniment proche de~$1$ inférieurement ; quant au corps résiduel de~$|.|_b$, c'est celui de~$\sch P$.

\trois{casbij3} {\em Supposons que~$x\in X\typ 3$}. En vertu de~\ref{valxitrois}, le théorème~\ref{bijvalgrad} se traduit alors comme suit : soit~$r$ un élément de~$|\kappa(x)\ti|$ n'appartenant pas à~$\sqrt{|k\ti|}$ et soit~$f\in {\sch O}_{X,x}$ telle que~$|f(x)|=r$ ; il y a au plus une branche de~$X$ issue de~$x$ le long de laquelle~$|f|<r$, au plus une le long de laquelle~$|f|>r$, et aucune le long de laquelle~$|f|=r$ ; il y a par conséquent au plus deux branches de~$X$ issues de~$x$, et l'égalité est atteinte si et seulement  si~$x$ est intérieur. 

Remarquons qu'il est possible que~$\br X x$ soit vide : cela se produit si et seulement si~$x$ est un point isolé de~$X$.

\deux{bijbr3} Notons une conséquence immédiate du~\ref{casbij3} ci-dessus : si~$Y\to X$ est un morphisme fini et plat entre courbes~$k$-analytiques et si~$y$ est un point de~$Y\typ 3$ dont on note~$x$ l'image sur~$X$, alors~$\br Y y \to \br X x$ est bijective. 

\deux{domcourbestop} Si~$X$ est une courbe~$k$-analytique et si~$Y$ est un domaine analytique fermé de~$X$, le bord topologique de~$Y$ dans~$X$ est un sous-ensemble fermé et discret de~$X$ contenu dans~$X\dtr$. On se propose, un tel ensemble étant donné, de décrire tous les domaines analytiques fermés de~$X$ dont il est le bord. 

\deux{propdomanferm} {\bf Proposition.} {\em Soit~$X$ une courbe~$k$-analytique et soit~$S$ un sous-ensemble fermé et discret de~$X$ contenu dans~$X\dtr$. Soit~$Y$ un sous-ensemble de~$X$. Les propositions suivantes sont équivalentes :

\medskip
1)~$Y$ est un domaine analytique fermé de~$X$ de bord égal à~$S$ ;

2)~$X-Y$ est une réunion~$\bigcup\limits _{i\in I} U_i$ de composantes connexes de~$X\setminus S$ telles que pour tout~$x\in S$ l'ensemble des indices~$i$ pour lesquels~$x\in \partial U_i$ soit fini et non vide.}

\medskip
{\em Démonstration.} L'implication 1)$\Rightarrow 2)$ découle de l'assertion iii) de la proposition~\ref{normsurb}. Supposons maintenant que 2) est vraie, et prouvons 1). L'hypothèse 2) implique que~$Y$ est fermé de bord égal à~$S$ ; comme la propriété d'être un domaine analytique fermé est locale, il suffit de vérifier que tout point~$x$ de~$X$ possède un voisinage~$U$ tel que~$V\cap Y$ soit un domaine analytique fermé de~$U$. C'est évident si~$x\notin Y$ ou si~$x$ appartient à l'intérieur topologique de~$Y$ ; il reste à traiter le cas où~$x\in S$. 

\medskip
Dans cette situation, 2) implique l'existence d'un voisinage ouvert~$V$ de~$x$ dans~$X$ qui est un arbre et est tel que~$V\cap Y$ soit de la forme~$V-\coprod_{j\in J} V_j$ où~$J$ est un ensemble {\em fini} et où les~$V_j$ sont des composantes connexes deux à deux disjointes de~$V\setminus\{x\}$. Chacune des~$V_j$ correspond à une branche~$b_j$ de~$X$ issue de~$x$, et donc à un point fermé~$\sch P_j$ de~$\red {(X,x)}$ en vertu du théorème~\ref{bijvalgrad} ci-dessus. 

\medskip
L'ouvert quasi-compact~$\red {(X,x)}\setminus\{\sch P_j\}_j$ de~$\red {(X,x)}$ s'identifie à~$\red{(Z,x)}$ pour un certain germe~$(Z,x)$ de domaine analytique fermé de~$(X,x)$ ; le théorème~\ref{bijvalgrad}  garantit que~$\br Z x= \br X x \setminus\{b_j\}_j$ ; l'assertion iii) de la proposition~\ref{normsurb} assure alors que~$Z$ et~$Y$ coïncident au voisinage de~$x$, ce qui achève la démonstration.~$\Box$

\subsection*{Image d'une branche par un morphisme de dimension relative nulle}

\deux{ouvcompphims} {\bf Lemme.} {\em Soit~$\phi : Y\to X$ une application continue entre espaces topologiques séparés, localement compacts et localement connexes. Soit~$V$ un ouvert connexe et non vide de~$Y$, et soit~$U$ un ouvert de~$X$ tel que~$\phi(V)\subset U$. Supposons que l'application~$V\to U$ induite par~$\phi$ est compacte ; l'ouvert~$V$ est alors une composante connexe de~$\phi\inv(U)$.}

\medskip
{\em Démonstration.} Il résulte de sa définition que~$V$ est une partie ouverte, connexe et non vide de~$\phi\inv(U)$ ; il suffit dès lors de vérifier qu'elle est également fermée dans~$\phi\inv(U)$. 

\medskip
Soit~$y\in \phi\inv(U)\setminus V$ et soit~$x$ l'image de~$y$ sur~$U$. Choisissons un voisinage~$U_0$ de~$x$ dans~$X$ tel que~$\overline{U_0}$ soit une partie compacte de~$U$. Soit~$W$ le compact~$(\phi_{|V})^{-1}(\overline{U_0})$. L'ouvert~$\phi\inv(U_0)-W$ de~$\phi\inv(U)$ contient~$y$ et ne rencontre pas~$V$ ; par conséquent,~$V$ est fermé dans~$\phi\inv(U)$.~$\Box$ 

\deux{phipasfin} {\bf Proposition.} {\em Soit~$\phi : Y\to X$ un morphisme entre deux courbes~$k$-analytiques, soit~$y\in Y$ et soit~$x$ son image sur~$X$ ; on suppose que~$\phi$ est de dimension relative nulle en~$y$.

\medskip
1) L'ensemble~$\br X x-\phi(\br Y y)$ est  fini. 

\medskip
2) Il existe~$U\in \arb X x$  et~$V\in \arb Y y$ possédant les propriétés suivantes : 

\medskip
 i)~$V$ est une composante connexe de~$\phi\inv(U)$ telle que~$\phi\inv(x)\cap V=\{y\}$ ; 
 
 ii)~$V\to U$ est compacte à fibres finies ; 
 
 iii)~$V\setminus\{y\}\to U\setminus\{x\}$ est finie, ouverte, et plate si~$X$ est génériquement réduite.}

\medskip
{\em Démonstration.} Notons pour commencer que l'on peut, pour montrer ces assertions, remplacer~$X$ par~$X_{\rm red}$ et~$Y$ par~$Y\times_X X_{\rm red}$ ; cela permet de se ramener au cas où la courbe~$X$ est réduite, et l'on distingue alors deux cas selon la nature du point~$x$.

\trois{brcasxrig} {\em Le cas où~$x$ est rigide.} L'ensemble~$\br X x$ est alors fini, d'où 1). Pour prouver 2), remarquons tout d'abord que~$y$ est également rigide,  ce qui entraîne que~$\phi$ est sans bord en~$y$ ; il induit de ce fait un morphisme fini d'un voisinage ouvert~$V$ de~$y$ sur un voisinage ouvert~$U$ de~$x$ ; on peut supposer que~$y$ est le seul antécédent de~$y$ sur~$V$, et que~$U$ est une composante connexe de~$\phi\inv(x)$ ; on peut également faire l'hypothèse que~$U$ et~$V$ sont des arbres. 

\medskip
Comme~$U$ est génériquement réduite, les anneaux locaux en ses points non rigides sont des corps, et il existe dès lors un ouvert de Zariski {\em dense} de~$U$ au-dessus duquel le morphisme~$V\to U$ est plat ; par conséquent il est loisible,  quitte à restreindre convenablement~$U$ et~$V$, de faire en sorte que la flèche~$V\setminus\{y\}\to U\setminus\{x\}$ soit finie et plate ; notons qu'elle est alors ouverte, ce qui achève de montrer 2), et partant la proposition, lorsque~$x$ est rigide. 

\trois{brcaspasrig} {\em  Le cas où~$x$ n'est pas rigide.} La version analytique du {\em Main Theorem} de Zariski assure l'existence d'un voisinage affinoïde~$Y_0$ de~$y$ dans~$Y$ et d'une factorisation de~$\phi$ sous la forme~$$Y_0\to W\hookrightarrow T\to X$$ où~$Y_0\to W$ est fini, où~$W$ s'identifie à un domaine affinoïde de~$T$, et où~$T\to X$ est étale ; soit~$t$ l'image de~$y$ sur~$W\subset T$. La courbe~$X$ étant génériquement réduite, il en va de même de~$T$ et de~$W$. Comme~$x$ n'est pas rigide,~$t$ ne l'est pas, et~$(Y_0,y)\to (W,t)$ est en conséquence fini et plat ; on peut donc restreindre~$W$ et~$Y_0$ de sorte que~$Y_0\to W$ soit fini et plat, et que~$y$ soit l'unique antécédent de~$t$ sur~$Y_0$. 

\medskip
Nous allons construire le voisinage~$U\in \arb X x$ par restrictions successives. On commence par le choisir de sorte qu'il existe~$Z\in \arb T t$ fini et étale sur~$U$ possédant les trois propriétés suivantes : 

\medskip
$\bullet$~$Z$ est une composante connexe de~$T\times_X U$ ; 

$\bullet$~$t$ est l'unique antécédent de~$x$ sur~$Z$ ; 

$\bullet$~$Z\cap \partial_T W\subset \{t\}$.  

\medskip
 Il existe~$W'\in \arb W t$ tel que~$Y_0\times_WW'$ ne rencontre pas~$\partial_Y Y_0$. Quitte à restreindre~$Z$ et~$U$, on peut supposer que~$Z\cap W\subset W'$. 
 
\medskip
Soit~$V$ l'image réciproque de~$Z\cap W$ dans~$Y_0$. C'est un ouvert connexe de~$Y_0$ qui ne rencontre pas~$\partial_Y Y_0$ ; c'est donc un ouvert de~$Y$. Il est inclus dans~$\phi\inv(U)$, et son intersection avec~$\phi\inv(x)$ est égale à~$\{y\}$. 

\medskip
La flèche~$V\to U$ est composée de~$V\to W\cap Z, W\cap Z\hookrightarrow Z$ et~$Z\to U$, qui sont toutes trois compactes et à fibres finies : la première et la troisième sont en effet finies, et la seconde est une immersion d'un domaine analytique fermé ; la flèche~$V\to U$ est dès lors elle-même compacte et à fibres finies ; il découle alors du lemme~\ref{ouvcompphims} que~$V$ est une composante connexe de~$\phi\inv U$. Par ailleurs,~$V\to W$ et~$Z\to U$ sont sans bord, et le bord de~$W\cap Z\hookrightarrow Z$ est contenu dans~$\{t\}$ ; par conséquent, le bord de~$V\to U$ est contenu dans~$\{y\}$. 

\medskip
On déduit de ce qui précède que~$V\setminus\{y\}\to U\setminus\{x\}$ est compacte, sans bord, et à fibres finies ; elle est donc finie. De plus, soit~$v\in V\setminus\{y\}$, soit~$w$ son image sur~$W$ et soit~$u$ son image sur~$U$. Comme le morphisme fini~$V\to Z\cap W$ est plat,~$\sch O_{V,v}$ est plat sur~$\sch O_{Z\cap W,w}$ ; comme~$w\neq t$, il est contenu dans l'intérieur topologique de~$W\cap Z$ dans~$Z$, et l'anneau local~$\sch O_{Z\cap W,w}$ est égal à~$\sch O_{Z,w}$ ; et comme~$Z\to U$ est étale,~$\sch O_{Z,w}$ est plat sur~$\sch O_{U,u}$. 

Il s'ensuit que le morphisme fini~$V\setminus\{y\}\to U\setminus\{x\}$ est plat ; cela entraîne qu'il est ouvert et achève de prouver 2). Pour établir 1), on remarque que~$\br V y \to \br {Z\cap W} t$ est surjectif, que~$\br {Z\cap W} t\to \br Z t$ est une injection dont l'image est de complémentaire fini (prop.~\ref{normsurb}), et enfin que~$\br Z t \to \br U x$ est surjective.~$\Box$ 

\deux{commentbrpasfp} {\em Commentaires.} Soit~$\phi : Y \to X$ un morphisme entre deux courbes~$k$-analytiques, soit~$y$ un point de~$Y$ en lequel~$\phi$ est de dimension nulle, et soit~$x$ son image sur~$X$ ; soient~$U\in \arb X x~$ et~$V\in \arb Y y$ satisfaisant les conditions i), ii), et iii) de la proposition~\ref{phipasfin} ci-dessus. 

\trois{uprimecommeu} Soit~$U'\in \arb U x$ et soit~$V'$ son image réciproque sur~$V$. La flèche~$V'\to U'$ est compacte et à fibres finies, et~$V'\setminus\{y\}\to U'\setminus\{x\}$ est finie, ouverte, et plate si~$X$ est réduite.  

Par ailleurs,~$V'$ est connexe. Pour le voir, on raisonne par l'absurde, en supposant qu'il possède une composante connexe~$V'_0$ qui ne contient pas~$y$. L'ouvert~$V'_0$ de~$V'$ est alors une composante connexe de~$V'\setminus\{y\}$, et son image~$U'_0$ sur~$U'$ est donc une composante connexe de~$U'\setminus\{x\}$. D'autre part,~$V'_0$ est fermé dans~$V'$, ce qui implique que la composante connexe~$U'_0$ de~$U'\setminus\{x\}$ est fermée dans l'arbre~$U'$, ce qui est absurde. 

On déduit de ce qui précède que~$U'$ et~$V'$ satisfont encore les conditions i), ii), et iii) de {\em loc. cit.}

\trois{degrebrmeme} Soit~$b\in \br Y y$ et soit~$a$ son image dans~$\br X x$. 

\medskip
L'image sur~$U$ de la composante connexe~$b(V)$ de~$V\setminus\{y\}$ est une composante connexe de~$U\setminus\{x\}$ qui est nécessairement~$a(U)$. Il s'ensuit, compte-tenu de l'assertion 1) de {\em loc. cit.}, que~$U-\phi(V)$ est une réunion finie de composantes connexes de~$U\setminus\{x\}$. 

\medskip
Si~$X$ est réduite, la flèche~$b(V)\to a(U)$ est finie et plate, et son degré ne dépend alors pas du choix de~$U$ et~$V$ (par exemple en vertu du~\ref{uprimecommeu} ci-dessus) ; on l'appellera le {\em degré de~$b$ sur~$a$} ;  cette définition coïncide avec la précédente dans le cas où~$(Y,y)\to (X,x)$ est fini et plat.

\subsection*{Fonctions sur une branche et hensélisé de la valuation associée}

\deux{introhensnorm} Soit~$X$ une courbe~$k$-analytique normale et soit~$x$ un point non rigide de~$X$. L'anneau local~${\sch O}_{X,x}$ est alors artinien et intègre ; il coïncide donc avec son corps résiduel~$\kappa(x)$. Soit~$b$ une branche de~$X$ issue de~$x$. Si~$Z$ est une section de~$b$ c'est un espace~$k$-analytique connexe, non vide et normal ; l'anneau~${\sch O}(Z)$ est par conséquent intègre ; il en résulte, par passage à la limite, que~${\sch O}(b)$ est intègre. Le corps~$\kappa(x)$ coïncidant avec~${\sch O}_{X,x}$, il se plonge dans~${\sch O}(b)$.

\deux{theohenbr}{\bf Théorème.} {\em Soit~$X$ une courbe~$k$-analytique normale et soit~$x$ appartenant à~$X-X\typ 0 -\brdan X$. Soit~$b$ une branche de~$X$ issue de~$x$ ; on note~${\sch O}_X(b)_{\rm alg}$ (resp.~${\sch O}_X(b)_{\rm sep}$) le corps égal à la fermeture algébrique (resp. séparable) de~$\kappa(x)$ dans~${\sch O}_X(b)$.  

\medskip
i) Si~$f\in {\sch O}_X(b)_{\rm alg}$ alors~$f$ est modérée. 

ii) Le sous-ensemble de~${\sch O}_X(b)_{\rm alg}$ formé des fonctions~$f$ telles que~$|f|\leq 1$ le long de~$b$ est l'anneau d'une valuation prolongeant la valuation~$|.|_b$ de~$\kappa(x)$, et notée encore~$|.|_b$. 

iii) Si~$x\in X\typ {23}$ le corps valué~$( {\sch O}_X(b)_{\rm sep}, |.|_b)$ est le hensélisé de~$(\kappa(x),|.|_b)$. En vertu de la
remarque~\ref{remvalbt3}, on a en particulier~$( {\sch O}_X(b)_{\rm sep}, |.|_b)=(\kappa(x),|.|)$ dès que~$x\in X\typ 3$.}

\medskip
{\em  Démonstration.} Commençons par deux remarques : l'assertion ii) est une conséquence triviale de i) ; et pour montrer i), on peut supposer, quitte à remplacer~$f$ par~$f^q$ pour une puissance convenable~$q$ de l'exposant caractéristique de~$k$, que~$f\in {\sch O}_X(b)_{\rm sep}$.

\medskip
Choisissons un système projectif filtrant~$((X_i,x_i))_i$ de germes d'espaces~$k$-analytiques galoisiens sur~$(X,x)$ tel que la limite inductive des~$\kappa(x_i)$ soit une clôture séparable de~$\kappa(x)$, que l'on notera~$\kappa(x)^s$ et que l'on verra comme la {\em réunion} des~$\kappa(x_i)$. Pour tout~$i$, l'ensemble~$\br {X_i}{x_i}$ est fini et possède au moins un élément situé au-dessus de~$b$ ; il existe donc un élément~$(b_i)_i$ de~$\lim\limits_{\leftarrow} \br {X_i}{x_i}$ tel que~$b_i$ soit situé pour tout~$i$ au-dessus de~$b$. 

\medskip
D'après le~\ref{brancheval} et la remarque~\ref{remtransnat}, la branche~$b_i$  induit pour tout~$i$ une valuation~$|.|_{b_i}$ de~$\kappa(x_i)$, et la famille des~$|.|_{b_i}$ se comporte bien vis-à-vis des restrictions ; elle définit donc une valuation~$|.|'$ sur~$\kappa(x)^s$ qui prolonge~$|.|_b$. Si~$f\in \kappa(x_i)$ pour un certain~$i$ alors~$|f|'\leq 1$ si et seulement si~$|f|_{b_i}\leq 1$, c'est-à-dire et seulement si~$|f|\leq 1$ le long de~$b_i$. 

\medskip
La famille~$({\sch O}_{X_i}(b_i))_i$ est un système inductif filtrant d'anneaux intègres dont les flèches de transition sont injectives (lemme~\ref{obgerm} iii) ). Sa limite inductive est une~$\kappa(x)$-algèbre intègre~$\mathsf A$, que l'on verra comme la {\em réunion} des~${\sch O}_{X_i}(b_i)$. La donnée pour tout~$i$ de l'injection de~$\kappa(x_i)$ dans ~${\sch O}_{X_i}(b_i)$ permet de voir~$\kappa(x)^s$ comme étant contenu dans~$\mathsf A$ ; il coïncide alors avec la fermeture séparable de~$\kappa(x)$ dans~$\mathsf A$.

\medskip
Soit~$f\in {\sch O}_X(b)_{\rm sep}$ ; son image dans~$\mathsf A$ vit dans~$\kappa(x)^s$, et partant dans~$\kappa(x_i)$ pour un certain~$i$ ; autrement dit, la fonction~$f$, {\em vue comme appartenant à~${\sch O}_{X_i}(b_i)$}, se prolonge en une fonction définie au voisinage de~$x_i$.

Ceci entraîne (lemme ~\ref{normsurb}) que~$f$ est modérée en tant qu'élément de ~${\sch O}_{X_i}(b_i)$ ; d'après le lemme ~\ref{obgerm} elle est alors modérée en tant qu'élément de~${\sch O}_X(b)$, ce qui achève de prouver i) , et ~$|f|\leq 1$ le long de~$b$ si et seulement si~$|f|\leq 1$ le long de~$b_i$, autrement dit si et seulement si~$|f|'\leq 1$ ($f$ est vue comme appartenant à~${\sch O}_X(b)$ pour la première condition, à~${\sch O}_{X_i}(b_i)$ pour la seconde, et à~$\kappa(x)^s$ pour la troisième). 

\medskip
Pour tout~$i$, notons~$\mathsf H_i$ le sous-groupe de~${\mathsf{Gal}} \;(\kappa(x_i)/\kappa(x))$ égal au stabilisateur de~$b_i$ {\em via} l'identification entre~${\mathsf{Gal}} \;(\kappa(x_i)/\kappa(x))$ et~${\mathsf{Gal}}\;((X_i,x_i)/(X,x))$ ; soit~$\mathsf H$ la limite projective des~$\mathsf H_i$. 

$\bullet$ Si~$x\in X\typ 2$  alors~$\mathsf H$ est le stabilisateur de~$|.|'$ dans~${\mathsf{Gal}}\;(\kappa(x)^s/\kappa(x))$ (\ref{casbij2}). 

$\bullet$ Si~$x\in X\typ 3$ alors~$\mathsf H={\mathsf{Gal}}\;(\kappa(x)^s/\kappa(x))$ (\ref{bijbr3}). Par ailleurs, on a dans ce cas~$|.|_b=|.|$ (rem.~\ref{remvalbt3}), 
et~$|.|'$ est donc l'{\em unique} prolongement de~$|.|_b$ à~$\kappa(x)^s$. Dès lors,~$\mathsf H$ est là encore
le stabilisateur de~$|.|'$ dans~${\mathsf{Gal}}\;(\kappa(x)^s/\kappa(x))$.

\medskip
Il découle du lemme~\ref{obgerm} que l'action de~$\mathsf H$ sur~$\kappa(x_s)$ s'étend à~$\mathsf A$ et que~$\mathsf A^{\mathsf H}={\sch O}_X(b)$. Le corps~$(\kappa(x)^s)^{\mathsf H}$ est donc égal à~$\kappa(x)^s\cap {\sch O}_X(b)$, c'est-à-dire à~${\sch O}_X(b)_{\rm sep}$. En conséquence~$\left({\sch O}_X(b)_{\rm sep},|.|'_{|{\sch O}_X(b)_{\rm sep}}\right)$ est le hensélisé de~$(\kappa(x),|.|_b)$. 
On a par ailleurs vu plus haut que l'anneau de~$|.|'_{|{\sch O}_X(b)_{\rm sep}}$ est constitué 
des fonctions majorées par~$1$ en valeur absolue le long de~$b$ ;  on a donc~$|.|'_{|{\sch O}_X(b)_{\rm sep}}=|.|_b$, ce qui achève la démonstration.~$\Box$

\trois{remmemeval} {\em Remarque.} Il découle de l'assertion iii) du théorème ci-dessus que lorsque~$x\in X\typ 2$, la valuation~$|.|_b$ de~${\sch O}_X(b)_{\rm sep}$ a même groupe des valeurs et même corps résiduel que sa restriction à~$\kappa(x)$ ; il nous arrivera pour cette raison dans ce cas de parler du groupe des valeurs et du corps résiduel de~$|.|_b$ sans expliciter le corps sur laquelle on la considère (lorsque~$x\in X\typ 3$ le problème ne se pose pas puisque les corps~${\sch O}_X(b)_{\rm sep}$ et~$\kappa(x)$ coïncident).

\trois{remextval}  {\em Remarque.} Soit~$(Y,y)\to(X,x)$ un morphisme fini et plat entre germes de courbes~$k$-analytiques normales et sans bord. Supposons que~$y\in Y\dtr$, soit~$b\in \br Y y$ et soit~$a$ son image dans~$\br X x$. Il est immédiat que~${\sch O}_X(a)_{\rm alg}$ (resp.~${\sch O}_X(a)_{\rm sep}$) s'envoie dans~${\sch O}_Y(b)_{\rm alg}$ (resp.~${\sch O}_Y(b)_{\rm sep}$) ; il résulte de l'assertion ii) du lemme~\ref{remtransnat} que~${\sch O}_Y(b)_{\rm alg}$ est une extension {\em de corps valués} de~${\sch O}_Y(b)_{\rm alg}$.

\deux{propbondeg} {\bf Proposition.} {\em Soit~$\phi : (Y,y)\to(X,x)$ un morphisme fini et plat entre germes de courbes~$k$-analytiques normales ; on suppose que~$y\in Y\dtr$. Soit~$b\in \br Y y$ et soit~$a$ son image dans~$\br X x$. L'extension~${\sch O}_{Y}(b)_{\rm sep}/{\sch O}_{X}(a)_{\rm sep}$ est finie de degré~$\deg \;(b\to a)$.}

\medskip 
{\em Démonstration.}  Les catégories des~$\kappa(x)$-algèbres étales, des~$\hres(x)$-algèbres étales, et des revêtements finis \'etales du germe~$(X,x)$ sont équivalentes. Soit~$F$ la fermeture séparable de~$\kappa(x)$ dans~$\kappa(y)$, et soit~$(X',x')$ le germe fini étale sur~$(X,x)$ qui lui correspond. Le corps~$\hres(x')\simeq F\otimes_{\kappa(x)}\hres(x)$ se plonge naturellement dans~$\hres(y)$ ; par conséquent,~$Y':=Y\times_XX'\to Y$ possède un~$\hres(y)$-point~$y'$ au-dessus de~$y$, et le morphisme fini étale~$(Y',y')\to (Y,y)$ est alors un isomorphisme ; la composée de sa réciproque et de la projection vers~$(X',x')$ définit une factorisation~$(Y,y)\to (X',x')\to (X,x)$ ; comme~$\kappa(x')\simeq F$, l'extension~$\kappa(x')\hookrightarrow \kappa(y)$ est purement inséparable. 

\medskip
On se ramène ainsi, pour démontrer la proposition, à traiter deux cas particuliers. 

\trois{deggermet} {\em Le cas où~$(Y,y)\to (X,x)$ est étale.} Soit~$\sch B$ le sous-ensemble de~$\br Y y$ formé des branches situées au-dessus de~$a$. Soit~$\beta\in \sch B$ et soit~$\delta$ le degré de~$\beta$ sur~$a$. L'on dispose d'un diagramme commutatif~$$\diagram \kappa(y)\rto & {\sch O}_Y(\beta)\\ \kappa(x)\uto \rto& {\sch O}_X(a)\uto \enddiagram\; ;$$ les anneaux ~$ {\sch O}_Y(\beta)$ et~${\sch O}_X(a)$ sont intègres, et~${\sch O}_Y(\beta)$ est de présentation finie et localement libre de rang~$\delta$ sur~${\sch O}_X(a)$ (lemme~\ref{obgerm}). Soit~$f\in {\sch O}_Y(\beta)_{\rm sep}$ et soit~$\chi$ son polynôme caractéristique sur~${\sch O}_X(a)$ ; il est de degré~$\delta$ et est égal à une puissance du polynôme minimal de~$f$ sur~${\rm Frac}\; {\sch O}_X(a)$. Étant séparable sur~$\kappa(y)$, l'élément~$f$ de~${\sch O}_Y(\beta)$ est séparable sur~$\kappa(x)$ et les coefficients de son polynôme minimal sur~${\rm Frac}\; {\sch O}_X(a)$ le sont alors aussi. Dès lors les coefficients de~$\chi$ sont des éléments de~${\sch O}_X(a)$ séparables sur~$\kappa(x)$ ; autrement dit, ils appartiennent à~${\sch O}_X(a)_{\rm sep}$. Ainsi, tout élément de~${\sch O}_Y(\beta)_{\rm sep}$ est annulé par un polynôme unitaire de degré~$\delta$ à coefficients dans~${\sch O}_X(a)_{\rm sep}$. 

\medskip
Il résulte du~\ref{casbij2} que~$\{|.|_\beta\}_{\beta\in \sch B}$ est l'ensemble des valuations de~$\kappa(y)$ qui prolongent~$|.|_a$ ; par ailleurs, le théorème~\ref{theohenbr} assure que~${\sch O}_X(a)_{\rm sep}$ est le hensélisé de~$(\kappa(x),|.|_a)$ et que~${\sch O}_Y(\beta)_{\rm sep}$ est le hensélisé de~$(\kappa(y),|.|_\beta)$ pour tout~$\beta$ dans~$\sch B$. Il en résulte que~$$\kappa(y)\otimes_{\kappa(x)}{\sch O}_X(a)_{\rm sep}=\prod_{\beta\in \sch B}{\sch O}_Y(\beta)_{\rm sep}.$$ En conséquence,~${\sch O}_Y(\beta)_{\rm sep}$ est pour tout~$\beta$ une extension finie séparable de~${\sch O}_X(a)_{\rm sep}$, et~$$\sum_{\beta\in \sch B}[{\sch O}_Y(\beta)_{\rm sep}:{\sch O}_X(a)_{\rm sep}]=[\kappa(y):\kappa(x)]=\deg^y \; \phi.$$

\medskip
On a vu plus haut que quelque soit~$\beta$ appartenant à~$\sch B$, chaque élément de~${\sch O}_Y(\beta)_{\rm sep}$ est de degré au plus égal à~$\deg\;(\beta\to a)$ sur~${\sch O}_X(a_{\rm sep})$ ; il en découle, grâce au théorème de l'élément primitif, que le degré de l'extension finie séparable~${\sch O}_Y(\beta)_{\rm sep}/{\sch O}_X(a)_{\rm sep}$ est au plus~$\deg\;(\beta\to a)$. Des deux égalités~$$\sum_{\beta\in \sch B} \deg\;(\beta\to a)=\deg^y \;\phi\;{\rm et}\; \sum_{\beta\in \sch B}[{\sch O}_Y(\beta)_{\rm sep}:{\sch O}_X(a)_{\rm sep}]=\deg^y \; \phi$$ l'on déduit alors que~$[{\sch O}_Y(\beta)_{\rm sep}:{\sch O}_X(a)_{\rm sep}]=\deg\;(\beta\to a)$ pour tout~$\beta\in \sch B$ et en particulier pour~$b$, ce qui achève la preuve dans le cas où~$\phi$ est étale.

\trois{deggerminsep} {\em Le cas où~$\kappa(y)$ est purement inséparable sur~$\kappa(x)$.} La valuation~$|.|_a$ admet alors un unique prolongement à~$\kappa(y)$, ce qui signifie, en vertu du~\ref{casbij2}, que~$b$ est la seule branche de~$(Y,y)$ située au-dessus de~$a$ ; par conséquent,~$\deg \;(b\to a)=\deg\;\phi$. Par ailleurs~${\sch O}_X(a)_{\rm sep}\otimes_{\kappa(x)}\kappa(y)$ étant un corps, c'est  le hensélisé de~$(\kappa(y),|.|_b)$, c'est-à-dire le corps~${\sch O}_Y(b)_{\rm sep}$ (th.~\ref{theohenbr}) ; il vient~$$[{\sch O}_Y(b)_{\rm sep}:{\sch O}_X(a)_{\rm sep}]=[\kappa(y):\kappa(x)]=\deg^y\;\phi=\deg \;(b\to a)\;.\;\Box$$

\section{Branches des courbes quasi-lisses}

\deux{unebrdeg} {\bf Lemme.} {\em Soit~$\phi : Y\to X$ un morphisme fini et plat entre courbes~$k$-analytiques, soit~$y\in Y$ et soit~$x$ son image sur~$X$.  Soit~$U$ une composante connexe de~$X\setminus\{x\}$ et soit~$V$ une composante connexe de~$\phi\inv(U)$. Supposons que~$\partial V=\{y\}$ et que~$\br Y y \ctd V$ est un singleton~$\{b\}$. Sous ces hypothèses,~$\deg\;(b\to \phi(b))=\deg \;(V\to U)$.}

\medskip
{\em Démonstration.} Posons~$a=\phi(b)$ ; comme~$b$ est contenue dans~$V$ son image~$a$ est contenue dans~$U$. Il existe alors une section~$Z$ de~$a$ qui est contenue dans~$U$ et une bijection~$\beta\mapsto Z_\beta$ entre~$\phi\inv(a)$ et~$\pi_0(\phi\inv(Z))$ telle que~$Z_\beta$ soit pour tout~$\beta$ une section de~$\beta$ dont le degré sur~$Z$ coïncide avec~$\deg\;(\beta\to a)$. Comme~$y$ est le seul antécédent de~$x$ qui adhère à~$V$, toute branche de~$Y$ située au-dessus de~$a$ et contenue dans~$V$ est issue de~$y$ ; comme~$\br Y y \ctd V=\{b\}$, la branche~$b$ est la seule branche de~$\phi\inv(a)$ qui soit contenue dans~$V$. Il s'ensuit que~$Z_b$ est la seule composante connexe de~$\phi\inv(Z)$ qui soit contenue dans~$V$ ; on a dès lors~$Z_b=V\times_UZ$, et partant~$$\deg\;(b\to a)=\deg\;(Z_b\to Z)=\deg\;(V\to U).\;\Box$$

\subsection*{Un lemme de Gabber et ses conséquences}

La suite de cet article repose de façon absolument cruciale sur la construction de fonctions méromorphes satisfaisant certaines conditions, elle-même fondée sur le lemme ci-dessous, communiqué à l'auteur par Gabber. 

\deux{lemgabber} {\bf Lemme (Gabber).} {\em Soit~$F$ un corps algébriquement clos et soit~${\sch X}$ une courbe algébrique projective, irréductible et lisse sur~$F$. Soient~$\mathsf V$ et~$\mathsf W$ deux sous-ensembles non vides de~${\sch X}(F)$ et soit~$N$ un entier strictement positif. Il existe deux diviseurs effectifs~$D$ et~$D'$ sur~${\sch X}$ tels que :

\medskip
1)~$D'$ est à support dans~$\mathsf W$ ;

2) le degré de~$D\cap \mathsf V$ est premier à~$N$ ;

3)~$D$ et~$D'$ sont linéairement équivalents.}

\medskip
{\em Démonstration.} Comme~$\mathsf V$ et~$\mathsf W$ sont non vides, on peut choisir un point~$P$ sur~$\mathsf V$ et un point~$Q$ sur~$\mathsf W$. La multiplication par~$N$ étant une isogénie de la variété abélienne~$\mathsf {Jac}\;{\sch X}$, il existe un fibré en droites~$\sch L$ sur~$\sch X$ tel que~${\sch L}^{\otimes N}\simeq {\sch O}(Q-P)$. Le théorème de Riemann-Roch assure qu'il existe un entier~$m>0$ et un diviseur effectif~$\Delta$ sur~$\sch X$ tel que~${\sch L}\otimes{\sch O}(mQ)\simeq {\sch O}(\Delta)$. Il vient :

$${\sch O}(N\Delta)\simeq {\sch L}^{\otimes N}\otimes{\sch O}(NmQ)\simeq {\sch O}(Q-P)\otimes{\sch O}(NmQ)\simeq {\sch O}((Nm+1)Q-P),~$$ et les diviseurs~$D:=P+N\Delta$ et~$D':=(Nm+1)Q$ satisfont les propriétés requises.~$\Box$

\deux{fmerop} {\bf Théorème.} {\em Soit~$X$ une courbe~$\KK$-analytique lisse et soit~$x\in X\dtr$. Soit~$b\in \br X x$. Il existe~$z\in \pkk$ et un morphisme fini étale de germes~$\phi:(X,x)\to (\pkk,z)$ tel que~$\deg\; (b\to \phi(b))$ soit premier à~$p$.}

\medskip
{\em Démonstration.} On peut appliquer la construction du théorème~\ref{theochir}, et ainsi supposer que~$X$ est l'analytifiée d'une~$k$-courbe algébrique projective, intègre et lisse ~${\sch X}$, telle que la propriété suivante soit satisfaite : si~$V$ désigne l'ouvert~$b(\sch X\an)$ alors~$\br {{\sch X}\an} x\ctd V=\{b\}$.

\medskip
Comme~$x\in X\dtr$, il résulte de~\ref{casbij2} et~\ref{casbij3} qu'il y a au moins deux branches de~${\sch X}\an$ issues de~$x$ ; il existe donc une composante connexe~$W$ de~${\sch X}\an\setminus\{x\}$ qui est différente de~$V$. Les sous-ensembles~$V(\KK)$ et~$W(\KK)$ de~${\sch X}(\KK)$ sont non vides (par le {\em Nullstellensatz} si~$|\KK\ti|\neq\{1\}$, et par la description explicite de~${\sch X}\an$ sinon). On peut donc appliquer le lemme~\ref{lemgabber} : il existe deux diviseurs effectifs linéairement équivalents~$D$ et~$D'$ sur~$\sch X$ tels que~$D'$ soit à support dans~$W(\KK)$ et tels que~$\deg\;(D\cap V(\KK))$ soit premier à~$p$. Soit~$f$ une fonction rationnelle sur~$\sch X$ de diviseur~$D'-D$. La fonction~$f$ ne s'annule pas sur~$V(\KK)$,  et la somme des ordres des pôles de~$f$ sur~$V(\KK)$ est première à~$p$ (en particulier,~$f$ a au moins un pôle sur~$V(\KK)$). 

\medskip
Soit~$\phi$ le morphisme fini et plat de~${\sch X}\an$ dans~$\pkk$ induit par~$f$ ; posons~$z=\phi(x)$. Le diviseur de~$f$ n'étant pas divisible par~$p$, le morphisme~$\phi$ est génériquement étale, et est donc étale en~$x$. L'ouvert~$U:=\phi(V)$ de~$\pkk$ ne contient pas~$0$ ; par conséquent,~$U$ est une composante connexe de~$\pkk\setminus\{z\}$ et~$V$ est une composante connexe de~$\phi\inv(U)$ (lemme~\ref{imtoutoucomp}). L'ouvert~$U$ contient~$\infty$, et le nombre de points comptés avec multiplicité de~$\phi_{|V}\inv(\infty)$ est premier à~$p$ ; par conséquent,~$\deg\;(V\to U)$ est premier à~$p$. 

\medskip
Nous sommes dans les conditions d'application du lemme~\ref{unebrdeg} ;  celui-ci garantit que~$\deg\;(b\to \phi(b))=\deg\;(V\to U)$ et est donc premier à~$p$, ce qu'il fallait démontrer.~$\Box$ 

\subsection*{Sections coronaires : le cas d'un point singulier isolé et «déployé»}

\deux{ptsingisol} Soit~$X$ une courbe~$k$-analytique génériquement quasi-lisse et soit~$x$ un point rigide de~$X$ ; soit~$X'$ la normalisée de~$X$ et soient~$x'_1,\ldots,x'_r
$ les antécédents de~$x$ sur~$X'$. {\em On fait l'hypothèse que les~$x'_i$ sont des~$k$-points lisses ; c'est par exemple toujours le cas si~$k$ est algébriquement clos.}

\trois{predesbrsingisol} Soit~$U\in \arb X x$ tel que~$U\setminus\{x\}$ soit lisse, et soit~$U'$ l'image réciproque de~$U$ sur~$X'$. Il existe une famille~$(V'_i)_{1\leq i\leq r}$ d'ouverts de~$U'$ telle que :

\medskip
$\bullet$ pour tout~$i$, l'ouvert~$V'_i$ est un~$k$-disque contenant~$x_i$ ; 

$\bullet$ les~$V'_i$ sont deux à deux disjoints. 

\medskip
Comme~$U\setminus\{x\}$ est normale, la réunion~$V'$ des~$V'_i$ est un ouvert saturé de~$U'$, qui est donc l'image réciproque d'un voisinage ouvert~$V$ de~$x$ dans~$U$ ; notons que l'image sur~$U$ de chacun des~$V'_i$ est connexe et contient~$x$ ; par conséquent,~$V$ est connexe et est dès lors un arbre. 

\medskip
La normalité de~$U\setminus\{x\}$ implique que~$V'\setminus\{x'_1,\ldots,x'_r\}\to V$ est un isomorphisme. Pour tout~$i$, l'ouvert~$V'_i\setminus\{x'_i\}$ est une~$k$-couronne de type~$]0,*[$ et l'intervalle ouvert~$\mathsf S(V'_i\setminus\{x_i\})$ aboutit proprement à~$x'_i$. 

\medskip
Il s'ensuit que~$V\setminus\{x\}$ s'écrit comme une réunion disjointe~$\coprod  W_i$, où chaque~$W_i$ est une~$k$-couronne de type~$]0,*[$ dont le squelette aboutit proprement à~$x$.

\trois{descbrsingisol} Soit~$b\in \br X x$ et soit~$\sch I$ le sous-ensemble de~$\intera X b$ formé des intervalles~$I$ tels que la section~$I^\flat$ de~$b$ soit une~$k$-couronne de type~$]0,*[$. Il résulte de ce qui précède que~$\sch I\neq \emptyset$. 

\medskip
Soit~$I\in \sch I$. On déduit de~\ref{extremcour} {\em et al.} que le bout de~$I^\flat$ qui correspond à~$x$ est son bout infini. Il s'ensuit que les sections de~$b$ de la forme~$J^\flat$, où~$J\in \inter X x\ctd I$ (qui constituent une base de sections de~$b$) sont encore des~$k$-couronnes de type~$]0,*[$.

\trois{sectinfcornalg} Soit~$X$ une~$\KK$-courbe génériquement quasi-lisse et soit~$x$ un~$\KK$-point de~$X$ ; soit~$b\in \br X x$. On appellera {\em section coronaire} de~$b$ toute section~$Z$ de~$b$ qui est une couronne de type~$]0,*[$ dont l'adhérence dans~$X$ est un arbre compact. Il résulte de ce qui précède : que si~$Z$ est une section coronaire de~$b$, le point~$x$ correspond au bout infini de~$Z$ ; et que~$b$ possède une base de sections coronaires. 

\subsection*{Sections et voisinages coronaires en un point de type 2 ou 3 : le cas d'un corps de base algébriquement clos}

\deux{sectioncour} {\bf Théorème.} {\em Soit~$X$ une courbe~$\KK$-analytique génériquement quasi-lisse et soit~$x$ appartenant à~$X\dtr$. 

\medskip
\begin{itemize}

\item[1)] Si~$x\in  X\typ 3$ il existe un voisinage ouvert~$Z$ de~$x$ dans~$X$ qui est une couronne dont l'adhérence dans~$X$ est un arbre compact, qui est telle que~$x\in \skelan Z$, et qui satisfait de surcroît les propriétés suivantes : 

\medskip
\begin{itemize}
\item [$\bullet$] si~$x\notin \partial \an X$ alors~$Z$ est de type~$]*,*[$ ; 

\item[$\bullet$] si~$x\in \partial \an X$ et si~$x$ n'est pas un point isolé de~$X$ alors~$Z$ est de type~$]*,*]$ et~$\partial \an Z=\{x\}$ ; 

\item[$\bullet$] si~$x$ est un point isolé de~$X$ alors~$Z=\{x\}=\partial\an Z$ (et~$Z$ est de type~$\{*\}$). 

\end{itemize}

\medskip
\item[2)] Si~$b\in \br X x$ il existe une section~$Z$ de~$b$ qui possède les propriétés suivantes :

\medskip
\begin{itemize}
\item[$i)$]~$Z$ est une couronne de type~$]*,*[$ dont l'adhérence dans~$X$ est un arbre compact ;

\item[$ii)$] il existe une fonction~$\theta\in \kappa(x)={\sch O}_{X,x}$ qui est définie sur~$Z$, en est une fonction coordonnée et est telle que~$\left|{\sch O}_X(b)_{\rm sep}\ti\right|_b=|\KK\ti|\oplus |\theta|_b^\ZZ$.

\end{itemize}
\end{itemize}}

\medskip
{\em Démonstration.} On prouve les deux assertions séparément, mais par des méthodes analogues. 

\trois{seccour3l} {\em Preuve de 1) dans le cas où~$x\notin \partial \an X$.} Soit~$b$ l'une des deux branches issues de~$x$ (\ref{casbij3}) ; le théorème~\ref{fmerop} fournit un point~$z$ appartenant à~$\pkk$ et un morphisme fini étale~$\phi : (X,x)\to (\pkk,z)$ entre germes de courbes~$k$-analytiques tel que le degré~$d$ de~$b$ au-dessus de~$\phi(b)$ soit premier à~$p$. En composant à gauche si nécessaire par (le germe d') une homographie, l'on se ramène au cas où~$z=\eta_{\KK,r}$ pour un certain~$r>0~$ n'appartenant pas à~$|\KK\ti|$.

\medskip
Le morphisme~$\phi$ induit une surjection~$\br X x\to \br \pkk {\eta_{\KK,r}}$. Or ces deux ensembles sont de cardinal~$2$ (\ref{casbij3}; notons que pour le second, il n'est pas besoin d'invoquer cette référence -- cela se voit directement) ; par conséquent,~$\phi$ induit une bijection~$\br X x\simeq \br \pkk {\eta_{\KK,r}}$. Il s'ensuit que~$\phi\inv(\phi(b))=\{b\}$ et donc que~$[\kappa(x):\kappa(\eta_{\KK,r})]=\deg\;\phi=d.$ 

\medskip
L'extension séparable~$\kappa(x)/\kappa(\eta_{\KK,r})$ de corps valués henséliens étant de degré premier à~$p$, elle est {\em modérément ramifiée.} On a ~$|\kappa(\eta_r)\ti|=|\KK\ti|\oplus |T(\eta_r)|^\ZZ$ ; comme~$|\KK\ti|$ est divisible et comme~$\red{\kappa(\eta_{\KK,r})}$ est algébriquement clos (il est égal à~$\kk$) la théorie de la ramification modérée assure que~$\kappa(x)$ est isomorphe à~$\kappa(\eta_{\KK,r})[\tau]/(\tau^d-T(\eta_{\KK,r}))$ ; il en résulte que le~$(\pkk,\eta_{\KK,r})$-germe fini étale~$(X,x)$ est isomorphe à~$(\pkk,\eta_{\KK,\sqrt[d]r})$, qui est vu comme~$(\pkk,\eta_{\KK,r})$-germe {\em via} l'élévation à la puissance~$d$ ; il existe donc un voisinage ouvert~$Z$ de~$x$ dans~$X$ qui est une couronne de type~$]*,*[$ telle que~$x\in \mathsf S(Z)$ ; quitte à remplacer~$Z$ par~$I^\flat$, où~$I$ est un intervalle ouvert relativement compact de~$\mathsf S(Z)$ contenant~$x$, on peut supposer que~$\overline Z$ est un arbre compact. Ainsi, 1) est démontré dans le cas où~$x\notin \partial\an X$. 

\trois{seccour3ql}  {\em Preuve de 3) dans le cas général.} Comme~$x$ est de type 3, la courbe génériquement quasi-lisse~$X$ est quasi-lisse en~$x$. Il existe donc un voisinage ouvert~$V$ de~$x$ dans~$X$ et une courbe~$k$-analytique lisse~$W$ tel que~$V$ s'identifie à un domaine analytique fermé de~$W$. En vertu du cas sans bord traité au~\ref{seccour3l} ci-dessus, il existe un voisinage ouvert~$T$ de~$x$ dans~$W$ qui est une couronne de type~$]*,*[$ dont~$x$ appartient au squelette ; notons~$I_1$ et~$I_2$ les deux composantes connexes de~$\mathsf S(Z)\{x\}$. Comme~$x$ est de type 3, il résulte de~\ref{voispkkt3} et de l'assertion~\ref{normsurb} que l'on peut supposer, quitte à restreindre~$T$, que~$T\cap V$ est de la forme~$\{x\}\cup\bigcup\limits_{i\in E} I_i^\flat$, où~$E$ est un sous-ensemble de~$\{1,2\}$ ; on peut également faire en sorte que~$\overline {T\cap V}^V$ soit un arbre compact. Le voisinage~$Z\cap V$ de~$x$ dans~$X$ répond alors aux conditions posées. 

\trois{seccour2} {\em Preuve de 2).} Comme~$X$ est quasi-lisse, il existe un voisinage de~$x$ dans~$X$ qui est isomorphe à un domaine analytique fermé d'une courbe lisse~$X'$ ; grâce à l'assertion iii) de la proposition~\ref{normsurb} et au~\ref{propdomanferm}, on peut remplacer~$X$ par~$X'$, c'est-à-dire supposer que~$X$ est lisse. Si~$x\in X\typ 3$, l'assertion 2) est une conséquence triviale de 1) ; on suppose à partir de maintenant que~$x\in X\typ 2$. Soit~$b\in \br X x$. Le théorème~\ref{fmerop} fournit un point~$z$ appartenant à~$\pkk$ et un morphisme fini étale~$\phi : (X,x)\to (\pkk,z)$ entre germes de courbes~$k$-analytiques tel que le degré~$d$ de~$b$ au-dessus de~$\phi(b)$ soit premier à~$p$. En composant à gauche si nécessaire par (le germe d') une homographie, l'on se ramène au cas où~$z=\eta_{\KK,1}$ et où~$\phi(b)$ est l'unique branche de~$\pkk$ issue de~$\eta_{\KK,1}$ le long de laquelle~$|T|<1$. 

\medskip
On a~$|T(\eta_{\KK,1})|=1$ et~$|T|_{\phi(b)}<1$ ; par conséquent,~$\red {T(\eta_{\KK,1})}$ s'annule au point de la courbe résiduelle en~$\eta_{\KK,1}$ qui correspond à~$\phi(b)$ ; cette courbe est isomorphe à~$ \PP^1_{\kk}$, et le point correspondant ne peut être que l'origine, en laquelle~$\red {T(\eta_{\KK,1})}$ s'annule à l'ordre~$1$. On déduit alors du~\ref{casbij2} et de l'égalité~$|\kappa(\eta_{\KK,1})\ti|=|\KK\ti|$ que le groupe des valeurs de~$|.|_{\phi(b)}$ est ~$|\KK\ti|\oplus |T|_{\phi(b)}^\ZZ$, et que son corps résiduel est ~$\kk$ et est en particulier algébriquement clos. 

\medskip
Par le théorème~\ref{theohenbr},~${\sch O}_X(b)_{\rm sep}$ (resp.~${\sch O}_{\pkk}(\phi(b))_{\rm sep}$) s'identifie naturellement au hensélisé de~$(\kappa(x),|.|_b)$ (resp.~$(\kappa(\eta_1),|.|_{\phi(b)})$. D'après la proposition~\ref{propbondeg}~${\sch O}_X(b)_{\rm sep}$ est de degré~$d$ sur~${\sch O}_{\pkk}(\phi(b))_{\rm sep}$.  L'extension {\em de corps valués henséliens} (rem.~\ref{remextval})~${\sch O}_X(b)_{\rm sep}/{\sch O}_{\pkk}(\phi(b))_{\rm sep}$ étant de degré premier à~$p$, elle est {\em modérément ramifiée}.

\medskip
Compte-tenu du caractère divisible de~$|\KK\ti|$, la théorie de la ramification modérée assure alors que ~${\sch O}_X(b)_{\rm sep}$ est engendré sur~${\sch O}_{\pkk}(\phi(b))_{\rm sep}$ par une racine~$d$-ième~$\tau$ de~$T$, et que l'on a~$\left|{\sch O}_X(b)_{\rm sep}\ti\right|_b=|\KK\ti|\oplus |\tau|_b^\ZZ$. Choisissons~$\theta$ dans~$\kappa(x)$ tel que~$|\theta|_b=|\tau|_b$ (rem.~\ref{remmemeval}) ; cela signifie que~$|\theta/\tau|=1$ le long de~$b$. 

\medskip
Nous pouvons supposer, quitte à restreindre~$X$, que le morphisme de germes~$\phi$ est induit par un « vrai » morphisme fini et plat, noté encore~$\phi$, de~$X$ sur un voisinage ouvert~$U$ de~$\eta_1$, que~$X$ est un arbre, et que~$x$ est le seul antécédent de~$\eta_1$ sur~$X$. 

Les ouverts de~$\pkk$ qui peuvent s'écrire~$I^\flat$, où~$I$ est un intervalle ouvert de la forme~$]\eta_{\KK,s};\eta_{\KK,1}[$ avec~$0<s<1$, forment une base de sections de~$\phi(b)$ ; il existe donc~$I$ de la forme évoquée tel que~$I^\flat\subset U$. 

\medskip
La section~$I^\flat$ de~$b$ est alors une couronne de type~$]*,*[$ contenue dans~$U$, dont~$T$ est une fonction coordonnée. Il existe une (unique) composante connexe~$Z$ de~$\phi\inv(I^\flat)$ qui soit une section de~$b$ ; la flèche~$Z\to I^\flat$ est finie étale de degré~$d$. Comme~$\{Z\times_{I^\flat} J^\flat\}_J$, où~$J$ parcourt~$\inter \pkk {\eta_{\KK,1}} \ctd I$, est une base de sections de~$b$ on peut supposer, quitte à restreindre~$I$ (et~$Z$), que~$\tau$ et~$\theta$ sont définies sur~$Z$, et que~$|\tau/\theta|=1$ identiquement sur~$Z$. 

\medskip
Soit~$Z'$ le revêtement de Kummer de~$I^\flat$ obtenu par adjonction d'une racine~$d$-ième de~$T$ ; la fonction~$\tau$ étant une racine~$d$-ième sur de~$T$ sur~$Z$, elle définit une section du revêtement étale ~$Z\times_{I^\flat}Z'\to Z$, c'est-à-dire une factorisation de~$Z\to I^\flat$ par~$Z'\to I^\flat$ ; par comparaison des degré,~$Z\to Z'$ est un isomorphisme ; par conséquent,~$Z$ est une couronne de type~$]*,*[$, la fonction~$\tau$ en est une fonction coordonnée et~$Z\to I^\flat$ est un revêtement de Kummer. Comme~$|\tau/\theta|$ est identiquement égale à~$1$ sur~$Z$ la fonction~$\theta$ est elle aussi une fonction coordonnée de~$Z$ ; en restreignant~$Z$ (il suffit de la remplacer par n'importe laquelle de ses sous-couronnes ouvertes strictes aboutissant à~$x$),  on peut faire en sorte que son adhérence dans~$X$ soit un arbre compact, ce qui achève la démonstration.~$\Box$

\deux{rembasecour} Faisons quelques commentaires à propos du théorème ci-dessus, en en conservant les notations. 

\trois{basevoist3d} Supposons que~$x\in X\typ3$ ; un voisinage ouvert~$Z$ de~$x$ satisfaisant les conditions de l'assertion 1) sera dit {\em coronaire}. Soit~$Z$ un voisinage coronaire de~$x$ ; les sous-couronnes de~$Z$ de la forme~$I^\flat$, où~$I$ est un voisinage ouvert de~$x$ dans~$\skelan Z$, sont encore des voisinages coronaires de~$x$ ; elles forment une base de voisinages de~$x$ dans~$Z$, et {\em a fortiori} dans~$X$ (\ref{voispkkt3}). Si~$f$ est une fonction analytique inversible sur~$Z$, alors~$f$ est une fonction coordonnée si et seulement si~$|f(x)|$ engendre~$|\hres(x)\ti|/|\KK\ti|$ : compte-tenu du fait que~$x$ est de type 3, cela découle de~\ref{oxmodulo} et~\ref{souscourtyp}. 

\trois{souscourmeme}
On ne suppose plus que~$x\in  X\typ3$. Nous dirons d'une section~$Z$ de~$b)$ qui satisfait la condition  i) de l'assertion 2) du théorème qu'elle est {\em coronaire} ; si~$Z$ est une section coronaire de~$b$ (resp. une section coronaire de~$b$ satisfaisant ii) ) et si~$Z'$ est une sous-couronne ouverte de~$Z$ aboutissant à~$x$, il est immédiat que~$Z'$ est coronaire (resp. est coronaire et satisfait ii) ).

\medskip
Par ailleurs, si~$Z$ est une section coronaire de~$b$, l'ensemble de ses sous-couronnes ouvertes aboutissant à~$x$ est exactement l'ensemble des~$I^\flat$, où~$I$ parcourt~$\inter X x \ctd {\mathsf S(Z)}$ ; c'est donc une base de sections de~$b$. 

\medskip
On déduit de ce qui précède que~$b$ possède une base de sections qui sont coronaires et satisfont ii). 

\trois{basecoron} Soit~$Z$ une section coronaire de~$b$. Par ce qui précède, il existe une section~$Z'$ de~$b$ qui est coronaire, contenue dans~$Z$ et satisfait ii). L'intersection~$\mathsf S(Z)\cap \mathsf S(Z')$ est un intervalle ouvert~$I$ aboutissant à~$x$. Comme~$I\subset \mathsf S(Z)\cap \mathsf S(Z')$, l'intervalle~$I$ est faiblement admissible dans~$Z$ aussi bien que dans~$Z'$, et~$I^\flat$ est une sous-couronne ouverte de~$Z$ aussi bien que de~$Z'$, qui aboutit à~$x$ puisque~$I$ aboutit à~$x$. En tant que sous-couronne de~$Z$ aboutissant à~$x$, l'ouvert~$I^\flat$ est une section coronaire de~$b$ ; comme la section~$I^\flat$ est aussi une sous-couronne de~$Z'$, elle satisfait ii). 

\medskip
Ainsi,~$Z$ possède une sous-couronne ouverte aboutissant à~$x$ et satisfaisant ii).

\trois{remcalcval} Soit~$f\in {\sch O}_X(b)_{\rm sep}\ti$. En vertu de~\ref{souscourmeme}, il existe une section coronaire~$Z$ de~$b$ satisfaisant ii) et sur laquelle~$f$ et~$f\inv$ sont définies. Soit~$\theta$ comme dans l'assertion 2), ii) ; comme~$f$ est une fonction inversible sur la couronne~$Z$ dont~$\theta$ est une fonction coordonnée, il existe~$\alpha\in \KK\ti$ et~$m\in \ZZ$ tel que~$|f|=|\alpha|.|\theta|^m$ identiquement sur~$Z$ , et~$f$ est une fonction coordonnée de~$Z$ si et seulement si~$m=1$ ou~$m=-1$ ; par définition même de~$|.|_b$, on a~$|f|_b=|\alpha|.|\theta|_b^m$. Il en résulte :

\begin{itemize}

\medskip
\item[$\alpha)$] qu'il existe un isomorphisme naturel~$|{\sch O}_X(Z)\ti|\simeq \left|{\sch O}_X(b)_{\rm sep}\ti\right|_b$ compatible aux plongements de~$|\KK\ti|$ dans ces deux groupes ; 

\item[$\beta)$] que~$f$ est une fonction coordonnée de la couronne~$Z$ si et seulement si~$\left|{\sch O}_X(b)_{\rm sep}\ti\right|_b=|\KK\ti|\oplus |f|_b^\ZZ$, ou encore si et seulement si~$|f|_b$ engendre~$\left|{\sch O}_X(b)_{\rm sep}\ti\right|_b$ modulo~$|\KK\ti|$. 
\end{itemize}

\medskip
Supposons que~$f$ soit une fonction coordonnée de~$Z$. La fonction~$|f|$ a une limite~$\lambda$ en~$x$, et de la définition de~$|.|_b$ découle la remarque suivante, qui permet de décrire entièrement~$\left|{\sch O}_X(b)_{\rm sep}\ti\right|_b=|\KK\ti|\oplus |f|_b^\ZZ$ comme groupe {\em ordonné} : si~$|f|$ est croissante (resp. décroissante) lorsqu'on oriente le squelette de~$Z$ dans la direction de~$x$, alors~$|f|_b$ est infiniment proche inférieurement (resp. supérieurement) de~$\lambda$. 

Notons que si~$\lambda\notin|\KK\ti|$ (on déduit aisément de ce qui précède que cela se produit si et seulement si~$x\in X\typ 3$), on obtient le même ordre sur~$|\KK\ti|\oplus |f|_b^\ZZ$ en décrétant que~$|f|_b$ est {\em égal} à~$\lambda$ ; on retrouve ainsi l'égalité attendue lorsque~$x\in X\typ 3$ entre~$|.|_b$ et~$|.|$ sur~$\kappa(x)={\sch O}_X(b)_{\rm sep}$ (rem.~\ref{remvalbt3} et assertion iv) du th.~\ref{theohenbr}).

\medskip
{\em Remarque.} L'isomorphisme naturel évoqué en~$\alpha)$ s'étend au cas où~$Z$ est coronaire sans satisfaire nécessairement ii) : en effet,~$Z$ contient  d'après~\ref{basecoron} une sous-couronne~$Z'$ qui aboutit à~$x$ et satisfait ii), et l'on conclut en utilisant le fait que la 
flèche naturelle~$|{\sch O}_X(Z)\ti|\to |{\sch O}_X(Z')\ti|$ est un isomorphisme.

\subsection*{Sections coronaires : le cas d'un point singulier isolé quelconque}

\deux{brsingisolkg} {\bf Lemme.} {\em Soit~$X$ une courbe~$k$-analytique génériquement quasi-lisse, soit~$x$ un point rigide de~$X$ et soit~$b$ une branche issue de~$X$. Il existe une section~$Z$ de~$b$ qui est une couronne virtuelle de type~$]0,*[$ dont l'adhérence dans~$X$ est un arbre compact.}

\medskip
{\em Démonstration.} Comme~$\sch O_{X,x}$ est hensélien, on peut restreindre~$X$ de sorte que~$\got s(x)\subset \sch O_X(X)$ ; on peut dès lors, quitte à remplacer~$k$ par~$\got s(x)$, supposer que~$\got s(x)=k$. Le point~$x$ admet sous cette hypothèse un unique antécédent~$y$ sur~$X_{\KK}$, lequel est fixe sous l'action de~$\mathsf G$ ; notons~$p$ la flèche~$X_{\KK}\to X$, qui s'identifie à~$X_{\KK}\to X_{\KK}/\mathsf G$.

\medskip
Choisissons une branche~$\beta$ de~$\br {X_{\KK}} y~$ située au-dessus de~$b$ ; comme~$y$ est rigide,~$\br {X_{\KK} y }$  est fini et le stabilisateur~$\mathsf H$ de~$\beta$ dans~$\mathsf G$ est donc ouvert. 

D'après~\ref{descbrsingisol}, il existe~$I\in  \intera {X_{\KK}} \beta$ tel que la~$\KK$-courbe~$I^\flat$ soit une couronne de type~$]0,*[$. En vertu de la proposition~\ref{corollstabbr}, il existe un intervalle ouvert~$J$ appartenant à~$\interac X x \ctd I$ tel que les propriétés suivantes soient satisfaites : 

\medskip

-~$p$ induit un homéomorphisme~$\overline J\simeq p(\overline J)$ ; 

- l'intervalle~$J$ est fixe point par point sous~$\mathsf H$ ; 

- pour tout~$g\in \mathsf G-\mathsf H$ on a~$g(J^\flat)\cap J^\flat=\emptyset$ ; 

\medskip
Remarquons que comme~$J$ aboutit à~$y$, l'ouvert~$J^\flat$ est lui aussi une couronne de type~$]0,*[$. Soit~$Z$ l'image de~$J^\flat$ sur~$X$ ; toujours d'après la proposition~\ref{corollstabbr}, l'intervalle ouvert~$p(J)$ de~$X$ aboutit proprement à~$x$, est faiblement admissible, et vérifie l'égalité~$Z=p(J^\flat)$ ; par conséquent,~$Z$ est un arbre à deux bouts dont l'adhérence dans~$X$ est un arbre compact. 

\medskip
Par construction,~$Z_{\KK}$ est la réunion disjointe des~$\mathsf g(J^\flat)$, pour~$g$ parcourant un système de représentants de~$\mathsf G/\mathsf H$ ; chacune des~$g(J^\flat)$ étant une couronne de type~$]0,*[$, il s'ensuit que~$Z$ est une~$\infty$-couronne virtuelle de type~$]0,*[$.~$\Box$ 

\deux{defsecinfcoron} Soit~$X$ une courbe~$k$-analytique génériquement quasi-lisse,  soit~$x$ un point rigide de~$X$ et soit~$b\in \br X x$. On appellera {\em section coronaire} de~$b$ toute section~$Z$ de~$b$ qui est une couronne virtuelle de type~$]0,*[$ dont l'adhérence dans~$X$ est un arbre compact ; cette définition est compatible avec la précédente lorsque~$\KK$ est algébriquement clos. Si~$Z$ est une section coronaire de~$b$, il en va de même de~$I^\flat$ pour tout intervalle ouvert de~$\mathsf S(Z)$ aboutissant à~$x$ ; par conséquent, comme le lemme~\ref{brsingisolkg} assure l'existence d'une section coronaire de~$b$, la branche~$b$ possède une base de telles sections. 

\medskip
Si~$Z$ est une section coronaire de~$b$, son bord est de la forme~$\{x,y\}$ avec~$y\neq x$. Il résulte de~\ref{extremcour} {\em et sq.} que le bout infini de~$Z$ est celui qui correspond à~$x$, et que~$y\in X\dtr$. 

\subsection*{Sections et voisinages coronaires en un point de type 2 ou 3 : le cas d'un corps de base quelconque} 

\deux{brcourvirt} {\bf Proposition.} {\em Soit~$X$ une courbe~$k$-analytique quasi-lisse et soit~$x\in X$. 

\begin{itemize}

\medskip
\item[1)] Si~$x\in  X\typ 3$ il existe un voisinage ouvert~$Z$ de~$x$ dans~$X$ qui est une couronne gentiment virtuelle dont l'adhérence dans~$X$ est un arbre compact, qui est telle que~$x\in \skelan Z$, et qui satisfait de surcroît les propriétés suivantes : 

\medskip
\begin{itemize}
\item [$\bullet$] si~$x\notin \partial \an X$ alors~$Z$ est de type~$]*,*[$ ; 

\item[$\bullet$] si~$x\in \partial \an X$ et si~$x$ n'est pas un point isolé de~$X$ alors~$Z$ est de type~$]*,*]$ et~$\partial \an Z=\{x\}$ ; 

\item[$\bullet$] si~$x$ est un point isolé de~$X$ alors~$Z=\{x\}=\partial\an Z$ (et~$Z$ est de type~$\{*\}$). 

\end{itemize}

\medskip
\item[2)] Si~$b\in \br X x$ il existe une section~$Z$ qui est une couronne gentiment virtuelle de type~$]*,*[$ dont l'adhérence dans~$X$ est un arbre compact. On peut de surcroît imposer à~$Z$ de vérifier la condition suivante : il existe une extension finie~$F$ de~$k$, et même finie séparable si~$|k\ti|\neq 1$, un antécédent~$x'$ de~$x$ sur~$X_F$, une composante connexe~$Z'$ de~$Z_F$ aboutissant à~$x'$ qui est une~$F$-couronne, et une fonction coordonnée sur~$Z'$ se prolongeant à un voisinage de~$x'$ dans~$X_F$. 
\end{itemize}}

\medskip
{\em Démonstration.} Comme~$\sch O_{X,x}=\kappa(x)$, on peut restreindre~$X$ de sorte que~$\got s(x)\subset \sch O_X(X)$ ; on peut dès lors, quitte à remplacer~$k$ par~$\got s(x)$, supposer que~$\got s(x)=k$. Le point~$x$ admet sous cette hypothèse un unique antécédent~$y$ sur~$X_{\KK}$, lequel est fixe sous l'action de~$\mathsf G$ ; notons~$p$ la flèche~$X_{\KK}\to X$, qui s'identifie à~$X_{\KK}\to X_{\KK}/\mathsf G$. 

\trois{courvirtt3}{\em Preuve de 1)}. Le point~$y$ étant lisse de type 3), il possède un voisinage coronaire~$Y$ dans~$X_{\KK}$ (th.~\ref{sectioncour} et~\ref{basevoist3d}). Le squelette analytique~$\skelan Y$ peut être vue comme une étoile de sommet~$y$ tracée sur~$X_{\KK}$, faiblement admissible,  dont la valence coïncide avec le cardinal de~$\br Y y$. On peut par conséquent appliquer la proposition~\ref{corollstabet} ; elle assure entre autres l'existence d'un intervalle~$I$ tracé sur~$\skelan Y$, ouvert dans ce dernier, contenant~$y$ et stable sous~$\mathsf G$. Quitte à choisir un tel~$I$ et à remplacer~$Y$ par sa sous-couronne~$I^\flat$, on peut supposer que~$Y$ elle-même est stable sous~$\mathsf G$ ; on note~$Z$ le quotient~$Y/\mathsf G$ ; c'est un voisinage ouvert de~$x$. Comme~$Y$ est un arbre ayant au plus deux bouts,~$Z$ est un arbre ayant au plus deux bouts ; son adhérence dans~$X$ est égale au quotient de~$\overline Y$ par~$\mathsf G$, et est donc un arbre compact.

\medskip
Nous allons montrer que~$\mathsf G$ agit par automorphismes directs sur~$Y$, en distinguant trois cas. 

Si~$Y$ est de type~$\{*\}$, autrement dit si~$Y=\{y\}$ cela découle du fait que~$\mathsf G$ opère par isométries sur~$\hres(y)$, et agit donc trivialement sur~$|\hres(y)\ti|/|k\ti|$. 

Si~$Y$ est de type~$]*,*]$, cela provient du fait que comme~$\skelan Y$ est un intervalle semi-ouvert,~$\mathsf G$ en préserve nécessairement les orientations.

Supposons maintenant que~$Y$ soit de type~$]*,*[$. C'est le cas où~$y\notin Y\an$, et donc où~$x\notin X\an$. Comme~$x$ est de type 3, il y a exactement deux branches de~$X$ issues de~$x$ (\ref{casbij3}), et~$Z\setminus\{x\}$ a de ce fait deux composantes connexes. Par conséquent, les images sur~$Z$ des deux composantes connexes de~$Y\setminus\{y\}$ sont disjointes, ce qui signifie que ces composantes ne sont pas échangées par~$\mathsf G$ ; autrement dit, l'action de ce dernier ne permute pas les orientations du squelette de~$Y$, ce qu'on souhaitait établir. 

\medskip
Soit~$r\in\RR\ti_+$ engendrant~$|\kappa(y)\ti|$ modulo~$|\KK\ti|$. La fermeture algébrique (et même séparable si~$|k\ti|\neq \{1\}$ ) de~$k$ dans~$\KK$ en est un sous-corps dense. Il existe par conséquent un voisinage ouvert~$U$ de~$x$ dans~$X$, une extension finie~$L$ de~$k$ (que l'on peut prendre séparable si~$|k\ti|\neq \{1\}$) et une fonction analytique inversible~$f$ sur~$U_L$ telle que~$|f(y)|=r$.

\medskip
Il découle de~\ref{basevoist3d}  qu'il est loisible de restreindre~$Y$ de sorte que~$Y\subset U_{\KK}$, et que~$f_{|Y}$ est alors une fonction coordonnée sur~$Y$.

\medskip
\medskip
On peut voir~$f$ comme une fonction analytique sur~$Z_L$, et elle induit par construction un isomorphisme entre~$(Z_L)\hotimes_L{\KK}$ et un ouvert de~$\Aff^{1,\rm an}_{\KK}$ défini par une condition de la forme~$R<|T|<R'$ ; il s'ensuit par descente (corollaire~\ref{desciso}) que~$f$ induit un isomorphisme  entre~$Z_L$ et l'ouvert de~$\Aff^{1,\rm an}_L$ défini par la condition~$R<|T|<R'$. Par conséquent,~$Z$ est une couronne {\em gentiment} virtuelle. Compte-tenu du fait que le type de~$Z$ est celui de~$Y$, que~$x\in \partial\an X$ si et seulement si~$y\in \partial \an Y$, et que~$x$ est isolé dans~$X$ si et seulement si~$y$ est isolé dans~$Y$, ceci achève la démonstration de i). 

\trois{courvirtt2}{\em Preuve de 2).} Si~$x\in X\typ 3$ alors 2) est une conséquence triviale de 1) ; supposons maintenant que~$x\in X\typ 2$, et notons~$\sch C$ la courbe résiduelle en~$y$ ; soit~$b$ une branche de~$X$ issue de~$x$. 

\medskip
 Choisissons une branche~$\beta$ de~$X_{\KK}$ issue de~$y$ et située au-dessus de~$b$ ; soit~$\mathsf H$ son stabilisateur ; c'est un sous-groupe ouvert de~$\mathsf G$ (prop.~\ref{corollstabbr}). Choisissons une section coronaire~$Y$ de~$\beta$. En vertu de la proposition~\ref{corollstabbr}, il existe un intervalle ouvert~$J$ appartenant à~$\interac X x \ctd {\skel Y}$ tel que les propriétés suivantes soient satisfaites : 

\medskip

-~$p$ induit un homéomorphisme~$\overline J\simeq p(\overline J)$ ; 

- l'intervalle~$J$ est fixe point par point sous~$\mathsf H$ ; 

- pour tout~$g\in \mathsf G-\mathsf H$ on a~$g(J^\flat)\cap J^\flat=\emptyset$ ; 

\medskip
Remarquons que~$J^\flat$ est une sous-couronne ouverte de~$Y$ aboutissant à~$y$, et est donc elle aussi une section coronaire de~$\beta$. Soit~$Z$ l'image de~$J^\flat$ sur~$X$ ; toujours d'après la proposition~\ref{corollstabbr}, l'intervalle ouvert~$p(J)$ de~$X$ aboutit proprement à~$x$, est faiblement admissible, et vérifie l'égalité~$Z=p(J^\flat)$ ; par conséquent,~$Z$ est un arbre à deux bouts dont l'adhérence dans~$X$ est un arbre compact. 

\medskip
Par construction,~$Z_{\KK}$ est la réunion disjointe des couronnes~$\mathsf g(J^\flat)$, pour~$g$ parcourant un système de représentants de~$\mathsf G/\mathsf H$ ; de ce fait,~$Z$ est une couronne virtuelle.

\medskip
La branche~$\beta$ correspond à un point fermé~$\sch P$ de~$\sch C$. La fermeture algébrique (et même séparable si~$|k\ti|\neq \{1\}$ ) de~$k$ dans~$\KK$ en est un sous-corps dense. Il existe par conséquent un voisinage ouvert~$U$ de~$x$ dans~$X$, une extension finie~$L$ de~$k$ (que l'on peut prendre séparable si~$|k\ti|\neq \{1\}$) et une fonction analytique~$f$ sur~$U_L$ telle que~$|f(y)|=1$ et telle que~$\red {f(y)}$ ait un zéro d'ordre~$1$ en~$\sch P$ (on note encore~$f$ la fonction analytique sur~$U_{\KK}$ induite par~$f$).  

\medskip
Il découle de~\ref{remcalcval} et~\ref{casbij2} qu'il est loisible de restreindre~$J$ de sorte que~$J^\flat\subset U_{\KK}$, et que~$f_{|J^\flat}$ soit une fonction coordonnée sur~$J^\flat$.

\medskip
Il existe une extension finie séparable~$F$ de~$L$ qui déploie, en tant qu'extension finie de~$k$, la~$k$-algèbre étale~$\got s(Z)$ ; toute composante connexe de~$Z_F$ est alors un espace~$F$-analytique géométriquement connexe ; il existe par conséquent une composante connexe~$Z'$ de~$Z_F$ et un plongement~$F\hookrightarrow \KK$ tels que~$Z'_{\KK}\simeq J^\flat$.

\medskip
On peut voir~$f$ comme une fonction analytique sur~$Z'$, et elle induit par construction un isomorphisme entre~$Z'_{\KK}$ et un ouvert de~$\Aff^{1,\rm an}_{\KK}$ défini par une condition de la forme~$R<|T|<R'$ ; il s'ensuit par descente (corollaire~\ref{desciso}) que~$f$ induit un isomorphisme  entre~$Z'$ et l'ouvert de~$\Aff^{1,\rm an}_F$ défini par la condition~$R<|T|<R'$. Par conséquent,~$Z'$ est une~$F$-couronne dont~$f$ est une fonction coordonnée, ce qui achève de prouver 2).~$\Box$

\deux{rembasecourkgen} Faisons quelques commentaires à propos du lemme ci-dessus, en en conservant les notations. 

\trois{basevoist3nondep} Supposons que~$x\in X\typ3$.  Un voisinage ouvert~$Z$ de~$x$ satisfaisant les conditions de l'assertion i) sera dit {\em coronaire} ;  cette terminologie est compatible avec celle déjà introduite dans le cas algébriquement clos. Si~$Z$ est un voisinage coronaire de~$x$, les sous-couronnes virtuelles de~$Z$ de la forme~$I^\flat$, où~$I$ est un intervalle de~$\skelan Z$ ouvert dans ce dernier et contenant~$x$, sont encore des voisinages coronaires de~$x$ ; elles forment une base de voisinages de~$x$ dans~$Z$, et {\em a fortiori} dans~$X$ ; cela découle, en vertu de la proposition~\ref{corollstabet}, du fait que~$\skelan Z$ peut être vu comme une étoile étoile tracée sur~$Z$ de sommet~$x$ et de valence égale au cardinal de~$\br X x$. 

\trois{comseccorpasalg} On ne suppose plus que~$x\in X\typ3$. Une section~$Z$ de la branche~$b$ qui satisfait la condition énoncée en ii) sera qualifiée de {\em coronaire} ; cette terminologie est compatible avec celle déjà introduite dans le cas algébriquement clos. 

\trois{basecoron-pasalg} Soit~$I\in \inter X b$ ; choisissons une section coronaire~$Z$ de~$b$. Comme~$\mathsf S(Z)$ et~$I$ définissent tous deux~$b$, leur intersection contient un intervalle ouvert~$J_0$ aboutissant à~$x$ ; si~$J$ est un intervalle ouvert aboutissant à~$x$ et contenu dans~$J_0$ il appartient à~$ \interac X x$, et~$J^\flat$ est une sous-couronne virtuelle ouverte de~$Z$, et partant une section coronaire de~$b$. Ainsi~$b$ possède-t-elle une base de sections coronaires de la forme~$J^\flat$, où~$J\in \interac X x \ctd I$. 

\trois{orientebranche} Supposons que~$b$ possède une section coronaire~$Z_0$ qui est une couronne ; 
elle possède alors une base de telles 
sections (par exemple, les sous-couronnes ouvertes de~$Z_0$ aboutissant
à~$x$). Soient~$Z_1$ et~$Z_2$ deux sections coronaires de~$b$ qui sont des
couronnes ; si~$Z_2$ est une sous-couronne de~$Z_1$ 
alors~$\sch Z(Z_1)\simeq \sch Z(Z_2)$. 
La famille~$\{\sch Z(Z)\}_Z$, où~$Z$ parcourt 
l'ensemble des sections coronaires de~$b$ qui sont des couronnes, 
apparaît ainsi comme un système filtrant dont toutes
 les flèches de transition sont des isomorphismes. On note~$\sch Z(b)$ sa limite (à la fois inductive et projective,
 selon le sens dans lequel on considère les isomorphismes) ; c'est un groupe
 libre de rang~$1$. Pour toute section coronaire~$Z$ de~$b$ qui est une couronne, 
 on dispose d'un isomorphisme naturel~$\sch Z(Z)\simeq \sch Z(b)$ ; on
 obtient par ce biais, 
 en vertu de~\ref{modcour}, une bijection entre l'ensemble des générateurs
 de~$\sch Z(b)$ et celui des orientations de la branche~$b$ (\ref{orbranchegr}).

 Lorsque~$k$ est algébriquement clos, on déduit de
 l'assertion~$\alpha)$ de~\ref{remcalcval} que~$\sch Z(b)$ est naturellement
 isomorphe à~$\left|\sch O_X(b)\ti_{\rm sep}\right|/|k\ti|$, et donc
 également à~$|\kappa(x)\ti|_b/|k\ti|$ (th.~\ref{theohenbr}).

\subsection*{Degré d'une branche sur la branche image et ramification résiduelle}

\deux{degetseccor} {\bf Proposition.} {\em Soit~$\phi:Y\to X$ un morphisme purement de dimension relative nulle entre courbes~$k$-analytiques génériquement quasi-lisses, et soit~$y$ un point de~$Y\geom$ ; posons~$x=\phi(y)$. 

\medskip
\begin{itemize}
\item[i)] Soit~$b\in \br Y y$ et soit~$a$ son image dans~$\br X x$. Il existe une section~$Z'$ de~$b$ et une section~$Z$ de~$a$ telles que : 

\medskip
\begin{itemize} 
\item[$\bullet$]~$Z'$ et~$Z$ sont coronaires  ; 

\item[$\bullet$]~$Z'$ est une composante connexe de~$\phi\inv(Z)$ ; 

\item[$\bullet$]~$Z'\to Z$ est fini et plat de degré~$\deg\;(b\to a)$. 
\end{itemize}

\medskip
\item[ii)] Supposons que~$y\in Y\typ 3$. Il existe un voisinage ouvert~$Z'$ de~$y$ dans~$Y$ et un voisinage ouvert~$Z$ de~$x$ dans~$X$ tels que : 

\medskip
\begin{itemize} 
\item[$\bullet$]~$Z'$ et~$Z$ sont coronaires ; 

\item[$\bullet$]~$Z'$ est une composante connexe de~$\phi\inv(Z)$ ; 

\item[$\bullet$] il existe un intervalle~$I$ de~$\skelan Z$ contenant~$x$, tel que~$\partial_Z I\subset \{x\}$, et tel que~$Z'\to Z$ induise un morphisme fini et plat~$Z'\to I^\sharp$. 
\end{itemize}

\end{itemize} }

\medskip
{\em Démonstration.} On prouve les deux assertions séparément, mais avec des arguments similaires.

\trois{casseccor} {\em Preuve de i)}. Commençons par remarquer que la proposition~\ref{phipasfin} permet de supposer, quitte à restreindre~$Y$ et~$X$, que tous deux sont des arbres, que~$y$ est le seul antécédent de~$x$ sur~$Y$, que~$Y\to X$ est compacte à fibre finie, et que le morphisme~$Y\setminus\{y\}\to X\setminus\{x\}$ est fini et plat. 
Soit~$V$ une section coronaire de~$a$ et soit~$V'$ la composante connexe de~$\phi\inv(V)$ qui correspond à~$b$. Il existe une section coronaire~$Z'$ de~$b$ contenue dans~$V'$ ; posons~$Z=\phi(Z')$. 

\medskip
Comme~$Z'$ est connexe et comme~$b\in \partial Z'$, l'ouvert~$Z$ de~$V$ est connexe et~$a\in \partial Z$ ; cela entraîne que~$Z\cap \mathsf S(V)\neq \emptyset$ ; le lemme~\ref{morcourvirt} assure alors que~$Z$ est de la forme~$I^\flat$ pour un certain intervalle ouvert non vide~$I$ de~$\mathsf S(V)$ et que~$Z'\to Z$ est fini et plat. Comme~$a\in \partial Z$, l'intervalle~$I$ aboutit nécessairement à~$x$ ; par conséquent,~$Z$ est une section coronaire de~$a$.

\medskip
La flèche finie (et plate)~$Z'\to Z$ est compacte ; il résulte alors du lemme~\ref{ouvcompphims} que~$Z'$ est une composante connexe de~$\phi\inv (Z)$. Comme~$Z'\subset V'$ on a nécessairement~$Z=V'\times_V Z$ ; il en résulte que~$\deg\;(Z'\to Z) =\deg\;(b\to a)~$.

\medskip
\trois{casvoiscor} {\em Preuve de ii).} Soit~$Z$ un voisinage coronaire de~$x$ dans~$X$ et soit~$Z'$ un voisinage coronaire de~$y$ dans~$\phi^{-1}(Z)$. Comme~$x\in \skelan Z$, le lemme~\ref{morcourvirt} assure~$Z'\to V$ induit un morphisme fini et plat de~$Z'$ sur~$I^\sharp$ pour un certain intervalle non vide~$I$ de~$\skelan Z$, qui contient évidemment~$x$ ; quitte à restreindre~$Z$, on peut supposer que~$\partial_Z I\subset \{x\}$. 

\medskip
La couronne virtuelle~$Z'$ est une partie ouverte, connexe et non vide de~$\phi\inv(Z)$. La flèche finie (et plate)~$Z'\to I^\sharp$ est compacte, et~$I^\sharp$ est fermée dans~$Z$  puisque~$x\in I$ et~$\partial_Z I\subset \{x\}$ ; par conséquent,~$I^\sharp\hookrightarrow Z$ est compacte, et~$Z'\to Z$ l'est aussi par composition. Le lemme~\ref{ouvcompphims} garantit alors que~$Z'$ est une composante connexe de~$\phi\inv (Z)$.~$\Box$

\deux{commcourt3} {\em Commentaires.} Dans la situation de l'assertion 2) (et en en conservant les notations),  on se trouve plus précisément dans l'un des cas suivants. 

\trois{typouvouv} {\em Le cas où~$Z'$ est de type~$]*,*[$, autrement dit où~$y\notin \partial\an Y$.} La couronne virtuelle~$I^\sharp$ est alors de type~$]*,*[$ (d'après le lemme~\ref{morcourvirt}), ce qui force~$Z$ à être égale à~$I^\sharp$ ; le morphisme~$\phi$ est alors sans bord en~$y$. 

\trois{typouvferm} {\em Le cas où~$Z'$ est de type~$]*,*]$, autrement dit où~$y\in \partial\an Y$ et où~$y$ n'est pas isolé dans~$Y$}. La couronne virtuelle~$I^\sharp$ est alors de type~$]*,*]$, ce qui laisse deux possibilités pour~$Z$ : ou bien~$Z=I^\sharp$, et~$\phi$ est sans bord en~$y$ ; ou bien~$Z$ est de type~$]*,*[$, et~$y\in \partial\an \phi$. 

\trois{typsingle} {\em Le cas où~$Z'$ est de type~$\{*\}$, autrement dit où~$y$ est un point isolé de~$Y$, et où~$Z'=\{y\}$}.  La couronne virtuelle~$I^\sharp$ est alors de type~$\{*\}$, ce qui veut dire que~$I=I^\sharp=\{x\}$ ;  cela laisse trois possibilités pour~$Z$ : ou bien~$Z=I^\sharp=\{x\}$, et~$\phi$ est sans bord en~$y$ ; ou bien~$Z$ est de type~$]*,*]$, et~$y\in \partial\an \phi$ ; ou bien~$Z$ est de type~$]*,*[$, et~$y\in \partial\an \phi$. 

\trois{commentdeg} Dans chacun des cas où~$\phi$ est sans bord en~$y$, c'est-à-dire où~$Z=I^\sharp$, le morphisme~$\phi$ est fini et plat en~$y$ ; comme~$y$ est le seul antécédent de~$x$ sur~$Z'$ (puisque~$\skelan {Z'}$ est égale à l'image réciproque de~$\skelan Z$ sur~$Z'$, et puisque~$\skelan {Z'}\to \skelan Z$ est un homéomorphisme), on a l'égalité~$$\deg^y \;\phi=\deg \;(Z'\to Z).$$

\deux{ramresi} {\bf Théorème.} {\em  Soit~$\phi : Y\to X$ un morphisme 
de dimension relative nulle entre deux courbes~$\KK$-analytiques quasi-lisses ; soit~$y\in Y\dtr$ et soit~$x$ son image sur~$X$. 

\medskip
i) Supposons que~$y\in Y\typ 2$, soit~$b\in \br Y y$ et soit~$a$ la branche~$\phi(b)$ ; posons~$d=\deg\;(b\to a)$. Notons~${\sch D}$ (resp.~${\sch C}$) la courbe résiduelle en~$y$ (resp.~$x$) et~$\sch Q$ (resp.~$\sch P$) le point fermé qui correspond à~$b$ (resp.~$a$) par la bijection du~\ref{casbij2} ; le point~$\sch P$ est l'image de~$\sch Q$ (rem.~\ref{remtransnat}) et l'on note~$e$ l'indice de ramification de~${\sch D}\to {\sch C}$ en~$\sch Q$. On a alors~$e=d$. 

\medskip
ii) Supposons que~$y\notin \partial\an \phi$ ; on alors ~$[\red{\kappa(y)}:\red{\kappa(x)}]=[\kappa(y):\kappa(x)]$ si~$y\in Y\typ 2$, et~$[|\kappa(y)\ti|:|\kappa(x)\ti|]=[\kappa(y):\kappa(x)]$ si~$y\in Y\typ 3$ ; par conséquent, l'extension~$\kappa(y)/\kappa(x)$ est sans défaut.

}

\medskip
{\em Démonstration.} Nous allons prouver les deux assertions séparément. 

\trois{ramresi2} {\em Preuve de i)}. La proposition~\ref{degetseccor} assure l'existence d'une section coronaire~$Z$ de~$a$ et d'une section coronaire~$Z'$ de~$b$ telle que~$\phi$ induise une flèche finie et plate de degré~$d$ de~$Z'$ vers~$Z$. D'après~\ref{morcourintersk} (et~\ref{oxmodulo}), la flèche~$\phi$ identifie~$|{\sch O}_X(Z)\ti|$ à un sous-groupe d'indice~$d$ de~$|{\sch O}_Y(Z')\ti|$. Il s'ensuit, en vertu du~\ref{remcalcval}, que~$d$ est l'indice de ramification de~$|.|_b$ sur~$|.|_a$. La valuation~$|.|_b$ (resp.~$|.|_a$) est la composée de~$|.|$ et de la valuation discrète d'anneau~$\sch O_{D, Q}$ (resp.~$\sch O_{C,P}$) ; comme~$|\kappa(y)\ti|=|\kappa(x)\ti|=|\KK\ti|$, l'indice de ramification de~$|.|_b$ sur~$|.|_a$ est égal à l'indice de ramification de~$\sch O_{D, Q}$ sur~$\sch O_{C,P}$. Ainsi,~$d=e$, d'où i).

\medskip
\trois{ramresi3} {\em Preuve de ii)}. Supposons tout d'abord que~$y\in Y\typ 2$. 
Le point~$x$ est alors également de type 2, et l'on note encore~$\sch C$ et~$\sch D$ les courbes
résiduelles respectives en~$x$ et~$y$. On identifie~$\red{(X,x)}$ (resp.~$\red{(Y,y)}$) à un ouvert dense de~$\sch C$ (resp.~$\sch D$). Choisissons
un point fermé~$\sch P\in \red{(X,x)}\subset \sch C$ ; il correspond à une branche~$a\in \br X x$. Notons~${\sch Q}_1,\ldots,{\sch Q}_r$ les antécédents de~${\sch P}$ sur~${\sch D}$ ; 
comme~$\phi$ est sans bord en~$y$, l'image réciproque de~$\red{(X,x)}$ sur~$\sch D$ est égale à~$\red{(Y,y)}$, ce qui implique que chacun des~$\sch Q_i$ appartient à~$\red{(Y,y)}$, et correspond
donc à une branche~$b_i$ issue de~$y$ et situé au-dessus de~$a$. 

\medskip
Étant sans bord et de dimension nulle en~$y$, le morphisme~$\phi$ est fini en~$y$ ; comme~$x$
n'est pas rigide et comme~$X$ est réduite par quasi-lissité~$(Y,y)\to (X,x)$ est fini et plat. Pour tout~$i$,
notons~$e_i$ l'indice de ramification de~${\sch Q}_i$ sur~$\sch P$ ; on a~$e_i=\deg \;(b_i\to a)$ pour tout~$i$ d'après l'assertion i) déjà établie. Il vient, en se souvenant que~$\sch D$ et~$\sch C$ sont des courbes projectives sur  le corps {\em algébriquement clos}~$\red k$ :~$$[\kappa(y):\kappa(x)]=\deg^y\phi=\sum_i\deg \;(b_i\to a)=\sum_i e_i$$~$$=\deg \;({\sch D}\to {\sch C})=[\red{\kappa(y)}:\red{\kappa(x)}].$$

\medskip
Supposons maintenant que~$y\in Y\typ 3$. En vertu de la proposition~\ref{degetseccor} et de~\ref{commentdeg}), il existe un voisinage coronaire~$Z$ de~$x$ dans~$X$ et un voisinage coronaire~$Z'$ de~$y$ tel que~$\phi$ induise un morphisme fini et plat~$Z'\to Z$, de degré égal à~$\deg^y\phi$. D'après~\ref{morcourintersk}, l'extension~$\hres(y)/\hres(x)$ est sans défaut et de degré~$\deg^y \phi$. Comme~$y$ et~$x$ sont de type 3, les corps résiduels~$\widetilde{\hres(y)}$ et~$\widetilde{\hres(x)}$ sont tous deux égaux à~$\kk$ ; par conséquent,~$[|\hres(y)\ti|:|\hres(x)\ti|]=\deg^y \phi$ ; par densité,~$[|\kappa(y)\ti|:|\kappa(x)\ti|]$ est aussi égal à~$\deg^y\phi$, ce qui achève la démonstration.~$\Box$

%

\subsection*{Stabilité des corps de degré de transcendance résiduel
gradué
maximal}

\deux{stable23} {\bf Théorème.} {\em Supposons que le corps~$k$ soit stable. Soit~$X$ un espace~$k$-analytique
et soit~$x$ un point d'Abhyankar de~$X$ dont on note~$d$ le rang ; le corps~$\hres(x)$ est alors stable}.

\medskip
{\em Démonstration.} Soit~$(Y,\phi)$ une présentation d'Abhyankar de~$x$, et soit~$\bf r$ le polyrayon 
tel que~$\phi(x)=\eta_{\bf r}$. Comme~$\hres(x)$ est une extension finie de~$\hres(\eta_{\bf r})$ (\ref{conclusion-presentation-abhyankar}), 
il suffit de démontrer que~$\hres(\eta_{\bf r})$ est stable. Par une récurrence immédiate, on se ramène finalement à montrer que~$\hres(\eta_r)$ est stable pour tout~$r>0$ ; fixons donc un tel~$r$, et écrivons~$\eta$ au lieu de~$\eta_r$, et~$\eta_F$ au lieu de~$\eta_{r,F}$ pour toute extension complète~$F$ de~$k$. Remarquons que si la valeur absolue de~$k$ est triviale, alors la valeur absolue de~$\hres(\eta)$ est ou bien triviale, ou bien discrète ; dans chacun de ces deux cas,~$\hres(\eta)$ est stable, et l'on peut donc supposer à partir de maintenant que la valeur absolue de~$k$ est non triviale. 

\medskip
Soit~$L$ une extension finie séparable de~$\hres(\eta)$ ; nous allons 
démontrer que~$[\red L:\red{\hres(\eta)}]=[L:\hres(\eta)]$, ce qui permettra de conclure, puisque la stabilité d'un corps complet se teste sur ses extensions finies séparables. La catégorie des~$\hres(\eta)$-algèbres finies étales étant naturellement équivalente à la catégorie des revêtements finis étales du germe~$(\Aff^{1,\rm an}_k,\eta)$, il existe un espace~$k$-analytique~$Z$, un point~$z$ de~$Z$ et un morphisme fini étale de~$Z$ sur un voisinage ouvert de~$\eta$ dans~$\Aff^{1,\rm an}_k$ qui envoie~$z$ sur~$\eta$ et est tel que l'extension~$\hres(z)$ de~$\hres(\eta)$ s'identifie à~$L$. Comme~$Z\to \Aff^{1,\rm an}_k$ est étale,~$Z$ est lisse. 

\medskip
Soient~$z_1,\ldots,z_m$ les antécédents de~$\eta_{\KK}$ sur~$Z_{\KK}$. Comme~$Z\to \Aff^{1,\rm an}$ est étale,~$\hres(z_i)\simeq \kappa(z_i)\otimes_{\kappa(z)}\hres(z)$ pour tout~$i$. En vertu du théorème~\ref{ramresi}, l'extension~$\kappa(z_i)/\kappa(\eta_{\KK})$ est sans défaut pour tout~$i$ ; par conséquent, l'extension~$\hres(z_i)/\hres(\eta_{\KK})$ est sans défaut pour tout~$i$ ; notons que~$\prod \hres(z_i)$ s'identifie à~$L\otimes_{\hres(\eta)}\hres(\eta_{\KK})$. 

\medskip
Comme la valeur absolue de~$k$ est non triviale,~$k^s$ est dense dans~$\KK$. Par le lemme de Krasner,~$L\otimes_{\hres(\eta)} (\hres(\eta)\otimes_k k^s)$ s'écrit sous la forme~$\prod L_i$, où~$L_i$ est pour tout~$i$ un sous-corps dense de~$\hres(z_i)$. Pour tout~$i$, l'extension finie séparable~$L_i/(\hres(\eta)\otimes_k k^s)$ est sans défaut ; on choisit un ensemble~$\sch B_i$ d'éléments de~$L_i\ti$ de cardinal~$\dim {\hres(\eta)\otimes_k k^s} L_i$ et tel que les~$\red b$ constituent, pour~$b$ parcourant~$\sch B_i$, une base de~$\red {L_i}$ sur~$\red{(\hres(\eta)\otimes_k k^s)}$. 

\medskip
Il existe une sous-extension finie~$F$ de~$k^s/k$ telle que les~$L_i$ et les~$\sch B_i$ soient «définis sur~$F$», c'est-à-dire plus précisément telle que~$L\otimes_{\hres(\eta)}\hres(\eta_F)$ s'écrive~$\prod L'_i$, où chaque~$L'_i$ est une extension finie séparable de~$\hres(\eta_F)$ dont le produit tensoriel avec~$\hres(\eta)\otimes_k k^s$ s'identifie à~$L_i$, et qui contient~$\sch B_i$ modulo cette identification. 

\medskip
Fixons~$i$. Les~$\red b$, pour~$b$ parcourant~$\sch B_i$, sont linéairement indépendants sur~$\red {(\hres(\eta)\otimes_k k^s)}$ ; ils le sont {\em a fortiori} sur~$\red {\hres(\eta_F)}$. Par conséquent, on a~$$[\red{L'_i}:\red {\hres(\eta_F)}]\geq {\rm card}\;(\sch B_i)=\dim {\hres(\eta)\otimes_k k^s} L_i=\dim {\hres(\eta_F)}L'_i,~$$ et l'on a finalement~$[\red{L'_i}:\red {\hres(\eta_F)}]=\dim {\hres(\eta_F)}L'_i$. 

\medskip
On déduit par ailleurs des descriptions explicites de~$\red{\hres(\eta_F)}$
et~$\red{\hres(\eta)}$ ({\em cf.}~\ref{remresidgauss})
que~$$[\red{\hres(\eta_F)}:\red {\hres(\eta)}]=[\red F:\red k]=[F:k],$$ cette dernière égalité provenant de l'hypothèse de stabilité de~$k$.

\medskip
On a alors d'une part~$$\sum [\red {L'_i}:\red {\hres(\eta)}]=[F:k]\sum [L'_i:\hres(\eta_F)]$$~$$=[F:k]\cdot[L:\hres(\eta)],~$$ et d'autre part~$$\sum [\red {L'_i}:\red {\hres(\eta)}]=[\red L:\red {\hres(\eta)}] \sum[ \red {L'_i}:\red L]$$ ~$$\leq [L:\hres(\eta)]\sum [L'_i:L]=[L:\hres(\eta)]\cdot[\hres(\eta_F):\hres(\eta)]=[L:\hres(\eta)]\cdot[F:k]$$~$$=\sum [\red {L'_i}:\red {\hres(\eta)}],$$ cette dernière égalité provenant de ce qui précède. 

Par conséquent, l'inégalité large ci-dessus est en réalité une égalité, et l'on a en particulier~$$[\red L:\red{\hres(\eta)}]=[L:\hres(\eta)],$$  ce qui achève la démonstration.~$\Box$

\deux{ant-canon-alclos}
{\bf Corollaire.}
{\em Soit~$X$ un espace~$\KK$-analytique, soit~$x$ un point d'Abhyankar
de~$X$
et soit~$L$ une extension complète de~$\KK$ ; on note~$\pi$
le morphisme canonique~$X_L\to X$. Le sous-ensemble
fini $\shil L x$ de~$\pi^{-1}(x)$ 
(\ref{ex-ant-can}) est alors un singleton~$\{y\}$, et~$\red{\hres(y)}$ s'identifie
au corpoïde des fractions de l'annéloïde~$\red{\hres(x)}\otimes_{\kk}\red L$ (lequel
est intègre d'après~\ref{geom-integre-oide}) .}

\medskip
{\em Démonstration.}
Soit~$(Y,\phi)$ une présentation d'Abhyankar de~$x$, 
et soit~$\bf r$ le polyrayon tel que~$\phi(x)=\eta_{{\bf r},\KK}$. 
On a~$\phi(x)=\eta_{{\bf r},\KK}$ et
$\phi^{-1}(\eta_{{\bf r},\KK})$ est de dimension nulle d'après~\ref{conclusion-presentation-abhyankar} ; 
soit~$d$ le degré de~$\hres(x)$ sur~$\hres(\eta_{{\bf r}, \KK})$. Comme~$\hres(\eta_{{\bf r}, \KK})$
est stable par le théorème~\ref{stable23}, on a $[\red{\hres(x)}:\red{\hres(\eta_{{\bf r}, \KK})}]=d$. 

\trois{majoration-dim-hresy}
En vertu de la proposition~\ref{prop-shilov-antec}, $\shil L x$
est
l'ensemble des antécédents de~$x$ sur~$\phi^{-1}(\eta_{{\bf r}, L})$. Soit~$y\in \shil L x$. 
Le corps~$\hres(y)$ est un quotient de~$\hres(x)\otimes_{\hres(\eta_{{\bf r}, \KK})}\hres(\eta_{{\bf r}, L})$. Par
conséquent, $[\hres(y):\hres(\eta_{{\bf r}, L})]\leq d$. 

\trois{hresx-tenseur-hresetarl}
Les corpoïdes~$\red{\hres(\eta_{{\bf r}, \KK})}$ et~$\red{\hres(\eta_{{\bf r}, L})}$ sont respectivement
égaux aux algèbres de polynômes en les~$\red{T_i(\eta_{{\bf r}, \KK})}$ à coefficients dans~$\kk$ et~$\red L$. 
On en déduit que~$\red{\hres(\eta_{{\bf r}, L})}$ s'identifie 
au corpoïde des fractions
de~$\red L \otimes_{\kk} \red{\hres(\eta_{{\bf r}, \KK})}$. 
Dès lors, la~$\red{\hres(\eta_{{\bf r},L})}$-algèbre
~$\red{\hres(x)}\otimes_{\red{\hres(\eta_{{\bf r},\KK})}}\red{\hres(\eta_{{\bf r},L})}$
est
le localisé de l'annéloïde intègre~$\red{\hres(x)}\otimes_{\kk} \red L$ par l'image
de~$(\red{\hres(\eta_{{\bf r}, \KK})}\otimes_{\kk} \red L)^{\neq 0}$, et est donc intègre.
Étant par ailleurs entière sur le corpoïde~$\red{\hres(\eta_{{\bf r}, L})}$, cette algèbre est un 
corpoïde, et c'est dès lors plus précisément
le corpoïde des fractions de~$\red{\hres(x)}\otimes_{\kk} \red L$
(puisque c'est un localisé de ce dernier). 

\trois{conclu-shilov-algclos}
Le morphisme
naturel du corpoïde~$\red{\hres(x)}\otimes_{\red{\hres(\eta_{{\bf r},\KK})}}\red{\hres(\eta_{{\bf r},L})}$
dans~$\red{\hres(y)}$ est nécessairement injectif. La dimension du~$\red{\hres(\eta_{{\bf r}, L})}$-espace
vectoriel~$\red{\hres(x)}\otimes_{\red{\hres(\eta_{{\bf r},\KK})}}\red{\hres(\eta_{{\bf r},L})}$ est égale à~$d$, 
et
$$[\red{\hres(y)}:\red{\hres(\eta_{{\bf r}, L})}]\leq [\hres(y):\hres(\eta_{{\bf r}, L})]\leq d$$

Il vient
$$\red{\hres(y)}=\red{\hres(x)}\otimes_{\red{\hres(\eta_{{\bf r},\KK})}}\red{\hres(\eta_{{\bf r},L})}$$
et
$$[\red{\hres(y)}:\red{\hres(\eta_{{\bf r}, L})}]=[\hres(y):\hres(\eta_{{\bf r}, L})]= d.$$
Cette dernière égalité assure
que
la flèche naturelle
$$\hres(x)\otimes_{\hres(\eta_{{\bf r}, \KK})}\hres(\eta_{{\bf r}, L})\to \hres(y)$$ est un isomorphisme,
ce qui entraîne que~$y$ est le seul antécédent de~$x$ sur~$\phi^{-1}(\eta_{{\bf r}, L})$. 

\medskip
On a vu ci-dessus que~$\red{\hres(y)}=\red{\hres(x)}\otimes_{\red{\hres(\eta_{{\bf r},\KK})}}\red{\hres(\eta_{{\bf r}, L})}$,
et 
ce dernier est en vertu de~\ref{hresx-tenseur-hresetarl}
le corpoïde des fractions
de~$\red{\hres(x)}\otimes_{\kk}\red L$ ; ceci achève la démonstration.~$\Box$

\deux{ant-canon-cardinal}
{\bf Corollaire.}
{\em Soit~$X$ un espace~$k$-analytique, soit~$x$ un point d'Abhyankar
de~$X$
et soit~$L$ une extension complète de~$k$. L'ensemble~$\shil L x$ est en bijection
naturelle avec celui des idéaux maximaux de~$L\otimes_k \got s(x)$.}

\medskip
{\em Démonstration}. Quitte à remplacer~$X$ par un voisinage convenable de~$x$, on peut
supposer que~$\got s(X)=\got s(x)$. Écrivons~$L\otimes_k \got s(x)=\prod L_i$, où chaque~$L_i$
est une extension finie séparable de~$L$. En vertu de~\ref{shil-l-x-connexe}, il suffit de montrer que
pour tout indice~$i$, le sous-ensemble~$\shil {L_i} x$
de~$X_{L_i}$ (défini en considérant~$X$ comme un espace~$\got s(x)$-analytique) est un singleton. 

\medskip
On se ramène ainsi au cas où~$\got s(x)=k$. Soit~$L'$ une extension complète
de~$k$ composée de~$\KK$
et~$L$. Comme~$\got s(x)=k$, le point~$x$ a un unique antécédent sur~$X_{\KK}$, qui coïncide avec
$\shil {\KK} x$. Il résulte du corollaire~\ref{ant-canon-alclos}
que~$\shil {L'} x=\shil {L'}{(\shil {\KK} x)}$ est un singleton. Puisque
$\shil {L'} x$ est la réunion des~$\shil {L'} y$ pour~$y$ parcourant
$\shil L x$, l'ensemble $\shil L x$ est lui-même un singleton.~$\Box$

\section{Voisinages des points d'une courbe quasi-lisse sur un corps algébriquement clos}

\deux{defgrenreres} Si~$x$ est un point de type 2 d'une courbe~$\KK$-analytique, on appellera {\em genre} de~$x$  le genre de la courbe résiduelle en~$x$ ; si~$x$ est un point de type 3, son genre sera {\em par définition} égal à~$0$ ; le genre d'un point~$x$ de type 2 ou 3 sera noté~$g(x)$. Par abus, nous dirons simplement \og point de genre tant\fg~au lieu de \og point de type 2 ou 3 et de genre tant\fg.

\medskip
Si~$x$ est un point de type 2 ou 3 de~$\pkk$ 
alors~$g(x)=0$ : si~$x$ est de type 3 c'est une conséquence
de la définition, et si~$x$ est de type~$2$ cela résulte de~\ref{extype-2}. 

\deux{fixation-ell}
On fixe pour toute cette section un entier~$\ell$ premier à~$p$. Dans ce qui suit, nous dirons «$\mu_\ell$-torseur»
pour «$\mu_\ell$-torseur étale».

\subsection*{Germes de
torseurs}

\deux{voisbiensect} Soit~$X$ une courbe~$\KK$-analytique quasi-lisse
et soit~$x\in X\dtr$. Nous dirons que $X$ est
{\em bien découpé autour de~$x$} si c'est un arbre, et si pour toute branche~$b$
de~$X$ issue de~$x$ l'une des deux conditions suivantes (exclusives l'une de l'autre)
est satisfaite :

\medskip
$\bullet$ $b(X)$ est relativement compacte dans~$X$ ; 

$\bullet$ $b(X)$ est une couronne. 

\medskip
Lorsque la seconde condition est satisfaite, nous dirons parfois que~$b(X)$ est une
{\em composante coronaire}
de~$X\setminus\{x\}$. 
 
 \trois{voisbiensectraff} Supposons que~$X$ soit bien découpé
 autour de~$x$, et soit~$U$ un voisinage ouvert de~$x$ dans~$X$,
 bien découpé autour de~$x$. 
 Soit~$Z$
 une composante coronaire de~$X\setminus\{x\}$.
L'intersection~$U\cap Z$ est une composante connexe de~$U\setminus\{x\}$. 
Si elle était relativement compacte dans~$U$, elle coïnciderait avec~$Z$, ce qui est absurde puisque cette dernière n'est pas
 relativement compacte dans~$X$. 
 En conséquence,~$U\cap Z$ est nécessairement une composante coronaire de~$U\setminus \{x\}$. 
  
 \medskip
 Nous dirons que le {\em couple}~$(U,X)$ est 
 bien découpé autour de~$x$ 
 si pour toute composante coronaire~$Z$ de~$U$,
 la composante coronaire
 de~$U\setminus\{x\}$ contenue dans~$U$ est une sous-couronne de~$Z$. 
 
 \trois{basebiendecoup} Il résulte de la proposition~\ref{xmoinssbr}
 et de~\ref{souscourmeme}
 que le point~$x$ possède une base de voisinages ouverts bien découpés
 autour de~$x$. 
 Si~$U$ est un tel voisinage, le point~$x$ possède même plus précisément une base 
 de voisinages ouverts~$V$ tels que le couple~$(V,U)$ soit bien découpé
 autour de~$x$. 
 
 \medskip
 En effet, soit~$W$ un voisinage ouvert de~$x$. Il existe un voisinage ouvert~$W_0$ de~$x$ 
 dans~$U\cap W$ qui est bien découpé autour de~$x$. Soit~$\Pi$ l'ensemble des composantes
 coronaires de~$U$ ; pour tout~$Z\in \Pi$, l'intersection~$Z\cap W_0$ 
 est une couronne contenue dans~$Z$ et aboutissant à~$x$, et il existe une
 sous-couronne~$C_Z$ de~$Z\cap W_0$ qui abouti à~$x$ et est aussi une sous-couronne
 de~$Z$ ; le complémentaire~$C'_Z$ de~$C_Z$ dans~$Z\cap W_0$ est fermé dans~$W_0$, par exemple
 en vertu de la proposition~\ref{xmoinssbr}. Le
 complémentaire
 $V$ de~$\coprod_{Z\in \Pi} C'_Z$ dans~$W_0$
 est alors un voisinage ouvert de~$x$ contenu dans~$W$, et
 le couple~$(V,U)$
 est bien découpé autour de~$x$.

\deux{trivetcompco} Soient~$X$ une courbe~$\KK$-analytique 
quasi-lisse et connexe, 
soit~$x\in X\typ 2$ et soit~$\phi : Y\to X$ un morphisme fini \'etale. 
Soit~$\sch C$
la courbe résiduelle en~$x$. Pour tout~$y\in Y_x$ on note~$\sch D_y$ la courbe résiduelle en~$y$. 

\trois{nonram-test-residuel}
Soit~$y\in Y_x$. Comme le corps~$\hres(x)$ est stable
(th.~\ref{stable23}) et comme~$|\hres(x)\ti|$ est divisible, 
on a l'égalité~$[\red{\hres(y)_1}:\red{\hres(x)}_1]=[\hres(y):\hres(x)]$ ; en conséquence,
$\hres(y)$ est non ramifiée sur~$\hres(x)$
si et seulement si~$\red{\hres(y)}_1$ est séparable
sur~$\red{\hres(x)}_1$, c'est-à-dire encore si et seulement si~$\sch D_y$ est génériquement
étale sur~$\sch C$. 

\trois{defcompPi}
Soit~$\Pi$ le sous-ensemble de~$\pi_0(Y\setminus\phi\inv(x))$ formé des composantes 
connexes~$V$ satisfaisant les conditions suivantes : 

\medskip
1) le bord de~$V$ est de la forme~$\{y\}$ pour un certain~$y\in \phi\inv(x)$ ; 

2) l'ensemble~$\br Y y\ctd V$ est un singleton~$\{b\}$ ; 

3) le morphisme~$\sch D_y\to \sch C$ est étale en le point fermé
de~$\sch D_y$ correspondant à~$b$. 

\medskip
On  désigne par~$\Pi'$ le sous-ensemble de~$\pi_0(X\setminus\{x\})$ formé
des composantes~$U$ telles que tout composante connexe de~$\phi\inv(U)$
appartienne à~$\Pi$.

\trois{vdanspitriv} Soit~$V\in \Pi$ et
soit~$y$ l'unique point
de son bord. 
L'image~$\phi(V)$ est une composante connexe
de~$X\setminus\{x\}$. L'ensemble~$\br Y y\ctd V$ est un singleton~$\{b\}$ ; 
posons~$a=\phi(b)$. 
Le lemme~\ref{unebrdeg} assure que~$\deg (V\to \phi(V))=\deg\;(b\to a)$. 
Par ailleurs, le morphisme
$\sch D_y\to \sch C$ est par 
hypothèse étale en le point
fermé de
~$\sch D_y$ qui correspond à
la branche~$b$ ; par le
théorème~\ref{ramresi}
ceci entraîne l'égalité~$\deg\;(b\to a)=1$. En conséquence,~$\deg (V\to \phi(V))=1$ ; 
autrement dit,~$\phi$ induit un
isomorphisme~$V\simeq \phi(V)$.

\trois{pprimenonram}
Supposons que~$\hres(y)$ est non ramifiée sur~$\hres(x)$ 
pour tout
antécédent~$y$
de~$x$, 
c'est-à-dire que~$\sch D_y$ est génériquement
étale sur~$\sch C$ pour tout~$y\in Y_x$. 
Cette hypothèse implique que
presque toutes les composantes connexes de~$Y\setminus\phi\inv(x)$
appartiennent à~$\Pi$, et partant 
que presque toutes les composantes
connexes de~$X\setminus\{x\}$ appartiennent à~$\Pi'$. 

\medskip
On déduit dès lors de~\ref{vdanspitriv}
que 
le revêtement~$Y\times_X U\to U$ est trivial pour presque toute composante connexe~$U$
de~$X\setminus \{x\}$. 

\deux{intromueltors} Soit~$X$ une courbe~$\KK$-analytique quasi-lisse
et soit~$x\in X$.  Le but de ce
qui suit
est d'étudier le groupe
$\H^1((X,x)\et,\mu_\ell)$, qui classifie les germes de~$\mu_\ell$-torseurs
en~$x$. Notons pour commencer l'existence d'isomorphismes naturels~$\H^1((X,x)\et,\mu_\ell)
\simeq \H^1(\hres(x),\mu_\ell)\simeq \H^1(\kappa(x),\mu_\ell)).$

\trois{explh1xxmu} Donnons une description explicite
de ces isomorphismes. 
Soit~$h$ une classe appartenant
à~$\H^1(\hres(x),\mu_\ell)$. Elle est de la forme~$(\lambda)$ pour
un certain~$\lambda$
appartenant à~$\hres(x)\ti$, bien déterminé modulo~$(\hres(x)\ti)^\ell$. Dire
que~$\H^1(\kappa(x),\mu_\ell))\simeq\H^1(\hres(x),\mu_\ell)$
signifie simplement
que le scalaire~$\lambda$ peut être choisi dans~$\kappa(x)\ti$,
et qu'un tel choix
est unique modulo~$(\kappa(x)\ti)^\ell$. 
Il existe donc
une fonction
$f\in \sch O_{X,x}\ti$ telle que~$h=(f(x))$, et la classe
de~$\H^1((X,x)\et,\mu_\ell)$ qui correspond à~$h$ est alors 
la classe de Kummer~$(f)$. En termes de torseurs, c'est la classe
du (germe en~$x$) de~$\sch M(\sch O_U[T]/(T^\ell-f))\to U$, où~$U$ est un voisinage
ouvert de~$x$ sur lequel~$f$ est définie et inversible.

 \trois{torslocalkum} On déduit de~\ref{explh1xxmu}
que tout~$\mu_\ell$-torseur 
sur~$(X,x)$
est de Kummer, c'est-à-dire encore que tout~$\mu_\ell$-torseur 
défini sur un ouvert contenant~$x$ est de Kummer au voisinage de~$x$. 

\trois{cardfibremutors} Soit~$h\in \H^1((X,x)\et,\mu_\ell)\simeq  \H^1(\hres(x),\mu_\ell)$,
et soit~$Y$ un~$\mu_\ell$-torseur sur~$(X,x)$
de classe~$h$.
Soit~$d$ l'ordre
de~$h$ dans~$\H^1(X,x)\et,\mu_\ell)$. 
On déduit de~\ref{muellegalzsurell}
que la fibre~$Y_x$ comporte~$\ell/d$ points. 

%

\deux{mueltorstriv} Supposons que~$x\in X\typ{14}$. Le
corps~$\hres(x)$ est complet et en particulier hensélien, 
$\red{\hres(x)}_1$ est algébriquement clos, et~$|\hres(x)\ti|$ est divisible. 
Il s'ensuit que~$\H^1(\hres(x),\mu_\ell)=0$ (\ref{ktilde1sepclos}
et~\ref{abskldiv}) ; en conséquence, $\H^1((X,x)\et,\mu_\ell)$
est trivial.  

\medskip
Autrement dit, tout~$\mu_\ell$-torseur défini sur un ouvert contenant~$x$
est trivial au voisinage de~$x$. 

\deux{mueltorsdtr} Supposons que~$x\in X\dtr$. Par convention, nous considérerons
dans ce qui suit que toute branche de~$X$ issue de~$x$ est orientée {\em vers~$x$}. 

\deux{branchekummer} Soit~$a\in \br X x$
et soit~$h\in \H^1((X,x)\et,\mu_\ell)$.


\trois{thetaakum}
La convention
d'orientation 
 adoptée
induit un système compatible
 d'isomorphismes
 $$\sigma_Z : \kum (Z)\simeq \ZZ/\ell \ZZ$$ pour~$Z$ parcourant l'ensemble
 des sections coronaires de~$a$. 
 
 Si~$Z$ est une section coronaire
de~$a$ et si~$Z'$ est un~$\mu_\ell$-torseur de Kummer sur~$Z$, l'image
par~$\sigma_Z$ de la classe de~$Z'$ dans~$\kum(Z)$ sera appelée
{\em l'invariant}
de~$Z'$. 
 
 \medskip
 Soit~$U$ un voisinage
 ouvert de~$x$ tel que la classe~$h$ provienne d'une classe
 $h'\in \H^1(U\et,\mu_\ell)$, et
 soit~$Z$ une section coronaire de~$a$ contenue dans
 l'ouvert~$U$ telle que~$h'_{|Z}$
soit de Kummer (une telle section existe toujours puisque~$h'$
est de Kummer au voisinage de~$x$, {\em cf.}
\ref{intromueltors}).
 L'élément~$\sigma_Z(h')$ de~$\ZZ/\ell \ZZ$
 ne dépend que de~$a$ (et pas du choix de~$U,Z$ et~$h')$ ; on le note~$\theta_a(h)$. 
 
 \trois{thetaahens} Soit~$|.|_a$ la valuation 
 de~$\kappa(x)$ associé à la variation en norme
 le long de~$a$.  Dans ce qui suit, nous allons utiliser
 implicitement le th.~\ref{sectioncour}
 ainsi que les commentaires qui en sont faits en~\ref{remcalcval}. 
 
 \medskip
 Soit~$t$ une fonction inversible sur un voisinage ouvert~$U_0$ de~$x$
 dans~$U$
 telle que pour toute section coronaire~$Z$ de~$a$
 dans~$U_0$, la fonction~$t_{|Z}$ soit une fonction coordonnée de~$Z$
 croissante en norme sur~$\skel Z$ (il suffit de vérifier ces conditions sur {\em une}
 telle section~$Z$). 
 
 \medskip
 Le groupe quotient~$|\kappa(x)\ti|_a/|\KK\ti|$ est libre de rang~$1$, et est engendré par la classe de~$|t(x)|_a$
 qui ne dépend pas du choix de~$t$. La formule~$|t(x)|_a^m\mapsto m$ définit donc sans ambiguïté un isomorphisme
 $|\kappa(x)\ti|_a/|\KK\ti|\simeq \ZZ$, puis par passage au quotient un isomorphisme
 $\iota : |\kappa(x)\ti|_a/|\kappa(x)\ti|_a^\ell\simeq \ZZ/\ell\ZZ$ (comme~$|\KK\ti|$ est divisible, il disparaît lors de
 cette dernière opération). 
 
 Soit~$\kappa(x)^h$ le hensélisé de~$\kappa(x)$ pour~$|.|_a$. Puisque le corps résiduel
 de~$|.|_a$ est algébriquement clos, 
 $\H^1(\kappa(x)^h,\mu_\ell)\simeq |\kappa(x)\ti|/|\kappa(x)\ti|^\ell$. On note~$\tau_a$ le morphisme
 composé
 $$\H^1((X,x)\et,\mu_\ell)\simeq \H^1(\kappa(x),\mu_\ell)\to \H^1(\kappa(x)^h,\mu_\ell)\simeq |\kappa(x)\ti|/|\kappa(x)\ti|^\ell\stackrel \iota\to \ZZ/\ell \ZZ
.$$

 \medskip
 Soit~$f$ une fonction inversible sur
 un voisinage~$U_1$ de~$x$ dans~$U_0$ telle que~$h'_{|U_1}=(f)$,
 et soit~$Z$
 une section coronaire de~$a$ contenue dans~$U_1$. 
 Il existe~$\lambda\in \KK\ti$
 et~$m\in \ZZ$
 tels que~$|f|=|\lambda|\cdot |t|^m$ identiquement sur~$Z$. Il résulte du choix de~$t$
 et de la définition de~$\theta_a$ que~$\theta_a(h)=\overline m\in \ZZ/\ell \ZZ$.

L'image de~$h$ dans~$\H^1(\kappa(x),\mu_\ell)$
est égale à~$(f(x))$. Comme~$|f|=|\lambda|\cdot |t|^m$ identiquement sur~$Z$, on a~$|f(x)|_a=|\lambda|\cdot |t(x)|_a^m$. 
Il vient~$$\tau_a(h)=\overline m=\theta_a(h).$$

\trois{thetares} Supposons maintenant
que~$x$ est de type 2, soit~$\sch C$ sa
courbe résiduelle et soit~$\sch P$
le point fermé de~$\sch C$ correspondant à~$a$. La
valuation~$|.|_a$
étant composée de
la valuation structurale de~$\kappa(x)$ et de la valuation discrète
associée à~$\sch P$, on peut choisir 
la fonction~$t$ de sorte que~$|t(x)|=1$ et que~$\red{t(x)}$ soit une uniformisante
de l'anneau de valuation discrète~$\sch  O_{\sch C,\sch P}$
(comme~$|t(x)|_a$ est alors infiniment proche {\em inférieurement}
de~$1$, la fonction~$|t|$ est bien croissante en norme le long de~$a$ à l'approche de~$x$).

Comme le point~$x$ est de type 2, le groupe~$|\kappa(x)\ti|$ est divisible, 
ce qui entraîne que~$\H^1(\kappa(x),\mu_\ell)\simeq\H^1(\kappa(x)\zero,\mu_\ell)\simeq \H^1(\red{\hres(x)}_1,\mu_\ell)$. 
Par conséquent, 
on peut supposer que~$|f(x)|=1$, et
la classe~$\red h$
de~$H^1(\red{\hres(x)}_1,\mu_\ell)$ qui correspond à~$h$ {\em via}
les isomorphismes
canoniques~$\H^1((X,x)\et,\mu_\ell)\simeq \H^1(\kappa(x),\mu_\ell)\simeq \H^1(\red{\hres(x)}_1,\mu_\ell)$
est alors précisément~$(\red{f(x)})$. 

\medskip
L'égalité~$|f(x)|_a=|\lambda|\cdot |t(x)|_a^m$ implique que~$|f(x)|=|\lambda|\cdot|t(x)|$ ; par conséquent, $|\lambda|=1$, et il vient
$|f(x)/t(x)^m|_a=1$. Cette dernière égalité signifie que~$\red{f(x)}/\red{t(x)}^m$ est de valuation~$\sch P$-adique (additive)
nulle. En conséquence, la valuation~$\sch P$-adique de~$\red{f(x)}$ est égale à~$m$, et le résidu~$\delta_{\sch P}(\red h)$ est égal
à~$\overline m$. Ainsi, $\theta_a(h)=\tau_a(h)=\delta_{\sch P}(\red h)$.

\medskip
Remarquons  pour terminer que comme
la~$\hres(x)$-algèbre étale
$$\hres(x)[T]/T^\ell-f(x)=\prod_{y\in Y_x}\hres(y)$$ est non ramifiée, 
il résulte de~\ref{pprimenonram}
que pour tout voisinage ouvert connexe~$V$ de~$x$ dans~$U$, 
la classe~$h'$
s'annule 
sur presque toutes les composantes connexes de~$V\setminus\{x\}$.

\deux{prepacondresidu} On revient aux hypothèses et notations précédant le~\ref{thetares} : le point~$x$
est supposé appartenir à~$X\dtr$, mais pas nécessairement à~$X\typ 2$. Si~$V$ est un voisinage ouvert connexe
de~$x$ dans~$U$, alors~$h'$  s'annule
au-dessus de presque toutes les composantes connexes de~$V\setminus\{x\}$ : en effet, 
si~$x\in X\typ 2$ cela découle de~\ref{thetares} ; et si~$x\in X\typ 3$ alors
comme~$\br X x$ est fini,~$\pi_0(V\setminus \{x\})$ l'est aussi,
et l'assertion requise est tautologique. 

\medskip
Il existe en conséquence une base de voisinages ouverts~$V$ de~$x$ dans~$U$ possédant les propriétés suivantes :
 
 \medskip
 $\bullet$ $V$ est bien découpé autour de~$x$ (\ref{voisbiensect}) ; 
 
 $\bullet$ la classe~$h'_{|V}$ est de Kummer ; 
 
 $\bullet$ la classe~$h'$ s'annule sur toute composante connexe de~$V\setminus\{x\}$
 qui est relativement compacte dans~$V$. 
 
 \medskip
 Dans ce cas, pour toute branche~$a\in \br X x$ telle que~$a(V)$ soit coronaire, la restriction~$h'_{|a(V)}$, est de Kummer, 
 et~$\sigma_{a(V)}(h'_{|a(V)})=\theta_a(h)$. 
 
 \deux{propcondresidutors}{ \bf Proposition.}
 {\em Soit~$X$ une courbe~$\KK$-analytique quasi-lisse et soit~$x$
 appartenant à~$X\dtr$. 
Soit~$X_0$ un voisinage
 ouvert de~$x$ dans~$X$
 bien découpé autour de~$X$. On note~$A$
 l'ensemble des branches~$a$ issues de~$x$
 telles que~$a(X_0)$ soit coronaire. 
 
 \medskip 
 \begin{itemize}
 \item[1)] 
 Soit~$(\lambda_a)_{a\in \br X x}$ une famille d'éléments de~$\ZZ/\ell\ZZ$
 tels que~$\lambda_a=0$ si~$a\notin A$. Les assertions suivantes sont
 équivalentes. 
 
 \medskip
 \begin{itemize}
 
 \item[i)] Il existe un~$\mu_\ell$-torseur~$Y\to X_0$ tel que : 
 
 \medskip
 \begin{itemize}
 \item[$\diamond$] $Y$ se trivialise au-dessus de toute
  composante connexe de~$X_0\setminus\{x\}$
 qui est relativement compacte dans~$X_0$ ; 
 
 \item[$\diamond$] pour tout~$a\in A$, le~$\mu_\ell$-torseur~$Y\times_X a(X_0)$ est de Kummer,
 d'invariant~$\lambda_a$.
 
 \end{itemize}
 
 \medskip
 
 \item[ii)] Il existe~$h\in \H^1((X,x)\et,\mu_\ell)$ telle que~$\theta_a(h)=\lambda_a$ pour
 toute branche~$a\in \br X x$. 
 \medskip
 
 \item[iii)] On a~$x\in \partial X$ ou~$\sum_{a\in A}\lambda_a=0$.
 \end{itemize}

\medskip
De plus :

\medskip
1a) si~i) est vraie, la classe~$h$ de ii) peut être prise égale à
la classe de~$Y\times_{X_0}(X,x)$ ; 

1b) si~ii) est vraie, on peut choisir~$Y$ dans~i) de sorte que~$Y\times_{X_0}(X,x)$ ait pour classe~$h$.

\medskip
\item[2)] 
Supposons que~$x\notin \partial X$ et soit~$\sch K$
le sous-groupe
de~$\H^1((X,x)\et,\mu_\ell)$ formé des classes~$h$ telles que
$\theta_a(h)=0$ pour tout~$a\in \br X x$. Si~$x\in X\typ 3$ alors~$\sch K$ est trivial. Si~$x\in X\typ 2$,
si~$\sch C$ désigne sa courbe résiduelle, et si~$\sch J$ désigne la jacobienne de~$\sch C$ alors

$$\sch K\simeq _\ell\sch J(\kk_1)\simeq (\ZZ/\ell\ZZ)^{2g(x)}$$ (le premier de ces
isomorphismes est canonique, le second ne l'est pas en général). 
 
 \end{itemize}}
 \medskip
 {\em Démonstration.}
 Nous allons montrer i)$\Rightarrow$ii), ii)$\Rightarrow$iii), iii)$\Rightarrow$ii)
 et~ii)$\Rightarrow$i), et les autres assertions de l'énoncé seront établies incidemment. 
 
 \trois{thetahzero}{\em Preuve de~i)$\Rightarrow$ii).}
 On suppose que~i) est vraie, et l'on note~$h$ la classe de~$Y\times_X (X,x)$
 dans~$\H^1((X,x)\et,\mu_\ell)$. On a par définition~$\lambda_a=\theta_a(h)$ pour tout~$a\in A$. 
 Soit~$a\in \br X x \setminus A$ et
 soit~$Z$ une section coronaire
 de~$a$. Le torseur~$Y$ se trivialise
au-dessus de~$a(X)$ puisque~$a\notin A$ ; il se trivialise
{\em a fortiori}
au-dessus de~$Z$, ce qui entraîne que~$\theta_a(h)=0$. Ainsi, ii) est vraie, et l'on a
par ailleurs établi~1a).

\trois{sommeinvzero}{\em Preuve de~ii)$\Rightarrow$iii).} On suppose
que~ii) est vraie. Pour montrer~iii),  on fait l'hypothèse
que~$x\notin \partial X$, et nous allons prouver que~$\sum_{a\in \br X x} \theta_a(h)=0$. 
Rappelons que pour toute branche~$a$ de~$X$ issue de~$x$ et toute section
 coronaire~$Z$ de~$a$, la couronne~$Z$ sera supposée orientée vers~$x$,
 et que~$\sigma_Z$
désigne
 l'isomorphisme : $\kum(Z)\simeq \ZZ/\ell\ZZ$
 déduit de ce choix d'orientation. 
 
\medskip

 {\em On suppose que~$x$ est de type 3}. L'ensemble~$\br X x$ possède alors
deux éléments~$a$ et~$b$. Soit~$U$ un voisinage coronaire de~$x$ dans~$X$ tel que~$h$
provienne d'une classe de Kummer~$h'$
sur~$U$. 
L'ouvert~$U\setminus\{x\}$ a deux composantes connexes, à savoir~$a(U)$ et~$b(U)$, qui sont deux sous-couronnes
ouvertes de~$U$. Munissons~$U$ de l'orientation induite
par celle de~$a(U)$. 
Ce choix induit un isomorphisme~$\rho : \kum (U)\simeq \ZZ/\ell\ZZ$. 

Comme l'orientation de~$U$ est induite par celle de~$a(U)$, on dispose d'un diagramme commutatif 
$$\diagram \kum (U) \rto^\rho \dto^\simeq &\ZZ/\ell \ZZ\dto^{\rm Id}\\
\kum(a(U))\rto^{\;\;\;\;\sigma_{a(U)}}&\ZZ/\ell \ZZ\enddiagram,$$ d'où l'on déduit que~$\rho(h')=\theta_a(h)$.

Comme  l'orientation de~$U$ est {\em opposée} à 
l'orientation induite par celle de~$b(U)$, on dispose d'un diagramme commutatif 
$$\diagram \kum (U) \rto^\rho \dto^\simeq &\ZZ/\ell \ZZ\dto^{-\rm Id}\\
\kum(b(U))\rto^{\;\;\;\;\sigma_{b(U)}}&\ZZ/\ell \ZZ\enddiagram,$$ d'où l'on déduit que~$\rho(h')=-\theta_b(h)$. 
En conséquence, $\theta_a(h)+\theta_b(h)=0$. 

\medskip
{\em On suppose que~$x$ est de type 2}.
Comme~$x\notin \partial X$ l'ensemble~$\br X x$ s'identifie à l'ensemble des 
points fermés de la courbe résiduelle~$\sch C$.  

On a vu en~\ref{thetares}
que~$\H^1((X,x)\et,\mu_\ell)\simeq\H^1(\red{\hres(x)}_1,\mu_\ell)$ ; 
on note~$\red h$
la classe de~$\H^1(\red{\hres(x)}_1,\mu_\ell)$
qui correspond à~$h$. Si~$a\in \br X x$, et si~$\sch P$
désigne le point fermé correspondant, on a 
d'après~{\em loc. cit.}
l'égalité $\theta_a(h)=\delta_{\sch P}(\red h)$. 

Comme $\sum_{\sch P \in \sch C(\kk_1)}\delta_{\sch P}(\red h)=0$
({\em cf.} \ref{suiteexacteresidus}),   
il vient~$\sum_{a\in \br X x}\theta_a(h)=0$, ce qu'il fallait démontrer. 

\trois{sommezeroinv3}
{\em Preuve de~iii)$\Rightarrow$ii) et de
l'assertion~2)
dans le cas d'un point de type~3.}
On fixe une famille~$(\lambda_a)_{a\in \br X x}$ satisfaisant les conditions de~iii), et l'on suppose que~$x$ est de type~3. 

\medskip
{\em Le cas où~$\br X x$ est de cardinal 2.}
C'est celui où~$x\notin \partial X$. Choisissons un voisinage coronaire~$U$
de~$x$, et notons~$a$ et~$b$ les deux branches issues de~$x$. Munissons~$U$ de l'orientation
induite par celle de sa sous-couronne~$a(U)$, et reprenons les notations du paragraphe de~\ref{sommeinvzero} consacré au cas du type 3. 
On déduit des deux diagrammes commutatifs
de~{\em loc. cit.}
que~$\sigma_{a(U)}(\rho^{-1}(\lambda_a))=\lambda_a$ 
et que~$\sigma_{b(U)}(\rho^{-1}(\lambda_a))=-\lambda_a$. 
Si~$h$ désigne l'image de~$\rho^{-1}(\lambda_a)$
dans~$\H^1((X,x)\et,\mu_\ell)$ on a donc
$\theta_a(h)=\lambda_a$ et~$\theta_b(h)=\lambda_b$, d'où~ii). 

Prouvons par ailleurs l'assertion~2), toujours en supposant que~$x\notin \partial X$. Soit~$h\in \sch K$, et soit~$U$ un voisinage
coronaire de~$x$ tel que~$h$ provienne d'une classe de Kummer~$h'$
sur~$U$. 
Comme~$\theta_a(h)=0$ la restriction
de~$h$ à~$a(U)$ est triviale ; comme~$\kum (U)\simeq \kum (a(U))$ on a~$h'=0$ , 
et donc~$h=0$ ; ainsi~$\sch K=0$, comme annoncé. 

\medskip
{\em Le cas où~$\br X x$ est de cardinal~$1$.}
Soit~$a$ l'unique branche de~$X$ issue de~$x$, et soit~$U$
un voisinage coronaire de~$x$
dans~$X$. On a~$U\setminus\{x\}=a(U)$. On munit~$U$ de l'orientation induite par celle
de~$a(U)$ ; ce choix induit un isomorphisme~$\rho : \kum (U)\simeq \ZZ/\ell \ZZ$. 
Soit~$V$ un~$\mu_\ell$-torseur sur~$U$ de classe~$\rho^{-1}(\lambda_a)$. Comme l'orientation
de~$U$ est induite par celle de~$a(U)$, la restriction de~$V$ à~$a(U)$ a pour invariant~$\lambda_a$. 
Si~$h$ désigne la classe de~$V\times_U(X,x)$ dans~$\H^1((X,x)\et,\mu_\ell)$ on a donc
$\theta_a(h)=\lambda_a$, d'où~ii). 

\medskip
{\em Le cas où~$\br X x=\emptyset$.} N'importe quelle classe~$h\in \H^1((X,x)\et,\mu_\ell)$ satisfait alors, 
pour des raisons tautologiques, la condition énoncée en~ii).

\trois{sommezeroinv2}
{\em Preuve de~iii)$\Rightarrow$ii)
et de l'assertion~2)
dans le cas d'un point de type~2.}
On fixe une famille~$(\lambda_a)_{a\in \br X x}$ satisfaisant les conditions de~iii), 
on suppose que~$x$ est de type~2, on note~$\sch C$ sa courbe résiduelle
et~$\sch J$ la jacobienne de~$\sch C$. 

Soit~$\sch C$ la courbe résiduelle en~$x$ ; on identifie~$\red{(X,x)}$ à un ouvert de Zariski non vide de~$\sch C$. Si~$\sch P$
est un point fermé de~$\red{(X,x)}$, il correspond à une branche~$a$
issue de~$x$, et l'on pose~$\mu_{\sch P}=\lambda_a$. 
Pour prouver~ii) il suffit, 
en vertu de~\ref{thetares}, de montrer l'existence d'une classe~$\red h\in \H^1(\red{\hres(x)}_1,\mu_\ell)$ telle que
$\delta_{\sch P}(\red h)=\mu_{\sch P}$ pour tout point fermé~$\sch P$ de~$\red{(X,x)}$. 

\medskip
{\em Supposons que~$x\notin \partial X$.}
On a alors~$\red{(X,x)}=\sch C$, et~$\sum\mu_{\sch P}=0$ 
puisque~$\sum \lambda_a=0$ en vertu de l'hypothèse~iii). L'existence d'une classe~$\red h$ satisfaisant les conditions
requises provient alors de la suite exacte des résidus
(\ref{suiteexacteresidus}). 

Prouvons par ailleurs l'assertion~2), toujours en supposant que~$x\notin \partial X$. 
Il résulte de~\ref{thetares}
que~$\sch K$ s'identifie à l'ensemble des
classes de~$\H^1(\red{\hres(x)}_1,\mu_\ell)$
dont tous les résidus en les points fermés de~$\sch C$ sont nuls. 
Il existe donc en vertu de~\ref{h1genreg}
un isomorphisme
canonique
entre~$\sch K$ et~$_\ell\sch J(\kk_1)$, lequel 
est lui-même isomorphe (non canoniquement
en général) à~$(\ZZ/\ell\ZZ)^{2g(x)}$, d'où~2). 

\medskip
{\em Supposons que~$x\in \partial X$.}
L'ouvert~$\red{(X,x)}$ de~$\sch C$ est alors strict. Choisissons un point fermé~$\sch Q$
de~$\sch C$ qui n'est pas situé sur~$\red{(X,x)}$, et posons~$\mu_\infty=-\sum \mu_{\sch P}$. La suite exacte des résidus assure l'existence d'une classe
$\red h\in \H^1((X,x)\et,\mu_\ell)$ telle que pour tout point fermé~$\sch P$ de~$\sch C$, 
le résidu~$\delta_{\sch P}(\red h)$ soit égal à : 

\medskip
$\bullet$ $\mu_{\sch P}$ si~$\sch P\in \red{(X,x)}$ ; 

$\bullet$ $\mu_\infty$ si~$\sch P=\sch Q$ ; 

$\bullet$ $0$ sinon. 

\medskip
La classe~$\red h$ satisfait alors les conditions requises. 

\trois{prolmuelletors}
{\em Preuve de~ii)$\Rightarrow$i)}. Supposons qu'il existe~$h$ comme dans~ii) ; choisissons
un voisinage ouvert~$U$ de~$x$ dans~$X_0$ et un~$\mu_\ell$-torseur~$V\to U$ tel que la classe
de~$V\times_U(X,x)$ dans~$\H^1((X,x)\et,\mu_\ell)$
soit égale à~$h$. Quitte à restreindre~$U$, on peut faire en sorte que les conditions suivantes soient
satisfaites : 

\medskip
$\bullet$ $V\to U$ est de Kummer ; 

$\bullet$ le couple~$(U,X_0)$ est bien découpé autour de~$x$. 

\medskip
Soit~$A'$ l'ensemble des branches~$a$
issues de~$x$ telles que~$a(U)$ soit coronaire. Si~$a\notin A'$ alors~$a(U)$ est relativement
compacte dans~$U$, et est donc égale à~$a(X_0)$. 
En conséquence, $X_0=U\coprod \left(\coprod_{a\in A'} a(X_0)\right)$. 

\medskip
Pour tout~$a\in A$, choisissons un~$\mu_\ell$-torseur
de Kummer~$W_a \to a(X_0)$ d'invariant~$\lambda_a$. Si~$a\in A'\setminus A$, on note~$W_a$
le~$\mu_\ell$-torseur trivial sur~$a(X_0)$. 

Soit~$a\in A'$. Le torseur de Kummer~$V\times_U a(U)\to a(U)$ a pour invariant~$\theta_a(h)=\lambda_a$ ; 
il est en particulier trivial si~$a\notin A$. 

\medskip
En conséquence, les restrictions à~$\coprod_{a\in A'}a(U)$ des torseurs~$V$ et~$\coprod_{a\in A'}W_a$
sont isomorphes ; en choisissant un isomorphisme arbitraire entre elles, on obtient par recollement un~$\mu_\ell$-torseur~$Y\to X_0$
qui satisfait les propriétés énonces en~i), et dont le germe ~$x$
a pour classe~$h$. Ceci achève la démonstration de l'implication~ii)$\Rightarrow$i)
et de l'assertion~1b), et partant celle de
la proposition.~$\Box$

\deux{resumtheotors} Soit~$X$ une courbe~$k$-analytique quasi-lisse, et soit~$x\in X\dtr\setminus \partial X$.
Il
résulte de
la proposition~\ref{propcondresidutors}
que l'on a une suite exacte 

$$\diagram
0\rto &(\ZZ/\ell\ZZ)^{2g(x)}\rto &\H^1((X,x)\et,\mu_\ell)\rrto^{\prod \theta_a}&&\bigoplus\limits_{a\in \br X x}
\ZZ/\ell\ZZ\rrto^\Sigma&&\ZZ/\ell\ZZ\rto&0\enddiagram.$$

\deux{mueltorstyp3}
Supposons~$x$ de type~$3$, et fixons un voisinage coronaire~$U$
de~$x$. On a~$g(x)=0$ et~$\br X x$ est de cardinal~$2$. On dispose donc
en vertu de la suite exacte ci-dessus d'un
isomorphisme 
entre~$\H^1((X,x)\et,\mu_\ell)$ et~$\ZZ/\ell \ZZ$, qui est non canonique si~$\ell >2$ : il
n'est déterminé
qu'au signe près et dépend du choix d'une des deux branches issues de~$x$, c'est-à-dire 
du choix d'une orientation de~$U$ (à la branche~$a$ est associée l'orientation de~$U$ qui induit
l'orientation de la sous-couronne~$a(U)$ vers le point~$x$). 

\deux{enfaitmueltype3}
En fait, l'existence d'un isomorphisme~$\H^1((X,x)\et,\mu_\ell)\simeq\ZZ/\ell \ZZ$
déterminé au signe près et dépendant du choix d'une orientation de~$U$ peut se montrer directement, 
sans référence à la suite exacte de~\ref{resumtheotors}, et sans même
supposer que
le point~$x\in X\typ 3$
est intérieur. 

En effet, supprimons donc l'hypothèse que~$x$ est intérieur, 
et désignons toujours par~$U$ un voisinage coronaire de~$x$. 
Puisque~$x$ est de type~3, 
le singleton~$\{x\}$
est une sous-couronne de~$U$. On dispose donc
d'isomorphismes 

$$\H^1((X,x)\et,\mu_\ell)\simeq\hres(x)\ti/(\hres(x)\ti)^\ell\simeq
\kum(\{x\})\simeq \kum(U)\simeq \ZZ/\ell\ZZ.$$
Tous ces isomorphismes sont canoniques, hormis le dernier qui n'est déterminé
qu'au signe près et dépend du choix d'une orientation sur~$U$, d'où notre assertion.

\medskip
Notons par ailleurs que
choisir une orientation sur~$U$ revient à en choisir une sur~$\{x\}$, 
c'est-à-dire encore à choisir un générateur du
groupe~$|\hres(x)\ti|/|(\KK)\ti|$ (qui est libre de rang~$1$). 

\deux{torsgratteciel}
{\bf Lemme.}
{\em Soit~$X$ une courbe~$\KK$-analytique
quasi-lisse, et soit~$x$
un point de~$X\typ 2$ de genre~$>0$. Il existe
un~$\mu_\ell$-torseur
sur~$X$ dont la fibre en~$x$ est connexe
et qui se trivialise au-dessus de~$Y\setminus \{x\}$.}

\medskip
{\em Démonstration.} 
Comme~$x$ est de type 2, le groupe~$|\hres(x)\ti|$ est divisible, 
et~$\H^1(\hres(x),\mu_\ell)\simeq \H^1(\red{\hres(x)_1}, \mu_\ell)$. Soit~$\sch C$
la courbe résiduelle de~$x$. Comme elle est de genre~$g(x)>0$,
on déduit de~\ref{h1genreg}
qu'il existe une classe~$\red h\in \H^1(\red{\hres(x)_1}, \mu_\ell)$
qui est d'ordre~$\ell$ et est telle que~$\delta_{\sch P}(\red h)=0$ pour
tout point fermé~$\sch P$ de~$\sch C$. 

Soit~$h$ la classe de~$\H^1(X,x)\et,\mu_\ell)\simeq \H^1(\hres(x),\mu_\ell))$
qui correspond à~$\red h$. Elle est d'ordre~$\ell$, et l'on a~$\theta_a(h)=0$
pour toute branche~$a\in \br X x$ en vertu de~\ref{thetares}. D'après~\ref{prepacondresidu}, 
il existe
un voisinage ouvert~$U$ de~$x$ bien découpé autour de~$x$ et un~$\mu_\ell$-torseur de Kummer~$V\to U$
dont le germe en~$x$ a pour classe~$h$, et qui est trivial au-dessus de toute
composante connexe de~$U\setminus \{x\}$ relativement compacte dans~$U$. Par ailleurs,
pour toute branche~$a\in \br X x$, le torseur~$V\times_U a(U)$ est de Kummer, d'invariant~$\theta_a(h)=0$. 
En conséquence, $V$ se trivialise également au-dessus de toute composante
coronaire de~$U\setminus\{x\}$, 
d'où finalement la trivialité de~$V\times_U (U\setminus \{x\})\to (U\setminus \{x\})$. En recollant~$V$
et~$\mu_\ell\times_{\KK} (X\setminus \{x\})$
par un isomorphisme arbitraire au-dessus de~$U\setminus \{x\}$ on obtient un~$\mu_\ell$-torseur~$Y\to X$, trivial
au-dessus de~$X\setminus\{x\}$, et dont le germe en~$x$ a pour classe~$h$. Comme l'ordre de~$h$ est
égal à~$\ell$, on déduit de~\ref{cardfibremutors}
que~$Y_x$ est connexe, ce qui achève la démonstration.~$\Box$ 
\medskip

\subsection*{Théorèmes de finitude
globale}

\deux{theopptdisc} {\bf Théorème.} {\em Soit~${\sch X}$ une~$\KK$-courbe algébrique projective, irréductible et lisse et soit~$x\in {\sch X}\an\typ 2$. Presque toutes les composantes connexes de~${\sch X}\an\setminus\{x\}$ sont des disques.}

\medskip
{\em Démonstration.} Soit~${\sch C}$ la courbe résiduelle en~$x$. Soit~$f$ une fonction rationnelle sur~${\sch X}$ telle que~$|f(x)|=1$ et telle que~$\red {f(x)}$ induise un morphisme génériquement étale de~${\sch C}$ sur~$\PP^1_{\kk}$ ; soit~$\phi : {\sch X}\an\to \pkk$ le morphisme fini et plat induit par~$f$. Soit~$\Pi$ l'ensemble des composantes connexes~$V$ de~${\sch X}\an\setminus\{x\}$ satisfaisant les conditions suivantes : 

\medskip
1)~$\br {{\sch X}\an}x \ctd V$ est un singleton~$\{b_V\}$ ;

2) l'application~${\sch C}\to \PP^1_{\kk}$ induite par~$\red {f(x)}$ est non ramifiée en le point fermé de~$\sch C$ qui correspond à~$b_V$ ; 

3)~$V$ ne contient aucun antécédent de~$\phi(x)$. 

\medskip
Il résulte du fait que~${\sch C}\to \PP^1_{\kk}$ est génériquement étale que presque toutes les composantes connexes de~${\sch X}\an\setminus\{x\}$ appartiennent à~$\Pi$. Soit~$V$ appartenant à~$\Pi$. La condition 3) implique, d'après le lemme~\ref{imtoutoucomp}, que~$V$ est une composante connexe de~$\phi\inv(U)$. On déduit alors de~\ref{vdanspitriv} que la flèche~$V\to U$ 
est un isomorphisme. En tant que composante connexe du complémentaire d'un point de type 2 dans~$\pkk$, l'ouvert~$U$ est un disque, et~$V$ est donc un disque, ce qui achève la démonstration.~$\Box$

\deux{lemrevg0} {\bf Lemme.} {\em Soit~${\sch X}$ une~$\KK$-courbe algébrique projective, irréductible et lisse, et soit~$x$ un point de~${\sch X}\an$ de genre zéro ; supposons que pour toute composante connexe~$V$ de~${\sch X}\an\setminus\{x\}$ l'ensemble~$\br{ {\sch X}\an} x \ctd V$ est un singleton.

\medskip
i) Si~$Y\to {\sch X}\an$ est un revêtement fini étale déployé au-dessus de~${\sch X}\an\setminus\{x\}$, il est trivial. 
 
 ii) Supposons qu'il existe un entier~$\ell$ premier à~$p$ et strictement supérieur à~$1$ tel que tout~$\mu_\ell$-torseur étale au-dessus de~${\sch X}\an\setminus\{x\}$ soit déployé ; la courbe~${\sch X}$ est alors isomorphe à~$\PP^1_{\KK}$.}
 
 \medskip
 {\em Démonstration.} Soit~$Y\to {\sch X}\an$ un revêtement fini étale déployé au-dessus de~${\sch X}\an\setminus\{x\}$ et soit~$y$ un antécédent de~$x$ sur~$Y$. Soit~$b\in \br Y y$ et soit~$a$ l'image de~$b$ dans~$\br X x$. Comme~$Y$ est déployé au-dessus de~${\sch X}\an\setminus\{x\}$, le degré de~$b$ au-dessus de~$a$ est égal à~$1$. 
 
 \medskip
 Si~$x$ est de type 3, il découle de l'assertion iii) du th.~\ref{theohenbr} et de la proposition~\ref{propbondeg} que~$\kappa(y)\simeq \kappa(x)$.
 
 \medskip
 Si~$x$ est de type 2, la courbe résiduelle en~$x$ est de genre zéro, c'est-à-dire isomorphe à~$\PP^1_{\kk}$. Il découle du théorème~\ref{ramresi} que la la courbe résiduelle en~$y$ est étale sur la courbe résiduelle en~$x$ ; comme le schéma~$\PP^1_{\kk}$ est simplement connexe,~$\red{\kappa(y)}\simeq \red{\kappa(x)}$ ; en utilisant encore la proposition~\ref{ramresi}, il vient~$\kappa(y)\simeq \kappa(x)$. 
 
\medskip
On a donc dans tous les cas~$\kappa(y)\simeq \kappa(x)$ pour tout antécédent~$y$ de~$x$. La fibre de~$Y$ au-dessus de~$x$ est par conséquent triviale, ce qui implique que~$Y$ est déployé sur un voisinage de~$x$ ; dès lors ,~$Y\to{\sch X}\an$ est un
revêtement
{\em topologique}.

\medskip
Comme ~$\br{ {\sch X}\an} x \ctd V$ est un singleton pour toute composante connexe~$V$ de~${\sch X}\an\setminus\{x\}$, le point~$x$ n'est situé sur aucune boucle du graphe~${\sch X}\an$. Le
revêtement topologique~$Y\to {\sch X}\an$ est donc déployé au-dessus de toutes les boucles de~$\sch X\an$, et partant trivial, ce qui montre i). 

 \medskip
 Supposons qu'il existe un entier~$\ell$ comme dans ii), et soit~${\sch Y}$ un~$\mu_\ell$-torseur étale sur~${\sch X}$. On déduit de i) et de l'hypothèse faite sur~$\ell$ que le revêtement~${\sch Y}\an\to {\sch X}\an$ est trivial ; par GAGA, il en va de même de~${\sch Y}\to {\sch X}$. Autrement dit,~$\H^1({\sch X}_{\tiny\mbox{ét}},\mu_\ell)=0$ ; par conséquent,~${\sch X}\simeq \PP^1_{\KK}$.~$\Box$

\deux{theogpos} {\bf Théorème.}
{\em Soit~$\sch X$ une~$\KK$-courbe algébrique projective, 
irréductible et lisse. L'ensemble des points de~${\sch X}\an$ de genre strictement positif est fini.}

\medskip
{\em Démonstration.} Soit~$\ell$ un entier premier à~$p$
et supérieur
ou égal à~$2$. Soit~$x\in \sch X\an$ un point
de genre~$>0$. Le lemme~\ref{torsgratteciel}
assure l'existence d'un~$\mu_\ell$-torseur~$Y(x)\to \sch X\an$ dont la fibre
en~$x$ est connexe, et qui est trivial au-dessus de~$\sch X\an\setminus \{x\}$. Si~$x$ et~$x'$ sont deux points
{\em distincts}
de~$\sch X\an$ 
de genre~$>0$, les torseurs~$Y(x)$ et~$Y(x')$ sont non isomorphes
(la fibre en~$x$ de~$Y(x)$ est connexe alors que la fibre en~$x$ de~$Y(x')$
est triviale).  Combinée à GAGA,
la flèche~$x\mapsto Y(x)$ induit
une
{\em injection} de l'ensemble des points 
de genre~$>0$ 
de~$\sch X\an$ dans le groupe fini~$\H^1({\sch X}_{\tiny\mbox{ét}},\mu_\ell)$.~$\Box$

\subsection*{Voisinages d'un point de type 4 et revêtements des disques et couronnes}

\deux{introvois4} Nous fixons pour tout ce paragraphe les notations suivantes : 

\medskip

$\bullet$~$\ell$ est un nombre premier différent de~$p$ ;

$\bullet$~${\sch X}$ est une~$\KK$-courbe algébrique projective, irréductible et lisse ; 

$\bullet$~$x$ est un point de~${\sch X}\an\typ 4$.

\deux{basevoist4} Comme~$x$ est de type 4, c'est d'après le théorème~\ref{lisseval} un point unibranche du graphe~$\sch X\an$ ; l'ensemble des points de genre~$>0$ de~${\sch X}\an$ étant par ailleurs fini (th.~\ref{theogpos}), le point~$x$ possède une base de voisinages~$D$ satisfaisant les deux conditions qui suivent :

\medskip
i)~$\partial D$ est un singleton~$\{\eta\}$ et~$\overline D= D\cup \{\eta\}$ est un arbre compact  ; 

ii)~$D$ ne contient aucun point de genre~$>0$. 

\deux{h1ettriv} {\bf Proposition.} {\em Soit~$D$ comme ci-dessus ; tout~$\mu_\ell$-torseur étale sur~$D$ est trivial}. 

\medskip
{\em Démonstration.} Soit~$Y\to D$ un~$\mu_\ell$-torseur étale, et soit~$U$ l'ensemble des points~$t$ de~$D$ tels que la fibre de~$Y$ en~$t$ soit triviale ; c'est un ouvert de~$D$. Le but de ce qui suit est de démontrer que~$U=D$ : cela assurera que~$Y\to D$ est un~$\mu_\ell $-torseur {\em topologique}, et partant trivial puisque~$D$ est un arbre. 

\medskip
On raisonne par l'absurde, en supposant que ce n'est pas le cas. On munit~$\overline D$ de l'ordre défini par~$\eta$. Soit~${\sch T}$ une chaîne non vide de points  de~$D-U$. Par compacité, le filtre de~$\overline D-U$ engendré par les sections commençantes de~$\sch T$ admet un point adhérent ; on vérifie aussitôt qu'il est différent de~$\eta$ et minore~$\sch T$ (c'en est plus précisément la borne inférieure). 

Il s'ensuit, par le lemme de Zorn, qu'il existe un point~$t$ de~$D-U$ qui est minimal. Comme~$t\in D-U$, la fibre de~$Y$ en~$t$ n'est pas déployée. Notons~$b$ la branche issue de~$t$ définie par l'intervalle~$[t;\eta]$. La fibre du~$\mu_\ell$-torseur~$Y$ en~$t$ n'étant pas déployée, c'est un singleton~$\{t'\}$, et~$\kappa(t')$ est une extension
cyclique de degré~$\ell$
de~$\kappa(t)$.

\medskip
{\em Le point~$t$ ne peut être de type 1 ou 4}. Cela résulte du fait que si~$F$ est une extension immédiate du corps algébriquement clos~$\KK$, toute extension de degré premier à~$p$ de~$\KK$ est triviale. 

\medskip
{\em Le point~$t$ ne peut être de type 3}. Supposons qu'il le soit. Il existe alors
un voisinage coronaire~$Z$ dans~$D$ au-dessus duquel~$Y$ est de Kummer. Comme~$Y_x$ est connexe, $Y\times_D Z$ est connexe
et est donc une couronne (\ref{revkummcor}
{\em et sq.}). Le morphisme~$Y\times_D Z\to Z$ est fini et plat ; il s'ensuit
en vertu de~\ref{morcourintersk}
et~\ref{caspartfinplat}
que tout point de~$\skel Z$ a un unique antécédent sur~$Y$. 

Par ailleurs, soit~$Z'$ la composante connexe de~$Z\setminus\{t\}$ qui ne contient pas~$b$. Par minimalité de~$t$, les fibres de~$Y$ au-dessus de~$Z'$ sont toutes triviales ; c'est en particulier le cas des fibres de~$Y$ au-dessus de~$\skel {Z'}\subset \skel Z$,
mais celles-ci sont des singletons
d'après ce qui précède ; on aboutit ainsi à une contradiction. 

\medskip
{\em Le point~$t$ ne peut être de type 2}. Supposons qu'il le soit. Le groupe~$|\kappa(t)\ti|$ est divisible ; le corps~$\red{\kappa(t')}$ est donc une extension cyclique de degré~$\ell$
de~$\red{\kappa(t)}$ (une extension modérément ramifiée étant toujours sans défaut, il n'y a pas lieu ici d'invoquer le th.~\ref{ramresi}). Soit~${\sch C}$ (resp.~${\sch C}'$) la courbe résiduelle en~$t$ (resp.~$t'$). 

\medskip
Soit~$\beta$ une branche de~$D$ issue de~$t$ et différente de~$b$ et soit~$Z$ une section de~$\beta$ ; la minimalité de~$t$ assure que les fibres de~$Y$ au-dessus de~$Z$ sont triviales. Par conséquent,~$Y\times_DZ\to Z$ est un~$\mu_\ell$-torseur {\em topologique} ; comme~$Z$ est un arbre, il est trivial. En vertu du théorème~\ref{ramresi}, ceci entraîne que~${\sch C}'\to {\sch C}$ est non ramifié au-dessus du point fermé de~$\sch C$ qui correspond à~$\beta$. 

\medskip
Soit~$\sch P$ le point fermé de~$\sch C$ correspondant à~$b$. D'après ce qui précède, la flèche~${\sch C}'\to \sch C$ induit un revêtement galoisien cyclique d'ordre~$\ell$
de
l'ouvert~${\sch C}\setminus\{\sch P\}$ ; mais~${\sch C}$ est de genre 0 (car~$D$ ne contient aucun point de genre strictement positif), et l'on aboutit ainsi à une contradiction, ce qui achève la démonstration.~$\Box$

\deux{commdisc} {\em Remarque.} La proposition précédente s'applique notamment lorsque~${\sch X}=\PP^1_{\KK}$ et lorsque~$D$ est un disque ; la proposition~\ref{h1ettriv} implique donc que tout~$\mu_\ell$-torseur étale sur un disque est trivial. Ce fait est dû à Berkovich, et notre démonstration est un décalque {\em mutatis mutandis} de sa preuve.

\deux{theovoisdiscalgclos} {\bf Théorème.} {\em Le point~$x$ possède un voisinage ouvert dans~$\sch X\an$ qui est un disque.} 

\medskip
{\em Démonstration.} Donnons-nous un voisinage~$D$ de~$x$ comme au~\ref{basevoist4}. Soit~$Z$ une section coronaire de la branche de~${\sch X}\an$ issue de~$\eta$ et contenue dans~$D$ ; choisissons un point~$\xi$ de~$]x;\eta[$ situé sur~$\mathsf S(Z)$, et appelons~$\Delta$ la composante connexe de~$D\setminus\{\xi\}$ qui contient~$x$. Remarquons que le voisinage~$\Delta$  de~$x$ vérifie lui aussi les conditions énoncées au~\ref{basevoist4} ; par conséquent, la proposition~\ref{h1ettriv} s'applique à~$\Delta$ et garantit que tout~$\mu_\ell$-torseur étale sur ce dernier est trivial. 

\medskip
Les composantes connexes de~$Z\setminus\{\xi\}$ sont :

\medskip
$\bullet$ la sous-couronne~$]\xi;\eta[^\flat$ de~$Z$ ;

$\bullet$ la sous-couronne~$(\mathsf S(Z)-[\xi;\eta[)^\flat$ de~$Z$, qui est contenue dans~$\Delta$ ;

$\bullet$ des disques aboutissant à~$\xi$. 

\medskip
Il en résulte que les composantes connexes de~$D\setminus\{\xi\}$ sont :

\medskip
$\bullet$ la sous-couronne~$]\xi;\eta[^\flat$ de~$Z$ ;

$\bullet$ l'ouvert~$\Delta$ ;

$\bullet$ des disques contenus dans~$Z$ et aboutissant à~$\xi$ ; 

\medskip
Comme~$D$ est un arbre, 
chacune de ces composantes contient une unique branche issue de~$\xi$. 

En prolongeant la couronne~$]\xi;\eta[^\flat$ en un disque, on immerge~$D$ dans une courbe~$\KK$-analytique propre, lisse et connexe qui est donc l'analytification d'une~$\KK$-courbe algébrique projective, lisse et intègre~${\sch Y}$. 

\medskip
Soit~$W$ une composante connexe de~${\sch Y}\an\setminus\{\xi\}$. Il découle de ce qui précède que~$W$ contient une et une seule branche issue de~$\xi$, et que~$W$ est ou bien égale à~$\Delta$, ou bien un disque ; ce dernier fait implique que tout~$\mu_\ell$-torseur étale sur~$W$ est trivial (prop.~\ref{h1ettriv} et rem.~\ref{commdisc}). Comme par ailleurs le point ~$\xi$ est de genre~$0$ (il appartient à~$\mathsf S(Z)$),  les hypothèses du lemme~\ref{lemrevg0} sont satisfaites ; il en résulte que~${\sch Y}\simeq \pkk$.

\medskip
L'ouvert~$\Delta$ s'identifie dès lors à une composante connexe du complémentaire d'un point de type 2 ou 3 dans~$\pkk$ ; par conséquent,~$\Delta$ est un disque.~$\Box$

\medskip
Terminons par une proposition (due à Berkovich) et un lemme qui résultent de la remarque~\ref{commdisc}. 

\deux{corollkum}{\bf Proposition (Berkovich).} {\em Soit~$X$ une~$\KK$-couronne et soit~$\ell$ un entier premier à~$p$. Tout~$\mu_\ell$-torseur étale sur~$X$ est de Kummer.}

\medskip
{\em Démonstration.}
Soit~$h\in \H^1(X_{\tiny\mbox{ét}},\mu_\ell)$. On sait que~$h$ 
est de Kummer au voisinage de tout point de~$X$ (\ref{torslocalkum}). Nous allons tout d'abord
préciser un peu cette assertion en montrant que
pour tout~$x\in \skel X$ il existe un entier~$r(x)$, unique modulo~$\ell$, tel que~$h=(t^{r(x)})$ au voisinage de~$x$. 

\trois{lockumsquel} {\em Existence de~$r(x)$.}

\medskip
{\em Le cas où~$x$ est de type 3}. Les sous-couronnes
ouvertes de~$X$ contenant~$x$ forment alors
une base de voisinages de~$x$ ; il existe donc une telle
sous-couronne~$Y$
telle que~$h_{|Y}$ soit de Kummer. Comme~$T_{|Y}$ est une
fonction coordonnée sur~$Y$, la classe~$h_{|Y}$ est bien de la forme~$(T^{r(x)})$
pour un certain entier~$x$. 

\medskip
{\em Le cas où~$x$ est de type 2.}
Le groupe~$|\hres(x)\ti|$ étant divisible, on dispose d'un isomorphisme
naturel~$\H^1(\hres(x),\mu_\ell)\simeq \H^1(\red{\hres(x)}_1,\mu_\ell)$ ; soit~$\red h$
la classe qui correspond par ce biais au germe de~$h$ en~$x$. 

\medskip
Comme~$|\hres(x)\ti|=|\KK\ti|$, on peut, quitte à multiplier~$T$ par
un élément convenable de~$\KK\ti$, supposer que~$|T(x)|=1$ ; la courbe résiduelle~$\sch C$
en~$x$ s'identifie alors à~$\PP^1_{\kk_1}$, et~$\red{T(x)}$ en est une fonction coordonnée. Notons 
son diviseur~$[0]-[\infty]$ ; les branches issues de~$x$
correspondant aux points fermés~$0$
et~$\infty$ sont précisément les deux branches issues de~$x$ et définies par~$\skel X$
(puisque~$|T|$ est localement constante en dehors de~$\skel X$ et strictement monotone sur~$\skel X$). 

\medskip
Soit~$a\in \br X x$ une branche qui n'est pas contenue dans~$\skel X$, soit~$U$ la composante
connexe de~$X\setminus \{x\}$ qui contient~$a$, et soit~$\sch P$
le point fermé de~$\sch C$ correspondant à~$a$. 
Comme~$U$ est un disque, on déduit de
la proposition~\ref{h1ettriv}
et de~\ref{commdisc}
que~$h_{|U}=0$ ; en particulier, la restriction de~$h$ à toute section coronaire
de~$a$ contenue dans~$U$ est nulle, ce qui signifie que~$\theta_a(h_{|(X,x)})=0$. 
Par conséquent, $\delta_{\sch P}(\red h) =0$ (\ref{thetares}). 

\medskip
Fixons un entier~$r(x)$ tel que~$\delta_0(\red h)=r(x)$ modulo~$\ell$.
D'après~\ref{suiteexacteresidus}, on a~$\delta_{\infty}(\red h)=-r(x)$ modulo~$\ell$. La classe
$\red h$ a ainsi les mêmes résidus que~$\left(\red{T(x)}^{r(x)}\right)$ en tout point fermé de
la courbe~$\sch C$.
Celle-ci étant
de genre~$0$, on en déduit que~$\red h=\left(\red{T(x)}^{r(x)}\right)$ (\ref{h1genreg}), 
et partant que~$h_{|(X,x)}=(T^{r(x)})$. 

\trois{rxmodellunique}
{\em Unicité de~$r(x)$
modulo~$\ell$.} 
Supposons qu'il existe un voisinage ouvert~$U$ de~$x$ et deux entiers~$r$ et~$s$
tels que~$h_{|U}=(T^r)=(T^s)$. Il existe une sous-couronne
ouverte~$Y$ de~$X$ contenue dans~$U$ (il suffit de prendre~$Y$ de la forme~$I^\flat$,
où~$I$ est un intervalle ouvert aboutissant à~$x$, tracé sur~$\skel X\cap U$, et suffisamment petit). 
L'égalité~$(T^r)=(T^s)$ dans~$\kum (Y)$ implique alors que~$r$ et~$s$ sont égaux modulo~$\ell$. 

\trois{conclusionkumcour}
{\em Conclusion.} Il résulte de la définition de~$r(x)$ et de son unicité
modulo~$\ell$ que~$x\mapsto r(x)$ définit sans ambiguité une application
localement constante de~$\skel X$ vers~$\ZZ/\ell\ZZ$.
Par connexité de~$\skel X$, cette application est constante ; soit~$r$ sa valeur.

Posons~$h'=h-(T^r)$, et soit~$Z\to X$ un~$\mu_\ell$-torseur de classe~$h'$. Par construction, $h'$ est nulle
au voisinage de tout point~$x$ de $\skel X$. Et comme toute composante connexe
de~$X\setminus \skel X$ est un disque, il résulte de la proposition~\ref{h1ettriv}
et de~\ref{commdisc}
que~$h'$ est triviale au-dessus de~$X\setminus \skel X$. En conséquence, $h'$ est triviale
au voisinage de tout point de~$X$, ce qui implique que~$Z$ est un~$\mu_\ell$-torseur topologique. 
La couronne~$X$ étant un arbre, $Z$ est trivial,
$h'=0$ et~$h=(T^r)$,
ce qui achève la démonstration.
~$\Box$

\deux{lemrevg0d} {\bf Lemme.} {\em Soit~$\sch Y$ une~$\KK$-courbe projective, irréductible et lisse. Supposons qu'il existe un point~$y$ de genre~$0$ sur~$\sch Y\an$ tel que~$\sch Y\an\setminus\{y\}$ soit réunion disjointe de disques. La courbe~$\sch Y$ est alors isomorphe à~$\PP^1_{\KK}$.}

\medskip
{\em Démonstration.} Fixons un entier~$\ell$ au moins égal à~$2$ et premier à~$p$. Si~$U$ est une composante connexe de~$\sch Y\an\setminus\{y\}$ c'est un disque, et donc un arbre à un bout. Par conséquent,~$\br {\sch Y\an} y \ctd U$ est un singleton, et la remarque ~\ref{commdisc} assure que~$\H^1(U,\mu_\ell)=0$ ; l'assertion requise découle alors de l'énoncé ii)  du lemme~\ref{lemrevg0}.~$\Box$ 

\section{Étude locale des courbes quasi-lisses sur un corps de base quelconque}

\subsection*{Étude locale des courbes quasi-lisses sur un corps de base quelconque : préliminaires}

Soit~$X$ une courbe~$k$-analytique. 

\deux{shilorb} Soit~$x\in X\dtr$. L'ensemble des antécédents de~$x$ sur~$X_{\KK}$ est alors fini.  en effet, on peut toujours supposer, quitte à le remplacer par le complété de sa clôture parfaite, que~$k$ est parfait, puis que~$X$ est réduite et par conséquent génériquement quasi-lisse, et en particulier quasi-lisse en~$x$. 

Il est alors loisible, en substituant à~$X$ un voisinage affinoïde convenable de~$x$, puis en immergeant ce dernier dans l'analytification d'une courbe affine lisse, de supposer que~$X$ est elle-même l'analytification d'une courbe affine. On peut maintenant procéder de deux façons différentes.

\trois{finorbgrth} {\em Première méthode.} Tout antécédent de~$x$ sur~$X_{\KK}$ est de type 2 ou 3, et est en particulier pluribranche (\ref{casbij2} et~\ref{casbij3}) ; il a donc a une orbite finie sous l'action de Galois. 
 
\trois{finorbshdsc} {\em Seconde méthode.} On se ramène, par normalisation de Noether, au cas où~$X=\Aff^{1,{\rm an}}_k$. Soit~$y$ un antécédent de~$x$ sur~$\Aff^{1,{\rm an}}_{\KK}$. Le point~$y$ est le bord de Shilov d'un disque fermé~$D$, qui possède un point rigide~$z$ ; si~$\mathsf H$ désigne le stabilisateur de~$z$ dans~$\mathsf G$ alors~$\mathsf H$ est un sous-groupe ouvert de~$\mathsf G$ qui stabilise~$D$, et partant~$y$.

\deux{finorbcons} Si~$x\in X\geom$ la fermeture séparable de~$k$ dans~$\kappa(x)$ sera notée~$\got s(x)$ ; c'est en vertu de ce qui précède une extension {\em finie} séparable de~$k$. Si~$x\in X\dtr$ et si~$y$ est l'un des antécédents de~$x$ sur~$X_{\KK}$, son genre ne dépend que de~$x$ ; nous l'appellerons encore le genre de~$x$ et le noterons~$g(x)$ ; {\em on prendra garde qu'il n'a aucune raison de coïncider avec le genre de la courbe résiduelle en~$x$, dont la formation ne commute pas à l'extension des scalaires en général.}

\deux{fermdiscgpos} L'ensemble~$\mathsf E$ des points de genre strictement positif de~$X$ est une partie fermée et discrète de~$X$. Pour le voir, on se ramène aussitôt au cas où~$X$ est compacte, et il s'agit alors de montrer que~$\mathsf E$ est fini. Il suffit de le faire après extension des scalaires à~$\KK$, ce qui autorise à supposer~$k$ algébriquement clos ; quitte à remplacer~$X$ par sa normalisée, on se ramène ensuite au cas où~$X$ est quasi-lisse, et donc réunion finie de domaines affinoïdes de courbes lisses ; on conclut alors à l'aide
du théorème~\ref{theogpos}.

\subsection*{Bases de voisinages sur une courbe quasi-lisse}

\deux{theovoisql}{\bf Théorème.} {\em Soit~$X$ une courbe~$k$-analytique quasi-lisse et soit~$x\in X$. Il possède un voisinage ouvert connexe~$V$ dans~$X$ qui satisfait les conditions énoncées ci-dessous.

\begin{itemize}

\medskip
\itb Le graphe~$\overline V$ est un arbre compact. 

\itb Si~$x\in X\typ{14}$ alors~$V$ est un disque virtuel. 

\itb Si~$x\in X\typ 3$ le voisinage~$V$ est coronaire.

\itb Si~$x\in X\typ 2$ alors~$V\setminus\{x\}$ est réunion disjointe de disques virtuels et d'un nombre fini de couronnes gentiment virtuelles.\end{itemize}}

\medskip
{\em Démonstration.} On traite séparément chacun des cas. 

\trois{voist1ql} {\em Le cas où~$x\in X\typ {14}$}. Choisissons un antécédent~$y$ de~$x$ sur~$X_{\KK}$. Comme~$y$ est de type 1 ou 4, c'est un point de l'intérieur analytique de~$X_{\KK}$ et il possède donc un voisinage ouvert dans~$X_{\KK}$ qui s'immerge dans l'analytification d'une~$\KK$-courbe projective et lisse ; en vertu du théorème~\ref{theovoisdiscalgclos}, il existe alors un voisinage~$Y$ de~$y$ dans~$X_{\KK}$ qui est un~$\KK$-disque. D'après la proposition~\ref{corollstabbr} on peut, quitte à le restreindre, supposer que~$Y$ possède les deux propriétés suivantes :

\medskip
$\bullet$ son adhérence dans~$X_{\KK}$ est un arbre compact ; 

$\bullet$ pour tout~$g\in \mathsf G$ l'on a~$g(Y)=Y$ ou~$g(Y)\cap Y=\emptyset$. 

\medskip
Soit~$V$ l'image de~$Y$ sur~$X$ ; c'est un voisinage de~$x$ dont l'adhérence dans~$X$ est un arbre compact. Comme~$\{g(Y)\}_{g\in \mathsf G}$ s'identifie à l'ensemble fini~$\pi_0(Z_{\KK})$, l'espace~$V_{\KK}$ est réunion disjointe de~$\KK$-disques, et~$V$ est de ce fait un disque virtuel. 

\trois{voist2l}  {\em Le cas où~$x\in X\typ 2-\brdan X$}. Le point~$x$ possède un voisinage~$V$ dans~$X$ tel que~$\overline V$ soit un arbre compact, et tel que~$V$ s'immerge dans le lieu lisse de~${\sch X}\an$ pour une certaine~$k$-courbe projective~${\sch X}$ ; dès lors~$V_{\KK}$ s'immerge dans~${\sch Y}\an$, où~$\sch Y$ est la normalisée de~${\sch X}_{\KK}$.

\medskip
Comme~$x$ est de type 2, l'ensemble de ses antécédents sur~$V_{\KK}$ est une partie finie~$S$ de~$V_{\KK}\subset {\sch Y}\an$, constituée de points de type 2. La courbe~${\sch Y}\an$ étant un graphe compact, presque toute composante connexe de~$V_{\KK}\setminus S$ est une composante connexe de~${\sch Y}\an\setminus\{y\}$ pour un certain~$y\in S$ ; il résulte alors du théorème~\ref{theopptdisc} que presque toutes les composantes connexes de~$V_{\KK}\setminus S$ sont des~$\KK$-disques ; par conséquent, pour presque toute composante connexe~$W$ de~$V\setminus\{x\}$, l'espace~$W_{\KK}$ est réunion disjointe de~$\KK$-disques, ce qui signifie que~$W$ est un disque virtuel. 

\medskip
Soit~$W$ une composante connexe de~$V\setminus\{x\}$ qui n'est pas un disque virtuel ; elle contient une unique branche de~$Z$ issue de~$x$ ; en vertu du lemme~\ref{brcourvirt}, cette branche possède une section qui est une couronne gentiment virtuelle. On peut donc, en rabotant un nombre fini de composantes connexes de~$V\setminus\{x\}$, faire en sorte que~$V$ possède la propriété requise. 

\trois{voist23ql} {\em Conclusion.} Si~$x\in X\typ 3$ l'existence d'un voisinage coronaire de~$x$ découle du théorème~\ref{sectioncour}. Supposons maintenant que~$x\in X\typ 2$. Il possède un voisinage ouvert dans~$X$ dont l'adhérence est un arbre compact, et qui s'identifie à un domaine analytique fermé d'une courbe lisse~$X'$. L'existence d'un voisinage~$V$ de~$x$ dans~$X$ qui soit de la forme voulue provient alors du cas sans bord traité au~\ref{voist2l}) ci-dessus et du corollaire~\ref{coroldomferm}.~$\Box$

\deux{commentbordsx} {\em Remarque}. Si~$x\in X\typ 3$ alors~$V$ est ou bien une couronne virtuelle dont~$x$ appartient au squelette  (si~$I=\{1,2\}$), ou bien une \og demi-couronne virtuelle de bord~$\{x\}$\fg~ (si~$I$ est un singleton), ou bien réduit au singleton~$\{x\}$  (si~$I$ est vide) ; le premier cas correspond à celui où~$x$ appartient à l'intérieur analytique de~$X$ (et~$V$ en est alors un voisinage coronaire), le dernier à celui où il en est un point isolé.

\subsection*{Toise canonique et paracompacité des courbes analytiques}

\deux{lemme-typ14-discret}
{\bf Lemme.}
{\em Soit~$X$ une courbe~$k$-analytique et soit~$\Gamma$ un sous-graphe localement fini de~$X$. Il existe un sous-ensemble
fermé et discret~$\Sigma$ de~$\Gamma$ tel que $\Gamma \setminus \Sigma \subset X\dtr$.}

\medskip
{\em Démonstration}. Soit~$S$ l'ensemble des points de~$X$ ayant au moins deux antécédents sur le normalisé de~$X$ ; c'est un sous-ensemble fermé et discret de~$X$ contenu dans~$X\typ 0$. Soit~$\Sigma$ la réunion de~$S\cap \Gamma$ et de l'ensemble des sommets de~$\Gamma$ ; c'est une partie fermée et discrète de~$\Gamma$. Si~$x\in \Gamma-\Sigma$ alors~$\br \Gamma x$ est de cardinal 2, et~$\br X x$ est donc de cardinal au moins égal à 2. Par conséquent,~$x\in X\geom$ ; comme par ailleurs~$x$ n'appartient pas à~$S$, il ne peut être rigide (il serait sinon unibranche), et appartient finalement à~$X\dtr$. Ainsi,~$\Gamma-\Sigma\subset X\dtr$.~$\Box$ 

\deux{metrinat} {\bf Proposition.} {\em Soit~$X$ une courbe~$k$-analytique. Il existe une toise canonique sur ~$X\dtr$.}

\medskip
{\em Démonstration.} La question est purement topologique. On peut donc, quitte à étendre les scalaires à~$\kparf$
et à remplacer~$X$ par la courbe réduite sous-jacente, supposer~$k$ parfait et~$X$ réduite ; chaque point de~$X\dtr$ est alors quasi-lisse. 

\medskip
Soit~$I$ un segment tracé sur~$X\dtr$ et soient~$x$ et~$y$ ses deux extrémités. En vertu de~\ref{basecoron}, il existe une suite~$x=x_0<x_1<\ldots<x_r=y$ d'éléments de~$I$ (orienté de~$x$ vers~$y$) telle que pour tout~$i$ compris entre~$1$ et~$r-1$ la condition suivante soit vérifiée : {\em  l'intervalle ouvert ~$]x_i;x_{i+1}[$ est faiblement admissible et~$]x_i;x_{i+1}[^\flat$ est une couronne (gentiment) virtuelle
de type~$]**[$.} 

La somme des logarithmes des modules des couronnes virtuelles ~$]x_i;x_{i+1}[^\flat$ne dépend pas du choix de la suite~$(x_i)$ : il suffit en effet de vérifier qu'elle est insensible à un raffinement de la subdivision initiale, ce qui est clair ; on la note~$l(I)$. Il résulte des propriétés élémentaires du module d'une couronne que~$I\mapsto l(I)$ définit une toise sur~$X\dtr$, qui est par construction invariante par tout automorphisme de~$X$.~$\Box$

\deux{corolltoiselocfin} {\bf Corollaire.} {\em Soit~$X$ une courbe~$k$-analytique et soit~$\Gamma$ un sous-graphe localement fini de~$X$. Le graphe~$\Gamma$ admet une toise.}

\medskip
{\em Démonstration.} D'après le lemme~\ref{lemme-typ14-discret}, il existe un sous-ensemble fermé et discret~$\Sigma$ 
de~$\Gamma$ tel que~$\Gamma \setminus \Sigma \subset X\dtr$. En vertu de la proposition~\ref{metrinat} ci-dessus, $\Gamma-\Sigma$ admet alors une toise,
et il s'ensuit d'après la proposition~\ref{longcaslocfin} que~$\Gamma$ admet lui-même une toise.~$\Box$  

\deux{prop-toise-presqualg}
{\bf Proposition.}
{\em
Soit
$X$ une courbe~$k$-analytique et soit~$\Gamma$ un 
un sous-graphe localement fini de~$X\dtr$. Soit~$F$ une extension presque algébrique de~$k$ et soit~$\Gamma_F$
l'image réciproque de~$\Gamma$ (d'après le cor.~\ref{corollimrecih}, $\Gamma_F$ est un sous-graphe localement fini
de~$X_F$, et~$\Gamma_F\to \Gamma$ est injective par morceaux). 
L'application~$\Gamma_F\to \Gamma$ commute aux toises canoniques de ses source et but}.

\medskip
{\em Démonstration}. 
Soit~$F'$ l'extension complète de~$k$ composée de~$F$ et~$\kparf$. Par construction, les toises canoniques sont invariantes par 
extension radicielle complétée des scalaires ; il suffit donc de démontrer que~$\Gamma_{F'}\to \Gamma_{\kparf}$ est une isométrie par morceaux, 
ce qui permet de supposer~$k$ parfait. Le problème considéré étant insensible aux phénomènes de nilpotence, on peut remplacer~$X$ par~$X_{\rm red}$, et donc
la supposer réduite, et partant génériquement quasi-lisse. Soit~$I$ un segment tracé sur~$\Gamma_F$, telle que~$X_F\to X$ induise un homéomorphisme de~$I$ sur son image~$J$. Nous allons démontrer que~$\ell(J)=\ell(I)$, ce qui permettra de conclure. 

Notons~$x$ et~$x'$ les deux extrémités de~$J$, et~$y$ et~$y'$ leurs antécédents
respectifs sur~$I$. En vertu de~\ref{basecoron}, il existe une suite~$x=x_0<x_1<\ldots<x_r=x'$ d'éléments de~$J$ (orienté de~$x$ vers~$y$) 
telle que pour tout~$i$ compris entre~$1$ et~$r-1$ la condition suivante soit vérifiée : {\em  l'intervalle ouvert ~$]x_i;x_{i+1}[$ est faiblement admissible et~$]x_i;x_{i+1}[^\flat$ est une couronne (gentiment) virtuelle de type~$]**[$}. Pour tout~$i$, notons~$y_i$ l'unique antécédent de~$x_i$ sur~$I$. Il suffit de vérifier que
l'on a~$\ell([y_i; y_{i+1}])=\ell([x_i;x_{i+1}])$ pour tout~$i$. On se ramène ainsi au cas où~$]x;x'[$ est faiblement admissible et où~$]x;'x[^\flat$ est une couronne
gentiment virtuelle de type~$]**[$. Notons
$L$ le corps~$\got s(]x; x'[^\flat)$ ; le produit
tensoriel~$F\otimes_k L$ est un produit fini~$\prod F_j$ d'extensions finies séparables de~$F$. 

\medskip
L'ouvert~$]x; x'[^\flat_F $ de~$X_F$  s'identifie à~$\coprod_j ]x; x'[^\flat\times_L F_j$. L'image réciproque de~$]x; x'[$ sur~$X_F$ est donc la réunion 
disjointe des~$\skelan{]x; x'[^\flat\times_L F_j}$ ; il s'ensuit que~$]y;y'[$ est égal à~$\skelan{ ]x; x'[^\flat\times_L F_j}$ pour un certain~$j$. Par définition de la toise canonique, 
$\ell([y;y'])$ est alors le module de la couronne virtuelle~$]x; x'[^\flat\times_L F_j$, qui est égal à celui de~$]x; x'[^\flat$, et partant à~$\ell([x;x'])$, ce qui achève la démonstration.~$\Box$

\deux{theofondmetr} {\bf Théorème.} {\em Soit~$X$ une courbe~$k$-analytique.

\medskip
1) Le graphe~$X$ admet une toise. 

2)  Soit~$\Gamma$ un sous-graphe admissible de~$X$ et soit~$r$ la rétraction canonique de~$X$ sur~$\Gamma$. L'inclusion de~$\Gamma$ dans~$X$ est une équivalence homotopique ; plus précisément il existe une application continue~$h: [0;1]\times X\to X$ telle que~$h(0,.)=r, h(1,.)=\mathsf{Id}_X$ et~$h(t,x)=x$ pour tout~$t\in [0;1]$ et tout~$x\in \Gamma$.

3) Le graphe~$X$ est paracompact, localement fortement contractile et a le type d'homotopie d'un graphe localement fini métrisable ; de plus,~$X$ est contractile si et seulement si~$X$ est un arbre non vide, et il est alors fortement contractile.}

\medskip
{\em Démonstration.} Tout sous-graphe compact de~$X$ admet une toise (\ref{theotoise}) ; comme par ailleurs tout sous-graphe localement fini de~$X$ admet une toise en vertu du corollaire~\ref{corolltoiselocfin}, le graphe~$X$ admet lui-même une toise d'après le lemme~\ref{xadmtoise}. Les assertions 2) et 3) découlent alors directement du théorème~\ref{theoadmtoise}.~$\Box$ 
\subsection*{Variation du corps des constantes sur une courbe}

\deux{conn23rig} Soit~$X$ une courbe~$k$-analytique et soit~$X'$ sa normalisation. Soit~$I$ un intervalle ouvert tracé sur~$X$ et soit~$x\in I$. Comme~$x$ est situé sur un intervalle ouvert,~$\br X x$ est de cardinal au moins égal à deux. Comme~$X\typ {1,4}\cap \brdan X=\emptyset$, il découle de~\ref{lisseval} que~$x\in X\geom$ ; si de plus~$X'\to X$ est un homéomorphisme ({\em e.g} si~$X$ est normale) alors~$\br X y$ est un singleton pour tout~$y\in X\typ 0$ (prop.~\ref{pointrigfinbr}), ce qui entraîne que~$x\in X\dtr$. 

\medskip
Il s'ensuit que~$X\geom$ est une partie convexe de~$X$, et qu'il en va de même de~$X\dtr$ si~$X'\to X$ est un homéomorphisme, et donc notamment si~$X$ est normale.

\deux{apropossx} Soit~$X$ une courbe~$k$-analytique, et soit~$p:X\to X_{\KK}$ la flèche naturelle. Soit~$x\in X\geom$ et soit~$b\in \br X x$. Comme~$x\in X\geom$ sa fibre sur~$X_{\KK}$ est finie ; le corollaire~\ref{corollimreci} assure l'existence d'un intervalle~$J_0\in \intera X b$ (nécessairement tracé sur~$X\geom$ en vertu du~\ref{conn23rig} ci-dessus) tel que pour tout~$J\in \inter X x \ctd{J_0}$ et tout~$y\in J$ le~$\mathsf G$-ensemble fini~$p\inv(y)$ s'identifie canoniquement à~$\pi_0(p\inv(J^\flat))$. 

\medskip
Par conséquent, pour tout~$J\in\inter X x \ctd{J_0}$ et tout~$y\in J$ l'inclusion naturelle~$\got s(J^\flat)\hookrightarrow \got s(y)$ est un isomorphisme ; cela entraîne immédiatement que pour tout couple~$(J,J')$ d'intervalles appartenant à~$\inter X x \ctd{J_0}$ et tels que~$J'\subset J$ l'inclusion naturelle~$\got s(J^\flat)\hookrightarrow \got s((J')^\flat)$ est un isomorphisme. 

\medskip
Posons~$\got s(b)=\lim\limits_{\stackrel\longrightarrow{Z\in \sbr b}}\got s(Z)$ ; en vertu de ce qui précède,~$\got s(b)$ est une extension finie séparable égale à~$\got s(J^\flat)$ pour n'importe quel~$J\in \inter X x \ctd{J_0}$, et si~$y\in J_0$ alors~$\got s(y)=\got s(b)$. 

\medskip
Les anneaux locaux de~$X$ étant henséliens,~$\got s(x)$ s'immerge dans~$\sch O_{X,x}$, et se plonge donc naturellement dans~$\got s(b)$. 

\deux{fontsxgr} Soit~$X$ une courbe~$k$-analytique et soit~$\Gamma$ un sous-graphe fermé et localement fini de~$X$ tracé sur~$X\geom$. 

\trois{imrecgamgeom} Chacun des points  de~$\Gamma$ a un nombre fini d'antécédents ; si~$F$ est une extension presque algébrique de~$k$, les assertions suivantes se déduisent alors de~\ref{corollimrecih} : l'image réciproque~$\Gamma_F$ de~$\Gamma$ sur~$X_F$ est un sous-graphe fermé et localement fini de~$X_F$ ; et la flèche~$\Gamma_F\to \Gamma$ est injective par morceaux. 

\trois{compsxsurgamma} Soit~$d$ l'application~$x\mapsto [\got s(x):k]$ de~$\Gamma$ dans~$\NN$. On  déduit de~\ref{apropossx} que~$d$ est constante par morceaux et semi-continue inférieurement, et que si~$x\in \Gamma$ les assertions suivantes sont équivalentes : 

\medskip
\medskip
i)~$d$ est continue en~$x$ ;

ii)~$d$ est constante au voisinage de~$x$ ;

iii) pour tout point~$y$ de~$\Gamma$ suffisamment proche de~$x$ l'on a~$\got s(x)\simeq \got s(y)$ ; 

iv) pour toute branche~$b\in \br \Gamma x$ la flèche naturelle~$\got s(x)\to \got s(b)$ est un isomorphisme.

\deux{introremseccorong} Soit~$X$ une courbe~$k$-analytique quasi-lisse et soit~$x\in X\dtr$.

\trois{remseccorong} Soit~$b$ une branche de~$X$ issue de~$x$. Si~$Z$ est  une section coronaire de~$b$ alors pour tout~$I\in \inter X x \ctd {\skel Z}$ l'ouvert~$I^\flat$ est une sous-couronne virtuelle de~$b$, et le plongement~$\got s(Z)\hookrightarrow \got s(I^\flat)$ est donc un isomorphisme. On déduit alors de~\ref{apropossx} que~$\got s(Z)$ s'identifie naturellement à~$\got s(b)$. 

\trois{remsxsv} Soit~$V$ un ouvert de~$X$ qui est un disque virtuel de bord~$\{x\}$ et soit~$b$ l'unique branche issue de~$x$ et contenue dans~$V$. On sait qu'il existe un ouvert~$Z$ de~$V$ qui est une couronne virtuelle aboutissant à~$x$ et telle que~$\got s(Z)=\got s(V)$ (\ref{discvirtcontcour} {\em et sq.}) ; la couronne virtuelle~$Z$ est section coronaire de~$b$, et l'on a par conséquent,~$\got s(b)=\got s(Z)$, d'où finalement l'égalité~$\got s(b)=\got s(V)$ ; ceci entraîne notamment que~$\got s(x)$ se plonge naturellement dans~$\got s(V)$. 

\medskip
Si~$y$ est un point de~$V\geom$ alors l'application de~$[x;y]$ dans~$\NN$ qui envoie~$z$ sur~$[\got s(z):k]$ est croissante lorsqu'on se dirige vers~$y$ : cela résulte de~\ref{anteccroissants}.  

\trois{remsxgen} Soit~$U$ un voisinage ouvert connexe de~$x$ qui est un arbre et qui est tel que~$U\setminus\{x\}$ soit réunion disjointe de disques virtuels et de couronnes virtuelles, et soit~$f\in \got s(x)$. 

On peut voir~$f$ comme une fonction définie sur un voisinage~$U'$ de~$x$, que l'on peut restreindre de sorte que la proposition suivante soit satisfaite : il existe une famille finie~$(U_i)$ de composantes connexes deux à deux disjointes de~$U\setminus\{x\}$, et pour tout~$i$ un fermé~$Z_i$ de~$U_i$, tels que~$U_i-Z_i$ soit pour tout~$i$ une section coronaire de l'unique branche issue de~$x$ et contenue dans~$U_i$, et tels que~$U'$ soit égal à~$U-\coprod Z_i$. Il résulte de~\ref{remseccorong} et de~\ref{remsxsv} que~$\got s(U_i)\to \got s(U_i-Z_i)$ est un isomorphisme pour tout~$i$ ; on en déduit que~$f$ admet un unique prolongement à~$U$. Ainsi, l'injection naturelle~$\got s(U)\to \got s(x)$ est un isomorphisme.

\subsection*{Modération des morphismes compacts et non constants entre courbes}

\deux{defdilat} Soient~$\Gamma$ et~$\Delta$ deux graphes localement finis, respectivement munis de deux toises~$\ell_\Gamma$ et~$\ell_\Delta$. On dira qu'une application continue~$\phi : \Delta\to \Gamma$ est une {\em dilatation} si~$\phi$ est un homéomorphisme et s'il existe un réel~$r\in \RR\ti_+$ tel que~$\ell_\Gamma(\phi(I))$ =$r\ell_\Delta(I)$ pour tout segment~$I$ tracé sur~$\Delta$ ; si~$\Delta$ n'est pas discret, le réel~$r$ en question est uniquement déterminé et est appelé le {\em rapport} de~$\phi$. Lorsqu'il vaut~$1$, on dit que~$\phi$ est une {\em isométrie}. 

\deux{defdilmorc} Soient~$\Gamma$ et~$\Delta$ deux graphes localement finis et soit~$\phi$ une application continue~$\Delta \to \Gamma$. Supposons qu'il existe un ensemble fermé et discret~$S$ de~$\Gamma$ tel que~$\phi\inv(S)$ soit discret dans~$\Delta$ ; et supposons que~$\Delta-\phi\inv(S)$ et~$\Gamma\setminus S$ soient chacun munis d'une toise. On dira que~$\phi$ est une {\em dilatation (resp. isométrie) par morceaux} s'il existe un sous-ensemble fermé et discret~$T$ de~$\Delta$ contenant~$\phi\inv(S)$ tel que pour toute composante connexe~$J$ de~$\Delta\setminus T$ l'application~$\phi$ induise une dilatation (resp. une isométrie)~$J\simeq \phi(J)$. 

\medskip
Toute dilatation par morceaux est injective par morceaux. 

\deux{toiseimpl} Soit~$\phi : Y\to X$ un morphisme de dimension relative nulle entre deux courbes~$k$-analytiques. Soit~$\Gamma$ un sous-graphe localement fini de~$X$ et soit~$\Delta$ un sous-graphe localement fini de~$Y$ tel que~$\phi(\Delta)\subset \Gamma$. Le sous-ensemble~$\Gamma\typ{14}$ de~$\Gamma$ en est une partie fermée et discrète, et~$\phi_{|\Delta}\inv(\Gamma\typ{14})=\Delta\typ {1 4}$. Lorsqu'on parlera de dilatation par morceaux dans ce contexte, ce sera toujours au sens du~\ref{defdilmorc} 
ci-dessus, avec les conventions implicites suivantes :~$S=\Gamma_{14}$  ; les graphes~$\Gamma\setminus S=\Gamma\dtr$ et~$\Delta- \phi_{|\Delta}\inv(S)=\Delta\dtr$ sont munis des restrictions respectives des toises canoniques de~$X\dtr$ et~$Y\dtr$ (\ref{metrinat}). 

\deux{exdilat} {\em Exemples.} 

\trois{gammaf-isom}
Soit~$X$ une courbe~$k$-analytique, soit~$F$ une extension presque algébrique de~$k$ et soit~$\Gamma$ un sous-graphe localement fini de~$X$. L'image réciproque~$\Gamma_F$ de
$\Gamma$ sur~$X_F$ est un sous-graphe localement fini de~$X_F$ (cor.~\ref{corollimrecih}), et il découle de la proposition~\ref{prop-toise-presqualg}
que~$\Gamma_F\to \Gamma$ est une isométrie par morceaux.

\trois{lemdilat}
Soient~$Y$ et~$X$ deux courbes~$k$-analytiques, chacune d'elles étant une couronne virtuelle. Si~$\phi : Y\to X$ est un morphisme fini et plat, il résulte du lemme~\ref{morcourvirt} que~$\phi$ induit une dilatation de rapport~$\displaystyle{\frac{\deg \;\phi}{[\got s(Y):\got s(X)]}}$ entre~$\skelan Y$ et~$\skelan X$. 

\medskip
Notons un cas particulier important :  si~$[\got s(Y):\got s(X)]=\deg\;\phi$, c'est-à-dire encore si~$Y$ s'identifie à~$X\times_{\got s(X)}\got s(Y)$ (ce qui est notamment le cas lorsque~$Y$ est une composante connexe de~$X_L$ pour une certaine extension finie~$L$ de~$k$), alors~$\skelan Y\to \skelan X$ est une dilatation de rapport 1, c'est-à-dire une isométrie ; on peut également voir ce résultat comme une déclinaison de l'exemple~\ref{gammaf-isom}
ci-dessus. Cela dit, le plus simple est de le vérifier directement : c'est une simple reformulation de l'invariance du module d'une couronne virtuelle par extension des scalaires
(cette invariance est d'ailleurs le cœur de la preuve de la proposition~\ref{prop-toise-presqualg}
sur laquelle se fonde l'exemple~\ref{gammaf-isom}). 

\deux{modnonconst} {\bf Théorème.} {\em Soit~$\phi : Y\to X$ un morphisme compact et de dimension relative nulle entre courbes~$k$-analytiques et soit~$\Gamma$ un sous-graphe fermé et localement fini de~$Y$. 

\medskip
\begin{itemize}
\item[1)] L'image de~$\phi(\Gamma)$ est un sous-graphe fermé et localement fini de~$X$. 

\item[2)] La flèche~$\Gamma\to X$ est injective par morceaux.

\item[3)] La flèche~$\Gamma\geom\to \phi(\Gamma\geom)$ est une dilatation par morceaux. 

\item[4)] On suppose de surcroît que~$Y$ et~$X$ sont génériquement quasi-lisses ; on appelle~$T$ le sous-ensemble fermé et discret de~$\Gamma$ égal à la réunion de~$\Gamma\cap \partial\an \phi$ et de l'ensemble des points rigides de~$\Gamma$ en lesquels~$\phi$ n'est pas plat. 

\medskip
\begin{itemize}

\item[$\alpha)$] si~$y\in \Gamma\typ 3\setminus T$, il existe un intervalle~$I$ qui est un voisinage ouvert de~$y$ dans~$\Gamma\dtr \setminus T$ tel que~$\phi$ induise une dilatation~$I\simeq \phi(I)$ de rapport~$\displaystyle{\frac {\deg^y\;\phi}{[\got s(y):\got s(\phi(y))]}}$, et tel que~$\deg^z\phi=\deg^y\;\phi$ pour tout~$z\in I$ ; 

\medskip
\item[$\beta)$] si~$y\in \Gamma\geom$ et si~$b\in \br \Gamma y$, il existe~$I\in \inter {\Gamma} b\ctd {\Gamma\dtr \setminus T}$, tel que~$\phi$ induise une dilatation~$I\simeq \phi(I)$ de rapport~$\displaystyle{\frac {\deg\;(b\to \phi(b))}{[\got s(b):\got s(\phi(b))]}}$, et tel que~$\deg^z\phi=\deg\;(b\to \phi(b))$ pour tout~$z\in I$ ; 

\medskip
\item[$\gamma)$] si~$y\in \Gamma$ et si~$b\in \br \Gamma y$, il existe~$I\in \inter {\Gamma} b\ctd{\Gamma\dtr \setminus T}$ tel que~$\deg^z\;\phi$ soit égal à~$\deg \;(b\to \phi(b))$ pour tout~$z\in I$ ; si de plus le point~$y$ appartient à~$\Gamma\typ {1,3,4}\setminus T$ et si~$Y$ est normale en~$y$ alors ce degré coïncide avec~$\deg^y\;\phi$.

\end{itemize}
\end{itemize}
}

\medskip
{\em Démonstration.} Les assertions à montrer sont locales sur le but, ce qui autorise à supposer que~$X$ est compacte ; la courbe~$Y$ l'est alors automatiquement, puisque~$\phi
$ est compact. On peut ensuite : étendre les scalaires au complété de la clôture parfaite de~$k$, c'est-à-dire supposer que~$k$ est parfait ; puis remplacer~$X$ et~$Y$ par les courbes réduites sous-jacentes, c'est-à-dire se ramener au cas où~$X$ et~$Y$ sont génériquement quasi-lisses. 

\medskip
Le graphe~$\Gamma$ est fini et compact ; l'assertion 1) sera dès lors une conséquence de 2). L'assertion 2) elle-même est locale sur~$\Gamma$, et découlera de 3) au voisinage des points de~$\Gamma \geom$ ; l'assertion 3) résultera quant à elle de 4),~$\beta)$. Il reste donc à prouver 2) au voisinage des points de type 1 ou 4 non rigides, et 4).

\trois{phiinjmorc14} {\em Preuve de 2) et de 4),~$\gamma)$ au voisinage des points de type 1 ou 4 non rigides.} Soit~$y\in \Gamma\typ{1,4}-\Gamma\typ 0$ et soit~$x$ son image sur~$X$ ; le point~$x$ appartient à~$X\typ{1,4}-X\typ 0$. Comme~$y$ n'appartient pas au bord analytique de~$Y$, comme il n'est pas rigide et comme~$X$ est réduit,~$\phi$ est fini et plat en~$y$. On déduit de~\ref{bondeviss} l'existence d'un voisinage affinoïde~$W$ de~$y$ dans~$Y$ et d'un voisinage affinoïde~$V$ de~$x$ dans~$X$ tel que~$\phi(W)\subset V$ et tel que~$W\to V$ admette une factorisation~$W\to Z\to V$ où~$W\to Z$ est fini, plat et radiciel et où~$Z\to V$ est fini étale ; on peut supposer que~$y$ est le seul antécédent de~$x$ sur~$W$. Soit~$\Delta$ l'image de~$\Gamma\cap W$ dans~$Z$ ; notons que~$(\Gamma\cap W)\to \Delta$ est un homéomorphisme ; par conséquence, le compact~$\Delta$ est un graphe fini ; soit~$\Delta_0$ l'image de~$\Delta$ dans~$V$, qui coïncide avec celle de~$\Gamma\cap W$. 

\medskip
Comme~$Z\to V$ est fini étale et comme~$Z$ et~$V$ sont connexes et non vides, il existe un morphisme fini étale~$Z'\to Z$ tel que~$Z'\to V$ soit un revêtement fini galoisien connexe dont on note~$\mathsf G$ le groupe. L'espace~$Z$ s'identifie à un quotient~$Z'/\mathsf H$ pour un certain sous-groupe~$\mathsf H$ de~$\mathsf G$. Il découle du corollaire~\ref{corollimrecih} que l'image réciproque~$\Delta_1$ de~$\Delta_0$ sur~$Z$ est un sous-graphe compact et fini, et que~$\Delta_1\to \Delta$ est injective par morceaux. L'image réciproque~$\Gamma_1$ de~$\Delta_1$ sur~$ W$ s'identifie à l'image réciproque de~$\Delta_0$ sur~$W$ ; comme elle est canoniquement homéomorphe à~$\Delta_1$, c'est un graphe compact et fini, et~$\Gamma_1\to \Delta_0$ est injective par morceaux. 

\medskip
Il s'ensuit que~$(\Gamma\cap Z)\to \Delta_0$ est injective par morceaux, ce qui prouve 2). Par ailleurs, comme~$y$ et~$x$ sont de type 1 ou 4 et non rigides, les ensembles~$\br X x$ et~$\br Y y$ sont des singletons ; notons~$b$ l'unique élément de~$\br Y y$. 

\medskip
Supposons que la branche~$b$ appartienne à~$\br \Gamma y$ (cela revient à demander que~$y$ ne soit pas isolé dans~$\Gamma$) ; elle appartient {\em a fortiori} à~$\br {\Gamma_1} y$, et l'on a dans ce cas~$\br {\Gamma_1}y=\br \Gamma y =\{b\}$ ; l'ensemble ~$\br  {\Delta_0}x$ est dès lors non vide, et coïncide nécessairement avec le singleton,~$\br X x$. 

On déduit de ce qui précède et du fait que~$y$ est le seul antécédent de~$x$ sur~$W$ qu'il existe un intervalle ouvert~$I\in \inter Y b \ctd {\Gamma_1}$ tel que pour tout 
$z\in I$, l'on ait~$\phi_{|W}\inv(\phi(z))=\{z\}$. On peut supposer, quitte à le restreindre, que~$I$ est contenu dans~$\Gamma\dtr$ et dans l'intérieur topologique de~$W$ dans~$Y$ ; cette dernière condition assure, la flèche~$W\to V$ étant finie et plate, que~$I\subset \Gamma\setminus T$.

Comme~$y$ est le seul antécédent de~$x$ sur~$W$ et comme~$b$ est la seule branche de~$Y$ issue de~$y$, on a~$\deg ^y\phi=\deg \;(b\to \phi(b))=\deg \;(W\to V)$ ; et si~$z\in I$ l'égalité~$\phi_{|W}\inv(\phi(z))=\{z\}$ assure que~$\deg ^z\phi=\deg \;(W\to V)=\deg \;(b\to \phi(b))$, ce qui achève de prouver 4),~$\gamma)$. 

\trois{phindilatmorc} {\em Preuve de 4).}  On démontre chacune des trois assertions séparément. 

\medskip
{\em Démonstration de~$\alpha)$.} Comme~$y\notin \partial\an \phi$, le morphisme~$\phi$ est fini et plat en~$y$. Il découle alors de l'énoncé ii) de la proposition~\ref{degetseccor} et de~\ref{commentdeg} qu'il existe un voisinage coronaire~$Z'$ de~$y$ dans~$Y$ et un voisinage coronaire~$Z$ de~$\phi(y)$ dans~$X$ tel que~$\phi$ induise un morphisme fini et plat~$Z'\to Z$ de degré~$\deg^y \phi$. Le point~$y$ étant de type 3 et situé sur~$\skelan {Z'}$, on a~$\br Y y =\br {Z'}y =\br {\skelan {Z'}} y$ et l'on peut donc restreindre~$Z'$ (et~$Z$) de sorte que~$I:=\Gamma\cap \skelan {Z'}$ soit un intervalle tracé sur~$\Gamma$ et ouvert dans ce dernier ; l'inclusion de~$I$ dans~$Z'$ garantit par ailleurs que~$I\subset \Gamma\dtr$, et que~$I$
évite~$T$, puisque~$Z'\to Z$ est fini et plat. L'assertion~$\alpha)$ découle alors de~\ref{lemdilat}, du fait que~$\got s(y)=\got s(Z')$ et~$\got s(\phi(y))=\got s(\phi(Z))$, et du fait que tout point~$z$ de~$I\subset \skelan Z'$ est l'unique antécédent de~$\phi(z)$ sur~$Z'$. 

\medskip
{\em Démonstration de~$\beta)$.} En vertu de l'énoncé i) de la proposition~\ref{degetseccor}, il existe une section coronaire~$Z'$ de~$b$ dans~$Y$ et une section coronaire~$Z$ de~$\phi(y)$ dans~$X$ telles
que~$\phi$ induise un morphisme fini et plat~$Z'\to Z$ de degré~$\deg\;(b\to \phi(b))$ ; on peut restreindre~$Z'$ (et~$Z$) de sorte que~$I:=\skelan {Z'}$ soit un contenu dans~$\Gamma$ ; l'inclusion de~$I$ dans~$Z'$ garantit par ailleurs que~$I\subset \Gamma\dtr$, et que~$I$ évite~$T$, puisque~$Z'\to Z$ est fini et plat. L'assertion~$\beta))$ découle alors de~\ref{lemdilat}, du fait que~$\got s(b)=\got s(Z')$ et~$\got s(\phi(b))=\got s(\phi(Z))$ (\ref{remseccorong}), et du fait que tout point~$z$ de~$I\subset \skelan Z'$ est l'unique antécédent de~$\phi(z)$ sur~$Z'$. 

\medskip
{\em Démonstration de~$\gamma)$.} Cette assertion a été établie au~\ref{phiinjmorc14} lorsque le point~$y$ appartient à~$\Gamma\typ{1,4}-\Gamma\typ 0$. Supposons maintenant que~$y\in \Gamma\geom$. L'existence d'un intervalle~$I$ ayant les propriétés requises découle alors immédiatement de l'assertion~$\beta)$ déjà prouvée ; si~$y\in \Gamma\typ 3\setminus T$ alors~$\phi \inv (\phi(b))=\{b\}$, d'où il découle que~$\deg \; (b\to \phi(b))=\deg^y\; \phi$ ; enfin, si~$y\in \Gamma\typ {1 4}\setminus T$, et si de plus~$Y$ est normale en~$y$ alors~$\br Y y =\{b\}$, et il s'ensuit là encore que~$\deg \;(b\to \phi(b))=\deg^y \phi$.~$\Box$ 

\deux{theomodimrecg} {\bf Théorème.} {\em Soit~$\phi : Y\to X$ un morphisme de dimension relative nulle entre courbes~$k$-analytiques et soit~$\Gamma$ un sous-graphe fermé et localement fini de~$X$. 

\medskip
1) L'image réciproque~$\phi^{-1}(\Gamma)$ est un sous-graphe fermé et localement fini de~$Y$. 

2) Si~$y\in \phi^{-1}(\Gamma)$, si~$x$ désigne son image sur~$X$ et si~$b\in \br Y y$ alors~$b\in \br {\phi^{-1}(\Gamma)}y$ si et seulement si~$\phi(b)\in \br \Gamma x$.} 

\medskip
{\em Démonstration.} La question étant purement topologique, on peut étendre les scalaires au complété de la clôture parfaite de~$k$, et donc supposer~$k$ parfait ; puis l'on peut remplacer~$Y$ et~$X$ par les courbes réduites associées, et donc supposer que~$Y$ et~$X$ sont réduites, et en partant génériquement quasi-lisses (le corps de base étant parfait). 

\medskip
ll est clair que~$\phi^{-1}(\Gamma)$ est fermé. Soit~$y\in \phi^{-1}(\Gamma)$, et soit~$x$ son image sur~$X$. Nous allons montrer que~$y$ possède un voisinage dans~$\phi^{-1}(\Gamma)$ qui est un arbre fini. Pour ce faire, on peut, quitte à restreindre~$Y$ et~$X$, supposer que ce sont des arbres compacts et que~$\phi^{-1}(x)=\{y\}$. 

Soit~$b\in \br \Gamma x$ et soient~$\beta_1,\ldots, \beta_r$ les antécédents de~$b$ dans~$\br Y y$. Il existe une section~$U$ de~$b$ telle que~$\phi^{-1}(U)$ soit de la forme~$\coprod\limits_{1\leq i\leq r} V_i$ où~$i$ est pour tout~$i$ une section de~$\beta_i$, finie et plate sur~$U$ de degré~$\deg (\beta_i\to b)$.  Fixons~$i$ et un point~$y_i$ de~$V_i$. En vertu de l'assertion 4),~$\gamma)$ du théorème~\ref{modnonconst} ci-dessus, il existe~$z_i\in ]x;y_i[$ tel que~$\deg^z\phi=\deg (\beta_i\to b)$ pour tout~$z\in ]y;z_i]$. Cette égalité, jointe au fait que~$\deg(V_i\to U)=\deg (\beta_i\to b)$, implique que si~$z\in  ]y;z_i]$ alors~$\phi_{|V_i}^{-1}(\phi(z))=\{z\}$ ; comme on a par ailleurs~$\phi^{-1}(x)=\{y\}$, il en résulte que~$\phi$ induit un homéomorphisme~$[y;z_i]\simeq [x;\phi(z_i)]$ et que~$$\phi_{|V_i}^{-1}(]x;\phi(z_i)]=]y;z_i].$$

Les intervalles~$]x;\phi(z_i)]$ sont tous tracés sur~$U$, et aboutissent tous à~$x$ ; de plus, la branche~$b$ appartient par hypothèse à~$\br \Gamma x$. Par conséquent, il existe~$t_b\in U$ tel que~$]x;t_b[$ soit contenu dans~$\Gamma\cap \bigcap ]x;z_i]$. Par construction, l'image réciproque de~$[x;t_b]$ est une réunion finie de segments dont les intersections deux à deux sont toutes égales à~$\{y\}$ ; c'est en particulier un arbre compact fini. 

\medskip
La réunion~$\Delta$ des segments~$[x;t_b]$ pour~$b$ parcourant~$\br \Gamma x$ est un voisinage de~$x$ dans~$\Gamma$ ; par conséquent,~$\phi^{-1}(\Delta)$ est un voisinage de~$y$ dans~$\phi^{-1}(\Gamma)$, et il résulte de ce qui précède que~$\phi^{-1}(\Delta)$ est un arbre compact et fini, ce qui achève la démonstration de 1).

\medskip
Montrons maintenant 2). Si~$b\in \br {\phi^{-1}(\Gamma)}y$il est clair que~$\phi(b)\in \br \Gamma x$ puisque~$\phi^{-1}(\Gamma)\to X$ se factorise par~$\Gamma$. Réciproquement, supposons que~$\phi(b)\in \br \Gamma x$. On déduit de l'assertion 2) du théorème~\ref{modnonconst} l'existence d'un intervalle~$I\in \intera Y b$ tel que~$\phi$ induise un homéomorphisme entre~$\overline I$ et~$\phi(\overline I)$ ; par conséquent,~$\phi(I)$ est un intervalle aboutissant proprement à à~$x$, et la branche qu'il définit est nécessairement~$\phi(b)$. Par hypothèse,~$\phi(b)\in \br \Gamma x$ ; il existe donc un intervalle ouvert~$J$ contenu dans~$\phi(I)$, aboutissant à~$x$ et inclus dans~$\Gamma$ ; son image réciproque~$\phi_{|I}^{-1}(J)$ est un intervalle qui est contenu dans~$I$ et aboutit à~$y$, et définit par conséquent la branche~$b$. Comme~$\phi{|I}^{-1}(J)\subset \phi^{-1}(\Gamma)$, on a bien~$b\in \br {\phi^{-1}(\Gamma)}y$.~$\Box$ 

\subsection*{Variation des fonctions holomorphes et loi des nœuds} 

\deux{defpentesholo} {\bf Graphe de variation.} Soit~$X$ une courbe~$k$-analytique et soit~$f$ une fonction analytique inversible sur~$X$ ; elle induit un morphisme~$\phi : X\to \gmk\an$ ; soit~$\Gamma$ l'image réciproque de~$\skel {\gmk\an}$ sur~$X$. Notons que comme~$\skel {\gmk\an}$ n'est constitué que de points de type 2 ou 3, le fermé~$\Gamma$ est contenu dans~$X\typ {2 3}$.

\trois{polvarfini} {\em Le fermé~$\Gamma$ de~$X$ en est un sous-graphe localement fini, et la flèche~$\Gamma\to \skel{\gmk \an}$ est une dilatation par morceaux.} En effet, soit~$x\in \Gamma$ ; comme~$x$ est de type 2 ou 3,  il est situé sur une et une seule composante irréductible~$X_0$ de~$X$. Comme l'ensemble~$f(X_0)$ contient un point de~$\skel {\gmk\an}$, il n'est pas réduit à un point rigide, ce qui signifie que la restriction de~$f$ à~$X_0$ est non constante, c'est-à-dire encore de dimension relative nulle. On déduit du théorème~\ref{theomodimrecg} que~$\Gamma\cap X_0$ est un graphe localement fini, et de l'assertion 3) du théorème~\ref{modnonconst} que~$\Gamma\cap X_0\to \skel {\gmk \an}$ est une dilatation par morceaux (pour appliquer en toute rigueur le théorème ~\ref{modnonconst}, on commence par se ramener en raisonnant localement sur~$\Gamma$ au cas d'une source compacte). Comme~$x$ n'appartient qu'à la composante  irréductible~$X_0$ de~$X$, celle-ci est un voisinage de~$x$ dans~$X$, et~$\Gamma\cap X_0$ est dès lors un voisinage de~$x$ dans~$\Gamma$, d'où notre assertion. 

\trois{polvarpolvar} La fonction~$|T|$ étant localement constante sur~$\gmk \an-\skel {\gmk \an}$, la fonction~$|f|$ est localement constante sur~$X-\Gamma$. Comme~$\Gamma\to \skel {\gmk \an}$ est une dilatation par morceaux, et comme~$|T|$ est strictement croissante sur~$\skel {\gmk \an}$ (orienté de~$0$ vers~$\infty$), la fonction~$|f|$ n'est constante sur aucun intervalle non vide et non singleton de~$\Gamma$. 

Soit~$x$ un point isolé de~$\Gamma$ et soit~$U$ la composante connexe de~$x$ dans l'ouvert~$X-(\Gamma\setminus\{x\})$. Si~$V$ est une composante connexe de~$U\setminus\{x\}$ alors~$|f|$ est localement constante, et donc constante, sur~$V$ ; et comme~$x\in \partial V$, la valeur constante de~$|f|~$sur~$V$ est égale à~$|f(x)|$. Ainsi,~$|f|$ est constante sur~$U$. Il s'ensuit que si~$\Gamma_0$ désigne l'ensemble des points non isolés de~$\Gamma$, l'ouvert~$X-\Gamma_0$ est le plus grand ouvert de~$X$ sur lequel~$|f|$ est localement constante ; on dira que le graphe fermé et localement fini~$\Gamma_0$ est le {\em graphe} de variation  de~$f$. 

\trois{rempasbord} Remarquons que~$\Gamma-\Gamma_0\subset \partial \an X$ (en particulier,~$\Gamma=\Gamma_0$ dès que~$X$ est sans bord). En effet, soit~$x\in \Gamma\setminus \partial \an X$, et soit~$\xi$ son image par~$\phi$. Comme~$x$ n'appartient pas au bord de~$X$ et comme~$\phi$ est de dimension nulle en~$x$, elle est finie en~$x$. Comme~$\xi$ est de type 2 ou 3 et comme~$\gmk \an$ est réduit,~$\sch O_{\gmk \an, \xi}$ est un corps, et~$\phi$ est donc fini et {\em plat} en~$x$. Par conséquent, si~$U$ est un voisinage ouvert de~$x$ dans~$X$ alors~$\phi(U)$ est un voisinage ouvert de~$\xi$ dans~$\gmk \an$, et contient en particulier un intervalle ouvert non vide de~$\skel {\gmk \an}$ ; il s'ensuit que~$U\setminus\{x\}$ rencontre~$\Gamma$, ce qui signifie que~$x$ n'est pas un point isolé de~$\Gamma$. 

\trois{defpentef} Soit~$x\in \Gamma$, soit~$\xi$ son image sur~$\gmk \an$ et soit~$b\in \br X x$. Soit~$b^-$ (resp.~$b^+$) la branche de~$\gmk \an$ issue de~$\xi$ définie par~$]0;\xi[$ (resp.~$]\xi;\infty[$). Comme~$\phi$ est de dimension nulle en~$x$, l'image~$\phi(b)$ est bien définie. Il résulte de l'assertion 2) de théorème~\ref{theomodimrecg} que~$b\in \br \Gamma x$ si et seulement si~$\phi(b)=b^+$ ou~$\phi(b)=b^-$ ; les branches~$b^+$ et~$b^-$ étant exactement les deux branches le long desquelles~$|T|$ n'est pas constante, on peut reformuler cette dernière assertion en disant que~$b\in \br \Gamma x$ si et seulement si~$|f|$ n'est pas constante le long de~$b$.

\medskip
Supposons que~$b\in \br \Gamma x$, et posons~$\epsilon=1$ (resp.~$-1$) si~$\phi(b)=b^+$ (resp.~$b^-$). 

Comme la flèche~$\Gamma\to \skel {\gmk \an}$ est une dilatation par morceaux, il existe~$I\in \inter \Gamma b$ et un réel~$r>0$ tel que~$\phi$ induise une dilatation~$I\simeq \phi(I)$ de rapport~$r$ ; le réel~$r$ ne dépend que de~$b$, et pas du choix de~$I$. Le réel~$\epsilon r[\got s(b):k]$, auquel on peut penser comme à la pente normalisée  de~$\log |f|$ le long de~$b$, sera noté~${\rm dlog}_b\; |f|$.

\medskip
{\em Le réel ~${\rm dlog}_b\; |f|$ est un entier, qui coïncide avec~$\epsilon \deg\;(b\to \phi(b))$ lorsque~$X$ est génériquement quasi-lisse.} En effet, le réel ~${\rm dlog}_b\; |f|$ n'est pas modifié si l'on remplace~$k$ par le complété de sa clôture parfaite, ni si l'on remplace~$X$ par~$X_{\rm red}$ ; on peut ainsi supposer que~$X$ est génériquement quasi-lisse, et il n'y a plus alors qu'à appliquer l'assertion 2) du théorème ~\ref{modnonconst}. 

\deux{theoloinoeuds} {\bf Théorème (harmonicité de~$\rm log |f|$ pour~$f$ holomorphe)}. {\em Soit~$X$ une courbe~$k$-analytique sans bord et soit~$f$ une fonction analytique inversible sur~$X$. Soit~$\Gamma$ son graphe de variation et soit~$x\in \Gamma$. On a alors la « loi des nœuds»~suivante :~$$ \sum _{b\in \br \Gamma x} {\rm dlog}_b |f|=0.$$}

\medskip
{\em Démonstration.} Pour montrer cette égalité, on peut étendre les scalaires au complété de la clôture parfaite de~$k$, puis remplacer~$X$ par~$X_{\rm red}$ ; autrement dit, on peut se ramener au cas où~$X$ est génériquement quasi-lisse. Soit~$\xi$ l'image de~$x$ par le morphisme~$\phi : X\to \gmk \an$ que définit~$f$. Comme~$X$ est sans bord, comme~$\phi$ est de dimension relative nulle en~$x$ et comme~$\sch O_{\gmk\an,\xi}$ est un corps (le point~$\xi$ étant de type 2 ou 3), le morphisme de germes~$(X,x)\to (\gmk\an,\xi)$ est fini et plat ; soit~$d$ son degré. Soit~$b^-$ (resp.~$b^+$) la branche de~$\gmk \an$ issue de~$\xi$ et définie par~$]0;\xi[$ (resp.~$]\xi;\infty[$) et soit~$\bnd B^-$ (resp.~$\bnd B^+$) l'ensemble des antécédents de~$b^-$ (resp.~$b^+$) dans~$\br X x$ ; d'après~\ref{defpentef}, on a~$\br \Gamma x=\bnd B^+\coprod \bnd B^-$. Il vient~$$ \sum _{b\in \br \Gamma x} {\rm dlog}_b |f|=\sum _{b\in \bnd B^+} {\rm dlog}_b |f|+\sum _{b\in \bnd B^-} {\rm dlog}_b |f|$$~$$=\sum _{b\in \bnd B^+}\deg(b\to \phi(b))-\sum _{b\in \bnd B^-}\deg(b\to \phi(b))=d-d=0.\;\;\Box$$

\deux{remharm} {\em Remarque.} La loi des nœuds énoncée ci-dessus vaut encore si l'on remplace~$\Gamma$ par n'importe quel sous-graphe fermé localement fini de~$X$ {\em contenant} le graphe de variation de~$f$, à condition de poser~${\rm dlog}_b|f|=0$ pour toute branche~$b$ le long de laquelle~$|f|$ est constante.

\chapter{Triangulations des courbes analytiques quasi-lisses}
\markboth{Triangulations}{Triangulations}

\section{Sous-graphes analytiquement admissibles, squelette analytique et triangulations}

\subsection*{Lemmes de fusion}

\deux{fusion} {\bf Lemme.} {\em Soit~$X$ une courbe~$k$-analytique, et soient~$U$ et~$V$ deux ouverts connexes et non vides de~$X$ possédant les propriétés suivantes :

\medskip
1)~$U$ est une couronne virtuelle ; 

2)~$V$ est un disque virtuel ; 

3)~$V$ rencontre le squelette de~$U$ et il existe un élément~$u$ de~$\partial U$ qui n'appartient pas à~$V$.

\medskip
Sous ces hypothèses :

\medskip
$\bullet$~$\got s(V)=\got s(U\cup V)=\got s(U)$ ; notons~$L$ ce corps ; 

$\bullet$~$U\cup V$ est un disque virtuel sur~$L$ dont l'adhérence est compacte, égale à~$U\cup V\cup\{u\}$, et est une composante connexe de~$X\setminus\{u\}$ ; 

$\bullet$ si~$F$ est une extension presque algébrique de~$L$ qui déploie~$U$ et si l'espace~$L$-analytique~$V$ a un~$F$-point, alors~$F$ déploie~$U\cup V$.}

\medskip
{\em Démonstration.} On procède en plusieurs étapes. 

\trois{propdiscv} Comme~$V$ est un disque virtuel, il existe un intervalle ouvert~$I$ tracé sur~$V$ et aboutissant à son unique bout~$\omega$ tel que~$\got s(x)=\got s(V)$ pour tout~$x\in I$ (\ref{existifixpp} et~\ref{uniqueantec})

\trois{ucupvarb} {\em Topologie de~$V\cup U$.} Comme le graphe~$V\cup U$ est réunion de deux ouverts connexes qui s'intersectent, il est connexe. Tout point de~$(V\cup U)-\skel U$ a dans~$V\cup U$ un voisinage qui est un arbre à un bout.
Le squelette de~$V\cup U$ est donc inclus dans~$\skel U$, et
il est connexe puisque~$V\cup U$ est connexe ; comme un squelette n'a jamais de point isolé ni unibranche,~$\mathsf S(V\cup U)$ ne peut-être qu'un intervalle ouvert de~$\skel U$, nécessairement strict puisque~$\skel U$ rencontre~$V$. Par ailleurs,~$\mathsf S(V\cup U)$ est fermé dans~$V\cup U$, et {\em a fortiori}
dans~$\skel U$ ; il en résulte finalement que~$\mathsf S(U\cup V)=\emptyset$, et donc que~$V\cup U$ est un arbre ayant au plus un bout. Le point~$u$ appartenant en vertu de 2) 
au bord de~$V\cup U$ dans le graphe~$X$, l'arbre~$V\cup U$ a exactement un bout ; son adhérence est égale à~$V\cup U\cup\{u\}$ et s'identifie à sa compactification arboricole. 

\medskip
Le sous-ensemble~$V\cup U$ de~$X$ est une partie connexe, non vide, ouverte et fermée de~$X\setminus\{u\}$ ; c'en est donc une composante connexe. 

\trois{sucupv} Comme~$V\cup U\cup\{u\}$ est un arbre compact, le bord de~$\skel U$ dans celui-ci est de la forme~$\{u,\xi\}$ pour un certain~$\xi\neq u$. Par convexité, l'ouvert non vide~$V\cap \skel U$ est un intervalle ouvert de~$\skel U$ ; soit~$\eta$ sa borne supérieure, lorsqu'on l'oriente en direction de~$u$. Le point~$\eta$ adhère à~$V$ mais n'appartient pas à~$V$ : cela résulte de sa définition et du fait que~$u\notin V$ par hypothèse. Comme~$\partial V$ ne peut contenir plus d'un élément ($V$ étant un arbre à un bout), la borne inférieure de~$V\cap \skel U$ appartient à~$V$, et par conséquent n'appartient pas à~$\skel U$ ; cela ne peut donc être que~$\xi$. 

\medskip
L'arbre à un bout~$V$ apparaît dès lors comme la composante connexe de~$V\cup U\cup\{u\}\setminus\{\eta\}$ contenant~$\xi$. Son intersection avec~$U$ est la sous-couronne virtuelle~$]\xi;\eta[^\flat$ de~$U$.

\medskip
L'intervalle~$I$ du~\ref{propdiscv} peut en conséquence être restreint de façon à être contenu dans~$]\xi;\eta[$ ; fixons~$x\in I$. Comme~$x$ est situé sur le squelette de~$U$, le plongement~$\got s(U)\hookrightarrow \got s(x)$ est un isomorphisme ; mais~$\got s(V)\hookrightarrow \got s(x)$ est aussi un isomorphisme par définition de~$I$. 

\medskip
Comme~$\got s(U)\hookrightarrow \got s(x)~$ et~$\got s(V)\hookrightarrow \got s(x)$ se factorisent par~$\got s(U\cap V)$ (ce dernier est bien défini puisque~$U\cap V$ est connexe et non vide), les injections canoniques~$\got s(U)\hookrightarrow \got s(U\cap V)$ et ~$\got s(V)\hookrightarrow \got s(U\cap V)$ sont des isomorphismes ; il s'ensuit que toute fonction appartenant à~$\got s(U)$ (resp.~$\got s(V)$ s'étend de manière naturelle en une fonction de~$\got s(U\cup V)$. 

\medskip
Par conséquent,~$\got s(U)=\got s(V)=\got s(U\cup V)$ ; notons~$L$ ce corps. 

\medskip
Pour toute extension presque algébrique~$F$ de~$L$, on désignera par~$x_F$ l'unique antécédent de~$x$ sur~$V\hotimes_LF$ ; pour tout intervalle~$J$ tracé sur~$U$, on désignera par~$J_F$ l'image réciproque de~$J$ sur~$U\hotimes_LF$ ; notons que~$\skel U_F=\mathsf S(U\hotimes_LF)$. 

\medskip
Soit~$F$ une extension presque algébrique de~$L$ qui déploie~$U$, et qui est telle que l'espace~$L$-analytique~$V$ ait un~$F$-point. On fixe un~$k$-plongement de~$F$ dans~$\KK$.

\trois{ucupvdisc} Les composantes connexes de~$(V\cup U)\hotimes_LF\setminus\{x_F\}$ sont :

\medskip
$\bullet$ la sous-couronne~$U':=]x;u[_F^\flat$ de~$U\hotimes_LF$ ; 

$\bullet$ la composante connexe~$W$ de~$V\hotimes_LF\setminus\{x_F\}$ qui contient~$]\xi;x[_F$ ; elle n'aboutit pas à~$\omega_F$,  et~$W\hotimes_F\KK$ est par conséquent réunion disjointe de composantes connexes de~$V\hotimes_L\KK\setminus\{x_{\KK}\}$ qui n'aboutissent pas à~$\omega_{\KK}$, et sont donc des disques ;

$\bullet$ les composantes connexes de~$U\hotimes_LF\setminus\{x_F\}$ qui ne rencontrent pas l'intervalle~$\skel U_F=S(U\hotimes_LF)$ ; chacune d'elle est un disque virtuel (\ref{compcourvirt}).

\trois{conclup1} Prolongeons la~$F$-couronne~$U'$ en un~$F$-disque~$Z$, et choisissons un~$F$-point~$z$ sur~$Z-U'$. En recollant~$Z$ et~$(V\cup U)\hotimes_LF$ le long de~$U'$, on obtient un espace~$F$-analytique lisse, topologiquement séparé, compact et connexe, qui est donc isomorphe à l'analytification d'une~$F$-courbe algébrique projective, irréductible et lisse~${\sch Y}$. On peut dès lors voir~$(V\cup U)\hotimes_LF$ comme une composante connexe de~${\sch Y}\an\setminus\{\zeta\}$, pour un certain point~$\zeta$ de~${\sch Y}\an\dtr$, et~$z$ comme un~$F$-point de~$\sch Y\an$ n'appartenant pas à ~$(V\cup U)\hotimes_LF$. 

\medskip
La courbe projective~${\sch Y}\otimes_F\KK$ est lisse, et irréductible puisque~$x_F$ n'a qu'un antécédent sur~$({\sch Y}\otimes_F\KK)\an$. Les composantes connexes de~$({\sch Y}\otimes_F\KK)\an\setminus\{x_{\KK}\}$ sont, par construction de~$\sch Y$ et en vertu du~\ref{ucupvdisc}, toutes des disques ; de plus le point~$x_{\KK}$, situé sur le squelette de la couronne~$U\otimes_L\KK$, est de genre 0. Il s'ensuit, d'après le lemme~\ref{lemrevg0d}, que~${\sch Y}\otimes_F\KK\simeq \PP^1_{\KK}$. 

\medskip
La~$F$-courbe projective, lisse et irréductible~${\sch Y}$ est ainsi de genre~$0$ ; comme~$\sch Y(F)\neq \emptyset$ (cet ensemble contient~$z$), la courbe~$\sch Y$ est isomorphe à~$\PP^1_F$. 

\medskip
Le point~$\zeta$ appartient à~${\sch Y}\an\dtr$, et l'ouvert~$(U\cup V)\hotimes_LF$ est une composante connexe de~${\sch Y}\an\setminus\{\zeta\}$ qui possède un~$F$-point par choix de~$F$. Par ailleurs, il existe un~$F$-point sur~$\sch Y\an$ qui n'appartient pas à~$(U\cup V)\hotimes_LF$, à savoir~$z$. On en déduit que~$(U\cup V)\hotimes_LF$  est un~$F$-disque (\ref{compp1msd}).~$\Box$

\deux{fuscour} {\bf Lemme.} {\em Soit~$X$ une courbe~$k$-analytique connexe et lisse. Supposons qu'il existe un point~$x$ de genre~$0$ sur~$X$ tel que~$X\setminus\{x\}$ soit réunion disjointe de deux~$k$-couronnes virtuelles~$X'$ et~$X''$ et de disques virtuels ; la courbe~$X$ est alors une~$k$-couronne virtuelle, qui est déployée par toute extension presque algébrique de~$k$ déployant~$X'$ et~$X''$.}

\medskip
{\em Démonstration.} Soit~$L$ une extension presque algébrique de~$k$ déployant~$X'$ et~$X''$. Par sa forme même, la courbe~$X$ est un arbre à deux bouts ; il suffit donc de montrer que~$X_L$ est une~$L$-couronne. 

\medskip
\trois{clotsepkappax} Tout élément de~$\got s(x)$ définit un élément de~$\got s(U)$ pour une certaine sous-couronne~$U$ de~$X'$ (par exemple) ; le corps des constantes d'une~$k$-couronne étant égal à~$k$, il vient~$\got s(x)=k$. Si~$F$ est une extension presque algébrique~$k$, le point~$x$ a donc un unique antécédent~$x_F$ sur~$X_F$. 

\medskip
\trois{prolong2cour} Prolongeons la~$L$-couronne~$X'_L$ (resp.~$X''_L$) en un~$L$-disque~$Z'$ (resp.~$Z''$), sur lequel on choisit un~$L$-point~$z'$ (resp.~$z''$) qui n'est pas situé sur~$X'_L$ (resp.~$Z''_L$). En recollant~$Z'\coprod Z''$
et~$X_L$ le long de~$X'_L\coprod X''_L$, on obtient une courbe~$L$-analytique lisse, connexe, et compacte, et partant isomorphe à~${\sch X}\an$ pour une certaine~$L$-courbe algébrique projective, irréductible et lisse~$\sch X$. Les composantes connexes de~${\sch X}\an_L\setminus\{x_L\}$ sont toutes des disques virtuels, et en particulier des arbres à un bout ; il en résulte que~${\sch X}\an_L$ est un arbre. Il existe par construction deux points~$\eta'$ et~$\eta''$ de~${\sch X}\an$, chacun étant de type 2 ou 3, tels que~$X_L$ s'identifie à la composante connexe de~${\sch X}\an\setminus\{\eta',\eta''\}$ contenant~$]\eta';\eta''[$, et tels qu'il existe deux composantes connexes~$T'$ et~$T''$ de ~${\sch X}\an\setminus\{\eta',\eta''\}$ possédant les propriétés suivantes : 

\medskip
$\bullet$~$T'(L)\neq \emptyset$ et~$T''(L)\neq \emptyset$ ; 

$\bullet$~$\partial T'=\{\eta'\}$ et~$\partial T''=\{\eta''\}$.

\trois{conclucourg0} Choisissons un~$k$-plongement de~$L$ dans~$\KK$. Comme~$x_L$ n'a qu'un antécédent sur~${\sch X}\an_{\KK}$, la~$\KK$-courbe projective et lisse~${\sch X}_{\KK}$ est irréductible. Par construction, toutes les composantes connexes de~${\sch X}\an_{\KK}\setminus\{x_{\KK}\}$ sont des disques ; d'autre part,~$x_{\KK}$ est par hypothèse de genre 0.  Il s'ensuit, d'après le lemme~\ref{lemrevg0d} que~${\sch X}_{\KK}\simeq \PP^1_{\KK}$.

\medskip
La~$L$-courbe~$\sch X$ est de genre 0 ; elle possède un~$L$-point par sa construction même (son analytifiée contient un~$L$-disque), et est donc isomorphe à~$\PP^1_L$ ; on déduit alors de la description de~$X_L$ comme ouvert de~${\sch X}\an$ donnée au~\ref{prolong2cour} et de~\ref{compp1msc} que~$X_L$ est une~$L$-couronne.~$\Box$

\subsection*{Sous-graphes analytiquement admissibles}

\deux{squelanaladm} Soit~$X$ une courbe~$k$-analytique et soit~$\Gamma$ un sous-graphe fermé de~$X$. Nous dirons que~$\Gamma$ est {\em analytiquement admissible} si toute composante connexe de~$X-\Gamma$ est un disque virtuel relativement compact dans~$X$ ; un sous-graphe analytiquement admissible de~$X$ est en particulier admissible, et partant convexe.

\deux{analadmchangeb} Soit~$X$ une courbe~$k$-analytique, soit~$\Gamma$ un sous-graphe analytiquement admissible de~$X$, soit~$F$ une extension complète de~$k$, et soit~$\Gamma_F$ l'image réciproque de~$\Gamma$ sur~$X_F$. 

\trois{imrecgenadm} Si~$U$ est une composante connexe de~$X-\Gamma$, son image réciproque sur~$X_F$ est une réunion finie disjointe de disques virtuels, relativement compacts puisque~$U$ est relativement compacte. On déduit alors de~\ref{remskgrad} que le fermé~$\Gamma_F$ de~$X_F$ en est un sous-graphe ; il est analytiquement admissible d'après ce qui précède. 

\trois{imrecpralgadm} Supposons de plus que~$F$ est une extension presque algébrique de~$k$ et que~$\Gamma$ est tracé sur~$X\geom$. Dans ce cas~$\Gamma_F$ est tracé sur~$X_{F,[0,2,3]}$ ; et si~$\Gamma$ est de plus localement fini, alors~$\Gamma_F$ est localement fini en vertu du corollaire~\ref{corollimrecih}.

\deux{exanaladm} {\em Exemple.} Soit~$X$ un disque virtuel sur~$k$, soit~$\omega$ son unique bout et soit~$\Gamma$ un sous-arbre fermé de~$X$ aboutissant à~$\omega$ ; nous allons montrer que~$\Gamma$ est un sous-graphe analytiquement admissible de~$X$. Soit~$\omega_{\KK}$ l'unique bout de~$X_{\KK}$, et soit~$\Gamma_{\KK}$ l'image réciproque de~$\Gamma$ sur~$X_{\KK}$.

\trois{gammachapk}{\em Montrons que~$\Gamma_{\KK}$ est un sous-arbre admissible de~$X_{\KK}$}. Comme~$\Gamma$ aboutit à~$\omega$, c'est un sous-arbre fermé et {\em non vide} de~$X$ ; par conséquent,~$\Gamma_{\KK}$ est un fermé non vide de~$X_{\KK}$.

\medskip
Soit~$x\in \Gamma_{\KK}$ et soit~$\xi$ son image sur~$X$. L'application~$X_{\KK}\to X$ induit un homéomorphisme de~$[x;\omega_{\KK}[$ sur~$[\xi;\omega[$ (\ref{unboutglobinv}), lequel est contenu dans~$\Gamma$ puisque ce dernier aboutit à~$\omega$. Par conséquent,~$[x;\omega_{\KK}[\subset \Gamma_{\KK}$. 

\medskip
Il s'ensuit que si~$x$ et~$y$ sont deux points de~$\Gamma_{\KK}$ alors~$[x;\omega_{\KK}[\cup [y;\omega_{\KK}[\subset \Gamma_{\KK}$ ; comme~$$[x;y]\subset X_{\KK}\cap( [x;\omega_{\KK}]\cup [y;\omega_{\KK}])=[x;\omega_{\KK}[\cup [y;\omega_{\KK}[,$$ il vient~$[x;y]\subset \Gamma_{\KK}$. Ainsi, le fermé non vide~$\Gamma_{\KK}$ de~$X_{\KK}$ est convexe ; c'est donc un sous-arbre fermé et non vide de~$X_{\KK}$. 

\medskip
Puisque~$[x;\omega_{\KK}[$ est , d'après ce qui précède, contenu dans~$\Gamma_{\KK}$ pour tout~$x\in \Gamma_{\KK}$, le sous-arbre fermé et non vide~$\Gamma_{\KK}$ de~$X_{\KK}$ aboutit à~$\omega_{\KK}$ ; c'est donc, en vertu de~\ref{arbradmarbr}, un sous-graphe admissible de~$X_{\KK}$.

\trois{gammakanad} Soit~$U$ une composante connexe de~$X_{\KK}-\Gamma_{\KK}$. Par ce qui précède,~$U$ est un arbre à un bout relativement compact dans~$X_{\KK}$ ; soit~$x$ l'unique point de~$\partial U$. Comme~$X_{\KK}$ n'est pas compact, il existe une composante connexe de~$X\setminus\{x\}$ qui n'est pas relativement compacte, et en particulier diffère de~$U$ ; le point~$x$ n'est dès lors pas unibranche, ce qui signifie qu'il est de type 2 ou 3. 

\medskip
Immergeons, au moyen d'une fonction coordonnée, le~$\KK$-disque~$X_{\KK}$ dans~$\pkk$. Étant relativement compact dans~$X_{\KK}$, l'ouvert~$U$ est une composante connexe de~$\pkk\setminus\{x\}$ ; il s'ensuit,~$x$ étant de type 2 ou 3, que~$U$ est un disque. 

\trois{concluadm} Ainsi toute composante connexe de~$X_{\KK}-\Gamma_{\KK}$ est-elle un disque relativement compact ; par conséquent, toute composante connexe de~$X-\Gamma$ est un disque virtuel relativement compact ; autrement dit,~$\Gamma$ est un sous-arbre analytiquement admissible de~$X$. 

\deux{graphanaladm} Soit~$X$ une courbe~$k$-analytique et soit~$\Gamma$ un sous-graphe analytiquement admissible de~$X$. 

\trois{bord23} Si~$V$ est une composante connexe de~$X-\Gamma$ son bord est un singleton~$\{x\}$, où~$x\in \Gamma$. Il résulte de
\ref{borddisct23}
que~$x\in X\dtr$
et que~$X$ est quasi-lisse en~$x$.

\trois{extanaladm} Si~$\Delta$ est un sous-graphe admissible de~$X$ contenant~$\Gamma$, alors~$\Delta$ est analytiquement admissible. En effet, soit~$V$ une composante connexe de~$X-\Delta$ ; nous allons montrer que c'est un disque virtuel, ce qui permettra de conclure. Soit~$x$ l'unique point de~$\partial V$. 

\medskip
Si~$x\in \Gamma$ alors~$V$ est une composante connexe de~$X-\Gamma$, et est donc un disque virtuel puisque~$\Gamma$ est analytiquement admissible. 

\medskip
Sinon, soit~$Y$ la composante connexe de~$X-\Gamma$ contenant~$x$ ; elle contient~$V$, qui apparaît donc comme une composante connexe de~$Y-\Delta$ ; soit~$y$ l'unique point de~$\partial Y$. 

L'intersection~$Y\cap \Delta$ est une partie fermée et convexe de~$Y$, c'est-à-dire un sous-arbre fermé de~$Y$, qui contient~$x$. Par ailleurs, comme~$y$ appartient à~$\Gamma$, il appartient à~$\Delta$ ; par convexité de~$\Delta\cap \overline Y$, l'intervalle~$[x;y[$ est contenu dans~$\Delta\cap Y$ ; ainsi,~$\Delta\cap Y$ aboutit à l'unique bout de~$Y$. 

\medskip
On déduit alors de l'exemple~\ref{exanaladm} que~$\Delta\cap Y$ est un sous-arbre analytiquement admissible de~$Y$ ; en conséquence, l'ouvert relativement compact~$V$ est un disque virtuel, ce qui achève de montrer que~$\Delta$ est analytiquement admissible.

\trois{corolextanal} Soit~$\Sigma$ un sous-graphe fermé de~$X$ et soit~$X'$ la normalisée de~$X$. Il existe un plus petit sous-graphe analytiquement admissible~$\Delta$ de~$X$ contenant~$\Gamma\cup \Sigma$, qui est localement fini si~$\Gamma$ et~$\Sigma$ sont eux-mêmes localement finis. De plus : 

\medskip
$\bullet$ si~$\Gamma$ et~$\Sigma$ sont tracés sur~$X\geom$ alors~$\Delta\subset X\geom$ ;  

$\bullet$ si~$\Gamma$ et~$\Sigma$ sont tracés sur~$X\dtr$ et si~$X'\to X$ est un homéomorphisme alors~$\Delta\subset X\dtr$. 

\medskip
En effet, la proposition ~\ref{pluspetitadm} assure l'existence d'un plus petit sous-graphe {\em admissible}~$\Delta$ de~$X$ contenant ~$\Gamma\cup \Sigma$, qui est localement fini si~$\Gamma$ et~$\Sigma$ le sont. Et elle fournit par ailleurs une description explicite de~$\Delta$, d'où il découle : que si~$\Gamma$ et~$\Sigma$ sont tracés sur~$X\geom$ il en va de même de~$\Delta$ par convexité de~$X\geom$  ; et que si~$\Gamma$ et~$\Sigma$ sont tracés sur~$X\dtr$ et si~$X'\to X$ est un homéomorphisme alors~$\Delta\subset X\dtr$, l'hypothèse faite sur~$X'\to X$ garantissant la convexité de~$X\dtr$. 

\medskip
En vertu de ~\ref{extanaladm}, le sous-graphe~$\Delta$ de~$X$ est analytiquement admissible ; comme tout sous-graphe analytiquement admissible de~$X$ est admissible,~$\Delta$ est bien le plus petit sous-graphe analytiquement admissible de~$X$ contenant~$\Gamma\cup \Sigma$.

\trois{deplofini} Soit~$V$ un ouvert de~$X$ qui est un disque virtuel relativement compact. Nous allons montrer que l'ouvert~$V$ est alors un disque {\em gentiment} virtuel. En effet, soit~$x$ l'unique point de~$\partial V$ ; en vertu de~\ref{borddisct23}, le point~$x$ appartient à~$X\dtr$ et au lieu quasi-lisse de~$X$. Comme~$V$ est un arbre à un bout, il contient une unique branche~$b$ de~$X$ issue de~$x$. Il existe une section~$Z$ de~$b$ contenue dans~$V$ qui est une couronne gentiment virtuelle. Par ailleurs,~$V$ possède un point rationnel sur une extension finie de~$\got s(V)$, et même finie séparable si~$|k\ti|\neq \{1\}$ : si~$|k\ti|=\{1\}$, cela est dû d'une part au fait que~$V\hotimes_{\got s(V)}\KK$ est un~$\KK$-disque, et donc possède un~$\KK$-point, et d'autre part au fait que~$\KK$ est alors une clôture algébrique de~$k$ ; si~$|k\ti|\neq\{1\}$, c'est une simple conséquence de la lissité de~$V$. 

Il découle dès lors du lemme~\ref{fusion} que~$\got s(Z)=\got s(V)$ et que le disque virtuel~$V$ (qui est égal à~$Z\cup V$) est déployé par une extension finie de~$\got s(V)$, et même finie séparable si~$|k\ti|\neq\{1\}$ ; autrement dit, c'est un disque gentiment virtuel.

\medskip
Notons que ce qui précède s'applique en particulier aux composantes connexes de~$X-\Gamma$ ; ce sont donc toutes dess disques gentiment virtuels.  

\deux{squelanalgenql} Soit~$X$ une courbe~$k$-analytique {\em génériquement} quasi-lisse et soit~$S$ le lieu singulier de~$X$ ; le sous-ensemble~$S$ de de~$X$ est fermé, discret, et contenu dans~$X\typ 0$. Soit~$F$ le complété d'une extension radicielle de~$k$, et soit~$X'$ la normalisée de~$X_F$ ; on note~$S'$ l'image réciproque de~$S$ sur~$X'$ ; c'est un sous-ensemble fermé et discret de~$X'$ contenu dans~$X'\typ 0$. La flèche~$X'\to X$ identifie topologiquement~$X$ au quotient de~$X'$ par la relation d'équivalence~$\sch R$ qu'elle définit ; elle induit un homéomorphisme~$X'\setminus S'\simeq X\setminus S$, et les fibres de~$S'\to S$ sont finies.

\trois{imgranad} Soit~$\Gamma'$ un sous-graphe analytiquement admissible de~$X'$ contenu dans~$X'\dtr$. Comme~$S'$ est un sous-graphe localement fini de~$X'$, il existe un plus petit sous-graphe analytiquement admissible~$\Delta'$ de~$X'$ qui contient~$S'\cup \Gamma'$, et~$\Delta'$ est localement fini dès que~$\Gamma'$ est localement fini (\ref{corolextanal}). Les fermés~$\Delta'$ et~$\Gamma'$ de~$X'$ sont saturés sous~$\sch R$ ; leurs images respectives~$\Delta$ et~$\Gamma$ sur~$X$ sont donc fermées, et s'identifient respectivement aux quotients de~$\Delta'$ et~$\Gamma'$ par la restriction de~$\sch R$ (notons que comme~$\Gamma'\subset X'\dtr$, il évite~$S'$, d'où il découle que~$\Gamma'\simeq \Gamma$) ; ce sont dès lors deux sous-graphes fermés de~$X$, qui sont localement finis si~$\Gamma'$ est localement fini ; nous allons montrer que~$\Delta$ est analytiquement admissible, contenu dans~$X\geom$, et satisfait l'égalité~$\Delta\cap X\typ 0=S$. 

\medskip
Par définition,~$\Gamma'$ est un sous-graphe admissible de~$X'$ qui évite~$S'$ ; soit~$r$ la rétraction canonique de~$X'$ sur~$\Gamma'$. Comme~$\Delta'$ est le plus petit sous-graphe admissible de~$X'$ contenant~$\Gamma'\cup S'$ (\ref{corolextanal}), il coïncide en vertu de la proposition~\ref{pluspetitadm}) avec~$\Gamma'\cup\bigcup\limits_{x\in S'}[x;r(x)[$. 

\medskip
Le graphe~$\Delta$ est égal à la réunion des sous-ensembles suivants de~$X$ : 

\medskip
$\bullet$ le graphe~$\Gamma$ qui, en tant qu'image de~$\Gamma'$, est contenu dans~$X\dtr$ ; 

$\bullet$ l'ensemble~$S$, contenu dans~$X\typ 0$ ; 

$\bullet$ la réunion~$\Sigma$ des images des intervalles~$]x;r(x)[$, où~$x$ parcourt~$S'$ ; comme~$X'$ est normale, chacun des~$]x;r(x)[~$ est contenu dans~$X'\dtr$ en vertu de~\ref{conn23rig} ; par conséquent,~$\Sigma\subset X\dtr$. 

\medskip
Il s'ensuit que~$\Delta\subset X\geom$, et que~$\Delta\cap X\typ 0=S$. Soit maintenant~$U$ une composante connexe de~$X-\Delta$ et soit~$U'$ son image réciproque sur~$X'$ ; c'est une réunion de composantes connexes de~$X'-\Delta'$. Comme~$S\subset \Delta$, l'ouvert~$U$ est contenu dans~$X\setminus S$ ; par conséquent,~$U'$ s'identifie à~$U_F$, et est en particulier connexe. On déduit de ce dernier fait que~$U'$ est une composante connexe de~$X'-\Delta'$ ; le sous-graphe~$\Delta'$ de~$X'$ étant analytiquement admissible,~$U'$ est un disque virtuel relativement compact dans~$X'$. La relative compacité de~$U'$ dans~$X'$ entraîne celle de~$U$ dans~$x$ ; et comme~$U'\simeq U_F$, la composante~$U$ est elle-même un disque virtuel, ce qui achève de montrer que~$\Delta$ est analytiquement admissible.

\subsection*{Le squelette analytique}

\deux{defsquelan} Si~$X$ est une courbe~$k$-analytique, l'ensemble des points de~$X$ possédant un voisinage dans~$X$ qui est un disque virtuel est un ouvert de~$X$ ; son fermé complémentaire sera appelé le {\em squelette analytique} de~$X$ et sera noté~$\skelan X$ ; c'est en vertu de~\ref{remskgrad} un sous-graphe fermé et convexe de~$X$. 

Il résulte des définitions que~$\skelan X \supset \skel X$ et que~$\skelan X$ est contenu dans tout sous-graphe analytiquement admissible de~$X$. 

\deux{squelanpseucour }{\em Exemple.} Soit~$X$ une couronne virtuelle sur~$k$. Par ce qui précède,~$\skelan X\supset \skel X$ ; et l'on déduit de~\ref{compcourvirt} que~$\skelan X\subset \skel X$ ; par conséquent,~$\skelan X=\skel X$. 

\deux{descgamma} Soit~$X$ une courbe~$k$-analytique {\em quasi-lisse} et soit~$x\in \skelan X$. On se donne un voisinage~$V$ de~$x$ de la forme (dépendant du type de~$x$) décrite par le théorème~\ref{theovoisql}. 

\trois{xpast14} Si~$x$ est de type 1 ou 4 alors~$V$ est un disque virtuel, ce qui contredit l'appartenance de~$x$ à~$\skelan X$ ; par conséquent,~$x$ est de type 2 ou 3. 

\trois{gamma23} Dans le cas du type 2 comme dans celui du type 3, l'adhérence de~$V$ dans~$X$ est un arbre compact, et~$\pi_0(V\setminus\{x\})$ peut s'écrire~${\sch C}\coprod \sch D$, où~$\sch C$ est fini et constitué de couronnes gentiment virtuelles, et où~$\sch D$ est constitué de disques virtuels ; notons que si~$x$ est de type 2 (resp. 3) alors~$\sch D$ est infini (resp. vide). Soit~$W\in\Pi_0(V\setminus\{x\})$. 
 
\medskip
{\em Supposons  que~$W\in \sch D$.} Dans ce cas,~$W\cap \skelan X=\emptyset$ par définition de~$\skelan X$.

\medskip
{\em Supposons que~$W\in \sch C$}. Les composantes connexes de~$W-\skel W$ sont des disques virtuels ; ceci entraîne que~$\skelan X\cap W$ est contenu dans~$\skel W$. 

Nous allons montrer que~$\skelan X\cap W$ est ou bien égal à~$\skel W$ tout entier, ou bien vide. Supposons donc que~$\skelan X\cap W$ est {\em strictement} contenu dans~$\skel W$. Il existe alors un point~$w$ de~$\skel W$ qui n'appartient pas à~$\skelan X$ ; cela signifie que~$w$ possède un voisinage~$W'$ dans~$X$ qui est un disque virtuel ; le point~$x$ étant situé sur~$\skelan X$, il n'appartient pas à~$W'$. Il s'ensuit, en vertu du lemme~\ref{fusion}, que~$W'\cup W$ est un disque virtuel, et une composante connexe de~$X\setminus\{x\}$ ; c'est donc {\em la} composante connexe de~$X\setminus\{x\}$ qui contient~$W$ ; celle-ci ne dépend ni du choix de~$w$ ni de celui de~$W'$ ; on la notera~$\varpi(W)$. Comme~$W$ est contenu dans le disque virtuel~$\varpi(W)$, on a~$\skelan X\cap W=\emptyset$, comme annoncé. 

\trois{propvprime} Soit~$V'$ le voisinage ouvert de~$x$ défini comme la réunion de~$V$ et des~$\varpi(W)$ où~$W$ parcourt l'ensemble des éléments de~$\sch C$ qui ne rencontrent pas~$\Gamma$. De ce qui précède découlent les faits suivants : 

$\bullet$~$V'$ est un arbre, et son adhérence dans~$X$ est un arbre compact ; 

$\bullet$ l'ensemble~$\pi_0(V'\setminus\{x\})$ est la réunion disjointe de deux ensembles~${\sch C}'$ et~${\sch D}'$, où~${\sch C}'$ est fini et constitué de couronnes virtuelles, et où~${\sch D}'$ est constitué de disques virtuels ; 

$\bullet$~$\skelan X\cap V'=\{x\}\cup \bigcup\limits_{W\in {\sch C}'}\skel W.$

$\bullet$ toute composante connexe de~$V'-\Gamma$ est un disque virtuel relativement compact dans~$V'$ ; c'est donc un disque virtuel relativement compact dans~$X$ dont l'unique point du bord est situé sur~$\skelan X \cap V'$.

\trois{propvprimebis} Soit~$U$ une composante connexe de~$X-\skelan X$ dont le bord contient~$x$. Son intersection avec~$V'$ est alors non vide ; comme elle est à la fois ouverte et fermée dans~$V'-\skelan X$, c'est une réunion non vide de composantes connexes de~$V'-\skelan X$. Il existe en conséquence au moins une composante connexe ~$U'$ de~$V'-\skelan X$ qui est contenue dans~$U$. D'après ce qui précède,~$\partial_X U'$ est un singleton contenu dans~$\skelan X$ ; il s'ensuit que~$U'$ est fermée dans~$U$, et donc égale à~$U$ par connexité de cette dernière. Ainsi,~$U$ est un disque virtuel relativement compact.

\deux{theoanaladm} {\bf Théorème.} {\em Soit~$X$ une courbe~$k$-analytique génériquement quasi-lisse et soit~$S$ son lieu singulier. 

\medskip

1) Le bord analytique de~$X$ est contenu dans~$\skelan X$. 

2) Le sous-graphe~$\skelan X$ de~$X$ est localement fini, contenu dans~$X\geom$, et~$\skelan X \cap X \typ 0 =S$ ; en particulier,~$\skelan X\subset X\dtr$ si et seulement si~$X$ est quasi-lisse. 

3) Si~$X$ est quasi-lisse et si~$\skelan X$ rencontre chaque composante connexe de~$X$, il est analytiquement admissible. 

4) La courbe~$X$ possède un sous-graphe localement fini et analytiquement admissible contenu dans~$X\geom$ et dont l'intersection avec~$X\typ 0$ est égale à~$S$.}

\medskip
{\em Démonstration.} Si~$x\in X-\skelan X$ alors~$x$ est contenu dans un disque virtuel, qui est une courbe sans bord ; par conséquent,~$x\notin \partial \an X$, d'où 1). 

\trois{analadmql} {\em Preuve de 2) dans le cas où~$X$ est quasi-lisse.} Il découle de~\ref{xpast14} que~$\skelan X\subset X\dtr$ ; on déduit par ailleurs de~\ref{propvprime} que~$\skelan X$ est localement fini. 

\medskip
\trois{skeladmql} {\em Preuve de 3) dans le cas où~$X$ est quasi-lisse.} Supposons que~$\skelan X$ rencontre toutes les composantes connexes de~$X$. Soit~$U$ une composante connexe de~$X-\skelan X$. En vertu de notre hypothèse,  le bord de~$U$ est non vide ; il s'ensuit, d'après~\ref{propvprimebis}, que~$U$ est un arbre virtuel relativement compact ; par conséquent,~$\skelan X$ est admissible. 

\trois{3casql}{\em Preuve de 4) dans le cas où~$X$ est quasi-lisse.} On se ramène en raisonnant composante par composante au cas où~$X$ est connexe et non vide. Si~$\skelan X$ est non vide, il répond au problème posé en vertu des  assertions 2) et 3) déjà établies. 

\medskip
Supposons maintenant que~$\skelan X=\emptyset$. Tout point de~$X$ a alors un voisinage qui est un disque virtuel. Comme~$X$ est non vide, il existe un ouvert~$U$ de~$X$ qui est un disque virtuel, et l'on peut toujours, quitte à restreindre~$U$, le supposer relativement compact (\ref{relcompdiscvirt}) ; soit~$u$ l'unique point de~$\partial U$. On sait (\ref{psdtop} et~\ref{existifixpp}-\ref{discvcourv}) qu'il existe un intervalle ouvert~$I$ tracé sur~$U\dtr$, aboutissant à~$u$ et qui possède la propriété suivante : {\em si~$x\in I$ et si~$V$ est une composante connexe de~$U\setminus\{x\}$ qui ne contient pas~$]x;u[$ (resp. qui contient~$]x;u[$) alors~$V$ est un disque virtuel (resp. une couronne virtuelle).} 

\medskip
Soit~$x\in I$ et soit~$\Pi$ l'ensemble des composantes connexes de~$U\setminus\{x\}$ ne contenant pas~$]x;u[$ ; on appelle~$W$ la composante connexe de~$U\setminus\{x\}$ contenant~$]x;u[$ ; c'est une couronne virtuelle. Si~$V\in \Pi$, son adhérence dans~$U$ est égale à~$V\cup\{x\}$ et est compacte ; cela implique que~$V$ est une composante connexe de~$X\setminus\{x\}$. Par ailleurs,~$x$ adhère à toute composante connexe de~$X\setminus\{x\}$ en vertu de la connexité de~$X$. Il en résulte que les composantes connexes de~$X\setminus\{x\}$ sont d'une part les~$V\in \Pi$ et d'autre part {\em la} composante connexe~$W'$ de~$X\setminus\{x\}$ qui contient~$W$. On distingue maintenant deux cas. 

\medskip
{\em Supposons qu'il existe~$y\in ]x;u[$ possédant un voisinage~$V'$ dans~$W'$ qui soit un disque virtuel}. Dans ce cas, le lemme~\ref{fusion} assure que~$W\cup V'$ est un disque virtuel, et une composante connexe de~$X\setminus\{x\}$, qui ne peut être que~$W'$. Dès lors, toutes les composantes connexes de~$X\setminus\{x\}$ sont des disques virtuels, et~$\{x\}$ est un sous-graphe localement fini et analytiquement admissible de~$X$, contenu dans~$X\dtr$. 

\medskip
{\em Supposons qu'aucun point~$y$ de~$]x;u[$ ne possède un voisinage dans~$W'$ qui soit un disque virtuel.} Dans ce cas, le squelette analytique~$\Gamma$ de~$W'$ contient~$]x;u[$ ; comme~$W'-]x;u[$ est réunion disjointe de disques virtuels,~$\Gamma\cap W=]x;u[$. Posons~$\Gamma'=\Gamma\cup\{x\}$ ; c'est un sous-graphe fermé et localement fini de~$X$. Il est contenu dans~$X\dtr$ : cela résulte du fait que~$x\in X\dtr$, et que~$\Gamma$ est, en tant que squelette analytique de l'ouvert quasi-lisse~$W'$, tracé sur~$X\dtr$ en vertu de l'assertion 2) déjà établie dans le cas quasi-lisse. 

Les composantes connexes de~$X-\Gamma'$ sont d'une part les composantes connexes de~$X\setminus\{x\}$ qui ne contiennent pas~$]x;u[$, d'autre part les composantes connexes de~$W'-\Gamma$ ; ce sont toutes des disques virtuels relativement compacts dans~$X$. Par conséquent,~$\Gamma'$ est un sous-graphe localement fini et analytiquement admissible de~$X$, contenu dans~$X\dtr$. 

\trois{analadmcagen} {\em Preuve de 2), 3) et 4) dans le cas général.} Soit~$F$ le complété de la clôture parfaite  de~$k$ et soit~$X'$ la normalisée de~$X_F$ ; soit~$S'$  l'image réciproque de~$S$ sur~$X'$. Comme~$F$ est parfait, la courbe~$X'$ est quasi-lisse. 

\medskip
D'après l'assertion 3) déjà établie,  la~$F$-courbe quasi-lisse~$X'$ possède un sous-graphe analytiquement admissible et localement fini~$\Gamma'$ tracé sur~$X'\dtr$. Soit~$\Delta'$ le plus petit sous-graphe analytiquement admissible de~$X'$ contenant~$\Gamma'\cup S'$, et soit~$\Delta$ l'image de~$\Gamma'$ sur~$X$. En vertu de~\ref{imgranad},~$\Delta$ est un sous-graphe localement fini et analytiquement admissible de~$X$, tracé sur~$X\geom$ et dont l'intersection avec~$X\typ 0$ est égale à~$S$ ; cela prouve 3). 

\medskip
Montrons maintenant 1). L'assertion 3) que l'on vient de prouver affirme l'existence d'un sous-graphe~$\Delta$ de~$X$, localement fini, analytiquement admissible, tracé sur~$X\geom$ et tel que~$\Delta\cap X\typ 0=S$. Le sous-graphe~$\skelan X$ étant contenu dans~$\Delta$, il est localement fini, tracé sur~$X\geom$, et ~$\skelan X\cap X\typ 0\subset S$. Par ailleurs si~$x\in S$ alors~$x$ est un point singulier de~$X$, et ne peut donc être contenu sur aucun disque virtuel ; il s'ensuit que~$x\in \skelan X$, et l'on a finalement~$\skelan X\cap X\typ 0=S$.~$\Box$ 

\subsection*{Les nœuds d'un sous-graphe localement fini et analytiquement admissible}

\deux{defnoeudsg} Soit~$X$ une courbe~$k$-analytique génériquement quasi-lisse et soit~$\Gamma$ un sous-graphe analytiquement admissible et localement fini de~$X$, tracé sur~$X\geom$. Nous dirons qu'un point~$x$ de~$\Gamma$ en est un {\em nœud} s'il satisfait {\em l'une au moins} des conditions suivantes : 

\medskip

1) le point~$x$ est rigide ; 

2) le point~$x$ appartient à~$\partial\an X$ ; 

3) le point~$x$ est un sommet topologique de~$\Gamma$, c'est-à-dire, rappelons-le, que la valence de~$(\Gamma,x)$ est différente de 2 ; 

4) il existe une branche~$\beta$ de~$\Gamma$ issue de~$x$ telle que~$\got s(\beta)$ soit une extension {\em stricte} de~$\got s(x)$ (autrement dit, les conditions équivalentes i), ii), iii) et iv) du~\ref{compsxsurgamma} ci-dessus ne sont pas vérifiées) ; 

5) le point~$x$ est de genre strictement positif.

\trois{remnoeudst2} Remarquons que 4) et 5) impliquent que~$x$ est de type 2 ; c'est évident pour 5), justifions-le pour 4). Si~$x$ est de type 3, il possède un voisinage~$V$ que l'on peut immerger dans une couronne virtuelle~$U$ de sorte que~$x\in \skel U$ et de sorte que~$V\setminus\{x\}$ soit réunion de~$0,1$ ou~$2$ composantes connexes de~$U\setminus\{x\}$ ; comme~$\br U x$ a pour cardinal 2, il existe un voisinage ouvert~$\Gamma_0$ de~$x$ dans~$\Gamma$ tel que~$\Gamma_0\subset \skel U$; mais on a alors~$\got s(y)\simeq \got s(U)\simeq \got s(x)$ pour tout~$y\in \Gamma_0$, et 4) n'est pas satisfaite.

\trois{remnoeudspasinclus} Il résulte de~\ref{apropossx} et~\ref{fermdiscgpos} que l'ensemble des nœuds de~$\Gamma$ en est une partie fermée et discrète. Si~$\Gamma'$ est un sous-arbre localement fini et analytiquement admissible de~$X$ contenant~$\Gamma$ et si~$x$ est un nœud de~$\Gamma$ il découle aussitôt des définitions que~$x$ est un nœud de~$\Gamma'$, sauf peut-être si les conditions suivantes sont satisfaites :~$x\in X\dtr\setminus \partial \an X, g(x)=0$ et~$x$ est un point isolé de~$\Gamma$ ou un point unibranche de~$\Gamma$ tel que~$\got s(b)=\got s(x)$ pour l'unique branche~$b$ de~$\br \Gamma x$. 

\trois{pretriang} Soit~$I$ un intervalle ouvert non vide tracé sur~$\Gamma$, relativement compact dans ce dernier et n'en contenant aucun nœud ; il n'en contient en particulier aucun sommet topologique, et est donc un ouvert de~$\Gamma$ ; supposons de plus que~$\partial I\subset X\dtr$. Nous allons montrer que~$I^\flat$ est une couronne gentiment virtuelle, évidemment relativement compacte.

\medskip
On déduit des hypothèses faites sur~$I$ l'existence d'une suite finie ~$$x_1<x_2<\ldots<x_n$$ de points de~$I$ (arbitrairement orienté) telle que l'on ait, si l'on note~$x_0$ et~$x_{n+1}$ les deux bouts de l'arbre~$I$ : 

\medskip
$\bullet$ pour tout~$i\in\{0,\ldots,n\}$ l'ouvert~$]x_i;x_{i+1}[^\flat$ est une couronne gentiment virtuelle ; 

$\bullet$ pour tout~$i\in \{1,\ldots,n\}$, les plongements~$$\got s(x_i)\hookrightarrow \got s(]x_{i-1};x_i[^\flat)\;{\rm et}\;\got s(x_i)\hookrightarrow \got s(]x_i;x_{i+1}[^\flat)$$ (\ref{remseccorong}) sont des isomorphismes. 

\medskip
Soit~$i\in\{1,\ldots,n\}$. Les composantes connexes de~$]x_{i-1};x_{i+1}[^\flat\setminus\{x_i\}$ autres que~$]x_{i-1};x_i[^\flat$ et~$]x_i;x_{i+1}[^\flat$ sont des disques virtuels ; il en découle, en vertu de~\ref{remsxgen}, que l'évaluation induit un isomorphisme entre~$\got s(]x_{i-1};x_{i+1}[^\flat)$ et~$\got s(x_i)$. On a donc pour tout~$i$ compris~$1$ et~$n$ des isomorphismes naturels~$$\got s(]x_{i-1},x_i[^\flat)\simeq  \got s(]x_{i-1};x_{i+1}[^\flat)\simeq \got s(]x_i;x_{i+1}[^\flat).$$ Il en résulte que pour tout élément~$i$ de~$ \{0,\ldots,n\}$, toute fonction appartenant à~$\got s(]x_i;x_{i+1}[^\flat)$ s'étend naturellement en une fonction appartenant à~$\got s(I^\flat)$ ; autrement dit, si~$F$ désigne le corps~$\got s(I^\flat)$ alors~$\got s(]x_i;x_{i+1}[^\flat)\simeq F$ pour tout~$i\in \{0,\ldots,n\}$.

\medskip
La courbe~$F$-analytique~$I^\flat$ est lisse (qu'elle soit sans bord résulte du fait que~$I$ ne contient aucun nœud de~$\Gamma$) ; pour tout~$i$ compris entre~$1$ et~$n$, les composantes connexes de~$I^\flat\setminus\{x_i\}$ sont d'une part~$]x_0;x_i[^\flat$ et~$]x_i;x_{n+1}[^\flat$, et d'autre part des disques virtuels ; par ailleurs, l'ouvert~$]x_i;x_{i+1}[\flat$ est pour tout~$i$ compris entre~$0$ et~$n$ une couronne gentiment virtuelle sur~$F$. 

\medskip
Une application répétée du lemme~\ref{fuscour} assure alors que~$I^\flat$ est une couronne gentiment virtuelle sur~$F$.

\subsection*{Triangulations}

\deux{deftri} Si~$X$ est une courbe~$k$-analytique quasi-lisse, nous appellerons {\em triangulation} de~$X$ la donnée d'un sous-ensemble~$S$ fermé et discret de~$X$, contenu dans~$X\dtr$, et  tel que toute composante connexe de~$X\setminus S$ soit un disque virtuel ou une couronne virtuelle et soit relativement compacte dans~$X$ ; on dira que l'ensemble~$S$ est l'ensemble des {\em sommets} de la triangulation en question. {\em Attention} : si~$U$ est une composante connexe de~$X\setminus S$ qui est une couronne virtuelle, on ne demande pas que son adhérence soit un arbre : elle peut être homotope à un cercle si~$\partial U$ est un singleton. 

\medskip
Remarquons que si~$S$ est une triangulation de~$X$, toute  composante connexe de~$X\setminus S$ a un bord non vide ; cela implique que~$S$ rencontre toutes les composantes connexes de~$X$. 

\deux{theotri}{\bf Théorème.} {\em Soit~$X$ une courbe~$k$-analytique quasi-lisse. 

\medskip
i) Si~$S$ est une triangulation de~$X$, la réunion~$\Gamma$ de~$S$ et des squelettes des composantes connexes de~$X\setminus S$ qui sont des couronnes virtuelles est un sous-graphe localement fini et analytiquement admissible de~$X$ tracé sur~$X\dtr$ ;  tous les nœuds de~$\Gamma$ appartiennent à~$S$, et~$\Gamma\setminus S$ est réunion disjointe d'intervalles ouverts d'adhérence compacte dans~$\Gamma$ ; on dira que~$\Gamma$ est le {\em squelette} de la triangulation~$S$. 

\medskip
ii) Réciproquement, soit~$\Gamma$ un sous-graphe localement fini et analytiquement admissible de~$X$ tracé sur~$X\dtr$ et soit~$S$ une partie fermée et discrète de~$\Gamma$ contenant ses nœuds et telle que~$\Gamma\setminus S$ soit réunion disjointe d'intervalles ouverts d'adhérence compacte dans~$\Gamma$. L'ensemble~$S$ est alors une triangulation de~$X$ dont~$\Gamma$ est le squelette, et toute composante connexe de~$X\setminus S$ est un disque ou une couronne gentiment virtuel(le). 

\medskip
iii) Si~$\Sigma$ est un sous-ensemble discret et fermé de~$X$ contenu dans~$X\dtr$, il existe une triangulation de~$X$ contenant~$\Sigma$ ; en particulier,~$X$ possède au moins une triangulation ; et si~$S$ est une triangulation de~$X$ alors toute composante connexe de~$X\setminus S$ est un disque ou une couronne {\em gentiment} virtuel(le).  

\medskip
iv) Si~$X$ est strictement~$k$-analytique, si~$|k\ti|\neq\{1\}$ et si~$\Sigma$ est un sous-ensemble discret et fermé de~$X$ constitué de points de type 2 alors il existe une triangulation de~$X$
contenant~$\Sigma$
et
dont tous les sommets sont de type 2 ; en particulier,~$X$ possède une triangulation dont tous les sommets sont de type 2.}

\medskip
{\em Démonstration}. Pour établir i), il suffit de vérifier que si~$x\in S$, il n'y a qu'un nombre fini de composantes connexes de~$X\setminus S$ aboutissant à~$x$ qui sont des couronnes virtuelles ; or c'est une conséquence du théorème~\ref{theovoisql} et du fait que tout voisinage de~$x$ dans le graphe~$X$ contient presque toutes les composantes connexes de~$X\setminus\{x\}$.

\medskip
L'assertion ii) est une reformulation du~\ref{deplofini} et du~\ref{pretriang} ; l'assertion iii) résulte de i), de ii), de l'existence d'un sous-graphe localement fini et analytiquement admissible de~$X$ tracé sur~$X\dtr$ et contenant~$\Sigma$ (th.~\ref{theoanaladm} et~\ref{corolextanal}), et de la métrisabilité d'un tel sous-graphe (prop.~\ref{metrinat}). 

\medskip 
Prouvons maintenant iv) ; on suppose donc que la courbe~$X$ est strictement~$k$-analytique ; il existe un sous-graphe localement fini et analytiquement admissible~$\Gamma$ de~$X$ tracé sur~$X\dtr$ et contenant~$\Sigma$ ({\em cf. supra}). Soit~$x$ un point de type 3 de~$\Gamma$. Comme~$X$ est strictement~$k$-analytique,~$x\notin \partial \an X$ ; par conséquent, il existe un voisinage ouvert~$V$ de~$x$ dans~$X$ qui est une couronne virtuelle telle que~$x\in \skel V$.

 Étant de type 3, le point~$x$ ne peut satisfaire les propriétés 3) et 4) du~\ref{defnoeudsg} ; l'on déduit de ce qui précède que~$x$ ne peut être un nœud de~$\Gamma$ que si le nombre de branches de~$\Gamma$ issues de~$x$ est inférieur ou égal à~$1$. Supposons que ce soit le cas ; il existe alors un intervalle ouvert~$I$ tracé sur~$\skel V$ aboutissant à~$x$ et qui ne rencontre pas~$\Gamma$. Comme~$|k\ti|\neq \{1\}$, les points de type 2 sont denses dans le squelette de~$V$, et {\em a fortiori} dans~$I$ ; choisissons un point~$y$ de type 2 sur~$I$ ; la réunion de~$\Gamma$ et de~$[y;x[$ est un sous-graphe localement fini admissible de~$X$ ; il est par conséquent analytiquement admissible (\ref{extanaladm}). Ainsi, l'on peut prolonger~$\Gamma$ au niveau de chacun de ses nœuds de type 3 de manière à obtenir un sous-arbre localement fini et analytiquement admissible~$\Delta$ de~$X$ contenant~$\Sigma$ et {\em dont tous les nœuds sont de type 2}. 

\medskip
Soit~$S_0$ la réunion de~$\Sigma$ et de l'ensemble des nœuds de~$\Delta$. La métrisabilité de~$\Delta$  (prop.~\ref{metrinat}) assure que le complémentaire de~$S_0$ dans~$\Delta$ est réunion disjointe de cercles et d'intervalles ouverts ; si~$J$ est l'un d'eux et si~$J$' est un intervalle ouvert tracé sur~$J$ dont l'adhérence est un segment, l'ouvert~$(J')\an$  est une couronne virtuelle (\ref{pretriang}) et les points de type 2 sont denses dans son squelette~$J'$ puisque la valeur absolue de~$k$ n'est pas triviale ; par conséquent, les points de type 2 sont denses dans~$\Delta\setminus S_0$. En conséquence, il existe un sous-ensemble fermé et discret~$S$ de~$\Delta$ contenant~$S_0$, constitué uniquement de points de type 2, et tel que~$\Delta\setminus S$ soit réunion disjointe d'intervalle ouverts relativement compacts dans~$\Delta$ ; en vertu de ii) (déjà établi)~$S$ est une triangulation de~$X$, dont tous les sommets sont par construction de type 2.~$\Box$

\deux{rayextshrx} {\bf Lemme.} {\em Soit~$X$ une courbe~$k$-analytique, soit~$U$ un ouvert de~$X$ et soit~$x\in \partial U$ ; supposons que~$x$ est de type 2 ou 3.

1) Si~$U$ est un~$k$-disque et si~$r$ désigne son rayon modulo~$|k\ti|$, on a l'inclusion~$|\hres(x)\ti|\subset |k\ti|.r^\ZZ$, avec égalité si~$k$ est algébriquement clos. Le point~$x$ est de type 2 (resp. 3) si et seulement si~$r\in \sqrt{|k\ti|}$ (resp.~$r\in \RR\ti_+- \sqrt{|k\ti|}$).

2) Si~$U$ est une~$k$-couronne, si~$\omega$ est un bout de~$U$ convergeant vers~$x$, et si~$r$ désigne le rayon extérieur de~$U$ modulo~$|k\ti|$ en~$\omega$, on a l'inclusion~$|\hres(x)\ti|\subset |k\ti|.r^\ZZ$, avec égalité si~$k$ est algébriquement clos. Le point~$x$ est de type 2 (resp. 3) si et seulement si~$r\in \sqrt{|k\ti|}$ (resp.~$r\in \RR\ti_+- \sqrt{|k\ti|}$).}

\medskip
{\em Démonstration.} Si~$U$ est un disque, il existe une couronne~$Z$ contenue dans~$U$ et aboutissant à l'unique bout~$\omega$ de~$U$. Le point~$x$ appartient à~$\overline Z$ ; lorsqu'on voit~$\omega$  comme bout de~$Z$, il converge vers~$x$ et le rayon extérieur modulo~$|k\ti|$ de~$Z$ en~$\omega$ est égal au rayon de~$U$ modulo~$|k\ti|$. Il suffit par conséquent de démontrer 2). 

\medskip
On se place donc sous les hypothèse de 2). L'inclusion~$|\hres(x)\ti|\subset |k\ti|.r^\ZZ$, la description du type de~$x$ en fonction de~$r$, et la quasi-lissité de~$X$ en~$x$ découlent de~\ref{extremcour} {\em et sq.}. 

\medskip
Supposons maintenant que~$k$ soit algébriquement clos, et soit~$b\in \br X x \ctd U$ la branche définie par~$\omega$.  Il existe une section coronaire~$Z$ de~$b$ et une fonction analytique~$f$ définie et inversible sur un voisinage ouvert de~$x$ contenant~$Z$ et telle que~$f_{|Z}$ soit une fonction coordonnée de~$Z$ ; quitte à remplacer~$f$ par~$f^{-1}$, on peut supposer que~$|f|$ croît dans la direction de~$x$. On peut par ailleurs restreindre~$Z$ de sorte qu'elle soit une sous-couronne de~$U$ aboutissant à~$\omega$. Le rayon extérieur de~$U$ modulo~$|k\ti|$ en~$\omega$ est alors égal à celui de~$Z$, qui coïncide avec la limite de~$|f|$ dans la direction de~$\omega$. Or cette limite vaut précisément~$|f(x)|$. Par conséquent~$r=|f(x)|$ modulo~$|k\ti|$, ce qui entraîne que~$r\in |\hres(x)\ti|$ et établit l'inclusion~$|k\ti|.r^\ZZ\subset |\hres(x)\ti|$, achevant ainsi la démonstration.~$\Box$ 
\subsection*{Ablation et adjonction de sommets} 

Soit~$X$ une courbe~$k$-analytique quasi-lisse et soit~$S$ une triangulation de~$X$ ; on note~$\Gamma$ le squelette de~$S$. 

\deux{condabl} Soit~$x\in S$ ; posons~$\Sigma=S\setminus\{x\}$, notons~$\Gamma_0$ la composante connexe de~$\Gamma-\Sigma$ contenant~$x$. Les composantes connexes de~$X-\Sigma$ sont d'une part~$\Gamma_0^\flat$, d'autre part les composantes connexes de~$X\setminus S$ dont le bord ne contient pas~$x$. Comme toute composante connexe de~$\Gamma\setminus S$ est relativement compacte dans~$\Gamma$, le graphe~$\Gamma_0$ est relativement compact dans~$\Gamma$, et~$\Gamma_0^\flat$ est par conséquent un ouvert relativement compact de~$X$. Il s'ensuit que~$\Sigma$ est une triangulation de~$X$ si et seulement si~$\Gamma_0^\flat$ est un disque virtuel ou une couronne virtuelle. 

\medskip
Supposons que ce soit le cas. Cela implique que~$\Gamma_0$ possède un ou deux bouts. Comme~$S$ contient tous les n\oe uds de~$\Gamma$, il en contient tous les sommets, ce qui implique que~$\Gamma_0$ ne  possède aucun sommet à l'exception éventuelle de~$x$. Il y a donc maintenant deux possibilités. 

\trois{casgammazerodeux} {\em Le cas où~$\Gamma_0$ possède deux bouts.} C'est alors un intervalle ouvert, dont l'adhérence dans~$\Gamma$ est ou bien un cercle se refermant sur un élément de~$S$, ou bien un intervalle compact joignant deux éléments de~$S$. Notons~$\beta_+$ et~$\beta_-$ les deux branches de~$\Gamma_0$ issues de~$x$

\medskip
Si~$\Gamma_0^\flat$ est une couronne virtuelle alors~$\Gamma_0$ est son squelette ; le point~$x$ appartient à l'intérieur analytique de~$X$, est de genre 0, et vérifie les égalités~$$\got s(\beta_+)=\got s(x)=\got s(\beta_-)=\got s(\Gamma_0\an).$$

Réciproquement, supposons que~$x$ appartienne à l'intérieur analytique de~$X$, qu'il soit de genre~$0$, et que~$\got s(\beta_+)=\got s(x)=\got (\beta_-)$ ; comme de plus le point~$x$ possède un voisinage dans~$\Gamma$ qui est un intervalle ouvert, ce n'est pas un n\oe ud de~$\Gamma$ ; comme~$\Gamma_0\cap S=\{x\}$, l'intervalle~$\Gamma_0$ ne contient aucun n\oe ud de~$\Gamma$. Par ailleurs,~$\Gamma_0$ est d'adhérence compacte dans~$\Gamma$ ; le~\ref{pretriang} assure alors que~$\Gamma_0^\flat$ est une couronne virtuelle. 

\trois{casgammazeroun} {\em Le cas où~$\Gamma_0$ possède un bout.} C'est alors un intervalle semi-ouvert dont~$x$ est l'extrémité, et dont le bord dans~$\Gamma$ consiste en un point~$y$ de~$S$ qui diffère de~$x$. Comme~$\Gamma_0\cap S=\{x\}$, l'intervalle~$]y;x[$ ne contient aucun n\oe ud de~$\Gamma$ ; par conséquent,~$]y;x[^\flat$ est une couronne virtuelle (\ref{pretriang}). On note~$b$ (resp.~$\beta$) la branche issue de~$y$ (resp. issue de de~$x$) et définie par~$]y;x[$. 

\medskip
Supposons que~$\Gamma_0$ soit un disque virtuel ; le point~$x$ appartient alors à l'intérieur analytique de~$X$ et est de genre~$0$, et l'on a~$\got s(b)=\got s(]y;x[^\flat)=\got s(\Gamma_0^\flat)$ (\ref{remsxsv}). L'on dispose par ailleurs de deux plongements naturels~$\got s(\Gamma_0^\flat)\hookrightarrow \got s(x)$ et~$\got s(x)\hookrightarrow \got s(\beta)=\got s(]y;x[^\flat)$ ; comme~$\got s(\Gamma_0^\flat)=\got s(]y;x[^\flat)$, on a~$\got s(x)=\got s(\beta)$. 

\medskip
Réciproquement, supposons que~$x$ appartienne à l'intérieur analytique de~$x$, soit de genre~$0$, et que~$\got s(x)=\got s(\beta)$. Nous allons montrer que~$\Gamma_0^\flat$ est un disque virtuel. 

Les composantes connexes de~$\Gamma_0^\flat\setminus\{x\}$ sont d'une part~$]y;x[^\flat$, d'autre part des composantes connexes de~$X-\Gamma$, et donc des disques virtuels. Il s'ensuit que~$\got s(\Gamma_0^\flat)=\got s(x)$ (\ref{remsxgen}) ; notons~$F$ ce corps et choisissons-en un~$k$-plongement dans~$\KK$. L'espace~$\Gamma_0^\flat$ est une courbe~$F$-analytique ; elle est sans bord, car le bord analytique de~$X$ est contenu dans~$S\setminus\{x\}$ qui ne rencontre pas~$\Gamma_0^\flat$. 

La courbe~$\KK$-analytique~$\Gamma_0^\flat\hotimes_F\KK$ est connexe, et~$x$ a un et un seul antécédent~$x_{\KK}$ sur celle-ci. Les composantes connexes de~$\Gamma_0^\flat\hotimes_F\KK\setminus\{x_{\KK}\}$ sont d'une part des~$\KK$-disques, d'autre part les composantes connexes de~$]y;x[^\flat\hotimes_F\KK$ ; comme ~$F=\got s(x)=\got s(\beta)=\got s(]y;x[^\flat)$, l'espace~$]y;x[^\flat\hotimes_F\KK$ est en fait connexe, et est donc une~$\KK$-couronne. Prolongeons-la en un~$\KK$-disque ; on obtient une courbe~$\KK$-analytique quasi-lisse, connexe, compacte et sans bord ; elle s'identifie dès lors à~${\sch X}\an$ pour une certaine~$\KK$-courbe projective, intègre et lisse~$\sch X$. Le point~$x_{\KK}$ est un point de genre~$0$ de ~${\sch X}\an$, et~${\sch X}\an\setminus\{x_{\KK}\}$ est réunion disjointe de disques. Il résulte alors du lemme~\ref{lemrevg0d} que~${\sch X}\simeq \PP^1_{\KK}$ ; par construction,~$\Gamma_0^\flat\hotimes_F\KK$ s'identifie à une composante connexe de~${\sch X}\an\setminus\{\eta\}$ pour un~$\eta\in{\sch X}\an\dtr$ ; c'est donc un disque, et~$\Gamma_0^\flat$ est bien un disque virtuel. 

\trois{recapitabl} {\em Récapitulation.} Le sous-ensemble~$S\setminus\{x\}$ est donc une triangulation si et seulement si toutes les conditions suivantes sont vérifiées :

\medskip
$\bullet$~$x\notin \partial X\an$ et est de genre~$0$ ; 

$\bullet$~$\Gamma_0$ est un intervalle ouvert, ou bien un intervalle semi-ouvert d'extrémité~$x$ ; 

$\bullet$~$\got s(x)=\got s(\beta)$ pour toute branche~$\beta$ de~$\Gamma_0$ issue de~$x$. 

\medskip
Remarquons que la dernière condition est automatiquement vérifiée lorsque~$k$ est algébriquement clos. 

\deux{introadj} Soit maintenant~$x$ un point de~$X\dtr\setminus S$ ; nous allons montrer qu'il existe une plus petite triangulation de~$X$ contenant~$S\cup \{x\}$, et la décrire. 

\trois{casxsurg} {\em Le cas où~$x\in \Gamma$}. Il est alors situé sur une composante connexe~$I$ de~$\Gamma\setminus S$, qui est un intervalle ouvert d'adhérence compacte dans~$X$. Soient~$I_+$ et~$I_-$ les deux composantes connexes de~$I\setminus\{x\}$. Les composantes connexes de~$X-(S\cup\{x\})$ sont d'une part les composantes connexes de~$X\setminus S$ ne contenant pas~$x$, d'autre part les ouverts~$I_+^\flat$ et~$I_-^\flat$ de~$X$. Chacun des intervalles~$I_-$ et~$I_+$ est d'adhérence compacte dans~$\Gamma$ (puisque c'est le cas de~$I$) et ne contient aucun n\oe ud de~$\Gamma$ (puisque ceux-ci appartiennent tous à~$S$). Il en résulte que~$I_+^\flat$ et~$I_-^\flat$ sont deux couronnes virtuelles (\ref{pretriang}) relativement compactes ; par conséquent,~$S\cup\{x\}$ est une triangulation de~$X$.

\trois{casxpassurg} {\em Le cas où~$x\notin \Gamma$.} Soit~$y$ l'image de~$x$ sur~$\Gamma$ par la rétraction canonique ; la réunion~$\Gamma'$ de~$\Gamma$ et de~$[x;y]$ est admissible, donc analytiquement admissible (\ref{extanaladm}). Soit~$S'$ la réunion de~$S$ et de l'ensemble des nœuds de~$\Gamma'$. 

\medskip
{\em Le point~$y$ appartient à~$S'$}. En effet, c'est évident s'il appartient à~$S$ ; sinon, ce n'est pas un nœud, et {\em a fortiori} pas un sommet de~$\Gamma$ ; il y a donc deux branches de~$\Gamma$, et partant trois branches de~$\Gamma'$,  issues de~$y$, ce qui implique que~$y$ est un sommet, et {\em a fortiori} un nœud, de~$\Gamma'$ ; il appartient en conséquence à~$S'$. Notons par ailleurs qu'un point de~$\Gamma$ différent de~$y$ est un nœud de~$\Gamma'$ si et seulement si c'est un nœud de~$\Gamma$ ; par conséquent,~$S'\cap \Gamma=S\cup\{y\}$. 

\medskip
Le point~$x$ est un sommet, et donc un nœud, de~$\Gamma'$ ; il appartient de ce fait à~$S'$. Si~$z\in ]x;y[$ alors~$z\notin \Gamma$ ; par conséquent,~$z$ est situé sur un disque virtuel ; c'est donc un point de genre zéro appartenant à l'intérieur analytique de~$X$ ; de plus, ce n'est pas un sommet de~$\Gamma'$. Il s'ensuit que~$z$ est un nœud de~$\Gamma'$ si et seulement si il vérifie la condition 3) du~\ref{defnoeudsg}. 

L'application~$z\mapsto [\got s(z):k]$ étant semi-continue inférieurement (\ref{defnoeudsg}), l'ensemble de ses points de discontinuité sur~$]x;y[$ est fini. Il s'ensuit, compte-tenu de ce qui précède, qu'il existe une suite~$y=x_0<x_1<\ldots<x_n=x$ de points de~$[y;x]$ tels que~$S'=S\cup\{x_i\}_{0\leq i\leq n}$.

\medskip
L'ensemble~$S'$ est une partie fermée et discrète de~$\Gamma'$ qui en contient tous les nœuds ; si~$J$ est une composante connexe de~$\Gamma'\setminus S'$ alors~$J$ est ou bien un intervalle de la forme~$]x_i;x_{i+1}[~$ avec~$0\leq i\leq n-1$, ou bien est contenue dans une composante connexe de~$\Gamma\setminus S$ ; dans les deux cas,~$J$ est un intervalle ouvert d'adhérence compacte dans~$\Gamma'$. On en conclut que~$S'$ est une triangulation de~$X$,  contenant par construction~$S\cup\{x\}$. 

Si~$\Sigma$ est une triangulation de~$X$ contenant~$S\cup\{x\}$ alors le graphe de~$\Sigma$ contient~$\Gamma'$ ; par conséquent,~$\Sigma$ contient~$S'$. Ainsi,~$S'$ apparaît comme la plus petite triangulation de~$X$ contenant~$S\cup\{x\}$. 

\trois{remadjpt} On conserve les hypothèses et notations du~\ref{casxpassurg} ; faisons quelques remarques. 

\medskip
{\em Première remarque.} Soit~$j\in \{0,\ldots,n\}$. Il découle de la construction ci-dessus si~$j\geq 1$ et du~\ref{casxsurg} si~$j=0$ que~$S\cup\{x_i\}_{0\leq i\leq j}$ est la plus petite triangulation de~$X$ contenant~$S\cup\{x_j\}$. 

On en déduit aisément que si~$\Sigma$ est une triangulation quelconque contenant~$S$ telle que~$\Sigma\setminus S$ soit fini, il existe une suite finie~$S=S_0\subset S_1\subset\ldots\subset S_m=\Sigma$ de triangulations de~$X$ telles que~$S_{i+1}\setminus S_i$ soit un singleton pour tout entier~$i\in \{0,\ldots,m-1\}$. 

\medskip
{\em Seconde remarque.} Si~$k$ est algébriquement clos la fonction~$z\mapsto [\got s(z):k]$ est constante (et égale à~$1$) sur~$\Gamma'$ ; par conséquent, la suite~$(x_i)$ est réduite à deux termes, à savoir~$x_0=y$ et~$x_1=x$. Ainsi,~$S'$ coïncide dans ce cas avec~$S\cup\{y,x\}$. 

\subsection*{Description du squelette analytique dans le cas algébriquement clos}

\deux{propcarintdc} {\bf Proposition.} {\em Soit~$X$ une courbe~$\KK$-analytique 
génériquement quasi-lisse et soit~$U$ un ouvert relativement compact de~$X$ dont le bord analytique est vide, dont le bord topologique est contenu dans~$X\typ{23}$, et qui ne contient aucun point singulier. L'ouvert~$U$ est un disque (resp. une couronne) si et seulement si~$U$ est un arbre à un bout (resp. à deux bouts) ne contenant aucun point de genre strictement positif.}

\medskip
{\em Démonstration.} Dans les deux cas, les implications directes sont évidentes ; il reste à montrer les réciproques. On suppose donc que~$U$ est un arbre ayant un ou deux bouts, sans point de genre strictement positif et sans point singulier. Quitte à remplacer~$X$ par son lieu quasi-lisse (qui contient~$U$ par hypothèse, ainsi
que~$\partial U$ car ce dernier est contenu dans~$X\dtr$), on peut alors
supposer que~$X$ elle-même est quasi-lisse. 

\medskip
Soit~$S$ le bord
topologique
de~$U$ dans~$X$
(il comprend un ou deux éléments). Fixons une
triangulation~$T$ de~$X$ qui contient~$S$, et soit~$\Gamma$ son squelette. 

\medskip
Supposons que~$T\cap U$ soit non vide ; dans ce cas~$\Gamma\cap U$ est non vide. Si~$x$ est un point unibranche de~$\Gamma\cap U$ alors~$x\in T$. En vertu des hypothèses faites sur~$U$ le point~$x$ est de genre zéro, et n'appartient pas à~$\partial \an X$. Il résulte alors de~\ref{recapitabl} que~$T\setminus\{x\}$ est encore une triangulation de~$X$, et elle contient~$S$ par construction ; en recommençant l'opération autant de fois  qu'il est nécessaire, on se ramène au cas où~$\Gamma\cap U$ ne contient plus aucun point unibranche. On distingue maintenant deux cas. 

\trois{carintdisc} {\em Supposons que~$\Gamma\cap U=\emptyset$.} Dans ce cas~$U$ est une composante connexe de~$X-\Gamma$ ; comme~$\Gamma$ est analytiquement admissible,~$U$ est disque. 

\trois{carintcour} {\em Supposons que~$\Gamma\cap U\neq \emptyset$.} Comme~$\Gamma$ contient~$S$ et rencontre~$U$, l'adhérence de~$\Gamma\cap U$ dans la compactification arboricole~$\wid U$ (adhérence qui s'identifie à~$\wid {(\Gamma\cap U)}$) n'est pas réduite à un singleton ; comme~$\wid {(\Gamma\cap U)}$ est par ailleurs un arbre compact et fini, il contient au moins deux points unibranches, dont aucun ne peut par hypothèse appartenir à~$U$. Comme~$\got d U$ possède au plus deux éléments, on déduit de ce qui précède qu'il en possède exactement deux, qui sont les points unibranches de~$\wid {(\Gamma\cap U)}$ ; autrement dit,~$U$ est un arbre à deux bouts et~$\Gamma\cap U$ est l'intervalle ouvert qui les joint, c'est-à-dire encore~$\skel U$. 

\medskip
L'intersection~$T\cap U$ est un ensemble fini. Si~$x\in T\cap U$, les hypothèses faites sur~$U$ assurent que~$x$ est de genre zéro, et n'appartient pas à~$\partial \an X$ ;  il résulte alors de~\ref{recapitabl} que~$T\setminus\{x\}$ est encore une triangulation de~$X$, et elle contient~$S$ par construction. En recommençant l'opération autant de fois qu'il est nécessaire, on se ramène au cas où~$T\cap U=\emptyset$. L'ouvert~$U$ apparaît alors comme une composante connexe de~$X\setminus T$ ; par définition d'une triangulation, et compte-tenu du fait que~$U$ est un arbre à deux bouts, il s'ensuit que~$U$ est une couronne.~$\Box$

\deux{skelanalskelgpos} {\bf Théorème.} {\em Soit~$X$ une courbe~$\KK$-analytique
génériquement quasi-lisse connexe et non vide et soit~$\Theta$ le sous-ensemble fermé et discret de~$X$ constitué des points~$x$ possédant l'une au moins des trois propriétés suivantes : 

$\bullet$ $x$ est singulier ;  

$\bullet$ $g(x)>0$ ; 

$\bullet$ $x\in \partial\an X$.

\medskip
1) Si~$\skel X$ est non vide et si~$r$ désigne la rétraction canonique de~$X$ sur~$\skel X$, alors~$\skelan X$ est la réunion de~$\skel X$ et des segments~$[x;r(x)]$ pour~$x$ parcourant~$ \Theta$.

2) Si~$X$ est un arbre à un bout~$\omega$ alors~$\skelan X$ est égal à la réunion des~$[x;\omega[$ pour~$x$ parcourant~$\Theta$.

3) Si~$X$ est un arbre compact,~$\skelan X$ est la réunion des segments~$[x;y]$ pour~$x$ et~$y$ parcourant~$\Theta$. }

\medskip
{\em Démonstration.} Dans chacune des situations 1), 2), et 3), on souhaite montrer l'égalité entre~$\skelan X$ et un certain sous-ensemble~$\sch E$ de~$X$ (dont la définition dépend du cas considéré), qui possède les propriétés suivantes : 

\medskip
$\bullet$ il est convexe (cela découle des définitions) ; 

$\bullet$ il est tracé sur~$X\geom$, en vertu de la convexité de ce dernier, du fait que~$\Theta\subset X\geom$, et du fait que~$\skel X\subset X\geom$ puisque tout point
non rigide de type~$1$ ou~$4$ de~$X$ est unibranche.

\medskip
De plus, l'ensemble~$\sch E$ est un sous-graphe fermé de~$X$. En effet, cela résulte dans le cas 1) de~\ref{convtestcomp}, dans le cas 2) de la proposition~\ref{propuniondrl}, et dans le cas 3) du fait que~$\sch E$ est alors une réunion finie de segments. Ce sous-graphe est admissible dès qu'il est non vide : dans le cas 1), c'est parce que~$\sch E$ est un sous-graphe convexe et fermé de~$X$ contenant son sous-graphe admissible~$\skel X$ ; dans les cas 2) et 3), c'est une conséquence de~\ref{arbradmarbr}.

\trois{envconvdiscferm} {\em Le squelette analytique~$\skelan X$ contient~$\sch E$.} C'est dû au fait qu'il contient~$\skel X\cup \Theta$, qu'il est convexe, et est admissible s'il est non vide (ce dernier point ne sert que dans le cas 2) -- il garantit alors que si~$\skelan X$ est non vide il aboutit à~$\omega$ d'après~\ref{arbradmarbr}). 

\trois{recenconvdiscferme} {\em Le sous-graphe convexe~$\sch E$ de~$X$ contient~$\skelan X$.} Soit~$x\in X\setminus \sch E$ ; nous allons montrer que~$x$ possède un voisinage ouvert dans~$X$ qui est un disque, ce qui assurera que~$x\notin \skelan X$ et achèvera la démonstration.

\medskip
Supposons tout d'abord que~$\sch E$ est non vide ; il  est dès lors admissible, et la composante connexe~$U$ de~$X\setminus \sch E$ est en conséquence un arbre à un bout relativement compact ; soit~$t$ l'unique point de son bord. Comme~$\sch E\supset \Theta$, l'ouvert~$U$ ne contient aucun point de~$\partial\an X$ ni aucun point de genre strictement positif. Soit~$b\in \br X t$ la branche définie par~$[x;t]$. Il existe~$s\in ]x;t[$ tel que~$]s,t[^\flat$ soit une section
coronaire de~$b$ ; soit~$s'\in ]s;t[$ et soit~$U'$ la composante connexe de~$X\setminus\{s'\}$ contenant~$s$. Par construction, $U'$ est un 
arbre à un bout relativement compact de bord~$\{s'\}$,
qui ne contient aucun point de~$\partial\an X$ ni aucun point de genre strictement positif. Comme~$s'$ est situé sur le squelette d'une couronne, il est de type 2 ou 3, et la proposition~\ref{propcarintdc} ci-dessus assure que~$U'$ est un disque (si~$t\in X\dtr$ on peut montrer directement par un argument analogue que~$U$ est un disque).

\medskip
Supposons maintenant que~$\sch E$ est vide. La courbe~$X$ est alors un arbre ou bien compact, ou bien à un bout. Si elle est compacte, elle ne peut être réduite à~$\{x\}$ (sinon~$x$ serait de type 3 et appartiendrait à~$\partial \an X$, et donc à~$\sch E$) et possède donc au moins un point~$y$ différent de~$x$ et non unibranche, donc de type 2 ou 3. 

Si~$X$ est un arbre à un bout~$\omega$, on choisit un point~$y$ sur~$]x;\omega[$ ; comme~$y$ est non unibranche, il est de type 2 ou 3. 

\medskip
Dans les deux cas, la composante connexe~$U$ de~$X\setminus\{y\}$ contenant~$x$ est un arbre à un bout relativement compact. Comme~$\sch E$ est vide,~$U$ ne contient aucun point appartenant à~$\partial\an X$, aucun point de genre strictement positif , ni aucun potin singulier ; la proposition~\ref{propcarintdc} ci-dessus assure alors que~$U$ est un disque.~$\Box$ 

\deux{corollunivalskelan} {\bf Corollaire.} {\em Soit~$X$ une courbe~$\KK$-analytique
génériquement quasi-lisse. Si~$x$ est un point unibranche de~$\skelan X$ alors ou bien~$x$ est de genre strictement positif, ou bien il
est singulier, 
ou bien il appartient à~$\partial\an X$.}

\medskip
{\em Démonstration.} Cela résulte immédiatement du
théorème~\ref{skelanalskelgpos} ci-dessus, et du fait que~$\skel X$ n'a pas de points unibranches.~$\Box$

\section{Triangulations et cohomologie étale}

\subsection*{Une suite exacte}

\deux{rapprpi}
Soit~$X$ une courbe~$k$-analytique, 
et soit~$\pi$ le morphisme du site
étale
~$X\et$ vers le site
topologique $X\top$. 

\trois{defpi-explicite}
Pour tout~$q$, le foncteur~${\rm R}^q\pi_*$ admet la description explicite suivante : 
il associe à un faisceau~$\sch F$ sur~$X\et$
le faisceau sur~$X\top$
associé au préfaisceau~$U\mapsto \H^q(U\et,\sch F)$.

\trois{suite-spectrale}
Soit~$\sch F$ un faisceau
sur~$X\et$. La dimension
topologique de~$X$ étant égale à~$1$
(prop.~\ref{propdimtop}), la suite spectrale
$$\H^p(X\top, {\rm R}^q\pi_*\sch F)\Rightarrow \H^{p+q}(X\et, \sch F)$$
induit pour tout entier~$p$ une suite exacte
$$0\to \H^{1}(X\top, {\rm R}^{p-1}\pi_*\sch F)\to \H^p(X\et,\sch F)\to \H^0(X\top, {\rm R}^p\pi_*\sch F)\to 0.$$

\trois{s-s-caspart1}
En prenant~$p=1$, on obtient la suite exacte
$$0\to \H^1(X\top, \pi_*\sch F)\to \H^1(X\et, \sch F)\to \H^0(X\top, {\rm R}^1\pi_*\sch F)\to 0,$$
qui est la seule que nous aurons à considérer ici. Nous
allons maintenant expliquer comment l'on peut,
grâce à l'interprétation explicite des trois groupes en jeu
et des flèches entre iceux, en démontrer l'existence directement,
sans référence à la machinerie cohomologique générale.

\trois{interp-gp} 
Un élément de~$\H^1(X\et, \sch F)$ peut s'interpréter
comme une classe d'isomorphie de~$\sch F$-torseurs
(sur~$X\et$).

\medskip
Un élément de~$\H^0(X\top, {\rm R}^1\pi_*\sch F)$
peut s'interpréter comme une classe de familles~$(\sch T_x)_{x\in X}$,
où~$\sch T_x$ est pour tout~$x$ un germe de~$\sch F$-torseur en~$x$,
et où il existe pour tout~$x\in X$ un voisinage ouvert~$U$ de~$x$ et un~$\sch F$-torseur~$\sch S$
sur~$U$ tel que~$\sch T_y=\sch S\times_U(X,y)$ pour tout~$y\in U$ ; 
le terme «classe»
fait référence à la relation suivante :  
on identifie~$(\sch T_x)_x$ et~$(\sch S_x)_x$
si les germes~$\sch T_x$ et~$\sch S_x$ sont isomorphes pour tout~$x$. 

\medskip
Un élément de~$ \H^1(X\top, \pi_*\sch F)$
peut s'interpréter
comme une classe d'isomorphie de~$\pi_*\sch F$-torseurs (sur~$X\top$) ;
 notons que~$\pi_*\sch F$ n'est autre que la restriction de~$\sch F$
à~$X\top$. 

\trois{interp-fl-1}  Soit~$\sch T$ un~$\pi_*\sch F$-torseur. Le produit contracté 
$\sch T':=\pi^*\sch T\times^{\pi^*\pi_*\sch F}\sch F$ est un~$\sch F$-torseur ; cette construction
induit un morphisme
de~$\H^1(X\top, \pi_*\sch F)$
dans~$\H^1(X\et, \sch F)$
qui est celui de~\ref{s-s-caspart1}.

Soit~$(U_i)$ un recouvrement ouvert de~$X$, muni d'une famille
$$(\sigma_i : \sch T_{|U_i}\simeq \pi_*\sch F_{|U_i})_i$$
de trivialisations. À cette donnée est associée un cocycle~$(h_{ij})$ à coefficients dans~$\pi_*\sch F$. La famille
$(\sigma_i : \sch T_{|U_i}\simeq \pi_*\sch F_{|U_i})$ induit une famille de trivialisations~$(\sigma'_i : \sch T'_{|U_i}\simeq \sch F_{|U_i})_i$,
et le cocycle correspondant est «égal»
à~$(h_{ij})$ (rappelons que~$\pi_*\sch F(U_i\cap U_j)=\sch F(U_i\cap U_j)$ pour tout~$(i,j)$). On en déduit que
pour tout ouvert~$U$ de~$X$, la flèche naturelle~$\sch T(U)\to \sch T'(U)$ est un isomorphisme (les deux ensembles admettent
en effet
la même description au moyen du cocycle~$(h_{ij})$ et des groupes~$\pi_*\sch F(U_i\cap U)=\sch F(U_i\cap U)$). 
On a en particulier~$\sch T(X)=\sch T'(X)$ ; par conséquent,~$\sch T(X)\neq \varnothing\iff \sch T'(X)\neq \varnothing$,
ce qui signifie que~$\H^1(X\top, \pi_*\sch F)\to\H^1(X\et, \sch F)$
est injective. 

\medskip
\trois{interp-fl-2}
Soit~$\sch S$ un~$\sch F$-torseur. La famille~$(\sch S\times_X(X,x))_x$ satisfait les conditions 
énoncées au~\ref{interp-gp} ; cette construction induit un morphisme
de~$\H^1(X\et, \sch F)$
dans~$\H^0(X\top, {\rm R}^1\pi_*\sch F)$,
qui est celui de~\ref{s-s-caspart1}. 

\medskip
Soit~$\sch T$ un~$\pi_*\sch F$-torseur ; on a vu au~\ref{interp-fl-1}
ci-dessus que
le~$\sch F$-torseur~$\sch T':=\pi^*\sch T\times^{\pi^*\pi_*\sch F}\sch F$ est 
trivialisé par un recouvrement ouvert de~$X$ (il suffit de prendre un recouvrement
ouvert trivialisant~$\sch T$) ; en conséquence, le germe~$\sch T'\times_X(X,x)$ est trivial
pour tout~$x\in X$. 

\medskip
Réciproquement, soit~$\sch S$ un~$\sch F$-torseur tel que~$\sch S\times_X(X,x)$
soit trivial pour tout~$x\in X$. Il existe alors un recouvrement ouvert~$(U_i)$ de~$X$, 
et pour tout~$i$ une trivialisation~$\sigma_i : \sch T_{|U_i}\simeq \sch F_{|U_i}$. Celle-ci
induit un cocycle~$(h_{ij})$ subordonné à~$(U_i)$ et à coefficients dans~$\sch F$, 
que l'on peut voir comme à coefficients dans~$\pi_*\sch F$ puisque~$\pi_*\sch F(U_i\cap U_j)=\sch F(U_i\cap U_j)$
pour tout~$(i,j)$.
Ce cocycle~$(h_{ij})$ définit donc un~$\pi_*\sch F$-torseur~$\sch T$, et
l'on a~$\sch S\simeq \pi^*\sch T\times^{\pi^*\pi_*\sch F}\sch F$ ({\em cf.}
\ref{interp-fl-1}). 

\medskip
En conséquence, la
suite de~\ref{s-s-caspart1}
est exacte en~$\H^1(X\et, \sch F)$. 

\trois{interp-fl3}
Soit~$(\sch T_x)_x$ une famille de germes de~$\sch F$-torseurs, 
comme au~\ref{interp-gp}. Choisissons un recouvrement ouvert~$(U_i)_{i\in I}$
de~$X$ et pour tout~$i$ un~$\sch F$-torseur~$\sch T_i$ sur~$U_i$ tel que
$\sch T_x\simeq \sch T_i\times_X(X,x)$ pour tout~$x\in U_i$. 

\medskip
Comme~$X$ est paracompact
de dimension topologique~$\leq 1$, on peut quitte à
raffiner~$(U_i)$ 
supposer qu'il est localement fini,
puis que
les intersections trois à trois des~$U_i$ sont vides. Par ailleurs, étant 
paracompact, l'espace~$X$ est normal ; il existe donc un recouvrement ouvert~$(V_i)$
de~$X$ tel que~$\overline{V_i}\subset U_i$ pour tout~$i$. On choisit un ordre total 
sur l'ensemble~$I$ des indices. 

\medskip
Soit~$x\in X$. Il y a au moins un, et au plus deux, éléments~$i\in I$
tels que~$x\in U_i$. Choisissons un voisinage
ouvert~$W_x$ de~$x$ dans~$X$ ayant les propriétés suivantes. 

\medskip
$\bullet$ On a~$W_x\subset V_i$ pour tout~$i$ tel que~$x\in V_i$, et~$W_x\subset U_i$
pour tout~$i$ tel que~$x\in U_i$. 

$\bullet$ Pour tout~$i$ tel que~$x\notin U_i$ l'intersection de~$W_x$ et~$\overline V_i$
est vide. 

$\bullet$ Si~$x\in U_i\cap U_j$ avec~$i< j$, il existe un isomorphisme~$\sigma_x : \sch T_{i|W_x}\simeq \sch T_{j|W_x}$. 

\medskip
Pour tout~$x\in X$ on choisit un indice~$i(x)$ tel que~$x\in V_{i(x)}$. Le but de ce qui suit est de montrer
l'existence d'un système compatible~$(\iota_{xy})$
d'isomorphismes autorisant le recollement des torseurs~$\sch T_{i(x)|W_x}$ lorsque~$x$
parcourt~$X$. 

Soient~$x$ et~$y$ deux points
de~$X$ tels que~$W_x\cap W_y\neq \emptyset$.  On distingue deux cas.

\medskip
$\bullet$ Si~$i(x)=i(y)$ on pose~$\iota_{xy}=\mathsf{Id}_{\sch T_{i(x)|W_x\cap W_y}}$. 

$\bullet$ Supposons maintenant~$i(x)\neq i(y)$. Comme~$W_y$ rencontre~$W_x\subset V_{i(x)}$, 
on a~$y\in U_{i(x)}$, et
de même~$x\in U_{i(y)}$. En conséquence, 
les points~$x$
et~$y$
appartiennent tous deux à~$U_{i(x)}\cap U_{i(y)}$. On pose alors
$\iota_{xy}=\sigma_{x_|W_x\cap W_y}$ si~$i(x)<i(y)$, et~$\iota_{xy}=\sigma^{-1}_{y|W_x\cap W_y}$
sinon. 

\medskip
Il reste à s'assurer que le système~$(\iota_{xy})$ ainsi construit satisfait les relations dites de cocycle. 
Soient donc~$x, y$ et~$z$ trois points
de~$X$ tels que~$W_x\cap W_y\cap W_z\neq \emptyset$. Par
ce qui précède, $x\in U_{i(x)}\cap U_{i(y)}\cap U_{i(z)}$ ; il s'ensuit que deux au moins des trois indices~$i(x), i(y)$
et~$i(z)$ sont égaux. La relation souhaitée découle alors immédiatement de la définition de~$\iota_{xy}, \iota_{yz}$
et~$\iota_{xz}$. 

\medskip
Soit~$\sch S$ le~$\sch F$-torseur
obtenu par recollement des~$\sch T_{i(x)|W(x)}$ le long du système
$(\iota_{xy})$. On a par construction~$\sch S\times_X(X,x)=\sch T_{i(x)}\times_{U_{i(x)}}(X,x)\simeq \sch T_x$ pour tout~$x\in X$. 

\medskip
On en déduit que~$\H^1(X\et, \sch F)\to \H^0(X\top, {\rm R}^1\pi_*\sch F)$
est surjective, ce qui achève la preuve «concrète»
de l'exactitude de la suite de~\ref{s-s-caspart1}. 

\subsection*{Calcul de~$\H^1(X\et,\mu_\ell)$
{\em via}
une triangulation
lorsque le corps de base est algébriquement clos}

{\em On fixe un entier~$\ell$ premier à~$p$}. 

\deux{r1pi-simple-con} Soit~$X$ une courbe~$\KK$-analytique. On a d'après~\ref{s-s-caspart1}
une suite exacte naturelle
$$0\to \H^1(X\top, \mu_\ell)\to \H^1(X\et,\mu_\ell)\to \H^0(X\top, {\rm R}^1\pi_*\mu_\ell)\to 0.$$

On en déduit que si~$\H^1(X\top, \mu_\ell)=0$ (ce qui équivaut à demander que
l'entier~$\ell$
soit égal à~$1$ ou que
toute composante connexe de~$X$ soit un arbre), 
la flèche naturelle de~$\H^1(X\et,\mu_\ell)$
vers~$\H^0(X\top, {\rm R}^1\pi_*\mu_\ell)$
est un isomorphisme. {\em C'est notamment le cas lorsque~$X$ 
est un disque ou une couronne}. 

\deux{pre-enonce-desch1}
Soit~$X$ une courbe~$\KK$-analytique 
{\em lisse}, soit~$S$ une triangulation de~$x$
et soit~$\Gamma$ le squelette de~$S$. 

Pour tout~$x\in S$, on note~$g(x)$
le genre de~$x$, et l'on 
définit le groupe~$\JJ(x)$ comme suit : 

\medskip
$\bullet$ si~$x$ est de type~$3$ alors~$\JJ (x)$ est trivial ; 

$\bullet$ si~$x$ est de type~$2$ alors~$\JJ(x)=_\ell\sch
J(\kk_1)$, où~$\sch J$
est la jacobienne de la courbe résiduelle en~$x$. 

\medskip
Dans tous les cas, $\JJ(x)$ est isomorphe
(non canoniquement en général)
à~$(\ZZ/\ell \ZZ)^{2g(x)}$. 

\deux{theo-desc-h1-triang}
{\bf Théorème.}
{\em On conserve les notations du~\ref{pre-enonce-desch1}
ci-dessus. On dispose d'une suite exacte naturelle

$$0\to \prod_{x \in S} \JJ(x)\to \H^0(X\top,{\rm R}^1\pi_*\mu_\ell)\to \mathsf{Harm}(\Gamma,\ZZ/\ell\ZZ)\to 0.$$}

\medskip
{\em Démonstration.}
Avant d'entamer la preuve proprement dite, nous
allons fixer quelques notations.

\trois{rem-r1-cour-disc}
Soit $U$ un ouvert de~$X$ qui est un arbre, et soit~$h$
une section de~${\rm R}^1\pi_*\mu_\ell)$ définie sur un ouvert
de~$X$ contenant~$U$. Il résulte de~\ref{r1pi-simple-con}
que~$h_{|U}$ peut être vue comme une
«vraie»
classe de cohomologie appartenant à~$\H^1(U\et,\mu_\ell)$. 
Considérons maintenant deux cas particuliers importants.  

\medskip
$\bullet$ Si~$U$ est un disque, $\H^1(U\et,\mu_\ell)=0$
(prop.~\ref{h1ettriv}
et remarque
\ref{commdisc}) ; 
en conséquence, $h_{|U}=0$. 

$\bullet$ Si~$U$ est une couronne, $\H^1(U\et,\mu_\ell)=
\kum(X)$ (th.~\ref{corollkum}). Fixons une orientation~$o$
sur~$U$ ; ce choix induit un isomorphisme
$$\sigma_{(U,o)}:\H^1(U\et,\mu_\ell)=\kum(X)\simeq \ZZ/\ell\ZZ.$$ On
écrira par abus~$\sigma_{(U,o)}(h)$
au lieu de~$\sigma_{(U,o)}(h_{|U})$.

\trois{rem-theta-hx}
Soit~$x\in X\dtr$ et soit~$a\in \br X x$. La
fibre
de~${\rm R}^1\pi_*\mu_\ell$
en~$x$ s'identifie à~$\H^1((X,x)\et,\mu_\ell)$. 

Si~$h$ est une section
de~${\rm R}^1\pi_*\mu_\ell$
de germe~$h_x$ en~$x$, 
l'élément~$\theta_a(h_x)\in \ZZ/\ell\ZZ$
qui a été défini au~\ref{thetaakum}
sera simplement noté~$\theta_a(h)$. 

Si~$Z$ est une section
coronaire de~$a$ contenue
dans l'ouvert de définition de~$h$, 
et si~$o$ désigne l'orientation 
de~$Z$
vers~$x$ on déduit de la définition de
l'invariant~$\theta_a$,
et du fait que~$\kum (Z)\simeq \kum(Z')$ pour toute sous-couronne
$Z'$ de~$Z$, que~$\theta_a(h)=\sigma_{(Z,o)}(h)$. 

\trois{def-cl-gratte}
Soit~$x\in X\dtr$ et soit~$U$ un voisinage ouvert de~$x$. 
On note~$G(x)$
le sous-groupe de~$\H^1((X,x)\et,\mu_\ell)$
formé des germes 
de classes 
qui admettent un représentant supporté par~$\{x\}$. 

Notons que si~$\alpha \in G(x)$ et si~$U$ est un voisinage
ouvert de~$x$ alors~$\alpha$ admet un représentant {\em sur~$U$}
qui est supporté par~$\{x\}$ : il suffit de choisir un représentant de~$\alpha$
sur un voisinage ouvert~$U_0$ de~$\{x\}$ qui est supporté par~$\{x\}$, et de le
prolonger par zéro. 

\trois{def-ux-harm}
Soit~$x\in S$. Notons~$\Pi_x$ l'ensemble
des composantes connexes de~$X\setminus \Gamma$
dont le bord est égal à~$\{x\}$. 
Pour toute branche~$a\in \br X x\ctd \Gamma$
on choisit un intervalle ouvert~$I_a$ contenu dans~$\Gamma$, 
aboutissant proprement à~$x$ et définissant~$a$ ; quitte à restreindre les~$I_a$, 
on peut les supposer deux à deux disjoints. La réunion

$$U_x:=\{x\}\coprod \left(\coprod_{V\in \Pi_x}V\right)\coprod \left(\coprod_{a\in \br X x 
\ctd \Gamma} I_a^\flat\right)$$ est un voisinage ouvert connexe de~$x$.

\trois{h0r1-vers-harm}
{\em Construction d'un
morphisme
de~$\H^0(X\top,{\rm R}^1\pi_*\mu_\ell)$
vers~$\mathsf{Harm}(\Gamma,\ZZ/\ell\ZZ)$}. Soit~$h\in \H^1(X\top, {\rm R}^1\pi_*\mu_\ell)$, 
et soit~$(I,o)$ un intervalle ouvert non vide
et orienté contenu dans~$\Gamma\setminus S$. 
Notons~$c_h(I,o)$ l'élément
$\sigma_{(I^\flat, o)}(h)$
de~$\ZZ/\ell\ZZ$. Par construction, $c_h$ est
une~$\ZZ/\ell \ZZ$-cochaîne sur~$(\Gamma,S)$

\medskip
{\em Montrons que~$c$ est harmonique.}
Soit~$x\in S$ et soit~$a\in \br X x$
une branche qui n'est pas contenue dans~$\Gamma$. Elle
est dans ce cas contenue dans une composante connexe~$V$
de~$X\setminus \Gamma$, 
qui est un disque. Choisissons une section coronaire~$Z$ de~$a$ contenue
dans~$V$. On a~$h_{|V}=0$ (\ref{rem-r1-cour-disc}), et
{\em a fortiori}
$h_{|Z}=0$. En conséquence, $\theta_a(h)=0$ (\ref{rem-theta-hx}). 

\medskip
Comme~$x\notin \partial X$
(la courbe~$X$ est lisse, et 
en particulier sans bord), 
il résulte de ce qui précède et de la 
proposition~\ref{propcondresidutors},
et plus précisément de l'implication~ii)$\Rightarrow$iii)
qui figure dans son énoncé, que
$$\sum_{a\in \br X x} \theta_a(h)=0.$$ Mais puisqu'on a
vu que
$\theta_a(h)=0$ dès que~$a$ n'est pas contenue dans~$\Gamma$,
il vient~$\sum_{a\in \br X x\ctd \Gamma} \theta_a(h)=0.$
On déduit alors
de la définition de~$c_h$ et de~\ref{rem-theta-hx}
que~$c_h$ est harmonique.

\trois{h01-vers-harm-surj}
{\em Surjectivité de la flèche~$\H^0(X\top, {\rm R}^1\pi_*\mu_\ell)\to \mathsf{Harm}(\Gamma,\ZZ/\ell\ZZ)$.}
Soit~$c$ une cochaîne harmonique sur~$\Gamma$
à coefficients dans~$\ZZ/\ell \ZZ$. 

\medskip
Soit~$I$ un intervalle ouvert non vide contenu
dans~$\Gamma\setminus S$. 
La classe~$\sigma_{(I^\flat, o)}^{-1}(c(I,o))$ ne dépend pas de
l'orientation~$o$ ; 
notons-la~$h_I$. Si~$J$ est un 
intervalle ouvert
non vide contenu dans~$I$ alors~$h_{I|J}=h_J$. 

\medskip
Comme la cochaîne~$c$ est harmonique, 
il résulte de la proposition~\ref{propcondresidutors},
et plus précisément de l'implication~iii)$\Rightarrow$i)
de son énoncé, qu'il existe une classe~$h^x$
appartenant à~$\H^1(U_{x,{\rm \acute{e}t}},\mu_\ell)$
ayant les propriétés suivantes : 

\medskip
$\bullet$
pour tout~$V\in \Pi$, on a~$h^x_{|V}=0$ ; 

$\bullet$ pour tout~$a\in \br X x\ctd \Gamma$, la restriction
de~$h^x$ à~$I_a^\flat$ est égale
à~$h_{I_a}$. 

\medskip
Par construction, les classes~$h^x$ (pour~$x$ parcourant~$S$)
et~$h_{I^\flat}$ (pour~$I$ parcourant~$\pi_0(\Gamma\setminus S)$) 
se recollent en une section
$h$ de~${\rm R}^1\pi_*\mu_\ell$ sur~$X$, 
dont l'image dans~$\mathsf{Harm}(\Gamma,\ZZ/\ell \ZZ)$ est égale à~$c$. 

\trois{h01-harm-noyau}
{\em Noyau de la flèche~$\H^0(X\top,{\rm R}^1\pi_*\mu_\ell)\to \mathsf{Harm}(\Gamma,\ZZ/\ell\ZZ)$.}
Soit~$h$
une classe appartenant à~$\H^0(X\top, {\rm R}^1\pi_*\mu_\ell)$
telle que~$c_h=0$ et soit~$V$ une composante connexe de~$\Gamma\setminus S$. 
On a~$h_{|V}=0$ : si~$V$
est un disque cela résulte de~\ref{rem-r1-cour-disc}, et si~$V$ est une couronne 
c'est une conséquence de la nullité de~$c_h$. Il s'ensuit : 

\medskip
$\bullet$ que la classe~$h$ est uniquement déterminée
par la famille de germes~$(h_x)_{x\in S}$ ; 

$\bullet$ que chacun des~$h_x$ appartient au groupe~$G(x)$ des classes gratte-ciel. 

\medskip
On a ainsi construit une injection
$${\rm Ker}\;[\H^0(X\top,{\rm R}^1\pi_*\mu_\ell)\to \mathsf{Harm}(\Gamma,\ZZ/\ell\ZZ)]
\hookrightarrow \prod_{x\in S} G(x).$$

Montrons qu'elle est surjective. Soit~$(\alpha_x)\in \prod_{x\in S}G(x)$.
Il existe pour tout~$x\in S$
une section
$\bar h^x$ de~${\rm R}^1\pi_*\mu_\ell$ sur l'ouvert~$(X\setminus S)\cup\{x\}$
de~$X$ de support contenu dans~$\{x\}$. Les classes~$\bar h^x$ pour~$x$ parcourant~$S$
se recollent en une section globale~$h\in \H^0(X\top,{\rm R}^1\pi_*\mu_\ell)$. On a par
construction d'une part ~$h_{|(X\setminus S)}=0$, ce qui implique immédiatement que~$c_h=0$ ; 
et d'autre part~$h_x=\alpha_x$ pour tout~$x$. Ainsi,

$${\rm Ker}\;[\H^0(X\top,{\rm R}^1\pi_*\mu_\ell\to \mathsf{Harm}(\Gamma,\ZZ/\ell\ZZ)]
\simeq \prod_{x\in S}G(x).$$

Soit~$x\in S$ et soit~$h\in \H^1((X,x)\et,\mu_\ell)$. Il résulte de~\ref{prepacondresidu}
que~$h\in G(x)$ si et seulement si~$\theta_a(h)=0$ pour tout~$a\in \br X x$. 
On en déduit, 
en vertu de l'assertion 2) de la proposition~\ref{propcondresidutors},
que le groupe~$G(x)$ est canoniquement isomorphe à~$\JJ(x)$, ce qui achève la démonstration.~$\Box$

\deux{desc-h1-etale}
On en déduit, compte-tenu de la suite exacte de~\ref{r1pi-simple-con}
et du fait que
l'inclusion $\Gamma\hookrightarrow X$ est une équivalence homotopique,
l'existence d'un diagramme commutatif
$$\diagram &&0\dto&&\\
& & \H^1(X\top,\mu_\ell)\dto\rto^\simeq&\H^1(\Gamma,\mu_\ell)&\\
&&\H^1(X\et,\mu_\ell)\dto&&\\
0\rto&\prod_{x\in S}\JJ(x)\rto&\H^0(X\top,{\rm R}^1\pi_*\mu_\ell)\dto
\rto&\mathsf{Harm}(\Gamma,\ZZ/\ell\ZZ)\rto &0\\
&&0&&\enddiagram$$

 dont la ligne du bas et la colonne sont exactes.
 
\deux{calcul-de-genre}
Supposons maintenant que~$X$ est non vide, propre et connexe. Elle 
s'identifie alors à l'analytifiée~$\sch X\an$ d'une~$\KK$-courbe
algébrique projective, irréductible et lisse~$\sch X$ ; soit~$g$ son genre. 

\medskip
Soit~$\beta$ le premier nombre de Betti du graphe~$\Gamma$, qui est aussi
celui de l'espace topologique~$X$. On déduit
du~\ref{desc-h1-etale}
ci-dessus que
le groupe de~$\ell$-torsion~$\H^1(X\et,\mu_\ell)$ est
(non canoniquement en général) isomorphe
à~$(\ZZ/\ell\ZZ)^{2\beta+2\sum_{x\in S}g(x)}$. 
On sait par ailleurs en vertu de GAGA qu'il est isomorphe à~$(\ZZ/\ell\ZZ)^{2g}$. 

\medskip
Ceci vaut pour tout entier~$\ell$
premier à~$p$ (et donc en particulier pour au moins
un entier~$\ell\neq 1$). Compte-tenu du fait que tout point
de~$X\dtr\setminus S$ est de genre~$0$, il vient

$$g=\beta+\sum_{x\in S} g(x)=\beta +\sum_{x\in X\typ 2} g(x).$$

\section{Extension des scalaires}

\subsection*{Premières propriétés}

\deux{intro-extscal-courbes}
Soit~$X$ une courbe~$k$-analytique, soit~$S$
un sous-ensemble de~$X\geom$ et soit~$L$ une extension
complète de~$k$. Nous noterons~$\shil L S$ le sous-ensemble de~$X_L$ 
égal à la réunion des ensembles
finis~$\shil L x$ pour~$x\in S$ 
(pour une définition de~$\shil L x$, {\em cf.}~\ref{ex-ant-can}) ;
l'image de~$\shil L S$
sur~$X$ est égale à~$S$. 
Nous allons énoncer quelques propriétés 
élémentaires de la formation de~$\shil L S$ que nous utiliserons librement
par la suite. 

\trois{extscal-courbes-presqualg}
Si~$L$ est presque algébrique
sur~$k$,
l'ensemble~$\shil L S$ est simplement l'image
réciproque de~$S$ sur~$X_L$ (\ref{shil-l-presquealg}). 

\trois{extscal-courbes-comp}
Si~$L'$ est une extension complète de~$L$ alors~$\shil {L'}S=\shil {L'}{(\shil L S)}$
(\ref{shil-l-trans}). 

\trois{extscal-courbes-sx}
Écrivons~$L\otimes_k \got s(X)$ comme un produit
fini~$\prod L_i$ d'extensions finies séparables de~$L$. On a alors
$X_L=\coprod X\times_{\got s(X)}L_i$, et il résulte de~\ref{shil-l-x-connexe}
que~$\shil L S$ est la réunion disjointe des~$\shil {L_i} S$
(où~$\shil {L_i} S\subset X\times_{\got s(X)}L_i$ est défini pour tout~$i$ en voyant~$X$
comme espace~$\got s(X)$-analytique). 

\trois{extscal-courbes-card}
Soit~$x\in X$. L'ensemble~$\shil L x$ est en bijection naturelle avec l'ensemble des idéaux maximaux
de~$L\otimes_k \got s(x)$ (corollaire~\ref{ant-canon-cardinal}). 

\deux{exemple-pseudo-cour}
Soit~$X$ une $k$-couronne virtuelle. Le sous-ensemble~$\shil L {\skelan X}$
de~$X_L$ est égal à~$\skelan {X_L}$. 

\medskip
En effet, soit~$F$ une extension complète de~$k$ composée
de~$L$ et~$\KK$. Comme~$\KK$ est une extension presque algébrique
de~$k$, le sous-ensemble~$\shil {\KK}{\skelan X}$ de~$X_{\KK}$ 
coïncide
avec l'image réciproque de~$\skelan X$ sur~$X_{\KK}$, à savoir~$\skelan {X_{\KK}}$. 

\medskip
Puisque la courbe
$X_{\KK}$ est une~$\KK$-couronne, il résulte de~\ref{shil-l-etar}
que
$\shil F {(\skelan {X_{\KK}} )}$
est égal à~$\skelan{X_F}$ ; son image sur~$X_L$ est donc égale à~$\skelan {X_L}$. 
Par ailleurs,
$$\shil F {(\skelan {X_{\KK}} )}=\shil F
{(\shil {\KK}{\skelan X}) }=\shil F {\skelan X}=\shil F{(\shil L {\skelan X})},$$
et son image sur~$X_L$ 
est donc égale à~$\shil L {\skelan X}$ ;
il vient~$\shil L {\skelan X}=\skelan {X_L}$, comme annoncé. 

\deux{genre-preserve-shil-l}
Soit~$x\in X\typ {23}$ ; on a alors~$g(y)=g(x)$ 
pour tout~$y\in \shil L x$. En effet, soit~$F$
le complété d'une clôture algébrique d'une extension complète
composée de~$L$ et~$\KK$,  
soit~$y\in \shil L x$, soit~$z\in \shil F y$, et soit~$t$ l'image de~$z$
sur~$X_{\KK}$ ; le point~$t$ appartient à~$\shil {\KK} x$. 

\medskip
Par définition, $g(t)=g(x)$.
D'après le corollaire~\ref{ant-canon-alclos}, le point~$z$ est l'unique
antécédent de~$t$ sur~$X_F$, et le corpoïde
$\red{\hres(z)}$ s'identifie au corpoïde des fractions
de~$\red{\hres(t)}\otimes_{\kk} \red F$. On distingue maintenant deux cas. 

\trois{genre-shil-l-2}
{\em Supposons que~$x$ est de type 2}. Dans ce cas~$t$ est de type 2 aussi, et
il résulte de ce qui précède que
$\red{\hres(z)}$ est le corps des fractions de~$\red{\hres(t)}_1\otimes_{\kk_1} \red F_1$ ; en conséquence, 
on a~$g(z)=g(t)=g(x)$. 

\trois{genre-shil-l-3}
{\em Supposons que~$x$ est de type 3}. Dans ce cas
$g(x)=0$, le point~$t$ est également de type 3, et
$\red{\hres(t)}$ est de la forme~$\kk(\tau)$ pour un certain élément~$\tau$ 
transcendant sur~$\kk$, dont le degré n'appartient pas à~$|\KK\ti|$. Le corpoïde
$\red{\hres(t)}$ est alors égal à~$\red F(\tau)$, l'élément~$\tau$ restant transcendant sur~$\red F$. 
Si le degré de~$\tau$ appartient à~$|F\ti|$, le point~$t$ est de type 2
et~$\red {\hres(t)}_1=(\red F(\tau))_1$ est transcendant pur sur~$\red F_1$, ce qui implique que~$g(t)=0$ ; si
le degré de~$\tau$ n'appartient à~$|F\ti|$, le point~$t$ est de type 3 et~$g(t)=0$ par définition. 

\trois{genre-shil-l-conclu}
Dans tous les cas, on a~$g(t)=g(x)$ ; comme~$F$ est le complété d'une clôture algébrique de~$L$ et comme
$y$ est un antécédent de~$\shil L x$ sur~$X_F$, il vient~$g(\shil L x)=g(t)=g(x)$,
ce qu'il fallait démontrer. 

\deux{prop-extension-branche}
{\bf Proposition.}
{\em Soit~$X$ une courbe~$k$-analytique et soit~$x\in X\geom$. Soit~$U$
une composante connexe de~$X\setminus \{x\}$. Le bord de~$U_L$ dans~$X_L$
est alors contenu 
dans~$\shil L x$.}

\medskip
{\em Démonstration}. Soit~$\pi \colon X_L\to X$
la flèche canonique. 
Nous allons commencer par un certain nombre de réductions. 

\trois{bonbord-shill-ql}
Soit~$k^{\rm parf}$ le complété
d'une clôture parfaite de~$k$, et soit~$L'$
une extension complète de~$k$ composée de~$L$ et~$k^{\rm parf}$. Les applications
continues
naturelles
$X_{k^{\rm parf }} \to X$ et~$X_{L'}\to X_L$ sont des homéomorphismes, et l'on a
$$\shil {L'}x=\shil {L'}{(\shil{k^{\rm parf}}x)}=\shil {L'}{(\shil L x)}.$$
On peut donc, en remplaçant~$k$ par~$k^{\rm parf}$, 
$x$ par~$\shil {k^{\rm parf}}x$ 
et~$L$ par~$L'$, se ramener ainsi au cas où~$k$ est parfait. 

L'assertion requise étant insensible aux phénomènes de nilpotence, on peut alors
remplacer
$X$ par~$X_{\rm red}$, et donc la supposer réduite ; elle est dès lors génériquement quasi-lisse. 

\trois{bonbord-shill-section}
Choisissons une famille~$(Z_b)_{b\in \br X x\ctd U}$ où chaque~$Z_b$ est une section coronaire de~$b$ contenue dans~$U$. 
Le fermé~$V:=U\setminus \coprod Z_b$ de~$U$ est fermé dans~$X$. Il s'ensuit
que le bord de~$U_L$ dans~$X_L$ est contenu dans~$\bigcup_b \partial Z_{b,L}\cap \pi^{-1}(x)$. 
Fixons~$b$ et écrivons~$Z$ au lieu de~$Z_b$. 
Par ce qui précède, il suffit pour conclure de prouver
que~$\partial Z_L\cap \pi^{-1}(x)\subset \shil L x$. 

\trois{bonbord-coronaire-rigide}
{\em Supposons que~$x$ est un point rigide.}
Dans ce cas~$\shil L x=\pi^{-1}(x)$ et il n'y a rien à démontrer.

\trois{bonbord-coronaire-23}
{\em Supposons que~$x\in X\typ {23}$.}
Choisissons un voisinage affinoïde~$V$ de~$x$
dans~$X$ et une fonction
analytique~$f$ sur~$V$ telle que~$|f(x)|\neq 0$ et telle
que~$\langle \red{f(x)}\rangle_ b <1$. Le morphisme~$\phi \colon V\to \Aff^{1,\rm an}_k$
définit une présentation d'Abhyankar de~$x$ ; on a~$\phi(x)=\eta_r$, où~$r=|f(x)|$. Par
construction, $\phi(b)$ est la branche issue de~$\eta_r$ et
contenue dans l'ouvert de~$\Aff^{1, \rm an}_k$
défini par la condition~$|T|<r$. 

D'après la proposition~\ref{degetseccor}, il existe une section coronaire~$Z'$ de~$b$ et une section
coronaire~$Z''$ de~$\phi(b)$ telles que~$\phi$ induise un morphisme fini et plat de~$Z'$ sur~$Z''$ ; quitte
à restreindre~$Z'$ et~$Z''$, on peut supposer que~$Z''$ est décrite par une inégalité de la forme~$s<|T|<r$ pour un
certain~$s\in ]0;r[$, et que~$Z'$ est une sous-couronne virtuelle de~$Z$. 

\medskip
Soit~$y$ un point de~$\overline{Z_L}$ situé au-dessus de~$x$ ; son image~$t$ sur~$\Aff^{1,\rm an}_L$
est située au-dessus de~$\eta_r$. 

Le point~$y$ appartient à l'adhérence d'une composante connexe~$Y$ de~$Z_L$, et partant
à~$\overline {\skelan Y}$. Soit~$Y'$ la sous-couronne virtuelle~$Y\times_Z Z'$ de~$Y$. Comme
le point~$x$
n'est pas adhérent à~$\skelan Z\setminus \skelan {Z'}$, le point~$y$ n'est pas adhérent à~$\skelan Y\setminus \skelan  {Y'}$, 
et il appartient en conséquence à~$\overline{\skelan {Y'}}$.
On en déduit que~$t\in \overline{\skelan {Z''_L}}=[\eta_{s,L};\eta_{r,L}]$. 
Mais puisque
$t$ est situé au-dessus
de~$\eta_r$, il
vient~$t=\eta_{r,L}$, et~$y$ appartient de ce fait à~$\shil L x$.~$\Box$

\subsection*{Les composantes qui apparaissent sont des disques virtuels}

\deux{theo-nouvelles-disques}
{\bf Théorème}. {\em Soit~$X$ une courbe~$k$-analytique génériquement quasi-lisse et soit~$x\in X\geom$. 
Soit~$L$ une extension complète de~$k$, et soit~$\pi \colon X_L\to X$ la flèche canonique. 
L'ouvert~$X_L\setminus \shil L x$
de~$X_L$
est réunion disjointe d'une part des composantes connexes de~$(X\setminus \{x\})_L$,
d'autre part de disques virtuels relativement compacts tous contenus dans~$\pi^{-1}(x)$.}

\medskip
{\em Démonstration}.
Si~$x$ est un point rigide $\shil L x=\pi^{-1}(x)$ et l'assertion 
est évidente. Supposons maintenant que~$x\in X\typ {23}$. 
Soit~$U$ une composante connexe de~$X\setminus \{x\}$. Il résulte
de la proposition~\ref{prop-extension-branche}
que le bord de~$U_L$ est contenu dans~$\shil L x$ ; par conséquent, $U_L$
est une réunion de composantes connexes de~$X_L\setminus \shil L x$. 

\medskip
On en déduit que~$X_L\setminus \shil L x$ est réunion disjointe des composantes connexes de~$(X\setminus \{x\})_L$,
et de ses
composantes connexes contenues dans~$\pi^{-1}(x)$. Il s'agit maintenant de montrer que ces dernières sont toutes des
disques virtuels relativement compacts. Soit donc~$V$ une telle composante. 

\trois{composantes-nouvelles-reduction}
{\em Réduction à un cas particulier plus simple.}
Choisissons un voisinage affinoïde~$W$ de~$x$ dans~$X$ tel que~$\got s(W)=\got s(x)$. 
Écrivons~$L\otimes_k \got s(W)=\prod L_i$, où les~$L_i$
sont des extensions finies séparables de~$L$ ; on a alors~$W_L=\coprod W\times_{\got s(W)}L_i$, et~$\shil L x$
possède un et un seul élément sur chaque~$W\times_{\got s(W)}L_i$. La composante connexe~$V$ étant 
contenue dans~$\pi^{-1}(x)$, elle est nécessairement incluse dans~$W\times_{\got s(W)}L_i$
pour un certain~$i$. En remplaçant~$X$ par~$W$, en considérant ce dernier comme un espace~$\got s(W)$-analytique
et en remplaçant~$k$ par~$\got s(W)$ et~$L$ par~$L_i$, 
on se ramène au cas où~$\got s(X)=\got s(x)=k$. Cette dernière condition assure que~$\shil F x$ est un singleton
pour toute extension complète~$F$ de~$k$. Par abus, son unique élément sera encore noté~$\shil F x$. 

\medskip
Soit~$L'$ le complété d'une clôture algébrique de~$L$ ; on fixe
un~$k$-plongement isométrique~$\KK\hookrightarrow L'$. 
Le point~$\shil {L'}x$ est l'unique antécédent de~$\shil L x$ sur~$X_{L'}$, et~$V_{L'}$ est donc une
union de composantes connexes
de~$X_{L'}\setminus \shil {L'}x$. Elle est contenue dans l'image réciproque de~$x$ par la flèche~$X_{L'}\to X$,
et donc dans l'image réciproque de~$\shil {\KK}x$ par la flèche~$X_{L'}\to X_{\KK}$, car~$\shil {\KK}x$
est l'unique
antécédent de~$x$ sur~$X_{\KK}$. Il suffit pour conclure de montrer que~$V_{L'}$ est une union disjointe de disques relativement compacts
de~$X_{L'}$. 

En vertu de ce qui précède on peut, en remplaçant~$k$ par~$\KK$, $x$ par~$\shil {\KK} x$, $L$ par~$L'$ et~$V$ par
une composante connexe de~$V_{L'}$, se ramener
au cas où~$k$ et~$L$ sont
algébriquement clos. 

\medskip
Comme~$V$ est contenue dans~$\pi^{-1}(x)$, on peut si besoin
restreindre ou agrandir~$X$ autour de~$x$. Puisque~$X$ est génériquement quasi-lisse, 
elle est quasi-lisse au voisinage de~$x$, ce qui autorise à la restreindre
de sorte que~$X\setminus \{x\}$ soit une union disjointe de disques ouverts relativement compacts et d'un nombre fini
de couronnes ouverte de bord~$\{x\}$. 

En prolongeant chacune de ces dernières en un disque, on obtient une courbe analytique propre, lisse, connexe et non vide ; on 
peut ainsi finalement supposer qu'il existe une~$k$-courbe projective, irréductible et lisse~$\sch X$ telle que~$X=X\an$ 
et telle que~$\sch X\an\setminus x$ soit réunion disjointe de disques. 

\trois{conclusion-nouvelles-composantes}
Soit~$g$ le genre de~$\sch X$. Par construction, $\{x\}$ est une triangulation de~$\sch X\an$. 
On déduit de~\ref{calcul-de-genre}
que~$g(x)=g$. 

\medskip
D'après~\ref{genre-preserve-shil-l}, on a~$g(\shil L x)=g(x)=g$. Comme
le genre de~$\sch X_L$ est égal à~$g$, il résulte de~\ref{calcul-de-genre}
que~$\sch X\an_L$ ne possède aucune boucle (comme~$\sch X_L\an$ est par ailleurs compacte,
connexe
et non vide pusique~$\sch X_L$ est projective et irréductible, $\sch X_L\an$ est un arbre compact), 
et que~$g(y)=0$ 
pour tout~$y\in (\sch X_L\an)\typ{23}\setminus \shil L x$. 

Puisque~$\sch X_L\an$ est un arbre compact, $V$ est un arbre à un bout relativement compact ; puisque~$\shil L x\notin V$, la composante~$V$
ne contient aucun point de genre~$>0$. L'unique point~$\shil L x$
de~$\partial V$ étant de type 2 ou 3, il résulte de la proposition~\ref{propcarintdc}
que~$V$ est un disque, ce qui achève la démonstration.~$\Box$ 

\deux{coro-nouvelles-composantes}
{\bf Corollaire}. 
{\em Soit~$X$ une courbe~$k$-analytique et soit~$x\in X\geom$. 
Soit~$L$ une extension complète de~$k$, et soit~$\pi \colon X_L\to X$ la flèche canonique. 
L'ouvert~$X_L\setminus \shil L x$
de~$X_L$
est réunion disjointe d'une part des composantes connexes de~$(X\setminus \{x\})_L$,
d'autre part d'arbres à un bout relativement compacts tous contenus dans~$\pi^{-1}(x)$.}

\medskip
{\em Démonstration}. Soit~$k^{\rm parf}$ le composé d'une clôture parfaite de~$k$, et soit~$L'$
une extension complète de~$k$ composée de~$L$ et~$k^{\rm parf}$. Les flèches
naturelles
$X_{k^{\rm parf }} \to X$ et~$X_{L'}\to X_L$ sont des homéomorphismes 
de réciproques respectives~$x\mapsto \shil {k^{\rm parf}}x$ et~$x\mapsto \shil{L'} x$, et l'on a
on a~$\shil {L'}x=\shil{L'}{(\shil {k^{\rm parf}}x)}=\shil {L'}{(\shil L x)}$ pour tout~$x\in X$. 
On peut dès lors, en remplaçant~$k$ par~$k^{\rm parf}$, $x$ par
son unique antécédent sur~$X_{k^{\rm parf}}$, 
et~$L$ par~$L'$, se ramener au cas où~$k$ est parfait. 
L'assertion requise étant insensible aux phénomènes de nilpotence, on peut remplacer
$X$ par~$X_{\rm red}$, et donc la supposer réduite, et partant génériquement quasi-lisse. 
Le corollaire est alors une conséquence immédiate
du théorème~\ref{theo-nouvelles-disques}.~$\Box$

\subsection*{Effets sur les sous-graphes et sur les triangulations}

\deux{theo-ferme-shil}
{\bf Théorème.}
{\em Soit~$X$ une courbe~$k$-analytique, soit~$\Gamma$ un fermé de~$X\geom$
et soit~$L$ une extension complète de~$k$. Soit~$\pi \colon X_L\to X$ la flèche canonique. 

\begin{itemize}

\medskip
\item[i)] Le sous-ensemble~$\shil L \Gamma$
de~$X_L$ est fermé. 

\medskip
\item[ii)] Les composantes connexes de~$X_L\setminus \shil L \Gamma$ sont
exactement les ouverts~$U$ de~$X$ satisfaisant 
l'une des deux assertions suivantes, exclusives l'une de l'autre. 

\begin{itemize}
\item[a)] L'ouvert~$U$ est une composante connexe de~$(X\setminus \Gamma)_L$. 

\item[b)] L'ouvert~$U$ est contenu dans~$\pi^{-1}(x)$ pour un certain
$x\in \Gamma\typ{23}$ et est une composante connexe de~$X_L\setminus \shil L x$ ; c'est alors un arbre à un bout
relativement compact, et un disque virtuel si~$X$ est génériquement quasi-lisse. 

\end{itemize}

\medskip
\item[iii)] Si~$\Gamma$ est un sous-graphe (resp. un sous-graphe localement fini)
de~$X$ alors~$\shil L\Gamma$ est un sous-graphe (resp. un sous-graphe localement fini)
de~$X_L$ ; la surjection $\shil L \Gamma \to \Gamma$ est ouverte, et est un homéomorphisme
si~$k$ est algébriquement clos.

\medskip
\item[iv)] Si~$\Gamma$ est un sous-graphe localement fini
de~$X$, l'application continue~$\shil L \Gamma \to \Gamma$ est une isométrie par morceaux. 

\medskip
\item[v)] Si~$X$ est génériquement quasi-lisse et si~$\Gamma$ est un sous-graphe de~$X$, il 
est analytiquement admissible si et seulement si le sous-graphe~$\shil L \Gamma$ de~$X_L$ est analytiquement
admissible.

\end{itemize}
}

\medskip
{\em Démonstration}. Nous allons commencer par démontrer à la fois les assertions~i) et~ii). 

\trois{memecoup-i-ii}
{\em Preuve de~i) et~ii)}.
Par définition, {\em l'ensemble}
$X_L\setminus \shil L \Gamma$ est la réunion disjointe de l'ouvert~$(X\setminus \Gamma)_L$
et des~$\pi^{-1}(x)\setminus \shil L x$ pour~$x$ parcourant~$\Gamma$. Or il résulte du théorème~\ref{theo-nouvelles-disques}
et de son corollaire~\ref{coro-nouvelles-composantes}
que pour tout~$x\in X\geom$,  l'ensemble~$\pi^{-1}(x)\setminus \shil L x$ est la 
réunion disjointe de composantes connexes de~$X_L\setminus \shil L x$, chacune
d'elle étant un arbre à un bout relativement compact, et un disque virtuel si~$X$ est génériquement quasi-lisse. 

\medskip
On déduit de ce qui précède que~$X_L\setminus \shil L \Gamma$ est un ouvert de~$X_L$, dont les composantes connexes sont 
précisément les ouverts décrits en~ii) ; on a ainsi démontré~i) et~ii). 

\medskip
\trois{shil-l-bon-top}
{\em Preuve de~iii)}. Commençons par une remarque. Si~$k$ est algébriquement clos, 
$\shil L x$ est un singleton pour tout~$x\in X \geom$, et la surjection continue~$\shil L \Gamma \to \Gamma$
est dès lors injective. 

\medskip
L'assertion~iii) étant locale sur~$X$, on peut 
supposer ce dernier compact. Soit~$L'$ le complété d'une clôture algébrique de~$L$ ; 
fixons un~$k$-plongement isométrique de~$\KK$ dans~$L'$. 
Supposons
que~$\Gamma$ soit un sous-graphe (resp. un sous-graphe fini-
de~$X$. 
Le compact~$\shil {\KK}\Gamma$ de~$X_{\KK}$ est l'image réciproque de~$\Gamma$ sur~$X_{\KK}$, 
et est donc en vertu du corollaire~\ref{corollimrecih}
un sous-graphe (resp. un sous-graphe fini)
de~$X_{\KK}$. 

\medskip
D'après l'assertion~i)
déjà établie, $\shil {L'}{(\shil {\KK}\Gamma)}$ est un compact de~$X_{L'}$. Puisque~$\KK$
est algébriquement clos,
la
surjection continue de~$\shil {L'}\Gamma=\shil {L'}{(\shil {\KK}\Gamma)}$ vers~$\shil {\KK}\Gamma$ est 
bijective par la remarque ci-dessus, et elle induit donc par compacité un homéomorphisme entre ses source et but ;
le compact~$\shil {L'}\Gamma$ est en conséquence un sous-graphe (resp. un sous-graphe fini)
de~$X_{L'}$. Comme~$L'$ est une extension presque algébrique de~$L$, le 
compact $\shil {L'}\Gamma=\shil {L'}{(\shil L \Gamma)}$ est l'image
réciproque de~$\shil L \Gamma$ sur~$X_{L'}$, et~$\shil L \Gamma$ s'identifie au quotient de~$\shil {L'}\Gamma$
par l'action de Galois. On déduit alors du théorème~\ref{theoquot}
que~$\shil L \Gamma$ est un sous-graphe (resp. un sous-graphe fini)
de~$X_L$. 

\medskip
Soit~$\Delta$ un ouvert de~$\shil L \Gamma$. Son image réciproque~$\shil {L'}\Delta$ sur~$\shil {L'}\Gamma$ est un ouvert
de ce dernier ; 
comme la flèche
$\shil {L'}\Gamma\to \shil {\KK}\Gamma$ est un homéomorphisme, elle identifie~$\shil {L'}\Delta$ à un ouvert~$\Omega$
de~$\shil {\KK}\Gamma$ ; l'application~$\shil {\KK}\Gamma\to \Gamma$ étant un quotient, elle induit
une surjection de~$\Omega$ sur un ouvert~$\Omega'$
de~$\Gamma$. Par surjectivité de~$\shil {L'}\Gamma \to \shil L \Gamma$, l'image de~$\Delta$
sur~$\Gamma$ est égale à~$\Omega'$. 

\medskip
Ainsi, la surjection $\shil L \Gamma \to \Gamma$ est ouverte.
Si~$k$ est algébriquement clos, on a signalé plus haut qu'elle est bijective, et c'est donc un homéomorphisme
(soit parce qu'elle est ouverte comme on vient de le voir, soit par compacité). 

\trois{shil-l-isom}
{\em Preuve de~iv).}
Soit~$F$
l'extension complète de~$k$ composée de~$L$ et~$k^{\rm parf}$. Les flèches
naturelles
$X_{k^{\rm parf }} \to X$ et~$X_F\to X_L$ sont des homéomorphismes 
de réciproques respectives~$x\mapsto \shil {k^{\rm parf}}x$ et~$x\mapsto \shil{L'} x$
qui préservent les toises canoniques, et l'on
a~$\shil {L'}x=\shil{L'}{(\shil {k^{\rm parf}}x)}=\shil {L'}{(\shil L x)}$ pour tout~$x\in X$. 
On peut dès lors, en remplaçant~$k$ par~$k^{\rm parf}$, $\Gamma$ par~$\shil {k^{\rm parf}}\Gamma$
et~$L$ par~$F$, se ramener au cas où~$k$ est parfait. 
L'assertion requise étant insensible aux phénomènes de nilpotence, on peut remplacer
$X$ par~$X_{\rm red}$, et donc la supposer réduite, et partant génériquement quasi-lisse. 

\medskip
Choisissons un ensemble fermé et discret~$S$ de points de~$\Gamma$ tel
que~$\Gamma \setminus S$ soit réunion disjointe d'intervalles ouverts, et soit~$S'$ une 
triangulation de~$X$ contenant~$S$. Soit~$\Gamma'$ le squelette de~$S$. Il contient~$\Gamma$, et il
suffit pour conclure
de montrer que~$\shil  L {\Gamma'} \to \Gamma'$ est une isométrie
par morceaux. 
Soit~$I$ une composante connexe de~$\Gamma'\setminus S$ ; c'est un intervalle
ouvert et~$I^\flat$ est une couronne virtuelle. On peut alors écrire~$I^\flat_L$ comme une union disjointe
finie~$\coprod Y_i$ de couronnes virtuelles. On a~$\shil L I=\coprod \skelan{Y_i}$, et
l'invariance du module des couronnes virtuelles entraîne que~$\skelan{Y_i}\to \skelan {I^\flat}=I$
est une isométrie pour tout~$i$, ce qui achève de montrer~iv).

\trois{shil-l-equivadm}
{\em Preuve de~v)}. Supposons~$\Gamma$ analytiquement admissible et soit~$U$ une composante connexe de~$X\setminus \shil L \Gamma$. 
D'après l'assertion~ii) déjà établie, on est dans l'un des deux cas suivants. 

\medskip
$\bullet$ {\em Premier cas.} L'ouvert~$U$ est une composante connexe de~$(X\setminus \Gamma)_L$. C'est alors une
composante connexe de~$V_L$ pour une
certaine composante connexe~$V$ de~$X\setminus \Gamma$ ; comme~$\Gamma$ est analytiquement admissible, 
$V$ est un disque virtuel relativement compact, et~$U$ est donc également un disque virtuel relativement compact. 

$\bullet$ {\em Second cas}. L'ouvert~$U$ est contenu dans l'image réciproque de~$x$ pour un certain~$x\in \Gamma$, et est un disque
virtuel relativement compact. 

\medskip
Ainsi, $U$ est dans tous les cas un disque virtuel relativement compact, et~$\shil L \Gamma$ est analytiquement admissible.

\medskip
Réciproquement, supposons que~$\shil L \Gamma$ est analytiquement admissible, et soit~$V$ une composante connexe de~$X\setminus \Gamma$. 
Toujours d'après~ii), $V_L$ est une union disjointe de composantes connexes de~$X_L\setminus \shil L \Gamma$, qui sont des disques virtuels relativement 
compacts puisque~$\shil L \Gamma$ est analytiquement admissible. On en déduit que~$V$ est un disque virtuel relativement compact, et~$\Gamma$ est dès lors
analytiquement admissible.~$\Box$ 

~

\deux{theo-invariance-topologie}
{\bf Théorème}.
{\em Soit~$X$ une courbe~$\KK$-analytique et soit~$L$ une extension complète de~$K$. L'application continue
compacte canonique~$X_L\to X$ est une équivalence homotopique, et~$\H^i_{\rm c}(X\tp, \Lambda)\to \H^i_{\rm c}(X_{L,\rm top}, \Lambda)$ 
est un isomorphisme pour tout~$i\in \{0,1\}$ et tout groupe abélien~$\Lambda$.}

\medskip
{\em Démonstration}.
Fixons un sous-graphe localement fini et analytiquement admissible~$\Gamma$ de~$X$. 
Le théorème~\ref{theo-ferme-shil}
assure que~$\shil L \Gamma$ est un sous-graphe analytiquement admissible de~$X_L$ 
et que~$\shil L \Gamma \to \Gamma$ est un homéomorphisme. Le théorème provient alors du fait 
que~$X$ se rétracte par déformation sur~$\Gamma$, que~$X_L$ se rétracte par déformation sur~$\shil L \Gamma$, 
et que les rétractions~$X\to \Gamma$ et~$X_L\to \shil L \Gamma$ sont compactes.~$\Box$

\deux{theo-invariance-triangulations}
{\bf Théorème.}
{\em Soit~$X$ une courbe~$k$-analytique génériquement quasi-lisse, soit~$S$ une partie fermée et discrète de~$X\geom$, et soit~$L$
une extension complète de~$k$. Les assertions suivantes sont équivalentes. 

\medskip
i) $S$ est une triangulation de~$X$. 

ii) $\shil L S$ est une triangulation de~$X_L$, et pour toute composante connexe~$U$ de~$X_L\setminus \shil L S$ qui est une
couronne virtuelle, la restriction de~$X_L\to X$ à~$\skelan U$ est injective.} 

\medskip
{\em Démonstration.}
Soit~$\pi \colon X_L\to X$ la flèche canonique. On procède par double implication. 

\trois{tri-shil-l-direct}
{\em Prouvons que~i)$\Rightarrow$ii)}. On suppose donc que~$S$ est une triangulation de~$X$. Il est immédiat que~$\shil L S$ est fermé
et discret dans~$X_L$ (les fibres de~$\shil L S\to S$ étant finies).
Soit~$U$ une composante connexe de~$X_L\setminus \shil L S$. D'après le théorème~\ref{theo-ferme-shil}, on est dans l'un des deux cas suivants : 

\medskip
$\bullet$ $U$ est contenu dans~$\pi^{-1}(S)$, et $U$ est un disque virtuel relativement compact ; 

$\bullet$ $U$ est une composante connexe de~$(X\setminus S)_L$.

\medskip
Supposons qu'on soit dans le second cas. Dans ce cas, $U$ est une composante connexe de~$V_L$ pour une certaine composante
connexe~$V$ de~$X\setminus S$. Comme~$S$ est une triangulation, on distingue à nouveau deux cas : 

\medskip
$\bullet$ $V$ est un disque virtuel relativement compact, auquel cas~$U$ est un disque virtuel relativement compact ; 

$\bullet$ $V$ est une couronne virtuelle relativement compacte ; dans ce cas, $U$ est une couronne virtuelle relativement compacte,
et~$\skelan U \to \skelan V$ est un homéomorphisme. 

\medskip
Ainsi, ii) est vérifiée. 

\trois{tri-shil-l-indirect}
{\em Prouvons que~ii)$\Rightarrow$i)}. On suppose donc que~$\shil L S$ est une triangulation de~$X_L$. Soit~$F$ le complété d'une clôture
algébrique de~$L$ ; fixons un~$k$-plongement isométrique de~$\KK$ dans~$L$. 

\medskip
Soit~$V$ une composante connexe de~$X\setminus S$, et soit~$U$ une composante
connexe de~$V_{\KK}$. D'après l'implication~i)$\Rightarrow$ii) déjà établie, 
$\shil F S$ est une triangulation de~$X_F$. Le théorème~\ref{theo-ferme-shil}
assure que~$V_F$ est une union disjointe de composantes connexes de~$X_F\setminus \shil F S$. 
Par conséquent, l'ouvert~$U_F$ est une composante connexe de~$X_F\setminus \shil F S$, et il y a deux cas à distinguer. 

\medskip
{\em Premier cas : $U_F$ est un disque virtuel relativement compact dans~$X_F$.}. Dans ce cas, $U$ est relativement compact
dans~$X_{\KK}$, 
et le théorème~\ref{theo-invariance-topologie}
assure que~$U$ est un arbre à un bout. Soit~$x$ l'unique point de~$\partial U_F$. Son image~$y$ sur~$X_F$ est l'unique point de~$\partial U$, et on a nécessairement
$x=\shil F y$. Comme~$U_F$ est un disque, $x$ est de type 2 ou 3 ; puisque~$t\mapsto \shil F t$ préserve le rang des points
d'Abhyankar, $y$ est de rang 1, donc de type 2 ou 3.  Enfin, puisque~$U_F$ est un disque virtuel, il ne contient aucun potin de genre~$>0$ ; par conséquent, $U$ ne contient aucun 
point de genre~$>0$. Il résulte alors de la proposition~\ref{propcarintdc}
que~$U$ est un disque. Par conséquent, $V$ est un disque virtuel, relativement compact dans~$X$
puisque~$U$ est relativement compact dans~$X_{\KK}$. 

\medskip
{\em Second cas : $U_F$ est une couronne virtuelle relativement compacte dans~$X_F$.} Il existe 
une composante connexe~$U'$ de~$V_L$ telle que~$U_F$ soit une composante connexe
de~$U'_F$ ; d'après le théorème~\ref{theo-ferme-shil}, $U'$ est une composante connexe de~$X_L\setminus \shil L S$. 
Comme~$U_F$ est une couronne virtuelle, $U'$ est une couronne virtuelle.

Comme~$U_F$ est une couronne virtuelle, 
c'est un arbre à deux bouts relativement compact dans~$X_F$, dont tout point du bord (qui en comprend au plus deux) est de type 2 ou 3, et qui ne contient
aucun point de genre~$>0$. Par un raisonnement analogue
à celui suivi dans le premier cas, $U$ est un arbre à deux bouts relativement compact
dans~$X_{\KK}$, dont tout point du bord est de type 2 ou 3, et qui ne contient aucun point de genre~$>0$. 
Il résulte alors de la proposition~\ref{propcarintdc}
que~$U$ est une couronne. Pour montrer que~$V$ est une couronne virtuelle
relativement compacte dans~$X$, il suffit de vérifier que~$\skelan U\to V$ est injective. 

D'après l'hypothèse~ii), $\skelan{U'}\to V$ est injective ; par ailleurs, la flèche 
$\skelan{U_F}\to \skelan{U'}$ est un homéomorphisme
puisque~$U'$ est une couronne virtuelle, et~$\skelan {U_F}\to \skelan U$ est un homéomorphisme car~$U$ est une couronne. 
On en déduit que~$\skelan U\to V$ est injective.~$\Box$

\section{Triangulations minimales d'une courbe analytique quasi-lisse compacte}

\subsection*{Triangulations des courbes de genre zéro}

\deux{gzeroptinv} Soit~${\sch X}$ une~$k$-courbe algébrique projective, lisse, géométriquement intègre de genre zéro. Comme~${\sch X}_{\KK}\simeq \pkk$ la courbe analytique~${\sch X}\an_{\KK}$ est un arbre, et~${\sch X}\an$ est par conséquent un arbre. 

\trois{ptinvtrimin} Soit~$x\in{\sch X}\an\dtr$ tel que~$\got s(x)=k$ ; le point~$x$ possède alors un unique antécédent~$x_{\KK}$ sur~${\sch X}_{\KK}\an$. La~$\KK$-courbe~${\sch X }\an_{\KK}$ étant isomorphe à~$\pkk$, son ouvert~${\sch X}\an_{\KK}\setminus\{x_{\KK}\}$ est réunion disjointe de disques ; par conséquent,~$\{x_{\KK}\}$ est une triangulation de~${\sch X}\an_{\KK}$ ; il s'ensuit que~$\{x\}$ est une triangulation de~${\sch X}\an$, évidemment minimale. 

\trois{triminptinv} Soit~$S$ une triangulation de~${\sch X}\an$, et soit~$\Gamma$ le squelette de la triangulation~$S_{\KK}$ de~${\sch X}\an_{\KK}$ ; c'est un sous-graphe compact de~${\sch X}_{\KK}\an$ qui est admissible ; c'est donc un arbre non vide. Si~$I$ est une arête de~$\Gamma$, l'action de~$\mathsf G$ n'échange pas les deux orientations de~$I$ (théorème~\ref{theo-invariance-triangulations}) ; il en découle qu'il y a au moins un sommet de~$\Gamma$ (et donc au moins un point de~$S_{\KK}$) qui est fixé par~$\mathsf G$ (\ref{quotarbrefin}). De ce fait, il existe~$x\in S$ tel que~$\got s(x)=k$. Ainsi, les triangulations minimales de~${\sch X}\an$ sont exactement les triangulations de la forme~$\{x\}$, où~$x$ est un point de type 2 ou 3 tel que~$\got s(x)=k$. 

\trois{discvirtk} Supposons qu'il existe ouvert~$V$ de~$X$ qui est un disque virtuel sur~$k$ (autrement dit, un disque virtuel géométriquement connexe). Il existe alors un ouvert ~$W$ de~$V$ qui est une couronne virtuelle sur~$k$ (\ref{discvcourv}). Si~$y\in \skel W$ le point~$y$ appartient à~$\sch X\an\dtr$ et~$\got s(y)=k$ ; par conséquent,~${\sch X}\an$ possède une infinité de triangulations minimales dont l'unique sommet est situé sur~$V$ ; si~$|k\ti|\neq\{1\}$, elle en possède même une infinité dont l'unique sommet est situé sur~$V$ et est de type 2. 

\deux{casdeuxtrimin} Supposons que~${\sch X}\an$ possède au moins deux triangulations minimales ; sous cette hypothèse, il existe deux points~$x$ et~$x'$ distincts sur~$\sch X\an\dtr$ et tels que~$\got s(x)=\got s(x')=k$. Le squelette de toute triangulation de~${\sch X}\an$ contenant le squelette analytique de~${\sch X}\an$, celui-ci est inclus dans~$\{x\}\cap\{x'\}=\emptyset$ et est donc vide. 

\medskip
Réciproquement, supposons que le squelette analytique de~${\sch X}\an$ soit vide. Choisissons une triangulation minimale de~${\sch X}\an$, c'est-à-dire un point~$x$ de~${\sch X}\an$ de type 2 ou 3 de~${\sch X}\an$ tel que~$\got s(x)=k$. Le squelette analytique de~${\sch X}\an$ étant vide,~$x$ a un voisinage~$V$ qui est un disque virtuel. Comme~$\got s(x)=k$, on a nécessairement~$\got s(V)=k$ ; il résulte alors de~\ref{discvirtk} que~${\sch X}\an$ possède une infinité de triangulations minimales, et même,  si~$|k\ti|\neq\{1\}$, une infinité de triangulations minimales dont l'unique sommet est de type 2. 

\deux{casunetrimin} Supposons maintenant que~${\sch X}\an$ possède une unique triangulation minimale (ou encore une plus petite triangulation), c'est-à-dire un unique point~$x$ de type 2 ou 3 tel que~$\got s(x)=k$. En vertu du~\ref{casdeuxtrimin} ci-dessus, le squelette analytique de~${\sch X}\an$ est non vide ; étant contenu dans le squelette de la triangulation~$\{x\}$, il est égal à~$\{x\}$. 

\medskip
Le point~$x$ est de type 2 : en effet s'il était de type~$3$ il possèderait un voisinage~$Z$ qui serait une couronne virtuelle telle que~$x\in \mathsf  S(Z)$ ; comme~$\got s(x)=k$ on aurait~$\got s(Z)=k$ et donc~$\got s(y)=k$ pour tout~$y\in \mathsf S(Z)$, ce qui fournirait une infinité de triangulations minimales de~${\sch X}\an$, contrairement à notre hypothèse. 

\medskip
Si~$V$ est une composante connexe de~${\sch X}\an\setminus\{x\}$, c'est un arbre à un bout (puisque~${\sch X}\an$ est un arbre) et donc un disque virtuel ; comme~${\sch X}\an$ ne possède qu'une triangulation minimale,~$V$ n'est pas géométriquement connexe (\ref{discvirtk}). 

\medskip
\deux{unetrigeco} Inversement, soit~$x$ un point de~${\sch X}\an$ de type~$2$ ou~$3$ tel que~$\got s(x)=k$ et tel qu'aucune composante connexe de~${\sch X}\an\setminus\{x\}$ ne soit géométriquement connexe ; dans ce cas, tout point~$y$ de~${\sch X}\an$ distinct de~$x$ appartient à un ouvert connexe de corps des constantes non trivial, ce qui entraîne que~$\got s(y)\neq k$ ; il en résulte que~$\{x\}$ est l'unique triangulation minimale de~${\sch X}\an$. 

\deux{extrigzero} {\bf Quelques exemples.} 

\trois{puninfin} La courbe~$\pk$ possède une infinité de triangulations minimales ; par exemple si~$a\in k$ et si~$r>0$ alors~$\{\eta_{a,r}\}$ est une triangulation minimale de~$\pk$. 

\trois{conpimpair} Supposons que~$p\neq 2$ et que~$\red k$ est parfait. Soit~$\sch X$ une~$k$-conique non déployée. Il existe~$\lambda\in k\ti$,~$a\in k$ et un morphisme fini et plat~${\sch X}\to \PP^1_k$, dont on notera~$\phi$ l'analytifié,  tel que~${\sch X}\times_{\PP^1_k}\Aff^1_k$ s'identifie à~$\spec k[T,X]/(X^2-\lambda(T^2-a))$. Comme~${\sch X}(k)=\emptyset$, l'élément~$a$ de~$k$ n'est pas un carré. Posons~$r=\sqrt {|a|}$. 

\medskip
{\em Le point~$\eta_{r,\KK}$ de~$\pkk$ a un et un seul antécédent sur~${\sch X}\an_{\KK}$}.  En effet, posons~$\tau=T(\eta_{r,\KK})$ et~$\alpha=\lambda(\tau^2-a)$. On a~$|\tau^2-a|=|\tau^2|=|a|=r^2$. Par conséquent,~$$\gred \alpha=\gred \lambda ((\red \tau^r)^2-\red a ^{r^2}).$$ Comme~$\red \tau^r$ est transcendant sur~$\gred {\KK}$  et engendre~$\gred {\hres(\eta_{r,\KK})}$ (en tant qu'extension de~$\gred {\KK}$), comme~$\red a^{r^2}\neq 0$ et comme~$p$ est impair,~$\gred \alpha$ n'est pas un carré dans~$\gred {\hres(\eta_{r,\KK})}$ ; par conséquent,~$\alpha$ n'est pas un carré dans~$\hres(\eta_{r,\KK})$, d'où notre assertion. 

\medskip
Il s'ensuit que le point~$\eta_r$ {\em de~$\pk$} a un et un seul antécédent~$x$ sur~${\sch X}\an$, et que~$\got s(x)=k$.

\medskip
{\em Aucune composante connexe de~${\sch X}\an\setminus\{x\}$ n'est géométriquement connexe.} En effet, soit~$V$ une composante connexe de~${\sch X}\an\setminus\{x\}$ ; nous allons montrer qu'elle n'est pas géométriquement connexe. Comme~$x$ est le seul antécédent de~$\eta_r$ par~$\phi$, l'ouvert~$V$ est une composante connexe de~$\phi\inv(U)$ pour une certaine composante connexe~$U$ de~$\pk\setminus\{\eta_r\}$. Si~$U$ n'est pas géométriquement connexe, l'assertion est claire. 

Supposons que~$U$ soit géométriquement connexe, et montrons que~$U$ est un~$k$-disque. Soit~$\psi :\pk\to \pk$ la flèche induite par la fonction~$T^2/a$. Soit~$S$ la fonction coordonnée sur le but de~$\psi$  ; posons~$\omega=\psi(\eta_r)$ et~$\sigma=S(\omega)$. On a~$|\sigma|=1$ et~$\gred \sigma=(\gred \tau)^2/\gred a$. Il en résulte : d'une part que~$\gred \sigma$ est transcendant sur~$\gred k$, et donc que~$\omega=\eta_1$ ; d'autre part que~$[\gred{\hres(\eta_r)}:\gred{\hres(\omega)}]\geq 2$. Comme~$\psi$ est lui-même de degré~$2$, on a en réalité~$[\gred{\hres(\eta_r)}:\gred{\hres(\omega)}]=2$ et donc~$[\hres(\eta_r):\hres(\omega)]=2$.

Par conséquent,~$\eta_r$ est le seul antécédent de~$\eta_1$ par~$\psi$, et~$U$ est donc une composante connexe de~$\psi\inv(U')$ pour une certaine composante connexe~$U'$ de~$\pk\setminus\{\eta_1\}$. Comme~$U$ est géométriquement connexe,~$U'$ est géométriquement connexe. 

En tant que composante connexe de~$\pk\setminus\{\eta_1\}$, l'ouvert~$U'$ est l'image réciproque d'un point fermé~$\bf u$ de~$\PP^1_{\red k}$ par la flèche de réduction. Le corps~$\red k$ étant parfait,~$\kappa({\bf u})$ est une extension finie séparable de~$\red k$ ; si~$L$ désigne l'extension non ramifiée de~$k$ correspondant à une clôture galoisienne de~$\kappa({\bf u})$ alors la courbe~$U'\times_kL$ possède~$[\kappa({\bf u}):\red k]$ composantes connexes. Comme~$U'$ est géométriquement connexe,  il s'ensuit que~$\kappa({\bf u})=\red k$, et~$U'$ est de ce fait un~$k$-disque. 

Si~$U$ contient~$0$ ou~$\infty$ alors~$U$ est un~$k$-disque. Sinon,~$\psi$ est non ramifiée sur~$U$ et~$U\to U'$ est donc un revêtement fini étale, de degré 1 ou 2. S'il est de degré 2 alors~$U_{\KK}\to U'_{\KK}$ est un revêtement fini étale de degré~$2$ ; d'après la remarque~\ref{commdisc}), il est trivial, ce qui contredit la connexité géométrique de~$U$ ; par conséquent~$U\to U'$ est de degré 1, et~$U$ est un~$k$-disque. 

\medskip
Le seul point au-dessus duquel~$\phi$ est ramifié est le point rigide de~$\pk$ d'équation~$T^2-a=0$. S'il appartenait à~$U$ alors~$U$ serait égal à l'ouvert de~$\pk$ défini par la condition~$|T^2-a|<r$, sur lequel~$a$ est un carré en vertu du lemme de Hensel, ce qui contredirait le caractère géométriquement connexe de~$U$. Par conséquent,~$V\to U$ est fini étale de degré 1 ou 2. Si ce degré valait 1 alors~$V$ serait un disque, et posséderait en particulier un~$k$-point, ce qui est impossible puisque~${\sch X}(k)=\emptyset$. Le degré de~$V\to U$ vaut donc 2 ; comme~$U$ est un disque, le~$\ZZ/2\ZZ$-torseur étale~$V_{\KK}\to U_{\KK}$ est trivial (rem.~\ref{commdisc}), ce qui montre que~$V$ n'est pas géométriquement connexe.

\medskip
{\em Conclusion.} On déduit de ce qui précède et du~\ref{unetrigeco} que~$\{x\}$ est l'unique triangulation minimale, ou encore la plus petite triangulation, de la courbe~${\sch X}\an$.

\trois{casq2} Soit~${\sch X}$ la~$\QQ_2$-conique non triviale ; il existe un morphisme fini et plat~${\sch X}\to \PP^1_{\QQ_2}$ tel que~${\sch X}\times_{\PP^1_{\QQ_2}}\Aff^1_{\QQ_2}$ s'identifie à~$\spec \QQ_2[T,X]/(X^2+T^2+1)$. Pour tout~$r>0$ notons~$\DD_r$ le~$\CC_2$-disque centré en l'origine et de rayon~$r$. Le rayon~$2$-adique de la série~$\sqrt {1+T^2}$ coïncide avec~$\epsilon:=|2|$ ; par conséquent, si~$r\in ]\epsilon;1[$, le revêtement étale de~$\DD_r$ obtenu par extraction d'une racine carrée de~$1+T^2$ est connexe ; il en va de même de celui obtenu par extraction d'une racine carrée de~$-1\setminus T^2$, puisque~$(-1)$ est un carré dans~$\CC_2$ ; autrement dit,~${\sch X}\an_{\CC_2}\times_{\PP^{1,{\rm an}}_{\CC_2}}\DD_r$ est connexe. 

\medskip
Soit~$r\in ]\epsilon ;1[$ ; l'image réciproque du point~$\eta_{r,\CC_2}$ de~$\PP^{1,{\rm an}}_{\CC_2}$ sur~${\sch X}\an_{\CC_2}$ est alors un singleton, par un argument dû à de Jong. Supposons en effet que ce ne soit pas le cas ; le revêtement~${\sch X}\an_{\CC_2}\to \PP^{1,{\rm an}}_{\CC_2}$ serait dès lors trivial au-dessus de~$\eta_{r,\CC_2}$, donc au voisinage de~$\eta_{r,\CC_2}$, et {\em a fortiori} au-dessus d'une couronne décrite par une inégalité de la forme~$r<|T|<s$ pour un certain~$s>r$. 

En recollant~${\sch X}\an_{\CC_2}\times_{ \PP^{1,{\rm an}}_{\CC_2}}\DD_s$ avec le revêtement trivial de l'ouvert de~$\PP^{1,{\rm an}}_{\CC_2}$ décrit par l'inégalité~$|T|>r$, on obtiendrait un revêtement fini étale de~$\PP^{2,{\rm an}}_{\CC_2}$, nécessairement trivial en vertu de GAGA  ; mais l'image réciproque de~$\DD_r$ sur ce revêtement s'identifierait à~${\sch X}\an_{\CC_2}\times_{\PP^{1,{\rm an}}_{\CC_2}}\DD_r$ qui est connexe, d'où une contradiction.

\medskip
Par conséquent, si~$r\in  ]\epsilon ;1[$, le point~$\eta_r$ de~$\PP^{1,{\rm an}}_{\QQ_2}$ a un et un seul antécédent~$x_r$ sur~${\sch X}\an$, et~$\got s(x_r)=k$ ; ainsi~${\sch X}\an$ possède-t-elle une infinité de triangulations minimales. 

\subsection*{Courbes compactes de squelette analytique vide}

\deux{gammavide} Soit~$X$ une courbe~$k$-analytique compacte, non vide et géométriquement connexe ; on fait de plus l'hypothèse que le squelette analytique de~$X$ est vide. On déduit des assertions i) et ii) du théorème~\ref{theoanaladm} que~$X$ est lisse ; étant par ailleurs compacte,~$X$ est propre et lisse et s'identifie donc à~$\sch X\an$ pour une certaine~$k$-courbe algébrique projective, lisse et géométriquement intègre~$\sch X$. 

\deux{gammvidealgclos} {\em Étude de~${\sch X}\an_{\KK}$.} Comme tout point de~${\sch X}\an$ a un voisinage qui est un disque virtuel, tout point de~${\sch X}\an_{\KK}$ a un voisinage qui est un disque, et en particulier un arbre à un bout ; il en résulte que~${\sch X}\an_{\KK}$ est un arbre compact.  Soit~$\ell$ un entier premier à~$p$ et soit~$Y\to {\sch X}\an_{\KK}$ un~$\ZZ/\ell \ZZ$-torseur étale. Si~$x\in {\sch X}\an_{\KK}$, il possède un voisinage qui est un disque, et par conséquent déploie~$Y$ (rem.~\ref{commdisc}) ; étant déployé au voisinage de tout point de~${\sch X}\an_{\KK}$, le~$\ZZ/\ell \ZZ$-torseur étale~$Y$ est en réalité un~$\ZZ/\ell \ZZ$-torseur {\em topologique} ; comme~${\sch X}_{\KK}\an$ est un arbre,~$Y$ est trivial. Combiné à GAGA, cela implique que~$\H^1({\sch X}_{\KK, \tiny\mbox{ét}},\ZZ/\ell \ZZ)=0$ ; par conséquent,~${\sch X}\an_{\KK}\simeq \pkk$ et~$\sch X$ est de genre~$0$. 

\medskip
\deux{equivsquelvide} Il en résulte, compte-tenu de ~\ref{gzeroptinv},~\ref{ptinvtrimin}-\ref{discvirtk},~\ref{casdeuxtrimin} et~\ref{conpimpair}, que si~$X$ est  une courbe~$k$-analytique compacte, non vide et géométriquement connexe les propriétés suivantes sont équivalentes : 

\medskip
i) le squelette analytique de~$X$ est vide ; 

ii) la courbe~$X$ est isomorphe à l'analytification d'une courbe algébrique projective, lisse et géométriquement intègre de genre~$0$, et il existe deux points distincts~$x$ et~$x$' sur~$X\dtr$ tels que~$\got s(x)=\got s(x')=k$ ;

iii) la courbe~$X$ est isomorphe à l'analytification d'une courbe algébrique projective, lisse et géométriquement intègre de genre~$0$ et possède un ouvert qui est un disque virtuel géométriquement connexe.

\trois{recapinftri} Si ces propriétés sont satisfaites, on déduit de~\ref{triminptinv} et~\ref{discvirtk} : que les triangulations minimales de~$X$ sont exactement les singletons de la forme~$\{x\}$, où~$x$ est de type 2 ou 3 et vérifie l'égalité~$\got s(x)=k$ ; et qu'il existe une infinité de telles triangulations et même, si~$|k\ti|\neq \{1\}$, une infinité de telles triangulations dont l'unique sommet est de type 2.

\trois{inftrip1} Lorsque~$\red k$ est parfait et de caractéristique différente de~$2$, il découle de l'exemple traité au~\ref{conpimpair} que les trois propriétés équivalentes ci-dessus sont satisfaites si et seulement si~$X\simeq \pk$ ; c'est également le cas, de façon évidente, lorsque que~$k$ est algébriquement clos. 

\trois{inftriconq2} L'exemple~\ref{casq2} assure que si~$k=\QQ_2$ et si~$X$ est l'analytifiée de la~$k$-conique projective non triviale (d'équation homogène~$T_0^2+T_1^2+T_2^2=0$) alors~$X$ satisfait  les trois propriétés équivalentes ci-dessus.

\deux{remgranaladmp} {\em Remarque.} Si~$\Gamma$ est un sous-graphe de~$\pkk$ alors~$\Gamma$ est analytiquement admissible si et seulement si~$\Gamma$ est un sous-arbre non vide ou bien non réduit à un point, ou bien de la forme~$\{x\}$ avec~$x$ de type 2 ou 3. 

\medskip
En effet, supposons~$\Gamma$ analytiquement admissible ; il est alors admissible et est donc un sous-arbre non vide de~$\pkk$ ; et s'il est réduit à un point~$x$, il résulte de~\ref{borddisct23} {\em et sq.} que~$x$ est de type 2 ou 3. 

\medskip
Réciproquement, si~$\Gamma$ est du type requis, alors il est admissible et contient au moins un point~$x$ de type 2 ou 3 (lorsque~$\Gamma$ n'est pas réduit à un point, il suffit de prendre pour~$x$ n'importe quel point pluribranche de~$\Gamma$). Mais~$\{x\}$ est un sous-arbre analytiquement admissible de~$\pkk$, ce qui implique que~$\Gamma$ lui-même est analytiquement admissible.

\subsection*{Les nœuds du squelette analytique} 

\deux{unibrskelan} Soit~$X$ une courbe~$k$-analytique quasi-lisse, soit~$x$ un point unibranche de~$\skelan X$ et soit~$b$ l'unique élément de~$\br \Gamma x$. Supposons que~$x$ appartient à l'intérieur analytique de~$X$ et qu'il est de genre zéro ; nous allons montrer que~$\got s(b)$ contient {\em strictement}~$\got s(x)$. Pour le voir, on choisit une triangulation~$S$ de~$X$ de squelette~$\skelan X$ et telle que~$S$ contienne strictement~$x$. La composante connexe~$\Gamma_0$ de~$x$ dans~$\Gamma-(S\setminus\{x\})$ est alors une demi-droite issue de~$x$. Si l'on avait~$\got s(b)=\got s(x)$, alors~$\Gamma_0^\flat$ serait un disque virtuel en vertu de~\ref{casgammazeroun}, contredisant ainsi l'appartenance de~$x$ à~$\skelan X$.

\deux{noeudsinclus} Soit~$X$ une courbe~$k$-analytique génériquement quasi-lisse et soit~$\Sigma$ l'ensemble des nœuds de~$\skelan X$ ; soit~$\Gamma$ un sous-graphe localement fini et analytiquement admissible de~$X$, tracé sur~$X\geom$. Nous allons montrer que~$\Sigma$ est contenu dans l'ensemble des nœuds de~$\Gamma$. Pour le voir, on peut supposer~$X$ connexe, auquel cas~$\Gamma$ et~$\skelan X$ le sont aussi, ce qui entraîne que~$\skelan X$ n'a aucun point isolé sauf si lui-même est un singleton. On déduit alors de~\ref{remnoeudspasinclus} et du~\ref{unibrskelan} ci)-dessus que l'assertion à démontrer est évidente sauf le cas où le graphe~$\skelan X$ est un singleton~$\{x\}$ avec~$x$ appartenant à~$X\dtr \setminus \partial  X\an$ et de genre zéro (notons que~$x$ est alors l'unique nœud de~$\skelan X$) ; nous nous plaçons donc sous cette hypothèse. La notion de nœud ne dépendant pas du fait que l'on voie~$X$ comme une courbe sur~$k$ ou sur~$\got s(x)$, on peut remplacer~$k$ par~$\got s(x)$, c'est-à-dire faire l'hypothèse que~$X$ est {\em géométriquement} connexe.

Comme~$X$ est connexe et comme~$\skelan X\neq \emptyset$, ce dernier est un sous-graphe analytiquement admissible de~$X$ ; sa compacité entraîne alors celle de~$X$. Comme~$\skelan X\subset X\dtr\setminus \partial \an X$, il résulte des assertions 1) et 2) du théorème~\ref{theoanaladm} que~$\partial \an X=\emptyset$ et que~$X$ est quasi-lisse ; la courbe compacte~$X$ est dès lors propre et lisse,  et s'identifie donc à l'analytification~$\sch X\an$ d'une~$k$-courbe algébrique projective, lisse et géométriquement intègre~$\sch X$. 

\medskip
L'image réciproque de~$x$ sur~$\sch X\an_{\KK}$ est un sous-graphe analytiquement admissible de~$\sch X\an_{\KK}$, et est en particulier connexe puisque~$\sch X\an_{\KK}$ est connexe. Cela signifie que~$\got s(x)=k$, et que l'ensemble des antécédents de~$x$ sur~$\sch X\an _{\KK}$ est un singleton~$x_{\KK}$ ; le sous-graphe~$\{x_{\KK}\}$ de~$\sch X\an _{\KK}$ étant analytiquement admissible, toutes les composantes connexes de~$\sch X\an_{\KK}\setminus\{x_{\KK}\}$ sont des disques. Par ailleurs,~$x$ est par hypothèse de genre~$0$, ce qui signifie {\em par définition} que~$x_{\KK}$ est de genre~$0$. Il s'ensuit (lemme~\ref{lemrevg0}) que~$\sch X_{\KK}\simeq \PP^1_{\KK}$, c'est-à-dire que~$\sch X$ est de genre 0. 

\medskip
Le squelette de~$\sch X\an$ étant le singleton~$\{x\}$, il résulte de~\ref{casdeuxtrimin} et~\ref{casunetrimin} que toute composante connexe de~$\sch X\an\setminus\{x\}$ est un disque virtuel {\em qui n'est pas géométriquement connexe.} 

Soit~$b\in \br {\sch X\an}x$ et soit~$V$ la composante connexe de~$\sch X\an\setminus\{x\}$ qui contient~$b$. Par ce qui précède,~$V$ n'est pas géométriquement connexe. En vertu de~\ref{remsxsv}, on a~$\got s(b)=\got s(V)$ ; par conséquent,~$\got s(b)$ contient strictement~$\got s(x)=k$. 

\medskip
Si~$\Gamma=\{x\}$ alors~$x$ est un nœud de~$\Gamma$. Sinon, il existe~$b\in \br \Gamma x$ (par connexité de~$\Gamma$). Par ce qui précède~$\got s(b)$ contient strictement~$\got s(x)$ et~$x$ est là encore un nœud de~$\Gamma$, ce qui achève la démonstration. 

\subsection*{Courbes compactes de squelette analytique non vide}

\deux{gammapasvide} Soit~$X$ une courbe~$k$-analytique que l'on suppose compacte, non vide, quasi-lisse et géométriquement connexe ; on fait de plus l'hypothèse que le squelette analytique~$\Gamma$ de~$X$ est non vide ; soit~$\Sigma$ l'ensemble des n\oe uds de~$\Gamma$. 

\trois{noeudsnonvide} {\em Supposons que~$\Sigma\neq \emptyset$.} Le graphe~$\Gamma$ est connexe, fini et compact, et l'ensemble~$\Sigma$ est non vide et contient tous les sommets de~$\Gamma$ ; le complémentaire de~$\Sigma$ dans~$\Gamma$ est donc réunion disjointe d'intervalles ouverts ; par conséquent,~$\Sigma$ est une triangulation de~$X$. Comme~$\Sigma$ est contenu d'après le~\ref{noeudsinclus} ci-dessus dans l'ensemble des nœuds de n'importe quel sous-graphe localement fini et analytiquement admissible de~$X$ tracé sur~$X\geom$, il est contenu dans toute triangulation de~$X$ et est donc la plus petite triangulation de~$X$. Remarquons que si~$|k\ti|\neq\{1\}$ et si~$X$ est strictement~$k$-analytique alors~$\Sigma$ est constitué uniquement de points de type 2. 

\trois{noeudsvide} {\em Supposons que~$\Sigma=\emptyset$.} Compte-tenu de la définition des n\oe uds de~$\Gamma$, cela équivaut à la conjonction des propriétés suivantes :

\medskip

$\bullet$~$X$ est sans bord, ce qui revient à demander qu'elle soit isomorphe à l'analytification d'une~$k$-courbe algébrique projective, lisse et géométriquement intègre~$\sch X$ ; 

$\bullet$~$\Gamma$ est sans sommet, ce qui revient à demander qu'il soit un cercle (et implique qu'il coïncide avec le squelette de~$X$) ; 

$\bullet$ tous les points de~$\Gamma$ sont de genre~$0$ ; 

$\bullet$ si~$x\in \Gamma$ et si~$b_0$ et~$b_1$ sont les deux branches de~$\Gamma$ issues de~$x$ alors~$\got s(b_0)=\got s(x)=\got s(b_1)$ ; il revient au même de demander que~$x\mapsto [\got s(x):k]$ soit constante sur~$\Gamma$, ou que~$\got s(x)\simeq \got s(y)$ pour tout~$(x,y)\in \Gamma^2$ ({\em cf.}~\ref{defnoeudsg}). 

\medskip
Notons maintenant quelques conséquences de ces hypothèses.  

\medskip
{\em La valeur absolue de~$k$ est non triviale.} En effet si elle l'était le squelette analytique de~${\sch X}\an$ serait, en vertu de la description explicite de cette dernière, réduit à un singleton. 

\medskip
{\em Les triangulations minimales de~$X$ sont exactement les singletons~$\{x\}$ avec~$x\in \Gamma$}. En effet si~$x\in \Gamma$ alors~$\{x\}$ contient l'ensemble des n\oe uds de~$\Gamma$ (lequel est vide !), et~$\Gamma\setminus\{x\}$ est un intervalle ouvert relativement compact dans~$\Gamma$ ; par conséquent,~$\{x\}$ est une triangulation de~$X$, évidemment minimale. 

Réciproquement, soit~$S$ une triangulation de~$X$. Toute composante connexe de~$X\setminus S$ est alors un arbre, ce qui exclut qu'une telle composante puisse contenir~$\Gamma$ ; par conséquent,~$S$ rencontre~$\Gamma$, d'où notre assertion. 

Remarquons que comme~$|k\ti|\neq \{1\}$ il y a une infinité de points de type 2 sur~$\Gamma$, et donc une infinité de triangulations minimales de~$X$ dont l'unique sommet est de type 2. 

\medskip
{\em Si~$I$ est un intervalle ouvert non vide tracé sur~$\Gamma$ alors~$I\an$ est une couronne gentiment virtuelle}. Cela résulte simplement du~\ref{pretriang}. 

\medskip
{\em Si~$x$ et~$y$ sont deux points distincts de~$\Gamma$ le choix d'une orientation sur~$\Gamma$ induit un~$k$-isomorphisme~$\got s(x)\simeq \got s(y)$.} En effet, choisissons un intervalle ouvert~$I$ de~$\Gamma$ contenant~$x$ et~$y$ et tel que le sous-intervalle~$[x;y]$ de~$I$ soit orienté dans le sens direct. Le corps~$\got s(I\an)$ est indépendant du choix de~$I$  à isomorphisme canonique près, et comme~$I$ est le squelette de la couronne virtuelle~$I\an$ les flèches~$\got s(I\an)\hookrightarrow \got s(x)$ et~$\got s(I\an)\hookrightarrow \got s(y)$ sont des isomorphismes.

\medskip
{\em Extensions des scalaires.} Soit~$F$ une extension presque algébrique de~$k$. Comme~$\Gamma$ est le squelette analytique de~$X$, il est tracé sur~$X\dtr$ et chacun de ses points a un nombre fini d'antécédents sur~$X_{\KK}$. Il en découle (corollaire~\ref{corollimrecih}) que si~$\Gamma_F$ désigne l'image réciproque de~$\Gamma$ sur~$X_F$ alors~$\Gamma_F$ est un sous-graphe compact et fini de~$X_F$, et que~$\Gamma_F\to \Gamma$ est injective par morceaux.

Si~$x\in\Gamma$, le nombre~$N$ d'idéaux maximaux de~$F\hotimes_k\hres(x)$  est égal au nombre d'idéaux maximaux de~$F\hotimes_k\got s(x)$ et ne dépend donc pas de~$x$. Il s'ensuit que~$\Gamma_F\to \Gamma$ est un revêtement topologique à~$N$ feuillets, ce qui entraîne que~$\Gamma_F$ est une réunion disjointe de cercles. Comme toutes les composantes connexes de~$X-\Gamma$ sont des disques virtuels, toutes les composantes connexes de~$X_F-\Gamma_F$ sont des disques virtuels, ce qui montre que~$\Gamma_F$ est un sous-graphe analytiquement admissible de~$X_F$. Il est en particulier admissible ; par conséquent~$X_F$ se rétracte sur~$\Gamma_F$ et celui-ci est donc connexe, et partant homéomorphe à un cercle ; le graphe~$\Gamma_F$ est dès lors le squelette de~$X_F$, et même son squelette analytique puisqu'il est analytiquement admissible. 

\medskip
Soit maintenant~$\mathsf H$ le noyau de l'opération naturelle de~$\mathsf G$ sur~$\Gamma_{\KK}$, qui est {\em a priori} un sous-groupe fermé de~$\mathsf G$. Soit~$\mathsf C$ le quotient~$\mathsf G/\mathsf H$, quotient que l'on voit comme un groupe de~$\Gamma$-homéomorphismes de~$\Gamma_{\KK}$. Comme~$\Gamma_{\KK}\to \Gamma$ est un revêtement fini entre cercles,~$\mathsf C$ est un groupe fini de «rotations» (et~$\mathsf H$ est ouvert) ; comme~$\mathsf C$ agit transitivement sur les fibres de~$\Gamma_{\KK} \to \Gamma$, et comme~$\mathsf{Aut}\;\Gamma_{\KK}/\Gamma$ est cyclique,~$\mathsf C=\mathsf{Aut}\;\Gamma_{\KK}/\Gamma$. Si~$y\in \Gamma_{\KK}$ son stabilisateur dans~$\mathsf C$ est trivial ; son stabilisateur dans~$\mathsf G$ est donc~$\mathsf H$. 

\medskip
Soit~$L$ la sous-extension finie de~$k^s$ correspondant à~$\mathsf H$. Il résulte de ce qui précède : que pour tout~$x\in \Gamma$ la~$k$-extension~$\got s(x)$ est (non canoniquement) isomorphe à~$L$ ; que son groupe de Galois est canoniquement isomorphe à~${\mathsf{Gal}}\;(L/k)$ ; que~$\Gamma_{\KK}\to \Gamma_L$ est un homéomorphisme et que~$\Gamma_L\to \Gamma$ est un revêtement topologique dont~${\mathsf{Gal}}\;L/k$ s'identifie naturellement au groupe d'automorphismes. Notons que de ce fait,~${\mathsf{Gal}}\;L/k$ possède une paire de générateurs distingués (qui correspondent aux deux rotations d'angle minimal sur~$\Gamma_L$) et même un générateur distingué si~$\Gamma$ a été préalablement orienté ; il en va de même de~${\mathsf{Gal}}\;(\got s(x)/k)$ pour tout~$x\in \Gamma$, en vertu de l'isomorphisme canonique ~${\mathsf{Gal}}\;(\got s(x)/k)\simeq {\mathsf{Gal}}\;L/k$.

\medskip
Comme~$\Gamma_{\KK}$ ne comporte que des points de genre~$0$ (puisque c'est le cas de~$\Gamma$), il ne possède aucun n\oe ud. Recouvrons-le par deux intervalles ouverts~$I$ et~$I'$ dont l'intersection est réunion disjointe de deux intervalles ouverts~$J$ et~$J'$. Les ouverts~$I\an, (I')\an, J\an$ et~$(J')\an$ de~$X_{\KK}$ sont des~$\KK$-couronnes (\ref{pretriang}) ; et~$J\an$ et~$(J')\an$ apparaissent tous deux comme des sous-couronnes de~$I\an$ aussi bien que de~$(I')\an$. La suite exacte de Mayer-Vietoris combinée à la proposition~\ref{corollkum} assure alors que~$\H^1(X_{\KK},\ZZ/\ell \ZZ)\simeq (\ZZ/\ell \ZZ)^2$ pour n'importe quel entier~$\ell$ premier à~$p$. 

Par conséquent, et en vertu de GAGA,~$\sch X$ est de genre~$1$ et est donc un espace principal homogène sous sa jacobienne~$\sch Y$, laquelle est une courbe elliptique. Le squelette de~$X_{\KK}$ étant un cercle,~${\sch Y}\an_{\KK}$ est une courbe de Tate, c'est-à-dire qu'elle est isomorphe, comme groupe~$\KK$-analytique, à un quotient de la forme~${\mathbb G}_{m,\KK}\an/q^{\ZZ}$ pour un certain~$q$ non nul de~$\KK$ tel que~$|q|<1$ ; cela revient à dire que le~$j$-invariant de~$\sch Y$ n'est pas entier. 

\medskip
 {\em Remarque.} Soit~$V$ une composante connexe de~$X-\Gamma$ ; c'est un disque virtuel dont le bord~$\{x\}$ est situé sur~$\Gamma$. Comme~$\got s(x)$ se plonge dans~$\got s(V)$ (\ref{remsxsv}), l'extension~$L$ se plonge dans~$\got s(V)$. Par conséquent, pour tout~$y\in X$ l'extension~$L$ se plonge dans~$\got s(y)$. Il s'ensuit notamment que si~$L\neq k$, alors~$X(k)$ est vide ; le~$\sch Y$-torseur~$\sch X$ est dans ce cas non trivial. Il s'ensuit aussi que~$X_L\to X$ est un revêtement topologique ; comme~$\Gamma\hookrightarrow X$ est une équivalence d'homotopie, ce revêtement est cyclique, et son groupe d'automorphismes s'identifie à~${\mathsf{Gal}}\;L/k$.

\deux{reciptwtate} Réciproquement, soit~$\sch X$ une~$k$-courbe algébrique projective, lisse et géométriquement intègre de genre~$1$, dont le~$j$-invariant de la jacobienne~$\sch Y$ n'est pas entier. Dans ce cas le squelette analytique de~${\sch X}\an_{\KK}$ coïncide avec son squelette, et est un cercle dont tous les points sont de genre~$0$ ; si~$\mathsf G$ opère sur le cercle en question par rotations, on vérifie immédiatement que la courbe~$X:={\sch X}\an$ vérifie les conditions décrites au~\ref{noeudsvide}.

\subsection*{Deux exemples en genre 1} 

\deux{extateprecis} Soit~$n$ un entier impair et premier à~$p$ ; supposons que~$k$ contient toutes les racines~$n$-ièmes de l'unité, et qu'il existe un élément~$q$ de~$k\ti$ tel que~$|q|<1$ et tel que~$T^n-q$ soit irréductible dans~$k[X]$ ({\em e.g}.~$k=\QQ_7, n=3, q=7$). Soit~$L$ l'extension cyclique de degré~$n$ de~$k$ contenue dans~$\KK$ et engendrée par une racine~$n$-ième~$\lambda$ de~$q$. Soit~$Y$ le quotient de~${\mathbb G}_{m,k}\an$ par~$q^\ZZ$ ; la courbe~$Y$ est isomorphe à l'analytification d'une~$k$-courbe elliptique~$\sch Y$ ; le~$j$-invariant de~$\sch Y$ n'est pas entier. 

\medskip
La classe~$\overline \lambda$ de~$\lambda$ dans~$Y(L)\simeq L\ti/q^\ZZ$ est de~$n$-torsion ; comme~$n$ est impair, le produit de toutes les racines~$n$-ièmes de l'unité de~$L$ est égal à~$1$ ; il s'ensuit que~$\overline \lambda$ est de norme~$1$. Donnons-nous un générateur~$g$ du groupe de Galois de~$L/k$ ; il existe alors un unique~$1$-cocycle~$\mathsf c$ de~${\mathsf{Gal}}\;(L/k)$ à valeurs dans~$Y(L)={\sch Y}(L)$ qui envoie~$g$ sur~$\overline \lambda$. La classe de~$\H^1(k,{\sch Y})$ définie par~$\mathsf c$ est celle d'un~$\sch Y$-torseur~$\sch X$ qui est déployé par~$L$. 

\medskip
L'espace analytique~${\sch X}\an$ s'identifie au quotient de~${\sch X}\an_{\KK}$ sous l'action de ~$\mathsf G$, c'est-à-dire encore au quotient de~${\sch Y}\an_{\KK}\simeq {\mathbb G}_{m,\KK}\an/q^\ZZ$ sous l'action de Galois tordue par la restriction de~$\mathsf c$ à~$\mathsf G$. Soit~$g'\in \mathsf G$ et soit~$N$ tel que la restriction de~$g'$ à~$L$ coïncide avec~$g^N$. On déduit de la description explicite de~$\mathsf c$ que la restriction au squelette de~${\sch Y}\an$, identifié à~$\RR\ti_+/|q|^\ZZ$, de l'action {\em tordue} de~$g'$, est la translation par~$\sqrt [n]{|q|^N}$ (qui est triviale si et seulement si~$n|N$, c'est-à-dire si et seulement si~$g'$ agit trivialement sur~$L$). Par conséquent,~${\sch X}$ est de la forme décrite au~\ref{reciptwtate}, et les résultats du~\ref{noeudsvide} s'appliquent ; on en déduit les faits qui suivent.

\medskip
Le squelette de~${\sch X}\an$ coïncide avec son squelette analytique, et est un cercle~$\Gamma$. Si~$y\in X$ il existe un~$k$-plongement (non canonique)~$L\hookrightarrow \got s(y)$ qui est un isomorphisme si~$y\in \Gamma$. 

Si~$F$ est une extension presque algébrique de~$k$ l'image réciproque~$\Gamma_F$ de~$\Gamma$ sur~${\sch X}\an_F$ est un cercle, qui n'est autre que le squelette, et même le squelette analytique, de~${\sch X}\an_F$ ; et la flèche ~$\Gamma_F\to \Gamma$ est un revêtement topologique. 

Le groupe~${\mathsf{Gal}}\;(L/k)$ s'identifie naturellement au groupe d'automorphismes du revêtement~$\Gamma_L\to \Gamma$, et~$X_L\to X$ est lui-même un revêtement topologique cyclique de groupe~${\mathsf{Gal}}\;(L/k)$ ; par ailleurs,~$\Gamma_{\KK}\to \Gamma_L$ est un homéomorphisme. 

\deux{tatetwist} Supposons que~$p\neq 2$, que~$|k\ti|\neq\{1\}$ et qu'il existe un élément~$a$ de~$k^0$ tel que~$\red a$ ne soit pas un carré dans~$\red k$ ; soit~$q\in k\ti$ tel que~$|q|<1$ ({\em e.g}~$k=\QQ_3, a=-1, q=3$). Soit~$\sch X$ la~$k$-courbe elliptique donnée par l'équation affine~$T_1^2=T_0(T_0-q)(T_0+a)$ ; la fonction~$T_0$ définit un morphisme fini et plat de degré 2 de~$\sch X$ vers~$\PP^1_k$. Soit~$\mathsf H$ le sous-groupe de~$\mathsf G$ fixant les racines carrées de~$a$. 

\medskip
Le~$j$-invariant de~$\sch X$ n'est pas entier ; la courbe~${\sch X}\an_{\KK}$ est donc une courbe de Tate. On vérifie aisément que si~$|q|<r<1$ alors~$\eta_{r,\KK}$ a deux antécédents sur~${\sch X}\an_{\KK}$ qui sont fixés par~$\mathsf H$, et que si~$r=|q|$ ou~$r=|1|$ alors~$\eta_{r,\KK}$ a un unique antécédent sur~${\sch X}\an_{\KK}$, nécessairement invariant par~$\mathsf G$. Il en résulte que l'image réciproque~$\Gamma$ de l'intervalle~$[\eta_{|q|,\KK};\eta_{1,\KK}]$ sur~${\sch X}\an_{\KK}$ est un cercle ; c'est le squelette, et même le squelette analytique, de~${\sch X}\an_{\KK}$. Notons~$x$ et~$y$ les antécédents respectifs de~$\eta_{|q|,\KK}$ et~$\eta_{1,\KK}$ sur~$\Gamma$, et~$\xi$ et~$\eta$ les images respectives de~$x$ et~$y$ sur~${\sch X}\an$. On peut également définir~$\xi$ (resp.~$\eta$) comme l'unique antécédent du point~$\eta_{|q|}$  sur~${\sch X}\an$. 

\medskip
Comme l'ensemble des nœuds de~$\Gamma$ est vide,~$\{x,y\}$ est une triangulation de~${\sch X}\an_{\KK}$ de squelette~$\Gamma$. Le groupe~$\mathsf G$ agit sur~$\Gamma$ {\em via} son quotient~$\mathsf G/\mathsf H$ en permutant les deux composantes connexes de~$\Gamma\setminus\{x,y\}$ et en fixant~$x$ et~$y$. Il s'ensuit que~${\sch X}\an$ est un arbre, et l'on déduit du théorème~\ref{theo-invariance-triangulations} que~$\{\xi,\eta\}$ est une triangulation de~${\sch X}\an$ de squelette~$[\xi;\eta]$. La considération de l'action de Galois sur~$\Gamma$ montre que~$\got s(\xi)=\got s(\eta)=k$ et que~$\got s(\zeta)\simeq k[\tau]/(\tau^2-a)$ pour tout~$\zeta\in ]\xi;\eta[$. 

Il découle alors de~\ref{recapitabl} que ni~$\{\xi\}$, ni~$\{\eta\}$ ne sont des triangulations de~${\sch X}\an$ ; par conséquent,~$\{\xi,\eta\}$ est une triangulation minimale de~${\sch X}\an$. Comme celle-ci comprend deux sommets, on est dans le cas traité au~\ref{noeudsnonvide} : l'ensemble~$\{\xi,\eta\}$ est la plus petite triangulation de~${\sch X}\an$ ; son squelette~$[\xi;\eta]$ est le squelette analytique de~${\sch X}\an$ ; et~$\{\xi,\eta\}$ est exactement l'ensemble des nœuds de~$[\xi;\eta]$.

\subsection*{Récapitulation} 

\deux{introrecaptriptt} On se donne une courbe~$k$-analytique géométriquement connexe, non vide, quasi-lisse et compacte~$X$. Soit~$\Sigma$ l'ensemble des nœuds de~$\skelan X$. On est dans l'un des trois cas suivants, exclusifs l'un de l'autre. 

\trois{sigmapasvd} {\em Le cas où~$\Sigma$ est non vide.} Remarquons que comme tout point de~$\partial \an X$ appartient ~$\skelan X$, et en est dès lors un nœud par définition,~$\Sigma$ est notamment non vide dès que~$\partial \an X\neq \emptyset$. 

\trois{sigmavidmais} {\em Le cas où~$\skelan X$ est non vide mais où~$\Sigma$ est vide.} En vertu de~\ref{noeudsvide}, il se produit si et seulement si~$X$ est isomorphe à l'analytification d'une~$k$-courbe algébrique, projective et lisse de genre 1 et possède la propriété suivante :~$\skelan X$ est un cercle (qui coïncide forcément avec~$\skel X$), et il existe une extension finie cyclique~$L$ de~$k$ telle que~$\got s(x)\simeq L$ pour tout~$x\in \skelan X$. Notons que si~$k$ est algébriquement clos cela signifie simplement que~$X$ est une courbe de Tate ; pour un exemple non trivial sur un corps non algébriquement clos, {\em cf.}~\ref{extateprecis}.

\trois{skelanmemevid} {\em Le cas où~$\skelan X$ est lui-même vide.} En vertu de~\ref{equivsquelvide}, il se produit si et seulement si~$X$ est isomorphe à l'analytification d'une~$k$-courbe algébrique, projective et lisse de genre 0 et possède la propriété suivante : il existe deux points distincts~$x$ et~$x'$ de~$X\typ{23}$ tels que~$\got s(x)=\got s(x')=k$. Si~$k$ est algébriquement clos ou si~$\red k$ est parfait de caractéristique différente de 2 cela signifie simplement que~$X\simeq \pk$ (\ref{inftrip1}) ; mais mentionnons que la conique non triviale sur~$\QQ_2$ satisfait également ces hypothèses (\ref{casq2}).  

\deux{recaptriplatt} Toute triangulation de $X$ ayant un nombre fini de sommets, 
il est évident \textit{a priori}
que $X$ possède des triangulations minimales. 

\medskip
Sous l'hypothèse de~\ref{sigmapasvd}, elle possède une plus petite triangulation, à savoir~$\Sigma$ (\ref{noeudsnonvide}). 

\medskip
Sous l'hypothèse de~\ref{sigmavidmais}, qui ne peut être satisfaite que si~$|k\ti|\neq \{1\}$, les triangulations minimales de~$X$ sont les singletons~$\{x\}$ avec~$x\in \skelan X$ (\ref{noeudsvide}) ; elle en possède donc une infinité, et même une infinité dont l'unique sommet est de type 2. Elle ne possède donc pas de plus petite triangulation. 

\medskip
Sous l'hypothèse de~\ref{skelanmemevid},  les triangulations minimales de~$X$ sont les singletons~$\{x\}$ où~$x$ est un point de~$X\typ{23}$ tel que~$\got s(x)=k$ ; la courbe~$X$ en possède une infinité, et même une infinité  dont l'unique sommet est de type 2 si~$|k\ti|\neq\{1\}$. Elle ne possède donc pas de plus petite triangulation. 

\deux{equivneoudsppt} Il découle de ce qui précède que~$X$ possède une plus petite triangulation {\em si et seulement si}~$\Sigma$ est non vide, et dans ce cas cette plus petite triangulation coïncide avec~$\Sigma$. 
\section{Raffinement : le cas des courbes marquées}

\subsection*{Triangulations adaptées}

\deux{gammasigmanonv} Soit~$X$ une courbe~$k$-analytique {\em quasi-lisse, connexe et compacte}. Soit~$\Gamma$ un sous-graphe localement fini et analytiquement admissible de~$X$ tracé sur~$X\geom$, soit~$\Sigma$ l'ensemble des nœuds de~$\Gamma$ et soit~$\Sigma\dtr$ le sous-ensemble de~$\Sigma$ formé des points  de type 2 ou 3. {\em On suppose que~$\Sigma\dtr$ est non vide}.  Nous allons montrer que sous cette hypothèse,~$\Sigma\dtr$ est une triangulation de~$X$, et que pour tout~$x\in \Sigma-\Sigma\dtr$, la composante connexe de~$X-\Sigma\dtr$ qui contient~$x$ est un disque virtuel. 

Pour ce faire, on fixe une composante connexe~$U$ de~$X-\Sigma\dtr$, dont nous allons établir qu'elle est un disque virtuel ou une couronne virtuelle (la relative compacité sera ici automatique,~$X$ étant compact), et que si elle rencontre~$\Sigma$ c'est obligatoirement un disque virtuel. Deux cas sont à distinguer. 

\trois{siunerenc} {\em Le cas où~$U$ ne rencontre pas~$\Gamma$.} C'est alors une composante connexe de~$X-\Gamma$, et partant un disque virtuel. 

\trois{siurenc} {\em Le cas où~$U$ rencontre~$\Gamma$}. La composante~$U$ est alors de la forme~$I^\flat$, où~$I$ est une composante connexe de~$\Gamma-\Sigma\dtr$. Comme~$X$ est compacte,~$I$ est relativement compacte ; comme~$X$ est connexe,~$\Gamma$ est connexe ; il s'ensuit,~$\Sigma\dtr$ étant non vide, que~$\partial I$ est lui-même une partie non vide de~$\Sigma\dtr$. Comme~$U\cap \Gamma=I$, la composante~$U$ rencontre~$\Sigma$ si et seulement si~$I$ rencontre~$\Sigma$. On doit maintenant distinguer deux sous-cas. 

\medskip
{\em Supposons que~$I$, et donc~$U$, ne rencontrent pas~$\Sigma$.} Dans ce cas,~$I$ ne contient en particulier aucun sommet (topologique) de~$\Gamma$, et est donc de valence 2 en chacun de ses points ; puisque~$\partial I$ est non vide,~$I$ est un intervalle ouvert. En vertu de~\ref{pretriang}, l'ouvert~$U=I^\flat$ est alors une couronne virtuelle. 

\medskip
{\em Supposons que~$I$, et donc~$U$, rencontrent~$\Sigma$.} La courbe~$X$ étant quasi-lisse, elle est normale et chacun de ses points rigides est donc unibranche. 

Puisque~$I$ ne rencontre pas~$\Sigma\dtr$, les seuls sommets de~$\Gamma$ éventuellement situés sur~$I$ sont rigides, et sont par conséquent des points unibranches de~$I$ ; il en découle que la valence de~$I$ en chacun de ses points vaut~$1$ ou~$2$, ce qui implique,~$\partial I$ étant non vide, que le graphe connexe~$I$ est un intervalle ouvert ou semi-ouvert. Or~$I$ rencontre par hypothèse~$\Sigma-\Sigma\dtr$, et possède de ce fait un point~$x$ qui est rigide, et dès lors unibranche ; par conséquent,~$I$ est un intervalle semi-ouvert issu de~$x$ et ~$\partial I$ est un singleton~$\{y\}$. On a~$\partial U=\{y\}$ et~$U$ est de ce fait une composante connexe de~$X\setminus\{y\}$.

En tant que point rigide d'une courbe quasi-lisse, le point~$x$ possède un voisinage ouvert~$V$ dans~$I^\flat$ qui est un disque virtuel ; l'intersection~$V\cap I$ est un voisinage de~$x$ dans~$I$ et n'est donc pas réduite à~$\{x\}$ ; choisissons~$x'$ sur~$(V\cap I)\setminus\{x\}$, et~$x''$ sur~$]x;x'[$. 

Comme~$x$ est le seul point unibranche de~$I$, c'en est le seul point rigide, ce qui signifie que~$I\cap \Sigma=\{x\}$ puisque~$I\cap \Sigma\dtr=\emptyset$. Par conséquent,~$]x'';y[$ ne rencontre pas~$\Sigma$. Par ailleurs, les points~$x''$ et~$y$ de~$\Gamma$ appartiennent tous deux à~$X\dtr$ : pour~$x''$ cela résulte par exemple de l'égalité~$I\cap \Sigma=\{x\}$ ; et pour~$y$, de l'inclusion~$\partial I\subset \Sigma\dtr$. 

En vertu de~\ref{pretriang}, l'ouvert~$]x'';y[^\flat$ de~$X$ est une couronne virtuelle ; et le voisinage ouvert~$V$ du point~$x'$ de~$]x'';y[$ est un disque virtuel, qui ne contient pas~$y$ (il est inclus dans~$I^\flat$). On déduit alors du lemme~\ref{fusion} que la réunion de~$V$ et~$]x'';y[^\flat$ est un disque virtuel, et est par ailleurs une composante connexe de~$X\setminus\{y\}$ ; mais comme~$V\cup ]x'';y[^\flat\subset I^\flat=U$, lequel est lui-même une composante connexe de~$X\setminus\{y\}$, on a nécessairement~$V\cup ]x'';y[^\flat=U$. Ainsi,~$U$ est-il un disque virtuel, ce qu'on souhaitait démontrer. 

\deux{triadapt} Soit~$X$ une courbe~$k$-analytique quasi-lisse et soit~$\mathsf E$ un sous-ensemble fermé et discret de~$X$ contenu dans~$X\typ 0$. Soit~$S$ une triangulation de~$X$. On dira que~$S$ est {\em adaptée} à~$\mathsf E$ si pour tout~$x\in \mathsf E$, la composante connexe~$V$ de~$X\setminus S$ contenant~$x$ est un disque virtuel tel que~$V\cap \mathsf E=\{x\}$ et tel que~$\got s(V)=\got s(x)$. 

\trois{triadaptbase} Soit~$F$ une extension presque algébrique de~$k$ et soit~$S$ une triangulation de~$X$ adaptée à~$\mathsf E$ ; la triangulation~$S_F$ de~$X_F$ est adaptée à l'image réciproque~$\mathsf E_F$ de~$\mathsf E$ sur~$X_F$ si et seulement si~$S$ est adaptée à~$\mathsf E$. En effet, il suffit de le démontrer lorsque~$F=\KK$ (et d'appliquer ensuite le résultat établi à l'extension~$F\hookrightarrow \KK$). Or c'est alors immédiat, en remarquant que l'égalité~$\got s(x)=\got s(V)$ équivaut à demander que~$x$ ait un et un seul antécédent sur chacune des composantes connexes de~$V_{\KK}$.

\trois{triadaptgr} Soit~$S$ une triangulation de~$X$ adaptée à~$\mathsf E$, soit~$\Gamma$ son squelette et soit~$r$ la rétraction canonique de~$X$ sur~$\Gamma$. Il résulte aisément de la définition d'une triangulation adaptée à~$\mathsf E$ que~$r(x)\in S$ pour tout~$x\in E$, que~$[x;r(x)[\cap [x';r(x')[=\emptyset$ pour tout couple~$(x,x')$ de points {\em distincts} de ~$\mathsf E$ et que si~$x\in \mathsf E$ alors~$\got s(y)$ s'identifie canoniquement à~$\got s(x)$ pour tout~$y\in [x;r(x)[$. 

Seul ce dernier point mérite peut-être une explication : si~$V$ désigne la composante connexe de~$X\setminus\{r(x)\}$ contenant~$x$ alors~$V$ est un disque virtuel, et contient donc une unique branche~$b$ issue de~$x$. On a~$\got s(b)=\got s(V)$ (d'après~\ref{remsxsv}), qui est lui-même égal à~$\got s(x)$ par hypothèse. Soit~$d$ l'application de~$[x;r(x)[$ dans~$\NN$ qui envoie~$z$ sur~$[\got s(z):k]$. Cette application est constante par morceaux sur~$[x;r(x)[$ et égale à~$[\got s(b):k]$ au voisinage de~$r(x)$ (\ref{apropossx}) ; elle est par ailleurs croissante lorsqu'on se dirige vers~$x$ (\ref{remsxsv}). Il s'ensuit que~$d(x)=d(y)$ pour tout~$y\in [x;r(x)[$ et l'assertion requise (plus précise) s'en déduit aussitôt grâce aux équivalences énoncées au~\ref{compsxsurgamma}.

\trois{grtriadapt} Réciproquement, soit~$\Gamma$ un sous-graphe localement fini et analytiquement admissible de~$X$ tracé sur~$X\dtr$ et soit~$r$ la rétraction canonique de~$X$ sur~$\Gamma$. Supposons que pour tout couple~$(x,x')$ d'éléments distincts de~$\mathsf E$ l'on ait~$[x;r(x)[\cap [x';r(x')[=\emptyset$, et que pour tout~$x\in \mathsf E$ la fonction~$y\mapsto [\got s(y):k]$ soit constante sur~$]x;r(x)[$ ; si~$S$ est une triangulation de~$X$ de squelette~$\Gamma$ et contenant l'ensemble des~$r(x)$ pour~$x$ parcourant~$\mathsf E$, alors~$S$ est adaptée à~$\mathsf E$. 

Là encore, seul un point mérite une vérification : le fait que si~$x\in \mathsf E$ et si~$V$ désigne la composante connexe de~$X\setminus\{r(x)\}$ contenant~$x$, alors~$\got s(x)=\got s(V)$. Reprenons les notations~$b$ et~$d$ du~\ref{triadapt} ci-dessus, et soit~$\delta$ la valeur constante de~$d$ sur~$]x;r(x)[$ ; par semi-continuité inférieure de~$d$, on a~$d(x)\leq \delta$. Comme~$d(y)=\delta$ pour tout~$y\in ]x;r(x)[$, c'est en particulier le cas au voisinage de~$r(x)$ et l'on a donc~$[\got s(b):k]=\delta$ ; par ailleurs~$\got s(b)=\got s(V)$ qui se plonge dans~$\got s(x)$, d'où l'inégalité ~$d(x)\geq \delta$, et finalement~$d(x)=\delta$, ce qui entraîne que~$\got s(V)\hookrightarrow \got s(x)$ est un isomorphisme. 

\deux{triminimarqueintro} Soit~$X$ une courbe~$k$-analytique quasi-lisse compacte et connexe ; on suppose que~$\skelan X$ est non vide ; c'est donc un sous-graphe analytiquement admissible de~$X$, et l'on note~$r$ la rétraction canonique de~$X$ sur~$\skelan X$. Soit~$\mathsf E$ un sous-ensemble fini de~$X$ contenu dans~$X\typ 0$ et soit~$S$ une triangulation de~$X$ adaptée à~$\mathsf E$. 

\medskip
Posons~$\Gamma=\skelan X\cup \bigcup\limits_{x\in \mathsf E} [x;r(x)].$ C'est un sous-graphe fermé localement fini de~$X$ tracé sur~$X\geom$, qui est le plus petit sous-graphe analytiquement admissible de~$X$ contenant~$\skelan X$ et~$\mathsf E$. Soit~$\Sigma$ l'ensemble des nœuds non rigides de~$\Gamma$ ; on suppose qu'il est {\em non vide}. C'est alors d'après~\ref{gammasigmanonv} une triangulation de~$X$, dont nous allons montrer qu'elle est adaptée à~$\mathsf E$ et contenue dans~$S$. 

\trois{szeroad} {\em La triangulation~$\Sigma$ est adaptée à~$\mathsf E$.} Soit~$x\in \mathsf E$ et soit~$V$ la composante connexe de~$X-\Sigma$ contenant~$x$. Il résulte de~\ref{gammasigmanonv} que~$V$ est un disque virtuel ; soit~$y$ l'unique point de son bord ; il appartient à~$\Sigma$. L'intervalle~$]x;y[$ tracé sur~$V$ est contenu dans~$\Gamma$ et n'est constitué que de points pluribranches, donc de type 2 ou 3. Il ne contient par hypothèse aucun nœud de~$\Gamma$, ce qui implique en particulier que~$z\mapsto [\got s(z):k]$ est constante sur~$]x;y[$. 

Pour montrer que~$\Sigma$ est adaptée à~$\mathsf E$, il suffit maintenant de vérifier que~$V\cap \mathsf E=\{x\}$. On raisonne par l'absurde, en supposant qu'il existe un point~$x'$ sur~$V\cap \mathsf E$ qui diffère de~$x$. L'intersection~$]x;y[\cap ]x';y[$ est alors de la forme~$]t;y[$, et la valence de~$(\Gamma,t)$ vaut au moins 3. Ainsi,~$t$ est un sommet topologique pluribranche de~$\Gamma$, et en particulier un nœud de ce dernier situé sur~$X\dtr$ ; il appartient donc à~$\Sigma$, ce qui est contradictoire. 

\trois{szerosubs} {\em On a~$\Sigma\subset S$.} Désignons par~$\Delta$ le squelette de la triangulation~$S$. Soit~$y\in \Sigma$. 

\medskip
Si~$y$ est un nœud de~$\skelan X$, il résulte de~\ref{noeudsinclus} que~$y$ est un nœud de~$\Delta$, et donc un élément de~$S$. 

\medskip
Supposons maintenant que~$y$ n'est pas un nœud de~$\skelan X$. Comme c'est un nœud de~$\Gamma=\skelan X\cup \bigcup\limits_{x\in \mathsf E} [x;r(x)]$ et comme~$\bigcup [x;r(x)]$ est un sous-graphe compact de~$\Gamma$, le point~$y$ appartient nécessairement à~$\bigcup [x;r(x)]$. 

Si~$x\in \mathsf E$, il résulte du~\ref{szeroad} ci-dessus : 

\medskip
$\bullet$ que~$]x;r(x)[$ ne rencontre aucun des~$]x';r(x')[$ pour~$x'\neq x'$, et partant ne contient aucun sommet topologique de~$\Gamma$ ; 

$\bullet$ que~$z\mapsto [\got s(z):k]$ est constante sur~$]x;r(x)[$. 

\medskip
Compte-tenu du fait que~$]x;r(x)[$ est tracé sur un disque virtuel, qui est sans bord, et que chacun de ses points est pluribranche, donc non rigide, il s'ensuit que~$]x;r(x)[$ ne contient aucun nœud de~$\Gamma$. 

On déduit de ce qui précède que~$\bigcup[x;r(x)[$ ne contient aucun nœud non rigide de~$\Gamma$ ; dès lors,~$y$ est nécessairement de la forme~$r(x)$ pour un certain~$x\in \mathsf E$ ; on note~$V$ la composante connexe de~$X-\skelan X$ qui contient~$x$ ; c'est une composante connexe de~$X\setminus\{y\}$. Comme on a fait l'hypothèse que~$y$ n'est pas un nœud de~$\skelan X$, la valence de~$(\skelan X,y)$ est 2. Comme~$\Delta\supset \skelan X$, la valence de~$(\Delta, y)$ est supérieure ou égale à~$2$. 

\medskip
{\em Supposons que la valence de~$(\Delta,y)$ est au moins égale à~$3$.} Dans ce cas~$y$ est un sommet topologique, et {\em a fortiori} un nœud, de~$\Delta$ ; c'est donc un élément de~$S$. 

\medskip
{\em Supposons que la valence de~$(\Delta,y)$ est égale à~$2$.} Dans ce cas les {\em germes}~$(\Gamma,y)$ et~$(\Delta,y)$ sont égaux. Par conséquent,~$\Delta$ ne contient aucun intervalle ouvert tracé sur~$V$ et aboutissant à~$y$ ; étant convexe,~$\Delta$ ne rencontre donc pas~$V$, qui apparaît ainsi comme la composante connexe de~$x$ dans~$X-\Delta$. Mais comme~$S$ est une triangulation adaptée à~$\mathsf E$, ceci entraîne que l'unique point du bord de~$V$, qui n'est autre que~$y$, appartient à~$S$. 

\subsection*{Triangulations des courbes marquées}

\deux{defcourbemar} On appellera {\em courbe~$k$-analytique marquée}  un couple~$(X,\mathsf E)$ où~$X$ est une courbe~$k$-analytique et où~$\mathsf E$ est un sous-ensemble fermé et discret de~$X$ contenu dans~$X\typ 0$. Une {\em triangulation} d'une courbe marquée~$(X,\mathsf E)$ sera par définition une triangulation de~$X$ adaptée à~$\mathsf E$. 

\deux{cmarquecomp} Soit~$(X,\mathsf E)$ une courbe~$k$-analytique marquée. On suppose que~$X$ est non vide, géométriquement connexe, quasi-lisse et compacte (cette dernière condition entraîne la finitude de~$\mathsf E$). Le but de ce qui suit est de déterminer sous quelles hypothèses supplémentaires il existe une plus petite triangulation de~$(X,\mathsf E)$ ; on note~$\Theta$ l'ensemble des nœuds de~$\skelan X$ et l'on distingue trois cas. 

\trois{marquegen} {\em Le cas où~$\skelan X$ et~$\Theta$ sont non vides.} On sait alors que~$\skelan X$ est un sous-graphe analytiquement admissible de~$X$, contenu dans~$X\dtr$. Soit~$\Gamma$ le plus petit sous-graphe analytiquement admissible de~$\mathsf X$ contenant~$\skelan X$ et~$\mathsf E$. L'ensemble~$\Sigma$ des nœuds non rigides de~$\Gamma$ contient~$\Theta$ en vertu de~\ref{noeudsinclus}, et est en particulier non vide ; il découle alors de~\ref{triminimarqueintro},~\ref{szeroad} et~\ref{szerosubs} que~$\Sigma$ est la plus petite triangulation de~$(X,\mathsf E)$. 

\trois{marqueunpeupart} {\em Le cas où ~$\skelan X$ est non vide et où~$\Theta$ est vide.} La courbe~$X$ est alors isomorphe à l'analytifiée d'une courbe projective, lisse et géométriquement intègre de genre 1, du type décrit au~\ref{noeudsvide} ; son squelette analytique~$\skelan X$ coïncide avec son squelette, et est un cercle. 

\medskip
Si~$\mathsf E$ est vide, on a vu au~\ref{noeudsvide} que~$X$ n'admet pas de plus petite triangulation, mais admet une infinité de triangulations minimales, à savoir les triangulations de la forme~$\{x\}$ avec~$x\in \skelan X$. 

\medskip
Si~$\mathsf E$ est non vide, appelons~$\Gamma$ le plus petit sous-graphe analytiquement admissible de~$X$ contenant~$\mathsf E$ et~$\skelan X$. Comme~$\skelan X$ est un cercle et comme~$\mathsf E$ est non vide,~$\Gamma$ possède au moins un sommet topologique pluribranche, et l'ensemble~$\Sigma$ de ses nœuds non rigides est dès lors non vide. Il découle alors de~\ref{triminimarqueintro},~\ref{szeroad} et~\ref{szerosubs} que~$\Sigma$ est la plus petite triangulation de~$(X,\mathsf E)$. 

\trois{marquetrespart} {\em Le cas où~$\skelan X$ est vide.} La courbe~$X$ est alors isomorphe à l'analytification d'une courbe projective, lisse et géométriquement intègre de genre zéro du type décrit au~\ref{equivsquelvide}. 

\medskip
Soit~$\mathsf E_{\KK}$ l'image réciproque de~$\mathsf E$ sur ~$X_{\KK}\simeq \pkk$. Le but de ce qui suit est de montrer que~$(X,\mathsf E)$ admet une plus petite triangulation si et seulement si la condition suivante est satisfaite : {\em~$\mathsf E_{\KK}$ compte au moins trois éléments, ou en compte exactement deux qui sont permutés par l'action de Galois.} Elle peut se retraduire comme suit : {\em~$\sum \limits_{x\in \mathsf E} [\got s(x):k]\geq 3$ ou~$\mathsf E$ est un singleton~$\{x\}$ tel que~$[\got s(x):k]=2$.} 

\medskip
Pour établir cette assertion, l'on distingue quatre cas. 

\medskip
\noindent
{\bf Premier cas :~$\mathsf E_{\KK}$ compte au moins trois éléments.} Soit~$\Gamma$ l'enveloppe convexe de~$\mathsf E_{\KK}$ et soit~$S$ l'ensemble des sommets pluribranches de~$\Gamma$. Comme ~$\mathsf E_{\KK}$ compte au moins trois éléments,~$S$ est non vide, et il coïncide avec l'ensemble des nœuds non rigides de~$\Gamma$. Le sous-arbre~$\Gamma$ de~$X_{\KK}\simeq \pkk$ est analytiquement admissible (remarque~\ref{remgranaladmp}) ; on déduit alors de~\ref{gammasigmanonv} que~$S$ est une triangulation de~$X_{\KK}$ adaptée à~$\mathsf E_{\KK}$. 

\medskip
L'arbre~$\Gamma$ et l'ensemble~$S$ sont stables sous~$\mathsf G$. Si~$I$ est une arête de~$\Gamma$ joignant deux points de~$S$ et dont~$\mathsf G$ renverse les orientations, il existe un unique point~$\xi(I)$ de~$I$ qui est fixé par le stabilisateur de~$I$ ; on note~$S'$ la réunion de~$S$ et des~$\xi(I)$, où~$I$ parcourt l'ensemble des arêtes de~$\Gamma$ joignant deux points de~$S$ dont~$\mathsf G$ renverse les orientations. L'ensemble~$S'$ est stable sous~$\mathsf G$ ; il constitue d'après~\ref{casxsurg} une triangulation de~$X_{\KK}$. 

\medskip
De plus, si un élément de~$S'$ est de la forme~$\xi(I)$, alors il appartient à la composante connexe~$I^\flat$ de~$X_{\KK}\setminus S$, qui est une couronne et qui ne contient donc aucun point de~$\mathsf E$, puisque~$S$ est adaptée à~$\mathsf E$ ; il s'ensuit que~$S'$ est encore adaptée à~$\mathsf E$. Par ailleurs, l'appartenance de chacun des~$\xi(I)$ à~$S'$ garantit que si~$J$ est un intervalle ouvert joignant deux points de~$S'$ alors~$\mathsf G$ ne permute pas les orientations de~$J$. 

\medskip
On déduit de ce qui précède, à l'aide du théorème \ref{theo-invariance-triangulations} et de~\ref{triadaptbase}, que~$S'$ est l'image réciproque sur~$X_{\KK}$ d'une triangulation~$T$ de~$X$ qui est adaptée à~$\mathsf E$. 

\medskip
Nous allons montrer que~$T$ est la plus petite triangulation de~$(X,\mathsf E)$. Soit donc~$\Theta$ une triangulation de~$(X,\mathsf E)$ ; il suffit de vérifier que l'image réciproque~$\Theta_{\KK}$ de~$\Theta$ sur~$X_{\KK}$ contient~$S'$ ; notons que~$\Theta_{\KK}$ est adaptée à~$\mathsf E_{\KK}$ (\ref{triadaptbase}). 

\medskip
Soit~$\Delta$ le squelette de~$\Theta_{\KK}$, soit~$r$ la rétraction canonique de~$X_{\KK}$ sur~$\Delta$, et soit~$t\in S$. Par définition de~$S$ il existe trois points deux à deux distincts~$x,y,$ et~$z$ de~$\mathsf E_{\KK}$ tels que~$[x;z]\cap [y;z]=[t;z]$. On a~$[x;z]=[x;r(x)]\cup I\cup [r(z);z]$ et~$[y;z]=[y;r(y)]\cup I'\cup [r(z);z]$ où~$I$ et~$I'$ sont deux segments tracés sur~$\Delta$. Comme~$\Theta_{\KK}$ est adaptée à~$\mathsf E_{\KK}$, le point~$t$ ne peut appartenir à~$[x;r(x)[\cup [y;r(y)[\cup [z;r(z)[$. 

Si~$t$ est égal à~$r(x),r(y)$ ou~$r(z)$ alors il appartient à~$\Theta_{\KK}$, là encore parce que~$\Theta_{\KK}$ est adaptée à~$\mathsf E_{\KK}$. 

Sinon~$t$ appartient à~$[x;z]\cup [y;z]-([x;r(x)]\cup [y;r(y)]\cup [z;r(z)])$, qui est un ouvert de~$[x;z]\cup [y;z]$ et un sous-graphe de~$\Delta$. Comme la valence de~$[x;z]\cup [y;z]$ en~$t$ vaut 3, la valence de~$\Delta$ en~$t$ vaut au moins 3 ; par conséquent,~$t$ est un nœud de~$\Delta$ et appartient de ce fait à~$\Theta_{\KK}$. 

\medskip
Ainsi,~$\Theta_{\KK}$ contient~$S$. Par ailleurs, soit~$I$ un intervalle ouvert joignant deux points de~$S$ et dont~$\mathsf G$ permute les orientations ; nous allons montrer par l'absurde que~$\Theta_{\KK}$ contient~$\xi(I)$. Supposons que ce ne soit pas le cas et soit~$I_0$ la composante connexe de~$\xi(I)$ dans~$I-\Theta_{\KK}$. Par hypothèse, il existe~$g\in \mathsf G$ stabilisant~$I$ et induisant sur ce dernier une symétrie de centre~$\xi(I)$ ; il s'ensuit que~$g$ stabilise~$I_0$ et induit sur celui-ci une symétrie de centre~$\xi(I)$, ce qui contredit le fait que~$\mathsf G$ ne permute pas les orientations des composantes connexes de~$\Delta- \Theta_{\KK}$ (théorème~\ref{theo-invariance-triangulations}). 

\medskip
On en conclut que~$\Theta_{\KK}$ contient~$S'$, ce qu'on souhaitait établir. 

\medskip
\noindent
{\bf Second cas :~$\mathsf E_{\KK}$ compte deux points échangés par l'action de~$\mathsf G$.} Soient~$x$ et~$y$ les deux points en question. Comme~$\mathsf G$ échange~$x$ et~$y$, il stabilise~$[x;y]$ et son action sur ce dernier se factorise par une surjection~$\mathsf G\to \{\mathsf{Id}, g\}$ où~$g$ est une symétrie ayant pour centre un certain~$z\in ]x;y[$. Le point~$z$ est stable sous~$\mathsf G$, c'est donc l'unique antécédent d'un point~$t$ de~$X$. Le singleton~$\{z\}$ est une triangulation de~$X_{\KK}$ (le point~$z$ est pluribranche, donc de type 2 ou 3, et~$X_{\KK}\simeq \pkk$) visiblement adaptée à~$\mathsf E_{\KK}=\{x,y\}$ ; on déduit alors du théorème~\ref{theo-invariance-triangulations} et du~\ref{triadaptbase} que~$\{t\}$ est une triangulation de~$X$ adaptée à~$\mathsf E$. 

\medskip
Nous allons montrer que~$\{t\}$ est la plus petite triangulation de~$(X,\mathsf E)$. Soit donc~$\Theta$ une triangulation de~$(X,\mathsf E)$ ; il suffit de vérifier que l'image réciproque~$\Theta_{\KK}$ de~$\Theta$ sur~$X_{\KK}$ contient~$z$. 

\medskip
Soit~$\Delta$ le squelette de~$\Theta_{\KK}$ et soit~$r$ la rétraction canonique de~$X_{\KK}$ sur~$\Delta$. Comme~$\Theta_{\KK}$ est adaptée à~$\mathsf E_{\KK}=\{x,y\}$ (\ref{triadaptbase}), les points~$r(x)$ et~$r(y)$ appartiennent à~$\Theta_{\KK}$. Par convexité,~$\Delta\cap [x;y]=[r(x);r(y)]$. Comme~$g$ stabilise~$\Delta$, il échange~$r(x)$ et~$r(y)$ et~$z$ appartient dès lors à~$[r(x);r(y)]$. 

Supposons que~$z$ n'appartienne pas à~$\Theta_{\KK}$. La composante connexe de~$z$ dans~$\Delta-\Theta_{\KK}$ serait alors un intervalle ouvert stable sous~$g$ et contenu dans~$]r(x);r(y)[$, contredisant ainsi le fait que~$\mathsf G$ agit sans permuter les orientations des composantes connexes de~$\Delta-\Theta_{\KK}$ (théorème~\ref{theo-invariance-triangulations}). 

\medskip
\noindent
{\bf Troisième cas :~$\mathsf E_{\KK}$ compte deux points distincts, fixes sous~$\mathsf G$.} Soient~$x$ et~$x'$ les deux points en question ; l'ensemble~$\mathsf E$ consiste alors lui-même en deux points~$y$ et~$y'$, qui sont les images respectives de~$x$ et~$x'$ sur~$X$, et~$\mathsf G$ fixe point par point l'intervalle~$[x;x']$ ; celui-ci s'envoie donc homéomorphiquement sur son image qui est de ce fait égale à~$[y;y']$. 

Soit~$z\in ]x;x'[$ et soit~$t$ son image sur~$X$. Le singleton~$\{z\}$ est une triangulation de~$X_{\KK}$ (le point~$z$ est pluribranche, donc de type 2 ou 3, et~$X_{\KK}\simeq \pkk$), visiblement adaptée à~$\mathsf E_{\KK}=\{x,x'\}$ ; on déduit alors du théorème
\ref{theo-invariance-triangulations}
et du~\ref{triadaptbase} que~$\{t\}$ est une triangulation de~$X$ adaptée à~$\mathsf E$. Elle est évidemment minimale ; ainsi,~$(X,\mathsf E)$ possède une infinité de triangulations minimales , à savoir toutes celles de la forme~$\{t\}$ avec~$t\in ]y;y'[$ ; il s'ensuit que~$(X,\mathsf E)$ ne possède pas de plus petite triangulation. 

Faisons une remarque : si~$\Theta$ est une triangulation de~$(X,\mathsf E)$ alors son image réciproque~$\Theta_{\KK}$ sur~$X_{\KK}$ est une triangulation adaptée à~$\mathsf E_{\KK}=\{x,x'\}$. Par conséquent, il y a au moins un sommet de~$\Theta_{\KK}$ qui est situé sur~$]x;x'[$. Autrement dit,~$\Theta$ contient une triangulation de la forme~$\{t\}$ pour un certain~$t\in ]y;y'[$.  Ainsi les triangulations minimales de~$(X,\mathsf E)$ sont-elles exactement les triangulations de la forme~$\{t\}$ avec~$t\in ]y;y'[$.  

\medskip
\noindent
{\bf Quatrième cas :~$\mathsf E_{\KK}$ compte au plus un point.} On sait que~$X$ possède une infinité de triangulations minimales, et que les triangulations minimales de~$X$ sont les triangulations de la forme~$\{x\}$ où~$x\in X\dtr$ et où~$\got s(x)=k$. Nous allons montrer que cette assertion reste vraie si l'on remplace~$X$ par~$(X,\mathsf E)$ (en particulier,~$(X,\mathsf E)$ n'a pas de plus petite triangulation). Pour le voir, on peut supposer~$\mathsf E$ non vide, et donc constitué d'un point~$y$ tel que~$\got s(y)=k$. Il suffit de vérifier que si~$x$ est un point de~$X\dtr$ tel que~$\got s(x)=k$, alors~$\{x\}$ est adaptée à~$\mathsf E$. Soit~$V$ la composante connexe de~$X\setminus\{x\}$ contenant~$y$ ; comme toute composante connexe de~$X\setminus\{x\}$, c'est un disque virtuel ; et comme~$\got s(V)\hookrightarrow \got s(y)=k$, on a~$\got s(V)=\got s(Y) \;(=k)$, ce qui termine la preuve. 

\subsection*{Exemples}

\deux{tricourbell} {\em Le cas d'une courbe elliptique.} Soit~$\sch X$ une courbe elliptique sur~$k$ et soit~$e\in \sch X(k)$ son élément neutre. Soit~$\Gamma$ le sous-graphe de~$\sch X\an$ égal à l'enveloppe convexe de~$\{e\}\cup \skelan {\sch X\an}$ et soit~$\Sigma$ l'ensemble des nœuds non rigides de~$\Gamma$. Il résulte de~\ref{marquegen} et~\ref{marqueunpeupart} que~$\Sigma$ est une triangulation de~$\sch X\an$, et que c'est plus précisément la plus petite triangulation de la courbe marquée~$(\sch X\an, e)$ ; ainsi, l'analytifiée d'une courbe elliptique possède une triangulation canonique. 

\deux{p1quad} {\em Le cas de~$\pk$ muni d'un point rigide quadratique.} Supposons que~$k$ est de caractéristique différente de 2, soit~$a$ un élément de~$k$ qui n'est pas un carré, et soit~$x$ le point rigide de~$\Aff^{1,{\rm an}}_k$ défini par l'équation~$T^2=a$. Soit~$\alpha$ une racine carrée de~$a$ dans~$\KK$ et soit~$\eta$ l'unique point du bord de Shilov du disque fermé~$V$ de~$\Aff^{1,{\rm an}}_{\KK}$ centré en~$\alpha$ et de rayon~$|2\alpha|$ ; le point~$\eta$  est invariant sous l'action de Galois (car c'est le cas de~$V$, qui est le plus petit disque fermé de~$\Aff^{1,{\rm an}}_{\KK}$ contenant~$\alpha$ et~$-\alpha$) ; soit~$\xi$ son image sur~$\pk$. Il résulte alors du second cas examiné au~\ref{marquetrespart} que~$\{\xi\}$ est la plus petite triangulation de~$(\pk,x)$. 

Nous allons maintenant décrire un peu plus précisément le point~$\xi$ en distinguant deux cas. 

\trois{etacardiff2} {\em Supposons que la caractéristique de~$\red k$ soit différente de 2.} Dans ce cas~$|2\alpha|=|\alpha|$ et~$\eta=\eta_{|\alpha|, \KK}$ ; par conséquent,~$\xi=\eta_{|\alpha|}=\eta_{\sqrt {|a|}}$. 

\trois{etacar2} {\em Supposons que la caractéristique de~$\red k$ soit égale à~$2$}. On a alors pour toute extension complète~$L$ de~$\KK$ et pour tout~$\lambda\in L$ équivalence entre les inégalités~$|\lambda-\alpha|\leq |2\alpha|$ et~$|\lambda^2-a|\leq |4\alpha^2|=|4|.|a|$. Il s'ensuit que~$V=U_{\KK}$, où~$U$ désigne le domaine analytique de~$\Aff^{1,{\rm an}}_k$ défini par l'inégalité~$|T^2-a|\leq |4|.|a|$ ; par conséquent,~$\xi$ peut être caractérisé comme l'unique point du bord de Shilov de~$U$.

\chapter{Modèles formels d'une courbe analytique en théorie de Berkovich}
\markboth{Modèles formels}{Modèles formels}
\section{Une alternative pour les courbes irréductibles compactes}

\deux{cortoutoucomp}{ \bf Lemme.} {\em Soit~$\sch X$ une courbe algébrique projective irréductible, munie d'un morphisme fini et dominant sur~$\PP^1_k$ dont on note~$\phi$ l'analytifié. Soit~$x$ un point de~$\sch X\an$ dont l'image par~$\phi$ est égale à~$\eta_r$ pour un certain~$r>0$, et soit~$\Omega$ une composante connexe de~$\sch X\an\setminus\{x\}$ telle que~$\br {\sch X\an}x \ctd \Omega$ soit un singleton et telle que~$\phi\inv(\infty)\subset\Omega$ ; soit~$\DD$ le disque fermé de~$\pk$ de centre~$0$ et de rayon~$r$. Le compact~$\sch X\an-\Omega$ est alors une composante connexe de~$\phi\inv(\DD)$.}

\medskip
{\em Démonstration.} La question étant purement topologique, on peut remplacer~$\sch X$ par~$\sch X_{\rm red}$, auquel cas le morphisme~$\phi$ est plat. Si~$V$ est une composante connexe de~$\sch X\an\setminus\{x\}$ qui n'est pas égale à~$\Omega$, il résulte de nos hypothèses que~$\infty\notin \phi(V)$ ; on déduit alors du lemme~\ref{imtoutoucomp} que~$\phi(V)$ est une composante connexe de~$\pk\setminus\{\eta_r\}$ ; comme elle ne contient pas~$\infty$, cette composante est incluse dans~$\DD$. Il s'ensuit, compte-tenu du fait que~$\phi(x)=\eta_r\in \DD$, que le compact~$\sch X\an-\Omega$ est contenu dans~$\phi\inv(\DD)$. 

\medskip
Posons~$U=\pk-\DD$. Par ce qui précède,~$\phi\inv(U)\subset \Omega$ ; si~$b$ désigne l'unique branche de~$\pk$ issue de~$\eta_r$ et contenue dans~$U$, chacune des branches de~$\sch X\an$ située au-dessus de~$b$ est contenue dans~$\Omega$ ; comme il existe par ailleurs au moins une branche de~$\sch X\an$ issue de~$x$ et située au-dessus de~$b$, on en déduit que l'unique branche~$\beta$ de~$\sch X\an$ issue de~$x$ et contenue dans~$\Omega$ vérifie l'égalité~$\phi(\beta)=b$. Par conséquent, il existe une section~$V$ de~$\beta$ telle que~$\phi(V)$ soit une section de~$b$ ; on a en particulier~$\phi(V)\subset U$. 

Puisque~$\beta$ est la seule branche de~$\sch X\an$ qui soit issue de~$x$ et contenue dans~$\Omega$, le fermé~$\Omega-V$ de~$\Omega$ est également fermé dans~$\sch X\an$ ; on en déduit que le compact connexe~$\sch X\an-\Omega$ est ouvert dans~$\sch X\an-V$ ; comme~$\phi(V)\subset U$, l'on a~$\phi\inv(\DD)\subset \sch X\an-V$ ; par conséquent,~$\sch X\an-\Omega$ est ouvert et fermé dans~$\phi\inv(\DD)$ ; étant de surcroît connexe et non vide, ~$\sch X\an-\Omega$ est une composante connexe de~$\phi\inv(\DD)$.~$\Box$ 

\deux{ablouvnorm} {\bf Proposition.} {\em Soit~$\sch X$ une~$k$-courbe algébrique projective et soit~$(\Omega_i)_{i\in I}$ une famille finie d'ouverts deux à deux disjoints de~$\sch X\an$ possédant les propriétés suivantes : 

\medskip
i) chacun des~$\Omega_i$ est connexe, normal et non vide ;

ii) pour tout~$i$, le bord topologique de~$\Omega_i$ dans~$\sch X\an$ est un singleton~$\{x_i\}$, où le point~$x_i$ est de type 2 ou 3, et~$\br {\sch X\an}{x_i}\ctd {\Omega_i}$ est également un singleton ;

iii) pour toute composante irréductible~$\sch Y$ de~$\sch X$, il existe un indice~$i$ tel que~$\Omega_i\cap \sch Y\an\neq \emptyset$.

\medskip
Sous ces hypothèses~$\sch X\an-\coprod \Omega_i$ est un domaine affinoïde de~$\sch X\an$.}

\medskip
{\em Démonstration.} Commençons par quelques remarques. Pour tout~$i$, l'ouvert~$\Omega_i$ est normal, et chacun de ses points appartient donc à l'analytification d'une {\em unique} composante irréductible de~$\sch X$ ; comme~$\Omega_i$ est par ailleurs connexe, il s'ensuit qu'il existe une composante irréductible~$\sch Y$ de~$\sch X$ telle que~$\Omega_i\subset \sch Y\an$, et telle que~$\Omega_i$ ne rencontre l'analytification d'aucune autre composante irréductible de~$\sch X$. 

Il en résulte notamment que dans l'énoncé de la condition iii), l'on peut remplacer «~$\Omega_i\cap \sch Y\an\neq \emptyset$ » par «~$\Omega_i\subset \sch Y\an$». 

\medskip
Par ailleurs,~$\sch X\an-\coprod \Omega_i$ est, d'après la proposition~\ref{propdomanferm}, un domaine analytique de~$\sch X\an$ ; il reste à s'assurer qu'il est affinoïde.

\trois{passkr}{\em Réduction au cas où la valeur absolue de~$k$ n'est pas triviale.} Soit~$\bf r$ un polyrayon~$k$-libre ; après extension des scalaires à~$k_{\bf r}$, toute composante irréductible de~$\sch X$ reste irréductible et~$\sch X$ reste réduite. Fixons un indice~$i$. L'ouvert~$ \Omega_{i, \bf r}$ de~$\sch X_{\bf r}\an$ est connexe et normal ; par ailleurs, l'hypothèse ii) assure l'existence d'un ensemble cofinal~$\sch K$ de compacts de~$\Omega_i$ tels que~$\Omega_i-\cal K$ soit connexe pour tout~${\cal K}\in \sch K$. Si~${\cal K}\in \sch K$, son image réciproque~${\cal K}_{\bf r}$ sur~$ \Omega_{i, \bf r}$ est un compact, et~$\{{\cal K}_{\bf r}\}_{{\cal K}\in \sch K}$ est un ensemble cofinal de compacts de~$\Omega_{i,\bf r}$ ; pour tout ~${\cal K}\in \sch K$ l'ouvert ~$\Omega_{i,\bf r}-{\cal K}_{\bf r}$ est égal à~$(\Omega_i-{\cal K})_{\bf r}$ et est de ce fait connexe ; il s'ensuit que le bord de~$\Omega_{i,\bf r}$ dans le graphe compact~$\sch X_{\bf r}\an$ est un singleton~$\{y_i\}$, et que~$\Omega_{i,\bf r}$ contient une et une seule branche issue de~$y_i$. 

Soit~$\sigma$ la section de Shilov de~$\sch X\an_{\bf r}\to \sch X\an$ ; le point ~$x_i$ adhère à~$\Omega_i$ ; il en résulte que~$\sigma(x_i)$ adhère à~$\Omega_{i,\bf r}$, et appartient de ce fait à son bord (puisque~$x_i\notin \Omega_i$). Par conséquent,~$y_i=x_i$, et le point~$y_i$ est donc de type 2 ou 3. 

\medskip
Par ailleurs, en vertu de la proposition~\ref{affkr},~$\sch X\an -\coprod \Omega_i$ est un domaine affinoïde de~$\sch X\an$ si et seulement si~$\sch X\an_{\bf r} -\coprod \Omega_{i,\bf r}$ est un domaine affinoïde de~$\sch X\an_{\bf r}$.  

\medskip
Les faits qui précèdent autorisent, quitte à étendre les scalaires à~$k_{\bf r}$ pour n'importe quel polyrayon~$k$-libre non vide~$\bf r$, à supposer que~$|k\ti|\neq \{1\}$. 

\trois{unseulomegai} Soit~$\sch J$ l'ensemble des parties~$J$ de~$I$ satisfaisant la condition suivante : {\em pour toute composante irréductible~$\sch Y$ de~$\sch X$ il existe un {\em et un seul}~$j\in J$ tel que~$\Omega_j\subset \sch Y\an$} ; on a~$$\sch X\an-\coprod_{i\in I} \Omega_i=\bigcap_{J\in \sch J}\left(\sch X\an-\coprod_{j\in J}\Omega_j\right).$$ Comme l'intersection d'une famille finie de domaines affinoïdes de l'espace séparé~$\sch X\an$ est encore un domaine affinoïde de~$\sch X\an$, on se ramène à traiter le cas où  pour toute composante irréductible~$\sch Y$ de~$\sch X$ il existe un {\em et un seul}~$i\in I$ tel que~$\Omega_i\subset \sch Y\an$. 

\medskip
Si~$i\in I$ on désignera désormais par~$\sch X_i$ la composante irréductible de~$\sch X$ dont l'analytification contient~$\Omega_i$ ; les~$\sch X_i$ sont deux à deux distinctes et toute composante irréductible de~$\sch X$ est l'une d'entre elles.

\trois{foncmer} Pour tout~$i$, il existe, en vertu du {\em Nullstellensatz}, un point rigide~$P_i$ sur~$\Omega_i$. Chaque~$P_i$ définit (par normalité de~$\Omega_i$) un diviseur de Cartier,  et donc un fibré en droites~$\sch O(P_i)$, sur~$\sch X$. Si~$\sch L$ désigne le fibré~$\bigotimes \sch O(P_i)$ l'on a pour tout entier~$n$ une suite exacte~$$0\to \sch L^{\otimes n-1}\to \sch L^{\otimes n}\to \bigoplus_i\sch L_{|P_i}\to 0.$$ Comme~$\{P_i\}$ rencontre toutes les composantes irréductibles de ~$\sch X$, le fibré~$\sch L$ est ample, et~$\H^1(\sch X, \sch L^{\otimes n-1})$ est donc nul pour~$n$ suffisamment grand ; il s'ensuit que si~$n$ est suffisamment grand, il existe une section globale de~$\sch L^{\otimes n}$ qui ne s'annule en aucun des~$P_i$ ; celle-ci peut s'interpréter comme une fonction~$f$ sur le schéma affine~$\sch X':=\sch X\setminus\{P_i\}_i$ qui possède un pôle en chacun des~$P_i$ ; pour tout~$i$, l'on pose~$\sch  X'_i=\sch X'\cap \sch X_i=\sch X_i\setminus\{P_i\}$. 

\medskip
Fixons~$j\in I$ ; soit~$g$ une fonction sur~$\sch X'$ nulle sur chacune des~$\sch X'_i$ pour~$i\neq j$ et dont la restriction à~$\sch X'_j$ est génériquement inversible. Cette dernière condition implique,~$x_j$ étant de type~$2$ ou~$3$, que~$g(x_j)\neq 0$. Comme~$|k\ti|\neq \{1\}$, il existe~$\lambda\in k\ti$ tel que~$|\lambda|<|g(x_j)|$ 
 ; si~$N$ est un entier suffisamment grand alors~$|f(x_i)|.|\lambda|^N< |f(x_j)|.|g(x_j)|^N$ pour tout~$i\neq j$ ; choisissons un tel~$N$ et posons~$h=f(g+\lambda)^N$. Par construction,~$h$ est une fonction sur~$\sch X'$ qui possède un pôle en chacun des~$P_i$ et qui est telle que~$|h(x_i)|<|h(x_j)|$ pour tout~$j\neq i$. 
 
 \medskip
La fonction~$h$ induit un morphisme fini de~$\sch X$ vers~$\PP^1_k$ dont la restriction à chacune des~$\sch X_i$ est dominante, et dont on note~$\phi$ l'analytification ; on a~$\phi\inv(\infty)=\{P_i\}_{i\in I}$. Soit~$T$ la fonction coordonnée sur~$\PP^1_k$ ; par choix de~$h$, l'on a~$|T(\phi(x_i))|<|T(\phi(x_j))|$ pour tout~$i\neq j$. Choisissons un ensemble fini~$\sch E$ de points rigides de~$\pk$ tel que~$\{0,\infty\}\subset \sch E$ et tel que l'enveloppe convexe~$\Gamma$ de~$\sch E$ contienne chacun des~$\phi(x_i)$, et soit~$\mathsf m$ un polynôme de~$k[T]$ dont~$\sch E$ est le lieu ensembliste des zéros. Le polynôme~$\mathsf m$ induit un morphisme fini et plat~$\mu : \pk\to \pk$ tel que~$\mu\inv(\infty)=\{\infty\}$ ; comme~$|\mathsf m|$ est localement constante en dehors de~$\Gamma$ et strictement croissante sur~$\Gamma$ lorsqu'on oriente celui-ci vers~$\infty$, l'image par~$\mu$ de~$\Gamma$ est l'intervalle~$[0;\infty]$. Dès lors,~$\mu(\phi(x_i))$ est pour tout~$i$ de la forme~$\eta_{r_i}$ pour un certain~$r_i>0$ ; de plus, l'inégalité~$|T(\phi(x_i))|<|T(\phi(x_j))|$, valable pour tout~$i\neq j$, garantit que quitte à augmenter suffisamment l'exposant de~$T$ dans le polynôme~$\mathsf m$ (sans toucher  à ses autres facteurs irréductibles), on peut faire en sorte que~$|\mathsf m(\phi(x_i))|<|\mathsf m(\phi(x_j))|$ pour tout~$i\neq j$, c'est-à-dire que~$r_i<r_j$ pour tout~$i\neq j$. 

\medskip
La considération du morphisme composé~$\mu\circ\phi$ conduit à l'énoncé suivant, dans lequel on réintroduit l'indice~$j$ dans les notations -- comme nous l'avions fixé, nous l'en avions provisoirement banni : {\em pour tout indice~$j$ il existe une application~$r_j: I\to \RR\ti_+$ dont le maximum est atteint en~$j$ et un  morphisme fini~$\psi_j : \sch X\an \to \pk$, dont la restriction à chaque~$\sch X_i$ est dominante, et qui possède les deux propriétés suivantes : 

i)~$\psi_j\inv(\infty)=\{P_i\}_{i\in I}$ ;

ii)~$\psi_j(x_i)=\eta_{r_j(i)}$ pour tout~$i$.}

\medskip
Pour tout réel strictement positif~$r$, on note~$\DD_r$ le disque fermé de centre~$0$ et de rayon~$r$.

\trois{conclupsijaff} Fixons~$j$. En vertu du corollaire ~\ref{cortoutoucomp}, le compact~$\sch X_i\an -\Omega_i$ est pour tout~$i$ une composante connexe de~$\psi_{j|\sch X_i\an}\inv(\DD_{r_j(i)})$ ; comme~$r_j(i)\leq r_j(j)$, cela implique que~$\sch X_i\an -\Omega_i\subset \psi_{j|\sch X_i\an}\inv(\DD_{r_j(j)})$.

\medskip
Si~$i\in I$, l'on désigne par~$V_{i,j}$ la composante connexe de~$\psi_j\inv(\DD_{r_j(j)})$ qui contient~$\sch X_i\an-\Omega_i$.On a vu ci-dessus que~$\sch X_j\an -\Omega_j$ est une composante connexe de~$\psi_{j|\sch X_j\an}\inv(\DD_{r_j(j)})$ ; les composantes connexes de~$\psi_j\inv(\DD_{r_j(j)})$ qui rencontrent~$\sch X_j\an$ sont donc d'une part~$V_{j,j}$, et d'autre part les composantes connexes de~$\psi_{j|\sch X_j\an}\inv(\DD_{r_j(j)})$ qui sont contenues dans~$\Omega_j$, et qui à ce titre ne peuvent coïncider avec l'une des~$V_{i,j}$. 

\medskip
Par conséquent, si l'on appelle~$W_j$ la réunion des~$V_{i,j}$ pour~$i$ dans~$I$, alors~$W_j$ est un domaine affinoïde de~$\sch X\an$ qui vérifie les propriétés suivantes :

\medskip

$\bullet$~$W_j\cap \sch X_i\an \supset \sch X_i\an-\Omega_i$ pour tout~$i$ ;

$\bullet$~$W_j\cap \sch X_j\an =\sch X_j\an-\Omega_j$. 

\medskip
On en déduit que~$\bigcap\limits_{j\in I}W_j=\sch X\an-\coprod\limits_{j\in I} \Omega_j$ ; comme chaque~$W_j$ est un domaine affinoïde de l'espace séparé~$\sch X\an$, le compact~$\sch X-\coprod\limits_{j\in I} \Omega_j$ est un domaine affinoïde de~$\sch X\an$.~$\Box$

\deux{altcompgen}{\bf Théorème.} {\em Soit~$X$ une courbe~$k$-analytique compacte et soit~$S$ son bord analytique. 

\medskip
\begin{itemize}
\item[A)] Si~$X$ est génériquement quasi-lisse, il existe une~$k$-courbe algébrique projective~$\sch X$ et un isomorphisme entre~$X$ et un domaine analytique de~$\sch X\an$. 

\medskip
\item[B)] Les propositions suivantes sont équivalentes :

\medskip
\begin{itemize}
\item[1)]~$X$ est affinoïde ;

\item[2)]~$S$ rencontre toutes les composantes irréductibles de~$X$.

\end{itemize}
\end{itemize}}

\medskip
{\em Démonstration}. L'implication 1)$\Rightarrow$2) de l'assertion B) est triviale ; il suffit donc de démontrer A) et l'implication 2)$\Rightarrow$1).

\trois{prolxprop} {\em Preuve de A).} On suppose donc~$X$ génériquement quasi-lisse. Soit~$x\in S$.  Comme~$X$ est quasi-lisse en~$x$, il existe une courbe~$k$-analytique lisse~$Y_x$ et un voisinage~$V_x$ de~$x$ dans~$X$ tel que~$V_x$ s'identifie à un domaine analytique fermé de~$Y_x$ ; fixons un voisinage affinoïde~$Z_x$ de~$x$ dans~$Y_x$ ; l'intersection~$W_x:=V_x\cap Z_x$ s'identifie à un voisinage analytique compact de~$x$ dans~$X$. 

\medskip
Comme~$S$ est fini l'on peut, en restreignant suffisamment chacun des~$Z_x$, supposer que les~$W_x$ (vus comme domaines analytiques de~$X$) sont deux à deux disjoints. Soit~$X'$ la courbe obtenue en recollant~$\coprod Z_x$ et~$X$ le long de~$\coprod W_x$ ; c'est une courbe~$k$-analytique dont~$X$ est naturellement un domaine analytique fermé, et dont le bord analytique ne rencontre pas~$X$. 

\medskip
Comme le bord analytique de~$X'$ ne rencontre pas~$X$, le bord topologique de~$X$ dans~$X'$ coïncide avec le bord analytique de~$X$, c'est-à-dire avec~$S$. On peut restreindre~$X'$ de manière à ce que son bord analytique soit vide, à ce que toutes ses composantes connexes rencontrent~$S$, et,  en vertu de l'assertion iii) de la proposition~\ref{normsurb}, à ce que l'ouvert~$X'-X$ ait un nombre fini de composantes connexes (qui sont alors nécessairement des composantes connexes de~$X'\setminus S$). Sous ces hypothèses, le théorème~\ref{theochir} fournit une identification entre~$X$ et une partie compacte de l'analytification~${\sch X}\an$ d'une~$k$-courbe projective et génériquement lisse~$\sch X$, modulo laquelle sont satisfaites les conditions suivantes : 

\medskip
$\bullet$ il existe un sous-ensemble fini et non vide~$\sch B$ de~$\br {{\sch X}\an}S$ tel que l'ouvert~${\sch X}\an-X$ s'écrive comme une réunion finie disjointe~$\coprod\limits_{b\in \sch B} \Omega_b$, où chaque~$\Omega_b$ est une composante connexe normale de~${\sch X}\an\setminus S$ ;

$\bullet$ pour tout~$b$ l'ensemble~$\br {{\sch X}\an}S\ctd{ \Omega_b}$ est le singleton~$\{b\}$ ; 

$\bullet$ il existe un voisinage ouvert de~$X$ dans~${\sch X}\an$ qui est isomorphe à un voisinage ouvert de~$X$ dans~$X'$ (ce qui entraîne que le bord topologique de~$X$ dans~$\sch X\an$ est égal à~$S$. )

\medskip
Il découle alors de la proposition~\ref{propdomanferm} que~$X$ est un domaine analytique compact de~$\sch X\an$, ce qui achève de montrer A). 

\trois{altaffcasql}{\em Preuve de 2)$\Rightarrow$1) dans le cas génériquement quasi-lisse.} On suppose toujours que~$X$ est génériquement quasi-lisse, on conserve les notations~$\sch X$ et~$\Omega_b$ du~\ref{prolxprop} ci-dessus, et l'on suppose de plus que 2) est satisfaite ; nous allons montrer que sous cette dernière hypothèse,~$X$ est un domaine {\em affinoïde} de~$\sch X\an$. 

Soit~$\sch Y$ une composante irréductible de~$\sch X$ ; vérifions qu'il existe au moins une branche~$b$ dans~$\sch B$ telle que~$\Omega_b\cap \sch Y\an\neq \emptyset$. On suppose que ce n'est pas le cas. Le compact~$\sch Y\an$ ne rencontrant alors aucun des~$\Omega_b$, il est contenu dans~$X$ ; c'en est un fermé de Zariski irréductible et de dimension~$1$, et partant une composante irréductible ; mais ceci implique, en vertu de l'hypothèse 2), que~$S\cap \sch Y\an\neq \emptyset$. Choisissons~$x\in \sch Y\an \cap S$. Appartenant à~$S$, le point~$x$ est un point lisse de~$\sch X\an$ ; par conséquent,~$\sch Y\an$ est un voisinage de~$x$ dans~$\sch X\an$. Comme~$S$ est exactement le bord topologique de~$X$ dans~$\sch X\an$, il existe~$b\in \sch B$ telle que~$x$ adhère à~$\Omega_b$ ; dès lors, le voisinage~$\sch Y\an$ de~$x$ rencontre~$\Omega_b$, ce qui est absurde. 

\medskip
Les hypothèses de la propostion~\ref{ablouvnorm} sont donc vérifiées ; par conséquent,~$X=\sch X\an-\coprod \Omega_b$ est un domaine affinoïde de~$\sch X\an$.

\trois{redqlassb} {\em Preuve de 2)$\Rightarrow$1) dans le cas général.} On ne suppose plus que~$X$ est génériquement quasi-lisse, et l'on se place sous l'hypothèse 2). 

\medskip
Comme~$X$ est compacte, il existe une extension complète~$F$ de~$k$, qui est égale à~$k$ si celui-ci est de caractéristique nulle, et de la forme~$k^{1/p^n}$ pour un certain~$n\in\NN$ sinon, telle que~$X_{F,{\rm red}}$ soit géométriquement réduite. 

\medskip
Si~$x\in X$, il possède un unique antécédent~$x_F$ sur~$X_F$ ; les extensions~$\gred F/\gred k$ et~$\gred {\kappa(x_F)}/\gred {\kappa(x)}$ étant radicielles,~$\PP_{\gred{\kappa(x_F)}/\gred F}\to \PP_{\gred{\kappa(x)}/\gred k}$ est un homéomorphisme ; combiné au critère de Temkin, ce fait assure que~$x$ appartient au bord analytique de~$X$ si et seulement si~$x_F$ appartient au bord analytique de~$X_{F,\rm red}$.

\medskip
Il s'ensuit que l'espace~$F$-analytique~$X_{F,\rm red}$ satisfait lui aussi 2) ; en vertu du cas génériquement quasi-lisse déjà traité au~\ref{altaffcasql},~$X_{F,\rm red}$ est affinoïde. On déduit alors du théorème~\ref{affred} que~$X_F$ est affinoïde, puis du théorème~\ref{affrad} que~$X$ lui-même est affinoïde.~$\Box$

\deux{corollaletrn} {\bf Corollaire.} {\em Si~$X$ est une courbe~$k$-analytique compacte et irréductible, alors~$X$ est ou bien projective, ou bien affinoïde.}~$\Box$

\section{Quelques compléments}

\subsection*{Branches non discales, nœuds du squelette analytique et morphismes finis}

\deux{defbranchedisc} Soit~$X$ une courbe~$\KK$-analytique quasi-lisse, soit~$x\in X\typ{23}$ et soit~$b\in \br X x$. On dira que la branche~$b$ est {\em discale (relativement à~$X$)} si la composante connexe~$U$ de~$X\setminus\{x\}$ contenant~$b$ est un disque. Notons que comme~$b$ est issue de~$x$, le point~$x$ appartient à~$\partial U$. Il s'ensuit, en vertu de la proposition~\ref{propcarintdc}, que~$b$ est discale sauf si l'une au moins des trois propriétés suivantes est satisfaite : 

\medskip
$\bullet$~$\overline U=U\cup\{x\}$ n'est pas compacte ; 

$\bullet$~$\overline U$ n'est pas un arbre, c'est-à-dire encore contient une boucle ; 

$\bullet$~$U$ contient un point de genre strictement positif.

\medskip

\deux{exbrnondisc} Soit~$X$ une courbe~$\KK$-analytique quasi-lisse, soit~$x\in \skelan X$ et soit~$b\in \br {\skelan X} x$. La branche~$b$ est alors non discale. En effet, soit~$U$ la composante connexe de~$X\setminus\{x\}$ contenant~$b$. Comme~$b\in \br {\skelan X} x$, il existe un intervalle ouvert~$I$ contenu dans~$\skelan X$, aboutissant proprement à~$x$ et définissant~$b$. Cet intervalle est nécessairement contenu dans~$U$, et l'inclusion~$I\subset \skelan X$ assure que~$U$ ne peut être un disque. 

\deux{carnoeudsintr} {\bf Lemme.} {\em Soit~$X$ une courbe~$\KK$-analytique quasi-lisse et soit~$x\in X\typ{23}$. Les propositions suivantes sont équivalentes : 

\medskip
i) le point~$x$ est un nœud de~$\skelan X$ ; 

ii) ou bien le point~$x$ est de genre strictement positif, ou bien il appartient à~$\partial \an X$, ou bien il existe trois branches deux à deux distinctes dans~$\br X x$ qui sont non discales. }

\medskip
{\em Démonstration.} Supposons que~$x$ soit un nœud de~$\skelan X$, qu'il n'appartienne pas à~$\partial\an X$ et soit de genre zéro. C'est alors un sommet topologique de~$\skelan X$. En vertu du corollaire~\ref{corollunivalskelan}, la valence de~$(\skelan X,x)$ vaut au moins 2 ; comme~$x$ est un sommet, elle vaut au moins 3, et il existe donc au moins trois éléments distincts dans~$\br{\skelan X}x$, qui peuvent être vues comme trois branches deux à deux distinctes issues de~$x$. En vertu de~\ref{exbrnondisc}, elles sont non discales, ce qui achève de montrer que i)$\Rightarrow$ii). 

\medskip
Supposons maintenant que ii) soit vraie. Si~$x$ est de genre strictement positif ou appartient à~$\partial \an X$, il appartient à~$\skelan X$ puisqu'il ne peut être situé sur un disque, et est alors par définition un nœud de~$\skelan X$. Supposons maintenant qu'il existe trois branches deux à deux distinctes dans~$\br X x$ qui sont non discales. 

\medskip
{\em Le point~$x$ appartient à~$\skelan X$.} Supposons en effet le contraire. Le point~$x$ est alors contenu dans un ouvert~$U$ de~$X$ qui est un disque, et si~$\omega$ désigne l'unique bout de ce dernier alors toutes les composantes connexes de~$U\setminus\{x\}$ à l'exception de celle contenant~$]x;\omega[$ sont des disques ; comme elles sont relativement compactes dans~$U$, ce sont encore des composantes connexes de~$X\setminus\{x\}$, et~$\br X x=\br U x$ compte ainsi au plus une branche qui n'est pas discale relativement à~$X$, à savoir éventuellement celle définie par~$]x;\omega[$ ; on aboutit donc à une contradiction.

\medskip
{\em Le point~$x$ est un sommet topologique, et en particulier un nœud, de~$\skelan X$.} En effet, dans le cas contraire, il existerait un intervalle ouvert~$I$ tracé sur~$\skelan X$, ouvert et relativement compact dans ce dernier, et n'en contenant aucun nœud. Si~$X'$ désigne la composante connexe de~$X$ contenant~$x$ alors comme~$\skelan {X'}=\skelan X\cap X'$ est non vide,~$\skelan X'$ est un sous-graphe analytiquement admissible de~$X'$, et il découle de~\ref{pretriang} que~$I^\flat$ est une couronne. Si~$I'$ et~$I''$ désignent les deux composantes connexes de~$I\setminus\{x\}$ alors toutes les composantes connexes de~$I^\flat\setminus\{x\}$ à l'exception de celle contenant~$I'$ et~$I''$ sont des disques ; comme elles sont relativement compactes dans~$I^\flat$, ce sont encore des composantes connexes de~$X\setminus\{x\}$, et~$\br X x=\br {I^\flat} x$ compte ainsi au plus deux branches qui ne sont pas discales relativement à~$X$, à savoir éventuellement celles définies par~$I'$ et~$I''$. On aboutit donc à une contradiction.~$\Box$  

\deux{finiplatnondisc} {\bf Lemme.} {\em Soit~$\phi : Y\to X$ un morphisme fini et plat entre courbes~$\KK$-analytiques quasi-lisses, soit~$y\in Y\typ{23}$ et soit~$x$ son image sur~$X$. Soit~$b$ appartenant à~$\br Y y$. Si~$b$ est discale,~$\phi(b)$ l'est également.}

\medskip
{\em Démonstration.} On prouve la contraposée ; on suppose donc que~$\phi(b)$ n'est pas discale, et nous allons montrer que~$b$ ne l'est pas non plus. On désigne par~$U$ la composante connexe de~$X\setminus\{x\}$ qui contient~$\phi(b)$, par~$V$ la composante connexe de~$\phi^{-1}(U)$ qui contient~$b$, et par~$V'$ la composante connexe de~$Y\setminus\{y\}$ qui contient~$b$ (ou~$V$, ce qui revient au même). Comme~$\phi(b)$ n'est pas discale, plusieurs cas peuvent se présenter, que l'on traite séparément. 

\medskip
{\em Premier cas : la composante~$U$ contient un point~$z\in \partial \an X$.} Comme~$\phi$ est fini et plat,~$V$ se surjecte sur~$U$, et il existe donc un antécédent~$t$ de~$z$ sur~$V$ ; la finitude de~$\phi$ assure que~$t\in \partial\an Y$. Ainsi,~$V'$ rencontre~$\partial \an Y$ et~$b$ n'est pas discale. 

\medskip
{\em Deuxième cas : la composante~$U$ contient un point~$z$ de genre strictement positif.} Là encore,~$z$ possède un antécédent~$t$ sur~$V$. Comme~$z$ est de type 2,~$t$ est de type 2 ; comme le genre du corps résiduel~$\red{\hres(x)}$ est strictement positif, celui de~$\red{\hres(t)}$ l'est aussi par le théorème de Lüroth. Ainsi,~$V'$ contient un point de genre strictement positif, et~$b$ n'est pas discale. 

\medskip
{\em Troisième cas : il existe une boucle~$C$ tracée sur~$\overline U$.} L'image réciproque~$C'$ de~$C$ sur~$Y$ est un graphe localement fini et la surjection~$C'\to C$ est injective par morceaux (th.~\ref{theomodimrecg} et th.~\ref{modnonconst}) ; par compacité de~$\phi$, le graphe~$C'$ est compact et fini. Comme~$C'\to C$ est ouverte et comme~$C$ est un cercle, le graphe~$C'$ ne comprend que des points pluribranches. Par conséquent,~$C'\cap (V'\cup\{y\})$ est un sous-graphe compact et fini du graphe connexe~$V'\cup\{y\}$ qui rencontre~$V'$ (puisque~$C$ rencontre~$U$) et dont tous les points situés sur~$V'$ sont pluribranches. Cela implique que~$C'$ contient une boucle, et~$b$ n'est pas discale. 

\medskip
{\em Quatrième cas :~$\overline U$ est non compacte.} Le fermé~$\overline U$ de~$X$ est contenu dans~$\phi(\overline{V'})$, qui ne peut dès lors pas être compact. Par conséquent~$\overline{V'}$ n'est pas compact, et~$b$ n'est pas discale.~$\Box$ 

\deux{transfernoeuds} {\bf Proposition.} {\em Soit~$\phi : Y\to X$ un morphisme fini et plat entre courbes~$\KK$-analytiques quasi-lisses, soit~$y\in Y$ et soit~$x$ son image sur~$X$. Supposons que~$x$ est un nœud de~$\skelan X$. Dans ce cas~$y$ est un nœud de~$\skelan Y$.} 

\medskip
{\em Démonstration.} Comme~$x$ est un nœud de~$\skelan X$ on se trouve nécessairement dans l'un des trois cas suivants. 

\medskip
{\em Premier cas : le point~$x$ appartient à~$\partial\an X$.} Dans ce cas~$y\in \partial\an Y$ et~$y$ est un nœud de~$\skelan Y$. 

\medskip
{\em Deuxième cas : le point~$x$ est de genre strictement positif.} Dans ce cas il est de type 2, et~$y$ aussi ; comme le genre du corps résiduel~$\red{\hres(x)}$ est strictement positif, celui de~$\red{\hres(y)}$ l'est aussi par le théorème de Lüroth. Ainsi,~$y$ est de genre strictement positif et est donc un nœud de~$\skelan Y$. 

\medskip
{\em Troisième cas : il existe au moins trois branches de~$X$ issues de~$x$ qui sont non discales.} Si~$b$ est une branche issue de~$x$, elle est de la forme~$\phi(\beta)$ pour une certaine~$\beta\in \br Y y$ et si~$b$ est non discale alors~$\beta$ est non discale en vertu du lemme~\ref{finiplatnondisc}. Il s'ensuit qu'il existe au moins trois branches de~$Y$  issues de~$y$ qui sont non discales, et~$y$ est un nœud de~$\skelan Y$.~$\Box$ 

\section{Espaces analytiques formels}

\deux{defdomform} Soit~$X$ un espace~$k$-affinoïde ; nous dirons qu'un domaine analytique~$Y$ de~$X$ est {\em spécial} si~$Y$ est réunion finie de domaines affinoïdes de~$X$ de la forme~$${\sch D}(f):=\{\xi\in X, |f(\xi)|=||f||_\infty\},$$ où~$f$ est une fonction analytique sur~$X$ et où~$||f||_\infty$ désigne sa semi-norme spectrale. 

\medskip
Il revient au même de demander que~$Y$ soit l'image réciproque par la flèche de spécialisation d'un ouvert de Zariski de la réduction graduée~$\gred X$ de~$X$. 

\medskip
Tout domaine analytique spécial de~$X$ est compact. Toute réunion (resp. intersection finie) de domaines spéciaux de~$X$ est un domaine spécial de~$X$.

\deux{defatlform} Soit~$X$ un espace~$k$-analytique. Un {\em atlas affinoïde formel} sur~$X$ est un ensemble~$\sch E$ de domaines affinoïdes de~$X$ telle que :

\medskip
i)~$\sch E$ constitue un G-recouvrement de~$X$ ; 

ii) pour tout couple~$(V,W)$ d'éléments de~$\sch E$ l'intersection~$V\cap W$ est un domaine analytique spécial de~$V$ et de~$W$.

\medskip
On dit que deux atlas affinoïdes formels~${\sch E}$ et~${\sch E}'$ sur~$X$ sont {\em équivalents} si leur réunion est un atlas affinoïde formel ; il revient au même de demander que pour tout couple~$(V,W)$ formé d'un élément de~$\sch E$ et d'un élément de~${\sch E}'$, l'intersection~$V\cap W$ soit un domaine formel de~$V$ et de~$W$. On définit ainsi une relation d'équivalence sur l'ensemble des atlas affinoïdes formels de~$X$. 

\medskip
Un atlas affinoïde formel sur~$X$ est dit {\em maximal} s'il est maximal pour l'inclusion (en tant que sous-ensemble de~${\sch P}(X)$) parmi tous les atlas affinoïdes formels de~$X$ ; toute classe d'équivalence d'atlas affinoïdes formels possède un unique élément maximal, à savoir la réunion de tous les atlas qu'elle contient. 

\deux{defstructform} Un {\em espace~$k$-analytique formel} est un couple~$(X,\got A)$ où~$X$ est un espace~$k$-analytique et~$\got A$ un atlas affinoïde formel maximal de~$X$. 

\deux{exaffform} Soit~$X$ un espace affinoïde et soit~$\got A$ l'atlas affinoïde formel maximal équivalent à~$\{X\}$ ; il est immédiat que~$\got A$ est l'ensemble des domaines affinoïdes spéciaux de~$X$. On dira que~$\got A$ est l'atlas affinoïde formel {\em canonique} de~$X$, et que~$(X,\got A)$ est l'espace~$k$-analytique formel {\em canoniquement associé à~$X$.}

\deux{structext} Si~$(X, \got A)$ est un espace~$k$-analytique formel et si~$L$ est une extension complète de~$k$, l'espace~$X_L$ hérite d'une structure naturelle d'espace~$L$-analytique formel : c'est celle définie par l'atlas maximal équivalent à~$\{V_L\}_V$, où~$V$ parcourt~$\got A$ ; cet atlas sera noté~$\got A_L$. 

\deux{shilovindep} Soit~$X$ un espace~$k$-affinoïde et soit~$f$ une fonction analytique sur~$X$. La flèche induite~$\gred {{\sch D}(f)}\to \gred X$ entre les réductions {\em graduées} des espaces en jeu identifie~$\gred {{\sch D}(f)}$ à l'ouvert de Zariski~$D(\gred f)$ de~$\gred X$ ; il en résulte que si~$x\in {\sch D}(f)$ alors~$x$ appartient au bord de Shilov de~${\sch D}(f)$ si et seulement si il appartient au bord de Shilov de~$X$.

\deux{sommatlform} Soit~$X$ un espace~$k$-analytique, soit~$\got A$ un atlas affinoïde formel sur~$X$ (non nécessairement maximal), soit~$x\in X$ et soient~$V$ et~$W$ deux éléments de~$\got A$. Supposons que~$x\in V\cap W$. Les assertions suivantes sont alors équivalentes :

\medskip
i) le point~$x$ appartient au bord de Shilov de~$V$ ;

ii) le point~$x$ appartient au bord de Shilov de~$W$. 

\medskip
Pour le voir il suffit, par symétrie, de montrer que~${\rm i)}\Rightarrow {\rm ii)}$. Supposons donc que i) soit vraie ; choisissons deux domaines affinoïdes~$V'$ et~$W'$ de~$V\cap W$ qui contiennent~$x$, le premier (resp. le second) étant de la forme~${\sch D}(f)$ (resp.~${\sch D}(g)$) où~$f$ (resp.~$g$) est une fonction analytique sur ~$V$ (resp.~$W$). 

\medskip
On déduit successivement du~\ref{shilovindep} : que~$x$ appartient au bord de Shilov de~$V'$ ; qu'il appartient au bord de Shilov de~$V'\cap W'={\sch D}(g_{|V'})$ ; qu'il appartient au bord de Shilov de~$W'$, puisque~$V'\cap W'={\sch D}(f_{|W'})$ ; et enfin, qu'il appartient au bord de Shilov de~$W$, ce qu'on souhaitait établir. 

Ainsi, si~$x$ est un point de~$X$, les propriétés suivantes sont équivalentes : 

\medskip
a) il existe~$V\in \got A$ tel que~$x$ appartienne au bord de Shilov de~$V$ ; 

b) pour tout~$V\in \got A$ tel que ~$x\in V$, le point~$x$ appartient au bord de Shilov de~$V$. 

\medskip
Si elles sont satisfaites, nous dirons que~$x$ est un {\em sommet} de~$\got A$. L'ensemble des sommets de~$\got A$ est une partie fermée et discrète de~$X$ : cela provient du fait que~$X$ est G-recouvert par les domaines affinoïdes appartenant à~$\got A$, et de la finitude du bord de Shilov d'un tel domaine. 

\medskip
Si~$\got B$ est un atlas affinoïde formel sur~$X$ qui est équivalent à~$\got A$, l'ensemble des sommets de~$\got A$ est égal à celui des sommets de~$\got B$. Pour le voir il suffit, par symétrie, de s'assurer que tout sommet de~$\got A$ est un sommet de~$\got B$. Soit donc~$x$ un sommet de~$\got A$, et soit~$V$ (resp.~$W$) un élément de~$\got A$ (resp.~$\got B$) contenant~$x$. Comme~$x$ est un sommet de~$\got A$, il appartient au bord de Shilov de~$V$. En appliquant l'équivalence entre les assertions a) et b) ci-dessus à l'atlas~$\got A\cup \got B$, on voit que~$x$ appartient au bord de Shilov de~$W$ ; par conséquent, c'est un sommet de~$\got B$.

\medskip
Si~$X$ est un espace affinoïde et si~$\got A$ désigne son atlas affinoïde formel canonique, on déduit de ce qui précède que l'ensemble des sommets de~$\got A$ coïncide avec l'ensemble des sommets de l'atlas~$\{X\}$, autrement dit avec le bord de Shilov de~$X$.

\deux{strictaffform} Soit~$(X,\got A)$ un espace~$k$-analytique formel. Comme un espace~$k$-affinoïde est strictement~$k$-affinoïde si et seulement si la semi-norme spectrale de son algèbre de fonctions est à valeurs dans le monoïde~$\{0\}\cup\sqrt {|k\ti|}$, les propriétés suivantes sont équivalentes : 

\medskip
i) il existe un atlas affinoïde formel sur~$X$ contenu dans~$\got A$ et dont les éléments sont strictement~$k$-affinoïdes ;

ii) pour tout sommet~$x$ de~$\got A$, on a~$|\hres(x)\ti|\subset \sqrt {|k\ti|}$ ;

iii) tout élément de~$\got A$ est strictement~$k$-affinoïde. 

\medskip
Si ces conditions sont satisfaites, on dira que~$(X,\got A)$ est un espace {\em strictement~$k$-analytique formel}. 

\deux{introdefdist} Soit~$(X,\got A)$ un espace strictement~$k$-analytique formel. On peut lui associer un schéma formel~$\form X A$, à savoir le~$k\zero$-schéma formel plat obtenu par recollement des~$\spf {\sch O}_X(V)\zero$ où~$V$ parcourt~$\got A$. 

\trois{condfinitudespf} Le~$k\zero$-schéma formel~$\form X A$ n'est pas en général localement topologiquement de présentation finie ; mais il l'est toutefois dans chacune des deux situations suivantes : 

\medskip
$\bullet$~$|k\ti|=\{1\}$ ; 

$\bullet$~$X$ est réduit et~$|k\ti|$ est libre de rang 1 ; 
 
$\bullet$~$X$ est réduit et~$k$ est algébriquement clos.

\medskip
Si~$\form X A$ est localement topologiquement de présentation finie alors~$X$ s'identifie à sa fibre générique~$\form X A_\eta$.

\trois{defspfs} On définit par recollement des~$\spec(\sch O(V)\zero/\sch O(V)\zeroo)$ la {\em fibre spéciale}~$\form X A _s$ de~$\form X A$ ; c'est un schéma localement de type fini et réduit sur~$\red k$. 

On construit, également par recollement, une flèche dite {\em de spécialisation}~$\rho : X\to\form X A_s$ qui est anticontinue.

\trois{spfpascommut} On prendra garde que si~$L$ est une extension complète de~$k$, la flèche naturelle~$\fort {X_L}{\got A_L}\to \form X A _{L\zero}$ n'est pas un isomorphisme en général.

\trois{defdist} Lorsque~$|k\ti|\neq \{1\}$, on dira que l'espace strictement~$k$-analytique formel~$(X,\got A)$ est {\em distingué} si chacun des éléments de~$\got A$ est distingué ;  on peut se contenter de le tester sur les éléments de n'importe quel atlas affinoïde formel de~$X$ contenu dans~$\got A$.  Si~$(X,\got A)$ est distingué, le schéma formel~$\form X A$ est localement topologiquement de présentation finie, et~$\form X A_s$ s'identifie à~$\form X A \otimes_{k\zero}\red k$ ; si de plus le~$\red k$ -schéma réduit~$\form X A  _s$ est {\em  géométriquement} réduit (c'est automatiquement le cas si~$\red k$ est parfait), alors pour toute extension complète~$L$ de~$k$ l'espace~$L$-analytique formel~$(X_L,\got A_L)$ est distingué, et~$\fort {X_L}{\got A_L}\simeq \form X A _{L\zero}$. 

\trois{suffdist} Soit~$S$ l'ensemble des sommets de~$\got A$ ; considérons les trois conditions suivantes : 

\medskip
$\alpha)$~$X$ est réduit ;

$\beta)$~$|\hres(x)\ti|=|k\ti|$ pour tout~$x\in S$. 

$\gamma)$ toute extension finie de~$k$ est sans défaut.

\medskip
Remarquons que~$\beta)$ est satisfaite si~$|k\ti|$ est divisible, et que~$\gamma)$ l'est si~$k$ est algébriquement clos ou de valuation discrète.

\medskip
Supposons que~$|k\ti|\neq \{1\}$. Alors pour que~$(X,\got A)$ soit distingué, il faut que~$\alpha)$ et~$\beta)$ soient satisfaites, et il suffit que~$\alpha)$,~$\beta)$ et~$\gamma)$ le soient. En particulier, si~$k$ est algébriquement clos et si~$X$ est réduit alors~$(X,\got A)$ est distingué. 

\deux{relnormal} Soit~$\got X$ un~$k\zero$-schéma formel localement topologiquement de présentation finie. Les assertions suivantes sont équivalentes : 

\medskip
i)~$\got X$ est de la forme~$\form{\got X_\eta} A$ pour un certain atlas strictement~$k$-analytique formel~$\got A$ sur~$\got X_\eta$ ; 

ii)~$\got X$ est plat et relativement normal dans~$\got X_\eta$ (si~$\got X_\eta$ est normal, cela revient à demander que~$\got X$ soit plat et normal). 

\medskip
On suppose qu'elles sont satisfaites. 

\trois{uniquatlas} L'atlas~$\got A$ de i) est alors unique : c'est nécessairement l'ensemble~$\mathsf{Aff}({\got X})$ des~$\got U_\eta$, pour~$\got U$ parcourant l'ensemble des ouverts formels affines de~$\got X$. 

\trois{xsxredk} Posons~$\got X_{\red k}=\got X\otimes_{k\zero}\red k$ ; on dispose d'une flèche de réduction anti-continue~$\got X_\eta\to \got X_{\red k}$. Comme le faisceau~$k\zeroo \sch O_{\got X}$ est un idéal de définition de~$\got X$, la fibre spéciale~$\got X_s$ est isomorphe à~$(\got X_{\red k})_{\rm red}$ ; on dispose en particulier d'une identification canonique entre les espaces topologiques~$\got X_s$ et~$\got X_{\rm red}$ modulo laquelle les flèches de réduction de~$\got X_\eta\to \got X_s$ et~$\got X_\eta \to \got X_{\red k}$ coïncident. 

\trois{normimplred} Si~$|k\ti|\neq \{1\}$ alors l'espace~$\got X_\eta$ est réduit. Pour le voir, on se ramène aussitôt au cas où~$\got X$ est affine, donc égal à~${\rm Spf}\; \mathsf A$ pour une certaine~$k\zero$-algèbre~$\mathsf A$ plate, topologiquement de présentation finie, et relativement normale dans son localisé~$\mathsf A\otimes_{k\zero} k$. Choisissons un épimorphisme~$k\zero\{T_1,\ldots, T_n\}\to \mathsf A$. Il induit un épimorphisme~$k\{T_1,\ldots, T_n\}\to \mathsf A\otimes_{k\zero}k$, et l'on note~$||.||$ la norme quotient correspondante ; par définition de~$||.||$, on a~$||a||\leq 1$ pour tout~$a\in \mathsf A$. Soit~$f$ un élément nilpotent de~$\mathsf A\otimes_{k\zero} k$. Pour tout~$\lambda\in k$ l'élément~$\lambda f$ de~$\mathsf A\otimes_{k\zero}k$ est nilpotent, et partant entier sur~$\mathsf A$ ; il appartient dès lors à~$\mathsf A$. Ceci implique que~$||\lambda.f||\leq 1$ pour tout~$\lambda\in k$, et donc que~$||f||=0$ puisque~$|k\ti|\neq \{1\}$. Ainsi,~$f=0$ et~$\mathsf A\otimes_{k\zero}k$ est réduit, ce qu'il fallait démontrer. 

\trois{cadisting} On fait toujours l'hypothèse que~$|k\ti|\neq\{1\}$. L'espace strictement~$k$-analytique formel~$(\got X_\eta,\mathsf{Aff}({\got X}))$ est alors distingué si et seulement si~$\got X_{\red k}$ est réduit, c'est-à-dire égal à~$\got X_s$.

\deux{typesomm} Soit~$(X,\got A)$ un espace~$k$-analytique formel. Si~$X$ est purement de dimension~$d$, alors pour tout sommet~$x$ de~$\got A$ l'on a~$${\rm deg.}\;{\rm tr.} \;\red{\hres(x)}/\red k+\dim \QQ \QQ\otimes_\ZZ\left(|\hres(x)\ti|/|k\ti|\right)=d.$$ Dans le cas particulier où~$d=1$, cela signifie que~$x$ est de type 2 ou 3 ; dans cette situation,~$|\hres(x)\ti|\subset \sqrt {|k\ti|}$ si et seulement si~$x$ est de type 2. 

\medskip
Autrement dit, une courbe~$k$-analytique formelle est {\em strictement}~$k$-analytique formelle si et seulement si tous ses sommets sont de type 2.   

\deux{domanaform} Soit~$(X,\got A)$ un espace~$k$-analytique formel et soit~$S$ l'ensemble des sommets de~$\got A$. Si~$Y$ est un sous-ensemble de~$X$ qui est réunion d'éléments de~$\got A$, c'est un domaine analytique fermé de~$X$ : en effet, si~$V\in \got A$ l'intersection~$Y\cap V$ est alors réunion de domaines analytiques spéciaux de~$V$, et est dès lors un domaine analytique spécial de~$V$ ; c'est en particulier un domaine analytique compact de~$V$, d'où notre assertion. 

\medskip
Un domaine analytique~$Y$ de~$X$ qui est réunion d'éléments de~$\got A$ sera dit {\em~$\got A$-formel.} Si~$(X,\got A)$ est strictement~$k$-analytique formel, il revient au même de demander que~$Y$ soit l"image réciproque par la flèche de spécialisation d'un ouvert de Zariski de~$\form X A_s$.

\medskip
Si~$Y$  est un domaine analytique~$\got A$-formel de~$X$ l'ensemble~$\got A_{|Y}$ des domaines affinoïdes appartenant à~$\got A$ et contenus dans~$Y$ est un atlas affinoïde formel maximal sur~$Y$ ;  l'ensemble des sommets de~$\got A_{|Y}$ est précisément~$S\cap Y$. 

\medskip
Un domaine affinoïde~$V$ de~$X$ est~$\got A$-formel si et seulement si il appartient à~$\got A$ ; si c'est le cas,~$\got A_{|V}$ est l'atlas affinoïde formel canonique de~$V$. 

\deux{fonctanaform} Soient~$(Y,\got B)$ et~$(X,\got A)$ deux espaces~$k$-analytiques formels. Un {\em morphisme} d'espaces~$k$-analytiques formels de~$(Y,\got B)$ vers~$(X,\got A)$ est un morphisme~$\phi: Y\to X$ d'espaces~$k$-analytiques tel que pour tout domaine analytique~$\got A$-formel~$V$ de~$X$, l'image réciproque~$\phi\inv(V)$ soit un domaine analytique~$\got B$-formel de~$Y$. 

Il suffit de vérifier cette propriété sur un G-recouvrement de~$X$ par des domaines analytiques~$\got A$-formels, et donc par exemple sur un atlas affinoïde formel de~$X$ contenu dans~$\got A$. 

\medskip
Si~$(X, \got A)$ et~$(Y,\got B)$ sont deux espaces {\em strictement}~$k$-analytiques formels, un morphisme de~$(Y,\got B)$ vers~$(X,\got A)$ sera simplement un morphisme d'espaces~$k$-analytiques formels de~$(Y,\got B)$ vers~$(X,\got A)$. 

\trois{morstrictaff} La flèche~$(X,\got A)\mapsto \form X A$ définit un foncteur de la catégorie des espaces strictement ~$k$-analytiques formels vers celle des~$k\zero$-schémas formels plats, foncteur qui est pleinement fidèle. 

\medskip
Avant de dire quelques mots à propos de son image essentielle dans certains cas particuliers, introduisons une convention et une notation : on dira qu'une courbe~$k$-analytique formelle~$(X,\got A)$ est réduite (resp. propre) si~$X$ est réduite (resp. propre) ; et l'on fixe un quasi-inverse~$\got X\mapsto \got X^{\rm alg}$ du foncteur~$\sch X\mapsto \widehat{\sch X}$, vu comme allant de la catégorie des~$k\zero$-courbes propres vers les~$k\zero$-courbes formelles propres.

\trois{equivvariees} Supposonsque~$|k\ti|\neq\{1\}$ et que~$k$ est ou bien algébriquement clos, ou bien de valuation discrète ({\em i.e.}~$|k\ti|$ est libre de rang 1). Il résulte alors de~\ref{condfinitudespf} et~\ref{relnormal} que~$(X,\got A)\mapsto \form X A$ induit une équivalence entre la catégorie des courbes strictement~$k$-analytiques formelles réduites et celle des~$k\zero$-schémas formels plats, localement de présentation finie, et relativement normaux dans leur fibre générique. Un quasi-inverse est fourni par~$\got X\mapsto (\got X_\eta,\mathsf {Aff}({\got X}))$. 

En vertu de GAGA et de l'algébrisation des courbes formelles propres (th.~\ref{algebrformcurv}, rem.~\ref{remgotxalg}), il s'ensuit que~$(X,\got A)\mapsto \form X A^{\rm alg}$ établit une équivalence entre la catégorie des courbes strictement~$k$-analytiques formelles réduites et propre et celle des~$k\zero$-schémas propres et plats relativement normaux dans leur fibre générique. Un quasi-inverse est fourni par~$\sch X\mapsto (\sch X_\eta\an=\widehat{\sch X}_\eta,\mathsf{Aff}(\widehat{\sch X}))$.

\subsection*{Reconstitution d'une structure formelle sur une courbe à partir de ses sommets}

\deux{remcourb} Soit~$X$ une courbe~$k$-affinoïde et soit~$S$ son bord de Shilov. Soit~$\rho$ la flèche de spécialisation~$X\to \gred X$. Soit~$\xi\in \gred X$. On est dans l'un des deux cas suivants :  

\medskip
$\bullet$~$\xi$ est un point générique de~$\gred X$ ; dans ce cas ~$\rho\inv(\xi)$ est de la forme~$\{x\}$ avec~$x\in S$, et tout élément de~$S$ s'obtient de cette manière (\ref{genshil}) ; 

$\bullet$ le corps résiduel gradué de~$\xi$ est fini sur~$\gred k$ ; dans ce cas~$\rho\inv(\xi)$ est un ouvert connexe et non vide de~$X$ (\ref{tubecon} dont le bord consiste exactement en les antécédents des points génériques des composantes irréductibles contenant~$\xi$ (\ref{tubeadh}). 

\medskip
Il découle de ce qui précède qu'un compact~$W$ de~$X$ en est un domaine formel si et seulement si il existe un sous-ensemble~$S'$ de~$S$ et un sous-ensemble~$\Pi$ de~$\pi_0(X\setminus S)$ tels que : 

\medskip
i)~$W=S'\cup \bigcup\limits_{V\in \Pi}V$ ;

ii) pour tout~$V\in \Pi$ l'on a~$\partial V\subset S'$ ; 

iii) presque toutes les composantes connexes de~$X\setminus S$ de bord contenu dans~$S'$ appartiennent à~$\Pi$.

\deux{structsomm} {\bf Théorème.} {\em Soit~$X$ une courbe~$k$-analytique. Soit~$S$ un ensemble fermé et discret de~$X$ contenu dans~$X\dtr$. Les propositions suivantes sont équivalentes :

\medskip
\begin{itemize}
\item[i)] il existe un atlas~$k$-affinoïde formel sur~$X$ dont~$S$ est l'ensemble des sommets ;  

\medskip
\item[ii)]~$S$ satisfait les conditions suivantes : 

\medskip
\begin{itemize} 
\item[$\alpha)$]~$\partial\an X\subset S$ ;

\item[$\beta)$] toute composante irréductible de~$X$ rencontre~$S$ ;

\item[$\gamma)$] toute composante connexe de~$X\setminus S$ est relativement compacte ; 

\item[$\delta)$] pour tout~$x\in S\setminus \partial  \an X$ et pour toute composante connexe~$U$ de~$X\setminus S$ telle que~$x\in \partial U$ il existe une branche de~$X$ issue de~$x$ et non contenue dans~$U$. 

\end{itemize}
\end{itemize}
\medskip

De plus si ces conditions sont satisfaites la classe d'équivalence de l'atlas évoquée au i) est unique ; il existe en particulier un et un seul atlas affinoïde formel maximal sur~$X$ dont~$S$ est l'ensemble des sommets.} 

\medskip
{\em Démonstration.} On procède par double implication. 

\trois{sommrelcomp} {\em Preuve de~$i)\Rightarrow ii)$.} On suppose que i) est vérifiée ; soit~$\got A$ un atlas affinoïde formel sur~$X$ dont~$S$ est l'ensemble des sommets.

\medskip
Soit~$x\in \partial\an X$ ; il est situé sur un domaine affinoïde~$V$ de~$X$ qui appartient à~$\got A$ ; il appartient nécessairement au bord analytique de~$V$, c'est-à-dire à son bord de Shilov (puisqu'on est en dimension 1) ; ainsi,~$x\in S$, et~$\alpha)$ est vraie.

\medskip
Soit~$Y$ une composante irréductible de~$X$. Il existe un élément~$V$ de~$\got A$ rencontrant~$Y$ ; l'intersection~$V\cap Y$ est une réunion non vide de composantes irréductibles de~$V$ ; le bord de Shilov de chacune d'elles est non vide, et est contenu dans le bord de Shilov de~$V$, et partant dans~$S$ ; ainsi,~$Y$ rencontre~$S$, et~$\beta)$ est vraie. 

\medskip
Soit~$U$ une composante connexe de~$X\setminus S$ ; il existe un élément~$V$ de~$\got A$ rencontrant~$U$. Comme le bord topologique de~$V$ dans~$X$ est contenu dans le bord de Shilov de~$V$, et donc dans~$S$, l'intersection~$U\cap V$ est contenue dans l'intérieur de~$V$ ; c'est donc une partie non vide, ouverte et fermée de~$U$ ; par connexité de~$U$, il s'ensuit que~$U\cap V=U$, c'est-à-dire que~$U\subset V$ ; ainsi,~$U$ est relativement compacte et~$\gamma)$ est vraie. 

\medskip
Soit maintenant~$x\in (S\setminus \partial \an X)\cap \overline U$. Comme~$x\in V\cap S$, il appartient au bord de Shilov de~$V$ ; comme il n'appartient pas à~$\partial\an X$, il appartient à~$\partial V$. Il existe donc une branche issue de~$x$ qui est contenue dans~$X-V$, et n'est dès lors pas contenue dans~$U$ . Ainsi,~$\delta)$ est vraie.

\trois{relcompsomm} {\em Preuve de~$ii)\Rightarrow i)$.} Supposons que ii) soit vérifiée. Soit~$\got A$ l'ensemble des parties compactes~$W$ de~$X$ telles qu'il existe un sous-ensemble fini~$S'$ de~$S$ et une famille~$(U_i)$ de composantes connexes deux à deux disjointes de~$X\setminus S'$ possédant les propriétés suivantes :

\medskip
a)~$W=X-\coprod U_i$ ;

b) pour tout~$x\in S'$ l'ensemble des indices~$i$ tels que~$x\in \partial U_i$ est fini, et non vide si~$x\notin \partial\an X$ ; 

c) ~$(S\setminus S')\subset \coprod U_i$ ; 

d)  toute composante irréductible de~$X$ qui ne rencontre pas~$S'$ est contenue dans~$\coprod U_i$. 

\medskip
Nous allons démontrer : que~$\got A$ est un atlas affinoïde formel sur~$X$ ; que~$S$ est l'ensemble de ses sommets ; et que tout atlas affinoïde formel sur~$X$ dont~$S$ est l'ensemble des sommets est contenu dans~$\got A$. L'assertion d'unicité requise en résultera -- ainsi par ailleurs que le caractère maximal de~$\got A$. 

\medskip
\noindent
{\bf Première étape :~$\got A$ constitue un  G-recouvrement de~$X$.}

\medskip
{\bf \em Un procédé général de constructions de parties appartenant à~$\got A$.} Soit~$U$ une composante connexe de~$X\setminus S$. Comme~$U$ est relativement compacte et~$S$ discrète, le bord de~$U$ est un sous-ensemble fini~$S'$ de~$S$. Soit~$Z$ la réunion des composantes irréductibles de~$X$ qui ne rencontrent pas~$S'$ ; c'est un fermé de~$X$. Soit~$\Pi$ l'ensemble des composantes connexes de~$X\setminus S'$ qui rencontrent~$Z$ ou~$S\setminus S'$. Si~$x\in S'$ alors comme~$S$ est discrète et~$Z$ fermé, il existe un voisinage ouvert connexe de~$x$ dans~$X$ ne rencontrant pas~$(S\setminus S')\cup Z$ ; ceci entraîne que l'ensemble des~$V\in \Pi$ telles que~$S'\cap \partial V\neq \emptyset$ est fini. Soit~$\Pi'$ un ensemble fini de composantes connexes de~$X\setminus S'$ différentes de~$U$, tel que pour tout élément~$x$ de~$S'$ qui n'appartient pas à~$\partial\an X$, il existe~$V\in \Pi\cup \Pi'$ avec~$x\in \partial V$. Posons~$W=X-\coprod\limits_{V\in \Pi\cup\Pi'}V$ ; le fermé~$W$ de~$X$ est un graphe qui vérifie par construction les conditions a), b), c) et d) énoncées au début de~\ref{relcompsomm}. Par ailleurs, si~$V$ est une composante connexe de~$W\setminus S'$ alors~$V$ ne rencontre pas~$S$ ; c'est donc une composante connexe de~$X\setminus S$, et elle est de ce fait relativement compacte dans~$X$, et partant dans~$W$ ; on en déduit que~$W$ est compact, et finalement qu'il appartient à~$\got A$. 

Montrons que~$W$ contient~$U$. Cela revient à prouver que~$U\notin \Pi\cup\Pi'$. La définition même de~$\Pi'$ assure que~$U\notin \Pi'$. Supposons que~$U\in \Pi$ ; comme~$U$ ne rencontre pas~$S$ (c'est une composante connexe de~$X\setminus S$), elle rencontre nécessairement~$Z$ ; il existe donc une composante irréductible~$T$ de~$Z$ telle que~$T\cap U\neq \emptyset$. L'ouvert~$U\cap T$ de~$T$ ne peut être fermé dans~$T$ : sinon, par connexité de~$T$, ce serait~$T$ elle-même et l'on aurait ainsi~$T\subset U$, ce qui contredirait le fait que~$T$ rencontre~$S$ ; dès lors il y au moins un point de~$S '$ qui est situé sur~$T$, ce qui est absurde puisque~$T$ est une composante irréductible de~$Z$, et ne rencontre donc pas~$S'$. 

\medskip
{ \bf \em Application du procédé décrit ci-dessus.} Soit~$x\in X$.

\medskip
{\em Supposons que~$x\notin S$}.  Soit~$U$ la composante connexe de~$X\setminus S$ contenant~$x$. Soient~$S'$  et~$\Pi$ comme ci-dessus. Pour tout point~$y$ de~$S'\setminus \partial  \an X$, il existe (propriété~$\delta)$ de l'énoncé du théorème) une branche~$\beta_y$ issue de~$y$ et non contenue dans~$U$. Soit~$\Pi'$ l'ensemble des composantes connexes de~$X\setminus S'$ contenant les différentes~$\beta_y$ ; il satisfait les conditions énoncées ci-dessus ; il s'ensuit que~$W:=X-\coprod\limits_{V\in \Pi\cup\Pi'}V$ est un élément de~$\got A$ qui contient~$U$, et est en particulier un voisinage de~$x$.

\medskip
{\em Supposons que~$x\in S$ et que ce n'est pas un point isolé de~$X$}. Soit~$b$ une branche issue de~$x$ et soit~$U$ la composante connexe de~$X\setminus S$ contenant~$b$. Soient~$S'$ et~$\Pi$ comme ci-dessus. Pour tout point~$y$ de~$S'\setminus \partial  \an X$, il existe (propriété~$\delta)$  de l'énoncé du théorème) une branche~$\beta_y$ issue de~$y$ et non contenue dans~$U$. Soit~$\Pi'$ l'ensemble des composantes connexes de~$X-\setminus S'$ contenant les différentes~$\beta_y$ ; il satisfait les conditions énoncées ci-dessus ; il s'ensuit que~$W:=X-\coprod\limits_{V\in \Pi\cup\Pi'}V$ est un élément de~$\got A$ qui contient~$U$, laquelle contient~$b$ ; il contient par ailleurs presque toutes les composantes connexes de~$X\setminus\{x\}$. 

\medskip
{\em Supposons que~$x$ est un point isolé de~$X$.} C'est alors à la fois une composante connexe, une composante irréductible et un domaine affinoïde de~$X$ ; il s'ensuit immédiatement que~$\{x\}$ est un élément de~$\got A$. 

\medskip
On déduit de ce qui précède que l'ensemble~$\got A$ constitue un G-recouvrement de~$X$.

\medskip
\noindent
{\bf Seconde étape : les éléments de~$\got A$ sont des domaines affinoïdes de~$X$.} 

\medskip
Soit~$W$ un élément de~$\got A$ ; il résulte de la définition de~$\got A$ et de la proposition~\ref{propdomanferm} que~$W$ est un domaine analytique fermé ; nous allons montrer qu'il est affinoïde de bord de Shilov~$W\cap S$. Remarquons déjà que si l'on se donne~$S'$ comme dans la définition de~$\got A$ alors~$S'=W\cap S$, et~$S'$ est la réunion de~$\partial W$ et de~$W\cap \partial \an X$ ; par conséquent,~$S'=W\cap S$ est exactement le bord analytique de~$W$. 

\medskip
Pour établir le caractère affinoïde de~$W$ il suffit, grâce au théorème~\ref{altcompgen}, de démontrer que toute composante irréductible de~$W$ a un bord analytique non vide. Nous allons procéder par l'absurde ; on se donne donc une composante irréductible~$Z$ de~$W$ dont on suppose que le bord analytique est vide. Dans ce cas,~$Z$ est propre, et est dès lors une composante irréductible de~$X$. Par conséquent,~$Z$ rencontre~$S$ ; comme elle est contenue dans~$W$, cela signifie qu'elle rencontre~$W\cap S$,  qui est précisément le bord analytique de~$W$. On aboutit ainsi à une contradiction, et ~$W$ est donc bien domaine affinoïde de~$\sch X\an$, comme annoncé ; son bord de Shilov est son bord analytique, à savoir~$W\cap S$.

\medskip
\noindent
{\bf Preuve du fait que~$\got A$ est un atlas affinoïde formel.}

\medskip
Soient~$V$ et~$W$ deux domaines affinoïdes de~$X$ appartenant à~$\got A$. Nous allons vérifier que~$V\cap W$ est un domaine analytique spécial de~$V$, ce qui permettra de conclure (comme~$V$ et~$W$ jouent le même rôle, ce sera aussi un domaine analytique spécial de~$W$).

Les bords de Shilov respectifs de~$V$ et~$W$ sont, en vertu de ce qui précède,~$S\cap V$ et~$S\cap W$. Par définition de~$\got A$, on peut écrire~$V$ (resp.~$W$) comme la réunion de~$S\cap V$ (resp.~$S\cap W$) et de composantes connexes de~$X\setminus S$ ; de plus, pour tout~$x\in S\cap V$ (resp.~$S\cap W$), presque toutes les composantes connexes de~$X\setminus S$ qui aboutissent à~$x$ sont contenues dans~$V$ (resp.~$W$). 

\medskip
Par conséquent, le compact~$V\cap W$ est de la forme~$(S\cap V\cap W)\cup \bigcup\limits_{U\in \Pi} U$, où~$\Pi$ est un ensemble de composantes connexes de~$V-(S\cap V)$ (ou, ce qui revient au même, de composantes connexes de~$X\setminus S$ contenues dans~$V$) possédant les deux propriétés suivantes : pour tout~$U\in \Pi$, le bord de~$U$ est contenu dans~$S\cap V\cap W$ ; presque toute composante connexe de~$V-(S\cap V)$ dont le bord est contenu dans~$S\cap V\cap W$ appartient à~$\Pi$. 

\medskip
En vertu du~\ref{remcourb},~$V\cap W$ est un domaine spécial de~$V$.

\medskip
\noindent
{\bf Détermination des sommets de~$\got A$.}

\medskip
Soit~$V\in \got A$. On a vu plus haut que le bord de Shilov de~$V$ est~$V\cap S$ ; comme tout point de~$S$ appartient par ailleurs à au moins un domaine~$V$ de l'atlas~$\got A$, l'ensemble des sommets de~$\got A$ est exactement~$S$ ; ceci achève de démontrer l'équivalence entre i) et ii).

\medskip
\noindent
{\bf Tout atlas affinoïde formel sur~$X$ ayant~$S$ pour ensemble des sommets est contenu dans~$\got A$}. 

\medskip
Soit~$\got B~$ un atlas~$k$-affinoïde formel sur~$X$ dont~$S$ est l'ensemble des sommets et soit~$W$ un domaine appartenant à~$\got B~$ ; nous allons montrer que~$W\in \got A$. 

\medskip
Le bord de Shilov~$S'$ de~$W$ est la réunion de~$\partial W$ et de~$W\cap \partial \an X$ ; il coïncide par ailleurs avec~$S\cap W$, puisque~$S$ est l'ensemble des sommets de~$\got B~$. Par conséquent,~$X-W$ s'écrit sous la forme~$\coprod U_i$, où~$(U_i)$ est une famille de composantes connexes de~$X\setminus S'$. 

\medskip
Soit~$Y$ une composante irréductible de~$X$ rencontrant~$W$. L'intersection de~$Y$ et~$W$ est une réunion non vide de composantes irréductibles de~$W$, et elle contient donc au moins un point du bord de Shilov de~$W$ ; par conséquent,~$Y\cap S'\neq \emptyset$. Par contraposition, si~$Y$ est une composante irréductible de~$X$ qui ne rencontrent pas~$S'$ alors~$Y\subset \coprod U_i$. 

\medskip
Soit~$x\in S'$ ; si~$x\notin \partial W$ alors l'ensemble des indices~$i$ tels que~$x\in \partial U_i$ est vide ; sinon, la proposition~\ref{propdomanferm} assure que cet ensemble est fini ; compte-tenu du fait que~$S\setminus \partial \an X\subset \partial W$, on en déduit que~$W\in \got A$.~$\Box$

\deux{commentdelta} {\em Commentaire.} On conserve les notations du théorème. Si~$x$ est un point de type 2 de~$S$ et si~$U$ est une composante connexe de~$X\setminus S$ telle que~$x\in \partial U$ alors comme~$\br X x \ctd U$ est fini et comme~$\br X x$ est infini, il existe une branche~$b$ issue de~$x$ qui n'est pas contenue dans~$U$. Il suffit par conséquent de vérifier la condition~$\delta)$ de l'assertion ii) sur les points de type 3 de~$S$ ; en particulier si~$S$ ne contient que des points de type 2 alors~$\delta)$ est automatiquement satisfaite. 

Remarquons que si~$X$ est propre, si son squelette est un cercle et si~$x$ est un point de type 3 sur celui-ci alors~$S:=\{x\}$ satisfait les propriétés~$\alpha),\beta)$ et~$\gamma)$ de l'assertion ii), mais pas~$\delta)$ : le bord analytique de~$X$ est vide, mais l'ouvert~$X\setminus\{x\}$ est connexe et contient les deux branches de~$X$ issues de~$x$.

\deux{defsommit} Si~$X$ est une courbe~$k$-analytique, on appellera {\em sous-ensemble sommital de~$X$} tout sous-ensemble~$S$ de~$X$ qui est l'ensemble des sommets d'un atlas affinoïde formel sur~$X$ ; en vertu du théorème~\ref{structsomm} ci-dessus un sous-ensemble~$S$ de~$X$ est sommital si et seulement si il est fermé, discret, contenu dans~$X\dtr$, et satisfait les conditions a), b), c), d) de l'assertion ii) de {\em loc. cit.} Il y a alors, toujours en vertu de {\em loc. cit.}, un unique atlas affinoïde formel maximal dont l'ensemble des sommets est~$S$, et il contient (et est donc équivalent) à tout atlas analytique formel sur~$X$ dont l'ensemble des sommets est~$S$ : on le notera~$\mathsf{Atl}(S)$. Il est strictement analytique si et seulement si~$S$ est constitué de points de type 2. 

\medskip
Si~$\got A$ est un atlas analytique formel sur~$X$, on notera~$\mathsf{Somm}(\got A)$ l'ensemble de ses sommets. Les faits suivantes résultent des définitions : si~$S$ est un sous-ensemble sommital de~$X$ alors~$\mathsf{Somm}(\mathsf{Atl}(S))=S$ ; si~$\got A$ est un atlas~$k$-analytique formel sur~$X$ alors~$\mathsf{Atl}(\mathsf{Somm}(\got A))$ est l'unique atlas~$k$-analytique formel maximal contenant~$\got A$. 

\deux{basestruct} {\em Extension des scalaires}. Soit~$X$ une courbe~$k$-analytique et soit~$S$ un sous-ensemble sommital de~$X$. 

\medskip
Soit~$F$ une extension presque algébrique de~$k$. L'ensemble~$\mathsf{Somm}(\mathsf {Atl}(S)_F)$ est égal à l'image réciproque~$S_F$ de~$S$ sur~$X_F$. Par conséquent,~$S_F$ est une partie sommitale de~$X_F$ et~$\mathsf{Atl}(S_F)=\mathsf{Atl}(S)_F$. 

\medskip
Soit~$V$ un domaine analytique~$\mathsf{Atl}(S)$-formel de~$X$. L'ensemble des sommets de~$\mathsf{Atl}(S)_{|V}$ étant égal à~$S\cap V$, l'ensemble~$S\cap V$ est une partie sommitale de~$V$ et l'on a~$\mathsf{Atl}(S\cap V)=\mathsf{Atl}(S)_{|V}$.  

\subsection*{Fonctorialité} 

\deux{cardomform} {\bf Proposition.} {\em Soit~$X$ une courbe~$k$-analytique et soit~$S$ un sous-ensemble sommital de~$X$. Soit~$V$ un domaine analytique fermé de~$X$. Les propositions suivantes sont équivalentes : 

\medskip
i)~$V$ est un domaine analytique~$\mathsf{Atl}(S)$-formel de~$X$ ; 

ii)~$\partial \an V\subset S$.}

\medskip
{\em Démonstration.} Supposons que i) soit vraie. L'ensemble des sommets de~$\mathsf{Atl}(S)_{|V}$ est alors égal à~$V\cap S$, et l'on sait d'après le théorème~\ref{structsomm} qu'il contient~$\partial \an V$, d'où ii). 

\medskip
Supposons maintenant que ii) soit vraie. Le sous-ensemble~$S\cap V$ de~$V$ en est une partie fermée et discrète, constituée de points de type 2 ou 3. 

\trois{vintersok} {\em Le couple~$(V,S\cap V)$ satisfait les conditions~$\alpha), \beta), \gamma)$ et~$\delta)$ introduites dans l'énoncé du théorème ~\ref{structsomm} pour le couple~$(X,S)$.} Notons pour commencer qu'en vertu justement du théorème~\ref{structsomm}, le couple~$(X,S)$ satisfait~$\alpha), \beta), \gamma)$ et~$\delta)$. 

\medskip
L'hypothèse faite sur~$V$ signifie exactement que~$(V,S\cap V)$ satisfait~$\alpha)$. 

\medskip
Vérifions que~$(V,S\cap V)$ satisfait~$\beta)$ ; soit~$Z$ une composante irréductible de~$V$. Si~$\partial \an Z\neq \emptyset$ alors~$\partial \an Z\subset \partial \an V\subset S\cap V$ ; et si ~$\partial \an Z= \emptyset$ la courbe~$Z$ est propre, et~$Z$  est dès lors une composante irréductible de~$X$, qui rencontre dès lors~$S$ puisque~$(X,S)$ satisfait~$\beta)$, et plus précisément~$S\cap V$ puisque~$Z\subset V$. 

\medskip
Vérifions que~$(V,S\cap V)$ satisfait~$\gamma)$. Soit~$U$ une composante connexe de~$V\setminus S$. Elle est contenue dans une composante connexe de~$X\setminus S$ et est dès lors relativement compacte dans~$X$ puisque~$(X,S)$ satisfait~$\gamma)$ ; le domaine~$V$ étant fermé dans~$X$, il s'ensuit que~$U$ est relativement compacte dans~$V$.

\medskip
Vérifions que~$(V,S\cap V)$ satisfait~$\delta)$. Soit~$x\in (S\cap V)\setminus \partial \an V$ et soit~$U$ une composante connexe de~$V\setminus S$ telle que~$x\in \partial U$. Soit~$U'$ la composante connexe de~$X\setminus S$ contenant~$U$. Comme~$x\notin \partial \an V, x\notin \partial \an X$. Puisque~$(X,S)$ satisfait~$\delta)$, il existe une branche~$b$ de~$X$ issue de~$x$ et non contenue dans~$U'$. 

Le point~$x$ n'appartenant pas au bord analytique de~$V$, il n'appartient en particulier pas à son bord topologique dans~$X$ ; autrement dit,~$V$ est un voisinage de~$x$ dans~$X$. Ainsi,~$b\in \br V x$ et comme~$b$ n'est pas contenue dans~$U'$, elle ne l'est {\em a fortiori} pas dans~$U$. 

\trois{structvinters} Ainsi~$S\cap V$ est-il un sous-ensemble sommital de~$V$. Soit~$W$ un élément de l'atlas~$\mathsf{Atl}(S\cap V)$ de~$V$ ; nous allons montrer qu'il appartient à~$\mathsf{Atl}(S)$, ce qui montrera que~$V$ est réunion d'éléments de~$\mathsf{Atl}(S)$ et achèvera la démonstration. 

\medskip
D'après ce qu'on a vu au cours de la preuve du théorème~\ref{structsomm}, il existe un sous-ensemble fini~$S'$ de~$S\cap V$ et une famille~$(V_j)$ de composantes connexes de~$V\setminus S'$ tels que les conditions a), b), c) et d) du~\ref{relcompsomm}, énoncées pour~$(X,S,S',W, (U_i))$, soient satisfaites par~$(V,S\cap V, S',W,(V_j))$. 

\medskip
{\em On a~$\partial_X W\subset S'$.} En effet,~$\partial_X W\subset \partial_V W\cup (\partial_X V\cap W)$ ; or~$\partial_XW\subset S'$ par la propriété a), et~$\partial_X V\subset S\cap V$ ; mais les propriétés a) et c) garantissent que~$(S\cap V)\cap W=S'$, d'où notre assertion. 

Par conséquent,~$W$ s'écrit~$X-\coprod U_i$ pour une certaine famille~$(U_i)$ de composantes connexes de~$X\setminus S'$. Nous allons montrer que~$(X,S,S',W, (U_i))$ satisfait également les propriétés a), b), c), et d), ce qui permettra de conclure. 

\medskip
La propriété a) résulte de la définition des~$U_i$. 

Vérifions  b). Soit~$x\in S'$ et soit~$I$ l'ensemble des indices tels que~$x\in \partial U_i$. Si~$x\notin \partial_X W$ alors~$I$ est vide, et si~$x\in \partial_X W$ l'ensemble~$I$ est fini en vertu de la proposition~\ref{propdomanferm}. Supposons de plus que~$x\notin \partial \an X$. Deux cas peuvent se présenter : ou bien~$x\notin \partial \an V$, auquel cas il existe~$j_0$ tel que~$x$ adhère à~$V_{j_0}$ (puisque~$(V,S\cap V, S',W,(V_j))$ satisfait b) ), et si~$i_0$ désigne l'indice tel que~$V_{j_0}\subset U_{i_0}$, alors~$i_0\in I$ ; ou bien~$x\in \partial \an V$, auquel cas~$x$ appartient à~$\partial_X V$, et {\em a fortiori} à~$\partial_X W$, ce qui entraîne que~$I$ est non vide. 

Vérifions c). Comme~$(V,S\cap V, S',W,(V_j))$ satisfait c), l'ensemble~$(S\cap V)\setminus S'~$ ne rencontre pas~$W$ ; on déduit alors de l'inclusion~$W\subset V$ que~$(S\setminus S')\cap W=\emptyset$. 

Vérifions d). Soit~$Z$ une composante irréductible de~$X$ ne rencontrant pas~$S'$. L'intersection~$Z\cap V$ s'écrit~$\bigcup Z_l$, où chaque~$Z_l$ est une composante irréductible de~$V$ ne rencontrant pas~$S'$. Comme~$(V,S\cap V, S',W,(V_j))$ satisfait d), on a~$Z_l\cap W=\emptyset$ pour tout~$l$ ; on déduit alors de l'inclusion~$W\subset V$ que~$Z\cap W=\emptyset$, ce qui achève la démonstration.~$\Box$ 

\deux{defcourbepingl} On définit comme suit la catégorie des  {\em courbes~$k$-analytiques épinglées} : ses objets sont les couples~$(X,S)$, où~$X$ est une courbe~$k$-analytique et~$S$ un ensemble sommital de~$X$ ; un morphisme~$\phi:(Y,T)\to (X,S)$ est un morphisme~$f:Y\to X$ telle que~$\phi^{-1}(S)\subset T$. On appelle catégorie des courbes {\em strictement}~$k$-analytiques épinglées la sous-catéorie pleine de la précédente dont les objets sont les couples~$(X,S)$ avec~$X$ strictement~$k$-analytique et~$S$ inclus dans~$X\typ 2$. 

\medskip
Remarquons qu'une courbe strictement~$k$-analytique épinglée compacte est simplement une courbe strictement~$k$-analytique compacte munie d'un ensemble fini de points de type 2 rencontrant toutes ses composantes irréductibles. 

\deux{fonctstruct} {\bf Proposition.} {\em Soit~$\phi: Y\to X$ un morphisme de courbes~$k$-analytiques, soit~$S$ un ensemble sommital de~$X$ et soit~$T$ un ensemble sommital de~$Y$. Les propositions suivantes sont équivalentes : 

1)~$\phi$ est un morphisme de courbes~$k$-analytiques formelles de~$(Y,\mathsf{Atl}(T))$ vers~$(X,\mathsf{Atl}(S))$ ; 

2)~$\phi$ est un morphisme de courbes~$k$-analytiques épinglées de~$(Y,T)$ vers~$(X,S)$.}

\medskip
{\em Démonstration.} Supposons que 1) soit vraie, soit~$x\in S$, soit~$y$ un antécédent de~$x$ et soit~$U$ un élément de~$\mathsf{Atl}(S)$ contenant~$x$. D'après l'hypothèse 1), l'image réciproque de~$U$ est réunion d'éléments de~$\mathsf{Atl}(T)$ ; il existe donc~$V\in \mathsf{Atl}(T)$ tel que~$\phi(V)\subset U$ et tel que~$y\in V$. Comme~$x$ appartient à~$U\cap S$, il appartient au bord de Shilov de~$U$ ; il s'ensuit que~$y$ appartient au bord de Shilov de~$V$, donc à~$V\cap T$ et en particulier à~$T$. Ainsi,~$\phi^{-1}(S)\subset T$ et 2) est vraie. 

\medskip
Supposons maintenant que 2) soit vraie et soit~$U\in \mathsf{Atl}(S)$ ; l'image réciproque~$V$ de~$U$ sur~$Y$ est un domaine analytique fermé de~$Y$, dont le bord analytique est contenu dans~$\phi\inv(\partial\an U)$, donc dans~$\phi\inv(S)$ et par conséquent dans~$T$ d'après l'hypothèse 2).  Il découle alors de la proposition~\ref{cardomform} ci-dessus que~$V\in \mathsf{Atl}(T)$, ce qui achève la démonstration.~$\Box$ 

\deux{reform} La proposition ci-dessus peut se reformuler comme suit : les flèches~$(X,S)\mapsto (X,\mathsf{Atl}(S))$ et~$(X,\got A)\mapsto (X,\mathsf {Somm}(\got A))$ établissent une équivalence entre la catégorie des courbes~$k$-analytiques épinglées (resp. des courbes strictement~$k$-analytiques épinglées) et celle des courbes~$k$-analytiques formelles (resp. des courbes strictement~$k$-analytiques formelles).

\deux{reformvdoualgc} Supposons de plus que~$|k\ti|\neq\{1\}$ et que~$k$ est ou bien algébriquement clos, ou bien de valuation discrète. Si~$\got X$ est une~$k\zero$-courbes formelle relativement normale dans sa fibre générique alors~$\mathsf{Somm}(\mathsf{Aff}(\got X))$ n'est autre que~$\rho^{-1}(\got X_{\red k}^{(0)})$, où~$\rho : \got X\to \got X_{\red k}$ est la flèche de réduction et~$\got X_{\red k}^{(0)}$ l'ensemble des points génériques des composantes irréductibles de~$\got X_{\red k}$. Bien entendu, on aurait tout aussi bien pu dans ce qui précède remplacer~$\got X_{\red k}$ par~$\got X_s$, puisque ces deux espaces topologiques coïncident. 

Il s'ensuit, en vertu de~\ref{reform} et~\ref{equivvariees} que les flèches~$(X,S)\mapsto  \widehat{(X,\mathsf{Atl}(S))}$ et~$\got X\mapsto (\got X_\eta, \rho^{-1}(\got X_{\red k}^{(0)}))$ établissent une équivalence entre la catégorie des courbes strictement~$k$-analytiques épinglées réduites et celle des~$k\zero$-courbes formelles relativement normales dans leur fibre générique. 

\medskip
Compte-tenu de GAGA et de l'algébrisation des~$k\zero$-courbes formelles propres, l'équivalence de catégories évoquée ci-dessus se décline comme suit dans le cas propre :~$(X,S)\mapsto \widehat{(X,\mathsf{Atl}(S))}^{\rm alg}$ et~$\sch X\mapsto (\sch X_\eta\an,  \rho^{-1}(\sch X_{\red k}^{(0)}))$ établissent une équivalence entre la catégorie des courbes strictement~$k$-analytiques épinglées réduites et propres, et celle des~$k\zero$-courbes algébriques propres relativement normales dans leur fibre générique.

\deux{modstabcan} {\bf Théorème.} {\em Soit~$X$ un courbe~$\KK$-analytique quasi-lisse, compacte, connexe et non vide admettant une plus petite triangulation~$S$ et soit~$Y$ une courbe~$\kk$-analytique quasi-lisse, compacte, connexe et non vide munie d'un morphisme fini et plat~$\phi : Y\to X$.

\medskip
1) La courbe~$Y$ admet une plus petite triangulation~$T$. 

2) Le morphisme~$\phi$ est un morphisme de courbes~$\KK$-analytiques formelles de~$(Y,T)$ vers~$(X,S)$.}

\medskip
{\em Démonstration.} Comme~$X$ possède une plus petite triangulation~$S$, l'ensemble des nœuds de~$\skelan X$ est non vide et coïncide avec~$S$. Soit~$T$ l'ensemble des nœuds de~$\skelan Y$. Il résulte de la proposition~\ref{transfernoeuds} que~$\phi^{-1}(S)\subset T$. Par conséquent,~$T$ est non vide, et~$Y$ possède ainsi une plus petite triangulation, à savoir~$T$. Comme~$\phi^{-1}(S)\subset T$, le morphisme~$\phi$ est bien, en vertu de la proposition~\ref{fonctstruct} ci-dessus, un morphisme de courbes~$\kk$-analytiques formelles de~$(Y,T)$ vers~$(X,S)$.~$\Box$ 

\deux{ideeacl} {\bf Proposition.}  {\em Soit~$X$ une courbe~$k$-analytique compacte, quasi-lisse, géométriquement connexe et non vide, et soit~$\phi : X\to X$ un morphisme fini et plat de degré strictement supérieur à 1. Supposons qu'il existe un ensemble sommital~$S$ de~$X$ tel que~$\phi^{-1}(S)\subset S$, c'est-à-dire tel que~$\phi$ soit un morphisme de courbes~$k$-analytiques formelles~$(X,S)\to (X,S)$ (prop.~\ref{fonctstruct}). L'ensemble~$S$ est alors un singleton~$\{x\}$, le point~$x$ est de type 2, et~$\got s(x)=k$.}

\medskip
{\em Démonstration.} On se ramène immédiatement par extension des scalaires à~$\KK$ au cas où~$k$ est algébriquement clos : le fait que~$S_{\KK}$ soit un singleton assurera que~$S$ est un singleton~$\{x\}$ et que~$\got s(x)=k$.

\trois{ideeaclssingle} Nous allons tout d'abord démontrer que~$S$ est un singleton. Pour cela, on raisonne par l'absurde. On suppose donc que~$S$ compte au moins deux éléments, et l'on va prouver que le degré de~$\phi$ est égal à~$1$. Comme~$\phi^{-1}(S)\subset S$ et comme~$\phi$ est surjective on a~$\phi^{-1}(S)=S$ pour des raisons de cardinal. Ainsi,~$\phi$ induit ainsi une permutation de~$S$. Quitte à remplacer~$\phi$ par~$\phi^{\circ n}$ pour~$n$ convenable, on peut supposer que~$\phi$ fixe~$S$ point par point. 

\medskip
Soit~$\Pi$ l'ensemble des composantes connexes de~$X\setminus S$ dont le bord contient au moins deux points ; c'est un ensemble fini, et non vide puisque~$S$ n'est pas un singleton. Si~$U$ est une composante connexe de~$X\setminus S$ alors~$\phi(U)$ est une composante connexe de~$X\setminus S$ de bord égal à~$\phi(\partial U)$, et donc à~$\partial U$. Par conséquent,~$\phi(U)\in \Pi$ si et seulement si~$U\in \Pi$. Ainsi,~$\phi$ induit une permutation de~$\Pi$ ; quitte, une fois encore, à remplacer~$\phi$ par ~$\phi^{\circ n}$ pour~$n$ convenable, on peut supposer que cette permutation est égale à l'identité. Choisissons~$U\in \Pi$ ; on a~$\phi(U)=U$ et donc~$\phi^{-1}(U)=U$. 

\medskip
{\em Le squelette analytique~$\skelan U$ est un graphe fini.} En effet, il existe un sous-graphe analytiquement admissible~$\Gamma$ de~$X$ qui contient~$S$ et est localement fini ; par compacité de~$X$, le graphe~$\Gamma$ est fini. Les composantes connexes de~$X\setminus S$ sont alors, d'une part les ouverts de la forme~$I^\flat$, où~$I$ est une composante connexe de~$\Gamma\setminus S$, et d'autre part les composantes connexes de~$X-\Gamma$ dont l'unique point du bord appartient à~$S$. Comme~$\partial U$ n'est pas un singleton,~$U$ est nécessairement de la forme~$I^\flat$, où~$I$ est une composante connexe de~$X\setminus S$ ; par conséquent,~$U$ possède un sous-graphe analytiquement admissible qui est fini, à savoir~$I$, et son squelette analytique est {\em a fortiori} fini. 

Soit~$T$ l'ensemble des nœuds de~$\skelan U$. Il est fini, car égal à la réunion : de l'ensemble des sommets topologiques de~$\skelan U$, qui est fini par finitude de ce dernier ; de l'ensemble~$\partial \an U$, qui est fini puisque contenu dans~$\partial \an X$ ; et de l'ensemble des points de~$U$ de genre strictement positif, qui est fini puisque contenu dans l'ensemble des points de genre strictement positif de ~$X$. La proposition~\ref{transfernoeuds} assure que~$\phi^{-1}(T)\subset T$. Par le même raisonnement que ci-dessus on peut, en remplaçant une dernière fois~$\phi$ par ~$\phi^{\circ n}$ pour~$n$ convenable, supposer que~$\phi$ fixe~$T$ point par point et que~$\phi(V)=\phi^{-1}(V)=V$ pour toute composante connexe~$V$ de~$U\setminus T$ dont le bord dans~$X$ compte au moins deux éléments. 

\medskip
{\em Il existe une composante connexe~$V$ de~$U\setminus T$ qui est un arbre à deux bouts.} Pour le voir, notons pour commencer que comme~$\partial U$ compte au moins deux points, l'ouvert~$U$ a au moins deux bouts ; dès lors,~$\skel U$ est non vide, et est donc un sous-graphe admissible de~$U$. Son bord dans~$X$ est dès lors égal à~$\partial U$, et comprend en particulier au moins deux points distincts~$x$ et~$y$. Il existe un intervalle ouvert contenu dans~$\skel U$ et joignant~$x$ à~$y$. 

\medskip
Si~$I$ rencontre~$T$, on pose~$J=]x;z[$ où~$z$ est le point de~$T$ le plus  proche de~$x$ ; sinon, l'on pose~$J=I$. Comme~$T$ contient tous les sommets topologiques de~$\skelan U$, l'intervalle~$J$ est à la fois ouvert et fermé dans~$\skelan U\setminus T$, et c'en est donc une composante connexe. Comme~$\skelan U$ est non vide, c'est un sous-graphe admissible de~$U$ et~$V:=J^\flat$ est une composante connexe de~$U\setminus T$. 

\medskip
La composante~$V$ est un arbre à deux bouts, évidemment relativement compact dans~$X$. Son bord est égal à~$\partial J$, et consiste donc en deux points de type 2 et 3. Enfin,~$V$ ne contient aucun point de~$T$, et donc aucun point de genre strictement positif ni aucun point de~$\partial \an X$. La proposition~\ref{propcarintdc} assure alors que~$V$ est une couronne. 

\medskip
Comme le bord de~$V$ contient au moins deux points, on a~$\phi^{-1}(V)=V$. Ainsi,~$\phi$ induit une flèche~$V\to V$ finie et plate de degré~$\deg\phi$. Il vient~$$\mathsf{Mod}(V)=\mathsf{Mod}(V)^{\deg \phi},$$ et partant~$\deg \phi=1$, ce qui est contradictoire avec notre hypothèse de départ, et montre ainsi que~$S$ est un singleton~$\{x\}$. 

\trois{ideeaclxt2} {\em Le point~$x$ est de type 2.} L'hypothèse~$\phi^{-1}(S)=S$ devient donc~$\phi^{-1}(x)=\{x\}$ ; par conséquent,~$\phi$ induit une~$k$-extension~$\hres(x)\hookrightarrow \hres(x)$ de degré~$\deg\phi$. Supposons~$x$ de type 3. L'égalité~$\red{\hres(x)}=\red k$ (rappelons que ~$k$ est supposé algébriquement clos) assure alors que la~$\red k$-extension~$\red{\hres(x)}\hookrightarrow \red{\hres(x)}$ induite par~$\phi$ est triviale ; comme le plongement~$|\hres(x)\ti|\hookrightarrow |\hres(x)\ti|$ induit par~$\phi$ est quant à lui évidemment trivial, la stabilité du corps~$\hres(x)$ (th.~\ref{stable23}) assure que~$\deg \phi=1$, ce qui est absurde.~$\Box$

\section{Application : la réduction semi-stable}

{\em Dans cette section, l'on suppose que~$|k\ti|\neq\{1\}$.}

\subsection*{Le cas analytique} 
On fixe pour ce paragraphe une courbe {\em strictement}~$k$-analytique quasi-lisse~$X$. 

\deux{trisommitale} En vertu de l'assertion iv) du théorème~\ref{theotri}, la courbe~$X$ admet une triangulation dont l'ensemble~$S$ des sommets est de type 2. Le couple~$(X, S)$ satisfait les conditions~$\alpha),\beta), \gamma)$ et~$\delta)$ de l'assertion  ii) du théorème~\ref{structsomm}. C'est en effet clair pour~$\alpha)$ et~$\gamma)$ ; pour~$\beta)$, cela résulte du fait que~$S$ rencontre toutes les composantes connexes de~$X$ (qui sont aussi ses composantes irréductibles par quasi-lissité) ; et pour~$\delta)$, c'est une conséquence de~\ref{commentdelta}.

\medskip
L'ensemble~$S$ est donc sommital, et définit de ce fait une courbe~$k$-analytique formelle~$(X,S)$. Les points de~$S$ étant tous de type 2, la courbe~$(X,S)$ est même  {\em strictement}~$k$-analytique formelle ; ceci permet de lui associer un~$k\zero$-schéma formel plat~$\fort X S$.

\deux{casdist} {\em Supposons que la courbe analytique formelle~$(X,S)$ est distinguée et que~$\fort X S _s$ est géométriquement réduit}. Dans ce cas le~$k\zero$-schéma formel~$\fort X S$ est localement topologiquement de présentation finie, et~$$\fort X S\hotimes_{k\zero}(\KK)\zero\simeq \fort {X_{\KK}}{S_{\KK}}.$$ Soit~$\bf x$ un point fermé de~$\fort {X_{\KK}}{S_{\KK}}_s$. Son image réciproque~$\rho\inv({\bf x})$ par la flèche de spécialisation~$\rho$ est une composante connexe de~$X_{\KK}\setminus S_{\KK}$. Le sous-ensemble~$S_{\KK}$ de~$X_{\KK}$ en est une triangulation, puisque~$S$ est une triangulation de~$X$ ; par conséquent,~$\rho\inv({\bf x})$ est un disque ou une couronne ; les sommets de~$S_{\KK}$ étant de type 2,  le lemme~\ref{rayextshrx} garantit  de plus que si~$\rho\inv({\bf x})$ est un disque (resp. une couronne) alors son rayon modulo~$|(\KK)\ti|$ est trivial (resp. ses deux rayons extérieurs modulo~$|(\KK)\ti|$ sont triviaux). Or en vertu des résultats de Bosch et Lütkebohmert, l'ouvert~$\rho\inv({\bf x})$ est un disque (resp. une couronne) de rayon modulo~$|(\KK)\ti|$ trivial (resp. de rayons extérieurs modulo~$|(\KK)\ti|$ triviaux) si et seulement si~$\bf x$ est un point lisse (resp. un  point double ordinaire) de~$\fort {X_{\KK}}{S_{\KK}}_s$ ; il s'ensuit que le~$(\KK)\zero$-schéma formel plat~$\fort {X_{\KK}}{S_{\KK}}=\fort X S\hotimes_{k\zero}(\KK)\zero$ est semi-stable ; par conséquent, le~$k\zero$-schéma formel~$\fort X S$ est semi-stable. 

\medskip
Donnons deux exemples, l'un général et l'autre extrêmement particulier, relevant de ce qui précède. 

\trois{casdistplus} Si~$k$ est algébriquement clos alors~$(X, S)$ est distinguée (\ref{suffdist}) et~$\fort X S$ est géométriquement réduit ; par conséquent,~$\fort X S$ est semi-stable.

\trois{poincare} Supposons : que~$k$ est une extension finie de~$\QQ_p$ ; que~$X$ est le {\em demi-plan de Poincaré sur~$k$}, c'est-à-dire l'ouvert~$\pk-\PP^{1}(k)$ de~$\pk$ (comme~$k$ est localement compact,~$\PP^1(k)$ est compact) ; et que~$S$ est le sous-ensemble de~$X$ formé des points~$\eta_{a,r}$ où~$a$ parcourt~$k$ et où~$r$ parcourt~$|k\ti|$ (on vérifie sans difficultés, par passage à~$\KK$, que c'est une triangulation de~$X$ ; elle est par construction constituée de points de type 2). Les conditions~$\alpha)$,~$\beta)$ et~$\gamma)$ du~\ref{suffdist} sont satisfaites ; par conséquent,~$S$ est distinguée. Comme~$\red k$ est fini, il est parfait et~$\fort X S_s$ est donc géométriquement réduit. Il s'ensuit que~$\fort X S$ est semi-stable. 

\deux{semstabpot} {\em Supposons que~$X$ est compacte, sans hypothèse particulière sur~$S$.} Dans ce cas, on peut recouvrir~$X$ par un nombre fini de domaines affinoïdes formels. Comme chacun d'eux devient distingué à fibre spéciale géométriquement réduite après une extension finie séparable convenable de~$k$, il existe une extension finie séparable~$L$ de~$k$ telle que~$(X_L,S_L)$  soit distinguée et telle que~$\fort {X_L}{S_L}_s$ soit géométriquement réduite ; le~$L\zero$-schéma formel~$\fort {X_L}{S_L}$ est alors en vertu du~\ref{casdist} un modèle semi-stable de~$X$ ; on vient ainsi de démontrer que {\em toute courbe strictement~$k$-analytique quasi-lisse et compacte admet un modèle formel semi-stable après extension finie séparable de~$k$.}

\subsection*{Le cas algébrique}

\deux{rssalg} Soit~$\sch X$ une~$k$-courbe projective et lisse. D'après le~\ref{semstabpot} ci-dessus, appliqué à~$X=\sch X\an$, il existe une extension finie séparable~$L$ de~$k$ et une~$L\zero$-courbe formelle semi-stable~$\got X$ tel que~$\got X_\eta\simeq (\sch X_L)\an$ ; en vertu du théorème ~\ref{algebrformcurv}, le schéma formel~$\got X$ s'identifie au complété formel d'un~$L\zero$-schéma propre et plat~$\bnd X$ le long de sa fibre spéciale ; on a alors~$$(\bnd X_\eta)\an\simeq \got X_\eta\simeq (\sch X_L)\an\;,$$ ce qui, en vertu de GAGA, entraîne que~$\bnd X_\eta\simeq \sch X_L$. Ainsi,~$\sch X$ admet un modèle semi-stable {\em algébrique} après extension finie séparable du corps de base. 

\deux{rsscdn} Supposons que~$|k\ti|$ est libre de rang 1, et soit~$k'$ un sous-corps dense de~$k$. Soit~$\sch X$ une~$k'$-courbe projective et lisse. D'après le~\ref{semstabpot} ci-dessus, appliqué à~$X=(\sch X_k)\an$, il existe une extension finie séparable~$L$ de~$k$ et une~$L\zero$-courbe formelle semi-stable~$\got X$ telle que~$\got X_\eta\simeq (\sch X_L)\an$. En vertu du lemme de Krasner, il existe une extension finie séparable~$L'$ de~$k'$ telle que~$L\simeq L'\otimes_{k'}k$ ; le corps~$L'$ s'identifie naturellement à un sous-corps dense de~$L$, et~$\sch X_L\simeq (\sch X_{L'})\otimes_{L'}L$. 

On déduit alors du théorème~\ref{algebrformcurv} l'existence d'un~$(L')\zero$-schéma projectif et plat~$\bnd Y$, purement de dimension relative 1, tel que~$\bnd Y_\eta\simeq \sch X_{L'}$ et tel que~$\got X$ s'identifie au complété formel de~$\bnd Y\otimes_{(L')\zero}L\zero$ le long de sa fibre spéciale. Ainsi, ~$\sch X$ admet un modèle semi-stable {\em algébrique} après extension finie séparable de son corps de base~$k'$.

\section{Le cas particulier de la valuation discrète}

\deux{introvd} On suppose pour toute cette section que~$|k\ti|$ est libre de rang 1. Les anneaux de la forme~$k\zero\{T_1,\ldots, T_n\}$ sont alors noethériens ; tout~$k\zero$-schéma formel topologiquement de type fini est donc noethérien. 

\medskip
Rappelons quelques faits établis plus haut. Soit~$\got X$ un~$k\zero$-schéma formel localement topologiquement de présentation finie et relativement normal dans sa fibre générique. Sa fibre générique~$\got x_\eta$ est réduite, l'ensemble~$\mathsf{Aff}(\got X)$ des domaines affinoïdes de~$\got X_\eta$ qui sont de la forme~$\got U_\eta$ où~$\got U$ est un ouvert affine de~$\got X$ est un atlas strictement~$k$-analytique formel sur~$\got X_\eta$, et~$\widehat{( \got X_\eta,\mathsf{Aff}(\got X))}\simeq \got X$. La fibre spéciale «absolue»~$\got X_s$ est naturellement isomorphe à~$(\got X_{\red k})_{\rm red}$. On dispose donc d'un homéomorphisme canonique~$\got X_s\simeq \got X_{\red k}$, modulo lequel les applications de réduction~$\got X_\eta\to \got X_s$ et~$\got X_\eta \to \got X_{\red k}$ coïncident. Nous nous permettrons de les désigner indistinctement par~$\rho$, sans référence explicite à~$\got X$. On désignera par~$\got X_{\red k}^{(0)}$ l'ensemble des points de codimension nulle de~$\got X_{\red k}$. L'ensemble~$\rho^{-1}(\got X_{\red k}^{(0)})$ est l'ensemble des sommets de l'atlas~$\mathsf{Aff}(\got X)$. Ce dernier est distingué si et seulement si~$|\hres(x)\ti|=|k\ti|$ pour tout~$x\in \rho^{-1}(\got X_{\red k}^{(0)})$, ou encore si et seulement si~$\got X_{\red k}$ est réduit, c'est-à-dire égal à~$\got X_s$. 

\deux{normcorpsres} {\bf Proposition.} {\em Soit~$\got X$ un~$k\zero$-schéma formel localement topologiquement de présentation finie et relativement normal dans sa fibre générique, soit~$\xi$ appartenant à~$\got X_{\red k}^{(0)}$ et soit~$x$ son unique antécédent par~$\rho$. 

\medskip
1) La multiplicité de la composante irréductible~$\overline{ \{\xi\}}$ de~$\got X_{\red k}$, c'est-à-dire la longueur de l'anneau local artinien~$\sch O_{\got X_{\red k},\xi}$, est égale à l'indice~$[|\hres(x)\ti:|k\ti|]$. 

2) Le corps~$\red{\hres(x)}$ s'identifie à~$\kappa(\xi)$.}

\medskip
{\em Démonstration.} La question est locale sur~$\got X$, ce qui autorise à le supposer affine, donc de la forme~${\rm Spf}\; \mathsf A$ où~$\mathsf A$ est une~$k\zero$-algèbre plate topologiquement de présentation finie intégralement close dans~$\mathsf A\otimes_{k\zero}k$ ; on peut également supposer qu'il est purement de dimension~$d$ pour un certain~$d$. Il en va alors de même de l'espace affinoïde~$\got X_\eta$ ; comme celui-ci est réduit, son lieu~$Y$ de non-normalité est de dimension strictement inférieure à~$d$. Le corps~$\red{\hres(x)}$ étant de degré de transcendance~$d$ sur~$\red k$, le fermé de Zariski~$Y$ de~$X$ ne contient pas~$x$ ; son image~$\rho(Y)$ évite donc~$\xi$. Comme elle est constructible (et même fermée), on peut, quitte à encore restreindre~$\got X$, supposer que~$\got X_\eta$ est normal et même intègre. Puisque~$\mathsf A$ est intégralement close dans~$\mathsf A\otimes_{k\zero}k$, l'algèbre intègre~$\mathsf A$ est normale. Soit~$\mathsf K$ son corps des fractions. 

\medskip
Le schéma~$\got X_{\red k}$ s'identifie à~$\spec \mathsf A\otimes_{k\zero}\red k$ ; le point~$\xi$ correspond par ce biais à un point de~$\spec \mathsf A$ que l'on note encore~$\xi$ et qui est de hauteur 1 dans le schéma noethérien normal~$\spec \mathsf A$. L'anneau local~$\sch O_{\spec \mathsf A,\xi}$ est donc l'anneau d'une valuation discrète~$v$ de~$\mathsf K$ qui prolonge la valuation de~$k\zero$, normalisée pour que~$v(\mathsf k\ti)=\ZZ$. Soit~$\pi$ une uniformisante de~$k\zero$ et soit~$e$ l'entier~$v(\pi)$. La longueur de~$\sch O_{\got X_{\red k},\xi}$ est égale à celle de~$\sch O_{\spec \mathsf A,\xi}/\pi$, et donc à~$e$.

\medskip
Soit~$\phi : \mathsf K\to \RR_+$ l'application qui envoie~$0$ sur~$0$ et~$a$ sur~$|\pi|^{v(a)/e}$ si~$a\neq 0$. C'est une valeur absolue dont la restriction à~$k$ coïncide avec la valeur absolue structurale de ce dernier, et qui est bornée par~$1$ sur~$\mathsf A$. Elle définit donc un point~$x'$ de~$\sch M(\mathsf A\otimes_{k\zero}k)$. Si~$a\in \mathsf A$ on a par construction~$a(\xi)=0$ si et seulement si~$\phi(a)<1$, ce qui signifie que~$\rho(x')=\xi$ ; par conséquent,~$x'=x$. 

\medskip
Par construction,~$(\mathsf K,\phi)$ est un sous-corps valué dense de~$\hres(x)$. Le corps résiduel de~$\hres(x)$ est dès lors égal à celui de~$(\mathsf K,\phi)$, c'est-à-dire à celui de~$v$ et donc à~$\kappa(\xi)$ ; et son groupe des valeurs n'est autre que~$\phi(\mathsf k\ti)=|\pi|^{(1/e)\ZZ}$. Par conséquent,~$[|\hres(x)\ti:|k\ti|]=e$, ce qui achève la démonstration.~$\Box$

\trois{precisepingl} La flèche~$\got X\mapsto (\got X_\eta, \rho\inv(\got X_{\red k}^{(0)}))$ établit une équivalence entre la catégorie des~$k\zero$-courbes formelles relativement normales dans leur fibre générique et celle des courbes strictement~$k$-analytiques épinglées réduites. Si~$(X,S)$ est une courbe strictement~$k$-analytique épinglée réduite, si~$\got X$ désigne la courbe formelle associée, et si~$x\in S$, la proposition~\ref{normcorpsres} ci-dessus assure que la composante irréductible de~$\got X_{\red k}$ qui correspond à~$x$ a pour multiplicité~$[|\hres(x)\ti|:|k\ti|]$ et pour corps des fonctions~$\red{\hres(x)}$.

\trois{precisepinglpro} Dans le cas propre on peut donner, grâce à GAGA et à l'algébrisation des courbes formelles propres, une variante algébrique (et non plus formelle) de l'équivalence évoquée ci-dessus, qui s'énonce comme suit. 

La flèche~$\bnd X\mapsto (\bnd X_\eta\an, \rho\inv(\bnd  X_{\red k}^{(0)}))$ établit une équivalence entre la catégorie des~$k\zero$-courbes algébriques projectives relativement normales dans leur fibre générique et celle des courbes strictement~$k$-analytiques épinglées réduites et propres. Si~$(X,S)$ est une courbe strictement~$k$-analytique épinglée réduite et propre, si~$\bnd X$ désigne la courbe algébrique sur~$k\zero$ associée, et si~$x\in S$, la proposition~\ref{normcorpsres} ci-dessus assure que la composante irréductible de~$\bnd X_{\red k}$ qui correspond à~$x$ a pour multiplicité~$[|\hres(x)\ti|:|k\ti|]$ et pour corps des fonctions~$\red{\hres(x)}$. 

\medskip
Rappelons qu'une courbe strictement~$k$-analytique épinglée propre est simplement une courbe~$k$-analytique propre munie d'un ensemble fini de points de type 2 rencontrant chacune de ses composantes irréductibles.

\deux{segmdense} Nous allons terminer ce paragraphe en mentionnant une propriété de rareté des points~$x$ d'une courbe~$k$-analytique qui sont de type 2 et satisfont la condition~$|\hres(x)\ti|=|k\ti|$. Soit donc~$X$ une courbe~$k$-analytique et soit~$I$ un segment tracé sur~$X\typ{23}$. L'ensemble~$\sch E$ des points~$x$ de~$I$ qui sont de type 2 satisfont la condition~$|\hres(x)\ti|=|k\ti|$ est alors fini. En effet, soit~$x\in I$ et soit~$J$ une composante connexe de~$I\setminus\{x\}$ ; elle définit une branche~$b$ issue de~$x$ et correspond donc à une valuation graduée non triviale~$\langle.\rangle$ du corps résiduel gradué~$\gred{\hres(x)}$. Choisissons~$f\in \kappa(x)\ti$ tel que~$\langle \gred{f(x)}\rangle\neq 1$ ; la fonction~$|f|$ est alors non constante le long de~$b$. Il découle dès lors de~\ref{defpentesholo} {\em et sq.} qu'il existe un intervalle ouvert~$J'$ contenu dans~$J$ et aboutissant à~$x$ tel que~$|f|$ induise un homéomorphisme de~$J'$ sur un intervalle ouvert relativement compact~$J_0$ de~$\RR\ti_+$. Comme~$|k\ti|$ est un sous-groupe discret de~$\RR\ti_+$, l'intersection~$J_0\cap |k\ti|$ est finie. Si~$x\in \sch E\cap J'$ alors~$|f(x)|\in J_0\cap |k\ti|$ ; il s'ensuit que~$J'$ ne contient qu'un nombre fini de points de~$\sch E$, ce qui permet de conclure par compacité de~$I$.

\subsection*{Disques unité maximaux sur une courbe quasi-lisse}

\deux{defdisqun} Commençons par un peu de vocabulaire. On appellera {\em~$k$-disque unité} tout~$k$-disque dont le rayon modulo~$|k\ti|$ est égal à~$1$. Et si~$X$ est une courbe~$k$-analytique et si~$x$ est un point de~$X\typ{23}$ on dira que~$x$ est {\em non ramifié} si~$|\hres(x)\ti|=|k\ti|$ (ce qui entraîne qu'il est de type 2) et si~$\red{\hres(x)}$ est une extension séparable de~$\red k$ -- c'est-à-dire si~$\Omega_{\red{\hres(x)}/\red k}$ est de dimension 1, ou encore si la~$\red k$-algèbre~$\red{\hres(x)}$ est géométriquement réduite.  

\deux{dvrdiscun} {\bf Proposition}. {\em Soit~$X$ une courbe~$k$-analytique et soit~$U$ un ouvert de~$X$ qui est un disque virtuel relativement compact contenant un~$k$-point. Soit~$x$ l'unique point de~$\partial U$, qui appartient à~$X\dtr$ en vertu de~\ref{borddisct23}. Les assertions suivantes sont équivalentes : 

\medskip
i)~$U$ est un disque unité ; 

ii)~$x$ est non ramifié.}

\medskip
{\em Démonstration.} D'après,~\ref{borddisct23} le point~$x$ appartient au lieu quasi-lisse de~$X$. On peut donc, quitte à remplacer~$X$ par son lieu quasi-lisse (qui contient le disque virtuel~$U$),  la supposer elle-même quasi-lisse. La réunion de~$\{x\}$ et de toutes les composantes connexes de~$X\setminus\{x\}$ qui sont relativement compactes est un domaine analytique compact~$X'$ de~$X$. En restreignant encore~$X'$ par ablation des composantes connexes de~$X'\setminus\{x\}$ qui rencontrent~$\partial \an X'$ (ce n'est pas le cas du disque virtuel~$U$), puis en substituant à~$X$ la composante connexe de~$x$ dans~$X'$, on se ramène finalement au cas où~$X$ est quasi-lisse, compacte, irréductible, et où~$\partial\an X\subset \{x\}$.  Le singleton~$\{x\}$ est alors un sous-ensemble sommital de~$X$. 

\medskip
\trois{discunimplnonram} {\em Supposons que i) soit vérifiée.} Le lemme~\ref{rayextshrx} assure alors que~$|\hres(x)\ti|\subset |k\ti|$ ; par conséquent,~$|\hres(x)\ti|=|k\ti|$, et~$x$ est de type 2. L'ensemble sommital~$\{x\}$ définit donc un modèle formel~$\got X$ de~$X$ sur~$k\zero$ qui est topologiquement de type fini, plat, et normal. L'ensemble~$\{x\}$ étant un singleton,~$\got X_{\red k}$ est irréductible. Soit~$\bf x$ le point fermé de~$\got X_{\red k}$ qui correspond à la composante connexe~$U$ de~$X\setminus\{x\}$ ; puisque~$U(k)$ est non vide,~$\bf x$ est un~$\red k$-point. 

Comme~$|\hres(x)\ti|=|k\ti|$ l'espace strictement~$k$-analytique formel~$(X,\mathsf{Atl}(\{x\}))$ est distingué, ce qui signifie que~$\got X_{\red k}$ est réduite. Le fait que l'image réciproque~$U$ de~$\bf x$ sur~$X$ soit un disque unité assure alors que~$\bf x$ est un point lisse de~$\got X_{\red k}$. La~$\red k$-courbe intègre~$\got X_{\red k}$ est dès lors génériquement lisse. Son corps des fonctions est de ce fait séparable sur~$\red k$ ; or il s'identifie à~$\red{\hres(x)}$ (prop.~\ref{normcorpsres}), d'où ii). 

\trois{nonramimpldiscun} {\em Supposons que ii) soit vérifiée.} L'égalité~$|\hres(x)\ti|=|k\ti|$ implique que~$(X,\mathsf{Atl}(\{x\}))$ est un espace strictement~$k$-analytique formel distingué. Il lui correspond donc un~$k\zero$-schéma formel topologiquement de type fini, plat, et normal~$\got X$ dont la fibre spéciale~$\got X_{\red k}$ est réduite ; comme le corps des fonctions de~$\got X_{\red k}$ s'identifie à~$\red{\hres(x)}$ (prop.~\ref{normcorpsres}), il est séparable sur~$\red k$, ce qui signifie que~$\got X_{\red k}$ est génériquement lisse. Soit~$F$ une extension finie purement inséparable de~$\red k$. Par ce qui précède, la~$F$-courbe~$\got X_F$ est génériquement réduite, et satisfait donc la propriété~$R_0$. Étant par ailleurs finie et plate sur la courbe~$\got X_{\red k}$, laquelle est réduite et satisfait en particulier~$S_1$, elle satisfait elle-même~$S_1$ ; par conséquent,~$\got X_F$ est réduite. Il en résulte que~$\got X_{\red k}$ est géométriquement réduite. 

Par conséquent,~$\form {X_{\KK}} {\{x\}_{\KK}}$ s'identifie à~$\got X_{(\KK)\zero}$. Soit~$\bf y$ l'unique antécédent du~$\red k$-point~$\bf x$ sur~$\got X_{\kk}$. L'image réciproque de~$\bf y$ dans~$X_{\KK}$ est égale à~$U_{\KK}$, qui est par hypothèse un disque ; l'unique point de~$\partial U_{\KK}$ étant situé au-dessus de~$x$, il est de type 2, ce qui force, en vertu du lemme~\ref{rayextshrx}, le rayon modulo~$|\KK\ti|$ de~$U_{\KK}$ à être égal à~$1$. Mais~$\bf y$ est alors un point lisse de~$\got X_{\kk}$ ; il s'ensuit que le~$\red k$-point~$\bf x$ est situé sur le lieu lisse de~$\got X_{\red k}$, et partant que~$U$ est un disque unité.~$\Box$

\bibliographystyle{plain}
\bibliography{aducros}

 \end{document}